\documentclass[11pt]{article}
\usepackage{bbm}
 \usepackage{amssymb}
\usepackage{amssymb, amsthm, amsmath, amscd}
\usepackage[usenames,dvipsnames]{color}

\newtheorem{thm}{Theorem}[section]
\newtheorem{lem}[thm]{Lemma}
\newtheorem{prop}[thm]{Proposition}
\newtheorem{cor}[thm]{Corollary}
\theoremstyle{definition}
\newtheorem{NN}[thm]{}
\theoremstyle{definition}
\newtheorem{df}[thm]{Definition}
\theoremstyle{definition}
\newtheorem{rem}[thm]{Remark}
\newtheorem{nota}[thm]{Notation}
\theoremstyle{definition}

\newcommand{\blue}{\color{black}}

\renewcommand{\phi}{\varphi}

\newcommand{\A}{\mathbb{A}}

\newcommand{\N}{\mathbb{N}}
\newcommand{\Z}{\mathbb{Z}}
\newcommand{\Q}{\mathbb{Q}}
\newcommand{\R}{\mathbb{R}}
\newcommand{\C}{\mathbb{C}}
\newcommand{\T}{\mathbb{T}}

\numberwithin{equation}{section}

\newcommand{\Aff}{\operatorname{Aff}}
\newcommand{\Inf}{\operatorname{Inf}}

\newcommand{\id}{\operatorname{id}}

\newcommand{\aff}{\rm aff}

\newcommand{\morp}{contractive completely positive linear map}
\newcommand{\cp}{completely positive linear map}

\newcommand{\hm}{homomorphism}
\newcommand{\dt}{\delta}
\newcommand{\ep}{\varepsilon}
\newcommand{\et}{\eta}
\newcommand{\sg}{\sigma}
\newcommand{\td}{\tilde}

\newcommand{\lr}{\longrightarrow}
\newcommand{\ld}{\lambda}
\newcommand{\cd}{\cdots}

\newcommand{\qq}{{\quad \quad}}
\newcommand{\sbs}{{\subset}}

\newcommand{\tht}{\theta}

\newcommand{\bb}{\mathfrak{b}}
\newcommand{\cc}{\mathfrak{c}}
\newcommand{\dd}{\mathfrak{d}}

\newcommand{\LD}{\Lambda}
\newcommand{\GM}{\Gamma}
\newcommand{\gm}{\gamma}
\newcommand{\sm}{\sigma}

\newcommand{\DT}{\Delta}

\newcommand{\spd}{(\diamondsuit)}
\newcommand{\spdd}{(\diamondsuit\diamondsuit)}
\newcommand{\spddd}{(\diamondsuit\diamondsuit\diamondsuit)}

\newcommand{\Tht}{\Theta}

\newcommand{\e}{\mbox{\large \bf 1}}
\newcommand{\0}{\mbox{\large \bf 0}}

\newcommand{\la}{\langle}
\newcommand{\ra}{\rangle}
\newcommand{\andeqn}{\,\,\,{\rm and}\,\,\,}
\newcommand{\rforal}{\,\,\,{\rm for\,\,\,all}\,\,\,}
\newcommand{\CA}{$C^*$-algebra}
\newcommand{\SCA}{$C^*$-subalgebra}

\newcommand{\af}{{\alpha}}
\newcommand{\bt}{{\beta}}
\newcommand{\dist}{{\rm dist}}

\newcommand{\one}{{\bf 1}}

\newcommand{\diag}{{\rm diag}}

\newcommand{\wilog}{without loss of generality}
\newcommand{\Wlog}{Without loss of generality}

\newcommand{\beq}{\begin{eqnarray}}
\newcommand{\eneq}{\end{eqnarray}}
\newcommand{\tforal}{\,\,\,\text{for\,\,\,all}\,\,\,}
\newcommand{\tand}{\,\,\,\text{and}\,\,\,}

\newcommand{\Blue}{\color{black}}
\newcommand{\Green}{\color{black}}

\usepackage{amsfonts}
\usepackage{mathrsfs}
\usepackage{textcomp}
\usepackage[all]{xy}

\usepackage{makeidx}

\makeindex

\title{A classification of finite simple amenable ${\cal Z}$-stable $C^*$-algebras, I:
$C^*$-algebras with generalized tracial rank one}
\author{Guihua Gong,  Huaxin Lin and Zhuang Niu
 }
\date{
}

\providecommand{\keywords}[1]
{\textbf{\textit{2010 AMS MathematicsSubject Classification:}} #1}
\begin{document}

\maketitle

\begin{abstract}

{{A class of \CA s, to be called those of generalized tracial
rank one, is introduced.
A
second class of unital simple separable amenable \CA s, those
whose tensor products with UHF-algebras of infinite type are in the
first class, to be referred to as those of rational generalized
tracial rank one, is proved to exhaust all possible values of the
Elliott invariant for unital finite simple separable amenable
${\cal Z}$-stable \CA s.
A number of results toward the classification of the second class are presented
including an
isomorphism theorem for a special sub-class of the first class, {{leading}}
 to the  general classification of
all  unital simple \CA s with rational generalized tracial
rank one in Part II.}}


\end{abstract}

\keywords{
 Primary 46L35, Secondary 46L80 and 46L05,\\
 \hspace{0.2in}{\bf {\em Keywords}}: Classification of Simple \CA s}

\tableofcontents

\section{Introduction}
The concept of a \CA\, exists harmoniously in many areas of mathematics. The
{{abstract}} definition of {{a \CA}}  axiomatized the norm closed self-adjoint subalgebras of $B(H),$ the algebra of all bounded linear operators on a Hilbert space $H.$ Thus, \CA s are operator algebras. The study of \CA s may also be viewed as the study of a non-commutative analogue of topology. This is because every unital commutative \CA\, is isomorphic to $C(X),$ the algebra of continuous functions on a compact Hausdorff space $X$ (by {{means}} of the Gelfand transform).
If we take our space $X$ and equip it with a group action via homeomorphisms, we enter the realm of topological dynamical systems, where remarkable progress has been made by considering the transformation \CA\, $C(X)\rtimes G$ arising {{via}}  the crossed product construction. Analogously, we can consider the study of general crossed products as the study of non-commutative topological dynamical systems.
{{The most special case is}} the group \CA\, of $G$, which is fundamental in the study of abstract harmonic analysis.
\CA\, theory is also basic
to the non-commutative geometry of Connes.
{{There are deep
and extremely important interactions between the theory of
\CA s and the theory of the very special concrete
\CA s (weak operator closed) called von Neumann
algebras.}}
{{One might continue in
this vein.}}
Therefore, naturally,  it would be of great interest to classify \CA s.
{{There has already
been noteworthy progress in this direction.}}

Early classification theorems start with the work of Glimm in the {{late}} 1950s who classified
{{infinite tensor products of matrix algebras, which he called
uniformly hyperfinite algebras}} (UHF-algebras) by supernatural numbers.
{{Dixmier, a few years later, classified
non-unital inductive limits of UHF-algebras, noting
that these were more complex. Bratteli,
in 1972, generalized the results of Glimm and Dixmier
to arbitrary inductive limits of finite-dimensional
\CA s (AF-algebras), using an analogue of
the combinatorial data of Glimm which is now called
a Bratteli diagram. (For instance, Pascal's triangle
is a Bratteli diagram.)}}
Bratteli's classification
of AF-algebras was reformulated
in a striking way by Elliott in 1976. He showed in
effect that Bratteli's equivalence relation on Bratteli
diagrams, which corresponded to isomorphism of the AF-algebras
they described, was in fact the same as
isomorphism of certain ordered groups arising in
a natural way from these diagrams: in a sense,
just their inductive limits. These ordered groups
turned out to be just the K-groups of the algebras
(generated by their Murray-von Neumann semigroups).

{{By 1989, Elliott
had begun  his classification program
by classifying AT-algebras  of
real rank zero by scaled ordered K-theory.
Notably, this used an approximate version of
an (exact) intertwining argument used by
both Bratteli and Elliott in the AF case.}}
Since then there has been rapid progress in the program to classify separable amenable
\CA s now known as the Elliott program.  Elliott {{and}} Gong (\cite{EG-RR0AH}) and Elliott, Gong, and Li (\cite{EGL-AH}) (together with a {{deep}} reduction theorem by Gong (\cite{Gong-AH})) classified
simple AH-algebras with no dimension growth by  {{means of}} the Elliott invariant (see {{Definition}} \ref{DEll} below).

The Elliott intertwining
argument
{{provided}} a framework for further classification proofs
 and
focused attention on {{the invariants of and the}}
maps between certain building block algebras
(for AF-algebras these building blocks would be finite dimensional algebras and for AH-algebras these building blocks would be  certain homogeneous \CA s). In particular, one wants to know when maps between invariants are induced by maps between building blocks (sometimes referred to as an existence theorem) and to know when maps between the building blocks are approximately unitarily equivalent (often called an uniqueness theorem).   To classify \CA s without assuming some particular inductive limit structure, one would like to
establish abstract existence theorems
and uniqueness theorems.   These efforts became the engine for these
rapid
developments
(\cite{Lnjotuni}, \cite{LinTAF2}, \cite{LnAUCT} and  \cite{DE} for example). Both existence theorems and uniqueness theorems
used $KL$-theory  (\cite{Ror-KL-I}) and the total $K$-theory ($\underline{K}(A)$) developed by
Dadarlat and Loring (\cite{DL}).   These existence and uniqueness theorems provide not only the technical tools for the classification program but also
{{the}} foundation
{{for understanding}} the morphisms in the category of \CA s.

 The rapid  developments mentioned above {{include}}  the Kirchberg-Phillips  {{classification}} (\cite{Kirch-Infty}, \cite{KP0}, and \cite{Ph1})  of purely infinite simple separable amenable \CA s which satisfy the UCT, by {{means of}} their $K$-theory.  There is also the classification of unital simple amenable \CA s in the UCT class which have tracial rank zero or one
(\cite{LinTAF2}, \cite{Lnduke}, and  \cite{LinTAI}).

On the other hand,  it had been suggested in  \cite{DNNP-AH} and \cite{BDR} that unital simple AH-algebras without a dimension growth condition might behave differently.  It was Villadsen (\cite{Vill-perf} and \cite{Vill-sr}) who showed that unital simple AH-algebras may have perforated $K_0$-groups and may have stable rank
{{equal to any non-zero natural number.}}
R{\o}rdam exhibited an  {{amenable}} separable simple $C^*$-algebra which is finite but not stably finite
 (\cite{Ror-infproj}).
{{It was  shown by}} Toms (\cite{Toms-Ann})  that there are unital simple AH-algebras of stable rank one with the same
Elliott invariant that are not isomorphic.
Before that,  Jiang {{and}} Su (\cite{JS}) constructed a unital simple ASH-algebra ${\cal Z}$ of stable rank one which has the same Elliott invariant as that of $\C.$ In particular,  ${\cal Z}$ has no non-trivial projections.

If $A$ is a {{simple}} separable
amenable  \CA\, with weakly unperforated $K_0(A)$ which belongs to a {{(reasonable)}} classifiable class, then
one would expect  that $A$ must be  isomorphic to $A\otimes {\cal Z},$
 since
these two algebras have the same Elliott invariant  {{(Theorem 1 of \cite{GJS}).}}
{{If $A$ is isomorphic to $A\otimes {\cal Z},$ then  $A$ is}} called ${\cal Z}$-stable. The existence of non-elementary simple \CA s which are not ${\cal Z}$-stable was first proved  by
Gong, Jiang {{and}} Su (see \cite{GJS}).  {{Toms's}} counterexample is in particular not ${\cal Z}$-stable. Thus, ${\cal Z}$-stability
should be added {{to}}
 the
hypotheses if one uses the conventional Elliott invariant.  {{(The class of simple AH-algebras of \cite{EGL-AH}
are known to be ${\cal Z}$-stable.  In fact,
all unital separable simple amenable \CA s with finite tracial rank are
${\cal Z}$-stable; see Corollary 8.4 of \cite{Lin-LAH}.)}}


The next development in this direction came from a new approach {{due to}}  Winter, who made use of the assumption
{{of}} ${\cal Z}$-stability  in a  remarkably  innovative way (\cite{Winter-Z}). His idea was to view $A\otimes {\cal Z}$ as an
inductive limit of algebras of paths in $A\otimes Q$ with endpoints in $A\otimes M_{\bf p}$ and $A\otimes M_{\bf q}$ (where $Q$ is the UHF-algebra {{with}}  $K_0(Q)=\Q,$ ${\bf p}$ and ${\bf q}$ are coprime supernatural numbers, and $M_{\bf p}$ and $M_{\bf q}$ their associated UHF algebras).
Suppose that the endpoint algebras are classifiable.
Winter  showed that  if somehow there is a continuous path of  isomorphisms
from one endpoint to the other,
then   the algebra $A$ itself is
also classifiable.

Winter's procedure  provided a new framework to carry out classification. {{H}}owever, to actually execute the continuation from endpoint to endpoint alluded to above, one
{{needs}} new types of uniqueness and existence theorems.
 In other words,  just like before, in Elliott's intertwining argument,  the new procedure  ultimately, but not surprisingly,
 depends on {{certain}} existence and uniqueness theorems. However, this time we need the uniqueness and existence theorems
  with respect to {{(one-parameter)}} asymptotic
unitary equivalence of the maps involved rather {{than  just (sequential)}} approximate
 unitary equivalence.  {{This}} is significantly more demanding. For example, in the case of the existence theorem, we need to construct  {{a map}} which {{lifts}} a prescribed $KK$-element
rather than {{just}} a $KL$-element.  {{This}} was once thought  to be out of reach for general stably finite algebras since the $KK$-functor does not preserve inductive limits {{(as the $KL$-functor does; see \cite{DL}).}} It was an unexpected usage of the Basic Homotopy Lemma that made this possible.  Moreover, the existence theorem also needs to respect
{{a}} prescribed rotation related map.  The existence
 theorems
are very different  {{from
those developed}} in the early study
{{of the subject.}}
Inevitably,  the uniqueness theorem also
becomes more complicated (again the Basic Homotopy Lemma plays the key role).

Once  we overcame these new hurdles {{arising  in following the Winter approach,}}
{{we were able to  classify}} the class of all {{unital separable simple}} ${\cal Z}$-stable \CA s $A$
{{which satisfy the UCT}} and whose tensor products with all UHF-algebras of infinite type are of tracial rank zero,
by means of the Elliott invariant,  in
 \cite{L-N}.
 We were then able in \cite{Lnclasn} to extend this result to the class ${\cal A}$
of all unital simple separable amenable
 ${\cal Z}$-stable  \CA s which satisfy the UCT and
 have  tracial rank  one (not just zero) after tensoring {{with (just)}} some
infinite dimensional UHF-algebra.
The class ${\cal A},$ and already for that
matter just the subclass mentioned above, the \CA s
rationally of tracial rank zero, carry classification
substantially beyond earlier results.
For  instance,
the Jiang-Su algebra ${\cal Z}$ is {{unital}} projectionless, but ${\cal Z}\otimes U\cong U$ for any infinite dimensional UHF-algebra $U.$
In fact, the class ${\cal A}$ exhausts all  Elliott invariants with simple weakly unperforated
 rationally Riesz groups as
ordered $K_0$-group and pairing with the trace simplex
taking extreme traces to extreme $K_0$-states
(\cite{LNjfa}).

The class ${\cal A}$ not only contains all unital simple separable amenable
\CA s with  tracial rank one  in the UCT class,  and  the Jiang-Su algebra, but also contains many other
simple \CA s.
In fact it unifies the previously classified classes such as the
simple limits of dimension drop {{interval}} algebras  {{and}} dimension
drop circle algebras which, {{like the Jiang-Su algebra, are not}}
AH-algebras (\cite{Lnjfa2010}).
However, the restriction on the pairing between traces
and
$K_0$
 prevents the class ${\cal A}$ from including
 inductive limits of
``point--line" algebras,
 which we called Elliott-Thomsen building blocks.
 {{This}} brings  us to the main goal of our {{work}}.

The goal of the present article, {{and its sequel (Part II),}} 
is to give a classification of a {{new}} class of unital simple separable amenable \CA s satisfying the UCT, by means of the Elliott invariant. This class is significant because it exhausts all possible values of the Elliott invariant for all unital, simple,
 finite, ${\cal Z}$-stable, separable amenable \CA s. { (}{{It strictly}} contains the class ${\cal A}$  {{mentioned above.}}{ )}

{{First}, {{we shall}} introduce a class of {{unital}} simple  separable \CA s
which we shall refer to as the
\CA s of
\emph{generalized} tracial rank one.
The definition {{is in}} the same spirit as that of tracial rank one, but, instead of using only
matrix algebras of continuous functions on a one-dimensional finite CW complex, it uses
{{the point--line algebras.}}  These were first introduced into the Elliott program by Elliott and Thomsen (\cite{ET-PL}),
{{in connection with determining
the range of the Elliott invariant.}}  Some time later,  these \CA s were also called
one dimensional non-commutative  CW complexes (NCCW).
This class of unital simple separable \CA s (of
generalized tracial rank one) will be denoted by ${\cal B}_1.$
If we insist that the point--line
algebras used in the definition
have trivial $K_1,$ then the resulted sub-class will be denoted by ${\cal B}_0.$
{{Amenable \CA s in the class ${\cal B}_1$ are proved here to be Jiang-Su stable.}}

Denote by ${\cal N}_1$  the family of unital simple separable amenable  \CA s  $A$ which satisfy the UCT such that $A\otimes Q\in {\cal B}_1,$  and by ${\cal N}_0$ the subclass of those \CA s $A$ such that
$A\otimes Q\in {\cal B}_0,$ where $Q$ is the UHF-algebra with $K_0(Q)=\Q$. 


As earlier for \CA s of tracial rank one,
we shall expand the new class of unital simple separable
\CA s of generalized tracial rank one (the class ${\cal B}_1$)
to include those ${\cal Z}$-stable
\CA s such that the tensor product with some
infinite-dimensional UHF algebra belongs to this class (${\cal B}_1$).
We shall show, in the present Part I
of this work, that this expanded new class, the class of ${\cal Z}$-stable unital
simple separable amenable \CA s rationally of generalized tracial
rank one,
exhausts the Elliott
invariant for finite ${\cal Z}$-stable unital simple separable
\CA s.
We shall also prove, in the present Part I of this work,
that, if $A$ and $B$ are amenable \CA s in ${\cal B}_0$
satisfying the UCT, then $A\otimes U\cong B\otimes U$
for some UHF-algebra $U$ of infinite type if, and only if, $A\otimes U$ and $B\otimes U$
have isomorphic Elliott invariants.
{{In other words, we classify a certain sub-class of \CA s of generalized tracial rank one.
In Part II, we shall classify the class of all \CA s rationally of generalized  tracial rank one.}}

The present part of the paper, Part I, is organized as
follows. Section 2 serves as preliminaries and
establishes some conventions.  In
Section 3, we study the class of unital
Elliott-Thomsen building blocks, denoted by ${\cal C}$ (see [26]
and \cite{point-line}). Elliott-Thomsen building blocks are also called point--line algebras, or one dimensional non-commutative CW complexes (NCCW complexes,
studied in \cite{ELP1} and \cite{ELP2}).
Sections 4 and 5 discuss the uniqueness theorem for maps from \CA s in ${\cal C}$ to finite dimensional
\CA s.  Section 8 presents a uniqueness theorem for maps from a $C^*$-algebra in ${\cal C}$ to another $C^*$-algebra in ${\cal C}.$ This is done by using
a homotopy lemma established in  Section 6 and existence theorems established in Section 7 to bridge the uniqueness theorems of  Sections 4 and 5
{{with those of}} Section 8.  In Section 9, the classes ${\cal B}_1$  {{(as above)}} and ${\cal B}_0$ {{(as above)}} are introduced.
Properties of \CA s in the
class ${\cal B}_1$ are discussed in Sections 9, 10,
and 11.
These
\CA s, unital separable simple \CA s of
generalized tracial rank (at most) one, can also be characterized as
being  tracially
 approximable by
(arbitrary)
subhomogeneous \CA s with one-dimensional
spectrum.
{{For example, we show, in Section 9, that ${\cal B}_1$ and ${\cal B}_0$ are not the same (unlike
the previous case, in which unital simple separable \CA s of tracial rank one are TAI), and}}
in Section 10,   {{ that amenable}} \CA s in ${\cal B}_1$ are ${\cal Z}$-stable.
Section 12 is dedicated to the main
uniqueness theorem used in the isomorphism theorem of
Section 21.


Sections 13 and 14 are devoted to
the range theorem, one of the main results: given any
possible Elliott invariant {{sextuple
for a ${\cal Z}$-stable  unital simple separable amenable \CA,
namely, any countable weakly unperforated simple order-unit
abelian group, paired with an arbitrary metrizable Choquet
simplex mapping onto the state space of the order-unit group,
together with an arbitrary countable abelian group,  there is
a \CA\, in the class ${\cal N}_1$
realizing this sextuple as its Elliott invariant.
(The construction of this model, an inductive limit of
subhomogeneous \CA s, is similar to that in \cite{point-line},
but additional work is needed to show that the inductive limits
in question belong to the class ${\cal N}_1$.)}}
{{This is shown in Section 13;
Section 14 gives a similar construction, for a restricted class
of invariants, and yielding a correspondingly restricted class
of (inductive limit) algebras. For reasons that will become
clear later, this second, restricted, model construction
is very important.}}

Sections  15 to 19 could all be
described as different stages in the development of the
existence theorem, to be used in the later sections as well as final isomorphism
theorem (in Part II).
These deal with
the issue of existence  for maps from \CA s in ${\cal C}$ to finite dimensional \CA s and then to
\CA s in ${\cal C}$ {{that}}
match prescribed $K_0$-maps and tracial information. {{The}} ordered {{$K_0$-}}structure and combined simplex information
of these \CA s
become complicated.  We also need to consider maps
from homogeneous \CA s to \CA s in ${\cal C}.$ The mixture with higher dimensional
{{noncommutative}} CW complexes does not ease the difficulties. However, in Section 18, we show that,  at least under certain restrictions, any given compatible triple
which
consists of a strictly positive $KL$-element,
a map on the tracial state space, and a \hm\, on a quotient of the unitary group,  it is possible to construct a
\hm\,  {{from a separable amenable \CA\, $A$ satisfying the UCT of the form
$B\otimes U$ for some $B\in {\cal B}_0$ and some UHF-algebra $U$ of infinite type to another
separable \CA\, $C$ of the form $D\otimes V,$ where $D\in {\cal B}_0$ and $V$ is a  UHF-algebra
of infinite type,}}
 which matches the triple. Variations of this are also discussed.
In Section 19, we show that ${\cal N}_1={\cal N}_0,$ {{even though ${\cal B}_1\not={\cal B}_0,$ as indicated in
Section 9.}}
In Section 20, we continue to study the existence theorem.
In Section 21, we show that any  unital simple \CA, which satisfies the UCT,
in the class {{${\cal N}_1$ absorbing tensorially a UHF algebra
of infinite type}} is isomorphic to an inductive limit
$C^*$-algebra  as constructed in Section 14, and
any two such \CA s are isomorphic if they have the
same Elliott invariant.
This isomorphism theorem is special but  is also    the foundation of our main
isomorphism theorem in  Part II of the paper. 



{\bf Acknowledgements}:
A large part of this article {{was}} written during the summers  of 2012, 2013, and 2014
when all three authors visited The Research Center for Operator Algebras in East China Normal University
 which is in part supported  by NNSF of China (11531003)  and Shanghai Science and Technology
 Commission (13dz2260400)
 and  Shanghai Key Laboratory of PMMP.
They were
partially
supported by the {{Center}} (also during summer 2017 when
some of revision were made).  Both the first named author and  the second named author were supported partially by NSF grants {(DMS 1665183 and DMS 1954600).}
 The work of the third named author was partially supported by a NSERC Discovery Grant, a Start-Up Grant from the University of Wyoming, a Simons Foundation Collaboration Grant, and a NSF grant (DMS-1800882).
 {{The}} authors would like to thank Michael Yuan Sun {{for reading}}
  part of
an earlier version of this {{article}}
 and
for his comments.
{\Green{Since the article was first posted at the beginning of 2015,
it went through many revisions. Over a dozen of referees made serious checking,
numerous corrections, suggestions, and comments over the years. The authors
made several rounds of corrections following those valuable corrections, suggestions and comments.
The authors would like to take this opportunity to express their most sincere gratitude to  the referees for their efforts.}}
The author would also like to thank the editor for checking lengthy referee's reports and corrections which
leads to the publication.

 \section{Notation and  Preliminaries}

This section includes a list of notations and definitions most of which are standard. We list them  here for
the reader's convenience. However,  we also recommend {skipping}
this section until some of these
notions appear.

\begin{df}\label{Du}
{\rm
Let $A$ be a unital \CA.
{{Denote by $A_{s.a.}$ the self-adjoint part of $A$
and $A_+$ the set of all positive elements of $A.$}}
Denote by $U(A)$\index{$U(A)$} the unitary group of $A$, and denote by  $U_0(A)$ the normal subgroup of $U(A)$ consisting of those unitaries which are in the  connected component  of $U(A)$ containing $1_A.$  Denote by
$DU(A)$ the commutator subgroup of $U_0(A)$\index{$U_0(A)$} and $CU(A)$\index{$CU(A)$} the closure of $DU(A)$ in $U(A).$}
\end{df}

\begin{df}\label{Aq}
{\rm Let $A$ be a unital \CA\, and let $T(A)$ denote the simplex of  tracial states of $A,$
{\blue{a compact subset of $A^*,$ the dual of $A,$ with the weak* topology (see also  \cite{Thoma} and II.4.4 of \cite{BH}).}}
Let $\tau\in T(A).$ We say that $\tau$ is faithful if $\tau(a)>0$ for all $a\in A_+\setminus\{0\}.$
Denote by $T_f(A)$ the set of all faithful tracial states.\index{$T_f(A)$}

{For each integer $n\ge 1$ and $a\in M_n(A),$
write $\tau(a)=(\tau\otimes \mathrm{Tr})(a),$ where $\mathrm{Tr}$ is the (non-normalized) standard trace on $M_n.$}

Let $S$ be a compact convex set.
Denote by $\Aff(S)$ the space of all real continuous affine functions on $S$ and
denote by ${\rm LAff}_b(S))$  the set of all bounded lower semi-continuous real affine functions on $S.$
Denote by $\Aff(S)_+$ the set of those non-negative valued functions in $\Aff(S)$ and
$\Aff(S)^+=\Aff(S)_+\setminus \{0\}.$\index{$\Aff(S)$}\index{$\Aff(S)_+$} {\blue{Also define $\Aff(S)^{++}=\{f\in \Aff(A)^{+}: f(s)>0,\,s\in S\}.$
Define ${\rm LAff}_b(S)_+$ to be the set of those non-negative valued functions
in ${\rm LAff}_b(S),$ and ${\rm LAff}_b(S)^{++}$ to be the set $\{f\in {\rm LAff}_b(S)^+: f(s)>0,\, s\in S\}.$}}
\index{${\rm LAff}_b(S)_+$}

Suppose that $T(A)\not=\emptyset.$ There is a linear   map
$r_{\aff}: A_{s.a.}\to \Aff(T(A))$ defined by
$$
r_{\aff}(a)(\tau)=\hat{a}(\tau)=\tau(a)\tforal \tau\in T(A)
$$
and for all $a\in A_{s.a.}.$

{{Let $A_0$ denote the  closure of the set of all self--adjoint elements of $A$  of the form $\sum x_ix_i^*-\sum x_i^*x_i$. Then by Theorem 2.9  of \cite{CP} (also see the proof of Lemma 3.1 of \cite{Thomsen-rims} for further explanation), we know that $A_0=\ker (r_{\aff}).$ {\index{$A_0$}}  As in \cite{CP}, denote by $A^q$ the quotient space $A_{s.a.}/A_0.$
It is a real Banach space. Denote by
$A_+^q$  the image of $A_+$ in $A/A_0.$ \index{$A_+^q$} }}
{{
Denote by $q: A_{s.a.}\to A^q$ the quotient map.
It follows from Proposition 2.7 of \cite{CP} that $T(A)$ is precisely the set of those real bounded linear functionals $f$ on $ A^q$
such that $f(x)\ge 0$ for all $x\in A_+^q$ and $f(q(1_A))=1.$ Moreover, the topology on $T(A)$ is
the weak* topology of $T(A)$ as a subset of the dual space of $A^q.$}}
}
{{For a given 
element $g\in \Aff(T(A))$, 
by  Proposition  2.8 of \cite{CP}, $g$ {{can be}} uniquely extended to a bounded linear functional on
$(A^q)^*,$ the dual space of $A^q.$ Since $g$ is continuous on $T(A),$ $g$ is weak*-continuous.
Therefore this  gives an element $\Gamma(g)$ in $A^q$ (this also shows that the map $r_{\aff}$ is surjective).
Let $\tau\in T(A).$ Then we view $\tau$ as an element in $(A^q)^*$ as above.
Then $\tau(\Gamma(g))=g(\tau)$ for all $\tau\in T(A).$ It follows that  $\Gamma$ is a  linear map
from $\Aff(T(A))$ to $A^q.$  It is clear that $\Gamma$ is injective. To see  it is surjective,
let $x\in A^q.$ Then, viewing $T(A)$ as a subset of the dual of $A^q,$ $\hat{x}(\tau)=\tau(x)$  (for all
$\tau\in T(A)$) defines an element {{of}} $\Aff(T(A)).$ It is clear that $\Gamma(\hat{x})=x.$
That the  map $\Gamma$ is an isometry now follows from the equation (3.1) of the proof of Lemma 3.1 of \cite{Thomsen-rims}
(see also Theorem 2.9 of \cite{CP}). (Hence $\Gamma$ is also positive.)}}

{{Suppose that $A$ and $B$ {{are}} two unital separable \CA s such that
$\Aff(T(A))\cong \Aff(T(B)),$ i.e., there is an  isometric  order isomorphism  $\gamma$  from
the real Banach space $\Aff(T(A))$ onto the real Banach space $\Aff(T(B))$  which
preserves the constant function 1. Then there is an isometric isomorphism from $A^q$ onto $B^q$  which maps
$A_+^q$ into $B_+^q$ and $q(1_A)$ to $q(1_B),$ and the inverse maps $B_+^q$ into $A_+^q$ and
maps $q(1_B)$ to $q(1_A).$ Since we have identified $T(A)$ and $T(B)$ with
the subset of $(A^q)^*$ which preserves the order and has value 1 on $q(1_A)$ and the subset
of $(B^q)^*$ which preserves the order  and has value 1 on $q(1_B),$ respectively,
the map $\gamma$ induces an affine homeomorphism from $T(B)$ onto  $T(A).$}}


\end{df}

\begin{df}\label{Drho}
Let $A$ be a unital stably finite \CA\, with $T(A)\not=\emptyset.$ Denote by $\rho_A: K_0(A)\to \Aff(T(A))$ the  order preserving \hm\,  defined by $\rho_A([p]){{(\tau)}}=\tau(p)$ for any projection $p\in M_n(A),$
$n=1,2,...$ (see the convention above).\index{$\rho_A$}

{{A}} map $s: K_0(A)\to \R$ is said to be a state if $s$ is an order preserving \hm\, such that
$s([1_A])=1.$ The set of states on $K_0(A)$ is denoted by $S_{[1_A]}(K_0(A)).$

{Denote by $r_A: T(A)\to S_{[1_A]}(K_0(A))$ the map defined by $r_A(\tau)([p])=\tau(p)$ for all projections
$p\in M_n(A)$ (for all integers $n$) and for all $\tau\in T(A).$ } \index{$r_A$}
\end{df}

{
\begin{df}\label{DEll}
Let $A$ be a unital simple \CA. The Elliott invariant of $A$, {{denoted}} by ${\rm Ell}(A),$ is the
sextuple\index{${\rm Ell}(A)$}
$$
{\rm Ell}(A)=(K_0(A), K_0(A)_+, [1_A], K_1(A), T(A), r_A).
$$
Suppose that $B$ is another unital simple \CA.
{{We say that $\Gamma: {\rm Ell}(A)\to {\rm Ell}(B)$ is a \hm\,   if
there are an order preserving  \hm\, $\kappa_0: K_0(A)\to K_0(B)$
such that $\kappa_0([1_A])=[1_B],$ a \hm\,
$\kappa_1:
K_1(A)\to K_1(B),$ and {{a
continuous affine map}}
$\kappa_T: T(B)\to T(A)$ such that
\beq\label{gong18-241}
r_A(\kappa_T(t))(x)=r_B(t)(\kappa_0(x))\rforal x\in K_0(A)\andeqn \rforal t\in T(B),
\eneq
and we write $\Gamma=(\kappa_0, \kappa_1,\kappa_T).$}}

We write  ${\rm Ell}(A)\cong {\rm Ell}(B)$ if there is $\Gamma$ as above such that
$\kappa_0$ is an order isomorphism
such that $\kappa_0([1_A])=[1_B],$
$\kappa_1$ is an isomorphism, and
$\kappa_T$ is an affine homeomorphism.
{{If, in addition, $A$ is separable and satisfies the UCT, then there exists an element
$\af\in KL(A, B)$ such that $\af|_{K_i(A)}=\kappa_i,$ $i=0,1$ (recall
that, by \cite{Ror-KL-I}, in this case, $KL(A,B)=KK(A,B)/Pext,$ where
${{Pext}}$ is the subgroup corresponding to the pure extensions of  $K_*(A)$ by $K_*(B)$).
If $B$ is also separable and satisfies the UCT,
then there also exists $\af^{-1}\in KL(B,A)$ with $\af^{-1}\times \af=[{\rm id}_A]$ and $\af\times \af^{-1}=[{\rm id}_B]$
(see 23.10.1  of \cite{Bla-Ktheory}).}}

{\blue{Any continuous affine map $\kappa_T: T(B)\to T(A)$ induces a map $\kappa^{\sharp}_T: \Aff(T(A)) \to \Aff(T(B))$ defined by $\kappa^{\sharp}_T(l)(\tau)=l(\kappa(\tau))$ for all $\tau \in T(B)$ and $l\in \Aff(T(A))$. Furthermore,  $\kappa_T$ is compatible with $\kappa_0$ in the sense of \eqref{gong18-241}
if and only if $\kappa^{\sharp}_T$ and $\kappa_0$ are compatible in the following sense:
$$\rho_B(\kappa_0(x))= \kappa^{\sharp}_T(\rho_A(x))\rforal x\in K_0(A).$$
Note  that any unital homomorphism $\phi: A\to B$ induces maps $\phi_{*,0}: K_0(A) \to K_0(B)$ and $\phi^\#:\Aff(T(A))\to \Aff(T(B))$, which are compatible. }}

\end{df}
}

\begin{df}\label{Dball}
Let $X$ be a compact metric space, let $x\in X$ be a point, and let $r>0.$
Denote by $B(x, r)$ the open unit ball $\{y\in X: {\rm dist}(x, y)<r\}.$\index{$B(x,r)$}

Let $\ep>0.$ Define $f_\ep\in C_0((0,\infty))$ to be the function with
 $f_\ep(t)=0$ if $t\in [0,\ep/2],$ $f_\ep(t)=1$ if $t\in [\ep,\infty),$ and
$f_\ep(t)=(2-\ep)/\ep$
if $t\in [\ep/2, \ep].$ Note that $0\le f\le 1$ and $f_\ep f_{\ep/2}=f_\ep.$

{{Denote by $t_+$ the continuous function
 $g(t):=(1/2)(t+|t|)$ for all $t\in \R.$ Then, if $A$ is a \CA\, and $a\in A_+,$ the element $a_+:=g(a)=(1/2)(a+|a|)$  is
 the positive part of $a.$ }}
\end{df}

\begin{df}\label{DW(A)}
Let $A$ be a \CA. Let $a, b\in M_n(A)_+$. Following Cuntz (\cite{Cu1}), we write
$a\lesssim b$\index{$\lesssim$} if there exists a sequence $\{x_n\}$ in $M_n(A)$ such that
$\lim_{n\to\infty} x_n^*bx_n=a.$ If $a\lesssim b$ and $b\lesssim a,$ then
we write $a\sim b.$ The relation ``$\sim$" is an equivalence relation.
Denote by $W(A)$\index{$W(A)$} the Cuntz {{semigroup,}} consisting of the equivalence classes of positive\index{Cuntz semigroup}
elements in $\bigcup_{m=1}^{\infty} M_m(A)$ with orthogonal
addition (i.e., $[a+b]=[a\oplus b]$).

If $p, q\in M_n(A)$ are projections, then $p\lesssim q$ if and only if $p$ is
Murray-von Neumann equivalent to a subprojection of $q$. In particular,  when $A$ is stably finite,  $p\sim q$ if and only if
$p$ and $q$  are {{Murray}}-von Neumann equivalent.

{\blue{Recall (see  II. 1.1 of \cite{BH}) that a  (normalized) 2-quasi-trace of a $C^*$-algebra $A$ is a function $\tau: A \to \mathbb{C}$ satisfting
\begin{enumerate}
\item $\tau(1)=1$,
\item $0\leq \tau(xx^*) = \tau(x^*x)$, $x\in A$,
\item $\tau(a+ib) = \tau(a) + i\tau(b)$, $a, b \in A_{s.a.}$,
\item $\tau$ is linear on abelian \SCA\,
of $A$, and
\item $\tau$ extends to a function from $M_2(A)$ to $\mathbb C$ satisfying the  conditions above.
\end{enumerate}}}

Denote by $QT(A)$ the set of normalized {{$2$-quasi-traces}} on $A.$\index{$QT(A)$}
{{It follows from II 4.1 of \cite{BH} that every $2$-quasi-trace $\tau$  extends  to a quasi-trace on $M_n(A)$
(for all $n\ge 1$).}}
For $a\in A_+$ and $\tau\in QT(A),$ define
$$
d_{\tau}(a)=\lim_{\ep\to 0} \tau(f_{\ep}(a)).
$$\index{$d_\tau$}
Suppose that $QT(A)\not=\emptyset.$
We say $A$ has strict comparison for positive elements if, for any $a, b\in M_n(A)$ (for all {{integers}} $n\ge 1$),
$d_\tau(a)<d_\tau(b)$ for $\tau\in QT(M_n(A))$  implies $a\lesssim b.$
\end{df}


\begin{df}\label{Aq1}
Let $A$ be a \CA.  Denote by $A^{\bf 1}$ the unit ball of $A$, {{and by}}
$A_+^{q, {\bf 1}}$ the image of the intersection of $A_+\cap A^{{\bf 1}}$ in $A_+^q.$\index{$A^{\bf 1}$}\index{$A_+^{q, {\bf 1}}$}  {(Recall that $A_+^q=r_{\aff}(A_+)$; see Definition \ref{Aq}.)}
\end{df}

\begin{df}\label{dInn}
Let $A$ be a unital \CA\,  and let $u\in U(A).$ We write
${\rm Ad}\, u$ for the automorphism $a\mapsto u^*au$ for all $a\in A.$\index{${\rm Ad}\,u$}
Suppose $B\subset A$ is a unital $C^*$-subalgebra.
Denote by $\overline{{\rm Inn}}(B,A)$ the set of all those monomorphisms $\phi: B\to A$ such that
there exists a sequence of unitaries $\{u_n\}\subset {{A}}$ with
$\phi(b)=\lim_{n\to\infty} u_n^*bu_n$ for all $b\in B.$
\end{df}

\begin{df}\label{dfcalN}
Denote by ${\cal N}$ the class of separable amenable \CA s which satisfy the Universal Coefficient Theorem (UCT).\index{${\cal N}$}

Denote by ${\cal Z}$ the Jiang-Su algebra (\cite{JS}). 
Note that $\mathcal Z$ has a unique trace {{state}} and $K_i({\cal Z})=K_i(\C)$ ($i=0,1$).\index{${\cal Z}$}\index{Jiang-Su algebra}
A \CA\, $A$ is said to be ${\cal Z}$-stable if $A\cong A\otimes {\cal Z}.$\index{${\cal Z}$-stable}
\end{df}

\begin{df}\label{DKLtriple}
{\rm

Let $A$ be a unital \CA\,. Recall that, following D\u{a}d\u{a}rlat and Loring (\cite{DL}), one defines\index{$\underline{K}(A)$}
\begin{equation}\label{Dbeta-5}
\underline{K}(A)=\bigoplus_{i=0,1}K_i(A)\oplus\bigoplus_{i=0,1}
\bigoplus_{k\ge 2}K_i(A,\Z/k\Z).
\end{equation}
There is a commutative \CA\, $C_k$ such that one may identify
$K_i(A\otimes C_k)$ with $K_i(A, \Z/k\Z).$
Let $A$ be a unital separable amenable \CA\,, and let $B$ be a $\sigma$-unital \CA.
Following R\o rdam (\cite{Ror-KL-I}),  $KL(A,B)$ is the quotient of $KK(A,B)$ by those elements
represented by limits of trivial extensions (see \cite{LnAUCT}). In the case that $A$ satisfies the UCT,
R\o rdam defines $KL(A,B)=KK(A,B)/{{Pext}},$ where
${{Pext}}$ is the subgroup corresponding to the pure extensions of  $K_*(A)$ by $K_*(B).$
In \cite{DL}, D\u{a}d\u{a}rlat and Loring  proved that\index{$KL(A,B)$} {{(if $A$  satisfies the UCT)}}
\beq\label{DKL-2}
KL(A,B)={\rm Hom}_{\Lambda}(\underline{K}(A), \underline{K}(B))
\eneq\index{${\rm Hom}_{\Lambda}(\underline{K}(A), \underline{K}(B))$}
{{(see page 362 of \cite{DL} for the definition of the group of module \hm s ${\rm Hom}_{\Lambda}(-,-)$).}}

{
Now suppose that $A$ is stably finite. Denote by
$KK(A,B)^{++}$\index{$KK(A,B)^{++}$} the set of those elements\index{$KK(A,B)^{++}$}
$\kappa\in KK(A,B)$
such that ${\Green{\kappa(K_0(A)_+)}}\setminus\{0\})\subset K_0(B)_+\setminus\{0\}.$ {\Green{(Warning: the notation here may be  different from other papers.)}}
{{In the absence of the UCT, we denote by ${\rm Hom}_{\Lambda}(\underline{K}(A), \underline{K}(B))^{++}$
\index{${\rm Hom}_{\Lambda}(\underline{K}(A), \underline{K}(B))^{++}$}
the set of those elements $\kappa\in {\rm Hom}_{\Lambda}(\underline{K}(A), \underline{K}(B))$
such that $\kappa(K_0(A)\setminus \{0\})\subset K_0(B)_+\setminus \{0\}.$}}
Suppose further that both $A$ and $B$ are unital. Denote by $KK_e(A,B)^{++}$\index{$KK_e(A,B)^{++}$} the subset of
those $\kappa\in KK(A,B)^{++}$ such that $\kappa([1_A])=[1_B].$
Denote by $KL_e(A,B)^{++}$\index{$KL_e(A,B)^{++}$} the image of $KK_e(A,B)^{++}$ in $KL(A,B).$ }

}
\end{df}

\begin{df}\label{Dmap}
Let $A$  and $B$ be  \CA s and  $\phi: A\to B$ be a linear map.
{{We will {{sometimes}}, without notice, continue to use $\phi$ for the induced
map $\phi\otimes {\rm id}_{M_n}: A\otimes M_n\to B\otimes M_n.$
Also, $\phi\otimes 1_{M_n}: A\to B\otimes M_n$ is  used for the amplification
which maps $a$ to ${{\phi(a)}}\otimes 1_{M_n},$ the diagonal element with ${{\phi(a)}}$ repeated $n$ times.}}
Throughout the paper,  if $\phi$ is a \hm, we will use $\phi_{*i}: K_i(A)\to K_i(B)$, $i=0, 1$, for the induced \hm.\index{$\phi_{*i}$}
We will use $[\phi]$ for the element of $KL(A,B)$ (or $KK(A,B)$ if there is no confusion) which is induced by $\phi.$
Suppose that $A$ and $B$ are unital and $\phi(1_A)=1_B.$
Then $\phi$ induces
an affine map $\phi_T: T(B)\to T(A)$\index{$\phi_T$} defined by $\phi_T(\tau)(a)=\tau(\phi(a))$ for all $\tau\in T(B)$ and
$a\in A_{s.a.}.$ Denote by $\phi^{\sharp}: \Aff(T(A))\to \Aff(T(B))$\index{$\phi^{\sharp}$} the affine continuous map defined by
$\phi^{\sharp}(f)(\tau)=f(\phi_T(\tau))$ for all $f\in \Aff(T(A))$ and $\tau\in T(B).$

\end{df}


\begin{df}\label{KLtriple}
Let $A$ be a unital separable amenable \CA\, and let  $x\in A.$ Suppose
that
$\|xx^*-1\|<1$ and $\|x^*x-1\|<1.$ Then $x|x|^{-1}$ is a unitary.
Let us use $\langle x \rangle $ to denote $x|x|^{-1}.$\index{$\la x\ra$}

Let ${\cal F}\subset A$ be a finite subset and $\ep>0$ be a positive number.
We may assume that $1_A\in {\cal F}.$
We say a map $L: A\to B$ is ${\cal F}$-$\ep$-multiplicative if
$$
\|L(xy)-L(x)L(y)\|<\ep\rforal x,\, y\in {\cal F}.
$$

Let ${\cal
P}\subset \underline{K}(A)$ be a finite subset.   {Let us first assume that $\mathcal P\subset K_0(A) \oplus K_1(A)$
}
{{Assume also that $\{p_i, p'_i: 1\le i\le m_0\}\subset M_N(A)$ {\blue{is}} a finite subset of projections and
$\{u_j: 1\le j\le m_1\}\subset M_N(A)$ {\blue{is}} a finite subset of unitaries
such that $\{[p_i]-[p_i'], [u_j]: 1\le i\le m_0, 1\le j\le m_1\}={\cal P}.$}}
{{Then there}} is $\ep>0$ and a finite subset
${\cal F}$ of $A$ satisfying the following condition: for any unital \CA\, $B$ and
any unital ${\cal F}$-$\ep$-multiplicative \cp\,
$L : A\to B,$ the map $L$
induces a \hm\, $[L]$ defined on  $G({\cal P}),$ where $G({\cal
P})$ is the subgroup generated by ${\cal P},$ to $\underline{K}(B)$
such that {{there are projections $q_i, q_i'\in M_N(B)$ with
$[q_i]=[L]([p_i]),$ $[q_i']=[L]([p_i'])$ in $K_0(B)$  ($1\le i\le m_0$) and  unitaries $u_j\in M_N(B)$
($1\le j\le m_1$)  with $[v_j]=[L]([u_j])$}}  such that \index{$[L]$}
\beq\label{KLtriple-1}
&&{{\|L(p_i)-q_i\|<1/2, \,\,\|L(p'_i)-q_i'\|<1/2\,\, (1\le i\le m_0) \andeqn}}\\
&& {{\|\langle L(u_j)\rangle -v_j\|<1/2\,\, (1\le j\le m_1).}}
\eneq
In general, ${\cal P}\cap K_i(A, \Z/k\Z)\not=\emptyset.$
Then the above  also applies to ${\cal P}\cap
K_i(A, \Z/k\Z)$ with a necessary modification, by replacing $L,$ by
$L\otimes {\rm id}_{C_k},$ where $C_k$ is the commutative \CA\,
referred to {{in}} \ref{DKLtriple}.
{\blue{Suppose that the triple $(\ep, {\cal F}, {\cal P})$ also has the following property:
if $L'$ is another such map
with the property
that  $\|L(a)-L'(a)\|<2\ep\rforal a\in {\cal F},$ then $[L]|_{\cal P}=[L']|_{\cal P}.$
Then}} such a triple $(\ep, {\cal F}, {\cal P})$ may be  called a
$KL$-triple for $A$ {{(see, for  example,  1.2  of \cite{Lnamj98} and  3.3 of \cite{DE}).}}
{{Note that these considerations, in particular, imply that, if $u_j\in U_0(A),$
then $[L]([u_j])\in U_0(B).$}}

{\blue{Suppose that $A$ is unital  and $L$ is  a \morp\,  which is  not unital.  We may always assume that $1_A\in {\cal F}.$
When $\ep<1/4,$ let $p=\chi(L(1_A)),$ where $0\le \chi(t)\le 1$ is a function  in  $C([0,1])$ which is
zero on $[0,1/3]$ and $1$ on $[1/2, 1].$ For small $\ep,$ $pL(1_A)p$ is close to $p$ and so invertible in $pBp.$
Let $b\in pBp$ be the inverse of $pL(1_A)p$ in $B.$  Then $\|b^{1/2}-p\|<4\ep.$
Define $L': A\to pBp$  by $L'(a)=b^{1/2}pL(a)pb^{1/2}$
for $a\in A.$ Note that, for any $\eta>0,$  $\|L-L'\|<\eta$  if $\ep$ is sufficiently small.
The convention of this article, as usual, will be that we may always assume that $L(1_A)$ is a projection when
we mention an ${\cal F}$-$\ep$-multiplicative map $L.$

Then, if $u$ is a unitary, with sufficiently large ${\cal G}$ and small $\ep,$ as above,  then $\la L(u)\ra$ is a unitary in $pBp.$
By  $[L]([u]),$ in $K_1(B),$ we then mean $[\la L(u)\ra +(1-p)].$}}

Suppose that $K_i(A)$ is finitely generated.
Then,
by  {{Proposition 2.4 of \cite{LnHomtp},}}
for some large ${\cal P},$  if $(\ep, {\cal F}, {\cal P})$ (with sufficiently small
$\ep,$ and sufficiently large ${\cal F}$) is a $KL$-triple for $A,$\index{$KL$-triple}
 then $[L]$ defines an element in $KL(A,B)=KK(A,B).$ In this case, we say $(\ep, {\cal F})$ is
 a $KK$-pair.\index{$KK$-pair}
 \end{df}
 \begin{lem}[{\rm Lemma 2.8 of \cite{LnHomtp}}]\label{oldnuclearity}
{{Let $A$ be a unital  {{ separable amenable}} \CA. Let $\ep>0,$ let
${\cal F}_0\subset A$ be a finite subset and let ${\cal F}\subset A\otimes C(\T)$ be a finite subset. There exist a finite subset ${\cal G}\subset A$ and $\dt>0$ satisfying the following condition:
For any {{unital}} ${\cal G}$-$\dt$-multiplicative \cp\, $\phi: A\to B$ (for some unital \CA\, $B$) and any unitary $u\in  B$ such that
\beq\label{oldnuc-1}
\|\phi(g)u-u\phi(g)\|<\dt\tforal g\in {\cal G},
\eneq
there exists a unital ${\cal F}$-$\ep$-multiplicative \cp\, $L: A\otimes C(\T)\to B$ such that
\beq\label{oldnuc-2}
\|\phi(f)-L(f\otimes 1)\|<\ep\andeqn
\|L(1\otimes z)-u\|<\ep
\eneq
for all $f\in {\cal F}_0,$ where $z\in C(\T)$ is the identity function
on the unit circle.}}
\end{lem}

\begin{proof}
{{This is well known.
The proof follows the same lines as that of 2.1 of \cite{LP1}.
We sketch the proof here.
Claim:
Suppose that there exists a sequence of unital \cp s $\phi_n: A\to B_n$
for some sequence of unital \CA s $\{B_n\}$ and a sequence
of unitaries $u_n\in B_n$ such that
$\lim_{n\to\infty} \|\phi_n(ab)-\phi_n(a)\phi_n(b)\|=0$ for all $a, b\in A$ and
$\lim_{n\to\infty}\|\phi_n(a)u_n-u_n\phi_n(a)=0$ for all $a\in A.$
Then there exists a sequence of unital \cp s $L_n: A\otimes C(\T)\to B_n$
such that $\lim_{n\to\infty}\|L_n(a\otimes 1_{C(\T)})-\phi_n(a)\|=0$ and
$\lim_{n\to\infty}\|L_n(1\otimes g)-g(u_n)\|=0$ for all $g\in C(\T).$
As in the proof of 2.1 of \cite{LP1}, the lemma follows from the claim.}}

{{Now we prove the claim with exactly the same argument as that of 2.1 of \cite{LP1}.
Consider the \CA s $C=\prod_{n=1}^{\infty} B_n$ and $C_0=\bigoplus_{n=1}^{\infty} B_n.$
Let $\pi: C\to C/C_0$ be the quotient map.
Define   $\Phi_A: A\to C$ by $\Phi_A(a)=\{\phi(a)\}$ for all $a\in A$ and
$\Phi_T: C(\T)\to C$ by $\Phi_T(g)=\{g(u_n)\}$ for all $g\in C(\T).$
Then $\pi\circ \Phi_A$ and $\pi\circ \Phi_T$  are unital \hm s such that  $\pi\circ \Phi_A(a)$
and $\pi\circ \Phi_T(g)$ commute for each $a\in A$ and $g\in C(\T).$
Thus, there is a unital \hm\, $\Psi: A\otimes C(\T)\to C/C_0$ such that
$\Psi(a\otimes g)=\Phi_A(a)\Phi_T(g)$ for all $a\in A$ and $g\in C(\T).$
Since $A\otimes C(\T)$ is amenable, by a result of Choi and Effros, just as in the proof
of 2.1 of \cite{LP1}, we obtain  a unital
\cp\, $L: A\otimes C(\T)\to C$ such that
$\pi\circ L=\Psi.$ Write $L(a)=\{L_n(a)\},$ where each $L_n: A\to B_n$ is a unital  \cp.
Then
\beq
&&\lim_{n\to\infty}\|L_n(a\otimes 1_{C(\T)})-\phi_n(a)\|=0\rforal a\in A
\andeqn\\
&&\lim_{n\to\infty}\|L_n(1\otimes g)-g(u_n)\|=0\rforal g\in C(\T),
\eneq
 as desired.}}
\end{proof}

\begin{df}\label{Dbeta}
{\rm Let $A$ be a unital \CA. Consider the tensor product
$A\otimes C(\T).$
By the K\"{u}nneth {{Formula}} {{(note that $K_*(C(\T))$ is finitely generated)}}, the tensor product induces two canonical injective \hm s
\begin{equation}\label{Dbeta-1}
\bt^{(0)}: K_0(A)\to K_1(A\otimes C(\T))\quad\mathrm{and}\quad
\bt^{(1)}: K_1(A)\to K_0(A\otimes C(\T)).
\end{equation}\index{$\bt^{(i)}$}
In this way (with further application of the  K\"{u}nneth
Formula), one may write
\begin{equation}\label{Dbeta-2}
K_i(A\otimes C(\T))=K_i(A)\oplus \bt^{(i-1)}(K_{i-1}(A)),\,\,\,i=0,1.
\end{equation}
For each $i\ge 2,$ one also obtains the following injective \hm s
\begin{equation}\label{Dbeta-3}
\bt_k^{(i)}: K_i(A, \Z/k\Z)\to K_{i-1}(A\otimes C(\T),\Z/k\Z),\,\,\,i=0, 1.
\end{equation}
Moreover, one may write
\begin{equation}\label{Dbeta-4}
K_i(A\otimes C(\T),\Z/k\Z)=K_i(A,\Z/k\Z)\oplus
\bt_k^{(i-1)}(K_{i-1}(A,\Z/k\Z)),\,\,\,i=0,1.
\end{equation}


If $x\in \underline{K}(A),$ {{let us write}}  ${\boldsymbol{\bt}}(x)$ for $\bt^{(i)}(x)$
if $x\in K_i(A)$ and for
$\bt_k^{(i)}(x)$ if $x\in K_i(A,\Z/k\Z).$ So {{we have}} an injective \hm\,
\beq\label{Dbeta-6}
{\boldsymbol{\bt}}: \underline{K}(A)\to \underline{K}(A\otimes C(\T)),
\eneq
and \index{${\boldsymbol{\bt}}$}
\beq\label{Dbeta-7}
\underline{K}(A\otimes C(\T))=\underline{K}(A)\oplus {\boldsymbol{\bt}}(\underline{K}(A)).
\eneq

Let $h: A\otimes C(\T)\to B$ be a unital \hm. Then
$h$ induces a \hm\, $h_{*i,k}: K_i(A\otimes C(\T),\Z/k\Z)\to K_i(B,\Z/k\Z),$
$k=0,2,3,...$ and $i=0,1.$
Suppose that $\phi: A\to B$ is a unital \hm\,  and
$v\in U(B)$
is a unitary
such that $\phi(a)v=v\phi(a)$ for all $a\in A.$
Then $\phi$ and $v$ determine a unital \hm\,
$h: A\otimes C(\T)\to B$ by $h(a\otimes z)=\phi(a)v$ for all $a\in A,$ where
$z\in C(\T)$ is the identity function on the unit circle $\T,$  {{and every unital \hm\, $A\otimes C(\T)\to B$ arises
in this way.}}
We use {${\rm Bott}(\phi,\, v): \underline{K}(A) \to \underline{K}(B)$ to denote the collection of}
all \hm s $h_{*i-1,k}\circ \bt_k^{(i)},$ {where $h: A\otimes C(\T)\to B$ is the homomorphism determined by $(\phi, v)$,}
and we write\index{${\rm Bott}$}
\beq\label{Dbeta-8}
{\rm Bott}(\phi, \,v)=0
\eneq
if $h_{{{*i-1}},k}\circ \bt_k^{(i)}=0$ for all $k$ and $i.$  
In particular, since $A$ is unital, \eqref{Dbeta-8} implies that $[v]=0$ in $K_1(B)$.
We also use ${\rm bott}_i(\phi, \, v)$ for
$h_{*i-1}\circ \bt^{(i)},$ $i=0,1.$

Suppose that  $A$ is a unital separable amenable \CA.
Let ${\cal Q}\subset \underline{K}(A\otimes C(\T))$ be a finite subset, let
${\cal F}_0\subset A$ be a finite subset, {{and let ${\cal F}_1\subset A\otimes C(\T)$ { also be a finite subset.}}}   Suppose that $(\ep, {\cal F}_{\blue{0}}, {\cal Q})$ is a $KL$-triple.
Then, {\blue{by \ref{oldnuclearity}}}
 there {{exist}} a finite subset ${\cal G}\subset A$ and $\dt>0$ satisfying the following condition:
For any unital ${\cal G}$-$\dt$-multiplicative \cp\, $\phi: A\to B$ (where $B$ is a unital \CA\,) and
{\blue{any}}  unitary
$v\in B$
such that
\beq\label{Dbeta-9}
\|[\phi(g), \, v]\|<\dt\rforal g\in {\cal G},
\eneq
there exists a unital ${\cal F}_1$-$\ep$-multiplicative
\cp\, $L: A\otimes C(\T)\to B$ such that
\beq\label{Dbeta-10}
\|L({{f\otimes 1}})-\phi(f)\|<\ep\rforal f\in {\cal F}_0\andeqn
\|L(1\otimes z)-v\|<\ep,
\eneq
where $z\in C(\T)$ is the standard unitary generator of $C(\T).$
In particular, $[L]|_{\cal Q}$ is well defined {{(see \ref{KLtriple}).}}
Let ${\cal P}\subset \underline{K}(A)$ be a finite subset.
There are $\dt_{\cal P}>0$ and a finite subset ${\cal F}_{\cal P}$ satisfying the following condition: if
$\phi: A\to B$ is a unital ${\cal F}_{\cal P}$-$\dt_{\cal P}$-multiplicative {\blue{\cp\,}} and
(\ref{Dbeta-9}) holds for $\dt_{\cal P}$ (in place of $\dt$) and ${\cal F}_{\cal P}$ (in place of ${\cal G}$),
then {{there exists a unital \cp\,   $L: A\otimes C(\T)\to B$ which satisfies \eqref{Dbeta-10} such that}}
$[L]|_{{\boldsymbol{\bt}}({\cal P})}$ is well defined,
{\blue{and $[L']|_{{\boldsymbol{\bt}}({\cal P})}=[L]|_{{\boldsymbol{\bt}}({\cal P})}$
if $L'$ also satisfies \eqref{Dbeta-10} {{(for the same $\phi$ and $v$)}} (see \ref{KLtriple}).}}  In this case, we will write
\beq\label{Dbeta-11}
{\rm Bott}(\phi,\, v)|_{\cal P}(x)=[L]|_{{\boldsymbol{\bt}}({\cal P})}(x)
\eneq
for all $x\in {\cal P}.$  In particular,
when
$
[L]|_{{\boldsymbol{\bt}}({\cal P})}=0,
$
we will write
\beq\label{Dbeta-12}
{\rm Bott}(\phi,\, v)|_{\cal P}=0.
\eneq
When $K_*(A)$ is finitely generated,
$\mathrm{Hom}_{\Lambda}(\underline{K}(A), \underline{K}(B))$ is determined
by a finitely generated subgroup of $\underline{K}(A)$ (see \cite{DL}).
Let ${\cal P}$ be a finite subset which generates this subgroup.
Then, in this case, instead of (\ref{Dbeta-12}), we may write
\beq\label{Dbeta-13}
{\rm Bott}(\phi, \, v)=0.
\eneq
In general, if ${\cal P}\subset K_0(A),$
we will write
\beq\label{Dbeta-14}
{\rm bott}_0(\phi, \,v)|_{\cal P}={\rm Bott}(\phi,\, v)|_{\cal P},
\eneq
and if ${\cal P}\subset K_1(A),$
we will write
\beq\label{Dbeta-15}
{\rm bott}_1(\phi,\, v)|_{\cal P}={\rm Bott}(\phi, \, v)|_{\cal P}.
\eneq
}
\end{df}

\begin{df}\label{expLR}
{\rm
Let $A$ be a unital \CA.  Each element $u\in U_0(A)$ can be written as  $u= e^{ih_1}e^{ih_2}\cdots e^{ih_k}$ for $h_1, h_2,..., h_k\in A_{{s.a.}}$.
We write  ${\rm cer}(u)\le k$ if $u=e^{ih_1}e^{ih_2}\cdots e^{ih_k}$ for   selfadjoint elements $h_1, h_2, ..., h_k.$
We  write  ${\rm cer}(u)=k$ if ${{\rm cer}}(u)\le k$ and
$u$ is {\it not} a norm limit of unitaries $\{u_n\}$ with ${\rm cer}(u_n)\le k-1.$
We write ${\rm cer}(u)=k+\ep$ if {{${\rm cer}(u)\not\le k$}} and there exists a sequence of unitaries
$\{u_n\}\subset A$ such that $u_n\in U_0(A)$ with ${\rm cer}(u_n)\le k.$

{{Let $u=u(t)\in C([0,1], U(A))$ be a unitary.
Let ${\cal Q}=\{0=t_0<t_1<\cdots t_m=b\}$ be a partition of $[0,1].$
Define $L_{\cal Q}((u(t))_{0\le t\le 1})=\sum_{i=1}^m \|u(t_i)-u(t_{i-1})\|,$ and
\beq
{\rm length}(u(t))_{0\le t\le 1}=\sup\{L_{\cal P}((u(t))_{0\le t\le 1}: {\cal P}\},
\eneq
where the supremum is taken among all possible partitions ${\cal Q}.$}}
Define
\index{${\rm cel}$}
$${\rm cel}(u)=\inf{\big{\{}}\mbox{length of } (u(t))_{0\leq t\leq 1}~|~~ u(t)\in {\blue{C([0,1],U_0(A))}},~ {{u(0))=u,~ u(1)}}=1{\big{\}}}.$$
Obviously, if $u= e^{ih_1}e^{ih_2}\cdots e^{ih_k}$, then
${\rm cel}(u) \leq \|h_1\|+\|h_2\|+\cdots + \|h_k\|~.$  {\blue{In fact (see \cite{Ringrose-cel}),
$$
{\rm cel}(u)=\inf\{\sum_{j=1}^n\|h_j\|: u=\prod_{i=1}^n e^{i h_j},\, h_j\in A_{s.a.}\}.
$$
}}}
\end{df}

\begin{df}\label{Dcu}
Suppose that $A$ is a unital \CA\, with $T(A)\not=\emptyset.$
Recall that $CU(A)$ is the closure of {{the}} commutator subgroup of $U_0(A).$
Let $u\in U(A).$  We {{shall}} use
${\bar u}$ {{to denote}} the image in $U(A)/CU(A).$  It was proved in
\cite{Thomsen-rims} that there is
a splitting short exact sequence
\begin{equation}\label{Dcu-1}
0\to \Aff(T(A))/\overline{\rho_A(K_0(A))}\to \bigcup_{n=1}^{\infty}U(M_n(A))/ \bigcup_{n=1}^{\infty}CU(M_n(A)){{\stackrel{\kappa_1^A}{\to}}}
K_1(A)\to 0.
\end{equation}
{\blue{In what follows, we will use $U(M_{\infty}(A))$ for $\bigcup_{n=1}^{\infty}U(M_n(A)),$
$U_0(M_{\infty}(A))$ for $\bigcup_{n=1}^{\infty}U_0(M_n(A))$, and
$CU(M_{\infty}(A))$ for $\bigcup_{n=1}^{\infty}CU(M_n(A)).$}}
Let $J_c$ {{(or $J_c^A$)}}\index{$J_c,$ $J_c^A$} be a  fixed  splitting map.  Then one may write
\begin{equation}\label{Duc-2}
 U(M_{\infty}(A))
 /CU(M_{\infty}(A))=\Aff(T(A))/\overline{\rho_A(K_0(A))}\oplus J_c(K_1(A)).
\end{equation}

{{As in \cite{GLX-ER}, denote by $P_n(A)$ the subgroup of $K_0(A)$ generated by the projections in $M_n(A)$. Denote by $\rho_A^n
(K_0(A))$ the subgroup $\rho_A
(P_n(A))$ of $\rho_A
(K_0(A))$. In particular,
$\rho_A^1
(K_0(A))$  is the subgroup of $\rho_A
(K_0(A))$ generated by the images of the projections
in A under the map $\rho_A$. Let $\beta: \pi_1(U(M_{\infty}(A))) \to K_0(A)$ denote the inverse  of the  Bott periodicity isomorphism. Then it is well known that the image of $\pi_1(U(M_{n}(A)))$, under $\beta$, contains $P_n(A)$. (Namely, for any $[p]\in K_0(A)$ represented by a projection $p\in M_n(A)$, we have $\beta([u])=[p]$, for $u\in C(\T, U_n(A))$ defined by $u(z)=z\cdot p+(1_n-p)$, for all $z\in \T$.) By Theorem 3.2 of \cite{Thomsen-rims} (see pertinent notation at the beginning of \S 3 of \cite{Thomsen-rims}), if $\rho_A^n
(K_0(A))=\rho_A
(K_0(A))$, then $U_0(M_n(A))/CU(M_n(A))\cong \Aff T(A)/{\blue{\overline{\rho_A
(K_0(A))}}}\cong U_0(M_{\infty}(A))/CU(M_{\infty}(A))$ (cf.
\cite{GLX-ER}).}}

{\blue{In general, let $\widetilde{r}_{\aff}^k: A_{s.a.}\to  U_0(M_k(A))/CU(M_k(A))$ be defined \index{$\widetilde{r}_{\aff}^k$}
by the composition  of $r_{\aff}: A_{s.a.}\to \Aff(T(A))$ and the quotient map
$\mathfrak{q}_k: \Aff(T(A))\to \Aff(T(A))/\overline{\rho_A(\pi_1(U_0(M_k(A))))}\cong U_0(M_k(A))/CU(M_k(A)),$
where the last isomorphism is given by Theorem 3.2 of \cite{Thomsen-rims}.
Denote by ${\widetilde{r}_{\aff}}: A_{s.a.}\to U_0(M_{\infty}(A))/CU(M_{\infty}(A))$
the composition of $r_{\aff}$ and the quotient map $\mathfrak{q}: \Aff(T(A))\to \Aff(T(A))/\overline{\rho_A(K_0(A))}.$ }}

{\blue{Suppose that}} {{$\Aff(T(A))/\overline{\rho_A(K_0(A))}\cong  U_0(M_k(A))/CU(M_k(A)).$}}
{\blue{Then the  map $\widetilde{r}_{\aff}=\widetilde{r}_{\aff}^k$ defined above}} {{can be defined concretely as
\beq\label{Dcu-f}
\widetilde{r}_{\aff}(h)=[\diag(\exp(2\pi ih),1_{k-1})]\in U_0(M_k(A))/CU(M_k(A))~~\mbox{for any}~~h\in A_{s.a}
\eneq
(see Cor 2.12 of \cite{GLX-ER}). The map $\widetilde{r}_{\aff}$ is surjective, since $r_{\aff}$ is surjective.}}

If $A$ has stable rank $k,$ then $K_1(A)=U(M_k(A))/U_0(M_k(A)).$
Note that
$${\rm csr}(C(\T, A))\le {\rm tsr}(A)+1=k+1.$$
It follows from Theorem 3.10  of \cite{GLX-ER} that
\begin{equation}\label{Dcu-3}
{\blue{U_0(M_{\infty}(A))/CU(M_{\infty}(A))}}=U_0(M_k(A))/CU(M_k(A)),
\end{equation}
{{whence it follows that this holds with $U$ in place of $U_0.$}}
Then, {{combining these facts,}}  one has the split short exact sequence
\begin{equation}\label{Dcu-4}
0\to \Aff(T(A))/\overline{\rho_A(K_0(A))}\to U(M_k(A))/CU(M_k(A))\,{{\stackrel{\kappa_1^A}{\to}}}\, U(M_k(A))/U_0(M_k(A))\to \,0,
\end{equation}
and one may write\index{$\kappa_1^A$}
\begin{eqnarray}\label{Dcu-5}
 && U(M_k(A))/CU(M_k(A))=
\Aff(T(A))/\overline{\rho_A(K_0(A))}\oplus J_c(K_1(A))\\
 &&=\Aff(T(A))/\overline{\rho_A(K_0(A))}\oplus J_c(U(M_k(A))/U_0(M_k(A))). \label{Dcu-6}
\end{eqnarray}
{{Note that $\kappa_1^A\circ J_c={\rm id}_{K_1(A)}.$}}
For each {{continuous and piecewise smooth path}} and  $\{u(t): t\in [0,1]\}\subset U(M_k(A)),$ define 
$$
D_A(\{u(t)\})(\tau)={1\over{2\pi i}}\int_0^1 \tau({du(t)\over{dt}}u^*(t))dt,\quad \tau\in T(A).
$$
For each $\{u(t)\},$ the map $D_A(\{u\})$ is a real continuous affine function on $T(A).$\index{$D_A,$ $\overline{D}_A$}
Let $$\overline{D}_A: U_0(M_k(A))/CU(M_k(A))\to \Aff(T(A))/\overline{\rho_A(K_0(A))}$$ denote the de la Harpe and
Skandalis determinant (\cite{HS})  given by
$$
\overline{D}_A(\bar u)=D_A(\{u\})+\overline{\rho_A(K_0(A))},\quad  u\in U_0(M_k(A)),
$$
where $\{u(t): t\in [0,1]\}\subset M_k(A)$ is a continuous and piecewise
smooth path of unitaries with $u(0)=1$ and  $u(1)=u.$  It is known that the de la {{Harpe}} and Skandalis determinant is independent
of the choice of representative for $\bar u$ and the choice of  path $\{u(t)\}.$
Define
\begin{equation}\label{July17-1}
\|\overline{D}_A({\bar u})\|=\inf\{\|D_A(\{v\})\|: v(0)=1,\,\,\, v(1)=v\andeqn {\bar v}={\bar u}\},
\end{equation}
where $\|D_A(\{v\})\|=\sup_{\tau\in T(A)}\|D_A(\{v\})(\tau)\|.$

Suppose that  $u, v\in U(M_k(A)).$ Define {\index{${\rm dist}(\overline{u},\overline{v})$}}
{\blue{
\beq\label{dddcu}
{\rm dist}(\overline{u}, \overline{v})=\inf\{\|uv^*-c\|: c\in CU(M_k(A))\}.
\eneq
It is a metric.}}  Note that ${\rm dist}(\overline{uv^{-1}},\overline{1})={\rm dist}(\overline{u}, \overline{v})\le {\rm dist}(\overline{u}, \overline{1})+{\rm dist}(\overline{v}, \overline{1})={\rm dist}(\overline{u}, \overline{1})+{\rm dist}(\overline{v^{-1}}, \overline{1}).$
 Define
\begin{equation}\label{July16-1}
{\rm d}({\bar u}, {\bar v})=
\left\{
\begin{array}{cl}
2, & \textrm{if $uv^*\not\in U_0(M_k(A))$,\,\,{{or $\|\overline{D}_A(\overline{uv^*})\|\ge 1/2$}}},\\
 \|{e^{2\pi i \|\overline{D}_A(\overline{uv^*})\|}-1}\|, & \textrm{otherwise}.
\end{array}
\right.
\end{equation}
This is also a metric {{(see the lines preceding Theorem 6.4 of \cite{KTm}).}}

Note that, if $u, v\in U_0(M_k(A)),$ then
${\rm d}(\overline{uv^*},\overline{1_k})={\rm dist}(\overline{u}, \overline{v}).$
Now suppose that $A$ has the property that  $\overline{\rho_A^k(K_0(A))}\supset \rho_A(K_0(A)).$
This means $\overline{\rho_A^1(P_k(A))}\supset \rho_A(P_k(A)),$ {{where
$P_k(A)$ is the subgroup of $K_0(A)$ which is generated by the elements in $K_0(A)$
represented by projections in $M_k(A)$.}}
By 3.6 of \cite{GLX-ER}, $U_0(M_k(A))/CU(M_k(A))=U_0(M_m(A))/CU(M_m(A))$ for all $m\ge k.$
{\blue{It follows from  the proof of Theorem  3.1 of  \cite{NT}
 that
 for any unitaries $u, v\in U_0(M_k(A)),$
 $\mathrm{d}(\overline{u}, \overline{v})={\rm dist}(\overline{u}, \overline{v})=\inf\{\|uv^*-w\|: w\in CU(M_k(A))\| \}.$}}
{{On the other hand, if $\|\overline{D}_A(\overline{uv^*})\|=\dt<1/2,$ then
\beq\label{9-23-2018}
\inf\{\|uv^*-w\|: w\in CU(M_k(A))\|\}<2\pi \dt.
\eneq}}
{\blue{See Proposition \ref{CUdist} below for further  discussion.}}
\end{df}


\begin{df}\label{DLddag}
Let $A$ be a unital separable amenable \CA\, and  $B$ be another unital \CA.
{{If $\phi: A\to B$ is a unital \hm,
then $\phi$ induces a continuous \hm\, $\phi^{\ddag}: U(M_m(A))/CU(M_m(A))\to U(M_m(B))/CU(M_m(B))$\index{$\phi^{\ddag}$}
which maps $U_0(M_m(A))/CU(M_m(A))$ to $U_0(M_m(B))/CU(M_m(B)))$  for each $m.$ Moreover,
$\kappa_1^B\circ \phi^{\ddag}=\phi_{*1}\circ \kappa_1^A.$}}

For any finite subset ${\cal U}\subset U(A),$ there exists $1>\dt>0$ and  a finite subset ${\cal G}\subset  A$ with
the following property:
If
$L: A\to B$ is a ${\cal G}$-$\dt$-multiplicative \morp,
then
$\overline{\langle L(u)\rangle}$ is a well-defined element in $U(B)/CU(B)$  for all $u\in {\cal U}.$
{{Recall that we have assumed that $L(1_A)=p$ is a projection in $B.$ Here $\la L(u)\ra$
is originally defined as a unitary in $pBp.$ But we will also use $\la L(u)\ra$ for   $\la L(u)\ra+(1-p),$
whenever it is convenient, and $\overline{\langle L(u)\rangle}$ is defined to be
$\overline{\la L(u)\ra+(1-p)}.$}}

 Let $G({\cal U})$ denote  the subgroup generated by {{ ${\cal U}$ and let $1/4\pi>\ep>0.$
 Denote by $\overline{G({\cal U})}$ the image of $G({\cal U})$ in $U(A)/CU(A).$
 Let ${\cal V}\subset G({\cal U})$ be another  finite subset.}}
 By the Appendix of \cite{Lin-LAH}, there is a  \hm\, $L^{\ddag}: \overline{G({\cal U})}\to U(B)/CU(B)$
 such that
 ${\rm dist}(L^{\ddag}(\overline{u}), \overline{\langle L(u)\rangle})<\ep$ for all $u\in {\cal V },$
 if ${\cal G}$ is sufficiently large and $\dt$ is sufficiently small.
 {{In particular, ${\rm dist}(\overline{\langle L(u)\rangle}, CU(B)))<\ep$ if $u\in CU(A)\cap {\cal V}.$}}
 {{Suppose that $u, v\in U(A)$ and ${\rm dist}(\overline{u}, \overline{v})<\ep.$
 Then ${\rm dist}(uv^*, CU(A))<\ep.$
  With sufficiently large ${\cal G}$ and small
 $\dt,$ we may assume that ${\rm dist}(\overline{\la L(uv^*)\ra}, {\bar 1_B})<\ep.$
 It follows that  we may also assume that ${\rm dist}(L^{\ddag}(\overline{u}), L^{\ddag}(\overline{v}))
 ={\rm dist}(L^{\ddag}(\overline{uv^{-1}}), \overline{1_B})<2\ep.$}}

 {{Note also that we may assume {{$L^{\ddag}((G({\cal U})\cap U_0(A))/CU(A))\subset U_0(B)/CU(B).$}}  To see this, let ${\cal U}_0$ be a finite subset of $U_0(A)\cap G({\cal U})$ which generates $U_0(A)\cap G({\cal U}).$
 (With sufficiently small $\dt$ and large ${\cal G}$---see \ref{KLtriple}), we may assume
 $\langle L(u)\rangle \in U_0(B)$ for each $u\in {\cal U}_0,$
 and
 $[L]$ is well defined on the subgroup generated by the image of ${\cal U}$ in $K_1(A).$
 Let $z\in U(B)$ be such that $\overline{z}=L^{\ddag}(\overline{u}).$ Then
 there exists $\zeta\in CU(B)$ such that $\|z(\la L(u)\ra)^{-1}-\zeta\|< \ep<1/2.$ It follows
 that there exists $y\in U_0(B)$ such that $z(\la L(u)\ra)^{-1}=y\zeta.$  Since $\la L(u)\ra\in U_0(B),$
 $\overline{z}\in U_0(B)/CU(B).$ This proves the assertion above.
 It follows that
 $\kappa_1^B\circ L^{\ddag}(\overline{u})=[L]\circ \kappa_1^A([u])$ for all $u\in G({\cal U}),$
 where $\kappa_1^C: \bigcup_{n=1}^{\infty} U(M_n(C))/\bigcup_{n=1}^{\infty} CU(M_n(C))\to K_1(C)$ is the quotient map
 for a unital \CA\, $C$ (see \ref{Dcu}).}}

 In what follows, whenever we write\index{$L^{\ddag}$}
 $L^{\ddag}$  {{(associated with ${\cal U}$ and $\ep$),
 ${\cal U}$ is specified, and   $1/2>\ep>0$ is given,
 we mean that $\dt$ is small enough and ${\cal G}$ is large enough  that
 we may choose $L^{\ddag}|_{\overline{G({\cal U})}}$ as above to be a \hm\, such that, if $u, v\in {\cal U},$ then}}
 \beq
&& {{{\rm dist}(L^{\ddag}(\overline{u}), \overline{\la L(u)\ra})<\ep/2 \andeqn}}\\
&&{{ {\rm dist}(L^{\ddag}(\overline{u}), L^{\ddag}(\overline{v}))<\ep,\,\,\,{\rm if}\,\,\, {\rm dist}(\overline{u}, \overline{v})<\ep/2,}}
 \eneq
{{and $\kappa_1^B \circ L^{\ddag}(\overline{u})=[\la L(u)\ra]=[L]\circ \kappa_1^B(\overline{u})$ for all $u\in G({\cal U}).$
The latter equation implies that  $L^{\ddag}(\overline{u})^{-1}(\overline{\la L(u)\ra})\in U_0(B)/CU(B).$ Note that such a choice is not unique.
However, if $L^{\dag}$ is another choice which satisfies the requirements above,
then ${\rm dist}(L^{\dag}(\overline{u}), L^{\ddag}(\overline{u}))<\ep$
and $\kappa_1^B\circ L^{\dag}(\overline{u})=\kappa_1^B\circ L^{\ddag}(\overline{u})$ for all $u\in {\cal U}.$ }}
 Moreover, for an integer $k\ge 1,$ we will also use $L^{\ddag}$ for the map on a subgroup of $U(M_k(A))/CU(M_k(A))$ induced by $L\otimes {\rm id}_{M_k}.$
{{Recall also, if $e\in B$ is a projection, and $w\in eBe$ is a unitary, then, we also write ${\bar w}$
for $\overline{w+(1-e)}$ in $U(B)/CU(B).$}}
 {{Finally, if $e_0, e_1\in A$ are two mutually orthogonal  projections and
 $\phi_i: C\to e_iAe_i$ ($i=0,1$) are two ${\cal G}$-$\dt$-multiplicative \morp s,
 then we shall write  $(\phi_0\oplus \phi_1)^{\ddag}(\overline{u}):=
 \phi_0^{\ddag}(\overline{u})\phi_1^{\ddag}(\overline{u})$ for $\overline{u}\in G(\overline{\cal U}).$}}
\end{df}

\begin{lem}\label{HvsU-lem-2018}
{{Let $A$ be a unital {\blue{\CA\,  and let}}
${\cal U}\in U_0(M_k(A))$ be a finite {\blue{subset.}} There {\blue{exists}} a finite {\blue{subset}} ${\cal H}\in A_+$ with
 $\widetilde{r_{\aff}^k}({\cal H})\supset \overline{\cal U}:=\{\bar{u}:~ u\in {\cal U}\}{\blue{\subset U_0(M_k(A))/CU(M_k(A))}}$ {\blue{with the following property:}}
for any $\ep>0$, there are a
{\blue{finite subset}} ${\cal G}\subset A$ ,  $\delta>0,$ and $\eta>0$ such that, for any two unital {\blue{${\cal G}$-$\dt$-multiplicative}}
completely positive
 linear maps $\phi, \psi: A\to B$ {\blue{(for any unital \CA\, $B$) with the property}}
$|\tau(\phi(h)-\psi(h))|<\eta$ for all $h\in {\cal H}$ and $\tau\in {\blue{T(B)}}$, {\blue{we have,}}  in $U_0(M_k(B))/CU(M_k(B))$,  $\dist(\phi^{\ddag}(\bar{u}),\psi^{\ddag}(\bar{u}))<\ep$ for all $u\in {\cal U}.$}}

\end{lem}

\begin{proof}
{{Assume that $\ep<1/4$. For the finite set ${\cal U}$,  since $\widetilde{r^k_{\aff}}$ is surjective, there is a finite set ${\cal H}\subset A_{s.a.}$ such that $\widetilde{r^k_{\aff}}({\cal H})\supset \overline{\cal U}$. We {\blue{may}} assume that $h\in A_+$ for all
$h\in {\cal H};$ otherwise replace $h$ by $h+2{\blue{m}}\pi 1_A$, for a large enough  positive integer ${\blue{m}}$, which has {\blue{the}} same image in $\Aff(T(A))/{\blue{\overline{\rho_A(\pi_1(U_0(M_k(A))))}}}\cong  U_0(M_k(A))/CU(M_k(A)).$
For each $u\in {\cal U}$, there are an $h\in {\cal H}$ and finitely many unitaries $u_i, v_i\in U_0(M_k(A))$ such that $\|u^*\diag(\exp(2\pi ih),1_{k-1})-\prod_{j}u_jv_ju^*_jv^*_j\|<\ep/8$.    {\blue{Choose a}}  finite set ${\cal G}\subset A$ {\blue{which}} contains all the
{\blue{elements}} $h$, $1_A,$ and all entries of $u, u^*, u_j,u^*_j, v_j, v^*_j$ {\blue{(as matrices in $M_k(A)$)}} for $u\in {\cal U}$.
{\blue{ Let $\eta=\ep/4\pi.$}}
If $\delta>0$ is small enough, then for any unital ${\cal G}${\blue{-$\delta$}}-multiplicative completely positive map $\phi: A\to B$, $\phi^{\ddag}|_{\overline{\cal U}}$ can be defined {\blue{so that}}
$$
{\rm dist}(\phi^{\ddag}(\overline{u}), \la \phi(u)\ra)<\ep/32,
$$
$$
\|\phi(u^*\diag(\exp(2\pi ih),1_{k-1}))-\langle \phi(u)\rangle^*\diag(\exp(2\pi i\phi(h)),1_{k-1})\|<\ep/32,~~\mbox{and}
$$
 $$
 \|\phi(\prod_{j}u_jv_ju^*_jv^*_j)-\prod_j\langle \phi(u_j)\rangle\langle \phi(v_j)\rangle\langle \phi(u_j)\rangle^*\langle \phi(v_j)\rangle^*\|<\ep/32.
 $$
 Consequently,
 $$
 \dist({\blue{\phi^{\ddag}(\bar{u})}}, \overline{\diag(\exp(2\pi i\phi(h)),1_{k-1})})<3\cdot(\ep/32)+\ep/8<\ep/4 \in U_0(M_k(B))/CU(M_k(B)).
 $$
 If $\psi: A\to B$ is another unital ${\cal G}${\blue{-$\delta$}}-multiplicative completely positive map such
 that\\
 ${\blue{|}}\tau(\phi(h)-\psi(h))|<\eta$ for all $h\in {\cal H}$ and $\tau \in {\blue{T(B)}}$, then by (\ref{9-23-2018}),
 $$\dist(\overline{\diag(\exp(2\pi i\phi(h)),1_{k-1})},\overline{\diag(\exp(2\pi i\psi(h)),1_{k-1})})<2\pi\eta=\ep/2.$$
 Hence $\dist(\phi^{\ddag}(\bar{u}),\psi^{\ddag}(\bar{u}))<\ep$ for all $u\in {\cal U}$.}}

\end{proof}

{\blue{The following  lemma is  known  and has appeared implicitly
in some of the proofs earlier. We present  it here for convenience.}}

\begin{lem}\label{approx-Aug-14-1}
{\blue{ For any $0<\dt<1/2^8,$
any pair of projections $p$ and $q,$ and any element $x$ in a \CA\, $C,$
if $\|p-x^*x\|<\dt$ and $\|q-xx^*\|<\dt,$
then there exists $w\in C$ such that $w^*w=p,$ $ww^*=q,$ and $\|w-x\|<{(4/3)(1+\dt)(\dt)^{1/4}+4(\dt)^{1/2}}.$}}

{{For any positive number  $\ep<1/4$, and positive integer $K\in \N$, there is a positive number  $\dt<\ep$ such that the following statements are true:}}

{{{\blue{\rm (a)}} For any unital
\CA\, $C$ and
{\blue{any \CA\,}} $C_1\subset C$ with $1_{C_1}=p$, if $e\in C$ is a projection with $\|ep-pe\|<\dt$ and $pep\in_{\dt} C_1$, then there exist projections $q\in C_1\subset pCp$ and $q_0\in (1-p)C(1-p)$ such that
$$ \|q-pep\|<\ep~~~~\mbox{and}~~~~\|q_0-(1-p)e(1-p)\|<\ep.$$
}}
{{{\blue{\rm (b)}} Under {\blue{the}} assumptions of (a), if  unitaries $v\in C$,  $u\in eCe,$ and $w=u\oplus (1-e)$
are such that  $\|vp-pv\|<\dt$, $\|up-pu\|<\dt$, $\|wp-pw\|<\dt,$ and $pvp\in_{\dt} C_1, ~pup\in_{\dt} C_1$, then there exist  unitaries $v'\in C_1\subset C$,  $v''\in (1-p)C(1-p),$ and $z\in qC_1q$ {\blue{($q$ {{as}} in (a))}} such that
$$\|v'-pvp\|<\ep,~~\|v''-(1-p)v(1-p)\|<\ep,~~\|z-quq\|<\ep,~~\mbox{and}~~\|z\oplus(p-q)-pwp\|<\ep.$$
}}
{{{\blue{\rm (c)}} Under {\blue{the}} assumptions of (a), if mutually orthogonal projections $ e_1, e_2,..., e_K\in C$ and partial isometries $s_1, s_2, ..., s_K\in C$ satisfy   $s^*_is_i=e$ and $s_is^*_i=e_i$ for all $i=1,2,...,K,$  where  $1-e=\sum_{i=1}^Ke_i$, and also satisfy
$$\|pe_i-e_ip\|<\dt,~~pe_i p\in_{\dt} C_1,~~\|ps_i-s_ip\|<\dt,~~~\mbox{and} ~~~~ps_i p\in_{\dt} C_1,$$
 for all $i=1,2,...,K$, then there are mutually orthogonal projections ${\blue{q}}, e'_1, e'_2, ..., e'_K\in C_1\subset pCp$, ${\blue{q_0}}, e''_1, e''_2, ..., e''_K\in  (1-p)C(1-p)$, and  partial isometries $s'_1, s'_2,..., s'_K\in C_1\subset pCp$, $s''_1, s''_2, ..., s''_k\in (1-p)C(1-p),$ such that $\|q-pep\|<\ep,$
 $$\|(p-q)-\sum_{i=1}^K e'_i\|<\ep,~~~\|{\blue{(s_i')^*}}s'_i-q\|<\ep,~ ~~~~~ s'_i{\blue{(s'_i)^*}}=e'_i,~~~\mbox{for all} ~~i=1,2,..., K,~~\mbox{and}$$
 $$\|((1-p)-q_0)-\sum_{i=1}^K e''_i\|<\ep,~~\|(s_i'')^*s''_i-q_0\|<\ep,~\tand s_i''(s'')_i^*=e''_i~~\mbox{for all} ~~i=1,2,..., K.$$
 Consequently, $[p-q]=K[q]\in K_0(C_1)$ and $[(1-p)-q_0]=K[q_0] \in K_0(C)$.}}

{\blue{\blue{{\rm (d)}}\,
Let $1/2>\ep>0.$  Under the assumptions of (a), if $\{e_1,e_2,...,e_m\}$ and $\{e_1', e_2',...,e_m'\}$ are two sets of
projections in $C$ and $\{v_1, v_2,...,v_m\}$ is a set of partial isometries such
that $v_j^*v_j=e_j$ and $v_jv_j^*=e_j',$  and if $u_1, u_2,...,u_m,\, u_1', u_2',...,u_m'\in C$ are unitaries
such that $u_j=z_ju_j',$ with $z_j=\prod_{k=1}^{n_j}\exp(i h_{k,j}),$ where $h_{k,j}\in C_{s.a.},$
$j=1,2,...,m,$
then there exist $\ep'>0$ and a finite subset ${\cal F}'\subset C$  with the following properties:

If $\|px-xp\|<\ep'$
and ${\rm dist}(pxp, C_1)<2\ep'$ for all $x\in {\cal F}',$
then there are projections $e_{j,0}, e_{j,0}'\in (1-p)C(1-p)$ and $v_{j,0}\in (1-p)A(1-p),$ and
there are projections $e_{j,1}, e_{j,1}'\in C_1$ and $v_{j,1}\in C_1,$ such that
$\|e_{j,0}-(1-p)e_j(1-p)\|<\ep,$ $\|e_{j,1}-(1-p)e_j'(1-p)\|<\ep,$ $\|e_{j,1}-pe_jp\|<\ep,$
$\|e_{j,1}'-pe_j'p\|<\ep,$ $v_{j,0}^*v_{j,0}=e_{j,0}$  and $v_{j,0}v_{j,0}^*=e_{j,0}',$ and
$v_{j,1}^*v_{j,1}=e_{j,1}$ and $v_{j,1}v_{j,1}^*=e_{j,1}',$
and there are unitaries $u_{j,0}, u_{j,0}'\in (1-p)C(1-p),$ and $h_{k,j,0}\in (1-p)A(1-p)_{s.a.},$
and unitaries $u_{j,1}, u_{j,1}'\in C_1,$ and $h_{k,j,1}\in (C_1)_{s.a.},$ $k=1,2,...,n_j+1,$
such that
$\|u_{j,0}-(1-p)u_j(1-p)]\|<\ep,$ $\|u_{j,0}'-(1-p)u_j'(1-p)\|<\ep,$
$\|u_{j,1}-pu_jp\|<\ep,$ $\|u_{j,1}'-pu_j'p\|<\ep,$   and $u_{j,0}=z_{j,0}u_{j,0}'$ and $u_{j,1}=z_{j,1}u_{j,1}',$
where $z_{j,0}=\prod_{k=1}^{n_j+1}\exp(i h_{k,j, 0})\in (1-p)A(1-p)$ and $z_{j,1}=\prod_{k=1}^{n_j+1}\exp(i h_{k,j,1})\in C_1,$
with $\|h_{n_j+1, k, 0}\|,\, \|h_{n_j+1, k, 1}\|\le \ep$ {{and}} $\|h_{j,k,0}\|, \|h_{k,j,1}\|\le \|h_{k,j}\|,$ $1\le k\le n_j$, $j=1,2,...,m.$
}}

\end{lem}

\begin{proof}
{\blue{Variations of the first part of the statement  are known.  In fact, $\||x|-p\|<\sqrt{\dt}$ (see Lemma 2.3 of \cite{CES}).
Write $x=u|x|$ as the polar decomposition of $x$ in $C^{**}.$
Then $x=u|x|\approx_{\sqrt{\dt}} up=up^2\approx_{\sqrt{\dt}} u|x|p=xp.$
Consider the polar decomposition $x^*=v(xx^*)^{1/2}$ in $C^{**}.$ Then, similarly, $\|qx-x\|<2\sqrt{\dt}.$
Put $y=qxp$ and $\eta:=(1+\dt)\sqrt{\dt}+\dt<1/16.$
 Then $\|y-x\|<4\sqrt{\dt}$ and $\|p-y^*y\|<\|x^*x\|\sqrt{\dt}+\dt\le (1+\dt)\sqrt{\dt}+\dt=\eta.$
 Also, $\|q-yy^*\|<\eta.$
 Thus, $y^*y$ is invertible in $pCp$ and $yy^*$ is invertible in $qCq.$
 One computes that $\|p-|y|^{-1}\|<{\sqrt{\eta}\over{1-\sqrt{\eta}}}<(4/3)\sqrt{\eta}.$
 Set $w=y|y|^{-1}.$ Then $\|w-y\|\le \|y\|\|p-|y|^{-1}\|\le (1+\dt)(4/3)\sqrt{\eta}.$
 It follows that $\|w-x\|<(1+\dt)(4/3)\sqrt{\eta}+4\sqrt{\dt}.$
 As in Lemma 2.5.3 of \cite{Lnbok}, one checks that $w^*w=p$ and $ww^*=q.$}}

{{ {\blue{For the second part of the statement,
one notes that it}} is straightforward to prove  part (a) and part (b) by  standard {\blue{perturbation}}
arguments
(see 2.5 of \cite{Lnbok}).}}

{\blue{For part (c), let $e=e_0.$
Let $\{e_{i,j}: 1\le i, j\le K+1\}$ be a system of matrix units for $M_{K+1}.$
There is a unital \hm\, $\phi: M_{K+1}\to C$ defined by $\phi(e_{i,i})=e_{i-1},$ $i=1,2,...,K+1,$ and
$\phi(e_{1,j})=s_{j-1},$ $j=2,3,...,K.$    Define $\phi_1: M_{K+1}\to pCp$ by
$\phi_1(a)=p\phi(a)p,$ and define $\phi_2: M_{K+1}\to (1-p)C(1-p)$ by $\phi_2(a)=(1-p)\phi(a)(1-p)$
for all $a\in M_{K+1},$ respectively.
Fix $\eta>0.$
Then, by semiprojectivity of $M_{K+1},$ if $\dt$ is sufficiently small,
there is a unital \hm\, $\psi_1:M_{K+1}\to pCp$ and a unital \hm\, $\psi_2:M_{K+1}\to (1-p)C(1-p)$
such that (note that the unital ball of $M_{K+1}$ is compact)
\beq
\max\{\|\phi_1(a)-\psi_1(a)\|: \|a\|\le 1\} <\eta\andeqn \max\{\|\phi_2(a)-\psi_2(a)\|: \|a\|=1\}<\eta.
\eneq
Since $pep, pe_ip, ps_ip\in_{\dt} C_1,$ $i=1,2,...,K,$ for sufficiently small $\dt$ we may assume that there are $a_{i,j}\in C_1$
such that
\beq
\|\psi(e_{i,j})-a_{i,j}\|<2\eta,\,\,\,1\le i,j\le K+1.
\eneq
Then, by semiprojectivity of $M_{K+1}$ again (see, for example, 2.5.9 of \cite{Lnbok}),
with sufficiently small $\eta$ (in other words, with sufficiently small $\dt$),
there is a unital \hm\, $\psi_3: M_{K+1}\to C_1$ such that
\beq
\sup\{\|\psi_3(a)-\psi_1(a)\|: \|a\|\le 1\}<\ep.
\eneq}}
{\blue{ Now let
$q=\psi_3(e_{1,1}),$ $e_i'=\psi_3(e_{i+1,i+1}),$ $s_i'=\psi_3(e_{1,i+1}),$
$q_0=\phi_2(e_{1,1}),$ $e_i''=\phi_2(e_{i+1}, e_{i+1}),$  and $s_i''=\phi_2(e_{1,i+1}),$
$i=1,2,...,K.$  One then verifies that part (c) of the lemma follows.
}}

{\blue{For part (d),
the statement for projections follows a  standard perturbation argument as above.
In fact, it is a direct consequence of (a) and the first part of the statement of the lemma.

To see the second part of (d), let $1/2>\ep>0$ be fixed, and  let
$u_j, u_j',$ and $h_{k,j}$ be given.  With sufficiently small $\ep'>0$ and  {{sufficiently}} large finite subset ${\cal F}',$
we have
\beq
\|(1-p)u_j(1-p)-z_{j,0}'(1-p)u_j'(1-p)\|<\ep/16\andeqn \|pu_jp-z_{j,1}'pu_j'p\|<\ep/16,
\eneq
with {{the}} unitaries
 $z_{j,0}'=\prod_{k=1}^{n_j}\exp(i(1-p)h_{k,j}(1-p))\in (1-p)C(1-p)$
and $z_{j,1}'=\prod_{k=1}^{n_j}\exp(i h_{k,j,1})\in C_1,$ where $h_{k,j,1}\in (C_1)_{s.a.}$  with $\|h_{k,j,1}-ph_{k,j}p\|<\ep'$
and $\|h_{k,j,1}\|\le \|h_{k,j}\|,$
 $j=1,2,...,m.$
Moreover, there are unitaries $u_{j,0}, u_{j,0}'\in (1-p)C(1-p)$ and $u_{j,1}, u_{j,1}'\in C_1$
such that
\beq
&&\|(1-p)u_j(1-p)-u_{j,0}\|<\ep/16, \,\,\|(1-p)u_j'(1-p)-u_{j,0}'\|<\ep/16,\\
 &&\|pu_jp-u_{j,1}\|<\ep/16,\andeqn \|pu_j'p-u_{j,1}'\|<\ep/16,
\eneq
$j=1,2,...,m.$
It follows that
\beq
\|u_{j,0}-z_{j,0}'u_{j,0}'\|<\ep/4\andeqn \|u_{j,1}-z_{j,1}'u_{j,1}'\|<\ep/4.
\eneq
Then, there are  $h_{n_j+1, j, 0}\in (1-p)A(1-p)_{s.a.}$ and $h_{n_j+1,j,1}\in (C_1)_{s.a.}$
such that $\|h_{n_j+1,j, i}\|\le 2\arcsin(\ep/4),$ $i=0,1,$ and
$u_{j,0}=z_{j,0}u_{j,0}'$ and $u_{j,1}=z_{j,1}u_{j,1}',$ where
$z_{j,0}=\exp(i h_{n_j+1,j,0})z_{j,0}$ and $z_{j,1}=\exp(i h_{n_j+1,j, 1} )z_{j,1}',$ $j=1,2,...,m.$
This proves the second part {{of}} (d).
}}

\end{proof}

\begin{df}\label{Mappingtorus}
Let $C$ and $B$ be  unital \CA s and let ${{\phi,\psi}}: C\to B$ be two monomorphisms.
{\blue{Consider the mapping torus}}
\begin{equation}\label{Maptorus-1}
M_{{{\phi, \psi}}}=\{(f,c): C([0,1], B)\oplus C: f(0)=\phi(c)\andeqn f(1)=\psi(c)\}.
\end{equation}
Denote by  $\pi_t: M_{{{\phi, \psi}}}\to B$  the point evaluation at $t\in [0,1].$
One has the
short exact sequence
\begin{equation*}
0\to SB\stackrel{\imath}{\to}M_{{{\phi, \psi}}} \stackrel{\pi_e}{\to} C\to 0,
\end{equation*}
where $\imath: SB\to M_{{{\phi, \psi}}}$ is the embedding and $\pi_e$ is the
quotient map from $M_{{{\phi, \psi}}}$ to $C.$
Denote by $\pi_0, \pi_1: M_{\phi, \psi}\to C$
the point evaluations at $0$ and $1,$ respectively. Since  both ${{\phi}}$ and ${{\psi}}$ are injective, one may identify  $\pi_e$ {{with}}
the point evaluation at $0$ for convenience.

Suppose that $[\phi]=[\psi]$ in $KL(C,B).$ Then $M_{\phi, \psi}$ corresponds {{to}}
{\blue{the zero}} element of $KL({{C}},B).$ In particular, the corresponding extensions
$$
0\to K_i(B) \stackrel{\imath_*}{\to} K_i(M_{\phi, \psi})\stackrel{\pi_e}{\to}K_i( C)\to 0\,\,\,\,\,\,\,{\rm (}i=0,1{\rm )}
$$
are pure {\blue{(see  Lemma 4.3 of \cite{LnAUCT}).}}
\end{df}

\begin{df}\label{DRphipsi}
{Suppose that $T(B)\not=\emptyset.$ Let $u\in M_l(M_{\phi, \psi})$ (for some integer $l\ge 1$) be a unitary
which is a {{piecewise smooth continuous}} function on $[0,1].$
{\blue{Recall from \ref{Dcu}}}
$$
D_B(\{u(t)\})(\tau)={1\over{2\pi} i}\int_0^1 \tau({du(t)\over{dt}}u^*(t))dt\, \tforal \tau\in T(B).
$$
(see \ref{Aq} for the extension of $\tau$ to $M_l(B)$).
}
{Suppose that $\tau\circ \phi=\tau\circ \psi$ for all $\tau\in T(B).$
Then there exists a \hm\,
$$
R_{\phi, \psi}: K_1(M_{\phi, \psi})\to \Aff(T(B)),
$$
defined by $R_{\phi, \psi}([u])(\tau)=D_B(\{u(t)\}){{(\tau)}}$ as above, which is independent  of the choice of the piecewise smooth
path $u$ in $[u].$
We have the following commutative diagram:}
$$
\begin{array}{ccc}
K_0(B) & \stackrel{\imath_*}{\longrightarrow} & K_1(M_{\phi, \psi})\\
 \rho_B\searrow && \swarrow R_{\phi, \psi}\\
  & \Aff(T(B))
  \end{array}.
  $$
{Suppose, in addition,  that $[{{\phi}}]=[{{\psi}}]$ in $KK(C,B).$  Then the following exact sequence splits:
\begin{equation}\label{Aug2-2}
0\to \underline{K}(SB)\to \underline{K}(M_{{\phi, \psi}})\,{\overset{[\pi_e]}{\underset{\theta}{\rightleftharpoons}}} \,\underline{K}(C)\to 0.
\end{equation}
We may assume that $[\pi_0]\circ [\theta]=[{{\phi}}]$ and
$[\pi_1]\circ [\theta]=[{{\psi}}].$
In particular, one may write $K_1(M_{\phi, \psi})=K_0(B)\oplus K_1(C).$ Then we obtain a \hm\,
$$
R_{\phi, \psi}\circ \theta|_{K_1(C)}: K_1(C)\to \Aff(T(B)).
$$
We shall say ``the rotation map vanishes" if there exists  {\blue{a}} splitting  map
$\theta,$ {{as above, such}} that $R_{\phi, \psi}\circ \theta|_{K_1(C)}=0.$}

{Denote by ${\cal R}_0$ the set of those elements $\lambda\in {\rm Hom}(K_1(C), \Aff(T(B)))$ for which there is a \hm\,
$h: K_1(C)\to K_0(B)$ such that $\lambda=\rho_B\circ h.$ It is a subgroup of ${\rm Hom}(K_1(C), \Aff(T(B))).$
If $[\phi]=[\psi]$ in $KK(C,B)$ and $\tau\circ \phi=\tau\circ \psi$
for all $\tau\in T(B),$ one has a well-defined element ${\overline{R_{\phi, \psi}}}\in {\rm Hom}(K_1(C), \Aff(T(B)))/{\cal R}_0$
(which is independent of the choice of $\theta$).}

{Under the assumptions that $[\phi] = [\psi]$ in $KK(C, B)$,  $\tau\circ \phi=\tau\circ \psi$ for all $\tau\in T(B)$,  and $C$ satisfies the UCT},  there exists a \hm\, $\theta'_1: K_1(C)\to K_1(M_{\phi, \psi})$ such that
$(\pi_e)_{*1}\circ \theta_1'={\rm id}_{K_1(C)}$ and $R_{\phi, \psi}\circ \theta'_1\in {\cal R}_0$ if, and only
if, there is $\Theta\in {\rm Hom}_{\Lambda}(\underline{K}(C),\underline{K}(M_{\phi, \psi}))$ such that
$$
[\pi_e]\circ \Theta=[{\rm id}_{C}]\,\,\,{\rm in}\,\,\,KK(C,B)\andeqn R_{\phi, \psi}\circ \Theta|_{K_1(C)}=0.
$$
{\blue{(See the proof of 4.5 of \cite{L-N}.)}}
In other words, $\overline{R_{\phi, \psi}}=0$ if, and only if, there is $\Theta$ as described above
such that $R_{\phi, \psi}\circ \Theta|_{K_1(C)}=0.$  When $\overline{R_{\phi, \psi}}=0,$ one has that
$\theta(K_1(C))\subset {\rm ker}R_{\phi, \psi}$ for some $\theta$ such that (\ref{Aug2-2}) holds. In this case
$\theta$ also gives the following {{decomposition:}}
$$
{\rm ker}R_{\phi, \psi}={\rm ker}\rho_B\oplus K_1(C).
$$

\end{df}



\begin{df}\label{appep}
Let $C$ be a  \CA, let $a, b\in C$ be two elements, and let $\ep>0.$
We write $a\approx_\ep b$\index{$a\approx_\ep b$} if $\|a-b\|<\ep.$
Suppose that $A$ is another \CA,  $L_1, L_2: C\to A$ are two maps, and ${\cal F}\subset C$ is a subset (usually finite).
We write
\beq\label{appep-1}
L_1\approx_{\ep} L_2\,\,\,{\rm on}\,\,\, {\cal F},
\eneq\index{$L_1\approx_{\ep} L_2$}
if $\|L_1(c)-L_2(c)\|<\ep$ for all $c\in {\cal F}.$

\end{df}

\begin{df}\label{Dfull}
Let $A$ and $B$ be \CA s, and assume that $B$ is unital. Let $\mathcal H\subset  A_+\setminus \{0\}$ be a finite subset, and  let
{\blue{functions}} $T: A_+\setminus\{0\}\to \R_+\setminus\{0\}$ and $N: A_+\setminus\{0\}\to \N$ {\blue{be given.}}
{\blue{A}} map $L: A\to B$ is said to be $T\times N$-$\mathcal H$-full\index{$T\times N$-$\mathcal H$-full} if for any $h\in \mathcal H$, {{there}} are $b_1, b_2, ..., b_{N(h)}\in B$ such that $\|b_i\|\leq T(h)$ and
$$\sum_{i=1}^{N(h)} b_i^*L(h)b_i=1_B.$$
{{We say $L$ is $T\times N$-full, if it is $T\times N$-${\cal H}$-full for every finite ${\cal H}\subset A_+\setminus \{0\}.$}}
\end{df}

\begin{prop}\label{full-2018-sept} {{ Let $A, T, N$ be as in Definition \ref{Dfull}. Let  ${{{\cal H}_1\subset {\cal H}_2\subset \dots}} \subset A_+^{{\bf 1}}\setminus \{0\}$ and ${{{\cal G}_1\subset {\cal G}_2\subset \cdots}} \subset {{A^{\bf 1}}}$
be two {{increasing}} sequences of finite subsets such that
$\bigcup_{n=1}^{\infty}{{{\cal H}_n}}$ is dense in the unit ball $A_+^{{\bf 1}}$ of $A_+$ and $\bigcup_{n=1}^{\infty}{\cal G}_n$ is dense in the unit ball $A^{\bf 1}$ of $A$.
{{Suppose also that, if $h\in {\cal H}_n$ then $f_{1/2}(h)\in {\cal H}_{n+1}\cup \{0\}.$}} Let $\{\dt_n\}$ be a decreasing
sequence of positive numbers with
{{$\lim_{n\to\infty}\dt_n=0,$ let}} $\{B_n\}$ be a sequence of unital \CA s, and let
 $\{\sigma_n: A\to B_n\}$ be a sequence of ${\cal G}_n${{-$\delta_n$-}}multiplicative {{and}} $T\times N$-$\mathcal H_n$-full unital completely positive maps. Then $\{\sigma_n\}$ induces a full homomorphism $\Sigma: A\to Q(B)$, where $B=\prod_{n=1}^{\infty} B_n$ and $B^0={{\bigoplus}}_{n=1}^{\infty} B_n$, and $Q(B)=B/B^0$.  }}
\end{prop}

 \begin{proof}
 Denote by $S: A\to \prod_{n=1}^{\infty}B_n$ the map defined by
 $S(a)=\{\sigma_n(a)\}$ for all $a\in A.$ Note that $\Sigma=\pi\circ S,$
 where $\pi: B\to Q(B)$ is the quotient map.
 By the properties of ${\cal G}_n$ and $\delta_n$, the induced  map $\Sigma: A\to Q(B)$
 is a homomorphism. To prove $\Sigma$ is full, we need to prove that for any $h\in {{A_+\setminus\{0\}}}$, $\Sigma(h)$ is full in $Q(B)$. For any $h\in {{\mathcal H}}_n$,  and for any $m\geq n$, there are $b^m_1, b^m_2, ..., b^m_{N(h)}\in B_m$ such that $\|b^m_i\|\leq T(h)$ and
$$\sum_{i=1}^{N(h)} (b^m_i)^*\sigma_m(h)b^m_i=1_{B_m}.$$
Set
$b_i={{(\underbrace{0,\cdots,0}_{n-1}, b^n_i, b^{n+1}_i,  \cdots)}}.$
Since $\|b_i^n\|\le T(h)$ for all $n,$
$b_i\in {{B}},$ $i=1,2,...,N(h).$ Then $\sum_{i=1}^{N(h)} {{\pi(b_i)}}^*\Sigma(h){{\pi(b_i)}}=1_{Q(B)}.$ That is, $\Sigma(h)$ is full for any $h\in {\cal H}_n$.
Now let $a\in A_+$  with $\|a\|=1.$
Then there exists $h\in  {\cal H}_n$ for some $n$ such that $\|a-h^2\|<1/16.$
It follows from Proposition 2.2 and part (a) of Lemma 2.3 of \cite{RorUHF2}  that there exists $r\in A$ such that $r^*ar=f_{1/2}(h).$
Since $\|a\|=1,$ $\|h\|>15/16.$ Thus, $f_{1/2}(h)\not=0.$ It follows that $f_{1/2}(h)\in {\cal H}_{n+1}.$
Since $\Sigma(f_{1/2}(h))$ is full in $Q(B),$ in other words, $\Sigma(r)^* \Sigma(a) \Sigma(r)$ is full,
one concludes that $\Sigma(a)$ is full.  Thus,  $\Sigma(a)$ is full in $Q(B)$ for any $a\in A_+\setminus \{0\}.$

 \end{proof}

\begin{NN}\label{UHFdense}
{{Let $A$ be a unital separable \CA\, and $U$ an infinite dimensional {{UHF-algebra}}.
Write $U=\overline{\bigcup_{n=1}^{\infty} B_n},$ where $B_n$ is a full matrix algebra and
$B_{n+1}=B_n\otimes M_{r_n}$ for some integer $r_n\to \infty$ and
where $B_n$ is identified with  $B_n\otimes 1_{r_n}$ as a \SCA\, of $B_{n+1}.$ }}

{{Let  $p\in A\otimes U$ be a projection, $r\in (0,1),$ and $\ep>0.$
Choose $n_0\ge 1$ such that $r_{n_0}>1/2\ep.$ There is a projection $q\in A\otimes B_n,$
for some $n\ge n_0,$ such that $\|q-p\|<1/2.$ Then there is $k\in \N$ such that
$|r-k/r_n|<\ep$ and a projection $e=q\otimes q'$  with $q'\in M_{r_n}$ and
$t(q')=k/r_n,$ where ${\rm tr}$ is the unique tracial state of $M_{r_n}.$ Moreover,
if $\tau\in T(A)$ then $\tau(e)={\rm tr}(q')\tau(p).$  In other words,
$|\tau(e)-r\tau(p)|<\ep$ for all $\tau\in T(A).$
We will use this fact later.}}
\end{NN}

\begin{NN}\label{Ddiag}
{\rm
{{Let $A$ be a unital \CA\, and let $e_1, e_2,...,e_n$ be mutually orthogonal
and mutually equivalent projections in $A$ and set  $p=\sum_{i=1}^n e_i.$
Let $v_i\in pAp$ be partial isometries such that $v_iv_i^*=e_1$ and $v^*_iv_i=e_i,$ $i=1,2,...,n.$
Then one may identify $pAp$ with $M_n(e_1Ae_1).$
Let $C$ be another \CA\, and let $\phi: C\to e_1Ae_1$ be a map.
Define $\phi_i: C\to e_iAe_i$ by $\phi_i(c)=v_i^*\phi(c)v_i$ for all $c\in C.$
In this paper, we shall often write
\beq\label{dddiag}
\sum_{i=1}^n\phi_i(c)=\diag(\overbrace{\phi(c),\phi(c),...,\phi(c)}^n)=\phi(c)\otimes 1_n\rforal c\in C.
\eneq}}
If $a,\, b\in A$ and $ab=ba=0,$ we will write $a+b=a\oplus b.$
}
\end{NN}

\section{The Elliott-Thomsen  building blocks}

 To generalize the class of \CA s of tracial rank at most one,
 {\blue{it would be  natural to}} consider
 all  subhomogeneous \CA s with one dimensional spectrum  which, in particular,
 include circle algebras as well as dimension drop interval algebras. We begin, however, with the following
 special {\blue{class}}:

\begin{df}[See \cite{ET-PL} and  \cite{point-line}]\label{DfC1}
{\rm
Let $F_1$ and $F_2$ be two finite dimensional \CA s.
Suppose that there are two unital \hm s
$\phi_0, \phi_1: F_1\to F_2.$
Consider  the mapping torus \index{$A(F_1, F_2,\phi_0, \phi_1)$} {\blue{$M_{\phi_0, \phi_1}$ (see \ref{Mappingtorus}):}}
$$
A=A(F_1, F_2,\phi_0, \phi_1)
=\{(f,g)\in  C([0,1], F_2)\oplus F_1: f(0)=\phi_0(g)\andeqn f(1)=\phi_1(g)\}.
$$
These \CA s {{were}} introduced into the Elliott program by Elliott and Thomsen (\cite{ET-PL}), and in \cite{point-line}, Elliott used this class of \CA s {{and some other building blocks with 2-dimensional spectra}} to realize any weakly unperforated {{simple}} ordered group {{with order unit}} as the $K_0$-group of a simple ASH \CA. Denote by ${\cal C}$ the class of all unital \CA s of the form $A=A(F_1, F_2, \phi_0, \phi_1)$ (which includes  all finite dimensional \CA s).
These \CA s will be called Elliott-Thomsen building blocks.

A unital \CA\, $C\in {\cal C}$ \index{${\cal C}$} is said to be {\it minimal} if
it is not the direct sum of  two non-zero
of  \CA s.
If $A\in {\cal C}$ is minimal and  is not finite dimensional,
{{then}} ${\rm ker}\phi_0\cap {\rm ker}\phi_1=\{0\}.$
In general, if $A\in {\cal C},$ and ${\rm ker}\phi_0\cap {\rm ker}\phi_1\not=\{0\}{{,}}$ {{then}}  $A=A_1 \oplus ({\rm ker}\phi_0\cap {\rm ker}\phi_1)$, where $A_1=A(F'_1, F_2,\phi'_0, \phi'_1)$,  $F_1=F'_1\oplus ({\rm ker}\phi_0\cap {\rm ker}\phi_1)$ and $\phi'_i=\phi_i|_{F'_1}$ for $i=1,2$. (Note that $A_1=A(F'_1, F_2,\phi'_0, \phi'_1)$ satisfies the condition  ${\rm ker}\phi'_0\cap {\rm ker}\phi'_1 =\{0\},$ and that ${\rm ker}\phi_0\cap {\rm ker}\phi_1$ is a finite dimensional \CA.)

{\blue{
Let $\lambda: A\to C([0,1], F_2)$ be defined by $\lambda((f,a))=f.$ Note that if $A$ is infinite dimensional and minimal,
then $\lambda$ is injective, since ${\rm ker}\phi_0\cap {\rm ker}\phi_1=\{0\}.$}}
{{As in Definition of \ref{Mappingtorus},}} for $t\in (0,1),$ define $\pi_t: A\to F_2$ by $\pi_t((f,g))=f(t)$ for all $(f,g)\in A.$
For  $t=0,$ define $\pi_0: A\to \phi_0(F_1)\subset F_2$ by $\pi_0((f, g))=\phi_0(g)$ for all $(f,g)\in A.$
For $t=1,$ define $\pi_1: A\to \phi_1(F_1)\subset F_2$ by $\pi_1((f,g))=\phi_1(g))$ for all $(f,g)\in A.$
In what follows, we will call $\pi_t$  a {{point evaluation}} of $A$ at $t.$
There is a canonical map $\pi_e: A \to F_1$ defined by $\pi_e(f,g)=g$ for all
pair $(f, g)\in A.$  It is a surjective map.\index{$\pi_e$}
{\it The notation $\pi_e$ will be used for this map throughout  this paper.}

If $A\in {\cal C}$, then $A$ is the pull-back  {{corresponding to the diagram}}
\begin{equation}\label{pull-back}
\xymatrix{
A \ar@{-->}[rr]^{\lambda} \ar@{-->}[d]^-{\pi_e}  && C([0,1], F_2) \ar[d]^-{(\pi_0, \pi_1)} \\
F_1 \ar[rr]^-{(\phi_0, \phi_1)} & & F_2 \oplus F_2\,.
}
\end{equation}
%
{{Conversely,}} every such pull-back is an algebra in ${\cal C}.$
Infinite dimensional  \CA s in ${\cal C}$ are also called {\it one-dimensional non-commutative finite CW complexes} (NCCW)
(see \cite{ELP1} and \cite{ELP2}). 

{{We would like to mention that the \CA s $C([0,1], F_2)$ and $C(\T, F_2)$ are in ${\cal C}.$
Suppose that  $F_1,$ $F_2,$ $\phi_0$ are as mentioned  above, and fix a point $t_0\in [0,1].$ Define
$$
B=\{(f, g)\in C([0,1], F_2)\oplus F_1: f(t_0)=\phi_0(g)\}.
$$
Then one easily verifies   that $B$ is also in ${\cal C}.$}}

Denote by ${\cal C}_0$ the sub-class of those \CA s $A$ in ${\cal C}$ such that $K_1(A)=\{0\}.$
\index{${\cal C}_0$}

{\CA s in ${\cal C}$ are finitely generated (Lemma 2.4.3 of  \cite{ELP1}) and semiprojective (Theorem 6.22 of \cite{ELP1}). We will use these important features later without further mention.}
}

\end{df}

\begin{lem}\label{2Lg1}
Let $f\in C([0,1], M_k)$ and let $a_0,a_1\in M_k$ be invertible elements with
$$\|a_0-f(0)\|<\ep\andeqn  \|a_1-f(1)\|<\ep.$$
Then there exists an invertible element  $g\in C([0,1], M_k)$ such that $g(0)=a_0,~ g(1)=a_1$, and
$$\|f(t)-g(t)\|<\ep \qq \tforal t\in [0,1].$$
\end{lem}
\begin{proof}
Let $S\sbs M_k$ denote the set consisting of all singular matrices.  Then $M_k$ is a $2k^2$-dimensional differential manifold (diffeomorphic to $\R^{2k^2}$), and  $S$ is a finite union of closed submanifolds of codimension at least two. Since each continuous map between two differential manifolds (perhaps with boundary) can be approximated
arbitrarily well by smooth map{s}, we can find $f_1\in C^{\infty}([0,1],M_k)$ with $f_1(0)=a_0$ and $f_1(1)=a_1$ and $\|f_1(t)-f(t)\|<\ep'<\ep$.  Apply  the relative version of the transversality theorem---the corollary on  page 73 of \cite{Guillemin-Pollack} and its proof (see pages 70 and 68 of \cite{Guillemin-Pollack}), for example, with $Z=S,$ $Y=M_k,$ $X=[0,1]$ and with $\partial X=\{0,1\}${{---}}to obtain $g\in C^{\infty}([0,1],M_k)$ with $g(0)=f_1(0)$,
 $g(1)=f_1(1)$, $g(t)\notin S,$ and $\|g(t)-f_1(t)\|< \ep-\ep'$. {Hence, $\|f(t) - g(t)\| <\ep$, $t\in[0, 1]$.}
\end{proof}

\begin{prop}\label{2pg3}
If  $A\in {\cal C}$, then $A$ has stable rank one.\index{stable rank one}
\end{prop}
\begin{proof}
Let $(f,a)$ be in $A$ with $f\in C([0,1], F_2)$ and $a\in F_1$  with $f(0)=\phi_0(a)$ and $f(1)=\phi_1(a)$.  For any $\ep>0$, since $F_1$ is a finite dimensional \CA, there is an invertible element $b\in F_1$ such that $\|b-a\|<\ep$.  Since $\phi_0$ and $\phi_1$ are unital, $\phi_0(b)$ and $\phi_1(b)$ are invertible. {\blue{Also}},
$$
\|\phi_0(b)-f(0)\|<\ep\andeqn \|\phi_1(b)-f(1)\|<\ep.
$$
By Lemma \ref{2Lg1} (applied to each direct summand of $F_2$),
there exists an invertible element $ g\in C([0,1], F_2)$ such that $g(0)=\phi_0(b),~~ g(1)=\phi_1(b),$ and
$$\|g-f\|<\ep.$$
This is what was desired.
\end{proof}

\begin{NN}\label{2Rg11}
Let $F_1=M_{R{{(1)}} }\oplus M_{R{{(2)}} }\oplus \cdots \oplus M_{R{{(l)}}},$  let
 $F_2=M_{r{{(1)}}}\oplus M_{r{(2)} }\oplus \cdots \oplus M_{r{(k)} }$  and  let $\phi_0,\,\,
  \phi_1:~ F_1\to F_2$ be unital homomorphisms, where $R(j)$ and $r(i)$ are positive integers.  Then $\phi_0$ and $\phi_1$ induce homomorphisms $$\phi_{0*}, \phi_{1*}: ~ K_0(F_1)=\Z^l \longrightarrow K_0(F_2)=\Z^k$$
{{represented}} by matrices $(a_{ij})_{k\times l}$ and $(b_{ij})_{k\times l}$, respectively, where
 $r{{(i)}}=\sum_{j=1}^l a_{ij}R{{(j)}}=\sum_{j=1}^l b_{ij}R(j)$ for $i=1,2,...,k.$
 {{Note that $a_{ij}$ and $b_{ij}$ are non-negative integers. The matrices
 $(a_{ij})_{k\times l}$
 and $(b_{ij})_{k\times l}$ {\blue{will be}} called the multiplicities of $\phi_0$ and
 $\phi_1,$ respectively.}}

\end{NN}

\begin{prop}[{{see also 2.1 of \cite{KT98}}}]\label{2Lg13}
{{Let $F_1$
and $F_2$ be as in \ref{2Rg11},}}
{\blue{and
let}} $A=A(F_1, F_2, \phi_0, \phi_1)$, {where $\phi_1, \phi_2: F_1\to F_2$ are two unital homomorphisms}.  Then $K_1(A)=\Z^k/\mathrm{Im}({\phi_0}_{*0}-{\phi_1}_{*0})$ and


\beq\label{2Rg12-n1}
K_0(A)\cong \left\{
\left(\!\!
         \begin{array}{c}
           v_1 \\
           v_2 \\
           \vdots \\
           v_l
         \end{array}\!\!
       \right) \in \Z^l~ , ~~  \phi_{~\!\!\!\!_{0*}}\!\!\left(\!\!
         \begin{array}{c}
           v_1 \\
           v_2 \\
           \vdots \\
           v_l
         \end{array}\!\!
       \right)=   \phi_{~\!\!\!\!_{1\!*}}\!\!\left(\!\!
         \begin{array}{c}
           v_1 \\
           v_2 \\
           \vdots \\
           v_l
         \end{array}\!\!
       \right)
\right\},
\eneq
with positive cone $K_0(A) \cap \Z^l_+$, and {scale} $\left(\!\!
         \begin{array}{c}
           \blue{R(1)} \\
          \blue{R(2)}\\
           \vdots \\
           \blue{R(l)}
         \end{array}\!\!
       \right)\in \Z^l$, where
$\Z_+^l=\left\{
\left(\!\!
         \begin{array}{c}
           v_1 \\
           v_2 \\
           \vdots \\
           v_l
         \end{array}\!\!
       \right)~;~~ v_i\geq 0
       \right\} \sbs\, \Z^l~.
$
Moreover, the map $\pi_e: A\to F_1$ induces the natural order embedding
$(\pi_e)_{*0}: K_0(A)\to K_0(F_1)=\Z^l;$ in particular,  ${\rm ker}\rho_A=\{0\}$
{{{\rm{(see Definition \ref{Drho} for $\rho_A$)}}}}.

Furthermore, if $K_1(A)=\{0\},$ then $\Z^l/K_0(A)\cong K_1(C_0((0,1), F_2))$ is torsion free,
and in this case, {{
$[\pi_e]|_{K_i(A, \Z/k\Z)}$ is injective for all $k\ge 2$ and $i=0,1.$}}

\end{prop}

\begin{proof}
Most of these {statements} are known. We sketch the proof here.
Consider the short exact sequence
\begin{displaymath}\label{sixterm-1}
\xymatrix{
0\ar[r] & \mathrm{C}_0((0,1), F_2) \ar[r] & A \ar[r]^{\pi_e} & F_1\ar[r] & 0.}
\end{displaymath}
We obtain  the  exact sequence
\begin{equation}\label{sixterm}
\xymatrix{
0 \ar[r] & K_0(A) \ar[r]^{{\pi_e}_*} & K_0(F_1) \ar[r] & K_0(F_2) \ar[r] & K_1(A)\ar[r] & 0,
}
\end{equation}
where the map $K_0(F_1)\to K_0(F_2)$ is given by ${\phi_0}_{*0}-{\phi_1}_{*0}$.
In particular, $(\pi_e)_{*0}$ is injective.
{{If $x\in K_0(A)$ is such that $(\pi_e)_{*0}(x)=[p],$
where $p\in \pi_e(M_m(A))$ is a projection,}} then, {{by \eqref{sixterm},}}
{{$((\phi_0)_{*0}-(\phi_1)_{*0})([p])=0,$ or}} $(\phi_0)_{*0}([p])=(\phi_1)_{*0}([p]).$
Therefore, $\phi_0(p)$ and $\phi_1(p)$ have the same rank. It follows  that there is a projection
$q\in M_m(A)$ such that $\pi_e(q)=p.$ {{Thus, $(\pi_e)_{*0}(x-[q])=0.$ Since
$(\pi_e)_{*0}$ is injective, $x=[q]\in K_0(A)_+.$}} This implies that $(\pi_e)_{*0}$ is an order embedding.
This also implies that ${\rm ker}\rho_A=\{0\}.$  The descriptions of $K_i(A)$ ($i=0,1$) also follow.
The quotient group
$K_0(F_1)/{{(\pi_e)_{*0}(K_0(A))}}$ is always torsion free, since $K_0(F_1)/(\pi_e)_{*0}(K_0(A))$ is a subgroup of $K_0(F_2)$. It happens that $K_0(F_1)/{{(\pi_e)_{*0}}}(K_0(A))=K_0(F_2)$ holds if and only if $K_1(A)=\{0\}.$

 In the case that $K_1(A)=\{0\},$
one also computes that $K_0(A, \Z/k\Z)$ may be identified with
$K_0(A)/kK_0(A)$ and $K_1(A,\Z/k\Z)=\{0\}$ for all $k\ge 2.$

{{To see that $[\pi_e]|_{K_0(A,\Z/k\Z)}$ is injective for $k\ge 2,$  let
${\bar x}\in K_0(A,\Z/k\Z)=K_0(A)/kK_0(A)$ and let
$x\in K_0(A)$ be such that its image in $K_0(A)/kK_0(A)$ is ${\bar x}.$
If $[\pi_e]({\bar x})=0,$ then $(\pi_e)_{*0}(x)\in kK_0(F_1).$ Let
$y\in K_0(F_1)$ be such that $ky=(\pi_e)_{*0}(x).$ Then $k{\bar y}={\overline{(\pi_e)_{*0}(x)}}=0$ in\\
$K_0(F_1)/(\pi_e)_{*0}(K_0(A)).$ This implies that $K_0(F_2)$ has torsion, a contradiction.
Therefore $[\pi_e]|_{K_0(A/\Z/k\Z)}$ is injective for $k\ge 2.$ }}
  Then {{ the rest of}} the proposition also  follows from \eqref{sixterm}.
\end{proof}

\begin{prop}\label{2Pg12}
For fixed finite dimensional \CA s $F_1, F_2$, the \CA\, $A=A(F_1, F_2, \phi_0, \phi_1)$ is completely determined
(up to isomorphism) by  the maps $\phi_{0*}, \phi_{1*}: \Z^l \longrightarrow \Z^k$.
\end{prop}

\begin{proof} Let $B=A(F_1, F_2, \phi'_0, \phi'_1)$ with $\phi'_{0*}=\phi_{0*}, \phi'_{1*}=\phi_{1*}$.  It is well known 
that there exist two unitaries $u_0, u_1\in F_2$ such that
$$u_0 \phi_0(a) u_0^*=\phi'_0(a), \quad a\in F_1,\andeqn
u_1 \phi_1(a) u_1^*=\phi'_1(a), \quad a\in F_1.$$
Since $U(F_2)$ is path connected, there is a unitary path $u:~ [0,1]\to U(F_2)$ with $u(0)=u_0$ and $u(1)=u_1$.  Define $\phi:~ A \to B$ by
$$\phi(f,a)=(g,c),$$
where $g(t)=u(t) f(t) u(t)^*$.
 Then a straightforward calculation shows that the map $\phi$ is a *-isomorphism.
\end{proof}

\begin{NN}\label{2Rg10}
Let $F_1=M_{R(1)}\oplus M_{R(2)}\oplus \cdots \oplus M_{R(l)}$ and ${{F}}_2=M_{r(1)}\oplus M_{r(2)}\oplus \cdots \oplus M_{r(k)}$ and let  $A=(F_1, F_2,\phi_0, \phi_1).$
Denote by $m$ the greatest common {{divisor}} of
${\big \{}R(1),R(2), \cdots, R(l){\big \}}$.  Then each $r(i)$ is also a multiple of $m$.
 Let $\widetilde{F}_1=M_{\frac{R(1)}{m} }\oplus M_{\frac{R(2)}{m} }\oplus \cdots \oplus M_{\frac{R(l)}{m} }$
 and
 $\widetilde{F}_2=M_{\frac{r(1)}{m} }\oplus M_{\frac{r(2)}{m} }\oplus \cdots \oplus M_{\frac{r(k)}{m} }$.  Let ${\widetilde{\phi}_0},\widetilde{\phi}_1:~ \widetilde{F}_1\to \widetilde{F}_2$ be maps such that  the $K_0$-maps
\vspace{-0.1in} \beq\nonumber
 \widetilde{\phi}_{0*}, \widetilde{\phi}_{1*}: ~ K_0(\widetilde{F}_1)=\Z^l \longrightarrow \Z^k
 \eneq
 satisfy $\widetilde{\phi}_{0*}=(a_{ij})_{k\times l}$ and $\widetilde{\phi}_{1*}=(b_{ij})_{k\times l}$.  That is, $\widetilde{\phi}_{0*}$
 and $\widetilde{\phi}_{1*}$ are the same as
 $\phi_{0*}$ and $ \phi_{1*}$.
By \ref{2Pg12},
$$A(F_1, F_2, \phi_0, \phi_1)\cong M_m(A(\tilde{F}_1, \widetilde{F}_2, \widetilde{\phi}_0, \widetilde{\phi}_1)).$$
\end{NN}

\begin{NN}\label{LgN1889}
 It is well known {{that}} the extreme points of $T(A)$  are in canonical one-to-one
 correspondence with the  
irreducible representations of $A$, which are given by
$$\coprod_{j=1}^k (0,1)_j \cup \{ \rho_1, \rho_2, ..., \rho_l\} = \mathrm{Irr}(A),$$
where  each $(0,1)_j$ is the open interval $(0,1).$ (We use the subscript $j$ to indicate the $j$-th copy.)

{{It follows from Lemma 2.3 of \cite{KT98} that the}} affine function space $\mathrm{Aff}(\mathrm{T}(A))$ can be identified with the subset of
$$\bigoplus_{j=1}^k C([0,1]_j,\R)\oplus \underbrace{(\R\oplus\R\oplus\cd\oplus \R)}_{l~~ copies}$$
consisting of {{the}} elements $(f_1,f_2,..., f_k, g_1, g_2,...,g_l)$ satisfying the conditions
$$f_i(0)=\frac1{R(i)}\sum_{j=1}^l a_{ij}g_j\cdot r_j \qq\qq \mbox{and} \qq\qq
f_i(1)=\frac1{R(i)}\sum_{j=1}^l b_{ij}g_j\cdot r_j,$$
where $(a_{ij})_{k\times l}=\phi_{0*}$ and $(b_{ij})_{k\times l}=\phi_{1*}$ as in \ref{2Rg11}.

{{Denote by $\tau_{t, j}\in T(A)$ the tracial state defined by
$\tau_{t,j}((g,a))={\rm tr}_j(\pi_j(g(t)))$ for all $(g, a)\in A,$ where
$t\in (0,1),$ ${\rm tr}_j$ is the tracial state of $M_{r(j)},$
and $\pi_j: F_2\to M_{r(j)}$ is the quotient map, $j=1,2,...,k.$
Denote  by $\tau_{0,j}\in T(A)$ the tracial state defined by
$\tau_{0,j}((g,a))={\rm tr}_j\circ \pi_j(g(0))$  and $\tau_{1,j}((g,a))={\rm tr}_j\circ \pi_j(g(1))$ for all $(g,a)\in A.$
If $t\to 0$ in $[0,1],$ then $\tau_{t,j}\circ \pi_j(g(t))\to \tau_{0,j}\circ \pi_j(g(0))$ for all $(g,a)\in A,$
$j=1,2,...,k,$ It follows that $\tau_{t,j}\to \tau_{0,j}$ in (the weak* topology of) $T(A).$
For exactly same reason $\tau_{t,j}\to \tau_{1,j}$ in $T(A)$ if $t\to 1$ in $[0,1].$
}}
{{Let $f=(f_1,f_2,...,f_k, g_1,g_2,...,g_l)\in \Aff(T(A)).$
Then $f(\tau_{0,j})=f_j(0)$ and $f(\tau_{1,j})=f_j(1),$ $j=1,2,...,k.$
Denote by $\tau_{e,i}\in T(A)$ the tracial state defined by
$\tau_{e,i}((g,a))={\rm tr}_{e,i}(\psi_{e,i}(a))$ for all $(g,a)\in A,$
where ${\rm tr}_{e,i}$ is the tracial state on $M_{R(i)}$ and ${\blue{\psi_{e,i}}}: F_1\to M_{R(i)}$
is the quotient map, $i=1,2,...,l.$}}

{{If  $h, h'\in (F_1)_+$ and
$\tau(h)=\tau(h')$ for all $\tau\in T(F_1),$
then, ${\rm Tr}_j\circ \pi\circ \phi_0$  and ${\rm Tr}_j\circ \pi_j\circ \phi_1$ are  traces of $F_1.$ It follows that
$$
{\rm Tr}_j(\pi_j(\phi_0(h)))={\rm Tr}_j(\pi_j(\phi_0(h'))\andeqn {\rm Tr}_j(\pi_j(\phi_1(h)))={\rm Tr}_j(\pi_j(\phi_1(h')),
$$
where $\pi_j: F_2\to M_{r(j)}$ is the quotient map and
${\rm Tr}_j$ is the standard trace on $M_{r(j)},$ $j=1,2,...,k.$
Furthermore, if
$h_1, h_2,..,h_m,, h_1', h_2',..., h'_{m'}\in (F_1)_+,$ and
$\sum_{n=1}^m\tau(h_n)=\sum_{s=1}^{m'} \tau(h_s')$ for all $\tau\in T(F_1),$ then
\beq\label{LgN1889-1}
\hspace{-0.2in}{\rm Tr}_j(\pi_j(\sum_{n=1}^m\phi_0(h_n))){\,=\,}{\rm Tr}_j(\pi_j\circ \phi_0(\sum_{n=1}^m(h_n)))
{\,=\,}{\rm Tr}_j(\pi_j\circ \phi_0(\sum_{s=1}^{m'}(h_s'))
{\,=\,}{\rm Tr}_j(\pi_j(\sum_{s=1}^{m'}\phi_0(h_s'))),
\eneq
$j=1,2,...,k.$}}
{{For exactly the same reason, we have
\beq\label{LgN1889-1+}
{\rm Tr}_j(\pi_j(\sum_{n=1}^m\phi_1(h_n)))
{\,=\,}{\rm Tr}_j(\pi_j(\sum_{s=1}^{m'}\phi_1(h_s'))),\,\,j=1,2,...,k.
\eneq
}}

{{Let $f\in {\rm LAff}_b(T(A))_+$ (see \ref{Aq}) and let $f_n\in \Aff(T(A))_+$ be such that
$f_n\nearrow f.$}}
{{Then
\beq\label{18815-s3-1}
\hspace{-0.2in}f(\tau_{0,i})&=&\lim_{n\to\infty}f_n(\tau_{0,i})=\lim_{n\to\infty}\frac1{R(i)}\sum_{j=1}^l a_{ij}f_n(\tau_{e,j})\cdot r(j)
=\frac1{R(i)}\sum_{j=1}^l a_{ij}f(\tau_{e,j})\cdot r(j)\\
\label{18815-s3-2}
\hspace{-0.2in}
&&\andeqn
f(\tau_{1,i})=\frac1{R(i)}\sum_{j=1}^l b_{ij}f(\tau_{e,j})\cdot r(j),\,\, i=1,2,...,k.
\eneq
}}
\end{NN}

\begin{prop}\label{2Rg13}
{\blue{Let $A$ be as in {Proposition} \ref{2Lg13}  and
let $u=(f,a)$ be a unitary in $A,$  where $a=(a_1,a_2,...,a_l)=e^{ih}\in \bigoplus_{j=1}^l M_{R(j)}=F_1$ and where
$h=(h_1, h_2, ..., h_l)\in (F_1)_{s.a.}.$

{\rm (1)} Then $uv^*\in U_0(A)$ (hence $[u]=[v]$ in $K_1(A)$),
where $v=(fe^{-i{g}}, 1_{F_1})$ and where
$$
g(t)=\left\{\begin{array}{lll}
             (1-2t) \phi_0(h),
            &  & 0\leq t\leq \frac{1}{2}, \\
                &   &  \\
               (2t-1) \phi_1(h), &  & \frac{1}{2}< t\leq 1.
            \end{array}
\right.
$$}}

{\blue{{\rm (2)} If ${\rm det}(u)(\psi)=1$ for all irreducible representations  $\psi$ of $A,$
then ${\rm det}(v)(\psi)=1$ for all $\psi.$}}

{\blue{{\rm (3)} Write $fe^{-ig}=(v_1, v_2,..., v_k),$ where $v_i\in U(C_0((0,1), F_2)^{\sim}).$
Then,
$[u]=[(s_1,s_2,...,s_k)]\in \Z^k/((\phi_0)_{*0}-(\phi_1)_{(0})(\Z^k),$
where $s_j$ is  the winding number of the map
$${[0,1] \ni t \mapsto} \mathrm{det}(v_j(t))\in \mathbb T\subset \C \,(j=1,2,...,k),$$
}}

{\blue{{\rm (4)} If ${\rm det}(u)(\psi)=1$ for all irreducible representations  $\psi$ of $A,$
then $u\in U_0(A).$}}
\end{prop}
\begin{proof}
Let $u=(f,a)\in U(A).$  Then $a\in U(F_1).$ Therefore, $a=e^{ih}$
for some
$h=(h_1, h_2, ..., h_l)\in \bigoplus_{j=1}^l M_{R(j)}=F_1.$   Define
\beq\label{1887-n1}
g(t)=\left\{\begin{array}{lll}
             (1-2t) \phi_0(h), &  & 0\leq t\leq \frac{1}{2}, \\
                &   &  \\
               (2t-1) \phi_1(h), &  & \frac{1}{2}< t\leq 1.
            \end{array}
\right.
\eneq
Then $(g, h)\in A.$  Let $v=(f_1, 1_{F_1})$, where $f_1(t)=f(t)e^{-ig(t)}$.
Let $U(s)=(f(t)e^{-ig(t)s}, e^{i h(1-s)})$ for $s\in [0,1].$ Then $\{U(s): s\in [0,1]\}$ is a continuous path
of unitaries in $A$ with $U(0)=u$ and $U(1)=(f_1, 1_{F_1}).$ Therefore $uv^*\in U_0(A).$
In particular, $[u]=[v] \in K_1(A).$  This proves (1).

To prove (2), suppose that ${\rm det}(u(\psi))=1$ for all irreducible representations $\psi$ of $A.$
Then ${\rm det}(e^{ih})=1.$ So
we may choose $h$ above so that ${\rm tr}_{e,j}(h_j)=0,$ where ${\rm tr}_{e,j}$ is the normalized tracial
state on $M_{R(j)},$ $j=1,2,...,l.$
Write $g(t)=(g_1(t),g_2(t),...,g_k(t))\in C([0,1], F_2).$ Let $\pi_m: F_2\to M_{r(m)}$ be
the quotient map.
Then $\pi_m\circ \phi_i: F_1\to M_{r(m)}$ is a \hm.
Write $\phi_i(h)=(g_1(i), g_2(i),...,g_k(i)),$ $i=0,1,$
where $g_m(i)\in M_{r(m)},$ $m=1,2,...,k.$
Then  ${\rm tr}_m(g_m(i))={\rm tr}_m(\phi_i(h)),$ $i=0,1,$  where ${\rm tr}_m$ is the normalized trace
on  $M_{r(m)},$ $m=1,2,...,k.$
Since $\phi_i$ is unital, ${\rm tr}_m\circ \phi_i$ is a trace on $F_1,$ $i=0,1$ and
$m=1,2,...,k.$ Thus ${\rm tr}_m(g_m(i))=0,$ $i=0,1$ and $m=1,2,...,k.$
It follows that (see \eqref{1887-n1}) ${\rm tr}_m(g(t))=0$ for all $t\in [0,1]$ and $1\le m\le k.$
In other words, ${\rm det}(e^{-g(t)})=1$ for all $t\in [0,1],$ as well as ${\rm det}(e^{ih})=1.$
This implies ${\rm det}(v(\psi))=1$ for all irreducible representations $\psi$ of $A.$

 For (3),
we may write  ${{fe^{-ig}}}=(v_1,v_2,...,v_k)$ with $v_m(0)=v_m(1)=1_{{M_{r(m)}}},$
where $v_m\in C([0,1], M_{r(m)}),$ $m=1,2,...,k.$
Thus, we may view   $v$ as an element {{of}} $C_0((0,1), F_2)^{\sim},$ the unitalization of the ideal $C_0((0,1), F_2)$
of $A,$  and write $v=(v_1, v_2,...,v_k)\in U(C_0((0,1), F_2)^{\sim})$
(as $v=(fe^{-ig}, 1_{F_1})$).
 Put $w'=({\rm det}(v_1(t)), {\rm det}(v_2(t)), ....,{\rm det}(v_k(t)))\in U(C([0,1])$
 and $w=\diag(d_1, d_2,...,d_k)\in  U(C_0([0,1], F_2)^{\sim}),$
 where $d_m=\diag(1_{r(m)-1}, {\rm det}(v_m(t)))\in U(C_0((0,1), M_{r(m)})^{\sim}),$ $m=1,2,...,k.$
 Then $w\in  U(C_0([0,1], F_2)^{\sim}).$ Note ${\rm det}(vw^*)(t)=1$ for all $t\in (0,1).$
 Then in $C_0([0,1], F_2)^{\sim},$ $vw^*\in U_0(C_0((0,1), F_2)^{\sim})\subset U_0(A).$

Let $(s_1,s_2,..., s_k)\in \Z^k=K_1(C_0((0,1), F_2))$, where $s_j$ is {{the}} winding number of the map
$${[0,1] \ni t \mapsto} \mathrm{det}(v_j(t))\in \mathbb T\subset \C.$$
In $K_1(C_0(0,1), F_2),$ $[v]=(s_1,s_2,..., s_k).$
Note  that (as $uv^*\in U_0(A),$ $[u]$ is the image of $[v]$ under the map from $K_0(F_2)=K_1(C_0((0,1), F_2))$ onto $K_1(A)$ given by
\eqref{sixterm}. Thus
$$
[u]=(s_1,s_2,..., s_k)\in \Z^k/(\phi_{1*}-\phi_{0*})(\Z^l)~.
$$

 Finally (4) follows from (3) and the facts that $A$ has stable rank one (by
 \ref{2pg3})  and stable rank one \CA s are $K_1$-injective by {\blue{\cite{Rf}}}.
\end{proof}

\begin{lem}\label{2Lg8}
Let $A=A(F_1, F_2, \phi_0, \phi_1)$ be as in Definition \ref{DfC1}. A unitary $u\in U(A)$ is in $CU(A)$ {{(see the definition \ref{Dcu})}} if, and only if, for each irreducible representation $\psi$ of $A$, one has $\mathrm{det}(\psi(u))=1$.
\end{lem}

\begin{proof}
One direction is obvious. {It} remains to show that the condition is also sufficient.
From
{{Proposition}} \ref{2Rg13} above, if $u\in U(A)$ with ${\rm det}(\psi(u))=1$ for all
irreducible representations $\psi,$ then  $u\in U_0(A).$

Write $F_1=M_{R(1)}\oplus {M_{R(2)}} \oplus \cdots \oplus M_{R(l)}$
and $F_2=M_{r(1)}\oplus M_{r(2)}\oplus \cdots \oplus M_{r(k)}.$

Then, since $u\in U_0(A),$ we may write $u=\prod_{n=1}^m \exp(i 2\pi h_n)$ for some  $h_n\in A_{s.a.},$ $n=1,2,...,m.$
{\blue{ We may write $h_n=(h_{nI},h_{nq})\in A,$
where $h_{nI}\in C([0,1], F_2)_{s.a.}$ and $h_{nq}\in (F_1)_{s.a.}$ with $\phi_i(h_{nq})=h_{nI}(i),$ $i=0,1.$}}
Write $h_{nq}=(h_{nq,1}, h_{nq,2},...,h_{nq,l}),$ where $h_{nq,j}\in (M_{R(j)})_{s.a.},$
$j=1,2,...,l.$
Let  $\pi_j^{F_1}: F_1\to M_{R(j)}$ denote  the quotient map and set
$\pi_{e,j}=\pi_j^{F_1}\circ \pi_e,$  $j=1,2,...,l.$

For any irreducible representation $\psi$ of $A$, {\blue{as ${\rm{det}}(\psi(u))=1,$ $\sum_{n=1}^m {\rm Tr}_\psi(\psi(h_n))\in \Z,$}} where
${\rm Tr}_\psi$ is the standard (unnormalized) trace on $\psi(A)\cong M_{n(\psi)}$ for some
integer $n(\psi).$ {By replacing each $h_n$ with $h_n+{\blue{J}}$ for a large enough integer ${\blue{J}}$, one may assume that { $h_n\in A_+$ and
$\sum_{n=1}^m {\rm Tr}_\psi(\psi(h_n))$} is positive.}
 Put $H_j(t)={\rm Tr}_j(\sum_{n=1}^m \pi_j(h_{nI}(t)))$ and $H(t)=({H_1(t), H_2(t),...,H_k(t)})$ for $t\in [0,1],$ where ${\rm Tr}_j$ is
 the standard trace on $M_{r(j)}$ and $\pi_j: F_2\to M_{r(j)}$ is the projection map, {for}  $j=1,2,...,k.$  {{Note that $H_j\in C([0,1])$ and $H_j(t)\in \Z$ for all $t\in [0,1],$
$j=1,2,...,k.$ It follows that
\beq\label{1888-n10}
H_j(t)=H_j(0)=H_j(1)\rforal t\in [0,1],\,\,\, j=1,2,...,k.
\eneq}}
 Put ${{b}}(\psi)={\blue{\sum_{n=1}^m {\rm Tr}_\psi(\psi(h_n))}}$ for {{$\psi=\pi_{e,j},$
 $j=1,2,...,l.$}}

Since ${\blue{\sum_{n=1}^m {\rm Tr}_\psi(\psi(h_n))}}\in \Z$ {{and is positive,}} there is a projection $p\in M_N(F_1)$ such that ${\blue{\mathrm{Tr}_{\pi_{e,j}}'{((\mathrm{id}_{M_N} \otimes \pi_j^{F_1})(p))}}}={{b}}(\pi_{e,j})$ for some $N\ge 1,$ {{where ${\rm Tr}_{\pi_{e,j}}'=Tr'_N\otimes {\rm Tr}_{\pi_{e,j}},$ and  where ${\rm Tr}_N'$ is the standard trace on $M_N.$}}
{{Write $p=(p_1, p_2,...,p_l),$ where $p_s\in M_N(M_{R(s)})$ is a projection, $s=1,2,...,l.$
Then ${\rm rank}(p_s)=\sum_{n=1}^m{\rm Tr}_{e,s}(h_{nq,s}),$  where
$h_{nq}=(h_{nq,1},h_{nq,2},...,h_{nq,l})\in F_1,$
where $h_{nq,s}\in (M_{r(s)})_+,$  and where ${\rm Tr}_{e,s}$ is
the standard trace on $M_{r(s)},$
$s=1,2,...,l.$}}
{{It follows that  from \eqref{LgN1889-1} that
\beq
&&{\rm Tr}_j'(\pi_j'\circ {(\mathrm{id}_{M_{\blue{N}}}\otimes\phi_0)}(p))=
\sum_{n=1}^m{\rm Tr}_j(\pi_j\circ \phi_0(h_{nq}))=H_j(0)\andeqn\\
&&{\rm Tr}_j'(\pi_j'\circ {(\mathrm{id}_{M_{\blue{N}}}\otimes\phi_1)}(p))=
\sum_{n=1}^m{\rm Tr}_j(\pi_j\circ \phi_1(h_{nq}))=H_j(1),
\eneq
$j=1,2,...,k,$ where
${\rm Tr}_j'$ is the standard trace on $M_N(M_{r(j)})$ and $\pi_j': M_N(F_2)\to M_N(M_{r(j)})$ is the projection map.}}
Then, {{by \eqref{1888-n10}, }}
$$
{\rm Tr}_j'(\pi_j'\circ {(\mathrm{id}_{M_{\blue{N}}}\otimes\phi_0)}(p))=
H_j(0)=H_j(1)={\rm Tr}_j'({\pi_j'(\mathrm{id}_{M_{\blue{N}}}\otimes\phi_1)}(p)),
$$
$j=1,2,...,k.$
{{It follows that $\pi_j'(\mathrm{id}_{M_{\blue{N}}}\otimes\phi_0)(p))$ and
$\pi_j'(\mathrm{id}_{M_{\blue{N}}}\otimes\phi_1)(p))$ have the same rank.
Therefore,}} there is a projection $P_j\in M_N(C([0,1], M_{r(j)}))$ such that $P_j(0)={\pi'_j}\circ {(\mathrm{id}_{M_n}\otimes\phi_0)}(p)$ and
$P_j(1)={\pi_j'}\circ {(\mathrm{id}_{M_n}\otimes\phi_1)}(p).$ Choose $P\in M_N(C([0,1], F_2))$ such that $\pi_j'(P)=P_j$ and {{put}}
$e=(P, p).$ Then $e\in M_N(A).$ Note that
since $P_j(t)$ has the same rank as $P_j(0)$ (and $P_j(1)$),
${\rm Tr}_j'(P_j(t))=H_j(t)=H_j(0)$ for all $t\in [0,1],$  $j=1,2,...,k.$
Consider
the continuous path $u(t)=\prod_{n=1}^m \exp(i 2\pi h_nt)$ for $t\in [0,1].$
Then  $u(0)=1,$  $u(1)=u,$ and
$$
\tau({du(t)\over{d{{t}}}}u^*(t))=i 2\pi \sum_{n=1}^m\tau(h_n)\rforal \tau\in T(A).
$$
But {{(see  Lemma 2.6 of \cite{Lnbirr}),}}  for  all $a\in A$ and $\tau\in T(A),$
$$
\tau(a)=\sum_{s=1}^l \af_s {{\rm tr}_{e,s} (\pi_{e,s}(a))}+\sum_{j=1}^k \int_{(0,1)} {{\rm tr}_j(\pi_j(a(t)))} d\mu_j(t),
$$
where $\mu_j$ is a Borel measure on $(0,1),$ {${\rm tr}_{e,s}$ is the tracial state on $M_{R(s)}$} and
{${\rm tr}_j$} is {the} tracial state on $M_{r(j)},$  and $\af_i\ge 0$ and
$\sum_{s=1}^l\af_s+\sum_{j=1}^k \|\mu_j\|=1.$
{{Note that
\beq
{\rm tr}_{e,s} (\pi_{e,s}(\sum_{n=1}^m h_n))&=&(1/R(s))({\rm Tr}_{\pi_{e,s}}(\pi_{e,s}(\sum_{n=1}^m h_n)))\\
&=&(1/R(s))b(\pi_{e,s})=
(1/R(s))({\rm Tr}_{\pi_{e,s}}'((\id_{M_N}\otimes \pi_s^{F_1})(p)).
\eneq
Thus (recall that we write $\tau(q)=(\tau\otimes Tr_N)(q)$ for $q\in M_N(A)$ as in \ref{Aq})
\beq\nonumber
&&\sum_{n=1}^m\tau(h_n)=\tau(\sum_{n=1}^m h_n)=
\sum_{s=1}^l \af_s {\rm tr}_{e,s} (\pi_{e,s}(\sum_{n=1}^m h_n))+\sum_{j=1}^k \int_{(0,1)} {\rm tr}_j(\pi_j(\sum_{n=1}^mh_{nI}(t))) d\mu_j(t)\\
&&=\sum_{s=1}^l \af_s {\rm tr}_{e,s}((\id_{M_N}\otimes \pi_s^{F_1})(p))+\sum_{j=1}^k \int_{(0,1)} (1/r(j))H_j(t)d\mu_j(t)\\
&&=\sum_{s=1}^l \af_s {\rm tr}_{e,s}((\id_{M_N}\otimes \pi_{e,s})(e))+\sum_{j=1}^k \int_{(0,1)} (1/r(j)){\rm Tr}_j'(P_j(t))d\mu_j(t)\\
&&=\sum_{s=1}^l \af_s {\rm tr}_{e,s}((\id_{M_N}\otimes \pi_{e,s})(e))+\sum_{j=1}^k \int_{(0,1)} {\rm tr}_j((\id_{M_N}\otimes \pi_j)(P(t)))d\mu_j(t)=\tau(e).
\eneq
}}
It follows that
$$
{1\over{2\pi i}}\int_0^1 \tau({du(t)\over{dt}}u(t)^*)dt=\tau(e)\tforal \tau\in T(A).
$$
In other words, $D_A(\{u(t)\})\in \rho_A(K_0(A))$ (see \ref{Dcu}).
It follows from {Lemma 3.1} of \cite{Thomsen-rims} and the fact that $A$ has stable rank one that $u\in CU(A)$ (see also
{{Corollary 3.11}} of \cite{GLX-ER}).
\end{proof}

The following  statement is known (see  \cite{Ell-RR0},  \cite{CE-AI},
and  \cite{Phil-cer}).

\begin{lem}\label{2Lgsphillips}
Let $u$ be a unitary in $C([0,1], M_n).$
Then, for any $\ep>0,$ there exist continuous functions
$h_j\in C([0,1])_{s.a.},$ {{$j=1,2,...,n,$}} such that
$$
\|u-u_1\|<\ep,
$$
where $u_1=\exp({{2}} i\pi H),$ $H=\sum_{j=1}^n h_jp_j$ {{where}} $\{p_1, p_2,...,p_n\}$ is a set of mutually orthogonal rank one projections in $C([0,1], M_n),$
and $\exp({{2}}i \pi h_j(t))\not=\exp({{2}} i\pi h_k(t))$ if $j\not=k$ for all
{{$t\in (0,1),$   and  $u_1(0)=u(0)$  and $u_1(1)=u(1).$}}

 Furthermore,
if ${\rm det}(u(t))=1$ for all $t\in [0,1],$ then
{\blue{$u_1$ can be chosen so}} that ${\rm det}(u_1(t))=1$ for all $t\in [0,1].$

\end{lem}
\begin{proof}
{{The statement}} follows from Lemma 1.5 of \cite{Phil-cer} with $m=2$ and $d=1,$
{{whose proof was inspired by the proof of Theorem 4 of \cite{CE-AI} for self adjoint elements (rather than unitaries).}}

\end{proof}

\begin{lem}\label{2Lg6}
Let $A=A(F_1, F_2, \phi_0, \phi_1)$ be as in Definition \ref{DfC1}. For any unitary $(f,a)\in U(A), ~\ep>0$, there is a unitary $(g,a)\in U(A)$ such that $\|g-f\|<\ep$ and, for each block $M_{r(j)} \sbs F_2= \bigoplus_{j=1}^k M_{r(j)}$, there are real valued functions $h^j_1,h^j_2,..., h^j_{{r(j)}}: [0,1]\to \R$  {{ such that
$g^j={{\sum_{i=1}^{r(j)} h^{j}_ip_i}}$ and $\{p_1, p_2,...,p_n\}$ is a set of mutually orthogonal rank one projections in $C([0,1], M_{r(j)})$ and $\exp({{2i \pi h^j_s}}(t))\not=\exp({{2i\pi h^j_{s'}}}(t))$ if ${{s\not=s'}}$ for all
$t\in (0,1).$ Moreover, if $(f,a)\in CU(A),$  one can choose $(g,a)\in CU(A).$}}
%
%
\end{lem}

\begin{proof}
{{By {{Lemma}} \ref{2Lgsphillips}, f}}or each unitary $f^j\in C([0,1], M_{r(j)})$, one can approximate $f^j$ by $g_1^j$ to within
{{$\ep$}} such that $g^j(0)=f^j(0)$, $g^j(1)=f^j(1),$ and for each $t$ in the open interval $(0,1)$, $g^j(t)$ has distinct eigenvalues.   {{If $(f,a)\in U(A),$ then  $(g,a)\in U(A),$ too.}} {{On combining {{Lemma}} \ref{2Lgsphillips} with  {{Lemma}} \ref{2Lg8}, the last statement also follows.}}
\end{proof}

\begin{rem}\label{2Rg7}
In 
{{Lemma \ref{2Lg6},}}
one may  ensure
that $h^j_1(0),h^j_2(0), ..., h^j_{{r(j)}}(0)\in [0,1),$  and that, for some $\dt\in (0,1),$ for all $t\in
(0,\dt)$,
\begin{equation}\label{eqn-G-ast}
\max \{ h^j_i(t);~ 1\leq i\leq {{r(j)}} \}-\min\{ h^j_i(t);~ 1\leq i\leq {{r(j)}} \}<1,
\end{equation}
and $h^j_{i_1}(t)\not= h^j_{i_2}(t)$ for $i_1\not= i_2$.  From the choice of $g^j$, we know that for any $t\in (0,1)$, $e^{2\pi i h^j_{i_1}(t)}\not= e^{2\pi i h^j_{i_2}(t)}$.  That is, $h^j_{i_1}(t) - h^j_{i_2}(t)\not\in \Z$.  This implies that \eqref{eqn-G-ast}
in fact holds for all $t\in (0,1)$. Hence, {{also,}} 
\begin{equation}\label{eqn-G-ast-01}
\max \{ h^j_i(1);~ 1\leq i\leq {{r(j)}} \}-\min\{ h^j_i(1);~ 1\leq i\leq {{r(j)}} \}\leq1.
\end{equation}
\end{rem}


\begin{lem}\label{2Lg9}
{{Let $A=A(F_1, F_2, \phi_0, \phi_1)$ be as in Definition \ref{DfC1}.}}
For any $u\in CU(A)$, one has $\mathrm{cer}(u)\leq 2+\ep$ and $\mathrm{cel}(u)\leq 4\pi$.
{Moreover, there exists a continuous path of unitaries $\{u(t): t\in [0,1]\}\subset CU(A)$  with length
at most $4\pi$ such that $u(0)=1_A$ and $u(1)=u.$}
\end{lem}

\begin{proof}
{{ Case (i): the case that}}
$u=(f,a)$ with $a=1_{{F_1}}$.
By
{{Lemma}} \ref{2Lg6},
 up to approximation to within an arbitrarily small  pre-specified
 number $\ep>0$, $u$ is unitarily equivalent to $v=(g,a)\in CU(A)$ with
{{$g=(g_1, g_2,..., g_k)\in C([0,1], F_2),$  where, for each $t\in (0,1),$ the unitary}}
$$g_j(t) =\diag{\big (}e^{2\pi i h^j_1(t)}, e^{2\pi i h^j_2(t)},..., e^{2\pi i h^j_{r(j)}(t)}{\big )}$$
has distinct eigenvalues.
(Note that since $f(0)=f(1)=1=g(0)=g(1)$, the unitary  intertwining the approximation of $u$ and $v$ can be chosen to be $1$ at $t=0, 1$, and therefore, the unitary is in $A=A(F_1, F_2, \phi_0, \phi_1)$.)  {{ Since $v\in CU(A),$ ${\rm det}(\psi(v))=1$ for every irreducible representation $\psi$ of $A$
(see Lemma \ref{2Lg8}).
Let $N\ge 1$ be an integer
such that every irreducible representation $\psi$ of $A$ has rank no more than $N.$ Then, we may assume that
$\|u-v\|<\ep<1/4N\pi.$ So $\|uv^*-1_A\|<\ep.$ Write $uv^*=e^{i h}$ for some $h\in A_{s.a.}$ with $\|h\|<2\arcsin(\ep/2)<1/2N.$
Note $uv^*\in CU(A).$ It follows that ${\rm det}(\psi(uv^*))=1$ for every irreducible representation $\psi$ of $A$ (see  Lemma \ref{2Lg8}).
Then, ${\rm Tr}_{\psi}(\psi(h))\in 2\pi \Z$ for every irreducible representation $\psi$ of $A,$
where ${\rm Tr}_{\psi}$ is the standard trace on $\psi(A)=M_{n(\psi)}$ (for some integer $n(\psi)\le N$).
Since $\|h\|<1/2N,$  $|{\rm Tr}_{\psi}(\psi(h))|<1/2.$ It follows that
${\rm Tr}_{\psi}(\psi(h))=0$ for every irreducible representation $\psi$ of $A.$
Define  $w(t)=e^{iht}v$ ($t\in [0,1]$). Since ${\rm Tr}_{\psi}(\psi(ht))=0$
for every irreducible representation $\psi$ of $A,$  ${\rm det}(\psi(w(t)))=1$ for each irreducible
representation $\psi$ of $A.$ It follows from Lemma \ref{2Lg8} that $w(t)\in CU(A)$ for all  $t\in [0,1].$
Note that  $w(0)=v$  and $w(1)=u.$
Moreover, the length of $\{w(t): t\in [0,1]\}$ is no more than $2\arcsin(\ep/2).$
Therefore, \wilog, we may assume that $u=v=(g,a)$ with $g$ as above.}}

Furthermore, one {{may}} assume {{that}}
$$h^j_1(0)={h^j_2(0)}=\cdots = h^j_{{r(j)}}(0)=0.$$
Since $\mathrm{det}(g_j(t))=1$ for all $t\in [0,1]$, one has  $h^j_1(t)+h^j_2(t)+ \cdots + h^j_{{r(j)}}(t)\in \Z$.

By the continuity of each $h^j_s(t)$,  {{It follows that}}
\begin{equation}\label{eq-G-0819-001}
\sum_{s=1}^{{r(j)}} h^j_s(t) = 0.
\end{equation}
Furthermore,  as $h^j_s(1)\in \Z$ (since $g_j(1)={1_{F_2}}$), we know that $h^j_s(1)=0$ for all $s\in \{1, 2, { ...}, {r(j)}\}$.  Otherwise, ${\displaystyle \min_s } {\big \{ }h^j_s(1) {\big \} }\leq -1$ and ${\displaystyle \max_s } {\big \{ } h^j_s(1){\big \} } \geq 1$ which implies that
\begin{eqnarray*}
{\displaystyle \max_s } {\big \{ } h^j_s(1){\big \} }-
{\displaystyle \min_s } {\big \{ }h^j_s(1) {\big \} } \geq 2,
\end{eqnarray*}
and this contradicts Remark \ref{2Rg7}.  That is, one has proved that $h={\big (} (h^1, h^2, \cdots, h^k), 0{\big )}$, where $h^j(t)=\diag {\big (} h^j_1(t), h^j_2(t), { ...}, h^j_{{r(j)}}(t) {\big )}$ is an element of $A=A(F_1, F_2, \phi_0, \phi_1)$ with $h(0)=h(1)={\blue{0}}$.
 As $g=e^{2\pi ih}$, we have $\mathrm{cer}(u)\leq 1+\ep$.
{We also have ${\rm tr}(h(t))=0$ for all $t\in [0,1].$}

It follows from \eqref{eq-G-0819-001} above and ${\displaystyle \max_s } {\big \{ } h^j_s(t){\big \} }-
{\displaystyle \min_s } {\big \{ }h^j_s(t) {\big \} } \leq 1~$ (see \eqref{eqn-G-ast-01} in Remark \ref{2Rg7}) that
\begin{equation*}
h^j_s(t)~ \sbs ~(-1, 1),\quad t\in [0,1],\ s=1,2, { ...}, {r(j)}~.
\end{equation*}
Hence $\|2\pi h\|\leq 2\pi,$ which implies $\mathrm{cel}(u) \leq 2\pi$.
Moreover, let $u(s)=\exp(i s2\pi h).$ Then $u(0)={1_{A}}$ and $u(1)=u.$  Since ${\rm tr}(s2\pi h(t))=2s\pi\cdot{\rm tr}(h(t))=0$
for all $t\in [0,1],$ {by Lemma \ref{2Lg8},} one has  $u(s)\in CU(A)$ for all $s\in [0,1].$

{{Case (ii):}} The general case.  {In this case} $a=(a^1, a^2, ..., a^l)$ with $\mathrm{det}(a^j)=1$ for $a^j\in F^j_1$.  So $a^j=\exp({2\pi i h^j})$ for $h^j\in F^j_1$ with $\mathrm{tr}(h^j)=0$ and $\|h^j\|<1$.  Define $H\in A(F_1, F_2, \phi_0, \phi_1)$ by
\begin{displaymath}
H(t)=
\left\{
\begin{array}{lcl}
  \phi_0(h^1, h^2, ... , h^l)\cdot (1-2t), & \textrm{if $0\leq t\leq \frac{1}{2}$}, \\
   & & \\
  \phi_1(h^1, h^2, ... , h^l)\cdot(2t-1),   & \textrm{if $\frac{1}{2}< t\leq 1$}.
\end{array}
\right.
\end{displaymath}
Note $H(\frac12)=0,~ H(0)= \phi_0(h^1, h^2, ..., h^l),$ and $ H(1)= \phi_1(h^1, h^2, ..., h^l),$ and therefore $H\in A(F_1, F_2, \phi_0, \phi_1)$.
{Moreover, ${\rm tr}(H(t))=0$ for all $t.$} Then $u'=u \cdot \exp({-2\pi i H})\in A(F_1, F_2, \phi_0, \phi_1)$ with
$u'(0)=u'(1)={1_{F_2}}$.  By {{Case (i),}} $\mathrm{cer}(u')\leq 1+\ep$ and $\mathrm{cel}(u')\leq 2\pi$,
and so  $\mathrm{cer}(u)\leq 2+\ep$ and $\mathrm{cel}(u)\leq 2\pi+2\pi \|H\| \leq 4\pi.$ Furthermore, we
note that $\exp(-2\pi s H)\in CU(A)$ for all $s$ as in {{Case (i)}}.
\end{proof}

\begin{thm}\label{FG-Ratn}
{\blue{Let $G$ be a subgroup of the ordered group $\Z^l$ (with the usual positive cone $\Z_+^l$).
Then  the semigroup $G_+=G\cap \Z^l_+$ is finitely generated.
In particular,  one has the following special case:}}

Let $A=A(F_1, F_2, \phi_1, \phi_2)$ be in ${\cal C}.$  Then $K_0(A)_+$ is finitely generated (by
its minimal elements); in other words, there are an integer $m\ge 1$ and finitely many
minimal projections of $M_m(A)$ such that these minimal projections generate the positive cone
$K_0(A)_+.$
\end{thm}
\begin{proof}
{\blue{Recall that an element $e\in G_+\setminus \{0\}$  is minimal,
if $x\in G_+\setminus\{0\}$ and $x\le e,$ then $x=e.$}}
We first show that
{\blue{$G_+\setminus\{0\}$}} has only finitely many minimal elements.

  Suppose otherwise that $\{q_n\}$ is an infinite set of minimal elements  of
  {\blue{$G_+\setminus\{0\}.$}}
Write $q_n=(m(1,n), m(2,n), ..., m({{l}}, n))\in \Z^l_+,$ where $m(i,n)$ are non-negative integers,
$i=1,2,...,l,$ and $n=1,2,....$ If there is an integer $M\ge 1$ such that
$m(i,n)\le M$ for all $i$ and $n,$ then $\{q_n\}$ is a finite set. So we may assume
that $\{m(i,n)\}$ is unbounded for some $1\le i\le l.$
There is a subsequence of $\{n_k\}$ such that
$\lim_{k\to\infty} m(i,n_k)=\infty.$ To simplify the notation, without loss of generality, we may assume
that  $\lim_{n \to\infty}m(i,n)=\infty.$ We may assume that, for some $j,$ {{the set}} $\{m(j,n)\}$ is bounded.
Otherwise, by passing to a subsequence,  we may assume that $\lim_{n\to\infty}m(i,n)=\infty$ for
all $i\in \{1,2,...,l\}.$
Therefore  $\lim_{n\to\infty}m(i,n)-m(i,1)=\infty.$
It follows that, for some $n\ge 1,$
$m(i,n)>m(i,1)$ for all $i\in\{1,2,...,l\}.$ Therefore $q_n\ge q_1,$ which contradicts the fact
that $q_n$ is minimal.
By passing to a subsequence, we may write $\{1,2,...,l\}=N\sqcup B$ in such {{a way}} that
$\lim_{n\to\infty} m(i,n)=\infty$ if $i\in N$ and
$\{m(i,n)\}$ is bounded if $i\in B.$ Therefore $\{m(j,n)\}$ has only finitely many different values {{for fixed}} $j\in B.$
Thus, by passing to a subsequence  again, we may assume that
$m(j,n)=m(j,1)$ for fixed  $j\in B.$ Therefore, for some $n>1,$
$m(i,n)>m(i,1)$ for all $n$ if $i\in N$ and $m(j,n)=m(j,1)$ for all $n$ if $j\in B.$
It follows that $q_n\ge q_1.$ This is impossible since $q_n$ is minimal.
This shows that ${\blue{G_+}}$ has only finitely many minimal elements.

To show that ${\blue{G_+}}$ is generated by these minimal elements,
fix an element $q\in {\blue{G_+\setminus\{0\}.}}$
If $q$ is not minimal, consider
the set of all elements of ${\blue{G_+\setminus\{0\}}}$
which are {{strictly}} smaller than $q.$
This set is finite.  Choose one which is minimal among them, say $p_1.$ Then $p_1$ is {{a}} minimal  element
of ${\blue{G_+\setminus\{0\}}},$
 {{as}} otherwise there is one smaller than $p_1.$  Since $q$ is not minimal, $q\not=p_1.$
Consider $q-p_1\in
{\blue{G_+}}\setminus\{0\}.$ If $q-p_1$ is minimal, then
$q=p_1+(q-p_1).$ Otherwise, we repeat the same argument to obtain a minimal element
$p_2\le q-p_1.$ If $q-p_1-p_2$ is minimal, then {{we have the decomposition}} $q=p_1+p_2+(q-p_1-p_2).$
Otherwise we repeat the argument again. This process is finite. Therefore $q$ is a finite sum of
minimal elements of $
{\blue{G_+}}\setminus \{0\}.$
%
%
%
%
\end{proof}



\begin{thm}\label{2Tg14}
The exponential rank of $A=A(F_1, F_2, \phi_0, \phi_1)$ is at most $3+\ep$.
\end{thm}

\begin{proof}

For each unitary $u\in U_0(A)$,
 {{as in Proposition \ref{2Rg13},
$u=(f, a)\in A$ where $a=e^{ih}$ for some $h\in (F_2)_{s.a.}.$
Therefore there is $x\in A_{s.a.}$ such that $\pi_e(x)=h$ and hence
$\pi_e(e^{ix})=e^{ih}=a.$}} Therefore one may
write $u=ve^{ix}$
for some $x\in A_{s.a.},$ where $v=(g,{1_{F_1}})$ with ${g}(0)={g}(1) ={1_{F_2}}\in F_2$.
  So we only need to prove the exponential rank of $v$ is at most $2+\ep$.  Consider $v$ as an element in $C_0((0,1), F_2)^\sim
$
which defines an element $(s_1,s_2,..., s_k)\in \Z^k=K_1(C_0((0,1), F_2))$.  Since $[v]=0$ in $K_1(A)$, there are $(m_1,m_2, ..., m_l)\in \Z^l$ such that
$$
(s_1, s_2,...,s_k)=((\phi_1)_{*0}-(\phi_0)_{*0})((m_1,m_2,...,m_l)).
$$
Note that
$$
(\phi_0)_{*0}(({R{(1)}, R{(2)},...,R{(l)}}))=(\phi_1)_{*0}(({R(1),R(2),...,R(l)}))=({r(1),r(2),...,r(k)})=
[{1}_{F_2}]\in K_0(F_2).
$$

Increasing $(m_1,m_2, ..., m_l)$ by adding a positive multiple of ${\big (}{R(1),R(2),...,R(l)}{\big )}$, we can assume $m_j\geq 0$ for all $j\in \{ 1, 2, ..., l\}$.  Let $a=(m_1P_1, m_2P_2, ..., m_lP_l)$, where
$$P_j=\left(
        \begin{array}{cccc}
          1 &   &   &   \\
            & 0 &  &  \\
           &  & \ddots &  \\
           &  &  & 0
        \end{array}
      \right)
      \in M_{ R(j)}\subset F_1.$$
Let $h$ be defined by
$$
h(t)=\left\{\begin{array}{lll}
              \phi_0(a)(1-2t), &  & 0\leq t\leq \frac12, \\
                &   &  \\
                \phi_1(a)(2t-1), &  & \frac12< t\leq 1.
            \end{array}
\right.
$$
Then $(h,a)$ defines a self-adjoint element in $A$.  One also has $e^{2\pi ih}\in C_0((0,1), F_2)^{\widetilde{}}$, since $e^{2\pi ih(0)}=e^{2\pi ih(1)}={{1_{F_2}}}$.  Furthermore, $e^{2\pi ih}$ defines
$(\phi_{1*}-\phi_{0*})((m_1,m_2,...,m_l))\in \Z^k$
as an element in $K_1(C_0((0,1), F_2))=\Z^k$.  Let $w=v e^{-2\pi ih}$.  Then $w$ satisfies $w(0)=w(1)={1_{F_2}}$ and $w\in {\widetilde{C_0((0,1), F_2)}}$ defines the element
$$(0,0,..., 0)\in K_1(C_0((0,1), F_2)).$$
{Up} to an approximation {{to}} within  a sufficiently small  $\ep$, {by Lemma \ref{2Lg6},} one {{may}} assume that $w=(w_1,w_2,..., w_k)$ with,  for all $j=1,2, ..., k$,

\begin{enumerate}
\item $w_j(t)=\diag{\big(}e^{2\pi ih^j_1(t)}, e^{2\pi ih^j_2(t)},...,e^{2\pi ih^j_{\blue{r(j)}}(t)}{\big )}$,
\item the numbers $e^{2\pi ih^j_1(t)}, e^{2\pi ih^j_2(t)},... ,e^{2\pi ih^j_{\blue{r(j)}}(t)}$ are distinct for all $t\in (0,1)$, and
%
%
\item $h_1^j(0)=h_2^j(0)=\cd= h_{\blue{r(j)}}^j(0)=0$.
\end{enumerate}
Since $w_j(1)={1_{F_2}}$, one has that  $h_i^j(1)\in \Z.$ 

On the other hand, the unitary $w$ defines
$$ h^j_1(1)+h^j_2(1)+\cd +h^j_{{r(j)}}(1)\in \Z\cong K_1(C_0((0,1), F_2^j))$$
which is zero {since $h_1^j(0)=h_2^j(0)=\cd= h_{\blue{r(j)}}^j(0)=0$}.  From \eqref{eqn-G-ast-01}, 
one has $h^j_{i_1}(1)-h^j_{i_2}(1)\leq 1$ {for $i_1 \neq i_2$}. This implies
$$ h^j_1(1)=h^j_2(1)=\cd=h^j_{\blue{r(j)}}(1)=0.$$ Hence, $h={\big(} (h^1, h^2, ..., h^k), 0 {\big )}$ defines a {{selfadjoint}}
element {{of}} $A$ and $w=e^{2\pi ih}$.
\end{proof}

\begin{NN}\label{2Rg15}
Let $A=A(F_1, F_2, \phi_0, \phi_1)\in {\cal C},$
where $F_1=M_{r(1)}\oplus M_{r(2)}\oplus\cdots \oplus M_{r(l)}$ and
$F_2=M_{R(1)}\oplus M_{R(2)}\oplus \cdots \oplus M_{R(k)}.$
{{By  \ref{2Rg11}, we may write
$(\phi_0)_{*0}=(a_{ij})_{k\times l}$ and $(\phi_1)_{*0}=(b_{ij})_{k\times l}.$
Denote by $\pi_i: F_2\to M_{r(i)}$ the projection map.}}
Let us calculate the Cuntz semigroup of $A$.

For any
$h\in (A\otimes {\cal K})_+$, {{consider the}} map $D_h:~ \mathrm{Irr}(A) \to \Z_+\cup\{ \infty\}$, for any $\pi\in \mathrm{Irr}(A)$
$$D_h(\pi)=\lim\limits_{n\to\infty}\mathrm{Tr}({{(\pi\otimes\id_{\cal K})}}(h^{\frac1n}))~,$$
where Tr is the {standard} unnormalized trace.  Then
we write $$D_h=(D_h^1,D_h^2,...,D_h^k,D_h(\rho_1),D_h(\rho_2),...,D_h(\rho_l)),$$
where $D_h^i=D_h|_{(0,1)_i}$ {{satisfy}} the following conditions:
\begin{enumerate}
\item $D_h^j$ is lower semi-continuous on each $(0,1)_j$,
\item $\liminf_{t\to 0}D_h^i(t) \geq \sum_j a_{ij} D_h(\rho_j){{={\rm rank}(\pi_i(\phi_0(h)))}}$ and\\
 $\liminf_{t\to 1}D_h^i(t) \geq \sum_j b_{ij} D_h(\rho_l){{={\rm rank}(\pi_i(\phi_1(h)))}}$.
\end{enumerate}
It is straightforward to verify that the image of the map $h\in (A\otimes {\cal K})_+ \to D_h$ is the subset of $\mbox{Map} (\mathrm{Irr}(A), \Z_+\cup\{ \infty\})$ consisting elements satisfying the above two conditions.

Note that $D_h(\pi)=\mbox{rank}((\pi\otimes \id_{\cal K})(h))$ for each $\pi\in \mathrm{Irr}(A)$ and $h\in (A\otimes {\cal K})_+$.

The following result is well known to experts (for example, see \cite{CEI-CuntzSG}).
{{We also keep the notation above in the following statement.}}
\end{NN}

\begin{thm}\label{2Tg16}
 Let $A=A(F_1, F_2, \phi_0, \phi_1)\in {\cal C}$ and let $n\ge 1$ be an integer.

\noindent {\rm (a)} The following statements are equivalent:
\begin{enumerate}
\item $h\in M_n(A)_+$ is Cuntz equivalent to a projection;
\item $0$ is an isolate{{d}} point in the spectrum of $h$;
\item  $D_h^{{i}}$ {{is}} continuous on each $(0,1)_{{i}}$, $\lim_{t\to 0}D_h^i(t) = \sum_j a_{ij} {{D_h(\rho_j)}},$
and $\lim_{t\to 1}D_h^i(t)= \sum_j b_{ij} {{D_h(\rho_j)}}$.
\end{enumerate}

\noindent {\rm (b)} For $h_1, h_2 \in M_n(A)_+$,  $h_1\lesssim h_2$
if and only if $D_{h_1}(\pi) \leq  D_{h_2}(\pi)$ for each $\pi\in \mathrm{Irr}(A)$. In particular, $A$ has {{strict}} comparison for positive elements.


\end{thm}

\begin{proof}

For part (b), \wilog, {by Subsection \ref{2Rg10}, we may assume that $h_1, h_2\in A$.}
Obviously, $h_1 {\lesssim} h_2$ implies $D_{h_1}(\phi)\leq D_{h_2}(\phi)$ for each $\phi\in {\rm Irr}(A)$. Conversely, assume that $h_1=(f,a)$ and $h_2=(g,b)$ satisfy $D_{h_1}(\phi)\leq D_{h_2}(\phi)$ for each $\phi\in {\rm Irr}(A)$.

First, there are strictly positive  {{functions}}
$s_1,s_2\in C_0((0,1])$
such that $s_1(a)$ and $s_2(b)$ are projections in $F_1.$
Note that $s_i(h_i)$ are Cuntz equivalent to $h_i$ ($i=1,2$).
By replacing $h_i$ by $s_i(h_i),$ \wilog\, we may assume that $a$ and $b$ are projections.

{{Let $\pi_i: F_2\to M_{R_i}$ be the projection map.}}
Fix $1/4>\ep>0.$ There exists $\dt_1>0$ such that
\beq\nonumber
\|f(t)-\phi_0(a)\|<\ep/64\,\,{\rm for\,\,all}\,\, t\in [0,2\dt_1]\andeqn \hspace{-0.05in}\|f(t)-\phi_1(a)\|<\ep/16\,\,{\rm for\,\,all}\,\, \hspace{-0.02in}t\in [1-2\dt_1,{{1].}}
\eneq
It follows that $f_{\ep/8}(f(t))$ is a projection in $[0,2\dt_1]$ and $[1-2\dt_1,1].$
Put $h_0=f_{\ep/8}(h_1).$  Note that $h_0=(f_{\ep/8}(f), a).$
 Then $D^i_{h_0}(\pi)$ is constant in $(0, 2\dt_1]\subset (0,1)_i$ and
$[1-2\dt_1,1]\subset (0,1)_i,$ $i=1,2,...,k.$
Choose $\dt_2>0$ such that
\beq\nonumber
&&D_{h_2}^i(t)\ge
{{{\rm rank}(\pi_i(\phi_0(b)))}}\rforal t\in (0,2\dt_2]_i\andeqn\\
&& D_{h_2}(t)\ge{{{\rm  rank}(\pi_i(\phi_1(b)))}}\rforal t\in [2\dt_2,1)_i,
\eneq
$i=1,2,...,k.$
Choose $\dt=\min\{\dt_1, \dt_2\}.$ Since $a\lesssim b$ in $F_1,$  there is a  unitary $u_e\in F_1$ such that
$u_e^*au_e=q\le b,$ where $q\le b$ is a {{subprojection}}.
Let $u_0=\phi_0(u_e)$ and $u_1=\phi_1(u_e)\in F_2.$ Then one can find  a unitary $u\in A$ such that
$u(0)=u_0$ and $u(1)=u_1.$
\Wlog, replacing $h_0$ by $u^*h_0u,$ we may assume that $a=q{{\le b}}.$

 Let $\bt: [0,1] \to [0,1]$ be a continuous function which is $1$ on the boundary and $0$ on $[\dt, 1-\dt]$. Let $h_3=(f_2, b-a)$ with $f_2(t)=\bt(t)\phi_0(b-a)$ for $t\in [0,\dt]$, $f_2(t)=\bt(t)\phi_1(b-a)$ for $t\in [1-\dt,\dt]$, and $f_2(t)=0$ for  $t \in[\dt, 1-\dt]$. Define $h_1'=h_0+h_3.$  Note that $h_1'$ has the form
 $(f',b)$ for some $f'\in C([0,1], F_2).$

  {Then,  for any  $\pi\in (\dt, 1-\dt)\subset (0,1)_i$,
 $$
 h_0\le h_1',\,\,\, D_{h_1'}^i(\pi)=D_{h_0}^i(\pi),$$
 for any $\pi\in (0,\dt)\subset (0,1)_i $,
 \begin{eqnarray*}
  {\rm rank}( h_1'(\pi)) & \le & {\rm rank} (h_0(t))+{\rm rank}(h_3(t)) \\
   & \le &  {\rm rank}(\pi_i(\phi_0(a)))+{\rm rank}(\pi_i(\phi_0(b-a))) \\
  & = & {\rm rank}(\pi_i(\phi_0(b))) \le  D_{h_2}^i(\pi),
  \end{eqnarray*}
  and for any $\pi\in (1-\dt,1)\subset (0,1)_i,$
  \begin{eqnarray*}
{\rm rank}( h_1'(\pi)) & \le & {\rm rank} (h_0(t))+{\rm rank}(h_3(t)) \\
& \le & {\rm rank}(\pi_i(\phi_1(a)))+{\rm rank}(\pi_i(\phi_1(b-a)))\\
 & = & {\rm rank}(\pi_i(\phi_1(b)))\le D_{h_2}^i(\pi).
 \end{eqnarray*} }
It follows that
\beq\label{N318-4}
D_{h_1'}(\pi)\le D_{h_2}(\pi)\rforal \pi\in {\rm Irr}(A).
\eneq
It follows from  (\ref{N318-4}) and  Theorem 1.1 of \cite{RL0} that
$h_1'\lesssim h_2$  in $C([0,1], F_2).$ Since $C([0,1], F_2)$ has stable rank one, there is a unitary
$w\in C([0,1], F_2)$ such that $w^*h_1'w=h_4\in \overline{h_2C([0,1], F_2)h_2}.$
Since $h_1'=(f', b),$ $h_2=(g,b),$ and  $b$ is a projection, $h_4(0)=\phi_0(b)$ and $h_4=\phi_1(b).$
In other words, $h_4\in A.$ In particular, $h_4\in \overline{h_2Ah_2}$ and $h_4\lesssim h_2.$
Note that $w^*\phi_i(b)w=\phi_i(b),$ $i=0,1.$ {{By}} using a continuous path of unitaries
which commutes {{with}} $\phi_i(b)$ ($i=0,1$) and connects to the identity, it is easy to find a sequence of unitaries $u_n\in C([0,1], F_2)$
{{with}}
$u_n(0)=u_n(1)=1_{F_2}$ such that
\beq\label{N318-5}
\lim_{n\to\infty}u_n^*h_1'u_n=h_4.
\eneq
{Notice that $u_n\in A$}.  We also have that $u_nh_4u_n^*\to h_1'.$ Thus, $h_1'\lesssim h_4.$ It follows that
\beq\label{N318-5-1}
f_{\ep}(h_1)\lesssim h_0\lesssim h_1'\lesssim h_4\lesssim h_2
\eneq
for all $\ep>0.$  This implies that  $h_1\lesssim h_2$ and part (b) follows.

To prove part (a), we note that (1) and (2) are obviously equivalent and both {imply} (3).
That  (3) implies (1) follows from the computation {{of}} $K_0(A)$ in \ref{2Lg13} and part (b).

\end{proof}

\begin{lem}\label{cut-full-pj}
Let $C\in \mathcal C$, and let $p\in C$ be a projection. Then $pCp\in \mathcal C.$
Moreover, if $p$ is full and $C\in\mathcal C_0$, then $pCp\in\mathcal C_0.$
\end{lem}
\begin{proof}
We may assume that $C$ is finite dimensional.
Write $C=C(F_1, F_2, \phi_0, \phi_1)$. Denote by $p_e=\pi_e(p)$, where $\pi_e: C\to F_1$ is the map defined in 
\ref{DfC1}.

For each $t\in [0,1],$ write
$\pi_t(p)=p(t)$ and ${\widetilde p}\in C([0,1], F_2)$ such that $\pi_t({\widetilde p})=p(t)$ for all $t\in [0,1].$
Then $\phi_0(p_e)=p(0)$, $\phi_1(p_e)=p(1),$ and
\begin{equation}\label{cfull-1}
pCp=\{(f,g)\in C: f(t)\in p(t)F_2p(t),
\andeqn g\in p_eF_1p_e\}.
\end{equation}

Put $p_0=p(0).$
There is a unitary $W\in C([0,1], F_2)$ such that
${{(W^*{\widetilde p}W)(t)}}=p_0$ for all $t\in [0,1].$
Define $\Phi: {\widetilde p}C([0,1], F_2){\widetilde p}\to
C([0,1], p_0F_2p_0)$ by $\Phi(f)=W^*fW$ for all $f\in  {\widetilde p}C([0,1], F_2){\widetilde p}.$
 Put
$F_1'=p_1F_1p_1$ and $F_2'=p_0F_2p_0.$
Define $\psi_0={\rm Ad}\, W(0) \circ \phi_0|_{F_1'}$ and $\psi_1={\rm Ad}\, W(1) \circ\phi_1|_{F_1'}.$ Put
$$
C_1=\{(f,g)\in C([0,1], F_2')\oplus F_1': f(0)=\psi_0(g)\andeqn f(1)=\psi_1(g)\},
$$
and note that $C_1\in\mathcal C$.
Define $\Psi: pCp\to C_1$ by
\beq\label{cfull-2}
\Psi((f,g))=(\Phi(f),g)\rforal f\in {\widetilde p}C([0,1], F_2){\widetilde p}\andeqn g\in F_1'.
\eneq
It is {{readily verified}} that $\Psi$ is an isomorphism.

If $p$ is full and $C\in {\cal C}_0,$  then, by Brown's theorem (\cite{Brown-Hereditary}), the hereditary {{\SCA\,}} $pCp$ is stably isomorphic to $C,$ and hence $K_1(pCp)=K_1(C)=\{0\}$; that is, $pCp\in \mathcal C_0$.
\end{proof}

The {{classes}} ${\cal C}$ and ${\cal C}_0$ are not closed under passing to quotients. However, we have the following
approximation rersult:

\begin{lem}\label{subapprox}
Any quotient of a \CA\, in ${\mathcal C}$ (or in ${\cal C}_0$)  can be locally approximated by \CA s in ${\cal C}$ (or in ${\cal C}_0).$
More precisely,
let  $A\in {\cal C}$ (or $A\in {\cal C}_0$), let $B$ be a quotient of $A,$ let ${\cal F}\subset B$ be a finite set, and let
$\ep>0.$ There exists a unital \SCA\, $B_0\subset B$ with $B_0\in {\cal C}$ (or $B_0\in {\cal C}_0$)  such that
\beq\nonumber
{\rm dist}(x, B_0)<\ep \tforal x\in {\cal F}.
\eneq
\end{lem}

\begin{proof}
Let $A\in\mathcal C.$ We may consider only those $A$'s which are minimal and are not finite dimensional.
 Let $I$ be an ideal of $A$. Write $A=A(E, F, \phi_0, \phi_1)$, where $E=E_1\oplus\cdots\oplus E_l$, ${{E_i}}\cong M_{k_i}$, and
{{where}} $F=F_1\oplus\cdots\oplus F_s$ with $F_j\cong M_{m_j}.$
Let $J=\{f\in C([0,1], F): f(0)=f(1)=0\}\subset A.$
As before,
we may write $[0,1]_j$ for the the spectrum of the $j$-th summand of $C([0,1],F_j),$ whenever it is convenient.
Put $\phi_{i,j}=\pi_j\circ \phi_i: E\to F_j,$ where $\pi_j: F\to F_j$ is the quotient map,
$i=0,1.$
Then  $A/I$ may  be written  (with a re-indexing) as
$$
\{(f,a): a\in {\widetilde E}, f\in \bigoplus_{1\le j\le s'} C({\tilde I}_j, F_j), f(0_j)={\widetilde \phi}_{0,j}(a), \,{\rm if}\,\, 0_j\in {\tilde I}_j,
f(1_j)={\widetilde \phi}_{1,j}(a), {\rm if}\,\, 1_j\in {\tilde I}_j\},\,\,\,{\bf(*)}
$$
where $s'\le s,$ ${{l'\le l}},$ ${\widetilde E}=\bigoplus_{i=1}^{l'} E_i$  and ${\widetilde \phi}_{i,j}=\phi_{i,j}|_{\widetilde{E}}$
{{and}} ${\tilde I}_j\subset [0,1]_j$ is a compact subset.

It follows from \cite{Freu} that,  {{for any $j,$}} there is a sequence of spaces $X_{n,j}$ which are
 finite disjoint unions of closed intervals
(including points) such that ${\tilde I}_j$ is the inverse limit of a sequence ${{(X_{n,j}, s_{n, n+1,j})}}$ and each  map
$s_{n,n+1,j}: X_{n+1, j}\to X_{n,j}$  is surjective and continuous.  Moreover, $X_{n,j}$ can be identified  with a disjoint union of closed subintervals of $[0,1].$
Let
$s_{n,j}: {\tilde I}_j\to X_{n,j}$  be the surjective  continuous map induced by the inverse limit system.

{{We may assume that ${\tilde I}_j=[0,1]_j,$ $j=1,2,...,t'\le s',$ and ${\tilde I}_j={\tilde I}_j^-\sqcup {\tilde I}_j^+,$
where ${\tilde I}_j^-\subset [0, t_j^-]$ and ${\tilde I}_j^+\subset [t_j^+,1]$ are compact subsets  for
some $0\le t_j^-<t_j^+\le 1,$ $t'<j\le s'.$
\Wlog\, (by applying \cite{Freu} to ${\tilde I}_j^-$ and to ${\tilde I}_j^+$), we may assume that,  for $j>t',$ each of the
disjoint closed interval in $X_{n,j}$ contains at most one of $s_{n,j}(0_j)$  and $s_{n,j}(1_j).$}}
Let ${{s_n^*}}: \bigoplus_{j=1}^{s'} C(X_{n,j}, F_j)\to \bigoplus_{j=1}^{s'}C({\tilde I}_j, F_j)$  be
the map induced by $s_{n,j}.$  Put $C_j=C(X_{n,j}, F_j).$
We may write $C({\tilde I}_j, F_j)=\overline{\bigcup_{n=1}^{\infty}s_{n,j}^*(C(X_{n,j}, F_j))}.$
Then, for all sufficiently large $n,$  for each $f\in {\cal F},$ there is $g\in \bigoplus_{j=1}^{s'} C_j$ such that
$g(s_{n,j}(0_j))=f(0_j),$ if $0_j\in {\tilde I}_j$ and $g(s_{n,j}(1_j))=f(1_j),$ if
$1_j\in {\tilde I}_j,$ and
\beq\label{1511-n2}
\|f|_{{\tilde I}_j}-s_{n,j}^*(g|_{X_{n,j}})\|<\ep/4.
\eneq
Note that $C_j$ is a unital \SCA\, of $C({\tilde I}_j, F_j),$ $j=1,2,...,s.$
Define
$$
C=\{(f, a):  f\in \bigoplus_{j=1}^{s'}C_j,\, a\in {\widetilde E}, f(s_{n,j}(0_j))={\widetilde \phi}_{0,j}(a),\,
{\rm if}\,\, 0_j\in {\tilde I}_j \andeqn f(s_{n,j}(1_j))={\widetilde \phi}_{1,j}(a), {\rm if}\,\,1_j\in {\tilde I}_j\}.
$$
Then $g$ in  (\ref{1511-n2}) is in $C$ and ${\cal F}\subset_{\ep} C.$
{{ Since each $C_j$ is a \SCA\, of $C({\tilde I}_j,F_j),$ one sees that
$C$ is a \SCA\, of $A/I.$
{{Note, $X_{n,j}$ is either equal to $[0,1]_j,$ or, $s_{n,j}(0_j)$ and $s_{n,j}(1_j)$ are in  different disjoint intervals of $X_{n,j}.$
It is then easy to check that $C\in {\cal C}$ (see \ref{DfC1}).}}
}}
This proves the lemma in the case that $A\in {\cal C}.$

Now suppose that $A\in {\cal C}_0.$
We will show that $K_1(A/I)=\{0\}.$
Consider the following six-term exact sequence:
$$
\xymatrix{
K_0(I) \ar[r] & K_0(A) \ar[r] & K_0(A/I) \ar[d] \\
K_1(A/I) \ar[u]^{\mathrm{\dt_1}} & K_1(A) \ar[l]  & K_1(I) . \ar[l]
}
$$
By {{Lemma}} \ref{2pg3}, $A$ and $A/I$ have stable rank one.  It follows from Proposition 4 of \cite{LR} that
$\dt_1=0.$ Since $K_1(A)=\{0\},$ it follows that  $K_1(A/I)=\{0\}.$

Denote by $\Pi: A\to A/I$ the quotient map.
Set ${\widetilde{F}}=\bigoplus_{i=1}^{s'} F_j$ and ${\widetilde{J}}=\Pi(J)=\Pi(C_0((0,1), F)).$
Note, since $K_1(C_0({\tilde I}_j^-\setminus \{s_{n,j}(0_j)\}, F_j))=\{0\}$ and
$K_1(C_0({\tilde I}_j^+\setminus \{s_{n,j}(1_j)\}, F_j))=\{0\}$ (for $t'<j\le s'$), one obtains
$K_1(J)=K_1(\bigoplus_{i=1}^{t'} C_0((0,1)_j, F_j))\cong K_0(\bigoplus_{i=1}^{t'} F_j).$
Consider the short exact sequence $0\to J\to A/I\to {\widetilde E}\to 0.$
Note that $K_1({\widetilde E})=0.$  Thus we have
$$
0\to K_0({\widetilde{J}})\to K_0(A/I)\to K_0({\widetilde E})\to K_1({\widetilde{J}})\to K_1(A/I)\to K_1({\widetilde E})=0.
$$
The fact that $K_1(A/I)=\{0\}$
implies that the map
$\oplus_{j=1}^{t'}(({\widetilde{\phi}}_{0,j})_{*0}-({\widetilde{\phi}}_{1,j})_{*1})$
from $K_0({\widetilde E})$ to $K_1({\widetilde{J}})\cong K_0(\bigoplus_{j=1}^{t'} F_j)$
 is surjective.

Now let $J_0={\widetilde{J}}\cap C.$ Then $K_1(J_0)=K_1({\widetilde{J}})\cong K_0(\bigoplus_j^{t'}F_j).$ The short exact sequence
$0\to J_0\to C\to {\tilde E}\to 0$ gives  the exact sequence
$$
0\to K_0(J_0)\to K_0(C)\to K_0({\widetilde E})\to K_1(J_0)\to K_1(C)\to K_1({\widetilde E})=0.
$$
The map from $K_0({\widetilde E})$ to $K_1(J_0)$ in the diagram above
is the same as $\oplus_{j=1}^{t'}(({\widetilde{\phi}}_{0,j})_{*0}-({\widetilde{\phi}}_{1,j})_{*1})$ which is surjective.
It follows that $K_1(C)=0.$
The lemma follows.
\end{proof}

We would like to return {{briefly}} to the beginning of this section by stating the following proposition which will not be used.

\begin{prop}[Theorem 2.15 of \cite{ENSTW}]\label{ASCAs}
Let $A$ be a unital \CA\, which is a subhomogeneous \CA\, with one dimensional spectrum.
Then, for any finite subset ${\cal F}\subset A$ and any $\ep>0,$ there exists  a unital  \SCA\, $B$ of $A$
 in the class ${\cal C}$ such that
\begin{equation*}
{\rm dist}(x, B)<\ep\tforal x\in {\cal F}.
\end{equation*}
\end{prop}

\begin{proof}
We use the fact that $A$ is an inductive limit of \CA s in ${\cal C}$ ({\blue{Theorem B and Definition 1.3 of}} \cite{ENSTW}).
Therefore, there is {{a}} \CA\, $C\in {\cal C}$ and a unital \hm\, $\phi: C\to A$ such that
\begin{equation*}
{\rm dist}(x, \phi(C))<\ep/2\rforal x\in {\cal F}.
\end{equation*}
Then we apply Lemma \ref{subapprox}.
\end{proof}

In \cite{Mg},  the following type of  unital  \CA s (dimension drop circle algebra) is studied:
\beq\label{dimdropcircle-d}
A=\{f\in C(\T, M_n): f(x_i)\in M_{d_i},\, i=1,2,...,N\},
\eneq
where $\{x_1, x_2,...,x_N\}$ are distinct points in $\T,$ and $d_i\in \N$ such that $m_id_i=n$ for some  integer $m_i>1.$

 {{The following fact will be used later in the paper.}}

\begin{prop}\label{dimcircle}
Every dimension drop circle algebra is in ${\cal C}.$
\end{prop}

\begin{proof}
Let $A$ be as in \eqref{dimdropcircle-d}. We identify $\T$ with the unit circle of the plane.
\Wlog, we may assume that $x_i\not=1\in \C,$ $i=1,2,...,N.$
Then,
$$A\cong \{f\in C([0,1], M_n): f(0)=f(1), f(x_i)\in M_{d_i},\,i=1,2,...,N\}.$$
We write $\bigoplus_{j=0}^{N}C([0,1], M_n)=\bigoplus_{j=0}^{N} C(I_j, M_n),$
where $I_j=[0,1]$ and we identify the end points of $I_j$ as $0_j$ and $1_j,$ $j=0,1,...,N.$
We may further write
\beq
A\cong \{f\in \bigoplus_{j=1}^{N+1} C(I_j, M_n): f(0_0)=f(1_{N}),\,\, f(1_j)=f(0_{j+1})\in M_{d_j},\,j=1,2,...,N\}.
\eneq
Let $h_i: M_{d_i}\to M_n$ be defined by $h_i(a)=a\otimes 1_{M_{m_i}}$ for all $a\in M_{d_i},$ $i=1,2,...,N.$
Consider $F_1=M_n\oplus (\bigoplus_{j=1}^NM_{d_i})$
and
$F_2=\bigoplus_{i=0}^{N+1}S_i,$ where $S_i\cong M_n.$
Let $\pi_0: F_1\to M_n$ and  $\pi_i: F_1\to M_{d_i},$ $i=1,2,...,N,$
{{denote}} the quotient maps. {{Also}} let  $\pi^i: F_2\to S_i$ {{denote the}}  quotient map, $i=0,1,...,N.$
Define $\phi_0: F_1\to F_2$ by $(\pi^0\circ \phi_0)|_{M_n}={\rm id}_{M_n},$
$(\pi^i\circ \phi_0)|_{M_n}=0,$  if $i\not=0,$ $(\pi^i\circ \phi_0)|_{M_{d_i}}=h_i,$ and $\pi^j\circ \phi_0|_{M_{d_i}}=0,$
$j\not=i,$ $i=1,2,...,N,$
and
define $\phi_1: F_1\to F_2$ by $(\pi^{N}\circ \phi_1)|_{M_n}={\rm id}_{M_n},$ where we identify
$S_N$ with $M_n,$ $(\pi^j\circ \phi_1)|_{M_n}=0,$ $j\not=N,$ and
$(\pi^{i-1}\circ \phi_1)|_{M_{d_i}}=h_i$ and $(\phi^j\circ \phi_1)|_{M_{d_i}}=0,$
if $j\not=i-1,$ $i=1,2,...,N.$ One checks that both $\phi_0$ and $\phi_1$ are unital.
Define
$$
B=\{(f,a)\in C([0,1], F_2)\oplus F_1: f(0)=\phi_0(a)\andeqn f(1)=\phi_1(a)\}.
$$
Let  $(f,a)\in B,$ where $f=(f_0, f_1, ...,f_N),$ $f_i\in C([0,1], M_n),$ $i=0,1,...,N,$
 and $a=(a_0, a_1,...,a_N),$ $a_0\in M_n,$ $a_i\in M_{d_i},$ $i=1,2,...,N.$
Then
$f_0(0)=\pi^0(f(0))=a_0$ and  $f_N(1)=\pi^N(f(1))=\pi^N(\phi_1(a))=a_0.$
So $f_0(0)=f_N(1),$
$f_0(1)=\pi^0(f(1))=\pi^0\circ \phi_1(a)=h_1(a_1),$ and  $f_1(0)=\pi^1(f(0))=\pi^1\circ \phi_0(a)=h_1(a_1).$
So $f_0(1)=f_1(0)=h_1(a_1);$ for $j=1,2,..., N,$
$f_j(1)=\pi^j(f(1))=\pi^j\circ \phi_1(a)=h_{j+1}(a_j)$ and $f_{j+1}(0)=\pi^{j+1}\circ \phi_0(a)=h_{j+1}(a_j).$
So $f_j(1)=f_{j+1}(0)=h_{j+1}(a_j),$ $j=1,2,...,N.$
Note that we may write $C([0,1], F_2)=\bigoplus_{i=0}^N C(I_j, M_n).$
Then,
$$
B\cong \{(f,a)\in (\bigoplus_{i=0}^N C(I_j, M_n)): f(0_0)=f(1_N),\, f(1_j)=f(0_{j+1}),\,\, j=1,2,...,N\}.
$$
Thus, $A\cong B.$ Since $B\in {\cal C},$ $A\in {\cal C}.$
\end{proof}

\begin{prop}\label{CUdist}
Let $A$ be a unital \CA\, with $T(A)\not=\emptyset.$
Suppose that $\overline{\rho_A^k(K_0(A))}\supset \rho_A(K_0(A))$ (see  Definition \ref{Dcu} for
the definition of $\rho_A^k(K_0(A))$). Then, for any pair $u, v\in U(M_k(A))$ with $uv^*\in U_0(M_k(A)),$
\beq
\label{Cu181226-l1}
{\rm d}(\overline{u},\overline{v})={\rm dist}(\overline{u}, \overline{v})
\eneq
(see \ref{Dcu}). Moreover, if $A$ also has stable rank one, then \eqref{Cu181226-l1} holds for
any
$u, v\in U(A).$
In particular, this is true
if $A\in {\cal C}.$
\end{prop}

\begin{proof}
If $u, v\in U(M_k(A))$ and $uv^*\in U_0(M_k(A)),$ then by \ref{Dcu},
\beq
{\rm d}(\overline{u}, \overline{v})={\rm d}(\overline{uv^*},\overline{1})={\rm dist}(\overline{uv^*}, \overline{1})=
{\rm dist}(\overline{u}, \overline{v}).
\eneq
Suppose $A$ also has stable rank one.
Let $u, v\in U(A)$ with  $uv^*\in U_0(A).$
Let ${\rm d}(\overline{u}, \overline{v})=\dt<2.$
Let $u_1=\diag(u, 1_{k-1})$ and $v_1=\diag(v,1_{k-1}).$
By Corollary 2.11 of \cite{GLX-ER}, for example, there is $a\in A_{s.a.}$ such that
$\tau(a)=D_A(\zeta)(\tau)$ for all $\tau\in T(A)$ for
some {{piecewise}} smooth path  $\zeta$ in $U_0(M_k(A))$
with $\zeta(0)=u_1v_1^*$ and $\zeta(1)=1_k,$ and $u_1v_1^*=\diag(\exp(i 2\pi a), 1_{k-1})z$ for some $z\in CU(M_k(A)).$
We may assume that $\omega(a)=\sup\{|\tau(a)|: \tau\in T(A)\}<1/2$ as $\dt<2.$
Since $A$ has stable rank one, by  Corollary 3.11 of \cite{GLX-ER}, $uv^*\exp(-i2\pi a)\in CU(A).$
This implies that
\beq
{\rm dist}(\overline{u}, \overline{v})={\rm dist}(\overline{uv^*}, \overline{1})={\rm dist}(\overline{\exp(-i2\pi a)}, \overline{1})=
{\rm dist}(\overline{\exp(i2\pi a)}, \overline{1}).
\eneq
By the first paragraph of the proof of Lemma 3.1 of \cite{KTm}, for any $\ep>0,$
there is $h\in A_0$ such that $
\omega(a)\le \|a-h\|\le \omega(a)+\ep<1/2.$  Let $b=a-h.$  Then, since $\exp(ih)\in CU(A)$ (see  Lemma 3.1 of \cite{Thomsen-rims}),
${\rm dist}(\overline{\exp(i2\pi a)}, \overline{1})={\rm dist}(\overline{\exp(i 2\pi b)},\overline{1})\le \|\exp(i 2\pi (\omega(a)+\ep))-1\|.$
Since this holds for any $\ep>0,$ we conclude that
\beq
{\rm dist}(\overline{u}, \overline{v})={\rm dist}(\overline{\exp(i2\pi a)}, \overline{1})\le  \|\exp(i 2\pi \omega(a))-1\|.
\eneq
This holds for any choice of $\zeta$ as above,  and therefore,
\beq
{\rm dist}(\overline{u}, \overline{v})\le \|\exp(i 2\pi  \|\overline{D}_A(uv^*)\|)-1\|={\rm d}(\overline{u}, \overline{v}).
\eneq
On the other hand,  by definition,
${\rm d}(\overline{u},\overline{v})={\rm d}(\overline{\exp(i2\pi a)}, \overline{1})=\dt.$
If ${\rm dist}((\overline{u},\overline{v})<\dt,$ then there is $z_1\in CU(A)$
such that $\|uv^*z_1-1\|<\dt<2.$  There exists $b_1\in A_{s.a.}$ such that $\|b_1\|<2\arcsin(\dt/2)$
and  $uv^*z_1=\exp(i2\pi b_1).$  It follows that $\tau(b)\le 2\arcsin(\dt/2)$ for all $\tau\in T(A).$
This implies that ${\rm d}(\overline{\exp(i2\pi b_1)}, \overline{1})<\dt.$ Note that $[z_1]=0$ in $K_1(A).$
Since $A$ has stable rank one, $z_1\in U_0(A).$ Therefore, by Lemma 3.1 of \cite{Thomsen-rims},
$\overline{D}_A(\overline{z_1})=0.$
It follows that
\beq
\overline{D}_A(\overline{\exp(i2\pi b_1}))=\overline{D}(\overline{uv^*z_1})=\overline{D}_A(\overline{uv^*})+\overline{D}_A(\overline{z_1})=\overline{D}_A(\overline{uv^*}).
\eneq
Hence  ${\rm d}(\overline{uv^*}, \overline{1})=\|\exp(i 2\pi \|\overline{D}_A(\overline{\exp(i2\pi b_1})\|)-1\|\le \|\exp(i2\pi\|b_1\|)-1\|<\dt.$
This contradicts the fact that ${\rm d}(\overline{u},\overline{v})=\dt.$
Therefore,
\beq
{\rm d}(\overline{u}, \overline{v})={\rm dist}(\overline{u}, \overline{v}).
\eneq

{{Now,}} $uv^*\not\in U_0(A),$ then, {{by Definition (see \ref{Dcu})}}, ${\rm d}(\overline{u}, \overline{v})=2.$
Suppose that $A$ has stable rank one. Since each element  of  $CU(A)$ gives the zero element
of $K_1(A),$  $CU(A)\subset U_0(A).$
Suppose that  $\|uv^*-z_2\|<2$ for some $z_2\in  CU(A).$   Then
$\|uv^*z_2^*-1\|<2.$ It follows that $uv^*z_2^*\in U_0(A).$ Therefore $uv^*\in U_0(A),$
 a contradiction. This shows that
${\rm dist}(\overline{u}, \overline{v})={\rm dist}(\overline{uv^*},1)=2.$

To see the last part of the statement, let $A\in {\cal C}.$
We note that, by Proposition \ref{2pg3}, $A$ has stable rank one.
Moreover, by Theorem \ref{FG-Ratn}, there exists $k\ge 1$ such that $M_k(A)$ contains
a set of projections  whose images in $K_0(A)$ generate
$K_0(A).$
It follows
that $\overline{\rho_A^k(K_0(A))}=\overline{\rho_A(K_0(A))}.$
Thus, the lemma {{applies in this case}}.

\end{proof}



\section{Maps to finite dimensional \CA s}

\begin{lem}\label{8-N-0}
Let $z_1, z_2,..., z_n$ {be} positive integers which may not be distinct.
{\blue{Set}} $T=n\cdot \max\{z_iz_j: 1\le i,j\le n\}.$
{\blue{Then,}}
for any two nonnegative integer linear {combinations}  $a=\sum_{i=1}^n a_i\cdot z_i$ and $b=\sum_{i=1}^n b_i\cdot z_i$, there are two {\blue{integer}} combinations $a'=\sum_{i=1}^n a'_i\cdot z_i$ and $b'=\sum_{i=1}^n b'_i\cdot z_i$ with $a'=b'$, $0\leq a'_i\leq a_i$, $0\leq b'_i\leq b_i$, and ${\rm min}\{a-a', b-b'\}\leq T$.

{Consequently, if $\dt>0$ and $|a-b|<\dt,$ we also have   $\max\{a-a', b-b'\}<\dt+T.$}
\end{lem}

\begin{proof}
To prove the first part,
{\blue{note first that
if $\min\{a, b\}\le T,$ we may choose $a'=b'=0.$ Therefore}}  it is enough to prove that if  ${\blue{\min\{a,b\}}} >T,$
then there are {{non-zero}} $0<a'=\sum_{i=1}^n a'_i\cdot z_i=b'=\sum_{i=1}^n b'_i\cdot z_i$ with $0\leq a'_i\leq a_i$, $0\leq b'_i\leq b_i,$
{\blue{and $\min\{a-a', b-b'\}\le T.$}}

{\blue{Suppose that $a,\, b>T.$
Then there is $i'$ such that $a_{i'}z_{i'}> \max_{i,j}z_iz_j,$ and there is $j'$ such that $b_{j'}z_{j'}\ge \max_{i,j}z_iz_j.$
Thus, $a_{i'}z_{i'}\ge z_{i'}z_{j'}$ and $b_{j'}z_{j'}\ge z_{i'}z_{j'}.$ It follows
that $a_{i'}\ge z_{j'}$ and $b_{j'}\ge z_{i'}.$}}
{Then choose} $a^{(1)'}=a^{(1)}_{i'}z_{i'}$ and $b^{(1)'}=b^{(1)}_{j'}z_{j'}$ {with $a^{(1)}_{i'}=z_{j'}$ and $b^{(1)}_{j'}=z_{i'}.$}
{\blue{Then $a^{(1)'}=b^{(1)'}\ge 1$ and $a_{i'}^{(1)}\le a_{i'}$ and $b_{j'}^{(1)}\le b_{j'}.$
Moreover, $\min\{a-a^{(1)'}, b-b^{(1)'}\}\le \min\{a, b\}-1.$
}}
{If ${\rm min}\{a-a^{(1)'},b-b^{(1)'}\}\le T,$ {\blue{then we are done, by choosing
$a^{(1)}_i=0$ if $i\not=i'$ and $b^{(1)}_j=0,$ if $j\not=j'.$ In particular, the lemma holds if $T-\min\{a, b\}\le 1.$ }}
If {\blue{$T-{\rm min}\{a-a^{(1)'},b-b^{(1)'}\}> 0,$}}  we repeat this on ${\blue{a-{a^{(1)'}}}}$ and $b-{\blue{b^{(1)'}}}$ and obtain
$a^{(2)}\le a-a^{(1)}$ and $b^{(2)}\le b-b^{(1)}$ such that $a^{(2)}=b^{(2)}{\blue{\ge 1}}.$  Put $a^{(2)'}=a^{(1)'}+a^{(2)}$ and
$b^{(2)'}=b^{(1)'}+b^{(2)}.$ Note we also have $a^{(2)'}=\sum_ia_i^{(2)}z_i$ and $b^{(2)'}=\sum_ib^{(2)}_iz_i$ with
$0\le a_i^{(2)}\le a_i$ and $0\le b_i^{(2)}\le b_i$ for all $i.$  Moreover,
${\blue{a^{(2)'}=b^{(2)'}\ge a^{(1)'}+1\ge 2}}.$
{\blue{It follows that  $\min\{a-a^{(2)'}, b-b^{(2)'}\}\le \min\{a, b\}-2.$}}
If $\min\{a-a^{(2)'}, b-b^{(2)'}\}\le T,$ then we are done. {\blue{ In particular,  the lemma holds for
$T-\min\{a, b\}\le 2.$}}
{\blue{If $T-\min\{a-a^{(2)'}, b-b^{(2)'}\}>0,$}} we continue.  An inductive argument
{\blue{{{establishes}}}} the first part of  the lemma.}

To see the second part, assume that $a-a'\le T.$ Then $b-b'<|a-b|+T.$
\end{proof}

\begin{thm}
[2.10 of \cite{Lin-AU11}; see also Theorem 4.6 of \cite{Lntams07}  and 2.15 of \cite{Li-interval}]\label{Oldthm201408}
Let {\blue{$A=PC(X, F)P,$}} where $X$ is  a compact metric space,
$F$ is a finite dimensional \CA\,, and
$P\in C(X,F)$ is a projection,  and let $\Delta:  A_+^{q, \bf 1}\setminus \{0\}\to (0,1)$ be an order preserving map.

For any $\ep>0,$ any finite set ${\cal F}\subset A,$   there exist a finite set
${\cal H}_1\subset A_+^{\bf 1}\setminus \{0\},$  a finite subset
${\cal P}\subset \underline{K}(A),$ a finite set ${\cal H}_2\subset A_{s.a.},$ and  $\dt>0$
with the following {{property}}:
If $\phi_1, \phi_2: A\to M_n$ (for some integer $n\ge 1$) are two unital \hm s such that
\begin{eqnarray*}
&&[\phi_1]|_{\cal P}=[\phi_2]|_{\cal P},\,\,\,
{\rm tr}\circ \phi_1 (h)\ge \Delta(\hat{h})\tforal h\in {\cal H}_1, {{\tand}}\\
&&|{\rm tr}\circ \phi_1(g)-{\rm tr}\circ \phi_2(g)|<\dt\tforal g\in {\cal H}_2,
\end{eqnarray*}
{\blue{where ${\rm tr}$ is the tracial state of $M_n,$}}
then there {\blue{exists}} a unitary $u\in M_n$ such that
\beq\label{8-NN-1}
\|{\rm Ad}\, u\circ \phi_1(f)-\phi_2(f)\|<\ep\tforal f\in {\cal F}.
\eneq
\end{thm}

\begin{proof}

{\blue{We first prove this for the case $A=C(X).$ This actually follows (as we shall show) from Theorem 2.10 of \cite{Lin-AU11}.
Fix $\ep>0$ and a finite subset ${\cal F}.$  Let $\eta>0$ be {\blue{as}}
in Theorem 2.10 of \cite{Lin-AU11}.

Let $x_1, x_2,...,x_K\in X$ be a finite subset such that
$X\subset \bigcup_{i=1}^KB(x_i, \eta/8).$  Set ${\cal H}_1=\{f_i: 1\le i\le K\},$
where $f_i\in C(X),$ $0\le f_i\le 1$ and the support of $f_i$ lies in $B(x_i, \eta/4)$ and
$f_i(x)=1$ if $x\in B(x_i, \eta/8),$ $i=1,2,...,K.$
Choose $\sigma=\min\{\Delta(\hat{f_i}): 1\le i\le K\}>0.$
Suppose that ${\rm tr}(\phi_1(f_i))\ge \Delta(\hat{f_i}),$ $i=1,2,...,K.$
Let $O_r$ be an open ball of $X$ with radius at least $r\ge \eta$ and with center $x\in X.$
Then $x\in B(x_i, \eta/8)$ for some $i.$  It follows that $B(x_i, \eta/4)\subset O_r.$
Therefore,
\beq
\mu_{{\rm tr}\circ \phi_1}(O_r)\ge {\rm tr}\circ \phi_1(f_i)\ge \sigma.
\eneq
Now let $\gamma>0,$ ${\cal P}\subset \underline{K}(A),$ and ${\cal H}_2\subset A_{s.a.}$ (in place of ${\cal H}$)  be as given by
2.10 of \cite{Lin-AU11} for the above $\ep,$ ${\cal F},$ $\eta$ and $\sigma$ (note we do not need to mention ${\cal G}$ and
$\dt$ since $\phi_1$ and $\phi_2$ are \hm s).

Choose $\dt=\gamma,$ ${\cal P},$ and ${\cal H}_2$ as above.   Then Theorem 2.10 of \cite{Lin-AU11} applies.
This proves that the theorem holds for $A=C(X).$
The case that $A=M_m(C(X))$  follows easily.
By considering each summand separately,  it is also easy to see that  the theorem also holds for $A=C(X, F).$

Now let $A=PC(X, F)P.$ Again, by considering each summand separately, we may reduce the general
case to the case that $A=PC(X, M_r)P$ for some integer $r\ge 1.$

 Note that $\{{\rm rank}(P(x)): x\in X\}$ is a finite set of positive integers.
Therefore the set $Y=\{x\in X: {\rm rank}(P(x))>0\}$ is compact and open.  Then we may write  $A=PC(Y,M_r)P.$
Thus, \wilog,  we may assume that $P(x)>0$ for all $x\in X.$ {\blue{With this assumption, we may assume that}} $P$ is a full projection {{in}} $C(X, M_r).$
Then,  by \cite{Brown-Hereditary}, $A\otimes {\cal K}\cong C(X,F)\otimes {\cal K}.$
It follows that there is an integer $m\ge r$ and a full projection $e\in M_m(A)$ such that
$eM_m(A)e\cong C(X,M_r).$  In particular, $e$ has rank $r$ everywhere.   Moreover, there {{exist}} an integer $m_1\ge 1,$ a full projection $e_A$ in  $M_{m_1}(eM_m(A)e)$, and a unitary $u\in M_{m_1m}(A)$
such that $u^*e_Au=1_A$ (we identify $1_A$ with $1_A\otimes e_{11}$ in $M_{m_1m}(A)$).
In particular, we may write $u^*e_A(M_{m_1}(eM_m(A)e))e_Au=A$ (where we identify  $A$ with $A\otimes e_{11}$).

Put $B=M_{m_1m}(A).$   Since $e\in M_m(A),$ $e\in B.$ Moreover, $eBe=eM_m(A)e$
and $M_{m_1}(eBe)\subset M_{m_1}(eM_m(A)e)\subset M_{m_1m}(A)=B.$
Put $e_1=\diag(\overbrace{e,e,...,e}^{m_1}).$ Then we may identify $M_{m_1}(eBe)$ with $e_1M_{m_1m}(A)e_1.$
Set $B_1=M_{m_1}(eBe)=e_1M_{m_1m}(A)e_1\subset B.$  Note $B_1\cong M_{m_1}(C(X, M_r)).$
Put $B_2=e_AM_{m_1}(eM_m(A)e)e_A=e_AB_1e_A.$ So $e_A\le e_1$ and
$B_2=e_ABe_A.$  Note $e_A$ is full in $B_1.$
Define $\Delta_1: B_+^{q, {\bf 1}}\setminus \{0\}\to (0,1)$ by $\Delta_1(\widehat{(a_{i,j})})=(1/m_1m)\sum_{i=1}^{m_1m}\Delta(\hat{a_{ii}})$
for
$0\le (a_{i,j})\le 1$ in $B_1^{\bf 1}\setminus\{0\}.$

Fix $\ep>0$ and a finite subset ${\cal F}\subset A.$
\Wlog, we may assume that $\|a\|\le 1$ for all $a\in {\cal F}.$
Put ${\cal F}_1=\{ufu^*: f\in {\cal F}\}\cup \{e_A\}.$
Let $\ep_1=\min\{\ep/15, 1/3\cdot 64\}.$
Let ${\cal P}_1\subset \underline{K}(B_1)$ (in place of ${\cal P}$),
${\cal H}_1'\subset B_1\setminus \{0\}$ (in place of ${\cal H}_1$),
and
let ${\cal H}_2'\subset (B_1)_{s.a.}$ (in place of ${\cal H}_2$), and $\dt'>0$ (in place of $\dt$) {{denote the finite subsets and constant provided}} by
the theorem for the case $A=C(X, M_r)$
for $\Delta_1,$ $\ep_1,$ and ${\cal F}_1$ (and for $B_1$ instead of $A$).
\Wlog, we may assume that $[e_A]$ and $[e_1]\in {\cal P}_1.$

Notice that $B_1\subset M_{m_1m}(A).$ Let ${\cal H}_1\subset A_+^{\bf 1}\setminus \{0\}$ be a finite subset such that
a finite subset of the  set $\{ (b_{i,j})_{m_1m\times m_1m}\in B_1\setminus\{0\}: b_{ii}\in {\cal H}\}$
contains ${\cal H}_1.$  We may view ${\cal P}_1$ as a finite subset of $\underline{K}(A)$ since $A\otimes {\cal K}\cong
B_1\otimes {\cal K}.$  Let ${\cal H}_2$ be a finite subset of $A$ such
that $\{(a_{ij}): a_{ii}\in {\cal H}_2''\}\supset {\cal H}_2'.$ Note if $(a_{i,j})\in M_{m_1m}(A)_{s.a.},$ then
$a_{ii}\in A_{s.a.}$ for all $i.$
Choose $\dt=\dt'.$

Now suppose that $\phi_1, \phi_2: A
\to M_n$ (for some integer $n\ge 1$) are two unital \hm s which satisfies the assumption for the above ${\cal H}_1, {\cal H}_2,$
${\cal P}_1$ (in place of ${\cal P}$) and $\dt.$

Set $\phi_i^\sim=(\phi_i\otimes {\id}_{m_1m}).$
Define  ${\tilde \phi}_i=\phi_i^\sim |_{B_1}: B_1\to M_{m_1mn},$  $i=1,2.$
Note that
\beq\label{s418109-00}
\phi_i(a)=\phi_i^\sim (u^*){\tilde \phi}_i(uau^*)\phi_i^\sim (u)  \rforal a\in A,\,\,i=1,2.
\eneq
Then,
\beq\label{s418109-0}
[{\tilde \phi}_1]|_{{\cal P}_1}=[{\tilde \phi}_2]|_{{\cal P}_1}.
\eneq
In particular, $[{\tilde \phi}_1(e_1)]=[{\tilde \phi}_2(e_1)].$
Therefore,  replacing ${\tilde \phi}_1$ by ${\rm Ad}\, w\circ {\tilde \phi}_1$ with
a unitary $w$ in $M_{m_1mn},$  \wilog, we may assume  ${\tilde \phi}_1(e_1)={\tilde \phi}_2(e_1).$
Put $F={\tilde \phi}_1(e_1)M_{m_1mn}{\tilde \phi}_1(e_1)\cong M_{m_1rn}.$

Then, for any $h=(h_{i,j})\in {\cal H}_1',$
\beq\label{s418109-1}
t({\tilde \phi}_1(h))=(1/m_1m)\sum_{i=1}^{m_1m}\tau(\phi_1(h_{i,i}))\ge \Delta_1(\hat{h}),
\eneq
where $\tau$ is the tracial state of $M_n$ and $t$ is the tracial state of $M_{m_1m}(M_n).$
Denote by ${\bar t}$ the tracial state of $F.$
Then, by \eqref{s418109-1},
\beq\label{s418109-2}
{\bar t}({\tilde \phi}_1(h)\ge \Delta_1(\hat{h})\rforal h\in {\cal H}_1'.
\eneq
If $h=(h_{i,j})\in {\cal H}_2',$
then
\beq
|t({\tilde \phi}_1(h))-t({\tilde \phi}_2(h))|&=&(1/m_1m)\sum_{i=1}^{m_1m}|\tau(\phi_1(h_{i,i})-\tau(\phi_2(h_{ii}))|\\
&\le & (1/m_1m)\sum_{i=1}^{m_1m}\dt=\dt.
\eneq
It follows that
\beq\label{s418109-3}
|{\bar t}({\tilde \phi}_1(h))-{\bar t}({\tilde \phi}_2(h))|<\dt\rforal h\in {\cal H}_2'.
\eneq
Since (as shown above) the theorem holds for $M_{m_1}(C(X,M_r))\cong B_1,$  by  \eqref{s418109-0},
\eqref{s418109-2}, and
\eqref{s418109-3}, we conclude
that there is a unitary $w_1\in F$ such
that
\beq\label{s418109-4}
\|{\rm Ad}\, w_1\circ {\tilde \phi}_1(a)-{\tilde \phi}_2(a)\|<\ep_1\rforal a\in {\cal F}_1.
\eneq
Note that $uau^*\in {\cal F}_1$ if $a\in {\cal F}.$ Set $w_2=w_1^*\phi_1^\sim(u^*)w_1\phi_2^\sim(u).$
Then,  by \eqref{s418109-00} and \eqref{s418109-4},
\beq\nonumber
{\rm Ad}\, w_1w_2\circ \phi_1(a)&=&w_2^*w_1^*{\phi}_1^\sim(u)^*{\tilde \phi}_1(uau^*)\phi_1^\sim(u)w_1w_2^*\\
&=&w_2^*w_1^* \phi_1^\sim (u)^*w_1w_1^*{\tilde \phi}_1(uau^*)w_1w_1^*\phi_1^\sim (u)w_1w_2\\
&\approx_{2\ep_1}&\phi_2^\sim (u^*)({\rm Ad}\, w_1\circ {\tilde \phi}_1(uau^*) )\phi^\sim_2(u)\\
&\approx_{\ep_1}& {\tilde \phi}_2(u^*){\tilde \phi}_2(uau^*){\tilde \phi}_2(u)
=\phi_2(a).
\eneq
Denote by $e_0$ the identity of $M_n.$ We also view $e_0$ as an element of $M_{m_1mn}.$
Note that $\phi_1(1_A)=\phi_2(1_A)=e_0.$
The above estimate implies that $\|e_0w_2^*w_1^*e_0w_1w_2e_0-e_0\|<3\ep_1.$ Thus there exists a unitary $w_3\in e_0M_{m_1mn}e_0=M_n$ such
that, as an element {{of}} $M_{m_1mn},$ $\|w_3-e_0w_1w_2e_0\|<6\ep_1.$
Thus we have, for all $a\in A,$
\beq
\|{\rm Ad}\, w_2\circ \phi_1(a)-\phi_2(a)\|<(6+3+6)\ep_1<\ep\rforal a\in {\cal F}.
\eneq
}}
\end{proof}

{\blue{Recall that $K_0(M_n)=\Z$ and $K_0(M_n)=\{0\}.$ So, if $X$ is connected, then, for any
two unital \hm s $\phi_1, \phi_2: C(X)\to M_n,$ $[\phi_1]=[\phi_2]$ in $KL(C(X), M_n).$

We state the following version of 2.10 of \cite{Lin-AU11} for convenience.

\begin{thm}[see 2.10 of \cite{Lin-AU11}, Theorem 4.6 of \cite{Lntams07}  and 2.15 of \cite{Li-interval}]\label{uniCMn}
 Let $X$ be a connected compact metric space, and {let $C= C(X)$}. Let $\mathcal F\subset C$ be a finite set, and let $\epsilon>0$ be a constant. There is a finite set $\mathcal H_1\subset C_+{\blue{\setminus\{0\}}}$ such that, for any $\sigma_1>0,$ there are a finite subset $\mathcal H_2\subset C$ and $\sigma_2>0$ such that for any  {unital }homomorphisms $\phi, \psi: C\to M_n$ for a matrix algebra $M_n$, satisfying
\begin{enumerate}
\item $\phi(h)>\sigma_1$ and $\psi(h)>\sigma_1$ for any $h\in\mathcal H_1$, and
\item $|{\rm tr}\circ \phi(h)-{\rm tr}\circ \psi(h)|<\sigma_2$ for any $h\in \mathcal H_2$,
\end{enumerate}
there is a unitary $u\in M_n$ such that
$$\|\phi(f)-u^*\psi(f)u\|<\epsilon\quad \textrm{for any $f\in\mathcal F$}.$$
\end{thm}

\begin{rem}\label{remark831KL}
{\rm
Let $X$ be a compact metric space and let $A=C(X).$
{\blue{Then $C(X)=\lim_{n\to\infty} C(Y_n),$ where each $Y_n$ is a finite CW complex.
It follows that $C(X, M_r)=\lim_{n\to\infty}C(Y_n, M_r)$ for any integer $r\ge 1.$
Let $P\in C(X, M_r)$ be a projection. Consider $A=PC(X, M_r)P.$ Note, since the rank function
is continuous and has integer values,
that $\{x\in X: P(x)>0\}$ is a clopen subset of $X.$ Thus, \wilog, we may assume
that $P(x)>0$ for all $x\in X.$  In particular, $P$ is a full projection {{in $C(X, M_r).$}}
We also have $A=\lim_{n\to\infty} P_nC(Y_n, M_r)P_n,$
where $P_n\in C(Y_n, M_r)$ is a projection.  Let $h_n: P_nC(Y_n, M_r)P\to A$
{{denote}} the \hm\, induced by the inductive limit. \Wlog, we may assume
that $P_n(y)>0$ for all $y\in Y_n,$ $n=1,2,....$}}

Suppose that $\phi_1, \phi_2: A\to M_n$ are  two  \hm s.
Then $(\phi_i)_{*1}=0,$ $i=1,2.$  Let ${\cal P}\subset \underline{K}(A)$ be a finite subset and {{denote by}}
 $G$ be the subgroup generated by ${\cal P}.$
There exists an integer $n\ge 1$ such that
$G\subset [h_n](\underline{K}(C(Y_n))).$ Define ${\bar G}=[h_n](\underline{K}(C(Y_n))).$
{\blue{Suppose that}} $\{p_1,p_2,...,p_k\}$ in $C(Y_n)$
{\blue{are}} mutually orthogonal projections which correspond to $k$ different
path connected components $Y_1, Y_2,..., Y_k$ of $Y$ {\blue{with $\bigsqcup_{i=1}^k Y_i=Y.$}}  Fix  $\xi_i\in Y_i,$
Let $C_i=C_0(Y_i\setminus \{\xi_i\}),$ $i=1,2,...,k.$
 Since $Y_i$ is path connected, by
considering the {{point evaluation}} at $\xi_i,$
it is easy to see that, for any {\blue{\hm}}\, $\phi: PC(Y_n, M_r)P\to M_n,$ $[\phi]|_{\underline{K}(C_i)}=0.$
Let $P_{n,i}=P_n|_{Y_i},$ $i=1,2,...,k.$
We may assume that $h_n(P_{n, i})(x)$ has two rank values $r_i\ge 1$ or zero.
Suppose that
$\tau\circ \phi_1(P_{n,i})=\tau\circ \phi_2(P_{n,i}),$ $i=1,2,...,k.$
{\blue{Then $[\phi_1]([p_i])=[\phi_2]([p_i]),$ $i=1,2,...,k.$
It follows that  $[\phi_1\circ h]=[\phi_2\circ h]$ in $KL(C(Y), M_n).$}}
One then computes that
\beq\label{KL-1408}
[\phi_1]|_{\cal P}=[\phi_2]|_{\cal P}.
\eneq

We will use this fact in the next proof.
}
\end{rem}

\begin{lem}\label{Aug-N-1}
Let $X$ be a compact metric space, let $F$ be a finite dimensional \CA\, and let $A=PC(X,F)P,$
where $P\in C(X,F)$ is a projection.
Let $\Delta: A_+^{q,{\bf 1}}\setminus \{0\}\to (0,1)$ be an order preserving map.

For any $\ep>0,$ any finite subset ${\cal F}\subset A$  and any $\sigma>0,$  there exists a finite subset
${\cal H}_1\subset A_+^{\bf 1}\setminus \{0\},$
a finite subset ${\cal H}_2\subset A_{s.a.}$, and  $\dt>0$  satisfying the following {{condition}}:
If $\phi_1, \phi_2: A\to M_n$ (for some integer $n\ge 1$) are two unital \hm s such that
\begin{eqnarray*}
\tau\circ \phi_1 (h)\ge  \Delta(\hat{h})\tforal h\in {\cal H}_1, {{\tand}}\\
|\tau\circ \phi_1(g)-\tau\circ \phi_2(g)|<\dt\tforal g\in {\cal H}_2,
\end{eqnarray*}
then, there exist  a projection $p\in M_n,$  a unital \hm\, $H: A\to pM_np, $
unital \hm s $h_1, h_2: A\to (1-p)M_n(1-p)$, and a unitary $u\in M_n$ such that
\begin{eqnarray*}
&&\|{\rm Ad}\, u\circ \phi_1(f)-(h_1(f)+H(f))\|<\ep,\\
&& \|\phi_2(f)-(h_2(f)+H(f)\|<\ep\tforal f\in {\cal F},\\
&&\andeqn
\tau(1-p)<\sigma,
\end{eqnarray*}
where $\tau$ is the tracial state of $M_n.$

\end{lem}

\begin{proof}
{\blue{Note that we may  rewrite $A=PC(X, M_r)P$ for some possibly large $r.$}}


Let $\Delta_1=(1/2)\Delta.$ Let ${\cal P}\subset \underline{K}(A)$ be a finite set,
${\cal H}_1'\subset A_+^{\bf 1}\setminus \{0\}$
(in place of ${\cal H}_1$)
be a finite set, ${\cal H}_2'\subset  A_{s.a.}$
(in place of ${\cal H}_2$)  be a finite set
and $\dt_1>0$ (in place of $\dt$) {\blue{be}} required by \ref{Oldthm201408} for
$\ep/2$ (in place of $\ep$), ${\cal F}$ and  $\Delta_1.$

%

Without loss of generality, we may assume that $1_A\in {\cal F},$ $1_A\in {\cal H}_1'\subset {\cal H}_2'$ and
${\cal H}_2'\subset A_+^{\bf 1}\setminus \{0\}.$
So, in what follows, ${\cal H}_2'\subset A_+^{\bf 1}\setminus \{0\}.$
Put
\beq\label{8-31-1}
\sigma_0=\min\{\Delta_1(\hat{g}): g\in {\cal H}_2'\}.
\eneq

Let $G$ be the subgroup generated by ${\cal P}$ and let ${\bar G}$ be {{as}} defined in  Remark \ref{remark831KL}.
{\blue{We keep the notation of \ref{remark831KL}.
Set $Q_i=h_n(P_{n,i}),$ $i=1,2,...,k.$  Denote by $h_{n*}: X\to Y_n$ the continuous map
induced by $h_n.$
Let ${\cal P}_0=\{[Q_1], [Q_2],..., [Q_k]\}$ and let $r_i\ge 1$ {{denote}} the rank of $Q_i(x)$ when $Q_k(x)\not=0$ (see
\ref{remark831KL}).}}






Let ${\cal H}_1={\cal H}_1'\cup\{Q_i: 1\le i\le k\}$
and ${\cal H}_2={\cal H}_2'\cup {\cal H}_1.$
{{Set}} $\sigma_1=\min\{\Delta(\hat{g}): g\in {\cal H}_2\}.$ Let  $r=\max\{r_1,r_2,...,r_k\}.$
Choose $\dt=\min\{\sigma_0\cdot \sigma/4kr, \sigma_0\cdot \dt_1/4kr, \sigma_1/16kr\}.$

Suppose now that $\phi_1, \phi_2: A\to M_n$ are two unital \hm s described in the lemma
for the above ${\cal H}_1,$ ${\cal H}_2$ and $\Delta.$

{\blue{If $x\in X,$ denote by $\pi_x: A\to M_{r(x)}$ the point evaluation at $x.$
We may write $\phi_j(f)=\sum_{s=1}^k (\sum_{i=1}^{n(s,j)} \psi_{s,i,j}(\pi_{x_{s,i,j}}(f)))$
for all $f\in P(C(X, M_r)P$ ($j=1,2$),
where $h_{n*}(x_{s,i,j})\in Y_{n,i},$  $\psi_{s,i,j}: M_{r_i}\to M_n$ is a \hm\, such
that $\psi_{s,i,j}(1_{M_{r_i}})$ has rank $r_i.$  Note that $x_{s,i,j}$ may be repeated, and
$\phi_j(Q_s)=\sum_{i=1}^{n(s,j)} \psi_{s,i,j}(\pi_{x_{s,i,j}}(Q_s)),$ $s=1,2,...,k,$ $j=1,2.$
Note also that $\psi_{s,i,j}(\pi_{x_{s,i,j}}(Q_s))$ has rank $r_s,$  and
the rank of $\phi_j(Q_s)$ {{is}} $n(s,j)r_s,$ $1\le s\le k,$ $j=1,2,$
}}

We have
\beq\label{8-31-2}
{\blue{({r_i\over{n}})|n(i,1)-n(i,2)|}}=|\tau\circ \phi_1(Q_i)-\tau\circ \phi_2(Q_i)|<\dt,\,\,\,i=1,2,...,k,
\eneq
where $\tau$ is the tracial state on $M_n.$
Therefore, {\blue{by comparing the {{ranks}} of $\phi_j(Q_i),$ $1\le i\le k$ ($j=1,2$), one finds}} a projection $P_{0,j}\in M_n$ such that
\beq\label{8-31-3}
\tau(P_{0,j})<k\dt<\sigma_0\cdot \sigma,\,\,\, j=1,2,
\eneq
{{and}}
${\rm rank}(P_{0,1})={\rm rank}(P_{0,2}),$
{{and}}
unital \hm s $\phi_{1,0}: A\to  P_{0,1}M_nP_{0,1}, $ $\phi_{2,0}: A\to P_{0,2}M_nP_{0,2},$
$\phi_{1,1}: A\to (1-P_{0,1})M_n(1-P_{0,1})$ and $\phi_{1,2}: A\to (1-P_{0,2})M_n(1-P_{0,2})$
such that
\beq\label{8-31-4}
\phi_1=\phi_{1,0}\oplus \phi_{1,1},\,\,\, \phi_2=\phi_{2,0}\oplus \phi_{2,1},\\\label{n151110-1}
\tau\circ \phi_{1,1}(Q_i)=\tau\circ \phi_{1,2}(Q_i),\,\,\,i=1,2,...,k.
\eneq
{{Replacing}} $\phi_1$ by ${\rm Ad}\, v\circ \phi_1,$ simplifying the notation, without loss of generality,
we may assume that $P_{0,1}=P_{0,2}.$
It follows   from (\ref{n151110-1}) that (see \ref{remark831KL})
\beq\label{8-31-5}
[\phi_{1,1}]|_{\cal P}=[\phi_{2,1}]|_{\cal P}.
\eneq
By (\ref{8-31-3}) and the choice of $\sigma_0,$ we also have
\beq\label{8-31-6}
\tau \circ \phi_{1,1}(g) \ge \Delta_1(\hat{g})\rforal g\in {\cal H}_1'\andeqn\\
|\tau \circ \phi_{1,1}(g)-\tau\circ \phi_{1,2}(g)|<\sigma_0\cdot \dt_1\rforal g\in {\cal H}_2'.
\eneq
Therefore,
\beq\label{8-31-6+}
t \circ \phi_{1,1}(g) \ge \Delta_1(\hat{g})\rforal g\in {\cal H}_1'\andeqn\\
|t \circ \phi_{1,1}(g)-t \circ \phi_{1,2}(g)|<\dt_1\rforal g\in {\cal H}_2',
\eneq
where $t$ is the tracial state on $(1-P_{1,0})M_n(1-P_{1,0}).$
By applying \ref{Oldthm201408},
there exists a unitary  $v_1\in (1-P_{1,0})M_n(1-P_{1,0})$ such that
\beq\label{8-31-7}
\|{\rm Ad}\, v_1\circ \phi_{1,1}(f)-\phi_{2,1}(f)\|<\ep/16\rforal f\in {\cal F}.
\eneq
Put $H=\phi_{2,1}$ and $p=P_{1,0}.$
The lemma follows.
\end{proof}

\begin{cor}\label{Aug-N-2}
Let $X$ be a compact metric space, let $F$ be a finite dimensional \CA\,, and let $A=PC(X,F)P,$
where $P\in C(X,F)$ is a projection.
Let $\Delta: A_+^{q,{\bf 1}}\setminus\{0\}\to (0,1)$ be an order preserving map.  Let $1>\af>1/2.$

For any $\ep>0,$ any finite subset ${\cal F}\subset A,$  any finite subset
${\cal H}_0\subset A_+^{\bf 1}\setminus \{0\},$ and any integer $K\ge 1,$   there {{are}} an integer $N\ge 1,$  a finite subset
${\cal H}_1\subset A_+^{\bf 1}\setminus \{0\},$
a finite subset ${\cal H}_2\subset A_{s.a.},$ and $\dt>0$
satisfying the following {{condition}}:
If $\phi_1, \phi_2: A\to M_n$ (for any integer $n\ge N$) are two unital \hm s such that
\begin{eqnarray*}
\tau\circ \phi_1 (h)\ge \Delta(\hat{h})\tforal h\in {\cal H}_1, {\tand}\\
|\tau\circ \phi_1(g)-\tau\circ \phi_2(g)|<\dt\tforal g\in {\cal H}_2,
\end{eqnarray*}
then there exist  mutually orthogonal non-zero projections $e_0, e_1,e_2,...,e_K\in M_n$ such that
$e_1, e_2,...,e_K$ are equivalent, $e_0\lesssim e_1$,   and $e_0+\sum_{i=1}^Ke_i=1_{M_n},$ and there are unital \hm s
$h_1, h_2: A\to e_0M_ne_0,$  $\psi: A\to e_1M_ne_1$
 and a unitary $u\in M_n$ such that
\begin{eqnarray*}
&&\|{\rm Ad}\, u\circ \phi_1(f)-(h_1(f)+{\rm diag}(\overbrace{\psi(f), \psi(f),...,\psi(f)}^K))\|<\ep,\\
 &&\|\phi_2(f)-(h_2(f)+{\rm diag}(\overbrace{\psi(f),\psi(f),...,\psi(f)}^K))\|<\ep\tforal f\in {\cal F},\\
&&\andeqn \tau\circ \psi(g)\ge \af{\Delta(\hat{g})\over{K}}\tforal g\in {\cal H}_0,
\end{eqnarray*}
where $\tau$ is the tracial state of $M_n.$

\end{cor}

\begin{proof}

By applying \ref{Aug-N-1}, it is easy to see that it suffices to prove the following statement:

Let $X,F,P$  $A$  and  $\af$ be as in the corollary.

Let $\ep>0,$  let ${\cal F}\subset A$ be a finite subset, let ${\cal H}_0\subset A_+^{\bf 1}\setminus \{0\}$ and let $K\ge 1.$  There are an integer $N\ge 1$  and a finite subset
${\cal H}_1\subset A_+^{\bf 1}\setminus \{0\}$ {{with}} the following {\blue{property}}:

Suppose that $H: A\to M_n$ (for some $n\ge N$) is a unital \hm\, such that
\beq\label{8-N-1-3}
\tau\circ H(g)\ge \Delta(\hat{g})\rforal g\in {\cal H}_0.
\eneq
Then  there  are mutually orthogonal projections $e_0, e_1,e_2,...,e_{2K}\in M_n,$ a unital \hm\, $\phi: A\to e_0M_ne_0,$
and a unital \hm\, $\psi: A\to e_1M_ne_1$ such that
\beq\label{8-N-1-4}
\|{\rm Ad}\, U\circ H(f)-(\phi(f)\oplus {\rm diag}(\overbrace{\psi(f), \psi(f),...,\psi(f)}^{2K}))\|<\ep\rforal f\in {\cal F},\\
\tau\circ \psi(g)\ge \af\Delta(\hat{g})/2K\rforal g\in {\cal H}_0
\eneq
{\blue{for some unitary $U\in M_n.$
(Note $2K$ is used since we will have $h_i(1_A)\oplus e_0\lesssim e_1+e_2$ ($i=1,2.$)}}

{\blue{We now prove the above statement.
We may rewrite $A=P(X, M_r)P$ for some  large integer $r.$
Note that there are finitely many values of the rank of $P(x).$
Thus, we may write $A=A_1\oplus A_2\oplus \cdots \oplus A_m,$
where each $A_i=P_iC(X_i, M_r)P_i,$ $X_i$ is a clopen subset of $X$ and $P_i(x)$ has {{constant}} rank for $x\in X_i,$ $i=1,2,...,m.$ Considering each summand separately, \wilog, we may
assume that $A=PC(X, M_r)P$ and $P$ has {{constant}} rank, say $r_0.$}}

Put
\beq\label{8-N-1-5}
\sigma_0= ((1-\af)/4)\min\{\Delta(\hat{g}): g\in {\cal H}_0\}>0.
\eneq
Let $\ep_1=\min\{\ep/16, \sigma_0, {\blue{1/2}}\}$ and let ${\cal F}_1={\cal F}\cup {\cal H}_0.$
Choose $d_0>0$ such that
\beq\label{8-N-1-6}
|f(x)-f(x')|<\ep_1\rforal f\in {\cal F}_1,
\eneq
provided that $x,x'\in X$ and ${\rm dist}(x,x')<d_0.$

Choose $\xi_1,\xi_2,...,\xi_m\in X$ such that $\bigcup_{j=1}^m B(\xi_j,d_0/2)\supset X,$
where $B(\xi, r)=\{x\in X: {\rm dist}(x, \xi)<r\}.$
There is $d_1>0$ such that $d_1<d_0/2$ and
\beq\label{8-N-1-7}
B(\xi_j, d_1)\cap B(\xi_i, d_1)=\emptyset,
\eneq
if $i\not=j.$
There is, for each $j,$  a function $h_j\in C(X)$ with $0\le h_j\le 1,$
$h_j(x)=1$ if $x\in B(\xi_j,d_1/2)$ and $h_j(x)=0$ if $x\not\in B(\xi_j,d_1).$
Define ${\cal H}_1={\cal H}_0\cup \{h_j: 1\le j\le m\}$ and put
\beq\label{8-N-1-8}
\sigma_1=\min\{\Delta(\hat{g}): g\in {\cal H}_1\}.
\eneq
Choose an integer $N_0\ge 1$ such that $1/N_0<\sigma_1 \cdot (1-\af)/4$ and $N=4{\blue{r_0}}m(N_0+1)^2(2K+1)^2.$

Now let $H: PC(X,M_r)P\to M_n$ be a unital \hm\, with $n\ge  N$ satisfying the assumption (\ref{8-N-1-3}).
Let $Y_1=\overline{B(\xi_1, d_0/2)}\setminus \bigcup_{i=2}^m B(\xi_i, d_1),$
$Y_2=\overline{B(\xi_2,d_0/2)}\setminus (Y_1\cup\bigcup_{i=3}^m B(\xi_i, d_1),$
$Y_j=\overline{B(\xi_j, d_0/2)}\setminus (\bigcup_{i=1}^{j-1} Y_i\cup \bigcup_{i=j+1}^m B(\xi_i, d_1)),$
$j=1,2,...,m.$ Note that $Y_j\cap Y_i=\emptyset$ if $i\not=j$ and
$B(\xi_j,d_1)\subset Y_j.$
We write that
\beq\label{8-N-1-9}
H(f)={\blue{\sum_{i=1}^{k_0} \psi_i(\pi_{x_i}(f))=\sum_{j=1}^m \sum_{x_i\in Y_j}\psi_i(\pi_{x_i}(f))}}
\rforal f\in PC(X, M_r)P,
\eneq
where {\blue{$\pi_{x_i}: A\to M_{r_0}$ is a point evaluation and
$\psi_i: M_{r_0}\to M_n$ is a \hm\, such that $\psi_i: M_{r_0}\to p_iM_np_i,$
is a unital \hm\, with multiplicity $1,$}} and
where $\{p_1,p_2,...,p_{k_0}\}$ is a set of mutually orthogonal rank  {\blue{$r_0$}}  projections in $M_n,$ {\blue{$n=k_0r_0,$ and }}
$x_1,x_2,...,x_{k_0}$ are in $X$ {\blue{(some of {{the}} $x_i$ could be repeated)}}.
Let $R_j$ be  the cardinality of $\{{\blue{x_i}}: x_i\in Y_j\},$ {\blue{counting multiplicities.}}
 Then, by (\ref{8-N-1-3}), for $j=1,2,...,m,$
\beq\label{8-N-1-10}
R_j\ge N\tau\circ H(h_j)\ge N\Delta(\hat{h_j})\ge {\blue{mr_0(N_0+1)^2 (2K+1)^2}}\sigma_1\ge mr_0(N_0+1)(2K+1)^2.
\eneq
Write $R_j=S_j2K+r_j,$ where $S_j\ge N_02Kmr_0$ and $0\le r_j<2K,$ $j=1,2,...,m.$
Choose $x_{j,1},x_{j,2},...,x_{j,r_j}$ in  $\{x_i: x_i\in Y_j\}$ and {{set}}
$Z_j=\{x_{j,1},x_{j,2},...,x_{j,r_j}\},$ $j=1,2,...,m,$ {\blue{counting multiplicities.}}

Then, we may write
\beq\label{8-N-1-11}
H(f)=\sum_{j=1}^m (\sum_{x_i\in Y_j\setminus Z_j} {\blue{\psi_i(\pi_{x_{i,j}}(f)))}}+\sum_{j=1}^m(\sum_{i=1}^{r_j}{\blue{\psi_i(\pi_{x_{i,j}}(f)))}}
\eneq
for $f\in C(X).$ Note that the cardinality of $\{x_i: x_i\in Y_j\setminus Z_j\}$
{\blue{(counting multiplicities)}} is $2KS_j,$ $j=1,2,...,m.$
{\blue{We write, counting multiplicities,  $\{x_i: x_i\in Y_j\setminus Z_j\}=\bigsqcup_{k=1}^{2K} \Omega_{k,j},$
where each $\Omega_{k,j}$ has exactly $S_j$ points in $\{x_i: x_i\in Y_j\setminus Z_j\},$ counting multiplicities.
Put $e_k=\sum_{x_i\in \Omega_{k,j}}p_i,$ $k=1,2,...,2K.$
Then each $e_k$ has rank $(\sum_{j=1}^m S_j)r_0,$ $k=1,2,...,2K.$}}
Define
{\blue{
\beq\label{8-N-1-12}
\phi(f)&=&\sum_{j=1}^m(\sum_{i=1}^{r_j}\psi_i(\pi_{x_{i,j}}(f)))\andeqn\\
\Psi(f)&=&\sum_{j=1}^m (\sum_{x_i\in Y_j\setminus Z_j} \psi_i(\pi_{\xi_j}(f)))=
\sum_{j=1}^m (\sum_{k=1}^{2K} (\sum_{x_i\in \Omega_{k,j}}\psi_i(\pi_{\xi_j}(f)))\\
&=&\sum_{k=1}^{2K}(\sum_{j=1}^m(\sum_{x_i\in \Omega_{k,j}}\psi_i(\pi_{\xi_j}(f)))
%
=\sum_{k=1}^{2K}(\sum_{j=1}^{m}\Psi_{k,j}(\pi_{\xi_j}(f))),
\eneq}}
where $\Psi_{k,j}$ is a direct sum of $S_j$ \hm s $\psi_i$ with $x_i\in \Omega_{k,j}.$
We estimate that
\beq\label{181020n1}
\|H(f)-(\phi(f)\oplus \Psi(f))\|<\ep_1\rforal f\in {\cal F}_1.
\eneq
{\blue{Note that each  $\psi_i$ is unitarily equivalent to $\psi_1.$  So each $\Psi_{k,j}$ is unitarily
equivalent to $S_j$ copies of $\psi_1.$
It follows that, for each $k,$
$\sum_{j=1}^m \Psi_{k,j}\circ \psi_{\xi_j}$ is unitarily equivalent to $\psi:=\sum_{j=1}^m \Psi_{1,j}\circ \pi_{\xi_j}.$
Thus there is a unitary $U_1\in M_n$
such that
\beq\label{8-N-1-14}
{\rm{Ad}}\, U_1\circ \Psi(f)={\rm diag}(\overbrace{\psi(f), \psi(f),...,\psi(f)}^{2K})\rforal f\in A.
\eneq
Put $e_0=\sum_{j=1}^m (\sum_{x_i\in Z_j})p_i.$ Then $e_0$ has rank $\sum_{j=1}^m r_jr_0.$
Moreover, $\phi$ is a unital \hm\, from $A$ into $e_0M_ne_0.$}}
Note that
\beq\label{8-N-1-13}
{\rm rank}(e_0)=\sum_{j=1}^mr_jr_0<mr_02K\andeqn
S_j\ge N_02Kmr_0>mr_02K,\,\,\,j=1,2,...,2K.
\eneq
It follows that $e_0\lesssim e_1$ and $e_i$ is equivalent to $e_1.$
Thus, by \eqref{181020n1}, for some unitary $U\in M_n,$
\beq\label{8-N-1-15}
\|{\rm Ad}\, U\circ H(f)-(\phi(f)\oplus {\rm diag}(\overbrace{\psi(f), \psi(f),...,\psi(f)}^{2K}))\|<\ep_1\rforal f\in {\cal F}_1.
\eneq
We also compute that, for all $g\in {\cal H}_0,$
\beq\label{8-N-1-16}
\tau\circ \psi(g)\ge (1/2K)(\Delta(\hat{g})-\ep_1-{\blue{{mr_02K\over{N_02Kmr_0}} )\ge \af{\Delta(\hat{g})\over{2K}}.}}
\eneq
\end{proof}

\begin{rem}\label{8-N-2-r}
{\rm
If we also assume that $X$ has infinitely many points, then Lemma \ref{Aug-N-2} holds without
mentioning the integer $N.$ This can be seen by taking larger ${\cal H}_1$ which {\blue{contains
at least $N$ mutually orthogonal non-zero elements. This}}  will force the integer
$n$  to be larger than $N.$
}
\end{rem}

\begin{df}\label{8-N-3}
Denote by ${\blue{{\cal {\bar D}}_0}}$ the class of all \CA s  {\blue{with the form
$PC(X, F)P,$ where $X$ is a compact metric space, $F$ is a finite dimensional \CA, and $P\in C(X,F)$
is a projection.}}
For $k\ge 1,$ denote by ${\blue{{\cal {\bar D}}_k}}$ the class of all \CA s with the form:
$$
A=\{(f,a)\in {\blue{PC(X, F)P\oplus B}}: f{|_{Z}}=\Gamma(a)\},
$$
where $X$ is a compact  metric space,  $F$ is a finite dimensional \CA, $P\in C(X,F)$ is a projection, { $Z\subset X$ is a {{non-empty}} proper {closed} subset of $X,$
$B\in {\blue{{\cal {\bar D}}_{m}}}$ {\blue{for some $0\le m<k,$}}}
and $\Gamma: B\to P|_ZC({ Z},F)P|_Z$ is a unital \hm,
where we assume that there is $d_{X,{ Z}}>0$
such that, for any $0<d\le d_{X,Z},$
there exists $s_*^d: \overline{X^d}\to { Z}$ such that
\beq\label{Amnotation-1}
s_*^d(x)=x\rforal x\in Z \andeqn \lim_{d\to 0} \|f|_{ Z}\circ s_*^d-f|_{\overline{X^d}}\|=0\rforal f\in C(X,F),
\eneq
where $X^d=\{x\in X: {\rm dist}(x, { Z})< d\}.$
We also assume that,  for any $0<d<{\blue{d_{X,Z}/2}}$ and for any $d>\dt>0,$
there is { a homeomorphism}  $r: X\setminus X^{d-\dt}\to X\setminus X^{d}$ such that
\beq\label{Amnotation-2+}
{\rm dist}(r(x), x)<\dt\rforal x\in X\setminus X^{d-\dt}.
\eneq

{\blue{In what follows, {{as in \ref{DfC1}},} we will use $\lambda: A\to PC(X, F)P$ for the \hm\, defined
by $\lambda((f,b))=f$ for all $(f,b)\in A.$}}

{\blue{ Note that, $\bar{D}_{k-1}\subset \bar{D}_k.$
Suppose that $A, B\in \bar{D_k}.$ Then, from the definition above, it is routine to check that $A\oplus B\in \bar{D}_k.$
 }}




Let $A_m$ be a unital \CA\, in ${\cal {\bar D}}_m.$
For $0\le k<m,$ let  $A_k\in {\cal {\bar D}}_k$ such that
$A_{k+1}=\{(f, a)\in P_{k+1}C(X_{k+1}, F_{k+1})P_{k+1}\oplus A_k: f|_{ Z_{k+1}}=\Gamma_{k+1}(a)\},$
where $F_{k+1}$ is a finite dimensional \CA, $P_{k+1}\in C(X_{k+1}, F_{k+1})$ is a projection,
$\Gamma_{k+1}: A_k\to P_{k+1}C({ Z_{k+1}}, F_{k+1})P_{k+1}$ is a unital \hm,
$k=0,1,2,...,m-1.$ Denote by $\partial_{k+1}: P_{k+1}C(X_{k+1}, F_{k+1})P_{k+1}\to
Q_{k+1} C(Z_{k+1}, F_{k+1})Q_{k+1}$  the map defined by
$f\mapsto f|_{Z_{k+1}}$ ($Q_{k+1}=P_{k+1}|_{Z_{k+1}}$).
We use $\pi_e^{(k+1)}: A_{k+1}\to A_k$ for the quotient map and
$\lambda_{k+1}: A_{k+1}\to P_{k+1}C(X_{k+1}, F_{k+1})P_{k+1}$ for the map  
$(f,a)\mapsto f.$

For each $k,$ one has the following commutative diagram:
\begin{equation}\label{pull-back-k}
\xymatrix{
A_{k} \ar@{-->}[rr]^{\lambda_k} \ar@{-->}[d]^-{\pi_e^{({\blue{k}})}}  && P_kC(X_k, F_k)P_k \ar[d]^-{\partial_k} \\
A_{k-1}\ar[rr]^-{\Gamma_k} & & Q_kC({ Z_k}, F_k)Q_k.
}
\end{equation}

In general, suppose that $A=A_m\in{\cal {\bar D}}_m$ is constructed as in the following sequence
$$A_0{\blue{\in \overline{{\cal D}}_0,}}
~~ A_1=P_1C(X_1, F_1)P_1\oplus_{Q_1C(Z_1, F_1)Q_1}A_0,~~A_2=P_2C(X_2, F_2)P_2\oplus_{Q_2C(Z_2, F_2)Q_2} A_1, \cdots,~~~~~~~~~~~~$$
$$ ~~~~~~~~~~~~~~~~~~~~~~ A_m=P_mC(X_m, F_m)P_m \oplus_{Q_mC(Z_m, F_m)Q_m} A_{m-1},$$
where $Q_i=P_i|_{Z_i},$ $i=1,2,...,m.$
With $\pi_e^{(k+1)}$ and $\lambda_{k}$ above we can define the quotient map $\Pi_k: A=A_m \to A_k$ and {{the}} homomorphism $\LD_k: A =A_m \to P_kC(X_k, F_k)P_k$ as follows:
$$\Pi_k=\pi_e^{(k+1)}\circ \pi_e^{(k+2)}\circ \cdots \circ \pi_e^{(m-1)}\circ \pi_e^{(m)}~~~~~~~~~~{\mbox {and}}~~~~~~~~~ \LD_k=\lambda_k\circ \Pi_k. $$
Combining all $\LD_k$ we get the inclusion homomorphism
$$\LD: A \to \bigoplus_{k=0}^mP_kC(X_k, F_k)P_k$$
with $X_0$ being {{the}} single point set. {\blue{In particular, $A$ is a subhomogeneous \CA.}}
For each $k\ge 1,$ we may write
$$
P_k(C(X_k, F_k))P_k=P_{k,1}C(X_k, M_{s(k,1)})P_{k,1}\oplus P_{k,2}C(X_k, M_{s(k,2)})P_{k,2}\oplus\cdots
{{\oplus}} P_{k,t_k}C(X_k, M_{s(k,t_k)})P_{t_k},
$$
where $P_{k,j}\in C(X_k, M_{s(k,j)})$ is a projection of rank $r(k,j)$ at each $x\in X_k.$
For each $x\in X_k$ and $j\leq t_k$ and $f\in A$, denote by $\pi_{(x,j)}(f)\in M_{r(k,j)}$ the evaluation of {{the}}
$j^{th}$ component of $\LD_k(f)$ at {{the}} point $x.$
 Then for each pair $(x,j)$, $\pi_{(x,j)}$ is a finite dimensional representation of $A$, and, furthermore if we assume $x\in X_k\setminus Z_k$ (and $P_{k,j}(x)\not=0$) then $\pi_{(x,j)}$ is an irreducible representation.


{{\blue{In}} the definition
of ${\cal {\bar D}}_m$ above,  if, in addition, $\blue{X_k}$ is path connected,}
${\blue{Z_k}}$  has finitely
many path connected components, and ${\blue{X_k}}\setminus {\blue{Z_k}}$ is path connected,
then we will use {${\cal D}_m$} for the resulting class of \CA s.
Note that {{${\cal D}_m\subset {\cal {\bar D}}_m$}} and ${\cal C}\subset {\cal D}_1.$

Note that if $X_k$ is  a
simplicial complex, $Z_k\subset X_k$  is sub-complex,
$k=1, 2, ..., m$,  then any iterated pull-back
$$
(P_m C(X_m, F_m) P_m\oplus_{Q_kC(Z_k, F_k)Q_k} P_{k-1} C(X_{k-1}, F_{k-1}) P_k)\oplus \cdots \oplus_{Q_1C(Z_2, F_1)Q_1} P_0 C(X_0, F_0) P_0
$$
{\blue{is in ${\cal D}_m.$}}

\end{df}


%

\begin{rem}\label{8-n-rn}
Let $A\in {\cal {\bar D}}_k$ (or $A\in {\cal D}_k$). It is easy to check that $C(\T)\otimes A\in {\cal {\bar D}}_{ k+1}$
(or ${\blue{C(\T)\otimes A\in {\cal D}_{k+1}}}$). First,  if $F_0$ is a finite dimensional $C^*$ algebra, then $C(\T)\otimes F_0\in {\cal {\bar D}}_1$ by putting $F_1=F_0$ and $X_1=\T$ with $Z_1=\{1\}\subset \T$ and $\GM_1: F_0\to C(Z_1, F_1)\cong F_0$ to be the identity map. And if a  pair of spaces $(X_k, Z_k)$ satisfies   the conditions described for the  pair $(X, Z)$ in the definition above, in particular, the existence of the retraction $s_*^d$ and homeomorphism $r$ as in {{\eqref{Amnotation-1} and \eqref{Amnotation-2+}}}, then the  pair  $(X_k \times \T , Z_k\times \T)$ also satisfies the same conditions.

\end{rem}

\begin{lem}\label{Lmeasdiv}
Let $X$ be a  measurable space with {{infinitely}} many points  with a {{specified}} set ${\cal M}$ of $k$ probability measures.
  Suppose $Y_1, Y_2,...,Y_m$ are disjoint measurable sets with $m\ge kN,$ where
  $N\ge 1$ is an integer.
  Then, for some $j,$ $\mu(Y_j)<1/N$ for all $\mu\in {\cal M}.$
\end{lem}

\begin{proof}
Write ${\cal M}=\{\mu_1, \mu_2,...,\mu_m\}.$
Since
$\sum_{i=1}^m \mu_1(Y_i)\le 1,$ there are at least $(k-1)N$ many $Y_i$'s such that
$\mu_1(Y_i)<1/N.$
We may assume that $\mu_1(Y_i)<1/N,$ $i=1,2,...,(k-1)N.$
Then among $\{Y_1, Y_2,...,Y_{(k-1)N}\},$ there are at least $(k-2)N$ $Y_i$'s such
that $\mu_2(Y_i)<1/N.$ By induction, one finds at least one $Y_i$ with the property
that $\mu_j(Y_i)<1/N$ for all $j.$

\end{proof}

\begin{lem}\label{8-N-4}
Let $A\in {\cal {\bar D}}_k$ be a unital \CA.
 Let $\Delta: A_+^{q, {\bf 1}}\setminus \{0\}\to (0,1)$ be an order preserving map.  Let $1>\af>1/2.$

Let $\ep>0,$ ${\cal F}\subset A$ be a finite subset, ${\cal H}_0\subset A_+^{\bf 1}\setminus \{0\}$ be a finite subset,
and $K\ge 1$ be an integer. There exist  an integer $N\ge 1,$ $\dt>0,$   a finite subset ${\cal H}_1\subset A_+^{\bf 1}\setminus \{0\},$
and  a finite
subset ${\cal H}_2\subset A_{s.a.}$  satisfying the following {{condition}}:

If $\phi_1, \phi_2: A\to M_n$(for some integer $n\ge N$) are two unital \hm s such that
\beq\label{8-N-4-1}
{\rm tr}\circ \phi_1(g)\ge \Delta(\hat{g})\tforal g\in {\cal H}_1\tand\\\label{8-N-4-1nn}
|{\rm tr}\circ \phi_1(g)-{\rm tr}\circ \phi_2(g)|<\dt\tforal g\in {\cal H}_2,
\eneq
where ${\rm tr}$ is the tracial state on $M_n,$ then
there exist mutually orthogonal projections\\
$e_0, e_1, e_2,...,e_K\in M_n$ such that
 $e_1,e_2,...,e_K$ are mutually equivalent, $e_0\lesssim e_1,$  and $e_0+\sum_{i=1}^K e_i=1_{M_n},$
unital \hm s $h_1, h_2: A\to e_0M_ne_0,$ a unital \hm\, $\psi: A\to e_1M_ne_1$,  and a unitary
$u\in M_n$ such that
\beq\label{8-N-4-3}
&&\|{\rm Ad}\, u \circ \phi_1(f)-(h_1(f)\oplus {\rm diag}(\overbrace{\psi(f), \psi(f),...,\psi(f)}^K))\|<\ep,\\
&&\|\phi_2(f)-(h_2(f)\oplus {\rm diag}(\overbrace{\psi(f), \psi(f),...,\psi(f)}^K))\|<\ep  \tforal f\in {\cal F}, \\
&&\tand {\rm tr}\circ \psi(g)\ge \af{\Delta(\hat{g})\over{K}}\tforal g\in {\cal H}_0,
\eneq
where ${\rm tr}$ is the tracial state on $M_n.$
\end{lem}

\begin{proof}
We will use {{the}} induction on the integer $k\ge 0.$ The case $k=0$ follows from Corollary \ref{Aug-N-2}. Assume
that the conclusion of the lemma holds for integers $0\le k\le  m.$

  We assume that $A\in {\blue{{\cal {\bar D}}_{m+1}}}.$
  We will retain the notation for $A$ as an algebra in ${\cal {\bar D}}_{m+1}$  in the  later part of  Definition \ref{8-N-3}.
Put $X_{m+1}=X,$ ${ Z_{m+1}= Z}$, $Y=X\setminus { {Z}},$
$X^0={ Z}=X\setminus Y,$ and
$I= PC_0(Y,F)P\subset A.$
We will write
\beq\label{adf}
A=\{(f,b)\in {\blue{PC(X, F)P\oplus B}}: f|_{X^0}=\Gamma(b)\},
\eneq
where $B\in {\cal {\bar D}}_m$  is a unital \CA\,  and will be identified with
$A/I.$
We also keep the notation $\lambda: A\to PC(X, F)P$ in the pull-back {\blue{of}}
 \ref{8-N-3}.
We will write $f|_S$ for $\lambda(f)|_S$ for $f\in A$  and $S\subset X$ in the proof when there is no confusion.
Let $d_{X,Z}> 0$ be {{as}} given in \ref{8-N-3}.
Denote by $\pi_I: A\to A/I$ the quotient map.
We may write
\beq\label{8-N-40-150825-1}
PC(X,F)P=\bigoplus_{j=1}^{k_2}P_{j}C({\blue{X_j}}, M_{s(k,j)})P_{j},
\eneq
where $P_{j}\in C(X, M_{s(k,j)})$ is a projection of rank $r(j)$ at each $x\in X_j.$
{{Let us say that the dimensions of the irreducible representations of $A/I$ are}} $l_1, l_2,...,l_{k_1}.$
Set
\beq\label{181127-412-T}
&&T=(k_1k_2)\cdot \max_{i,j}\{z_iz_j: z_i, z_j\in \{l_1,l_2,...,l_{k_1}, r(1),r(2),...,r(k_2)\}\}\andeqn\\
\label{8-N-4-3+}
&&\dt_{00}=\min\{\Delta(\hat{g})/2: g\in {\cal H}_0\}>0.
\eneq
Let $\bt=\sqrt{1-(1-\af)/8}=\sqrt{(7+\af)/8}.$  Note  that $1>\bt^2>\af.$
Fix $N_{00}\ge 4$ such that
\beq\label{181128-412-bt}
1/N_{00}<{(1-\bt)\dt_{00}\over{64}}.
\eneq
Fix $\ep>0,$ a finite subset ${\cal F}\subset A,$ {{and}} a finite subset ${\cal H}_0\subset A_+^{\bf 1}\setminus \{0\},$
{\blue{and let $K>0$ be an
integer. Let $K_0=N_{00}K.$}}
We may assume that $1_A\in {\cal H}_0\subset {\cal F}.$ Without loss of generality, we may also assume that
${\cal F}\subset A_{s.a.}$ and $\|f\|\le 1$ for all $f\in {\cal F}.$
Write $I=\{f\in PC(X, F)P: f|_{X^0}=0\}.$
There is $d>0$ such that
\beq\label{8-N-4-4}
\|\pi_{x,j}(f)-\pi_{x',j}(f)\|<\min\{\ep, \dt_{00}\}/256KN_{00}\rforal f\in {\cal F},
\eneq
provided that ${\rm dist}(x, x')<d$ for any pair $x, x'\in X$
(here we  identify $\pi_{x,j}(f)$ with $\pi_{x,j}(\LD(f))$
---see {{Definition}} { \ref{8-N-3}}).
Put $\ep_0=\min\{\ep, \dt_{00}\}/16KN_{00}.$

We also assume that, {\blue{for any $x\in X^d=\{x\in X: {\rm dist}(x, X^0)<d\},$}} choosing a smaller $d$ if necessary,
\beq\label{8-N-4-4+}
\|\pi_{x,j}\circ s^d\circ { (\lambda(f)|_Z)}-\pi_{x,j}(f)\|<\ep_0/16{\blue{\rforal f\in {\cal F},}}
\eneq
where $s^d: QC(Z, F)Q\to P|_{X^d}C(\overline{X^d}, F)P|_{X^d}$ is induced by $s_*^d: \overline{X^d} \to Z$ (see \ref{8-N-3}). Note that $s^d$ also induces a map (still denoted by $s^d$)
\beq\label{112618-4-1}
s^d: B \to {\blue{P|_{\overline{X^d}}C(\overline{X^d}, F)P|_{\overline{X^d}}\oplus_{QC(Z,F)Q} B,}}
\eneq
where $Q=P|_Z, $   {\blue{by $s^d(a)=(s^d(\Gamma(a)), a)$ for all $a\in B.$}}
To simplify notation, let us  assume that $d<d_{X, Z}/2.$


For any $b>0,$ as in \ref{8-N-3}, we will continue to use $X^b$ for
$\{x\in X: {\rm dist}(x, X^0)<b\}.$

Let  ${\blue{Y_{0,d/2}}}=X\setminus X^{d/2}.$ {\blue{Note $Y_{0,d/2}$ is closed.}}
Put ${\blue{C_{I,0}}}=PC(Y_{0,d/2}, { F})P.$
Let ${\cal F}_{I,0}=\{f|_{Y_{0,d/2}}: f\in {\cal F}\}$ and {{set}} ${\cal H}_{0,I,0}=\{h|_{Y_{0,d/2}}: h\in {\cal H}_0\}.$
%
{Let $f_{0,0}\in C_0(Y)_+$ be such that $0\le  f_{0,0}\le 1,$ $f_{0,0}(x)=1$ if $x\in X\setminus X^d,$
$f_{0,0}(x)=0$ if $x\notin Y_{0,d/2},$ and $f_{0,0}(x)>0$ if ${\rm dist}(x, X^0)>d/2.$}

Let $\Delta_{I,0}: ({\blue{C_{I,0}}})_+^{q, {\bf 1}}\setminus \{0\}\to (0,1)$ be defined by
\beq\label{8-N-4-5}
\Delta_{I,0}(\hat{g})=\bt \Delta(\widehat{g'})\tforal g\in ({\blue{C_{I,0}}})_+^{\bf 1}\setminus\{0\},
\eneq
where $g'={\blue{(f_{0,0}\cdot P)}}\cdot g$
is viewed as an element in $I_+^{\bf 1}.$
Note that if $g\in C_{I,0}$ and $g\not=0,$ then $(f_{0,0}\cdot P)\cdot g\not=0.$ So $\Delta_{I,0}: ({\blue{C_{I,0}}})_+^{q, {\bf 1}}\to (0,1)$ is an order preserving
map.  Let $N^I\ge 1$ be an integer (in place of $N$) as
provided by \ref{Aug-N-2} for ${\blue{C_{I,0}}}$ (in place of $A$),
$\Delta_{I,0}$ (in place of $\Delta$), $\ep_0/16$ (in place of $\ep$), ${\cal F}_{I,0}$ (in place of ${\cal F}$),
${\cal H}_{0,I,0}$ (in place of ${\cal H}_0$), {and $2K_0$ (in place of $K$)}.

Let ${\cal H}_{1, I,0}\subset ({\blue{C_{I,0}}})_+^{\bf 1}\setminus \{0\}$ (in place of ${\cal H}_1$),
${\cal H}_{2,I,0}\subset ({\blue{C_{I,0}}})_{s.a.}$ (in place of ${\cal H}_2$), and $\dt_{1}>0$ (in place of $\dt$)
{{be}} {{the finite subsets and constant}}
provided
by \ref{Aug-N-2} for $\ep_0/16$ (in place of $\ep$), ${\cal F}_{I,0}$ (in place of ${\cal F}$), $2K_0$ { (in place of $K$),}  ${\cal H}_{0,I,0}$ associated
with ${\blue{C_{I,0}}}$ (in place of $A$),
$\Delta_{I,0}$ (in place of $\Delta$), {and $\beta$ (in place of $\alpha$)}.
Without loss of generality, we may assume that
$\|g\|\le 1$ for all $g\in {\cal H}_{2,I,0}.$ {\blue{We may assume
that $1_{{\blue{C_{I,0}}}}\in {\cal H}_{1,I,0}$ and $1_{C_{I,0}}\in {\cal H}_{2,I,0}.$}}

Let ${\cal F}_{\pi}=\pi_I({\cal F}).$
Let $g_0'\in C(X)_+$ with $0\le g_0'\le 1$ such that
$g_0'(x)=0$ if ${\rm dist}(x, X^0)<d/256$ and $g_0'(x)=1$ if ${{{\rm dist}(x,Y_{0,d/2})\le d/16}}.$
Define $g_0=1_A-g_0'\cdot P{\blue{=((1-g_0') P, 1_B)}}.$ Since ${\blue{g_0'\cdot P}}\in I,$ we view $g_0$ as an element {{of}} $A.$
Hence, for $g\in A/I=B,$
$g_0\cdot s^d(g)=((1-g_0')P\cdot s^d({\blue{\Gamma(g)}}), g)\in A$ {\blue{(see \eqref{112618-4-1}).}}
Define
\beq\label{8-N-4-9}
\Delta_\pi(\hat{g})=\bt\Delta(\widehat{g_0\cdot s^d(g)})\rforal g\in (A/I)_+^{\bf 1}.
\eneq
{{We will {{later}} use the fact that $g_0(x)=0$ if ${\rm dist}(x, Y_{0,d/2})\le d/16.$}}

 Note if $g$ is non-zero, so is $s^d(g).$  {{Since $g_0|_{X^0}=1,$}} we have   $g_0 \cdot s^d(g)\not=0.$
It follows that $\Delta_\pi: (A/I)_+^{q, {\bf 1}}\setminus \{0\}\to (0,1)$ is an order preserving map.

Put ${\cal H}_{0,\pi}=\pi_I({\cal H}_0).$
Let $N^{{{\pi}}}\ge 1$ be the integer associated  with $A/I\, (=B),$ $\Delta_\pi,$ $\ep_0/16,$ ${\cal F}_{\pi}$ and
${\cal H}_{0,\pi}$ (as required by the inductive assumption that the lemma holds for integer $m$).

 Let ${\cal H}_{1,\pi}\subset (A/I)_+^{q, {\bf 1}}{\blue{\setminus\{0\}}}$
(in place of ${\cal H}_1$),
${\cal H}_{2,\pi}\subset A/I_{s.a.}$
(in place of ${\cal H}_2$), and
$\dt_2>0$ (in place of $\dt$)
 {{denote the finite subsets and constant}}  provided
{{by the inductive assumption that the  lemma {\blue{holds}} for the case  that $k=m$}}  for $\ep_0/16$ (in place of $\ep$), ${\cal F}_{\pi}$ (in place of ${\cal F}$), ${\cal H}_{0, \pi}$ (in place of ${\cal H}_0$),
$2K_0$ associated with $A/I$ (in place of $A$), $\Delta_{\pi}$ (in place of $\Delta$), and $\bt$ (in place of $\af$).
{{\Wlog,  we may assume that $\|h\|\le 1$ for all $h\in {\cal H}_{2,\pi}.$}}

Set ${{{\dt_{000}}}}=\min\{{{ \dt_1,\dt_2, \ep_0}}\}.$ There is an integer $N_0\ge 256$ such that
\beq\label{8-N-4-5+}
&&\hspace{-0.8in}1/N_0<\Delta(\widehat{f_{0,0}\cdot P})\cdot{{\dt_{000}}}^2\cdot \min\{\Delta_{I,0}(\hat{g}): g\in {\cal H}_{1,I,0}\}\cdot
\min\{\Delta_{ \pi}(\hat{g}): g\in {\cal H}_{1,\pi}\}/64K_0N_{00}.
\eneq
{\blue{Choose $0<d_0<d$ such that, if ${\rm dist}(x, x')<d_0,$ {{then}}
\beq\label{181128-412-d0}
\|g(x)-g(x')\|<1/16N_0^2\andeqn \|f_{00}Pg(x)-f_{00}P g(x')\|<1/16N_0^2\rforal  g\in {\cal H}_{1,I,0}.
\eneq
}}
Define $Y_k$ to be the closure of $\{y\in Y: {\rm dist}(y, {\blue{Y_{0,d/2}}})<kd_0/64N_{ 0}^2\},$
$k=1,2,..., 4N_0^2.$

Let ${\cal F}_{I,k}=\{f|_{Y_k}: f\in {\cal F}\}$ and let ${\cal H}_{0,I,k}=\{h|_{Y_k}: h\in {\cal H}_0\}.$
Put ${\blue{C_{I,k}}}=P|_{Y_k}C(Y_k, F)P|_{Y_k}.$
Let $f_{0,k}\in C_0(Y)_+$  be such that $0\le f_{0,{ k}}\le 1,$ $f_{0,k}(x)=1$ if $x\in Y_{k-1},$
$f_{0,k}(x)=0$ if $x\not\in Y_k$ and $f_{0,k}(x)>0$ if ${\rm dist}(x, Y_{0,d/2})<kd_0/64N_0^2,$ $k=1,2,...,4N_0^2.$


Let $r_k: Y_k\to Y_{0,d/2}$ be a homeomorphism such that
\beq\label{8-N-4-6}
{\rm dist}(r_k(x),x)<d_0/16\tforal x\in Y_k,\,\,\,k=1,2,..., 4N_0^2
\eneq
(see \ref{8-N-3}).
%
Let ${\cal F}_{I,k}'=\{f\circ r_k: f\in {\cal F}_{I,0}\}$ and ${\cal H}_{0,I,k}'=\{g\circ r_k: g\in {\cal H}_{0,I,0}\},$
$k=1,2,...,4N_0^2.$

Any unital \hm\, $\Phi: {\blue{C_{I,k}}}\to D$ (for a unital \CA\, $D$)
 induces  a unital \hm\, $\Psi: {\blue{C_{I,0}}}\to D$ defined by $\Psi(f)=\Phi(f\circ r_k)$ for all $f\in {\blue{C_{I,0}}}.$ Note
that $f\mapsto f\circ r_k$ is an isomorphism from ${\blue{C_{I,0}}}$ onto ${\blue{C_{I,k}}}.$
{\blue{Therefore,  applying  Corollary \ref{Aug-N-2} $4N_0$ times,  }}
for $\ep_0/16$ (in place of $\ep$), ${\cal F}_{I,k}'$ (in place of ${\cal F}$), $2K_0$ (in place of $K$),  and ${\cal H}_{0,I,k}'$
(in place of ${\cal H}_0$) associated
with ${\blue{C_{I,k}}}$ (in place of $A$), $\Delta_{I,0}$ (in place of $\Delta$), and $\bt$ (in place of $\af$),
we obntain
${\cal H}_{1,I,k}$ (in place of ${\cal H}_1$), which  we may suppose equal to
${\blue{{\cal H}}}_{1,{ I},0}\circ r_k,$
${\cal H}_{2,I,k}$ (in place of ${\cal H}_2$), which we may suppose equal to
${\cal H}_{2,{ I},0}\circ r_k,$ and $\dt_1$ (in place of $\dt$).
{\blue{(Note that $\Delta_{I,0}$  is the same as above.)}}

We also note that
\beq\label{8-N-4-8}
\|f-f|_{{{Y_{0,d/2}}}}\circ r_k\|<\min\{\ep, \dt_{00}\}/64K_0N_{00}\tforal f\in {\cal F}_{I,k}.
\eneq
Define
\beq\label{181127-412-sm}
\sigma_0=\min\{\min_{0\le k\le 4N_0^2} \{\min\{\Delta_{I,0}(\widehat{g\circ r_k}): g\in {\cal H}_{1, I,k}\}\}, \min\{\Delta_{\pi}(\hat{g}): g\in {\cal H}_{1,\pi}\}\}.
\eneq
{\blue{Then $\Delta(\widehat{f_{00}\cdot P})\ge \sigma_0.$}}
Choose an integer $N\ge (N^\pi+N^I)$ such that
\beq\label{8-N-4-13-}
T/N< \sigma_0\cdot \min\{\dt_1/64, \dt_2/64, \ep_0/64K_0\}/N_{00}{\blue{(N^l+N^\pi)}}.
\eneq
Put
\beq\label{8-N-4-10}
{\cal H}_1=\bigcup_{k=0}^{4N_0^2}\{f_{0,k}\cdot P\circ g\circ { r_k}: g\in {\cal H}_{1,I ,0}\}\cup
\{(g_0\cdot P\cdot s^{ d}(g),g): g\in {\cal H}_{1,\pi}\}.
\eneq
With  the convention that $r_0: Y_{0,d/2}\to Y_{0,d/2}$ is the identity map, put
\beq\label{8-N-11}
{\cal H}_2'=\bigcup_{k=0}^{4N_0^2}\{f_{0,k}\cdot P\cdot g\circ r_k: g\in {\cal H}_{2, I,0}\}\cup\{f_{0,k}\cdot P: 0\le k\le 4N_0^2\}.
\eneq
 Put $g_{0.k}=1_A-{\blue{f_{0,k}}}\cdot P.$
  Note, since $f_{0,k} \cdot P\in I,$
$g_{0,k}{\blue{=((1-f_{0,k})P, 1_B)}}\in A.$
Define
\beq\label{8-N-12}
{\cal H}_2''=\bigcup_{k=1}^{4N_0^2}\{(g_{0,k}\cdot  s^{ d}({\blue{\Gamma(g)}}),g)\in A: g\in {\cal H}_{2, \pi}\}\cup {\cal F}.
\eneq
Put
${\cal H}_2={\cal H}_2'\cup{\cal H}_2''.$
Let
\beq\label{181127-412-dt}
\dt={\sigma_0\cdot \min\{\dt_1/64, \dt_2/64, \ep_0/64K\}\over{4K_0{\blue{N_0}}N_{00}}}.
\eneq


Now let $\phi_1,\phi_2: A\to M_n$ (for some integer $n\ge N$) be two unital \hm s  {{satisfying the assumptions}} for the above ${\cal H}_1,$ ${\cal H}_2,$ and $\dt.$

Consider the two finite  Borel  measures on $Y$ defined by
\beq\label{8-N-4-15}
\int_Y f\mu_i={\rm tr}\circ \phi_i(f\cdot P)\rforal f\in C_0(Y),\,\,\,i=1,2.
\eneq

Note that $\{Y_k\setminus Y_{k-1}: k=1,2, { ...},4N_0^2\}$ is a family of $4N_0^2$ disjoint Borel sets.
{\blue{By \ref{Lmeasdiv},}} there exists $k$ such that
\beq\label{8-N-4-17}
\mu_i(Y_k\setminus Y_{k-1})<1/N_0,
,\,\,\, i=1,2.
\eneq
{\blue{We fix this $k.$}} We may write
\beq\label{8-N-4-18}
\phi_1=\Sigma^1_\pi\oplus \Sigma_b^1\oplus \Sigma^1_s\oplus \Sigma^1_I\andeqn \phi_2=\Sigma^2_\pi\oplus
\Sigma_b^2\oplus \Sigma^2_s\oplus \Sigma^2_I,
\eneq
where $\Sigma^1_I$  and $\Sigma^2_I$ are finite direct sums of terms of the form $\pi_{x,j}$ for $x\in Y_{k-1},$
$\Sigma_s^1$ and $\Sigma_s^2$ are finite direct sums {{of terms of the form}} $\pi_{x,j}$ for $x\in Y_k\setminus Y_{k-1},$
$\Sigma_b^1$ and $\Sigma_b^2$ are    finite direct sums {{of terms}} of the form
$\pi_{x,j}$ for $x\in Y\setminus Y_k,$
and $\Sigma^1_{\pi}$ and $\Sigma^2_{\pi}$ are finite direct sums {{of terms}} of the form
 ${\bar\pi}_{x,i}$
 given by irreducible representations of $A/I$ {{(note that these $\pi_{x,j}$ or ${\bar \pi}_{x,i}$ can be repeated)}}.

Define $\psi_I^{1,0}, \psi_I^{2,0}: C_{I,k}\to M_n$ by
\beq\label{8-N-4-19}
\psi_I^{i,0}(f)=\Sigma_I^i( f)\tforal f\in {\blue{C_{I,k}}},\,\,\,i=1,2.
\eneq
By {\blue{\eqref{8-N-4-17},}}
the choice of ${\cal H}_2,$ {\blue{\eqref{8-N-4-1nn},  and \eqref{8-N-4-5+},}}
we estimate that
\beq\label{8-N-4-20}
&&|{\rm tr}\circ \psi_I^{1,0}(1_{{\blue{{\blue{C_{I,k}}}}}})-{\rm tr}\circ \psi_I^{2,0}(1_{{\blue{C_{I,k}}}})|\\
&\le& |{\rm tr}\circ \Sigma_I^1(f_{0,k}\cdot P)-{\rm tr}\circ \phi_1(f_{0,k}\cdot P)|+|{\rm tr}\circ \phi_1(f_{0,k}\cdot P)-{\rm tr}\circ \phi_2(f_{0,k}\cdot P)|
\\\label{8-N-4-20-0}
&+&|{\rm tr}\circ \phi_2(f_{0,k}\cdot P)-\Sigma_I^2(f_{0,k}\cdot P)|
<1/N_0+\dt+1/N_0\\
&\le & \dt+\Delta(\widehat{f_{00}}\cdot P)\dt_1\min\{\Delta_{I,0}(\hat{g}): g\in {\cal H}_{1,I,0}\}/32K_0N_{00}.
\eneq
{\blue{Note that $n{\rm tr}\circ \psi_I^{1,0}(1_{{\blue{C_{I,k}}}})=\sum a_i z_i$ and $n{\rm tr}\circ \psi_I^{2,0}(1_{{\blue{C_{I,k}}}})=\sum_j b_jz_j$
(integer combinations), where  $z_i\in \{r(1), r(2),...,r(k_2)\}\subset \{l_1,l_2,...,l_{k_1}, r(1),r(2),...,r(k_2)\}$
(see \eqref{181127-412-T}).}}
It follows from  Lemma \ref{8-N-0} that there are two equivalent projections {\blue{$p_{1,0},\,p_{2,0}\in M_n$}} such that
$p_{i,0}$ commutes with $\psi_I^{i,0}(f)$ for all $f\in {\blue{C_{I,k}}}$ and $p_{i,0}\psi_I^{i,0}(1_{{\blue{C_{I,k}}}})={\blue{p_{i,0},}}$ $i=1,2,$
and {{so}} ({{by}} \eqref{8-N-4-20-0}), for $i=1,2,$
\beq\label{8-N-4-21}
&&\hspace{-0.2in}0\le {\rm tr}\circ \psi_I^{i,0}(1_{{\blue{C_{I,k}}}})-{\rm tr}(p_{i,0})<(1/N_0+\dt+1/N_0)+T/n\\\label{8-N-4-21+}
&&<\dt+\Delta(\widehat{f_{00}\cdot P})\dt_1\min\{\Delta_{I,0}(\hat{g}): g\in {\cal H}_{1,I,0}\}/32{{K_0}}N_{00}+T/n {\blue{ (<1/2)}}.
\eneq
Since $Y_{0,d/2}\subset Y_{k-1},$  ${\rm supp}(f_{00}){ =Y_{0,d/2}}\subset Y_{k-1}.$
Therefore ({{by}} \eqref{8-N-4-13-}),
\beq\label{10-5-N1}
{\rm tr}\circ \psi_I^{1,0}(1_{{\blue{C_{I,k}}}})\ge {\rm tr}\circ \psi_I^{1,0}(f_{00}\cdot P)\ge \Delta(\widehat{f_{00}\cdot P})\ge \sigma_0{\blue{>4N^I/n}}.
\eneq
{\blue{This, in particular, shows that $\psi_I^{1,0}(1_{{\blue{C_{I,k}}}})$ has rank at least $4N^I.$
Then, by \eqref{8-N-4-21+},  $p_{1,0}$ has rank at least $N^l.$}}
Moreover, {\blue{by  \eqref{10-5-N1}, \eqref{8-N-4-21+}, \eqref{8-N-4-5+}, \eqref{181127-412-sm}, \eqref{181127-412-dt}, and \eqref{8-N-4-13-},}}
\beq\label{10-5-N2}
\hspace{-0.4in}{\rm tr}(p_{2,0})&>&{\blue{\Delta(\widehat{f_{00}\cdot P})-(\dt+\Delta(\widehat{f_{00}\cdot P})\dt_1\min\{\Delta_{I,0}(\hat{g}): g\in {\cal H}_{1,I,0}\}/32{{K_0}}N_{00}+T/n)}}\\
&\ge& {\blue{31\Delta(\widehat{f_{00}\cdot P})/32-(\dt+T/N)
\ge \max\{31\sigma_0/32, 64/(N_0\dt_1)\}-(\dt+T/N)}}\\\label{10-5-N22}
&\ge& \max\{\sigma_0, 64/(N_0\dt_1)\}/2.
\eneq
Put $q_{i,0}=\psi_I^{i,0}(1_{{\blue{C_{I,k}}}})-p_{i,0},$ $i=1,2.$
There is a unitary $U_1\in M_n$ such that $U_1^*p_{1,0}U_1=p_{2,0}.$
Define $\psi_I^1: {\blue{C_{I,k}}}\to p_{2,0}M_np_{2,0}$ by
$\psi_I^1(f)=U_1^*p_{1,0}\psi_I^{1,0}(f)U_1$ for all $f\in {\blue{C_{I,k}}}$ and define
$\psi_I^2: {\blue{C_{I,k}}}\to p_{2,0}M_np_{2,0}$ by $\psi_I^2(f)=p_{2,0}\psi_I^{2,0}(f)$ for all $f\in {\blue{C_{I,k}}}.$
We compute
(using {\blue{\eqref{181128-412-d0},  (\ref{8-N-4-21}), \eqref{181127-412-dt},  \eqref{8-N-4-5+}, \eqref{8-N-4-13-},
\eqref{8-N-4-1}, and \eqref{181128-412-bt}}}) that
\beq\label{8-N-4-22}
{\rm tr}\circ \psi_I^1(g\circ r_{k}) &\ge & {\rm tr}\circ \psi_I^{1,0}((f_{0,0} Pg)\circ r_{k})-{\rm tr}(q_{1,0})\\
&>&  {\rm tr}\circ \psi_I^{1,0}(f_{0,0} Pg)-
1/16N_0^2-(\dt+2/N_0+T/N)\\
&>&  {\rm tr}\circ \psi_I^{1,0}(f_{0,0} Pg)-
5\min\{\Delta_{I,0}{\blue{(\hat{g})}}: g\in {\cal H}_{1,I,0}\}/64N_{00}\\
&>&(1-(1-\bt)/64) {\rm tr}\circ \psi_I^{1,0}(f_{0,0} Pg)\\
&>&\bt\Delta(\widehat{f_{0,0}Pg})
=\Delta_{I,0}(g)
\eneq
for all $g\in {\cal H}_{1,I,0}.$
Therefore,
\beq\label{8-N-4-23}
t\circ \psi_I^1(g)\ge \Delta_{I,0}(g)\tforal g\in {\cal H}_{1,I,k},
\eneq
where $t$ is the tracial state on $p_{2,0}M_np_{2,0}.$
We also estimate that {\blue{(using \eqref{8-N-4-17}, \eqref{8-N-4-1nn}, \eqref{10-5-N22}, \eqref{8-N-4-21}, \eqref{181127-412-dt},
and \eqref{8-N-4-13-})}},
\beq\label{8-N-4-24}
\hspace{-0.5in}|t\circ \psi_I^1(g)-t\circ \psi_I^2(g)| &=& (1/{\rm tr}({\blue{p_{1,0})}})|{\rm tr}\circ \psi_I^1(g)-{\rm tr}\circ \psi_I^2(g)|\\
&\le&  (1/{\rm tr}(p_{1,0}))|{\rm tr}\circ \psi_I^1(g)-{\rm tr}\circ \psi_I^{1,0}(g)|\\
&+& (1/{\rm tr}(p_{1,0}))|{\rm tr}\circ \psi_I^{1,0}(g)-
{\rm tr}(\phi_1(f_{0,k}\cdot 1_A\cdot  g))|\\
&+& (1/{\rm tr}(p_{2,0}))|{\rm tr}(\phi_1(f_{0,k}\cdot 1_A\cdot g))-{\rm tr}(\phi_2(f_{0,k}\cdot 1_A\cdot g))|\\
&+& (1/{\rm tr}(p_{2,0}))|{\rm tr}(\phi_2(f_{0,k}\cdot 1_A\cdot g)-{\rm tr}\circ \psi_I^{2,0}(g)|\\
&+&|{\rm tr}\circ \psi_I^{2,0}(g)-{\rm tr}\circ \psi_I^2(g)|\\\label{8-N-4-24+}
&<& (1/{\rm tr}(p_{2,0}))({\rm tr}(q_{1,0})+1/N_0+ \dt +1/N_0+{\rm tr}(q_{2,0}))\\\label{8-N-4-24++}
&<&{\blue{(N_0\dt_1/32)((4/N_0+2 \dt +2T/n)}}<\dt_1
\eneq
for all $g\in {\cal H}_{2,I,k}.$
{\blue{Recall}} that  $p_{2,0}$ has rank at least $N^I.$
It follows ({\blue{by \eqref{8-N-4-23} and \eqref{8-N-4-24++},}} by {{Corollary}} \ref{Aug-N-2}, {\blue{and  by the choice of ${\cal H}_{1,I,k},$ ${\cal H}_{2,I,k}$,  and $\dt_1$}}) that there are mutually orthogonal projections $e_0^I,e_1^I,e_2^I,...,e_{2K_0}^l\in
p_{2,0}M_np_{2,0}$ such that $e_0^I+\sum_{i=1}^{2K_0}e_i^I=p_{2,0},$ $e_0^I\lesssim e_1^I,$ and {{all}}
$e_j^I$ are equivalent to $e_1^I,$ two  unital \hm s $\psi_{1,I,0}, \psi_{2,I,0}: {\blue{C_{I,k}}}{\blue{\cong {\blue{C_{I,0}}}}}\to e_0^IM_ne_0^I,$ a
unital \hm\, $\psi_{I}: {{{\blue{C_{I,k}}}}}\to e_1^IM_ne_1^I$, and a
unitary $u_1\in p_{2,0}M_np_{2,0}$
such that
\beq\label{8-N-4-25}
\|{\rm Ad}\, u_1\circ \psi_I^1(f)-(\psi_{1,I,0}(f)\oplus {\rm diag}(\overbrace{\psi_{I}(f), \psi_I(f),...,\psi_I(f)}^{{\blue{2K_0}}}))\|<\ep_0/16\\
\andeqn\,\,\,
\|\psi_I^2(f)-(\psi_{2,I,0}(f)\oplus {\rm diag}(\overbrace{\psi_{I}(f), \psi_I(f),...,\psi_I(f)}^{2K_0}))\|<\ep_0/16
\eneq
for all $ f\in {{{\cal F}_{I,k}'}}.$
By (\ref{8-N-4-8}), {{it follows}} that, for all $ f\in {{{\cal F}_{I,k}}},$
\beq\label{8-N-4-26}
\|{\rm Ad}\, u_1\circ \psi_I^1(f)-(\psi_{1,I,0}(f)\oplus {\rm diag}(\overbrace{\psi_{I}(f), \psi_I(f),...,\psi_I(f)}^{2K_0}))\|<\ep_0/8\\
\label{8-N-4-26+}
\andeqn\,\,\,
\|\psi_I^2(f)-(\psi_{2,I,0}(f)\oplus {\rm diag}(\overbrace{\psi_{I}(f), \psi_I(f),...,\psi_I(f)}^{2K_0}))\|<\ep_0/8.
\eneq

For each $x\in X\setminus Y_k$ such that $\pi_{x,j}$ appeares in ${{\Sigma_b^1}},$ or ${{\Sigma_b^2}},$
by {\blue{\eqref{8-N-4-4+}}},
\beq\label{8-N-4-27}
\|\pi_{x,j}(f)-\pi_{x,j}\circ s^d\circ \pi_I(f))\|<\ep_0/16\tforal f\in {\cal F}.
\eneq
Define $\Sigma_{\pi,b,i}:=\Sigma_b^i\circ s^d: A/I\to M_n,$ $i=1,2.$

Define $\Phi_1: A/I\to (1-p_{2,0})M_n(1-p_{2,0})$ by
\beq\label{8-N-4-28}
\Phi_1(f)={\rm Ad}\, U_1\circ (\Sigma_\pi^1\oplus \Sigma_{\pi, b,1})(f)\rforal f\in A/I.
\eneq
Define $\Phi_2: A/I\to (1-p_{2,0})M_n(1-p_{2,0})$ by
\beq\label{8-N-4-29}
\Phi_2(f)=(\Sigma_\pi^1\oplus \Sigma_{\pi, b,2})(f)\rforal f\in A/I.
\eneq
Note that
\beq\label{8-N-4-29+}
\Phi_1(1_{A/I})=\Sigma_\pi^1(g_{0,k})\oplus \Sigma_b^1(g_{0,k})\andeqn
\Phi_2(1_{A/I})=\Sigma_\pi^1(g_{0,k})\oplus \Sigma_b^2(g_{0,k}).
\eneq
We {{compute}} that
\beq\label{8-N-4-30}
\hspace{-0.8in}|{\rm tr}\circ \Phi_1(1_{A/I})-{\rm tr}\circ \Phi_2(1_{A/I})|
&\le & |{\rm tr}\circ \Phi_1(1_{A/I})-{\rm tr}\circ \phi_1(g_{0,k})|\\
&&\hspace{-1.3in}+|{\rm tr}\circ \phi_1(g_{0,k})-{\rm tr}\circ \phi_2(g_{0,k})|+
|{\rm tr}\circ \phi_2(g_{0,k})-{\rm tr}\circ \Phi_2(g_{0,k})|\\
&<& 1/N_0+\dt+1/N_0\\
&<&\dt+\Delta(\widehat{f_{00}\cdot P})\dt_{ 2}\min\{\Delta_{ \pi}(\hat{g}): g\in {\cal H}_{1,{ \pi}}\}/32N_{00}.\hspace{-0.4in}.
\eneq
It follows from Lemma \ref{8-N-0} that there are two mutually equivalent projections $p_{1,1}$ and
$p_{2,1}\in (1-p_{2,0})M_n(1-p_{2,0})$ such that
$p_{i,1}$ commutes with $\Phi_i(f)$ for all $f\in A/I$ and $p_{i,1}\Phi_i(1_{A/I})=p_{i,1}.$ $i=1,2,$
and, for $i=1,2,${\rm tr}
\beq\label{8-N-4-32}
&&\hspace{-0.3in}0\le {\rm tr}\circ \Phi_i(1_{A/I})-{\rm tr}(p_{i,1})<{\blue{\dt+2/N_0+T/n}}\\\label{8-N-4-32-+}
&&<\dt+
\Delta(\widehat{f_{00}})\dt_{ 2}\min\{\Delta_{\pi}(\hat{g}): g\in {\cal H}_{1,\pi}\}/{ (32N_{00})}+T/n {\blue{(<1/2)}}.
\eneq
{{Since $g_0(x)=0,$ if ${\rm dist}(x, Y_{0,d/2})\le d/16,$}} we have, {\blue{by \eqref{8-N-4-5+},}}
$${\rm tr}\circ \Phi_1(1_{A/I}) > \DT(\widehat{g_0\cdot s^d(1)})> \DT_{\pi}(\hat{1})\geq {{\max}}\{\sigma_0, {\blue{64KN_{00}}}/(N_0\dt_2)\},$$
and  (by \eqref{8-N-4-32-+}) ${\rm tr}(p_{2,1})\geq \max\{\sigma_0, 64K_0N_{00}/(N_0\dt_2)\}/2.$
{\blue{Since $\sigma_0>4N^\pi/N\ge 4N^\pi/n$ (see \eqref{8-N-4-13-}), $\Phi_i(1_{A/I})$ has rank at least $4N^\pi.$
Then, by \eqref{8-N-4-32-+}, {{it follows}} that $p_{i,1}$ has rank at least $N^\pi,$ $i=1,2.$}}

Put $q_{i,1}=\Phi_i(1_{A/I})-p_{i,1},$ $i=1,2.$
There is a unitary $U_2\in (1-p_{2,0})M_n(1-p_{2,0})$ such that $U_2^*p_{1,1}U_2=p_{2,1}.$
Define $\Phi_\pi^1: A/I\to p_{2,1}M_np_{2,1}$ by
$\Phi_\pi^1(f)=U_2^*p_{1,1}\Phi_1(f)U_2$ for all {\blue{$f\in  A/I$}}
and define
$\Phi_\pi^2: A/I\to p_{2,1}M_np_{2,1}$ by $\Phi_\pi ^2(f)=p_{2,1}\Phi_2(f)$ for all $f\in A/I.$

{\blue{Since $g_0(x)=0,$ if ${\rm dist}(x, Y_{0,d/2})\le d/16,$ and $\{g_0s^d(g): g\in {\cal H}_{1,\pi}\}\subset {\cal H}_1,$}} we compute (using {\blue{ \eqref{8-N-4-1}, \eqref{8-N-4-32}, \eqref{181127-412-dt},  \eqref{8-N-4-13-},
 \eqref{181128-412-bt}, and \eqref{8-N-4-9}}}) {{that}}
\beq\label{8-N-4-33}
{\rm tr}\circ \Phi_\pi^1(g) &\ge & {\blue{\phi_1(g_0 s^d(g))-{\rm tr}(q_{1,1})}}>
 \Delta(\widehat{g_0 s^{ d}(g)})-\sigma_0/N_{00}\\
&>& \bt \Delta(\widehat{g_0 s^{ d}(g)})=\Delta_{\pi}(g)
\eneq
for all $g\in {\cal H}_{1,\pi}.$
Therefore, for the tracial state $t_1$ of $p_{2,1}M_np_{2,1},$
\beq\label{8-N-4-34}
t_1\circ \Phi_\pi^1(g)\ge \Delta_{\pi}(g)\tforal g\in {\cal H}_{1,\pi}.
\eneq

We also estimate  ({{in a way}} similar to the estimate of (\ref{8-N-4-24++})) that, for all $g\in {\cal H}_{2,\pi},$
\beq\label{8-N-4-35}
\hspace{-0.6in}|t_1\circ \Phi_\pi^1(g)-t_1\circ \Phi_\pi^2(g)| &=& (1/{\rm tr}(p_{2,1}))|{\rm tr}\circ \Phi_\pi^1(g)-{\rm tr}\circ \Phi_\pi^2(g)|\\
&\le&  (1/{\rm tr}(p_{2,1}))|{\rm tr}\circ \Phi_\pi^1(g)-{\rm tr}\circ \Phi_1(g)|\\
&+& (1/{\rm tr}(p_{2,1}))|{\rm tr}\circ \Phi_1(g)-{\rm tr}\circ \phi_1(g_{0,k} s^d(g))|\\
&&\hspace{-0.3in}+ (1/{\rm tr}(p_{2,1}))|{\rm tr}\circ \phi_1(g_{0,k} s^d(g))-{\rm tr}\circ \phi_2(g_{0,k} s^d(g))|\\
&+& (1/{\rm tr}(p_{2,1}))|{\rm tr}\circ \phi_2(g_{0,k} s^d(g))-{\rm tr}\circ \Phi_2(g)|\\
&+&(1/{\rm tr}(p_{2,1}))|{\rm tr}\circ \Phi_2(g)-{\rm tr}\circ \Phi_\pi^2(g)|\\\label{8-N-4-35+}
&<&  (1/{\rm tr}(p_{2,1}))(1/N_0+ \dt +1/N_0)<\dt_{ 2}.
\eneq
Recall that $p_{2,1}$ has rank at least $N^\pi.$ It follows from {\blue{the induction assumption that
the theorem holds {{for}} $A/I$
(and from \eqref{8-N-4-34} and \eqref{8-N-4-35+})}} that there are mutually orthogonal projections $e_0^\pi,e_1^\pi,e_2^\pi,...,e_{2K_0}^\pi\in
p_{2,1}M_np_{2,1}$ such that $e_0^\pi\lesssim e_1^\pi$ and {{all}}
$e_j^\pi$ are equivalent to $e_1^\pi,$ two  unital \hm s $\psi_{1,\pi,0}, \psi_{2,\pi,0}: A/I\to e_0^\pi M_ne_0^\pi,$ a
unital \hm\, $\psi_{\pi}: A/I\to e_1^\pi M_ne_1^\pi,$ and a
unitary $u_2\in p_{2,1}M_np_{2,1},$
such that
\beq\label{8-N-4-36}
\|{\rm Ad}\, u_2\circ \Phi_\pi^1(f)-(\psi_{1,\pi,0}(f)\oplus {\rm diag}(\overbrace{\psi_{\pi}(f), \psi_\pi(f),...,\psi_\pi(f)}^{2K_0}))\|<\ep_0/16\\
\andeqn\,\,\,
\|\Phi_\pi^2(f)-(\psi_{2,\pi,0}(f)\oplus {\rm diag}(\overbrace{\psi_{\pi}(f), \psi_\pi(f),...,\psi_\pi(f)}^{2K_0}))\|<\ep_0/16
\eneq
for all $ f\in {\cal F}_{\pi}.$
Let ${{\psi_\pi^{1}}}: A\to  p_{2,1}M_np_{2,1}$ {{be defined}} by
$\psi_\pi^{1}(f)={\rm Ad}\, u_2\circ{\rm Ad}\, U_2 (p_{2,1}(\Sigma_\pi^1\oplus \Sigma_b^1)(f))$ and
define $\psi_\pi^{2}: A\to p_{2,1}M_np_{2,1}$ by
$\psi_\pi^{2}(f)=p_{2,1}(\Sigma_\pi^2\oplus \Sigma_b^2)(f)$
 for all $f\in A.$
Then, by (\ref{8-N-4-27}),
\beq\label{8-N-4-37}
\|\psi_\pi^{i}(f)-(\psi_{i,\pi,0}\circ\pi_I(f)\oplus {\rm diag}(\overbrace{\psi_{\pi}(\pi_I(f)), \psi_\pi(\pi_I(f)),...,\psi_\pi(\pi_I(f))}^{2K_0}))\|<\ep_0/8
\eneq
for all $f\in {\cal F},$ $i=1,2.$

Put $e_i={\blue{\sum_{j=1}^{2N_{00}} e_{2(i-1)N_{00}+j}^I\oplus \sum_{j=1}^{2N_{00}}e_{2(i-1)N_{00}+j}^\pi}},$
$i=1,2,...,K.$  {\blue{Denote by $\bar{\psi}_I$  and ${\bar \psi_\pi}$  the direct sums of
$2N_{00}$ copies of $\psi_I$ and $\psi^\pi,$ respectively.}}
Define $\psi: A\to  e_1M_ne_1$ by
$$
\psi(f)={\rm diag}({\bar \psi_I}(f|_{\blue{Y_k}}),{\bar \psi}_\pi({\blue{\pi_I}}(f)))
$$
 for all $f\in A.$
 By (\ref{8-N-4-21}),  (\ref{8-N-4-32}),
 {{and}}
{\blue{\eqref{8-N-4-17}, for $i=1,2,$}}
 \beq\label{8-N-4-38}
 {\rm tr}(q_{i,0})+{\rm tr}(q_{i,1})+{\rm tr}(\Sigma_s^i(1_A))<  1/64K_0+ 1/64K_0+1/N_0<1/16K.
 \eneq
We have, for $f\in A,$
\beq\label{8-N-4-39}
\phi_2(f)={{\psi_\pi^{2}}}(f)\oplus q_{2,1}(\Sigma_\pi^2+\Sigma_b^2)(f)\oplus \Sigma_s^2(f)\oplus \psi_I^2(f|_{{{Y_{k}}}})\oplus q_{2,0}{{\psi_I^{2,0}}}(f{{|_{Y_{k}}}}).
\eneq
Put $e_0=e_0^I\oplus e_0^\pi+q_{2,1}+\Sigma_s^2(1_A)+q_{2,0}.$
Then
\beq\label{8-N-4-50}
{\rm tr}(e_0)<{\rm tr}(e_0^I)+{\rm tr}(e_0^\pi)+1/16K\le {\rm tr}(e_1).
\eneq
In other words, $e_0\lesssim e_1.$ Moreover $e_1$ is equivalent to each $e_i,$ $i=1,2,...,K.$
{{Define $h_2: A\to e_0M_ne_0$ by,  for each $f\in A,$
\beq\label{181128-412h2}
h_2(f)=\psi_{2,I,0}(f|_{Y_{k}})\oplus \psi_{2,\pi,0}(\pi_I(f)) \oplus q_{2,1}(\Sigma_\pi^2+\Sigma_b^2)(f)\oplus \Sigma_s^2(f)\oplus q_{2,0}\psi_I^{2,0}(f|_{Y_{k}}).
\eneq
}}
It follows from (\ref{8-N-4-26+}), (\ref{8-N-4-37}), and {\blue{\eqref{181128-412h2}}} that
\beq\label{8-N-4-51}
\|\phi_2(f)-(h_{{2}}(f)\oplus {\rm diag}(\overbrace{{\bar\psi}(f), {\bar \psi}(f),...,{\bar \psi}(f)}^K)\|<\ep_0/8\rforal f\in {\cal F}.
\eneq
Similarly, there {{exist}} a unitary $U\in M_n$ and a unital \hm\, ${\blue{h_1}}: {{A\to}} e_0M_ne_0$ such that
\beq\label{8-N-4-52}
\|{\rm Ad}\, U\circ \phi_1(f)-(h_{{1}}(f)\oplus{\rm diag}(\overbrace{{\bar \psi}(f), {\bar \psi}(f),...,{\bar \psi}(f)}^K)\|<\ep_0/8\rforal f\in {\cal F}.
\eneq
Since we  assume that ${\cal H}_0\subset {\cal F},$  {\blue{by \eqref{8-N-4-1},}}
 the choice of $\ep_0$, {\blue{and \eqref{181128-412-bt},   we also have}}
that
\beq\label{8-N-4-53}
\hspace{-0.4in}{\rm tr}\circ {\bar \psi}(g) &\ge & {\blue{(1/K)({\rm tr}\circ \phi_1(g)-\ep_0/8-{\rm tr}(h_1(g)))}}\\
&>&{\blue{(1/K)(\Delta(\hat{g})-\dt_{00}/16KN_{00}-1/KN_{00})>}}
\af{\Delta({\blue{\hat{g}}})\over{K}}\tforal g\in {\cal H}_0.
\eneq
Thus the {{conclusion of the}} theorem holds for $m.$

\end{proof}


\begin{rem}\label{8-N-4-r}
{\rm
If we assume that $A$ is infinite dimensional, then Lemma \ref{8-N-4} still holds  without the assumption about the integer $N.$  This could be easily seen by taking a larger ${\cal H}_1$ which contains {\blue{at least $N$ mutually orthogonal non-zero elements}} as we remarked in \ref{8-N-2-r}.

}
\end{rem}

{\blue{ 
\begin{cor}\label{repcor}
Let  $A_0\in {\overline {\cal D}}_s$ be a unital \CA,
let $\ep>0$, and
let ${\cal F}\subset A_0$ be a finite subset.
Let $\Delta: (A_0)_+^{q, {\bf 1}}\setminus \{0\}\to (0,1)$ be an order preserving  map.

Suppose that ${\cal H}_1\subset (A_0)_+^{\bf 1}\setminus \{0\}$ is a finite subset,  $\sigma>0$
is {{a}} positive number and $n\ge 1$ is an integer.
There exists a finite subset ${\cal H}_2\subset (A_0)_+^{\bf 1}\setminus\{0\}$ satisfying the following {{condition}}:
Suppose that $\phi: A=A_0\otimes C(\T)\to M_k$ (for some integer $k\ge 1$) is a unital \hm\, and
\beq\label{repcor-1}
tr\circ \phi(h\otimes 1)\ge \Delta(\hat{h})\tforal h\in {\cal H}_2.
\eneq
Then there exist mutually orthogonal  projections $e_0,e_1,e_2,...,e_n\in
M_k$ such that
$e_1, e_2,...,e_n$ are equivalent and $\sum_{i=0}^n e_i=1,$ and there {{exist}}
unital \hm s $\psi_0: A= A_0\otimes C(\T)\to e_0M_ke_0$ and $\psi: A=A_0\otimes C(\T)\to e_1M_ke_1$ such that
\beq\label{repcor-2}
&&\|\phi(f)-{\rm diag}( \psi_0(f),\overbrace{\psi(f),\psi(f),...,\psi(f)}^n)\|<\ep\\
&&\tand {\rm tr}(e_0)<\sigma
\eneq
for all $f\in {\cal F},$ where $tr$ is the tracial state on $M_k.$
Moreover, {{$\psi$ can be chosen such that }}
\beq\label{repcor-3}
tr(\psi(g\otimes 1))\ge {\Delta(\hat{g})\over{{\blue{2n}}}}\tforal g\in {\cal H}_1.
\eneq
\end{cor}
}}

{
\begin{proof}
The statement {\blue{follows directly}} from Lemma \ref{8-N-4} with $\phi_1 = \phi_2 = \phi$, $K=n,$ and $\alpha=\frac{1}{2}$.
\end{proof}
}

The following is known and is taken from Theorem 3.9 of \cite{LnAUCT}.

\begin{thm}\label{Lnuct}
Let $A$ be a unital separable amenable \CA\,  which satisfies the UCT and let $B$ be a unital \CA.  Suppose that
$h_1,\, h_2: A\to B$ are  \hm s such that
$$
[h_1]=[h_2]\,\,\,{\text in}\,\,\, KL(A,B).
$$
Suppose that $h_0: A\to B$ is a unital full monomorphism. Then, for any $\ep>0$ and   any
finite subset ${\cal F}\subset A,$ there {{exist}} an integer $n\ge 1$ and {a partial  isometry} $W\in M_{n+1}(B)$
such that
$$
\|W^*{\rm diag}(h_1(a),h_0(a),...,h_0(a))W-{\rm diag}(h_2(a),h_0(a),...,h_0(a))\|<\ep
$$
for all $a\in {\cal F},$  {{and
$W^*pW=q,$}} where
\beq\nonumber
p={\rm diag}(h_1(1_A),h_0(1_A),...,h_0(1_A))\tand q={\rm diag}(h_2(1_A),h_0(1_A),...,h_0(1_A)).
\eneq
In particular,  if $h_1(1_A)=h_2(1_A),$ {{we may choose}} $W\in U(pM_{n+1}(B)p).$
\end{thm}

\begin{proof}
This is a slight variation of Theorem 3.9 of \cite{LnAUCT}.  If $h_1$ and $h_2$ are both unital, then
it is exactly the same as Theorem 3.9 of \cite{LnAUCT}.  So suppose that $h_1$ is not unital. Let $A'=\C\oplus A.$
Choose $p_0=1_B-h_1(1_A)$ and $p_1={\rm diag}(p_0, 1_B).$ Put $B'=p_1M_2(B)p_1.$
Define $h_1': A'\to B'$ by $h_1'(\lambda \oplus a)=\lambda\cdot {\rm diag}(p_0,p_0)\oplus h_1(a)$ for all
$\lambda\in \C$ and $a\in A,$  and define $h_2': A'\to B'$ by
$h_2'(\lambda \oplus a)=\lambda \cdot {\rm diag}(p_0, 1_B-h_2(1_A))\oplus h_2(a)$ for all $\lambda\in \C$ and $a\in A.$
Then $h_1'$ and $h_2'$ are unital and $[h_1']=[h_2']$ in $KL(A', B').$ Define
$h_0': A'\to B'$ by $h_0'(\lambda\oplus a)=\lambda\cdot p_0\oplus h_0(a)$ for all $\lambda\in \C$ and $a\in A.$
Note that $h_0'$ is full in $B'.$
So, Theorem 3.9 of \cite{LnAUCT} applies. It follows that there {{are}} an integer $n\ge 1$ and a unitary
$W'\in M_{n+1}(B')$ such that
\beq\nonumber
\|(W')^*{\rm diag}(h_1'(a),h_0'(a),...,h_0'(a))W'-{\rm diag}(h_2'(a),h_0'(a),...,h_0'(a))\|<\min\{1/2, \ep/2\}
\eneq
for all $a\in {\cal F}\cup\{1_A\}.$  In particular,
\beq\label{lnauct-3}
\|(W')^*pW'-q\|<\min\{1/2,\ep/2\}.
\eneq
There is a unitary $W_1\in M_{n+1}(B')$ such that
\beq\label{linauct-3}
\|W_1-1_{M_{n+1}(B')}\|<\ep/2\andeqn W_1^*(W')^*pW'W_1=q.
\eneq
Put ${\blue{W=pW'W_1q}}.$ Then
\beq\label{linatct-4}
\|W^*{\rm diag}(h_1(a),h_0(a),...,h_0(a))W-{\rm diag}(h_2(a),h_0(a),...,h_0(a))\|<\ep
\eneq
 for all $a\in {\cal F}$, {{as desired.}} (The last statement is clear.)
\end{proof}

\begin{lem}\label{Lauct2} {\rm (cf.{{ 5.3 of \cite{Lnjotuni}, Theorem 3.1 of}} \cite{GL-almost-map}, \cite{Da1}, 5.9 of \cite{LnAUCT}, and
{{Theorem 7.1 of \cite{Lin-hmtp}}})}
Let $A$ be a unital separable amenable  \CA\, which satisfies the UCT and let $\Delta: A_+^{q, {\bf 1}}\setminus\{0\} \to (0,1)$ be an order preserving
map.  For any $\ep>0$ and any finite subset ${\cal F}\subset A,$ there {{exist}}
$\dt>0,$ a finite subset ${\cal G}\subset A,$
a finite subset ${\cal P}\subset \underline{K}(A),$ a finite subset ${\cal H}\subset A_+^{\bf 1}\setminus \{0\}$, and an integer $K\ge 1$ satisfying the following {{condition}}:
For any two unital ${\cal G}$-$\dt$-multiplicative
{\blue{\cp s}} $\phi_1, \phi_2: A\to M_n$ (for some integer $n$) and any
unital ${\cal G}$-$\dt$-multiplicative {\blue{\cp}} $\psi: A\to M_m$ with $m\ge n$ such that
\beq\label{Lauct-1}
{\rm tr}\circ \psi(g)\ge \Delta(\hat{g})\tforal g\in {\cal H}\tand\,
[\phi_1]|_{\cal P}=[\phi_2]|_{\cal P},
\eneq
{\blue{where ${\rm tr}$ is the tracial state of $M_m,$}}
there exists a unitary $U\in M_{Km+n}$ such that
\beq\label{Lauct-2}
\|{\rm Ad}\, U\circ (\phi_1\oplus \Psi)(f)-(\phi_2\oplus \Psi)(f)\|<\ep\tforal f\in A,
\eneq
where
$$
\Psi(f)={\rm diag}(\overbrace{\psi(f), \psi(f),...,\psi(f)}^K)\tforal f\in A.
$$
\end{lem}

\begin{proof}
This follows, {{as we shall show}}, from {Theorem} \ref{Lnuct}.

Fix $\Delta$ as given.
Suppose {{that the conclusion}} is false. Then there exist  $\ep_0>0$ and a finite subset
${\cal F}_0\subset A,$  an increasing sequence of finite subsets
$\{{\cal P}_n\}$ of $\underline{K}(A)$ {{with}} union $\underline{K}(A),$  an increasing sequence
of finite subsets $\{{\cal H}_{n}\}\subset A_+^{\bf 1}\setminus \{0\}$ {{with}} union  dense
in $A_+^{\bf 1}$ {{and {{such that}} if $a\in {\cal H}_n$ and $f_{1/2}(a)\not=0,$ then
$f_{1/2}(a)\in {\cal H}_{n+1},$}}
three increasing sequences of    integers  $\{R(n)\}$,  $\{r(n)\}$, {{and}} $\{s(n)\}$ (with $s(n)\ge r(n)$),
 two sequences of {\blue{unital \cp s}}
 $\phi_{1,n}, \phi_{2,n}, : A\to M_{r(n)}$
 with the properties that
 \beq\label{Lauct-3-}
 [\phi_{1,n}]|_{{\cal P}_n}&=&[\phi_{2,n}]|_{{\cal P}_n}\andeqn\\
 \lim_{n\to\infty}\|\phi_{i,n}(ab)-\phi_{i,n}(a)\phi_{i,n}(b)\|&=&0\rforal a,b \in A,\,\,\,i=1,2,
 \eneq
 and a sequence of unital {\blue{\cp s}}
 $\psi_n: A\to M_{s(n)}$ with the properties
 that
 \beq\label{Lauct-3}
&& {\rm tr}_n\circ \psi_n(g)\ge \Delta(\hat{g})\rforal g\in {\cal H}_{n}\andeqn\\
 && \lim_{n\to\infty}\|\psi_{n}(ab)-\psi_{n}(a)\psi_{n}(b)\|=0\rforal a,b \in A
 \eneq
such that
\beq\label{Lauct-4}
\inf\{ \sup\{\|{\rm Ad}\, U_n\circ (\phi_{1,n}(f)\oplus {\tilde \psi}^{R(n)}_n(f))-(\phi_{2,n}(f)\oplus {\tilde \psi}^{R(n)}_n(f)\|: f\in {\cal F}\}\}\ge \ep_0,
\eneq
where ${\rm tr}_n$ is the normalized trace  on ${\blue{M_{s(n)}}},$
${\tilde \psi}^{(R(n))}(f)={\rm diag}(\overbrace{\psi_n(f),\psi_n(f),...,\psi_n(f)}^{R(n)})$ for all $f\in A,$
and the infimum  is taken among all unitaries $U_n\in {\blue{M_{r(n)+R(n)s(n)}}}.$

Note that, by (\ref{Lauct-3}), since $\{{\cal H}_n\}$ is increasing, for any $g\in {\cal H}_{n}\subset A_+^{\bf 1},$  we compute that
\beq\label{Lauct-5}
{\rm tr}_m(p_m)\ge \Delta(\hat{g})/2,
\eneq
where $p_m$ is the spectral projection of $\psi_m(g)$ corresponding to the subset $\{\lambda>\Delta(\hat{g})/2\}$
for all $m\ge n.$
It follows that (for all sufficiently large $m$)  there are {{elements}} $x_{g,i,m}\in M_{s(m)}$ with $\|x_{g, i,m}\|\le 1/\Delta(\hat{g}),$
$i=1,2,...,N(g)$, such that
\beq\label{Lauct-6}
\sum_{i=1}^{N(g)} x_{g,i,m}^*\psi_m(g)x_{g,i,m}=1_{s(m)},
\eneq
where $1\le N(g)\le 1/\Delta(\hat{g})+1.$
Define  $X_{g,i}=\{x_{g,i,m}\},$ $i=1,2,...,N(g).$ Then $X_{g,i}\in \prod_{n=1}^{\infty}M_{r(n)}.$
Let $Q(\{M_{r(n)}\})=\prod_{n=1}^{\infty}M_{r(n)}/\bigoplus_{n=1}^{\infty}M_{r(n)},$
$Q(\{M_{s(n)}\})=\prod_{n=1}^{\infty}M_{s(n)}/\bigoplus_{n=1}^{\infty}M_{s(n)},$ and let
$\Pi_1: \prod_{n=1}^{\infty}M_{r(n)}\to \prod_{n=1}^{\infty}M_{r(n)}/\bigoplus_{n=1}^{\infty}M_{r(n)},$
$\Pi_2:\prod_{n=1}^{\infty}M_{s(n)}\to \prod_{n=1}^{\infty}M_{s(n)}/\bigoplus_{n=1}^{\infty}M_{s(n)}$ be  the quotient maps.
Denote by $\Phi_i: A\to Q(\{M_{r(n)}\})$ the \hm s
$\Pi_1\circ \{\phi_{i,n}\}$ and denote by ${\bar \psi}: A\to Q(\{M_{s(n)}\})$ the \hm\, $\Pi_2\circ \{\psi_n\}.$
For each $g\in \bigcup_{n=1}^{\infty} {\cal H}_n,$
\beq\label{Lauct-7}
\sum_{i=1}^{N(g)}\Pi_2(X_{i,g})^*{\bar \psi}(g)\Pi_2(X_{i,g})=1_{Q(\{M_{s(n)}\})}.
\eneq
{{Note that if $g\in \bigcup_{n=1}{\cal H}_n$ and $f_{1/2}(g)\not=0,$ then $g_{1/2}(g)\in \bigcup_{n=1}{\cal H}_n,$ and
$\bigcup_{n=1}{\cal H}_n$ is dense in $A_+^{\bf 1}.$}} This implies  that ${\bar \psi}$ is full {\blue{(see the proof of
Proposition \ref{full-2018-sept})}}.
Note that both $\prod_{n=1}^{\infty} M_{r(n)}$ and $Q(\{M_{r(n)}\})$ have stable rank one and real rank zero.
One then computes {\blue{(see {{Corollary}} 2.1 of \cite{GL-almost-map})}} that
\beq\label{Lauct-10}
[\Phi_1]=[\Phi_2]\,\,\,{\rm in}\,\,\, KL(A, Q(\{M_{r(n)}\}).
\eneq
By applying Theorem \ref{Lnuct} (Theorem 3.9 of \cite{LnAUCT}),
one obtains an integer $K\ge 1$ and a unitary $U\in PM_{K+1}(Q(\{M_{s(n)}\})P,$ where $P={\rm diag}(1_{Q(\{M_{r(n)}\}}, 1_{M_K(Q(\{M_{s(n)}\}))}),$
such that
\beq\label{Lauct-11}
\|{\rm Ad}\, U\circ (\Phi_1(f)\oplus {\rm diag}(\overbrace{{\bar \psi}(f),...,{\bar \psi}(f)}^K)
-(\Phi_2(f)\oplus {\rm diag}(\overbrace{{\bar \psi}(f),...,{\bar \psi}(f)}^K)\|<\ep_0/2
\eneq
for all $f\in {\cal F}_0.$
It follows that there are unitaries
$$
\{U_n\}\in \prod_{n=1}^{\infty}M_{r(n) +Ks(n)}
$$
such that, for all
large $n,$
\beq\label{Lauct-12}
\hspace{-0.4in}\|{\rm Ad}\, U_n\circ \phi_{1,n}(f)\oplus {\rm diag}(\overbrace{\psi_n(f),...,\psi_n(f)}^K)-
\phi_{2,n}(f)\oplus {\rm diag}(\overbrace{\psi_n(f),...,\psi_n(f)}^K)\|<\ep_0/2
\eneq
for all $f\in {\cal F}_0.$  This {{is in direct}} contradiction with (\ref{Lauct-4}) when we choose $n$ with $R(n)\ge K.$
\end{proof}

\begin{rem}\label{Re414}
{{The preceding lemma}} holds in a much more general  setting and variations of it have appeared. We state this version here for  our immediate purpose (see \ref{Suni} and part (1) of \ref{Rsuni} for more comments).
\end{rem}

It should be noted that, in the following statement the integer $L$ and {{the map}} $\Psi$ depend not only on $\ep,$ ${\cal F},$ and
${\cal G},$  but also on $B,$ as well as $\phi_1$ and $\phi_2.$

\begin{lem}\label{Newstableuniq}
Let $C$ be a unital amenable separable residually finite dimensional \CA\, which satisfies the UCT. For any $\ep>0$ and  any finite subset
${\cal F}\subset C,$ there {{exist}} a finite subset ${\cal G}\subset C,$
$\dt>0,$ {{and}} a finite subset ${\cal P}\subset \underline{K}(C)$
satisfying the following {{condition:}}
For any unital ${\cal G}$-$\dt$-multiplicative {\blue{\cp s}} $\phi_1, \phi_2:
C\to A$ (for any unital \CA\, $A$)
such that
\beq\label{New-1}
[\phi_1]|_{\cal P}=[\phi_2]|_{\cal P},
\eneq
there exist an integer $L\ge 1,$ a  unital \hm~ $\Psi: C\to M_L\subset M_L(A)$, and  a unitary $U\in U(M_{L+1}(A))$   such that,
for all $f\in {\cal F},$
\beq\label{New-2}
\|{\rm Ad}\, U\circ {\rm diag}(\phi_1(f), \Psi(f))-{\rm diag}(\phi_2(f),
\Psi(f)\|<\ep.
\eneq
\end{lem}

\begin{proof}
The proof is almost the same as that of Theorem 9.2 of \cite{Linajm}.
Suppose that the conclusion is false.
We then obtain a  positive number $\ep_0>0,$ a finite {\blue{subset}} ${\cal F}_0\subset C,$ a sequence  of finite
subsets ${\cal P}_n\subset \underline{K}(C)$ with ${\cal P}_n\subset {\cal P}_{n+1}$ and
$\bigcup_n {\cal P}_{n+1}=\underline{K}(C),$ a sequence of unital \CA s $\{A_n\},$ {{and}}  sequences
of unital {\blue{\cp s}} $\{L_n^{(1)}\}$ and $\{L_n^{(2)}\}$, from $C$ to $A_n$,
such that
\beq\label{12/31/24-1}
&&\hspace{-0.3in}\lim_{n\to\infty}\|L_n^{(i)}(ab)-L_n^{(i)}(a)L_n^{(i)}(b)\|=0\rforal a, b\in C,\\
&&\hspace{-0.3in}{[}L_n^{(1)}{]}|_{{\cal P}_n}={[}L_n^{(2)}{]}|_{{\cal P}_n}, \tand \\\label{12/31/24-1+}
&&\hspace{-0.3in}\inf\{\sup\{\|u_n^*{\rm diag}(L_n^{(1)}(a), \Psi_n(a)){{u_n}}-{\rm diag}(L_n^{(2)}(a), \Psi_n(a)\|: a\in {\cal F}_0{\blue{\}\}}}\ge \ep_0,
\eneq
where the infimum
is taken among all integers $k >1,$  all possible unital homomorphisms
$\Psi_n: C\to  M_k,$ and all possible unitaries ${\blue{u_n}} \in M_{k+1}(A_n).$  We may assume
that $1_C\in {\cal F}.$
Define $B_n=A_n\otimes {\cal K},$ $B=\prod_{n=1}^{\infty} B_n,$ and $Q_1=B/\bigoplus_{n=1}^{\infty} B_n.$
Let $\pi: B\to Q_1$ be the quotient map. Define $\phi_j:C\to B$ by $\phi_j(a)=\{L_n^{(j)}(a)\}$ and
let  $\bar{\phi_j}=\pi\circ \phi_j,$ $j=0,1.$  Note that $\bar{\phi_j}: C\to Q_1$ is a \hm.
As in the proof of 9.2 of \cite{Linajm},
we have
$$
[{\bar{ \phi_1}}]=[\bar{\phi_2}]\,\,\,{\rm in}\,\,\, KL(C, Q_1).
$$
%

Fix an irreducible representation  $\phi_0': C\to M_{r}.$
Denote by $p_n$ the unit of  the unitization ${\tilde B_n}$ of $B_n,$ $n=1,2,....$
Define a  \hm\,
$\phi_0^{(n)}: C\to M_r({\tilde B}_n)=M_r\otimes {\tilde B}_n$ by
$\phi_0^{(n)}(c)=\phi_0'(c)\otimes {\blue{1_{{\tilde B_n}}}}$ for all $c\in C.$ Put
$$
e_A=\{1_{A_n}\},\,\,  P=\{1_{M_r({\tilde B}_n)}\}+e_A.
$$
Put also  $Q_2=\pi(P)M_{r+1}({\tilde Q_1})\pi(P)$ and
define
${\blue{{\bar \phi}_j'}}={\blue{{\bar \phi}_j}}\oplus \pi\circ \{\phi_0^{(n)}\},$ $j=1,2.$
Then
\beq\label{Jukly31-1}
[{\bar \phi}_1']=[{\bar \phi}_2']\,\,\,{\rm in}\,\,\, KL(C, Q_2).
\eneq
The {{reason for adding}} $\pi\circ \{\phi_0^{(n)}\}$ is that, now, ${\bar \phi}_1'$ and ${\bar \phi}_2'$ are unital.
It follows from Theorem 4.3 of \cite{Da1} that there exist an integer $K>0,$ a unitary $u\in M_{1+K}(Q_2),$ and
a unital \hm\, $\psi: C\to M_{K}\subset M_{K}(Q_2)$
($M_{K}$ {{idenitified}} with the {{natural}} unital subalgebra of
$M_K(Q_2)$) such that
\begin{equation*}
{\rm Ad}\, u\circ {\rm diag}({\blue{{\bar \phi}_1'}}, \psi)\approx_{\ep_0/4} {\rm diag}({\blue{{\bar \phi}_2'}}, \psi) \,\,\,{\rm on}\,\,\, {\cal F}_0.
\end{equation*}
There exists a unitary $V=\{V_n\}\in M_{1+K}(PM_{r+1}({\tilde B})P)$ such that
$\pi(V)=u.$ It follows (on identifying $M_K$ with $M_K\otimes 1_{Q_2}$ {{as above}}) that for all sufficiently large $n,$
\begin{equation*}
{\rm Ad}\, V_n\circ {\rm diag}(L_1^{(n)}\oplus \phi_n^{(n)}, \psi)\approx_{\ep_0/3} {\rm diag}(L_2^{(n)}\oplus \phi_n^{(n)}, \psi)
\,\,\,{\rm on}\,\,\, {\cal F}_0.
\end{equation*}
{{For}} each integer $k\ge 1,$ {{write}} $e_{n,k,0}={\rm diag}(\overbrace{1_{A_n}, 1_{A_n},...,1_{A_n}}^k)\in A_n\otimes {\cal K}=B_n,$
\beq\label{July31-3}
e_{n,k}'&=&{\rm diag}(1_{A_n},\overbrace{e_{n,k,0},e_{n,k,0},...,e_{n,k,0}}^r)\in PM_{1+r}(B_n)P, \andeqn\\
e_{n,k}''&=& {\rm diag}(\overbrace{e_{n,k}',e_{n,k}',...,e_{n,k}'}^K)\in M_{K}(PM_{1+r}(B_n)P).
\eneq
It should be noted that
$e_{n,k}''$ commutes with $\psi$ and $e_{n,k}'$ commutes with $\phi_n^{(0)}.$ Put $e_{n,k}=e_{n,k}'\oplus e_{n,k}''$
in $M_{1+K}(PM_{1+r}(B_n)P).$ Then $\{e_{n,k}\}$ {{is}} an approximate identity for
$M_{1+K}(PM_{1+r}(B_n)P).$ Note that $V_n\in M_{1+K}(PM_{r+1}({\tilde B})P).$ It is easy to check that
\beq\label{Aug1st-4}
\lim_{k\to\infty}\|[V_n, e_{n,k}]\|=0.
\eneq
It follows that there exists a unitary $U_{n,k}\in e_{n,k}M_{1+K}(PM_{1+r}(B_n)P)e_{n,k}$  for each $n$ and $k$ such that
\beq\label{Aug1st-5}
\lim_{k\to\infty}\|e_{n,k}V_ne_{n,k}-U_{n,k}\|=0.
\eneq
For each $k,$ there is $N(k)=rk+K(rk+1)$ such that
\beq\label{Aug1st-6}
M_{N(k)}(A_n)=((e_{n,k}'-1_{A_n})\oplus e_{n,k}'')M_{1+K}(PM_{r+1}(B_n)P)((e_{n,k}'-1_{A_n})\oplus e_{n,k}'').
\eneq
Moreover, $e_{n,k}M_{1+K}(PM_{1+r}(B_n)P)e_{n,k}=M_{N(k)+1}(A_n).$
Define
$\Psi_n(c)=(e_{n,k}'-1_{A_n})\phi_0^{(n)}(c)(e_{n,k}'-1_{A_n})\oplus e_{n,k}''\psi(c)e_{n,k}$ for $c\in C.$ Then, for  large
$k$ and large $n,$
\beq\label{Aug1st-7}
{\rm Ad}\, U_n\circ {\rm diag}(L_1^{(n)}, \Psi_n)\approx_{\ep_0/2}{\rm diag}(L_2^{(n)}, \Psi_n)\,\,\,{\rm on}\,\,\, {\cal F}_0.
\eneq
This {{is in}} contradiction {\blue{{{with}} \eqref{12/31/24-1+}.}}
%
%
\end{proof}

\begin{thm}\label{UniqAtoM}
Let $A\in {\cal {\bar D}}_s$
and let $\Delta: A_+^{q,{\bf 1}}\setminus \{0\}\to (0,1)$ be an order preserving map.

For any $\ep>0$ and {\blue{any}} finite subset ${\cal F}\subset A,$ there {{exist}} $\dt>0,$ a finite subset
${\cal P}\subset \underline{K}(A),$ a finite subset ${\cal H}_1\subset A_+^{\bf 1}\setminus \{0\}$, and
a finite subset ${\cal H}_2\subset A_{s.a.}$ satisfying the following {{condition}}:

If $\phi_1, \phi_2: A\to M_n$ are  unital \hm s such that
\beq\label{UnAM-1}
&&[\phi_1]|_{\cal P}=[\phi_2]|_{\cal P},\\\label{UNAM-1+}
&&{\rm tr}\circ \phi_1(g)\ge \Delta(\hat{g})\tforal g\in {\cal H}_1,\tand\\
&&|{\rm tr}\circ \phi_1(h)-{\rm tr}\circ \phi_2(h)|<\dt\tforal h\in {\cal H}_2,
\eneq
where ${\rm tr}\in T(M_n),$
then there exists a unitary $u\in M_n$ such that
\beq\label{UnAM-2}
\|{\rm Ad}\, u\circ \phi_1(f)-\phi_2(f)\|<\ep\tforal f\in {\cal F}.
\eneq
\end{thm}

\begin{proof}
If $A$ has finite {{dimension}}, the lemma is known. So, in what follows, we will assume
that $A$ is infinite dimensional.

Define $\Delta_0: A_+^{q, {\bf 1}}\setminus \{0\}\to (0,1)$ by
$\Delta_0=(3/4)\Delta.$
Fix $\ep>0$ and a finite subset ${\cal F}\subset A.$
Let ${\cal P}\subset \underline{K}(A)$, ${\cal H}_0\subset A_+^{q, {\bf 1}}\setminus \{0\}$
(in place of ${\cal H}$), and $K\ge 1$ be {{the finite subsets and integer provided}} by  Lemma \ref{Lauct2} for $\ep/2$ (in place of $\ep$),
${\cal F},$ and $\Delta_0.$

Choose $\ep_0>0$  and a finite subset ${\cal G}\subset A$ such that $\ep_0<\ep$ and
\beq\label{UnAM-5}
[\Phi_1']|_{\cal P}=[\Phi_2']|_{\cal P}
\eneq
for any pair of unital \hm s from $A$  satisfying
\beq\label{UnAM-6}
\|\Phi_1'(g)-\Phi_2'(g)\|<\ep_0\rforal g\in {\cal G}.
\eneq
We may assume that ${\cal F}\subset {\cal G}$ and $\ep_0<\ep/2.$

Let $\af=3/4.$  Let $N\ge 1$, $\dt_1>0$ (in place of $\dt$), ${\cal H}_1\subset A_+^{\bf 1}\setminus \{0\}$, and ${\cal H}_2\subset A_{s.a.}$ {{be the constant and finite subsets provided}} by \ref{8-N-4}
for $\ep_0/2$ (in place $\ep$), ${\cal G}$ (in place of ${\cal F}$), ${\cal H}_0,$ $K,$ and $\Delta_0$ (in place of $\Delta$).
Choosing a larger ${\cal H}_1,$ since $A$ has infinite dimension, we may assume
that ${\cal H}_1$ contains at least $N$ mutually orthogonal non-zero positive elements.

Now suppose that $\phi_1, \phi_2$ are  unital \hm s satisfying the {{assumptions}}
for  the  ${\cal P},$ ${\cal H}_1,$ and ${\cal H}_2$ above.
The assumption (\ref{UNAM-1+}) implies that $n\ge N.$
{{Applying}} \ref{8-N-4}, we obtain a unitary $u_1\in M_n,$ mutually orthogonal non-zero projections
$e_0, e_1, e_2,...,e_K\in M_n$ with $\sum_{i=0}^K e_i=1_{M_n},$ $e_0\lesssim e_1,$
$e_1$ equivalent to $e_i,$ $i=1,2,...,K,$ unital \hm s $\Phi_1,\Phi_2: A\to e_0M_ne_0$, and
a unital \hm\, $\psi: A\to e_1M_ne_1$ such that
\beq\label{UnAM-7}
\|{\rm Ad}\, u_1\circ \phi_1(f)-(\Phi_1(f)\oplus \Psi(f))\|<\ep_0/2\rforal f\in {\cal G},\\\label{UnAM-7+}
\|\phi_2(f)-(\Phi_2(f)\oplus \Psi(f))\|<\ep_0/2\rforal f\in {\cal G}, \andeqn\\
{\rm tr}\circ \psi(g)\ge (3/4)\Delta{\blue{(\hat{g})}}/K\rforal g\in {\cal H}_0,
\eneq
where $\Psi(a)={\rm diag}(\overbrace{\psi(a), \psi(a),...,\psi(a)}^K)$ for all $a\in A$ and
${\rm tr}$ is the tracial state on $M_n.$

Since $[\phi_1]|_{\cal P}=[\phi_2]|_{\cal P},$ by the choice of $\ep_0$ and ${\cal G},$ we compute that
\beq\label{UNAM-8}
[\Phi_1]|_{\cal P}=[\Phi_2]|_{\cal P}.
\eneq
Moreover,
\beq\label{UnAM-9}
t\circ \psi(g)\ge (3/4)\Delta(\hat{g})\rforal g\in {\cal H}_0,
\eneq
{\blue{where}} $t$ is the tracial state of $e_1M_ne_1.$
 By \ref{Lauct2}, there is a unitary $u_2\in M_n$ such that
 \beq\label{UnAM-10}
 \|{\rm Ad}\, u_2\circ (\Phi_1\oplus \Psi)(f)-{\blue{(\Phi_1\oplus \Psi)}}(f))\|<\ep/2\tforal f\in {\cal F}.
 \eneq
 Put $U=u_2u_1.$ Then, by (\ref{UnAM-7}), (\ref{UnAM-7+}), and (\ref{UnAM-10}),
 \beq\label{UnAM-11}
 \|{\rm Ad}\, U\circ \phi_1(f)-\phi_2(f)\|<\ep_0/2+\ep/2+\ep_0/2<\ep\rforal f\in {\cal F}.
 \eneq
\end{proof}

\begin{lem}\label{Lfullab}
Let $A\in \bar{{\cal D}}_s$ be a unital \CA\, and let $\Delta: A_+^{q, \bf 1}\setminus \{0\}\to (0,1)$
{{be an order preserving map.}}
Let ${\cal P}_0\subset K_0(A)$ be a finite subset.
Then there {{exist}} an integer $N({\cal P}_0){{\ge 1}}$ and a finite subset ${\cal H}\subset A_+^{\bf 1}\setminus \{0\}$ satisfying the following {{condition}}:
For any unital \hm\, $\phi: A\to M_k$ (for some $k\ge 1$) and any unital \hm\,
$\psi: A\to M_R$ for some integer $R\ge N({\cal P}_0)k$
such that {\blue{(with ${\rm tr}$  the tracial state of $M_R$)}}
\beq\label{Lfullab-0}
{\rm tr}\circ \psi(g)\ge \Delta(\hat{g})\tforal g\in {\cal H},
\eneq
there exists a unital \hm\,
$h_0: A\to M_{R-k}$ such that
\beq\label{Lfullab-1}
[\phi\oplus h_0]|_{{\cal P}_0}=[\psi]|_{{\cal P}_0}.
\eneq
\end{lem}

\begin{proof}
{{Denote by}} $G_0$  the subgroup of $K_0(A)$ generated by ${\cal P}_0.$   We may assume, without loss of generality, that
${\cal P}_0=\{[p_1],[p_2],...,[p_{m_1}]\}\cup \{z_1,z_2,...,z_{m_2}\},$
where $p_1, p_2,...,p_{m_1}\in M_l(A)$ are  projections {\blue{(for some integer $l\ge 1$)}}
 and $z_j\in {\rm ker}\rho_A,$ $j=1,2,...,m_2.$

We prove the lemma  by induction.
Assume first that
$A=PC(X,F)P{,}$ {{w}}here $X$ is  {\blue{a compact metric space}}.
This, of course, includes the case that $X$ is a single point.
There is $d>0$ such that
\beq\label{Lfullab-5}
\|\pi_{x,j}\circ p_i-\pi_{x',j}\circ p_i\|<1/2,\,\,\, i=1,2,...,m_1,
\eneq
provided that ${\rm dist}(x,x')<d,$ where $\pi_{x,j}$ is identified with $\pi_{x,j}\otimes {\rm id}_{M_l}.$
Since $X$ is compact, we may assume that $\{x_1,x_2,...,x_{m_3}\}$ is a $d/2$-dense set.
Write $P_{x_i}FP_{x_i}=M_{r(i,1)}\oplus M_{r(i,2)}\oplus \cdots \oplus M_{r(i,k(x_i))},$ $i=1,2,...,m_3.$

There are $h_{i,j}\in C(X)$ with $0\le h_{i,j}\le 1,$ $h_{i,j}(x_i)=1_{M_{r(i,j)}}$, and
$h_{i,j}h_{i',j'}=0$ if $(i,j)\not=(i',j').$  Moreover we {{may}} assume
that $h_{i,j}(x)=0$ if ${\rm dist}(x, x_i)\ge d.$

Put $g_{i,j}=h_{i,j}\cdot P\in A,$ $j=1,2,...,k(x_i), i=1,2,...,m_3.$
Let
\beq\label{Lfullab-6}
\sigma_0=\min\{\Delta(\hat{h_{i,j}}): 1\le j\le k(x_i),\,1\le i\le m_3\}
\eneq
{{and choose an integer $N(\mathcal P_0) \geq 2/\sigma_0$}}. Put ${\cal H}=\{h_{i,j}: 1\le j\le k(x_i),\, 1\le i \le m_3\}.$

Now suppose that {{maps}} $\phi: A\to M_k$ and $\psi: A\to M_R$ {{are given}} with $R\ge N({\cal P}_0)k$ and
\beq\label{Lfullab-7}
{\rm tr}\circ \psi(g)\ge \Delta(\hat{g})\rforal g\in {\cal H}.
\eneq
Write $\phi=\bigoplus_{i,j}^{m_3}\Pi_{y_i,j},$ where $\Pi_{y_i,j(i)}$ is $T_{i,j}$ copies of $\pi_{y_i,j}.$
 Note {\blue{that}} $k-T_{\blue{i,j}}>0$ for all {\blue{$i,\,j.$}}
Since $R\ge N({\cal P}_0)k,$ (\ref{Lfullab-7}) implies
that $\psi$ {{is the}} direct sum of at least
\beq\label{Lfullab-8}
\Delta(\widehat{h_{j,i}})\cdot (2k/\sigma_0)>2k
\eneq
copies of $\pi_{x,j}$ with ${\rm dist}(x, x_i)<d,$ $i=1,2,...,m_3.$
Rewrite
$\psi$ {{as}} $\Sigma_1\oplus \Sigma_2,$ where
$\Sigma_1$ contains exactly $T_{i,j}$ copies of $\pi_{x,j}$ with ${\rm dist}(x, x_i)<d$ for each $i$ and $j.$
Then
\beq\label{Lfullab-9}
{\rm rank} \, \Sigma_1(p_i)={\rm rank} \phi(p_i),\,\,\,i=1,2,...,m_1.
\eneq
Put $h_0=\Sigma_2.$
Note for any unital \hm\, $h: A\to M_n,$ $[h(z)]=0$ for all $z\in {\rm ker}\rho_A.$
{{So $[\phi\oplus h_0]|_{{\cal P}_0}=[\psi]|_{{\cal P}_0}.$}}
This proves the case that $A=PC(X,F)P$ as above, in particular,  the case that $A\in  {\bar{{\cal D}}}_0.$

Now assume the conclusion of the lemma holds for any \CA\, $A\in \bar{{\cal D}}_m.$

Let $A$ be a \CA\, in ${\bar{{\cal D}}}_{m+1}.$
We assume that $A\subset PC(X,F)P\oplus B$ is a unital \SCA\,
and $I=\{f\in PC(X,F)P: f|_{X^0}=0\},$ where $X^0=X\setminus Y$ and $Y$ is an open subset of $X$
and $B\in {\cal D}_m'$ and $A/I\cong B.$
We assume that, {{if ${\rm dist}(x,x')<2d,$}} {{then}}
\beq\label{Lfullab-10}
\|\pi_{x, j}(p_i)-\pi_{x',j}(p_i)\|<1/2\andeqn
\|\pi_{x,j}\circ s\circ  \pi_I(p_i)-\pi_{x,j}(p_i)\|<1/2,
\eneq
where $s: A/I\to A^{d}=\{{\blue{(f|_{{\overline{X^d}}}, b): (f,b)}}\in A\}$ is {{the}} injective \hm\, given by \ref{8-N-3},
{\blue{and $X^d=\{x\in X: {\rm dist}(x, X^0)<d\}.$}}  We also assume
that $2d<d_{\blue{X, X^0}}.$
Define $\Delta_\pi: (A/I)_+^{q, {\bf 1}}\setminus \{0\}\to (0,1)$ by
\beq\label{Lfullab-11}
\Delta_\pi(\hat{g})=\Delta(\widehat{g_0\cdot P\cdot s(g)})\rforal g\in (A/I)_+^{\bf 1}\setminus \{0\},
\eneq
where $g_0\in C(X^d)_+$ with $0\le g\le 1,$ $g_0(x)=1$ if $x\in X^0,$
$g_0(x)>0$ if ${\rm dist}(x, X^0)<d/2,$ and
$g_0(x)=0$ if ${\rm dist}(x, X^0)\ge d/2\}.$
Note that $g_0\cdot P s(g)>0$ if $g\in (A/I)_+\setminus \{0\}.$ Therefore, $\Delta_\pi$ is indeed an order preserving map
from $(A/I)_+^{q, {\bf 1}}\setminus \{0\}$ into $(0,1).$

Note that $A/I\in {\cal D}_m$ {{(see the later part of Definition \ref{8-N-3}).}}
By the inductive assumption,
there is an integer $N_\pi({\cal P}_0)\ge 1,$ a finite subset ${\cal H}_\pi\subset (A/I)_+^{\bf 1}\setminus \{0\}$ satisfying the following {{condition}}:
if $\phi': A/I\to M_{k'}$ is a unital \hm\, and $\psi': A/I\to M_{R'}$ is a unital \hm\, for some $R'>N_\pi({\cal P}_0)k'$ such that
\begin{equation*}
t\circ \psi'(\hat{g})\ge \Delta_\pi(\hat{g})\rforal g\in {\cal H}_\pi,
\end{equation*}
where $t$ is the tracial state of $N_{R'},$ then there exists a unital \hm\, $h_\pi: A/I\to M_{R'-k'}$ such that
\begin{equation*}
(\phi'\oplus h_\pi)_{*0}|_{\bar{\cal P}_0}=(\psi_\pi)_{*0}|_{\bar{\cal P}_0},
\end{equation*}
where $\bar{\cal P}_0=\{(\pi_I)_{*0}(p): p\in {\cal P}_0\}.$

{\blue{For $\dt>0,$ define $Y^\dt=X\setminus X^\dt.$}}
Let $r: Y^{d/2}\to Y^d$ be a homeomorphism such that ${\rm dist}(r(x), x)<d$ for all $x\in Y^{d/2}$ {\blue{(see  Definition \ref{8-N-3}).}}
Set $C=\{f|_{Y^d}: f\in I\}.$
{\blue{Let $\pi_C: A\to C$ be defined by
$\pi_C(a)=\lambda(a)|_{Y^d}$ for all $a\in A$ (see \ref{8-N-3} for $\lambda: A\to PC(X, F)P$).}}
Define $\Delta_I: C_+^{q, {\bf 1}}\setminus \{0\}\to (0,1)$ by
\begin{equation}
\Delta_I(\hat{g})=\Delta(\widehat{(f_0P g)\circ r})\rforal g\in  C_+^{q,{\bf 1}}\setminus \{0\},
\end{equation}
where $f_0\in C_0(Y)_+$ with $0\le f_0\le 1,$ $f_0(x)=1$ if $x\in Y ^d,$ $f_0(x)=0$ if ${\rm dist}(x, X^0)\le d/2$
and $f_0(x)>0$ if ${\rm dist}(x, X^0)>d/2\}.$

Note that $C=P|_{Y^d} C({\blue{Y^d}}, F)P_{Y^d}.$  By what has been proved, there {{are}} an integer $N_I({\cal P}_0)\ge 1$ and a finite subset ${\cal H}_I\subset C_+^{\bf 1}\setminus \{0\}$ satisfying the following {{condition}}:
if $\phi'': C\to M_{k''}$ is a unital \hm\, and $\psi'': C\to M_{R''}$ (for some $R''\ge N_I({\cal P}_0)k''$) is another  unital \hm\, such that
{\blue{($t$ is the tracial state on $M_{R''}$)}}
\begin{equation*}
t\circ \psi''(\hat{g})\ge \Delta_I(\hat{g})\tforal g\in {\cal H}_I,
\end{equation*}
then there exists a unital \hm\, $h'': C\to M_{R''-k''}$ such that
\beq\label{Lfullab-14}
(\phi''\oplus h'')_{*0}|_{{\blue{[\pi_C]}}({\cal P}_0)}=(\psi''|_{*0})|_{{\blue{[\pi_C]}}({\cal P}_0)}.
\eneq
Put
\beq\label{Lfullab-13+}
\sigma=\min\{\min\{\Delta_\pi(\hat{g}): g\in {\cal H}_\pi\},\min\{\Delta_I(\hat{g}): g\in {\cal H}_\pi\}\}>0.
\eneq

Let $N=(N_\pi({\cal P}_0)+N_I({\cal P}_0))/\sigma$ and let
\beq\label{Lfullab-15}
{\cal H}=\{g_0\circ s(g): g\in {\cal H}_\pi\}\cup \{f_0\cdot g\circ r\}.
\eneq

Now suppose that $\phi: A\to M_k$ and $\psi: A\to M_R$ satisfy the {{assumptions}} for $N=N({\cal P}_0)$ and
${\cal H}$ as above {\blue{(and $R>Nk$).}}
We may write $\phi=\Sigma_{\phi, \pi}\oplus \Sigma_{\phi, I},$
where $\Sigma_{\phi, \pi}$
is
the  (finite) direct sum of irreducible representations
of $A$ which {{factor}} through $A/I$ and $\Sigma_{\phi, I}$ is
the (finite) direct sum
of irreducible representations of $I.$
We {{may}} also write
\beq\label{Lfullab-16}
\psi=\Sigma_{\psi, \pi}\oplus \Sigma_{\psi, b}\oplus \Sigma_{\psi, I'},
\eneq
where $\Sigma_{\psi, \pi}$ {{is}} the direct sum of irreducible representations
of $A$ which {{factor}} through $A/I,$ $\Sigma_{\psi, b}$
is the (finite) direct sum
of irreducible representations  which factor through point evaluations at $x\in Y$ with ${\rm dist}(x, X^0)<d/2$
and $\Sigma_{\psi, I'}$ {{is the direct sum}}  of irreducible representations which
{{factor through}} point evaluations at $x\in Y$ with ${\rm dist}(x, X^0)\ge d/2$.

Put $q_\pi=(\Sigma_{\psi, \pi}\oplus \Sigma_{\psi,b})(1_A)$ and
$k'={\rm rank} \Sigma_{\phi,\pi}(1_A).$
Define $\psi_\pi: A/I\to M_{{\rm rank}(q_\pi)}$ by
$\psi_\pi(a)=(\Sigma_{\psi,\pi}\oplus \Sigma_{\psi,b})\circ {\blue{s(a)}}\rforal a\in A/I.$
Then, {\blue{by \eqref{Lfullab-0}
and the choice of ${\cal H},$}}
\beq\label{Lfullab-17-}
t_1\circ \psi_\pi(g)\ge {\rm tr}\circ {\blue{\psi( g_0 P \cdot s(g))}}\ge \Delta(\widehat{g_0 P\cdot s(g)})=\Delta_\pi(\hat{g})\tforal g\in {\cal H}_\pi,
\eneq
where $t_1$ is the tracial state on $M_{{\rm rank}(q_\pi)}.$
Note that
\beq\label{Lfullab-17}
{\rm tr}\circ \psi(g_0{\blue{ P}})\ge \Delta(\widehat{\blue{g_0 P}})=\Delta_\pi(\widehat{1_{A/I}}),
\eneq
Therefore,
\beq\label{Lfullab-18}
{\rm rank} (q_\pi)\ge R\Delta_\pi(1_{A/I})\ge N_\pi({\cal P}_0)k'.
\eneq
By the inductive assumption, there is a unital \hm\, $h_\pi: A/I\to M_{{\rm rank}q_\pi-k'}$ such that
\beq\label{Lfullab-19}
(\Sigma_{\phi,\pi}\oplus h_\pi)_{*0}|_{\blue{\bar{{\cal P}}_0}}=(\psi_\pi)_{*0}|_{\bar{{\cal P}}_0}.
\eneq

Put $q_I=\Sigma_{\psi, I'}(1_A)$ and $k''={\rm rank}(\Sigma_{\phi, I}(1_A)).$
Define $\psi_I: C\to M_{{\rm rank}q_I}$ by
\beq\label{Lfullab-20}
\psi_I(a)=\Sigma_{\psi,I'}(a\circ r)\rforal a\in C.
\eneq
Then
\beq\label{Lfullab-21}
&&\hspace{-0.8in}t_2\circ \psi_I(g)= {\rm tr}\circ \Sigma_{\psi, I'}({\blue{g}}\circ r)\ge \psi(f_0{\blue{P g}}\circ r)\ge
\Delta(\widehat{(f_0P g)\circ r})=\Delta_I(\hat{g})\rforal g\in {\cal H}_I,
\eneq
where $t_2$ is the tracial state on $M_{{\rm rank}(q_I)}.$
Note that
\beq\label{Lfullab-22}
{\rm tr}\circ \psi({\blue{f_0 P}})\ge \Delta(\widehat{\blue{f_0 P}})=\Delta_I(\hat{1_A}).
\eneq
Therefore,
\beq\label{Lfullab-23}
{\rm rank}(q_I)\ge R\Delta_I(1_A)\ge N_I({\cal P}_0)k''.
\eneq
There are $0<d_1<d<d_{\blue{X,X^0}}$ such that all irreducible representations appearing in $\Sigma_{\phi,I}$
factor through point evaluations at $x$ with ${\rm dist}(x, X^d)\ge d_1.$
{{Choose a \hm}}\, $r': Y^{d_1}\to Y^d$ {{as  in \ref{8-N-3}.}}
Define $\phi_I: C\to M_{k''}$ by $\phi_I(f)=\Sigma_{\phi, I}(f\circ r').$

By  what has been proved, the choice of  $N_I({\cal P}_0),$ and by \eqref{Lfullab-23},
there is a unital \hm\, $h_I: C\to M_{{\rm rank}(q_I)-k''}$ such that
\beq\label{Lfullab-24}
(\Sigma_{\phi, I}\oplus h_I)_{*0}|_{{\blue{[\pi_C]({\cal P}_0)}}}=(\psi_I)_{*0}|_{[\pi_C]({\cal P}_0)}.
\eneq
Define $h: A\to M_{R-k}$ by $h(a)=h_\pi(\pi_I(a))\oplus h_I(a|_{Y^d})\oplus \Sigma_{\psi,b}(a) $ for all $a\in A.$
Then, for each $i,$  {\blue{by  \eqref{Lfullab-19}  and \eqref{Lfullab-24}}},
\beq\label{Lfullab-25}
\hspace{-0.7in}{\rm rank}\phi(p_i)+{\rm rank}h(p_i)&=&{\rm rank}(\Sigma_{\phi, \pi}(p_i))+
{\rm rank}(\Sigma_{\phi, I}(p_i))\\
&& { +} {\rm rank}(\Sigma_{\psi, b}(p_i))+{\rm rank}h_\pi(p_i)+{\rm rank}h_I(p_i)\\
&=& {\rm rank}\psi_\pi(p_i)+{\rm rank}(\Sigma_{\psi, b}(p_i))+{\rm rank}\psi_I(p_i)\\
&=&{\rm rank} \psi(p_i),\,\,\, i=1,2,...,m_1.
\eneq
Since $(\phi)_{*0}(z_j)=h_{*0}(z_j)=\psi_{*0}(z_j){\blue{=0}},$ $j=1,2,...,m_2,$ we conclude that
\beq\label{Lfullab-26}
(\phi\oplus h)_{*0}|_{{\cal P}_0}=\psi_{*0}|_{{\cal P}_0}.
\eneq
This completes the induction process.
\end{proof}

\begin{lem}\label{fullabs}
Let $A\in \bar{\cal D}_s$ be a unital \CA\, and let
$\Delta: A^{q,{\bf 1}}_+\setminus \{0\}\to (0,1)$ be
an order preserving map.
For any $\ep>0$ and any finite subset ${\cal F}\subset A$ there exist a finite
subset ${\cal H}\subset A_+^{\bf 1}\setminus \{0\}$ and an integer
$L\ge 1$ satisfying the following {{condition}}:
For any unital \hm\, $\phi: A\to M_{k}$ and
any unital \hm\, $\psi: A\to M_{R}$ for some $R\ge Lk$ such that
\beq\label{fullabs-1}
{\rm tr}\circ \psi(h)\ge \Delta(\hat{h})\tforal h\in {\cal H}\,\,\,{\blue{({\rm tr}\in T(M_R))}},
\eneq
there exist a unital \hm\, $\phi_0: A\to M_{R-k}$ and a unitary $u\in M_{R}$ such that
\beq\label{fullabs-2}
\|{\rm Ad}\, u\circ {\rm diag}(\phi(f), \phi_0(f))-
\psi(f)\|<\ep
\eneq
for all $f\in {\cal F}.$

\end{lem}

\begin{proof}
Let $\dt>0,$ ${\cal P}\subset \underline{K}(A)$ be a finite subset, ${\cal H}_1\subset A_+^{\bf 1}\setminus \{0\}$
be a finite
subset, ${\cal H}_2\subset A_{s.a.}$ be a finite subset,  and $N_0$ be an integer as
provided by Theorem \ref{UniqAtoM} for $\ep/4$ (in place of $\ep$),
${\cal F},$ $(1/2)\Delta,$ and $A.$
Without loss of generality, we may assume that ${\cal H}_2\subset A_+^{\bf 1}\setminus \{0\}.$
Let {$\sigma_0=\min\{\Delta(\hat{g}): g\in {\cal H}_1 \cup {\cal H}_2\}$.}

Let $G$ denote the subgroup of $\underline{K}(A)$ generated by ${\cal P}.$ Put
${\cal P}_0={\cal P}\cap K_0(A).$  We may also assume, without loss of generality, that
${\cal P}_0=\{[p_1],[p_2],...,[p_{m_1}]\}\cup \{z_1,z_2,...,z_{m_2}\},$
where $p_1, p_2,...,p_{m_1}$ are projections in $M_l(A)$ {\blue{(for some integer $l$)}}
 and $z_j\in {\rm ker}\rho_A,$ $j=1,2,...,m_2.$
Let $j\ge 1$ be an integer such that $K_0(A, \Z/j'\Z)\cap G=\emptyset$ for all $j'\ge j.$
Put $J=j!.$

Let $N({\cal P}_0)\ge 1$ {{denote the}} integer and ${\cal H}_3\subset A_+^{\bf 1}\setminus \{0\}$ {{the}}  finite subset {{provided}} by Theorem \ref{Lfullab} for ${\cal P}_0.$
Let $p_s=(a^{(s)}_{i,j})_{l\times l},$ $s=1,2,...,m_1$, and
choose $\ep_0>0$ and a finite subset ${\cal F}_0$ such that
\beq\label{fullabs-9}
[\psi']|_{\cal P}=[\psi'']|_{\cal P}{{,}}
\eneq
{{whenever}} $\|\psi'(a)-\psi''(a)\|<\ep_0$ for all $a\in {\cal F}_0.$

Put ${\cal F}_2={\cal F}\cup {\cal F}_1\cup {\cal H}_2$ and put $\ep_1=\min\{\ep/16, {\blue{\ep_0/2,\dt/2}}\}.$
Let $K>8((N({\cal P}_0)+1)(J+1)/\dt{\blue{\sigma_0}})$ be an integer.
Let ${\cal H}_0={\cal H}_1\cup {\cal H}_3.$
Let $N_1\ge 1$ (in place of $N$), $\dt_1>0$ (in place of $\dt$), ${\cal H}_4\subset
A_+^{\bf 1}\setminus \{0\}$ (in place of ${\cal H}_1$), and ${\cal H}_5\subset A_{s.a.}$
(in place of ${\cal H}_2$) {{be the constants and finite subsets provided}} by Lemma \ref{8-N-4} for $\ep_1$ (in place of $\ep$),
${\cal F}_2$ (in place of ${\cal F}$),  ${\cal H}_0,$ and $K.$
Let $L=K(K+1),$  let ${\cal H}={\cal H}_4\cup {\cal H}_0$ and let  $\af=15/16.$
Suppose that $\phi$ and $\psi$ {{satisfy}}  the assumption \eqref{fullabs-1} for the above $L$ and ${\cal H}.$

Then, by Lemma \ref{8-N-4}, there are mutually orthogonal projections
$e_0, e_1, e_2,...,e_K\in M_R$ such that $e_0\lesssim e_1$  and $e_i$ is equivalent to $e_1,$ $i=1,2,...,K,$
a unital \hm\, $\psi_0: A\to e_0M_Re_0,$ and a unital \hm\, $\psi_1: A\to e_1M_Re_1$ such that
\beq\label{fullabs-10}
\|\psi(a)-(\psi_0(a)\oplus {\rm diag}(\overbrace{\psi_1(a),\psi_1(a),...,\psi_1(a)}^{K}))\|<\ep_1\tforal a\in {\cal F}_2\\\label{fullabs-10+0}
\andeqn
{\rm tr}\circ \psi_1(g)\ge (15/16){\Delta(\hat{g})\over{K}}\rforal g\in {\cal H}_0.
\eneq
Put $\Psi=\psi_0\oplus {\rm diag}(\overbrace{\psi_1(a),\psi_1(a),...,\psi_1(a)}^{K})).$
We compute that
\beq\label{fullabs-10+}
[\Psi]|_{\cal P}=[\psi]|_{\cal P}.
\eneq
Let $R_0={\rm rank}(e_1).$ Then, {\blue{by \eqref{fullabs-10+0}}},
\beq\label{fullabs-11}
R_0=R{\rm tr}\circ \psi_1(1_A)&\ge& Lk(15/16) {\Delta(\widehat{1_A})\over{K}}\ge k(K+1)(15/16)\Delta(\widehat{1_A})\\
&\ge & k(15/16)8N({\cal P}_0)(J+1)/\dt.
\eneq
Moreover,
\beq\label{fullabs-12}
{\rm tr}'\circ \psi_1(\hat{g})\ge (15/16)\Delta(\hat{g})\tforal g\in {\cal H}_3,
\eneq
where ${\rm  tr}'$ is the tracial state of $M_{R_1}.$
It follows from Lemma \ref{Lfullab} that there exists a unital \hm\, $h_0: A\to M_{R_0-k}$ such that
\beq\label{fullabs-13}
(\phi\oplus h_0)_{*0}|_{{\cal P}_0}=(\psi_1)_{{\cal P}_0}.
\eneq
Put
\beq\label{fullabs-14}
h_1=h_0\oplus {\rm diag}(\overbrace{\phi\oplus h_0,\phi\oplus h_0,...,\phi\oplus h_0}^{J-1})\andeqn
\psi_2={\rm diag}(\overbrace{\psi_1, \psi_1,...,\psi_1}^J).
\eneq
Then
\beq\label{fullabs-15}
[\phi\oplus h_1]|_{\cal P}=[\psi_2]|_{\cal P}.
\eneq
Put $\Psi'={\rm diag}(\overbrace{\psi_1,\psi_1,...,\psi_1}^{K-J})=\psi_I\otimes 1_{K-J}.$
Let $\phi_0=h_1\oplus \psi_0\oplus \Psi'.$ Then {\blue{(see \eqref{fullabs-10})}}
\beq\label{fullabs-16}
[\phi\oplus \phi_0]|_{\cal P}=[\psi_0\oplus \psi_2\oplus \Psi']|_{\cal P}=[\psi]|_{\cal P}.
\eneq
Since $J/K<{\blue{\dt/4}},$  by \eqref{fullabs-10},
\beq\label{fullabs-16+1}
|{\rm tr}\circ (\phi(g)\oplus \phi_0(g))-{\rm tr}\circ \psi(g)|<2J/K+\ep_1<\dt  \rforal g\in {\cal H}_2.
\eneq
Then,  by (\ref{fullabs-1}), applying  Theorem \ref{UniqAtoM}, {{we obtain}} a unitary $u\in M_R$ such that
\beq\label{fullabs-17}
\|{\rm Ad}\, u\circ (\phi(f)\oplus \phi_0(f))-\psi(f)\|<\ep\tforal f\in {\cal F}.
\eneq
\end{proof}




\section{Almost multiplicative maps to finite dimensional \CA s}





In the following statement,  $n$ is given and $(\mathcal G, \delta)$ depends on $n:$ {{the proof
is a standard  compactness argument.}}

\begin{lem}\label{LtoMN}
Let $n\ge 1$ be an integer and let $A$ be a unital separable  \CA.
For any $\ep>0$ and any finite subset ${\cal F}\subset A,$ there exist
$\dt>0$ and a finite subset ${\cal G}\subset A$ such that, for any unital ${\cal G}$-$\dt$-multiplicative {{\cp}} $\phi: A\to M_n,$ there exists a
unital \hm\, $\psi: A\to M_n$ such that
$$
\|\phi(a)-\psi(a)\|<\ep\tforal a\in {\cal F}.
$$
\end{lem}

\begin{proof}

Suppose that {{the conclusion}} is not true for {{a}}  certain finite set ${\cal F} \sbs A$ and $\ep_0{{>0}}$. Let $\{{\cal G}_k\}_{k=1}^{\infty}$ be a sequence of finite subsets of $A$ with ${\cal G}_k\subset  {\cal G}_{k+1}$ and $\overline{\bigcup_k {\cal G}_k}=A$
and let $\{\dt_k\}$ be a decreasing sequence of positive numbers with $\dt_k \rightarrow 0$.
Since the {{conclusion}} is assumed not to be true, there are unital ${\cal G}_k$-$\dt_k$-multiplicative
{{\cp s}} $\phi_k: A\to M_n$ such that
\beq\label{150103-sec5-1}
 \mbox{inf}\{\mbox{max}_{a\in {\cal F}}\|\phi_k(a)-\psi(a)\|:~~~\psi: A\to M_n ~\mbox{{{a}} homomorphisms} \} \geq \ep_0.
\eneq
For each pair $(i,j)$ with $1\leq ~i,j~\leq n$, let $l^{i,j}: M_n \to \C$ be the map defined by taking the matrix $a\in M_n$ to the entry of {{the}} $i^{th}$ row and $j^{th}$ column of $a$. Let $\phi_k^{i,j} =l^{i,j}\circ \phi_k: A\to \C$. Note that the unit ball of the dual space of $A$ (as {{a}} Banach space) is {{weak$^*$ compact}}.  Since $A$ is separable, there is a subsequence (instead of subnet) of $\{\phi_k\}$ (still denoted by $\phi_k$) such that $\{\phi_k^{i,j}\}$ is weak$^*$ convergent for all $i,j$. In other words,  $\{\phi_k\}$ {{converges}}  pointwise. Let $\psi_0$ be the the limit. Then $\psi_0$ is a homomorphism and for $k$ large enough, we have
$$\|\phi_k(a)-\psi_0(a)\|< \ep_0, ~~\rforal a\in {{\cal F}}.$$
This is {{in}} contradiction {{with}} (\ref{150103-sec5-1}) above.
\end{proof}

\begin{lem}[cf. Lemma 4.5 of \cite{Lncrell}]\label{0dig}
Let $A$ be a unital \CA\, arising from a locally trivial continuous field of \CA s isomorphic to $M_n$ over a compact metric space $X.$ Let $T$ be a finite subset of tracial states on $A.$ For any finite subset ${\cal F}\subset A$ and for any $\ep>0$ and any $\sigma>0,$ there {{are}} an ideal $J\subset A$ such that
$\|\tau|_J\|<\sigma$ for all $\tau\in T,$ a finite dimensional \SCA\, $C\subset A/J$, and a unital \hm\, $\pi_0$
{{from $A/J$}}
such that
\beq\label{0dig-1}
{\rm dist}(\pi(x), C)<\ep\tforal x\in {\cal F} \tand
\pi_0(A/J)=\pi_0(C)\cong C,
\eneq
where $\pi: A\to A/J$ is the quotient map.
\end{lem}
\begin{proof}
This  follows from
Lemma 4.5 of \cite{Lncrell}.  In fact, the only difference is the existence of $\pi_0.$
{\blue{We will keep the notation {{of}} the proof of Lemma 4.5 of \cite{Lncrell}.
Let $B_j$ and $F_j$ be as in the proof of Lemma 4.5 of \cite{Lncrell}.
Note that in the proof Lemma 4.5 of \cite{Lncrell}, $\phi(\{g_{ij}\})\cong D_\xi\cong M_n.$
In other words,  {{in the proof of Lemma 4.5 of \cite{Lncrell},}}  $(B_j)|_{F_j}\cong M_n.$}}
{\blue{Recall that
$$
J=\{f\in A: f(\zeta)=0\rforal \zeta\in F\},\, F=\bigsqcup_{i=1}^k F_i, \andeqn F_j\subset B(\xi_j, \dt_j),\,\,j=1,2,...,k,
$$
 as in the proof of 4.5 of \cite{Lncrell}.
Recall also  $D=\bigoplus_{i=1}^k f_j B_j$ and $\pi(D)=C,$
where $f_j|_{F_j}=1$ and $f_j|_{F_i}=0,$ if $j\not=i$.
Choose  $x_j\in F_j,$ }}
%
$j=1,2,...,k.$ One  defines
${{\pi_0'}}=\bigoplus_{j=1}^k \pi_{x_j}.$ Then
$${{\pi'_0(D) = \bigoplus_{j=1}^l\pi_{x_j}(B_j)\cong \bigoplus_{j=1}^l M_n{\blue{\cong}} C.}}$$
Note that $\xi_j\in F_j\subset F,$ $j=1,2,...,k.$ Therefore, each $\pi_{\xi_j}$ induces
a \hm\, $\psi_j$ of $A/J$ such that $\psi_j(\pi(a))=\pi_{\xi_j}(a)$ for all $a\in A.$
{\blue{In particular, $\psi_j(A/J)\cong M_n,$ $j=1,2,...,k.$
Put $\pi_0=\bigoplus_{j=1}^k \psi_j.$
 Then $\pi_0(A/J)=\pi_0'(A)=\bigoplus_{j=1}^k \pi_{\xi_j}(A)=\bigoplus_{j=1}^k\pi_{\xi_j}(D)=\pi_0'(D)=
 \pi_0(\pi(D))=\pi_0(C)=\pi_0'(D)\cong C.$}}
\end{proof}

\begin{lem}[cf. Lemma 4.7 of \cite{Lncrell}]\label{abdig}
Let $A$ be a unital  separable subhomogeneous \CA. 
Let $T\subset T(A)$ be a finite subset. For any finite subset ${\cal F}\subset A,$ $\ep>0$ and $\sigma>0,$ there are an ideal $J\subset A$ such that
$\|\tau|_J\|<\sigma$ for all $\tau\in T,$ a finite dimensional \SCA\, $C\subset A/J$, and  a unital \hm\, $\pi_0$
{{from $A/J$}}
such that
\beq\label{abdig-1}
{\rm dist}(\pi(x), C)<\ep\tforal x\in {\cal F} {{\andeqn}}
\pi_0(A/J)=\pi_0(C)\cong C,
\eneq
where $\pi: A\to A/J$ is the quotient map.
\end{lem}

\begin{proof}
The proof is in fact contained in that of Lemma 4.7 of \cite{Lncrell}. Each time Lemma 4.5 of \cite{Lncrell} {{is}} applied, one can apply
 Lemma \ref{0dig} above.
 {\blue{Let us assume $A$  and $m$ to be as in Lemma 4.7 of \cite{Lncrell} and keep the notation {{of}}
  the proof of Lemma 4.7 of \cite{Lncrell}.
 Let us point out what the map $\pi_0$ is.  When $m=1,$ $\pi_0$ is the {{map of}}  Lemma  \ref{0dig}.}}

 {\blue{Otherwise we proceed {{as in}} the proof of Lemma 4.7 of \cite{Lncrell} to where $\pi_1$ and $\pi_2$ are constructed.
 Instead of applying Lemma 4.5 of \cite{Lncrell}, we apply Lemma \ref{0dig} {{above}} to obtain  a \hm\,
 $\pi_0^{(1)}$ (in place of $\pi_0$) from $\pi_2\circ \pi_1(A)$ and a finite
 dimensional \CA\, $C_1=\pi_2\circ \pi_1(D_1)=\pi_2\circ \pi_1(D_1')$ such that $\pi_0^{(1)}(\pi_2\circ \pi_1(A))=\pi_0^{(1)}(C_1)\cong C_1.$

 In the proof of Lemma 4.7 of \cite{Lncrell}, the second time Lemma 4.5 of \cite{Lncrell} is applied, we again,
 instead, {{apply}} Lemma \ref{0dig} to obtain a \hm\, $\pi_0^{(2)}$ (in place of $\pi_0$) from $\pi_4\circ \pi_3(A)$
 and $C_2=\pi_4\circ \pi_3(D_2)=\pi_4\circ \pi_3(D_2')$ such that
 $\pi_0^{(2)}(\pi_4\circ \pi_3(A))=\pi_0^{(2)}(C_2)\cong C_2.$
 Note, in the proof of Lemma 4.7 of \cite{Lncrell},  $\pi(A)=\pi_2\circ \pi_1(A)\oplus \pi_4\circ \pi_3(A)$
 and
 $C=\pi(D_1'\oplus D_2')=\pi(D_1')\oplus \pi(D_2')
 =\pi_2\circ \pi_1(D_1')\oplus \pi_4\circ \pi_3(D_2')=C_1\oplus C_2.$
 Define $\pi_0=\pi_0^{(1)}\oplus \pi_0^{(2)}$ from $\pi(A).$ Then
 $\pi_0(\pi(A))=\pi_0^{(1)}(\pi_2\circ \pi_1(A))\oplus \pi_0^{(2)}(\pi_4\circ \pi_3(A))
 =\pi_0^{(1)}(C_1)\oplus \pi_0^{(2)}(C_2)=\pi_0(\pi(C))=\pi_0^{(1)}(C_1)\oplus \pi_0^{(2)}(C_2)\cong C_1\oplus C_2=C.$
 }}
\end{proof}


\begin{lem}\label{LtoMn}
Let $A$ be a unital {{separable}} subhomogeneous \CA. 
Let $\ep>0,$ let ${\cal F}\subset A$ be a finite subset, and let $\sigma_0>0.$
There exist
 $\dt>0$ and a finite subset ${\cal G}\subset A$ satisfying the following {{condition}}:
Suppose that $\phi: A\to M_n$ (for some integer $n\ge 1$) is a {{unital}}
${\cal G}$-$\dt$-multiplicative {{\cp.}}
Then,
there {{exist}} a projection $p\in M_n$  and a unital \hm\,
$\phi_0: A\to pM_np$ such that
\beq\label{LtoMn-2}
\|p\phi(a)-\phi(a)p\|&<&\ep\tforal a\in {\cal F},\\
\|\phi(a)-[(1-p)\phi(a)(1-p)+\phi_0(a)]\|&<&\ep\tforal  a\in {\cal F}, \tand\\
{\rm tr}(1-p)&<&\sigma_0,
\eneq
where ${\rm tr}$ is the normalized trace on $M_n.$
\end{lem}

\begin{proof}
{{Assume}} that the {{conclusion}} is false.
Then there {exist} $\ep_0>0,$ a finite subset ${\cal F}_0,$ a positive number $\sigma_0>0,$
an increasing sequence
of finite subsets ${\cal G}_n\subset A$ such that ${\cal G}_n\subset  {\cal G}_{n+1}$ and such that
$\bigcup_{n=1}{\cal G}_n$ is dense in $A,$ a {{decreasing}} sequence of positive
numbers $\{\dt_n\}$ with $\sum_{n=1}^{\infty}\dt_n<\infty,$
a sequence of integers $\{m(n)\},$ and a sequence
of unital ${\cal G}_n$-$\dt_n$-multiplicative
{{\cp s}}
$\phi_n: A\to M_{m(n)}$ satisfying the following condition:
\beq\label{LtoMn-6}
\inf \{\max\{\|\phi_n(a)-[(1-p)\phi_n(a)(1-p)+\phi_0(a)]\|: a\in {\cal F}_0\}\}\ge \ep_0,
\eneq
where the infimum is taken among all projections $p\in M_{m(n)}$
with ${\blue{\rm{tr}}}_n(1-p)<\sigma_0,$
 and $${{\|p\phi_n(a) - \phi_n(a) p\| < \ep_0,}}$$
where ${\blue{\rm{tr}}}_n$ is the normalized
trace on $M_{m(n)}$, and all possible  unital \hm s $\phi_0: A\to pM_{m(n)}p.$
By virtue of \ref{LtoMN},
one may also assume that $m(n)\to \infty$ as $n\to\infty.$

Note  that  $\{{\blue{\rm{tr}}}_n\circ \phi_n\}$ is a sequence of (not {{necessarily}} tracial) states of $A.$ Let
$t_0$ be a weak{ *} limit of $\{{\blue{{\rm{tr}}}}_n\circ \phi_n\}.$ Since $A$ is separable, there is a subsequence (instead of a subnet) of
 $\{{\blue{\rm{tr}}}_n\circ \phi_n\}$ converging to $t_0$.
Passing to a subsequence,
we may assume that
${\blue{\rm{tr}}}_n\circ \phi_n$ converges to $t_0.$ By the {{${\cal G}_n$-$\dt_n$}}-multiplicativity of $\phi_n$,
{{the limit}}
$t_0$ is a tracial state on $A$.

Consider the ideal $\bigoplus_{n=1}^{\infty}M_{m(n)},$ where
$$
\bigoplus_{n=1}^{\infty}M_{m(n)}=\{\{a_n\}: a_n\in M_{m(n)}\andeqn \lim_{n\to \infty}\|a_n\|=0\}.
$$
Denote by $Q$ the quotient
$\prod_{n=1}^{\infty}M_{m(n)}/\bigoplus_{n=1}^{\infty}M_{m(n)}.$
Let
$\pi_\omega: \prod_{n=1}^{\infty}M_{m(n)}\to  Q$ be the quotient map.
Let $A_0=\{\pi_{\omega}(\{\phi_n(f)\}): f\in A\}$ which is a  subalgebra of
$Q.$
Denote by $\Psi$ the canonical
unital \hm~ from $A$ to  $Q$
with $\Psi(A)=A_0$. If $a\in A$ has  zero image in $\pi_\omega(A_0)$, i.e., $\phi_n(a) \to 0$, then $t_0(a)=\lim_{n\to \infty} {\blue{\rm{tr}}}_n(\phi_n(a))=0$.
  So we  may view $t_0$ as a state on $A_0=\Psi(A).$


It follows from Lemma \ref{abdig} that there {{are}} an
ideal $I\subset \Psi(A)$ and a finite dimensional \SCA\, $B\subset \Psi(A)/I$
and a unital \hm\, $\pi_{00}:\Psi(A)/I \to B$ such that
\beq\label{LtoMn-10}
{\rm dist}(\pi_I\circ \Psi(f), B)&<&\ep_0/16\tforal f\in {\cal F}_0,\\\label{LtoMn-11}
\|(t_0)|_I\|&<&\sigma_0/2\andeqn
\pi_{00}|_B=\mbox{id}_B.
\eneq
Note that $\pi_{00}$ can be regarded as map from ${{\Psi(A)}}$ to $B$ {{with}} ${\rm ker} \pi_{00}\supset I$.
There is, for each $f\in {\cal F}_0,$ an element $b_f\in B$ such that
\beq\label{LtoMn-12}
\|\pi_I\circ \Psi(f)-b_f\|<\ep_0/16.
\eneq
Put $C'=B+I$ and $I_0=\Psi^{-1}(I)$ and $C_1=\Psi^{-1}(C').$ For each $f\in {\cal F}_0,$ there exists
$a_f\in C_1\subset A$ such that
\beq\label{LtoMn13}
\|f-a_f\|<\ep_0/16\andeqn \pi_I\circ \Psi(a_f)=b_f.
\eneq
Let $a\in (I_0)_+$ be a strictly positive element and let
$J=\overline{\Psi(a)Q\Psi(a)}$ {{denote}} the hereditary \SCA\, of $Q$ generated by $\Psi(a).$
{{Since $Q$ has real rank zero, so does $J$ (see   \cite{BP}).}}
Put $C_2=\Psi(C_1)+J.$ Then $J$ is a ($\sigma$-unital)  ideal of $C_2.$
Denote by $\pi_J: C_2\to B$ the quotient map.
Since $Q$
has real rank zero, {{$J$ is a hereditary \SCA\, of $Q,$}}  and
$C_2/J=B$ has finite {{dimension}}, by Lemma 5.2  of \cite{LnTAF}, $C_2$ has real rank zero
{{and projections in $B$ lifts to a projection in $C_2.$}} It follows {{(see Theorem 9.8 of \cite{Eff-dim})}} that the extension
$$
0\to J\to C_2\to B\to 0
$$
splits and is a quasidiagonal {{(see also the proof of Theorem 5.3 of \cite{LnTAF})}}.  As in Lemma 4.9 of \cite{Lncrell}, there are a projection
$P\in J$ and a unital \hm\, $\psi_0: B\to (1-P)C_2(1-P)$ such that
\beq\label{LtoMn-15}
\|P\Psi(a_f)-\Psi(a_f)P\|<\ep_0/8\andeqn \|\Psi(a_f)-[P\Psi(a_f)P+\psi_0\circ \pi_J\circ\Psi(a_f)]\|<\ep_0/8
\eneq
for all $f\in {\cal F}_0.$ Let $H: A\to \psi_0(B)$ be defined by
$H=\psi_0\circ \pi_{00}\circ \pi_I\circ \Psi$.
One estimates that
\beq\label{LtoMn-16}
\|P\Psi(f)-\Psi(f)P\|<\ep_0/2\andeqn\\\label{LtoMn-17}
\|\Psi(f)-[P\Psi(f)P+H(f)]\|<\ep_0/2
\eneq
for all $f\in {\cal F}_0.$ Note that ${\rm dim}H(A)<\infty,$ and that $H(A)\subset Q.$
There is a
\hm\, $H_1: H(A)\to \prod_{n=1}^{\infty}M_{m(n)}$ such that
$\pi_{{\omega}}\circ H_1\circ H=H.$ One may write $H_1=\{h_n\},$ where
each $h_n: H(A)\to M_{m(n)}$ is a (not {{necessarily}} unital) \hm, $n=1,2,....$
There is also a sequence of projections $q_n\in M_{m(n)}$ such that
$\pi_{{\omega}}(\{q_n\})=P.$ Let $p_n=1-q_n,$ $n=1,2,....$ Then, for sufficiently large $n,$
by (\ref{LtoMn-16}) and (\ref{LtoMn-17}),
\beq\label{LtoMn-18}
&&{{\|(1-p_n)\phi_n(f)-\phi_n(f)(1-p_n)\|}}<\ep_0\andeqn\\\label{LtoMn-18+}
&&\|\phi_n(f)-[(1-p_n)\phi_n(f)(1-p_n)+h_n\circ H(f)]\|<\ep_0
\eneq
for all $f\in {\cal F}_0.$ Moreover, since $P\in J,$ for any $\eta>0,$ there is $b\in I_0$ with
$0\le b\le 1$ such that
\begin{equation*}
\|\Psi(b)P-P\|<\eta.
\end{equation*}
However, by (\ref{LtoMn-11}),
\beq\label{LtoMn-20}
0\,{\blue{\le}}\,t_0(\Psi(b))<\sigma_0/2\tforal b\in I_0\,\,\, {\rm with}\,\,\, 0\le b\le 1.
\eneq
By choosing sufficiently small $\eta,$
for all sufficiently large $n,$ {{we have}}
\begin{equation*}
{\blue{\rm{tr}}}_n(1-p_n)<\sigma_0.
\end{equation*}
This, {{together with \eqref{LtoMn-18}}} and \eqref{LtoMn-18+}, contradicts (\ref{LtoMn-6}).
\end{proof}



%
%

\begin{cor}\label{LtoMncor}
Let $A$ be a unital subhomogeneous \CA.  
Let $\eta>0,$ let ${\cal E}\subset A$ be a finite subset, and let $\eta_0>0.$
There exist
 $\dt>0$ and a finite subset ${\cal G}\subset A$ satisfying the following condition:
Suppose that $\phi,~ \psi: A\to M_n$ (for some integer $n\ge 1$) are two {{unital}}
${\cal G}$-$\dt$-multiplicative {{\cp s.}}
Then,
there exist  projections $p, q\in M_n$ with ${\rm rank}(p)={\rm rank}(q)$ and  unital \hm s\,
$\phi_0: A\to pM_np$ and $\psi_0: A\to qM_nq$ such that
$$
\|p\phi(a)-\phi(a)p\|<\eta,\quad \|q\psi(a)-\psi(a)q\|<\eta,\quad a\in {\cal E},$$
$$\|\phi(a)-[(1-p)\phi(a)(1-p)+\phi_0(a)]\|<\eta,\quad \|\psi(a)-[(1-q)\psi(a)(1-q)+\psi_0(a)]\|<\eta,\quad  a\in {\cal E}$$
$$\tand {\rm tr}(1-p)={\rm tr}(1-q)<\eta_0,$$
where ${\rm tr}$ is the normalized trace on $M_n.$
\end{cor}

For convenience  {{in}} future use, we  {{have}} used $\eta$, $\eta_0$ and ${\cal E}$ to replace  {{the}} $\ep$, $\sg_0,$ and ${\cal F}$  {{of}} \ref{LtoMn}.

\begin{proof}

By Lemma \ref{LtoMn}, we can {{obtain}} such decompositions for $\phi$ and $\psi$ separately.  {{Thus}} the only missing part  is that ${\rm rank}(p)={\rm rank}(q)$.
Let $\{z_1,z_2,...,z_m\}$ be the set of ranks of irreducible representations of $A$ and let
$T$ be the number given by \ref{8-N-0} corresponding to  $\{z_1,z_2,...,z_m\}.$
We apply \ref{LtoMn} to $\eta_0/2$ instead of $\sg_0$ (and, $\eta$ and ${\cal E}$ in places of  $\ep$ and ${\cal F}$). By
Lemma \ref{LtoMN}, we can assume the size $n$ of the matrix algebra $M_n$ is  large enough that $T/n<\eta_0/2$. By Lemma \ref{8-N-0}, we may take sub-representations out of $\phi_0$ and $\psi_0$ (one of them has size at most $T$) so that the remainder of $\phi_0$ and $\psi_0$ have same size---that is for ${\rm rank}({\rm new}~ p)={\rm rank}({\rm new}~ q)$, and $\mathrm{tr}(1-({\rm new}~p))=\mathrm{tr}(1-({\rm new}~q))<\eta_0/2+T/n<\eta_0$.
\end{proof}

\begin{lem}\label{Combinerep}
Let $A\in \overline{ {\cal D}}_s$ be an infinite dimensional unital \CA, let $\ep>0$   and let
${\cal F}\subset A$ be  a finite subset.
 {{Let}} $\ep_0>0$ and let ${\cal G}_0\subset A$ be a finite subset.
Let $\Delta: A_+^{q, {\bf 1}}\setminus\{0\}\to (0,1)$ be
{{an order preserving}} map.

Suppose that
${\cal H}_1\subset A_+^{\bf 1}\setminus \{0\}$ is a finite subset,
$\ep_1>0$ is a
positive number and
$K\ge 1$ is an integer.
There  {{exist}} $\dt>0,$  $\sg>0$, and  a finite subset
${\cal G}\subset A$ and a finite subset ${\cal H}_2\subset A_+^{\bf 1}\setminus\{0\}$ satisfying the following  {{condition}}:
Suppose that $L_1, L_2: A\to M_n$ (for some integer $n\ge 1$) are unital ${\cal G}$-$\dt$-multiplicative
{{\cp s}} such that
\beq\label{Comb-1}
tr\circ L_1(h)\ge \Delta({\hat h})
~~\tforal h\in {\cal H}_2,\ and
\eneq
\beq\label{Comb-2}
|tr\circ L_1(h)-tr\circ L_2(h)|<\sg \rforal ~h\in {\cal H}_2.
\eneq
Then there exist mutually orthogonal  projections $e_0, e_1,e_2,...,e_K\in
M_n$ such that
$e_1, e_2,...,e_K$ are pairwise equivalent, $e_0\lesssim e_1,$
$tr(e_0)<\ep_1$, and $e_0+\sum_{i=1}^K e_i=1,$ and there exist  unital
${\cal G}_0$-$\ep_0$-multiplicative {{\cp s}}
 $\psi_1,\psi_2: A\to e_0{{M_n}}e_0$,  a unital \hm\, $\psi: A\to e_1{{M_n}}e_1$, and  {{a}} unitary $u\in M_n$ such that
\beq\label{Com-3}
\|L_1(f)-{\rm diag}(\psi_1(f), \overbrace{\psi(f),\psi(f),...,\psi(f)}^K)\|<\ep \tand
\eneq
\beq\label{Com-4}
\|{\blue{u^*}}L_2(f){\blue{u}}-{\rm diag}(\psi_2(f), \overbrace{\psi(f),\psi(f),...,\psi(f)}^K)\|<\ep
\eneq
for all $f\in {\cal F},$ where ${{\mathrm{tr}}}$ is the tracial state of $M_n.$
Moreover,
\beq\label{Com-5}
{{\mathrm{tr}}}(\psi(g))\ge {\Delta({\hat g})\over{{\blue{2}}K}}\tforal g\in {\cal H}_1.
\eneq
\end{lem}

\begin{proof}
{{Let $\ep>0,$ $\ep_0>0,$ ${\cal F}$ and ${\cal G}_0$  and $\Delta$ be given
as stated.}}
First note that the following statement is evident. For any $C^*$-algebra $A$, {{any}} finite subset ${\cal G}_0\subset A$  and  $\ep_0>0$, there are a
{{finite subset}}  ${\cal F}'\subset A$ which contains ${\cal F}$ and $\ep'>0$
{{which is smaller than $\min\{\ep/2, \ep_0/2\}$}}
 satisfying the following condition.  If $L: A\to B$ is a unital
{\blue{${\cal F}'$-$\ep'$}}-multiplicative {{\cp}}\,,  $p_0,p_1\in B$ are projections with $p_0+p_1= 1_B$, and $L'_0: A\to p_0Bp_0$, $L'_1: A\to p_1Bp_1$ are unital \cp s with
$$\|L(f)-{\rm diag}(L'_0(f),L'_1(f))\|< \ep'\rforal f\in {\cal F}',$$
then  both $L'_0$ and $L'_1$ are {\blue{${\cal G}_0$-$\ep_0$}}-multiplicative.

{{Let $\ep_1>0,$ ${\cal H}_1$ and $K$ be given as in the statement of the lemma.}}
{\blue{Choose an integer $k_0\ge 1$ such that $1/k_0<\ep_1$
and let $K_1=k_0K.$}}
Put
\beq\label{8-Nsec4-1-1}
{\blue{\ep_2}}=\min\{\ep/16, \ep'/16,{\blue{\ep_1/2, 1/2}}\}.
\eneq
Let $\Delta_1=(3/4)\Delta.$  Let $\dt_1>0$ (in place of $\dt$), ${\cal H}_{1,0}\subset A_+^{\bf 1}\setminus \{0\}$ (in place of ${\cal H}_1$), 
and ${\cal H}_{2,0}\subset A_{s.a.}$ (in place of ${\cal H}_2$) be  {{the constant and finite subsets provided by
Lemma}} \ref{8-N-4} (see also Remark \ref{8-N-4-r}) for ${\blue{\ep_2}}$ (in place of $\ep$),  ${\cal F}'$
{{(in place of ${\cal F}$),}}
$2{\blue{K_1}}$ (in place of $K$), {{${\cal H}_1$ (in place of ${\cal H}_0$),}} $\Delta_1,$ and $A,$ as well as $\af=3/4.$  {\blue{ It is clear that, \wilog, we may assume that ${\cal H}_{2,0}$ is in the self-adjoint part of the unit ball of $A.$}}
{\blue{Since every $h\in {\cal H}_{2,0}$ can be written as $h=(|h|+h)/2+(h-|h|)/2,$
choosing  an even smaller $\dt_1$ (half  the size), \wilog,}} we may assume that
${\cal H}_{{2,0}}\subset A_+^{\bf 1}\setminus \{0\}.$

Let $\eta_0=\min\{\dt_1/16, {\blue{\ep_2}}/16, \min\{{{\Delta(\hat{h})}}: h\in {\blue{{\cal H}_{1,0}/16}}\}\}.$
{{Let}} $\dt_2>0$ (in place of $\dt$) and let  ${\cal G}_1\subset A$ (in place of ${\cal G}$) be  {{the constant and}} the finite subset  {{provided}} by \ref{LtoMncor} for
$\eta=\eta_0\cdot \min\{{\blue{\ep_2}},\dt_1/4\},$  $\eta_0,$ and  ${\cal E}={\blue{{\cal F}'\cup}} {\cal H}_{1,0}\cup {\cal H}_{2,0}\cup {\cal H}_1.$
Let  $\dt=\eta_0\cdot \min\{\dt_2/2, \dt_1/2, {\blue{\ep_0/2}}\},$ let $\sigma=\min\{\eta_0/2, {{\eta}}/2\},$
let ${\cal G}={\blue{{\cal G}_0}}\cup {\cal G}_1\cup {\cal F}\cup {\cal F}'\cup {\cal E},$ and
let ${\cal H}_2={\cal H}_{1,0}\cup {\cal H}_{2,0}\cup {\cal H}_1.$
{{Recall that we have assumed that ${\cal F}'\supset {\cal F}$ and $\ep'<\min\{\ep/2, \ep_0/2\}.$}}

Now suppose that $L_1$ and $L_2$ satisfy the  {{assumptions}} of the lemma with respect to the $\dt,$ $\sigma$ and
${\cal G}$, ${\cal H}_2$  {{above.}}
It follows from Corollary \ref{LtoMncor} that there exist a projection $p\in M_n,$
two unital \hm s $\phi_1, \phi_2: A\to pM_np,$
and a unitary $u_1\in M_n$ such that
\beq\label{8-Nsec4-1-2}
&&\|{\blue{u_1^*L_1(a)u_1}}-((1-p)u_1^*L_1(a)u_1(1-p)+\phi_1(a))\|<\eta,\\
&&\|L_2(a)-((1-p)L_2(a)(1-p)+\phi_2(a)\|<\eta
\eneq
for all $a\in {\cal E},$ and
\beq\label{8-Nsec4-1-3}
{\blue{{\rm tr}}}(1-p)<\eta_0,
\eneq
where ${\blue{{\rm tr}}}$ is the tracial state on $M_n.$

We compute that
\beq\label{8-Nsec4-1-4}
&&{\blue{{\rm tr}}}\circ \phi_1(g)\ge \Delta(\hat{g})-\eta-\eta_0\ge (3/4)\Delta(\hat{g}) \tforal g\in {\cal H}_{1,0}\andeqn\\
&&|{\blue{{\rm tr}}}\circ \phi_1(g)-\tau\circ \phi_2(g)|<2\eta+2\eta_0+{\blue{\sigma}}<\dt_1\tforal g\in {{{\cal H}_{2,0}.}}
\eneq

It follows from Lemma \ref{8-N-4} (and  {{Remark}} \ref{8-N-4-r}) that  there {{exist}}
mutually orthogonal projections $q_0, q_1,...,q_{2K}\in pM_np$ such that
$q_0\lesssim q_1$ and $q_i$ is equivalent to $q_1$ for all $i=1,2,...,2{\blue{K_1}},$ two unital \hm s
$\phi_{1,0}, \phi_{2,0}: A\to q_0M_nq_0,$ a unital \hm\, $\psi': A\to {{q}}_1M_n{{q_1}}$, and a unitary
$u_2\in pM_np$ such that
\beq\label{8-Nsec4-1-5}
\|{\blue{u_2^* \phi_1(a)u_2}}-(\phi_{1,0}{\blue{(a)}}\oplus {\rm diag}(\overbrace{\psi'(a),\psi'(a), ...,\psi'(a)}^{\blue{2K_1}}))\|<{\blue{\ep_2}}\\
\andeqn\|\phi_2(a)-(\phi_{2,0}{\blue{(a)}}\oplus {\rm diag}(\overbrace{\psi'(a),\psi'(a), ...,\psi'(a)}^{{\blue{2K_1}}}))\|<
{\blue{\ep_2}}
\eneq
 for all $a\in {\cal F}\cup {\cal F}'.$ Moreover,
 \beq\label{8-Nsec4-1-6}
 \tau\circ{\blue{\psi'}}(g)\ge (3/4)^2{\Delta(\hat{g})/{2{\blue{K_1}}}}\rforal g\in {\cal H}_1,
 \eneq
 {\blue{where $\tau$ is the tracial state of $pM_np.$}}
Let $u={\blue{u_1}}((1-p)+u_2),$ $e_0=(1-p)\oplus q_0,$\linebreak
$e_i={\blue{\sum_{j=1}^{2k_0}q_{2k_0(i-1)+j}}},$
$i=1,2,...,K,$ let $\psi_1(a)=(1-p)u_1^*L_1(a)u_1(1-p)\oplus {{\phi}}_{1,0},$
$\psi_2(a)=(1-p)L_2(a)(1-p)\oplus {{\phi}}_{2,0}$ and $\psi(a)={\rm diag}({\blue{\overbrace{\psi'(a),\psi'(a),...,\psi'(a)}^{2k_0}}})$ for $a\in A.$
{{Then}}
\beq\label{8-Nsec4-1-7}
\|u^*L_1(f)u-(\psi_1(f)\oplus {\rm diag}(\overbrace{\psi(f), \psi(f),...,\psi(f)}^K))\|<{\blue{\eta+\ep_2<\ep'<}}\ep\\
\andeqn \|L_2(f)-(\psi_2(f)\oplus {\rm diag}(\overbrace{\psi(f), \psi(f),...,\psi(f)}^K))\|<{\blue{\eta+\ep_2<\ep'<}}\ep
\eneq
for all $f\in {\cal F}{\blue{\cup{\cal F}'}}.$
{\blue{It follows  (by the choice of $\ep'$ and ${\cal F}'$)  that}} $\psi_1$ and $\psi_{{2}}$ are ${\cal G}_0$-$\ep_0$-multiplicative.
{\blue{Moreover,
$$
{\rm tr}(e_0)={\rm tr}(1-p)+{\rm tr}(q_0)<\eta_0+1/(2k_0K+1)<
\ep_1.$$
}}
Further,
\beq\label{8-Nsec4-1-8}
{\blue{\rm{tr}}}\circ \psi(g)&=&{\blue{\rm{tr}(p)\tau\circ \psi(g)\ge (1-\eta_0)\tau\circ \psi(g)}}\\
&&{\blue{\ge (1-{1\over{32}})(9/16)(2k_0){\Delta(\hat{g})\over{2k_0K}}\ge {\Delta(\hat{g})\over{2K}}}}\rforal g\in {\cal H}_1.
\eneq
\end{proof}

\begin{lem}[{{9.4 of \cite{LinTAI}}}]\label{measureexistence}
Let $A$ be a unital separable \CA. For any $\ep>0$ and  any  finite subset ${\cal H}\subset A_{s.a.},$ there  {{exist}} a finite subset ${\cal G}\subset A$ and $\dt>0$  satisfying the following  {{condition}}:
Suppose that $\phi: A\to B$ {{(for some unital  \CA\, $B$)}} is a  {{unital}} ${\cal G}$-$\dt$-multiplicative
\cp\, and $t\in T(B)$ is a tracial state of $B.$ Then, there exists
a tracial state $\tau\in T(A)$ such that
\beq\label{mext-1}
|t\circ\phi(h)-\tau(h)|<\ep\tforal h\in {\cal H}.
\eneq
\end{lem}

\begin{proof}
{\color{Green} This follows from the same proof of Lemma 9.4 of \cite{LinTAI}.}
\end{proof}

\begin{thm}\label{Newunique1}
Let $A\in {\overline{\cal D}}_s$ be a unital \CA.
Let $ \Delta: A_+^{ q, {\bf 1}}\setminus \{0\}\to (0,1)$ be an order preserving map.

Let  $\ep>0$ and let ${\cal F}\subset A$ be a finite subset. There exist a finite subset
${\cal H}_1\subset A_+^{{ {\bf 1}}}\setminus \{0\},$
a finite subset ${\cal G}\subset A,$
$\dt>0,$ a finite subset ${\cal P}\subset \underline{K}(A),$
a finite subset ${\cal H}_2\subset A_{s.a.},$ and $\sigma>0$ satisfying the following  {{condition}}:
Suppose that $L_1, L_2: A\to M_k$ (for some integer $k\ge 1$) are two unital
${\cal G}$-$\dt$-multiplicative
\cp s such that
\beq\label{Newu-1}
[L_1]|_{\cal P}&=&[L_2]|_{\cal P},\\\label{Newu-1+}
{\rm tr}\circ L_1(h)&\ge& \Delta({\hat{h}})\,\,\,
 \tforal h\in {\cal H}_1, \tand\\
|{\rm tr}\circ L_1(h)-\mathrm{tr}\circ L_2(h)|&<&\sigma\tforal h\in {\cal H}_2.
\eneq
 {{Then}} there exists a unitary $u\in M_k$ such that
\beq\label{Newu-4}
\|{\rm Ad}\, u\circ L_1(f)-L_2(f){\|}<\ep\tforal f\in {\cal F}.
\eneq
\end{thm}

\begin{proof}
The proof is {{almost}}  the same as that of \ref{UniqAtoM}.

{\blue{If $A$ is finite dimensional, then $A$ is semiprojective. Therefore, it is easy to see
that the general case can be reduced to the case that both $L_1$ and $L_2$ are unital \hm s.
Then we return to the situation  {{of}} \ref{UniqAtoM}.

So we now assume that $A$ is infinite dimensional.

Let $\Delta_1=(1/3)\Delta.$ Fix  {{$\ep>0$}} and a finite subset ${\cal F}\subset A.$
Let $\dt_1>0$ (in place of $\dt$), ${\cal G}_1\subset A$ (in place of ${\cal G}$),
 ${\cal P}\subset \underline{K}(A)$,  ${\cal H}_1'\subset A_+^{\bf 1}\setminus \{0\}$
 {{(in place of ${\cal H}$)}} be the finite subsets
and $K\ge 1$  {{be the integer provided}} by Lemma \ref{Lauct2}  (see also Remark \ref{Re414}) for $\ep/2$ (in place of $\ep$),
${\cal F},$ $\Delta_1,$ and $A.$

Let ${\cal G}_0={\cal G}_1$ and $\ep_0=\dt_1/2.$
\Wlog, we may assume the following:
for any two  ${\cal G}_0$-$\ep_0$-multiplicative \morp s $\Phi_1, \Phi_2: A\to C$ (for any unital \CA\, $C$),
$[\Phi_1]|_{\cal P}$ and $[\Phi_2]|_{\cal P}$ are well defined, and,
if $\|\Phi_1(a)-\Phi_2(a)\|<\ep_0$ for all $a\in {\cal G}_0,$ then
\beq
[\Phi_1]|_{\cal P}=[\Phi_2]|_{\cal P}.
\eneq

Let $\dt_2>0$ (in place  {{of}} $\dt$), $\sigma_1>0$ (in place of $\sigma$), a finite subset ${\cal G}_2\subset A$
(in place of ${\cal G}$)  and  a finite subset ${\cal H}_2'\subset A_+^{\bf 1}\setminus \{0\}$ (in place of  ${\cal H}_2$) be as given
by Lemma \ref{Combinerep} for ${\cal H}_1'$ (in place of ${\cal H}_1$), $K,$  $\ep_0,$ ${\cal G}_0$  $\dt_1/4$ (in place of $\ep$),
${\cal G}_1={\cal G}_0$ (in place of ${\cal F}$), $\Delta,$ and $A.$

Let $\dt=\min\{\dt_1/4, \dt_2, \ep/4\},$ ${\cal G}={\cal G}_0\cup {\cal F}\cup {\cal G}_2,$
${\cal H}_1={\cal H}_2'$ and ${\cal H}_2={\cal H}_2'.$

Now let $L_1, \, L_2: A\to M_k$ be two unital ${\cal G}$-$\dt$-multiplicative {{\cp s}} which satisfy
the assumption for the above ${\cal H}_1,$ ${\cal G},$ $\dt,$ $\sigma,$ ${\cal P},$ and ${\cal H}_2.$

By  \ref{Combinerep}, we obtain a unitary $u_1\in M_k,$ mutually orthogonal non-zero projections
$e_0, e_1, e_2,...,e_K\in M_k$ with $\sum_{i=0}^K e_i=1_{M_k},$ $e_0\lesssim e_1,$
$e_i$   equivalent to $e_1,$ $i=1,2,...,K,$ unital ${\cal G}_0$-$\ep_0$-multiplicative
\cp s $\Phi_1,\Phi_2: A\to e_0M_ke_0,$ and
a unital \hm\, $\psi: A\to e_1M_ke_1$ such that
\beq\label{nUnAM-7}
\|u_1^*\circ L_1(f)u_1-(\Phi_1(f)\oplus \Psi(f))\|<\dt_1/4\rforal f\in {\cal G}_1,\\\label{nUnAM-7+}
\|L_2(f)-(\Phi_2(f)\oplus \Psi(f))\|<\dt_1/4\rforal f\in {\cal G}_1,\andeqn\\\label{nUnAM-7++}
\tau\circ \psi(g)\ge \Delta(\hat{g})/2K\rforal g\in {\cal H}_1',
\eneq
where $\Psi(a)={\rm diag}(\overbrace{\psi(a), \psi(a),...,\psi(a)}^K)$ for all $a\in A$ and
$\tau$ is the tracial state on $M_n.$

Let $\tau_1$ be the tracial state of $e_1M_ke_1.$  Then \eqref{nUnAM-7++} implies
\beq\label{18831-n1}
\tau_1\circ \Psi(g)\ge \Delta(\hat{g})/2\ge \Delta_1(\hat{g})\rforal g\in {\cal H}_1.
\eneq

By the choice of $\ep_0$ and ${\cal G}_0,$  {{by \eqref{nUnAM-7} and \eqref{nUnAM-7+},}}
one has, for all $x\in {\cal P},$
\beq\label{18831-n1}
[\Phi_1](x)+[\Psi](x)=[\Phi_2](x)+[\Psi](x)
\eneq
in  {{the}} group $\underline{K}(A).$ It follows that, for all $x\in {\cal P},$
\beq\label{18831-n2}
[\Phi_1](x)=[\Phi_2](x).
\eneq
Now,  by \eqref{18831-n2}, \eqref{18831-n1}, and the choice  {{of}} $K,$    on applying  Lemma \ref{Lauct2},
{{one obtains}} a unitary $u_2\in M_k$ such
that
\beq\label{18831-n3}
\|u_2^*(\Phi_1(a)\oplus \Psi(a))u_2-(\Phi_2(a)\oplus {{\Psi}}(a))\|<\ep/2\tforal a\in {\cal F}.
\eneq
Choose $u=u_1u_2.$ Then,  by \eqref{nUnAM-7} and \eqref{18831-n3}, for all $a\in {\cal F},$
\beq\nonumber
\hspace{-0.2in}\|u^*L_1(a)u- L_2(a)\| &\le& \|u_2^*(u_1^*L_1(a){{u_2}}-(\Phi_1(a)\oplus {{\Psi}}(a))\|\\\nonumber
&&\hspace{-0.3in}+\|u_2^*(\Phi_1(a)\oplus \Psi(a))u_2-(\Phi_1(a)\oplus {{\Psi}}(a))\|<\ep/2+\ep/2=\ep.
\eneq
}}
\end{proof}

%







\section{Homotopy Lemma in finite dimensional \CA s}

{{
\begin{lem}\label{homhom}
Let $S$ be a subset of $M_k$ (for some integer $k\ge 1$), and let $u\in M_k$ {\blue{be}}  a unitary
such that
\beq\label{homhom-1}
ua=au\tforal a\in S.
\eneq
Then there exists a continuous path of unitaries
$\{u_t: t\in [0,1]\}\subset M_k$ such that
\beq\label{homhom-2}
u_0=u,\,\,\, u_1=1,\,\,\,
u_ta=au_t\tforal a\in {\blue{S}}
\eneq
and for all $t\in [0,1],$  {{and}} moreover,
\beq\label{homhom-3}
{\rm length}(\{u_t\})\le \pi.
\eneq
\end{lem}
}}

\begin{proof}
{{There is a continuous function}} $h$  from $sp(u)$ to
$[-\pi, \pi]$ such that
\beq\label{homhom-4}
\exp(i h(u))=u.
\eneq
{{Because $h$ is a continuous function of $u,$}}
\beq\label{homhom-5}
{\blue{a}}h(u)=h(u){\blue{a}}\tforal a\in {\blue{S}}.
\eneq
Note that $h(u)\in (M_k)_{s.a.}$ and $\|h(u)\|\le \pi.$
Define $u_t=\exp(i(1-t)h(u))$ ($t\in [0,1]$). Then
$
u_0=u\andeqn u_1=1.
$
Also,
$$
u_t{\blue{a}}={\blue{a}}u_t
$$
for all $a\in {\blue{S}}$ and $t\in [0,1].$
Moreover, one has
$$
{\rm  {{length}}}(\{u_t\})\le \pi,
$$
as desired.
\end{proof}

%

\begin{lem}\label{EXTMM}
Let $A\in {\bar{\cal D}}_s$ be a unital \CA, let ${\cal H}\subset (A\otimes C(\T))_{s.a.}$ be a finite subset,
  let $1> \sigma>0$ be a positive number, and let $\Delta: {{(A\otimes C(\T))}}_+^{q,\bf 1}\setminus \{0\}\to (0,1)$ be an o{{r}}der preserving map.  Let  $\ep>0,$  {{let}} ${\cal G}_0\subset A\otimes C(\T)$ be a finite subset, let ${\cal P}_0, {\cal P}_1\subset \underline{K}(A)$ be finite subsets, and  {{write}}
  ${\cal P}={\cal P}_0\cup \bt({\cal P}_1)\subset \underline{K}(A\otimes C(\T)).$
There exist $\dt>0,$ a finite subset ${\cal G}\subset A\otimes C(\T),$
and a finite subset ${\cal H}_1\subset (A\otimes C(\T))_+^{\bf 1}\setminus \{0\}$ satisfying the following condition:
Suppose that $L: A\otimes C(\T)\to M_k$ (for some integer ${{k}}\ge 1$) is a {{unital}}
${\cal G}$-$\dt$-multiplicative  \cp\,
such that
\beq\label{EXTMM-1}
{\rm tr}\circ L(h) &\ge& \Delta(\hat{h})\tforal h\in {\cal H}_1,
\tand\\\label{EXTMM-1+}
{[L]}|_{{\boldsymbol{\bt}}({\cal  P}_1)}&=&0,
\eneq
{{where ${\mathrm{tr}}$ is the tracial state of $M_n.$}}
Then there exists a unital ${\cal G}_0$-${{\ep}}$-multiplicative \cp\,
$\psi: A\otimes C(\T)\to M_k$ such that
$u=\psi(1\otimes z)$ is a unitary,
\beq\label{EXTMM-2}
u\psi(a\otimes 1)&=&\psi(a\otimes 1)u\tforal a\in A\\
{[L]}|_{\cal P}&=&[\psi]|_{\cal P}\,\,\text{and},\\\label{EXTMM-3}
|{\rm tr}\circ L(h)-{\rm tr}\circ \psi(h)|&<&\sigma\tforal h\in {\cal H}.
\eneq

\end{lem}

\begin{proof}
Let  ${\cal H}$ and ${{\sigma}},$ $\ep$ and ${\cal G}_0$ be given.
{\blue{It is clear that, \wilog,  we may assume that ${{{\cal H}\subset (A\otimes C(\T))_{s.a.}^{\bf 1}}}.$}}
{\blue{ By writing $h=h_+-h_-,$ where $h_+=(|h|+h)/2$ and $h_-=(|h|-h)/2,$ and {{choose}} a smaller
$\sigma,$ to simplify notation,
\wilog, we may assume that ${\cal H}\subset (A\otimes C(\T))_+^{\bf 1}\setminus \{0\}.$}}
We may also assume that
$$
{\cal G}_0=\{g\otimes f: g\in {\cal G}_{0A}\andeqn f\in {\cal G}_{1T}\},
$$
where ${{1_A}}\subset {\cal G}_{0A}\subset A$ and ${\cal G}_{1T}\subset C(\T)$ are finite
subsets. To simplify  {{matters}} further, we may assume, without loss of generality, that
${\cal G}_{1T}=\{1_{C(\T)}, z\},$ where $z\in C(\T)$ is the standard unitary generator.

{We may assume that $\mathcal G_{0A}$ is sufficiently large and $\ep$ is sufficiently small
that for any unital $\mathcal G_0$-$\ep$-multiplicative \cp\,
$L$, $[L]|_{\mathcal P}$ is well defined,
{\blue{any unital ${\cal G}_{0A}$-$\ep$-multiplicative  \cp\,
$\Phi$ from $A,$
$[\Phi]|_{{\cal P}_0}$ is well defined,}} and for any unital $\mathcal G_{0}$-$\ep$-multiplicative
\cp s
$L_1$ and $L_2$ with
$$
L_1\approx_\ep L_2\quad \textrm{on $\mathcal G_0$},
$$
we have $$[L_1]|_{\mathcal P} = [L_2]|_{\mathcal P},$$}
{\blue{and,  {{furthermore}}, for any
unital ${\cal G}_{0A}$-$\ep$-multiplicative \cp s
$\Phi_1$ and $\Phi_2$ (from $A$)
with
$$\Phi_1\approx_{\ep} \Phi_2\,\,\,{\rm on}\,\,\, {\cal G}_{0A},$$
we have
$$
[\Phi_1]|_{{\cal P}_0} =[\Phi_2]|_{{\cal P}_0}.
$$
}}

{{There is $\dt_0>0$ satisfying the following condition:
if $L', L'': A\otimes C(\T)\to C$ (for any unital \CA\, $C$) are unital ${\cal G}_0$-$\dt_0$-multiplicative \cp s
such that
$\|L'(a)-L''(a)\|<\dt_0$ for all $a\in {\cal G}_0,$ then
$\|L'(h)-L''(h)\|<\sigma/4\rforal h\in {\cal H}.$}}

Let $n$ be an integer such that $1/n<\sigma/2.$
Note that $A\otimes C(\T)\in {\overline {\cal D}}_s.$
{\blue{Put $\ep_0=\min\{\sigma/2, \ep/2, \dt_0/2\}.$}}

Let $\dt>0,$ {{$\sigma_0>0$ (in place of $\sigma$)}} ${\cal G}\subset A\otimes C(\T),$ and ${\cal H}_1\subset
A\otimes C(\T)^{\bf 1}_+\setminus \{0\}$ (in place of ${\cal H}_2$) be the constants  and finite subsets
provided by
{{Lemma}} \ref{Combinerep} for $A\otimes C(\T)$ (in place of $A$), {\blue{$\ep_0$}} (in place of $\ep$), ${\cal G}_0$ (in place of ${\cal F}$),
${\cal H}$ (in place of ${\cal H}_1$), and $\Delta.$
{\blue{\Wlog,  we may assume that $\dt<\ep,$ and by choosing larger a ${\cal G}$ if necessary, we may assume that ${\cal G}_0\subset {\cal G}.$ }}
Now suppose that $L: A\otimes C(\T)\to M_k$ satisfies the assumption
for the above $\dt,$ ${\cal G},$ and ${\cal H}_1.$  {{In particular, $[L]|_{\cal P}$ is well defined.}} It follows from {{Lemma}} \ref{Combinerep} (for $L_1=L_2=L$) that there is a projection $e_0\in M_k$ and a unital ${\cal G}_0$-${\blue{\ep_0}}$-multiplicative 
\cp\, $\psi_0: A\otimes C(\T)\to e_0M_ke_0$ and a unital \hm\, $\psi_1: A\otimes C(\T)\to (1-e_0)M_k(1-e_0)$
such that
\beq\label{Exmm-4}
&&{\rm tr}(e_0) <1/n<\sigma,\\\label{Exmm-5}
&&\|L(a)-\psi_0(a)\oplus \psi_1(a)\|<\ep_0\tforal a\in {\cal G}_0.
\eneq
 Define $\psi: A\otimes C(\T)\to M_k$ by
 $\psi(a)=\psi_0(a)\oplus \psi_1(a)$ for all $a\in A$ {{(see \ref{Ddiag})}} and
 $\psi(1\otimes z)=e_0\oplus \psi_1(1\otimes z).$
Put $u=\psi(1\otimes z).$
{\blue{Consider $\Phi_1(a)=L(a\otimes 1)$ and $\Phi_2(a)=\psi(a\otimes 1).$
Then, by \eqref{Exmm-5} and the choices of ${\cal G}_0$ and $\ep,$
\beq\label{Exmm-10}
[L]|_{{\cal P}_0}=[\psi]|_{{\cal P}_0}.
\eneq
On the other hand,  define $\Psi_0: A\otimes C(\T)\to e_0M_ke_0$
by $\Psi_0(a\otimes f)=\psi_0(a)(f(1)e_0)$ for all $a\in A$ and $f\in C(\T)$ (and where $1$ is the point on
the unit circle).   Then
$$
[\psi]|_{{\boldsymbol{\bt}}({\cal  P}_1)}=[\Psi_0]|_{{\boldsymbol{\bt}}({\cal  P}_1)}+[\psi_1]|_{{\boldsymbol{\bt}}({\cal  P}_1)}.
$$
Since ${{\Psi_0(1\otimes z)=e_0}},$  one concludes that
$[\Psi_0]|_{{\boldsymbol{\bt}}({\cal  P}_1)}=0.$ On the other hand,
$\psi_1$ is a \hm\, from $A\otimes C(\T)$ into $M_k,$  {{and so}}  $[\psi_1]|_{{\boldsymbol{\bt}}({\cal  P}_1)}=0$
(this also follows from the first part of Lemma \ref{homhom}).
Thus, by \eqref{EXTMM-1+} and  \eqref{Exmm-10},
$$
[L]|_{\cal P}=[L]|_{{\cal P}_0\cup {{\boldsymbol{\bt}}({\cal  P}_1)}}=[\psi]|_{\cal P}.
$$
}}
{\blue{Since ${\cal H}\subset {\cal G}_0$  and $\ep_0<\sigma,$ by \eqref{Exmm-5} and \eqref{Exmm-4},
the inequality
{{\eqref{EXTMM-3}}} also holds. }}
\end{proof}


{\blue{\begin{df}\label{DDel1}
Let $A$ be a unital \CA, let $X$ be a compact metric space, and  let $B=A\otimes C(X).$
We identify $B$ with $C(X, A).$
Let $b\in C(X, A)_+^{\bf 1}\setminus \{0\}.$ Choose $x_0\in X$ such that
$\|b(x_0)\|>0$ and choose $\ep=\|b(x_0)\|/4.$ Then $b(x_0)>b(x_0)-\ep>0.$
Since $b\in C(X,A)_+,$ there exists a neighborhood $N(x_0)$ of $X$ such that
$b(x)>b(x_0)-\ep.$ Choose  $f\in C(X)_+$ with $0\le f\le 1$ such that its support
lies in $N(x_0).$ Then $b(x)>(b(x_0)-\ep)\otimes f.$ Note
$(b(x_0)-\ep)\otimes f\ge 0.$

Let $\Delta: (A\otimes C(X))_+^{q,\bf 1}\to (0,1)$ be an order preserving map.
Then,
$$
\Delta_0(\hat{h})=\sup \{\Delta(\widehat{h_1\otimes h_2}): h_1\otimes h_2\le h, h_1\in A_+^{\bf 1}\setminus \{0\},
h_2\in C(X)_+^{\bf 1}\setminus \{0\}\}>0
$$
for any $h\in (A\otimes C(X))^{\bf 1}_+\setminus \{0\}.$    Note if $h\ge h',$ then
$$
\Delta_0(\hat{h})\ge \Delta_0(\hat{h'}).
$$
In other words, $\Delta_0: (A\otimes C(X))_+^{q, {\bf 1}}\to (0,1)$ is an order preserving map.
Moreover, $\Delta_0(\hat{h})\le \Delta(\hat{h})$ for all $h\in (A\otimes C(X))_+^{\bf 1}\setminus \{0\}.$
If $h=h_1\otimes h_2$ for some $h_1\in A_+^{\bf 1}\setminus \{0\}$ and
$h_2\in C(X)_+^{\bf 1}\setminus \{0\},$ then $\Delta_0(\hat{h})=\Delta(\hat{h}).$
\end{df}
}}

\begin{lem}\label{homfull}
Let $A\in {\overline{\cal D}}_s$ be a unital \CA\, and
let $\Delta: (A\otimes C(\T))_+^{q, {\bf 1}}\setminus\{0\}\to (0,1)$ be an order  preserving map.
Let $\ep>0$ and let
${\cal F}\subset A$ be a finite subset. There exist a finite subset
${\cal H}_1\subset A_+^{\bf 1}\setminus \{0\},$ a finite subset
${\cal H}_2\subset C(\T)_+^{\bf 1}\setminus \{0\},$
a finite subset ${\cal G}\subset A,$ $\dt>0$
and a finite subset ${\cal P}\subset \underline{K}(A)$ such that,
if $L: A\otimes C(\T)\to M_k$ (for some integer $k\ge 1$) is {{a}}  unital ${\cal G}'$-$\dt$-multiplicative \cp\,
where ${\cal G}'=\{g\otimes f: g\in {\cal G}, f=\{1, z,z^*\}\}$, and $u\in M_k$ is a unitary
such that
\beq\label{homfull-1}
{{\|L(1\otimes z)-u\|}}&<&\dt,
\\\label{homfull-2}
[L]|_{\bt(\cal P)}&=&0\tand\\\label{homfull-3}
{\rm tr}\circ L(h_1\otimes h_2)&\ge &\Delta(\widehat{h_1\otimes h_2})
\eneq
for all $h_1\in {\cal H}_1$ and $h_2\in {\cal H}_2,$ {{where
${\rm tr}$ is the tracial state of $M_k,$}} then
there exists a continuous path of unitaries $\{u_t: t\in [0,1]\}\subset M_k$
with $u_0=u$ and $u_1=1$ such that
\beq\label{homfull-3+}
\|L(f\otimes 1)u_t-u_tL(f\otimes 1)\|<\ep\tforal f\in {\cal F}
\eneq
and $t\in [0,1].$
Moreover, $\{u_t\}$ can be chosen  such that
\beq\label{homfull-4}
{\rm length}(\{u_t\})\le \pi+\ep.
\eneq
\end{lem}

\begin{proof}
{\blue{Let $\Delta_0$ be  {{as}} 
  associated with $\Delta$  in Definition \ref{DDel1}}}. Let $\Delta_1=(1/2)\Delta_{\blue{0}},$ {{let}}
${\cal F}_0=\{f\otimes 1: 1\otimes z: f\in {\cal F}\}$ and
let $B=A\otimes C(\T).$  Then $B\in {\bar{\cal D}}_s.$ Let ${\cal H}'\subset B_+^{\bf 1}\setminus\{0\}$ (in place of ${\cal H}_1$),  {{ ${\cal H}_0\subset {{B}}_{s.a.}$}} {{(in place of ${\cal H}_2$),}} {\blue{$\sigma_0>0$
(in place of $\sigma$),}} ${\cal G}_1\subset A\otimes C(\T)$ (in place of ${\cal G}$), $\dt_1>0$ (in place of $\dt$), {{and}} ${\cal P}'\subset \underline{K}(B)$ (in place of ${\cal P}$) be the finite sets and constants {{provided}} by
Theorem \ref{Newunique1} (for $B$ instead of $A$) for $\ep/16$ (in place of
$\ep$), ${\cal F}_0$ (in place of ${\cal F}$), and $\Delta_{{1}}.$
{\blue{ It is clear that, \wilog, we may assume that every element of ${\cal H}_0$ has norm no more than 1.
If $h\in {\cal H}_0,$ then one may write $h=h_+-h_-,$ where $h_+, h_-\in B_+$ and $\|h_+\|, \, \|h_-\|\le 1.$
Therefore, {{choosing $\sigma_0$ even smaller}},
\wilog,}}
 we may assume that ${\cal H}_0\subset B_+^{\bf 1}\setminus \{0\}.$
 {\blue{Since the elements of the form $\sum_{i=1}^m \af_i a_i,$
 where $\af_i\ge 0,$ and $a_i=a^{(i)}\otimes b^{(i)}$ with $a_i\in A_+^{\bf 1}$ and
 $b_i\in C(\T)_+^{\bf 1},$ are dense in $B_+^{\bf 1},$  choosing an even smaller $\sigma_0,$ \wilog, we may further
 assume that ${\cal H}_0=\{h_1\otimes h_2: h_1\in {\cal H}_1'\andeqn h_2\in {\cal H}_2'\},$
 where}} ${\cal H}_1'\subset A_+^{\bf 1}\setminus \{0\}$ and ${\cal  H}_2'\subset C(\T)_+^{\bf 1}\setminus \{0\}$
 are finite subsets.  {\blue{Similarly, {{choosing even smaller $\sigma_0$ and $\dt_1,$}} }}  we may assume that ${\cal G}_1=\{g\otimes f: g\in {\cal G}_1'\andeqn f\in \{1,z,z^*\}\}$
{\blue{for some}} finite subset ${\cal G}'\subset A.$

{\blue{By the definition of $\Delta_0$ (see \ref{DDel1}),
for each $h\in {\cal H}',$ there exist $a_h\in A_+^{\bf 1}\setminus\{0\}$ and $b_h\in C(X)_+^{\bf 1}\setminus \{0\}$
such that $\hat{h}\ge \widehat{a_h\otimes b_h}$ and
\beq\label{homfull-8-18-1}
\Delta_0(\hat{h})\le (16/15)\Delta(\widehat{a_h\otimes b_h}).
\eneq
Choose finite subsets  ${\cal H}_A\subset A_+^{\bf 1}\setminus \{0\}$ and ${\cal H}_T\subset C(\T)_+^{\bf 1}\setminus \{0\}$
such that, for each $h\in {\cal H}',$ there are $a_h\in {\cal H}_A$ and $b_h\in {\cal H}_T$
such that \eqref{homfull-8-18-1} holds for the triple $h,$ $a_h$ and $b_h.$
Put ${\cal H}''=\{a\otimes b: a\in {\cal H}_A\andeqn b\in {\cal H}_T\}.$
 Replacing both ${\cal H}_1'$ and ${\cal H}_A$  by ${\cal H}_1'\cup {\cal H}_A,$ and
 both ${\cal H}'_2$ and ${\cal H}_T$ by ${\cal H}_2'\cup {\cal H}_A,$
 \wilog, we may assume that ${\cal H}_0={\cal H}'' .$
Let
\beq\label{homfull-8+1}
{\blue{\sigma}}={\blue{(1/4)\min\{}}\min\{\Delta_1(\hat{h}): h\in {\cal H}'\}{\blue{, \sigma_0\}.}}
\eneq
}}
Without loss of generality, {{we}} may  assume that
\beq\label{homfull-8}
{\cal P}'={\cal P}_0\sqcup {\cal P}_1,
\eneq
where ${\cal P}_0\subset \underline{K}(A)$ and
${\cal P}_1\subset {\boldsymbol{\bt}}(\underline{K}(A))$ are finite subsets.
Let ${\cal P}\subset \underline{K}(A)$ be a finite subset such that
${\boldsymbol{\bt}}({\cal P})={\cal P}_1.$

{{Let}} $\dt_2>0$ (in place of $\dt$) with $\dt_2<\ep/16,$ {{the}} finite subset ${\cal G}_2\subset A\otimes C(\T)$ (in place of ${\cal G}$), and
the  finite subset ${\cal H}_3\subset (A\otimes C(\T))_+^{\bf 1}\setminus \{0\}$ (in place of ${\cal H}_1$)
{{be as provided}} by Lemma \ref{EXTMM} for ${\blue{\sigma}},$
{\blue{$(3/4)\Delta_0$ (in place of $\Delta$),}}  ${\blue{\cal H}''}$ (in place of
${\cal H}$), $\min\{\ep/16,\dt_1/2\}$ (in place of $\ep$),
${\cal G}_1$ (in place of ${\cal G}_0$),  {{and}} ${\cal P}_0$ and ${\cal P}$
(in place of ${\cal P}_0$ and ${\cal P}_1$).
{{Choosing a}} smaller $\dt_2,$ we may also assume that
$$
{\cal G}_2=\{g\otimes f: g\in {\cal G}_2'\andeqn f\in \{1, z,z^*\}{{\}}}
$$
for a finite set ${\cal G}_2'\subset A$.
Let
$$
{\cal H}_3'=\{h_1\otimes h_2: h_1\in {\cal H}_4\andeqn h_2\in {\cal H}_5\}
$$
for 
finite {\blue{subsets}} ${\cal H}_4\subset A_+^{\bf 1}\setminus \{0\}$ and
${\cal H}_5\subset C(\T)_+^{\bf 1}\setminus \{0\}$ {{be}}  {\blue{such that,
if $h\in {\cal H}_3,$ then there are $a_h\in {\cal H}_4$ and $b_h\in {\cal H}_5$
such that $\widehat{h}\ge {\widehat{a_h\otimes b_h}}$ and
\beq\label{6418827-1}
\Delta_0(\widehat{h})\le (16/15)\Delta(\widehat{a_h\otimes b_h}).
\eneq}}
Let ${\cal G}={\cal F}\cup {\cal G}_1'\cup {\cal G}_2',$ $\dt=\min\{\dt_1/2,\dt_2/2, \ep/16\},$
${\cal H}_1={\cal H}_{\blue{A}}\cup {\cal H}_4,$ and
${\cal H}_2={\cal H}_{\blue{T}}\cup {\cal H}_5.$

Now suppose that one has a {\blue{unital \cp\,}} $L: A\otimes C(\T)\to M_k$ and a unitary $u\in M_k$ satisfying the assumptions
(\eqref{homfull-1} to \eqref{homfull-3}) with
the above ${\cal H}_1,$ ${\cal H}_2,$ ${\cal G},$ ${\cal P},$  $\dt$ and $\sigma.$
{\blue{In particular,  by \eqref{6418827-1} and by the assumption \eqref{homfull-3}, if $h\in {\cal H}_3,$  there
are $a_h\in {\cal H}_4\subset {\cal H}_1$ and $b_h\in {\cal H}_5\subset {\cal H}_2$ such that
\beq
{\rm tr}(L(\hat{h})))\ge {\rm tr}(L(\widehat{a_h\otimes b_h}))\ge \Delta(\widehat{a_h\otimes b_h})
\ge (15/16)\Delta_0(\widehat{h})\ge (3/4)\Delta_0(\widehat{h}).
\eneq}}
It follows from Lemma \ref{EXTMM} {{(on using also \eqref{homfull-2})}} that there is a unital  {{${\cal G}_1$-$\min\{\ep/16, \dt_1/2\}$}}-multiplicative \cp\,
$\psi: A\otimes C(\T)\to M_k$ such that $w=\psi(1\otimes z)$ is a unitary,
\beq\label{homfull-9}
w\psi(g\otimes 1)&=&\psi(g\otimes 1)w\tforal g\in A,\\\label{homfull-9n}
[\psi]|_{\cal P'}&=&[L]|_{\cal P'},\tand\\\label{homfull-9+}
|{\rm tr}\circ L(g)-{\rm tr}\circ \psi(g)|&<&\sigma \tforal g\in {\cal H}''.
\eneq
It follows that, {\blue{for $h\in {\cal H}'',$ there are $a_h\in {\cal H}_A$ and $b_h\in {\cal H}_T$
such that $h=a_h\otimes b_h$ and (by  \eqref{homfull-3}),}}
\beq\label{homfull-10}
\mathrm{tr}\circ \psi(h)&\ge& \mathrm{tr}\circ L(h)-\sigma={\rm tr}(L(a_h\otimes b_h))-\sigma\\
&\ge &
\Delta(\widehat{a_h\otimes b_h})-\sigma\ge {\blue{(3/4)\Delta(\widehat{a_h\otimes b_h})=(3/4)\Delta(\hat{h}).}}
\eneq
{\blue{If $h\in {\cal H}',$ there is $h'\in {\cal H}''$ such that $h\ge h'$ and
\beq
(15/16)\Delta_0(\hat{h})\le \Delta(\widehat{h'})
\eneq
(see \eqref{6418827-1}).
It follows that
\beq\label{homfull-18-n2}
\mathrm{tr}\circ \psi(h) &\ge&  \mathrm{tr}\circ \psi(h')\ge (3/4)\Delta(\hat{h'})\\\label{homfull-18-n3}
&\ge & (3/4)(15/16)\Delta_0(\hat{h})\ge \Delta_1(\hat{h})\rforal h\in {\cal H}'.
\eneq}}

{\blue Note that we have assumed that ${\cal H}_0={\cal H}''.$}
Combining (\ref{homfull-9n}), {\blue{\eqref{homfull-18-n3}}},
and (\ref{homfull-9+}),
{{and}} applying  Theorem \ref{Newunique1},
one obtains a unitary $U\in M_k$ such that
\beq\label{homfull-12}
\|{\rm Ad}\,U\circ \psi(f)-L(f)\|<\ep/16\tforal f\in {\cal F}_0.
\eneq
Let $w_1={\rm Ad}\,U\circ \phi(1\otimes z).$ Then
\beq\label{homfull-13}
\|u-w_1\| &\le & \|u-L(1\otimes z)\|+\|L(1\otimes z)-{\rm Ad}\, U\circ \psi(1\otimes z)\|\\
&<&\dt+\ep/16<\ep/8.
\eneq
{{Thus there is $h\in (M_k)_{s.a.}$ with $\|h\|<\ep\pi/8$ such
that  $uw_1^*=\exp(ih).$}}
It follows that  there is a continuous path of unitaries $\{u_t\in [0,1/2]\}\subset M_k$
such that
\beq\label{homfull-13+}
\|u_t-u\|<\ep/8,\,\,\,\|u_t-w_1\|<\ep/8,\,\,\, u_0=u,\,\,\, u_{1/2}=w_1,\\
\andeqn
{\rm length}(\{u_t: t\in [0,1/2]\})<\ep\pi/8.
\eneq
It follows from Lemma \ref{homhom} that there exists a continuous path of unitaries $\{u_t: t\in [1/2,1]\}\subset M_k$ such that
\beq\label{homfull-14}
u_{1/2}=w_1,\,\,\, u_1=1\andeqn u_t({\rm Ad}\, U\circ {{\psi}}(f\otimes 1))=({\rm Ad}\, U\circ {{\psi}}(f\otimes 1))u_t
\eneq
for all $t\in [1/2,1]$ and $f\in A\otimes 1.$
Moreover, we may assume that
\beq\label{homfull-15}
{\rm length}(\{u_t: t\in [1/2,1]\})\le \pi.
\eneq
It follows that
\beq\label{homfull-16}
{\rm length}(\{u_t: t\in [0,1]\}) \le \pi+\ep\pi/6.
\eneq
Furthermore,
\beq\label{homfull-17}
\|u_tL(f\otimes 1)-L(f\otimes 1)u_t\|<\ep\tforal f\in {\cal F}
\eneq
and $t\in [0,1].$
%
\end{proof}

\begin{lem}\label{changespectrum}
Let $A\in {\bar {\cal D}}_s$ be a unital \CA,
{\blue{l}}et  $\ep>0$  and let ${\cal F}\subset A$ be a finite subset.
 Let ${\cal H}_1\subset A_+^{\bf 1}\setminus\{0\}$ and let  ${\cal H}_2\subset C(\T)_+^{\bf 1}\setminus \{0\}$  be  finite subsets.
For any order preserving map $\Delta: A_+^{q, {\bf 1}}\setminus \{0\}\to (0,1),$ there exists a finite subset ${\cal G}\subset A,$
a finite
subset ${\cal H}_1'\subset A_+^{\bf 1}\setminus \{0\},$
and $\dt>0$ satisfying the following condition: for any unital ${\cal G}$-$\dt$-multiplicative
\cp\, $\phi: A\to M_k$ (for some integer $k\ge 1$) and any unitary $u\in M_k$ such that
\beq\label{change-1}
\|u\phi(g)-\phi(g)u\|<\dt\tforal g\in {\cal G} \tand \\ \label{change-1+}
{\blue{\rm{tr}}}\circ \phi(h)\ge \Delta(\hat{h})\tforal h\in {\cal H}_1',
\eneq
there exists a continuous path of unitaries
$\{u_t: t\in [0,1]\}\subset M_k$  such that
\beq\label{change-2}
&&u_0=u,\,\,\, u_1=w,\,\,\,\|u_t\phi(f)-\phi(f)u_t\|<\ep\tforal f\in {\blue{{\cal F}}}\tand t\in [0,1], {{\text{and}}}\\
&&\hspace{0.4in}\mathrm{tr}\circ L(h_1\otimes h_2)\ge  \Delta(\widehat{h_1}) \tau_m(h_2)/4
\eneq
for all $h_1\in {\cal H}_1$ and $h_2\in {\cal H}_2,$
where $L: A\otimes C(\T)\to M_k$ is a unital   {\blue{${\cal F}_1$-$\ep$-multiplicative}}  \cp\,
such that
\beq\label{changes-4}
\|L(f\otimes 1_{\blue{C(\T)}})-\phi(f)\|<\ep \tforal f\in {\cal F}, \tand
\|L(1\otimes z)-w\|<\ep,
\eneq
and
$\tau_m$ is the tracial state on $C(\T)$ induced by the Lebesgue
measure
on the circle, {\blue{where ${\cal F}_1=\{f\otimes g: f\in {\cal F}\cup\{1\}, g\in \{1_{C(\T)}, z, z^*\}\}.$}}
Moreover, {{$\{u_t\}$ can be chosen such that}}
\beq\label{changes-5}
{\rm length}(\{u_t\})\le 2\pi+\ep.
\eneq

\end{lem}


\begin{proof}
{\blue{\Wlog, we may assume that ${\cal F}$ is in the unit ball of $A$ and
if $f\in {\cal F},$ then $f^*\in {\cal F}.$}}

There exists an integer $n\ge 1$ such that
\beq\label{changes-6}
 (1/n)\sum_{j=1}^n f(e^{\theta+j 2\pi i/n})\ge (63/64)\tau_m(f)
\eneq
for all $f\in {\cal H}_2$ and for any $\theta\in [-\pi, \pi].$
We may also assume that $16\pi/n<\ep.$


Let
$$
\sigma_1=(1/2^{10})\inf\{{\blue{\Delta}}({{\hat{h}}}): h\in {\cal H}_1\}\cdot  \inf\{\tau_m(g): g\in {\cal H}_2\}.
$$
Let ${\blue{{\cal H}_{1a}}}\subset A_+^{\bf 1}\setminus \{0\}$
(in place of ${\cal H}_2$)
be a finite subset {{ as provided}} by  Corollary
\ref{repcor} for $\min\{\ep/32,\sigma_1/16\}$ (in place of $\ep$), ${\cal F}\cup {\cal H}_1$ (in place of ${\cal F}$), ${\cal H}_1$ (in place of ${\cal H}$), $(1-1/2^{12})\Delta$ (in place of $\Delta$),  $\sigma_1/16$ (in place of $\sigma$), and {{the}} integer
$n$ {\blue{(and for $A$).}}

{\blue{Put ${\cal H}_1'={\cal H}_{1a}\cup {\cal H}_1.$
Put
$\sigma=(1/2^{10})\inf\{{\blue{\Delta}}({{\hat{h}}}): h\in {\cal H}_1'\}\cdot  \inf\{\tau_m(g): g\in {\cal H}_2\}.$ }}
{\blue{Note that $A\otimes C(\T)$ is a {{subhomogeneous}} \CA.}}
Let ${\cal F}'=\{f\otimes 1_{\blue{C(\T)}}, f\otimes z: f\in {\cal F}\cup{\blue{{\cal H}_1'}}\}.$
Let $\dt_1>0$ (in place of $\dt$) and ${\cal G}_1\subset A\otimes C(\T)$ (in place of ${\cal G}$) be {{as provided}} by
Lemma \ref{LtoMn} for ${\blue{\ep_0:=\min\{\ep/64, \sigma/16\}}}$
(in place of $\ep$), ${\cal F}'$ (in place of ${\cal F}$), and ${\blue{\sigma}}/16$ (in place of $\sigma_0$).
{\blue{Choosing a smaller $\dt_1$ if necessary,}} without loss of generality, one may assume that, for a finite set ${\cal G}_2\subset A,$
$$
{\cal G}_1=\{g\otimes 1, 1\otimes z: g\in {\cal G}_2\}.
$$
{\blue{We may also assume that if $g\in {\cal G}_2,$ then $g^*\in {\cal G}_2.$}}
Put
$$
{\cal H}'=\{h_1\otimes h_2, h_1\otimes 1, 1\otimes h_2: h_1\in {\blue{{\cal H}_1'}}\andeqn h_2\in {\cal H}_2\}.
$$

Let ${\cal G}_3={\cal G}_2\cup  {\cal H}_1'.$ To simplify  notation, without loss of generality, let us assume that
${\cal G}_3$  {\blue{is}} in the unit ball of {\blue{$A$}} and  ${\cal F}'$ {\blue{is in the unit ball
of}}  $A\otimes C(\T),$ {\blue{respectively}}. Let
$\dt_2=\min\{\ep/64,\dt_1/2, {\blue{\sigma}}/16\}.$

Let ${\cal G}_4\subset A$ be a finite subset (in place of ${\cal G}$) and let
$\dt_3$ (in place of $\dt$) be {\blue{a}} positive {\blue{number}} as provided by \ref{oldnuclearity} for ${\cal G}_3$ (in place of ${\cal F}_0$), ${\cal F}'$ (in place of ${\cal F}$),  and $\dt_2$ (in place of $\ep$).

Let ${\cal G}={\cal G}_4\cup {\cal G}_3\cup {\cal F}$ and
$\dt=\min\{\dt_1/4, \dt_2/2, \dt_3/2\}.$
Now let $\phi: A\to M_k$ be a unital ${\cal G}$-$\dt$-multiplicative \cp, and let $u\in M_k$ be a unitary
such that (\ref{change-1}) and (\ref{change-1+}) hold for the above $\dt,$
${\cal G},$ and
${\cal H}_1'.$

{{By}} Lemma \ref{oldnuclearity}, there
exists a {{unital}} ${\cal G}_3$-$\dt_2$-multiplicative \cp\,
$L_1: A\otimes C(\T)\to M_k$ such that
\beq\label{changes-8}
\|L_1(g\otimes 1_{\blue{C(\T)}})-\phi(g)\|<\dt_2\tforal g\in {\cal G}_{\blue 3}\andeqn
\|L_1(1\otimes z)-u\|<\dt_2.
\eneq
Then
\beq\label{changes-8+1}
{\rm tr}\circ L_1(h\otimes 1)&\ge &{\rm tr}\circ \phi(h)-\dt_2\\\label{changes-8+2}
&\ge & \Delta(\hat{h})-\sigma/16\ge (1-1/2^{14})\Delta(\hat{h})
\eneq
for all $h\in {\cal H}_1'.$
It follows from Lemma   \ref{LtoMn} that there exist a projection
$p\in M_k$ and a unital \hm\, $\psi: A\otimes C(\T)\to pM_kp$ such that
\beq\label{changs-9}
&&\|pL_1(f)-L_1(f)p\|<{\blue{\ep_0}}
\tforal f\in {\cal F}',\\\label{changes-9+}
&&\|L_1(f)-((1-p)L_1(f)(1-p)+\psi(f))\|<{\blue{\ep_0}}
\tforal f\in {\cal F}' ,\\\label{changes-10}
&&\andeqn {\rm tr}(1-p)<{\blue{\sigma}}/16.
\eneq
Note that $pM_kp\cong M_m$ for some $m\le k.$ It follows from (\ref{changes-8+2}),
(\ref{changes-9+}),  and (\ref{changes-10}) that
\beq\label{changes-11}
{\blue{{\rm tr}'\circ \psi(h)}}\ge {\rm tr}\circ \psi(h)\ge (1-1/2^{14})\Delta(\hat{h})-{\blue{\sigma}}/16-{\blue{\sigma}}/16 \ge (1-1/2^{12})\Delta(\hat{h})
\eneq
{\blue{for all $h\in {\cal H}_1',$ where ${\rm tr}'$ is the normalized trace on $pM_kp\cong M_m.$}}
It follows from \ref{homhom} that  there is a continuous path of unitaries $\{u_t: t\in [1/4, 1/2]\}\subset pM_kp$ such that $u_{1/4}=\psi(1\otimes z),$  $u_{1/2}=p$, and
\beq\label{add1601-1}
u_t\psi(f\otimes 1_{\blue{C(\T)}})=\psi(f\otimes 1_{\blue{C(\T)}})u_t\rforal  f\in A\andeqn \rforal t\in [{\blue{1/4}},1/2],
\eneq
and ${\rm length}(\{u_t: t\in [1/4,1/2]\})\le \pi.$

By  Corollary \ref{repcor}, there are mutually orthogonal
projections $e_0, e_1, e_2,...,e_n\in pM_kp$  such that  $e_1, e_2,...,e_n$ are {\blue{mutually}} equivalent and
{\blue{$\sum_{i=0}^ne_i=p,$}} and there are unital \hm s $\psi_0: A
\to e_0M_ke_0$ and
$\psi_1: A
\to e_1M_ke_1$ such that
\beq\label{changes-12}
&&{\blue{{\rm tr}'}}(e_0)<\sigma_1/16\quad {{\textrm{and}}}\\\label{chnages-12+}
&&\hspace{-0.6in}\|\psi(f\otimes 1_{C(\T)})-(\psi_0(f)\oplus \diag(\overbrace{\psi_1(f\otimes 1_{C(\T)}),...,\psi_1(f\otimes 1_{C(\T)}))}^n)\|<\min\{\ep/32, \sigma_1/16\}
\eneq
{\blue{for all $f\in {\cal F}\cup {\cal H}_1,$ where we identify $(\sum_{i=1}^ne_i)M_k(\sum_{i=1}^ne_i)$ with
$M_n(e_1M_ke_1)$ {\blue{(using the convention introduced in  \ref{Ddiag}).}} }}
Moreover,
\beq\label{changes-13}
{\blue{{\rm tr}'}}(\psi_1(h\otimes 1))\ge (1-1/2^{12})\Delta(\hat{h})/{\blue{2n}}\tforal h\in {\cal H}_1.
\eneq
{\blue{Note this implies that, by \eqref{changes-10},}}
\beq\label{cahnges-188-n1}
{\blue{{\rm tr}(\psi_1(h\otimes 1))\ge (1-\sigma/16)(1-1/2^{12})\Delta(\hat{h})/2n\tforal h\in {\cal H}_1.}}
\eneq
Let  $w_{0,j}=\exp(i(2\pi j/n))e_j,$ $j=1,2,...,n.$  Define {\blue{$w_{00}:=\sum_{j=1}^n w_{0,j}=\diag(w_{0,1},w_{0,2},...,w_{0,n}),$
and define}}
\beq\label{1205-65-1}
w_0'= e_0\oplus {\blue{w_{00}}}=e_0\oplus { {\rm diag}(w_{0,1},w_{0,2},...,w_{0,n})}.
\eneq
Then  $w_0'$ {{commutes}} with $\psi_0(f{\blue{\otimes 1_{C(\T)}}})\oplus {\rm diag}(\psi_1(f{\blue{\otimes 1_{C(\T)}}}),\psi_1(f{\blue{\otimes 1_{C(\T)}}}),...,\psi_1(f{\blue{\otimes 1_{C(\T)}}}))$
for all $f\in A.$
As in \ref{homhom}, there exists a continuous path $\{u_t: t\in [1/2, 1]\}\subset pM_kp$
such that $u_{1/2}=p,$ $u_{1}=w_0'$ and $u_t$ commutes with
$\psi_0(f)\oplus {\rm diag}(\psi_1(f),\psi_1(f),...,\psi_1(f))$ for all $f\in A,$ and moreover,
${\rm length}(\{u_t: t\in [1/2, 1]\}\le \pi.$

There is a unitary $w_0''\in (1-p)M_k(1-p)$ such that (see \eqref{changes-8} and \eqref{changs-9})
\beq\label{changes-15}
\|w_0''-(1-p)L_1(1\otimes z)(1-p)\|<\ep/16.
\eneq
{\blue{For $f\in {\cal F},$  by \eqref{changes-15}, \eqref{changs-9},  and
the fact that $L_1$ is ${\cal G}_3$-$\dt_2$ multiplicative,
\beq
&&\hspace{-0.6in}w_0''(1-p)L_1(f\otimes 1_{C(\T)})(1-p)\approx_{\ep/16} (1-p)L_1(1\otimes z)(1-p)L_1(f\otimes 1_{C(\T)})(1-p)\\
&&\approx_{\ep_0}(1-p)L_1(1\otimes z)L_1(f\otimes 1_{C(\T)})(1-p)\\
&&\approx_{\dt_2} (1-p)L_1((f\otimes 1_{C(\T)})(1\otimes z))(1-p)\\
&&\approx_{\dt_2} (1-p)L_1(f\otimes 1_{C(\T)})L_1(1\otimes z))(1-p)\\
&&\approx_{\ep_0}(1-p)L_1(f\otimes 1_{C(\T)})(1-p)(1-p)L_1(1\otimes z))(1-p)\\\label{Change-18826-1}
&&\approx_{\ep/16} (1-p)L_1(f\otimes 1_{C(\T)})(1-p) w_0''.
\eneq
In other words,  for all $f\in {\cal F},$
\beq\nonumber
&&\hspace{-0.6in}\|w_0''(1-p)L_1(f\otimes 1_{C(\T)})(1-p)-(1-p)L_1(f\otimes 1_{C(\T)})(1-p)w_0''\|\\\label{Change-18826-2}
&&<
\ep/16+2\ep_0+2\dt_2+\ep/16<\ep/16+\ep/32+\ep/32+\ep/16=3\ep/16.
\eneq}}
{\blue{Note also,  by \eqref{changes-8} and \eqref{changes-9+},
\beq\label{Changes-18826-3}
\hspace{-0.2in}\|\phi(f)-((1-p)L_1(f\otimes 1_{C(\T)})(1-p)+\psi(f\otimes 1_{C(\T)}))\|<\dt_2+\ep_0<\ep/32\tforal f\in {\cal F}.
\eneq
}}
Put
$u_0'=w_0''\oplus \psi(1\otimes z).$
Then $u_0'$ is a unitary and
\beq\label{changes-17}
\|u-u_0'\| &\le & \|u-L_1(1\otimes z)\|+\|L_1(1\otimes z)-u_0'\|\\
&\le & \dt_2+\ep/16<\ep/8.
\eneq
{\blue{{{As above,}} we obtain a continuous path of unitaries
$\{w_t\in [0,1/4]\}\subset M_k$ such that
$w_0=u$ and $w_{1/4}=u_0'$ and ${\rm cel}(\{w_t: t\in [0,1/4]\})<2\arcsin(\ep/16).$}}
Define $w_t=w_0''\oplus u_t$ for $t\in [1/4, 1].$
Then $w_{1/4}=w_0''\oplus \psi(1\otimes z)=u_0',$
$w_1=w_0''\oplus w_0'.$
{\blue{From the construction, $\{w_t: t\in [0,1]\}$ is a continuous path of unitaries in $M_k$
such that
\beq
w_0=u,\,\,\, w_1=w_0''\oplus w_0' \andeqn
{\rm cel}(\{w_t: t\in [0,1]\})\le 2\pi +\ep.
\eneq
}}
{\blue{By \eqref{change-1} and the choice of $\dt,$
for $t\in [0,1/4]$ and $f\in {\cal F},$ }}
\beq\label{changes-1825-1}
\|w_t\phi(f)-\phi(f)w_t\|<\ep\rforal  f\in {\cal F}.
\eneq
{\blue{For $t\in [1/4, 1/2],$ by \eqref{Changes-18826-3} and by \eqref{Change-18826-2},
for all $f\in {\cal F},$
\beq\nonumber
&&\hspace{-0.5in}w_t\phi(f)=(w_0''\oplus u_t)\phi(f)\approx_{\ep/32} (w_0''\oplus u_t)((1-p)L_1(f\otimes 1_{C(\T)})(1-p)\oplus
\psi(f\otimes 1_{C(\T)}))\\\nonumber
&&\approx_{3\ep/16} ((1-p)L_1(f\otimes 1_{C(\T)})(1-p)\oplus
\psi(f\otimes 1_{C(\T)}))(w_0''\oplus u_t) \approx_{\ep/32} \phi(f)w_t.
\eneq
}}
{\blue{ Similarly, since $u_t$ commutes with $\psi_0(f)\oplus {\rm diag}(\psi_1(f),\psi_1(f),...,\psi_1(f))$ for all $f\in A,$
one also has, for $t\in [1/2, 1],$
\beq\label{changes-1825-2}
\|w_t\phi(f)-\phi(f)w_t\|<\ep\rforal  f\in {\cal F}.
\eneq
Therefore, \eqref{changes-1825-2} holds for all $t\in [0,1].$}}
\
Define $L: A\otimes C(\T)\to M_k$ by
\vspace{-0.1in}\beq\nonumber
L(a\otimes f)=(1-p)L_1(a\otimes f)(1-p)\oplus
({\rm diag}(\psi_0(a{\blue{\otimes 1_{C(\T)}}}), \overbrace{\psi_1(a{\blue{\otimes 1_{C(\T)}}}),...,\psi_1(a{\blue{\otimes 1_{C(\T)}}})}^n)f({\blue{w_0'}}))
\eneq
for all $a\in A$ and $f\in C(\T).$
{\blue{Then}}
\beq\label{changes-19}
\|L(f\otimes 1)-\phi(f)\|<\ep\tforal f\in {\cal F}.
\eneq
{\blue{Moreover {{(see \ref{Ddiag})}},
\beq
L(1\otimes z)=(1-p)L_1(1\otimes z)(1-p)\oplus ({\rm diag}(\psi_0(1), \overbrace{\psi_1(1),...,\psi_1(1)}^n)(w_0'))\\
=(1-p)L_1(1\otimes z)(1-p)\oplus pw_0'\approx_{\ep/16} w_0''\oplus w_0'=w_1.
\eneq
One then {{verifies}} that {{$L$ is}} ${\cal F}_1$-$\ep$-multiplicative, where ${\cal F}_1=\{f\otimes g: f\in {\cal F}\cup\{1\}, g\in \{1_{C(\T)}, z, z^*\}\}.$}}

{\blue{In the next few lines, for $h_1\in A_+$ and $h_2\in C(\T)_+,$ we write  $h_1':=h_1\otimes 1_{C(\T)}$ and $h_2':=1\otimes h_2.$ Also, we view $h_2$ as a positive function on $\T.$
{{Recall}} that $w_{0,j}=\exp(i(2\pi j/n))e_j$ is a scalar multiple of $e_j.$
Therefore, the element $h_2(w_{0,j})=h_2(e^{i2\pi j/n})e_j$ is also a scalar multiple of $e_j.$
Moreover,
$h_2(w_{00})$ may be written as ${{\diag(h_2(e^{i2\pi /n}), h_2(e^{i2\pi 2/n}),...,h_2(e^{i2\pi n/n}))}}.$}}
{\blue{Finally, we  estimate,}} {\blue{for for all $h_1\in {\cal H}_1$ and $h_2\in {\cal H}_2,$
keeping in mind of \ref{Ddiag}, by \eqref{1205-65-1},}}
\beq\nonumber
{\rm tr}\circ L(h_1\otimes h_2) &\ge & {\rm tr}(\psi_0({\blue{h_1'}}))+{\rm  tr}(\diag(\overbrace{\psi_1({\blue{h_1'}}),...,\psi_1({\blue{h_1'}})}^n)h_2(w_{\blue{00}}))\\\nonumber
&\ge &{\blue{{\rm  tr}(\diag(\overbrace{\psi_1({\blue{h_1'}}),...,\psi_1({\blue{h_1'}})}^n)h_2(w_{00}))}}\\\nonumber
&=&{\blue{\sum_{j=1}^n{\rm  tr}(\psi_1(h_1')h_2(e^{i2\pi j/n}))=\sum_{j=1}^n{\rm  tr}(\psi_1(h_1'))(h_2(e^{i2\pi j/n}))}}\\\nonumber
&\ge & {\blue{((1-\sigma_2/16)(1-1/2^{12})\Delta(\hat{h})/3n)(\sum_{j=1}^n h_2(e^{i2\pi j/n}))}}\,\,\,\,\,\,\,\,\,\,\,\,\,\,\,\,
~~~~{\rm {\blue{(by\,\,\,\eqref{cahnges-188-n1})}}}\\\nonumber
&\ge & (1-1/2^{10})\Delta(\widehat{h_1})(63/64)\tau_m(h_2)/3
\ge \Delta(\widehat{h_1})\cdot \tau_m(h_2)/4.\,\,\,\,\,\,
~~~~{\blue{\rm (by\,\,\,\eqref{changes-6})}}
\eneq
\end{proof}

\begin{df}\label{Ddel}
Let $A$ be a unital \CA\, with $T(A)\not=\emptyset$ and let $\Delta: A_+^{q,{\bf 1}}\setminus \{0\}\to (0,1)$ be an order preserving
map. {{Let}} $\tau_m: C(\T)\to \C$ {{denote}} the tracial state given by normalized Lebesgue
 measure.
Define
$\Delta_1: (A\otimes C(\T))_+^{q, {\bf 1}}\setminus\{0\}\to (0,1)$ by
\beq\label{Ddel-1}
\hspace{-0.3in}\Delta_1(\hat{h})=\sup\{ {\Delta(h_1)\tau_m(h_2)\over{4}}:\,
\hat{h}\ge \widehat{h_1\otimes h_2}\andeqn h_1\in A_+\setminus \{0\},\,\,\, h_2\in C(\T)_+\setminus\{0\}\}.
\eneq
\end{df}

\begin{lem}\label{homotopy1}
Let $A\in  {\bar{\cal D}}_s$ be a unital \CA.
Let $\Delta: A_+^{q, {\bf 1}}\setminus \{0\}\to (0,1)$ be an order preserving  map.
For any $\ep>0$ and any finite subset ${\cal F}\subset A,$ there exists a finite subset ${\cal H}\subset A_+^{\bf 1}\setminus\{0\},$
 $\dt>0,$  a finite subset
${\cal G}\subset A$,  and a finite subset ${\cal P}\subset \underline{K}(A)$ satisfying the following {{condition}}:
For any unital  ${\cal G}$-$\dt$-multiplicative \cp\,
$\phi: A\to M_k$ (for some integer $k\ge 1$), and any unitary $v\in M_k$, such that
\beq\label{homot-1}
&&tr\circ \phi(h)\ge \Delta(\hat{h}) \tforal h\in {\cal H},\\\label{homot-2}
&&\|\phi(g)v-v\phi(g)\|<\dt\tforal g\in {\cal G}, \tand\\\label{homot-2+}
&&{\rm  Bott}(\phi,v)|_{\cal P}={{0}},
\eneq
there exists a continuous path of unitaries
$\{u_t: t\in [0,1]\}\subset M_k$ such that
\beq\label{homot-3}
u_0=v,\,\,u_1=1,\tand \|\phi(f)u_t-u_t\phi(f)\|<\ep
\eneq
for all $t\in [0,1]$ and $f\in {\cal F},$ {{and, moreover,}}
\beq\label{homot-4}
{\rm length}(\{u_t\})\le 2\pi+\ep.
\eneq

\end{lem}

\begin{proof}
Let $\Delta_1$ be as in Definition \ref{Ddel} {\blue{(associated with the given $\Delta$).}}
Let ${\cal H}_1\subset A_+^{\bf 1}\setminus \{0\}$ and ${\cal H}_2\subset C(\T)_+^{\bf 1}\setminus \{0\}$ {{and}}
${\cal G}_1\subset A$ (in place of ${\cal G}$), and  ${\cal P}\subset \underline{K}(A),$ be  finite subsets {{and $\delta_1>0$ (in place of $\delta$) the constant, as provided }} by \ref{homfull} for $\ep/4$ (in place of $\ep$),
 ${\cal F}$ and $\Delta_1.$ {\blue{We may assume that $1\in {\cal F}.$
 \Wlog, we may assume that $(\dt_1, {\cal G}_1, {\cal P})$ is a $KL$-triple (see \ref{KLtriple})
 Moreover, we may assume that $\dt_1<\dt_{\cal P}$ and ${\cal G}_1\supset {\cal F}_{\cal P}$
 (see \ref{Dbeta}).}}


Let ${\cal G}_2\subset A$ (in place of ${\cal G}$) and ${\cal H}_1'\subset A_+^{{\bf 1}}\setminus \{0\}$ {{denote the}} finite subsets, {{and}} $\dt_2>0$ (in place of $\dt$) {{the constant, provided}} by Lemma
\ref{changespectrum} for $\min\{\ep/16, \dt_1/2\}$ (in place of $\ep$),
${\cal G}_1\cup {\cal F}$ (in place of ${\cal F}$), and ${\cal H}_1$ and ${\cal H}_2$ {\blue{as well as $\Delta.$}}
Let ${\cal G}={\cal G}_2\cup {\cal G}_1{\blue{\cup}} {\cal F},$ {{let ${\cal H}={{{\cal H}_1'}},$ and}}
let $\dt=\min\{\dt_2, \ep/16\}.$

Now suppose that $\phi: A\to M_k$ is a unital ${\cal G}$-$\dt$-multiplicative \cp\,
and $u\in M_k$ is a unitary which satisfy the assumptions  {\blue{(in particular, \eqref{homot-1}, \eqref{homot-2}, and \eqref{homot-2+} hold)}} for the above ${\cal H},$ $\dt,$
${\cal G}$, and  ${\cal P}.$

{{Applying}} Lemma \ref{changespectrum}, one obtains a continuous path of unitaries $\{u_t: t\in [0,1/2]\}\subset M_k$ such that
\beq\label{homot-5}
u_0=u,\,\,\, u_1=w,\,\,\,\|u_t\phi(g)-\phi(g)u_t\|<\min\{\dt_1, \ep/4\}
\eneq
for all $g\in {\cal G}_1\cup {\cal F}$ and $t\in [0,1/2].$
Moreover, there is a unital {\blue{${\cal G}'$-$\min\{\dt_1/2, \ep/16\}$-multiplicative}} \cp\,
$L: A\otimes C(\T)\to M_k,$
{\blue{where ${\cal G}'=\{g\otimes f: g\in {\cal G}_1\cup {\cal F}\andeqn g\in \{1_{C(\T)}, z,z^*\}\}$,}} such that
\beq\label{homot-6}
&&\|L(g\otimes 1)-\phi(g)\|<\min\{\dt_1{\blue{/2}}, \ep/{\blue{16}}\}\tforal g\in {\cal G}_1\cup {\cal F},\\\label{homot-7}
&&\|L(1\otimes z)-w\|<\min\{\dt_1{\blue{/2}}, \ep/{\blue{16}}\},
\\\label{homot-8}
&&\andeqn {\rm tr}\circ L(h_1\otimes h_2)\ge \Delta(h_1)\tau_m(h_2)/4
\eneq
for all $h_1\in {\cal H}_1$ and $h_2\in {\cal H}_2.$ Furthermore,
\beq\label{homot-12}
{\rm length}(\{u_t:t\in [0,1/2]\})\le \pi+{{\ep/16}}.
\eneq

Note that {\blue{(see \ref{Dbeta})}}
\beq\label{homot-13}
[L]|_{\bt({\cal P})}={\rm Bott}(\phi, w)|_{\cal P}=
{\rm Bott}(\phi, u)|_{\cal P}=0.
\eneq

By (\ref{homot-6}), (\ref{homot-7}), (\ref{homot-13}), and (\ref{homot-8}),
applying Lemma \ref{homfull},
one obtains a continuous path of unitaries
$\{u_t\in [1/2, 1]\}\subset M_k$ such that
\beq\label{homot-14}
u_{1/2}=w,\,\,\,u_1=1,\,\,\,
\|u_t\phi(f)-\phi(f)u_t\|<\ep/4\tforal f\in {\cal F},\\
\andeqn {\rm length}(\{u_t: t\in [1/2, 1]\})\le \pi+\ep/4.
\eneq
Therefore, $\{u_t: t\in [0,1]\}\subset M_k$ is a continuous path of unitaries in $M_k$ with $u_0=u$ and $u_1=1$ such that
\beq\label{homot-15}
\|u_t\phi(f)-\phi(f)u_t\|<\ep\tforal f\in {\cal F}
\andeqn
{\rm length}(\{u_t: t\in [0,1]\})\le 2\pi+\ep.
\eneq
\end{proof}

\section{An Existence Theorem for  Bott maps }

\begin{lem}\label{Ext1}
Let $A$ be a unital {{amenable separable}} residually finite dimensional  \CA\, which satisfies the UCT, 
let
$G=\Z^r\oplus {\rm Tor}(G)\subset K_0(A)$ be a finitely generated subgroup with $[1_A]\in G,$ and
let
$J_0, J_1\ge 0$ be integers.

For any $\dt>0,$ any finite subset ${\cal G}\subset A,$ and any finite subset ${\cal P}\subset \underline{K}(A)$ with
{{$[1_A]\in {\cal P}$ and}} ${\cal P}\cap K_0(A)\subset G,$
there exist integers $N_0, N_1,...,N_k$ and
unital \hm s $h_j: A\to M_{N_j},$ $j=1,2,...,k,$ satisfying the following {{condition}}:

For any $\kappa\in Hom_{\Lambda}(\underline{K}(A), \underline{K}({\cal K})),$
with $|\kappa([1_A])|\le J_1$ {{(note that $K_0(\mathcal K) = \mathbb Z$)}} and
\beq\label{Ext1-1}
J_0=\max\{|\kappa(g_i)|: g_i=(\overbrace{0,...,0}^{i-1}, 1, 0,...,0)\in \Z^r: 1\le i\le r\},
\eneq
there exists a {{unital}} ${\cal G}$-$\dt$-multiplicative \cp\,
$\Phi: A\to M_{N_0+\kappa([1_A])}$
such that
\beq\label{Ext1-2}
[\Phi]|_{\cal P}=(\kappa+[h_1]+[h_2]+\cdots +[h_k])|_{\cal P}.
\eneq

(Note that, as $\Phi$ is unital, $N_0=\sum_{i=1}^kN_i.$ The notation $J_0$ is for later use.)

\end{lem}

\begin{proof}
It follows from 6.1.11 of \cite{Lnbok} (see also \cite{LinTAF2} and \cite{DE}) that, for each such $\kappa,$
there is a unital ${\cal G}$-$\dt$-multiplicative
\cp\, $L_{\kappa}: A\to M_{n(\kappa)}$ (for {{some}} integer $n(\kappa)\ge 1$) such that
\beq\label{Ext1-3}
[L_{\kappa}]|_{\cal P}=(\kappa+[h_{\kappa}])|_{\cal P},
\eneq
where $h_{\kappa}: A\to M_{N_{\kappa}}$ is a unital \hm.
There are only finitely many different $\kappa|_{\cal P}$ such that
(\ref{Ext1-1}) holds {{and $|\kappa([1_A])|\le J_1,$}} say $\kappa_1, \kappa_2,...,\kappa_k.$
Set $h_i=h_{\kappa_i},$ $i=1,2,...,k.$ Let $N_i=N_{\kappa_i},$ $i=1,2,...{{,}}k.$
Note that $N_i=\kappa_i([1_A])+n(\kappa_i),$ $i=1,2,...,k.$
Define
$$
N_0=\sum_{i=1}^k N_i.
$$
If $\kappa=\kappa_i,$ the map  $\Phi: A\to M_{N_0+\kappa([1_A])}$ defined by
$$
\Phi=L_{\kappa_i}+\sum_{j\not=i}h_j
$$
satisfies the requirements.
\end{proof}

\begin{lem}\label{Ext2}
Let $A$ be a unital \CA\, as in \ref{Ext1} and let {{$G=\Z^r\oplus {\rm Tor}(G)$ with}}
$[1_A]\in G$ {{be also as in \ref{Ext1}.}}
There exist $\Lambda_i\ge 0,$ $i=1,2,...,r,$ {{such that}} the following {{statement holds}}:
For any $\dt>0,$ any finite subset ${\cal G}\subset A,$ and any finite subset ${\cal P}\subset \underline{K}(A)$ with
{{$[1_A]\in {\cal P}$ and}}
${\cal P}\cap K_0(A)\subset G,$  there exist integers
$N(\dt, {\cal G}, {\cal P}, i)\ge 1,$  $i=1,2,...,r,$ satisfying the following {{condition}}:

Let $\kappa\in Hom_{\Lambda}(\underline{K}(A), \underline{K}({\cal K}))$
and
$
S_i=\kappa(g_i),
$
where $g_i=(\overbrace{0,...,0}^{i-1},1,0,...,0)\in \Z^r.$ {{There exist}}
a unital ${\cal G}$-$\dt$-multiplicative \cp\,
$L: A\to M_{N_1}$ and
a \hm\, $h: A\to M_{N_1}$ such that
\beq\label{EXt2-1}
[L]|_{\cal P}=(\kappa+[h])|_{\cal P},
\eneq
where $N_1=\sum_{i=1}^r (N(\dt, {\cal G},{\cal P}, i)+{{{\rm sign}(S_i)}}\cdot\Lambda_i)\cdot |S_i|$.
\end{lem}

\begin{proof}
Let $\psi_i^{+}: G\to \Z$ be {{the}} \hm\, defined by
$\psi_i^{+}(g_i)=1,$ $\psi_i^{+}(g_j)=0$ if $j\not=i,$ and
 $\psi_i^{+}|_{{\rm Tor}(G)}=0,$ and similarly  let $\psi_i^{-}(g_i)=-1,$  $\psi_i^{-}(g_j)=0$ if $j\not=i,$  and $\psi_i^{-}|_{{\rm Tor}(G)}=0,$  $i=1,2,...,r.$
 Note that $\psi_i^-=-\psi_i^+,$ $i=1,2,...,r.$
 Let $\Lambda_i=|\psi_i^+([1_A])|,$ $i=1,2,...,r.$

Let $\kappa_i^{+}, \kappa_i^-\in \mathrm{Hom}_{\Lambda}(\underline{K}(A), \underline{K}({\cal K}){{)}}$ be such that
$\kappa_i^+|_G=\psi_i^+$ and $\kappa_i^-|_G=\psi_i^-,$ $i=1,2,...,r.$
Let $N_0(i)\ge 1$ (in place of $N_0$) be {{as provided by}} \ref{Ext1} for
$\dt,$ ${\cal G},$ $J_0=1,$ and $J_1=\Lambda_i.$
Define
$N(\dt, {\cal G},{\cal P}, i)=N_0(i),$ $i=1,2,...,r.$

Let $\kappa\in \mathrm{Hom}_{\Lambda}(\underline{K}(A), \underline{K}({\cal K})).$
Then
$\kappa|_G=\sum_{i=1}^rS_i\psi_i^+,$
where $S_i=\kappa(g_i),$ $i=1,2,...,r.$
Applying  Lemma \ref{Ext1}, one obtains  unital ${\cal G}$-$\dt$-multiplicative
\cp s $L_i^{\pm}: A\to M_{N_0(i)+\kappa_i^{\pm}([1_A])}$  and \hm s $h_i^{\pm}: A\to M_{N_0(i)}$ such that
\beq\label{Ext2-3}
[L_i^{\pm}]|_{\cal P}=(\kappa_i^{\pm}+[h_i^{\pm}])|_{\cal P},\,\,\,i=1,2,...,r.
\eneq
Define
$L=\sum_{i=1}^r L_i^{{{{\rm sign}(S_i)}},|S_i|},$
where $L_{{i}}^{{{{\rm sign}(S_i)}}, |S_i|}: A\to M_{|S_i|N_0(i)}$ is defined by
$$L_{{i}}^{{{{\rm sign}(S_i)}}, |S_i|}(a)={\rm diag}(\overbrace{L_i^{{{\rm sign}(S_i)}}(a), ...,L_i^{{\rm sign}(S_i)}(a)}^{|S_i|})=L_i^{{\rm sign}(S_i)}(a)\otimes 1_{|S_i|}$$ for
all $a\in A.$
One checks that {{the map}}
$L: A\to M_{N_1},$ where $N_1=\sum_{i=1}^r|S_i|(\Lambda_i'+N(\dt, {\cal G}, {\cal P}, i))$ with $\Lambda_i'=\psi^+_i([1_A])$ if $S_i>0,$ or $\Lambda_i'=-\psi_i^+([1_A])$ if $S_i<0,$
is a
unital ${\cal G}$-$\dt$-multiplicative
\cp\, and
$$
[L]|_{\cal P}=(\kappa+[h])|_{\cal P}
$$
for some \hm\, $h: A\to M_{N_1}.$

\end{proof}

\begin{lem}\label{Ext3}
Let $A\in {\bar{{\cal D}}}_s$ be a unital \CA\, and let ${\cal P}\subset \underline{K}(A)$ be a finite subset.
{{Denote by}} $G\subset \underline{K}(A)$  the group generated by ${\cal P}$, and write
$G_1=G\cap K_1(A)=
\Z^r\oplus ({\rm Tor}(K_1(A)\cap G).$
Let ${\cal F}\subset A$ be a finite subset, let
$\ep>0$, and let $\Delta: A_+^{q, {\bf 1}}\setminus \{0\}\to (0,1)$ be an order preserving map.

There exist
$\dt>0,$  a finite subset ${\cal G}\subset A,$ a
 finite subset ${\cal H}\subset A_+^{{{\bf 1}}}\setminus \{0\}$, and an integer $N\ge 1$ satisfying the following {{condition}}: Let $\kappa\in KK(A\otimes C(\T), \mathbb C)$ and put
\begin{equation}\label{Ext3-2}
K=\max\{|\kappa(\boldsymbol{\bt}(g_i))|:1\le i\le r\}\,\,\,
{{{\text{(see\,\, \ref{Dbeta}\, for\,
the\, definition\, of}}\,\, {\boldsymbol{\bt}})}},
\end{equation}
where $g_i=(\overbrace{0,...,0}^{i-1},1,0,...,0)\in \Z^r.$ Then for any unital ${\cal G}$-$\dt$-multiplicative
{{\cp}}\,
$\phi: A\to M_R$ such that  $R\ge N(K+1)$ and
\begin{equation}\label{Ext3-1}
{\rm tr}\circ \phi(h)\ge \Delta(\hat{h})\tforal h\in {\cal H},
\end{equation}
{{where ${\rm tr}$ is the tracial state of $M_R,$}} there exists a
unitary $u\in M_R$ such that
\beq\label{EXt3-3}
&&\|[\phi(f),\, u]\|<\ep\tforal f\in {\cal F}
\tand\\
&&{\rm Bott}(\phi,\, u)|_{\cal P}=\kappa\circ {\boldsymbol{\bt}}|_{\cal P}.
\eneq
\end{lem}

\begin{proof}


To simplify notation, without loss of generality,  we may assume that ${\cal F}$ is a subset of the unit ball.
Let $\Delta_1=(1/8)\Delta$ and $\Delta_2=(1/16)\Delta.$

Let $\ep_0>0$ and let ${\cal G}_0\subset A$ be a finite subset {{satisfying}} the following condition:
If $\phi': A\to B$ (for any unital \CA\, $B$) is a unital ${\cal G}_0$-$\ep_0$-multiplicative
{{\cp}} and $u'\in B$ is a unitary such that
\beq\label{Ext3-n+1}
\|\phi'(g)u'-u'\phi'(g)\|<4\ep_0\tforal g\in {\cal G}_0,
\eneq
then ${\rm Bott}(\phi',\, u')|_{\cal P}$ is well defined {{(see \ref{Dbeta}).}} Moreover,
if $\phi': A\to B$ is another unital ${\cal G}_0$-$\ep_0$-multiplicative
{{\cp\, such that}}
\beq\label{Ext3-n-3}
\|\phi'(g)-\phi''(g)\|<4\ep_0\andeqn \|u'-u''\|<4\ep_0\tforal g\in {\cal G}_0,
\eneq
then
${\rm Bott}(\phi', \, u')|_{\cal P}={\rm Bott}(\phi'', \, u'')|_{\cal P}.$
We may assume that $1_A\in {\cal G}_0.$
Let
$$
{\cal G}_0'=\{g\otimes f: g\in {\cal G}_0\} \andeqn f=\{1_{C(\T)}, z, z^*\},
$$
where $z$ is the identity function on the unit circle $\T.$
We also assume that if $\Psi': A\otimes C(\T)\to C$ {{(}}for a unital \CA\, $C${{)}}
is a  unital ${\cal G}_0'$-$\ep_0$-multiplicative \cp,
then there exists a unitary $u'\in C$ such that
\beq\label{Ext3-n-4}
\|\Psi'(1\otimes z)-u'\|<4\ep_0.
\eneq

Without loss of generality,  we may assume that
${\cal G}_0$  is  contained in the unit ball of $A.$
Let $\ep_1=\min\{\ep/64, \ep_0/512\}$ and ${\cal F}_1={\cal F}\cup {\cal G}_0.$

Let ${\cal H}_0\subset A_+^{\bf 1}\setminus \{0\}$ (in place of ${\cal H}$) be the finite subset and  $L\ge 1$ the integer  provided by Lemma \ref{fullabs} for $\ep_1$ (in place of $\ep$) and ${\cal F}_1$ (in place of ${\cal F}$) as well as $\Delta_2$ (in place of $\Delta$).

Let ${\cal H}_1\subset A_+^{\bf 1}\setminus\{0\}$,
${\cal G}_1\subset A$ (in place of ${\cal G}$),
$\dt_1>0$ (in place of $\dt$), ${\cal P}_1\subset  \underline{K}(A)$
(in place of ${\cal P}$), ${\cal H}_2\subset A_{s.a.}$, and $1>\sigma>0$ be {{as}} provided by Theorem \ref{Newunique1}
for $\ep_1$ (in place of $\ep$),  ${\cal F}_1$ (in place of ${\cal F}$),
and $\Delta_1.$  We may assume that $[1_A]\in {\cal P}_2,$ ${\cal H}_2$ is in the unit ball of $A,$
and ${\cal H}_0\subset {\cal H}_1.$



Without loss of generality, we may assume that
$\dt_1<\ep_1/16, \sigma<\ep_1/16$, and ${\cal F}_1\subset {\cal G}_1.$
Put ${\cal P}_2={\cal P}\cup {\cal P}_1.$

{{Denote by $\{r_1, r_2,...,r_k\}$  the set of all ranks of
irreducible representations of $A.$}}
Fix an irreducible representation
$\pi_0: A\to M_{r_1}.$
Let $N(p)\ge 1$ (in place of $N({\cal P}_0))$ and ${\cal H}_0'\subset A_+^{\bf 1}\setminus \{0\}$ (in place of ${\cal H}$) {{denote the integer and}} finite subset {{provided}} by  Lemma \ref{Lfullab} for  $\{1_A\}$ (in place of ${\cal P}_0$) and $(1/16)\Delta.$ Let ${\cal H}_1'={\cal H}_1\cup{\cal H}_0'.$

 Let $G_0=G\cap K_0(A)$ and
write $G_0=\Z^{s_1}\oplus \Z^{s_2}\oplus {\rm Tor}(G_0),$ where $\Z^{s_2}\oplus {\rm Tor}(G_0)
\subset {\rm ker}\rho_A.$
Let $x_j=(\overbrace{0,...,0}^{j-1},1,0,...,0)\in \Z^{s_1}\oplus \Z^{s_2},$ $j=1,2,...,s_1+s_2.$ Note that
$A\otimes C(\T)\in {\overline{\cal D}}_s$ and
$A\otimes C(\T)$ {{has}} irreducible representations of ranks
$r_1, r_2,...,r_k.$
Let
$$
{\bar r}=\max \{ |(\pi_0)_{*0}(x_j)|: 0\le j\le s_1+s_2\}.
$$

Let ${\cal P}_3\subset \underline{K}(A\otimes C(\T))$ be a finite {{subset}}
containing ${\cal P}_2,$ $\{{\boldsymbol{\bt}}(g_j): 1\le j\le r\}$, and a finite subset which generates
${\boldsymbol{\bt}}({\rm Tor}(G_1)).$  Choose $\dt_2>0$ and a finite subset
$$
\overline{{\cal G}}=\{g\otimes f: g\in {\cal G}_2,\,\,\, f\in \{1, z, z^*\}\}
$$
in $A\otimes C(\T),$ where
${\cal G}_2\subset A$ is a finite subset
such that, for any unital $\overline{{\cal G}}$-$\dt_2$-multiplicative \cp\,
$\Phi': A\otimes C(\T)\to C$ {{(for any unital \CA\, $C$ with ${\rm Tor}(K_0(C))={\rm Tor}(K_1(C))=\{0\}$),}}
$[\Phi']|_{{\cal P}_3}$ is well defined and
\beq\label{Ext3-6}
[\Phi']|_{{\rm Tor}(G_0)\oplus \boldsymbol{\bt}({\rm Tor}(G_1))}=0.
\eneq
We may assume ${\cal G}_2\supset {\cal G}_1\cup {\cal F}_1.$

Let $\sigma_1=\min\{\Delta_2(\hat{h}): h\in {\cal H}_1'\}.$
Note $K_0(A\otimes C(\T))=K_0(A)\oplus \boldsymbol{\bt}(K_1(A))$ and
$\underline{K}(A\otimes C(\T))=\underline{K}(A)\oplus \boldsymbol{\bt}(\underline{K}(A)).$
Consider the subgroup of
$K_0(A\otimes C(\T)){{= K_0(A)\oplus {\boldsymbol{\bt}}(K_1(A))}}$ given by
$$
\Z^{s_1}\oplus \Z^{s_2}\oplus \Z^r\oplus {\rm Tor}(G_0)\oplus {\boldsymbol{\bt}}({\rm Tor}(G_1)).
$$

Let $\dt_3=\min\{\dt_1, \dt_2\}.$ Let $N(\dt_3, \overline{{\cal G}}, {\cal P}_3,i)$ and $\Lambda_i,$
$i=1,2,...,s_1+s_2+r,$ be as provided by Lemma \ref{Ext2} for $A\otimes C(\T)$.
Choose an integer $n_1 {{>}} N(p)$ such that
\beq\label{Ext3-7}
{(\sum_{i=1}^{s_1+s_2+r} N(\dt_3, \overline{{\cal G}}, {\cal P}_3,i)+1+\Lambda_i)N(p)\over{n_1-1}}<\min\{\sigma/16, \sigma_1/2\}.
\eneq
Choose $n>n_1$ such that
\beq\label{Ext3-7+1}
{n_1+2\over{n}}<\min\{\sigma/16, \sigma_1/2, 1/(L+1)\}.
\eneq

Let $\ep_2>0$ and let ${\cal F}_2\subset A$ be a finite subset
such that
$[\Psi]|_{{\cal P}_2}$ is well defined {{for any ${\cal F}_2$-$\ep_2$-multiplicative \morp\, $\Psi:
A\to B$ (for any \CA\, $B$).}}
Let $\ep_3=\min\{\ep_2/2, \ep_1\}$ and ${\cal F}_3={\cal F}_1\cup {\cal F}_2.$

Denote by $\dt_4>0$ (in place of $\dt$), ${\cal G}_3\subset A$
(in place of ${\cal G}$), ${\cal H}_3\subset A_+^{\bf 1}\setminus\{0\}$ (in place of ${\cal H}_2$) the {{constant and}} finite sets {{provided}}  by Lemma \ref{Combinerep} for $\ep_3$ (in place of $\ep$),
${\cal F}_3\cup {\cal H}_1'$ (in place of ${\cal F}$), $\dt_3/2$ (in place of $\ep_0$),
${\cal G}_2$ (in place of ${\cal G}_0$),
$\Delta,$ ${\cal H}_1'$ (in place of ${\cal H}$),
$\min\{\sigma/16, \sigma_1/2\}$  (in place of {{$\ep_1$}}), and $n^2$ (in place of $K$)
(with $L_1=L_2$ {{so no $\sigma$ is needed in \ref{Combinerep}}}).

Set ${\cal G}={\cal F}_3\cup {\cal G}_1\cup {\cal G}_2\cup {\cal G}_3$, set  $\dt=\min\{\ep_3/16, \dt_4, \dt_3/16\}{{,}}$
and
set ${\cal G}_5=\{g\otimes f: g\in {\cal G}_4,\,\,\,f\in\{1, z, z^*\}\}.$

Let ${\cal H}={\cal H}_1'\cup {\cal H}_3.$  Define $N_0=(n+1)^2N(p)(\sum_{i=1}^{s_1+s_2+r}N(\dt_3, {\cal G}_0, {\cal P}_3,i)+\Lambda_i+1)$ and define $N=N_0+N_0{\bar r}.$
Fix any $\kappa\in KK(A\otimes C(\T), \C)$
with
$$
K=\max\{|\kappa({\boldsymbol{\bt}}(g_i)|: 1\le j\le r\}.
$$
Let $R>N(K+1).$
Suppose
that $\phi:A\to M_R$ is a unital ${\cal G}$-$\dt$-multiplicative \cp\,
such that
\beq\label{Ext3-8}
{\rm tr}\circ \phi(h)\ge \Delta(\hat{h})\tforal h\in {\cal H}.
\eneq
Then, by \ref{Combinerep}, there exist mutually orthogonal projections
$e_0, e_1,e_2,...,e_{n^2}\in M_R$ such that $e_1, e_2,...,e_{n^2}$ are equivalent, {{$e_0\lesssim e_1,$}}
${\rm tr}(e_0)<\min\{\sigma/64, \sigma_1/4\}$  and $e_0+\sum_{i=1}^{n^2}e_i=1_{M_R},$ and there {{exist}} a unital ${\cal G}_2$-$\dt_3/2$-multiplicative \cp\,
$\psi_0: A\to e_0M_Re_0$ and a unital \hm\, $\psi: A\to e_1M_Re_1$ such that
\beq\label{Ext3-9}
&&\|\phi(f)-(\psi_0(f)\oplus \overbrace{\psi(f),\psi(f),...,\psi(f)}^{n^2})\|<\ep_3\tforal f\in {\cal F}_3\andeqn\\\label{Ext3-9+}
&&{\rm tr}\circ \psi(h)\ge \Delta(\hat{h})/{{2}}n^2\tforal h\in {{{\cal H}_1'.}}
\eneq
Let $\alpha\in \mathrm{Hom}_{\Lambda}(\underline{K}(A\otimes C(\T)), \underline{K}(M_r))$ be defined as follows:
$\alpha|_{\underline{K}(A)}=[\pi_0]$ and
$\alpha|_{\boldsymbol{\bt}(\underline{K}(A))}
=\kappa|_{\boldsymbol{\bt}(\underline{K}(A))}.$
Note that
$$
\max{{\{\max}}\{|\kappa\circ \boldsymbol{\bt}(g_i)|: {{1\le i\le r\}}},\,{{\max\{}}  |\pi_0(x_j)|:1\le j\le s_1+s_2\}\}
\le \max\{K, {\bar r}\}.
$$

Applying  Lemma \ref{Ext2}, we obtain a unital ${\cal G}$-$\dt_3$-multiplicative \cp\,
$\Psi: A\otimes C(\T)\to M_{N_1'},$ where $N_1'\le N_1:=\sum_{j=1}^{s_1+s_2+r} (N(\dt_3, {\cal G}_0, {\cal P}_3,j)+\Lambda_i)\max\{K, {\bar r}\},$
and a \hm\, $H_0: A\otimes C(\T)\to H_0(1_A)M_{N_1'}H_0(1_A)$ such that
\beq\label{Ext3-10-}
[\Psi]|_{{\cal P}_3}=(\af+[H_0])|_{{\cal P}_3}.
\eneq
{{Note that (since $H_0$ is a \hm\, with finite dimensional range)
\beq\label{Ext3-10+++}
[\Psi]|_{\boldsymbol{\bt}({\cal P})}=\kappa|_{\boldsymbol{\bt}(\cal P)}.
\eneq}}
In particular, since $[1_A]\in {\cal P}_2\subset {\cal P}_3,$
${\rm rank} {{(}}\Psi(1_A){{)}}=r_1+{\rm rank} (H_0).$
Note that
\beq\label{Ext3-10}
{N_1'+N(p) \over{R}}\le {N_1+N(p)\over{N(K+1)}}<1/(n+1)^2.
\eneq
Let $R_1={\rm rank}\,e_1.$
Then $R_1\ge R/(n+1)^2.$
Hence by (\ref{Ext3-10}),
$R_1\ge N_1+N(p).$
In other words, $R_1-N_1'\ge N(p)>0.$
Note that, {{by \eqref{Ext3-9+},}}
\begin{equation*}
{\rm tr}'\circ  \psi(\hat{g})\ge (1/3)\Delta(\hat{g})\ge \Delta_2(\hat{g})\rforal g\in {\cal H}_0',
\end{equation*}
where ${\rm tr}'$ is the tracial state on $M_{R_1}.$ Note that $n\ge N(p).$
Applying Lemma \ref{Lfullab} to the
{{pair}} $\pi_0\oplus H_0|_{A\otimes 1_{C(\T)}}$ (in place {{of}} $\phi$) {{and}} ${\tilde \psi}$ {{(in place of $\psi$),}}
where ${\tilde \psi}$ is an amplification of $\psi$ with $\psi$ repeated $n$ times, and ${\cal P}_0=\{[1_A]\},$
we obtain  a unital \hm\, $h_0: A\otimes C(\T)\to M_{nR_1-N_1'} $
such that $h_0(1\otimes {{1}})=1_{M_{nR_1-N_1'}}.$
Define $\psi_0': A\otimes C(\T)\to e_0M_Re_0$ {{by}}
$\psi_0'(a\otimes f)=\psi_0(a)\cdot (f(1)\cdot e_0)$ for all $a\in A$ and $f\in C(\T),$ where $1\in  \T.$
Define $\psi': A\otimes C(\T)\to e_1M_Re_1$ by
$\psi'(a\otimes f)=\psi(a)\cdot (f(1)\cdot e_1)$ for all $a\in A$ and $f\in C(\T).$
{{Note that
\beq\label{Extnn-1}
[\psi']|_{\boldsymbol{\bt}({\cal P})}=[\psi_0']|_{\boldsymbol{\bt}({\cal P})}=\{0\}.
\eneq}}
Put $E_1={\rm diag}(e_1,e_2,...,e_{nn_1}).$
Define $L_1: A\to E_1M_RE_1$ by
$L_1(a)=\pi_0(a)\oplus H_0|_{A}(a)\oplus h_0(a\otimes 1)\oplus (\overbrace{\psi(a),...,\psi(a)}^{n(n_1-1)})$ for $a\in A$, and define
$L_2: A\to E_1M_RE_1$ by
$L_2(a) =\Psi(a\otimes 1)\oplus h_0(a\otimes 1)\oplus (\overbrace{\psi(a),...,\psi(a)}^{n(n_1-1)})$ for $a\in A.$
Note that
\beq\label{Ext3-12}
&&[L_1]|_{{\cal P}_1}=[L_2]|_{{\cal P}_1},\\
&&{\rm tr}\circ L_1(h)\ge \Delta_1(\hat{h}),\,\,\,{\rm tr}\circ L_2(h)\ge \Delta_1(\hat{h})\tforal h\in {\cal H}_1, \andeqn\\
&&|{\rm tr}\circ L_1(g)-{\rm tr}\circ L_2(g)|<\sigma \tforal g\in {\cal H}_2.
\eneq
It follows from Theorem \ref{Newunique1} that
there exists a unitary $w_1\in E_1M_RE_1$ such that
\beq\label{Ext3-13}
\|{\rm {{Ad}} }\, w_1\circ L_2(a)-L_1(a)\|<\ep_1\tforal a\in {\cal F}_1.
\eneq
Define $E_2=(e_1+e_2+\cdots +e_{n^2})$ and define
$\Phi: A\to E_2M_RE_2$ by
\beq\label{Ext3-12+n-1}
\Phi(f)(a)={\rm diag}(\overbrace{\psi(a),\psi(a),...,\psi(a)}^{n^2})\tforal a\in A.
\eneq
Then
\beq\label{EXt3-12+n-2}
{\rm tr}\circ \Phi(h)\ge \Delta_2(\hat{h})\tforal h\in {\cal H}_0.
\eneq
By (\ref{Ext3-7+1}), one has $n/(n_1+2)>L+1.$ Applying Lemma \ref{fullabs}, we obtain a unitary $w_2\in E_2M_RE_2$ and
a unital \hm\, $H_1: A\to (E_2-E_1)M_R(E_2-E_1)$ such that
\beq\label{Ext3-14}
\|{\rm ad}\, w_2\circ {\rm diag}(L_1(a), H_1(a))-\Phi(a)\|<\ep_1\tforal a\in {\cal F}_1.
\eneq
Put
$$
w=(e_0\oplus w_1\oplus (E_2-E_1))(e_0\oplus w_2)\in M_R.
$$
Define $H_1': A\otimes C(\T)\to (E_2-E_1)M_R(E_2-E_1)$ by
$H_1'(a\otimes f)=H_1(a)\cdot f(1)\cdot (E_2-E_1)$ for all $a\in A$ and
$f\in C(\T).$
Define $\Psi_1: A\to M_R$ by
\beq\label{Ext3-11}
\Psi_1(f)=\psi_0'(f)\oplus \Psi(f)\oplus h_0(f)\oplus { (}\overbrace{\psi'(f),...,\psi'(f)}^{n(n_1-1)})\oplus H_1'(f)\tforal f\in A\otimes C(\T).
\eneq
{{By  \eqref{Ext3-10+++}, \eqref{Extnn-1},
\beq\label{Ext-RR}
[\Psi_1]|_{\boldsymbol{\bt}({\cal P})}=\kappa|_{\boldsymbol{\bt}(\cal P)}.
\eneq}}
It follows from (\ref{Ext3-13}), (\ref{Ext3-14}), and (\ref{Ext3-9}) that
\beq\label{Ext3-15}
\|\phi(a)-w^*\Psi_1(a\otimes 1)w\|<\ep_1+\ep_1+\ep_3\tforal a\in {\cal F}.
\eneq
Now pick a unitary $v\in M_R$ such that
\beq\label{Ext3-15+}
\|\Psi_1(1\otimes z)-v\|<4\ep_1.
\eneq
Put $u=w^*vw.$
Then, we estimate  that
\beq\label{Ext3-16}
\|[\phi(a), \, u]\|<\min\{\ep, \ep_0{{\}}}\tforal a\in {\cal F}_1.
\eneq
Moreover, by {{\eqref{Ext-RR} and by the choice of $\ep_0,$}}
one has
\beq\label{Ext3-17}
{\rm Bott}(\phi,\, u)|_{{\cal P}}={\rm Bott}({\Psi_1}|_A, \Psi_1(1\otimes z))
|_{\cal P}=[\Psi_1]\circ \boldsymbol{\bt}|_{\cal P}=
\kappa\circ \boldsymbol{\bt }|_{\cal P}.
\eneq

\end{proof}

\section{A Uniqueness Theorem for \CA s
 in $\mathcal D_s$}

The main goal of this section is to prove Theorem \ref{UniqN1}.
\begin{df}\label{1endpoints}
Let $A$ be a unital \CA,
$C=C(F_1, F_2, \phi_0, \phi_1)\in {\cal C},$  {{$\pi_i: C\to F_2$ ($i=0,1,$), and $\pi_e: C\to F_1$}} be as in Definition \ref{DfC1}.
Suppose that $L: A\to C$ is a \morp.
Define $L_e=\pi_e\circ L.$
Then
 $L_e: A\to F_1$ is a \morp\, such that
\beq\label{1end-1}
\phi_0\circ L_e=\pi_0\circ L\andeqn
\phi_1\circ L_e=\pi_1\circ L.
\eneq
Moreover, if $\dt>0$ and ${\cal G}\subset A$ and $L$ is ${\cal G}$-$\dt$-multiplicative, then $L_e$ is also ${\cal G}$-$\dt$-multiplicative.
\end{df}

\begin{lem}\label{endpoints}
Let $A$ be a unital \CA\, and let
$C=C(F_1, F_2, \phi_0, \phi_1) \in {\cal C}$ be as in  {\rm \ref{DfC1}}.
Let  $L_1, L_2: A\to C$ be   unital completely positive linear maps, let  $\ep>0,$ and let ${\cal F}\subset A$ be a {{finite}} subset.
 Suppose that there are unitaries $w_0\in \pi_0(C)\subset F_2$  and $w_1\in \pi_1(C)\subset F_2$ such that
 \beq\label{endpoints-1}
 &&\|w_0^*\pi_0\circ L_1(a)w_0-\pi_0\circ L_2(a)\|<\ep\quad {\textrm{and}}\\
 \label{endpoints-1+1}
 && \|w_1^*\pi_1\circ L_1(a)w_1-\pi_1\circ L_2(a)\|<\ep \tforal a\in {\cal F}.
 \eneq
 Then there exists a unitary $u\in F_1$ such that
 \beq\label{endpoints-2}
 &&\|\phi_0(u)^*\pi_0\circ L_1(a)\phi_0(u)-\pi_0\circ L_2(a)\|<\ep \quad {\textrm{and}}
\\\label{endpoints-2+1}
&&\|\phi_1(u)^*\pi_1\circ L_1(a)\phi_1(u)-\pi_1\circ L_2(a)\|<\ep\tforal a\in {\cal F}.
\eneq
\end{lem}

\begin{proof}
Write $F_1
=M_{n(1)}\oplus M_{n(2)}\oplus\cdots \oplus M_{n(k)}$ and
$F_2=M_{r(1)}\oplus M_{r(2)}\oplus\cdots \oplus M_{r(l)}.$
We may assume that ${\rm ker}\phi_0\cap {\rm ker}\phi_1=\{0\}$ {\blue{(see \ref{DfC1}).}}

We may assume that there are $k(0)$ and $k(1)$ such that $\phi_0|_{M_{{n(j)}}}$ is injective, ${j}=1,2,...,k(0),$ with $k(0)\le k,$
$\phi_0|_{M_{{n(j)}}}=0$ if $\blue{j}> k(0),$ and
$\phi_1|_{M_{{n(j)}}}$ is injective, ${j}=k(1), k(1)+1,...,k,$ with $k(1)\le k,$
$\phi_1|_{M_{{n(j)}}}=0,$ if ${j}<k(1).$
Write $F_{1,0}=\bigoplus_{{j}=1}^{k(0)}M_{n(j)}$ and $F_{1,1}=\bigoplus_{j=k(1)}^{k}M_{n(j)}.$
Note that $k(1)\le k(0)+1,$ {{and}} $\phi_0|_{F_{1,0}}$ and $\phi_1|_{F_{1,1}}$ are injective.
Note $\phi_0(F_{1,0})=\phi_0(F_1)=\pi_0(C)$ and $\phi_1(F_{1,1})=\phi_1(F_1)=\pi_1(C).$

For each fixed $a\in A$, since $L_i(a)\in C~(i={1, 2})$, there are elements
$$
g_{a,i}=g_{a,i,1}\oplus g_{a,i,2}\oplus \cd \oplus g_{a,i,k(0)}\oplus\cdots \oplus g_{a,i,k}\in F_1,
$$
such that $\phi_0(g_{a,i})=\pi_0\circ L_i(a)$ and $\phi_1(g_{a,i})=\pi_1\circ L_i(a),$ $i=1,2$,
 where $g_{a,i,j}\in M_{n(j)},$ $j=1,2,...,k$ and $i=1,2.$ Note that such $g_{a,i}$ is unique since ${\rm ker}\phi_0\cap {\rm ker}\phi_1=\{0\}.$
 Since $w_0\in \pi_0(C)=\phi_0(F_1)$, there is a unitary
 $$u_0=u_{0,1}\oplus u_{0,2} \oplus \cd \oplus u_{0, k(0)}\oplus \cd \oplus u_{0,k}$$
 such that $\phi_0(u_0)=w_0$. Note that the first $k(0)$ components of $u_0$ are uniquely determined by $w_0$ (since $\phi_0$ is injective on this part) and the components after {{the}} $k(0)$'th component can be {chosen} arbitrarily (since $\phi_0=0$ on this part). Similarly there exists
 $$u_1=u_{1,1}\oplus u_{1,2} \oplus \cd \oplus u_{1, k(1)}\oplus \cd \oplus u_{1,k}$$
 such that  $\phi_1(u_1)=w_1.$

 Now by  (\ref{endpoints-1}) and (\ref{endpoints-1+1}), we have
 \beq\label{endpoints-3}
 \|\phi_0(u_0)^*\phi_0(g_{a,1})\phi_0(u_0)-\phi_0(g_{a,2})\|<\ep\andeqn\\
\|\phi_1(u_1)^*\phi_1(g_{a,1})\phi_1(u_1)-\phi_1(g_{a,2}))\|<\ep\tforal a\in {\cal F}.
\eneq

Since $\phi_0$ is injective on $M_{{n_j}}$ for  ${j}\leq k(0)$   and $\phi_1$ is injective on $M_{{n(j)}}$ for ${j}>k(0)$ (note that we use $k(1)\leq k(0)+1$),
we have
\beq\label{endpoints-3+1}
 \|(u_{0,{j}})^*(g_{a,1,{j}})u_{0,{j}}-(g_{a,2,{j}})\|<\ep\quad \tforal {j}\leq k(0), \andeqn\\
\label{endpoints-3+2}
\|(u_{1,{j}})^*(g_{a,1,{j}})u_{1,{j}}-(g_{a,2,{j}})\|<\ep\quad \tforal {j}>k(0), \tforal a\in {\cal F}.
\eneq
Let $u=u_{0,1}\oplus \cd \oplus u_{0,k(0)}\oplus u_{1,k(0)+1} \oplus \cd \oplus u_{1,k} \in F_1$---that is, for the first $k(0)$ components of $u$, we  use $u_0$'s corresponding components,  and for the last $k-k(0)$ components of $u$, we use $u_1$'s. From (\ref{endpoints-3+1}) and  (\ref{endpoints-3+2}){,}
we have
$$\|u^*g_{a,1}u-g_{a,2}\|<\ep~~\tforal a\in {\cal F}.$$
{Once we apply} $\phi_0$ and $\phi_1$ to the above inequality, we get (\ref{endpoints-2}) and (\ref{endpoints-2+1}) as desired.
\end{proof}

{\blue{Let us very briefly describe the proof of  Theorem \ref{UniqN1}.
The key ingredients are Theorem \ref{Newunique1},  Lemma \ref{homotopy1}, and Lemma \ref{Ext3}.
First fix $\ep>0$ and a finite subset ${\cal F}\subset A.$
We note that the result has been established {{when}} the target algebra is finite dimensional (Theorem \ref{Newunique1}).
So we may reduce the general case to the case that the target algebra  {{is}} infinite dimensional, {and}} has only a single direct summand
(minimal) in ${\cal C}.$}}  {\blue{  So we write the target algebra $C$ {{as}} $A(F_1, F_2, h_0, h_1)$ and {{note}}
$\lambda: C\to C([0,1], F_2)$ {{as}}  given in \ref{DfC1} is injective {{(as $C$ is minimal)}}.

The first idea is to consider a partition $0=t_0<t_1<\cdots t_n=1$ such that \\
$\|\pi_t\circ \phi(f)(t)-\pi_{t_i}\circ \phi(f)(t_i)\|$ is very small if $t\in [t_i, t_{i+1}],$  where $\pi_{t_i}: C\to F_2$ is the  point evaluation of $C$ at $t_i.$
We will then consider {{the}}  pair $\phi_i=\pi_{t_i}\circ \phi$ and
$\psi_i=\pi_{t_i}\circ \psi.$}}
 {\blue{At each point $t_i,$ we apply \ref{Newunique1} to obtain {{the}}  unitary $w_i.$ We then connect these $w_i$ to obtain
the unitary in $C$ that we need to find.  In other words, the continuous path from $w_i$ to $w_{i+1}$ given by $u$
should change $\{\pi_t\circ \phi(f)(t): f\in {\cal F}, t\in [t_i, t_{i+1}]\}$  very little from
$w_i(\pi_{t}\circ \phi(f))w_i^*.$  One then observes that $w_iw_{i+1}^*$ almost commutes
with $\{\pi_{t_i}\circ \phi(f): f\in {\cal F}\}$ (and hence $\{\phi(f)(t): f\in {\cal F}, t\in [t_i, t_{i+1}]\}$).
 A basic homotopy lemma (such as Lemma \ref{homotopy1}) would provide a path $v(t)$ from
$w_iw_{i+1}^*$  to $1$ which also almost commutes with the set. Then one {{considers}}  $v(t)w_{i+1}$ which starts as $w_i$ and
ends at $w_{i+1}$ which would be a {{path as}} desired. However, the basic homotopy lemma may have  an obstacle, {{the class
${\boldsymbol{\rm Bott}}(\phi_i, w_iw_{i+1}^*$)}}. We will then consider $w_iz_i,$ where $z_i$ is as provided by Lemma \ref{Ext3},
so that $z_i$ almost commutes with $\{\phi(f)(t_i): f\in {\cal F}\}$ and {{the class}}  ${\boldsymbol{\rm Bott}}(\phi_i, w_iz_i(w_{i+1}z_{i+1})^*)$ will be zero.
This is possible because {{of}} the condition \eqref{Uni1-3+1}.  To simplify the process, we kill torsion elements in $\underline{K}(A)$
by repeating $\phi$ (and $\psi$) $N$ times.
Much of the proof is to make sure the idea can actually be carried out. One also needs to {{exercise}}   special care at the endpoints, applying Lemma \ref{endpoints}.}}

\begin{thm}\label{UniqN1}
Let $A\in {{{\bar{{\cal D}}}}}_s$ be a unital \CA\, with finitely generated $K_i(A)$ ($i=0,1$).
Let ${\cal F}\subset A$ be a finite subset, let
$\ep>0$ be a positive number, and let $\Delta: A_+^{q, {\bf 1}}\setminus \{0\}\to (0,1)$  be an order preserving map. There {{exist}}  a finite subset ${\cal H}_1\subset A_+^{\bf 1}\setminus \{0\},$
$\gamma_1>0,$ $\gamma_2>0,$ $\dt>0,$ a finite subset
${\cal G}\subset A$, a finite subset ${\cal P}\subset \underline{K}(A),$ a finite subset ${\cal H}_2\subset A$, a finite subset ${\cal U}\subset J_c(K_1(A))$ {\rm (see {{Definition}} \ref{Dcu})} for which $[{\cal U}]\subset {\cal P}$, and $N\in \mathbb N$ satisfying the following {{condition}} :
For any {pair of} unital ${\cal G}$-$\dt$-multiplicative  \cp s
$\phi, \psi: A\to C$,
for some $C\in {\cal C}$, such that
\beq\label{Uni1-1}
&&[\phi]|_{\cal P}=[\psi]|_{\cal P},
\\\label{Uni1-2}
&&\tau(\phi(a))\ge \Delta(a),\,\,\, \tau(\psi(a))\ge \Delta(a),\quad \textrm{for all $\tau\in T(C) \tand  a\in {\cal H}_1$},\\
\label{Uni1-3}
&&|\tau\circ \phi(a)-\tau\circ \psi(a)|<\gamma_1 \tforal a\in {\cal H}_2, \andeqn\\
\label{Uni1-3+1}
&&{\rm dist}(\phi^{\ddag}(u), \psi^{\ddag}(u))<\gamma_2 \tforal u\in {\cal U},
\eneq
there exists a unitary $W\in C\otimes M_{N}$ such that
\begin{equation}\label{Uni1-4}
\|W(\phi(f)\otimes 1_{M_{{N}}})W^*-(\psi(f)\otimes 1_{M_{N}})\|<\ep,\tforal f\in {\cal F}.
\end{equation}
\end{thm}

\begin{proof}
{{If $A$ is finite dimensional,}} {\blue{then it is semiprojective. So, with sufficiently large ${\cal G}$  and sufficiently small
$\dt,$ both $\phi$ and $\psi$ are close to \hm s to within $\ep/2$ and ${\cal F}.$ Therefore the  general case is reduced
to the case that both $\phi$ and $\psi$ are \hm s.}}
{{Since  two homomorphisms from $A$ to a \CA\, with stable rank one are unitarily equivalent if and only if they induce the same map on the ordered $K_0$ groups (see, for instance, Lemma 7.3.2 (ii) of \cite{RLL-Ktheory}), the theorem holds (with $N=1$ and no requirements on $\mathcal H_1$, $\mathcal H_2$, $\mathcal U$, $\gamma_1$, and $\gamma_2$).}}
So, {{it remains to consider the case}}  that $A$ is infinite dimensional.


Since $K_*(A)$ is finitely generated,  there is $n_0$ such that $\kappa\in\mathrm{Hom}_\Lambda(\underline{K}(A), \underline{K}(C))$ is determined by its restriction to $K(A, \Z/n\Z)$, $n=0,..., n_0$.
Set $N=n_0!$.


{Let ${\cal H}_1'\subset A_+\setminus \{0\}$ (in place of ${\cal H}$),}
 $\dt_1>0$ (in place of $\dt$), ${\cal G}_1\subset A$ (in place of ${\cal G}$),
 and ${\cal P}_0\subset \underline{K}(A)$ (in place of ${\cal P}$) be the finite subsets and constant be as provided by Lemma \ref{homotopy1}
for $\ep/32$ (in place of $\ep$), ${\cal F}$ and  $\Delta.$
We may assume that $\dt_1<{\min\{\ep/32, 1/64\}}$ and $(2\dt_1, {\cal G}_1)$
is a $KK$-pair (see the end of Definition \ref{KLtriple}).

Moreover, we may assume that $\dt_1$ is sufficiently small  that
if $\|uv-vu\|<3\dt_1,$ then the Exel formula
$$
\tau({\rm bott}_1(u,v))={1\over{2\pi\sqrt{-1}}}(\tau(\log(u^*vuv^*))
$$
holds for any pair of unitaries $u$ and $v$ in any unital \CA\, $C$ {{and}}
any $\tau\in T(C)$ {\blue{(see   Lemma 3.1 and 3.2 of \cite{Ex}, also, Theorem 3.7 of \cite{HL}).}}
{\blue{Furthermore
(see 2.11 of \cite{Linajm}
and 2.2 of \cite{ER}), we may assume that,
if $\|uv_i-v_iu\|<3\dt_1,$ $i=1,2,3,$ {{then}}
\beq\label{bott}
{\rm bott}_1(u, v_1v_2v_3)=\sum_{i=1}^3{\rm bott}_1(u,v_i).
\eneq
}}

{Fix a decomposition $K_1(A) = \mathbb Z^{k(A)} \oplus \mathrm{Tor}(K_1(A))$, where $k(A)$ is a positive integer. Choose $g_1, g_2,...,g_{k(A)}\in U(M_{m(A)}(A))$ (for the some integer $m(A)\ge 1$) such that
$\{\bar{g_1}, \bar{g_2},..., \bar{g}_{k(A)}\}\subset J_c(K_1(A))$ and
$[g_i]=(\overbrace{0,...,0}^{i-1},1,0,...,0)\in \Z^{k(A)}.$ {{Set}}  ${\cal U}=\{\bar{g_1}, \bar{g_2},..., \bar{g}_{k(A)}\}\subset
{\blue{J_c(K_1(A)).}}$
}


Let ${\cal U}_0\subset A$ be a finite subset
such that
$$
\{g_1, g_2,...,g_{k(A)}\}\subset\{(a_{i,j}): a_{i,j}\in {\cal U}_0\}.
$$
Let $\dt_u=\min\{1/256m(A)^2, \dt_1/{24}m(A)^2\},$
${\cal G}_u={\cal F}\, \cup\, {\cal G}_1\,\cup\, {\cal U}_0$ and
let ${\cal P}_u={\cal P}_0\,\cup\, \{[g_1],[g_2],...,[g_{k(A)}]\}.$
Let $\dt_2>0$ (in place of $\dt$), ${\cal G}_2\subset A$ (in place of ${\cal G}$), ${\cal H}_2'\subset A_+\setminus \{0\}$ (in place of ${\cal H}$), and $N_1\ge 1$ (in place of $N$) be the finite subsets and  constants {{provided}} by Lemma \ref{Ext3} 
for the data
$\dt_u$ (in place of $\ep$), ${\cal G}_u$ (in place of ${\cal F}$), ${\cal P}_u$ (in place of ${\cal P}$), and
$\Delta$, and with ${[g_j]}$ (in place of $g_j$), $j=1,2,...,k(A)$ (with $k(A)=r$).



  Let {${\blue{\ep_1}}=\min\{1/192N_1m(A)^2, \dt_u/2,  \dt_2/2m(A)^2\}$.}  






Let $\mathcal H_3'\subset A_+^{\bf 1}\setminus \{0\}$ (in place of $\mathcal H_1$), $\dt_4>0$ (in place of $\dt$),  ${\cal G}_3\subset A$  (in place of ${\cal G}$), ${\cal H}_4'\subset A_{s.a.}$ (in place of $\mathcal H_2$), ${\cal P}_1\subset \underline{K}(A)$ (in place of ${\cal P}$), and $\sigma>0$
be the finite subsets and constants {{provided}}  by Theorem \ref{Newunique1} with respect to $\ep_1/4$ (in place {{of}} $\ep$),
${\cal G}_u$ (in place of ${\cal F}$), and $\Delta$.


 Choose ${\cal H}_5'\subset A_+^{\bf 1} \setminus \{0\}$ and $\dt_5>0$ and a finite subset ${\cal G}_4\subset A$ such that, for any {$m\in\mathbb N$} and any unital ${\cal G}_4$-$\dt_5$-multiplicative \cp\,
 $L': A\to M_m,$ if
 ${\rm tr}\circ L'(h)>{\Delta(\hat{h})}\tforal h\in {\cal H}_5',$ then $m\ge {\blue{4N_1}}.$
 {\blue{This is  possible because we can apply  Lemma \ref{Combinerep} (taking $K=4N_1,$ $L_1=L_2=L'$), since we now assume
 that $A$ is infinite dimensional.}}

Put $\dt=\min\{\ep_1/16,  \dt_4/4m(A)^2, \dt_5/4m(A)^2\},$
${\cal G}={\mathcal G_2\, \cup}\, {\cal G}_u\,\cup {\cal G}_3\cup {\cal G}_4$, and  ${\cal P}={\cal P}_u\cup {\cal P}_1.$
%
Put
$${\cal H}_1={\cal H}_1'\cup {\cal H}_2'\cup {\cal H}_3'\cup {\cal H}_4'\cup {\cal H}_{{5}}'$$
and let ${\cal H}_2={\cal H}_4'.$
Let $\gamma_1=\sigma$ and let
$0<\gamma_2<\min\{{1/64m(A)^2N_1}, \dt_u/9m(A)^2, 1/256m(A)^2\}.$

{{We assume that $\dt$ is sufficiently small and $\mathcal G$ is sufficiently large that for {{any}} unital $\mathcal G$-$\delta$-multiplicative completely positive linear map ${\blue{\Phi}}: A\to B$, where $B$ is a unital \CA\,}}
{\blue{(so that $\Phi\otimes {\rm id}_{M_{m(A)}}$ is approximately multiplicative),  one has that}}
\begin{equation}\label{amm-unitary}
\|({\blue{\Phi}}\otimes\mathrm{id}_{M_{m(A)}})(g_j) - \langle({\blue{\Phi}}\otimes\mathrm{id}_{M_{m(A)}})(g_j)\rangle\| < \ep_1, \quad j=1, 2, ..., k(A).
\end{equation}

Now suppose that $C\in {\cal C}$ and $\phi, \psi: A\to C$ are two unital
${\cal G}$-$\dt$-multiplicative \cp s
satisfying the condition of the theorem for the given $\Delta,$ ${\cal H}_1,$ $\dt,$ ${\cal G},$ ${\cal P},$ ${\cal H}_2,$ $\gamma_1,$ $\gamma_2$, and ${\cal U}.$

We write $C=A(F_1, F_2, h_0, h_1),$
$F_1=M_{m_1}\oplus M_{m_2}\oplus \cdots \oplus M_{m_{F(1)}},$ and
$F_2=M_{n_1}\oplus M_{n_2}\oplus \cdots \oplus M_{n_{F(2)}}.$
For each $t\in [0,1],$ we will write $\pi_t: C\to C([0,1], F_2)$  {{for}} the point evaluation at $t$ as defined in \ref{DfC1}.
{\blue{Note that, when $C$ is finite dimensional, the theorem holds by Theorem \ref{Newunique1}.
So, we may assume that $C$ is infinite dimensional.  It is also clear that the  general case
can be reduced to the case that $C$ is minimal (see \ref{DfC1}).
As in Definition \ref{DfC1}, {{then}}, we may assume
that ${\rm ker}h_0\cap {\rm ker}h_1=\{0\}$, {{and}}
$\lambda: C\to C([0,1], F_2)$ defined by $\lambda(f,a)=f$ is unital and injective.
}}

Let
$$
0=t_0<t_1<\cdots <t_n=1
$$
be a partition of $[0,1]$ {{such}}  that
\beq\label{Uni-11}
\|\pi_{t}\circ \phi(g)-\pi_{t'}\circ \phi(g)\|<\ep_1/16\andeqn
\|\pi_{t}\circ \psi(g)-\pi_{t'}\circ \psi(g)\|<\ep_1/16
\eneq
for all $g\in {\cal G},$ provided $t, t'\in [t_{i-1}, t_i],$ $i=1,2,...,n.$

{{{{Set}}  $V_{i,j}=\langle \pi_{t_i}\circ \phi\otimes {\rm id}_{M_{m(A)}}(g_j) \rangle,$ $j=1,2,...,k(A)$ and $i=0,1,2,...,n.$}}

Applying Theorem \ref{Newunique1}, one obtains a unitary
$w_i\in F_2,$ if $0<i<n,$ $w_0\in h_0(F_1),$ if $i=0,$
and $w_n\in h_1(F_1),$ if $i=n,$
such that
\beq\label{Uni-12}
\|w_i\pi_{t_i}\circ \phi(g)w_i^*-\pi_{t_i}\circ \psi(g)\|<\ep_1/16\rforal g\in {\cal G}_u.
\eneq

It follows from  Lemma \ref{endpoints} that  we may assume that there is a unitary  $w_e\in F_1$ such that $h_0(w_e)=w_0$ and $h_1(w_e)=w_n.$ Since we assume that $\lambda$ is injective, we also have
\beq\label{Uni-12++}
\|w_e(\pi_e\circ \phi(g))w_e^*-\pi_e\circ\psi(g)\|<\ep_1/16\rforal g\in {\cal G}_u.
\eneq

By \eqref{Uni1-3+1}, {one has (see Definition \ref{Dcu}; {\blue{see also Proposition \ref{CUdist})}}
$$\mathrm{dist}(\overline{\langle\phi\otimes {\rm id}_{M_{m(A)}}(g_j^*)\rangle \langle (\psi\otimes {\rm id}_{M_{m(A)}})(g_j)\rangle}, \overline{1}) < {\blue{\gamma_2.}}$$}
{\blue{It follows that}} there is a unitary  ${\theta}_j\in M_{m(A)}(C)$
such that ${\theta}_j\in CU(M_{m(A)}(C))$ {(see Definition \ref{DLddag})} and
\beq\label{Uni-13}
\|\langle (\phi\otimes {\rm id}_{M_{m(A)}}(g_j^*)\rangle \langle (\psi\otimes {\rm id}_{M_{m(A)}})(g_j)\rangle  -{\theta}_j\|<\gamma_2,\,\,\,j=1,2,...,k(A).
\eneq
{By Theorem \ref{2Tg14}, {{we can}} write}
$$
{\theta}_j=\prod_{l=1}^{{4}}\exp(\sqrt{-1}a_j^{(l)})
$$
for self-adjoint {elements} $a_j^{(l)}\in M_{m(A)}(C),$
$l=1,2,...,{4},$ $j=1,2,...,k(A).$
Write
\beq\nonumber
a_j^{(l)}=(a_j^{(l,1)}, a_j^{(l,2)},...,a_j^{(l,{\blue{F(2)}})})\andeqn {\theta}_j=({\theta}_{j,1}, {\theta}_{j,2},...,{\theta}_{j,{\blue{F(2)}}})
\eneq
in $C([0,1], {M_{m(A)}} \otimes F_2)=C([0,1],M_{{m(A)}n_1})\oplus \cdots \oplus C([0,1], M_{{m(A)}n_{{F(2)}}}),$
where ${\theta}_{j,s}=\prod_{l=1}^{{4}}\exp(\sqrt{-1}a_j^{(l,s)}),$ $s=1,2,...,F(2).$
Then {{(see Lemma \ref{2Lg8})}}
$$
\sum_{l=1}^{{4}}{n_s({\rm tr}_{n_s}\otimes {\rm Tr}_{m(A)})(a_j^{(l,s)}(t))\over{2\pi}}\in \Z,\quad t\in [0,1],
$$
where ${\rm tr}_{n_s}$ is the normalized trace on $M_{n_s},$
$s=1,2,...,F(2).$
In particular,
\beq\label{Uni-15}
\sum_{l=1}^{{4}}n_s({\rm tr}_{n_s}\otimes {\rm Tr}_{m(A)})(a_j^{(l,s)}(t))=\sum_{l=1}^{{4}}n_s({\rm tr}_{n_s}\otimes {\rm Tr}_{m(A)})(a_j^{(l,s)}(t'))
\rforal t, t''\in [0,1].
\eneq

Let $W_i=w_i\otimes { 1}_{M_{m(A)}},$ $i={0, 1, ..., n}$
and $W_e=w_e\otimes { 1}_{\blue{M_{m(A)}}}.$
Then, {\blue{by \eqref{amm-unitary},  \eqref{Uni-12}, the choice of $\mathcal G_5,$  \eqref{amm-unitary} again,  \eqref{Uni-13}, and the choices of $\ep_1$ and $\gamma_2,$}}
\begin{eqnarray}
&&\hspace{-0.4in}{{\|\pi_{t_i}(\langle{(} \phi\otimes {\rm id}_{M_{m(A)}})(g_j^*)\rangle) W_i(\pi_{t_i}(\langle{(} \phi\otimes {\rm id}_{M_{m(A)}})(g_j)\rangle)W_i^*-{\theta}_j(t_i)\|}} \nonumber \\\nonumber
&<& {\blue{\ep_1+\|\pi_{t_i}(\langle(\phi\otimes {\rm id}_{M_{m(A)}})(g_j^*)\rangle) W_i(\pi_{t_i}( \phi\otimes {\rm id}_{M_{m(A)}})(g_j))W_i^*-{\theta}_j(t_i)\|}}\\\nonumber
&<& {\blue{\ep_1+m(A)^2\ep_1/16+\|\pi_{t_i}(\la (\phi\otimes {\rm id}_{M_{m(A)}})(g_j^*)\ra) \pi_{t_i}(\psi\otimes\mathrm{id}_{M_{m(A)}})(g_j)
-{\theta}_j(t_i)\|}}\\
&<& {\blue{\ep_1+m(A)^2\ep_1/16+\ep_1+ \|\pi_{t_i}(\la (\phi\otimes {\rm id}_{M_{m(A)}}(g_j^*)\ra)\pi_{t_i}( \langle (\psi\otimes {\rm id}_{M_{m(A)}})(g_j)\rangle)  -{\theta}_j(t_i) \|}} \nonumber \\
&<&{\ep_1 + m(A)^2\ep_1/16 + \ep_1 + \gamma_2 }
<{3m(A)^2\ep_1+\gamma_2}\  (\leq 1/32). \label{Ui-16}
\end{eqnarray}
{\blue{Similarly (but using \eqref{Uni-12++} instead of \eqref{Uni-12},}}  we also have (with $\phi_e=\pi_e\circ \phi$)
\beq\label{Uni-17}
\hspace{-0.3in}\|\langle (\phi_e\otimes {\rm id}_{M_{m(A)}})(g_j^*)\rangle W_e(\langle {(}\phi_e\otimes {\rm id}_{M_{m(A)}})(g_j)\rangle)W_e^*-\pi_e({\theta}_j)\|<3m(A)^2\ep_1+{\blue{\gamma_2}}\ {(\leq 1/32)}.
\eneq

It follows from {\eqref{Ui-16}} that, {for $j=1, ..., k(A)$ and $i=0, ..., n$}, there exist self-adjoint elements $b_{i,j}\in M_{m(A)}(F_2)$ and $b_{e,j}\in M_{m(A)}(F_1)$ such that
\beq\label{Uni-18}
&&\hspace{-0.4in}\exp(\sqrt{-1}b_{i,j})={\theta}_j(t_i)^*(\pi_{{t_i}}(\langle \phi\otimes {\rm id}_{M_{m(A)}})(g_j^*)\rangle)W_i(\pi_{{t_i}}(\langle \phi\otimes {\rm id}_{M_{m(A)}})(g_j)\rangle) W_i^*,
\\\label{Uni-19}
&&\hspace{-0.4in}\exp(\sqrt{-1}b_{e,j})=\pi_e({\theta}_j)^*(\pi_e(\langle \phi\otimes {\rm id}_{M_{m(A)}})(g_j^*)\rangle)W_e(\pi_e(\langle \phi\otimes {\rm id}_{M_{m(A)}})(g_j)\rangle) W_e^*,\text{and}
\\\label{Uni-20}
&&\|b_{i,j}\|<2\arcsin (3m(A)^2\ep_1/2+ {\gamma_2/2}),\,\,\,j=1,2,...,k(A),\,i={0, 1, ..., n}, e.
\eneq
{\blue{Note that (recall $h_0(w_e)=w_0$ and $h_1(w_e)=w_n$)}}
\beq\label{Uni-20+}
{{(h_0\otimes {\rm id}_{M(A)})}}(b_{e,j})=b_{0,j}\andeqn {{(h_1\otimes {\rm id}_{M(A)})}}(b_{e,j})=b_{n,j}.
\eneq
Write
\beq\nonumber
b_{i,j}=(b_{i,j}^{(1)},b_{i,j}^{(2)},...,b_{i,j}^{F(2)})\in {{M_{m(A)}}}(F_2)\andeqn
b_{e,j}=(b_{e,j}^{(1)}, b_{e,j}^{(2)},...,b_{e,j}^{(F(1))})\in {{M_{m(A)}}}(F_1).
\eneq
Note {\blue{also}} that, {\blue{for $i=0,1,...,n,e,$}}
\beq\label{Uni-21}
(\pi_{t_i}(\langle \phi\otimes {\rm id}_{M_{m(A)}}(g_j^*)\rangle))W_i(\pi_{t_i}(\langle \phi\otimes {\rm id}_{M_{m(A)}})(g_j)\rangle) W_i^*
=\pi_{t_i}({\theta}_j)\exp(\sqrt{-1}b_{i,j}),
\eneq
{\blue{$1\le j\le k(A).$}}
Then, {\blue{for $s=1,2,...,F(2),$ $j=1,2,...,k(A),$ and $i={0, 1, ..., n},$}}
\beq\label{Uni-22}
{n_s\over{2\pi}}({\rm tr}_{n_s}\otimes {\rm Tr}_{M_{m(A)}})(b_{i,j}^{(s)})\in \Z,
\eneq
where ${\rm tr}_{n_s}$ is the normalized trace on $M_{n_s}.$
We also have
\beq\label{Uni-23}
{m_s\over{2\pi}}({\rm tr}_{m_s}\otimes {\rm Tr}_{M_{m(A)}})(b_{e,j}^{(s)})\in \Z,
\eneq
where ${\rm tr}_{m_s}$ is the normalized trace on $M_{m_s},$ $s=1,2,...,F(1),$
$j=1,2,...,k(A).$
Put
\begin{equation}\label{defn-L1}
\lambda_{i,j}^{(s)}={n_s\over{2\pi}}({\rm tr}_{n_s}\otimes {\rm Tr}_{M_{m(A)}})(b_{i,j}^{(s)})\in \Z,
\end{equation}
$s=1, 2, ...,F(2)$, $j=1,2,...,k(A),$ and  $i=0, 1, ..., n$,
and
put
\begin{equation}\label{defn-L2}
\lambda_{e,j}^{(s)}={m_s\over{2\pi}}({\rm tr}_{m_s}\otimes {\rm Tr}_{M_{m(A)}})(b_{e,j}^{(s)})\in \Z,
\end{equation}
$s=1,2,...,F(1)$ and $j=1,2,...,k(A).$
{{Set}}
\begin{equation}\label{defn-Lambdas}
\lambda_{i,j}=(\lambda_{i,j}^{(1)}, \lambda_{i,j}^{(2)},...,\lambda_{i,j}^{(F(2))})\in \Z^{F(2)}\quad\mathrm{and}\quad
\lambda_{e,j}=(\lambda_{e,j}^{(1)},\lambda_{e,j}^{(2)},...,\lambda_{e,j}^{(F(1))})\in \Z^{F(1)}.
\end{equation}
{{
We have, by (\ref{Uni-20}) {\blue{and \eqref{defn-L1},}}
that, for $1\le s\le F(2)$, ${\blue{1\le j}}\le k(A),$ $0\le i\le n,$
\begin{eqnarray}\label{Uni-24}
|{\lambda_{i,j}^{(s)}\over{n_s}}|&<& \frac{1}{2\pi}m(A) \|b_{i, j}^{(s)}\|
< \frac{1}{\pi}m(A) \arcsin (3m(A)^2\ep_1/2+ {\gamma_2/2})  \nonumber \\
& < & \frac{1}{2\pi\cos{(3m(A)^2\ep_1+ \gamma_2)}}m(A) (3m(A)^2\ep_1+ {\gamma_2}) \nonumber \\
& < &\frac{1}{2\pi\cos{(1/32)}}(1/64N_1 + 1/64N_1) < 1/4N_1.
\end{eqnarray}}}
{\blue{Similarly (using  (\ref{Uni-20}) and \eqref{defn-L2}),}}
\begin{equation}\label{Uni-24+}
|{\lambda_{e,j}^{(s)}\over{m_s}}| < 1/4N_1,\,\,\,s=1,2,...,F(1), \,\,\,{\blue{j=1,2,...,k(A).}}
\end{equation}

{\blue{Note that $K_1(A)= \Z^{k(A)} \oplus \mathrm{Tor}(K_1(A)).$}}
Define $\af_i^{(0,1)}: K_1(A)\to \Z^{F(2)}$ {\blue{by
$\af_i^{(0,1)}|_{\mathrm{Tor}(K_1(A))}=0$ and by sending}} $[g_j]$ to $\lambda_{i,j},$ $j=1,2,...,k(A)$, $i=0,1,2,...,n,$ and
define
$\af_e^{(0,1)}: K_1(A)\to \Z^{F(1)}$ {\blue{by $\af_e^{(0,1)}|_{\mathrm{Tor}(K_1(A))}=0$ and}} by
 sending $[g_j]$ to $\lambda_{e,j},$ $j=1,2,...,k(A).$
We write $K_0(A\otimes C(\T))=K_0(A)\oplus {\boldsymbol{\bt}}(K_1(A)))$
(see   \ref{Dbeta}
for the definition
of ${\boldsymbol{\bt}}$).
Define $\af_i: K_*(A\otimes C(\T))\to K_*(F_2)$ as follows:
On $K_0(A\otimes C(\T)),$ define
\beq\label{Uni-25}
\af_i|_{K_0(A)}=[\pi_i\circ \phi]|_{K_0(A)},\,\,\,
\af_i|_{{\boldsymbol{\bt}}(K_1(A))}=\af_i\circ {\boldsymbol{\bt}}|_{K_1(A)}=\af_i^{(0,1)}, \nonumber
\eneq
and on $K_1(A\otimes C(\T)),$ define
\beq\label{Uni-25+}
\af_i|_{K_1(A\otimes C(\T))}=0,\,\,\,i=0,1,2,...,n. \nonumber
\eneq
Also define  $\af_e\in {\rm Hom}(K_*(A\otimes C(\T)), K_*(F_1)),$ by
\beq\label{Uni-26}
\af_e|_{K_0(A)}=[\pi_e\circ \phi]|_{K_0(A)},\,\,\,
\af_e|_{{\boldsymbol{\bt}}(K_1(A))}=\af_e\circ {\boldsymbol{\bt}}|_{K_1(A)}=\af_e^{(0,1)}
\eneq
on $K_0(A\otimes C(\T))$, and  $(\af_e)|_{K_1(A\otimes C(\T))}=0.$
Note that {\blue{(see \eqref{Uni-20+})}}
\beq\label{Uni-26+0}
(h_0)_{*}\circ \af_e=\af_0\andeqn (h_1)_{*}\circ\af_e=\af_n.
\eneq
Since $A\otimes C(\mathbb{T})$ satisfies the UCT,  the map $\alpha_e$ can be lifted to an element of $KK(A\otimes C(\mathbb T), F_1)$ which {{will}}  still {{be}}  denoted by $\alpha_e$. Then define
\beq\label{Uni-26+}
\af_0=\af_e\times [h_0] \andeqn \af_n=\af_e\times [h_1]
\eneq
in $KK(A\otimes C(\mathbb T), F_2)$.
For $i=1, ..., n-1$, also pick a lifting of $\alpha_i$ in $KK(A\otimes C(\mathbb T), F_2)$, and still denote it by $\alpha_i$.
{\blue{Combining \eqref{Uni-11} and \eqref{Uni-12},}} we {{compute}} that
\begin{equation}\label{pre-amc}
\|(w_{i}^*w_{i+1})\pi_{t_i}\circ \phi(g)-\pi_{t_i}\circ \phi(g)(w_{i}^*w_{i+1})\|<\ep_1/4\tforal g\in {\cal G}_u,
\end{equation}
$i=0,1,...,n-1.$

{Recall that} $V_{i,j}=\langle \pi_{t_i}\circ \phi\otimes {\rm id}_{M_{m(A)}}(g_j) \rangle,$ $j=1,2,...,k(A)$ and $i=0,1,2,...,n,$ {{and}}
$$
V_{i,j}=(V_{i,j, 1}, V_{i,j,2},..., V_{i,j,F(2)})\in M_{m(A)}(F_2),\,\,\, j=1,2,...,k(A),\,\,\, i=0,1,2,...,n.
$$
Similarly, write
$$
W_i=(W_{i,1},W_{i,2},...,W_{i,F(2)})\in M_{m(A)}(F_2),\,\,\, i=0,1,2,...,n.
$$
{By \eqref{pre-amc} and \eqref{amm-unitary}, we have,} {{for}}  {\blue{$i=0,1,...,n,$}}
\begin{equation}\label{Uni-27}
{\|W_i^*W_{i+1} V_{i, j} - V_{i, j} W_i^*W_{i+1}\|< m(A)^2\ep/4+2\ep_1}\le {{(9/4)m(A)^2\ep_1}}< \delta_1,
\end{equation}
and
{{by \eqref{Uni-11} and \eqref{amm-unitary}, for any $j=1,2,...,k(A)$, $i=0, 1, ..., n-1,$}}
\begin{equation}\label{cut-V}
{{ \|V_{i,j}-V_{i+1,j}\|< \ep_1 + \ep_1m(A)^2/16 + \ep_1\leq 3m^2(A)\delta_u/2 < \dt_1/12}}.
\end{equation}
{{Thus,}} {\blue{combining with \eqref{Uni-27},}}
{{\beq\nonumber
&&\hspace{-0.8in}\|W_{i}V_{i,j}^*W_{i}^*V_{i,j}V_{i+1,j}^*W_{i+1} V_{i+1,j}W_{i+1}^*-1\|\\\label{8418nW}
&\approx_{\delta_1/6}&  \|W_{i}V_{i,j}^*W_{i}^*V_{i,j}V_{i,j}^*W_{i+1} V_{i,j}W_{i+1}^*-1\|  \\\label{8418nW-1}
&=& \|W_{i}V_{i,j}^*W_{i}^*W_{i+1} V_{i,j}W_{i+1}^*-1\|
\approx_{\blue{\delta_1}} 0.
\eneq}}
{By \eqref{cut-V}}, there is a continuous path $Z(t)$ of unitaries  such that $Z(0)=V_{i,j}$, $Z(1)=V_{i+1,j},$  and
\begin{equation}\label{cont-discrete-1}
\|Z(t)-Z(1)\|<\dt_1/6,\quad t\in [0,1].
\end{equation}
We also write
$$
Z(t)=(Z_1(t), Z_2(t), ...,Z_{F(2)}(t))\in F_2\andeqn t\in [0,1].
$$
We obtain a continuous path
$$
W_{i}V_{i,j}^*W_{i}^*V_{i,j}Z(t)^*W_{i+1} Z(t)W_{i+1}^*
$$
which is in {$CU(M_{m(A)}(F_2))$}
for all $t\in  [0,1],$ and {\blue{by \eqref{cont-discrete-1}, and combining \eqref{8418nW} and \eqref{8418nW-1},}}
$$
\|W_{i}V_{i,j}^*W_{i}^*V_{i,j}Z(t)^*W_{i+1} Z(t)W_{i+1}^*-1\|< {\blue{2\dt_1/6+\dt_1/16+3\dt_1}} <1/8\tforal t\in [0,1].
$$
It follows that {{the integer}}
$$
(1/2\pi\sqrt{-1})({\rm tr}_{n_s}\otimes {\rm Tr}_{M_{m(A)}})[\log(W_{i,s}V_{i,j,s}^*W_{i,s}^*V_{i,j,s}Z_s(t)^*W_{i+1,s}Z_s(t)W_{i+1,s}^*)]
$$
{{is independent of $t\in [0, 1],$
 $s=1, 2, ..., F(2).$}}
In particular,
\beq\label{Uni-28}
&&\hspace{-0.6in}(1/2\pi\sqrt{-1})({\rm tr}_{n_s}\otimes {\rm Tr}_{M_{m(A)}})(\log(W_{i,s}V_{i,j,s}^*W_{i,s}^*W_{i+1,s} V_{i,j,s}W_{i+1,s}^*)) \nonumber \\\label{Uni-28+}
&&\hspace{-0.3in}=(1/2\pi\sqrt{-1})({\rm tr}_{n_s}\otimes {\rm Tr}_{M_{m(A)}})(\log(W_{i,s}V_{i,j,s}^*W_{i,s}^*V_{i, j}V_{i+1,j,s}^*W_{i +1,s}{V_{i+1,j,s}}W_{i+1,s}^*))
\eneq
{for $s=1, 2, ..., F(2)$}.
One also has {(by \eqref{Uni-18})}
\begin{eqnarray}\label{Uni-29}
&&W_{i}V_{i,j}^*W_{i}^*V_{i, j}V_{i+1,j}^*W_{i+1} V_{i+1,j}W_{i+1}^* \nonumber \\
&&=({\theta}_j(t_{i})\exp(\sqrt{-1}b_{i,j}))^*{\theta}_j(t_{i+1})
\exp(\sqrt{-1}b_{i+1,j}) \nonumber \\ \label{Uni-30}
&&=\exp(-\sqrt{-1}b_{i,j}){\theta}_j(t_{i})^*{\theta}_j(t_{i+1})
\exp(\sqrt{-1}b_{i+1,j}).
\end{eqnarray}
Note that, {\blue{by \eqref{Uni-13},
{\eqref{amm-unitary}},  \eqref{Uni-11}, and \eqref{Uni-13},}}  for $t\in [t_i, t_{i+1}],$ $1\le i\le n-1,$ $1\le j\le k(A),$
{\blue{
\beq\label{theta}\nonumber
\|{\theta}_j(t_{i})^*{\theta}_j(t)-1\|\approx_{\gm_2}\|\theta_j(t_i)^*(\pi_{t}(\langle (\phi\otimes {\rm id}_{M_{m(A)}}(g_j^*)\rangle \langle (\psi\otimes {\rm id}_{M_{m(A)}})(g_j)\rangle))-1\|\\\nonumber
\approx_{2\ep_1}\|\theta_j(t_i)^*\pi_{t}((\phi\otimes {\rm id}_{M_{m(A)}}(g_j^*)(\psi\otimes {\rm id}_{M_{m(A)}})(g_j))-1\|\\\nonumber
\approx_{m(A)^2\ep_1/8}\|\theta_j(t_i)^* \pi_{t_i}((\phi\otimes {\rm id}_{M_{m(A)}}(g_j^*)(\psi\otimes {\rm id}_{M_{m(A)}})(g_j))-1\|\\\nonumber
\approx_{2\ep_1}\|\theta_j(t_i)^* \pi_{t_i}(\la(\phi\otimes {\rm id}_{M_{m(A)}}(g_j^*)\ra\la(\psi\otimes {\rm id}_{M_{m(A)}})(g_j))\ra-1\|\\
\approx_{\gamma_2} \|\theta_j(t_i)^*\theta_j(t_i)-1\|=0.
\eneq}}
{\blue{Note that $m(A)^2\ep_1/8+2\ep_1+2\gamma_2<\dt_1<1/32.$}}
By Lemma 3.5 of \cite{Lin-AU11},
\begin{equation}\label{LeasyApp}
({\rm tr}_{n_s}\otimes {\rm Tr}_{m(A)})(\log({\theta}_{j,s}(t_{i})^*{\theta}_{j,s}(t_{i+1})))=0.
\end{equation}
It follows that (by the Exel formula (see  3.1 and 3.2 of \cite{Ex} and 3.7 of
{{\cite{HL}}}),
using (\ref{Uni-27}),
\eqref{Uni-28+}, (\ref{Uni-30}), \eqref{theta}, Lemma 2.11  of \cite{Linajm},
(\ref{LeasyApp}), {{and the choice of $\dt_1$}}), {{we have}}
\beq\label{Uni-31}
&&\hspace{-0.6in}(t\otimes {\rm Tr}_{m(A)})({\rm bott}_1(V_{i,j}, W_{i}^*W_{i+1})) \nonumber \\
 \hspace{-0.2in}&=&
({1\over{2\pi \sqrt{-1}}})(t\otimes {\rm Tr}_{m(A)})(\log(V_{i,j}^*W_{i}^*W_{i+1}V_{i,j}W_{i+1}^*W_{i})) \nonumber \\
 \hspace{-0.2in}&=&({1\over{2\pi \sqrt{-1}}})(t\otimes {\rm Tr}_{m(A)})(\log(W_{i}V_{i,j}^*W_{i}^*{W_{i+1}}V_{i,j}W_{i+1}^*)) \nonumber \\
&=&({1\over{2\pi \sqrt{-1}}})(t\otimes {\rm Tr}_{m(A)})(\log(W_{i}V_{i,j}^*W_{i}^*V_{i,j}V_{i+1,j}^*
W_{i+1}V_{i+1,j}W_{i+1}^*)) \nonumber \\
&=& ({1\over{2\pi \sqrt{-1}}})(t\otimes {\rm Tr}_{m(A)})(\log(\exp(-\sqrt{-1}b_{i,j}){\theta}_j(t_{i})^*
{\theta}_j(t_{i+1})\exp(\sqrt{-1}b_{i+1,j})) \nonumber \\
&=& ({1\over{2\pi \sqrt{-1}}})[(t\otimes {\rm Tr}_{m(A)})(-\sqrt{-1}b_{i,j})+(t\otimes {\rm Tr}_{m(A)})(\log({\theta}_j(t_{i})^*{\theta}_j(t_{i+1})) \nonumber \\
&&\hspace{1.4in}+(t\otimes {\rm Tr}_{m(A)})(\sqrt{-1}{b_{i+1,j}})] \nonumber \\
&=&{1\over{2\pi}}(t\otimes {\rm Tr}_{m(A)})(-b_{i,j}+b_{i+1,j}) \nonumber
\eneq
for all $t\in T(F_2).$
{By \eqref{defn-L1}, \eqref{defn-L2}, and \eqref{defn-Lambdas}, one has}
\beq\label{Uni-32}
{\rm bott}_1(V_{i,j}, W_{i}^*W_{i+1}))=-\lambda_{i,j}+\lambda_{i+1,j}
\eneq
$j=1,2,...,m(A),$  $i=0,1,...,  n-1.$

{{Note that $\phi$ is $\mathcal G_2$-$\delta_2$-multiplicative,}}
{\blue{and,  by \eqref{Uni1-2},  for $h\in\mathcal H_2' \cup {\cal H}_5'\subset {\cal H}_1,$
\beq\label{tr>D1}
&&(\mathrm{tr}\circ\pi_{t_i})(\phi(h)) \geq \Delta(\hat{h})\rforal {\rm tr}\in T(F_2)\\\label{tr>D2}
&& \andeqn  (\mathrm{tr}_e\circ\pi_e)\circ\phi(h))\ge \Delta(\hat{h}) \rforal {\rm tr}_e\in T(F_1).
\eneq}}
{\blue{Recall that these  inequalities {{imply}}  that $n_s\ge 4N_1$ {{and}} $m_s\ge 4N_1$ (by the choice of ${\cal H}_5'$).}}
{\blue{ Then, by
\eqref{Uni-24},
\beq\nonumber
n_s/N_1\ge n_s/2N_1+n_s/2N_1\ge 1+\max\{|\lambda_{i, j}^{(s)}|:  j=1, 2, ..., k(A)\}\\\label{Uni-32+++}
\ge (1+\max\{|\alpha_e(\boldsymbol{\bt}(g_j))|: j=1, 2, ..., k(A)\}).
\eneq}}
{\blue{Similarly (by \eqref{Uni-24+}),
\beq\label{Uni-32++}
m_s/N_1\ge  (1+\max\{|\alpha_e(\boldsymbol{\bt}(g_j))|: j=1, 2, ..., k(A)\}).
\eneq}}
Applying \ref{Ext3}
(using {\blue{\eqref{Uni-32+++}, \eqref{Uni-32++}, \eqref{tr>D1} and  \eqref{tr>D1}),}}  {{we obtain}} unitaries
$z_i\in F_2$, $i=1,2,...,n-1$, and $z_e\in F_1$ such that
\begin{eqnarray}
&&\|[z_i,\,\pi_{t_i}\circ\phi(g)]\|<\dt_u {{\tand \|[z_e,\, \pi_e\circ \phi(g)]\|<\dt_u}} \tforal g\in {\cal G}_u {,} \label{Uni-33} \\
&& {\rm Bott}(z_i, \pi_{t_i}\circ \phi)=\alpha_i,\andeqn
 {\rm Bott}(z_e,\pi_e\circ \phi)=\alpha_e.
\end{eqnarray}
Put $z_0=h_0(z_e)$ and $z_n=h_1(z_e).$ {{Recall that $\pi_{t_0}\circ \phi=h_0(\pi_e\circ \phi)$ and
$\pi_{t_n}\circ \phi=h_1(\pi_e\circ \phi).$}}
One verifies (by \eqref{Uni-26+}) that
\begin{eqnarray}\label{Uni-34}
{\rm Bott}(z_0, \pi_{t_0}\circ \phi)=\alpha_0\andeqn
{\rm Bott}(z_n, \pi_{t_n}\circ \phi)=\alpha_n.
\end{eqnarray}
Let $U_{i,i+1}={z_{i}w_{i}^*w_{i+1}z_{i+1}^*},$
$i=0,1,2,...,n-1.$
Then, {by \eqref{Uni-33}, \eqref{pre-amc}, and \eqref{Uni-11},
one has}
\beq\label{NT-8}
\|[U_{i,i+1}, \pi_{t_i}\circ \phi(g)]\|<2\ep_1+2\dt_u<\dt_1/2
\rforal g\in {\cal G}_u,\ i=0,1,2,...,n-1.
\eneq
Moreover, for $0\le i\le n-1$, $1\le j\le k(A)$, by \eqref{Uni-32} {\blue{(and by \eqref{bott}, and \eqref{Uni-33}, \eqref{Uni-27}),}}
\begin{eqnarray*}
{\rm bott}_1(U_{i,i+1}, \pi_{t_i}\circ \phi){([g_j])}&=&
{\rm bott}_1(z_{i},\, \pi_{t_i}\circ \phi){([g_j])}+{\rm bott}_1(w_{i}^*w_{i+1}, \pi_{t_i}\circ \phi){([g_j])}\\
&&+{\rm bott}_1(z_{i+1}^*,\, \pi_{t_i}\circ\phi){([g_j])}\\
&=&(\lambda_{i,j})+(-\lambda_{i,j}+\lambda_{i+1,j})+(-\lambda_{i+1,j})=0.
\end{eqnarray*}
Note that for any $x\in\bigoplus_{*=0, 1}\bigoplus_{k=1}^{n_0} K_*(A\otimes C(\mathbb T),\mathbb Z/k\mathbb Z)$, one has $N x=0$.
Therefore,
\begin{equation}\label{Uni-37}
{\rm Bott}((\underbrace{U_{i,i+1}, ..., U_{i, i+1}}_N), (\underbrace{\pi_{t_i}\circ \phi, ..., \pi_{t_i}\circ\phi}_N))|_{\cal P}=N{\rm Bott}(U_{i,i+1}, \pi_{t_i}\circ \phi)|_{\cal P} =0,
\end{equation}
$i=0, 1, 2,...,n-1.$
Note that, by the assumption \eqref{Uni1-2},
\beq\label{Uni-38}
t_{n_s}\circ \pi_t\circ \phi(h)\ge \Delta(\hat{h})\tforal h\in {\cal H}_1',
\eneq
where $t_{n_s}$ is the normalized trace on $M_{n_s},$ $1\le s\le F(2).$

By Lemma \ref{homotopy1} {and our choice of $\mathcal H_1'$, $\delta_1$, $\mathcal G_1$, $\mathcal P_0$,}  {{in view of}} \eqref{Uni-38}, \eqref{NT-8}, and \eqref{Uni-37}, there exists
a continuous path of unitaries, $\{\tilde{U}_{i,i+1}(t): t\in [t_i, t_{i+1}]\}\subset F_2\otimes M_N(\mathbb C),$ such that
\begin{equation}\label{Uni-39}
\tilde{U}_{i,i+1}(t_i)={1}_{F_2\otimes M_N(\mathbb C)},\,\,\, \tilde{U}_{i, i+1}(t_{i+1})=(z_iw_i^*w_{i+1}z^*_{i+1})\otimes 1_{M_N(\mathbb C)},
\end{equation}
and
\begin{equation}\label{Uni-39+}
\|\tilde{U}_{i, i+1}(t)(\underbrace{\pi_{t_i}\circ \phi(f), ..., {{\pi}}_{t_i}\circ\phi(f)}_N)\tilde{U}_{i, i+1}(t)^*-(\underbrace{\pi_{t_i}\circ \phi(f), ..., {{\pi}}_{t_i}\circ\phi(f)}_N)\|<\ep/32
\end{equation}
for all $f\in {\cal F}$ and for all $t\in [t_i, t_{i+1}].$ Define
$W=(W(t), \pi_e(W))\in C\otimes M_N$ by
\beq\label{Uni-40}
W(t)=(w_iz_i^*\otimes 1_{M_N})\tilde{U}_{i, i+1}(t)\tforal t\in [t_i, t_{i+1}],
\eneq
$i=0,1,...,n-1,$ {\blue{and $\pi_e(W)=w_ez_e^*\otimes 1_{M_N}.$}} Note that $W(t_i)=w_iz_i^*\otimes 1_{M_N},$ $i=0,1,...,n.$
Note also that
$$W(0)=w_0z_0^*\otimes 1_{M_N} =h_0(w_ez_e^*)\otimes 1_{M_N}$$ and
$$W(1)=w_nz_n^*\otimes 1_{M_N}=h_1(w_ez_e^*)\otimes 1_{M_N}.$$
So $W\in C\otimes M_N.$
One then checks that, by  (\ref{Uni-11}), (\ref{Uni-39+}), (\ref{Uni-33}){\blue{, (\ref{Uni-12}), and \eqref{Uni-11} again,}}
\begin{eqnarray*}
&&\hspace{-0.4in}\|W(t)((\pi_t\circ \phi)(f)\otimes 1_{M_N})W(t)^*-(\pi_t\circ \psi)(f)\otimes 1_{M_N}\|\\
&&<\|W(t)((\pi_t\circ \phi)(f)\otimes 1_{M_N})W(t)^*-W(t)((\pi_{t_i}\circ \phi)(f)\otimes 1_{M_N})W^*(t)\|\\
 &&\hspace{0.2in}+\|W(t)(\pi_{t_i}\circ \phi)(f)W(t)^*-W(t_i)(\pi_{t_i}\circ \phi)(f) W(t_i)^*\|\\
 &&\hspace{0.3in}+\|W(t_i)((\pi_{t_i}\circ \phi)(f)\otimes 1_{M_N}) W(t_i)^*-(w_i(\pi_{t_i}\circ \phi)(f)w_i^*)\otimes 1_{M_N}\|\\
 &&\hspace{0.4in}+\|w_i(\pi_{t_i}\circ \phi)(f)w_i^*-\pi_{t_i}\circ \psi(f)\|
+\|\pi_{t_i}\circ \psi(f)-\pi_t\circ \phi(f)\|\\
 &&<\ep_1/16+\ep/32+\dt_u+\ep_1/16+\ep_1/16<\ep
\end{eqnarray*}
 for all $f\in {\cal F}$ and for $t\in [t_i, t_{i+1}],$ {\blue{$i=0,1,...,n-1.$}}

 {\blue{Since $\lambda: C\to C([0,1], F_2)$ is assumed to be injective,
  {{this}} implies that
 \beq
 \|W(\phi(f)\otimes 1_{M_N})W^*-(\psi(f)\otimes 1_{M_N})\|<\ep\rforal f\in {\cal F}.
 \eneq}}
\end{proof}

\begin{rem}\label{Rem-on-N}
{Although it will not be needed in this paper, it is perhaps worth  {{pointing}} out that with} some  modification, the proof also works without assuming that $K_*(A)$ is  finitely generated.
In Theorem \ref{UniqN1}, the multiplicity $N$ only depends on $\underline{K}(A)$ as $\underline{K}(A)$ is finitely generated. However, if $K_*(A)$ is not finitely generated, the multiplicity $N$ {\blue{then will depend}} on $\mathcal F$ and $\ep$. {\blue{On the other hand,}}
 if $K_1(A)$ is torsion free,
or if $K_1(C)=0$, then the multiplicity $N$ can {always} be chosen to be $1$. {\blue{This  also will  not be needed here.}}
\end{rem}

\begin{cor}\label{CUniN1}
{ The statement of} Theorem \ref{UniqN1} holds if $A$ is replaced by $M_m(A)$ for any
integer $m\ge 1.$

\end{cor}







\section{\CA s in ${\cal B}_1$}

\begin{df}\label{DB1}
Let $A$ be a unital
simple \CA. We say $A\in {\cal B}_1$\index{${\cal B}_1$} if the following property holds:
Let $\ep>0,$ let $a\in A_+\setminus \{0\},$ and let ${\cal F}\subset A$ be a
finite subset.
There exist a non-zero projection $p\in A$ and a \SCA\, $C\in {\cal C}$
with $1_C=p$ such that
\begin{eqnarray}
&&\|xp-px\|<\ep\tforal x\in {\cal F},\\
&&{\rm dist}(pxp, C)<\ep\tforal x\in {\cal F},\andeqn \\
&&1-p\lesssim a.
\end{eqnarray}

If $C$ as above can always be chosen in $\mathcal C_0$, that is, with $K_{{1}}(C)=\{0\},$ then we say that $A\in {\cal B}_0.$
\index{${\cal B}_0$}
\end{df}

\begin{df}\label{DgTR}
Let $A$ be a unital simple \CA. We say $A$ has  generalized  tracial rank at most one,
\index{generalized tracial rank at most one}  if the following  {{property holds}}:

Let $\ep>0,$ let $a\in A_+\setminus \{0\}$ and let ${\cal F}\subset A$ be a
finite subset.
There exist a non-zero projection $p\in A$ and a  unital \SCA\, $C$ which is a subhomogeneous \CA, with  {{at most}}
one dimensional
spectrum, {{ in particular, }}
 a finite dimensional \CA\,
with $1_C=p$ such that
\begin{eqnarray}\label{Dgtr-1}
&&\|xp-px\|<\ep\tforal x\in {\cal F},\\\label{Dgtr-2}
&&{\rm dist}(pxp, C)<\ep\tforal x\in {\cal F},\andeqn \\\label{Dgtr-3}
&&1-p\lesssim a.
\end{eqnarray}

In this case, we write $gTR(A)\le 1.$\index{$gTR(A)\le 1$}
\end{df}

\begin{rem} It follows from  {\rm \ref{ASCAs}} that $gTR(A)\le 1$ if and only if $A\in {\cal B}_1.$
\end{rem}



Let ${\cal D}$ be a class of unital \CA s.
\begin{df}\label{DATD}
Let $A$ be a unital simple \CA. We say $A$ is tracially approximately in ${\cal D},$ denoted by $A\in\mathrm{TA}{\mathcal D}$, if
 the following property holds:

 For any $\ep>0,$ any $a\in A_+\setminus \{0\}$ and any finite subset ${\cal F}\subset A$,
there exist a non-zero projection $p\in A$ and a unital \SCA\, $C\in {\cal D},$
with $1_C=p$ such that
\begin{eqnarray}\label{DtrD-1}
&&\|xp-px\|<\ep\tforal x\in {\cal F}, \\\label{DtrD-2}
&&{\rm dist}(pxp, C)<\ep\tforal x\in {\cal F},\  \textrm{and} \\
&&1-p\lesssim a.
\end{eqnarray}
(see Definition 2.2 of \cite{EN-Tapprox}).
Note that $\mathcal B_0=\mathrm{TA}\mathcal C_0$ and $\mathcal B_1=\mathrm{TA}\mathcal C$.
If in the above definition, only (\ref{DtrD-1}) and (\ref{DtrD-2}) hold, then we say $A$ has the property
($L_{\cal D}$). 

{{ The property ($L_{\cal D}$) is a generalization of Popa's property in Theorem 1.2 of \cite{Popa} (also see
Definition 1.2 of \cite{Brown-Mem}  and
Definition 3.2 of \cite{LnTR}).}}

\end{df}


{{An earlier version of}} the following proposition  first appeared in an unpublished paper
of the second named author
distributed in 1998 {{(see  {{Corollary 6.4}} of \cite{LnTR} and 5.1 of \cite{LinTAI})}}.

\begin{prop}\label{Adset9201411}
Let $A$ be a unital simple \CA\, which has the property  {\rm ($L_{\cal D}$)}.
Then, for any $\ep>0$ and any finite subset ${\cal F}\subset A,$ there exist a projection
$p\in A$ and a \SCA\, $C\in {\cal D}$ with $1_C=p$ such that
\beq\label{Addb-1}
&& \|[x,\, p]\|<\ep\tforal f\in {\cal F}, \\
&& {\rm dist}(pxp, C)<\ep,\ \textrm{and} \\
&& \|pxp\|\ge \|x\|-\ep\tforal x\in {\cal F}.
\eneq
\end{prop}

\begin{proof}
Fix $\ep\in (0,1)$ and a finite subset ${\cal F}\subset A.$
\Wlog, we may assume  {{$1_A
\in {\cal F}.$ It follows from Proposition 2.2 of \cite{Bl-IND} that there is a unital separable simple \SCA\, $B\subset A$
which contains ${\cal F}.$  By Definition \ref{DATD}, there exist a sequence of \SCA\, $C_n\in {\cal D}$ and
a sequence of non-zero projections $p_n=1_{C_n}$ such that
\beq
&&\lim_{n\to\infty}\|p_nb-bp_n\|=0,\,\,\, \lim_{n\to\infty}{\rm dist}(p_nbp_n, C_n)=0,\andeqn\\
&& \lim_{n\to\infty}\|L_n(ab)-L_n(a)L_n(b)\|=0\rforal a,\in b\in B,
\eneq
where $L_n: B\to p_nBp_n$ is defined by $L_n(b)=p_nbp_n$ for all $b\in B.$
Consider  the map $L: B\to \prod_{n=1}^{\infty} p_nBp_n$ which is a unital \cp.
Let $\pi: \prod_{n=1}^{\infty} p_nBp_n\to \prod_{n=1}^{\infty} p_nBp_n/\bigoplus_{n=1}^{\infty}p_nBp_n$ be the quotient map.
Set $\phi=\pi\circ L.$ Then $\phi$ is a unital \hm. Since $B$ is a simple,
$\phi$ is an isometry.  It follows that there exists a subsequence $\{n_k\}$  such
that
\beq
\|p_{n_k}xp_{n_k}\|\ge \|x\|-\ep\rforal x\in {\cal F}.
\eneq
Note that we have $\lim_{k\to\infty}\|p_{n_k}x-xp_{n_k}\|=0$ and $\lim_{k\to\infty}{\rm dist}(p_{n_k}xp_{n_k}, C_{n_k})=0$
for all $x\in {\cal F}.$  Choosing $p:=p_{n_k}$ and $C=C_{n_k}$ for some sufficiently large $k,$ the conclusion of the lemma holds.}}

\end{proof}

\begin{thm}\label{B1sp}
Let $A$ be a unital simple separable \CA\, in {\blue{${\rm TA}{\cal D},$ where
${\cal D}$ is a class of unital \CA s.}}
Then either
$A$ {{is locally approximated by subalgebras in ${\cal D},$}}
or $A$  has the property (SP).
{\blue{In  {{the}} case that ${\cal D}$ is a class of semiprojective \CA s, then,
when $A$ does not have (SP), $A$ is an inductive limit of \CA s in ${\cal D}$ (with not necessarily injective maps).}}

\end{thm}

\begin{proof}
This follows from Definition \ref{DB1} immediately.  
Let ${\cal F}_1, {\cal F}_2,...,{\cal F}_n,...$ be a sequence of increasing
finite subsets of the unit ball of $A$ whose union is dense in the unit ball.   If $A$ does not have property (SP), then there is a non-zero positive element $a\in A$ such that $\overline{aAa}\not=A$ and
$\overline{aAa}$ has no non-zero projection.  Then, for each $n\ge 1,$  there is a projection
$1_A-p_n\lesssim a$ and a \SCA\, $C_n\in {{{\cal D}}}$
such that
$1_{C_n}=p_n$ and
\beq\label{B1sp-1}
\|p_nx-xp_n\|<{{1/n}}\andeqn {\rm dist}(p_nxp_n, C_n)<1/n \rforal x\in {\cal F}_n.
\eneq
Since $\overline{aAa}$ does not have any non-zero projection, one has $1_A-p_n=0.$
In other words,
$1_A=p_n$ and
\beq\label{B1sp-2}
{\rm dist}(x, C_n)<1/n\rforal x\in {\cal F}_n, \,n=1,2,...,
\eneq
as asserted.
{\blue{In  {{the}} case that  {{the}} \CA s in ${\cal D}$ are}} semiprojective,
$A$ is in fact an inductive limit of \CA s in {{$\mathcal D$ (with not necessarily injective maps).}}
\end{proof}

\begin{thm}\label{B1stablerk}
Let $A\in {\cal B}_1.$ Then $A$ has stable rank one.\index{stable rank one}
\end{thm}
\begin{proof}
This follows from {{Proposition \ref{2pg3}}} and Theorem 3.3 of \cite{Fan-sr1} {\blue{(see}} also 4.3 of \cite{EN-Tapprox}).
\end{proof}

\begin{lem}\label{MF}
Let ${\cal D}$ be a family of unital separable  \CA s which are  residually finite dimensional.
Any unital separable simple \CA\,
with property ($L_{\cal D}$)\index{property $L_{\cal D}$} can be embedded in $\prod M_{r(n)}/\bigoplus M_{r(n)}$ for some sequence of integers $\{r(n)\}.$
\end{lem}
\begin{proof}
Let $A$ be a unital separable simple \CA\,  with property ($L_{\cal D}$).
Let $\mathcal F_1\subset\mathcal F_2\subset \cdots\subset\mathcal F_i\subset\cdots$ be an increasing sequence of finite subsets of $A$ with  union dense in  $A$.
Since
$A$ has property ($L_{\cal D}$),
for each $n,$ there are a projection $p_n\in A$ and $C_n\subset A$ with $1_{C_n}=p_n$ and $C_n\in {\mathcal D}$ such that
\beq\label{MF-1}
\hspace{-0.1in}\| p_nf-fp_n \|<1/2^{n+2},\,\, \|p_nfp_n\|\ge \|f\|-1/2^{n+2}, \andeqn  p_nfp_n\in_{1/2^{n+2}} C_n
\eneq
{{  for all $f\in \mathcal F_n.$}}
For each $a\in {\cal F}_n,$ there exists $c(a)\in C_n$ such that
$\|p_nap_n-c(a)\|<1/2^{n+{{2}}}.$
  There are unital \hm s $\pi_n': C_n\to B_n,$ where $B_n$ is a finite dimensional \SCA\, such that
\beq\label{MF-4-}
&&\|\pi_n'(c(a))\|=\|{{c(a)}}\|\ge
{{\|p_nap_n\|-1/2^{n+2}}}
\\\label{MF-4}
&&{{
\ge \|a\|-(1/2^{n+2}+1/2^{n+2})=}} \|a\|-1/2^{n+1} \rforal \in {\cal F}_n,\,\,\,\,n=1,2,....
\eneq
There is an integer $r(n)\ge 1$ such that $B_n$ is unitally embedded into
$M_{r(n)}.$ Denote by $\pi_n: C_n\to M_{r(n)}$ the composition of $\pi_n'$ and the embedding.
Note $C_n\subset p_nAp_n.$  Then there is  {{a}}  unital
\cp\, $\Phi_n': p_nAp_n\to M_{r(n)}$ such that
\beq\label{MF-4+}
\Phi_n'|_{C_n}=\pi_n.
\eneq
Define $\Phi_n: A\to M_{r(n)}$ by $\Phi_n(a)=\Phi_n'(p_nap_n)$ for all $a\in A.$  It is a unital  \cp.
Moreover,
\beq\label{MF-2}
\|\Phi_n(p_nap_n)-\Phi_n(c(a))\|<1/2^{n+1}\rforal a\in {\cal F}_n,
\eneq
$n=1,2,...,.$
Combining with (\ref{MF-1}), we obtain that
\beq\label{MF-2ad1411}
 \|\Phi_n(f)\|\ge \|f\|-1/2^{n}\tforal f\in {\cal F}_n,\,\,\,\, n=1,2,....
\eneq
 Define $\Phi: A\to \prod_{n=1}^{\infty} M_{r(n)}$ by
$\Phi(a)=\{{\Phi_n(a)}\}$ for all $a\in A.$  Let
$$
\Pi: \prod_{n=1}^{\infty} M_{r(n)}\to \prod_{n=1}^{\infty} M_{r(n)}/\bigoplus_{n=1}^{\infty}M_{r(n)}
$$
be the quotient map.
Put $\Psi=\Pi\circ \Phi.$ One easily checks that $\Psi$ is in fact a unital \hm.
{{Since $A$ is simple,}} $\Psi$ is {{ a  monomorphism}}.
%
\end{proof}

\begin{thm}\label{B1hered}
Let $A\in {\cal B}_1$ (or $A\in {\cal B}_0$). Then, for any projection $p\in A,$ one has  $pAp\in {\cal B}_1$ (or $pAp\in {\cal B}_0$).
\end{thm}

\begin{proof}
{{Let us assume $p\not=0$.}}
Let $1/4>\ep>0,$ let $a\in (pAp)_+\setminus \{0\},$ and let ${\cal F}\subset pAp$ be a finite subset.
{{Without loss of generality, we may assume that $\|x\|\leq 1$ for all $x\in {\cal F}.$}}
Since $A$ is unital and simple, there are $x_1,x_2, ...,x_m\in A$ such that
\beq\label{B1hered-1-1}
\sum_{i=1}^m x_i^*px_i=1_A.
\eneq
{\blue{Since $x_i^*px_i\le 1_A,$ $\|px_i\|\le 1.$
By replacing $x_i$ by $px_i,$ we may assume that $\|x_i\|\le 1.$}}
Put ${\cal F}_1=\{p, x_1,x_2,...,x_m, x_1^*,x_2^*,...,x_m^*\}\cup {\cal F}.$
  Since $A\in {\cal B}_1,$ there is a projection $e\in A$ and a unital \SCA\, $C_1\in {\cal C}$ (or $C_1\in {\cal C}_0$)
with $1_{C_1}=e$ such that
\beq\label{B1hered-1}
\|xe-ex\|&<&\ep/64({{m}}+1)\tforal x\in {\cal F}_1,\\\label{B1here-1+1}
{\rm dist}(exe, C_1)&<&\ep/64({{m}}+1)\tforal x\in  {\cal F}_1,\andeqn\\
1-e &\lesssim & a.
\eneq
Put $\eta=\ep/64({{m}}+1).$ Then $0<\eta<1/2^8.$
Since $p\in {\cal F}_1,$ {{we have $\|epe-(epe)^2\|<\eta$}} {\blue{ and
$\|pep-(pep)^2\|<\eta.$  Moreover, there is $c(p)\in (C_1)_{s.a.}$ such that
\beq\label{B1hered-1+n2}
\|epe-c(p)\|<\eta.
\eneq}}
{\blue{One estimates that, since $0<\eta<1/2^8,$
\beq\nonumber
\sqrt{1-4\eta}>1-2\eta-4\eta^2=1-(2+4\eta)\eta \andeqn \sqrt{1-4\eta}<1-2\eta,\\
{\rm or}\,\, (1-\sqrt{1-4\eta})/2<(1+2\eta)\eta \andeqn (1+\sqrt{1-4\eta})/2>1-\eta.
\eneq}}
{\blue{One then  computes ${\rm sp}(epe),\, {\rm sp}(pep)\subset [0, (1+2\eta)\eta]\sqcup [1-\eta, 1].$
By the functional calculus, one obtains a projection $q_1\in pAp$ such that
\beq\label{B1hered-1+n1}
\|q_1-pep\|<(1+2\eta)\eta.
\eneq
If $\lambda\in ((1+2\eta)\eta, 1-\eta),$ then
\beq\label{B1here-1+2}
\|e-(\lambda e-epe)^{-1}(\lambda e-c(p))\| &\le& \|(\lambda e-epe)^{-1}\|\|(\lambda e-epe)-(\lambda e-c(p)\|\\\label{B1here-1+3}
&<&{\eta\over{\min\{\lambda-(1+2\eta)\eta, (1-\eta)-\lambda\}}}.
\eneq
It follows that $(\lambda e-epe)^{-1}(\lambda e-c(p)$ is invertible in $eAe$ when  the expression  in
\eqref{B1here-1+3}
is less than $1$ in which case $\lambda e-c(p)$ is invertible in $eAe$ (or in $C_1$).
Thus,
\beq
{\rm sp}(c(p))\subset [-\eta, (2+2\eta)\eta]\sqcup [1-2\eta, 1+\eta].
\eneq}}
{{Consequently}} {\blue{(by the functional calculus),}}  there is a projection
$q\in C_1$ such that
\beq\label{B1hered-18-1}
\|epe-q\|<{{(2+2\eta)\eta+\eta=(3+2\eta)\eta}}.
\eneq
Moreover, there are $y_1,y_2,...,y_m\in C_1$ {{such that $\|y_i-ex_ie\|<\eta$. Then}}
{\blue{
\beq
y_i^*qy_i\approx_{\eta(1+\eta)} y_i^*qex_ie\approx_{\eta} ex_i^*eqex_ie
\approx_{(3+2\eta)\eta} ex_i^*epex_ie
\approx_{2\eta} ex_i^*px_ie.
\eneq
}}
Therefore {{(by also \eqref{B1hered-1-1}),}}
\beq\label{B1hered-2+1}
\|\sum_{i=1}^m y_i^*qy_i-e\|<{{m(\eta(1+\eta)+(3+2\eta)\eta+3\eta)}}
{{=m\eta(7+3\eta)}}<\ep<1/4.
\eneq
{{Therefore}} $q$ is full in $C_1.$
It follows from \ref{cut-full-pj} that
$qC_1q\in {\cal C}$ (or $qC_1q\in {\cal C}_0$).
{{Note that
\beq\label{B1hered-18-2}
\|pep- epe\|\leq\|pep-epep\|+\|epep-epe\|< 2\eta.
\eneq}}
From \eqref{B1hered-1+n1}, \eqref{B1hered-18-2}, and \eqref{B1hered-18-1},
\beq\nonumber
\|q_1-q\|&<&\|q_1-pep\|+\|pep-epe\|+\|epe-q\|\\
 &<& (1+2\eta)\eta+2\eta+(3+2\eta)\eta=(6+4\eta)\eta<1.
\eneq
Hence, there is a unitary $u\in A$ such that
$
u^*qu=q_1\le p.
$
Put
$C=u^*qC_1qu.$ Then $C\in {\cal C}$ (or $C\in {\cal C}_0$) and $1_C=q_1.$
We also have, {{by  \eqref{B1hered-1+n1} and \eqref{B1hered-18-1}),}}
\beq\label{B1hered-3=}
\|epe-q_1\|<{{(1+2\eta)\eta+(3+2\eta)\eta=(4+4\eta)\eta}}.
\eneq
If $x\in {\cal F},$
then
\beq\label{B1hered-5}
\|q_1x-xq_1\| &\le& 2\|(q_1-epe)x\|+\|epex-xepe\|\\
&<& {{2(4+4\eta)\eta+3\eta=(11+8\eta)\eta}}<\ep
\tforal x\in {\cal F}.
\eneq
Similarly, we estimate that
\beq\label{B1hered-4}
{\rm dist}(q_1xq_1, C)<\ep\tforal x\in {\cal F}.
\eneq
We also have {{(by \eqref{B1hered-1+n1})}}
\beq\label{B1hered-5}
\|(p-q_1)-(p-pep)\|=\|q_1-pep\|<{{(1+2\eta)\eta.}}
\eneq
Put $\sigma=(1+2\eta)\eta<1/16.$
Let $f_{\sigma}(t)\in C_0((0,\infty))$ be as in \ref{Dball}.
Then, by 2.2 of \cite{RorUHF2},
\beq\label{B1hered-6}
p-q_1=f_{\sigma}(p-q_1)\lesssim p-pep\lesssim 1-e\lesssim a.
\eneq
This shows that $pAp\in {\cal B}_1.$
\end{proof}

\begin{prop}\label{Ltrace}
{{Let ${\cal D}$ denote the class of unital separable amenable  \CA s with faithful tracial states.
Let $A$ be a non-zero unital  simple separable  \CA\, which is $TA{\cal D}.$
Then $QT(A)=T(A)\not=\emptyset.$}}
\end{prop}

\begin{proof}
{{One may assume that $A$ is infinite dimensional.
Since $A$ is a unital infinite dimensional simple \CA,
there are
$n$ mutually orthogonal non-zero positive elements, for any integer $n\ge 1$
(see, for example, 1.11.45 of \cite{Lnbok}).
By repeatedly applying Lemma 3.5.4 of \cite{Lnbok} (see also Lemma 2.3 of \cite{Lncons}),  one finds  a sequence of positive elements
$\{b_n\}$ which has the following property:
$b_{n+1}\lesssim b_{n,1},$
where $b_{n,1}, b_{n,2},...,b_{n,n}$ are mutually orthogonal non-zero positive elements
in $\overline{b_nAb_n}$ such that $b_nb_{n,i}=b_{n,i}b_n=b_{n,i},$
$i=1,2,...,n,$ and $ b_{n,i}\sim  b_{n,1},$ $1\le i\le n,$ $n=1,2,....$
Note that
\beq\label{99-+}
\lim_{n\to\infty}\sup\{\tau(b_n): \tau\in QT(A)\}=0.
\eneq
}}
{{One obtains two sequences of unital \SCA s $A_{0,n}:=e_nAe_n,$ $D_n$ of $A,$ where
$D_n\in {\cal D}$ with $1_{D_n}=(1-e_n),$   two sequences
of unital \cp s
$\phi_{0,n}: A\to A_{0,n}$ (defined by $\phi_{0,n}(a)=e_nae_n$ for all $a\in A$)
and
$\phi_{1,n}: A\to D_n$  with $A_{0,n}\perp D_n$ satisfying the following conditions:
\beq\label{TD0qc-1}
\lim_{n\to\infty}\|\phi_{i,n}(ab)-\phi_{i,n}(a)\phi_{i,n}(b)\|=0\rforal a,\, b\in A,\\\label{TD0qc-1+}
\lim_{n\to\infty}\|a-(\phi_{0,n}+\phi_{1,n})(a)\|=0\rforal a\in A,\\\label{TD0qc-1++}
(1-e_n)\lesssim b_n,\,\text{and}\\\label{TD0qc-4}
\lim_{n\to\infty}\|\phi_{1,n}(x)\|=\|x\|\rforal x\in A.
\eneq
}}
{{
 Since quasitraces are norm continuous (Corollary II 2.5 of \cite{BH}),
 by \eqref{TD0qc-1+},
 \beq\label{TD0qc-4+3}
\lim_{n\to\infty}(\sup|\tau(a)-\tau((\phi_{0,n}+\phi_{1,n})(a))|:\tau\in QT(A)\})=0\rforal a\in A.
\eneq
Since $\phi_{0,n}(a)\phi_{1,n}(a)=\phi_{1,n}(a)\phi_{0,n}(a)=0,$ for any $\tau\in QT(A)$  we have
\beq
\tau((\phi_{0,n}+\phi_{1,n})(a))=\tau(\phi_{0,n}(a))+\tau(\phi_{1,n}(a))\rforal a\in A.
\eneq
Note that, by \eqref{TD0qc-1++} and \eqref{99-+},
\beq\label{TD0qc-4++}
\lim_{n\to\infty}\sup\{\tau(\phi_{0,n}(a)): \tau\in QT(A)\}=0\rforal a\in A.
\eneq
Therefore
\beq\label{TD0qc-10}
\lim_{n\to\infty}(\sup\{|\tau(a)-\tau\circ \phi_{1,n}(a)|: \tau\in {\mathrm{QT}(A)}\})=0\rforal a\in A.
\eneq
{{It follows}}  $\lim_{n\to\infty}\|\tau\circ \phi_{1,n}\|=\|\tau\|=1$ for all $\tau\in {\mathrm{QT}(A)}.$}}

{{For any $\tau\in QT(A),$  let
$t_n=(\|\tau\circ \phi_{1,n}\|^{-1})\tau\circ\phi_{1,n}.$  Note  that $t_n|_{D_n}$ is a tracial state.
Therefore  $t_n$ is a state of $A.$
Let $t_0$ be a weak* limit of $\{t_n\}.$
Then, as $A$ is unital,  $t_0$ is a state of $A.$}}

{{It follows from \eqref{TD0qc-1}  and the fact that $t_n|_{D_n}$ is a tracial state that $t_0$ is a trace. Then, by
\eqref{TD0qc-10},  for every $\tau\in {\mathrm{QT}(A)},$
\beq\nonumber
\tau(a+b)=\tau(a)+\tau(b)\rforal a, b\in A.
\eneq
It follows that $\tau$ is a trace.  Therefore $QT(A)=T(A).$}}

{{To see that $T(A)\not=\emptyset,$  let $\tau_n\in T(D_n)$ be a faithful tracial state, $n=1,2,....$
Define $t_n=(\|\tau_n\circ \phi_{1,n}\|^{-1})\tau_n\circ \phi_{1,n}.$ Let $t_0$ be a weak * limit of $\{t_n\}.$
As above, $t_0$ is a tracial state of $A.$}}

\end{proof}

\begin{thm}\label{Comparison}
Let ${\cal D}$ be a class of {{unital}} \CA s  which is closed under tensor products with a finite dimensional \CA\, and which  has
 the strict comparison property for positive elements (see \ref{DW(A)}).
 Let $A$ be a unital simple {{separable}} \CA\, in the class $\mathrm{TA}{\cal D}.$  Then $A$ has strict comparison for positive elements. In particular, if $A\in\mathcal B_1$, then $A$ has strict comparison for positive elements and
$K_0(A)$ is weakly unperforated.
\end{thm}

%

\begin{proof}

{{Note that from \ref{Ltrace}, we have $QT(A)=T(A)\not=\emptyset$.}}
By a result of R\o rdam (see, for example,  Corollary 4.6 of \cite{Ror-srZ}; {{note also  that the exactness is only used to get $QT(A)=T(A)\not=\emptyset$ there}}), to show that $A$ has strict comparison for positive elements, it is enough to show that $W(A)$ is almost unperforated, i.e., for any positive elements $a, b$ in a matrix algebra over $A$, if $(n+1) [a] \leq n [b]$ for some $n\in\mathbb N$, then $[a]\le [b]$.

Let $a, b$ be such positive elements. Since any matrix algebra over $A$ is still in $\mathrm{TA}{\mathcal D}$, let us assume that $a, b\in A$.

First we consider the case that $A$ does not have (SP) property. In this case, by the proof of  \ref{B1sp},
$A=\overline{\bigcup_{n=1}^{\infty} A_n},$ where $A_n\in {\cal D}$ ($\{A_n\}$ {{may}} not be increasing).

Without loss of generality, we may assume that $0\le a, \, b\le 1.$
 Let $\ep>0.$ It follows from an argument of R\o rdam (see Lemma 5.6 of \cite{Niu-MD}) that there
 {{exist}} an integer $m\ge 1,$ {{and positive elements}}
 $a', b'\in A_m$ such that
 \beq\label{Compr-n-1}
 \|a'-a\|<\ep/2,\,\,\, \|b'-b\|<\ep/2,\,\,\, b'\lesssim b \andeqn\\
 {\rm diag}(\overbrace{f_{\ep/2}(a'), f_{\ep/2}(a'),...,f_{\ep/2}(a')}^{n+1})\lesssim
 {\rm diag}(\overbrace{b',b',...,b'}^n)\,\,\,{\rm in}\,\,\, A_m.
 \eneq
Since $A_m$ has strict comparison (see part (b) of  Theorem \ref{2Tg16}), one has
\beq\label{Compr-n-2}
f_{\ep/2}(a')\lesssim b'\,\,\,{\rm in}\,\,\, A_m.
\eneq
{\blue{Note that $\|a-(a'-\ep/2)_+\|\le \|a-a'\|+\|a'-(a'-\ep/2)\|<\ep.$}}
It follows,  on using 2.2 of \cite{RorUHF2}, that
\beq\label{Compri-n-3}
f_{2\ep}(a)\lesssim {{(a'-\ep/2)_+}}\lesssim f_{\ep/2}(a')\lesssim b'\lesssim b
\eneq
for every $\ep>0.$
It follows that $a\lesssim b.$

Now we assume that $A$ has (SP).
Let $1/4>\ep>0.$  We may further assume that $\|b\|=1.$
Since $A$ has (SP) and is simple, { { by Lemma 3.5.6 and Lemma 3.5.7 of \cite{Lnbok}
(also see Theorem I of \cite{Zh}),}}
 there are mutually orthogonal and mutually equivalent non-zero projections
$e_1,e_2,...,e_{n+1} \in \overline{f_{3/4}(b)Af_{3/4}(b)}.$
Put $E=e_1+e_2+\cdots +e_{n+1}.$
By 2.4 of \cite{RorUHF2}, we also have that
\beq\label{Compr-n2}
(n+1)[f_{{\ep/4}}(a)]\le n[f_{\dt}(b)]
\eneq
for some $\ep>\dt>0.$  Put $0<\eta<\min\{\ep/4, \dt/4, 1/8\}.$
It follows from  Definition \ref{DB1} that there are a \SCA\, $C=pAp\oplus S$ with $S\in\mathcal D$ and
$a', b', E', e_i'\in C$ ($i=1,2,...,n+1$) such that
$0\le a', b'\le 1$ and $E', e_i'$ are projections in $C$,
\beq\label{Compri-n3}
||a-a'||<\eta,\  b'\lesssim f_{\dt}(b), \|f_{1/2}(b') E'-E'\|<\eta,\\
E'=\sum_{i=1}^{n+1}e_i',\,\,\,\|e_i-e_i'\|<\eta \andeqn  \|E-E'\|<\eta<1,
\eneq
and
\begin{equation}\label{Compri-n4}
{\rm diag}(\overbrace{f_{\ep/2}(a'),f_{\ep/2}(a'),...,f_{\ep/2}(a')}^{n+1}) \lesssim
{\rm diag}(\overbrace{b',b',...,b'}^n)\,\,\,
 \textrm{in $C$}
\end{equation}
(see Lemma 5.6 of \cite{Niu-MD}).
Moreover, the projection $p$ can be chosen so that $p\lesssim e_1$.
From (\ref{Compri-n3}), there is a projection ${e_i''}$, ${E}''\in \overline{f_{1/2}(b')Cf_{1/2}(b')}$
($i=1,2,...,n+1$) such that
$\|E'-E''\|<2\eta,$  $\|e_i''-e_i'\|<2\eta,$ $i=1,2,...,n+1,$ and
$E''=\sum_{i=1}^{n+1}e_i''$ (we also assume that $e_1'',e_2'',...,e_{n+1}''$ are mutually orthogonal).
Note that $e_i'$ and
$e_i''$ are equivalent.
Choose a function $g\in C_0((0,1])_+$ with $g\le 1$ such that $g(b')f_{1/2}(b')=f_{1/2}(b')$
and $[g(b')]=[b']$ in $W(C).$ In particular,
$g(b')E''=E''.$

Write
$$a'=a'_0 \oplus a'_1,\quad g(b')=b'_0 \oplus b'_1,\, e_i''=e_{i,0}\oplus e_{i,1}, \,\,\textrm{and}\quad E''=E'_0 \oplus E'_1$$
with $a'_0, b'_0, E_0', e_{i,0}\in pAp$ and $a'_1, b'_1, e_{i,1},\, E_1' \in S,\, i=1,2,...,n+1$.
Note that $E_1'b_1'=E_1'b_1'=E_1'.$
This, in particular, implies that
\begin{equation}\label{12/27/04-comp-1}
\tau(b_1')\ge (n+1)\tau(e_{1,1})\rforal \tau\in T(S).
\end{equation}
It follows from \eqref{Compri-n4} that
$$
d_\tau(f_{\ep/2}(a_1')) \leq \frac{n}{n+1} d_\tau (b_1'),
\rforal \tau\in\mathrm{T}(S).
$$
{{Note that  $(b_1'-e_{1,1})e_{1,1}=0$ and $b_1'=(b_1'-e_{1,1})+e_{1,1}.$}} {{F}}or all $\tau\in T(S),$
$$
d_\tau((b'_1-e_{1,1}))=d_\tau(b'_1)-\tau(e_{1,1})>d_\tau(b'_1)-\frac{1}{n+1}d_\tau(b'_1)\geq d_\tau(f_{\ep/2}(a_1'))).
$$
Since $S$ has the strict comparison, one has
$$
{{f_{\ep/2}(a_1')}}\lesssim  (b_1'-{\blue{e_{1,1}}}).
$$
{{
Just as in the calculation of (\ref{Compri-n-3}),}} {\blue{ $f_{\ep}(a)\lesssim (a'-\ep/4)_+\sim  f_{\ep/2}(a')$ as $\eta<\ep/4.$}} {{ Consequently,}}
\beq\label{Compari-n6}
&&f_{\ep}(a)\lesssim f_{\ep/2}(a')\lesssim  p\oplus f_{\ep/2}(a_1')\lesssim p\oplus(b'_1-e_{1,1})\\
&& \lesssim  e_1 \oplus(b_1'-e_{1,1})\lesssim e_1\oplus (b_1'-e_{1,1})+(b_0'-e_{1,0})\\
&&\sim e_1''\oplus (g(b')-e_1'')\sim g(b')\sim b'\lesssim b.
\eneq
Since $\ep$ is arbitrary, one has that $a\lesssim b.$

Hence one always has that $a\lesssim b$, and therefore $W(A)$ is almost unperforated.
\end{proof}

The following  fact is known to experts. We include it here for the reader's convenience.

{{
\begin{lem}\label{LCuP}
Let $A$ be a unital simple \CA\, with $T(A)\not=\emptyset.$ Let  $a,\, b\in A_+$ and let $1>\ep>0$ such that
$|\tau(a)-\tau(b)|<\ep$ for all $\tau\in T(A).$ Then
there are $x_1, x_2,...,x_n\in A$ such that
\beq
\|\sum_{i=1}^n x_i^*x_i-a\|<2\ep\andeqn \|\sum_{i=1}^n x_ix_i^*-b\|<\ep.
\eneq
\end{lem}
}}
\begin{proof}
This follows from results in \cite{CP}.
Let $\dt=\max\{|\tau(a)-\tau(b)|: \tau\in T(A)\}.$ Then $0\le \dt<\ep.$
Let $\eta=\ep-\dt.$
By Theorem 9.2 of \cite{LinTAI}, there exists $c\in A_{s.a.}$ with $\|c\|<\dt+\eta/4$
such that $\tau(c)=\tau(b)-\tau(a)$ for all $\tau\in T(A).$
Consider $a_1=a+c+\|c\|$ and $b_1=b+\|c\|.$
Note that $a_1\ge 0.$
Then $\tau(a_1)=\tau(b_1)$ for all $\tau\in T(A).$ It follows from (iii) of Theorem 2.9  of \cite{CP} that
$a_1-b_1\in A_0$ (in the notation of \cite{CP}). It follows from Theorem 5.2 of \cite{CP}
that $a_1\sim b_1$ (in the notation  of \cite{CP}; see the lower half of page136 of \cite{CP}). Thus, there are $x_1, x_2,...,x_n\in A$
such that
\beq\nonumber
\|a_1-\sum_{i=1}^nx_i^*x_i\|<\eta/4\andeqn \|b_1-\sum_{i=1}^nx_ix_i^*\|<\eta/4.
\eneq
It follows that
\beq\nonumber
\|a-\sum_{i=1}^nx_i^*x_i\|\le \eta/4+\|a+\|a\|\|\le \eta/4+2\dt<2\ep\andeqn\\
\|b-\sum_{i=1}^nx_ix_i^*\|\le \eta/4+\|a\|\le \eta/4+\dt<\ep.
\eneq
\end{proof}

\begin{lem}\label{Tapprox}
Let $\mathcal D$ be a class of unital amenable \CA s, let $A$ be a separable unital \CA\, which
is $TA{\cal D},$ and
let $C$ be a unital  {{amenable}}  \CA.

Let $\mathcal F, \mathcal G\subset C$ be finite subsets, let ${{\ep}}>0$ and $\delta>0$ be positive numbers. Let $\mathcal H\subset C_+^{\bf 1}$ be a finite subset, and let $T: C_+\setminus\{0\}\to \mathbb R_+\setminus\{0\}$ and $N: C_+\setminus\{0\}\to\mathbb N$ be maps. Let $\Delta: C_+^{q, {\bf 1}}\setminus\{0\}\to (0, 1)$ be an order preserving map. Let $\mathcal H_1\subset C_{+}^{\bf 1}\setminus \{0\}$, $\mathcal H_2\subset {{C_+}}$ and $\mathcal U\subset {\blue{U(M_k(C))/CU(M_k(C))}}$
{\blue{(for some $k\ge 1$)}} be finite subsets. Let $\sigma_1>0$ and $\sigma_2>0$ be constants.   Let $\phi, \psi: C\to A$ be  unital $\mathcal G$-$\dt$-multiplicative
{{completely positive}} linear maps
 such that
\begin{enumerate}
\item\label{cond-708-001} $\phi$ and $\psi$ are $T\times N$-$\mathcal H$-full  (see the definition \ref{Dfull}),
\item\label{cond-708-002} $\tau\circ\phi(c)>\Delta(\hat{c})$ and $\tau\circ\psi(c)>\Delta(\hat{c})$
{{ for any $\tau\in T(A)$ and}} for any $c\in\mathcal H_1$,
\item\label{cond-708-003} $|\tau\circ\phi(c)-\tau\circ\psi(c)|<\sigma_1$ for any $\tau\in T(A)$ and any $c\in \mathcal H_2$,
\item\label{cond-708-004} $\mathrm{dist}(\phi^{\ddagger}(u), \psi^{\ddagger}(u))<\sigma_2$ for any $u\in\mathcal U$.
\end{enumerate}

Then, for any finite subset $\mathcal F'\subset A$ and ${{\ep}}'>0$, there exists a $C^*$subalgebra $D\subset A$ with $D\in\mathcal D$ such that if $p=1_D$, then, for any $a\in\mathcal F'$,

\indent
{\rm (a)} $\|pa-ap\|<{{\ep}}'$,

\indent
{\rm (b)} $pap\in_{{{\ep}}'} D,$ and

\indent
{\rm (c)} $\tau(1-p)<{{\ep}}'$, for any $\tau\in T(A)$.

There are also a (completely positive) linear map $j_0: A \to (1-p)A(1-p)$ and
a unital
\cp\,.
$j_1: A\to D$ such that
\beq\nonumber
j_0(a)=(1-p)a(1-p) \rforal a\in A\andeqn\\\nonumber
\|j_1(a)-pap\|<3{{\ep}}' \rforal a\in {\mathcal F}.
\eneq

Moreover, define
$$
\phi_0=j_0\circ\phi,\,\, \psi_0=j_0\circ\psi,
\phi_1=j_1\circ\phi\quad\textrm{and}\quad \psi_1=j_1\circ\psi.
$$
With a sufficiently large $\mathcal F'$ and small enough ${{\ep}}'$, one has that $\phi_0$, $\psi_0$, $\phi_1$ and $\psi_1$ are $\mathcal G$-$2\dt$-multiplicative and
\begin{enumerate}\setcounter{enumi}{4}
\item\label{concl-708-001} $\|\phi(c)-(\phi_0(c)\oplus\phi_1(c))\|<{{\ep}}$ and $\|\psi(c)-(\psi_0(c)\oplus\psi_1(c))\|<{{\ep}}$, for any $c\in\mathcal F$,
\item\label{concl-708-002} $\phi_0, \psi_0$ and $\phi_1, \psi_1$ are $2T\times N$-$\mathcal H$-full,
\item\label{concl-708-003} $\tau\circ\phi_1(c)>\Delta(\hat{c})/2$ and $\tau\circ\psi_1(c)>\Delta(\hat{c})/2$ for any $c\in\mathcal H_1$ {\blue{and for any $\tau\in T(D),$}}
\item\label{concl-708-004} $|\tau\circ\phi_1(c)-\tau\circ\psi_1(c)|<3\sigma_1$ for any $\tau\in \mathrm{T}(D)$ and any $c\in \mathcal H_2$, and
\item\label{concl-708-005} $\mathrm{dist}(\phi_{{i}}^{\ddagger}(u), \psi_{{i}}^{\ddagger}(u))<2\sigma_2$ for any $u\in\mathcal U,$
{$i=0,1.$}
\end{enumerate}

{\blue{If, furthermore, ${\cal P}\subset \underline{K}(C)$ is a finite subset and
$[L]|_{\cal P}$ is well defined for
any ${\cal G}$-$2\dt$-multiplicative \morp\, $L,$  and $[\phi]|_{\cal P}=[\psi]|_{\cal P},$
then, we may require, with possibly smaller $\ep'$ and larger ${\cal F}',$  that
\begin{enumerate}\setcounter{enumi}{9}
\item\label{concl-708-006} $[\phi_i]|_{\cal P}=[\psi_i]|_{\cal P},$ $i=0,1.$
\end{enumerate}
}}
\end{lem}

\begin{proof}
Without loss of generality, one may assume that each element of $\mathcal F$, $\mathcal G,$
${\cal H}_2,$ or $\mathcal F'$ has norm at most one and that $1_A\in \mathcal F'$.

Since $\phi$ and $\psi$ are $T\times N$-$\mathcal H$-full, for each $h\in\mathcal H$, there are $a_{1,h}, ..., a_{N(h),h}$ and $b_{1,h}, ..., b_{N(h),h}$ in $A$ with $\|a_{i,h}\|, \|b_{i,h}\|\le T(h)$ such that
\beq\label{911-n1}
\sum_{i=1}^{N(h)} a^*_{i,h}\phi(h)a_{i,h}=1_A \quad\textrm{and}\quad \sum_{i=1}^{N(h)} b^*_{i,h}\psi(h)b_{i,h}=1_A.
\eneq
Put
$
d_0=\min\{\Delta(\hat{h}): h\in {\cal H}_1\}.
$
By  (\ref{cond-708-003}),
It follows from \ref{LCuP} that there are, for each $c\in {\cal H}_2,$
$x_{1,c},x_{2,c},...,x_{t(c), c}\in A$ such that
\beq\label{911-n2}
\|\sum_{i=1}^{t(c)} x_{i, c}^*x_{i,c}-\phi(c)\|\le {{\sigma_3}}
\andeqn
\|\sum_{i=1}^{t(c)} x_{i,c}x_{i,c}^*-\psi(c)\|\le {{\sigma_3}}
\eneq
{{for some $0<\sigma_3<6\sigma_1/5.$}}
Let
$$
t({\cal H}_2)=\max \{(\|x_{i,c}\|+1)(t(c)+1): 1\le i\le t(c): c\in {\cal H}_2\}.
$$
For the given finite subset $\mathcal F'\subset A,$
and given ${{\ep}}'>0$,
since $A$ can be tracially approximated by the \CA s
in the class $\mathcal D$, there exists a \SCA\, $D\subset A$ with $D\in\mathcal D$ such that if $p=1_D$, then, for any $a\in\mathcal F'$,

(i) $\|pa-ap\|<{{\ep}'}$,

(ii) $pap\in_{{{\ep}'}} D$, and

(iii)  $\tau(1-p)<{{\ep}'}$, for any $\tau\in T(A)$.

On the way to making ${\cal F}'$ large and ${{\ep}'}$ small,
we may assume that ${\cal F}'$  contains ${\cal F},$ ${\cal G},$ ${\cal H},$
$\phi({\cal G}\cup {\cal F}),$ $\psi({\cal G}\cup {\cal F}),$
$\phi({\cal H}),$ $\psi({\cal H}),$ ${\cal H}_1,$ ${\cal H}_2,$
$x_{i,c}, x_{i,c}^*,$ $i=1,2,...,t(c)$ and $c\in {\cal H}_2,$ as well as
$a_{i,h}, a_{i,h}^*, b_{i,h}, b_{i,h}^*,$
$i=1,2,...,N(h)$ and $h\in {\cal H},$ and
$$\ep'<\min\{\min\{1/64(T(h)+1)(N(h)+1): h\in {\cal H}\}, \ep, \dt,d_0, \sigma_1, \sigma_2, (2\sigma_1-\sigma_3)\}/64(t({\cal H}_2)+1)^2.$$

For each $a\in\mathcal F'$, choose $d_a\in D$ such that $\|pap-d_a\|<{{\ep'}}$ (choose $d_{1_A}=1_D$).  Consider the finite subset $\{d_ad_b:\ a, b\in\mathcal F'\}\subset D$. Since $D$ is an amenable \SCA\,  of $pAp$,
{\blue{by 2.3.13 of \cite{Lnbok} (with $D=B=C$ and $\phi={\rm id}_D$),}} there is unital completely positive linear map $L: pAp\to D$ such that
$$
\|L(d_ad_b)-d_ad_b\|<{{\ep'}},\quad a, b\in\mathcal F'.
$$
Define $j_1: A\to D$ by $j_1(a)=L(pap).$
Then, for any $a\in\mathcal F'$, one has
\begin{eqnarray}\label{dj1-18}
\|j_1(a)-pap\| & = &\| L(pap)-pap\|=
\|L(d_a)-d_a\| + 2{{\ep}'}={{2\ep'}}.
\end{eqnarray}
Note that $j_0$ and $j_1$ are  {{${\cal F'}$-$7{{\ep}}'$}}
-multiplicative, and
\beq\label{dj1-18+}
\|a-j_0(a)\oplus j_1(a)\|<4{{\ep}}' \rforal  a\in\mathcal F'.
\eneq
Define
$$
\phi_0=j_0\circ\phi,\,\, \psi_0=j_0\circ\psi,
\phi_1=j_1\circ\phi\quad\textrm{and}\quad \psi_1=j_1\circ\psi.
$$
Then (by the choices of ${\cal F}'$ and $\ep'$),
the maps $\phi_0$, $\psi_0$, $\phi_1$, and $\phi_1$ are $\mathcal G$-$2\dt$-multiplicative, and for any $c\in\mathcal F$,
\begin{equation}\label{app-709-0002}
\|\phi(c)-(\phi_0(c)\oplus\phi_1(c))\|<{{\ep}}\quad\textrm{and}\quad \|\psi(c)-(\psi_0(c)\oplus\psi_1(c))\|<{{\ep}}.
\end{equation}
So \eqref{concl-708-002} holds. Apply $j_1$ to both sides of both equations in (\ref{911-n1}).
One obtains two invertible elements  $e_h:=\sum_{i=1}^{N(h)} j_1(a_i^*)\phi_1(h)j_1(a_i)$ and $f_h:=\sum_{i=1}^{N(h)} j_1(b^*_i)\psi_1(h)j_1(b_i)$  such that
$ |\|e_h^{-\frac{1}{2}}\|-1|<1$ and
 $|\|f_h^{-\frac{1}{2}}\|-1|<1.$
Note that
$$\sum_{i=1}^{N(h)} e_h^{-\frac{1}{2}}j_1(a_i^*)\phi_1(h)j_1(a_i)e_h^{-\frac{1}{2}}=1_D,
\sum_{i=1}^{N(h)} f_h^{-\frac{1}{2}}j_1(b_i^*)\psi_1(h)j_1(b_i)f_h^{-\frac{1}{2}}=1_D,$$
$$\|j_1(a_i)e_h^{-\frac{1}{2}}\|<2 T(h), \quad\mathrm{and}\quad \|j_1(b_i)f_h^{-\frac{1}{2}}\|<2 T(h).$$
Therefore, $\phi_1$ and $\psi_1$ are $2T\times N$-$\mathcal H$-full, and this proves \eqref{concl-708-003}. The same calculation also shows that $\phi_0$ and $\psi_0$ are $2T\times N$-$\mathcal H$-full.
Note that we have shown
(\ref{concl-708-001})
holds.
Since $\ep'<d_0/16,$ it is also easy to check that  \eqref{concl-708-003} holds.

To see \eqref{concl-708-005}, one notes that
\beq\label{911-n4}
\|\sum_{i=1}^{t(c)} d_{x_{i,c}}^*d_{x_{i,d}}-\phi_1(c)\|<\sigma_3+t({\cal H}_2)\ep'<7\sigma_1/5\andeqn
\|\sum_{i=1}^{t(c)}d_{x_{i,c}}d_{x_{i,d}}^*-\psi_1(c)\|<7\sigma_1/5
\eneq
for all $c\in {\cal H}_2.$ Then \eqref{concl-708-005} also holds.


%
Let us show \eqref{concl-708-004} holds for sufficiently large ${\cal F}'$ and $\ep'.$
Since $A$ is separable, there are an increasing sequence of finite subsets $\mathcal F_1'\subset\mathcal F_2'\subset\cdots$ such that $\bigcup\mathcal F'_n$ is dense in the unit ball of $A$. Set $\epsilon_n'=\frac{1}{n}$. Suppose
\eqref{concl-708-004} were not true, for each $\mathcal F'_n$ and each $\epsilon'_n$, there are C*-subalgebra $D_n\in\mathcal D$ and $j_{1, n}: A\to D_n$ as constructed above (in place of $j_1$), and there is $\tau_n\in \mathrm{T}(D_n)$ such that there is $c\in\mathcal H_1$
$$\tau_n\circ\phi_{1,n}(c)\leq \Delta(\hat{c})/2\quad \textrm{or}\quad \tau_n\circ\psi_{1,n}(c)\leq\Delta(\hat{c})/2$$
(where $\phi_{1,n}$ and $\psi_{1,n}$ are as $\phi_1$ and $\psi_1$ corresponding to ${\cal F}'={\cal F}'_n$ and $\ep'=\ep_n'=1/n$
for all large $n$).
Passing to a subsequence, one may assume that
\begin{equation}\label{contr-709-0003}
\tau_n\circ\phi_{1,n}(c)\leq \Delta(\hat{c})/2.
\end{equation}
Consider $\tau_n\circ j_{1, n}: A\to \mathbb C$, and pick an accumulating point $\tau$ of $\{ \tau_n\circ j_{1, n}: n\in\mathbb N\}$. Since $j_{1, n}$ is $7\epsilon'_n$-$\mathcal F_n'$-multiplicative, it is straightforward to verify that $\tau$ is actually a tracial state of $A$.
Passing to a subsequence, we may assume that $\tau(a)=\lim_{n\to\infty}\tau_n\circ j_{1,n}(a)$ for all $a\in A.$
By \eqref{contr-709-0003}, for any $0<\eta<d_0/4,$ there exists $n_0\ge 1$ such that, for all $n\ge n_0,$
$$
\tau\circ \phi(c)\le \tau_n\circ j_{1,n}\circ \phi(c)+\tau_n\circ \phi_{2,n}(c)+\eta\le
\Delta(\hat{c})/2 +d_0/4+\eta<\Delta(\hat{c}),
$$
which contradicts to the assumption \eqref{cond-708-002}.


Let us show that  (9)
holds with sufficiently large $\mathcal F'$ and sufficiently small $\epsilon'$.

Choose unitaries $u_1, u_2, ..., u_n\in {\blue{M_k(C)}}$ such that $\mathcal U=\{\overline{u_1}, \overline{u_2}, ..., \overline{u_n}\}$. Pick unitaries $w_1, w_2, ..., w_n\in {\blue{M_k(A)}}$ such that
each $w_i$ is a commutator and
$$\mathrm{dist}(\left<\phi(u_i)\right>\left<\psi(u^*_i)\right>, w_i)<\sigma_2.$$
Choose $\mathcal F'$ sufficiently large and $\epsilon'$ sufficiently small such that there are commutators $w_1', w_2, ..., w_n'\in
CU({\blue{M_k(D)}})$  and commutators $w_1'', w_2'',...,w_n''\in (1-p)A(1-p)\otimes M_k$ satisfying
$$\|j_1(w_i)-w_i'\|<\sigma_2/2 \andeqn \|j_0(w_i)-w_i''\|<\sigma_2/2, \quad 1\leq i\leq n,$$
(see Appendix of \cite{Lin-LAH}) and
$$\| \left<\phi_{k}(u_i)\right>\left<\psi_{k}(u^*_i)\right> - j_{k}(\left<\phi(u_i)\right>\left<\psi(u^*_i)\right>)\|<\sigma_2/2, \quad 1\leq i\leq n \andeqn k=0,1.$$
{\blue{(Recall we use  $\phi, \psi, j_0,$ and  $j_1$ for $\phi\otimes {\rm id}_{M_k},$
$\psi\otimes {\rm id}_{M_k},$ $j_0\otimes {\rm id}_{M_k}$ and $j_1\otimes{\rm id}_{M_k},$ respectively.)}}
Then
\begin{eqnarray}
&&\hspace{-0.4in}\|  \left<\phi_{k}(u_i)\right>\left<\psi_{k}(u^*_i)\right> - w_i'  \|
\leq  \|  \left< \phi_{k}(u_i)\right>\left<\psi_{k}(u^*_i)\right> - j_{k}(w_i) \| + \| j_{k}(w_i)  -w_i'\| \\
&\leq & \|  \left< \phi_{k}(u_i)\right>\left<\psi_{k}(u^*_i)\right> - j_{k}(\left<\phi(u_i)\right>\left<\psi(u^*_i)\right>)\| + \\
&& \| j_{k}(\left<\phi(u_i)\right>\left<\psi(u^*_i)\right>)-j_{k}(w_i)\| + \sigma_2/2
\leq 2\sigma_2,\,\,\, k=0,1.
\end{eqnarray}
This proves \eqref{concl-708-005}.

{\blue{To see the last part of the lemma,
let ${\cal P}_{0}$ be a finite subset of projections of $C,$
let  $q_\phi, q_\psi\in A$ be projections, and let
$v_q\in A$ be a partial isometry (for each $q\in {\cal P}_0$) such that
$v_q^*v_q=q_\phi,$ $v_qv_q^*=q_\psi,$ $\|q_\phi-\phi(q)\|<\dt'/2<1/4,$ and $\|q_\psi-\psi(q)\|<\dt'/2<1/4$
for some $\dt'>0.$

As in part (d) of Lemma \ref{approx-Aug-14-1}, with sufficiently small  $\ep'$ and large ${\cal F}',$ one obtains projections $q_{\phi,0}, q_{\psi,0}\in (1-p)A(1-p)$ and
$q_{\phi, 1}, q_{\psi, 1}\in D,$ and
partial isometries $w_{q,0,}\in (1-p)A(1-p)$ and $w_{q,1}\in D$ such that
\beq
&&w_{q,i}^*w_{q,i}=q_{\phi,i}, w_{q,i}w_{q,i}^*=q_{\psi, i},\\
&&\|\phi_i(q)-q_{\phi,i}\|<\dt',\andeqn \|\psi_i(q)-q_{\psi,i}\|<\dt'
\eneq
for all $q\in {\cal P}_{0}.$ This implies that $[\phi_i(q)]=[\psi_i(q)]$ for all $q\in {\cal P}_0,$ $i=0,1.$

Suppose that ${\cal U}_0$ is a finite subset of $U(C)$ and $\la \phi(u)\ra=z_u\la \psi(u)\ra,$
where $z_u=\prod_{k=1}^{l(u)}\exp(i h_{u,k})$ and where $h_{u,i}\in A_{s.a.},$ for each $u\in {\cal U}_0.$
By virtue of part (d) of  Lemma \ref{approx-Aug-14-1}, there are $h_{u,k,0}\in (1-p)A(1-p)_{s.a.}$ and $h_{u,k,1}\in D_{s.a.}$ such
that $\la \phi_s(u)\ra =(\prod_{k=1}^{l(u)}\exp(i h_{u,k,s}))\la \psi_s(u)\ra,$ $u\in {\cal U}_0,$ $s=0,1.$
This implies that $[\phi_s]|_{{\cal U}_0}=[\psi_s]|_{{\cal U}_0},$ $s=0,1.$

If ${\cal P}_1$ is a finite subset of projections and unitaries  in
$M_m(A)$ for some integer $m\ge 1,$ by considering $\phi\otimes {\rm id}_{M_m}$ and
$\psi\otimes {\rm id}_{M_m},$ with sufficiently small $\ep'$ and large ${\cal F}',$ we conclude
that we can require that $[\phi_i]|_{{\cal P}_1}=[\psi_i]|_{{\cal P}_1},$ $i=0,1.$

In general, consider a finite subset ${\cal P}_k\subset K_0(A\otimes B_0)$ for some
$B_0={\tilde B},$ where $B$ is a commutative \CA\, such that $K_0(B)=\Z/k\Z$ and $K_1(B)=\{0\},$
$k=2,3,....$
Also consider ${\tilde \phi}=\phi\otimes {\rm id}_{B_0}$ and ${\tilde\psi}=\psi\otimes {\rm id}_{B_0}.$
We will replace $(1-p)A(1-p)$ by $A_0,$ where\\
 $A_0:=(1-p\otimes 1_{B_0})(A\otimes B_0)(1-p\otimes 1_{B_0}),$
and $D$ by $D\otimes B_0$ in the above argument.
We will also consider $(j_i\otimes {\rm id}_{B_0})\circ {\tilde \phi}$ and
$(j_i\otimes {\rm id}_{B_0})\circ {\tilde \psi}.$
Note that  $1-(p\otimes 1_{B_0})$ almost commutes with ${\tilde \phi}$  and ${\tilde \psi}$ on a given
finite subset provided that $\ep'$ is sufficiently small and ${\cal F}'$ is sufficiently large.
Thus, as above, with sufficiently small $\ep'$ and ${\cal F}',$
$[{\tilde \phi}]|_{{\cal P}_k}=[{\tilde \psi}]|_{{\cal P}_k}.$
This implies that the last part of the lemma holds. }}

\end{proof}

\begin{prop}\label{B0not=B1}
{{${\cal B}_1\not={\cal B}_0.$}}
\end{prop}

\begin{proof}
{{It follows from Theorem 1.4 of \cite{Mg} that there is a unital separable simple \CA\,
$A$ which is an inductive limit of dimension drop circle algebras such
that $A$ has a unique tracial state, $(K_0(A), K_0(A)_+, [1_A])=(\Z, \Z_+, 1),$ and  $K_1(A)=\Z/3\Z.$
Note that dimension drop circle algebras are in ${\cal C}$ (see \ref{dimcircle}).}}

{{Note also that $A$ is a unital projectionless \CA.
If $A$ were in ${\cal B}_0,$ since $A$ does not have (SP),  by Theorem \ref{B1sp},  $A$ would be
an inductive limit (with not necessarily injective maps) of \CA s in ${\cal C}_0.$ However, every \CA\, in ${\cal C}_0$ has
trivial $K_1.$ This would imply that $K_1(A)=\{0\},$ a contradiction.
Thus, $A\not\in {\cal B}_0.$}}
\end{proof}

\begin{rem}
{{The \CA\, $A$ in the proof of  Proposition \ref{B0not=B1} is rationally of tracial rank zero (see \cite{LNjfa}), since it has a
unique tracial state.  Later we will see that there are many other \CA s which are in ${\cal B}_1$ but not
in ${\cal B}_0.$}}
\end{rem}

{{\begin{prop}\label{Psp}
Let $A$ be a \CA\, in  ${\cal B}_1$ and let $U$ be an infinite dimensional  UHF-algebra.
Then $A\otimes U$ has the property (SP).
\end{prop}
}}
\begin{proof}
{{It suffices to show that, for any non-zero positive element $b\in A\otimes U,$ there
exists a non-zero projection $e\in A\otimes U$ such that
$e\lesssim b.$  \Wlog, we may assume that $0\le b\le 1.$
Let  $\sigma=\inf\{\tau(b): \tau\in T(A)\}>0.$  Then
there is a non-zero projection $e\in A\otimes U$ with the form
$1_A\otimes p,$ where $p\in U$ is a non-zero projection, such that
\beq
\tau(e)<\sigma.
\eneq
By strict comparison for positive elements (Theorem \ref{Comparison}), $e\lesssim b,$ which implies that $\overline{bAb}$ has a projection
equivalent to $e.$}}
\end{proof}

\section{${\cal Z}$-stability}

\begin{lem}\label{Affon1}
Let $A\in {\cal B}_1$ {\rm (}or ${\cal B}_0$ \rm{)} be a unital infinite dimensional simple \CA. Then,  for any $\ep>0,$ any $a\in A_+\setminus\{0\},$ any finite subset ${\cal F}\subset A$ and any
integer $N\ge 1,$ there exist a projection $p\in A$ and a \SCA\, $C\in {\cal C}$  (or $\in {\cal C}_0$)  with $1_C=p$ {{that}} satisfy the following conditions:

{\rm (1)} ${\rm dim}(\pi(C))\ge N^2$ for every irreducible representation $\pi$ of $C,$

{\rm (2)} $\|px-xp\|<\ep$ for all $x\in {\cal F},$

{\rm (3)} ${\rm dist}(pxp, C)<\ep$ for all $x\in {\cal F},$ and

{\rm (4)} $1-p\lesssim a.$

\end{lem}

\begin{proof}
Since $A$ is an  infinite dimensional simple \CA, there are
$N+1$ mutually orthogonal non-zero positive elements
$a_1, a_2,...,a_{N+1}$ in $A.$  Since $A$ is simple,
there are $x_{i,j}\in A,$ $j=1,2,...,k(i),$ $i=1,2,...,N+1,$ such that
$$
\sum_{j=1}^{k(i)}x_{i,j}^*a_ix_{i,j}=1_A.
$$
Let
$$
K=(N+1)\max\{\|x_{i,j}\|+1: 1\le j\le k(i),\,\,\, 1\le i\le N+1\}.
$$
Let $a_0\in A_+\setminus \{0\}$ be such that
$a_0\lesssim a_i$ for all $1\le i\le N+1,$  {{and furthermore $a_0\lesssim a$.}}
Since $\overline{a_0Aa_0}$ is also an infinite dimensional simple \CA, one obtains $a_{01}, a_{02}\in \overline{a_0Aa_0}$ which are mutually orthogonal and non-zero. One then
obtains a non-zero element ${{b}} \in \overline{a_{01}Aa_{01}}$ such that ${{b}}\lesssim a_{02}.$

Let
$$
{\cal F}=\{a_i: 1\le i\le N+1\}\cup\{x_{i,j}: 1\le j\le k(i),\,1\le i\le N+1\}\cup \{a\}.
$$
Now since $A\in {\cal B}_1,$ there are a projection $p\in A$ and $C\in {\cal C}$ with $1_C=p$ such that

(1) $\|xp-px\|<\min\{1/2,\ep\}/2K$ for all $x\in {\cal F},$

(2) ${\rm dist}(pxp, C)<\min\{1/2, \ep\}/2K$ for all $x\in {\cal F}$, and

(3) $1-p\lesssim {{b}}.$

Then, with a standard computation, we obtain mutually orthogonal
 non-zero positive elements $b_1,b_2,...,b_{N+1}\in C$ and
 $y_{i,j}, \in C$ ($1\le j\le k(i)),$ $i=1,2,...,N+1,$ such that
\beq\label{Affon1-2}
\|\sum_{j=1}^{k(i)}y_{i,j}^*b_iy_{i,j}-p\|<\min\{1/2, \ep/2\}.
\eneq
For each $i,$ we find another element $z_i\in C$ such that
\beq\label{Affon1-3}
\sum_{j=1}^{k(i)}z_i^*y_{i,j}b_iy_{i,j}z_i=p.
\eneq
Let $\pi$ be an irreducible representation of $C.$ Then
by (\ref{Affon1-3}),
\beq\label{Affon1-4}
\sum_{j=1}^{k(i)}\pi(z_i^*y_{i,j})\pi(b_i)\pi(y_{i,j}z_i)=\pi(p).
\eneq
Therefore, $\pi(b_1), \pi(b_2),...,\pi(b_{N+1})$ are mutually orthogonal
non-zero positive elements in $\pi(A).$ Then (\ref{Affon1-4}) implies
that $\pi(C)\cong M_n$ with $n\ge N+1.$ This proves the lemma.
\end{proof}

\begin{cor}\label{CTdense}
Let $A\in {\cal B}_1$ {{(or $A\in {\cal B}_0$)}} be  {{an infinite dimensional}} unital simple \CA. Then, for any $\ep>0$ and $f\in \Aff(T(A))^{{++}},$ there exist
a \SCA\, {{$C\in {\cal C}$ (or $C\in {\cal C}_0$ )}} in $A$ {{and}}  an element $c\in C_+$ such that
\beq\label{CTdense-1}
&&{\rm dim}\pi(C)\ge  (4/\ep)^2\,\,\,\text{for \,each \,irreducible\,representation} \,\,\pi\,\,\,of \,\,C,\\
&&0<\tau(f)-\tau(c)<\ep/2\rforal \tau\in T(A).
\eneq
\end{cor}

\begin{proof}
{\blue{By 9.3 of \cite{LinTAI},}}  {{there is an element $x\in A_+$ such that $\tau(x)=\tau(f)$ for all $\tau\in T(A)$.
Let $N\geq {\blue{16}}(\|x\|+1)/\ep$.}}  {\blue{As in the beginning of the proof of Proposition \ref{Ltrace}, one finds
a non-zero element $a\in A_+$
such that $N[a]\le 1.$}}
{{ Apply Lemma \ref{Affon1} to  $\ep/(16(\|x\|+1)),$ $N$ and $a\in A_+$ and ${\cal F}=\{x, 1\}$, to get $C$ and $p=1_C$ as
in that   lemma.}} {\blue{Then $C$ satisfies \eqref{CTdense-1} and, by (4) of  Lemma \ref{Affon1},  $\tau(1-p)<1/N$
for all $\tau\in T(A).$ It follows that
$0<\tau(x)-\tau(pxp)<2\|x\|/N<\ep/8$ for all $\tau\in T(A).$}}
  {{Then choose $c'\in C$ with $\|pxp-c'\|<\ep/{\blue{16(\|x\|+1)}}$.
  Replacing $c'$ by $(c'+(c')^*)/2,$ we may assume that $c'\in C_{s.a.}.$
Since $pxp\ge 0,$ we obtain a positive element $c\in C$
such that $\|pxp-c\|<\ep/8.$}}
 {{ We have
  \beq
  0<\tau(f)-\tau(c)=\tau(x)-\tau(c)<\ep/2\rforal \tau\in T(A),
  \eneq
  as desired.}}
\end{proof}

The following is  known. {{In the following statement, we identify $[0,1]$ with
the space of extreme points of $T(M_n(C[0,1])).$}}

\begin{lem}\label{Affon1+}
Let $C=M_n({{C([0,1])}})$ and $g\in {\rm LAff}_b(T(C))_+$ {{with $0\le g(t)\le 1.$}} Then there exists $a\in
C_+$ with $0\le a\le 1$ such that
$$
0\le g(t)-d_t(a)\le 1/n \tforal t\in [0,1],
$$
where $d_t(a)=\lim_{k\to\infty} {{{\rm tr}}}(a^{1/k}(t))$ for all $t\in [0,1],$ {{where ${\rm tr}$ is the tracial state of $M_n.$}}
\end{lem}

\begin{proof}
We will use the proof of Lemma 5.2 of \cite{BPT}.
For each $0\le i\le  n-1,$ define
$$
X_i=\{t\in [0,1]: g(t)>i/n\}.
$$
Since $g$ is lower semi-continuous, $X_i$ is open in $[0,1].$ There
is a continuous function $g_i\in C([0,1])_+$ with $0\le g_i\le 1$
such that
$$
\{t\in [0,1]: g_i(t)\not=0\}=X_i,\,\,\,i=0,1,...,n-1.
$$
Let $e_1,e_2,...,e_n$ be $n$ mutually orthogonal rank one
projections in $C=M_n(C([0,1]).$ Define \beq\label{Affon1++1}
a=\sum_{i=1}^{n-1} g_i e_i\in C. \eneq Then $0\le a\le 1.$ Put
$Y_i=\{t\in [0,1]: (i+1)/n\ge g_i(t)>i/n\}=X_i\setminus
\bigcup_{j>i} X_j,$ $i=0,1,2,...,n-1.$ These are mutually disjoint
sets. Note that
$$
[0,1]=([0,1]\setminus X_0)\cup \bigcup_{i=0}^{n-1}Y_i.
$$
If $x\in ([0,1]\setminus X_0)\cup Y_0,$ {{then}} $d_t(a)=0.$ So $0\le g(t)-d_t(a)(t)\le 1/n$ for all such $t.$
If $t\in Y_j,$
\beq\label{Affon1++2}
d_t(a)=j/n.
\eneq
Then
\beq\label{Affon1++3}
0\le g(t)-d_t(a)\le 1/n\tforal t\in Y_j.
\eneq
It follows that
\beq\label{Affon1++4}
0\le g(t)-d_t(a)\le 1/n\tforal t\in [0,1].
\eneq
\end{proof}

\begin{lem}\label{Affon2}
Let $F_1$ and $F_2$ be two finite dimensional \CA s such that each
simple {direct} summand of $F_1$ and $F_2$ has rank at least $k$, where $k\ge
1$ is an integer.  Let $\phi_0, \phi_1: F_1\to F_2$ be unital \hm s.
Let $C=A(\phi_0, \phi_1, F_1, F_2).$ Then, for any $f\in {\rm
LAff}_b(T(C))_+$ with $0\le f \le 1,$ there exists a positive
element $a\in M_2(C)$ such that
$$
\max_{\tau\in T(C)}|d_{\tau}(a)-f(\tau)|\le 2/k.
$$
\end{lem}

\begin{proof}
Let $I=\{{{(g,a)}}\in C{{\subset C([0,1], F_2)\oplus F_1}}: g(0)=g(1)=0{{=a}}\}.$
Note that $C/I{{\cong F_1.}}$
Let
$$
T=\{\tau\in T(C):\, {\rm ker}\tau\supset I\}.
$$
Then $T$ may be identified with $T(C/I){{=T(F_1)}}.$ Let $f\in  {\rm
LAff}_b(T(C))_+$ with $0\le f \le 1.$
{{T}}here exists $b\in (C/I)_+$
such that
\beq\label{Affon2-1}
0\le f(\tau)-d_{\tau}(b) \rforal \tau\in T\andeqn
\max\{f(\tau)-d_\tau(b): \tau\in T\}\le 1/k,
\eneq
 {and furthermore,
if $f(\tau)>0$, then $f(\tau)-d_\tau (b)>0.$}
{\blue{To see this,}} {{write $F_1=M_{R(1)}\oplus M_{R(2)}\oplus\cdots \oplus M_{R(l)}.$
Note that $R(s)\ge k,$ $s=1,2,....$ Denote by $\tau_{q(s)}$ the tracial state of $M_{R(s)},$
$s=1,2,...,l.$
One can find an integer $J_s\ge 0$
such that
\beq\label{18814-n10}
{J_s\over{R(s)}}\le f(\tau_{q(s)})\le {J_s+1\over{R(s)}},\,\,\, s=1,2,...,l.
\eneq
Moreover, if $f(\tau_{q(s)})>0,$ we may assume that $J_s/R(s)<f(\tau_{q(s)}).$
Let $b\in F_1=C/I$ be a projection with rank $J_s$ in $M_{R(s)},$ $s=1,2,...,l.$}}
For such a
choice,  {{since $b$ is a projection,}} we have $d_\tau (b)= \tau (b)$ for all $\tau\in T.$
Then, by \eqref{18814-n10}, \eqref{Affon2-1} holds. Moreover, if $f(\tau)>0,$
$f(\tau)-d_\tau(b)>0.$ {{In particular,
if $d_\tau(b)>0,$ then
\beq\label{18815-1}
f(\tau)-d_\tau(b)>0.
\eneq}}

Recall  that $b=(b_1, b_2,...,b_l)\in C/I=F_1.$
Let $b^0=\phi_0(b)\in F_2$ and
$b^1=\phi_1(b)\in F_2.$
Write $F_2=M_{r(1)}\oplus M_{r(2)}\oplus\cdots\oplus M_{r(m)}.$ Write
$b^0=b_{0,1}\oplus b_{0,2}\oplus \dots \oplus b_{0,r(m)}$ and
$b^1=b_{1,1}\oplus b_{1,2}\oplus\cdots \oplus b_{1,r(m)},$ where
$b_{0,j},\,b_{1,j}\in M_{r(j)},$ $j=1,2,...,m.$
 Let $\tau_{t,j}={\rm tr}_j\circ \Psi_j\circ \pi_t,$ where ${\rm tr}_j$ is the normalized trace on $M_{r(j)},$
$\Psi_j: F_2\to M_{r(j)}$ is the quotient map and
$\pi_t: A\to F_2$ is the point evaluation at $t\in (0,1).$

{{Denote by $\tau_{0,j}$ the tracial state of $C$ defined by
$\tau_{0,j}((g,a))={\rm tr}_j(g(0))={\rm tr}_j(\phi_0(a))$ for all $(g,a)\in C,$ and denote
by $\tau_{1,j}$ the tracial state of $C$ defined by
$\tau_{1,j}((g,a))={\rm tr}_j(g(1))={\rm tr}_j(\phi_1(a))$ for all $(g,a)\in C,$ respectively, $j=1,2,...,m.$
It follows from {{the third paragraph of  Subsection}} \ref{LgN1889}  that $\tau_{t,j}\to \tau_{0,j}$ in $T(A)$ if
$t\to 0$ in $[0,1]$ and
$\tau_{t,j}\to \tau_{1,j}$ in $T(A)$ if $t\to 1$ in $[0,1].$ }}

{{Recall from Section 3 (see {\blue{\ref{LgN1889}}} and
{\blue{\ref{2Rg11}}}),
 we have
\beq\label{18815-s10-1}
{\rm tr}_j(b_{0,j})=\frac{1}{r(j)}\sum_{s=1}^l a_{js} R(s)\cdot \tau_{q(s)}(b_{s})~~\mbox{and}~~{\rm tr}_j(b_{1,j})=\frac{1}{r(j)}\sum_{s=1}^l b_{js} R(s)\cdot \tau_{q(s)}(b_s),
\eneq
where $\{a_{js}\}$ and $\{b_{js}\}$ are   matrices of  {\blue{non-}}negative integers given by the maps\\
$(\phi_0)_{*0}: K_0(F_1)=\Z^l\to K_0(F_2)=\Z^m$ and  $(\phi_1)_{*0}: K_0(F_1)=\Z^l\to K_0(F_2)=\Z^m$.}}

{{If $g\in {\rm LAff}_b(T(A))_+,$ then, by \eqref{18815-s3-1} and \eqref{18815-s3-2} in Section 3,
$$
g(\tau_{0,j})=\sum_{s=1}^l a_{js} R(s)\cdot g(\tau_{q(s)})~~\mbox{and}~~g(\tau_{1,j})=\frac{1}{r(j)}\sum_{s=1}^l b_{js} R(s)\cdot  g(\tau_{q(s)}),
$$
}}

 {{Recall that $b$ is a projection,}} since $f$ is lower semicontinuous on $T(C),$  {{we have, by \eqref{Affon2-1}  and
 \eqref{18815-1},}}
\beq\nonumber
&&\hspace{-0.4in}\liminf_{t\to 0} f(\tau_{t,j})\ge  {{f(\tau_{0,j})}}=\sum_{s=1}^l a_{js} R(s)\cdot f(\tau_{q(s)})
\ge \frac{1}{r(j)}\sum_{s=1}^l a_{js} R(s)\tau_{q(s)}(b_{s})=
{\rm tr}_j(b_{0,j})\\\label{18815-s10-2}
&&\andeqn
 \liminf_{t\to
1} f(\tau_{t,j})\ge \frac{1}{r(j)}\sum_{s=1}^l b_{js} R(s)\cdot  f(\tau_{q(s)})\ge {\rm tr}_j(b_{1,j}).
\eneq
Note that,
{{since $\phi_0$ is unital,  $\frac{1}{r(j)}\sum_{s=1}^la_{js}R(s)=1.$ }}Therefore, if ${\rm tr}_j(b_{0,j})>0,$ then
\beq\nonumber
f(\tau_{0,j})=f(\sum_{s=1}^l {{\frac{1}{r(j)}a_{js}R(s)}} {{\tau}}_{q(s)})>\sum_{s=1}^l {{\frac{1}{r(j)}a_{js}R(s)}} {{\tau}}_{q(s)}(b_{s})={\rm tr}_j(b_{0,j}),\\
f(\tau_{1,j})=f(\sum_{s=1}^l {{\frac{1}{r(j)}a_{js}R(s)}} {{\tau}}_{q(s)})>\sum_{s=1}^l {{\frac{1}{r(j)}b_{js}R(s)}} {{\tau}}_{q(s)}(b_{s})={\rm tr}_j(b_{1,j}).
\eneq

{Hence, if ${\rm tr}_j(b_{0,j})>0$ (or ${\rm tr}_j(b_{1,j})>0$),
then $\liminf_{t\to 0} f(\tau_{t,j})> {\rm tr}_j(b_{0,j})$ (or
$\liminf_{t\to 1} f(\tau_{t,j})> {\rm tr}_j(b_{1,j})$).}
 Therefore, there exists $1/8>\dt>0$ such that
\beq\label{Affon2-2}
f(\tau_{t,j}) &\ge& {\rm tr}_j(b_{0,j})\tforal t\in (0,2\dt)\andeqn\\
f(\tau_{t,j}) &\ge&  {\rm tr}_j(b_{1,j})\tforal t\in (1-2\dt,1),
\,\,\,j=1,2,...,l. \eneq Let \beq\label{Affon2-3}
c(t)&=&({\dt-t\over{\dt}})b_0\,\,\,  \text{if}\,\,\,t\in [0,
\dt),\\\label{Affon2-3+1} c(t)&=& 0 \,\,\,\text{if}\,\,\, t\in [\dt,
1-\dt],\andeqn\\\label{Affon2-3+2}
c(t)&=&({t-1+\dt\over{\dt}})b_1\tforal t\in (1-\dt, 1].
\eneq
Note
that $c\in A.$ Define
\beq\label{Affon2-4} g_j(0)&=&
0,\\\label{Affon2-4+1}
g_j(t)&=&f(\tau_{t,j})-{\rm tr}_j(b_{0,j})\tforal
t\in (0, \dt],\\\label{Affon2-4+2}
g_j(t)&=& f(\tau_{t,j}) \tforal
t\in (\dt, 1-\dt),\\\label{Affon2-4+3}
g_j(t)&=& f(\tau_{t,j})-{\rm tr}_j(b_{1,j})\tforal  t\in [1-\dt, 1),\andeqn\\
g_j(1)&=&0.
\eneq
One verifies that $g_j$ is lower semicontinuous on
$[0,1].$ It follows {{from}} Lemma \ref{Affon1+} that there exists $a_1\in C([0,1],
F_2)_+$
such that
\begin{equation}\label{Affon2-5}
0\le g_j(t)-d_{{\rm tr}_{t,j}}(a_1)\le 1/r(j)\leq1/k \tforal t\in [0,1].
\end{equation}
Note that $a_1(0)=0$ and $a_1(1)=0.$
Therefore  $a_1\in I\subset C.$  Now let $a=c\oplus a_1\in M_2(C).$
 Note that
\begin{eqnarray*}
d_{\tau}(a)&=& d_{\tau}(c)+d_{\tau}(a_1)=d_{\tau}(b)\,\,\, \text{if}\,\,\, t\in T,\\
d_{{\rm tr}_{t,j}}(a)&=& d_{{\rm tr}_{t,j}}(c)+d_{{\rm tr}_{t,j}}(a_1)=d_{t,j}(b_0)+d_{{\rm tr}_{t,j}}(a_1)\tforal t\in (0, \dt),\\
d_{{\rm tr}_{t,j}}(a)&=& d_{{\rm tr}_{t,j}}(a_1)\tforal t\in [\dt, 1-\dt],\andeqn\\
d_{{\rm tr}_{t,j}}(a)&=&d_{{\rm tr}_{t,j}}(c)+d_{{\rm tr}_{t,j}}(a_1)=d_{t,j}(b_1)+d_{{\rm tr}_{t,j}}(a_1) \tforal t\in (1-\dt, 1).
\end{eqnarray*}
{{Hence, by}}
\eqref{Affon2-1}, \eqref{Affon2-3}, \eqref{Affon2-3+1}, \eqref{Affon2-3+2}, \eqref{Affon2-4},
\eqref{Affon2-4+1}, \eqref{Affon2-4+2}, \eqref{Affon2-4+3} and \eqref{Affon2-5},
\beq\label{Affon2-6}
0\le f(\tau)-d_{\tau}(a)\le 2/k\tforal \tau\in T \andeqn \tau={\rm tr}_{t,j}, j=1,2,...,l, \,\,\,t\in (0,1).
\eneq
Since $T\cup \{{\rm tr}_{t,j}: 1\le j\le l, \andeqn t\in (0,1)\}$ contains all  the extreme  points of $T(C),$ we conclude that
\beq\label{Affon2-7}
0\le f(\tau)-d_{\tau}(a)\le 2/k\tforal \tau\in T(C).
\eneq
\end{proof}


\begin{thm}\label{Aff=W-L}
Let $A\in {\cal B}_1$ be  {{ an infinite dimensional}}  unital simple \CA. Then
the map $W(A)\to V(A)\sqcup {\rm LAff}_b(A)_{++}$ is surjective.
\end{thm}

 \begin{proof}
The proof follows the lines of {{the proof of}} Theorem 5.2 of \cite{BPT}.
 It suffices to show that the map $a\mapsto d_\tau(a)$ is surjective from $W(A)$ onto ${\rm LAff}_b(T(A)).$
Let $f\in {\rm LAff}_b(A)_+$ with $f(\tau)>0$ for all  $\tau\in T(A).$
We may assume that $f(\tau)\le 1$ for all $\tau\in T(A).$
 As in the proof of 5.2 of \cite{BPT}, it suffices to find a sequence of elements $a_i\in M_2(A)_+$ such that
 $a_i\lesssim a_{i+1},$
 $[a_n]\not=[a_{n+1}]$ (in $W(A)$), and
 $$
 {\lim_{n\to\infty}d_{\tau}(a_n)=f(\tau)\rforal \tau\in T(A).}
 $$
Using the semicontinuity of $f,$  we find a sequence $f_n\in
\Aff(T(A))_{++}$ such that \beq\label{Aff=W-L-1}
f_n(\tau)< f_{n+1}(\tau)\rforal \tau\in T(A),\,\,\,n=1,2,...,\\
\lim_{n\to\infty} f_n(\tau)=f(\tau)\tforal \tau\in T(A).
 \eneq
Since $f_{n+1}-f_n$ is continuous and strictly positive on the
compact set $T,$ there is ${{1>}}\ep_n>0$ such that
$(f_{{n+1}}-f_{{n}})(\tau)>\ep_n$ for all $\tau\in T(A),$ $n=1,2,....$
{{By Corollary}} \ref{CTdense}, for each $n,$ there is a \SCA\,
$C_n$ of $A$ with $C_n\in {\cal C}$ and an element $b_n\in (C_n)_+$
such that \beq\label{Aff=W-L-2}
{\rm dim}\pi(C_n)\ge  (16/\ep_n)^2\,\,\,\text{for \,each \,irreducible\,representation} \,\,\pi\,\,\,of \,\,C_n,\\
0<\tau(f_n)-\tau(b_n)<\ep_n/8\rforal \tau\in T(A).
\eneq
Applying Lemma \ref{Affon2}, one obtains an element $a_n\in M_2(C_n)_+$
such that
\beq\label{Aff=W-L-3} 0<t(b_n)-d_t(a_n)<{{2\over{k}}}<\ep_n/4\tforal
t\in T(C_n),
\eneq
{{where $k$ is the minimal rank of all irreducible representations of $C_n$ and $k\geq {\blue{16/\ep_n}}$.}}
It follows that \beq\label{Aff=W-L-4}
0<\tau(f_n)-d_\tau(a_n)<\ep_n/2\tforal \tau\in T(A). \eneq One then
checks that $\lim_{n\to\infty}d_\tau(a_n)=f(\tau)$ for all $\tau\in
T(A).$ Moreover, $d_\tau(a_n)<d_\tau(a_{n+1})$ for all $\tau\in
T(A),$ $n=1,2,....$ It follows  from Theorem  \ref{Comparison} that
$a_n\lesssim a_{n+1},$ $[a_n]\not=[a_{n+1}],$ $n=1,2,....$  This
ends the proof.
\end{proof}


\begin{thm}\label{Aff=W}
Let $A\in {\cal B}_1$ be {{an infinite dimensional}} unital simple \CA. Then $W(A)$ is
$0$-almost divisible.
\end{thm}

\begin{proof}
Let $a\in M_n(A)_+\setminus \{0\}$ and $k\ge 1$ be an integer. We
need to show that there exists an element $x\in M_{m'}(A)_+$ for
some $m'\ge 1$ such that \beq\label{Aff=W-1} k [x] \le [a] \le (k+1)
[ x ] \eneq in $W(A).$ {
It follows from Theorem \ref{Aff=W-L} that, since $kd_\tau(a)/(k^2+1)\in {\rm
LAff}_b(T(A)),$ there is $x\in M_{2n}(A)_+$ such that
\beq\label{Aff=W-2} d_\tau(x)=kd_\tau(a)/(k^2+1)\rforal \tau\in
T(A). \eneq Then, \beq\label{Aff=W-3}
kd_\tau(x)<d_\tau(a)<(k+1)d_\tau(x)\tforal \tau\in T(A). \eneq It
follows from Theorem \ref{Comparison} that \beq\label{Aff=W-4} k[x]\le [a]\le
(k+1) [x]. \eneq }
\end{proof}

\begin{thm}\label{Zstable}
Let $A\in {\cal B}_1$ be {{an infinite dimensional}} unital separable simple amenable \CA. Then
$A\otimes {\cal Z}\cong A.$
\end{thm}

\begin{proof}
Since $A\in {\cal B}_1,$ $A$ has finite weak tracial nuclear dimension (see 8.1 of \cite{Lin-LAH}).
By Theorem \ref{Comparison}, $A$ has the strict comparison property for positive elements.
 Note that, by Theorem \ref{B1hered}, every unital hereditary \SCA\, of $A$ is in ${\cal B}_1.$
 Thus, by Theorem \ref{Aff=W},  its Cuntz semigroup also has $0$-almost divisibility. It follows from 8.3 of \cite{Lin-LAH} that $A$ is ${\cal Z}$-stable.
\end{proof}

\section{The unitary groups}

\begin{thm} [{{\rm cf.~Theorem 6.5 of \cite{LinTAI}}}]\label{UL1}
Let $K\in \N$ be an integer and let ${\cal B}$ be a class of unital \CA s which has the property
that ${\rm cer}(B)\le K$
for all $B\in {\cal  B}.$
Let $A$ be a unital simple \CA\, which is
{\blue{${\rm {TA}}{\cal B}$ (see Definition \ref{DATD})}} and let $u\in U_0(A).$ {\blue{Suppose  that
$A$ has the cancellation property for projections,  in particular, if $e, q\in A$ are projections and $v\in A$ such that
$v^*v=e$ and $vv^*=q,$ then  there is a $w\in A$ such that $w^*w=1-e$ and $ww^*=1-q.$}}
Then, for any $\ep>0,$  there exist unitaries
$u_1, u_2\in A$ such that $u_1$ has exponential length no more than $2\pi,$ $u_2$ has exponential rank $K,$ and
$$
\|u-u_1u_2\|<\ep.
$$
Moreover, ${\rm cer}(A)\le K+2+\ep.$
\end{thm}

\begin{proof}
{\blue{Suppose that  $A$ is locally approximated by \SCA s  $\{A_n\}$ in ${\cal B}.$
One may write (for some integer $m\ge 1$) $u=\exp(ih_1)\exp(i h_2)\cdots \exp(ih_m),$ where $h_i\in A_{s.a.},$ $i=1,2,...,m.$
\Wlog, we may assume that, for some large $n,$ $h_i\in A_n,$ $i=1,2,...,m.$
In other words, we may assume that $u\in U_0(A_n).$ Then the conclusion of the lemma follows with $u_1=1,$
since ${\rm cer}(A_n)\le K$ as assumed.  This also follows
from the proof of Theorem 6.5 of \cite{LinTAI} with $p=1.$}}

{\blue{Now as in the proof of  \ref{B1sp} we may assume that $A\in {\rm TA}{\cal B}$ and $A$ has property (SP).}}
{{Without loss of generality, we {\blue{may}} assume  {\blue{that}} $A$ is an infinite dimensional \CA. Let $n$  {\blue{
be a positive integer.}}}}

{\blue{As in the beginning of the proof of  \ref{Ltrace},}}
one  can find mutually orthogonal {\blue{non-zero}} positive elements
 $a, b_{1}, b_{2},...,b_{\blue{2(n+1)}}\in A$ such that $a\lesssim b_i$ for all $i\in \{1,2,\cdots, 2{\blue{(n+1)}}\}.$
 Hence, with  a non-zero projection ${\blue{e}}$ such that $e\lesssim a,$
 {\blue{there is a projection $q\in \overline{aAa}$ and $v\in A$ such that
 $v^*v=e$ and $vv^*=q.$
 Note that $b_1,b_2,..., b_{2(n+1)}\in (1-q)A(1-q).$ It follows
 that $2(n+1)[e]\le [1-q],$ or, $2(n+1)[1-p]\le [p],$
  where $p=1-e.$   {\blue{By the assumption that $A$ has cancellation of projections, $[1-q]=[p].$}}
  So, if we apply the property that $A$ is $TA{\cal B},$ as in the proof of Theorem  6.5 of \cite{LinTAI},
   (3) in that
   proof
   {\blue{holds}}
 {\blue{as $(1-p)\lesssim a.$}}  With this fact in mind, the rest of}} proof is exactly the same as that of Theorem 6.5 of \cite{LinTAI}.
\end{proof}

\begin{cor}\label{cerB1}
Any \CA\, in the class ${\cal B}_1$ has exponential rank at most $5+\epsilon$.
\end{cor}
\begin{proof}
By Theorem \ref{2Tg14},  \CA s in ${{\mathcal C}}$ have exponential rank at most $3+\epsilon$.
{\blue{Moreover, by \ref{B1stablerk}, $A$ has stable rank one. Therefore  projections of $A$  have cancellation.}} Therefore, by Theorem \ref{UL1}, any \CA\,  in $\mathcal B_1$ has exponential rank at most $5+\epsilon$.
\end{proof}

\begin{thm}\label{Ulength}
Let $L>0$ be a positive number and let ${\cal B}$ be a class of unital \CA s such that ${\rm cel}(v)\le L$ for
every  unitary $v$ {in their closure of commutator subgroups}.
Let $A$ be a unital simple \CA\,  which is tracially in ${\cal B}$  and let $u\in CU(A).$
{\blue{Suppose that $A$ has cancellation property for projections.}}
Then $u\in U_0(A)$ and ${\rm cel}(u)\le {\blue{2\pi+L+\ep}}.$
\end{thm}

 \begin{proof}
 Let $1>\ep>0.$
 There are $v_1, v_2,....,v_k\in U(A)$ such that
 \beq\label{Ulength-1}
 \|u-v_1v_2\cdots v_k\|<\ep/16
 \eneq
 and $v_i=a_ib_ia_i^*b_i^*,$ where $a_i, b_i\in U(A).$
 Let $N$ be an integer as in Lemma 6.4 of \cite{LinTAI} ({{with  $8\pi  k+\ep$ playing the role of $L$ there}}).
 {{We further assume that
$(8\pi k+1)\pi/N <\frac{\ep}4$.}}
Since $A$ is tracially in ${\cal B},$ there are a projection $p\in A$ and a unital \SCA\, {{$B$}} in ${\cal B}$ with $1_B=p$
 such that
 \beq\label{Ulnegth-2}
 \|a_i-(a_i'\oplus a_i'')\|<\ep/{32}k,\,\,\,\|b_i-(b_i'\oplus b_i'')\|<\ep/{32}k,\,\,\,i=1,2,...,k, \andeqn\\
 \|u-\prod_{i=1}^k (a_i'b_i'(a_i')^*(b_i')^*\oplus a_i''b_i''(a_i'')^*(b_i'')^*\|<\ep/8,
 \eneq
 where $a_i', b_i'\in U((1-p)A(1-p)),$ $a_i'', b_i''\in U_0(B)$ and
 $6N[1-p]\le [p]$ {\blue{(if $A$ is locally approximated by
  \CA s in ${\cal B},$ we can choose $p=1_A;$ otherwise, we can apply
 the argument in  the proof of Theorem \ref{UL1}).}}
Put
 \beq\label{Ulength-3}
 w=\prod_{i=1}^k a_i'{{b_i'}}(a_i')^*(b_i')^*
 \andeqn z=\prod_{i=1}^k a_i''b_i''(a_i''){{(b_i'')^*.}}
 \eneq
 Then $z\in  CU(B).$ Therefore ${\rm cel}_{B}(z)\le L$  in $B\subset pAp.$
 It is standard to show that then
 $$
 a_i'b_i'(a_i')^*(b_i{\blue{'}})^*\oplus (1-p)\oplus (1-p)
 $$
 is in $U_0(M_{{3}}((1-p)A({{1-p}})))$ and it has exponential length no more than $4(2 \pi)+2\ep/16k.$
 This implies
 $$
 {\rm cel}(w\oplus (1-p)\oplus (1-p))\le 8\pi k+\ep/4
 $$
 in $U(M_{{3}}((1-p)A(1-p))).$  {{Since $6N[1-p]\leq [p]$, there is a projection $q\in M_3(pAp)$ such that
 $(1-p)\oplus (1-p)\oplus q$ is Murray von Neumann equivalent to  $ p$ {\blue{and}}
 $N[(1-p)\oplus(1-p)\oplus(1-p)]\leq [q]$. {\blue{View}} $w\oplus(1-p)\oplus(1-p){\blue{\oplus q}}$ as a unitary in $$\big((1-p)\oplus(1-p)\oplus(1-p)\oplus q\big)M_3(A)\big((1-p)\oplus(1-p)\oplus(1-p)\oplus q\big).$$
 {\blue{It}} follows from Lemma 6.4 of \cite{LinTAI} that
 \beq\label{Ulength-2008}
 {\rm cel}(w\oplus 1-p)\oplus(1-p)\oplus q)\le 2\pi+\frac{8\pi k+\ep/4}{N}\pi< 2\pi+\ep/4.
 \eneq }}
 It follows that
 $$
 {\rm cel}((w\oplus p)((1-p)\oplus z))<2\pi+\ep/4+L+\ep/16.
 $$
 {{Therefore}} ${\rm cel}(u)\le 2\pi+L+\ep.$
  \end{proof}

{\blue{
\hspace{-0.2in}{\bf Remark}:
The cancellation property of projections in both Lemma \ref{UL1} and  Theorem \ref{Ulength} can be replaced by
requiring that every unital hereditary \SCA\, of $A$ is in ${\rm{TA}}{\cal B}.$ We do not need this though.
}}

\begin{cor}\label{celB1}
Let $A$ {{be a unital simple \CA\,}} in $\mathcal B_1$, and let $u\in CU(A)$. Then $u\in U_0{{(A)}}$ and $\mathrm{cel}(u)\leq 7\pi$.
\end{cor}
\begin{proof}
{{This}}  follows from Lemma \ref{2Lg9} and Theorem \ref{Ulength}.
\end{proof}

\begin{lem}\label{UCUdiv}
Let $A$ be a unital \CA. Then the group $U_0(A)/CU(A)$ is divisible. Let $U$ be a UHF-algebra {{of infinite type}},  and let $B=A\otimes U.$
Then the group $U_0(B)/CU(B)$ is torsion free.

\end{lem}

\begin{proof}
It is well known that $U_0(A)/CU(A)$ is always divisible. Indeed, pick $\bar{u}\in U_0(A)/CU(A)$ for some $u\in U_0(A)$, and pick $k\in\mathbb N$. Since $u\in U_0(A)$, there are self-adjoint elements $h_1, h_2, ..., h_n\in A$ such that
$u={{\exp(ih_1)\exp(ih_2)\cdots \exp(ih_n)}}).$
The{{n}} the unitary \\
$w:={{\exp(ih_1/k)\exp(ih_2/k)\cdots \exp(ih_n/k)}}$
satisfies $(\bar{w})^k = \bar{u}$ in $U_0(A)/CU(A)$, and thus \linebreak
$U_0(A)/CU(A)$ is divisible.


Now, consider $B=A\otimes U$, where $U$ is an infinite dimensional UHF-algebra.
{\blue{By Corollary. 6.6 and Theorem 6.8 of \cite{RorUHF1}, $B$ is either purely infinite,
or $B$ has stable rank one. If $B$ is purely infinite, by Corollary. 2.7 of \cite{GLX-ER},
$U_0(B)/CU(B)$ is zero, whence $U_0(B)/CU(B)$ is torsion free.  If $B$ has stable rank one, by Corollary. 3.11 of \cite{GLX-ER},
the map from $U_0(B)/CU(B)$ to $U_0(M_n(B))/CU(M_n(B))$ is an isomorphism for all $n\ge 1.$}}
{\blue{Therefore, by}} Theorem 3.2 of \cite{Thomsen-rims}, the group $U_0(B)/CU(B)$ is isomorphic to $\Aff(T(B))/\overline{\rho_B(K_0(B))}$.
 {\color{red}{Let $D=K_0(U).$
Then we may view $D$ as a dense subgroup of $\Q.$  It follows that, for any $x\in K_0(B)$ and $r\in D,$
$r\rho_B(x)\in \rho_B(K_0(B)).$ Therefore $\overline{\rho_B(K_0(B))}$ is divisible.
Hence $U_0(B)/CU(B)$ is torsion free.}}

\end{proof}

\begin{thm}\label{Utorlength}
Let $A$ be a unital \CA\, such that
there is a number $K>0$ such that
${\rm cel}(u)\le K$ for all $u\in CU(A).$
Suppose that  $U_0(A)/ CU(A)$ is torsion free and suppose that
 $u,\, v\in U(A)$ such that $u^*v\in U_0(A).$
Suppose also that there is $k\in \N$ such that ${\rm cel}((u^k)^*v^k){{\le}}L$ for some $L>0.$
Then
\beq\label{Utorl-1}
{\rm cel}(u^*v)\le K+L/k.
\eneq
\end{thm}

\begin{proof}
It follows from \cite{Ringrose-cel} that, for any $\ep>0,$   there are $a_1, a_2,...,a_N\in A_{s.a.}$ such that
\beq\label{Utorl-2}
(u^k)^*v^k=\prod_{j=1}^N \exp(\sqrt{-1}a_j)\andeqn \sum_{j=1}^N \|a_j\|\le L+\ep/2.
\eneq
Choose
$
w=\prod_{j=1}^N\exp(-\sqrt{-1}a_j/k).
$
Then
$
(u^*vw)^k\in CU(A).
$
 Since $U_0(A)/CU(A)$ is assumed to be torsion free, it follows  that
\beq\label{Utorl-3}
u^*vw\in CU(A).
\eneq
Thus,
$
{\rm cel}(u^*vw)\le K.
$
Note that
$
{\rm cel}(w)\le L/k +\ep/2k.
$
It follows that
$$
{\rm cel}(u^*v)\le K+L/k+\ep/2k.
$$
\end{proof}

\begin{cor}\label{Unotrosion}
Let $A$ be a unital simple \CA\, in ${\cal B}_1,$  and let $B=A\otimes U,$ where
$U$ is a UHF-algebra
{{of}}  infinite type.
Then
\begin{enumerate}
\item $U_0(B)/CU(B)$ is torsion free and divisible; and
\item if $u, v\in U(B)$ with ${\rm cel}((u^*)^kv^k)\le L$ for some integer $k>0,$  then
$$
{\rm cel}(u^*v)\le {7}\pi+L/k.
$$
\end{enumerate}
\end{cor}
\begin{proof}
The lemma follows from Lemma \ref{UCUdiv}, Corollary \ref{celB1}, and
Theorem \ref{Utorlength}.
\end{proof}

{\blue{The following  lemma  consists  of some  standard perturbation results  in the same spirit as that of
Lemma \ref{approx-Aug-14-1}. Some of these  have been used
in some of the proofs before. We present them here for convenience.}}

\begin{lem}\label{approx-Aug-14-2}
{{ Let $A$ be a unital $C^*$-algebra in ${\cal B}_1$. Let  $e\in A$ be a {\blue{projection}}  with $[1-e]=  K[e]\,\,\,{\rm in}\,\,\, K_0(A)$ for some positive integer $K$.  Suppose that  $u\in
U_0(eAe)$ and $w=u+(1-e)$ with ${\rm dist}({\bar w},{\bar 1})\leq \eta<2$. Suppose that ${\blue{{\cal F}}}\subset A$ is a finite set, $R$ is a positive integer, and $\ep>0.$
Then there are  a non{\blue{-zero}} projection $p\in A$  and a
\SCA\, $D\in {\cal C}$ with $1_D=p$ such that}}
\begin{enumerate}

\item {{$\|[p,\,x]\|<\ep$ for all $x\in {\cal F}\cup \{u, w, e, (1-e)\},$}}

\item {{$ pxp\in_{\ep} D$ {\blue{for}} $x\in {\cal F}\cup \{u, w, e, (1-e)\},$}}

\item {{there are a  projection $q\in D$, a unitary $z_1\in qDq,$ and $c_1\in CU(D)$
such that $\|q-pep\|<\ep,$ $\|z_1-quq\|<\ep,$ $\|z_1\oplus (p-
q)-pwp\|<\ep,$ and $\|z_1\oplus (p-q)-c_1\|<\ep+{\eta},$}}

\item  {{there are a  projection $q_0\in (1-p)A(1-p)$, a unitary
$z_0\in q_0Aq_0,$ and $c_0\in CU((1-p)A(1-p){\blue{)}}$ such that\\ $\|q_0-(1-p)e(1-p)\|<\ep,$
$\|z_0-(1-p)u(1-p)\|<\ep,$ $\|z_0\oplus
(1-p-q_0)-(1-p)w(1-p)\|<\ep,$  {\blue{and}} $\|z_0\oplus
(1-p-q_0)-c_0\|<\ep+{\eta},$}}

\item {{$[p-q]=K [q]$ in $K_0(D),$ $[(1-p)-q_0]=K[q_0]$ in $K_0(A),$}}

\item {{$R[1-p]<  [p]$ in $K_0(A),$ and }}

\item {{$ {\rm cel}_{(1-p)A(1-p)}(z_0\oplus (1-p-q_0))\le {\rm cel}_A(w)+{\blue{\ep.}}$}}

\end{enumerate}

\end{lem}

\begin{proof}
{{Without {\blue{loss}} of generality,   we may assume that $\ep<\frac{1}{2}{\rm min}{\blue{\{}}\eta, 2-\eta, {\blue{1/2\}}}$.}}

{{Since $[1-e]=  K[e]$ {\blue{and $A$ has stable rank one (see \ref{B1stablerk}),}}  there are mutually orthogonal projections $e_1, e_2,\cdots, e_K\in A$ and partial isometries $s_1, s_2, \cdots, s_K\in A$ satisfy{\blue{ing}}
$s^*_is_i=e,~s_is^*_i=e_i$ for all $i=1,2,\cdots,K$  and
$1-e=\sum_{i=1}^Ke_i$. }}

{{ Let ${\rm cel}_A(w)=L.$}} {\blue{ Then (see the last line of \ref{expLR})
there are $h_1,h_2,...,h_M\in A_{s.a.}$ such that
$
w=\prod_{j=1}^M \exp(i h_j)\andeqn \sum_{j=1}^M \|h_j\|<L+\ep/4.
$}}
{\blue{Choose $\dt_0>0$ such that, if $h'\in A_{s.a.}$ (for any unital \CA\, $A$) with $\|h'\|\le L+1,$
$p'\in A$ is a projection and
$\|p'h'-h'p'\|<\dt_0,$ then
\beq\label{18818-sec11-n2}
\|p'\exp(ih')-p'\exp(ip'h'p')p'\|<\ep/(64N(L+1)(M+1))).
\eneq
}}
{{Since ${\rm dist}({\overline{w}},{\bar 1})\leq \eta<2$, there exist unitaries $\{u_i, v_i\}_{i=1}^N\subset  {\blue{A}}$ such that \beq\label{August-16-2018-b}\|w-\Pi_{i=1}^N u_iv_iu^*_iv^*_i\|<\eta+\frac{\ep}{4}.\eneq}}
{{Let ${\cal F}'$ {\blue{be}} the set
$${\cal F}\cup\{u, w, e, 1-e;~~ e_k, {\blue{s}}_k: {\blue{1\le k\le K;
~~ h_j:  1\le j\le M}};~~ u_i, v_i:{\blue{1\le i\le N}}\}.$$
Let $\dt$ ($<\ep/16N(L+1)(M+1)$)
be the positive number $\dt$   in Lemma \ref{approx-Aug-14-1} for the positive integer $K$ and
${\blue{\min\{\dt_0,}}\ep/16N(L+1)(M+1)
{\blue{\}}}$
(in place of $\ep$).
{\blue{Recall that (by \ref{B1stablerk})
 $A$ has stable rank one.}}
As in the proof of Theorem \ref{UL1}, choose any positive element $a\in A_+$ such that $1-p\lesssim a$ implies $R[1-p]< [p]$ in $K_0(A)$.}}

{{Since $A\in {\cal B}_1$, there are a non{\blue{-zero}} projection $p\in A$ and a
\SCA\, $D\in {\cal C}$  with $1_D=p$ such that{\blue{,}} for all $x\in {\cal F}'${\blue{,}}
\beq\label{August-17-2018}(i)~~\|[p,\,x]\|<\dt,~~~~ ~~(ii)~~  pxp\in_{\dt} D,~~~ \mbox{and}~~~~~~(iii)~~
1-p\lesssim a.
\eneq}}
{{Then (1) follows from $(i)$, and (2) follows from $(ii)$. If we do not require the existence of $c_1$ in  (3) and $c_0$ in  (4) and the estimates involving $c_1$ and $c_0$, then  (3) (the existence of $q$, $z_1$ and all estimates not involving $c_1$) and  (4) (the existence of $q_0$, $z_0$ and all estimates not involving $c_0$)  follow from part (a) and part (b) of Lemma \ref{approx-Aug-14-1}, and (5) follows from part (c) of Lemma {\blue{\ref{approx-Aug-14-1}}}. Furthermore,  (6) follows from $(iii)$ above by the choice of $a$. We emphasize that, since $\dt$ is the number $\dt$ for $\frac{\ep}{16N(L+1)(M+1)}<\frac{\ep}{16N(L+1)}<\frac{\ep}{16N}$ (instead of $\ep$) in Lemma \ref{approx-Aug-14-1}, we
have \beq\label{August-16-2018-c}
\|z_1\oplus (p-
q)-pwp\|<\frac{\ep}{16N},~
\|z_0\oplus
(1-p-q_0)-(1-p)w(1-p)\|<\frac{\ep}{16N(L+1)}.~~~~
\eneq}}
{\blue{It remains to show that (7) holds and the existence of $c_0$ and $c_1.$}}
{\blue{Note, since $\dt<\dt_0,$
\beq
 \|ph_j-h_jp\|<\dt_0,\,\,\, j=1,2,...,M.
\eneq
 }}

{{By  part (a) and {\blue{part}} (b) of Lemma \ref{approx-Aug-14-1}, there are untaries
$\{u'_i,v'_i\}_{i=1}^N\subset D\subset pAp$,  $\{u''_i,v''_i\}_{i=1}^N\subset (1-p)A(1-p)$
such that
\beq\label{August-16-2018-1}
&&\|u'_i-pu_ip\|<\frac{\ep}{16N},~~~~\|v'_i-pv_ip\|<\frac{\ep}{16N},
\\\label{August-16-2018-1a}
&&\|u''_i-(1-p)u_i(1-p)\|<\frac{\ep}{16N},~~ \|v''_i-(1-p)v_i(1-p)\|<\frac{\ep}{16N},
\eneq
for all $i=1,2,\cdots, N$, and such that,}} {\blue{by the choice of $\dt_0$ and by \eqref{18818-sec11-n2},}}
\beq\label{18818-sec11-n1}
{\blue{\|(1-p)w(1-p)-(1-p)(\prod_{j=1}^M\exp(i(1-p)h_j(1-p)))(1-p)\|<\ep/(64N(L+1)).}}
\eneq
{\blue{Put $w_0=(1-p)\prod_{j=1}^M\exp(i(1-p)h_j(1-p))(1-p)$ and $h_j'=(1-p)h_j(1-p),$
$j=1,2,...,M.$ Then $w_0\in U_0((1-p)A(1-p)).$ In $(1-p)A(1-p),$
$w_0=\prod_{j=1}^M \exp(i h_j').$
Note that $\sum_{j=1}^M\|h_j'\|\le \sum_{j=1}^M \|h_j\|< L+\ep/4.$}}
{{From (\ref{August-16-2018-c}) and \eqref{18818-sec11-n1}
we have
$$
\|z_0\oplus
(1-p-q_0)-w_0\|\leq \|z_0\oplus
(1-p-q_0)-(1-p)w(1-p)\|+\|(1-p)w(1-p)-w_0\|<\frac{\ep}{8(L+1)}.
$$
}}
{\blue{There is $h_0'\in ((1-p)A(1-p))_{s.a.}$ such that
\beq
\|h_0\|<2\arcsin(\ep/16(L+1))\andeqn z_0\oplus(1-p-q_0)=w_0\exp(ih_0)=\prod_{j=0}^M\exp(ih_j').
\eneq
Then ${\rm cel}_{(1-p)A(1-p)}(z_0\oplus (1-p-q_0))<L+\ep/4+\ep/2<{\rm cel}(w)+\ep.$ So (7) holds.
}}

{{Let $c_1=\Pi_{i=1}^N u'_iv'_i(u'_i)^*(v'_i)^*$ and $c_0=\Pi_{i=1}^N u''_iv''_i(u''_i)^*(v''_i)^*$.
Then by (\ref{August-16-2018-c}), (\ref{August-16-2018-b}), $(i)$ of (\ref{August-17-2018}), and (\ref{August-16-2018-1}), one {\blue{gets}}
$$
\|z_1\oplus (p-q)-c_1\|<\|z_1\oplus (p-
q)-pwp\|+\|pwp-p(\Pi_{i=1}^N u_iv_iu^*_iv^*_i)p\|~~~~~~~~~~~~~~~~~~~~$$
$$~~~~~~~~~~~~~+\|p(\Pi_{i=1}^N u_iv_iu^*_iv^*_i)p-\Pi_{i=1}^N(pu_ipv_ipu^*_ipv^*_ip)\|+\|\Pi_{i=1}^N(pu_ipv_ipu^*_ipv^*_ip)-c_0\|
$$
$$<\frac{\ep}{16N}+(\eta+\frac{\ep}{4})+(4N\dt)+(4N\frac{\ep}{16N})<\eta+\ep.
$$
Similarly, $\|z_0\oplus
(1-p-q_0)-c_0\|<\ep+{\eta}$. That is, (3) and (4)}} {\blue{hold}}.

\end{proof}

\begin{lem}\label{ph}
Let $K\ge 1$ be an integer.  Let $A$ be a  unital  simple \CA\, in
${\cal B}_1.$  Let  $e\in A$ be a projection and let $u\in
U_0(eAe).$ Suppose that $w=u+(1-e)\in U_0(A)$ and let   $\eta\in (0,2].$ Suppose
that
\beq\label{PH-1} [1-e]\le  K[e]\,\,\,{\text in}\,\,\, K_0(A)
\tand {\rm dist}({\bar w},{\bar 1})\leq \eta. \eneq

Then, if $\eta<2,$ one has
$$
{\rm
dist}({\bar u}, {\bar e})<(K+{{9}}/8)\eta\tand
{\rm cel}_{eAe}(u)<({{{(K+1)}}\pi\over{2}}+1/16)\eta+6\pi,
$$
and if $\eta=2,$ one has
$$
{\rm cel}_{eAe}(u)<{\blue{({9K\over{8}}+1){\rm cel}(w)}}+1/16+6\pi.
$$
\end{lem}

\begin{proof}
We  assume that (\ref{PH-1}) holds. Note that $\eta\le 2.$ Put
$L={\rm cel}(w).$
We first consider the case that $\eta<2.$ There is a projection
$e'\in M_m(A)$ (for some integer $m$) such that
$$
[(1-e)+e']=K[e].
$$
{{Note $M_m(A)\in {\cal B}_1.$}}
Replacing $A$ by
$(1_A+e')M_m(A)(1_A+e')$ {{(which is in ${\cal B}_1$ by Theorem  \ref{B1hered})}} and $w$ by $w+e',$ without loss of
generality, we may now assume that
\beq\label{PH-10}
[1-e]=K[e]\andeqn {\rm dist}({\bar w},{\bar 1})<\eta.
\eneq
There is $R_1>1$ such that
$\max\{L/R_1,2/R_1,\eta\pi/R_1\}<\min\{\eta/64, 1/16\pi\}.$

For
any
$\eta/32(K+1)^2\pi>\ep>0$ with $\ep+\eta<2,$  since
$gTR(A)\le 1,$ {{by Lemma \ref{approx-Aug-14-2},}} there exist a {{non-zer{\blue{o}}}} projection $p\in A$ and a \SCA\, $D\in
{\cal C}$ with $1_D=p$ such that
\begin{enumerate}

\item $\|[p,\,x]\|<\ep$ for $x\in \{u, w, e, (1-e)\},$

\item $pwp, pup, pep, p(1-e)p\in_{\ep} D,$

\item there are a  projection $q\in D,$  a unitary $z_1\in qDq,$ and $c_1\in CU(D)$
such that $\|q-pep\|<\ep,$ $\|z_1-quq\|<\ep,$ $\|z_1\oplus (p-
q)-pwp\|<\ep,$ and $\|z_1\oplus (p-q)-c_1\|<\ep+{\eta},$

\item  there are a  projection $q_0\in (1-p)A(1-p)$ and a unitary
$z_0\in q_0Aq_0$ such that\\ $\|q_0-(1-p)e(1-p)\|<\ep,$
$\|z_0-(1-p)u(1-p)\|<\ep,$ $\|z_0\oplus
(1-p-q_0)-(1-p)w(1-p)\|<\ep$ {\blue{and}} $\|z_0\oplus
(1-p-q_0)-c_0\|<\ep+{\eta},$

\item $[p-q]=K [q]$ in $K_0(D),$ $[(1-p)-q_0]=K[q_0]$ in $K_0(A),$

\item $2(K+1)R_1[1-p]<  [p]$ in $K_0(A),$ and

\item $ {\rm cel}_{(1-p)A(1-p)}(z_0\oplus (1-p-q_0))\le L+\ep,$

\end{enumerate}
where
 $c_1\in CU(D)$ and
$c_0\in CU((1-p)A(1-p)).$

By  Lemma \ref{2Lg8}, one has that ${\rm det}({{\psi_D(c_1)}})=1$
{{for every irreducible representation $\psi_D$ of $D.$}} Since $\ep+\eta<2,$
there is $h\in D_{s.a}$ with $\|h\| <
2\arcsin((\ep+\eta)/2)$ such that (by (3) above)
\beq\label{ph-5}
(z_1\oplus (p-q))\exp(ih)=c_1.
\eneq
{\blue{By (5) above, since $D$ has stable rank one, we may write $D=M_{K+1}(D_0),$
where $D_0\cong qDq.$ In particular, $D_0\in {\cal C}$ (see \ref{cut-full-pj}).
Let $\{e_{ij}\}$ be a system of matrix units for $M_{K+1}.$
Define $h_i=e_{ii}he_{ii},$ $i=1,2,...,K+1,$ and $h_0=\sum_{j=1}^{K+1}e_{1i}he_{i1}.$
Then $\|h_0\|<2(K+1)\arcsin((\ep+\eta)/2).$
Note $\tau(h_0)=\tau(\sum_{j=1}^{K+1}h_j)=\tau(h)$ for all $\tau\in T(D).$
We may identify   $h_0$ with  an element in $(qDq)_{s.a.}.$ Put
$\omega_0=\exp(i h_0)\in  U_0(qDq).$ Then, ${\rm det}(\psi_D((\omega_0\oplus (p-q))\exp(-ih))=1$ for every irreducible
representation
$\psi_D$ of $D.$  It follows by \eqref{ph-5} that,  for every irreducible representation $\psi_D$ of $D,$
\beq
{\rm det}(\psi_D(z_1\omega_0\oplus (p-q)))=1.
\eneq
 By Lemma \ref{2Lg8}, this implies
$z_1w_0\oplus (p-q)\in CU(M_{K+1}(qDq)).$ Since $qDq$ has stable rank one (by Proposition \ref{2pg3}), by Corollary 3.11 of \cite{GLX-ER}, $zw_0\in CU(qDq).$}}
It follows that, if $2{{(K+1)}}\arcsin({\ep+\eta\over{2}})< \pi,$
\beq\label{ph-7}
{\blue{{\rm dist}(\overline{z_1}, \overline{q})={\rm dist}(\overline{\omega_0}, \overline{q})<2\sin((K+1)\arcsin({\ep+\eta\over{2}}))
\le (K+1)(\ep+\eta)}}
\eneq
{\blue{(see \eqref{dddcu}).}} If $2{{(K+1)}}\arcsin({\ep+\eta\over{2}})\ge \pi,$ then
$
2{{(K+1)}}({\ep+\eta\over{2}}){\pi\over{2}}\ge \pi.
$
It follows that
\beq\label{ph-8-}
{{(K+1)}}(\ep+\eta)\ge 2\ge {\rm dist}(\overline{z_1}, {\overline{q}}).
\eneq
By combining both (\ref{ph-8-}) and (\ref{ph-8-}), one obtains that
\beq\label{ph-8+}
{\rm dist}(\overline{z_1}, {\overline{q}})\le {{(K+1)}}(\ep+\eta) \le
{{(K+1)}}\eta+{\eta\over{32(K+1)\pi}}.
\eneq
{{It}} follows from  Lemma \ref{2Lg9} that
\beq\label{ph-8+1}
{\rm cel}_{qDq}(z_1)\le
{{(K+1)}}(\ep+\eta){\pi\over{2}}+4\pi
\le
({{(K+1)}}{\pi\over{2}}+{1\over{64(K+1)}})\eta +4\pi .
\eneq
By (5) and (6) above,
$$
(K+1)[q]=[p-q]+[q]=[p]>2(K+1)R_1[1-p].
$$
 {{By Theorem \ref{Comparison},}} $K_0(A)$ is weakly unperforated. {{Hence,}}
\beq\label{ph-9}
2R_1[1-p]<[q].
\eneq

There is a unitary $v\in A$ such that
$v^*(1-p-q_0)v\le q.$
Put $v_1=q_0\oplus (1-p-q_0)v.$ Then
\beq\label{ph-11}
v_1^*(z_0\oplus (1-p-q_0))v_1=z_0\oplus v^*(1-p-q_0)v.
\eneq
Note that {\blue{(by (4))}}
\beq\label{ph-10+n}
\|(z_0\oplus v^*(1-p-q_0)v)v_1^*c_0^*v_1-q_0\oplus
v^*(1-p-q_0)v\|<\ep+\eta.
\eneq
Moreover, by (7) above,
\beq\label{ph-10+}
{\rm cel}(z_0\oplus v^*(1-p-q_0)v)\le L+\ep.
\eneq
It follows from (\ref{ph-9}) and Lemma 6.4 of \cite{LinTAI} that
\beq\label{ph-10+1}
{\rm cel}_{(q_0+q)A(q_0+q)}(z_0\oplus q)\le 2\pi+(L+\ep)/R_1.
\eneq
Therefore, on combining {{this with}} (\ref{ph-8+1}),
\beq\label{ph-10+2}
{\rm cel}_{(q_0+q)A(q_0+q)}(z_0+z_1)\le
2\pi+(L+\ep)/R_1+({{(K+1)}}{\pi\over{2}}+{1\over{64(K+1)}})\eta +{{4}}\pi.
\eneq
By (\ref{ph-10+}), (\ref{ph-9}), and Lemma 3.1 of \cite{Lin-hmtp},
in $U_0((q_0+q)A(q_0+q))/CU((q_0+q)A(q_0+q)),$
\beq\label{ph-13}
{\rm dist}(\overline{z_0+q},
{\overline{q_0+q}})<{(L+\ep)\over{R_1}}.
\eneq
Therefore, by {\blue{(\ref{ph-13}), (\ref{ph-8+})}} {\blue{(and by the line below \eqref{dddcu}),}}
\beq\label{ph-14}
{\rm dist}({\overline{z_0\oplus z_1}},
{\overline{q_0+q}})<{(L+\ep)\over{R_1}}+{{(K+1)}}\eta+{\eta\over{32(K+1)\pi}}<
(K+{{\frac{17}{16}}})\eta.
\eneq
We note that
\beq
\|e-(q_0+q)\|<2\ep\andeqn \|u-(z_0+z_1)\|<2\ep.
\eneq
It follows that
\beq\label{ph-15}
{\rm dist}({\bar u}, {\bar e})<4\ep +(K+{{\frac{17}{16}}})\eta<(K+{{\frac{9}{8}}})\eta.
\eneq

Similarly, by (\ref{ph-10+2}),
\beq\label{ph-15+1}
{\rm cel}_{eAe}(u)&\le& 4\ep\pi+2\pi+(L+\ep)/R_1+({{(K+1)}}{\pi\over{2}}+{1\over{64(K+1)}})\eta
+4\pi\\
&<&({{(K+1)}}{\pi\over{2}}+1/16)\eta+6\pi.
\eneq

This proves the case that $\eta<2.$

Now suppose that $\eta=2.$
Define {\blue{an integer}} $R=[{\rm cel}(w)+1]{{\ge 4}}.$  Note that ${\rm cel}(w)/{R}<1.$
There is a projection $e'\in M_{{R'}}(A)$ {{(for some integer $R'$)}} such that
$$
{\blue{[e']=2R[e]+2R[1-e]=2R[1_A].}}
$$
It follows from Lemma 3.1 of \cite{Lin-hmtp} that
\beq\label{PH-111}
{\rm dist}(\overline{w\oplus e'},{\overline{ 1_A+e'}})<{{\rm
cel}(w)\over{{R}}}.
\eneq
Put $K_1{\blue{=2R(K+1)+K.}}$
{\blue{Then
\beq
([1-e]+[e'])\le K[e]+2R[e]+2R[1-e]\le (K+2R+2RK)[e]=K_1[e].
\eneq
}} 
It follows from the first part of the lemma that
\beq\nonumber
\hspace{-0.3in}&&{\rm cel}_{eAe}(u)<
({{{(K_1+1)}}\pi\over{2}}+{1\over{16}}){{\rm cel}(w)\over{R}}+6\pi
{\blue{=(K+1+K/2R+1/16R){\rm cel}(w)+6\pi}}\\\nonumber
&&\hspace{-0.2in}\le {\blue{(K+ 1+K/8){\rm cel}(w)+(1/16)+6\pi=(9K/8+1){\rm cel}(w)+1/16+6\pi.}}
\eneq

\end{proof}

%

%

\begin{thm}(Theorem 4.6 of \cite{GLX-ER})\label{UCUiso}
{Let $A$ be a unital  simple \CA\, of stable rank one} and let
$e\in A$ be a non-zero projection. Then the map $u\mapsto u+(1-e)$
induces an isomorphism  from $U(eAe)/CU(eAe)$ onto $U(A)/CU(A).$
\end{thm}

\begin{proof}
{{This is }} Theorem 4.6 of \cite{GLX-ER}.
\end{proof}

\begin{cor}\label{c1}
Let {$A$ be a unital  simple \CA\, of stable rank one.} Then
the map $j: a\to {\rm diag}(a, \overbrace{1,1,...,1}^{{n-1}})$ from $A$ to
$M_n(A)$ induces an isomorphism from $U(A)/CU(A)$ onto
$U(M_n(A))/CU(M_n(A))$ for any integer $n\ge 1.$
\end{cor}

\begin{proof}
This follows from \ref{UCUiso} but also follows from 3.11 of \cite{GLX-ER}.
\end{proof}

\section{A Uniqueness Theorem for \CA s in ${\cal B}_1$ }

{{The following is {\blue{taken from}} Definition 2.1 of \cite{GL-almost-map}.}}

\begin{df}\label{K-rank-2018}
{{Let $ {\bf r}: \N\to \N\cup\{0\},$ ${\bf T}: \N^2\to \N,$  and ${\bf E}: \R_+\times \N\to \R_+$ be  maps, and ${\bf k},{\bf R}\in \N$ be two positive integers. Let $A$ be a unital $C^*$-algebra.}}

{{(1)~We say that $A$ has $K_0$-${\bf r}$-cancellation if
$p\oplus {\bf 1}_{M_{{\bf r}(n)}(A)}\sim q\oplus {\bf 1}_{M_{{\bf r}(n)}(A)}$ for any two projections $p, q\in M_n(A)$ with
${\blue{[p]=[q]}}$ in $K_0(A)$.  }}

{{(2)~We say that $A$ has $K_1$-${\bf r}$-cancellation if $u\oplus {\bf 1}_{M_{{\bf r}(n)}(A)}$ and $v\oplus {\bf 1}_{M_{{\bf r}(n)}(A)}$ are in the same connected component of $U(M_{{\bf r}(n)}(A))$  for any pair $u, v\in M_n(A)$ with $[u]=[v]$ in $K_1(A)$. }}

{{(3)~We say that $A$ has $K_1$-stable rank at most ${\bf k}$ if $U(M_{\bf k}(A))$ is mapped surjectively to $K_1(A)$.}}

{{(4)~We say that $A$ has stable exponential  rank at most ${\bf R}$ if ${\blue{{\rm cer}(M_m(A))}}\leq {\bf R}$ for all $m$.}}

{{(5)~We say that $A$ has $K_0$-divisible rank ${\bf T}$ if for any $x\in K_0(A)$ {\blue{and any pair}}
$(n,k)\in \N\times \N,$
$$-n[{\bf 1}_A]\leq kx \leq n[{\bf 1}_A]$$ implies that
$$-{\bf T}(n, k)[{\bf 1}_A]\leq x \leq {\bf T}(n, k)[{\bf 1}_A].$$}}

{{(6)~We say that $A$ has $K_1$-divisible rank ${\bf T}$ if $k[x]=[u]$ in $K_1(A) $ for some unitary $u\in M_n(A)$ implies that $[x]=[v]$ for some unitary $v\in M_{{\bf T}(n, k)}(A)$. }}

{{(7)~  We say that $A$ has exponential length divisible rank ${\bf E}$ if $u\in U_0(M_n(A))$ with ${\rm cel}(u^k)\leq L$ implies that ${\rm cel}(u)\leq {\bf E}(L,k)$.}}

{{(8)~Let ${\bf C}_{{\bf R, r, T, E}}$ be the class of unital $C^*$ algebras which have  $K_i$-${\bf r}$-cancellation, $K_i$-divisible rank ${\bf T}$,
exponential length divisible rank ${\bf E},$  and stable exponential  rank  at most ${\bf k}$.}}

\end{df}

\begin{rem}\label{K-rank-rem}

   {{If $A$ has stable rank one, then,  by a standard result of Rieffel (\cite{Rf}),   $A$ has $K_i$-$0$-cancellation,
   $K_1$-stable rank $1$ and $K_1$-divisible rank ${\bf T}$ for any ${\bf T}$ (since any element in $K_1(A)$ can be realized by a unitary $u\in A$).  If $K_0(A)$ is weakly unperforated, then by Proposition 2.2 (5) of \cite{GL-almost-map}, $A$ has $K_0$-divisible rank ${\bf T}$ for ${\bf T}(n, k)=n+1$. }}

\end{rem}

The following  theorem follows from Theorem 7.1 of \cite{Lin-hmtp}.
{\blue{We refer to \ref{K-rank-2018} for some of the notation used below.
Recall that,  as before, if $L: A\to B$ is a map, we continue to use $L$ for
the extension $L\otimes {\rm id}_{M_n}: A\otimes M_n\to B\otimes M_n$ (see also \ref{Dmap}).
Recall also  that, when $A$ is unital,  in the following statement and its proof, we always assume
that $\phi(1_A)$ and $\psi(1_A)$ are projections (see the 5th paragraph of \ref{KLtriple}).
We also use the convention $\la \phi(v)\ra$ for $\la \phi(v)\ra+(1-\phi(1_A))$ when
$\phi$ is a approximately multiplicative map and $v$ is a unitary.}}

\begin{thm}\label{Suni} {{\rm  (cf.~Theorem 7.1 of \cite{Lin-hmtp})}}
Let  $A$ be a unital separable amenable \CA\, which satisfies the UCT, {{let}} $T\times N: A_+\setminus \{0\}
\to \R_+\setminus \{0\}\times \N$ be a map  and let ${\bf L}: U(M_{\infty}(A))\to
\R_+$ be another map. For any $\ep>0$ and any finite subset ${\cal F}\subset A,$ there exist $\dt>0,$ a finite subset
${\cal G}\subset A,$ a finite subset ${\cal H}\subset A_+\setminus \{0\},$ a finite subset ${\cal U}\subset
{\bigcup_{m=1}U(M_m(A))},$ a finite subset ${\cal P}\subset \underline{K}(A),$ and an integer $n>0$ satisfying the following
condition:
for any unital separable simple \CA\, $C$ in ${\cal B}_1,$ if $\phi, \psi, \sigma: A\to B
=C\otimes U,$ where
$U$ is a UHF-algebra of infinite type,
are three ${\cal G}$-$\dt$-multiplicative \cp s
and  $\sigma$ is  unital and $T\times N$-${\cal H}$-full ({{see}} \ref{Dfull}) with
{the properties that}
\beq\label{Suni-1}
[\phi]|_{\cal P}=[\psi]|_{\cal P}\tand
{\rm cel}(\langle \psi ({{v}})\rangle^*\langle \phi(v)\rangle)\le {\bf L}(v)
\eneq
for all $v\in {\cal U},$
{{then}} there {{exists}}
a unitary $u\in M_{n+1}(B)$ such that
$$
\|u^*{\rm diag}(\phi(a), \overline{\sigma}(a))u-{\rm diag}(\psi(a), \overline{\sigma}(a))\|<\ep
$$
for all $a\in {\cal F},$ where
$$
\overline{\sigma}(a)={\rm diag}(\overbrace{\sigma(a),\sigma(a),...,\sigma(a)}^n)=\sigma(a)\otimes 1_n \tforal a\in A.
$$

{{{\blue{Moreover, the conclusion  holds if one only assumes}}
$B\in {\bf C}_{{\bf k, r, T, E}}$ for certain ${\bf k, r, T, E}$ as in Definition \ref{K-rank-2018}.}}

{{Furthermore,  if $K_1(A)$ is finitely generated, and  ${\cal U}_0\in U(M_m(A))$ is any finite set such that $\{[u]\in K_1(A):u\in {\cal U}_0\}$ generates the whole group $K_1(A)$, then, in the above statement,  we can always choose the set ${\cal U}$ to be this fixed set ${\cal U}_0$ for any ${\cal F}$ and $\ep$.  }}

\end{thm}

\begin{proof}
{{By Theorem \ref{B1stablerk},}} $B\,{{(=C\otimes U\in {\cal B}_1)}}$  has stable rank one{{. By Corollary \ref{cerB1}, }}  ${\rm cer}(M_m(B))\le  6$ {{ for all $m\in \N$. By part (2) of Corollary of \ref{Unotrosion}, $A$ has}} exponential length divisible rank
${{\bf E}}(L,k)=7\pi +L/k$ {{(see (7) of \ref{K-rank-2018}). By Theorem \ref{Comparison}, }}  $K_0(B)$ {{is}} weakly unperforated (in particular, {{by Remark \ref{K-rank-rem},
$B$}} has $K_0$-divisible rank ${{\bf T}}(L,k)=L+1$). {\blue{In other words, $B=C\otimes U\in {\bf C}_{{\bf R, r, T, E}}$ for ${\bf R}=6$, ${\bf r}(n)=0$, ${\bf T}(L,k)=L+1$ and ${\bf E}(L,k)=7\pi +L/k.$ Therefore, it suffices to prove the ``Moreover" part
(and ``Furthermore" part) of
the statement.}}
{{Note the case that both $\phi$ and $\psi$ are unital follows from
Theorem 7.1 of \cite{Lin-hmtp}.}} {\blue{However, the following  argument also works for the non-unital case.
Denote by ${\bf C}$  the class ${\bf C}_{{\bf R, r, T, E}}$ of unital \CA s  for ${\bf R}=6$, ${\bf r}(n)=0$, ${\bf T}(L,k)=L+1,$ and ${\bf E}(L,k)=7\pi +L/k.$
}}


Suppose that the conclusion of the theorem is false.  Then there exists  $\ep_0>0$ and a finite subset ${\cal F}\subset A$
such that there are a sequence of positive numbers $\{\dt_n\}$ with $\dt_n\searrow 0,$ an increasing sequence
$\{{\cal G}_n\}$ of finite subsets of $A$ such that $\bigcup_n {\cal G}_n$ is dense in $A,$
an increasing sequence $\{{\cal P}_n\}$ of finite subsets of $\underline{K}(A)$ with
$\bigcup_{n=1}{\cal P}_n=\underline{K}(A),$  an increasing sequence
$\{{\cal U}_n\}$of finite subsets  of $U(M_{\infty}(A))$ such that {{$\{[u]\in K_1(A): u\in \bigcup_{n=1}^{\infty}{\cal U}_n\}$ generates $K_1(A)$}}
{{(for the case that $K_1(A)$  is  generated by the single  finite set ${\cal U}_0\in M_m(A)$, we set ${\cal U}_n={\cal U}_0$)}}, an increasing sequence $\{{\cal H}_n\}$ of finite subsets
of  $A_+^{\bf 1}\setminus \{0\}$ such that if $a\in {\cal H}_n$
and $f_{1/2}(a)\not=0,$ 
then
$f_{1/2}(a)\in {\cal H}_{n+1},$ and $\bigcup_{n=1}{\cal H}_n$ is dense
in $A_+^{\bf 1},$  a {{sequence}} of positive integers $\{k(n)\}$ with
$\lim_{n\to\infty} k(n)=\infty$,  a sequence of unital \CA s $B_n\in {\bf C},$
 sequences of ${\cal G}_n$-$\dt_n$-multiplicative
\morp s $\phi_{n}, \psi_{n}: A\to B_n,$
such that
\beq\label{stableun2-1}
[\phi_{n}]|_{{\cal P}_n}=[\psi_{n}]|_{{\cal P}_n}\andeqn
{\rm cel}(\langle \phi_{n}(u)\rangle \langle \psi_{n}(u)\rangle^{{*}})\le {\bf L}(u)
\eneq
for all $u\in {\cal U}_n,$
and a sequence of unital ${\cal G}_n$-$\dt_n$-multiplicative 
{{ completely}} positive linear map{{s}}, $\sigma_n: A\to
B_n,$ which {{are}}  ${{T\times N}}$-${\cal H}_n$-full, {{such that for any $n\in \N$,}}
\beq\label{stableun2-2}
&&\hspace{-0.3in} \inf{{\Big\{}}}\{\sup\|{\blue{v_n}}^*{\diag}( \phi_{n}(a), S_n(a)){\blue{v_n}}-{\diag}(\psi_n(a), S_n(a))\|\!: a\!\in\! {\cal F}\}\!: {\blue{v_n}}\!\in\! {{U(}}M_{{k(n)+1}}(B_n){{{)}}{{\Big\}}} \ge \ep_0, ~~~~~~~~~~~~~~
\eneq
where
$$
S_n(a)={\diag}(\overbrace{\sigma_n(a),\sigma_n(a),...,\sigma_n(a)}^{k(n)})=\sigma_n(a)\otimes 1_{k(n)} \rforal a\in A.
$$

Let $C_0=\bigoplus_{n=1}^{\infty}B_n,$ $C=\prod_{n=1}^{\infty}B_n,$ $Q(C)=C/C_0,$ and
$\pi: C\to Q(C)$ {{be}} the quotient map.  Define $\Phi, \Psi, S: A\to C$ by
$\Phi(a)=\{\phi_n(a)\},$ $\Psi(a)=\{\psi_n(a)\},$  and  $S(a)=\{\sigma_n(a)\}$ for
all $a\in A.$ Note that $\pi\circ \Phi,$ $\pi\circ \Psi$ and $\pi\circ S$ are \hm s.

As in the proof of 7.1 of \cite{Lin-hmtp}, since $B_n\in {\bf C},$ one computes that
\beq\label{stableun2-3}
[\pi\circ \Phi]=[\pi\circ \Psi]\,\,\,{\rm in}\,\,\, KL(A, Q(C))
\eneq
{{as follows.}}

{{For each fixed element $x\in K_1(A)$,  there are  integers $m$, $r$ such that $ {\cal U}_m\subset  M_{r}(A)$ and such that
$x$ is in the group generated by $\{[u]\in K_1(A): u\in {\cal U}_m\}$. {\blue{In the case that
$K_1(A)$ is finitely generated, $x$ is in the group generated by ${\cal U}_0.$}}
 Then, for $n\geq m$, we have that
${\rm cel}(\langle \phi_{n}(u)\rangle \langle \psi_{n}(u)\rangle^*)\le {\bf L}(u)$. Hence by Lemma 1.1 of \cite{GL-almost-map},
 there are a constant ${\bf K}=K({\bf L}(u))$ {\blue{and}}   $U_n(t)\in M_{r}(B_n)$
 such that $U_n(0)=\langle \phi_{n}(u)\rangle$,
 $U_n(1)=\langle \psi_{n}(u)\rangle,$ and
 $\|U_n(t)-U_n(t')\|\leq {\bf K}|t-t'|$. Therefore, $\{U_n\}_{n=m}^{\infty}\in C([0,1],M_{r}( \prod_{n=m}^{\infty}B_n))$ and consequently
 $$
 (\pi\circ \Phi)_{*1} ([u])= [\{\langle \phi_{n}(u)\rangle\}_n]=[\{\langle \psi_{n}(u)\rangle\}_n]=(\pi\circ \Psi)_{*1} ([u])\,\,\,{\rm in}\,\, K_1(Q(C)) ~~\mbox{for all}~~u\in {\cal U}_m.
 $$
 {\blue{Thus,}} $(\pi\circ \Phi)_{*1} (x)= (\pi\circ \Phi)_{*1} (x)$.}}  {\blue{Note that this includes the case that ${\cal U}_m={\cal U}_0$
 when $K_1(A)$ is generated by $\{[u]: u\in {\cal U}_0\}.$}}
 {{In other words, $[\pi\circ \Phi]_{K_1(A)}= [\pi\circ \Psi]_{K_1(A)}: K_1(A) \to K_1(Q(C))$.}}

 {{Since all $B_n$ have  $K_0$-${\bf r}$-cancellation, by (1)
 (see (5) also) of Proposition 2.1 of page 992 in \cite{GL-almost-map},
 $K_0(\prod B_n)=\prod_b K_0(B_n)$ (see page 990 of \cite{GL-almost-map} for the definition of $\prod_b $).
 Hence by 5.1 of \cite{Gong-AH},
$K_0(Q(C))=\prod_b K_0(B_n)/ {\blue{\bigoplus}} K_0(B_n)$.
 Since all $B_n$ have  $K_i$-${\bf r}$-cancellation, $K_i$-divisible rank ${\bf T}$,
exponential length divisible rank ${\bf E},$  and  stable exponential  rank  at most ${\bf R}$, by Theorem 2.1 and (3) of Proposition 2.1  on page 994 of \cite{GL-almost-map}, we have
$$K_i\Big(\prod_n B_n, \Z/k\Z\Big)\subset \prod_n K_i(B_n,\Z/k\Z), K_i(Q(C),\Z/k\Z)\subset \prod_n K_i(B_n,\Z/k\Z)/\bigoplus_n K_i(B_n,\Z/k\Z).$$
That is,  (2) of Proposition 2.2 of \cite{GL-almost-map} holds for $B_n$, {\blue{even though}} we do not assume that $B_n$ have stable rank one. Consequently, {\blue{by \eqref{stableun2-1},}}
we have $[\pi\circ \Phi]_{K_0(A)}= [\pi\circ \Psi]_{K_0(A)}: K_0(A) \to K_0(Q(C)),$ and, for $i=0,1,$ $[\pi\circ \Phi]_{K_i(A, \Z/k\Z)}= [\pi\circ \Psi]_{K_i(A, \Z/k\Z)}: K_i(A, \Z/k\Z) \to K_i(Q(C),\Z/k\Z)$, since
$\bigcup_{n=1}{\cal P}_n=\underline{K}(A)$ and ${\cal P}_n\subset {\cal P}_{n+1},$ $n=1,2,....$}}

{{Since $A$ satisfies UCT, we have  $[\pi\circ \Phi]=[\pi\circ \Psi]\,\,\,{\rm in}\,\,\, KL(A, Q(C)),$ {\blue{i.e., \eqref{stableun2-3} holds.}}}}

{{ By Proposition \ref{full-2018-sept}, }}
 $\pi\circ S$ is a full \hm.
It follows from Theorem \ref{Lnuct} that there exist an integer $K\ge 1$ and a unitary $U\in M_{K+1}(Q(C))$
such that
\beq\label{stableun2-4}
\|U^*{\diag}(\pi\circ \Phi(a), \Sigma(a))U-{\diag}(\pi\circ \Psi(a),{{\Sigma(a)}})\|<\ep_0/4\rforal a\in {\cal F}{{,}}
\eneq
{{where $\Sigma(a)={\diag}(\overbrace{(\pi\circ S)(a), (\pi\circ S)(a),...,(\pi\circ S)(a)}^K)=(\pi\circ S)(a)\otimes 1_K.$}} It follows that there exist a
{\blue{unitary}} $V=\{v_n\}\in C$  and an integer $N\ge 1$ such that, for any $n\ge N$ with
$k(n)\ge K,$
\beq\label{stableun2-5}
\|v_n^*{\diag}(\phi_n(a), {\overline{\sigma}_n}(a))v_n-{\diag}(\psi_n(a), {\overline{\sigma}_n}(a))\|<\ep_0/2
\eneq
for all $a\in {\cal F},$ where
$$
{\overline{\sigma}_n}(a)={\diag}(\overbrace{\sigma_n(a), \sigma_n(a),...,\sigma_n(a)}^K)\rforal a\in A.
$$
This contradicts  (\ref{stableun2-2}).
\end{proof}

\begin{rem}\label{Rsuni}

(1)
{{N}}ote
that $\phi$ and $\psi$ are not assumed to be unital. Thus Theorem \ref{Lauct2} can also  be  viewed
as a special case of Theorem \ref{Suni}.

(2) \, \,{{If $A$ has $K_1$-stable rank
${\blue{k}},$ then in the proof of Theorem \ref{Suni}, we can choose the matrix size $r={\blue k}$ (for $u\in {\cal U}_m\subset M_r(A)$ to represent an element $x\in K_1(A)$). }} Suppose that there exists an integer $n_0\ge 1$ such that
$U(M_{n_0}(A))/{U_0(M_{n_0}(A))}\to U(M_{n_0+k}(A))/U_0(M_{n_0+k}(A))$ is an isomorphism
for all $k\ge 1.$ Then $A$ {\blue{has}}  {{$K_1$-stable rank $n_0$, and the map}} ${\bf L}$ may be replaced by a map from $U(M_{n_0}(A))$ to $\R_+$,  and $\mathcal U$ can be chosen in {$U(M_{n_0}(A))$}.

{
 Moreover, {{for the case $B=C\otimes U$ for $C\in {\cal B}_1$ in the theorem,}}
 {\blue{by Corollary \ref{celB1},}}
 the condition that
${\rm cel}(\la \phi(u)\ra \la \psi(u){{\ra}}^*)\le {\bf L}(u)$ may,  in practice,  be replaced
by the  stronger condition that,  for all ${{\bar u}}\in \overline{{\cal U}},$
\beq\label{Rsuni-1}
{\rm dist}(\phi^{\ddag}({{\bar u}}), \psi^{\ddag}({{\bar u}}))<{{\bf L}},
\eneq
where $\overline{{\cal U}}\subset U(M_m(A))/CU(M_m(A))$ is a finite subset
and where ${{\bf L}}<2$ is a given constant,
{\blue{and where $\phi^{\ddag}$ and $\psi^{\ddag}$ are maps
from  $U(M_m(A))/CU(M_m(A))$ to $U(M_m(B))/CU(M_m(B))$}} (see also \ref{DLddag}).}

To see this, let ${\cal U}$ be a finite subset of $U(M_m(A))$  for
some large $m$ whose image in the group $U(M_m(A))/CU(M_m(A))$ is
$\overline{{\cal U}}.$ Then (\ref{Rsuni-1}) implies that
\beq\label{Rsuni-2} \|\la \phi(u)\ra\la \psi(u)\ra^{{*}} -v\|<{{2}}
\eneq
for some $v\in CU(M_m({\blue{B}})),$ provided that $\dt$ is
sufficiently small and ${\cal G}$ is sufficiently large. {{By \ref{celB1},}} ${\rm
cel}(v)\le {7}\pi.$ {{Also from (\ref{Rsuni-2}), one gets $(\la \phi(u)\ra\la \psi(u)\ra^*)v^*={\rm exp}(i h)$ with $\|h\|<\pi$. W}}e conclude that
\beq\label{Rsuni-3}
{\rm cel}(\la \phi(u)\ra \la \psi(u^*)\ra)\le \pi+{{{7}\pi}}
\eneq for
all $u\in {\cal U},$ {and take ${\bf L}: U_\infty(A) \to \R_+$ to
be constant equal to ${{8\pi}}$.}
Furthermore,  we may assume $$\overline{{\cal U}}\subset U(M_{n_0}(A))/CU(M_{n_0}(C)),$$ if
$K_1(C)=U(M_{n_0}(C))/U_0(M_{n_0}(C)).$

\end{rem}

\begin{lem}\label{Ldet}
Let $A$ be a unital separable simple \CA\, with $T(A)\not=\emptyset.$  {{Then t}}here exists an order preserving map
$\Delta_0: A_+^{q, {\bf 1}}\setminus \{0\}\to (0,1)$ {{with}} the following {{property}}:
For any finite subset ${\cal H}\subset A_+^{\bf 1}\setminus \{0\},$ there {{exist}} a finite subset ${\cal G}\subset
A$  and $\dt>0$ such that, for any unital \CA\, $B$ with $T(B)\not=\emptyset$ and any unital ${\cal G}$-$\dt$-multiplicative {{completely}} positive linear map
$\phi: A\to B,$ one has
\beq\label{Ldet-1}
\tau\circ \phi(h)\ge \Delta_0(\hat{h})/2\tforal h\in {\cal H}
\eneq
and for all $\tau\in T(B),$ and moreover, one may assume that $\Delta_0(\widehat{1_A})=3/4$.
\end{lem}

\begin{proof}

Define, for each $h\in A_+^{\bf 1}\setminus \{0\},$
\beq\label{Ldet-2}
\Delta_0(\hat{h})=\min\{3/4, \inf\{\tau(h): \tau\in T(A)\}\}.
\eneq

Let ${\cal H}\subset A_+^{\bf 1}\setminus\{0\}$ be a finite subset.
Define
\beq\label{Ldet-3}
d=\min\{\Delta_0(\hat{h})/4: h\in {\cal H}\}>0{.}
\eneq

Let $\dt>0$ and let ${\cal G}\subset A$ be a finite subset {{as}}
provided by \ref{measureexistence} for $\ep=d$ and ${\cal F}={\cal H}.$

Suppose that $\phi: A\to B$ is a unital ${\cal G}$-$\dt$-multiplicative {{completely}} positive linear map. 
Then,
by \ref{measureexistence}, for each $t\in T(B),$ there exists $\tau\in T(A)$ such that
\beq\label{Ldet-4}
|t\circ \phi(h)-\tau(h)|<d\tforal h\in {\cal H}.
\eneq
It follows that
$
t\circ \phi(h)>\tau(h)-d\tforal h\in {\cal H}.
$
Thus,
\beq\label{Ldet-6}
t\circ \phi(h)>\Delta_0(\hat{h})-d>\Delta_0(\hat{h})/2\tforal h\in {\cal H} {{\andeqn \rforal t\in T(B)}}.
\eneq
\end{proof}


\begin{lem}\label{Fullmeasure}
Let $C$ be a unital \CA, and let $\Delta: C_+^{q, {\bf 1}}\setminus \{0\}\to (0,1)$ be an order preserving
map.  There exists  a map  $T\times N: C_+\setminus \{0\}\to \R_+\setminus \{0\}\times \N$ {{with}}
the following {{property}}:
For any finite subset ${\cal H} \subset C_+^{\bf 1}\setminus \{0\}$ and any unital \CA\, $A$
 with {{strict}} comparison of positive elements, if $\phi: C \to A$ is a unital \cp\,
 satisfying
\beq\label{Fullm-1}
\tau\circ \phi(h)\ge \Delta(\hat{h})\tforal h\in {\cal H} \tforal \tau\in\mathrm T(A),
\eneq
then $\phi$ is $(T\times N)$-${\cal H}$-full.
\end{lem}


\begin{proof}
For each $\delta\in(0, 1)$, let $g_\delta:[0, 1] \to [0, +\infty)$ be the continuous function defined by
$$g_{\delta}(t)=\left\{
\begin{array}{ll}
0 & \textrm{if $t\in[0,
\dt/4]$},\\
f_{\dt/2}(t)/t & \textrm{otherwise},
\end{array}
\right.$$
where $f_{\dt/2}$ is as defined in \ref{Dball}. Note that
\begin{equation}\label{factor}
g_\delta(t) t  =  f_{\dt/2}(t)\rforal t\in[0, 1].
\end{equation}
Let $h\in C_+^{\bf 1}\setminus\{0\}$. Then define
$$
T(h)= \|(g_{\Delta(\widehat{h})})^{\frac{1}{2}}\|={{\sqrt{\frac {2}{\Delta(\widehat{h})}}}}\andeqn N(h)=\lceil\frac{2}{\Delta(\widehat{h})}\rceil.
$$
Then the function $T\times N$ {{has the property of}}  the lemma.

Indeed, let $\mathcal H\subset C^1_+\setminus \{0\}$ be a finite subset. Let $A$ be a unital
\CA\,  with
{{strict}} comparison for positive elements, and let $\phi: C\to A$ be a unital positive linear map satisfying
\begin{equation}\label{eq-lb-tr}
\tau\circ\phi(h)\geq\Delta(\hat{h})\rforal h\in \mathcal H \rforal \tau\in \mathrm{T}(A).
\end{equation}
Put $\phi(h)^{-}=(\phi(h)-\frac{\Delta(\hat{h})}{2})_+.$
Then, by  \eqref{eq-lb-tr}, one has that, since $0\le h\le 1,$
$$d_\tau(\phi(h)^{-})=d_\tau((\phi(h)-\frac{\Delta(\hat{h})}{2})_+)\ge \tau((\phi(h)-\frac{\Delta(\hat{h})}{2})_+)\geq \frac{\Delta(\hat{h})}{2}\rforal \tau\in\mathrm{T}(A).$$
 {{This}} shows {{in particular}} that $\phi(h)^{-}\not=0.$
Since $A$ has {{strict}} comparison for positive elements, one has
$K\left<(\phi(h)^{-}\right> > \left< 1_A \right>,$ where $K=\lceil \frac{2}{\Delta(\hat{h})}\rceil$
and where $\la x\ra$ denotes the class of $x$ in  $W(A).$

Therefore
there is a partial isometry $v=(v_{ij})_{K\times K}\in \mathrm{M}_{K}(A)$ such that
$$
vv^*=1_A\quad\textrm{and}\quad v^*v\in \overline{(\phi(h)^{-}\otimes 1_K) M_K(A)(\phi(h)^{-}\otimes 1_K)}.
$$
Note  that
$
c(f_{\dt/2}(\phi(h))\otimes 1_K)={{(f_{\dt/2}(\phi(h))\otimes 1_K)c}}=c\rforal c\in\overline{(\phi(h)^{-}\otimes 1_K)M_K(A)(\phi(h)^{-}\otimes 1_K)},
$
where $\dt=\Delta(\hat{h}),$
and therefore
$$
v(f_{\dt/2}(\phi(h))\otimes 1_K)v^*=vv^*=1_A.
$$
Considering the upper-left corner of $M_K(A),$ one has
$\sum_{i=1}^K v_{1,i}f_{\dt/2}(\phi(h))v_{1,i}^*=1_A,$
and therefore, by \eqref{factor}, one has
$$
\sum_{i=1}^K v_{1,i} (g_{\Delta(\hat{{h}})}(\phi(h)))^{\frac{1}{2}}\phi(h)(g_{\Delta(\hat{{h}})}(\phi(h)))^{\frac{1}{2}}v_{1,i}^*=1_A.
$$
Since $v$ is a partial isometry, one has $\|v_{i,j}\|\leq 1$, $i,j=1, ..., K$, and therefore
$$\|v_{1,i} (g_{\Delta(\hat{{h}})}(\phi({\blue{h}})))^{\frac{1}{2}}\|\leq \|(g_{\Delta(\hat{{h}})}(\phi({{h}})))^{\frac{1}{2}}\|\leq \|(g_{\Delta(\hat{{h}})})^{\frac{1}{2}}\|=T({{h}}).$$ Hence the map $\phi$ is $T\times N$-$\mathcal H$-full, as desired.
\end{proof}



\begin{thm}\label{UniCtoA}
Let $C$ be a unital \CA\, in ${\bar{\cal D}}_s$ {\rm (}see \ref{8-N-3}{\rm{)}} {{with finitely generated $K_i(C)$ ($i=0,1$)}}.
Let ${\cal F}\subset C$ be a finite subset, let $\ep>0$ be a
positive number and let $\Delta: C_+^{q, {\bf 1}}\setminus \{0\}\to
(0,1)$ be an {order preserving } map. There exist  a finite
subset ${\cal H}_1\subset C_+^{\bf 1}\setminus \{0\},$
constants
$\gamma_1>0,$ $1>\gamma_2>0,$ and $\dt>0,$ a finite subset ${\cal
G}\subset C,$ a finite subset ${\cal P}\subset \underline{K}(C)$, 
a finite subset ${\cal H}_2\subset C_{s.a.},$
and a finite subset
${\cal U}\subset {\blue{U(M_{n_0}(C))/CU(M_{n_0}(C))}}$ (for some $n_0\ge 1$)
for which $[{\cal U}]\subset {\cal P}$
satisfying the following condition: For any
unital ${\cal G}$-$\dt$-multiplicative \cp s
$\phi, \psi: C\to A$,
where $A=A_1\otimes U$ for some
$A_1\in {\mathcal B_1}$ and a UHF-algebra $U$ of infinite type,
satisfying
\beq\label{CUni-1}
&& [\phi]|_{\cal P}=[\psi]|_{\cal P},\\
\label{CUni-2}
&&\tau(\phi(a))\ge \Delta(\hat{a}),\,\,\,
\tau(\psi(a))\ge \Delta(\hat{a})
 \eneq
 for all $\tau\in T(A)$ and
for all $a\in {\cal H}_1,$
\beq\label{CUni-3}
 |\tau\circ\phi(a)-\tau\circ \psi(a)|&<&\gamma_1\tforal a\in {\cal
H}_2,\tand\\\label{CUni1-3+1}
{\rm dist}(\phi^{\ddag}({\blue{z}}),
\psi^{\ddag}({\blue{z}}))&<&\gamma_2\tforal {\blue{z}} \in {\cal U},
\eneq
there exists a
unitary $W\in A$ such that \beq\label{CUni-4}
\|W^*\phi(f)W-\psi(f)\|<\ep\tforal f\in {\cal F}. \eneq
\end{thm}

\begin{proof}
Let $T'\times N: C_+\setminus\{0\}\to\mathbb R_+\{0\}\times\mathbb N$ be the map of Lemma \ref{Fullmeasure} with respect to $C$ and $\Delta/4$. Let $T=2T'.$

Define
${\bf L}=1.$
Let $\delta_0>0$ (in place of $\delta$), $\mathcal G_0\subset C$ (in place of $\mathcal G$), $\mathcal H_0\subset C_+\setminus\{0\}$ (in place of $\mathcal H$), $\mathcal V_0\subset U(M_{n_0}(C))$
(in place of $\mathcal U$), and
$\mathcal P_0\subset \underline{K}(C)$ (in place of $\mathcal P$) {be finite subsets} and {$n_1$ (in place of $n$) be an integer as {{provided}} by} Theorem \ref{Suni} {{(with the modification as in   (2) of \ref{Rsuni},  and with
the second inequality of (\ref{stableun2-1}) replaced by (\ref{Rsuni-1}) )}} with respect to $C$ (in place of $A$), $T\times N$, ${\bf L}$, $\mathcal F,$ and $\epsilon/2$.
{\blue{Put ${\cal U}_0=\{\bar{v}\in U(M_{n_0}(C))/CU(M_{n_0}(C)): v\in {\cal V}_0\}.$}}

Let $\mathcal H_{1, 1}\subset C_+^1\setminus\{0\}$ (in place of $\mathcal H_1$), $\mathcal H_{1, 2}\subset A$ (in place of $\mathcal H_2$), $1>\gamma_{1, 1}>0$ (in place of $\gamma_1$), $1>\gamma_{1, 2}>0$ (in place of $\gamma_2$), $\delta_{1}>0$ (in place of $\delta$), $\mathcal G_1\subset C$ (in place of $\mathcal G$), $\mathcal P_1\subset\underline{K}(C)$ (in place of $\mathcal P$), $\mathcal U_1\subset J_c(K_1(C))$ (in place of $\mathcal U$),
{and  $n_2$ (in place of $N$)} be the finite subsets and constants provided by Theorem \ref{UniqN1} with respect to $C$ (in place of $A$), $\Delta/4$, $\mathcal F,$ and $\epsilon/4$.

Put $\mathcal G=\mathcal G_0\cup\mathcal G_1$, $\delta=\min\{\delta_0/4, \delta_1/4\}$, $\mathcal P=\mathcal P_0\cup\mathcal P_1$, $\mathcal H_1=\mathcal H_{1, 1}$, $\mathcal H_2=\mathcal H_{1, 2}$, $\mathcal U=\mathcal U_0\cup\mathcal U_1$, $\gamma_1=\gamma_{1, 1}{/2}$,  and $\gamma_2=\gamma_{1, 2}{/2}$. We assert {{that}} these are the desired finite subsets and constants (for $\mathcal F$ and $\epsilon$).  {We may assume that
$\gamma_2<1/4.$}

In fact, let $A=A_1\otimes U$, where $A\in\mathcal B_1$ and $U$ is a UHF-algebra
{of infinite type}. Let $\phi, \psi: C\to A$ be $\mathcal G$-$\dt$-multiplicative maps satisfying {{(\ref{CUni-1})
to (\ref{CUni1-3+1}) for the above-chosen ${\cal G},$ ${\cal H}_1,$ ${\cal P},$ ${\cal H}_2,$ ${\cal U},$
$\gamma_1,$ and $\gamma_2.$}} {{{\blue{Applying}} Lemma  \ref{Fullmeasure}, and {\blue{by}} the choice of $T'\times N$ at the beginning of the proof, we know that both $\phi$ and $ \psi$ are $T'\times N$-${\cal H}_1$-full.}}
%

Since $A=A_1\otimes U$, $A\cong A\otimes U.$ Moreover, the map
$j\circ \imath: A\to A$ is approximately inner, where
$\imath: A\to A\otimes U$ is defined by
$a\mapsto  a\otimes 1_U$  and $j: A\otimes U\to A$ is some isomorphism.
Thus, we may assume that $A=A_1\otimes U\otimes U=A_2\otimes U,$ where $A_2=A_1\otimes U.$
Moreover, without loss of generality, we may assume that the images of $\phi$ and $\psi$
are in $A_2.$
Since $A_2\in\mathcal B_1$,   for {{every}}  finite subset ${\cal G}''\subset A_2,$
$\delta'>0,$  and integer $m\ge 1,$  there are a projection $p\in A_2$ and a $C^*$-subalgebra
$D\in \mathcal C_1$
with $p=1_D$
such that
\begin{enumerate}
\item $\|pg-gp\|<\delta'$ for any $g\in\mathcal G'',$
\item $pgp\in_{\delta'} D$, and
\item $\tau(1-p)<\min\{\delta', \gamma_1/4, 1/8m\}$ for any $\tau\in\mathrm{T}(A)$.
\end{enumerate}
Define $j_0: A_2\to (1-p)A_2(1-p)$ by $j_0(a)=(1-p)a(1-p)$ for all $a\in A_2.$
{\blue{ With sufficiently small $\dt'$ and large ${\cal G}'',$ applying Lemma  \ref{Tapprox}, one obtains  a {{unital \cp\,}}
$j_1: A_2\to D$ with $\|j_1(a)-pap\|<\dt'$ for all $a\in {\cal F}$ such that
$j_i\circ \phi$ and $j_i\circ \psi$
are $\mathcal G$-$2\dt$-multiplicative, $i=0,1,$  and
\begin{enumerate}\setcounter{enumi}{3}
\item\label{concl-708-001} $\|\phi(c)-(j_0\circ \phi(c)\oplus j_1\circ \phi(c))\|<\ep/16$ and $\|\psi(c)-(\psi_0(c)\oplus\psi_1(c))\|<\ep/16$, for any $c\in\mathcal F$,
\item\label{concl-708-002} $j_0\circ \phi, j_0\circ \psi$ and $j_1\circ \phi, j_1\circ \psi_1$ are $2T'\times N$-$\mathcal H_1$-full,
\item\label{concl-708-003} $\tau\circ(j_1\circ \phi(c)>\Delta(\hat{c})/2$ and $\tau\circ j_1\circ \psi(c)>\Delta(\hat{c})/2$ for any $c\in\mathcal H_1$ and for any $\tau\in T(D),$
\item\label{concl-708-004} $|\tau\circ\phi_1(c)-\tau\circ\psi_1(c)|<2\gamma_1$ for any $\tau\in T(D)$ and any $c\in \mathcal H_2,$
\item\label{concl-708-005} $\mathrm{dist}((j_i\circ \phi)^{\ddagger}(u), (j_i\circ \psi)^{\ddagger}(u))<2\gamma_2$ for any $u\in\mathcal U,$
{$i=0,1,$} and
\item\label{concl-181028} $[j_0\circ \phi]|_{\cal P}=[j_0\circ \psi]|_{\cal P}\andeqn
[j_1\circ \phi]|_{\cal P}=[j_1\circ \psi]|_{\cal P}.$
\end{enumerate}
}}

{Choose an integer $m\ge 2(n_1+1)n_2$ and mutually orthogonal and
mutually equivalent projections $e_1,e_2...,e_m\in U$ with
$\sum_{i=1}^me_i=1_U.$} {Define $\phi_i', \psi_i': C\to A\otimes U$
by $\phi_i'(c)=\phi(c)\otimes e_i$ and $\psi_i'(c)=\psi(c)\otimes
e_i$ for all $c\in C,$ $i=1,2,...,m.$ Note that \beq\label{n107-1}
[\phi_1']|_{{\cal P}}=[\phi_i']|_{\cal P}=[\psi'_1]|_{\cal
P}=[\psi_i']|_{\cal P}, \eneq $i=1,2,...,m.$ } {Note also that $\phi_i',
\psi_i': C\to e_iAe_i$ are ${\cal G}$-$\dt$-multiplicative.}

Write $m=kn_2+r,$ where $k\ge n_1+1$ and $r<n_2$ are integers.
Define\\
$
{\tilde \phi}, {\tilde \psi}: C\to (1-p)A_2(1-p)\oplus \bigoplus_{i=kn_2+1}^m A_2\otimes e_i
$
by
\beq\label{nn107-1}
&&{\tilde\phi}(c)=j_0\circ \phi(c)\oplus \sum_{i=kn_2+1}^m  j_1\circ \phi(c)\otimes e_i\andeqn\\
&&{\tilde \psi}(c)=j_0\circ \psi(c)\oplus \sum_{i=kn_2+1}^m  j_1\circ
\psi(c)\otimes e_i \eneq for all $c\in C.$  With sufficiently large
${\cal G}''$ and small $\dt',$ we may assume that ${\tilde \phi}$
and ${\tilde \psi}$ are ${\cal G}$-$2\dt$-multiplicative
 and, by
(\ref{n107-1}),
\beq\label{nn107-2}
 [{\tilde \phi}]|_{{\cal
P}}=[{\tilde \psi}]|_{{\cal P}}.
\eneq
Moreover, {\blue{by  (8) above,}}
we {{have}}
{{$${\rm dist}({(j_i\circ \phi)}^{\ddag}({\bar{v}}), {(j_i\circ
\psi)}^{\ddag}({\bar{v}}))<2\gamma_2\le \gamma_{1,2}\le {\bf L}~~\mbox{for}~~ i=1,2, ~~\mbox{and }~~
{\bar{v}}\in {\cal U},$$}}
{{which implies}}
\beq\label{nn107-2+}
{\rm dist}({\tilde \phi}^{\ddag}({\bar v}), {\tilde
\psi}^{\ddag}({\bar v})){{<}} {\bf L} \eneq for all ${\bar
v}\in {\cal U}.$ Define  $\phi_i^1,\psi_i^1: C\to D\otimes e_i$ by
$\phi_i^1(c)=j_1\circ \phi(c)\otimes e_i$ and $\psi_i^1(c)=j_1\circ
\psi(c)\otimes e_i.$
By  {\blue{(7) and (8) above,}}
we {{have}}
\beq\label{nn107-10}
&&\tau\circ \phi_i^1(h)\ge \Delta(\hat{h})/2\andeqn \tau(\psi_i^1(h))\ge \Delta(\hat{h})/2\rforal h\in {\cal H}_1,\\\label{nn107-10+}
&&|\tau \circ \phi_i^1(c)-\tau\circ \psi_i^1(c)|<\gamma_{1,1}\rforal c\in {\cal H}_2
\eneq
and for all $\tau\in T(pAp\otimes e_i),$ and, {{furthermore,}}
\beq\label{nn107-10++}
{\rm dist}((\phi_i^1)^{\ddag}({\bar v}),(\psi_i^1)^{\ddag}({\bar v})) <\gamma_{1,2} \rforal {\bar v}\in {\cal U}.
 \eneq
{\blue{By (5) above, we have that}}  $\phi_i^1$ and $\psi_i^1$ are $T\times N$-${\cal H}_1$-full{{, since $T=2T'$.}}
Moreover, we may also assume that
$\phi_i^1$ and $\psi_i^1$ are ${\cal G}$-$2\dt$-multiplicative,
and, by \eqref{concl-181028},
\beq\label{nn107-10+++}
[\phi_i^1]|_{\cal P}=[\psi_i^1]|_{\cal P}.
\eneq

Define
$\Phi, \Psi: C\to  \bigoplus_{i=1}^{kn_2} D\otimes e_i$
by
\beq\label{nn107-3}
\Phi(c)=\bigoplus_{i=1}^{kn_2} \phi_i^1(c)\andeqn
\Psi(c)=\bigoplus_{i=1}^{kn_2} \psi_i^1(c)
 \eneq
for all $c\in C.$
By (\ref{nn107-10+++}), (\ref{nn107-10}), (\ref{nn107-10+}), (\ref{nn107-10++}), and
by {{Theorem}} \ref{UniqN1}, there exists a unitary $W_1\in (\sum_{i=1}^{kn_2}p\otimes e_i)(A_2\otimes U)(\sum_{i=1}^{kn_2}p\otimes e_i)$ such that
\beq\label{nn107-15}
\|W_1^*\Phi(c)W_1-\Psi(c)\|<\ep/4\rforal c\in {\cal F}.
\eneq
 Note  that
 \beq\label{nn107-16}
  \tau(1-p)+\sum_{kn_2+1}^m\tau(e_i)<(1/m)+(r/m)\le n_2/m
 \eneq
  for all $\tau\in T(A).$  Note also that $k\ge n_1.$
  By (\ref{nn107-2}) and (\ref{nn107-2+}), since $\psi_i^1$ is $T\times N$-${\cal H}_1$-{{full}}, on applying {{Theorem}} \ref{Suni}, 
  one obtains a unitary $W_2\in A$ such that
  \beq\label{nn107-17}
  \|W_2^*({\tilde \phi}(c)\oplus \Psi(c))W_1-({\tilde \psi}(c)\oplus \Psi(c))\|<\ep/2
  \rforal c\in {\cal F}.
  \eneq
  Set
  $$
  W=({\rm diag}(1-p, e_{kn_2+1}, e_{kn_2+2},...,e_m)\oplus W_1)W_2.
  $$
  Then
  we compute that
  \beq\label{nn107-18}
  \|W^*({\tilde \phi}(c)\oplus \Phi(c))W-({\tilde \psi}(c)\oplus \Psi(c))\|<\ep/2+\ep/4
  \eneq
  for all $c\in {\cal F}.$  By \eqref{concl-708-001},
  we have
  \beq\label{nn107-19}
  \|W^*\phi(c)W-\psi(c)\|<\ep\rforal c\in {\cal F}
  \eneq
  as desired.
\end{proof}

\vspace{0.2in}

Theorem \ref{UniCtoA} can be strengthened as follows.

\begin{cor}\label{RemUniCtoA}
{\blue{(1) Theorem \ref{UniCtoA} still holds, if, in \eqref{CUni-2}, only one  of the two
inequalities holds.}}

{{ (2) In Theorem \ref{UniCtoA}, one can choose ${\cal U}\subset U(M_{n_0}(C))/CU(M_{n_0}(C))$ to be}} {\blue{a finite}} subset of  {\blue{a torsion free subgroup of}} $J_c(K_1(C))$ {{(see \ref{Dcu})}}.
{\blue{Furthermore, if $C$ has stable rank $k,$ then ${{{\cal U}\subset}} J_c(K_1(C))$
${{(=J_c(U(M_k(C))/U_0(M_k(C))))}}$ may be {{chosen to be {\blue{a finite}} subset of a free subgroup of}}
{{(the abelian group)}} $U(M_k(C)/{{CU}}(M_k(C)).$
In the case that  $C$ has stable rank one,
then ${\cal U}$ may be assumed to be a subset of a free subgroup of ${{U}}(C)/{{CU}}(C).$ In the case that $C=C'\otimes C(\T)$ for
some $C'$ with stable rank one,  then the stable rank of $C$ is no more than $2.$ Therefore, in this case,
${\cal U}$ may be assumed to be  a finite subset of a free subgroup of
$U(M_2(C))/CU(M_2(C)).$}}

\end{cor}

\begin{proof}

{\blue{ For part (1), let $\Delta$ be given. Choose $\Delta_1=(1/2)\Delta.$
Suppose ${\cal G},$ ${\cal P},$ ${\cal H}_1$ and ${\cal H}_2,$
$\dt,$ $\gamma_1,$ and $\gamma_2$ are as provided by Theorem \ref{UniCtoA} for  the given $\ep,$ ${\cal F},$ and  $\Delta_1$ (instead of $\Delta$).

Set $\sigma=\min\{\Delta(\hat{a}): a\in {\cal H}_1\}$ and
 $\gamma_1'=\min\{\gamma_1, \sigma/2\}.$
Choose  ${\cal H}_3={\cal H}_1\cup {\cal H}_2.$
Now suppose that
\beq
\phi(a)\ge \Delta(\hat{a})\rforal {\cal H}_1\andeqn
|\tau(\phi(b))-\tau(\psi(b))|<\gamma_1'\rforal b\in {\cal H}_3.
\eneq
Then
\beq
\psi(a)\ge \Delta(\hat{a})-\gamma_1'\ge \Delta_1(\hat{a})\rforal a\in {\cal H}_1.
\eneq
This shows that,  replacing ${\cal H}_2$ by ${\cal H}_3$ (or by choosing ${\cal H}_2\supset {\cal H}_1$) and
replacing $\gamma_1$ by $\gamma_1',$
we only need one inequality in \eqref{CUni-2}.  This proves part (1).}}

{\blue{For part (2), let ${\cal U}'\in U(M_{n_0}(C))/CU(M_{n_0}(C))$ be a finite subset and $\gamma_2'>0$ be given.
Suppose that ${\cal U}\subset U(M_{n_0}(C))/CU(M_{n_0}(C))$ is another finite subset such that
${\cal U}'\subset G({\cal U}),$ the subgroup  generated by ${\cal U}.$
Then, it is routine to check}} that, there exist a sufficiently small $\gamma_2>0$ and
a sufficiently large finite subset ${\cal G}\subset C$ and small  $\delta>0$ that, for any ${\cal G}$-$\delta$-multiplicative completely positive contractions $\phi,~ \psi : C\to B$, if $${\rm dist}(\phi^{\ddag}(u),
\psi^{\ddag}(z))<\gamma_2\tforal z\in {\cal U},$$
then $${\rm dist}(\phi^{\ddag}(u),
\psi^{\ddag}(z))<\gamma'_2\tforal z\in {\cal U}'.$$
{\blue{This shows that, we may replace a finite subset ${\cal U}$ by any generating subset of
the subgroup $G({\cal U})$ (with possibly larger ${\cal G}$ and  smaller $\dt$ and $\gamma_2$).

Note that,  since $K_0(C)$ is finitely generated,
$\rho^{n_0}_C(K_0(C))=\rho_C(K_0(C))$  if $n_0$ is chosen  large enough
(see Definition \ref{Dcu}).
By  Definition \ref{Dcu}, this implies that
$U(M_{\infty}(C))/CU(M_{\infty}(C))\cong \Aff(T(A))/\overline{\rho_A(K_0(A))}\cong U(M_{n_0}(C)/CU(M_{n_0}(C)).$
Since $K_i(C)$ is finitely generated,  {{on}} choosing a larger $n_0$ if necessary,
$K_1(C)$ is generated by images of unitaries in $U(M_{n_0}(C)).$
Write
$$
U(M_{n_0}(C))/CU(M_{n_0}(C))=U(M_{\infty}(C))/CU(M_{\infty}(C))=\Aff(T(A))/\overline{\rho_A(K_0(A))}\bigoplus J_c(K_1(C))
$$
as in \ref{Dcu}.  We may write $K_1(C)=G_f\oplus G_t,$ where $G_f$ is a finitely generated free abelian group and $G_t$ is a finite
group.
By the previous paragraph, we may assume that
${\cal U}={\cal U}_0\sqcup {\cal U}_f\sqcup {\cal U}_t,$  where
${\cal U}_0\subset U_0(M_{n_0}(C))/CU(M_{n_0}(C)),$ ${\cal U}_f\subset J_c(G_f),$ and ${\cal U}_t\subset J_c(G_t).$

Suppose, in Theorem \ref{UniCtoA}, that ${\cal U}$ has been chosen as above. We then enlarge ${\cal P}$ so that ${\cal P}$
contains $G_t.$
Then, if $z\in {\cal U}_t,$ by the assumption $[\phi]|_{\cal P}=[\psi]|_{\cal P},$
$\kappa_1^C(\phi^{\ddag}(z)-\psi^{\ddag}(z))=0,$ where $\kappa_1^C : U(M_{\infty}(C))/CU(M_{\infty}(C))\to K_1(C)$ is the  quotient map defined in \ref{Dcu}.
 It follows that
$\phi^{\ddag}(z)-\psi^{\ddag}(z)\in U_0(M_{n_0}(A))/CU(M_{n_0}(A))=\Aff(T(A))/\overline{\rho_A(K_0(A))}$. Also, since $J_c(G_t)$ is a torsion group, there is an $m>0$ such that $m\big(\phi^{\ddag}(z)-\psi^{\ddag}(z)\big)=0.$
Since $A=A_1\otimes U$ and $A_1\in {\cal B}_1,$ by \ref{Unotrosion}, ${\rm Aff}(T(A))/\overline{\rho_A(K_0(A))}$ is torsion free. Therefore $\phi^{\ddag}(z)-\psi^{\ddag}(z)=0$. {{In other words, $[\phi]|_{\cal P}=[\psi]|_{\cal P}$ implies that
$\phi^{\ddag}|_{{\cal U}_t}=\psi^{\ddag}|_{{\cal U}_t}$ when $G_t\subset {\cal P},$
which means   we may assume  that ${\cal U}={\cal U}_0\sqcup {\cal U}_f$}}
(by choosing ${\cal P}\supset G_t$).}}

{\blue{Next,  fix the finite subset ${\cal  U}_0.$ If $\gamma_2>0$ is given, by Lemma \ref{HvsU-lem-2018},
if ${\cal H}_2$ and ${\cal G}$ are sufficiently large and $\dt$ and $\gamma_1$ are sufficiently small, then
\beq\label{C181029-1}
{\rm dist}(\phi^{\ddag}(z), \psi^{\ddag}(z))<\gamma_2\rforal z\in {\cal U}_0.
\eneq
In other words, \eqref{C181029-1} follows from \eqref{CUni-3}, provided that ${\cal H}_2$ and ${\cal G}$
are sufficiently large and $\dt$ and $\gamma_1$ are sufficiently small.
Thus the first part of (2) of the corollary follows.}}

{\blue{The case that $C$ has stable rank $k$ also follows since we have $K_1(C)=U(M_k(C))/U_0(M_k(C)).$
In particular, if $k=1,$ then $K_1(C)=U(C)/U_0(C).$
If $C=C'\otimes C(\T)$ for some $C'$ which has stable rank one, then $C$ has stable rank  at most two
(see Theorem 7.1 of \cite{Reff}).
So the last statement also follows.}}

\end{proof}

\begin{cor}\label{Unitaryuni}
Let
$\ep>0$ be a positive number and let $\Delta: C(\T)_+^{q, {\bf 1}}\setminus \{0\}\to (0,1)$ be {{an order
preserving}} map.  There exist a finite subset
${\cal H}_1\subset C(\T)_+^{\bf 1}\setminus \{0\},$
$\gamma_1>0,$ $1>\gamma_2>0,$ and {\blue{a}}  finite subset ${\cal H}_2\subset C(\T)_{s.a.}$
satisfying the following  condition:
For any two unitaries $u_1$ and $u_2$ in a unital separable simple
\CA\, {{ $A=A'\otimes U$ with $A'\in {\cal B}_1$ and $U$ a UHF-algebra of infinite type}}
such that
\beq\label{UUni-1}
[u_1]=[u_2]\in K_1(A),\,
\tau(f(u_1)),\,\tau(f(u_2))\ge \Delta(\hat{f})
\eneq
for all $\tau\in T(C)$ and for all $f\in {\cal H}_1,$ and
\beq\label{UUni-3}
|\tau(g(u_1))-\tau(g)(u_2)|<\gamma_1\tforal g\in {\cal H}_2\tand
{\rm dist}(\bar{u_1}, \bar{u_2})<\gamma_2,
\eneq
there exists a unitary $W\in C$ such that
\beq\label{UUni-4}
\|W^*u_1W-u_2\|<\ep.
\eneq
\end{cor}

\begin{lem}\label{tensorprod}
Let $A$ be a \CA\, and $X$ be a compact metric space.
Suppose that $y\in (A\otimes C(X))_+\setminus \{0\}.$
Then {{there exist}} $a(y)\in A_+\setminus \{0\},$ $f(y)\in C(X)_+\setminus \{0\},$ and
$r_y\in A\otimes C(X)$ such that
$\|a(y)\|\le \|y\|,$
$\|f(y)\|\le 1,$
and
$r_y^*yr_y=a(y)\otimes f(y).$
\end{lem}

\begin{proof}
Identify $A\otimes C(X)$ with $C(X,A).$
Let $x_0\in X$ be such that $\|y(x_0)\|=\|y\|.$
There is $\dt>0$  such that $\|y(x)-y(x_0)\|<\|y\|/16$ for   $x\in B(x_0, 2\dt).$
Let $Y=\overline{B(x_0, \dt)}\subset X.$  Let $z(x)=(y(x_0)-\|y\|/4)_+$
for all $x\in Y.$
Note $z\not=0.$
By 2.2 and  (iv) of 2.4  of \cite{RorUHF2}, there exist $r\in C(Y, A)$
such that
$r^*(y|_Y)r=z.$ Choose $g\in C(X)_+\setminus \{0\}$ such that
$0\le g\le 1,$ $g(x)=0$ if ${\rm dist}(x, x_0)\ge \dt,$ and $g(x)=1$ if ${\rm dist}(x, x_0)\le \dt/2.$
Since $g$ is a zero outside $Y,$ one may view $rg^{1/2}, zg\in C(X, A).$ Put $r_y=rg^{1/2},$ $a(y)=(y(x_0)-\|y\|/4)_+$ and
$f(y)=g.$ Then
$$
r_y^*yr_y=zg=a(y)\otimes f.
$$
\end{proof}

\begin{thm}\label{MUN1}{{Part (a).}} Let $A\in {\cal B}_1$
be a unital simple \CA\, which satisfies the UCT.
For any $\ep>0,$ {\Blue{and}} any finite subset ${\cal F}\subset A,$
there exist
$\dt>0,$ a finite subset ${\cal G}\subset A,$  $\sigma_1, \sigma_2>0,$ a finite subset ${\cal P}\subset \underline{K}(A),$ a finite subset ${\cal U}\subset {{U(A)/CU(A)}}$
{(see \ref{ReMUN1}),}
and a finite subset ${{\cal H}}\in A_{s.a}$ satisfying the following condition:

Let $B'\in {\cal B}_1,$ let $B=B'\otimes U$ for some  UHF-algebra $U$ of infinite type{\Blue{,}} and let $\phi, \psi: A \to B$ be two unital ${\cal G}$-$\dt$-multiplicative \cp s
such that
\beq\label{MUN1-1}
[\phi]|_{\cal P}&=&[\psi]|_{\cal P},\\
\label{MUN1-1++}
|\tau\circ \phi(a)-\tau\circ \psi(a)|&<&\sigma_1\tforal a\in {{\cal H}},\text{and} \\\label{MUN1-1+2}
{\rm dist}(\phi^{\ddag}({\bar{u}}), \psi^{\ddag}({\bar u}))&<&\sigma_2\tforal {\bar u}\in {\cal U}.
\eneq
Then there exists a unitary $u\in U(B)$ such that
\beq\label{ }
\|{\rm Ad}\, u\circ \phi(f)-\psi(f)\|<\ep\tforal f\in {\cal F}.
\eneq

{{Part (b).}} Let $A_1\in {\cal B}_1$
be a unital simple \CA\, which satisfies the UCT,
{\Blue{and let}} $A=A_1\otimes C({{\T}}).$
For any $\ep>0,$ any finite subset ${\cal F}\subset A,$  and any order preserving map $\Delta:
C({{\T}})_+^{{\bf 1}}\setminus \{0\}\to (0,1),$
there exist
$\dt>0,$ a finite subset ${\cal G}\subset A,$  $\sigma_1, \sigma_2>0,$ a finite subset ${\cal P}\subset \underline{K}(A),$ a finite subset ${\cal H}_1\subset C({{\T}})_+^{\bf 1}\setminus\{0\},$ a finite subset ${\cal U}\subset {U(M_2(A))/CU(M_2(A))}$
{(see \ref{ReMUN1}),}
and a finite subset ${\cal H}_2\in A_{s.a}$ satisfying the following condition:

Let $B'\in {\cal B}_1,$ let $B=B'\otimes U$ for some  UHF-algebra $U$ of infinite type and let $\phi, \psi: A \to B$ be two unital ${\cal G}$-$\dt$-multiplicative \cp s
such that
\beq\label{MUN1-1}
[\phi]|_{\cal P}&=&[\psi]|_{\cal P},\\\label{MUN1-1+}
\tau\circ \phi(1\otimes h) &\ge & \Delta({{\hat{h}}})\tforal  h\in {\cal H}_1\tand \tau\in T(B),\\\label{MUN1-1++}
|\tau\circ \phi(a)-\tau\circ \psi(a)|&<&\sigma_1\tforal a\in {\cal H}_2,\text{and} \\\label{MUN1-1+2}
{\rm dist}(\phi^{\ddag}({\bar{u}}), \psi^{\ddag}({\bar u}))&<&\sigma_2\tforal {\bar u}\in {\cal U}.
\eneq
Then there exists a unitary $u\in U(B)$ such that
\beq\label{ }
\|{\rm Ad}\, u\circ \phi(f)-\psi(f)\|<\ep\tforal f\in {\cal F}.
\eneq

\end{thm}

\begin{proof} {{We only prove part (b){\Blue{;}} part(a) is simpler.}}
Let $1>\ep>0$ and let ${\cal F}\subset A$ be a finite subset.
Without loss of generality, we  may assume that {{$1_A\in {\cal F},$}}
$$
{\cal F}=\{ a\otimes f: a\in {\cal F}_1\andeqn f\in {\cal F}_2\},
$$
where ${\cal F}_1\subset A$ is a finite subset and ${\cal F}_2\subset C({{\T}})$
is also a finite subset. We further assume that
${\cal F}_1$ and ${\cal F}_2$ are  in the unit ball  of $A$ and $C({{\T}}),$ respectively.

{{Let}} $\Delta: C({{\T}})_+^{\bf 1}\setminus \{0\}\to (0,1)$ be an order preserving map.
Let $T'\times N': C({{\T}})_+\setminus \{0\}\to \R_+\setminus \{0\}\times \N$ be the map as  given by \ref{Fullmeasure}  with respect to  $3\Delta/16.$
Since $A_1$ is a unital separable simple \CA, the identity map on $A_1$  is
$T''\times N''$-full for some $T''\times N'': (A_1)_+\setminus \{0\}\to \R_+\setminus \{0\}\times \N.$

Define  a map $T\times N: A_+\setminus \{0\}\to \R_+\setminus \{0\}\times \N$ as follows:
For any $y\in A_+\setminus \{0\},$ {{by \ref{tensorprod}}},
choose $a(y)\in (A_1)_+\setminus \{0\},$  $f(y)\in C({{\T}})_+\setminus \{0\},$
and $r_y\in A\otimes C({{\T}})$
 such that \beq\label{May-23-2019}{{r_y^*yr_y=a(y)\otimes f(y){\Blue{,}}
 ~~\|a(y)\|\le \|y\|{\Blue{,}}~~
 \mbox{and}~~ \|f(y)\|\le 1.}}\eneq

There are $x_{a(y),1},x_{a(y),2},...,x_{a(y), N''(a(y))}\in A_1$ with
$\max\{\|x_{a(y),i}\|: 1\le i\le N''(a(y))\}=T''(a(y))$ such that
\beq\label{12-fulleq}
\sum_{i=1}^{N''({a(y)})}x_{a(y),i}^*a(y)x_{a(y),i}=1_{A_1}.
\eneq
Then define
\beq\label{MUN1-6-3}
\hspace{-0.3in}(T\times N)(y)=(1+\max\{T''(a(y)), T'(f(y))\} \cdot 2\max\{1, \|y\|, {{\|r_y\|}}\}, N''(a(y)) \cdot N'(f(y))).
\eneq

 {{Let ${\bf L}=1.$}} Let $\ep/16>\dt_1>0$ (in place of $\dt$), ${\cal G}_1\subset A$ (in place of ${\cal G})$,
 ${\cal H}_0\subset A_+\setminus \{0\}$ (in place of $\mathcal H$),
${\cal U}_1\subset {U(M_2(A))/CU(M_2(A))}$ {{(in part (a), one can choose  ${\cal U}_1\subset U(A)/CU(A)$)}} (in place of ${\cal U}$---see  (2) of \ref{Rsuni}) {{(Note that $A\in {\cal B}_1$ has stable rank $1$ in part (a), $A=A_1\otimes C(\T)$ has stable rank at most $2$ in part (b)}}), ${\cal P}_1\subset \underline{K}(A)$ (in place of ${\cal P}$), and $n_0\ge 1$ (in place of $n$) be the finite subsets and constants  provided by \ref{Suni} for $A,$ ${\bf L}{{=1}},$ $\ep/16$ (in place of $\ep$), ${\cal F},$ and $T\times N.$
{{Note, here we refer to  the inequality \eqref{Rsuni-1} instead of the inequality \eqref{Suni-1}.
Moreover, this also implies that $[L']|_{{\cal P}_1}$ is well defined, for any
${\cal G}_1$-$\dt_1$-multiplicative \morp\, $L':A\to B'$ (for any \CA\, $B'$).}}

Without loss of generality,  we may {{assume that}} 
\beq\label{MUN1-6}
{\cal G}_1=\{a\otimes g: a\in {\cal G}_1'\andeqn g\in {\cal G}_1''\},
\eneq
where ${\cal G}_1'\subset A_1$  and
${\cal G}_1''\subset C({{\T}})$ are finite subsets. We may further assume
that ${\cal F}_1\subset {\cal G}_1'$ and ${{{\cal F}_2}}\subset {\cal G}_1'',$
and {\Blue{that}} ${\cal H}_0,{\cal G}_1'$, and ${\cal G}_1''$ are  all in the unit balls, respectively.  In particular,
${\cal F}\subset {\cal G}_1.$
Let
\beq\label{MUN1-6+}
\overline{{\cal H}_0}=\{a(y)\otimes f({{y}}): {{y}}\in {\cal H}_0\}=
\{a\otimes f: a\in {\cal H}_0'\andeqn f\in {\cal H}_0''\},
\eneq
where {{$a(y)$, $f(y)$ {\Blue{are}} defined in (\ref{May-23-2019}),}} ${\cal H}_0'\subset (A_1)^{\bf 1}_+\setminus \{0\}$ and ${\cal H}_0''\subset C({{\T}})^{\bf 1}_+\setminus \{0\}$ {{are}} finite
subsets.
{{We may assume that $1_{A_1}\in {\cal H}_0'$ and
$1_{C({{\T}})}\in {\cal H}_0''.$}} For convenience, {\Blue{let us}} further assume that,  for the above  integer $n_0\ge 1,$
\beq\label{n107-n1}
1/n_0<\inf\{\Delta(h): h\in {\cal H}_0''\}/16.
\eneq
{{Set
$
M_0=\sup\{((T'(h)+1)\cdot N'(h): h\in {\cal H}_0''\}
$
and choose $n\ge n_0$ such that $K_0(U)\subset \Q$ is divisible by $n,$ i.e.,  $r/n\in K_0(U)$ for all $r\in K_0(U).$}}


{Let ${\cal U}_1=\{\bar{v_1},\bar{v_2},...,\bar{v_K}\},$ where
$v_1,v_2,...,v_K\in U(M_2(A)).$ Put
${\cal U}_0=\{v_1,v_2,...,v_K\}.$}
Choose a finite subset
${{{\cal G}_u}}\subset A$ 
such that
\beq\label{MUN1-6+++}
{v_j\in \{(a_{i,j})_{1\le i, j\le 2}: a_{i,j}\in {{{\cal G}_u}}\tforal v_j\in {\cal U}_0.}
\eneq
Choose  {{a small enough}} $\dt_1'>0$ and a
large enough finite subset ${\cal G}_v\subset A_1$
{{that the following condition {\Blue{holds}}:}}
If $p\in A_1$ is a projection such that
$$
\|px-xp\|<\dt_1'\tforal x\in {\cal G}_v,
$$
then {there are unitaries $w_j\in ({\rm diag}(p,p)\otimes
1_{C({{\T}})})M_2(A)({\rm diag}(p,p)\otimes 1_{C({{\T}})})$ such that (for $1\le j\le K$)
\beq\label{n107-n2} \|({\rm diag}(p,p)\otimes 1_{C({{\T}})})v_j({\rm
diag}(p,p)\otimes 1_{C({{\T}})})-w_j\|<\dt_1/16n\tforal v_j\in {\cal U}_0.
\eneq }
Let
\beq\nonumber
{\cal G}_2'={\cal F}_1\cup {\cal G}_1'\cup {\cal H}_0'\cup {{{\cal G}_v}}\cup \{a(y),x_{a(y),j}, x_{a(y),j}^*, r_y, r_y^*: y\in {\cal H}_0\}{{\subset A_1}},
\andeqn\\
M_1=(\max\{\|x_{a(y), j}\|: y\in {\cal H}_0'\}+1) \cdot \max\{N(y): y\in {\cal H}_0\}.
\eneq
{{Put}}
\beq\label{May29-2019}
\dt_1''=\min\{\dt_1', \dt_1\}/(2^{16}(n+1)M_1^2M_0^2).
\eneq
Since $A_1\in  {\cal B}_1,$
there exist mutually orthogonal
projections $p_0', p_1'\in A_1,$
a \SCA\, {{(of $A_1$)}} $C\in {\cal C}$
with $1_{C}={{p_1'}}$, {\Blue{and}} unital {${\cal G}_2'$}-$\dt_1''/16$-multiplicative \cp s
$\imath_{00}': A_1\to p_0'A_1p_0'$ and $\imath_{01}': A_1\to
C$ such that
\beq\label{MUN1-7}
{{\|{{p_1'}}x-x{{p_1'}}\|<\dt_1''/2,}}
\,\,\,{{\rm diag}(\overbrace{p_0',
p_0',...,p_0'}^{n+1})\lesssim p_1',\andeqn \|x-\imath_{00}'(x)\oplus
\imath_{01}'(x)\|<\dt_1''}
\eneq
for all $x\in
{{\tilde{{\cal G}_2'}}},$ where $p_0'+p_1'=1_{A_1},$ $\imath_{00}'(a)=p_0'ap_0'$ for
all $a\in A_1,$ and $\imath_{01}'$ factors through the map $a\mapsto p_1'ap_1',$
and ${{\tilde{{\cal G}_2'}}}=\{xy: x, y\in {\cal G}_2'\}\cup {\cal G}_2'$ {{(see the lines around
\eqref{dj1-18} and \eqref{dj1-18+}).}}
 Define $p_0=p_0'\otimes 1_{C({{\T}})},$ $p_1=p_1'\otimes
1_{C({{\T}})}.$
\Wlog, we may assume that $p_0'\not=0.$  Since $A_1$  is simple, there is an integer $N_0>1$  such that
\beq\label{MUN1-n7}
{{(N_0-1)}}[p_0']\ge [p_1']\,\,\,{\rm in}\,\,\, W(A_1).
\eneq
This also implies that
\beq\label{MUN1-n8}
{{(N_0-1)}}[p_0]\ge [p_1].
\eneq

 Define
$\imath_{00}: A\to p_0Ap_0$ by
$\imath_{00}(a\otimes f)=\imath_{00}'(a)\otimes f$ and $\imath_{01}: A\to C\otimes C({{\T}})$ by
$\imath_{01}(a\otimes f)=\imath_{01}'(a)\otimes f$ for all $a\in A_1$ and
$f\in C({{\T}}).$
Define $L_0: A\to A$ by
$$
L_0(a)=\imath_{00}(a)\oplus \imath_{01}(a) \rforal a\in A.
$$
{For each $v_j\in {\cal U}_0,$ there exists a unitary $w_j\in M_2(p_0Ap_0)$ such that
\beq\label{n107-nn1}
\|{\rm diag}(p_0,p_0)v_j{\rm diag}(p_0,p_0)-w_j\|<\dt_1/16n, \,\,j=1,2,...,K.
\eneq
}
{{Since $C\subset A_1,$}} $C_{{\T}}:=C\otimes C({{\T}}){\subset  A_1\otimes C({{\T}})=A}.$

Let $\imath_0: C\to A_1$ be the natural embedding of $C$ as a unital \SCA\, of ${{p_1'A_1p_1'}}.$ Let $\imath_0^{\sharp}: C^q\to A_1^q$ be
defined by $\imath_0^{\sharp}(\hat{c})=\hat{c}$ for $c\in C.$
Let $\Delta_0: A_1^{q,{\bf 1}}{{\setminus \{0\}}}\to (0,1)$ be the map given by {{Lemma}} \ref{Ldet}
and
define
$$
\Delta_1(\hat{h})
=\sup\{\Delta_0(\imath_0^{\sharp}(\hat{h_1}))\Delta({{\hat h}}_2)/8: h\ge h_1\otimes h_2,\,\,\,
h_1\in C\setminus \{0\},\andeqn h_2\in C({{\T}})_+\setminus \{0\}\}
$$
for all $h\in (C_X)^{\bf 1}_+\setminus \{0\}.$

Let ${\cal G}_3'=\imath_{0,1}'({{\tilde{{\cal G}_2'}}}).$ {{(Note that ${\cal G}_1'\subset {\cal G}_2'\subset \tilde{{\cal G}_2'}$.)}}
Let ${\cal G}_3=\{a\otimes f: a\in {\cal G}_3'\andeqn f\in {\cal G}_1''\}.$
Let ${\cal H}_3\subset (C_{{\T}})_+\setminus \{0\}$  (in place of ${\cal H}_1$), $\gamma_1'>0$ (in place of $\gamma_1$), $\gamma_2'>0$ 
 (in place of $\gamma_2$), $\dt_2>0$
 (in place of $\dt$), ${\cal G}_4\subset C_{{\T}}$ (in place of ${\cal G}$), ${\cal P}_2\subset \underline{K}(C_{{\T}})$ (in place of ${\cal P}$), ${\cal H}_4'\subset (C_{{\T}})_{s.a.}$ (in place of ${\cal H}_2$), and
$\overline{{\cal U}_2}\subset {U(M_2(C_{{\T}}))/CU(M_2(C_{{\T}}))}$  {(in
 place of ${\cal U}$---{\Blue{see}} Corollary \ref{RemUniCtoA})}
 be the finite subsets and constants
 provided by {\Blue{Theorem}} {\ref{UniCtoA}} for $\dt_1/16$ (in place of $\ep$),
${\cal G}_3$ (in place of ${\cal F}$),
$\Delta_1/2$ (in place of $\Delta$), and $C_{{\T}}$ (in place of $A$).

Let ${\cal U}_2\subset U(M_2(C_{{\T}}))$ be a finite subset {which  has {{a}} one-to-one correspondence to its image
in $U(M_2(C_{{\T}}))/CU(M_2(C_{{\T}}))$ which is exactly  $\overline{{\cal U}_2}.$}  We also assume
that $\{[u]: u\in {\cal U}_2\}\subset {\cal P}_2.$

Without loss of generality, we may assume that
\beq\label{MUN-8}
{\cal H}_3=\{h_1\otimes h_2: h_1\in {\cal H}_3'\andeqn h_2\in {\cal H}_3''\},
\eneq
where ${\cal H}_3'\subset C^{\bf 1}_+\setminus \{0\}$ and ${\cal H}_3''\subset C({{\T}})^{\bf 1}_+\setminus \{0\}$ are
finite subsets,  $1_C\in {\cal H}_3'$ and $1_{C({{\T}})}\in {\cal H}_3'',$ and
\beq\label{MUN-8+1}
{\cal G}_4=\{a\otimes f: a\in {\cal G}_4' \andeqn f\in {\cal G}_4''\}{{\cup\{p_0, p_1\}}},
\eneq
where ${\cal G}_4'\subset C$ and ${\cal G}_4''\subset C({{\T}})$
are finite subsets.
{{Set}}
\beq\label{125-sigma0}
{{\sigma_0=}}{{\min\{}}{{\inf\{\Delta_1(\hat{h}): h\in {\cal H}_3\}}}, {{\inf\{\Delta(\hat{h}): h\in {\cal H}''_0\}\}}}{{  >0.}}
\eneq
{{Note that $\Delta_0$ is the map given by Lemma \ref{Ldet} for the simple $C^*$ algebra $A_1$. For ${\cal H}_3'$ (in place of ${\cal H}$), there are  }} $\dt_3>0$ (in place of $\dt$) and  ${\cal G}_5\subset A_1$ (in place of ${\cal G}$) {{as}} the constant and finite subset
provided by {{Lemma \ref{Ldet}.}}
Set
\beq\label{UNN1-9}
\dt_0={\min\{1/16, \ep/16, \dt_1,\dt_1'', \dt_2, \dt_3, {{\sigma_0, \gamma_1', \gamma_2'}}\}\over{128(N_0+2)(n+1)}}
\eneq
and
set
$$
{\cal G}_6=\{x, xy: x,y\in {{\tilde{{\cal G}_2'}}}\cup L_0({{\tilde{{\cal G}_2'}}})\cup \imath_{01}'({{\tilde{{\cal G}_2'}}})\}\cup
\imath_{01}'({{\tilde{{\cal G}_2'}}})\}\cup {\cal G}_4'\cup {\cal G}_5\andeqn
$$
\beq\label{dg0}
\hspace{-0.3in} {\cal G}_0= \{a\otimes f: a\in {\cal G}_6  \andeqn f\in {\cal G}_1''\cup {\cal H}_0''\cup {\cal G}_4''\}\cup\{p_j: 0\le j\le 1\}\cup\{v_j,w_j: 1\le j\le K\} .
\eneq
To simplify notation, without loss of generality, {\Blue{let us}} assume that ${\cal G}_0\subset A^{\bf 1}.$
{{Let}}
\beq\label{May30-2019}
{\cal P}={\cal P}_1\cup \{p_j: 0\le j\le 1\}\cup [\imath]({\cal P}_2)\cup [\imath_{00}]({\cal P}_1),\eneq
 where $\imath: C_{{\T}}\to A$ is the embedding.
Let ${\cal H}_1={\cal H}_0''\cup {\cal H}_3''.$

{{By \eqref{MUN1-n8}  and \eqref{MUN1-7}, we may assume that ${\cal G}_0$ is sufficiently large (with even smaller $\dt_0$,
  if necessary),
that any ${\cal G}_0$-$\dt_0$-multiplicative \morp\, $L'$ from $A$ (to any unital \CA\, $B'$) has the properties
\beq\label{1107-n1}
{{(N_0-1)}}([L'(p_0)])\ge [L'(p_1)],\,\, (n{{+1}})[L'(p_0)]\le [L'(p_1)],
\andeqn\\ \label{1107-n1+}
\,\tau(L'(1-p_1))<16/15n\rforal \tau\in T(B').
\eneq
}}
{Let ${\cal U}_2'=\{{\rm diag}(1-p_1,1-p_1)+w: w\in {\cal U}_2\}$ and let
${\cal U}_0''=\{ w_j+{\rm diag}(p_1, p_1): 1\le j\le K\}.$
Let ${\cal U}=\{\bar{v}: v\in {\cal U}_1\cup {\cal U}_2' {{\cup\, {\cal U}_0''}}\}$
and let } ${\cal H}_2={\cal H}_4'.$
Let $\sigma_1=\min\{{1\over{4n}},{\gamma_1'\over{16n{{(N_0+2)}}}}\}$ and $\sigma_2=\min\{{1\over{16n{{(N_0+2)}}}},{\gamma_2'\over{16n{{(N_0+2)}}}}\}.$

{{We then choose a finite ${\cal G}\supset {\cal G}_0$ and a positive number $0<\dt<\dt_0/64$ with
the following property: If $L':  A\to B'$ (for any unital \CA\, $B'$) is any unital ${\cal G}$-$\dt$-multiplicative \cp\,
then there exist a projection $q'\in B'$ such that
$\|L'(p_1)-q'\|<\dt_0/64$, {{an}} element $b_0'\in (1-q')B'(1-q')$ {{such that}}
 \beq\|b_0'-(1-q')\|<\dt_0/16~~ {{\mbox{and}~~~ b_0'L'(1-p_1)b_0'=1-q',}}\eneq
 and {{an element}} $b_1'\in q'Bq'$ such that
 \beq \|b_1'-q'\|<\dt_0/16 ~~{{\mbox{and}~~~ b_1'L'(p_1)b_1'=q'{\Blue{,}}}}\eneq
 {\Blue{with the elements $b_0'$ and $b_1'$, cutting $L'$ down approximately as follows:}}
$\|b_0'L'(x)b_0'-L'(x)\|<\dt_0/4$ for all $x\in (1-p_1)A(1-p_1)$ with $\|x\|\le 1,$
 and $\|b_1'L'(x)b_1'-L'(x)\|<\dt_0/4$  for all $x\in p_1Ap_1$ with $\|x\|\le 1.$}}



Now {\Blue{let us}} assume that $B$ is as in the statement {\Blue{of the theorem, and}} $\phi,\,\,\psi: A\to B$ are two unital ${\cal G}$-$\dt$-multiplicative \cp s
satisfying the assumption for $\dt,$ ${\cal G},$
${\cal P},$ ${\cal H}_1,$ ${\cal U},$ ${\cal H}_2,$ $\sigma_1,$ and $\sigma_2$ ({\Blue{as defined above}}).

Note that $B'$ is in  ${\mathcal B_1}$  {and $B=B'\otimes U.$
We may also write $B=B_1\otimes U,$ where $B_1=B'\otimes U,$ since $U$ is strongly self absorbing.
Without loss of generality, by the fact that $U$ is strongly self absorbing, we may assume
that the image of both $\phi$ and $\psi$ are in ${{B_1\otimes 1_U}}.$}

{{As mentioned two paragraphs above, there are two  projections $q_\phi, q_\psi\in B_1,$   unital
\cp s
$\phi_0': (1-p_1)A(1-p_1)\to ((1-q_\phi)B_1(1-q_\phi ))\otimes 1_U,$ $\psi_0': (1-p_1)A(1-p_1)\to (1-q_\psi)B_1(1-q_\psi)\otimes 1_U,$
and unital 
\cp s $\phi_1': pAp\to (q_\phi B_1q_\phi)\otimes 1_U$ and $\psi_1': p_1Ap_1\to (q_\psi Bq_\psi)\otimes 1_U$
such that
\beq
\|\phi_0'-\phi|_{(1-p_1)A(1-p_1)}\|<\dt_0/4,\,\, \|\psi_0'-\psi|_{(1-p_1)A(1-p_1)}\|<\dt_0/4,\\\label{1212-phi1}
\|\phi_1'-\phi|_{p_1Ap_1}\|<\dt_0/4\andeqn \|\psi_1'-\psi|_{p_1Ap_1}\|<\dt_0/4.
\eneq
}}
{{Note since $[\phi(p_1)]=[\psi(p_1)]$ and $B_1$ has stable rank one (see \ref{B1stablerk}), there exists a unitary $u_0\in B_1$
such that $u_0^*(q_\psi)u_0=q_\psi.$  Then $\|u^*\psi(p_1)u-q_\phi\|<\dt.$ Thus,
\wilog, by replacing $\psi$ by ${\rm Ad}\, (u_0\otimes 1_U)\circ \psi,$
we may assume that $q:=q_\phi=q_\psi.$ Note, by \eqref{1107-n1},
\beq\label{1212-tauq}
n\tau(1-q)<\tau(q)~~{{\mbox{and}~~N_0\tau(1-q)>\tau (q)}}~~\rforal \tau\in T(B).
\eneq
}}
{Since $K_0(U)$ is divisible by $n,$ there
are mutually orthogonal and
unitarily equivalent
projections $e_1,...,e_{n}\in U$ such that $\sum_{j=1}^{n}{{e_j}}=1_U.$}
{Define $\phi_{00}', \psi_{00}': A\to (1-q)B_1(1-q)\otimes 1_U$ by
\beq\label{n107-12}
\phi_{00}'(a)=\phi_0'\circ \imath_{00}(a)\otimes 1_U\andeqn
\psi_{00}'(a)=\psi_0'\circ \imath_{00}(a)\otimes 1_U
\eneq
for all $a\in A.$ Define $\Phi_A, \Psi_A: A\to qB_1q\otimes 1_U$ by
\beq\label{n107-13}
\Phi_A(a)=\phi_1'\circ \imath_{01}(a)\otimes 1_U\andeqn
\Psi_A(a)=\psi_1'\circ \imath_{01}(a)\otimes 1_U,
\eneq
and define $\Phi_{i,A}\hspace{-0.05in}: A\to qB_1q\otimes e_i$ by
$\Phi_{i,A}(a)=\Phi_A(a)e_i\rforal a\in A$ ($1\le i\le n$).
Define $\Phi_C, \Psi_C: C_{{\T}}\to (qB_1q){{\otimes 1_U}}$ by
 \beq\label{n107-14}
 \Phi_C=\phi_1'\circ \imath\andeqn \Psi_C=\psi_1'\circ \imath.
 \eneq
}
Note, by the choice of $\dt$ and ${\cal G},$ {{$\phi|_{A_1\otimes 1_{C({{\T}})}}$ and $\psi|_{A_1\otimes 1_{C({{\T}})}}$ are
${\cal G}_5$-$\dt_3$-multiplicative.  By \eqref{1212-phi1}, and  then by  Lemma \ref{Ldet} (see also \eqref{UNN1-9} and \eqref{125-sigma0}),
\beq
\tau(\Phi_C(h\otimes 1_{C({{\T}})}))\ge \tau(\phi(h\otimes 1_{C({{\T}})}))-\dt_0/4
\ge \Delta_0(\hat{h})/2-\dt_0/4\\
\ge  \Delta_0(\hat{h})/2-\sigma_0/16\ge \Delta_0(\hat{h})/3\rforal h\in
{\cal H}_3'.
\eneq
}}
{{Similarly, by the assumption \eqref{MUN1-1+}, we also have that
\beq
\tau(\Phi_C(1_{A_1}\otimes h))\ge \tau(\phi(\imath(1_{A_1}\otimes h))-\dt_0/4\ge
\Delta(\hat{h})-\sigma_0/16\\
\ge \Delta(\hat{h})/2\rforal h\in
{\cal H}_3''.
\eneq}}
It follows that {{(see the definition of $\Delta_1$)}}
\beq\label{n107-15}
\tau(\Phi_C(h))\ge \Delta_1(\hat{h})/2\rforal h\in
{\cal H}_3.
\eneq
By {the} assumption \eqref{MUN1-1++},
and  \eqref{1212-phi1}, for all $\tau\in T(B_1),$
\begin{equation}\label{n107-16}
|\tau(\Phi_C(c))-\tau(\Psi_C(c))|<\sigma_1-\dt_0/4\rforal c\in {\cal H}_4'={\cal H}_2.
\end{equation}
Therefore, {{by \eqref{1212-tauq},}}
for all $t\in T(qB_1q),$
\begin{equation}\label{n107-17}
|t(\Phi_C(c))-t(\Psi_C(c))|<{{{n\over{n+1}}(\sigma_1-\dt_0/4)<(1/2)(\sigma_1-\dt_0/4)}}<\gamma_1'
\rforal c\in{\cal H}_4'.
\end{equation}
Since $[\imath]({\cal P}_2)\subset {\cal P},$ by {the}
assumptions, one  also has
\beq\label{n107-18}
[\Phi_C]|_{{\cal
P}_2}=[\Psi_C]|_{{\cal P}_2}.
 \eneq
 By \eqref{MUN1-1+2},
  \eqref{1212-phi1},  {{and applying  Lemma \ref{ph} (with $K=1$), since ${\cal U}_2'\subset {\cal U},$}}  one has
 \beq\label{n107-19}
{\rm
dist}_{M_2(q{{B_1}}q)}(\Phi_C^{\ddag}({\bar v})\Psi_C^{\ddag}({\bar
v}^*), {\diag}(q,q))< (1+9/8)(\sigma_1+\dt_0/2)< \gamma_2'\rforal v\in {\cal U}_2.
\eneq
By the choices of
$\dt,$ ${\cal G}_4,$ $\gamma_1',$ $\gamma_2',$ ${\cal P}_2,$ ${\cal
H}_3,$ ${\cal H}_4',$ and $\overline{{\cal U}_2},$ and by applying Theorem
\ref{UniCtoA} (see also  Corollary \ref{RemUniCtoA}),  there exists {{a unitary}} $u_1\in qB_1q$ such that
\beq\label{n107-20}
\|u_1^*\Psi_C(c)u_1-\Phi_C(c)\|<\dt_1/16\rforal c\in
{\cal G}_3.
\eneq
Thus, {{since $\imath'_{01}({\cal G}_1'))\subset {\cal G}_3',$}} by \eqref{MUN1-6}, and by (\ref{MUN1-7}),  one obtains
\beq\label{n107-21}
\|u_1^*\Psi_A(a)u_1-\Phi_A(a)\|<\dt_1/16 +\dt'' \rforal a\in {\cal G}_1.
\eneq
{{Note that, by \eqref{UNN1-9} {{and (\ref{125-sigma0})}}, $\dt_0{{<}}{\sigma_0\over{128(N_0+2)}}{{<}} {\inf\{\Delta(\hat{h}): h\in {\cal H}_0''\}\over{128(N_0+2)}}.$}}
One has, {{by \eqref{1212-phi1},
 \eqref{1107-n1+}, the assumption}}
(\ref{MUN1-1+}), {{
(\ref{n107-n1}){\Blue{,}}
and  the choice of $\dt_0,$}}
\beq\label{n107-22}
&&\hspace{-0.7in}\tau(\Phi_A(1_C\otimes h))\ge {{\tau(\phi(p_1'\otimes h))-\dt_0/4=\tau(\phi(p_1(1\otimes h)))-\dt_0/4}}\\
&&\hspace{-0.3in}=\tau(\phi(1_{A_1}\otimes h))-\tau(\phi{{((1-p_1)\otimes h)}})-\dt_0/4\\
&&\hspace{-0.3in}{{\geq}}{{\tau(\phi(1_{A_1}\otimes h))-\tau(\phi(1-p_1))-\dt_0/4
\ge \tau(\phi(1_{A_1}\otimes h))-16/15n-\dt_0/4}}\\
&&\hspace{-0.3in}{{\ge \Delta(\hat{h})-16/15n-\dt_0/4}}\ge 13\Delta(\hat{h})/15
\rforal h\in {\cal H}_0'' \andeqn\hspace{-0.1in}\rforal
\tau\in T(B_1).
\eneq
By {\Blue{Lemma}} \ref{Fullmeasure}, it follows that  $\Phi_A{{|_{1_{A_1}\otimes C(\T)}:~C(T)\to qB_1q\otimes 1_U}}$
is $T'\times N'$-${\cal H}_0''$-full. {{In other words,  for any $h{{=f(y)}}\in {\cal H}_0''$ {{(where $y\in {\cal H}_0$)}}{\Blue{,}} there are $z_1, z_2,...,z_m\in qB_1q$
with $\|z_i\|\le T'(h)$ and $m\le N'(h)$ such that
$\sum_{j=1}^m z_i^*\Phi_C(1_{A_1}\otimes h)z_i=q\otimes 1_U.$}}
{{Note that $\imath_{01}'$
is ${{\tilde{{\cal G}_2'}}}$-$\dt''/16$-multiplicative, and
${\cal G}_0\subset {\cal G}$ (see \eqref{dg0}).
By\eqref{1212-phi1},  the fact that $\phi$ is ${\cal G}$-$\dt$-multiplicative {\Blue{(used several  times)}},
\eqref{1212-phi1} again,  the fact that $\phi$ is ${\cal G}$-$\dt$-multiplicative again,  the fact that
$\imath_{01}'$
is ${{\tilde{{\cal G}_2'}}}$-$\dt''/16$-multiplicative,  \eqref{12-fulleq}, and then using the linearity of the maps involved,
one has,
that for $y\in {\cal H}_0,$}}
{{\beq\nonumber
&&\sum_{i=1}^{N''(a)}\Phi_A(x_{a(y),i}^*\otimes 1_{C({{\T}})})\Phi_A(a(y)\otimes f(y))\Phi_A(x_{a(y),i}\otimes 1_{C({{\T}})})\\\nonumber
&&\hspace{0.1in}\approx_{(M_1^2+M_1){\dt_0\over{4}}}\sum_{i=1}^{N''(a)}\Phi_A(x_{a(y),i}^*\otimes 1_{C({{\T}})})\phi((\imath_{01}'(a(y))\otimes 1_{C({{\T}})})(1_C\otimes f(y)))\phi(\imath_{01}'(x_{a(y),i})\otimes 1_{C({{\T}})})\\\nonumber
&&\hspace{0.2in}\approx_{M_1\dt}\sum_{i=1}^{N''(a(y))}\Phi_A(x_{a(y),i}^*\otimes 1_{C({{\T}})})\phi\left((\imath_{01}'(a(y))\otimes 1_{C({{\T}})})(1_C\otimes f(y))(\imath_{01}'(x_{a(y),i})\otimes 1_{C({{\T}})})\right)\\\nonumber
&&\hspace{0.3in}\approx_{M_1\dt}\sum_{i=1}^{N''(a)}\Phi_A(x_{a(y),i}^*\otimes 1_{C({{\T}})})\phi(\imath_{01}'(a(y))\imath_{01}'(x_{a(y),i})\otimes 1_{C({{\T}})})\phi(1_C\otimes f(y))\\\nonumber
&&\hspace{0.3in}\approx_{(M_1^2+M_1){\dt_0\over{4}}}\sum_{i=1}^{N''(a(y))}\phi(\imath_{01}'(x_{a(y),i}^*)\otimes 1_{C({{\T}})}))\phi(\imath_{01}'(a(y))\imath_{01}'(x_{a(y),i})\otimes 1_{C({{\T}})}))\Phi_A(1_C\otimes f(y))\\\nonumber
&&\hspace{0.3in}\approx_{M_1\dt}\sum_{i=1}^{N''(a(y))}\phi((\imath_{01}'(x_{a(y),i}^*)\imath_{01}'(a(y))\imath_{01}'(x_{a(y),i}))\otimes 1_{C({{\T}})})\Phi_A(1_C\otimes f(y))\\\nonumber
&&\hspace{0.3in}\approx_{2M_1(\dt_1''/16)}\sum_{i=1}^{N''(a(y))}\phi(\imath_{01}'(x_{a(y),i}^*a(y)x_{a(y),i})\otimes 1_{C({{\T}})}))\Phi_A(1_C\otimes f(y))\\\nonumber
&&\hspace{0.3in}\approx_{\dt_0/4} (q\otimes 1_U)\Phi_A(1_C\otimes f(y))
=\Phi_A(1_C\otimes f(y)).
\eneq}}
{{From {\Blue{these}} estimates, it follows that, for any $y\in {\cal H}_0$ {{with $f(y)=h\in {\cal H}''_0$}} (note $M_1\le M_1^2$ and $\dt_0/4<\dt_1''/16$)}}
{{\beq\nonumber
&&\hspace{-0.3in}\|\sum_{j=1}^mz_j^*(\sum_{i=1}^{N''(a(y))}\Phi_A(x_{a(y),i}^*\otimes 1_{C({{\T}})})\Phi_A(r_y^*)\Phi_A(y)\Phi_A(r_y)\Phi_A(x_{a(y),i})\otimes 1_{C({{\T}})}))z_j-q\otimes 1_U\|\\\nonumber
&&<2mT'(h)^2T''(a(y))^2N''(a(y))(\dt_0/2+\dt_1''/16)+
\eneq
\beq\nonumber
&&\hspace{0.4in}\|\sum_{j=1}^mz_j^*(\sum_{i=1}^{N''(a)}\Phi_A(x_{a(y),i}^*\otimes 1_{C({{\T}})})\Phi_A(a(y)\otimes {{f(y)}})\Phi_A(x_{a(y),i}\otimes 1_{C({{\T}})})))z_j-q\otimes 1_U\|\\
\nonumber
&&<2M_0^2M_1^2(\dt_0/2+\dt_1''/16)+mT'(h)(10M_1^2)(\dt_1''/16)+\|\sum_{j=1}^mz_j^*\Phi_A(1_C\otimes h)z_j-q\otimes 1_U\|\\
&&=2M_0^2M_1^2(\dt_0/2+\dt_1''/16)+mT'(h)^2(10M_1^2)(\dt_1''/16)\\
&&\le  2M_0^2M_1^2(\dt_0/2)+12M_0^2M_1^2(\dt_1''/16)<1/64(n+1).
\eneq
}}
{{Note also the image of $\Phi_A$ is in $qB_1q\otimes 1_U.$
There is $c_y\in (qBq)_+\otimes 1_U$ with $\|c_y-q\otimes 1_U\|<1/64(n+1)$}}
such that
\beq
\sum_{j=1}^m c_yz_j^*(\sum_{i=1}^{N''(a(y))}\Phi_A(x_{a(y),i}^*)\Phi_A(r_y^*)\Phi_A(y)\Phi(r_y)\Phi_A(x_{a(y),i}))z_jc_y=q{{\otimes 1_U}}.
\eneq
{{By the definition of $T\times N$ (see \eqref{MUN1-6-3}),
one then computes}}
that $\Phi_A$ is $T\times N$-${\cal H}_0$-full. It follows that
$\Phi_{i,A}$ is $T\times N$-${\cal H}_0$-full (as a map to $qB_1q\otimes e_i$), $i=1,2,...,n.$
{{Since $[\imath_{00}]({\cal P}_1)\subset {\cal P},$}}
\beq\label{n107-24}
[\psi_{00}']|_{{\cal
P}_1}=[\phi_{00}']|_{{\cal P}_1}.
\eneq

By \eqref{MUN1-1+2},
 {{since ${\cal U}_0''\subset {\cal U},$}}
\beq\label{n107-25}
{\rm dist}(\phi^{\ddag}(\overline{w_j+{\rm
diag}(p_1, p_1)}), \psi^{\ddag}(\overline{w_j+{\rm diag}(p_1,
p_1)}))<\sigma_2,\,\, 1\le j\le K.
\eneq
{{Put $w_j^\sim=w_j+{\rm diag}(p_1,
p_1),$ $J=1,2,...,K.$
By \eqref{1212-phi1}, one has
\beq\nonumber
&&\hspace{-0.6in}\|{\diag}(1-q,1-q)(\phi\otimes {\rm id}_{M_2})(w_j^\sim){\diag}(1-q,1-q)-(\phi_{00}'\otimes {\rm id}_{M_2})(w_j)\|<\dt_0\andeqn\\\nonumber
&&\hspace{-0.6in}\|{\diag}(1-q,1-q)(\phi\otimes {\rm id}_{M_2})(w_j^\sim){\diag}(1-q,1-q)-(\phi_{00}'\otimes {\rm id}_{M_2})(w_j)\|<\dt_0.
\eneq}}
It follows from \eqref{n107-25} and the above inequalities, that for $1\le j\le K,$
\beq
{\rm dist}(\overline{(\phi_{00}'\otimes {\rm id}_{M_2})(w_j)+{\diag}(q,q)}, \overline{({{\psi}}_{00}'\otimes {\rm id}_{M_2})(w_j)+{\diag}(q,q)})<2\dt_0+\sigma_2.
\eneq
{{It follows  from (\ref{1212-tauq}) that $N_0[1-q]> [q]$.}} By (\ref{n107-nn1}) {{(with $N_0$ in place of $K$ and $1-q$ in place of $e$)}},
\eqref{MUN1-n8},  and by 
{{Lemma \ref{ph},
 in ${\diag}(1-q, 1-q)M_2(B){\diag}(1-q, 1-q),$}} one has
\beq\label{n107-26-}
{\rm dist}((\phi_{00}')^{\ddag}(\bar{v_j}),(\psi_{00}')^{\ddag}(\bar{v_j}))
<\dt_1/8n+{{{\rm dist}((\phi_{00}')^{\ddag}(\bar{w_j}), (\psi_{00}')^{\ddag}(\bar{w_j}))}}\\\label{n107-26}
{{<\dt_1/8n+(N_0+9/8)(2\dt_0+\sigma_2)}}<1/2<{\bf L},\,\,\
j=1,2,...,K.
\eneq
{{Note that, by \eqref{1212-tauq},  since $B$ has strict comparison, $[(1-q)\otimes 1_U]\le [q\otimes e_1]=[q\otimes e_i]$
($1\le i\le n$). Note also that $\Phi_{i,A}$ is unitarily equivalent to $\Phi_{1,A},$ $1\le i\le n.$}}
It follows from Theorem \ref{Suni} and Remark
\ref{Rsuni}, by \eqref{n107-24} and \eqref{n107-26},   that there exists a unitary $u_2\in B$ such that
\beq\nonumber
 \|u_2^*(\phi_{00}'(a)\oplus \Phi_{1,A}(a)\oplus\cdots
\oplus \Phi_{n,A}(a))u_2-(\psi_{00}'(a)\oplus \Phi_{1,A}(a)\oplus \cdots
\oplus \Phi_{n,A}(a))\| <\ep/16
\eneq
for all $a\in {\cal F}.$ In other
words,
\beq\label{n107-28}
\|u_2^*(\phi'_{00}(a)\oplus
\Phi_A(a))u_2-\psi_{00}'(a)\oplus \Phi_A(a)\|<\ep/16\rforal a\in {\cal
F}.
\eneq
Thus, by (\ref{n107-21}),
 \beq\label{n107-29}
&&\hspace{-0.9in}\|u_2^*(\phi_{00}'(a)\oplus \Phi_A(a))u_2-(\psi_{00}'(a)\oplus
u_1^*\Psi_A(a)u_1)\|<\ep/16+\dt_1/16+\dt''\rforal a\in {\cal F}.
\eneq Let
$u=u_2(q+u_1^*)\in U(B).$  {{Then,
for all $a\in {\cal F},$
\beq
&&\hspace{-0.7in}\|u^*(\phi_{00}'(a)\oplus \Phi_A(a))u-(\psi_{00}'(a)\oplus
\Psi_A(a))\|<\ep/16+\dt_1/8+\dt''.
\eneq
}}
{{It  then follows from   \eqref{1212-phi1},
\eqref{n107-12}, and \eqref{n107-13} that, for all $a\in {\cal F},$
\beq\nonumber
\|u^*(\phi\circ \imath_{00}(a)+\phi\circ \imath_{01}(a))u-(\psi\circ \imath_{00}(a)+\psi\circ \imath_{01}(a))\|<\ep/16+\dt_1/8+\dt''+4(\dt_0/4).
\eneq}}
Then, by
\eqref{MUN1-7}, finally,  one has
\begin{equation}\label{n107-30}
\|u^*\phi(a)u-\psi(a)\|<\ep\rforal a\in {\cal F},
\end{equation}
as desired.
\end{proof}

\begin{rem}\label{ReMUN1}

 {{As}} in  Corollary \ref{RemUniCtoA},  the condition  ${\cal U}\subset U(M_2(A))/CU(M_2(A))$ can be replaced
by ${\cal U}\subset J_c(U(M_2(A))/U_0(M_2(A))),$ which generates a torsion free subgroup.  Moreover, {{for part (a),}} or  equivalently, {{the case}}
$A\in {\cal B}_1,$ we may take ${\cal U}\subset J_c(U(A)/U_0(A)),$  {{so that}}  generates a torsion free subgroup. {{(Note that $A$ has stable rank one.)}}

In the case that $A=\overline{\bigcup_{n=1}^{\infty} A_n}~{{\in{\cal B}_1}},$ in the theorem above, one may choose ${\cal U}$ to be
in ${{U(A_n)/CU(A_n)}}$ for some sufficiently large $n.$
\end{rem}

\section{The range of the invariant}

\begin{nota}\label{ktimes}
 In this section we will  use the concept   of {{set}} with multiplicity.\index{set with multiplicity.}
 {\blue{Fix a set  ${\cal X}.$ A subset with multiplicities in ${\cal X}$ is a collection of
 elements in ${\cal X}$ which may be repeated finitely many times.
 Therefore
 $X_1=\{x,x,x,y\}$ is different from $X_2=\{x,y\}.$
 Let $x\in {\cal X}$ be a single element. Denote by
 $x^{\sim k}$  the subset with multiplicities}}
 $\{\underbrace{x,x,...,x}_k\}$ (see 1.1.7 of \cite{Gong-AH}).
  Let $X=\{x_1^{\sim i_1}, x_2^{\sim i_2}, ... , x_n^{\sim i_n}\}$ and $Y= \{x_1^{\sim j_1}, x_2^{\sim j_2}, ... , x_n^{\sim j_n}\}$ (some of the $i_k$'s (or $j_k$'s) may be zero which means the element $x_k$  does not appear in the set $X$ (or $Y$)).
 If $i_k\leq j_k$ for all $k=1,2, ... , n$, then we say that $X\subset Y$ (see 3.21 of \cite{Gong-AH}). We  define
 \beq\label{XcupY}
 &&X\cup Y:=\{x_1^{\sim \max(i_1,j_1)},x_2^{\sim \max(i_2,j_2)}, ... , x_n^{\sim \max(i_n,j_n)}\}\andeqn\\
 \label{XsqcupY}
 &&X\sqcup_{mult} Y=\{x_1^{\sim(i_1+j_1)}, x_2^{\sim(i_2+j_2)},\cdots, x_n^{\sim (i_n+j_n)}\}.
 \eneq
 We may also use the convention
$ \{X, Y\}:=X\sqcup_{mult} Y.$
 For example, $\{x^{\sim 2}, y^{\sim 3}\}=\{x,x,y,y,y\}$ and $\{x^{\sim 2}, y^{\sim 3}, z^{\sim 1}\}=\{x,x,y,y,y,z\}.$
   By $X^{\sim k}$ we mean the set $\{x_1^{\sim ki_1},x_2^{\sim ki_2}, ... , x_n^{\sim ki_n}\}$.
   Note
   $$
   \{X^{\sim k}, Y^{\sim n}\}=X^{\sim k}\sqcup_{mult} Y^{\sim n}.
   $$
 For example, $$\{x,\{x,y\}^{\sim 2},\{y,z\},x, w\}=\{x^{\sim 4},y^{\sim 3},z,w\}=\{x,x,x,x,y,y,y,z,w,\}.$$
\end{nota}

\begin{nota}\label{homrestr}
Let $A$ be  a unital subhomogeneous \CA; that is,  the maximal dimension of irreducible representations of $A$  is finite.

Let us use $Sp(A)$ \index{$Sp(A)$} to denote the set of {equivalence classes of} all irreducible representations  of $A$. {{The set $Sp(A)$ will serve as base set when we talk about a finite set with (finite) multiplicities.}}

Since $A$ is of type I,  the set $Sp(A)$ has {{a}} one-to-one correspondence  to the set of primitive  ideals of $A.$
Let $X\subset Sp(A)$ be a closed subset, then $X$ corresponds to the idea{{l}} $I_X=\bigcap_{\psi\in X} \ker \psi$. In this section, let us use  $A|_X$ to denote the quotient algebra $A/I_X$.
If $\phi: B \to A$ is a homomorphism, then we will use $\phi|_X: B \to A|_X$ to denote the composition $\pi\circ \phi$, where $\pi: A \to A|_X$ is the quotient map. As usual, if $B_1$ is a subset of $B$, we will also use $\phi|_{B_1}$ to denote the restriction of $\phi$ to $B_1$. These two notation{{s}} will not be confused, since it will be clear from {{the context}} which notation we refer to.

If $ \phi: A \to B$ is a homomorphism, then we write  $Sp(\phi)= \{x\in Sp(A): \ker \phi \subset \ker x\}$.\index{$Sp(\phi)$}
{\blue{Denote by}} $RF(A)$ \index{$RF(A)$}
the set of equivalence classes of all (not necessarily irreducible) finite dimensional
 representations. {{ An element $[\pi]\in RF(A)$ will be identified with the set with multiplicity $\{[\pi_1], [\pi_2], \cdots, [\pi_n]\}$ in $Sp(A)$, where $\pi_1, \pi_2,\cdots, \pi_n$ are irreducible representations and $\pi=\pi_1\oplus\pi_2\oplus\cdots\oplus \pi_n$.}}
{{ For $X, Y\in RF(A),$ as sets with multiplicities,}}  {\blue{we write}} $X\subset Y$ {{in the sense of  \ref{ktimes},}} if  {{and only if}} the representation corresponding to  $X$ is equivalent to a sub-representation of that corresponding to   $Y$.  Any finite subset of $RF(A)$ also defines an element of $RF(A)$ which is the equivalence class of  the direct sum of all corresponding representations in the set with the correct multiplicities. { {For example, the subset  $\{X, Y\}\subset RF(A)$ defines an element  $X\oplus Y\in RF(A)$ which is direct sum of two representations corresponding to $X$ and $Y$ (see the {{notation $\{X,Y\}$}} in \ref{ktimes}).}}

{{W}}hen we write $X=\{z_1^{\sim k_1}, z_2^{\sim k_2},..., z_m^{\sim k_m}\},$  we do not insist that $z_i$ be itself in $Sp(A)$.
{\blue{In other words, $z_i$ could be itself a subset {{with}} multiplicity of  $Sp(A).$}}
It {{might}} be a list of several elements {{of}} $Sp(A)$---that is, we do not insist that $z_i$  should  be irreducible (but as we know, it can always be decomposed into {{irreducibles}}). {{So,}} in this notation, we do not differentiate between $\{x\}$ and $x;$ both give the same element of  $RF(A)$ and {{the}} same set with multiplicities whose elements are in $Sp(A)$.

 If $\phi:A \to M_m$ is a homomorphism, let us use $SP(\phi)$  \index{$SP(\phi)$} to denote the corresponding equivalen{{ce}} class of $\phi$ in $RF(A)$.



 {{I}}f $\phi: A\to M_m$ is a homomorphism and if
 $SP(\phi)=\{x_1^{\sim k_1},x_2^{\sim k_2}, ... ,x_i^{\sim k_i}\}$, with $x_1, x_2, ... ,x_i$  {{irreducible}} representations {{ and $k_j>0,$  $j=1,2,\cdots, i$}}, then
 $Sp(\phi)=\{x_1,x_2,...,x_i\}\subset Sp(A).$ So $Sp(\phi)$ is an ordinary set which is a subset of $Sp(A)$, while $SP(\phi)$ is a set with multiplicities, whose elements are also elements of ${Sp(A)}$.
\end{nota}

\begin{nota}\label{density}
Let us recall some notation from Definition \ref{8-N-3}. Suppose that $A=A_m\in{\cal D}_m$ {{(see the end of  \ref{8-N-3})}} is as  constructed in the following {{sequence:}}
\beq\nonumber
&&
A_0={{P_0C(X_0, F_0)P_0}},~~
A_1= {{P_1C(X_1, F_1)P_1\oplus_{Q_1C(Z_1, F_1)Q_1}A_0}},\\
&&\hspace{0.4in}..., A_m=P_mC(X_m, F_m)P_m\oplus_{Q_mC(Z_m, F_m)Q_m} A_{m-1},
\eneq
{\blue{where $F_j=M_{s(j,1)}\oplus M_{s(j,2)}\oplus \cdots \oplus M_{s(j, t_j)}$ is a finite
dimensional \CA,  $P_j\in C(X_j, F_j)$ is a projection, $j=0,1,...,m,$ and}}
$$
{{P_k(C(X_k, F_k))P_k=\bigoplus_{i=1}^{t_k} P_{k,i}
C(X_k, M_{s(k,i)})P_{k,i},\, k=1,2,...,m,}}
$$
{{are}} as in Definition \ref{8-N-3}.
Let $\LD: A \to \bigoplus_{k=0}^mP_kC(X_k, F_k)P_k$ be the {{inclusion}} as in Definition \ref{8-N-3}.
{\blue{Note that  $Z_0=\emptyset.$}}  Let  $\pi_{(x,j)}$, where $x\in X_k$ and
the positive integer $j$ refers to  the
$j$-th block of
$\bigoplus_{i=1}^{{{t}}_k}P_{k,i}C(X_k, F_k)P_{k,i}$ (and
$P_{k,j}(x)\not=0$),  be the  finite dimensional representations of $A$ as in \ref{8-N-3}. According to \ref{homrestr}, one has that $\pi_{(x,j)}\in RF(A)$.
{{If $x\in X_k\setminus Z_k$, then $\pi_{(x, j)}\in Sp(A)$, i.e., it is irreducible. In fact, $Sp(A)=\{\pi_{(x,j)}: x\in X_k\setminus Z_k,
k=0,1,...,m,\,
\andeqn P_k(x)\not=0,\,j=1,2,\cdots, t_k\}$.}} {{Recall that}} all $X_k$ are {{compact}} metric spaces. {{For $k\leq m$ and $j\leq t_k$, we use  $X_{k,j}$
to denote all the (not necessary irreducible) non-zero representations ${\blue{\pi_{(x,j)}}}$ for $x\in X_k$ {\blue{(see \ref{8-N-3}).}}
{\blue{Note that, as a set, $X_{k,j}=\{x\in X: P_{k,j}(x)\not=0\}.$}}
{\blue{Set $Z_{k,j}=\{\pi_{(x,j)}\in X_{k,j}: x\in Z_k\}.$}}
{\blue{Note that}} $X_{k,j}$ {\blue{(and $Z_{k,j}$, respectively)}}
has a natural metric induced from $X_k$---{{i.e.,}} $\dist({\blue{\pi_{(x,j)},\pi_{(y,j)}}})=\dist(x,y)$ for  $x, y \in X_k$. }}
{\blue{In what follows, if $\theta\in RF(A)$ and $\theta=\pi_{(x,j)}$ for some $\pi_{(x, j)}\in X_{k,j},$ then
we will write $\theta\in X_{k,j}.$}}

  {{The {\blue{usual}} topology on $Sp(A_m)$ is
  in general  not Hausdorff.
  For each $\theta\in Sp(A)$ {\blue{and $\dt>0,$
  the subset $B_{\dt}(\theta)$ of
   $Sp(A),$}}
   called  the $\delta$-{\blue{neighborhood}} of $\theta,$  is defined as {{follows}}:}}

  {{An irreducible representation $\Theta\in Sp(A_m)$ is in $B_{\delta}(\theta)$ if there is a sequence of finite
  {\blue{dimensional}} representations $\theta=\sigma_0, \theta_0, \sigma_1, \theta_1, \sigma_2, \theta_2,\cdots, \sigma_k, \theta_k=\Theta$ such that\\
  1. For each  $i=1,2,\cdots, k$,
  $\theta_{i-1}\subset
  \sigma_i$ (see \ref{homrestr});\\
  2. For each  $i=0,1,\cdots, k$,
  there exists a pair $(l_i, j_i)$
  such that $\sigma_i\in X_{l_i, j_i}$ and $\theta_i\in X_{l_i, j_i}\setminus Z_{l_i, j_i}$;\\
  3. $\sum_{i=0}^{k}\dist_{X_{l_i, j_i}}(\sigma_i, \theta_i)<\delta.$\\
  {\blue{Note that, since $\theta_{i-1}\in X_{l_{i-1}, j_{i-1}}\setminus Z_{l_{i-1}, j_{i-1}},$   $\sigma_i\in X_{l_i, j_i},$
  and $\theta_{i-1}\subset \sigma_i,$ one must have $l_i\ge l_{i-1},$ $i=1,2,...,k.$}}
  Note  also that we borrow the concept  of $\delta$-neighborhood from metric space theory, but there is no metric on $Sp(A)$. In fact $\Theta \in B_{\delta}(\theta)$ does not imply that $\theta\in B_{\delta}(\Theta)$. Note that if $\Theta\in B_{\delta_1}(\Theta_1)$ and $\Theta_1\in B_{\delta_2}(\theta)$, then $\Theta\in B_{\delta_1+\delta_2}(\theta)$.}}


  Let $\phi: A \to M_{\bullet}$ be a homomorphism. We {{shall}} say that $SP(\phi)$
  is $\dt$-dense in $Sp(A)$ \index{$\dt$-dense in the spectrum of subhomogeneous algebra}
   if for each irreducible representation $\theta$ of $A$, {{there is $\Theta\in B_{\delta}(\theta)$ such that $\Theta\subset SP(\phi)$. For any algebra in ${\cal D}_1$ (see the end of  \ref{8-N-3}){{, including}} all  Elliott-Thomsen building blocks, $\delta$-density of $SP(\phi)$ means}} for each irreducible representation $\theta$ of A, {{either $\theta\subset SP(\phi)$}}{{, or}} there are two points  $x, y \in  {{ X_1}}$ {{(for the case of an Elliott-Thomsen building block, $ X_1=[0,1]$)}} and $j$ (a single  $j$) such that $\dist (x,y) <\dt$, and such that
$\tht \subset \pi_{(y,j)}$ and $\pi_{(x,j)}\subset SP(\phi)$. This will be used in this section and {\blue{the}} next.

If $X_k=[0,1]$, then we use the standard metric of the interval $[0,1]$.
\end{nota}

\begin{lem}\label{lem-simple}
{\blue{Let $A\in {\cal D}_m.$}}

{\blue{(1) If $f\in A$ and $\theta\in Sp(A)$ such that $\theta(f)\not=0,$  then there exist $\dt>0$ and $d_0>0$ such that
$\|\Theta(f)\|\ge d_0$ for all $ \Theta\in B_\dt(\theta).$}}

{{(2) 
For any $\dt>0$, there is a finite set ${\cal F}\subset A_+\setminus \{0\}$ satisfying the following condition. For any irreducible representation $\theta \in Sp(A)$, there is an element
$f\in {\cal F}$ such that if $\Theta\in Sp(A)$ satisfies  $\Theta(f)\not=0$, then $\Theta\in B_{\delta}(\theta)$. (We do not require that $\theta(f)\not=0$ for the element $f$ corresponding to $\theta$.)}}

\end{lem}

\begin{proof}
{\blue{In what follows we will keep the notation introduced in \ref{density}.}}

 {\blue{ For part (1), we
 will prove it by  induction.
  It is clear that part (1)  holds when $A\in {\cal D}_0.$
  Assume that part (1)  holds for all $A\in {\cal D}_k$ with $k< m.$}}

 {\blue{ Let $A\in {\cal D}_m$ and let
  $\theta(f)\not=0.$
    Write $A_m=P_m(C(X_m, F_m)P_m\oplus_{Q_mC(Z_k,F_m)Q_m}A_{m-1}$ and
  $f=(g, h),$ where $g\in P_m(C(X_m, F_m)P_m$ and $h\in A_{m-1}.$
  If $\theta\in \{\pi_{(x,j)}: x\in X_m\setminus Z_m\andeqn P_{m,j}(x)\not=0\},$
 for some $j,$ then, by continuity of $g$ at $x$ in the $j$th block,  one obtains  $\dt>0$
 such that $\|\pi_{(y,j)}(f)\|\ge d_{0}$ for all $(y,j)\in X_{m,j}$ such that ${\rm dist}(y, x)<\dt,$
 where $d_{0}>0.$
 In other words,  $\|\Theta(f)\|\ge d_{0}$ for all $\Theta\in B_{\dt}(\theta).$}}

 {\blue{ Suppose $\theta\in \{\pi_{(x,j)}: x\in X_k\setminus Z_k\andeqn P_{k,j}(x)\not=0\}$ for some $j$ and
  for some $0<k<m.$
  Thus, one may view $\theta$ as a point in
  $A_{m-1}$ and $\theta(h)\not=0.$ By the induction assumption, there are $\dt_0>0$  and $d_{00}>0$ such that, for
  any $\Theta'\in   B^{m-1}_{\dt_0}(\theta),$ $\|\Theta'(h)\|\ge d_{00}$ (where
  $B^{m-1}_{\eta}(\theta)$ is a $\eta$-neighborhood of $Sp(A_{m-1})$ for any $\eta>0$).
  Note that we may also view $Sp(A_{m-1})$ as a subset of $Sp(A).$
For each $j=1,2,..., t_m,$ let
$$
Y_{m,j}=\{\pi_{(z, j)}\subset Z_{m,j}:  \exists~\theta'\in B^{m-1}_{\dt_0/2}(\theta)\,\, \mbox{such that}~~\theta'\subset \pi_{(z,j)}\}.
$$
If $\pi_{(x,j)}\in Y_{m, j},$ then $\|\pi_{(x, j)}(f)\|\ge \|\theta'(f)\|\ge d_{00}.$
Note that $g|_{Z_m}=\Lambda_m(h).$  In the $j$th block, $\|\pi_{(x,j)}(f)\|\ge d_{00}$
for all $\pi_{(x,j)}\in \overline{Y_{m,j}},$ the closure of $Y_{m,j}.$
Since $\overline{Y_{ m,j}}$ is compact, there is $\dt_j>0$ such that
$\|\pi_{(x,j)}(f)\|\ge d_{00}/2$ for all $\pi_{(x, j)}\in X_{m,j}$ such that
${\rm dist}(x, Y_{m,j})<\dt_j.$
%
If $Y_{m,j_0}=\emptyset,$ put $\dt_{j_0}=1.$
Choose $\dt=\min\{1/2, \dt_0/2, \dt_j: 1\le j\le t_m\}.$}}

{\blue{We claim that $\|\Theta(f)\|\ge d_{00}/2$ for all $\Theta\in B_{\dt}(\theta).$
Let  $\Theta\in X_{k,j}$ for some $j.$ Then, by the definition of $B_{\dt}(\theta),$ there exist
finite dimensional  representations $\theta=\sigma_0, \theta_0, \sigma_1, \theta_1, \sigma_2, \theta_2,\cdots, \sigma_n, \theta_n=\Theta$  (for some integer $n\ge 1$) such that\\
  1. For each  $i=1,2,\cdots, n$,
  $\theta_{i-1}\subset
  \sigma_i;$\\
  2. For each  $i=0,1,\cdots,n$,
  there exist  a pair $(l_i, j_i)$
  such that $\sigma_i, \theta_i\in X_{l_i, j_i}$  $\theta_i\in X_{l_i, j_i}\setminus Z_{l_i,j_i};$ and\\
  3. $\sum_{i=1}^{n}\dist_{X_{l_i, j_i}}(\sigma_i, \theta_i)<\delta.$}}\\
  {\blue{Note, by  \ref{density}, if $\Theta\in B_{\dt}(\theta)\cap Sp(A_{m-1}),$
  then all $l_i<m.$ It follows that $\Theta\in B_{\dt_0}^{m-1}(\theta).$ Then, by
the choice of $\dt_0,$  $\|\Theta(f)\|\ge d_{00}.$}}

 {\blue{ Otherwise, $\Theta\in X_{m,j}.$ Choose  the largest $i$ such that  $l_i<m.$
  Then $\theta_i\in X_{l_i, j_i}\setminus Z_{l_i, j_i}$ and $\sigma_{i+1}\in X_{m, j_{i+1}}$ for some $j_{i+1}.$
  Therefore, by the definition of $B_{\dt}^{m-1}(\theta),$  $\theta_i\in B_{\dt}^{m-1}(\theta).$
  Note that, for any $i'>i,$ $l_{i'}=m.$ It follows that $\sigma_{i'},\, \theta_{i'}\in X_{m, j_{i'}}.$
 On the other hand,  since $\theta_n=\Theta\in X_{m,j},$
  $\sigma_n\in X_{m,j}.$ Since $\theta_{n-1}\subset \sigma_n,$
  either $l_{n-1}<m,$ in which case, $i=n-1,$ or $\theta_{n-1}=\sigma_n$ as $\theta_{n-1}\in X_{l_{n-1}, j_{n-1}}\setminus Z_{l_{n-1}, j_{n-1}},$ in which case $l_{n-1}=m$ and $j_{n-1}=j.$ By repeatedly using 1,2,3, above, one concludes that
  $j_{i'}=j$ for $i'>i.$
  In other words,  $\sigma_{i+1}\in X_{m,j}.$ Note that $\theta_i\subset \sigma_{i+1}$ and $\theta_i\in  B_{\dt}^{m-1}(\theta).$
  It follows that $\sigma_{i+1}\in Y_{m,j}.$
  One checks (by 3 above) that ${\rm dist}(\Theta, \sigma_{i+1})={\rm dist}(\theta_n, \sigma_{i+1})<\dt<\dt_0/2.$
  By the choice of $\dt,$ one obtains $\|\Theta(f)\|\ge d_{00}/2.$
  This completes the induction and part (1)  holds}}.

{{We will {\blue{also}} prove part (2)  by induction.}} {\blue{ If $A=A_0=\bigoplus_{j=1}^{t_0}P_{0,j}C(X_0,M_{s(0,j)})P_{0,j} \in {\cal D}_0,$
then it is easy to see that this reduces to the case that $A=A_0=P_0C(X_0, M_s)P_0,$ where $s\ge 1$ is an integer
and $P_0\in C(X_0, M_s)$ is a projection. Given $\dt>0,$ let $\{U_i: 1\le i\le m\}$ {{be}}  an open cover
of $X_0$ with the diameter of each $U_i$
smaller than $\dt/2.$ Consider
a partition of unity ${\cal F}=\{f_i\in C(X_0): 1\le i\le m\}$
subordinate to the open cover $\{U_i: 1\le i\le m\}.$  Then, clearly
${\cal F}$ satisfies the requirements.}}
{{ Let us assume that $A=A_m\in {\cal D}_m$ with $m\geq 1$ and that the conclusion of the lemma is true for algebras in ${\cal D}_{m-1}$.}}

{{Write $A_m=P_mC(X_m, F_m)P_m\oplus_{Q_mC(Z_m, F_m)Q_m}A_{m-1}$ with $\Gamma_{m}: A_{m-1}\to Q_{m}C({ Z_{m}}, F_{m})Q_{m}$, where
 $Q_m=P_m|_{Z_m}$ (see \ref{8-N-3}). By {\blue{the}} induction assumption, there is a {\blue{finite
 subset}}  ${\cal F}_1\subset (A_{m-1})_+\setminus \{0\}$ with the following
  property: For any irreducible representation $\theta \in Sp(A_{m-1})$, there is an element $f\in {\cal F}_1$ such that {\blue{if}} $\Theta\in Sp(A_{m-1})$ satisfies {{$\Theta(f)\not=0$}}, then $\Theta\in B^{m-1}_{\delta/2}(\theta)$ (notation from
  the proof of part (1)).}}
  For each $f\in {\cal F}_1$, by the Tietze Extension
   Theorem, there is an $$h=(h_1,h_2,\cdots, h_{t_m})\in P_mC(X_m, F_m)P_m=\bigoplus_{j=1}^{t_m}P_{m,j}C(X_m, M_{s(m,j)})P_{m,j}$$ such that $\Gamma_{m}(f)=h|_{Z_m}$. For each $j\leq t_m$, let $\Omega_{m,j}\subset Z_m~(\subset X_m)$ be the closure of the 
   set
  $$\{z\in Z_m:~~\exists~\theta' \in Sp(A_{m-1})~~\mbox{such that}~~\theta'\subset \pi_{(z,j)}~~\mbox{and}~~\theta'(f)\not = 0\},$$
 {\blue{where $\theta\subset \pi_{(x, j)}$ as subsets of $RF(A)$ with multiplicities.}}
{{
Choose, {\blue{in the case $\Omega_{m,j}\not=\emptyset,$}}  a continuous function $\chi_j:X_m\to [0,1]$ such that $$\chi_j(x)=1~~\mbox{if}~~x\in \Omega_{m,j}~~~\mbox {and}~~ \chi_j(x)=0~~\mbox{if} ~~
\dist(x, \Omega_{m,j})\geq {\blue{\delta/3}},$$
{\blue{and, in the case that $\Omega_{m,j}=\emptyset,$
let $\chi_j=0.$
Set}} $h'=(\chi_1\cdot h_1, \chi_2\cdot h_2,\cdots,\chi_{t_m}\cdot h_{t_m}).$
Then $h'|_{Z_m}={\blue{h|_{Z_m}}}=\Gamma_{m}(f).$}}
{\blue{Thus,}} {{${\tilde f}=(h',f)\in P_mC(X_m, F_m)P_m\oplus A_{m-1}$ defines an element of $A_m=P_mC(X_m, F_m)P_m\oplus_{Q_mC(Z_m, F_m)Q_m}A_{m-1}$.}}

{\blue{Now let $\Theta\in Sp(A_m)$ be such that $\Theta({\tilde f})\not=0.$
If $\Theta\in X_{m,j}\setminus Z_{m,j}$ for some $j,$
then $\Omega_{m,j}\not=\emptyset$ and ${\dist}(\Theta, \Omega_{m,j})<\dt/3.$
Therefore, there exists $\sigma\in Z_{m,j}$ such that ${\rm dist}(\sigma, \Theta)\le \dt/3$
and  $\Theta'\in Sp(A_{m-1})$ such that $\Theta'\subset \sigma$ and $\Theta'(f)\not=0.$
By the inductive assumption, $\Theta'\in B_{\dt/2}(\theta).$
By the definition of $B_{\dt}(\theta),$
this implies that $\Theta\in B_{\dt}(\theta).$
}}
{\blue{If $\Theta\in Sp(A_{m-1}),$ then $\Theta\in B_{\dt/2}^{m-1}(\theta),$ which also implies
that $\Theta\in B_{\dt/2}(\theta).$}}
Let ${\cal {\tilde F} }_1=\{\tilde{f}: f\in {\cal F}_1\}.$

{{{\blue{Choose}} a finite subset $\Xi:=\{x_1, x_2, \cdots, x_{\bullet}\}\subset X_m $ which is $\delta/4$-dense in $X_m$---that is, if $x\in X_m$, then there is $x_i$ such that $\dist(x, x_i)<\delta/4$. We need to modify the set $\Xi$ so that $\Xi\subset X_m\setminus Z_m$. If $x_i\in Z_m$ and
$$
W:=\{x\in (X_m\setminus Z_m): \dist (x, x_i)<\delta/4\}\not=\emptyset,
$$
then replace $x_i$ by any element of $W$;  if $W$ is the empty set, then simply delete $x_i$.  After the modification, we have $\Xi\subset X_m\setminus Z_m$ and $\Xi$ is $\delta/2$ (instead of $\delta/4$) dense in $X_m\setminus Z_m$ (instesd of $X_m$). For each $x_i\in \Xi$, choose an open set $U_i\ni x_i$ such that $U_i\subset X_m\setminus Z_m$ and ${\blue{\dist}}(x, x_i)<\delta/2$ for any $x\in U_i$. For each $x_i\in \Xi$, choose a function $g_{x_i}:X_m \to [0,1]$ such that $g_{x_i}(x_i)=1$ and $g_{x_i}(x)=0$ for $x\notin U_i$. For $j\leq t_m$ let
$$
g_{x_i, j}= (\underbrace{0,\cdots, 0}_{j-1}, g_{x_i}\cdot P_{m,j}, \underbrace{0,\cdots, 0}_{t_m-j})\oplus 0 \in \bigoplus_{j=1}^{t_m}P_{m,j}C(X_m, M_{s(m,j)})P_{m,j}\oplus A_{m-1},$$
which defines an element (still denoted by $g_{x_i, j}$) of $A_m$. Define ${\cal G}=\{g_{x_i, j}: x_i\in \Xi, j\leq t_m\}$. {\blue{Set}} ${\cal F}={\cal G}\cup {\cal {\tilde F} }_1.$ }}

{\blue{Now fix $\theta\in Sp(A).$ If $\theta\in Sp(A_{m-1}),$ then choose  $f\in {\cal F}_1$ with the property
$\Theta(f)\not=0$ implies  $\Theta\in B_{\dt/2}^{m-1}(\theta).$  Now consider ${\tilde f}.$
By what has been proved, if $\Theta({\tilde f})\not=0,$ then $\Theta\in B_{\dt}(\theta).$}}

{\blue{If $\theta\in X_{m,j}\setminus Z_{m,j},$   then there is $x_i\in \Xi,$ such that
${\rm dist}(\theta, x_i)<\dt/2.$  Then $g_{x_i, j}\in {\cal F}$ and
$\pi_{(x_i,j)}(g_{x_i, j})\not=0.$ Moreover, by the construction of $g_{x_i,j},$ if  $\Theta\in Sp(A)$ with
$\Theta(g_{x_i, j})\not=0,$  then $\Theta\in B_{\dt}(\theta).$  This ends the induction.}}



\end{proof}

{{The special case of the following lemma for {\blue{AH-algebras can be found in}}
 \cite{DNNP-AH}.}} {\blue{In the following statement, we use notation introduced in \ref{homrestr}.}}

\begin{prop}\label{simplelimit}

Let $A=\varinjlim (A_n,\phi_{n,m})$ be a
{{unital}} inductive limit {\blue{of \CA s, where
$A_n\in {\cal D}_{\tiny{l(n)}}$ (for some $l(n)$) and where each $\phi_{n,m}$ is injective.}}
{\blue{Then the}}  limit $C^*$-algebra is simple if {{and only if}} the inductive system satisfies the following condition: for any $n>0$ and $\dt>0$, there is an integer $m>n$ such that for any ${{\sigma}}\in Sp(A_m)$, $SP(\phi_{n,m}|_{{\sigma}})$ is $\dt$-dense
 {\blue{in $Sp(A_n)$}}---equivalently, for any $m'\ge m$ and any ${{\sigma'}}\in Sp(A_{m'})$, $SP(\phi_{n,m'}|_{{\sigma'}})$ is $\dt$-dense {{in $Sp(A_n).$}}
\end{prop}

\begin{proof}
{\blue{Note that, since $A$ is unital and each $A_n$ is unital, \wilog, we may assume
that all $\phi_{n,m}$ are unital.}}
 The proof of this {{proposition}}  is standard. {{Suppose that the condition holds.}}
  {{F}}or any non-zero element $f\in A_n$, {{there is  $\theta\in Sp(A_n)$ such that $\theta(f)\not=0$. Consequently,}}
  {\blue{by part (1) of \ref{lem-simple},}} {{there is $B_\dt(\theta)$}} (for some $\dt>0$) {\blue{such that $\Theta(f)\not=0$
  for any $\Theta\in B_{\delta}(\theta).$}}
  {{Then, by the given condition, there is an {\blue{integer}} $m$ such that for any $m'\ge m$ and any  irreducible representation}} ${{\sigma}}$ of $A_{m'}$, one has {\blue{$Sp(\phi_{n,m'}|_\sigma)\cap B_\dt(\theta)\not=\emptyset.$
  Let $\Theta\in Sp(\phi_{n,m'}|_\sigma)\cap B_\dt(\theta).$ Then, $\Theta\circ \phi_{n,m'}(f)\not=0.$}}
  {\blue{Consequently, $\|\sigma\circ \phi_{n,m'}(f)\|\ge \|\Theta\circ \phi_{n,m'}(f)\|>0.$}}
  It follows that the ideal $I$ generated by $\phi_{n,m'}(f))$ in $A_{m'}$ equals $A_{m'}$---otherwise, any irreducible representation ${{\sigma}}$ of $A_{m'}/I\not=0$ (which is also an irreducible representation of $A_{m'}$) satisfies ${{\sigma}}(\phi_{n,m'}(f))=0$, and this contradicts  the fact that ${{\sigma}}(\phi_{n,m'}(f))\not=0$ for every irreducible representation ${{\sigma}}$.
  {\blue{Fix $m'\ge m.$ Then $\phi_{n,m'}(f)$ is full.  Since each $\phi_{n,m}$ is unital,
  it follows that $\phi_{m', \infty}\circ \phi_{n,m'}(f)=\phi_{n,\infty}(f)$ is full.
  In other words (since  $f$ above is arbitrary), for any  proper ideal  $I$ of $A,$ $\phi_{n,\infty}(A_n)\cap I=\{0\}.$ It is standard to show
  that this implies that $A$ is simple.
 In fact, for any $a\in \phi_{n,\infty}(A_n),$ $\pi_I(\|a^*a\|)=\|a^*a\|,$ where $\pi_I: A\to A/I$
  is the quotient map. It follows that $\|\pi_I(a)\|=\|a\|.$  Note that $\bigcup_{n=1}^{\infty}(A_n)$ is dense
  in $A.$  {{The quotient map $\pi_I$ is thus isometric,  and so $I=\{0\}.$ (This proof is due to  Dixmier  (see the proof
  of Theorem 1.4 of \cite{Dix}).)}}
 Therefore $A$ is simple.}}


 {{Suppose the unital limit algebra $A$ is simple. For any $A_n$ and $\delta>0$, let ${\blue{{\cal F}_n\subset (A_n)_+}}\setminus \{0\}$ be as in the Lemma \ref{lem-simple}. Since $A$ is simple and $\phi_{n,m}$  is injective, for any $f\in {\blue{{\cal F}_n}}$, the ideal generated by $\phi_{n, \infty}(f)\in A$ contains ${\bf 1}_A$. Hence there is $m_f> n$ such that if $m'\geq m_f$ then the ideal generated by $\phi_{n,m'}(f)$ in $A_{m'}$ contains ${\bf 1}_{A_{m'}}$ and therefore is  the whole of $A_{m'}$. Let $m=\max\{m_f: f\in {\blue{{\cal F}_n}}\}$.   For any $\sigma \in Sp(A_{m'})$ (with $m'\geq m$) and $f\in {\blue{{\cal F}_n}}$, we have $\sigma(\phi_{n,m'}(f))\not= 0$.}}

 {{We are  going to verify that for any $\sigma \in Sp(A_{m'})$, {\blue{the}} set $SP(\phi_{n,m'}|_{\sigma})$ is $\dt$-dense
 {\blue{in $Sp(A_n).$}} For any $\theta\in Sp(A_n)$, there is an  $f\in {\blue{{\cal F}_n}}$ {\blue{such that, if $\Theta(f)\not=0,$
 then $\Theta\in B_{\dt}(\theta).$}}
 From $\sigma(\phi_{n,m'}(f))\not= 0$, one {\blue{then}} concludes that there is an irreducible representation $\Theta\subset SP(\phi_{n,m'}|_{\sigma})$ such that $\Theta(f)\not=0$. Hence $\Theta\in B_{\delta}(\theta)$ by Lemma \ref{lem-simple}. {\blue{This implies that}} $SP(\phi_{n,m'}|_{\sigma})$ is $\dt$-dense. }}
\end{proof}

\begin{df}\label{Class0}
Denote by ${\cal N}_0$ \index{${\cal N}_0$} the class of those unital  simple \CA s
$A$  in ${\cal N}$ for which $A\otimes U\in {\cal
N}\cap {\cal B}_0$ for any  UHF-algebra $U$ of
infinite type (see \ref{dfcalN} for the definition of {{the}}  class ${\cal N}$).

Denote by ${\cal N}_1$ \index{${\cal N}_1$} the class of those unital \CA s $A$ in ${\cal N}$ for which $A\otimes U\in {\cal N}\cap {\cal B}_1$ for any  UHF-algebra of
infinite type. In Section 19, we will show that ${\cal N}_1={\cal N}_0.$

Also denote by ${\cal N}_0^{\cal Z}$ \index{${\cal N}_0^{\cal Z}$}(respectively, ${\cal N}_1^{\cal Z}$ \index{${\cal N}_0^{\cal Z}$}) the class of all ${\cal Z}$-stable $C^*$-algebras in ${\cal N}_0$ (respectively, ${\cal N}_1$).

\end{df}


\begin{NN}\label{range 0.1}
Let $(G,G_+,u)$ be a scaled ordered abelian group $(G,G_+)$ with order
unit $u\in G_+\setminus \{0\},$ with
{the scale  given by
$\{g\in G_+: g\le u\}$. Sometimes we will also call $u$ the scale of the group.}
Let $S(G):=S_u(G)$ be the state space of $G.$
Suppose that
$((G,G_+,u),K,\Delta,r)$ is a weakly  unperforated Elliott invariant---that is, $(G,G_+,u)$ is a simple scaled ordered countable group,
$K$ is a countable
abelian group, $\DT$ is a metrizable Choquet simplex, and $r: \DT \to S(G)$
is a surjective affine map such that for any $x\in G$,
\beq\label{1508/star13-1}
 x\in G_+\setminus \{0\}\, \mbox{ if and only if } \,\, r(\tau)(x)>0
 \tforal \tau\in \DT.
\eneq
{{C}}ondition (\ref{1508/star13-1})  above is also called  weak  unperforation \index{weakly unperforated} for the
simple ordered group. (Note that this condition is equivalent to  the condition that
$x\in G_+\setminus \{0\}$ if and only if for any $f\in S(G)$,
$f(x)>0$.  The latter condition does not mention Choquet simplex
$\DT$.) In this paper, we only consider the Elliott invariant for stably finite simple {{unital nuclear}} $C^*$-algebras and therefore $\DT$ is not empty. {{Evidently, t}}he above weak unperforation condition {{implies the condition}} that $x>0$ if $nx>0$ for some positive integer $n$. {{The converse follows from   Proposition 3.2 of \cite{RorUHF2}.}}

In this section, we will show that for any weakly unperforated Elliott invariant \\
$((G,G_+,u),K,\Delta,r)$, there is a unital simple \CA\, $A$ in the class ${\cal N}_0^{\cal Z}$ such that $$ ((K_0(A),K_0(A)_+,
[{\mbox{\large \bf 1}}_A]),K_1(A),T(A),r_A)\cong ((G,G_+,u),K,
\Delta,r).
$$
{\blue{A similar general range theorem  was presented by Elliott
in \cite{point-line}.
To obtain our version, we will modify the construction given by Elliott in \cite{point-line}.
  Our modification  will  reveal more details in the   construction  and
  will also ensure that the algebras constructed are actually in the class ${\cal N}_0^{\cal Z}.$}}
One difference is that, at an  important step of the construction, we will use
a finite subset $Z_i$  of a certain space $X_i$  instead of a one-dimensional subspace of $X_i.$

\end{NN}

\begin{NN}\label{range 0.2}  Our construction will be a modification of the Elliott
construction mentioned above.  As {{a}} matter
of fact, for the case that $K=\{0\}$ and $G$ {{is}} torsion free, our
construction uses the same building blocks, in ${\cal C}_0,$ as in \cite{point-line}.  We will repeat a part of the
construction of Elliott for this case.  There are two steps in
Elliott{{'s}} construction:

Step 1.  Construct an inductive limit
$$ C_1\lr C_2 \lr \cd\lr C $$
with inductive limit of ideals
$$ I_1\lr I_2 \lr \cd\lr I $$
such that the non-simple limit $C$ has the prescribed Elliott invariant
and the quotient $C/I$ is a simple AF algebra. For the case  $K=\{0\}$ and $G$ torsion free, we will use the notation $C_n$ and $C$ for the construction, and reserve $A_n$ and $A$ for the general case.

Step 2.  Modify the above inductive limit to make $C$ (or $A$ in the general case) simple without
changing the Elliott invariant of $C$ (or $A$).

For {{the}} reader's convenience, we will repeat Step 1 of Elliott{{'s}} construction with minimum modification. For Step 2,
{\blue{we are not able to reconstruct
the
argument of \cite{point-line} and in particular we do not know  how to make the eigenvalue patterns given on pages 81--82
\cite{point-line}   satisfy the required boundary conditions. We}}
will use a  different way of modifying the inductive limit {{(see \ref{range 0.32} for the details)}}. {{This new way will}} be {{also}} more suitable for our purpose---{{i.e.,}} to construct an inductive limit $A\in {\cal
N}_0^{\cal Z}$ in  the general case, with  possibl{{y}} nontrivial $K_1$ and nontrivial
${\rm Tor}(K_0(A))$.

\end{NN}

\begin{NN}\label{range 0.3}
{\blue{Suppose that $((G, G_+, u), K, \DT, r)$ is a weakly unperforated Elliott invariant as defined in \ref{range 0.1}.}}
Let $\rho:~ G\to \Aff(\DT)$ be the
dual map of $r:\DT \to S(G)$.  That is, for every $g\in G, ~
\tau\in \DT,$
$$
\rho (g)(\tau)= r(\tau)(g)\in \R.
$$
Since $G$ is weakly unperforated, one has that $g\in G_+\setminus \{0\}$ if and only if $\rho (g)(\tau)>0$
for all $\tau\in \DT$.  Note that $\Aff(\DT)$ is an ordered vector
space with the strict (pointwise) order, i.e., $f\in \Aff(\DT)_+\setminus\{0\}$ if and only if $f(\tau)
>0$ for all $\tau\in\DT$.   {\blue{Note that $\DT$ is a metrizable compact convex set, and
$\Aff(\DT)$ is  a norm closed subspace of $C_\R(\DT),$ the real Banach space
of real continuous functions on $\DT.$ Consequently $\Aff(\DT)$  is separable.
{{We assume that $G\not=\{0\}$  and  $\Delta\not=\emptyset.$}}
Therefore}} there is a countable dense subgroup $G^1\subset \Aff(\DT).$   Put
$H=G\oplus G^1$ and define  { ${\tilde \rho}: H\to \Aff(\DT)$ by
${\tilde \rho}((g,f))(\tau)=\rho(g)(\tau)+f(\tau)$ for all $(g,f)\in G\oplus G^1$ and
$\tau\in \DT.$ Define
$H_+\setminus\{0\}$ to be the set of  elements
$(g,f)\in G\oplus G^1$ with ${\tilde \rho}((g,f))(\tau)>0$ for all $\tau\in \DT.$}
The order unit (or scale) $u\in G_+,$ regarded as $(u,0)\in G\oplus G^1=H,$
is an order unit for  $H_+$ (still denote it by $u$).

 Then $(H, H_+,u)$
is a simple  ordered group {\blue{with the Riesz
interpolation property. With the strict order,
$(H, H_+,u)$ is a simple ordered group.}}
{\blue{Since $\DT$ is  a simplex, by Corollary II.3.11 of \cite{Alf}, $\Aff(\DT)$ has the weak Riesz interpolation property.}}
Since ${\tilde \rho}(H)$ is  dense,
it is {{straightforward}} to prove that
$(H, H_+,u)$ is a  Riesz group. Let us  give a brief proof of this fact. Let $a_1, a_2, b_1, b_2\in H$ with $a_i<b_j$ for $i,j\in\{1,2\}.$
Then, {\blue{since $\Aff(\DT)$ has the weak Riesz interpolation property (see page 90 of \cite{Alf}),  there exists $f\in \Aff(\DT)$ such
that ${\tilde\rho}(a_i)<f<{\tilde \rho}(b_j),$ $i,j\in\{1,2\}.$ }}
Let
$$
{\blue{d=\min\{\min\{\tau(f)-\tau(a_i):\tau\in \DT,\,\,1\le i\le 2\}, \min\{\tau(b_j)-\tau(f): \tau\in \DT,\,\, 1\le j\le 2\}\}.}}
$$
 {\blue{Then $d>0.$ Since ${\tilde \rho}(H)$ is dense in $\Aff(DT),$ there exists
$h\in H$ such that $\|{\tilde \rho}(h)-f\|<d/2.$ Then
${\tilde \rho}(a_i)<{\tilde \rho}(h)<{\tilde \rho}(b_j),$ $i,j\in \{1,2\}.$ }}
Consequently, $a_i<h<b_j$ for $i,j\in\{1,2\}$.

As a
direct summand of $H$, the subgroup $G$ is {\it relatively divisible}
subgroup of $H;$ that is, if $g\in G,$  $m\in \N\setminus \{0\},$ and $h\in H$
such that $g=mh,$ then there is $g'\in G$ such that $g=mg'.$\index{relatively divisible subgroup}
Note that $G\subsetneqq H$ since $\Delta\not=\emptyset.$

{{Now we assume, until \ref{range 0.14}, that $G$ is torsion free and $K=0.$}}
Then  $H$ is also torsion free. Therefore, $H$ is a simple dimension group.

\begin{rem}\label{March-22-2019}

{{(1)}} {\blue{There is a unital simple AF-algebra $B$ such that $(K_0(B), K_0(B)_+, [1_B])=
(H, H_+, u)$ and $S_u(H)=\DT=T(B),$  where $S_u(H)$ is the state space of $H$  and
$T(B)$ is the tracial state space of $B.$}}
{\blue{ In fact, by \cite{EHS},  there is a unital simple
AF-algebra $B$ such that $(K_0(B), K_0(B)_+, [1_B])=(H, H_+, u).$}}
{\blue{
It follows from Theorem III.1.3 of \cite{BH}  that the state space $S_u(H)$ of $K_0(B)=H$ is $T(B)$ (with the
topology of pointwise convergence on $S_u(H)$ and the weak* topology of $T(B)$).}}
{\blue{On the other hand, evaluation at a point of $\DT$ also gives a state of $H.$
Therefore $\DT$ is a closed convex subset of $S_u(H).$
Since $x<y$ in $H$ if and only if ${\tilde \rho}(x)(s)<{\tilde \rho}(y) (s)$ for all
$s\in \DT,$ by Lemma 2.9 of \cite{Blatrace}, $\DT=S_u(H)=T(B).$ }}{\blue{Furthermore, the map $\rho_B: K_0(B)=H\to \Aff(T(B))=\Aff(\DT)$ is the same as ${\tilde \rho}: H\to \Aff(\DT)$.}}

{\blue{(2) Suppose that $A$ is a unital $C^*$-algebra and $\phi: A\to B$ is a unital homomorphism (where $B$ is as in  part (1)). Suppose that $(K_0(A), K_0(A)_+, [1_A])=
(G, G_+, u)$,  and $T(A)\cong T(B)$ and suppose that the induced maps $\phi_{*0}: K_0(A) \to K_0(B)  $  and $T(\phi): T(B) \to T(A)$ (of the homomorphism $\phi$) are the inclusion from $G$ to $H$  and the affine
homeomorphism between $T(B)$ and $T(A)$. Then {{the map}}  $\rho_A: K_0(A)=G\to\Aff(T(A))=\Aff(\DT)$ is the same as $\rho: G \to \Aff(\DT)$, under the identification of $T(A)=T(B)=\DT$. This is true because
$\rho={\tilde \rho}\circ \iota_{G, H}: G\to \Aff(\DT)$, where $\iota_{G, H}: G\to H$ is the inclusion.  }}

\end{rem}

\end{NN}


\begin{NN}\label{range 0.4} In \ref{range 0.3} we can choose the dense
subgroup $G^1\subset  \Aff(\DT)$ to contain at least three elements
$x,y,z\in \Aff(\DT) $ such that $x,y$ and $z$ are $\Q$-linearly
independent.
With this
choice, when we write $H$ as the inductive limit of a sequence
$$
H_1\lr H_2\lr \cd
$$
of finite  direct sums of copies of the ordered group $(\Z,{{\Z}}_+)$ as in Theorem 2.2
of \cite{EHS}, we can assume all $H_n$ have
at {{least}} three copies of $\Z$.

Note that  the homomorphism
$$
\gm_{n,n+1}:~ H_n=\Z^{p_n}\lr H_{n+1}=\Z^{p_{n+1}}
$$ is given by a $p_{n+1}\times p_n,$ matrix $\cc=(c_{ij})$ of nonnegative integers,
where $i=1,2,...,p_{n+1},~ j=1,2,...,p_n,$ and $c_{ij}\in
\Z_+:=\{0,1,2,...\}.$ For $M>0$, if all $c_{ij}\geq M$, then we will
say $\gm_{n,n+1}$ is at least $M$-large or has multiplicity at least
$M$.  Note that since $H$ is a simple ordered group, {{passing}} to a subsequence, we
{may assume  that} at each step $\gm_{n,n+1}$ is at least $M_n$-large for
{{an}} arbitrary choice of $M_n$ depending on our construction up to step
$n$.

\end{NN}

\begin{NN}\label{range 0.5}
Recall that with the dimension group $H$ as in \ref{range  0.3}  and \ref{range  0.4}, we have $G\subset H$ with
$G_+=H_+\cap G,$ and both $G$ and $H$ share the same order unit $u\in G\subset H$.  As in \ref{range 0.4}, write $H$ as {{the}} inductive limit of
$H_n$---{{a sequence of }} finite direct sum{{s}} of three or more copies of the ordered group $(\Z, \Z_+)$
(with connecting maps with large multiplicities).  Let $G_n={\gm_{n, \infty}^{-1}(\gm_{n, \infty}(H_n)\cap G)},$
 where $\gm_{n,\infty}:~ H_n\to H$ is {{the}}  canonical
 map to the limit.  {{There is a $k_0\in \N$ such that $u\in G_{k_0}$,
 and without {\blue{loss}} of generality,  we {\blue{may}}  assume that $k_0=1.$  {\blue{In other words, we}}}}
{{may}}
assume $u\in G_n\subset  H_n$ for each $n$. {{Since $G_+=H_+\cap G$, if we
{{confer an}}  order structure on $G_n$ by setting}} $(G_n)_+{{:}}=(H_n)_+\cap G_n$, {{then we have $G_+=\lim_{n\to \infty} (G_n)_+$.}}

{{We claim that $G_n$ is a relatively divisible subgroup of $H_n$. Let us suppose that $g\in G_n$ and $g=mh$ for $h\in H_n$ {\blue{and for some integer $m\ge 1.$}} Since $G$ is relatively divisible in $H$, there is a $g'\in G$ such that $\gamma_{n, \infty }(g)=m g'$. Noting that $\gamma_{n, \infty}(g)=m\gamma_{n, \infty}(h)$, we have $m(g'-\gamma_{n, \infty}(h))=0$. Hence $\gamma_{n, \infty}(h)=g'\in G$ which implies that $h\in \gamma_{n, \infty}^{-1}(\gamma_{n, \infty}(H_n)\cap G)$. That is, $h\in G_n$.  }}
Since $G_{{n}}$ is a relatively divisible subgroup of $H_{{n}}$ and $H_{{n}}$ is torsion free, the quotient
$H_n/G_n$ is {{a}}  torsion-free finitely generated abelian group and {{therefore}} a direct sum of copies
of $\Z$, denoted by $\Z^{l_n}$. Then we have the following commutative diagram:
\begin{displaymath}
    \xymatrix{G_{1_{\,}} \ar[r]^{\gm_{_{12}}|_{_{G_{1}}}} \ar@{^{(}->}[d] & G_{2_{\,}} \ar[r]\ar@{^{(}->}[d]&\cd \ar[r]&G_{\,} \ar@{^{(}->}[d]\\
         H_1 \ar[r]^{\gm_{_{12}}} \ar[d] & H_2 \ar[r]\ar[d]&\cd \ar[r]&H \ar[d]\\
         H_1/G_1 \ar[r]^{\td\gm_{_{12}}}  & H_2/G_2 \ar[r]&\cd \ar[r]&H/G.}
\end{displaymath}
Let $H_n=\left(\Z^{p_n}, \Z^{p_n}_+, u_n\right)$, where
$u_n{{:=}}([n,1],[n,2], ..., [n,p_n])\in (\Z_+\setminus \{0\})^{p_n}$.
Then $H_n$ can be realized as the $K_0$-group of
$F_n=\bigoplus_{i=1}^{p_n}M_{[n,i]}$---that is,
$$
(K_0(F_n),K_0(F_n)_+,[1_{F_n}])=(H_n,(H_n)_+,u_n).
$$

If there are infinitely many $n$ such that the inclusion maps
$G_n\rightarrow H_n$
are isomorphisms, then,  passing to a subsequence, we have that
${{G\rightarrow H}}$
is also an isomorphism which contradicts $G\subsetneqq H$ in \ref{range 0.3}. Hence, without loss of generality, we may assume that for all $n$, $G_n\subsetneqq H_n$, and therefore $H_n/G_n\not= 0$.

To construct a \CA\, with $K_0$ equal to  $(G_n,(G_n)_+,u_n)$, we
consider the {\blue{quotient}} map $\pi:~ H_n\to H_n/G_n$
as a map
(still denoted in the same way)
$$\pi:~ \Z^{p_n}\lr \Z^{l_n},$$
as  in \cite{point-line}. We emphasize that $l_n>0$ for all $n$ {{(as $H_n/G_n\not=\{0\}$).}}
Such a map can be realized as difference of two
maps
$$\bb_0,~ \bb_1:~ \Z^{p_n}\lr \Z^{l_n}$$
corresponding to two $l_n\times p_n$ matrices of strictly positive
integers $\bb_0=(b_{0,ij})$ and $\bb_1=(b_{1,ij})$. That is,
$$
\qq\qq\pi\left(
     \begin{array}{c}
       t_1 \\
       t_2 \\
       \vdots \\
       {{t_{p_n}}} \\
     \end{array}
   \right)
   =
   \big(\bb_1-\bb_0\big)
   \left(
     \begin{array}{c}
       t_1 \\
       t_2 \\
       \vdots \\
       {{t_{p_n}}} \\
     \end{array}
   \right)\in \Z^{l_n},\qq
\mbox{for any} \quad \left(
     \begin{array}{c}
       t_1 \\
       t_2 \\
       \vdots \\
       t_{{p_n}} \\
     \end{array}
   \right)\in \Z^{p_n}.
$$
Note that $u_n=([n,1],[n,2], ..., [n,p_n])\in G_n$ and hence
$\pi(u_n)=0$.  Consequently,
$$
 \bb_1
\left(
     \begin{array}{c}
       {[n,1]} \\
      {[n,2]} \\
       \vdots \\
       {[n,p_n]} \\
     \end{array}
   \right)
   =
\bb_0 \left(
     \begin{array}{c}
       [n,1] \\
       {[n,2]} \\
       \vdots \\
       {[n,p_n]} \\
     \end{array}
   \right)
{{=:}} \left(
     \begin{array}{c}
       \{n,1\} \\
       \{n,2\}\\
       \vdots \\
       \{n,p_n\} \\
     \end{array}
   \right),
$$
i.e.,
$$
\{n,i\}{{:=}}\sum_{j=1}^{p_n} b_{1,ij}[n,j]=\sum_{j=1}^{p_n}
b_{0,ij}[n,j].
$$

Let $E_n=\bigoplus_{i=1}^{l_n} M_{\{n,i\}}$.  Choose any two
homomorphisms $\bt_0,\bt_1:~ F_n\to E_n$ such that $(\bt_0)_{*0}=\bb_0$
and $(\bt_1)_{*0}=\bb_1$.  Then define
$$C_n:=\Big\{(f,a)\in C([0,1],E_n)\oplus F_n;~ f(0)=\bt_0(a), f(1)=\bt_1(a)\Big\}{{=C([0,1],E_n)\oplus_{\bt_0,\bt_1}F_n}},$$
which is $A(F_n, E_n, \bt_0,\bt_1)$ as in the definition of \ref{DfC1}.
{{By Proposition \ref{2Lg13}}} and the fact that  the map $\pi{{=({{\bt_1}})_{*0}-(\bt_0)_{*0}}}$ {{(playing  the role of ${\phi_1}_{*0}-{\phi_0}_{*0}$ there)}}  is surjective {{(see definition of $\pi$, $\bb_0,$ and $\bb_1$ above)}},
we have $K_1(C_n)=0$ and
\beq\label{1311k0}
\Big(K_0(C_n),K_0(C_n)_+, \e_{C_n} \Big)=\Big({G_n},(G_n)_+,u_n\Big).
\eneq
{\blue{As in (\ref{sixterm}), the map}}
$K_0(F_n)=\Z^{p_n} \to
{\blue{K_0(E_n)=}}K_1(C_0\big((0,1),E_n\big))=\Z^{l_n}$ is given by $\bb_1-\bb_0\in
M_{l_n\times p_n}(\Z)$, which is surjective, as {{$\pi$ is the}} quotient map $H_n
\,(=\Z^{p_n})\to H_n/G_n\, (=\Z^{l_n})$.
{\blue{{\it In particular, all $C_n\in {\cal C}_0.$}}}

As observed in \cite{point-line}, in the construction of $C_n$ {{with \eqref{1311k0},}}
we have the
freedom to choose the pair of the $K_0$-maps $(\bt_0)_{*0}=\bb_0$ and
$(\bt_1)_{*0}=\bb_1$ as long as the difference is the same  map
$\pi:~ H_n\, (=\Z^{p_n})\lr H_n/G_n\, (=\Z^{l_n})$.
For example, if
$(m_{ij})\in M_{l_n\times p_n}\big(\Z_+\setminus \{0\}\big)$ is any
$l_n\times p_n$ matrix of positive integers, then we can replace
$b_{0,ij}$ by $b_{0,ij}+m_{ij}$ and, at the same time, replace
$b_{1,ij}$ by $b_{1,ij}+m_{ij}$. That is, we can assume that each
entry of $\bb_0$ (and of $\bb_1$) is larger than any fixed integer
$M$ which {{may depend}} on $C_{n-1}$ and $\psi_{n-1,n}: F_{n-1}\to F_n$.
Also, we can make all the entries of one column (say, the third
column) of both $\bb_0$ and $\bb_1$  much larger than all the
entries of another column
(say, the second), by choosing  $m_{i3} \gg m_{j2}$ for all $i,~ j$.  \\

\end{NN}

\begin{NN}\label{range 0.5a}

{\blue{Let us consider a slightly more general case than in \ref{range 0.5}. Let}}
$G_n=\Z^{{p_n^0}}\oplus G_n'$ and $H_n=\Z^{{p_n^0}}\oplus H_n'$,
with the inclusion map being the identity for the first ${{p_n^0}}$ copies
of $\Z$. In this case, the quotient map $H_n\, (=\Z^{p_n}) \to
H_n/G_n\, (=\Z^{l_n})$ given by the matrix $\bb_1-\bb_0$ maps the first
${{p_n^0}}$ copies of $\Z$ to zero. For this case, it will be much
more convenient to assume that the first ${{p_n^0}}$ columns of  both the
matrices $\bb_0$ and $\bb_1$ are zero and each entry of the last
${{p_n^1:=}}~ p_n-{{p_n^0}}$ {\blue{columns}} of them are larger than any previously  given integer
$M$. Now we have that the entries of the matrices $\bb_0,\bb_1$ are strictly positive integers except the ones in the first ${{p_n^0}}$ columns which are zero.

Consider the following diagram:
\begin{displaymath}
    \xymatrix{\Z^{p_1^0}\oplus G'_1 \ar[r]^{\gm_{_{12}}|_{_{G_1}}} \ar@{^{(}->}[d] & \Z^{p_2^0}\oplus G'_2 \ar[r]\ar@{^{(}->}[d]&\cd \ar[r]&G_{\,}\, \ar@{^{(}->}[d]\\
       \Z^{p_1^0}\oplus  H'_1 \ar[r]^{\gm_{_{12}}} \ar[d] & \Z^{p_2^0}\oplus H'_2 \ar[r]\ar[d]&\cd \ar[r]&H \ar[d]\\
         H_1/G_1 \ar[r]^{\td\gm_{_{12}}}  & H_2/G_2 \ar[r]&\cd \ar[r]&H/G\,.}
\end{displaymath}
That is, $G_n=\Z^{p_n^0}\oplus G_n'$ and $H_n=\Z^{p_n^0}\oplus H_n'$,
with the inclusion map being  the identity for the first $p_n^0$ copies
of $\Z$.
We are now assuming that the entries of the matrices $\bb_0,\bb_1$ are strictly positive integers for the last $p_n^1$ {{c}}olumns and are zeros for  the first $p_n^0$ {{c}}olumns. Note that since $l_n>0$ (see \ref{range 0.5}), we have $p_n^1>0$.

{ {{The inductive limits $H=\varinjlim(H_n,\gm_{n,n+1})$ and $G=\varinjlim(G_n,\gm_{n,n+1}|_{G_n})$ (with $G_n\subset H_n$)  constructed in \ref{range 0.5} are in fact special cases of the {{present}} construction
when we assume that $p_n^0=0.$} }} One {notices} that, for the case $G_n=\Z^{p_n^0}\oplus G_n'$ and $H_n=\Z^{p_n^0}\oplus H_n'$,
with the inclusion map being {{the}} identity for the first $p_n^0$ copies
of $\Z$, one {{could}} still use the construction of \ref{range 0.5} to make all the entries of $\bb_0$ and $\bb_1$ (not only the entries of the last $p_n^1$) be strictly positive (of course with the first $p_n^0$ {{columns of the}} matrices $\bb_0$ and $\bb_1$ equal to each {{other}}).
{\blue{However,}}
for the {{algebras with the property (SP), it is possible to assume}}  that $p_n^0\not= 0$ for all $n$, {{and to get a certain}} decomposition {{property}} that we will discuss in the next section. In fact, for the case that $p_n^0\not= 0$ for all $n$, the construction is much simpler than the case that $p_n^0=0$ for all $n$.

If we write $F_n$ {{above}} as
$\bigoplus_{i=1}^{p_n^0}M_{[n,i]}\oplus F_n'$, where
$F_n'=\bigoplus_{i=p_n^0 +1}^{p_n}M_{[n,i]},$ then the maps
$\bt_0$ and $\bt_1$ are zero on the part
$\bigoplus_{i=1}^{p_n^0}M_{[n,i]}$. Moreover, the algebra
$C_n=\Big\{(f,a)\in C([0,1],E_n)\oplus F_n;~ f(0)=\bt_0(a),
f(1)=\bt_1(a)\Big\}$ {{can}} be written as
$\bigoplus_{i=1}^{p_n^0}M_{[n,i]} \oplus C'_{\blue{n}}$, where
$$C'_{\blue{n}}=\Big\{(f,a)\in C([0,1],E_n)\oplus F_n';~ f(0)=\bt_0(a),
f(1)=\bt_1(a)\Big\}{\blue{=C([0,1],E_n)\oplus_{\bt_0|_{F'_n},\bt_1|_{F'_n}}F'_n}}$$ as in {\blue{\ref{range 0.5}. }}
{\blue{It should be remembered
 that \ref{range 0.5} is a special case of \ref{range 0.5a}.}}

\end{NN}

\begin{NN}\label{conditions}
Let us {{emphasize}}  that once
$$(H_n,(H_n)_+,u_n)=\big(\Z^{p_n}, (\Z_+)^{p_n},
([n,1],[n,2], ..., [n,p_n])\big),\andeqn \bb_0,\bb_1: \Z^{p_n}
\to
\Z^{l_n},$$
are fixed, then the algebras $F_n=\bigoplus_{i=1}^{p_n} M_{[n,i]}$, $E_n=\bigoplus_{i=1}^{l_n} M_{\{n,i\}}$ {{are fixed, }} where \\
$\{n,i\}{{:=}}\sum_{j=1}^{p_n} b_{1,ij}[n,j]=\sum_{j=1}^{p_n}
b_{0,ij}[n,j]${{;}} and, {{by Proposition \ref{2Pg12}, the algebras}} $C_n$ \linebreak $=A(F_n,E_n,\bt_0,\bt_1)$, with
${{(\bt_i)_{*0}}}=\bb_i$ ($i=0,1$), are determined  up to isomorphism.

To construct the inductive limit, we not only need to construct $C_n$'s (later on $A_n$'s for the  general case) but also
need to construct $\phi_{n,n+1}$   which realizes the corresponding K-theory map---i.e., $(\phi_{n,n+1})_{*0}=\gm_{_{n,n+1}}|_{G_n}$. In {{addition}}, we need to make the limit algebras {{have}}  the desired tracial state space.
In order to do all these {{things}}, we need some extra conditions on the maps $\bb_0, \bb_1$ for $C_{n+1}$ (or for $A_{n+1}$)  depending on $C_n$ (or $A_n$).
We will divide the construction into several steps with gradually stronger conditions on $\bb_0, \bb_1$ (for $C_{n+1}$)---of course depending on $C_n$ and the map $\gm_{n,n+1}: H_n\to H_{n+1}$,  to guarantee the construction can go through.

Let $G_n \subset H_n=\Z^{p_n}$,  $G_{n+1} \subset H_{n+1}=\Z^{p_{n+1}}$, and $\gm_{_{n,n+1}}:~ H_n\to H_{n+1}$, with $\gm_{_{n,n+1}} (G_n) \subset G_{n+1},$ be as in \ref{range 0.4} and \ref{range 0.5} (also see \ref{range 0.5a}). Then $\gm_{_{n,n+1}}$ induces a map $\td\gm_{_{n,n+1}}: H_n/G_n\,(=\Z^{l_n}) \to H_{n+1}/G_{n+1}= (\Z^{l_{n+1}})$.
Let $\gm_{_{n,n+1}}:~H_n\left(=\Z^{p_n}\right)\to
H_{n+1}\left(=\Z^{p_{n+1}}\right)$ be given by the matrix
$\cc=(c_{ij})\in M_{_{p_{n+1}\times p_n}}(\Z_+\setminus \{0\})$ and
$\td\gm_{_{n,n+1}}:~\Z^{l_n} \to \Z^{l_{n+1}}$ (as a map from
$H_n/G_n \to H_{n+1}/G_{n+1}$) be given by the matrix $\dd=(d_{ij})$. Let $\pi_{n+1}: H_{n+1} (=\Z^{p_{n+1}}) \to H_{n+1}/G_{n+1}(=\Z^{l_{n+1}})$ be the quotient map.

Let us use $\bb'_0, \bb'_1:\Z^{p_{n+1}} \to \Z^{l_{n+1}}$ to denote the maps
%
required for
the construction of $C_{n+1}$, and reserve $\bb_0, \bb_1$ for $C_n$.  Of course, $\pi_{n+1}=\bb'_1-\bb'_0$.

We will prove that if $\bb'_0, \bb'_1$ satisfy:
\beq\label{13spd-1}
\spd \qq\qq\qq\qq  \td b_{0,ji},~\td b_{1,ji}~ > ~
\sum_{k=1}^{l_n}\big(|d_{jk}|+2\big)\max(b_{0,ki},~b_{1,ki})
\eneq
for all
$i\in\{1,2,...,p_n\}$ and for all  $j\in\{1,2,...,l_{n+1}\},$
where {{$\td b_{0,ji}$ and $\td b_{1,ji}$ are the entries of}} $\td \bb_0{{:=}}\bb'_0\cdot \cc$ and $\td
\bb_1{{:=}}\bb'_1\cdot \cc$, {{respectively,}} then one can construct the homomorphism $\phi_{n,n+1}: C_n \to C_{n+1}$ to {{realize}} the desired $K$-theory map (see \ref{range 0.6} below). If $\bb'_0, \bb'_1$ satisfy the stronger condition
 \beq\label{13spd-2}
\spdd\qq\qq\qq\qq \td b_{0,ji}, \td
 b_{1,ji}
 > 2^{2n}\!\left(\sum_{k=1}^{l_n}
 (|d_{jk}|+2)\{n,k\}\!\!\right)
 \eneq
 for all
$i\in\{1,2,...,p_n\}$ and for all  $j\in\{1,2,...,l_{n+1}\}$, then we can prove the limit algebra constructed has the desired tracial state space (see {{the corresponding calculation in}} \ref{range 0.8}--{{\ref{range 0.10}, which will be used in the proof of Theorem \ref{EllofA}}}). (It follows from {{the fact}} that $\{n,k\}{{:=}}\sum_{i=1}^{p_n} b_{1,ki}[n,i]=\sum_{i=1}^{p_n}b_{0,ki}[n,i]$
{{that}}
the inequality $\spdd$ is stronger than $\spd$.)

{{From the definition of $\td \bb_0$ and $\td \bb_1$, we have
\beq\label{13Dec9-2018}
\td b_{0,ji}=\sum_{k=0}^{p_{n+1}}b_{0,jk}'c_{ki}~~~~~~\mbox{and}~~~~~\td b_{1,ji}=\sum_{k=0}^{p_{n+1}}b_{1,jk}'c_{ki},
\eneq
where $b_{0,jk}'$ and $b_{1,jk}'$ are the entries of $\bb_0'$ $\bb_1'$, respectively. }}

Let us {\blue{also}} emphasi{{ze}} that when we modify inductive limit system to make it simple, we never change the algebras $C_n$ (or $A_n$ in the general case), what will be changed are the connecting homomorphisms.


\end{NN}

{\it Let $A$ and $B$ be  \CA s,  $\phi: A\to B$ be a \hm, and $\pi\in RF(B).$
{For the rest of} this  section, we will use $\phi|_\pi$ for the composition $\pi\circ \phi,$ in particular,
in the following statement and its proof. This notation is consistent with  {\rm \ref{homrestr}}.}

{In the following lemma we will give the construction of  $\phi_{n,n+1}$
{{in the case that}} $\bb'_0$ and $\bb'_1$ satisfy Condition $\spd$ {{of}} \ref{conditions}---of course, the condition depends on the previous  step.
So this lemma {provides the}
{{$(n+1)$st}} step of the  construction. Again, we first have $G,$ and then {{obtain}} $H,$
$H_n,$ and $G_n$ as constructed in \ref{range 0.3} and \ref{range 0.4}. }

\begin{lem}\label{range 0.6}
Let
\beq\nonumber
(H_n,(H_n)_+,u_n)=\big(\Z^{p_n}, (\Z_+)^{p_n},
([n,1],[n,2], ... , [n,p_n])\big),\\
\nonumber
F_n=\bigoplus_{i=1}^{p_n} M_{[n,i]},\,\,\, \bb_0,\bb_1: \Z^{p_n} \to
\Z^{l_n},\,\,\,
E_n=\bigoplus_{i=1}^{l_n} M_{\{n,i\}},~ \bt_0,\bt_1:
F_n \to E_n
\eneq
 with $(\bt_0)_{*0}={\blue{\bb_0=(b_{0,ij})_{l_n\times p_n}}}, (\bt_1)_{*0}=\bb_1={\blue{(b_{1,ij})_{l_n\times p_n}}}$, and
$C_n=A(F_n,E_n,\bt_0,\bt_1)$ with $K_0(C_n)=G_n$ be as in {\blue{{\rm \ref{range 0.5a}}, which includes the case of
{\rm \ref{range 0.5}} as the special case  $p_n^0=0$.}}
  Let
$$
(H_{n+1},H_{n+1}^+,u_{n+1})=\big(\Z^{p_{n+1}}, (\Z_+)^{p_{n+1}},
([{n+1},1],[{n+1},2], ..., [{n+1},p_{n+1}])\big),
$$
let  $\gm_{n,n+1}:~
H_n\to H_{n+1}$ be { an} {{order}} homomorphism with
${{\gm}}_{n,n+1}(u_n)=u_{n+1}$ {\rm (}as in {\rm \ref{range  0.5}} or {\rm \ref{range  0.5a}} {\rm )}, and let
$G_{n+1}\sbs H_{n+1}$ be a subgroup containing $u_{n+1}$ {\rm (}as in
{\rm \ref{range 0.4}}{\rm )} and satisfying $\gm_{n,n+1}(G_n)\sbs G_{n+1}$.  Let
$\pi_{n+1}: H_{n+1}~ (=\Z^{p_{n+1}}) \to H_{n+1}/G_{n+1} ~(=\Z^{l_{n+1}})$ denote the quotient map,
{\blue{and let $\gamma_{n,n+1}$ be represented by the $p_{n+1}\times p_n$-matrix $(c_{ij}).$}}

Suppose that the {{maps}}
$\bb'_0=(b'_{0,ij}),~\bb'_1=(b'_{1,ij}): \Z^{p_{n+1}} \to \Z^{l_{n+1}}$  satisfy $\bb'_1-\bb'_0=\pi_{n+1}$ and satisfy Condition $\spd$ of {\rm \ref{conditions}}. (As a convention, we assume that the entries of the first $p_{n+1}^0$ columns of the matrices $\bb'_0$ and $\bb'_1$ are zeros, and the entries of the last $p_{n+1}^1=p_{n+1}-p_{n+1}^0$ columns are strictly positive. Note that $p_{n+1}^0$ might be zero as in the special case {\rm \ref{range  0.5}}.) 

Put
$F_{n+1}=\bigoplus_{i=1}^{p_{n+1}} M_{[{n+1},i]}$ and
$E_{n+1}=\bigoplus_{i=1}^{l_{n+1}} M_{\{{n+1},i\}}$ with $$\{n+1,i\}=\sum_{j=1}^{p_{n+1}} b'_{1,ij}[n+1,j]=\sum_{j=1}^{p_{n+1}}
b'_{0,ij}[n+1,j],$$
and pick unital homomorphisms $\bt_0',\bt_1': F_{n+1}\to E_{n+1}$ with $$(\bt_0')_{*0}=\bb'_0\quad \textrm{and}\quad (\bt_1')_{*0}=\bb'_1.$$
Set
$$C_{n+1}=A(F_{n+1},E_{n+1},\bt_0',\bt_1').$$

Then there is a  \hm~ $\phi_{n, n+1}:~
C_n\to C_{n+1}$ satisfying the following conditions:
\begin{enumerate}
\item[\rm{(1)}] $K_0(C_{n+1})=G_{n+1}$ as {{ scaled ordered groups}} (as already verified in \ref{range 0.5}).
\item[\rm{(2)}] $(\phi_{n, n+1})_{*0}:~ K_0(C_n)=G_n \to K_0(C_{n+1})= G_{n+1}$
satisfies $(\phi_{n, n+1})_{*0}=\gm_{n, n+1}|_{G_n}$.
\item[\rm{(3)}] $\phi_{n, n+1}(C_0\big((0,1), E_n\big))\sbs C_0\big((0,1),
E_{n+1}\big)$.
\item[\rm{(4)}] Let $\td\phi_{n, n+1}:~ F_n\to F_{n+1}$ be the quotient map
induced by $\phi_{n, n+1}$ (note from (3), we know that this
quotient map exists); then
 $(\td\phi_{n, n+1})_{*0}=\gm_{n, n+1}:~ K_0(F_n)=H_n \to K_0(F_{n+1})=H_{n+1}$.
\item[\rm{(5)}] For each $y\in Sp(C_{n+1})$, $Sp(F_{n})\subset Sp(\phi_{n,n+1}|_y)$.
\item[\rm{(6)}] For each {{pair}} $j_0\in \{ 1, 2, ..., l_{n+1}\}, ~i_0\in \{ 1, 2, ... ,
l_{n}\}$, one of the following {{properties}}  holds:
\begin{enumerate}
\item[{\rm (i)}]  for each $t\in (0,1)_{j_0}\sbs
Sp(I_{n+1})=\bigcup_{j=1}^{l_{n+1}}(0,1)_j\sbs Sp(C_{n+1})$,
$Sp(\phi_{n, n+1}|_t)\cap (0,1)_{i_0}$ contains $t\in (0,1)_{i_0}\sbs
Sp(C_n)$; or
\item[{\rm (ii)}] for each $t\in (0,1)_{j_0}\sbs
Sp(I_{n+1})=\bigcup_{j=1}^{l_{n+1}}(0,1)_j\sbs Sp(C_{n+1})$,
$Sp(\phi_{n, n+1}|_t)\cap (0,1)_{i_0}$ contains $1-t\in
(0,1)_{i_0}\sbs Sp(C_n)$.
\end{enumerate}

\item[\rm{(7)}] {\blue{The map $\phi_{n, n+1}$ is injective.}}

\item[\rm{(8)}]  {\blue{If $X\sbs Sp(C_{n+1})$ is $\dt$-dense,
then $\bigcup_{x\in X}Sp(\phi_{n, n+1}|_x)$ is $\dt$-dense in $Sp(C_n)$ (see \ref{density}).}}
\end{enumerate}


\end{lem}


\begin{rem}\label{range 0.7}
Let $I_n=C_0\big((0,1),E_n\big)$  and
$I_{n+1}=C_0\big((0,1),E_{n+1}\big)$.  If $\phi_{n, n+1}:~ C_n\to
C_{n+1}$
 is as described in \ref{range  0.6}, then we have the following {{map}} between  exact sequences:
\begin{displaymath}
\xymatrix{0\ar[r] &
K_0(C_n)\ar[r]^{~}\ar@{->}[d]_{\gm_{_{n,n+1}}|_{{_{G_n}}}} &
K_0(C_n/I_n) \ar@{->}[d]_{\gm_{_{n,n+1}}} \ar[r]^{~} & K_1(I_n)
\ar@{->}[d]_{\td\gm_{_{n,n+1}}} \ar[r]&0
\\
0\ar[r] & K_0(C_{n+1})\ar[r]^{~}& K_0(C_{n+1}/I_{n+1})
 \ar[r]^{~} & K_1(I_{n+1}) \ar[r] & 0,}
\end{displaymath}
where $K_0(C_n)$ and $K_0(C_{n+1})$ are identified with $G_n$ and
$G_{n+1}$, $K_0(C_n/I_n)\,(=K_0(F_n))$ and
$K_0(C_{n+1}/I_{n+1})\,(=K_0(F_{n+1}))$ are identified with $H_n$ and
$H_{n+1}$, $K_1(I_n)$ is identified with $H_n/G_n$, and
$K_1(I_{n+1})$ is identified with $H_{n+1}/G_{n+1}$. Moreover,
$\td\gm_{_{n,n+1}}$ is induced by $\gm_{_{n,n+1}}:~ H_n\to H_{n+1}$.

Consider the matrix
$\cc=(c_{ij})_{p_{n+1}\times p_n}$ with $c_{ij}\in \Z\setminus \{0\}$
which is induced by the map
$\gm_{_{n,n+1}}: H_n\left(=\Z^{p_n}\right) \to H_{n+1}\left(=\Z^{p_{n+1}}\right)$
and consider the matrix $\dd=(d_{ij})$ which is induced by the map
$\td\gm_{_{n,n+1}}:~\Z^{l_n} \to \Z^{l_{n+1}}$ (as a {{map}}
$H_n/G_n \to H_{n+1}/G_{n+1}$).
Note that $\pi_n:~ H_n\,(=\Z^{p_n})\to  H_n/G_n~(=\Z^{l_n})$ is given
by $\bb_1-\bb_0,$ a $l_n\times p_n$-matrix with entries in $\Z.$
{Here,  in the situation of \ref{range 0.5a},
 we assume the first
$p_n^0$ columns of both $\bb_0$ and $\bb_1$ are zero
and the last $p_n^1$ columns are  strictly positive.}
Let $\pi_{n+1}:~ H_{n+1}\, (=\Z^{p_{n+1}})
\to H_{n+1}/G_{n+1}\, (=\Z^{l_{n+1}})$ be the quotient map.
Write $G_{n+1}=\Z^{p_{n+1}^0}\oplus G_{n+1}'$,
$H_{n+1}=\Z^{p_{n+1}^0}\oplus H_{n+1}'.$ Then we can choose  both $\bb_0'$ and $\bb_1'$, with the first $p_{n+1}^0$ columns  zero and the last $p_{n+1}^1=p_{n+1}-p_{n+1}^0$ columns  strictly positive, so that $\pi_{n+1}=\bb'_1-\bb'_0$ and
Condition $\spd$ is satisfied, {\blue{i.e.,}}
\beq\nonumber
 {{\td b_{0,ji}}},~\td b_{1,ji}~ > ~
\sum_{k=1}^{l_n}\big(|d_{jk}|+2\big)\max(b_{0,ki},~b_{1,ki})
\eneq
for all
$i\in\{1,2,...,p_n\}$ and for all  $j\in\{1,2,...,l_{n+1}\},$
where $\td \bb_0=\bb'_0\cdot \cc=(\td b_{0,ji})$ and $\td
\bb_1=\bb'_1\cdot \cc=(\td b_{1,ji})$.


Indeed, note that  $l_{n+1}>0$ and  $p_{n+1}^1>0$. So we can make $\spd$ hold by only increasing the last
$p_{n+1}^1$ columns of the the matrices $\bb_0'$ and
$\bb_1'$---that is, the first $p_{n+1}^0$ column{{s}} of the matrices
are still kept to be  zero, since all the entries in $\cc$ are strictly positive.  Note that, even though the first $p_{n+1}^0$ columns of $\bb_0'$ and
$\bb_1'$ (as $l_{n+1}\times p_{n+1}$ matrices) are zero,  all entries of $\td \bb_0$ and $\td \bb_1$ (as $l_{n+1}\times p_n$ {{matrices}}) {\blue{have}} been
made strictly positive.
Again note that  the case of \ref{range 0.5} is the special case of \ref{range 0.5a} for $p_{n+1}^0=0$, so one does not need to deal with this case separately.

 \end{rem}

  \begin{proof}[Proof of \ref{range 0.6}]
 Suppose that
 $\bb_0'$ and $\bb_1'$  satisfy $\bb'_1-\bb'_0=\pi_{n+1}$ and
 (\ref{13spd-1}) (and the first $p_{n+1}^0$ columns {{are}} zero).  Now let $E_{n+1},~\bt'_0,\bt'_1:~ F_{n+1}\to
 E_{n+1}$, and $C_{n+1}=A(F_{n+1},E_{n+1},\bt'_0,\bt'_1)$ be as
 constructed in \ref{range 0.5}.
  We will define $\phi_{n,n+1}:~C_n\to C_{n+1}$
to satisfy the conditions (2)--{{(8)}}  of \ref{range  0.6} (the condition (1) is a property of $C_{n+1}$ which is verified in \ref{range 0.5}).

As usual, let us use $F_n^i$ (or $E_n^i$) to denote the $i$-th block of $F_n$ (or $E_n$).

There exists a unital \hm\,
$\td\phi_{n, n+1}:~ F_n\to F_{n+1}$  such that
\beq\label{134-n1}
(\td\phi_{n,n+1})_{*0}=
\gm_{n, n+1}:~ K_0(F_n)=H_n \to K_0(F_{n+1})=H_{n+1},
\eneq
where
$\gm_{n,n+1}$ is {{as described in the hypotheses of {\blue{Lemma}} \ref{range 0.6} (see also}} \ref{range 0.4} {{ and \ref{conditions})}}.
Note that
$Sp(C_{n+1})=\bigsqcup_{j=1}^{l_{n+1}}(0,1)_j\cup Sp(F_{n+1})$ (see
\S 3). Write $Sp(F_{n+1})=(\tht'_1,\tht_2',...,\tht'_{p_{n+1}})$ {and}
 $Sp(F_{n})=(\tht_1,\tht_2, ... ,\tht_{p_{n}})$. To define
$\phi_{n,n+1}:~C_n\to C_{n+1}$,  we need to specify {{each map}} $\phi_{n,n+1}|_y=e_y\circ \phi_{n,n+1},$
{{i.e.,}}  for each
$y\in Sp(C_{n+1})$, the composed map
\begin{displaymath}
\xymatrix{\phi_{_{n,n+1}}|_y:~C_n\ar[r]&  C_{n+1}~ \ar[r]^{e_{_y}} &
C_{n+1}|_y}~,
\end{displaymath}
{{with}} $e_y$  the point evaluation at $y.$
%
%

{To actually construct $\phi_{n, n+1},$ we  first construct
a} homomorphism
$\psi:~ C([0,1],E_n)\to
C([0,1],E_{n+1})$.
This can be done by
{\blue{defining the}} map
\beq\label{13Dec21-2018-1}
\psi^j:~C([0,1],E_n)\to C([0,1],E^j_{n+1}),\,\,\mbox{for
each }~ j=1,2,..., l_{n+1},\eneq
{\blue{as follows.}}
Let $(f_1,f_2,..., f_{l_n})\in C([0,1],E_n)$.  For any $k\in
\{1,2,..., l_{n}\}$, if $d_{jk}>0$, then let
\beq\label{13Jan16-2019}
F_k(t)=\diag\big(\underbrace{f_k(t),f_k(t),...,
f_{k}(t)}_{d_{jk}}\big) \in C\left([0,1],M_{d_{jk}\cdot
\{n,k\}}\right);\eneq
 if $d_{jk}<0$, then let
\beq\label{13Jan16-2019-1}F_k(t)=\diag\big(\underbrace{f_k(1-t),f_k(1-t),...,
f_{k}(1-t)}_{|d_{jk}|}\big) \in C\left([0,1],M_{|d_{jk}|\cdot
\{n,k\}}\right);\eneq and if $d_{jk}=0$, then let
\beq\label{13Jan16-2019-2}
F_k(t)=\diag\big(f_k(t),f_k(1-t)\big) \in C\left([0,1],M_{2\cdot
\{n,k\}}\right).
\eneq
With the above notation, define
\beq\label{defpsi}
\hspace{-0.3in}\psi^j(f_1,f_2,..., f_{l_n})(t)\!=\!\diag(F_1(t),F_2(t),..., F_{l_n}(t))
\in C\!\left([0,1],M_{\left({\sum_{k=1}^{l_n}d_k'}\cdot
\{n,k\}\right)}\right)\!,~~
\eneq
where
$$
d_k'=\left\{\begin{array}{cc}
                |d_{jk}| & \mbox{ if } d_{jk}\not= 0, \\
                ~&\\
                2 & \mbox{ if } d_{jk} = 0.
              \end{array}\right.
$$
Note that
\beq\label{sizeE}
\{n+1,j\}= \sum_{l=1}^{p_{n+1}}b_{0,jl}'[n+1,l]
=\sum_{l=1}^{p_{n+1}}\sum_{i=1}^{p_{n}}b_{0,jl}'c_{li}[n,i]
=\sum_{i=1}^{p_{n}}\td b_{0,ji}[n,i].
\eneq
Recall that $\td \bb_0=\bb_0'\cdot \cc=\big(\td b_{0,ji}\big)$. From
 (\ref{13spd-1}),  {{(\ref{sizeE})}}, $d_k'\leq |d_{kj}|+2,$ and $\{n,k\}=\sum
b_{0,ki}[n,i]$, we  {{deduce}}
$$
\{n+1,j\} =\sum_{i=1}^{p_{n}}\td b_{0,ji}[n,i]> \sum_{i=1}^{p_{n}} \big( \sum_{k=1}^{l_n}\big(|d_{jk}|+2\big)b_{0,ki}\big)[n,i]\geq \sum_{k=1}^{l_n}d_k'\{n,k\}.
$$
Hence the \CA\,
$C\left([0,1],M_{\left({\sum_{k=1}^{l_n}d_k'}\cdot
\{n,k\}\right)}\right)$ can be regarded as a corner of the \CA\,
$C([0,1],E_{n+1}^j)=C\left([0,1],M_{ \{n+1,j\}}\right)$, and
consequently, $\psi^j$ can be regarded
{\blue{as a map from   $C([0,1], E_n)$}}  into
$C([0,1],E_{n+1}^j)$.  {{Putting}} all $\psi^j$ together we get {{a map}} $\psi:~
C([0,1],E_n)\to C([0,1],E_{n+1})$ {\blue{defined by $\psi(f)=(\psi^1(f), \psi^2(f),...,\psi^{l_n}(f))$ for all $f\in C([0,1], E_n).$}}

Define $\psi_0,\psi_1:~ C_n\to E_{n+1}$ to be
$$
\psi_0(f)=\psi(f)(0) \andeqn \psi_1(f)=\psi(f)(1)
$$
for any $f\in C_n\sbs C([0,1],E_n)$.  Since
$\psi_0(C_0((0,1),E_n))=0$ and $\psi_1(C_0((0,1),E_n))=0$, this
defines maps $\af_0,~\af_1:~ F_n\to E_{n+1}$.
 Note that for each $j\in\{1,2, ... , l_{n+1}\}$, the maps
 $\af_0^j,~ \af_1^j :~F_n\to E_{n+1} \to E_{n+1}^j$ have spectra
$$
SP(\af_0^j) =\left\{\tht_1^{\sim i_1},\tht_2^{\sim i_2},...,
\tht_{p_n}^{\sim i_{p_n}}\right\} \andeqn SP(\af_1^j)
=\left\{\tht_1^{\sim i_1'},\tht_2^{\sim i_2'},..., \tht_{p_n}^{\sim
i_{p_n}'}\right\},
$$
respectively
(see \ref{ktimes} for the notation used here), where
$$i_l=\sum_{{\{k:d_{jk}<0\}}}|d_{jk}|b_{1,kl}+\sum_{{\{k:d_{jk}>0\}}}
|d_{jk}|b_{0,kl}
+\sum_{{\{k:d_{jk}=0\}}}(b_{0,kl}+b_{1,kl}),
$$
$$\qq\mbox{ and }\qq
i_l'=\sum_{{\{k:d_{jk}>0\}}}|d_{jk}|b_{1,kl}+\sum_{{\{k:d_{jk}<0\}}}|d_{jk}|b_{0,kl}
+\sum_{{\{k:d_{jk}=0\}}}(b_{0,kl}+b_{1,kl}).\qq\qq
$$
Note that $i_l'-i_l=\sum_{k=1}^{l_n}d_{jk}(b_{1,kl}-b_{0,kl})$, and note
that, if $l\leq p_n^0$ in the case \ref{range 0.5a}, then
$b_{0,kl}=b_{1,kl}=0,$ and consequently, $i_l=i_l'=0$. Put
\begin{displaymath}
\xymatrix{ \td\af_0=\bt_0'\circ \td\phi_{{n,n+1}}:~F_n~
\ar[r]^{~~~~~~~~~~~~ \td\phi_{_{n, n+1}}} & F_{n+1} \ar[r]^{\tiny
\bt'_{0}}& E_{n+1}~}\andeqn\\\\
\xymatrix{ \td\af_1=\bt_1'\circ \td\phi_{n,n+1}:~F_n~
\ar[r]^{~~~~~~~~~~~~ \td\phi_{n, n+1}} & F_{n+1} \ar[r]^{\tiny
\bt'_{1}}& E_{n+1}~}.
\end{displaymath}
{\blue{Then,}} for each $j\in\{1,2,...,l_{n+1}\} $, the maps
$\td\af_0^j,~\td\af_1^j:~  F_n\to E_{n+1}\to E_{n+1}^j$ have spectra
$$
SP(\td\af_0^j) =\left\{\tht_1^{\sim \bar i_1},\tht_2^{\sim \bar
i_2},..., \tht_{p_n}^{\sim \bar i_{p_n}}\right\} \andeqn
SP(\td\af_1^j) =\left\{\tht_1^{\sim \bar i_1'},\tht_2^{\sim \bar
i_2'},..., \tht_{p_n}^{\sim \bar i_{p_n}'}\right\},
$$
where
$$\bar i_l=\sum_{k=1}^{p_{n+1}}b'_{0,jk}c_{kl}=\td b_{0,jl}
\andeqn \bar i'_l=\sum_{k=1}^{p_{n+1}}b'_{1,jk}c_{kl}=\td
b_{1,jl}.
$$
From  (\ref{13spd-1}), we have that $\bar i_l> i_l$ and $\bar
i_l'> i_l'$. Furthermore, $\bar i_l'-\bar
i_l=\sum_{k=1}^{p_{n+1}}(b'_{1,jk}-b'_{0,jk})c_{kl}$. Since
$(\bb_1'-\bb_0')\cc=\dd(\bb_1-\bb_0)$, we have that $\bar i_l'-\bar i_l=
i_l'- i_l$, and hence $\bar i_l'-i_l'= \bar i_l- i_l{{:=}} r_l>0$.
Note that these numbers are defined for the homomorphisms $\af_0^j,
\af_1^j, \td\af_0^j, \td\af_1^j :~F_n\to  E_{n+1}^j$. So, strictly
speaking,  $r_j>0$ means
$r_l^j> 0$.

Define a unital homomorphism
 $\Phi:~ C_n\to C([0,1], E_{n+1})=\bigoplus_{j=1}^{l_{n+1}}C([0,1],E_{n+1}^j)$ by
\beq\label{13Nov30-2018}
\Phi^j(f,  {{( a_1,a_2,... ,a_{p_n})}})=
\diag\big(\psi^j(f),
~ a_1^{\sim r_1^j},a_2^{\sim
r_2^j},... ,a_{p_n}^{\sim r_{p_n}^j} \big)
\eneq
{\blue{for all $f\in C([0,1], E_n)$ and $(a_1,a_2,...,a_{p_n})\in F_n.$}}
Again, define the maps $\Phi_0,~ \Phi_1:~ C_n\to E_{n+1}$ by
$$\Phi_0(F)=\Phi(F)(0)\andeqn\Phi_1(F)=\Phi(F)(1),$$
for $F=(f_1,f_2,...,f_{l_n}{{;}}~a_1,a_2,...,a_{p_n})\in C_n$.  These two
maps induce two quotient maps
$$ \td{\td\af}_0,~\td{\td\af}_1,~:F_n\to E_{n+1},$$
{\blue {as $\Phi_0(I_n)=0$ and $\Phi_1(I_n)=0$.}}

From our calculation, for each $j\in\{1,2,...,l_{n+1}\} $, the map $
\td{\td\af}_0^j$ (resp. $\td{\td\af}_1^j$) has {{the}}  same
spectrum {\blue{(with multiplicities)}} as ${\td\af}_0^j$ (resp. ${\td\af}_1^j$)
does.  That is,
 $(\td{\td\af}_0^j)_{*0}=({\td\af}_0^j)_{*0}$ and $(\td{\td\af}_1^j)_{*0}=({\td\af}_1^j)_{*0}$.
 There are unitaries $U_0,~ U_1\in E_{n+1}$ such that $\mathrm{Ad}\,U_0\circ \td{\td\af}_0^j
 ={\td\af}_0^j$ and $\mathrm{Ad}\,U_1\circ \td{\td\af}_1^j
 ={\td\af}_1^j$.  Choose a unitary path $U\in C([0,1],E_{n+1})$ such that $U(0)=U_0$
 and $U(1)=U_1$.  Finally, set $\phi_{n, n+1}:~ C_n\to C([0,1],E_{n+1})$ to
 be defined as
{{\beq\label{13Dec10-2018}
\phi_{n, n+1}=\mathrm{Ad}\,U\circ \Phi.\eneq}}
From the construction, we have that $\psi(C_0\big((0,1), E_n\big))\sbs C_0\big((0,1),
E_{n+1}\big)$ and consequently, $\Phi(C_0\big((0,1), E_n\big))\sbs C_0\big((0,1),
E_{n+1}\big)${{, i.e.,}} $\phi_{n, n+1}(I_n) \subset I_{n+1}$.
{\blue{So (3) follows.}}

Since $\mathrm{Ad}\,U(0)\circ
\td{\td\af}_0^j
 ={\td\af}_0^j=\bt_0'\circ \td\phi_{n, n+1}$ and $\mathrm{Ad}\,U(1)\circ \td{\td\af}_1^j
 ={\td\af}_1^j=\bt_1'\circ\td\phi_{n, n+1}$ we have that $\phi_{n, n+1} (C_n)\subset C_{n+1}$ and furthermore,  the quotient map from
$C_n/I_n\to C_{n+1}/I_{n+1}$ induced by $\phi_{n, n+1}$ is the
same as $\td\phi_{n, n+1}$ (see definition of ${\td\af}_0^j$ and
${\td\af}_1^j$).
{\blue{{{Condition}} (4) follows from  \eqref{134-n1}, and  {{condition}} (2)  follows from
{{condition}} (4) and the following commutative diagram
\begin{displaymath}
\xymatrix{
K_0(C_n)=G_n\ar[r]^{~}\ar@{->}[d]_{(\phi_{n,n+1})_{*0}
} &
K_0(C_n/I_n)=H_n \ar@{->}[d]_{(\td\phi_{n, n+1})_{*0}}
\\
 K_0(C_{n+1})=G_{n+1}\ar[r]^{~}& K_0(C_{n+1}/I_{n+1})=H_{n+1}.
  }
\end{displaymath}
}}

 If $x\in Sp(F_{n+1})\subset Sp(C_{n+1})$, then $Sp(F_n)\subset Sp(\phi_{n,n+1}|_x)\, (=Sp(\td \phi_{n,n+1}|_x)),$
{by}  the fact that all entries of $\cc$ are strictly positive.
If $x\in (0,1)_j=Sp(C_0((0,1), E_{n+1}^j)$, then each $\tht_i$, as the only element in $Sp(F_n^i)\, (\subset Sp(F_n))$,
 appears $r_i^j>0$ times in
$Sp(\phi_{n,n+1}|_x)$ {{(see \eqref{13Nov30-2018}.}}
 Consequently, we also have $Sp(F_n)\subset Sp(\phi_{n,n+1}|_x)$. Hence condition (5) holds.

 {{To see (6),  we will use \eqref{13Nov30-2018}  and \eqref{defpsi}.
 Note each $F_k(t)$ appears in  \eqref{defpsi}.}}
 {\blue{Therefore if $d_{jk}>0$, then, {{by (\ref{13Jan16-2019}),}}
 the point  $t\in (0,1)_{n,k}{{\subset}} Sp(C([0,1], E_n^k))$ is {{in}}
 $Sp(\phi_{n,n+1}|_{t})$ for $t\in (0,1)_{n+1,j} {{\subset}} Sp(C([0,1], E_{n+1}^j))$; if $d_{jk}<0$, then, {{by (\ref{13Jan16-2019-1}),}} the
  point  \\
  $1-t\in (0,1)_{n,k} {{\subset}} Sp(C([0,1], E_n^k))$ is in the
  $Sp(\phi_{n,n+1}|_{t})$ for $t\in (0,1)_{n+1,j}\subset Sp(C([0,1], E_{n+1}^j))$; and if $d_{jk}=0$, then,
  {{by (\ref{13Jan16-2019-2}),}} both points $t$ and  $1-t$ from
  $ (0,1)_{n,k}\subset Sp(C([0,1], E_n^k))$ are  in the $Sp(\phi_{n,n+1}|_{t})$ for $t\in (0,1)_{n+1,j}\subset Sp(C([0,1],
   E_{n+1}^j))$. Hence (6) follows.}}

 {\blue{From (5), we know $Sp(F_n)\subset Sp(\phi_{n,n+1}|_y)$ for any $y\in Sp(C_{n+1})$. From (6) and $l_{n+1}>0$, we know that ${{\bigcup}}_{i=1}^{l_n} (0,1)_{n,i}\subset {{\bigcup}}_{t\in (0,1)_{n+1, 1}} Sp(\phi_{n,n+1}|_t)$. {{Thus}} $Sp(C_n)\subset Sp(\phi_{n,n+1}).$
 {{By the definition of $Sp(\phi_{n,n+1}),$ this implies
 ${\rm ker}\phi_{n,n+1} \subset \bigcap_{x\in Sp(C_n)}\{{\rm ker}\, x\}=\{0\}.$}} It follows that $\phi$ is injective.
  Hence (7) holds.
   Note that  if $X\subset (0,1)_{n+1,1}$ is $\dt$-dense, then the sets $\{x: x\in  X\}=X$ and $\{1-x, x\in X\}$ are $\dt$-dense in $(0,1)$. Hence  $\bigcup_{t\in X} Sp(\phi_{n,n+1}|_t)$ is $\dt$-dense in $\bigcup_{i=1} (0,1)_{n,i}$. Since $Sp(C_n)=Sp(F_n)\cup\bigcup_{i=1} (0,1)_{n,i}$, combining with (5), we get (8). }}
This finishes the proof.
\end{proof}








\begin{NN}\label{range 0.8}  Let $\phi:~ C_n\to C_{n+1}$ be as in the proof above.
We will calculate the {\blue{contractive linear}}
map
$$\phi_{n,n+1}^{\sharp}:~ \Aff(T(C_n))\to \Aff(T(C_{n+1}))$$
{\blue{which also preserves the order (i.e., maps $\Aff(T(C_n))_+$ to $\Aff(T(C_{n+1}))_+$).}}
Recall from {{\ref{LgN1889} (see \ref{2Rg15} also)}} that
$\Aff(T(C_n))$ {{is the subset of}} $ \bigoplus_{i=1}^{l_n}C([0,1]_i,\R)\oplus \R^{p_n}$
consisting of {{the}} {{elements}} $(f_1,f_2,...,f_{l_n}{{;}}~h_1,h_2,...,h_{p_n})$ which  {{satisfy the}}
conditions
\beq\label{1508/13star-2}
 f_i(0)=\frac1{\{n,i\}}\sum b_{0,ij}h_j\cdot [n,j] \andeqn
f_i(1)=\frac1{\{n,i\}}\sum b_{1,ij}h_j\cdot [n,j],
\eneq
and $\Aff(T(C_{n+1}))$ is {{the}} subset of $\bigoplus_{i=1}^{l_{n+1}}C([0,1]_i,\R)\oplus
\R^{p_{n+1}}$ consisting of {\blue{the elements}}\\
$(f'_1,f'_2,...,f'_{l_{n+1}}{{;}}~h'_1,h'_2,...,h'_{p_{n+1}})$ which
{\blue{satisfy}}
\beq\label{13Dec9-2018-2}
\hspace{-0.2in}f'_i(0)=\frac1{\{n+1,i\}}\sum b'_{0,ij}h'_j\cdot [n+1,j] \,\,\,{\rm and}\,\,\,
 f'_i(1)=\frac1{\{n+1,i\}}\sum b'_{1,ij}h'_j\cdot
[n+1,j].
\eneq
{{Recalling}} that $(c_{ij})_{p_{n+1}\times p_n}$ is the matrix
corresponding to $(\td\phi_{n,n+1})_{*0}=\gm_{n,n+1}$ for
$\td\phi_{n,n+1}: F_n\to F_{n+1}$, and {{noting}} that since  $\td\phi_{n,n+1}$ is unital, one has
$$\sum_{j=1}^{p_n}c_{ij}[n,j]=[n+1,j].$$
{{Let}}
$$\phi_{n,n+1}^{\sharp}(f_1,f_2,...,f_{l_n}{{;}}~h_1,h_2,...,h_{p_n})
=(f'_1,f'_2,...,f'_{l_{n+1}}{{;}}~h'_1,h'_2,...,h'_{p_{n+1}}).$$ Then
$$h_i'=\frac1{[n+1,i]}\sum_{j=1}^{p_n}c_{ij}h_j[n,j].$$
{{Combining this with (\ref{13Dec9-2018}), we have
\beq\label{13Dec9-2018-1}
f'_i(0)=\frac1{\{n+1,i\}}\sum\td b_{0,il}h_l[n, l]~~\mbox{and}~~f'_i(1)=\frac1{\{n+1,i\}}\sum\td b_{1,il}h_l[n, l]
\eneq}}
Also note that
\beq\nonumber
f'_i(t)=\frac1{\{n+1,i\}}\left\{\sum_{d_{ik}>0}d_{ik}f_k(t)\{n,k\}+
\sum_{d_{ik}<0}|d_{ik}|f_k(1-t)\{n,k\}+\right.\qq\qq\qq\qq\\\label{13starss}
\qq\qq\qq\qq +\left.\sum_{d_{ik}=0}(f_k(t)+f_k(1-t))\{n,k\}+
\sum_{l=1}^{p_n}r_l^ih_l[n,l]\right\}, \qq\qq
\eneq
where
 \begin{eqnarray}\label{rli}
r_l^i & = & \sum_{k=1}^{p_{n+1}}b_{0,ik}'c_{kl}-
\left(\sum_{d_{ik}<0}|d_{ik}|b_{1,kl}+
\sum_{d_{ik}>0}|d_{ik}|b_{0,kl}+ \sum_{d_{ik}=0}(b_{0,kl}+b_{1,kl})
\right)\nonumber \\
& = & \sum_{k=1}^{p_{n+1}}b_{1,ik}'c_{kl}-
\left(\sum_{d_{ik}>0}|d_{ik}|b_{1,kl}+
\sum_{d_{ik}<0}|d_{ik}|b_{0,kl}+ \sum_{d_{ik}=0}(b_{0,kl}+b_{1,kl})
\right).
\end{eqnarray}
It follows from the last paragraph of \ref{range  0.5} and {{from}} \ref{range
0.5a} that, when we define $C_{n+1}$, we can always increase the entries
of the last $p_{n+1}^1=p_{n+1}- p_{n+1}^0$  columns of the matrices
$\bb_0'=(b_{0,ik}')$ and $\bb_1'=(b_{1,ik}')$ by adding an
{{arbitrary}} (but {{the}} same for $\bb_0'$ and $\bb_1'$)
 matrix $(m_{ik})_{l_{n+1}{{\times p}}_{n+1}^1}$, with each $m_{ik}>0$
 sufficiently large, to the last $p_{n+1}^1$ columns of the matrices.
 In particular we can strengthen the requirement $\spd$ (see (\ref{13spd-1})) to {condition $\spdd$, i.e.,}
 \beq\label{13spdd}
\td b_{0,il}\left(= \!\sum_{k=1}^{p_{n+1}}b_{0,ik}'c_{kl}\right),\
\td b_{1,il}\left(=\sum_{k=1}^{p_{n+1}}b_{1,ik}'c_{kl}\!\!\right)
 > 2^{2n}\!\left(\sum_{k=1}^{l_n}
 (|d_{ik}|+2)\{n,k\}\!\!\right)
 \eneq
 for all
 $i\!\in\!\{1,...,l_{n+1}\}.$
This condition and \eqref{rli} (and note that $b_{0,kl}\leq \{n,k\}$, $b_{1,kl}\leq \{n,k\}$ for any $k\leq l_n$) imply
\beq\label{13spdd'}
r_l^i \geq \frac{2^{2n}-1}{2^{2n}}\td b_{0,il}, \qq \mbox{or
equivalently} \qq 0\leq \td b_{0,il}-r_l^i< \frac{1}{2^{2n}}\td
b_{0,il}.
\eneq
Recall that $\phi_{n,n+1}^{\sharp}:~\hspace{-0.05in}\Aff(T(C_{n}))\to \Aff(T(C_{n+1}))$ and
$\td\phi_{n,n+1}^{\sharp}:~\hspace{-0.06in} \Aff(T(F_{n}))\to \Aff(T(F_{n+1}))$ are the maps induced by the
homomorphisms
$\phi_{n,n+1}$ and $\td\phi_{n,n+1}$; and the same for $\pi_n^{\sharp}:~\hspace{-0.07in} \Aff(T(C_{n}))\to \Aff(T(F_{n}))$ and
$\pi_{n+1}^{\sharp}: ~ \Aff(T(C_{n+1}))\to \Aff(T(F_{n+1}))$, which are the maps
induced by the quotient
 maps
  $\pi_n: C_n\to F_n$ and $\pi_{n+1}: C_{n+1}\to F_{n+1},$ respectively.
    Since $\td\phi_{n,n+1}\circ \pi_n=\pi_{n+1}\circ \phi_{n,n+1}$, we have
\beq\label{ee1316}
\pi_{n+1}^{\sharp}\circ\phi_{n,n+1}^{\sharp}=\td\phi_{n,n+1}^{\sharp}\circ\pi_n^{\sharp}:~ \Aff(T(C_{n+1}))\to \Aff(T(F_{n+1})).
\eneq

\end{NN}

\begin{NN}\label{range 0.9}
For each $n,$
we will now define
a map $\GM_n:~ \Aff(T(F_{n}))\to \Aff(T(C_{n}))$ which is a right inverse of
$\pi_n^{\sharp}:~ \Aff(T(C_{n}))\to \Aff(T(F_{n}))$---that is, $\pi_n^{\sharp}\circ
\GM_n=\id|_{_{\Aff(T(F_n))}}.$

Recall that $C_n=A(F_n,E_n,\bt_0,\bt_1)$ with unital homomorphisms
$\bt_0,\bt_1: F_n\to E_n$ whose K-theory maps satisfy
$(\bt_0)_{*0}=\bb_0=(b_{0,ij})$ and  $(\bt_1)_{*0}=\bb_1=(b_{1,ij})$.
{\blue{Let $\bt_i^{\sharp}: \Aff(T(F_n))\to \Aff(T(E_n))$ be the contractive linear order preserving
map induced by the \hm\, $\bt_i,$ $i=0,1.$}}
{\blue{For each $h\in \Aff(T(F_n)),$ consider the function}}
\beq\label{13Dec11-3-2018}
\GM_n'(h)(t)= t\cdot \bt_1^{\sharp}(h)+(1-t)\cdot \bt_0^{\sharp}(h),
\eneq
an element of  $C([0,1],\R^{l_n})=\bigoplus_{i=1}^{l_n}C([0,1]_i,\R)$.
Finally, define the map
\beq\nonumber
&&\hspace{-0.2in}\GM_n:~ \Aff(T(F_{n}))=\R^{p_n}\to \Aff(T(C_{n})) \sbs
\bigoplus_{j=1}^{l_n}C([0,1]_j,\R)\oplus \R^{p_n} \,\,\,{\rm by}\\ \label{13Dec11-4-2018}
&&\GM_n(h)=(\GM_n'(h), h)\in \bigoplus_{j=1}^{l_n}C([0,1]_j,\R)\oplus
\R^{p_n}.
\eneq
Note that $\GM_n(h)\in \Aff(T(C_{n}))$ (see (\ref{1508/13star-2})
in \ref{range 0.8}).
{\blue{One verifies that the map $\Gamma_n$ is a contractive linear and order preserving map
from $\Aff(T(F_n))$ to $\Aff(T(C_n)).$}}
 Evidently, $\pi_n^{\sharp}\circ \GM_n=\id|_{_{\Aff(T(F_n))}}$.\\

\end{NN}

\begin{lem}\label{range 0.10}
 If {Condition} $\spdd$ (see $(\ref{13spdd})$) holds, then for any $f\in
\Aff(T(C_{n}))$ with $\|f\|\leq 1$, and $f'{{:=}}\phi_{n,n+1}^{\sharp}(f)\in \Aff(T(C_{n+1}))$, we
have
$$\| \GM_{n+1}\circ\pi_{n+1}^{\sharp}(f')-f'\|<\frac2{2^{2n}}.$$
\end{lem}

\begin{proof} {{As in \ref{range 0.8} (see \ref{LgN1889} also), one can
w}}rite
\beq\nonumber
f=(f_1,f_2,...,f_{l_n}{{;}}~h_1,h_2,...,h_{p_n})\in \Aff(T(C_{n}))
\andeqn\\
f'=(f'_1,f'_2,...,f'_{l_{n+1}}{{;}}~h'_1,h'_2,...,h'_{p_{n+1}})\in
\Aff(T(C_{n+1})).
\eneq
{{(Note that the  element $f$ satisfies (\ref{1508/13star-2}) and (\ref{1508/13star-2}); and the element $f'$ satisfies (\ref{13Dec9-2018-2}) and  (\ref{13Dec9-2018-3}).)}} {{Also, we have $$\|f_i\|, \|h_j\|\leq \|f\|\leq 1~~ \mbox{for}~~ 1\leq i\leq l_n, 1\leq j\leq p_n,~~\mbox{and}$$
$$\|f'_i\|, \|h'_j\|\leq \|f'\|\leq 1~~ \mbox{for}~~ 1\leq i\leq l_{n+1}, 1\leq j\leq p_{n+1}.$$}}  Since
$\pi_{n+1}^{\sharp}\circ\GM_{n+1}=\id|_{_{\Aff(T(F_{n+1}))}}$, one has
$$\GM_{n+1}\circ\pi_{n+1}^{\sharp}(f'):=g'
:=(g'_1,g'_2,... ,g'_{l_{n+1}}{{;}}~ h'_1,h'_2,... ,h'_{p_{n+1}});$$
that is, $f'$ and $g'$ have the same boundary value
$(h'_1,h'_2,...,h'_{p_{n+1}})$.

Note that{{, from (\ref{13Dec9-2018-2}),}}
the evaluation of $f'$ at zero, $(f'_1(0),f'_2(0),...,f'_{l_{n+1}}(0))$, and the evaluation at {{one}}, $(f'_1({{1}}),f'_2({{1}}),...,f'_{l_{n+1}}({{1}})),$
are
completely determined by $h'_1,h'_2,...,h'_{p_{n+1}}$.
{{By}} {{(\ref{13Dec9-2018-1}) and}} (\ref{13starss}), we also have
$$
f'_i(t)-f'_i(0)=\frac1{\{n+1,i\}}\left(\sum_{d_{ik}>0}d_{ik}f_k(t)\{n,k\}+
\sum_{d_{ik}<0}|d_{ik}|f_k(1-t)\{n,k\}\right.+\qq\qq\qq\qq\qq$$
\beq\label{13Dec11-2-2018}
\qq\qq\left.+\sum_{d_{ik}=0}(f_k(t)+f_k(1-t))\{n,k\}-
\sum_{l=1}^{p_n}(\td b_{0,il}-r_l^i)h_l[n,l]\right).
\eneq
From {Condition} $\spdd$ (see \eqref{13spdd}) {{and $\|f\|\leq1$}}, one has
\beq\nonumber
&&\hspace{-0.2in}\left|\sum_{{\{k:d_{ik}>0\}}}d_{ik}f_k(t)\{n,k\}+
\sum_{{\{k:d_{ik}<0\}}}|d_{ik}|f_k(1-t)\{n,k\}
+\sum_{{\{k:d_{ik}=0\}}}(f_k(t)+f_k(1-t))\{n,k\}\right|\\\nonumber
&&\leq \sum_{k=1}^{l_n}(|d_{ik}|+2)\{n,k\}{{\|f_k\|}}\qq\qq\,\,\,\,\,{{\mbox{(since}~~ |d_{ik}|\leq |d_{ik}|+2,~ 2\leq |d_{ik}|+2)~ }}\qq\qq\\
\nonumber
&&\leq \frac1{2^{2n}} \td
b_{0,il}\qq\qq\qq\qq\qq\qq {{\mbox{(by}~ \eqref{13spdd}~~ \mbox{and}~ \|f_k\|\leq\|f\|\leq 1)}}\qq\qq\qq\qq\qq\\
\label{13Dec-11-2018}
&&< \frac1{2^{2n}}\cdot\{n+1,i\}\qq\qq\qq\hspace{0.4in}{{\mbox{(by (\ref{sizeE}) and } [n,j]\ge 1)}}.\hspace{-0.2in}
\eneq
{{By (\ref{13spdd'}), (\ref{sizeE}), and $\|h_l\|\leq \|f\|\leq 1$, one {\blue{obtains}}
\beq\label{13Dec-11-1-2018}
\left|\sum_{l=1}^{p_n}(\td b_{0,il}-r_l^i)h_l[n,l]\right|\leq\frac{1}{2^{2n}}\td b_{0,il}\|h_l\|[n,l]\leq \frac{1}{2^{2n}}\{n+1,i\}.
\eneq}}

Combining {{(\ref{13Dec11-2-2018}), (\ref{13Dec-11-2018}), and (\ref{13Dec-11-1-2018}),}} {\blue{one has}}
$$|f_i'(t)-f_i'(0)|<\frac2{2^{2n}}.$$
Similarly, {{one has}}
$$|f_i'(t)-f_i'(1)|<\frac2{2^{2n}}.$$
But, by the definition of $\GM_{n+1}$ {{(with $n+1$ in place of $n$) in (\ref{13Dec11-3-2018}) and (\ref{13Dec11-4-2018})}}, we have
$$g_i'(t)=tg_i'(1)+(1-t)g_i'(0).$$
Combining  {{this}} with $g_i'(0)=f_i'(0)$ and $g_i'(1)=f_i'(1)$, we have
$$|g_i'(t)-f_i'(t)|<\frac2{2^{2n}}\rforal i,$$
as desired.
\end{proof}

\begin{thm}\label{range 0.12} {{Let $((G, G_+, u), K, \DT,r)$ be {\blue{a}} weakly unperforated Elliott invariant {{as}}  defined in {\rm \ref{range 0.1}} with $G$ torsion free and $K=0$.}} {{Let $(H, H_+, u)$ be the simple ordered group with $H\supset G$ defined in {\rm \ref{range 0.3}}}}
{\blue{Let $C_n,$ $\phi_{n,n+1}: C_n\to C_{n+1},$ {{$F_n$,}} and ${\bar \phi}_{n,n+1}: F_n\to F_{n+1}$ be as {{in}} Lemma {\rm \ref{range  0.6}}, {{and}} let
$I_n$  be as {{in Remark}} {\rm{\ref{range 0.7}}}.
{{T}}hen the inductive limit $C=\lim
(C_n,\phi_n)$ has the property
$({{(}}K_0(C),
K_0(C)_+,\e_C{{), K_1(C)}})
=({{(}}G,G_+,{{u), K,}}).$
In particular, $K_1(C)=\{0\}.$
Moreover, $I=\lim_{n\to\infty}(I_n, \phi_{n,n+1}|_I)$ is an ideal of $C$ such that
$C/I=F=\lim_{n\to\infty}(F_n, {\bar \phi}_{n,n+1})$ with $K_0(F)=H$ and $T(C/I)=\Delta.$}} {{(If we further assume that
Condition $\spd$ in $(\ref{13spd-1})$
is replaced by the stronger Condition $\spdd$ in  $(\ref{13spdd})$, then we  have $T(C)=T(C/I)=\DT$. We will
not use this. However, it follows from the proof of Theorem
\ref{EllofA}.)}}
\end{thm}


\begin{proof} {{Recall, from \ref{range 0.5}, that $C_n=C([0,1],E_n)\oplus_{\bt_0,\bt_1}F_n$ with ideal $I_n=C_0((0,1), E_n)$ and quotient $C_n/I_n=F_n$. Hence, }}  we have the following infinite  commutative diagram:
{\small{\begin{displaymath}
    \xymatrix{ I_1 \ar[r]\ar[d] & I_2 \ar[r]\ar[d]& I_3 \ar[d]\ar[r]&\cd I \\
               C_1 \ar[d]\ar[r] & C_2 \ar[r]\ar[d] &C_3 \ar[d]\ar[r]&\cd C \\
               {{F_1}} \ar[r] & {{F_2}} \ar[r] &{{F_3}} \ar[r]&\cd C/I~{{.}}}
\end{displaymath}}}
  Also from the construction, we have the
following diagram:
{\small{
\begin{displaymath}
    \xymatrix{0 \ar[r]\ar[d] & 0 \ar[r]\ar[d]& 0 \ar[d]\ar[r]&\cd  \\ K_0(C_1)=G_1 \ar[r]\ar[d]^{\blue{\iota_1}} & K_0(C_2)=G_2 \ar[r]\ar[d]^{{\blue{\iota_2}}}& K_0(C_3)=G_3 \ar[d]^{{\blue{\iota_3}}}\ar[r]&\cd  \\
               K_0({{F_1}})=H_1 \ar[d]\ar[r] & K_0({{F_2}})=H_2 \ar[r]\ar[d] &K_0({{F_3}})=H_3 \ar[d]\ar[r]&\cd \\
               K_1(I_1)\ar[d] \ar[r] & K_1(I_2)\ar[d] \ar[r] &K_1(I_3)\ar[d] \ar[r]&\cd \\
               0 \ar[r] & 0 \ar[r]& 0 \ar[r]&\cd ~,}
\end{displaymath}
}}
$\\${\blue{where the inclusion map{{s}} $\iota_{{n}}: K_0(C_{{n}})= G_{{n}}\to K_0({{F_n}})=H_{{n}}$
are  induced by the \hm s $\pi_i: C_n\to F_n.$}}
So $K_0(C)=\lim(G_n,\gm_{n,n+1}|_{_{G_n}})=G$.
{\blue{Let $\iota: K_0(C)\to K_0(C/I)$ be the \hm\, given by the above diagram.}} {\blue{Note that the AF algebra $C/I$ has  $K_0(C/I)=\lim(H_n,\gm_{n,n+1})=H$ as scaled ordered group. Hence by the  {{part (1) of Remark \ref{March-22-2019}}}, we have $T(C/I)=\DT$.}}

\end{proof}

\begin{NN}\label{notation-Dec12-2018}
{{Let us  fix some notation. Recall that $C_n=C([0,1],E_n)\bigoplus_{(\bt_{n,0},\bt_{n,1})}F_n$, where $E_n={{\bigoplus}}_{i=1}^{l_n}E_n^i=\bigoplus_{i=1}^{l_n}M_{\{n,i\}}$, $F_n=\bigoplus_{i=1}^{p_n}F_n^i=\bigoplus_{i=1}^{p_n}M_{[n,i]}$, and
$\bt_{n,0}, \bt_{n,1}: F_n\to E_n$ are unital homomorphisms. (Here we use  $\bt_{n,0}, \bt_{n,1}$ instead of $\bt_0, \bt_1$ to distinguish  the maps for different $n$.)}}

{{In the rest of this section, we will use $t_{n,j}$ to denote the representation
$$C\ni (f_1, f_2, \cdots, f_{l_n};~ a_1, a_2, \cdots, a_{p_n})\mapsto f_j(t)\in E^j_n=M_{\{n,j\}},$$
and $\theta_{n, i}$ to denote the representation
$$C\ni (f_1, f_2, \cdots, f_{l_n};~ a_1, a_2, \cdots, a_{p_n})\mapsto a_i\in F^i_n=M_{[n,j]},$$
for $t\in [0,1]$,  $j=1,2, \cdots l_n,$ and  $i =1,2, \cdots, p_n$.}}

{{Note that $\theta_{n, i}\in Sp(C_n)$ and $t_{n,j}\in Sp(C_n)$ for $t\in (0,1)$, but in general, $0_{n, j}, 1_{n, j}\notin Sp(C_n)${{; they}} may not be irreducible. In $RF(C_n)$, in the notation of \ref{homrestr}, we have
$$0_{n, j}=\{\theta_{n, 1}^{\sim b_{0,j1}},\theta_{n, 2}^{\sim b_{0,j2}}, \cdots,\theta_{n, p_n}^{\sim b_{0,jp_n}}\}~~ \mbox{and}~~1_{n, j}=\{\theta_{n, 1}^{\sim b_{1,j1}},\theta_{n, 2}^{\sim b_{1,j2}}, \cdots,\theta_{n, p_n}^{\sim b_{1,jp_n}}\}, $$
where $(b_{0,ji})_{l_n\times p_n}$ and $(b_{1,ji})_{l_n\times p_n}$ are the matrices corresponding to $(\bt_{n,0})_*,(\bt_{n,1})_*: \Z^{p_n}\to \Z^{l_n}$.
{\blue{Let us denote}} the set of all points $t_{n,j}$ for $t\in [0,1]$ by $[0,1]_{n,j}$. Hence $Sp(A_n)=\{\theta_{n,1},\theta_{n,2},\cdots, \theta_{p_n}\}\cup \bigcup_{j=1}^{l_n}(0,1)_{n,j}$. }}

{{In the previous construction, we only describe how to construct $C_{n+1}$ from $C_n$}}
{\blue{for a fixed $n.$}}
{\blue{Now we need to let $n$ move.
Accordingly, we need to change  some notation.}}
{\blue{For example,}} {{$\cc=(c_{ij})$ of \ref{range 0.4} will be denoted by $\cc_{n,n+1}=(c^{n,n+1}_{ij})$ as it is the matrix of $
\gm_{n,n+1}:~ H_n=\Z^{p_n}\lr H_{n+1}=\Z^{p_{n+1}}
$; $\bb_0=(b_{0,ji})$ and $\bb_1=(b_{1,ji})$  of \ref{range 0.5} (and $\bb'_0=(b'_{0,ji})$ and $\bb'_1=(b'_{1,ji})$ of \ref{conditions})   will be denoted by $\bb_{n,0}=(b^n_{0,ji})$ and $\bb_{n,1}=(b^n_{1,ji})$ (and $\bb_{n+1,0}=(b^{n+1}_{0,ji})$ and $\bb_{n+1,1}=(b^{n+1}_{1,ji})$, respectively); $\dd=(d_{ij})$ of \ref{conditions} will be denoted by $\dd_{n,n+1}=(d^{n,n+1}_{ij})$ as it is the matrix of $\td\gm_{n,n+1}:~H_n/G_n=\Z^{l_n} \to H_{n+1}/G_{n+1}=\Z^{l_{n+1}}$;  $\td \bb_0= (\td b_{0,ji})$ {{and}} $\td \bb_1= (\td b_{1,ji})$  will be denoted by $\td \bb_{(n,n+1),0}= (\td b^{n,n+1}_{0,ji})$ {{and}} $\td \bb_{(n,n+1),1}= (\td b^{n,n+1}_{1,ji})$ (hence,
the relations $\td \bb_0=\bb'_0\cdot \cc$ and $\bb_1=\bb'_1\cdot \cc$  there will be $\td \bb_{(n,n+1),0}=\bb_{n+1,0}\cdot \cc_{n,n+1}$ {{and}}  $\td \bb_{(n,n+1),1}=\bb_{n+1,1}\cdot \cc_{n,n+1}$).}}


\end{NN}

\begin{NN}\label{1323proj}
{{In this paper, a projection $p\in M_l(C(X))$ is called trivial \index{trivial projection} if there is a unitary $u\in M_l(C(X))$ such that}} {\blue{$u^*(x)p(x)u(x)=\diag(1,1,\cdots,1, 0,\dots,0).$}}
{\blue{By part 2 of Remark 3.26 of \cite{EG-RR0AH},
if $X$ is a finite CW complex with dimension at most three,
then}}
{{the above statement is equivalent to {{the statement}} that $p$ is}} {{ Murray-von}}
{\blue{Neumann}} {{equivalent to}} {\blue{$\diag(1, 1, \cdots, 1,0,\cdots, 0).$}}
{{A projection $p\in QM_l(C(X))Q$ is trivial if it is trivial when regarded as a projection in $M_l(C(X))$.}}

{\blue{Let $X$ be a  connected finite CW complex with dimension at most three with  the base point $x_0.$
Then $K_0(C(X))=K_0(C_0(X\setminus \{x_0\}))\oplus K_0(\C)=K_0(C_0(X\setminus \{x_0\}))\oplus \Z.$}}

{\blue{We will use the fact  $K_0(C(X))_+\setminus \{0\}=\{(x, n)\in K_0(C_0(X\setminus \{x_0\}))\oplus \Z, \, n>0\}.$
It is clear that $K_0(C(X))_+\setminus \{0\}\subset \{(x, n)\in K_0(C_0(X\setminus \{x_0\}))\oplus \Z, \, n>0\}.$
Let $(x, n)\in  \{(x, n)\in K_0(C_0(X\setminus \{x_0\}))\oplus \Z, \, n>0\}.$
Choose $p\in M_m(C(X))$ such that $[p]-[1_k]=(x, n)$ for some  integers $0\le k\le m.$
Then $p$ has rank $n+k.$
Since $n\ge 1\ge {{(3+1)/2-1}},$ it follows from Theorem 9.12 of \cite{Hus} that $1_{k}$ is  unitarily equivalent to a subprojection of $p.$
Therefore there is a projection $q\le p$ such that $[q]=[1_{k}].$ Then $(x,n)=[p-q]=[p]-[q]\in K_0(A)_+\setminus \{0\}.$
But $[p]-[q]=(x,n).$
}}

\end{NN}

The following theorem is in Section 3 of \cite{EG-RR0AH}
{\blue{(see Proposition  3.16 and Theorem 3.10 there).}}


\begin{prop}\label{range 0.14}
Let $X$ and $Y$ be path-connected finite CW complexes of dimension
at most three, with base {{points}} $x_0\in X$ and $y_0\in Y$, such that
the cohomology groups $H^3(X)$ and $H^3(Y)$ are finite. Let $\af_0:~
K_0(C(X))\to K_0(C(Y))$ be a homomorphism
{\blue{such that}} $\af_0$ is at least 12-large {\blue{(see the definition in Section 3 of  {\rm{\cite{EG-RR0AH})}}}}
and
\beq\label{13Dec-14-2018}
\af_0(K_0(C(X))_+\setminus \{0\})\sbs K_0(C(Y))_+\setminus \{0\},\eneq and let $\af_1:~
K_1(C(X))\to K_1(C(Y))$ be any homomorphism.  Let $P\in M_\infty
(C(X))$ be any {non-zero} projection and $Q\in M_\infty (C(Y))$ be a
projection with $\af_0([P])=[Q]$ (such projections always exist{{; see the proof below}}).
Then there exists a unital homomorphism $\phi:~ PM_\infty (C(X))P\to
QM_\infty (C(Y))Q$ such that $\phi_{*0}=\af_0$ and
$\phi_{*1}=\af_1$, and such that
$$
\phi(PM_\infty (C_0(X\setminus \{x_0\}))P)\,\sbs\, QM_\infty
(C_0(Y\setminus \{y_0\}))Q.$$
That is, if $f\in PM_\infty (C(X))P$ satisfies $f(x_0)=0$, then $\phi(f)(y_0)=0$.
{\blue{Moreover, if $\af_0$ is at least {\rm 13}-large and $Y$ is not a single point,
 then one may further require that $\phi$ is injective.}}

{\blue{Suppose that $B=\bigoplus_{i=1}^m B_i$ and $D=\bigoplus_{j=1}^n D_j,$ where each $B_i$ has the form $P_iM_{\infty}(C(X_i))P_i$
and {{each}} $D_j$ {{has the form}} $Q_jM_{\infty}(C(Y_j))Q_j,$ where $X_i,$ $Y_j$ {{are connected finite CW complexes with at least one $Y_i$ not a single point, and}}  $P_i$ and $Q_j$ are as in the first part of the {{statement.}}
Fix the base points $x_i\in {{X_i}}$ and $y_j\in Y_j.$    Let $I_i=\{f\in B_i: f(x_i)=0\}$ and
$J_j=\{g\in D_j: g(x_j)=0\}.$
{{Let}} $\af_0: K_0({{B}})\to K_0({{D}})$ be an order preserving \hm\, with $\af_0([1_B])=[1_D]$ such that $\af_0|_{K_0(B_i)}$ is {{at least}} $13$-large {{in}} each component of $K_0(D_j),$ and {{let}} $\af_1: K_1(B)\to K_1(D)$ {{be}}
any \hm.}}
{\blue{Then there is a unital injective \hm\, $\phi: B\to D$ such that
$(\phi)_{*i}=\af_i,$ $i=0,1,$ and $\phi(I)\subset J,$
where $I=\bigoplus_{i=1}^m I_i$ and $J=\bigoplus_{j=1}^n J_j.$}}

\end{prop}

\begin{proof} {{
Let us first prove the first part  of the lemma.  Since $[P]\in K_0(C(X))_+\setminus \{0\}$, by the condition (\ref{13Dec-14-2018}), we have $\af_0([P])\in K_0(C(Y))_+\setminus \{0\}$. Hence there is a projection $Q\in M_\infty (C(Y))$ such that $\af_0([P])=[Q]$. }}
{\blue{Similarly, there is a projection $q\in M_k(C(Y))$  (for some integer $k\ge 1$) such that $\af_0([1_{C(X)}])=[q].$}}

{\blue{This  proposition is  a special case  of Theorem 9.1 of \cite{Lin-LAH}.
To see this, {{first recall}} $K_0(C(X))=K_0(C_0(X\setminus \{x_0\}))\oplus {{\Z}}$ and
$K_0(C(Y))=K_0(C_0(Y\setminus \{y_0\})\oplus {{\Z}}.$ 
{{By hypothesis,}}
${\rm rank}\,q\,\ge 12.$}}
{\blue{Note also $K_1(C(X))=K_1(C_0(X\setminus \{x_0\}))$ and $K_1(C(Y))=K_1(C_0(Y\setminus \{x_0\})).$
By \ref{1323proj}, $K_0(C(X))_+\setminus \{0\}=\{(z_1, z_2)\in K_0(C_0(X\setminus\{x_0\}))\oplus {{\Z}}: z_2>0\}.$
  One then computes}} {{that  the condition (\ref{13Dec-14-2018}) implies $\af_0(K_0 (C_0(X\setminus\{x_0\})))\subset K_0 (C_0(Y\setminus\{y_0\})).$}}
 {\blue{ Let $\kappa\in KK(C_0(X\setminus \{x_0\}), C_0(Y\setminus \{y_0\}))$ be such that
  $\kappa|_{K_0(C_0(X\setminus \{x_0\})}=(\af_0)|_{K_0(C_0(X\setminus \{x_0\}))}$ and
  $\kappa|_{K_1(C_0(X\setminus \{x_0\}))}=\af_1.$}}
{{By the definition of $F_3K_*(C(X))$ in 3.7 of \cite{EG-RR0AH} and by 3.4.4 of \cite{DN-Shape}, we know $F_3K_*(C(X))=H^3(X)$, for a finite CW complex $X$ of dimension {\blue{at}} most three. Note that by our assumption on $X$ and $Y$, we have  $H^3(X)={\rm Tor} (K_1(C(X))),$ and $H^3(Y)= {\rm Tor} (K_1(C(Y))),$ {{and so}} ${\blue{\kappa_*}}(F_3K_*(C(X)))\subset F_3K_*(C(Y)),$
  since $\kappa_*$ maps torsion elements to
   torsion elements.}}

{\blue{Since ${\rm rank}\, q\ge 12\ge 3({\rm dim}Y+1),$  we may write $q=q_0\oplus e,$
where $e$ is zero, or $e$ is a rank one trivial projection when ${\rm rank}\, q\ge 13$ (in the case $\af_0$ is
at least 13-large), by Theorem 9.1.2 of \cite{Hus}.
If $e\not=0,$ choose a surjective homotopy-trivial
continuous map $s:~ Y\to X$, which induces an injective homomorphism
$s^*:~ C(X)\to C(Y).$
Let $\af_0'(g,z)=\af_0(g)+{{z([q]-[e])}}$ {{(note that $\af_0(0,1)=[q]$ as $(0,1)=[1_C(X)]\in K_0(C(X))$)}} for all $(g,z)\in
K_0(C_0(X\setminus \{x_0\})\oplus \Z=K_0(C(X))$ as in \ref{1323proj}.  For $\af_0'$ and $\af_1,$
since $\af_0'$ is at least 12-large, there is a unital \hm\, $\psi_0: C(X)\to q_0M_{\infty}(C(Y))q_0$
such that $(\psi_0)_{*0}=\af_0'$ and $(\psi_0)_{*1}=\af_1.$
In both cases,
by  Theorem 9.1 of \cite{Lin-LAH}, we obtain a unital \hm\,
$\psi_0: C(X)\to q_0M_{{\infty}}(C(Y))q_0\subset C(Y)\otimes {\cal K}$ such that
$[\psi_0|_{C_0(X\setminus \{x_0\})}]=\kappa.$    In the case $e=0,$ $q=q_0.$ Let $\psi=\psi_0.$
If $e\not=0,$
define $\psi: C(X)\to qM_k(C(Y))q$ by $\psi(f)=\psi_0(f)\oplus s^*(f)$
for all $f\in C(X).$
Then, in this case, $\psi$ is a unital injective \hm\, and we still have $[\psi|_{C_0(X\setminus \{x_0\})}]=\kappa.$

{{We still need to define $\phi: PM_{\infty}(C(X))P\to QM_{\infty}(C(X))Q$.}}
\Wlog, we may assume that $P\in M_m(C(X))$ for some integer $m\ge 1.$ {{Then $Q$ (with $[Q]=\af_0([P])$ and ${\rm rank} (Q)\geq 12$) can be regarded as  a subprojection of $q\otimes 1_m$, since $[q\otimes 1_m]=\af_0([1_{M_m(C(X))}])\geq \af_0([P])$. Hence we may assume that $Q\in M_m(q(M_{\infty}(C(Y)))q)$. Furthermore, the homomorphism $\psi\otimes {\rm id}_m: M_m(C(X))\to M_m(q(M_{\infty}(C(Y)))q)$ satisfies  $[\psi\otimes {\rm id}_m(P)]=\af_0([P])=[Q].$}}
By 9.15 of \cite{Hus},  there is  a unitary $u\in {{M_m(q(M_{\infty}(C(Y)))q)}}$ such that $u^*{{(\psi\otimes {\rm id}_m)(P)}}u=Q.$
Define $\phi={\rm Ad}\, {{u}} \circ {{(\psi\otimes {\rm id}_m)|_{PM_m(C(X))P}}}$ (see \ref{dInn}).  One then checks that $\phi$ meets all requirements
including the  injectivity when $\af_0$ is  at least 13-large.}}


{{The second  part follows from the first part by considering each partial map $\af_0^{i,j}: K_0(B_i)\to K_0(D_j)$
and $\af_1^{i,j}: K_1(B_i)\to K_1(D_j)$ separately. Note that, since at least one of $Y_j$ is not a single point, say
$Y_{j_0}$ is not a single point, then, each partial map $\phi_{i,j_0}: B_i\to D_{j_0}$ can be chosen injective for each $i.$ It follows that $\phi$ is injective.}}

\end{proof}
The following lemma is elementary.

\begin{lem}\label{range 0.15a} Let $0\to {{D}}\to H \to H/{{D}}\to 0$ be a
short exact sequence of countable abelian groups with $H/{{D}}$ torsion
free{{, and}} let

 \begin{displaymath}
    \xymatrix{
        {{H'_1}}\ar[r]^{\gm'_{1,2}} & {{H'_2}}\ar[r]^{\gm'_{2,3}}&
        {{H'_3}}\ar[r]^{\gm'_{3,4}} &\cdots \ar[r]& H/{{D}} }
\end{displaymath}
be an inductive system with limit $H/{{D}}$ such that each $H^{{'}}_i$ is
{{a}} finitely generated free abelian group.  Then there are an increasing
sequence of finitely generated subgroups ${{D}}_1\subset
{{D}}_2\subset\cd\subset {{D}}_n \subset \cd \subset {{D}}$ with
${{D}}=\bigcup_{i=1}^{\infty}{{D}}_i$, and an inductive system
\begin{displaymath}
    \xymatrix{
       {{D}}_1 \oplus {{H'_1}}\ar[r]^{\gm_{1,2}} & {{D}}_2 \oplus {{H'_2}} \ar[r]^{\gm_{2,3}}&{{D}}_3 \oplus {{H'_3}}
         \ar[r]^{\gm_{3,4}} &\cdots \ar[r]& H }
\end{displaymath}
with limit $H$ {{satisfying}} the following conditions:

{\rm (i)} $\gm_{n,n+1} ({{D}}_n) \subset {{D}}_{n+1}$ and $\gm_{n,n+1}|_{{D}_n}$ is
the inclusion from ${{D}}_n$ to ${{D}}_{n+1}$.

{\rm (ii)} If $\pi_{n+1}: {{D}}_{n+1}\oplus {{H'_{n+1}}}\to {{H'_{n+1}}}$ is the
canonical projection, then $\pi_{n+1}\circ
\gm_{n,n+1}|_{{H'_{n}}}=\gm'_{n,n+1}.$

(Here, we do not assume $\gm'_{n,n+1}$ to be injective.)

\end{lem}

\begin{proof} {{Since $D$ is countable, one {\blue{can}} list all the elements of $D$ as}} ${{D}}=\{e_i\}_{i=1}^{\infty}$. We will
construct the system $({{D}}_n\oplus H_{n},\gm_{n, n+1})$, inductively.
Let us assume {{that}} we already have ${{D}}_n \subset {{D}}$ with
$\{e_1,e_2,\cd,e_n\}\subset {{D}}_n$ and a map $\gm_{n,\infty}:
{{D}}_n\oplus {{H'_{n}}}\to H$ such that $\gm_{n,\infty}|_{{{D}}_n}$ is the
inclusion and $\pi\circ \gm_{n, \infty}|_{{H'_{n}}} =\gm'_{n, \infty}$,
where $\pi: H \to H/{{D}}$ is the quotient map.
{{To begin  the induction process we  must start with $n=0,$  and $D_0=H_0'=\{0\}.$}}
Note that, since
$\gm'_{n+1,\infty}({{H'_{n+1}}})$ is a finitely generated free abelian
subgroup of $H/{{D}}$, one has a lifting map {$\gm_{n+1, \infty}: {{H'_{n+1}}} \to H$ such
that $\pi\circ \gm_{n+1, \infty} =\gm'_{n+1, \infty}$.}
For each
$h\in {{H'_{n}}}$, we have  $\gm (h)
:=\gm_{n,\infty}(h)-\gm_{n+1,\infty}(\gm'_{n,n+1}(h))\in {{D}}$. Let
${{D}}_{n+1}\subset {{D}}$ be
{the}
finitely generated subgroup generated {by}
${{D}}_n\cup\{e_{n+1}\}\cup \gm({{H'_{n}}})$ and extend the map $\gm_{n+1,
\infty} $ on ${{D}}_{n+1}\oplus {{H'_{n+1}}}$ by defining it to be inclusion
on ${{D}}_{n+1}$. And finally let $\gm_{n,n+1}: {{D}}_n\oplus {{H'_{n}}} \to
{{D}}_{n+1}\oplus {{H'_{n+1}}}$ be defined by $\gm_{n,n+1}(e,
h)=(e+\gm(h),\gm'_{n,n+1}(h))\in {{D}}_{n+1}\oplus {{H'_{n+1}}}$ for each
$(e,h) \in {{D}}_n\oplus {{H'_{n}}} $. Evidently, $\gm_{n, \infty} =\gm_{n+1,
\infty} \circ \gm_{n,n+1}$.
\end{proof}

\begin{NN}\label{range 0.16}

{\blue{We now fix {{an object}}
 $((G,G_+,u),K,\Delta,r)$  as described in}}
\ref{range  0.1}, {{which is to be shown to be in the range of Elliott invariant.}}   In general, $G$ may have torsion and $K$ may not be zero. Let $G^1\subset \Aff(\DT)$ be the  dense subgroup
with at least three $\Q$-linearly independent elements and
let $H=G\oplus G^1$ be as in \ref{range 0.3}. The order unit $u\in G_+$ {{when}}  regarded as $(u,0)\in G\oplus G^1=H$
{{is also an}} order
 unit of $H_+$. Note that
 ${\rm Tor}(H)={\rm Tor}(G).$
We have {{the}}  split short exact sequence
$$0\lr G\lr H\lr H/G\,(=G^{{1}})\lr 0$$
with $H/{\rm Tor}(H)$ a dimension group (see \ref{range 0.3}). {{Since $H/{\rm Tor}(H)$ is a {\blue{simple}} dimension group, as in \ref{range 0.4}, {\blue{
we may}}  write $H/{\rm Tor}(H)$ as  an inductive limit
\begin{displaymath}
    \xymatrix{
        H'_1 \ar[r]^{\gm'_{1,2}} & H'_2 \ar[r]^{\gm'_{2,3}}&H'_3
         \ar[r]^{\gm'_{3,4}} &\cdots \ar[r]& H/ {\rm Tor}(H),}
\end{displaymath}
where $H'_n=\Z^{p_n}$ with standard positive cone $(H'_n)_+=(\Z_+)^{p_n}$.}}
{\blue{Moreover, we may assume that $\gamma_{n,n+1}'$ is $2$-large.}}
Applying Lemma \ref{range 0.15a}
to the short exact sequence $0 \to {\rm Tor}(H) \to H {\blue{{\stackrel{\pi}{\longrightarrow}}}} H/{\rm Tor}(H) \to 0,$ {{with $D={\rm Tor}(H)$}},
{\blue{we may}} write $H$ as {{an}} inductive
limit of finitely generated abelian {{groups,}}
\begin{displaymath}\label{1326-n}
    \xymatrix{
        H_1 \ar[r]^{\gm_{1,2}} & H_2 \ar[r]^{\gm_{2,3}}&H_3 \ar[r]^{\gm_{3,4}} &\cd \ar[r] & H,}
\end{displaymath}
where {{$H_n=D_n\oplus H'_n= D_n\oplus \Z^{p_n},$}}  {\blue{and (i) and (ii) of \ref{range 0.15a} hold.}}
{\blue{It follows that  $\gamma_{n, n+1}(d,h)=(\gamma_{D,n, n+1}(d+h), \gamma_{n, n+1}'(h)),$
where $\gamma_{D,n,n+1}: D_n\oplus H_n'\to D_{n+1}$ is a \hm.
Note that $h\in H_+\setminus \{0\}$  if and only if}} {{$\pi(h) > 0.$}}
{\blue{We may write}}
$H_n=\bigoplus_{j=1}^{p_n} H_n^j$ with $H_n^1=\Z\oplus {\rm Tor}(H_n),$ {\blue{and}}
$H_n^i=\Z$ for all $i\geq 2.$
Define
$$(H_n)_+=\big((\Z_+\setminus \{0\})\oplus {\rm Tor}(H_n)\cup \{{\blue{(0,0)}}\}\big)\oplus \Z_+^{p_n-1}.$$
{\blue{Since $\gamma_{n, n+1}'$ is positive  and $2$-large, one checks easily that each $\gamma_{n, n+1}$
is also positive with respect to the positive cones $(H_n)_+$ {{and $(H_{n+1})_+.$}}
Let $H_{+'}=\bigcup_{n=1}^{\infty} \gamma_{n, \infty}((H_n)_+).$
If $x\in H_+\setminus \{0\},$ then, as mentioned above, $\pi(x)>0.$
Suppose that $y=(y_{_D}, y_{_H})\in H_n=D_n\oplus H_n'$ is such that $\gamma_{n, \infty}(y))=x.$
Then $\gamma_{n, m}(y)=(y_{_D}', \gamma_{n,m}'(y_{_H})),$ where $\gamma_{n,m}'=\gamma_{m-1, m}'\circ{{\gamma_{m-2, m-1}'\circ\cdots \circ}} \gamma_{n,n+1}'$
and $m>n.$ It follows that, for some large $m,$ $\gamma_{n,m}'(y_{_H})\in (H_n')_+\setminus \{0\}.$
Since each $\gamma_{n,n+1}'$ is $2$-large, there is $m_1>m$ such
that $\gamma_{m, m_1}'(\gamma_{n,m}(y_{_H}))\in \{(z_1,z_2,...,z_{p_{m_1}})\in \Z^{p_{m_1}}_+\andeqn z_1>0\}.$
In other words, $y:=\gamma_{n, m_1}(y)\in (H_n)_+\setminus \{0\}.$ Since $\gamma_{m_1, \infty}(y)=x,$
this implies $x\in H_{+'}.$}}

{\blue{Conversely, if $x\in H_{+'}\setminus \{0\},$  then there {{exist}} $n>0$ and $y_1\in  (H_n)_+$ such
that $\gamma_{n, \infty}(y_1)=x.$ It follows that $y=(z_1, \eta)\oplus (z_2, z_3,...,z_{p_n})$
with $z_1\in \Z_+\setminus \{0\},$ $\eta\in D_n,$ and $(z_2, z_3,...,z_{p_n})\in {{\Z_+^{p_n-1}}}.$
{{By}} {{ construction}}, $(z_1, z_2,..., z_{p_n})\in (H_n')_+\setminus \{0\}.$
It follows {{that}} $\gamma'_{n, \infty}((z_1,z_2,...,z_{p_n}))\in H_{+'}\setminus \{0\}.$ By  (ii) of  \ref{range 0.15a},
this implies that $\pi(x)>0.$ Therefore $x\in H_+.$
This shows that $H_+=H_{+'}.$}}

{\blue{{{The notation}} $H$ and $(H_n)_+$ above will be used later.}}

{\blue{Let $\gamma_{n,n+1}': H_n'\to H_{n+1}'$ be determined by the  $p_n\times p_{n+1}$ matrix of positive intergers $\cc_{n, n+1}$ which
we assume  at least $2$-large.
We now represent  $\gamma_{n,n+1}$}}  by a  $p_{n+1}\times p_n$ matrix of {\blue{homomorphisms}} $\td \cc_{{n, n+1}}=(\td c^{{n, n+1}}_{ij})$,
where $\td c^{{n, n+1}}_{ij}:~ H_n^j\to H_{n+1}^i$ {{are described  as follows.}}
{{Note $\td c^{{n, n+1}}_{ij}=P_i\circ \gamma_{n,n+1}|_{H_n^j},$ where $P_i: H_{n+1}\to H_{n+1}^i.$}}

 If $i>1,~ j{{>}} 1$, then
{\blue{define}} $\td c^{{n, n+1}}_{ij}=c^{{n, n+1}}_{ij}$
{\blue{(recall  $H_n^j=\Z$ and $H_n^i=\Z,$ and
$c^{n,n+1}_{ij}$ maps $m$ to $c^{n,n+1}_{ij}m$).}}

 {{If $i>1, j=1$,}}  {\blue{define $\td c^{n,n+1}_{i1}: H_n^1=\Z\oplus {\rm Tor}(H_n)\to H_{n+1}^i=\Z$
 by $\td c^{n,n+1}_{i1}(z,t)=c^{n,n+1}_{i1}z$ for all $(z, t)\in \Z\oplus {\rm Tor}(H_n)${{, and}}
{{ let us still
 denote  this}}  by $c^{n,n+1}_{i1}.$}}

 {{If $i=1$}},  let $Q_1\hspace{-0.05in}:\hspace{-0.02in} H_{n+1}^1=\Z\oplus {\rm Tor}(H_{n+1})\to  {\rm Tor}(H_{n+1})\subset H_{n+1}^1$ be a projection.
 Define $T^{n,n+1}_j=Q_1\circ P_1\circ \gamma_{n,n+1}|_{H_n^j}.$
 Viewing $c^{n,n+1}_{1j}$ as a \hm\, from
 $\Z$ to $\Z$ as mentioned above (see \ref{notation-Dec12-2018}),
  {\blue{we define}} ${\blue{\td{c}^{n, n+1}_{1j}}}=c^{{n, n+1}}_{1j}+T^{{n,n+1}}_j$ as follows.
If $j\not=1,$ $\td c^{{n, n+1}}_{ij} (m)=c^{{n, n+1}}_{1j}m+T^{n,n+1}_j(m)$ for $m\in H_n^j=\Z;$
if $j=1,$ $\td c^{{n, n+1}}_{ij} ((m,t))=c^{{n, n+1}}_{1j}m+T^{n,n+1}_j((m,t))$ for all $(m,t)\in \Z\oplus {\rm Tor}(H_n).$



Since
$\gm_{n,n+1}$ satisfies $\gm_{n,n+1}({\rm Tor}(H_n))\subset
{\rm Tor}(H_{n+1})$, it induces the  map
$$\gm'_{n,n+1}:~ H_n/{\rm Tor}(H_n)=\Z^{p_n}\lr H_{n+1}/{\rm Tor}(H_{n+1})=\Z^{p_{n+1}}.$$
{{(By (ii) of Lemma \ref{range 0.15a}, this {{map}}
is the same as {{the map}} $\gm'_{n,n+1}: H'_n~(=H_n/{\rm Tor}(H_n))\to H'_{n+1}~(=H_{n+1}/{\rm Tor}(H_{n+1}))$.)}}
{{Thus,}} $\gm'_{n,n+1}$ is given by the matrix $\cc_{{n,n+1}}=(c^{{n, n+1}}_{ij})$ {{with}} {{positive integer entries.}}
{\blue{Put $G_k'=G_k/{\rm Tor}(G_k)$ ($k=1,2,...$). Then $\gamma_{n,n+1}'(G_n')\subset G_{n+1}'.$
Note that,  passing to a subsequence,
we may always {{assume}}
$c^{{n, n+1}}_{ij}>2.$}}
\end{NN}

\begin{NN}\label{range 0.17}

Let $G_n=H_n\cap \gm^{-1}_{n,\infty}(G)$ with $(G_n)_+=(H_n)_+\cap
G_n$.  Then the order unit\\  ${\blue{u_n= (([n, 1], \tau_{n}), [n,2], ..., [n, p_n] )}}$ {{of}}  $(H_n)_+$ is also an order unit for $(G_n)_+.$


Recall  that ${\rm Tor}(G)={\rm Tor}(H)$, {{which implies that
$\gm_{n,\infty}( {\rm Tor}(H_n))\subset {\rm Tor}(H)\subset G$. {{Therefore,}}
$ {\rm Tor}(H_n)\subset H_n\cap \gm^{-1}_{n,\infty}(G) =G_n$.  Since $G_n$ is a subgroup of $H,$}} ${\rm Tor}(H_n)={\rm Tor}(G_n).$ Furthermore, we have the following commutative diagram:

\begin{displaymath}{\small{
    \xymatrix{0\ar[d] &0\ar[d]&&0\ar[d]\\
        G_1 \ar[r]^{\gm_{12}|_{_{G_1}}} \ar[d] & G_2 \ar[r]\ar[d]&\cd \ar[r]&G \ar[d]\\
         H_1 \ar[r]^{\gm_{12}} \ar[d] & H_2 \ar[r]\ar[d]&\cd \ar[r]&H \ar[d]\\
         H_1/G_1 \ar[r]^{\td\gm_{12}} \ar[d] & H_2/G_2 \ar[r]\ar[d]&\cd \ar[r]&H/G \ar[d]\\
        0 &0&&0 }}}
\end{displaymath}
where $\td\gm_{n,n+1}$ is induced by $\gm_{n,n+1}$.  Note that
the inductive limit of the quotient groups, $H_1/G_1\to H_2/G_2\to \cd
\to H/G,$ has no
{{obvious}}  order structure.
{\blue{Since ${\rm Tor}(H_n)={\rm Tor}(G_n),$ the quotient map $H_n\to H_n/G_n$
induces
a surjective map $\pi_n: H_n/{\rm Tor}(H_n)=\Z^{p_n}\to H_n/G_n=\Z^{l_n},$
which will  be used in  the our construction later. {{We will reserve the notation $\pi_n$ for
this map from $\Z^{p_n}$ to $\Z^{l_n}$. The map $\pi_n: H_n\to H'_n=H_n/{\rm Tor}(H_n)$ {{which}} appeared in \ref{range 0.15a} (see $\pi_{n+1}$ there) will be denoted by $\pi_{H_n,H'_n}$ from now on.}}
Note also $u_n'={{\pi_{H_n,H'_n}(u_n)=}}([n, 1], [n,2], ..., [n, p_n] ).$  Denote by $\pi_{G,H/{\rm Tor}(H)}$ the {{composed}}  map from $G$ to $H$ and then to $H/{\rm Tor}(H).$ }}

 Also write the group {{$K$ of
\ref{range 0.16}}} as {{an}} inductive limit,
\begin{displaymath}
    \xymatrix{
        K_1 \ar[r]^{\chi_{12}}  & K_2 \ar[r]^{\chi_{23}}&K_3 \ar[r]^{\chi_{34}}&\cd \ar[r]& K ~,}
\end{displaymath}
where each $K_n$ is finitely generated.\\
\end{NN}

\begin{NN}\label{range 0.18}
~~Recall from \cite{EG-RR0AH}  that the finite CW {{complex}} $T_{2,k}$  (
$T_{3,k}${{, respectively}}) is defined to be a 2-dimensional connected finite CW
complex with $H^2(T_{2,k})=\Z/k$ and
 $H^1(T_{2,k})=0$ (3-dimensional finite CW complex with $H^3(T_{3,k})=\Z/k$
 and $H^1(T_{3,k})=0=H^2(T_{3,k})${{, respectively}}). (In {\cite{EG-RR0AH}} the spaces are denoted
 by $T_{II,k}$ and $T_{III,k}$.)  For each $n$, {{write $$H^1_n=\Z\oplus {\rm Tor}(H_n):=\Z\oplus (\Z/k_1\Z)\oplus (\Z/k_2\Z)\oplus\cdots\oplus (\Z/k_i\Z)~~ \mbox{ and}$$
  $$K_n:=
  {\blue{\Z^l}}\oplus(\Z/m_1\Z)\oplus (\Z/m_2\Z)\oplus\cdots\oplus (\Z/m_j\Z).$$ {\blue{Set}}}}
$$X_n'=\overbrace{S^1\vee S^1\vee\cd \vee S^1}^{{l}}\vee T_{2,k_1} \vee T_{2,k_2} \vee \cd \vee T_{2,k_i} \vee T_{3,m_1} \vee T_{3,m_2}\vee \cd \vee T_{3,m_j}{{.}}$$
{{Then}} $K_0(C(X_{{n}}'))=H^1_n=\Z\oplus {\rm Tor}(H_n)$ and $K_1(C(X_{{n}}'))=K_n$. Let $x^{{1}}_n$ be the base point of $X_n'$ which is {{the}} common point of all {{copies of the}} spaces $S^1$, $T_{2,k}$, $T_{3,k}$ appearing above in the wedge $\vee$.
{\blue{By \ref{1323proj},}} there is a projection $P_n\in M_\infty (C(X{{'}}_n))$ such that
$$[P_n]=([n,1],\tau_n)\in K_0(C(X^{{'}}_n))=\Z\oplus {\rm Tor}(K_0(C(X^{{'}}_n))),$$
where $([n,1],\tau_n)$ is the first component of the unit $u_n\in H_n.$

{{ Suppose that  $P_n\in M_N(C({\blue{X_n'}})),$
where $N$ is a large enough integer. Since}}  ${\rm rank}(P_n)=[n,1]$, {{ there is a unitary $u\in M_N$ such that $uP_n(x^{{1}}_n)u^*=\diag(\e_{[n,1]},\underbrace{0,\cdots, 0}_{N-[n,1]})$. Replacing $P_n$ by $uP_nu^*$, }}
{{we may}}
assume  that $P_n(x^{{1}}_n)=\e_{M_{[n,1]}}$,
where $M_{[n,1]}$ is identified with the upper left corner of
$M_{{\blue{N}}}$. Define $X_n=[0,1]\vee X_n'$ with $1\in [0,1]$
identified with the base point $x^{{1}}_n\in X_n'$.
{{Let us label the point}}
$0\in [0,1]$
by the symbol {{$x^0_n$}}. So $[0,1]$ is identified with $[{{x^0_n}},{{x^1_n}}]${{---the convex combinations of $x^0_n$ and $x^1_n$}}. Under this identification,   we have
$X_n=[{{x^0_n}},{{x^1_n}}]\vee X_n'$. {{ It is convenient to write the point $(1-t)x^0_n+ t x^1_n\in [x^0_n, x^1_n]$ as $x^0_n+t$ for any $t\in [0,1]$. In particular, we have $x^1_n=x^0_n+1$.}}
{\blue{The projection}} $P_n\in {\blue{M_N(C(X_n'))}}$ can
be extended to a projection {\blue{in $M_N(C(X_n)),$}} still called $P_n,$
by $P_n({{x^0_n}}+t)=\e_{{{M_{[n,1]}}}}$ for each
$t\in(0,1)$.
{{Let us  choose}} ${{x^0_n}}$ {{as}} the base point of $X_n$. The (old) base point of $X_n'$ is ${{x^0_n}}+1{{=x^1_n}}$.

\end{NN}

\begin{NN}\label{range 0.19}
~~  {\blue{Note that the}} map
$${\blue{\pi_n:}} H_n/{\rm Tor}(H_n)~ (=\Z^{p_n})\lr H_n/G_n~ (=\Z^{l_n})$$
{\blue{may be written as}}
${\blue{\pi_n=\bb_{n,0}-\bb_{n,1}}},$ {\blue{where
$\bb_{{{n}},0},\bb_{{{n}},1}:~ \Z^{p_n}\to \Z^{l_n},$
are two order preserving \hm s.  Note also that we may assume that}}
 $\bb_{{{n}},0}=(b^{{n}}_{0,ij})$ and $\bb_{{{n}},0}=(b^{{n}}_{1,ij})$  are
two $l_n\times p_n$ matrices of  non-negative integers {{which satisfy the conditions of \ref{range 0.5a} and \ref{conditions}.
}}
%

 Exactly as in \ref{range 0.5} (in which we considered the special
case of torsion free $K_0$ and trivial $K_1$), we can define
\beq\label{150102-1}
\{n,i\}:= \sum_{j=1}^{p_n}b^{{n}}_{0,ij}[n,j]=\sum_{j=1}^{p_n}b^{{n}}_{1,ij}[n,j]~.
\eneq

{\blue{As in}} {{{{\ref{range 0.5a}}}, set $F_n=\bigoplus_{i=1}^{p_n} F^i_n=\bigoplus_{i=1}^{p_n} M_{[n,i]},$}}
$E_n=\bigoplus_{i=1}^{l_n}M_{\{n,j\}},$ and let $\bt_{{{n}},0},\bt_{{{n}},1}:~
{{F}}_n\to E_n$ be  homomorphisms
such that $(\bt_{{{n}},0})_{*0}=\bb_{{{n}},0}$
and
$(\bt_{{{n}},1})_{*0}=\bb_{{{n}},1}$.
{\blue{Since $P_n(x^0_n)$ has rank $[n,1],$ there is an isomorphism $j_n: F_n^1=M_{[n,1]}\to (P_nM_{N}(C(X))P_n)|_{x^0_n}.$ }}
{{Let $\beta_{X_n}: F_n\to (P_nM_{\blue{N}}(C(X))P_n)|_{x^0_n}$ be defined by $\beta_{X_n}(a_1, a_2, \cdots, a_{p_n})={\blue{j_n(a_1)}}\in (P_nM_{\blue{N}}(C(X))P_n)|_{x^0_n}.$}}

Now {{consider  the algebra}}
\beq\nonumber A_n:=\big\{(f,g, {{a}})\in C([0,1],E_n)\oplus {{P_n M_{\blue{N}}(C(X_n))P_n}}\oplus F_n;~\\ \label{13Dec15-2018} f(0)=\bt_{{{n}},0}({{a}}),
f(1)=\bt_{{{n}},1}({{a}}), {{g(x^0_n)=\bt_{X_n}(a)}}\big\}.\eneq
{{We may also write $A_n=\big(C([0,1],E_n)\oplus P_n M_{\blue{N}}(C(X_n))P_n\big)\oplus_{\bt_{n,0},\bt_{n,1},\bt_{X_n}} F_n.$
{{Denote by}} $\pi_{A_n}^e: A_n\to F_n$ the quotient map.}}


{{ As in
the construction of $C_n$ (see \ref{conditions}), {{once $E_n=\bigoplus_{i=1}^{l_n}M_{\{n,j\}}, F_n=\oplus_{i=1}^{p_n} M_{[n,i]},$
and
 $P_n M_{\blue{N}}(C(X_n))P_n$ with $P(x^0_n)=\diag(\one_{[n,1]},0,...,0)\in M_{N}$ are fixed, }} the algebra $A_n$ depends only on  $\bb_{{{n}},0},\bb_{{{n}},1}: K_0({{ F}}_n)=\Z^{p_n} \to \Z^{l_n}$ {{up to isomorphism}}.
{{Thus,}} once ${{E_n,}} F_n, \bb_{{{n}},0},
\bb_{{{n}},1},$  {{ and   $ P_n M_{\blue{N}}(C(X_n))P_n$}}
{{are specified,}}
the construction
 of the algebra $A_n$ is  {{complete.}}}}

Note that {{by}}  (\ref{150102-1}), $\bt_{{{n}},0}$ {{and }} $\bt_{{{n}},1}$ are
unital homomorphisms and therefore $A_n$ is a unital algebra. Note
that this algebra, in general,  is {{$not$}} the  direct sum of a homogeneous
algebra and an algebra in ${\cal C}_0.$
{{It is easy to verify that}} $A_n\in {\cal D}_{{1}}$ {{  in the sense of Definition \ref{8-N-3}, by writing $A_n=PC({\blue{X_n}}, F)P\oplus_{\Gamma} B $
with {\blue{$B=F_n,$ $X=[0,1]
\sqcup X_n$, $Z=\{0,1,x_n^0\}\subset X,$}}
{\blue{and $F=E_n\oplus M_N$,}}  
{\blue{with}} $P$ defined by $P=\one_{E_n}$ on $[0,1],$ and $P=P_n$ on $X_n$, and, finally, $\Gamma:
F\to P|_ZC(Z,F)P|_Z $   defined by $\bt_{n,0}$ for the point $0\in Z\subset X$, $\bt_{n,1}$ for the point $1\in Z\subset X,$ and $\bt_{X_n}$ for the point $x^0_n\in Z\subset X$.}}
Later, we will
deal with a nicer special case {{in which}}  $A_n$ is {{in fact}} a direct sum of a  homogeneous \CA\,  and
a \CA\, in ${\cal C}_0.$
\end{NN}


\begin{NN}\label{construction}
Let $ {{J}}_n=\{(f, g, {{a}})\in A_n: {{f=0,~a=0}}\}$
and $I_n=\{(f,g, {{a}})\in A_n: g=0,~ {{a=0}}\}.$
Denote the quotient
algebra $A_n/ {{J}}_n$ by $ A_{{C},n},$ and $A_n/I_n$ by $A_{X,n}.$
{\blue{Let $\pi_{J,n}: A_n\to A_n/J_n$ and $\pi_{I,n}: A_n\to A_n/I_n$ denote the quotient maps.  Set ${\bar I_n}=\pi_{J,n}(I_n)$ and ${\bar J_n}=\pi_{I,n}(J_n).$}} {{Since $I_n\cap J_n=\{0\}$, the map $\pi_{J,n}$ is injective on $I_n$, whence ${\bar I_n}\cong I_n$. {{Similarly,}} ${\bar J_n}\cong J_n$.}}
Note that
\beq\nonumber
\hspace{-0.05in}
{A_{{C},n}}&=&{{C([0,1],E_n)\oplus_{\bt_{n,0}, \bt_{n,1}} F_n}}
\\ \label{construction-1}
&=&\big\{(f,{{a}})\in C([0,1],E_n)\oplus {{ F}}_n;~ f(0)=\bt_{{n,0}}({{a}})~f(1)=\bt_{{n,1}}({{a}})\big\},\andeqn
\hspace{-0.1in}\\
A_{X,n}&=&P_nM_{\blue{N}}(C(X_n))P_n\oplus_{\bt_{X_n}}F_n.
\eneq
{\blue{Denote by $\pi_{C,n}^e: A_{C,n}\to F_n$ and $\pi_{X,n}^e: A_{X,n}\to F_n$ the quotient maps, and let
$\lambda_{C,n}: A_{C,n}\to C([0,1], E_n)$
be the \hm\,    defined in  \eqref{pull-back-k} (as $\lambda_k$  with $A_{k-1}=F_n$).
Then ${\bar I}_n=C_0((0,1),E_n)={\rm ker}\pi_{C,n}^e$ and ${\bar J}_n=\{g\in P_nM_N(C(X_n))P_n : g(\xi_n^0)=0\}={\rm ker}\pi_{X,n}^e.$  Note that $\pi_{A_n}^e=\pi_{C,n}^e\circ \pi_{J,n}=\pi_{X, n}^e\circ \pi_{I,n}.$}}

{{Evidently, we can write
\beq\label{13-March26-2019}
Sp(I_n)=\bigcup_{j=1}^{l_n} (0,1)_{n,j}, ~~Sp(J_n)=X_n\setminus\{x_n^0\}~\mbox{and} ~Sp(F_n)=\{\tht_{n,1},\tht_{n,2},\cdots, \tht_{n,p_n}\}.
\eneq
 We also have:
\beq\label{13-March26-2019-1}
Sp(A_{C,n})=Sp(I_n){{\sqcup}} Sp(F_n),~Sp(A_{X,n})=Sp(J_n){{\sqcup}} Sp(F_n)\, \mbox{{{(disjoint union);}}}\eneq
\beq\label{13-March26-2019-2}
Sp(A_n)=Sp(I_n)\sqcup Sp(A_{X,n})=Sp(J_n)\sqcup Sp(A_{C,n})\,\mbox{(disjoint union); and}
\eneq
\beq\label{13-March26-2019-3}
~Sp(A_n)=Sp(A_{C,n})\cup Sp(A_{X,n})~\mbox{with}~Sp(A_{C,n})\cap Sp(A_{X,n})=Sp(F_n).
\eneq
We can also see  that $Sp(A_{X,n})=X_n\cup Sp(F_n),$ with $\tht_{n,1}\in Sp(F_n)$ identified with $x_n^0\in X_n$. Also note that for any $\tht \in Sp(A_n)$, {{by}}  (\ref{13-March26-2019-2}), we have
\beq\label{13-March26-2019-4}
\tht\in Sp(A_{X,n})~~\mbox{if and only if}~~ \tht|_{I_n}=0.
\eneq
}}

The homomorphism $\phi_{n,n+1}: A_n\to A_{n+1}$ to be constructed should
 satisfy the conditions $\phi_{n,n+1}(I_n)\subset I_{n+1}$ and
 $\phi_{n,n+1}( {{J}}_n)\subset {{J}}_{n+1}$, and therefore {{should induce}} two homomorphism{{s}} $\psi_{n,n+1}:A_n/I_n={{A_{X,n}}}\to A_{n+1}/I_{n+1}={{
 A_{X,n+1}}}$ and
  $\bar \phi_{n,n+1}: A_n/{{J}}_n=A_{{C},n}\to A_{n+1}/{{J}}_{n+1}=A_{{C},n+1}$.
 Conversely, if two  homomorphisms $\psi_{n,n+1}: {{A_{X,n}}}\to {{A_{X,n+1}}}$
 (necessarily injective) and $\bar \phi_{n,n+1}: A_{{C},n}\to A_{{C},n+1}$ satisfy the {{two}} {{conditions}}

(a) $\psi_{n,n+1}({{\bar J}}_n)\subset {{\bar J}}_{n+1}$ and $\bar \phi_{n,n+1} (\bar I_n)\subset \bar I_{n+1}$, and

(b) the homomorphism ${{\psi^{q}_{n,n+1}:}}{{A_{X,n}}}/{{\bar J}}_n= {{F}}_n \to {{A_{X,n+1}}}/{{\bar J}}_{n+1}= {{F}}_{n+1} $ induced by $\psi_{n,n+1}$
and the homomorphism ${{\phi^{q}_{n,n+1}:}}A_{{C},n}/\bar I_n={{ F}}_n \to
    A_{{C},n+1}/\bar I_{n+1}={{ F}}_{n+1}$ induced by $\bar \phi_{n,n+1}$ are the same,

\noindent
then, there is a unique  {\blue{(necessarily injective)}} homomorphism $\phi_{n,n+1}: A_n\to A_{n+1}$ satisfying

(c)\,$\phi_{n,n+1}(I_n)\subset I_{n+1}$, $\phi_{n,n+1}( {{J}}_n)\subset {{J}}_{n+1}$, and $\phi_{n,n+1}$
 {\blue{induces}} the homomorphisms  $\psi_{n,n+1}$ and $\bar \phi_{n,n+1}$.

{\blue{To see this, 
define $\phi_{n, n+1}: A_n\to C([0,1], E_{n+1})\oplus P_{n+1}M_{\infty}(C(X_{n+1}))P_{n+1}\oplus F_{n+1}$ by
\beq\label{dofphi}
\hspace{-0.4in}\phi_{n, n+1}(x)=(\lambda_{C,n+1}( {{{\bar \phi}}}_{n,n+1}(\pi_{J,n}(x))),  \lambda_{X,n+1}(\psi_{n, n+1}(\pi_{I,n}(x))), \psi_{n,n+1}^q(\pi_{A_n}^e(x)))
\eneq
for all $x\in A_n.$ It is a unital \hm.  Since $\psi^q_{n,n+1}$ is induced by ${{\psi}}_{n, n+1},$
one sees that
\beq\label{1326nl-1}
\hspace{-0.4in} \lambda_{X,n+1}(\psi_{n, n+1}(\pi_{I,n}(x)))(x_{n+1}^0)=\bt_{X, n+1}(\psi_{n,n+1}^q(\pi_{X,n}^e(\pi_{I,n}(x)))
=
\bt_{X_{n+1}}(\psi_{n,n+1}^q(\pi_{A_n}^e(x)))
\eneq
for all $x\in A_n.$}} {\blue{On the other hand, since $\phi_{n, n+1}^q$ is induced by ${{\bar \phi}}_{n, n+1},$ also, by (b),  one has
\beq
\pi_{t_i}(\lambda_{C,n+1}(
{{\bar \phi}}_{n,n+1}(\pi_{J,n}(x))))=\bt_{n+1,i}(\phi_{n,n+1}^q(\pi_{C,n}^e(\pi_{J,n}(x))))\\\label{1326nl-2}
=\bt_{n+1,i}(\phi_{n,n+1}^q(\pi_{A_n}^e(x)))=\bt_{n+1,i}(\psi_{n,n+1}^q(\pi_{A_n}^e(x))),\,\, i=0,1,
\eneq
for all $x\in A_n,$   where $\pi_{t_i}: C([0,1], E_n)\to E_n$ is the point evaluation at $t_i,$ where $t_0=0$ and $t_1=1.$}}
{\blue{By \eqref{1326nl-1} and \eqref{1326nl-2}, one has $\phi_{n,n+1}(x)\in A_{n+1}$
for all $x\in A_n.$ Thus,  $\phi_{n,n+1}$ is a
 \hm\, from $A_n$ into $A_{n+1}.$   {{By}} definition, $\phi_{n,n+1}$ induces ${{\bar \phi}}_{n,n+1}$ as well as
$\psi_{n, n+1}.$  There is only one such map.  Note also, if both $\psi_{n, n+1}$ and ${\bar \phi}_{n,n+1}$ are
injective, then $\phi_{n,n+1}$ defined in \eqref{dofphi} is also injective, as $J_n\cap I_n=\{0\}.$}}

Note that $A_{{C},n}$ is the same as $C_n$ in \ref{range 0.5}--\ref{range 0.12},
with $ {{F}}_n$ and $\bar I_n={\blue{\pi_{J,n}}}(I_n)$ in place of
 $F_n$ and $I_n$ in \ref{range 0.5}--\ref{range 0.12}. Therefore the
  construction of $\bar \phi_{n,n+1}$ can be carried out as
  in \ref{range 0.5}--\ref{range 0.12},  with the map
  ${{F}}_n \to {{F}}_{n+1}$ being given by the matrix
  $\cc_{{n,n+1}}$ as in \ref{range 0.16} above---of course,
   we need to assume that the corresponding maps $\bb_{{{n}},0}$
    and $\bb_{{{n}},1}$ in this case (see \ref{range 0.19}) satisfy
     $\spdd$ {{(in place of $\bb_0$ and $\bb_1$)}}. So, in what
      follows, we will focus on the construction of $\psi_{n,n+1}$. But
      before the construction,  {{let us {{introduce}} the following notation.}}

 {{Note that $A_{X,n}=(P_nM_{\infty}(C(X_n))P_n\oplus_{(\bt_{X_n})|_{F_n^1}}F_n^1)\oplus \bigoplus_{i\ge 2}F_n^i.$
 Since $F_n^1=M_{[n,1]}$ and $P_n$ has rank $[n,1]$ (see \ref{range 0.19}),
 one verifies that  ${{P_nM_{\infty}(C(X_n))P_n\oplus_{(\bt_{X_n})|_{F_n^1}}F_n^1}}=P_nM_{\infty}(C(X_n))P_n.$
 Hence $
 A_{X,n}=P_nM_{\infty}(C(X_n))P_n\oplus \bigoplus_{i=2}^{p_n}F^i_n$. Let us emphasize that we will write $A_{X,n}=\bigoplus_{i=1}^{p_n} A_{X,n}^i $ with $A_{X,n}^1=P_nM_{\infty}(C(X_n))P_n$ and $A_{X,n}^i=F^i_n$ for $i\geq 2$.}}
 {\blue{Note that ${\bar J_n}\subset A_{X,n}^1.$
 Note that
 we may also write
 \beq\label{Adef}
 &&\hspace{-0.3in}A_n=C([0,1], E_n)\oplus_{\bt_{n,0}\circ \pi_{X,n}, \bt_{n,1}\circ \pi_{X,n}} A_{X,n}\\\label{Adef2}
&&\hspace{-0.9in} =\{(f,g): f\in C([0,1], E_n), \, g\in A_{X,n}, f(0)=\bt_{n,0}(\pi_{X,n}(g))\andeqn f(1)=\bt_{n,1}(\pi_{X,n}(g))\}.
 \eneq}}
 {{Note that from $\phi_{n,n+1}(I_n)\subset I_{n+1}$ for each $n$, we know that for any $m>n$, $\phi_{n,m}(I_n)\subset I_m$. By (\ref{13-March26-2019-4}), for any $m>n$ and  $y\in Sp(A_{X,m}),$ we have
 \beq\label{13-March26-2019-5}
 Sp(\phi_{n,m}|_y)\subset Sp(A_{X,n}) \subset Sp(A_n).
 \eneq}}

\end{NN}

\begin{lem}\label{13-Feb17-2019}
Let $(H_n, (H_n)_+, u_n)$, $G_n\subset H_n$, $(G_n)_+=(H_n)_+\cap G_n$, {{and}} $\pi_n:H_n/{\rm Tor} (H_n)=\Z^{p_n} \to H_n/G_n=\Z^{l_n}$ be as in \ref{range 0.17}.
Let $\bt_{n,0},\bt_{n,1}: H_n/{\rm Tor}(H_n)= \Z^{p_n}\to \Z^{l_n}$ satisfy $\bt_{n,1}-\bt_{n,0}=\pi_n$ as in \ref{range 0.19}. Let $F_n={{\bigoplus}}_{i=1}^{p_n}M_{[n,i]}, E_n={{\bigoplus}}_{j=1}^{l_n}M_{\{n,j\}}, A_{X,n}^1=P_nM_N(C(X_n))P_n,$
$A_{X,n}=A_{X,n}^1\oplus_{i=2}^{p_n}M_{[n,i]},$ and
$A_n=C([0,1], E_n)\oplus_{\bt_{n,0}\circ \pi_{X,n}, \bt_{n,1}\circ \pi_{X,n}} A_{X,n}$
be as in  \eqref{Adef}.

(1)~ Suppose that $$(K_0(F_n),  K_0(F_n)_+, [1_{F_n}])=(H_n/{\rm Tor} (H_n), (H_n/{\rm Tor} (H_n))_+, ([n,1],[n,2],\cdots [n,p_n]))$$ (recall $H_n/{\rm Tor} (H_n)=\Z^{p_n}$) and
$$(K_0(A_{X,n}^1),K_0(A_{X,n}^1)_+, [1_{A_{X,n}^1}], K_1(A_{X,n}^1))=(\Z\oplus
{\rm Tor}(H_n), (\Z_+\setminus\{0\})\oplus {\rm Tor}(H_n) \cup\{(0,0)\}, ([n,1],\tau_n), K_n).$$
Then $$(K_0(A_n), K_0(A_n)_+, [1_{A_n}], K_1(A_n))=(G_n, (G_n)_+, u_n, K_n),$$
$K_0({\bar J}_n)={\rm Tor}(G_n)={\rm Tor}(H_n)$, and $K_1({\bar J}_n)=K_n$.

(2)~Suppose that ${\phi}_{n,n+1}: A_{C,n}\to A_{C, n+1}$ and $\psi_{n,n+1}: A_{X,n}\to A_{X,n+1}$ satisfy the conditions (a) and (b) of
\ref{construction}. If ${\psi_{n,n+1}}_{*0}=\gm_{n,n+1}: K_0(A_{X,n})=H_n\to K_0(A_{X,n+1})=H_{n+1}$ and ${\psi_{n,n+1}}_{*1}=\chi_{n,n+1}: K_1(A_{X,n})=H_n\to K_1(A_{X,n+1})=H_{n+1}$, then ${\phi_{n,n+1}}_{*0}=\gm_{n,n+1}|_{G_n}: K_0(A_n)=G_n\to K_0(A_{n+1})=G_{n+1}$ and ${\phi_{n,n+1}}_{*1}=\chi_{n,n+1}: K_1(A_n)=K_n\to K_1(A_{n+1})=K_{n+1}$.

\begin{proof} Part (1): As ${\bar J}_n$ is an ideal of $A_{X,n}^1$ with quotient $A_{X,n}^1|_{x_n^0}\cong M_{[n,1]}$, evidently, $K_0({\bar J}_n)={\rm Tor}(G_n)={\rm Tor}(H_n)$, and $K_1({\bar J}_n)=K_n$. Since $A_{X,n}=A_{X,n}^1\oplus \bigoplus_{i=2}^{p_n}F_n^i$, {{and this differs from $F_n$ by  replacing $F_n^1$ by $A_{X,n}^1$, from the description of $K_0(A_{X,n}^1)_+=K_0(C(X_n))_+$ in \ref{1323proj} and the choice of $X_n$ in \ref{range 0.18}}},  we have $(K_0(A_{X,n}),K_0(A_{X,n})_+,[1_{A_{X,n}}], K_1(A_{X,n}))=(H_n, (H_n)_+, u_n, K_n).$ {{On c}}onsidering the six term exact sequence for the K-theory of  the short exact sequence
$$
0\to I_n \to A_n \to A_{X,n}\to 0,
$$
the proof of part (1) follows  the lines
of the proof of Proposition \ref{2Lg13},
with $F_n$ {{replaced}} by $A_{X,n}$.
(Note that $I_n\cong C_0((0,1), E_n),$ and the boundary map $K_0(A_{X,n})\to K_1(I_n)=K_0(E_n)$ is given by
$$
(\bt_{n,1}-\bt_{n,0})\circ \pi=\pi_n\circ \pi:~ K_0(A_{X,n})=H_n\stackrel{\pi}{\longrightarrow} K_0(F_n)=H_n/{\rm Tor}(H_n)=\Z^{p_n}\stackrel{\pi_n}{\longrightarrow}K_0(E_n)=\Z^{l_n},
$$
and this map (playing the role of the map ${\phi_1}_{*0}-{\phi_0}_{*0}$ in \ref{2Lg13}) is surjective, since both $\pi_n$ and $\pi$ are surjective.) In particular, $(\pi_{I,n})_{*0}: K_0(A_n)=G_n\to K_0(A_{X,n})=H_n$  is the inclusion map, and $(\pi_{I,n})_{*1}: K_1(A_n)=K_n\to K_1(A_{X,n})=K_n$ is the identity map.

{\blue{For part (2), we know that, in the
commutative diagram
\begin{displaymath}
\xymatrix{
K_{*i}(A_n)\ar[r]^{(\pi_{I,n})_{*i}}\ar@{->}[d]_{(\phi_{n,n+1})_{*i}}
 &
K_{*i}(A_{X,n})
\ar@{->}[d]_{(\psi_{n, n+1})_{*i}}
\\
 K_{*i}(A_{n+1})\ar[r]^{(\pi_{I,n+1})_*}& K_{*i}(A_{X,n+1})
  }
\end{displaymath}
($i=0,1$), the horizontal maps are injective. Thus, $(\phi_{n,n+1})_*$ is uniquely determined by $(\psi_{n,n+1})_*.$}}

\end{proof}

\end{lem}

\begin{NN}\label{condition2}
{\blue{Recall that $(G, G_+, u, K, \Delta, r)$ is
fixed, as in \ref{range 0.16} and \ref{range 0.17}.}}
The construction of $A_{n+1}$ and $\phi_{n,n+1}$ will be done by induction.
Suppose that we already have the first part of the inductive sequence:
\begin{displaymath}
    \xymatrix{
        A_1 \ar[r]^{\phi_{1,2}} & A_2 \ar[r]^{\phi_{2,3}}&A_3 \ar[r]^{\phi_{3,4}} &\cd \ar[r]^{\phi_{n-1,n}} & A_n,}
\end{displaymath}
satisfying the following four conditions: for each $m=1,2,\cd n-1$,

(a)\, {\blue{For  subsequences  $G_{{k_n}}$ and $H_{{k_n}}$ {{(of $G_n$ and $H_n$), we have}}
{{$$(K_0(A_{X,n}),K_0(A_{X,n})_+,[\one_{A_{X,n}}], K_1(A_{X,n}))=({{H}}_{k_n}, ({{H}}_{k_n})_+,u_{k_n}, K_{k_n}),$$
$$(K_0(F_n),K_0(F_n)_+,[\one_{F_n}])=(H_{k_n}/{\rm Tor}(H_{k_n}), (H_{k_n}/{\rm Tor}(H_{k_n}))_+,{{\pi_{H_{k_n}, H_{k_n}'}}}(u_{k_n})), \text{and}$$
$$(K_0(A_{n}),K_0(A_{n})_+,[\one_{A_{n}}], K_1(A_{n}))=({{G}}_{k_n}, ({{G}}_{k_n})_+,u_{k_n}, K_{k_n}),$$ where ${{\pi_{H_{k_n}, H_{k_n}'}}}$ is the projection from $H_{k_n}$ to ${{H_{k_n}'=}}H_{k_n}/{\rm Tor}(H_{k_n}).$ }}
 }}
({{Without loss of generality, we may relabel  $k_n$ by $n$, and then $\gm_{k_n,k_{n+1}}:=\gm_{k_{n+1}-1,k_{n+1}}\circ \gm_{k_{n+1}-2,k_{n+1}-1}\circ\cdots \gm_{k_{n},k_{n}+1}$ and $\chi_{k_n,k_{n+1}}:=\chi_{k_{n+1}-1,k_{n+1}}\circ \chi_{k_{n+1}-2,k_{n+1}-1}\circ\cdots \chi_{k_{n},k_{n}+1}$ will  become
$\gm_{n,n+1}$ {{and}}  $\chi_{n,n+1},$ respectively.}}) {\blue{Furthermore,}}
$\phi_{m,m+1}(I_m)\subset I_{m+1}$ and $\phi_{m,m+1}({{J}}_m)\subset {{J}}_{m+1},$ and therefore
the map $\phi_{{m,m+1}}$ induces {{three}} homomorphisms $\psi_{m,m+1}: {{A_{X,m}}}\,{{(=A_m/I_m)}}\to
    {{A_{X,m+1}}},$ $\bar \phi_{m,m+1}: A_{{C},m}\,{{(=A_m/J_m)}} \to A_{{C},m+1},$ {{and $\psi^q_{m,m+1}=\phi^q_{m,m+1}: F_m\to F_{m{{+}}1}$, where $F_m$  arises
    as a quotient algebra in three ways: $F_m=A_m/(I_m\oplus J_m)=A_{X,m}/{\bar J}_m=A_{C,m}/{\bar I}_m$ (see \ref{construction}  )}};

(b)\, all  the homomorphisms $\phi_{m,m+1}$, $\psi_{m,m+1},$ and $\bar \phi_{m,m+1}$ are injective; in particular, $\psi_{m,m+1}|_{{{\bar J}}_m}: {{\bar J}}_m\to {{\bar J}}_{m+1}$ is injective, {\blue{and
$(\phi_{m, m+1})_{*0}=\gamma_{m, m+1}|_{G_n},$ $(\phi_{m,m+1})_{*1}=\chi_{m,m+1},$
$(\psi_{m,m+1})_{*0}=\gamma_{m,m+1},$ and $(\psi_{m,m+1})_{*1}=\chi_{m,m+1};$}}


(c)\, the induced map $\bar \phi_{m,m+1}:A_{{C},m} \to A_{{C},m+1}$ satisfies the conditions (1)--(8) of \ref{range 0.6} with $m$ in place of $n$, with $A_{{C},m}$  in place of $C_n$ (and of course naturally with $A_{{C},m+1}$ in place of $C_{n+1}$), {{ with}} $G_m/\mathrm{Tor}(G_m)$ (or $G_{m+1}/ \mathrm{Tor}(G_{m+1})$) and $H_m/\mathrm{Tor}(H_m)$ (or $H_{m+1}/\mathrm{Tor}(H_{m+1})${{, all from \ref{range 0.16},}}  in place of $G_n$ (or $G_{n+1}$) and $H_n$ (or $H_{n+1}$),
{{and with
$$
{\blue{({{\phi}}^{q}_{m,m+1})_{*0}}}=\gm'_{m,m+1}{\blue{=(c_{i,j}^{m,m+1})}}:~ H_m/{\rm Tor}(H_m)=\Z^{p_m}\lr H_{m+1}/{\rm Tor}(H_{m+1})=\Z^{p_{m+1}}
$$
 (satisfying $\gm'_{m,m+1}(G_m/\mathrm{Tor}(G_m))\subset G_{m+1}/ \mathrm{Tor}(G_{m+1})$) from \ref{range 0.16} in place of \\ $\gm_{n,n+1}:~
H_n\to H_{n+1}$ (satisfying $\gm_{n,n+1}(G_n)\subset G_{n+1}$),}}  respectively; {{moreover,}}
 $$
{{\spdd_1}} \qq\qq\qq\qq\qq\qq\qq c^{{m,m+1}}_{ij}> 13 \cdot 2^{2m}\cdot {\blue{(M_m+1)}}L_m\qq\,\,\,
\rforal \,\,i,j,\qq\qq\qq\qq,
$$
where $M_m=\max\{ b^{{m}}_{0,ij}:~ i=1,2,..., p_m,~ j=1,2, ..., l_m\}$
and $L_m$ {{is}} {{specified below,}} 

(d)\, the matrices ${{\bb_{m+1,0}}}$ and ${{\bb_{m+1,1}}}$ for each $A_{m+1}$ {\blue{satisfy
the condition
$$\spdd~~~
\td  b_{0,ji}:=\sum_{k=1}^{p_{m+1}} b^{{m+1}}_{0,ik}\cdot
c^{{m,m+1}}_{kl}, \,\,\td
 b_{1,ji}:=\sum_{k=1}^{p_{m+1}} b^{{m+1}}_{1,ik}\cdot c^{{m,m+1}}_{kl}
 > 2^{2n}\!\left(\sum_{k=1}^{l_n}
 (|d_{jk}^{n, n+1}|+2)\{n,k\}\!\!\right),
$$
where $\td \gamma_{n,n+1}=(d_{ij}^{n,n+1}): H_n/G_n\to H_{n+1}/G_{n+1}.$}}

{\blue{The number $L_m$ above which was to be chosen
after the  $m$-th  step will now be
specified:}}

Choose a finite set $Y_m\subset X_m\setminus \{{{x^0_m}}\}$ (where ${{x^0_m}}$ is the base point of $X_m$) such that for each $i < m,~
\bigcup_{y\in Y_{{m}}}Sp(\phi_{{i,m}}|_{y})$ is $\frac1m$-dense in $X_i$.
This can be done since the corresponding map $\psi_{{i,i+1}}|_{{{{\bar J}}_i}}:~{{\bar J}}_i (\subset {{A_{X,i}^1}})\to {{\bar J}}_{i+1} (\subset {{A_{X,i+1}^1}})$ is injective for
each $i <m,$  by the induction assumption (see (b) above).
{{Recall from \ref{notation-Dec12-2018} that we}} denote $t\in(0,1)_j\subset
Sp\big(C([0,1],E_m^j)\big) $ by $t_{{m,j}}$ to distinguish the
spectra of
different direct summands of $C([0,1], E_m)$ {{and different $n$}}.
Let $T_m\subset Sp(A_m)$ be defined by
$$T_m=\left\{ (\frac km)_{{m,j}};~ j=1,2, ...,l_m;~ k=1,2,..., m-1 \right\}.$$
 Let $Y_m=\{y_1,y_2,... , y_{_{L_{m,Y}}}\}\subset X_m$ and
  let $L_m=l_m\cdot (m-1)+L_{m,Y}=\#(T_m\cup Y_m)$.

\vspace{0.1in}

{{We will}} construct $A_{n+1}$ now and {{later}} the homomorphism $\phi_{n,n+1}$.
{{As {\blue{part of the}} induction assumption, suppose that {{the algebra}} $A_n=\big(C([0,1],E_n)\oplus {{P_n M_{\infty}(C(X_n))P_n}}\big)\oplus_{\bt_{n,0},\bt_{n,1},\bt_{X_n}} F_n$ is already}} {{constructed, }} with $(\bt_{{n,0}})_{*0}=\bb_{{n,0}}=(b^{{n}}_{0,ij})$
and
$(\bt_{{n,1}})_{*0}=\bb_{{n,1}}=(b^{{n}}_{1,ij})$.
So $M_n$ and $L_n$ can be chosen as described above.

{\blue{As at the end of  \ref{range 0.16}, choose $m'>n$ such that
$\gamma_{n, m}'$ has multiplicity at least  $13 \cdot 2^{2m}\cdot (M_n+1)L_n.$   Then, by renaming $H_{m'}$ as $H_{n+1},$
\wilog,
we may assume that
 $$
c^{{n,n+1}}_{ij}> 13 \cdot 2^{2n}\cdot (M_n+1)L_n\,\,\,
\rforal \,\,i,j,$$ in order to satisfy ${{\spdd_1}}$ in  condition (c).
}}

{\blue{Note that we still have \ref{range 0.16} and
\ref{range 0.17} (passing to a subsequence).}}
{\blue{Since the scaled ordered group $(G_{n+1}, u_{n+1})\subset (H_{n+1},u_{n+1})$, with $u_{n+1}=\big(([n+1,1],\tau_{n+1}),[n+1,2],...,[n+1,p_n]\big)$ and ${\rm Tor}(G_{n+1})={\rm Tor}(H_{n+1}),$ and  {{the group}} $K_{n+1}$  are {{now}} chosen,
the algebra $F_{n+1}=\bigoplus_{i=1}^{p_{n+1}}M_{[n+1,i]}$ (with $K_0(F_{n+1})=H_{n+1}/{\rm Tor}(H_{n+1})$ as scaled ordered group) can now be defined.
Furthermore, the space $X_{n+1}$ (with base point $x^0_{n+1}$ and $K_0(X_{n+1})=\Z\oplus {\rm Tor}(G_{n+1})$, $K_1(X_{n+1})=K_{n+1}$) and {{the projection}} $P_{n+1}$ (with $[P_{n+1}]=([n+1,1], \tau_{n+1})\in \Z\oplus {\rm Tor}(G_{n+1})$), and the identification $\bt_{X_{n+1}}$ of $F^1_{n+1}$ with $\big(P_{n+1} M_{\infty}(C(X_{n+1}))P_{n+1}\big)|_{x^0_{n+1}}$  {{can now be defined}}
as in \ref{range 0.18} and \ref{range 0.19}.}}

{\blue{Let us modify $b_{n+1,0}$ and $b_{n+1,1}$ given by \ref{range 0.19}.}}
 {\blue {Let
 \beq\label{Lambdanbt}
 \Lambda_n=2^{2n}\!\left(\sum_{k=1}^{l_n}
 (|d_{jk}|+2)\{n,k\}\!\!\right)
\eneq
and let $k_0=\max\{p_{n+1}^0+1,3\}.$
{{We replace}}   $b^{n+1}_{1,ik_0}$ and $b^{n+1}_{1,ik_0}$ by
 $b^{n+1}_{0,ik_0}+\Lambda_n,$ and by $b^{n+1}_{1,j,k_0}+\Lambda_n,$   $i=1,2,....,l_{n+1}.$ Then, since $c_{ij}^{n,n+1}$ is at least $26,$  one easily sees that, with {{the}} new
 $b^{n+1}_{0,ik_0},$
 $\spdd$ holds for $n+1.$
 Moreover,  the first {{$p_{n+1}^0$}} columns of both  $\bb_{n+1,0}$ and $\bb_{n+1,1}$  are still zero, and the last {{$p_{n+1}-p_{n+1}^0$}}  columns of both  $\bb_{n+1,0}$ and $\bb_{n+1,1}$ are strictly positive (see \ref{range 0.5a}).}}

{\blue{With the choice of $m'$ (as $n+1$) {\blue{above,}} and the choice of $\bb_{n+1, 0}$ {{and}} $\bb_{n+1, 1}$ {\color{blue}{above}},
one  defines $A_{n+1}$ as  in \eqref{13Dec15-2018}  which satisfies}}
 Conditions ${{\spdd_1}}$ and $\spdd.$
 {\blue{Note that as $A_{n+1}$ is constructed, we also obtain $A_{C,n+1}.$}}

{\blue{Next,}} we will begin the construction of $\phi_{n,n+1}$ by constructing
{\blue{${\bar \phi}_{n,n+1}$}} and $\psi_{n,n+1}$ first. {{Note that we already have {{the}} map $\gm_{n,n+1}=\td \cc_{n,n+1}=(\td
c^{n,n+1}_{ij}): H_n\to H_{n+1}$ with $\gm_{n,n+1}(G_n)\subset G_{n+1}$ from \ref{range 0.16}.}}
{{Recall that $H_k'=H_k/{\rm Tor}(G_k)$ in  \ref{range 0.16}, $k=1,2,....$
Denote by $u_k'$ the image of $u_k$ in $H_k'.$
{{Therefore,}} $\gamma_{n,n+1}': H_n'\to H_{n+1}'$ is also defined (with new 
{{subscripts}}).}}
{\blue{Note {{also}} that {{the}} algebras $A_{C,n}$
and $A_{C, n+1}$
have {{the}} same property as $C_n$ and $C_{n+1}$ {{of}}  Lemma \ref{range 0.6}
with $(H_n, (H_n)_+, u_n)$ and $(H_{n+1}, (H_{n+1})_+, u_{n+1})$
replaced by $(H_n', (H_n')_+, u_{n}')$ and $(H_{n+1}', (H_{n+1}')_+, u_{n+1}'),$ {\color{blue}{respectively.}}
Moreover, $\gamma_{n, n+1}'(G_n')\subset G_{n+1}'.$
To apply \ref{range 0.6}, we also
replace $\gamma_{n,n+1}$ by $\gm_{n,n+1}',$
{{and}} $G_n$  and $G_{n+1}$ by $G_n'$  and $G_{n+1}'.$
Then, by  Lemma \ref{range 0.6},  there is an injective homomorphism ${\bar \phi}_{n,n+1}{{:A_{C,n}\to A_{C,n+1}}}$
(in place of the homomorphism $\phi_{n,n+1}$) which satisfies  (1)--(8) of Lemma  \ref{range 0.6}.
 In particular, the homomorphism ${\bar \phi}_{n,n+1}$ satisfies\\
(a'): ${\bar \phi}_{n,n+1}({\bar I}_n)\subset {\bar I}_{n+1}$ (this is (3) of Lemma \ref{range 0.6}) and\\
(b') the K-theory map $ {{(\phi^q_{n,n+1})_{*,0}}}: K_0(F_n)=H'_n \to K_0(F_{n+1})={{H'_{n+1},}}$ induced by {{the}} quotient map ${{ \phi^q_{n,n+1}}}: F_n\to F_{n+1}$,
 is the same as $\gm'_{n,n+1}=\cc_{n,n+1}: K_0(F_n)=H'_n=\Z^{p_n} \to K_0(F_{n+1})=H'_{n+1}=\Z^{p_{n+1}}$ (this is (4) of Lemma \ref{range 0.6}). }}

\end{NN}

\begin{NN}\label{range 0.20}
~~ Recall that ${\blue{(K_0(F_n), 1_{F_n})}}=(H_n, u_n)$,
${\blue{(K_0(F_{n+1}), [1_{F_{n+1}}])}}=(H_{n+1},u_{n+1}),$ and the map $\td \cc=(\td
c_{ij}):~ H_n\to H_{n+1}$ {{are}} as in \ref{range 0.16}.  Assume that
$c_{ij}>13$ for any $i$ and $j$ (which is a consequence of ${{\spdd_1}}$).  We shall define the unital  homomorphism
$\psi_{{n,n+1}}:~{{A_{X,n}}}\to {{A_{X,n+1}}}$ to satisfy the following
conditions:
 \beq\nonumber
&&\hspace{-1.1in}(1)\hspace{0.1in} (\psi_{{n,n+1}})_{*0}={\blue{\gm_{{n,n+1}}}}: K_0({{A_{X,n}}})(=H_n)\lr K_0({{A_{X,n+1}}})\,(=H_{n+1}),\andeqn\\
\nonumber
&& (\psi_{{n,n+1}})_{*1}=\chi_{{n,n+1}}:~ K_1({{A_{X,n}}})\,(=K_n)\lr K_1({{A_{X,n+1}}})\,(=K_{n+1});
 \eneq
(2)\hspace{0.2in} $\psi_{{n,n+1}}({{\bar J}}_n)\subset {{\bar J}}_{n+1}$, and the {{quotient}} map ${\psi}^{{q}}_{{n,n+1}}:~{ {F}}_n\to { {F}}_{n+1}$  induced by $\psi_{{n,n+1}}$  satisfies $$({\psi}^{{q}}_{{n,n+1}})_{*0}={\blue{\gm_{n,n+1}'}}=\cc_{{n,n+1}}=(c^{{n,n+1}}_{ij}):
 K_0({{ F}}_n)(=\Z^{p_n})\to ({{ F}}_{n+1})_{*0}(=\Z^{p_{n+1}}).$$


{\blue{Denote by $F_k'=\bigoplus_{{i\ge 2}} F_k^i=\bigoplus_{{i\ge 2}} A_{X,k}^i,$ $\pi_k'': F_k\to F_k',$
and $\pi_k^{-,1}: F_k\to F_k^1$ the projection maps, $k=1,2,....$
{\blue{Let $\pi_{x_n^0}: P_nM_{\infty}(C(X_n))P_n\to F_n^1$ be the point evaluation at $x_n^0.$
Define $\pi_{x_n^0}^\sim: A_{X, n}=A_{X,n}^1
\oplus F_n'\to F_n$ by
$\pi_{x_n^0}^\sim (a,b)=(\pi_{x_n^0}(a), b)$ for all $(a,b)\in A_{X,n}^1\oplus F_n'.$
It is well known (see \cite{Ell-76}) that there is a \hm\, ${\psi}^{q}_{n,n+1}{'}:~{{ F}}_n\to {{ F}}_{n+1} $
such that $(\psi_{n, n+1}^{q}{'})_{*0}=\cc_{n,n+1}.$}}
Then, define $\psi_{n, n+1}'': A_{X,n}\to F_{n+1}'$
by $\psi_{n, n+1}''=\pi_{n+1}''\circ \psi^q_{n,n+1}{'}\circ \pi_{x_n^0}^\sim.$}}

  {{Recall from \ref{range 0.16}, $\gm_{{n,n+1}}^{i,1}([\e_{
     A_{X,n}^i}])=(c^{n,n+1}_{{1,i}}\cdot[n,i]+T^{n,n+1}_i({\blue{[\e_{
     A_{X,n}^i}]}}))\in \Z\oplus {\rm Tor}(H_{n+1})=K_0(A_{X,n+1}^1)${\blue{, where $[1_{
     A_{X,n}^i}]=[n,i]$ if $i\geq 2$, or $[1_{
     A_{X,n}^i}]=[n,1]+\tau_n$ if $i=1$.}} }}
{\blue{By \ref{1323proj},}}
     one can find  projections $Q_1,Q_2, ... , Q_{p_n}$ such that
      $\gm_{{n,n+1}}^{i,1}([\e_{{{
     A_{X,n}^i}}}])=[Q_i] \in K_0(P_{n+1}M_{\infty}(C(X_{n+1}))P_{n+1})$. {{Since  ${{\bigoplus}}_{i=1}^{p_n}\gm_{{n,n+1}}^{i,1}([1_{
     {A_{X,n}^i}}])=\gm_{{n,n+1}}^{-,1}([1_{{A_{X,n}}}])=[1_{
     A_{X,n+1}^1}]=[P_{n+1}]$, by Remark 3.26 of \cite{EG-RR0AH}, one can make $\{Q_i\}_{i=1}^{p_n}$ mutually orthogonal}} {{with}}
     $$Q_1+Q_2+\cdots +Q_{p_n}=P_{n+1}.$$
    Since $c_{ij}>13$ for all $i,j$---{{i.e.,}} ${\rm rank}(Q_i)/{\rm rank} (\e_{{{A_{X,n}^i}}})=c_{1,i}>13${{---}}, by \ref{range 0.14}
    {\blue{(using the base point $x_{n+1}^0$ (as $y_0$)),}}
    there are unital homomorphisms $\psi_{{n,n+1}}^{i,1}:{{A_{X,n}^i}}\to
    Q_i{{A_{X,n+1}^1}}Q_i$ which realize the K-theory map $\gm_{{n,n+1}}^{i,1}
:K_0({{A_{X,n}^i}})\to K_0({{A_{X,n+1}^1}}),$ and $\chi_{{n,n+1}}:~
K_1({{A_{X,n}^1}})~(=K_n)\to K_1({{A_{X,n+1}^1}})~(=K_{n+1})$ (note that
$K_1({{A_{X,n}^i}})=0$ for $i\geq 2$).
{\blue{Moreover,  since, in {{the}}  application of \ref{range 0.14}, we used the base point $x_{n+1}^0$ (as $y_0$), we have
$\psi_{n,n+1}^{1,1}({\bar J}_n)\subset {\bar J}_{n+1}$ (see also the last line of \ref{construction}).
It follows that $\psi_{n,n+1}^{1,1}$  induces a {{(not necessarily}} unital{{)}} \hm\,
$\psi_{n, n+1}^{1,1,q}: A_{X,n}^1/{\bar J}_n\to A_{X, n+1}^{{1}}/{\bar J}_{n+1}.$
{{That is, $\psi_{n,n+1}^{1,1}$  induces a (not necessarily}} unital{{)}} \hm\, $\psi_{n, n+1}^{1,1, q}: F_n^1\to F_{n+1}^1.$
Note that $K_0(A_{X,n}^1)=\Z\oplus {\rm Tor}(H_n)$ and $K_0(A_{X,n+1}^1)=\Z\oplus {\rm Tor}(H_{n+1}).$
It is then clear {{that}} ${\psi_{n,n+1}^{1,1,q}}_{*0}=c_{1,1}^{n,n+1}$ (see  \ref{range 0.16}).
}}
{\blue{Define $\psi_{n, n+1}^{-,1}: A_{X,n}\to A_{X,n+1}^1$ by
$\psi_{n,n+1}^{-,1}|_{A_{X,n}^i}=\psi_{n,n+1}^{i,1},$ $i=1,2,...,p_n.$
It follows {\blue{that}} $\psi_{n,n+1}^{-,1}$ induces a unital \hm\,
$\psi_{n,n+1}^{-,1,q}: A_{X,n}/{\bar J}_n\to A_{X,n+1}^1/{\bar J}_{n+1}.$
One computes that $(\psi_{n,n+1}^{-,1,q})_{*0}=(\pi_{n+1}^{-,1})_{*0}\circ \gamma_{n,n+1}.$
Now define $\psi_{n,n+1}: A_{X,n}\to A_{X,n+1}{{=A_{X,n+1}^1\oplus F_n'}}$
by $\psi_{n,n+1}=\psi_{n,n+1}^{-,1}\oplus {{\psi_{n,n+1}''}}.$  One then checks
that (1) above holds.   Let
$\psi_{n,n+1}^q: F_n\to F_{n+1}$ be the induced unital \hm. Then, as above,
$(\pi_{n+1}^{-,1}\circ \psi_{n,n+1}^q)_{*0}=(\pi_{n+1}^{-,1})_{*0}\circ \gamma_{n,n+1}'$ and
$(\pi_{n+1}''\circ \psi_{n,n+1}^q)_{*0}=(\pi_{n+1}''\circ \psi_{n,n+1}^q{'})_{*0}=(\pi_{n+1}'')_{*0}\circ \gm_{n,n+1}'.$
{{Therefore,}}
$(\psi_{n,n+1}^q)_{*0} =\gm_{n,n+1}'.$
}}

\end{NN}

\begin{NN}\label{range 0.21}
{\blue{In this subsection, we will describe $\psi_{n,n+1}^{i,1}$ for $i\ge 2.$}}

{\blue{Recall  that $\psi_{n, n+1}^{i,1}$ is a unital \hm\, from $A_{X,n}^i$ to $Q_iA_{X, n+1}^1Q_i,$
$i=1, 2,...,p_n.$
Since $n$ is fixed, to simplify notation, in this subsection, let us use $\psi^i$ for $\psi_{n, n+1}^i.$
Note that $A_{X,n+1}^1|_{[x_{n+1}^0, x_{n+1}^0+1]}$ may be identified with $M_{[n+1, 1]}(C([0,1])).$
{{Without loss of generality,}} by fixing a system of matrix {{units}}, we may write
\beq\nonumber
&&\hspace{-0.3in} {\mbox{(a)}}~~~ Q_i|_{[{{x^0_{n+1}}},{{x^0_{n+1}}}+1]}=\diag (\0_{c_{11}[n,1]},\0_{c_{12}[n,2]},..., \0_{c_{1~i-1}[n,i-1]},\e_{c_{1i}[n,i]},\0_{c_{1~i+1}[n,i+1]},..., \0_{c_{1p_n}[n,p_n]}),
\eneq
where $P_{n+1}|_{[{{x^0_{n+1}}}, {{x^0_{n+1}}}+1]}$ is identified with
$\e_{[n+1,1]}\in M_{[n+1,1]}(C[{{x^0_{n+1}}}, {{x^0_{n+1}}}+1]).$}}

{\blue{Let
{$\{e_{{kl}}\}$} be the matrix units of $F_n^i=M_{[n,i]}$ and let
$q=\psi^{{i}}(e_{11}).$ The unital \hm\, $\phi^i: M_{[n,i]}\to Q_iA_{X,n+1}^1Q_i$ allows us
to write $Q_iA_{X,n+1}^1 Q_i=qA_{X,n+1}^1q\otimes M_{[n,i]}=qM_{\infty}(C(X_{n+1}))q\otimes M_{[n,i]},$
and   to write}}
\beq\label{0.21*}
\psi^i\big( (a_{ij})\big)= q\otimes (a_{ij}) \in qM_{\infty}(C(X_{n+1}))q\otimes M_{[n,i]}.
\eneq
{{By Remark 3.26 of \cite{EG-RR0AH} (see also  \ref{1323proj})
w}}e can write $q=q_1+q_2+\cd+q_d+p$, where $q_1,q_2,\cd, q_d$ are
mutually equivalent trivial rank $1$ projections and $p$ is a
({{possibly}} {{non-trivial}}) rank $1$ projection. Under the identification
$Q_i=q\otimes \e_{M_{[n,i]}}$, we {{write}}  $\hat{q}_j{{:=}}q_j\otimes
\e_{M_{[n,i]}}$ and $\hat{p}{{:=}}p\otimes \e_{M_{[n,i]}}$.

$\\$ (b)~~ Projections $\hat{q}_1,\hat{q}_2,\cd, \hat{q}_d,$ and $ \hat{p}$ for $\psi^i$ can be chosen such that $\hat{q}_1|_{[{{x^0_{n+1}}}, {{x^0_{n+1}}}+1]}$, $\hat{q}_2|_{[{{x^0_{n+1}}}, {{x^0_{n+1}}}+1]}$ ,..., $\hat{q}_d|_{[{{x^0_{n+1}}}, {{x^0_{n+1}}}+1]},$
   and $ \hat{p}|_{[{{x^0_{n+1}}}, {{x^0_{n+1}}}+1]}$ are diagonal
   {{matrices}} with $\e_{[n,i]}$ in the correct place {\blue{(see below)}} when
    $Q_iM_{\infty}(C[{{x^0_{n+1}}}, {{x^0_{n+1}}}+1])Q_i$  {{is}}  identified with
     $M_{c_{i1}[n,i]}(C[{{x^0_{n+1}}}, {{x^0_{n+1}}}+1])$. That {{is,}}
$$\hat{q}_j= \diag (\underbrace{\0_{[n,i]},
..., \0_{[n,i]}}_{j-1}, \e_{[n,i]},\0_{[n,i]},..., \0_{[n,i]})\quad
 \mbox{and}\quad  \hat{p}=\diag (\underbrace{\0_{[n,i]},..., \0_{[n,i]} }_{d},\e_{[n,i]}).$$

\end{NN}

\begin{lem}\label{range 0.22}
 Let $i\geq 2$ and $\psi^i: {{A_{X,n}^i}}\to
Q_i{{A_{X,n+1}^1}}Q_i$ be as in  $\ref{range 0.21}$ above.
Suppose $m\leq d
=c^{{n,n+1}}_{1i}-1$. Let $\LD:  Q_i{{A_{X,n+1}^1}}Q_i \to M_m(Q_i{{A_{X,n+1}^1}}Q_i)$ be the
amplification defined by $\LD(a)=a\otimes {\bf 1}_{m}.$
Then there
is a projection $R^i\in M_m(Q_i{{A_{X,n+1}^1}}Q_i)$ such that

\hspace{-0.15in}{\rm (i)} $R^i$ commutes with $\LD(\psi^i(F_n^i))$ {{(note that $F_n^i=A_{X,n}^i$ for $i\geq 2$), and }}

\hspace{-0.15in}{\rm (ii)} $R^i({{x^0_{n+1}}})= Q_i({{x^0_{n+1}}})\otimes \left(
                                        \begin{array}{cc}
                                          \e_{m-1} & 0 \\
                                          0 & 0 \\
                                        \end{array}
                                      \right)= \diag(\underbrace{Q_i({{x^0_{n+1}}}),..., Q_i({{x^0_{n+1}}})}_{m-1}, 0))\in M_m({{A_{X,n+1}^1}}|_{{x^0_{n+1}}}).$
Consequently, $rank (R^i)= c^{{n,n+1}}_{i1}(m-1)[n,i]=(d+1)(m-1)[n,i]$.

 Let $\pi{\blue{:=\pi_{X,n+1}^e|_{M_m(Q_iA_{X,n+1}^1Q_i)}\otimes {\rm id}_{M_m}}}:
 M_m(Q_i{{A_{X,n+1}^1}}Q_i) \to {\blue{M_m(\pi_{X,n+1}^e(Q_i A^1_{X, n+1}Q_i))}}.$
 Then $\pi$ {\blue{maps}} $R^iM_m(Q_i {{A_{X,n+1}^1}} Q_i)R^i$ onto \\ $M_{m-1}({\blue{\pi_{X,n+1}^e(Q_i A^1_{X, n+1}Q_i)}})
 \subset M_{m}({\blue{\pi_{X,n+1}^e(Q_i A^1_{X, n+1}Q_i)}}).$

{{Below we will use the same notation $\pi$ to denote {{the}}  restriction {{of $\pi$}} to $R^iM_m(Q_i {{A_{X,n+1}^1}} Q_i)R^i$, whose codomain is $M_{m-1}({\blue{\pi_{X,n+1}^e(Q_i A^1_{X, n+1}Q_i)}})$.}}

{\rm (iii)} There is {{a unital embedding }}
 $$ \iota: M_{m-1}(\pi_{X,n+1}^e(Q_i A^1_{X, n+1}Q_i))
 \hookrightarrow R^iM_m({\blue{Q_i {{A_{X,n+1}^1}} Q_i}})R^i$$
 such that $\pi\circ \iota=\id|_{M_{m-1}({\blue{\pi_{X,n+1}^e(Q_i A^1_{X, n+1}Q_i}}))}$
 and such that $R^i(\LD (\psi^i(F_n^i)))R^i \subset \mathrm{Image}(\iota)$.

\end{lem}

 \begin{proof} In the proof of this lemma, $i\geq 2$ is fixed.
%

 The homomorphism $\LD\circ \psi^i:
 {\blue{A_{X,n}^i=}}F_n^i= M_{[n,i]} \to M_m(Q_i{{A_{X,n+1}^1}}Q_i)$ can be regarded as ${{\LD_1\otimes \id_{F^i_n}=}}\LD_1\otimes \id_{[n,i]}$,
 where $\LD_1: \C \to M_m(q{{A_{X,n+1}^1}}q)$
 is the  unital
 homomorphism given by
 \beq\label{13Dec19-2018-1}\LD_1(c)= c\cdot(q\otimes \e_m),
 \eneq
 {{ and $q=\psi^i(e_{11})$ {\blue{is}} as in \ref{range 0.21}, by identifying $M_m(Q_iA_{X,n+1}^1Q_i)$ {{with}} $M_m(qA_{X,n+1}^1q)\otimes F^i_n$. }}

  Note that $q=q_1+q_2+\cd +q_d+p$ with {{$\{q_i\}$  mutually}} equivalent rank
  one projections. Furthermore, $p|_{[{{x^0_n}},{{x^0_n}}+1]}$ is also a {\blue{rank one trivial}}
  projection. Let $r\in M_m(q{{A_{X,n+1}^1}}q)=
  q{{A_{X,n+1}^1}}q\otimes M_m$ be defined
  as {{follows:}}
  \beq\label{13Dec19-2018}
  r({{x^0_n}})=q({{x^0_n}})\otimes \left(
                         \begin{array}{cc}
                           \e_{m-1} & 0 \\
                           0 & 0 \\
                         \end{array}
                       \right)\qq\qq\qq\qq\qq\qq\qq\qq\\
                       \nonumber
                       = (q_1({{x^0_n}})+q_2({{x^0_n}})+\cd +q_d({{x^0_n}})+p({{x^0_n}}))\otimes \left(
                         \begin{array}{cc}
                           \e_{m-1} & 0 \\
                           0 & 0 \\
                         \end{array}
                       \right),\eneq
$$r({{x^0_n}}+1)= (q_1({{x^0_n}}+1)+q_2({{x^0_n}}+1)+\cd +q_{m-1}({{x^0_n}}+1))\otimes \e_m +~~~~~~~~~~~~~~$$
$$~~~~~~~~~~~~~+ (q_m({{x^0_n}}+1)+\cd +q_{d}({{x^0_n}}+1)) \otimes \left(
                         \begin{array}{cc}
                           \e_{m-1} & 0 \\
                           0 & 0 \\
                         \end{array}
                       \right).$$
In {{this definition}}, between ${{x^0_n}}$ and ${{x^0_n}}+1$, $r(t)$ can be
defined to be any continuous path connecting the projections
$r({{x^0_n}}), $ and $r({{x^0_n}}+1)$, both of rank
$(d+1)(m-1)=(m-1)m+(d-m+1)(m-1)$. (Note that all $q_i(t)$ and $p(t)$
are constant on $[{{x^0_n}}, {{x^0_n}}+1]$.)
{{Finally,}}  for $x\in X_{n+1}'\subset X_{n+1}$, {{define}}
 $$ r(x)= (q_1(x)+q_2(x)+\cd +q_{m-1}(x))\otimes \e_m + (q_m(x)+\cd +q_{d}(x)) \otimes \left(
                         \begin{array}{cc}
                           \e_{m-1} & 0 \\
                           0 & 0 \\
                         \end{array}
                       \right).$$
{{
Since all $\{q_k\}_{k=1}^d$ are trivial projections, $r|_{X_{n+1}'}$ is also a trivial projection. Hence $r$ itself is also a trivial projection as the inclusion map $X_{n+1}'\to X_{n+1}$ is a homotopy equivalence. }}

Let $R^i=r\otimes \e_{[n,i]}$ under the identification of
$M_m(Q_i{{A_{X,n+1}^1}}Q_i)$ with $M_m(q
{{A_{X,n+1}^1}}q)\otimes \e_{[n,i]}$.
Since $\LD_1: \C \to M_m(q{{A_{X,n+1}^1}}q)$ sends $\C$ to the center of
 $M_m(q{{A_{X,n+1}^1}}q)$ {{(see (\ref{13Dec19-2018-1}))}}, we have that $r$ commutes with $\LD_1(\C),$ and consequently,
 $R^i=r\otimes \e_{[n,i]}$ commutes with $\LD(\psi^i(F_n^i))$ as
 $\LD\circ \psi^i= \LD_1\otimes \id_{[n,i]}$. That is, condition (i) holds.
 Condition (ii) follows from the definition of $r({{x^0_{n+1}}})$ {{(see (\ref{13Dec19-2018}))}} and
 $R^i({{x^0_{n+1}}})=r({{x^0_{n+1}}})\otimes \e_{[n,i]}$.

Note that $r\in M_m(q{{A_{X,n+1}^1}}q)= M_m(qM_{\infty}(C(X_{n+1}))q)$ is
a trivial projection of rank $(d+1)(m-1)$ and
$r({{x^0_{n+1}}})=q({{x^0_{n+1}}})\otimes {{\diag(}}\e_{m-1}{{,0)}}$.

Note also that
$$
{\blue{\pi(r)M_m(\pi_{X,n+1}^e(qA_{X,n+1}^1q))\pi(r)=
M_{m-1}(\pi_{X,n+1}^e(qA_{X,n+1}^1q))}}
\cong
M_{(d+1)(m-1)}.$$
Let $r_{ij}^0, 1\leq i,j\leq (d+1)(m-1)$ be a system of  matrix units
for $M_{(d+1)(m-1)}$.
Since $r$ is a trivial projection, one
can construct
$r_{ij}\in rM_m(q{{A_{X,n+1}^1}}q)r$, $1\leq i,j\leq (d+1)(m-1),$
 with ${\blue{\pi(r_{ij})}}
 =r_{ij}^0$ serving as a system of matrix units for $M_{(d+1)(m-1)}\subset
 rM_m(q{{A_{X,n+1}^1}}q)r \cong M_{(d+1)(m-1)}(C(X_{n+1}))$.
 Here, by matrix units, we mean $r_{ij}r_{kl}= \dt_{jk}r_{il}$
 and $r=\sum_{i=1}^{(d+1)(m-1)} r_{ii}$. We {\blue{define}}
 $$\iota_1:  M_{m-1}(\pi_{X,n+1}^e(qA_{X,n+1}^1q))~
 {{\big(}}\cong \pi(r)
 M_{m}(\pi_{X,n+1}^e(qA_{X,n+1}^1q)
 \pi(r)){{\big)}}
 \,\,\hookrightarrow rM_m(q {{A_{X,n+1}^1}} q)r
 $$
 by $\iota_1(r_{ij}^0)=r_{ij}$. {{Finally,}} we define $\iota=\iota_1\otimes \id_{[n,i]}$. {{Since
 $\pi(r_{ij})
 =r_{ij}^0$, we have}}\\ {\blue{$\pi\circ \iota=\id|_{M_{m-1}(\pi_{X,n+1}^e(Q_iA_{X,n+1}^1Q_i))}.$}}
 {{Note that $r(\LD_1(\C))r=\C\cdot r\subset \iota_1(M_{m-1}(\pi(qA_{X,n+1}^1q))).$
 On the other hand,}} using the identification {{$M_m(Q_i{ A}_{X,n+1}^1Q_i)=M_m(q{ A}_{X, n+1}^1q)\otimes F^i_n$ (recall that $F^i_n=M_{[n,i]}$), we have}} $R^i=r\otimes \e_{[n,i]}$ and {{$\LD\circ\psi^i=\LD_1\otimes \e_{[n,i]}$. Hence $R^i(\LD (\psi^i(F_n^i)))R^i \subset \mathrm{Image}(\iota)$.}}\\
{\blue{So}} (iii) follows.
\end{proof}


\begin{NN}\label{range 0.23}
~~Now  we would like to choose $\psi_{n,n+1}^{1,1}$ in a {{specially}} simple form described below.
 We know that ${\rm rank} (Q_1)= c^{{n,n+1}}_{11}[n,1]$,
 where $[n,1]={\rm rank} (P_n)$ for $
 {{A_{X,n}^1}}=P_nM_{\infty}(C(X_n))P_n$.
 Note that $c^{{n,n+1}}_{11} > 13$. {{Set}} $d=c^{{n,n+1}}_{11} - 13$.
 {\blue{We may write $Q_1=Q'\oplus {\tilde Q},$ where $\td{Q}$ is a trivial projection of rank $d.$}} {{Recall that as in \ref{range 0.21}, we use $\psi^i$ for $\psi_{n,n+1}^{i,1}$. }}
 Then we can choose {{$\psi^1:=\psi_{n,n+1}^{1,1}$}} satisfying
the following {{conditions:}}

(a)\, $Q_1=\e_{d{{[n,1]}}}\oplus \td{Q}:=Q'\oplus \td{Q} \in M_{\infty}(C(X_{n+1}))$
  such that $\td{Q}|_{[{{x^0_{n+1}}},{{x^0_{n+1}}}+1]}=\e_{13{{[n,1]}}}$ (but in the lower right
  corner of $Q_1|_{[{{x^0_{n+1}}},{{x^0_{n+1}}}+1]}$);

(b)\,
   ${{\psi^1}}: {{A_{X,n}^1}} \to Q_1 {{A_{X,n+1}^1}} Q_1$  can be
   {{ decomposed as}}
   ${{\psi^1}}=\psi_1\oplus \psi_2,$  where
   $\psi_1: {{A_{X,n}^1}}\to M_{d[n,i]}(C(X_{n+1}))= Q'M_{\infty}(C(X_{n+1}))Q'$ and  $\psi_2:{{A_{X,n}^1}}\to \td{Q}M_{\infty}(C(X_{n+1}))\td{Q}$  are as follows:
   \begin{enumerate}
    \item[(b1)]
   the unital  homomorphism  $\psi_1: {{A_{X,n}^1}}\to M_{d[n,i]}(C(X_{n+1}))= Q'M_{\infty}(C(X_{n+1}))Q'$ is
    defined by $\psi_1(f)= \diag(\underbrace{f({{x^0_n}}),f({{x^0_n}}),..., f({{x^0_n}})}_d)$ as
    a constant function on $X_{n+1}$;

    \item[(b2)] the unital {\blue{injective}} homomorphism
    $\psi_2:{{A_{X,n}^1}}\to \td{Q}M_{\infty}(C(X_{n+1}))\td{Q}$ is a
    homomorphism satisfying $(\psi_2)_{*0}= \td{c}^{{n,n+1}}_{11}-d= c^{{n,n+1}}_{11}-d+T^{{n,n+1}}_1$
    (where $T^{{n,n+1}}_1: H_n^1(=K_0({{A_{X,n}^1}}))\to {\rm Tor}(H_{n+1}) \subset H_{n+1}^1$ is as
    in \ref{range  0.16})
    and $(\psi_2)_{*1}=\chi_{n,n+1}: K_1({{A_{X,n}^1}})~(=K_n)\to K_1({{A_{X,n+1}^1}})~(=K_{n+1})$
    (such $\psi_2$ exists  {\blue{since ${{\tilde Q}}$ has rank $13{{[n,1]}}$,}} {{by}}  \ref{range 0.14}).
    \end{enumerate}
  {\blue{Moreover, we may write
  $\psi_2(f)(x_{n+1}^0+1)=\diag(\underbrace{f({{x^0_n+1}}),f({{x^0_n+1}}),...,
     f({{x^0_n+1}})}_{13})$ for a fixed system of matrix {{units}} for $M_{13[n,1]}(C([x_{n+1}^0, x_{n+1}^0+1])).$
    Then,  we may  {{change}} $\psi_2(f)|_{[x_{n+1}^0, x_{n+1}^0+1]}$ so {{that}} it also satisfies the {{following}} condition:}}
\begin{enumerate}
    \item[(b3)]
    For $t\in [0,\frac12]$,
     $$\psi_2(f)({{x^0_{n+1}}}+t)=\diag(\underbrace{f({{x^0_n}}),f({{x^0_n}}),...,
     f({{x^0_n}})}_{13}),$$
    and for $t\in [\frac12, 1]$,
    $$\psi_2(f)({{x^0_{n+1}}}+t)=\diag(\underbrace{f({{x^0_n}}+ 2t{{-1}}),f({{x^0_n}}+ 2t{{-1}}),...,
    f({{x^0_n}}+ 2t{{-1}})}_{13}).$$
    Here, $f({{x^0_n}}+s)\in P_n({{x^0_n}}+s)M_{\infty}P_n({{x^0_n}}+s)$ is
    regarded as an $[n,1]\times[n,1]$ matrix for each $s\in [0,1]$ by
    using the fact $P_n|_{[{{x^0_n}}, {{x^0_n}}+1]}=\e_{[n,1]}$.
    \end{enumerate}

 Let us remark that
$\psi_1: {{A_{X,n}^1}}\to M_{\blue{d[n,1]}}(C(X_{n+1}))(=Q'M_{\infty} (C(X_{n+1}))Q'{{)}}$
  factors through ${{F}}_n^1=M_{[n,1]}(\C)$, and the
  restriction ${{\psi^1}}|_{[{{x^0_{n+1}}},{{x^0_{n+1}}}+\frac12]}$
 also factors
  through ${{F}}_n^1$,
  as $${{\psi^1}}(f)(x)=\diag(\underbrace{f({{x^0_n}}),..., f({{x^0_n}})}_{d+13})\,\,\,{\rm for\,\,\, any}\,\,\,x\in [{{x^0_{n+1}}},{{x^0_{n+1}}}+\frac12].$$

  \end{NN}

\begin{lem}\label{range 0.24}  Suppose that $Q_1$ and ${{\psi^1}}:
{{A_{X,n}^1}} \to Q_1{{A_{X,n+1}^1}}Q_1$ satisfy conditions (a) and (b) {{above}} (including (b1), (b2), and (b3)).
 Suppose that $13m\leq d=c^{{n,n+1}}_{11}-13$. Let $\LD:  Q_1{{A_{X,n+1}^1}}Q_1 \to M_m(Q_1{{A_{X,n+1}^1}}Q_1)$
 be {{the}} amplification  defined by
$\LD(a)=a\otimes {\bf 1}_m.$
There is a projection $R^1\in
M_m(Q_1{{A_{X,n+1}^1}}Q_1)$ satisfying the following conditions:

{\rm (i)} $R^1$ commutes with $\LD(\psi^1({{A_{X,n}^1}}))$ and

{\rm (ii)} $R^1({{x^0_{n+1}}})= Q_1({{x^0_{n+1}}})\otimes \left(
                                        \begin{array}{cc}
                                          \e_{m-1} & 0 \\
                                          0 & 0 \\
                                        \end{array}
                                      \right)= \diag(\underbrace{Q_1({{x^0_{n+1}}}),\cd, Q_1({{x^0_{n+1}}})}_{m-1}, {{0}})\in M_m(
                                      {{A_{X,n+1}^1}}|_{{x^0_{n+1}}}).$
Consequently, $rank (R^1)= c^{{n,n+1}}_{11}(m-1)[n,1]=(d+13)(m-1)[n,1]$.

 Let $\pi:={\pi_{X,n+1}^e}|_{M_m(Q_1{{A_{X,n+1}^1}}Q_1)}\otimes {\rm id}_{M_m}: M_m(Q_1{{A_{X,n+1}^1}}Q_1) \to M_{m}
 ({\blue{\pi_{X,n+1}^e(Q_1A_{X,n+1}^1Q_1}})).$
{\blue{T}}hen $\pi$ {\blue{maps}} $R^1M_m(Q_1 {{A_{X,n+1}^1}} Q_1)R^1$ onto  $M_{m-1}({\blue{\pi_{X,n+1}^e(Q_1A_{X,n+1}^1Q_1}}))
 \subset M_{m}({\blue{\pi_{X,n+1}^e(Q_1A_{X,n+1}^1Q_1}})).$

{{Below we will use the same notation $\pi$ to denote {{the}} restriction {{of $\pi$}}  to $R^1M_m(Q_1 {{A_{X,n+1}^1}} Q_1)R^1$, whose codomain is $M_{m-1}({\blue{\pi_{X,n+1}^e(Q_1A_{X,n+1}^1Q_1}}))$.}}

{\rm (iii)} There is {{a unital embedding}}
 $$ \iota: M_{m-1}({\blue{\pi_{X,n+1}^e(Q_1A_{X,n+1}^1Q_1}}))
 \hookrightarrow R^1M_m(Q_1 {{A_{X,n+1}^1}} Q_1)R^1$$ such that $\pi\circ \iota=\id|_{M_{m-1}({\blue{\pi_{X,n+1}^e(Q_1A_{X,n+1}^1Q_1}}))}$
 and such that $R^1(\LD (\psi^1({{A_{X,n}^1}})))R^1 \subset
 {\rm Image} (\iota)$.

\end{lem}

The notation $\LD$, $d$, {{and}} $m$ in the lemma above, and  $q, q_1,q_2,..., q_d, p, {{r,}} $ and $\LD_1$ in the proof below, are also used in Lemma \ref{range 0.22} and
its proof for the case $i\geq 2$ (comparing with $i=1$ here). Since they are used for the same purpose, we choose the same notation.

 \begin{proof} The map
 \beq\label{13Dec-23-2018}
\xymatrix{\psi_1: {{A_{X,n}^1}}\ar[r]^{\pi}&  {{F}}_n^1 \ar[r] &
M_{d[n,1]}(C(X_{n+1}))=Q'M_{\infty}(C(X_{n+1}))Q'}
\eneq
(where $Q'=\e_{d[n,1]}$)
can be written as $(\LD_1\otimes
\id_{[n,1]})\circ \pi$, where $\LD_1: \C \to M_d(C(X_{n+1}))$ is
the map sending $c\in \C$ to $c\cdot\e_d$. We write
$\LD_1(1):=q'=q_1+q_2+\cd+q_d$, with each $q_i$ a trivial constant
projection of rank $1$. Here $q'$ is {\blue{a}} constant subprojection of $Q'$
with $Q'=q'\otimes \e_{[n,1]}$. Consider the map
$\td{\psi}_2:=\psi_2|_{[{{x^0_{n+1}}},{{x^0_{n+1}}}+\frac12]}: {{A_{X,n}^1}} \to
\td{Q}
{{A_{X,n+1}^1}}\td{Q}|_{[{{x^0_{n+1}}},{{x^0_{n+1}}}+\frac12]},$ and $\td{\psi}^1:=\psi^1|_{[{{x^0_{n+1}}},{{x^0_{n+1}}}+\frac12]}=
(\psi_1+\psi_2)|_{[{{x^0_{n+1}}},{{x^0_{n+1}}}+\frac12]}: {{A_{X,n}^1}} \to
{Q_1}{{A_{X,n+1}^1}}{Q_1}|_{[{{x^0_{n+1}}},{{x^0_{n+1}}}+\frac12]}$. As
pointed out in
\ref{range 0.23}, $\td{\psi}_2$ has the factorization
\begin{displaymath}
\xymatrix{{{A_{X,n}^1}}\ar[r]^{\pi} &  {{F}}_n^1\ar[r] &
M_{13[n,1]}(C[{{x^0_{n}}},{{x^0_{n+1}}}+\frac12]).}
\end{displaymath}
Hence $\td{\psi}^1$ has the factorization
\beq\label{13Dec-23-2018-1}
\xymatrix{{{A_{X,n}^1}}\ar[r]^{\pi} &  {{F}}_n^1 \ar[r] &
M_{(d+13)[n,1]}(C[{{x^0_{n+1}}},{{x^0_{n+1}}}+\frac12]).}
\eneq
The map $\td{\psi}^1$ can be written as $(\LD_2\otimes \id_{[n,1]})\circ \pi$,
where $\LD_2: \C \to M_{d+13}(C[{{x^0_{n+1}}},{{x^0_{n+1}}}+\frac12])$ is the map defined by
sending $c\in \C$ to $c\cdot\e_{d+13}$. We write
$\LD_2(1):=q=q_1+q_2+\cd+q_d+ p$ with each $q_i$ {{the}}
restriction of $q_i$ {{appearing}} in the definition of $\LD_1(1)$
on $[{{x^0_{n+1}}},{{x^0_{n+1}}}+\frac12]$, and {{$p$ a}} rank $13$ {{constant}}
projection. Here $q$ is a constant projection on
$[{{x^0_{n+1}}},{{x^0_{n+1}}}+\frac12]$ and
$Q_1|_{[{{x^0_{n+1}}},{{x^0_{n+1}}}+\frac12]}=q\otimes \e_{[n,1]}$. Let $r\in
M_m({\blue{Q_1A_{X,n+1}^1Q_1}})= {\blue{Q_1A_{X,n+1}^1Q_1}}\otimes M_m$ be defined as {{follows:}}
\beq\label{13Dec23-2018-2}
r({{x^0_{n+1}}})=q({{x^0_{n+1}}})\otimes \left(
                         \begin{array}{cc}
                           \e_{m-1} & 0 \\
                           0 & 0 \\
                         \end{array}
                       \right)\qq\qq\qq\qq\qq\qq\qq\qq\\
                       \nonumber= (q_1({{x^0_{n+1}}})+q_2({{x^0_{n+1}}})+\cd +q_d({{x^0_{n+1}}})+p({{x^0_{n+1}}}))\otimes \left(
                         \begin{array}{cc}
                           \e_{m-1} & 0 \\
                           0 & 0 \\
                         \end{array}
                       \right);\eneq
for $t\in [\frac12, 1]$,
$$r({{x^0_{n+1}}}+t)= \Big(q_1({{x^0_{n+1}}}+t)+q_2({{x^0_{n+1}}}+t)+\cd +q_{13(m-1)}({{x^0_{n+1}}}+t)\Big)\otimes \e_m + $$
$$+\Big(q_{13(m-1)+1}({{x^0_{n+1}}}+t)+\cd +q_{d}({{x^0_{n+1}}}+t)\Big) \otimes \left(
                         \begin{array}{cc}
                           \e_{m-1} & 0 \\
                           0 & 0 \\
                         \end{array}
                       \right);$$
and, for $ x\in X_{n+1}'\subset X_{n+1}$,
 $$ r(x)= (q_1(x)+q_2(x)+\cd +q_{13(m-1)}(x))\otimes \e_m + (q_{13(m-1)+1}(x)+\cd +q_{d}(x)) \otimes \left(
                         \begin{array}{cc}
                           \e_{m-1} & 0 \\
                           0 & 0 \\
                         \end{array}
                       \right).$$
In the {{definition}} above, between ${{x^0_{n+1}}}$ and ${{x^0_{n+1}}}+\frac12$, $r(t)$ can
 be defined to be any continuous path connecting the projections
 $r({{x^0_{n+1}}})$ and $r({{x^0_{n+1}}}+\frac12);$ {{note that}} both have  {{rank}}
 $(d+13)(m-1)=13(m-1)m+(d-13(m-1))(m-1)$. (Note that all
 $q_i(x) $ are constant on $x\in X_{n+1}=[{{x^0_{n+1}}},{{x^0_{n+1}}}+1]\vee X_{n+1}'$
   and $p(t)$ is constant for $t\in [{{x^0_{n+1}}}, {{x^0_{n+1}}}+1]$.)
   Note that for $x \in [{{x^0_{n+1}}}+\frac12,{{x^0_{n+1}}}+1]\vee X_{n+1}'$,
    $r(x)$ has the same form as $r({{x^0_{n+1}}}+\frac12)$ which is {\blue{a}} constant
    sub-projection
    of {{the}} constant projection $q'\otimes \e_m$. {{Hence $r$ is a trivial projection.}} We will  define
    $R^1$
    to be $r\otimes \e_{[n,i]}$
    under {{a certain}} identification described below.  Note that the
     projection $Q_1$ is identified with $q\otimes \e_{[n,1]}$
     only on {{the}} interval $[{{x^0_{n+1}}}, {{x^0_{n+1}}}+1]${{,
     so}} the definition of $R^1$ will be
     divided into two parts. For the part on $[{{x^0_{n+1}}}, {{x^0_{n+1}}}+\frac12]$, we
     use the identification of $Q_1$ with $q\otimes \e_{[n,1]}$, and for the
      part
      that $x\in [{{x^0_{n+1}}}+\frac12,{{x^0_{n+1}}}+1]\vee X_{n+1}'$, we use
       the identification of $Q'=\e_{d[n,1]}$ with
       $q'\otimes \e_m$ (of course, we use the fact
       that $r$ is {\blue{a}} sub-projection of $q'$
       on this part). This is the only difference between the proof of this lemma
       and that of Lemma \ref{range 0.22}. The definition of
       $ \iota: M_{m-1}({\blue{\pi_{X,n+1}^e(Q_1A_{X, n+1}^1Q_1)}})
       \hookrightarrow R^1M_m(Q_1 {{A_{X,n+1}^1}} Q_1)R^1$ and
      {{the}} verification
       that $\iota$ and ${\blue{R^1}}$ satisfy the conditions are {{exactly}}  the same
       as in the proof of \ref{range  0.22}, with $(d+1)(m-1)$ replaced by $(d+13)(m-1)$. {{We will now give  the details.}} 

       {{On $[x^0_{n+1}, x^0_{n+1}+\frac{1}{2}]$, $Q_1\otimes \one_m\in M_m(Q_1({ A}^1_{X,n+1}|_{[x^0_{n+1}, x^0_{n+1}+\frac{1}{2}]})Q_1)$ is identified with $(q\otimes \one_m)\otimes \one_{[n,i]}$
    (as $Q_1$ is identified with $q\otimes \one_{[n,i]}$).
     As $r$ is a sub-projection of $q\otimes \one_m$, we can define $R=r\otimes \one_{[n,i]}$ as  a sub-projection of $Q_1\otimes \one_m$ on $[x^0_{n+1}, x^0_{n+1}+\frac{1}{2}].$
     {\blue{N}}ote that on $[x^0_{n+1}, x^0_{n+1}+\frac{1}{2}]$, $R$ commutes with $\LD\circ {\tilde \psi^1}$
     (recall ${\tilde \psi^1}=\psi^1|_{[x^0_{n+1}, x^0_{n+1}+\frac{1}{2}]}$), since $\td{\psi}^1=(\LD_2\otimes \id_{[n,1]})\circ \pi,$ and $R=r\otimes \one_{[n,1]},$ and $r$ commutes {{with}} $\mbox{range} (\LD\circ \LD_2)=\C\cdot q\otimes \one_m$, as $r$ is a sub-projection of $q\otimes \one_m$. }}
     {{On $[x^0_{n+1}+\frac{1}{2}, x^0_{n+1}+1]\vee X_{n+1}'$, $Q'\otimes \one_m\in M_m(Q_1{ A}^1_{X, n+1}|_{[x^0_{n+1}+\frac{1}{2}, x^0_{n+1}+1]\vee X_{n+1}'})Q_1$ is identified with $(q'\otimes \one_m)\otimes \one_{[n,i]}$ (as $Q'$ is identified with $q'\otimes \one_{[n,i]}$). As $r$ is a sub-projection of $q'\otimes \one_m$, we can define $R^1=r\otimes \one_{[n,i]}$ as a sub-projection of $Q'\otimes \one_m$ on $[x^0_{n+1}+\frac{1}{2}, x^0_{n+1}+1]\vee X_{n+1}'$. On this part, $R^1$ commutes with $\LD\circ { \psi_1}$ since ${\psi}_1=(\LD_1\otimes \id_{[n,1]})\circ \pi$ and $R^1=r\otimes \one_{[n,1]}$ and $r$ {{commutes}} {\blue{with}}}}
       $\mbox{range} (\LD\circ \LD_1)=\C\cdot q'\otimes \one_m$, as $r$ is a sub-projection of $q'\otimes \one_m$. Note that $R^1$ is {\blue{a}} sub-projection of $Q'\otimes \one_m$ and therefore
       {\blue{is orthogonal}} to $\td{Q}\otimes \one_m$ and  {\blue{the range}}
       of $\LD\circ \psi_2$. Hence on this part,  $R^1$ also commutes with $\LD\circ \psi^1$ as $\psi^1=\psi_1+\psi_2$. {{On combining this}} with {\blue{the}} previous paragraph, (i) follows.

On the other hand, (ii) follows from  {\blue {the}} definition of $R^1$ and (\ref{13Dec23-2018-2}).

      {{From {\blue{the}} last paragraph, we know that $ { A}^1_{X,n} \ni f \mapsto R^1(\LD\circ \psi^1(f))R^1\in R^1(M_m(Q_1{ A}^1_{X,n+1}Q_1))R^1$ is a homomorphism{{.  We}} denote it by $\Xi$.
      From (\ref{13Dec-23-2018}) and (\ref{13Dec-23-2018-1}), we know that $\Xi$ {\blue{factors}} through $F^1_n={ A}^1_{X,n}|_{x^0_n}$ as $\Xi=\Xi'\circ \pi$. Also when we identify $R^1=r\otimes \one_{[n,i]}$, {\blue{the map}} $\Xi':M_{[n,1]}\to R^1(M_m(Q_1{ A}^1_{X,n+1}Q_1))R^1$ can be identified as $\xi\otimes \id_{[n,1]}$, where $\xi: \C \to r(M_m(Q_1{ A}^1_{X, n+1}Q_1))r$ is defined by $\xi(c)=c\cdot r$. }}

{{Let $r_{ij}^0, 1\leq i,j\leq (d+13)(m-1)${{, be}} the matrix units
for $M_{(d+13)(m-1)}$.
Since $r$ is a trivial projection of rank $(d+13)(m-1)$, one
can construct
$r_{ij}\in rM_m({\blue{Q_1({ A}_{X,n+1}^1Q_1))}}r$, $1\leq i,j\leq (d+13)(m-1)${{,
 with}} ${\blue{\pi(r_{ij})}}
 =r_{ij}^0$ serving as {\blue{a system of}} matrix units for
 $M_{(d+13)(m-1)}\subset
 rM_m({\blue{Q_1A_{X,n+1}^1Q_1}})r \cong M_{(d+13)(m-1)}(C(X_{n+1}))$.}}
 {\blue{Recall {{that}} $j_{n+1}: F_{n+1}^1\to \pi_{x_{n+1}^0}(A_{X,n+1}^1)$ is  {{an isomorphism}} (see \ref{range 0.19}).
Note that  $q(x_{n+1}^0)$ is a projection in $\pi_{x_{n+1}^0}(A_{X,n+1}^1).$  Set $\hat{q}=j_{n+1}^{-1}(q({x_{n+1}^0})).$
Then $\hat{q}$ has rank $d+13.$
Define}}
 $$\iota_1:  M_{m-1}(\hat{q}F_{n+1}^1\hat{q})
\to
  rM_m({\blue{Q_1A_{X,n+1}^1 Q_1}})r
  $$
  by $\iota_1(r_{ij}^0)=r_{ij}.$
  {\blue{Note that $M_{m-1}(\pi_{X,n+1}^e(Q_1 A_{X,n+1}^1Q_1))=M_{m-1}(\hat{q}F_{n+1}\hat{q})\otimes M_{[n,1]}.$}}
  Define $\iota=\iota_1\otimes {\blue{\id_{[n,1]}}}: {\blue{M_{m-1}(\pi_{X,n+1}^e(Q_1 A_{X,n+1}^1Q_1))\to R(M_m(Q_1A_{X,n+1}^1Q_1))R.}}$

  {{  Note that $\mbox{range } (\xi)=\C\cdot r\subset \iota_1(M_{m-1}(
  \hat{q}F_{n+1}^1\hat{q}))=
    \mathrm{range }(\iota_1)$. {\blue{S}}ince $\Xi'=\xi\otimes \id_{[n,i]}$ and  $\iota=\iota_1\otimes {\blue{\id_{[n,1]}}}$, {\blue{we}} have $\mbox{range } (\Xi')\subset \mathrm{range }(\iota).$  {\blue{H}}ence $R(\LD (\psi^1({ A}^1_{X,n})))R=\Xi({ A}^1_{X,n})=\Xi'(F^1_n) \subset \mathrm{range }(\iota)$.}}

 \end{proof}



\begin{NN}\label{13-Jan31-2019}

 {\blue{Recall that we have constructed $A_{X, n+1}$ and the unital injective \hm\, $\psi_{n,n+1}$
 (with {{the}} specific {{form described}}   in \ref{range  0.23}).  Note that,  {\blue{by}} the end of \ref{condition2},
 we {{know that}} $A_{C,n+1}$ is fixed and  the injective map ${\bar \phi}_{n,n+1}: A_{C, n}\to A_{C, n+1}$ is defined.
 One also has $A_{n+1}$ as defined {\blue{in}}  \eqref{13Dec15-2018}.
 As in \ref{range 0.20}, $(\phi_{n, n+1}^q)_{*0}=\gm_{n,n+1}'=(\psi_{n,n+1}^q)_{*0}.$ {{Therefore,}} there  exists a unitary
$U\in F_{n+1}$ such that ${\rm Ad}\, U\circ \psi_{n, n+1}^q={\phi}_{n,n+1}^q.$
In $F_{n+1},$ there exist $h\in  (F_{n,n+1})_{s.a.}$ such that $U=\exp(ih).$
Let $H\in (A_{X,n+1})_{s.a.}$ {{be}}  such that $\pi_{X,n+1}^e(H)=h.$ Define $V=\exp(iH).$
Note that $V^*{\bar J}_{n+1}V\subset {\bar J}_{n+1}.$
If we replace $\psi_{n, n+1}$ by ${\rm Ad}\, V\circ \psi_{n,n+1},$ then we still have (1) and (2) {{of}}  \ref{range 0.20}.
Moreover, \ref{range 0.22}  and \ref{range 0.24}   also hold (up to unitary equivalence in $A_{X,n+1}$). More importantly, ${{\psi_{n, n+1}^q}}=\phi_{n,n+1}^q.$  It follows from \ref{construction}    that there is
a unital  injective \hm\, $\phi_{n,n+1}: A_n\to A_{n+1}$ which satisfies (a), (b){{, and}} (c) {{of}} \ref{construction}.
Moreover, $A_{n+1}$ satisfies (d) of \ref{condition2}.
Therefore, by \ref{13-Feb17-2019},  one checks  {\blue{that}} $\phi_{n,n+1}$ also satisfies (a), (b){{, and}} (c)  {{of}} \ref{condition2}.
This ends the induction {{of}} \ref{condition2}.

 Now we have
 \begin{displaymath}
    \xymatrix{
        A_1 \ar[r]^{\phi_{{1,2}}} & A_2 \ar[r]^{\phi_{{2,3}}}&A_3 \ar[r]^{\phi_{{3,4}}} &\cd \ar[r]^{\phi_{i-1,i}} & A_i \ar[r]\cd & A,}
\end{displaymath} satisfying conditions (a)--(d) {{of}} \ref{condition2}.
From (a) of \ref{condition2}, we obtain {{the}}  two inductive limits
 \begin{displaymath}
    \xymatrix{
        A_{C,1} \ar[r]^{{\bar\phi}_{{1,2}}} & A_{C,2} \ar[r]^{{\bar\phi}_{{2,3}}}&A_{C,3} \ar[r]^{{\bar\phi}_{{3,4}}} &\cd \ar[r]^{\bar{\phi}_{i-1,i}} & A_{C,i} \ar[r]\cd & {{A_C}}~\mbox{and}}
        \end{displaymath}
        \begin{displaymath}
    \xymatrix{
        A_{X,1} \ar[r]^{\psi_{{1,2}}} & A_{X,2} \ar[r]^{\psi_{{2,3}}}&A_{X,3} \ar[r]^{\psi_{{3,4}}} &\cd \ar[r]^{\psi_{i-1,i}} & A_{X,i} \ar[r]\cd & A_X.}
\end{displaymath}}}

 Since
 $\psi_{n,n+1}({{\bar J}}_n)\subset {{\bar J}}_{n+1}$, this procedure also gives
 an
 inductive limit of quotient algebras ${{F}}=\lim ({{F}}_n, \psi^{{q}}_{{n,n+1}})$,
 where ${{F}}_n={{A_{X,n}}}/{{\bar J}}_n$. Evidently, ${{F}}$ is an AF algebra
 with $K_0({{F}})= H/{\rm Tor}(H)$. {\blue{(Note that $\psi^q_{n,n+1}={\bar \phi}^q_{n,n+1}: F_n\to F_{n+1}$.)}}
{\blue{Since $\phi_{n,n+1}(I_n)\subset I_{n+1},$  $\phi_{n,n+1}(J_n)\subset J_{n+1},$
$\{\pi_{I,n}\}$ {{and}}  $\{\pi_{J,n}\}$ induce the quotient maps $\pi_I: A\to A_X$ and $\pi_J: A\to  {{A_C}},$ respectively.
Moreover, $\{\pi_{X,n}^e\}$ and $\{\pi_{C,n}^e\}$ induce the quotient maps
$\pi_X^e: A_X\to F,$ and $\pi_C^{{e}}: {{A_C}}\to F,$ respectively.
{{Finally,}} $\pi_A^e=\pi_X^e\circ {{\pi}}_I=\pi_C^e\circ \pi_J: A\to F${{, and this map is}}  induced by $\{\pi_{A_n}^e\}.$}}

 Combining \ref{range  0.22} and \ref{range  0.24}, we have the following
 theorem which  will be
  used to conclude that the algebra $A$ ({{which}}  will be constructed later) {{has the property }} that
   $A\otimes U$ is in ${\cal B}_0.$

\end{NN}

\begin{thm}\label{range 0.25}
  Suppose that $1< m\leq {\min} \{{{(c_{11}-13)/13}}, c_{12}-1, c_{13}-1,..., c_{1p_n}-1\}$. Let $\psi: {{A_{X,n}}} \to
  {{A_{X,n+1}^1}}$ be the composition
 \begin{displaymath}
\xymatrix{{{A_{X,n}}}\ar[r]^{\psi_{n,n+1}}& {{A_{X,n+1}}} \ar[r]^{{\blue{\pi^1}}} &
{{A_{X,n+1}^1,}}}
\end{displaymath}
{\blue{where}} ${\blue{\pi^1}}$ is the quotient map to the first block. Let $\LD:
{{A_{X,n+1}^1}} \to M_m({{A_{X,n+1}^1}})$ be
the amplification {{defined }}  by
$\LD(a)=a\otimes {\bf1}_m.$
There is a projection $R\in
M_m({{A_{X,n+1}^1}})={{A_{X,n+1}^1}}\otimes M_m$ and there is a unital inclusion
 homomorphism
 $ \iota: M_{m-1}({{F}}_{n+1}^1) ={{F}}_{n+1}^1\otimes M_{m-1} \hookrightarrow RM_m( {{A_{X,n+1}^1}})R${{, satisfying}} the following {{three}} conditions:

{\rm (i)} $R$ commutes with $\LD(\psi({{A_{X,n}}}))${{, and}}

{\rm (ii)} $R({{x^0_{n+1}}})= \e_{{{F}}_{n+1}^1}\otimes \left(
                                        \begin{array}{cc}
                                          \e_{m-1} & 0 \\
                                          0 & 0 \\
                                        \end{array}
                                      \right).$

Consequently, the map $\pi: {{A_{X,n+1}^1}} \to {{F}}_{n+1}^1$ takes
$RM_m( {{A_{X,n+1}^1}} )R$ onto  $M_{m-1}({{F}}_{n+1}^1)$ .

{{Below we will use the same notation $\pi$ to denote {{the}} restriction {{of $\pi$}} to $RM_m(A_{X,n+1}^1)R$, whose codomain is $M_{m-1}(F_{n+1}^1)$.}}

{\rm (iii)}  $\pi\circ \iota=\id|_{M_{m-1}({{F}}_{n+1}^1)}$,  and $R(\LD
(\psi({{A_{X,n}}})))R \subset \mathrm{range} (\iota)$.

\end{thm}

\begin{proof} Choose $R=\bigoplus_{i=1}^{p_n} R^i\in
M_m({{A_{X,n+1}^1}})$, where $R^1\in
M_m(Q_1{{A_{X,n+1}^1}}Q_1)$ is {{as described }} in Lemma  \ref{range 0.24} and $R^i\in M_m(Q_i{{A_{X,n+1}^1}}Q_i)$ (for $i\geq 2$) is  {{as described }} in Lemma \ref{range 0.22}.
Then the theorem follows.
\end{proof}

\begin{thm}\label{range 0.27}
{\blue{Let
$(H, H_+, u)$  be as in \ref{range 0.16} and \ref{range 0.17}.}}
{\blue{Then
$$
((K_0({{A_{X}}}),K_0({{A_{X}}})_+, [{\mbox{\large \bf 1}}_{{A_{X}}}]),K_1({{A_{X}}}))
\cong
((H,H_+,u),K).
$$}}
\end{thm}
\begin{proof} Note that $\psi_{n,n+1}$ satisfies (1) in \ref{range 0.20}, {{and,}}  consequently,
$$ ((K_0({{A_{X}}}),K_0({{A_{X}}})_+, [{\mbox{\large \bf 1}}_{{A_{X}}}]),K_1({{A_{X}}}))\cong
((H,H_+,u),K, ).$$
\end{proof}
Let ${\blue{\pi_{X,n}^e}}: {{ A_{X,n}}}\to  {{F}}_n$ be
{{as}} in \ref{construction}. Then ${\blue{(\pi_{X,n}^e)^{\sharp}}}: \Aff(T({{ A_{X,n}}}))=C(X_n, \R)\oplus \R^{p_n-1} \to \Aff (T ({{F}}_n)) =\R^{p_n}$
is given by ${\blue{{(\pi_{X,n}^e)^{\sharp}}}}(g,h_2,h_3,...,h_{p_n})=(g(\tht_1), h_2,h_3,..., h_{p_n})$. Define $\GM_n: \Aff (T ({{F}}_n)) =\R^{p_n} \to \Aff (T ({{ A_{X,n}}}))=C(X_n, \R) \oplus \R^{p_n-1}$ to be the right inverse of $(\pi_{X,n}^e)^{\sharp}${{, given}} by $\GM_n(h_1, h_2,h_3,..., h_{p_n})=(g,h_2,h_3,...,h_{p_n})$ with $g$ {{the}} constant function $g(x)=h_1$ for all $x\in X_n$. Then with the condition $c_{ij}> 13\cdot 2^{2n}$, we have the following {\blue{lemma}}:

\begin{lem}\label{Lemaff}
For any $f\in \Aff (T({{ A_{X,n}}}))$ with $\|f\|\leq 1$ and $f'{{:=}}\psi_{n,n+1}^{\sharp} (f)\in \Aff (T({{ A_{X,n+1}}}))$, we have
$$\|\GM_{n+1}\circ{\blue{(\pi_{X,n+1}^e)^{\sharp}}}(f')-f'\|<\frac{2}{2^{2n}}.$$
\end{lem}

\begin{proof}
Write $f=(g, h_2,..., h_{p_n})$ and $f'=(g', h'_2,..., h'_{p_{n+1}})$.
Then $\GM_{n+1}\circ{\blue{(\pi_{X,n+1}^e)^{\sharp}}}(f')=(g'', h'_2,..., h'_{p_{n+1}})$ with \beq\label{13Dec21-2018}
g''(x)=g'({{x^0_{n+1}}})~~~~ \mbox{for all} ~~~~~x\in X_{n+1}.\eneq
Recall that $\psi_{n,n+1}^{i,1}$ is denoted by $\psi^i: {{ A_{X,n}}} \to Q_i{{ A_{X,n+1}^1}}Q_i$ and $\psi^1=\psi_1+\psi_2$ with $\psi_1: {{ A_{X,n}^1}} \to Q'{{ A_{X,n+1}^1}}Q' ,$  and $\psi_2:
{{ A_{X,n}^1}} \to {\td Q}{{ A_{X,n+1}^1}}{\td Q} $ as in \ref{range 0.21} and \ref{range 0.23}. Note that
$$\frac{{\rm rank}(Q_i)}{{\rm rank}(P_{n+1})}=\frac{c_{i,1}}{\sum_{j=1}^{p_n}c_{j,1}},~~~~\frac{{\rm rank}(Q')}{{\rm rank}(P_{n+1})}=\frac{c_{1,1}-13}{\sum_{j=1}^{p_n}c_{j,1}},~~~~{\mbox{and}}~~~\frac{{\rm rank}({\td Q})}{{\rm rank}(P_{n+1})}=\frac{13}{\sum_{j=1}^{p_n}c_{j,1}}.$$
Hence
$$g'=\frac{c_{1,1}-13}{\sum_{j=1}^{p_n}c_{j,1}}\psi_1^{\blue{\sharp}}(g)+
\frac{13}{\sum_{j=1}^{p_n}c_{j,1}}\psi_2^{\blue{\sharp}}(g)+
\sum_{i=2}^{p_n}\frac{c_{i,1}}{\sum_{j=1}^{p_n}c_{j,1}}(\psi^i)^{\blue{\sharp}}(h_i).$$
Also from the construction in \ref{range 0.21} and \ref{range 0.23}, we know that $\psi_1^{\blue{\sharp}}(g)$ and $(\psi^i)^{\blue{\sharp}}(h_i)$ ($i\geq 2$) are constant. So we have
$$|g'(x)-g'({{x^0_{n+1}}})|\leq \frac{2\times 13}{\sum_{j=1}^{p_n}c_{j,1}} < \frac{2}{2^{2n}}.$$
Then the lemma follows from
{{from (\ref{13Dec21-2018})}}.
\end{proof}

Using Lemma \ref{Lemaff} one can actually prove (see the proof of \ref{EllofA}) that \\
$((K_0({{A_{X}}}),K_0({{A_{X}}})_+, [{\mbox{\large \bf 1}}_{{A_{X}}}]),K_1({{A_{X}}}),T({{A_{X}}}), r_{{A_{X}}})$
is isomorphic to $((H,H_+,u),K, \Delta,r).$

\begin{NN}\label{range 0.28+Dec-2018}

{{We will also  {\blue{compute}} the tracial state space for the $C^*$-algebra $A$ in \ref{13-Jan31-2019}.
}}

{{As in}} \ref{range 0.8} (see \eqref{1508/13star-2}), {{the subspace}}
$$\Aff(T(A_n)) \sbs \bigoplus_{i=1}^{l_n}C([0,1]_i,\R)\oplus C(X_n, \R)\oplus \R^{p_n-1}$$
{consists} of {{ the elements}} $(f_1,f_2,... ,f_{l_n}{{;}}~g,h_2,... ,h_{p_n})$  {{(here we do not need $h_1$, since it is identified with $g(x^0_n)$)}} which {{satisfy}}
the conditions
\beq\label{1star'}
\qq\qq\qq\qq f_i(0)&=&\frac1{\{n,i\}}\big(b^{{n}}_{{0,i1}}g({{x^0_n}})[n,1]+\sum_{j=2}^{p_n} b^{{n}}_{0,ij}h_j\cdot [n,j]\big)\,\,\,\,\,\,\,\,\,\,\andeqn\qq\qq\,\,\,\,
\\\label{2star'}
%
\qq\qq\qq\qq f_i(1)&=&\frac1{\{n,i\}}\big(b^{{n}}_{{1,i1}}g({{x^0_n}})[n,1]+
\sum_{j=2}^{p_n}b^{{n}}_{1,ij}h_j\cdot [n,j].\big.\qq\qq\qq\qq\,
\eneq

For $h=(h_1,h_2,...,h_{p_n})\in \Aff T({{F}}_n)$, let
$\GM_n'(h)(t)= t\cdot \bt_{{n,1}}^{\sharp}(h)+(1-t)\cdot \bt_{{n,0}}^{\sharp}(h)$ {{(see \ref{range 0.9} and}} {{ \ref{range 0.10}),}} which
gives an element $C([0,1],\R^{l_n})=\bigoplus_{i=1}^{l_n}C([0,1]_i,\R)$.
{\blue{Let}}
$$
\GM_n:~ \Aff(T({{F}}_{n}))=\R^{p_n}\to \Aff(T(A_{n})) \subset
\bigoplus_{j=1}^{l_n}C([0,1]_j,\R)\oplus C(X_n, \R)\oplus \R^{p_n-1}
$$
be defined by
$
\GM_n(h_1,h_2,..., h_{p_n})=(\GM_n'(h_1,h_2,..., h_{p_n}), g,h_2,..., h_{p_n} )\in \Aff(T(A_n)),
$
where $g\in C(X_n,\R)$ is the constant function $g(x)=h_1$.

\end{NN}

Now we are ready to show that the Elliott invariant of $A$ is as desired.

\begin{thm}\label{EllofA}
Let $A$ be as constructed.
Then
$$\big((K_0(A),K_0(A)_+, [1_A]),\ K_1(A), T(A), r_{\!\!_A} \big) \cong \big((G,G_+,u),~ K, \DT, r \big).$$

\end{thm}

\begin{proof}



\noindent{\rm Claim:} For any $f\in \Aff(T (A_n))$ with $\|f\|\leq 1$ and $f'{{:=}}\phi_{n,n+1}^{\sharp} (f)\in \Aff (T(A_{n+1}))$, we have
$$\|\GM_{n+1}\circ{\blue{\pi_{A_{n+1}}^{e\,\,\sharp}}}(f')-f'\|<\frac{4}{2^{2n}},$$
where ${\blue{(\pi_{A_{n+1}}^e)^{\sharp}}}: \Aff(T (A_{n+1})) \to \Aff(T({{F}}_{n+1}))$ is induced by $\pi_{A_{n+1}}^e$
{\blue{(see \ref{construction}).}}

\noindent Proof of the claim:
For any $n\in\Z_+$, write
\begin{displaymath}
\xymatrix{\GM_n=\GM_n^1\circ\GM_n^2:
\Aff(T({{F}}_n))\ar[r]^{\qq\qq\GM_n^2}& \Aff(T(
{{A_{X,n}}}))  \ar[r]^{\GM_n^1} &
\Aff(T({ A}_n))}
\end{displaymath}
with $\GM_n^2: \Aff(T({{F}}_n))=\R^{p_n} \to
\Aff(T({{A_{X,n}}}))=C(X_{{n}},\R)\oplus \R^{p_{n-1}}$
defined by $$\GM_n^2(h_1,h_2,..., h_{p_n})=(g,h_2,..., h_{p_n}),$$
 where $g$ is the constant function $g(x)=h_1$, and with
 $\GM_n^1: \Aff(T({{A_{X,n}}}))=C(X_{{n}},\R)\oplus \R^{p_{n-1}} \to \Aff(T({ A}_n))$ defined by $$\GM_n^1(g,h_2,..., h_{p_n})=(\GM_n'(g({\blue{x_n^0}}),h_2,..., h_{p_n}), g,h_2,..., h_{p_n} ).$$
 {\blue{Define $\Gamma_{n,C}: \Aff(T(F_n))\to \Aff(T(A_{C,n}))$
 by
 $$
 \Gamma_{n,C}(h_1,h_2,...,h_{p_n})=(\Gamma_n'(h_1,h_2,...,h_{p_n}), h_1,h_2,...,h_{p_n}).
 $$
 }}
{\blue{Note also $\pi_{A_n}^{e\,\,\sharp}=\pi_{X,n}^{e\,\,\sharp}\circ \pi_{I,n}^{e\,\,\sharp}$ (see \ref{construction}).}}

{\blue{One has, for any $(f,g,h_2,h_3,...,h_{p_n})\in \Aff(T(A_n)),$
\beq\nonumber
&&\hspace{-0.4in}\pi_{J,n}^{\sharp}\circ \Gamma_{n}^1 \circ \pi_{I,n}^{\sharp}(f,g,h_2,h_3,...,h_{p_n})
=\pi_{J,n}^{\sharp}\circ \Gamma_{n}^1(g,h_2,...,h_{p_n})\\\nonumber
&&=\pi_{J,n}^{\sharp}\circ (\Gamma_n'(g({\blue{x_n^0}}), h_2,...,h_{p_n}), g, h_2,...,h_{p_n})\\\nonumber
&&=(\Gamma_n'(g({\blue{x_n^0}}), h_2,...,h_{p_n}), g(x_n^0), h_2,...,h_{p_n})
=\Gamma_{n,C}(g({\blue{x_n^0}}), h_2,...,h_{p_n})\\\nonumber
&&=\Gamma_{n,C}\circ \pi_{C,n}^{e\,\,\sharp}(f, g({\blue{x_n^0}}), h_2,h_3,...,h_{p_n})
=\Gamma_{n,C}\circ  \pi_{C,n}^{e\,\,\sharp}\circ \pi_{J,n}(f,g,h_2,...,h_{p_n}).
\eneq}}
{\blue{In other words,  for all $n\in \Z_+,$
\beq\label{Aff-1340}
\pi_{J,n+1}^{\,\,\sharp}\circ \Gamma_{n+1}^1\circ \pi_{I,n+1}^{\,\,\sharp}=\Gamma_{n+1,C}\circ \pi_{C,n+1}^{\,\,\sharp}\circ\pi_{J,n+1}^{\,\,\sharp}.
\eneq}}
{\blue{One also checks that
\beq\label{Aff-1340-2}
\pi_{I,n+1}^{\,\,\sharp}\circ \Gamma_{n+1}^1\circ \pi_{I, n+1}^{\,\,\sharp}=\pi_{I,n+1}^{\,\,\sharp}.
\eneq}}


For any $f\in \Aff(T({ A}_n))$ with $\|f\|\leq 1$, write $f_1=\pi_{I,n}^{\sharp}(f)$, $f'=\phi_{n,n+1}^{\sharp} (f)$ and $f'_1=\psi_{n,n+1}^{\sharp} (f_1)$.
By the condition $c_{ij}>13\cdot2^{2n}$ and  the claim in the proof of Theorem \ref{range 0.27}, we have
\beq\label{0.28d}
\|\GM_{n+1}^2\circ\pi_{X,n+1}^{e\,\,\sharp}(f'_1)-f'_1\|<\frac{2}{2^{2n}}.
\eneq
{{Using}} {{condition $\spdd$,}} 
 applying Lemma \ref{range 0.10} to {{the map}} $A_{{C},n} \to A_{{C},n+1}$ as $C_n \to C_{n+1}$ (note that ${\blue{\GM_{n,C}}}$ is the same as  ${\blue{\GM_n}}$ in \ref{range 0.9}  and $\bar \phi_{n,n+1}: A_{{C},n} \to A_{{C},n+1}$ is the same as $\phi_{n,n+1}:C_n \to C_{n+1}$), {\blue{one has
\beq\label{Aff-1340-3}
\|\Gamma_{n+1,C}\circ \pi_{C,n+1}^{e\,\,\sharp}\circ {\bar \phi}_{n,n+1}^{\,\,\sharp}(\pi_{J,n}^{\sharp}(f))-{\bar \phi}_{n,n+1}^{\,\,\sharp}(\pi_{J,n}^{\,\,\sharp}(f))\|<\frac{2}{2^{2n}}.
\eneq}}
{\blue{Note that ${\bar \phi}_{n,n+1}^{\,\,\sharp}\circ \pi_{J,n}^{\sharp}=\pi_{J,n+1}^{\,\,\sharp}\circ \phi_{n,n+1}^{\,\,\sharp}.$
By \eqref{Aff-1340},
\beq\label{Aff-1340-4}
\|\pi_{J,n+1}^{\,\,\sharp}\circ \GM_{n+1}^1\circ \pi_{I,n+1}^{\sharp}(f')-\pi_{J,n+1}^{\,\,\sharp}(f')\|<\frac{2}{2^{2n}}.
\eneq
Combining this with \eqref{Aff-1340-2}, we have}}
\begin{equation}
\nonumber
\left(\GM_{n+1}^1\circ \pi_{I,n+1}^{\,\,\sharp}(f')-f'\right)|_{\blue{T(A_{X,n+1})}}=0 \andeqn \|\left(\GM_{n+1}^1\circ\pi_{I,n+1}^{\,\,\sharp}(f')-f'\right)|_{\blue{T(A_{C,n+1})}}
\|<\frac{2}{2^{2n}}.
\end{equation}
{{Recall that $Sp(A_{X, n+1})\cup Sp(A_{C, n+1})= Sp(A_{n+1})$ {\blue{(see (\ref{13-March26-2019-3})).}}
{\blue{By Lemma  2.16 of \cite{Lncrell},}} we know that any extreme trace of $A_{n+1}$ {\blue{is induced by}} either {{an}} extreme  trace of $A_{X, n+1}$ or {{an}} extreme trace of $A_{C, n+1}$.}} It follows
that
\begin{equation}\label{0.28f}
\|\GM_{n+1}^1\circ\pi_{I,n+1}^{\sharp}(f')-f'\|<\frac{2}{2^{2n}}.
\end{equation}
Consequently (applying  (\ref{0.28d}) {{and}} (\ref{0.28f})), we {\blue{obtain}}
\begin{eqnarray*}
&&\hspace{-0.3in}\|\GM_{n+1}\circ{\blue{\pi_{A_{n+1}}^{e\,\,\sharp}}}(f')-f'\|
= \|\GM_{n+1}^1\circ\GM_{n+1}^2\circ{\blue{\pi_{X,n+1}^{e\,\,\sharp}}}
\circ\pi_{I,n+1}^{\,\sharp} (f')-f'\| \\
 & &= \|\GM_{n+1}^1\circ\GM_{n+1}^2\circ\pi_{X,n+1}^{e\,\,\sharp}\circ
 \psi_{n,n+1}^{\,\sharp}\circ\pi_{I,n}^{\,\sharp} (f)-f'\| \\
& = & \|\GM_{n+1}^1\circ\GM_{n+1}^2\circ\pi_{X,n+1}^{\,\sharp}\circ\psi_{n,n+1}^{\,\sharp}(f_1)-f'\|
=\|\GM_{n+1}^1\circ\GM_{n+1}^2\circ\pi_{X,n+1}^{e\,\,\sharp}(f'_1)-f'\| \\
& < & \|\GM_{n+1}^1(f'_1)-f'\|+\frac{2}{2^{2n}}\quad\quad(\mbox{by}~(\ref{0.28d})) \\
& = &\|\GM_{n+1}^1\circ\psi_{n,n+1}^{\sharp}\circ\pi_{I,n}^{\sharp}(f)-f'\|+
\frac{2}{2^{2n}}
=\|\GM_{n+1}^1\circ(\pi_{I,n+1})^{\sharp}\circ\phi_{n,n+1}^{\sharp}(f)-f'\|+\frac{2}{2^{2n}} \\
& = &\|\GM_{n+1}^1\circ(\pi_{I,n+1})^{\sharp}(f')-f'\|+\frac{2}{2^{2n}}
 < \frac{2}{2^{2n}}+\frac{2}{2^{2n}}.
\end{eqnarray*}
This proves the claim.

Using the claim, one {\blue{obtains}} the following {{approximate}} intertwining diagram:
\begin{displaymath}
    \xymatrix{
        \Aff(T(A_1)) \ar@/_/[d]_{{\blue{\pi^{e\, \sharp}_{A_1}}}}\ar[r]^{\phi_{1,2}^{\blue{\sharp}}} & \Aff(T(A_2)) \ar[r]^{\phi_{2,3}^{\sharp}} \ar@/_/[d]_{\blue{\pi_{A_2}^{e\,\sharp}}}& \Aff(T(A_3)) \ar[r] \ar@/_/[d]_{\blue{\pi_{A_3}^{e\,\sharp}}}& \cd \Aff(T({{A}} ))\\
        \Aff(T({{F}}_1)) \ar[r]^{{\psi}_{1,2}^{{q},\sharp}}\ar@/_/[u]_{\GM_1} &
         \Aff(T({{F}}_2)) \ar[r]^{{\psi}_{2,3}^{{q},\sharp}}\ar@/_/[u]_{\GM_2}&
         \Aff(T({{F}}_3)) \ar[r]\ar@/_/[u]_{\GM_3}& \cd
         \Aff(T({{F}})). }
\end{displaymath}
{\blue{Recall {{that}} $\pi_A^e=\pi_X^e\circ \phi_I=\pi_C^e\circ \pi_J: A\to F$ is induced by $\{\pi_{A_n}^e\}.$
{{Thus,}} the above approximately intertwining diagram  shows that $\pi_A^{e\, \sharp}: \Aff(T(A))\to \Aff(T(F))$
is an isometric isomorphism of Banach spaces which also preserves the order and the inverse also preserves the order.
By \ref{Aq}, this implies the inverse $(\pi_A^{e\,\sharp})^{-1}$ induces an affine homeomorphism from $\Delta$ onto
$T(A).$
Note that we also have
\beq\label{Aff-1340-6}
\pi_A^{e\, \sharp}\circ \rho_A=\rho_F\circ (\pi_A^e)_{*0}.
\eneq
By \ref{13-Jan31-2019}, $(\pi_A^e)_{*0}=\pi_{G, H/{\rm Tor}(H)}.$
Therefore $\pi_A^{e\,\,\sharp}\circ \rho_A=\rho: G\to \Aff(\Delta)$ is given by}} {{the part (2) of Remark \ref{March-22-2019}.}}
{\blue{It then follows from the fact that $\pi_A^{e\,\sharp}$ is an affine homeomorphism that
$$\big((K_0(A),K_0(A)_+, [1_A]),\ K_1(A), T(A), r_{\!\!_A} \big) \cong \big((G,G_+,u),~ K, \DT, r \big).$$}}
\end{proof}


\begin{NN}\label{range 0.29}
The algebra $A$ of \ref{13-Jan31-2019}
{{and \ref{range 0.28+Dec-2018} }} is not simple, {{and so}} we need to modify
the \hm s  $\phi_{{n,n+1}}$ to make the limit algebra simple. Let us
emphasize that every \hm\,  $\phi: A_n\to A_{n+1}$ is completely
determined by $\phi_x=\pi_{x}\circ \phi$ for each $x\in
Sp(A_{n+1})$, where the map $\pi_x:~ A_{n+1}\to A_{n+1}|_{x}$ is
the corresponding irreducible representation.

Note that from the definition of $\phi: A_n\to A_{n+1}$ and from the assumption that
$c_{ij}>13$ for each entry of $\cc=(c_{ij})$, we know that for any
$x\in Sp(A_{n+1})$,
\beq\label{150102-p1}
 Sp(\phi_{{n,n+1}}|_x)\supset
 Sp({{F}}_n)=(\tht_{{n,1}},\tht_{{n,1}},...,\tht_{{n,p_n}}).
 \eneq
(See \ref{homrestr} and \ref{ktimes} for {{notation}}.)

To make the limit algebra simple, we will {\blue{change $\phi_{n,n+1}$ to $\xi_{n,n+1}$ so that}} the set
$Sp(\xi_{n,m}|_x)$  {{is}}  sufficiently dense
in $Sp(A_n)$, for any $x\in
Sp(A_m)$, provided $m$ {{is}} large enough.

Write
$Sp(A_n)={{\bigcup}}_{j=1}^{l_n}(0,1)_{{n,j}}\cup X_n\cup S_n,$
where $S_n=Sp(\bigoplus_{i=2}^{p_n}M_{[n,i]}).$
Let $0<d<1/2$ and let $Z\subset Sp(A_n)$ be a subset.

{\it Recall from \ref{density}  {{that}}
$Z$ is $d$-dense in $Sp(A_n)$ if the following (sufficient) condition
holds:
$Z\cap (0,1)_{{n,j}}$ is $d$-dense in $(0,1)_{n,j}$ with the usual metric,
$Z\cap X_n$ is $d$-dense in $X_n$ with a given metric of {{$X_n$,}} and
{{$Sp(F_n)\subset  Z.$}} }

{\blue{For each {{fixed}} $n,$ let $\{{\cal F}_{n,k}:1\le k<\infty\}$ be an increasing sequence of finite subsets of {{the unit ball of }} $A_n$
such that ${{\bigcup}}_{k=1}^{\infty} {\cal F}_{n,k}$ is dense in {{the unit ball of}} $A_n.$}}
{{Fix a sequence of positive numbers $\et_1>\et_2>\cdots>\et_n\cdots>1$ such that $\prod_{n=1}^\infty \et_n< 2$. }}  Now we will change $\phi_{{n,n+1}}$ to a {\blue{unital injectve \hm\,}} $\xi_{{n,n+1}}:
A_n\to A_{n+1}$ satisfying:

(i) $(\xi_{{n,n+1}})_{*i}=(\phi_{{n,n+1}})_{*i}.$

(ii)\,
{\blue{$\|\phi_{{n,n+1}}^{\sharp}({{\hat{f}}}))-\xi_{{n,n+1}}^{\sharp}({{\hat{f}}}))\|\leq \frac1{2^{{2n-2}}}\rforal f\in {\cal G}_n,$
where {{$\hat{f}\in \Aff T(A_n)$ is defined by $\hat{f}(\tau)=\tau(f)$,}} ${\cal G}_1={\cal F}_{1,1},$
 ${\cal G}_k={\cal F}_{k,k}\cup (\bigcup_{i=1}^{{k-1}}\phi_{i,k}({{{\cal F}_{i,k}\cup}}{\cal G}_i))\cup
(\bigcup_{i=1}^{{k-1}}\xi_{i,k}({{{\cal F}_{i,k}\cup}}{\cal G}_i)),$ $k=2,3,...,n.$}}

(iii)\,
{\blue{for any $y\in Sp(A_{X,n+1})\,(\subset Sp(A_{n+1})),$ }} {{$Sp(\xi_{n,n+1}|_y)\supset Sp(\phi_{n,n+1}|_y)$;}}
{\blue{and for any $y\in Sp(A_{n+1}),$ $S{{p}}(\xi_{n,n+1}|_y)\supset Sp(F_{{n}}).$}}


(iv)\,{\blue{For any $\dt>0$ and
 any finite subset $S{{\subset Sp(A_{C, n+1})\subset Sp(A_{n+1})}},$ {{if $S\supset Sp(F_{n+1})$ and {{$S$}} is }}
$d$-dense {{in $Sp(A_{C, n+1})$}}, then
${{\big({{\bigcup}}_{s\in S}Sp(\xi_{n,n+1}|_s)\big)\cap Sp(A_{C,n}) }}$ is  $\eta_n\dt$-dense in  ${{ Sp(A_{C,n}).}}$
}}

((v)\, $Y_n\cup {{T}}_n {{\cup Sp(F_n)}}\subset Sp(\xi_{n,n+1}|_{\theta_{n+1,2}})$ {{(see \ref{condition2} for $Y_n$ and $T_n$)}} and {{consequently,}}
$Sp(\xi_{{{\blue{n}},n+1}}|_{{\tht_{{n+1,2}}}})\ \  \mbox{is  ${{1/n}}$\,-  dense in  } Sp(A_n),$
where $\tht_{{n+1,2}}\in Sp(F^2_{n+1})\sbs
Sp({ {F}}_{n+1})$\\
$=\{\tht_{{n+1,1}},\tht_{{n+1,2}},...,\tht_{{n+1,p_{n+1}}}\}$ is the second
point  of $\{\tht_{{n+1,1}},\tht_{{n+1,2}},...,\tht_{{n+1,p_{n+1}}}\}$
(note that $\tht_{{n+1,1}}$ is {{identified with }}the base point {{$x^0_{n+1}$}} of
$[{{x^0_{n+1}}},{{x^0_{n+1}}}+1]\vee X_{n+1}'=X_{n+1}=Sp({{A_{X,n}^1}})$,
and we do not
want to modify this one).

(vi)\,For any $x\in
Sp({{A_{X,n+1}}})$ satisfying $x\not= \tht_{{n+1,2}},$
\begin{equation}\label{150102-p2-1}
\phi_{{n,n+1}}|_{x} = \xi_{{n,n+1}}|_{x}.
\end{equation}
In
particular we have,
(vi') $\phi_{n,n+1}|_{Sp
({{A_{X,n+1}^1}})}=\xi_{n,n+1}|_{Sp({{A_{X,n+1}^1}})}$,  or equivalently, for any $x\in X_{n+1}=Sp({{A_{X,n+1}^1}}),$
\begin{equation}\label{0.29-1}
\phi_{{n,n+1}}|_{x}
= \xi_{{n,n+1}}|_{x}.
\end{equation}

{\blue{Moreover, we have the following remarks:}}

 {\blue{(vii)   Suppose that (iii)  holds, for $n=1,2,...,k.$
  Let $i<k${{. For any}} $y\in Sp(A_{X,k+1}),$
   { {we have}} $Sp(\xi_{i,k+1}|_y){{\supset Sp(\phi_{i,k+1}|_y)}}$.   {{Furthermore,}}  $Sp(F_i)\subset SP(\xi_{i,{{k}}+1}|_y)$ for all $y\in Sp(A_{{{k}}+1}).$}}

 {\blue{The proof of (vii):}} {{We will prove it by reverse induction. From the assumption that (iii) holds for $n=k$, we know {\blue{that}} both  conclusions
 above  are true for $i=k$. Let us assume both conclusions  hold  for $i=j\leq k$ (and $j\geq 2$), {{i.e.,}} for any $y\in Sp(A_{X,k+1})$, $Sp(\xi_{j, k+1}|_y)\supset Sp(\phi_{j, k+1}|_y)$; and for any $y\in Sp(A_{k+1})$,
 $Sp(F_j)\subset Sp(\xi_{j, k+1}|_y)$. We will prove both conclusions for $i=j-1$. Assume that $y\in Sp(A_{k+1})$. Choose any $z_0\in Sp(F_j)\subset Sp(\xi_{j,k+1}|_y)$ (the inclusion is true by  {\blue{the}} {{induction}} assumption). By the second part of (iii) for $n=j-1$, we have $Sp(F_{j-1})\subset Sp(\xi_{j-1,j}|_{z_0})$, and consequently $Sp(F_{j-1})\subset Sp(\xi_{j-1,j}|_{z_0})\subset \bigcup_{z\in Sp(\xi_{j, k+1}|_y)}Sp(\xi_{j-1,j}|_{z})=Sp(\xi_{j-1,k+1}|_y).$ That is, {\blue{the second part of (vii)}} is true for $i=j-1$.}}

 {{{\blue{For the first statement  of the claim for $i=j-1,$ l}}et $y\in Sp(A_{X,k+1})$. By {\blue{the}} induction assumption,
 \beq\label{13-March26-2019-6}
 Sp(\xi_{j,k+1}|_y)\supset Sp(\phi_{j,k+1}|_y).
 \eneq
 Hence we have
 $$Sp(\phi_{j-1,k+1}|_y)={{\bigcup}}_{z\in Sp(\phi_{j,k+1}|_y)} Sp(\phi_{j-1,j}|_z)\qq\qq\qq\qq\qq$$
 $$\subset {{\bigcup}}_{z\in Sp(\phi_{j,k+1}|_y)} Sp(\xi_{j-1,j}|_z)~~~\qq\mbox{(by (\ref{13-March26-2019-5}) and by (iii) for } n=j-1)$$
 $$\hspace{-0.8in}\subset {{\bigcup}}_{z\in Sp(\xi_{j,k+1}|_y)} Sp(\xi_{j-1,j}|_z)\qq~~~\mbox{(by (\ref{13-March26-2019-6}))}\qq\qq\qq$$
 $$\hspace{-0.65in}=Sp(\xi_{j-1,k+1}|_y).\qq\qq\qq\qq\qq\qq\qq\qq\qq\qq$$ }}

  {\blue{(viii) {{Suppose that (iii)  holds, for $n=1,2,...,k.$
  Let $i<k$.}} {{Then}} ${{\big(}}\bigcup_{y\in Y_n}S{{p}}(\xi_{{i,n}}|_{y}){{\big) \cap X_i}}$ is ${{1/n}}$-dense in $X_i$ for all $0<i\le n$ (see {{the construction of  $Y_n$ from}} \ref{condition2}). {{Furthermore, $\big(\bigcup_{y\in Y_n\cup Sp(F_n)}SP(\xi_{{i,n}}|_{y})\big)\cap Sp(A_{X,i})$ is ${{1/n}}$-dense in $Sp(A_{X,i})$.}}

  Note that $\bigcup_{y\in Y_n}S{{p}}(\phi_{{i,n}}|_{y}) {~~{\cap X_i}}$ is ${{1/n}}$-dense in $X_i${{, and $\bigcup_{y\in Y_n\cup Sp(F_n)}S{{p}}(\phi_{{i,n}}|_{y})(\subset Sp(A_{X,i}))$ is ${{1/n}}$-dense in $Sp(A_{X,i})$.}}  So (viii) follows
  from (vii).}}

{\blue{ (ix)  Suppose that {{(iii),}} (iv), (v), and (viii) hold for $n=1,2,...,k.$ Then, for any $1\le i\le k,$
 $${{Sp(\xi_{i,k+1}|_{{\tht_{{k+1,2}}}})}}
\ \ \mbox{{{contains $Sp(F_i)$ and }} is  ${{2/k}}$\,-  dense in  } Sp(A_i).$$
The proof of (ix):}}
  {\blue{First note that, by (v) (holds {{for}} $1\le n\le k$),  $Y_k\cup {{T}}_k{{\cup Sp(F_k)}} \subset S{{p}}(\xi_{k,k+1}|_{\theta_{k+1,2}}).$
      It follows that $S{{p}}(\xi_{k,k+1}|_{\theta_{{{k}}+1,2}})$ is $1/k$-dense in ${{Sp(A_k)}}.$ That is, the statement holds for $i=k$ as ${{1/k}}$-dense implies ${{2/k}}$-dense.}} {{By (viii), $\big(\bigcup_{y\in Y_k\cup Sp(F_k)} Sp(\xi_{i,k}|_y)\big)\cap Sp(A_{X,i})$ is ${{1/k}}$-dense in $Sp(A_{X,i})$. Hence
 \beq\label{13-March27}
 Sp(\xi_{i,k+1}|_{\tht_{k+1,2}})\cap Sp(A_{X,i})~\big(\supset\bigcup_{y\in Y_k\cup Sp(F_k)} Sp(\xi_{i,k}|_y)\cap Sp(A_{X,i})\big)~\mbox{is}~ {{1/k}}\mbox{-dense in}~Sp(A_{X,i}).\qq
 \eneq }}
 {{From the fact that $T_k\cup Sp(F_k)$ is ${{1/k}}$-dense in $Sp(A_{C,k})$ (see the construction of $T_k$ in \ref{condition2}), by applying (iv) to $S=T_k\cup Sp(F_k) \subset Sp(A_{C,k})$ (and (iii), respectively), we know that $\big(\bigcup_{s\in T_k\cup Sp(F_k)} Sp(\xi_{k-1,k}|_s)\big)\cap Sp(A_{C,k-1})$
 is ${{\eta_{k-1}/k}}$-dense in $Sp(A_{C,k-1})$ (and contains $Sp(F_{k-1})$, respectively).
 Hence $Sp(\xi_{k-1,k+1}|_{\tht_{k+1,2}}))\cap Sp(A_{C,k-1})$ is ${{\eta_{k-1}/k}}$-dense in $Sp(A_{C,k-1})$.
 Now applying (iv) to $S'=\big(\cup_{t\in T_k\cup Sp(F_k)} Sp(\xi_{k-1,k}|_t)\big)\cap Sp(A_{C,k-1})$, we have $\big(\bigcup_{s\in S'}Sp(\xi_{k-2,k-1}|_s)\big)\cap Sp(A_{C,k-1})$ is ${{(\eta_{k-2}\eta_{k-1})/k}}$-dense in $Sp(A_{C,k-2})$.
 Since $S'\subset Sp(\xi_{k-1,k+1}|_{\tht_{k+1,2}}))\cap Sp(A_{C,k-1})$, $Sp(\xi_{k-2,k+1}|_{\tht_{k+1,2}}))\cap Sp(A_{C,k-2})$ is ${{(\eta_{k-2}\eta_{k-1})/k}}$-dense in $Sp(A_{C,k-2})$. Similarly, by induction (reversely), one gets that $Sp(\xi_{i,k+1}|_{\tht_{k+1,2}}))\cap Sp(A_{C,i})$ is $ {{(\eta_i\cdots \eta_{k-2}\eta_{k-1})/k}}$-dense in $Sp(A_{C,i})$.  Note that $\prod_{n=i}^{k-1} \eta_n\leq \prod_{n=1}^\infty \et_n< 2$. So $Sp(\xi_{i,k+1}|_{\tht_{k+1,2}}))\cap Sp(A_{C,i})$ is ${{2/k}}$-dense in $Sp(A_{C,i})$. Combining with (\ref{13-March27}), we get the desired conclusion.  }}

{\blue{(x): Suppose that $\phi_{n,n+1}$ are constructed which satisfy (i)--(vi)
for $n=1,2,...,k+2.$  For $i<k,$
then, {{for any $y\in Sp(A_{k+2})$, we have}}
\beq
Sp(\xi_{i,k+2}|_y) \,\,\, {\rm is}\,\, {{2}}/k{\rm -dense}\,\, {\rm in}\,\,\, Sp(A_i).
\eneq
}}
{\blue{To see this, we note, by  (iii), ${{\tht_{k+1,2}\in }}Sp(F_{k+1})\subset Sp(\xi_{k+1, k+2}|_y)$ for {{any}} $y\in Sp(A_{k+2}).$
Therefore
$Sp(\xi_{i,k+2}|_y)\supset Sp(\xi_{i,k+1}|_{\theta_{k+1,2}}).$
Combining this  with (ix), we obtain (x).
}}

  Property  (v{{i}}')  will be used
(together with
\ref{range 0.25})  to prove
that the  limit algebra $A$ {has} the
property that
 $A\otimes U\in {\cal B}_0$
 for any UHF-algebra $U$.






\end{NN}

\begin{NN}\label{range 0.30}
Suppose that we
have constructed {{the finite sequence}}
\begin{displaymath}
\xymatrix{ A_1  \ar[r]^{\xi_{{1,2}}} &A_2  \ar[r]^{\xi_{{2,3}}} &
\cd \ar[r]^{\xi_{{n-1,n}}}& A_n   ,}
\end{displaymath}
{{in such a way}} that for every $i\leq n-1$, $\xi_{i,i+1}$ satisfies conditions (i)
\,--\,(vi) (with $i$ in place of $n$) in \ref{range 0.29}. We will
construct the map $\xi_{{n,n+1}}:~ A_n\to A_{n+1}$.
{{Write }}
\beq\nonumber
&&\vspace{-0.1in}A_n=\big\{(f,g)\in C([0,1],E_n)\oplus {{A_{X,n}}};~ f(0)=\bt_{{n,0}}(\pi_{X,n}^e(g)), f(1)=\bt_{{n,1}}(\pi_{X,n}^e(g))\big\}\,\,\,{\rm with}\\
&&\vspace{-0.1in}{{A_{X,n}}}=P_n M_{\infty}(C(X_n))P_n\oplus
\bigoplus_{i=2}^{p_n} M_{[n,i]},
\eneq
where ${\rm rank}(P_n)=[n,1]$ ({\blue{recall that $\pi_{X,n}^e: A_{X,n}\to F_n$ is from \ref{construction})}}.

Also recall {{that we}}  denote $t\in(0,1)_{{n,j}}\sbs
Sp\big(C([0,1],E_n^j)\big) $ by $t_{{n,j}}$ to distinguish
{{the spectra}}
from different direct {{summands}} of $C([0,1], E_n)$, and we denote  $0\in [0,1]_{{n,j}}$ by
$0_{{n,j}},$ and $1\in [0,1]_{{n,j}}$ by $1_{{n,j}}.$ Note that $0_{{n,j}}$ {{and}} $1_{{n,j}}$ do not
correspond to {{single}}  irreducible representations.  In fact, $0_{{n,j}}$
corresponds to the direct sum of irreducible
representations for the
set
$$\left\{ \tht_{{n,1}}^{\sim b^{{n}}_{0,j1}}, \tht_{{n,2}}^{\sim b^{{n}}_{0,j2}},...,~ \tht_{{n,p_n}}^{\sim b^{{n}}_{0,jp_n}}\right\},$$
and $1_{{n,j}}$ corresponds to the set
$$\left\{ \tht_{{n,1}}^{\sim b^{{n}}_{1,j1}}, \tht_{{n,2}}^{\sim b^{{n}}_{1,j2}},..., ~\tht_{{n,p_n}}^{\sim b^{{n}}_{1,jp_n}}\right\}.$$
Again recall from \ref{condition2},  $T_n\sbs Sp(A_n)$ is defined by
$$T=\left\{ (\frac kn)_{{{n,j}}};~ j=1,2,... ,l_n;~ k=1,2, ... , n-1 \right\}.$$
%

 Recall that   in {{the condition}} ${{\spdd_1}},$
 $L_n=l_n\cdot (n-1)+L_{{n,Y}}=\#(T_{{n}}\cup Y_{{n}})$ and  $M=\max\{ {{b^{{n}}_{0,ij}:}}~ i=1,2,...,p_n;~ j=1,2,..., l_n\}$, where we write $Y_n=\{y_1,y_2,..., y_{_{L_{n{{,Y}}}}}\}\subset X_n$.


\end{NN}

\begin{NN}\label{range 0.31}~
{\blue{First we define a unital \hm, $\xi_X: A_n\to A_{X, n+1}.$ Denote, only in this subsection, by
$\Pi': A_{X, n+1}\to \bigoplus_{i\not=2} A_{X,n+1}^i{{(=A_{X,n+1}^1\oplus \bigoplus_{i=3}^{p_{n+1}}F_{n+1}^i)}}$ and $\Pi^{(2)}: A_{X, n+1}\to A_{X,n+1}^2 {{(=F_{n+1}^2)}}$ the quotient maps.
Define $\xi_X': A_n\to \bigoplus_{i\not=2} A_{X,n+1}^i$ by $\xi_X'=\Pi'\circ \circ \pi_{I,n+1}\circ  \phi_{n,n+1}.$}}
{\blue{We note that, for $a\in A_n$ such that $\pi_{X,n}^1(a)\not=0,$ by the definition of $\psi_{n,n+1},$
$\xi_X'(a)\not=0.$}}

{\blue{Note that, by (a) of  \ref{condition2}, $\phi_{n, n+1}(I_n+J_n)\subset I_{n+1}+J_{n+1}.$   Therefore,
\beq\label{1345sp}
SP(\phi_{{n,n+1}}|_{{\tht_{{n+1,2}}}})=\left\{ \tht_{{n,1}}^{\sim c_{21}}, \tht_{{n,2}}^{\sim c_{22}},..., \tht_{{n,p_n}}^{\sim c_{2p_n}}\right\}.
\eneq
}}\\
{{(Here $c_{jk}$ means $c^{n,n+1}_{jk}$.)}}
{\blue{Let $Y_n=\{y_1,y_2,...,y_{_{L_{n,Y}}}\}$ be as in \ref{condition2}.  Since $X_n$ is path connected, for each $i,$
there is a continuous {{simple}} path $\{y_i(s): s\in [0, 1]\}\subset X_n$ such that $y_i(0)=y_i$  and $y_i(1)=x_n^0.$ {{Note that {\blue{$P_n(y_i)M_{\infty}P_n(y_i)$}} can
 be identified with $M_{[n,1]}=F_n^1$. We will also  identify
 $P_n(y_i(s))M_{\infty}P_n(y_i(s))$ with $M_{[n,1]}=P_n(x_n^0)M_{\infty}P_n(x_n^0)=A_n|_{\tht_{n,1}}=
 F_n^1$---such
 an identification could be chosen to be continuously depending on $s$. (Here,
 we only use the fact that any projection (or vector bundle) over
the  interval is trivial to make such {{an}} identification. Since  the
 projection $P_n$ itself may
 not be trivial, it is possible
 that the paths for different $y_i$ and $y_j$ ($i\not=j$) may intersect at
 $y_i(s_1)=y_j(s_2)$
 and we may use {{a}}  different identification of
 $P_n(y_i(s_1))M_{\infty}P_n(y_i(s_2))= P_n(y_j(s_2))M_{\infty}P_n(y_j(s_2))$
 with $M_{[n,1]}$ for $i$ and $j$.)
 When we talk about $f(y_i(s))$ later, we will {{consider}}  it to
 be an element {{of}}  $M_{[n,1]}$ (rather than {{of}}
 $P_n(y_i(s))M_{\infty}P_n(y_i(s))$).
 }}
Define  (recall $\pi_{I,n}: A_n\to A_{X,n}$ is the quotient map)
$$
\Omega_{Y,s}(f)=\diag(\pi_{I, n}(f)(y_1(s)),...,\pi_{I,n}(f)(y_{_{L_{n,Y}}}(s))~{{\in M_{L_{n,Y}[n,1]}}}~\rforal f\in A_n\andeqn s\in [0,1].
$$
Note that {{$\Omega_{Y,s}(\one_{A_n})=\one_{_{M_{L_{n,Y}[n,1]}}}$ is independent of $s$ and that}}
$\Omega_{Y,1}(f)=(\theta_{n,1}^{\sim L_{n,Y}})(f)$ for all $f\in A_n.$}}

{\blue{For each $(\frac{{h}}n)_{n,j}\in T_n$, {{where $1\leq h\leq n-1$ }} (see \ref{condition2}), there exists also a continuous path $\{g_{n,j{{,h}}}(s): s\in [0,1]\}\subset
[0,1]_{n,j}$ such that $g_{n,j{{,h}}}(0)=(\frac hn)_{n,j}$ and $g_{n,j{{,h}}}(1)=0_{n,j}.$
 Define, for each $f\in A_n$ and $s\in [0,1],$
 $$
  \Omega_{I,s}(f)=\bigoplus_{j=1}^{l_n}\diag(\pi_{I,n}(f)(g_{n,j{{,1}}}(s)), \pi_{I,n}(f)(g_{n,j{{,2}}}(s)),..., \pi_{I,n}(f)(g_{n,j{{,n-1}}}(s))){{\in M_{\sum_{j=1}^{l_n}(n-1)\{n,j\}}}}.
  $$
  {{(Recall that $\{n,j\}=\sum_{k=1}^{p_n}b^n_{0,jk}=\sum_{k=1}^{p_n}b^n_{1,jk}$ is the rank of the representation of $A_n$ {{corresponding}}  to any point $t\in [0,1]_{n,j}$ as $E_n^j=M_{\{n,j\}}$.)}}
  {{Note that  $\Omega_{I,s}{{(\one_{A_n})}}\in M_{\sum_{j=1}^{l_n}(n-1)\{n,j\}}$ is independent of $s$.  }}
 Then
  $$
 \Omega_{I,0}(f)=\bigoplus_{j=1}^{l_n}\diag(\pi_{I,n}(f)((\frac 1n)_{n,j}, \pi_{I,n}(f)(\frac 2n)_{n,j},..., \pi_{I,n}(f)((\frac {(n-1)}n)_{n,j})),
 $$
 and,  as $0_{{n,j}}=\left\{ \tht_{{n,1}}^{\sim b^{{n}}_{0,j1}}, \tht_{{n,2}}^{\sim b^{{n}}_{0,j2}},..., \tht_{{n,p_n}}^{\sim b^{{n}}_{0,jp_n}}\right\},$
$$SP(\Omega_{I,1})=\left\{ \tht_{{n,1}}^{\sim b_1}, \tht_{{n,2}}^{\sim b_2},..., \tht_{{n,p_n}}^{\sim b_{p_n}}\right\},$$
where $b_k=(n-1)(\sum_{j=1}^{l_n} b_{n,jk}^n),$ $k=1,2,...,p_n.$}}
{\blue{Put
\beq\nonumber
   a_1= c_{21}-L_{n,Y}-\left(\sum_{j=1}^{l_n}b^{{n}}_{0,j1} \right)(n-1), \,\,
   a_k = c_{2k}-\left(\sum_{j=1}^{l_n}b^{{n}}_{0,jk} \right)(n-1),\,\, k=2,3,...,p_n.
 \eneq
}}
{\blue{Let  $\xi_{X,s}'$  be {{the}}  finite dimensional representation of $A_n$
defined by, for each $f\in A$ and $s\in [0,1],$
\beq\label{Dxiexs}
\xi_{X,s}'(f)=\diag(\tht_{{n,1}}^{\sim a_1}(f), \tht_{{n,2}}^{\sim a_2}(f),..., \tht_{{n,p_n}}^{\sim a_{p_n}}(f))
\oplus\Omega_{X,s}(f)\oplus \Omega_{I,s}(f).
\eneq
In particular,
\beq
SP(\xi_{X,0}')&=& \left\{ \tht_{{n,1}}^{\sim a_1}, \tht_{{n,2}}^{\sim a_2},..., \tht_{{n,p_n}}^{\sim a_{p_n}}\right\} \cup T_n \cup Y_n,\andeqn\\
SP(\xi_{X,1}') &=&\left\{ \tht_{{n,1}}^{\sim c_{21}}, \tht_{{n,2}}^{\sim c_{22}},..., \tht_{{n,p_n}}^{\sim c_{2p_n}}\right\}=SP(\phi_{n,n+1}|_{\theta_{n+1,2}}).
\eneq
Since $\xi'_{X,0}$ is homotopic to $\xi'_{X,1},$  and $\xi'_{X,s}(1_{A_n})=\xi'_{X,1}(1_{A_n})$ for all $s\in [0,1],$
 up to unitary equivalence,
we may view  $\{\xi_{X,s}': s\in [0,1]\}$ as a  continuous path of unital \hm s from $A_n$ into $A_{X,n+1}^2=F_{n+1}^2. $}}
View
\beq
&&e_u:=\diag(\tht_{{n,1}}^{\sim a_1}(1_{A_n}), \tht_{{n,2}}^{\sim a_2}(1_{A_n}),..., \tht_{{n,p_n}}^{\sim a_{p_n}}(1_{A_n}))\andeqn\\
&&e_c:=\Omega_{{{Y}},s}(1_{A_n})\oplus \Omega_{I,s}(1_{A_n})
\eneq
as two projections in $F_{n+1}^2${{, which do not depend on $s$.}}
From ${{\spdd_1}}$ of {\blue{\ref{condition2},}}
 we know that
\beq\label{150102-e1}
a_i \geq \frac{2^{2n}-1}{2^{2n}}c_{2i}.
\eneq
{\blue{Then,  for the tracial state $\tau$ of $F_{n+1}^2,$
\beq\label{1344trace1}
\tau(e_c)<1/2^{2n}\andeqn \tau(e_u)>1-(1/2^{2n}).
\eneq}}
{\blue{Define $\xi_{X,s}: A_n\to A_{X, n+1}$ by
$\xi_{X,s}=(\xi_X'\oplus \xi_{X,s}').$  By replacing $\xi_{X,s}$ by ${\rm Ad}\, U\circ \xi_{X,s}$
for a suitable unitary {{path}} $U_{{s}}{{=\oplus U_s^j}}\in {{\oplus A_{X,n+1}^j}}$ {{(with $U_s^j=1$ if $j\not=2$)}},  we may assume that
\beq\label{13-March28}
\xi_{X,1}=\pi_{I, n+1}\circ \phi_{n,n+1}.\eneq
{{Since $U_s^j=1$ for $j \not=2$, from the definition of $\xi_{X}'$, we get}} $\xi_{X,s}|_x=\phi_{n, n+1}|_x$ for all $x\in S{{p}}(A_{X,n+1})$ {{with}} $x\not=\theta_{n+1,2}.$
It follows that
\beq\label{1344trace0}
{{\xi_{X,1}(f)-\xi_{X,s}(f)}}={{e_c(\xi_{X,1}(f)-\xi_{X,s}(f))}}\rforal {{s\in [0,1]~\mbox{and}}}~ f\in A_n.
\eneq
Define $\xi_X:=\xi_{X,0}.$
Therefore $\xi_X$ is homotopic to $\pi_{I,n+1}\circ \phi_{n,n+1}$
and $\xi_X|_x=\phi_{n, n+1}|_x$ for all $x\in S{{p}}(A_{X,n+1})$ {{with}} $x\not=\theta_{n+1,2}.$}}
{\blue{Since $\xi_X$ is homotopic to $\pi_{I,n}\circ \phi_{n,n+1},$
$(\xi_X)_{*i}=(\pi_{I,n}\circ \phi_{n,n+1})_{*i},$ $i=0,1.$
From \eqref{1344trace0} and \eqref{1344trace1}, we also have
\beq\label{1344trace2}
|\tau(\pi_{I,n+1}\circ \phi_{n,n+1}(f))-\tau(\xi_X(f))|<{{\tau(e_c)\|f\|\leq}}(1/2^{2n-1})\|f\|\tforal \tau\in T(A_{X,n+1}).\qq
\eneq
From the first paragraph of this section and by \eqref{Dxiexs},
if $a\in A_n\setminus I_n,$ then $\xi_X'(a)\not=0.$ }}

 \end{NN}

\begin{NN}\label{range 0.32}
~~   {\blue{In \ref{range 0.31}, we have defined a unital \hm\,
$\xi_X: A_n\to A_{X,n+1}.$
In this subsection, we define the map $\xi_{n,n+1}.$ We first define a unital \hm\, $\xi_E: A_n\to C([0,1], E_{n+1}).$
}}
Define
\beq\label{xiebn1}
&&{\xi_E}|_{{0_{{n+1,j}}}}=\pi^{{j}}\circ \bt_{{n+1,0}}\circ \pi_{X, n+1}^e\circ\xi_X:~ A_n\to E_{n+1}^j\andeqn\\\label{xiebn2}
&&{\xi|_E} |_{{1_{{n+1,j}}}} = \pi^{{j}}\circ \bt_{{n+1,1}}\circ \pi_{X, n+1}^e\circ \xi_X:~ A_n\to
E_{n+1}^j,
\eneq
  where $\pi^{{j}}:~ E_{n+1}\to E_{n+1}^j$.
Now we need to connect ${\xi_E}|_{{0_{{n+1,j}}}}$ and ${\xi_E}|_{{1_{{n+1,j}}}}$
to obtain  $\xi_E.$

{\blue{Fix a finite subset ${{\cal G}}_n\subset A_n$ {{in (ii) of \ref{range 0.29},}} {{There}} is $\dt>0$ such
that
\beq
\|{\bar \phi}_{n,n+1}\circ \pi_{J,n}(f)(s)-{\bar \phi}_{n,n+1}\circ \pi_{J,n}(f)(s')\|<1/2^{2n+1}\rforal f\in {\cal G}_n,
\eneq
{{provided}}  $|s-s'|<3\dt$ {{and}} $s,s'\in [0,1]${{, and such that ${{1/(1-2\dt)}}<\eta_n$.}}}}
{\blue{Let $h_0(t)=t/\dt$ for $t\in [0,\dt],$ $h_z(t)=(t-\dt)/(1-2\dt)$ for $t\in [\dt, 1-\dt]$ and $h_1(t)=(1-t)/\dt$ for
$t\in [1-\dt,1].$
Define, for each $f\in A_n,$
\beq
{\xi_E}(f)|_{[0,\dt]_{n+1,j}}(t)=\pi^j\circ \bt_{{n+1,0}}\circ \pi_{X,n+1}^e\circ\xi_{X, h_0(t)}(f)\tforal t\in [0,\dt],\\\label{dexizz1}
{\xi_E}(f)|_{(\dt, 1-\dt]_{n+1,j}}(t)=(\pi^j\circ {\bar \phi}_{n,n+1}\circ \pi_{J,n}(f))(h_z(t))\tforal t\in (\dt, 1-\dt]\andeqn\\
{\xi_E}(f)|_{(1-\dt, 1]_{n+1,j}}(t)=\pi^j\circ \bt_{{n+1,1}}\circ \pi_{X, n+1}^e\circ\xi_{X, h_1(t)}(f)\tforal t\in (1-\dt, 1],
\eneq
where ${\bar \phi}_{n,n+1}: A_{C,n}\to A_{C,n+1}\subset C([0,1], E_{n+1})$ is the injective \hm\, given by
\ref{13-Jan31-2019}.}} {\blue{In particular,  using the fact that $\phi_{n,n+1}$ is a map from
$A_n$ to $A_{n+1},$ {{and (\ref{13-March28}) (note that $h_0(\dt)=1=h_1(1-\dt)$),}} for all $f\in A_n,$ {{we have}}
\beq
&&\hspace{-0.5in}\pi^j\circ {\bar \phi}_{n,n+1}\circ \pi_{J,n}(f)(h_z(\dt))=
(\pi^j\circ {\bar \phi}_{n,n+1}\circ \pi_{J,n}(f))(0)\\
&&\hspace{-0.4in}=\pi^j\circ \bt_{{n+1,0}}\circ \pi_{X, n+1}^e\circ \phi_{n,n+1}(f)=
\pi^j\circ \bt_{{n+1,0}}\circ \pi_{X, n+1}^e\circ\xi_{X, h_0(\dt)}(f), \,\,{\rm{{and}}}\\\label{dexizz}
&&\hspace{-0.5in}(\pi^j\circ {\bar \phi}_{n,n+1}\circ \pi_{J,n}(f))(h_z(1-\dt))=(\pi^j\circ {\bar \phi}_{n,n+1}\circ \pi_{J,n}(f))(1)\\
&&\hspace{-0.2in}=\pi^j\circ \bt_{{n+1,1}}\circ \pi_{X, n+1}^e\circ \phi_{n,n+1}(f)
=\pi^j\circ \bt_{{n+1,1}}\circ \pi_{X, n+1}^e\circ\xi_{X, h_1(1-\dt)}(f).
\eneq}}
\noindent
Thus, {\blue{$\xi_E$ defines a unital   \hm\, from $A_n$ into $C([0,1], E_{n+1})$
which is injective on $I_n.$
Finally, define $\xi_{n, n+1}: A_n\to A_{n+1}=C([0,1], E_{n+1})\oplus_{\bt_{n{{+1}},0}\circ \pi_{X,n+1}^e, \bt_{n+1,1}\circ \pi_{X,n+1}^e} A_{X,n+1}$ by
$$
\xi_{n, n+1}(f)=(\xi_E(f), \xi_X(f))\rforal f\in A_n.
$$
By \eqref{xiebn1} and \eqref{xiebn2},  this is {{indeed}} a unital \hm\,
from $A_n$ into $A_{n+1}$  (see \eqref{Adef} and \eqref{Adef2}).
Since $\xi_E$ is injective, by the end of \ref{range 0.31}, $\xi_{n,n+1}$ is injective.
Note that, by  (the end of) \ref{range 0.31},  $(\pi_{I,n+1}\circ \xi_{n, n+1})_{*i}=(\xi_X)_{*i}=(\pi_{I,n+1}^e\circ \phi_{n,n+1})_{*i},$ $i=0,1.$
It follows from {{part (1) of \ref{13-Feb17-2019} (also see the last sentence of the proof of part (1))}}  that $(\pi_{i,n+1})_{*i}: K_i(A_{n+1})\to K_i(A_{X,n+1})$  is injective.
Therefore $(\xi_{n,n+1})_{*i}=(\phi_{n,n+1})_{*}.$   So (i) {{of \ref{range 0.29} }}holds. {{We {{have}} already {{mentioned}} that $\xi_X|_x=\phi_{n, n+1}|_x$ for all $x\in Sp(A_{X,n+1})$ with $x\not=\theta_{n+1,2}$ in \ref{range 0.31}, {{and}} so (vi) of \ref{range 0.29} holds.}} By the presence of the representations corresponding to $Y_n$, {{$T_n$,}} {{and $Sp(F_n)$ }}in
\eqref{Dxiexs}, we also see that (v) holds. }}
{\blue{{{For}} $y\in S{{p}}(A_{X,n+1})$  and $y\not=\theta_{n+1,2}$, by (v),
$S{{p}}(\phi_{n,n+1}|_y)=S{{p}}(\xi_{n,n+1}|_y).$ For $y=\theta_{n+1,2},$ by \eqref{1345sp}, \eqref{Dxiexs}, and
\eqref{150102-e1},
if $x\in Sp(A_n)$ and $x\in S{{p}}(\phi_{n,n+1}|_y),$ then $x\in S{{p}}(\xi_{n,n+1}|_y).$}}{{That is, the first part of (iii) of \ref{range 0.29} holds.}}
{\blue{Moreover,
$Sp(F_n)\subset S{{p}}(\xi_{n,n+1}|_y)$ for all $y\in SP(A_{X,n+1}).$
Therefore, $Sp(F_n)\subset S{{p}}(\xi_{n,n+1}|_y)$ for all $y\in [0, \dt)_{n+1,j}\cup (1-\dt, 1)_{n+1,j}.$ By  (c) of \ref{condition2}  and (5) of \ref{range 0.6},
$Sp(F_{n})\subset SP({\bar \phi}_{n,n+1}|_s)$ for each $s\in (0,1)_{n+1, j}.$  It follows from
\eqref{dexizz} that $Sp(F_n)\subset SP(\xi_{n,n+1}|_s)$ for all $s\in [\dt, 1-\dt]_{n+1,j}.$
{{Therefore,}} $Sp(F_n)\subset SP(\xi_{n,n+1}|_y)$ for all $y\in Sp(A_{n+1}).$ This {{proves the second part of}} (iii).}}\\

{{{\bf Claim:} {{If}} a finite subset $Z\subset (0,1)$ {{is}} such that $Z\cup \{0,1\}$ is $d$-dense in $[0,1]$, then the finite subset $h_z(Z\cap(\dt,1-\dt))\cup \{0,1\}$ is $\eta_n d$-dense in $[0,1]$}}

\begin{proof} {{Let us order the set $Z=\{z_j\}_{j=1}^k$ as
$0<z_1<z_2<\cdots<z_k<1.$
Then $Z\cup \{0,1\}$ is $d$-dense in $[0,1]$ if and only if
$$z_1<2d,~~1-z_k<2d,~~\mbox{and}~~ z_{j+1}-z_j<2d~ \mbox{for all}~j=1,2,\cdots,k-1.$$
For {\blue{convenience,}} let $z_0=0$ and $z_{k+1}=1$. Let $j_1$ be the smallest index such that $z_{j_1}>\dt$ and $j_2$ be the largest index such that $z_{j_2}<1-\dt$. Then $Z\cap (\dt,1-\dt)=\{z_j\}_{j=j_1}^{j_2}$ and the set $h_z(Z \cap (\dt,1-dt))\cup\{0,1\}$ can be {{listed}} as
$$h_z(z_{j_1-1})=0<h_z(z_{j_1})< h_z(z_{j_1+1})<\cdots <h_z(z_{j_2})<1=h_z(z_{j_2+1}).$$
The claim follows from the fact that}} $$h_z(z_{j+1})-h_z(z_j)\leq {{(z_{j+1}-z_j)/(1-2\dt)}}<{{2d/(1-2\dt)}}<2\eta_nd.$$
\end{proof}
{\blue{{{Let us verify that (iv) of \ref{range 0.29} holds.}}
For any $d>0,$ let $S\subset {{Sp(A_{C,n+1})}}$ be $d$-dense {{in $Sp(A_{C,n+1}),$ and satisfy}} that {{$S\supset Sp(F_{n+1})$ (note that $Sp(F_{n+1})$ is a subset of $Sp(A_{C,n+1})$).}}
Let ${{Z}}=S\cap {{(0,1)}}_{n+1,{{1}}}$ and {{$Z_0=S\cap (\dt,1-\dt)_{n+1,1}$.}} {{It follows from the $d$-density of $S$ in $Sp(A_{C,n+1})$, regarding $Z_0\subset Z$ as subsets of the open interval $(0,1)$, {{that}} $Z\cup \{0,1\}$ is $d$-dense in $[0,1]$. Hence by the claim, $h_z(Z_0)\cup\{0,1\}$ is {\blue{$\eta_n d$-}}dense in $[0,1]$.}}
Then, by (c) of \ref{condition2} and  {{(6) of \ref{range 0.6} (applied to all indices $i_0=1,2,\cdots l_n$ and $j_0=1$), we know that $\big({{\bigcup}}_{z\in h(Z_0)}Sp(\phi_{n,n+1}|_z)\cap Sp(A_{C,n})\big)\cup Sp(F_n)$}} is $\eta_nd$-dense in $Sp(A_{C,n})$. {{By {{(\ref{dexizz1})}}, ${{\bigcup}}_{z\in Z_0}Sp(\xi_{n,n+1}|_z)={{\bigcup}}_{z\in h(Z_0)}Sp(\phi_{n,n+1}|_z)$.}}
{{ By (iii) of \ref{range 0.29}, we know that
${{\bigcup}}_{s\in S}Sp(\xi_{n,n+1}|_{s})\supset Sp(F_n)\cup \big({{\bigcup}}_{z\in Z_0}Sp(\xi_{n,n+1}|_z)\big)$.}}
{{ It follows that ${{\bigcup}}_{s\in S}Sp(\xi_{n,n+1}|_{s})\cap Sp(A_{C,n})$ is $\eta_n d$-dense in $Sp(A_{C,n})$. Hence (iv) of \ref{range 0.29} hold.}}}}

{\blue{It remains to check (ii) {{of \ref{range 0.29}}} holds.
Note that,
by \ref{range 0.31},
$Sp(\xi_{n,n+1}|_y)=Sp(\phi_{n,n+1}|y)$ for all $y\in Sp(A_{X,n+1}^j)$ for $j\not=2.$
Note that, for the tracial state $t$ of $E_{n+1}^j,$  {{the map $$a\mapsto {t(\pi^j\circ \bt_{n+1, 0}(a))\over{t(\pi^j\circ \bt_{n+1,0}(1_{F_{n+1}^2}))}}~~\mbox{for}~~a\in F_{n+1}^2 $$}}
 is a tracial state of $F_{n+1}^2.$
It follows from  \eqref{1344trace1} that
\beq\label{1345trace1}
{t(\pi^j\circ \bt_{n+1, 0}(e_c))\over{t(\pi^j\circ \bt_{n+1,0}(1_{F_{n+1}^2}))}}<1/2^{2n}.
\eneq
}}
{\blue{Define $e_{c,j}'=\pi^j\circ \bt_{n+1, 0}(e_c)).$ Then $t(e_c')<1/2^{2n}$ for $t\in E_{n+1}^j.$
It follows that, for each $t\in [0,\dt]_{n+1,j},$ and $\tau\in T(E_{n+1}^j)$ (see also {{\eqref{1344trace0}),}}
\beq\nonumber
&&|\tau(\xi_{n,n+1}(f)(s))-\tau(\phi_{n,n+1}(f)({{0_{n+1,j}}}))|\\
&&=|\tau(e_{c,j}'((\xi_{n,n+1}(f)(s))-\phi_{n,n+1}(f)({{0_{n+1,j}}})))|<(1/2^{2n-1})\|f\|\rforal f\in A_n.\qq\label{1345trace2}
\eneq
Exactly the same computation shows that,  for all $s\in [1-\dt, 1]_{n+1,j}$ and  for all $\tau\in T(E_{n+1}^j),$
\beq\label{1345trace3}
|\tau(\xi_{n,n+1}(f)(s))-\tau(\phi_{n,n+1}(f)({{1_{n+1,j}}}))|<(1/2^{2n-1})\|f\|\rforal f\in A_n.
\eneq
}}
{\blue{Note {{that}}  $|s-h_z(s)| \leq \dt<3\dt$ for all $s\in [0,1].$  By the choice of $\dt,$
\beq\nonumber
&&\hspace{-0.6in}\|\xi_{n,n+1}(f)(s)-\phi_{n,n+1}(f)(s)\|=\|({\bar \phi}_{n,n+1}\circ \pi_{J,n}(f))(h_z(s))-({\bar \phi}_{n,n+1}\circ \pi_{J,n}(f))(s)\|\\\label{1345trace5}
&&\hspace{0.8in}<1/2^{2n+1}\rforal f\in {\cal G}_n
\eneq
for all $s\in [\dt, 1-\dt]_{n+1,j}.$ }}{{Again, by the choice of $\dt$, we also have,  for all {{$f\in {\cal G}_n$,}}
\beq\label{13-March29}
\|\phi_{n,n+1}(f)(s)-\phi_{n,n+1}(f)(0_{n+1,j})\|<1/2^{2n+1}~ \mbox{for all}~ s\in[0,\dt]_{n+1,j} ~~\mbox{and}
\eneq
\beq\label{13-March29-1}
\|\phi_{n,n+1}(f)(s)-\phi_{n,n+1}(f)(1_{n+1,j})\|<1/2^{2n+1}~ \mbox{for all}~ s\in[1-\dt, 1]_{n+1,j}{{.}}
\eneq}}

{\blue{Combining  \eqref{1345trace2}, {{\eqref{13-March29},}} \eqref{1345trace3}, {{\eqref{13-March29-1}}}{{, \eqref{1345trace5}, and}} \eqref{1344trace2},
\beq\label{1345trace6}
\|\xi_{n,n+1}^{\sharp}({{\hat f}})-\phi_{n,n+1}^{\sharp}({{\hat f}})\|<1/2^{2n-1}{{
+1/2^{2n+1}<1/2^{2n-2}}}~\rforal f\in {\cal G}_n.
\eneq
}}
{\blue{This proves (ii). By induction, this completes the construction of $\xi_n.$}}

\end{NN}


\begin{NN}\label{range 0.33}
   Let $ B=\lim(A_n, \xi_{{n,n+1}})$. Recall that $A=\lim(A_n, \phi_{{n,n+1}})$.
   {\blue{By (i) {{of \ref{range 0.29}}}, ${\xi_{n,n+1}}_{*i}={\phi_{n,n+1}}_{*i},$ $i=0,1.$ It follows that
   $$
   (K_0(B), K_0(B)_+, [1_B], K_1(B))=(K_0(A), K_0(A)_+, [1_A], K_1(A)).
   $$}}
    {\blue{For each $n$ and $\sigma>0,$ choose
    $m>n+1$ with $2/m-2<\sigma.$  Then, by (x) of \ref{range 0.29},
$Sp(\xi_{{n, m}}|_{x})$ is ${{2}}/(m-2)$-dense in $Sp(A_n)$ for
any $x\in Sp(A_m).$ It follows from  Proposition \ref{simplelimit} (see also  the end of \ref{density})
that $B$ is a simple \CA.}}

 {\blue{ We will show, in fact,
   \beq
  \hspace{-0.3in} (K_0(B), K_0(B)_+, [1_B], K_1(B), T(B), r_B)=(K_0(A), K_0(A)_+, [1_A], K_1(A), T(A), r_A).
   \eneq

Consider the following non-commutative diagram:
\begin{displaymath}
    \xymatrix{
        \Aff(T(A_1)) \ar@/_/[d]_{\id_1^\sharp}\ar[r]^{\phi_{1,2}^\sharp}  & \Aff(T(A_2)) \ar[r]^{\phi_{2,3}^\sharp} \ar@/_/[d]_{\id_2^\sharp}& \Aff(T(A_3)) \ar[r] \ar@/_/[d]_{\id_3^\sharp}& \cd \Aff(T(A ))\\
        \Aff(T(A_1)) \ar[r]^{{\xi}_{1,2}^{\sharp}}\ar@/_/[u]_{\iota_1^{\sharp}} &
         \Aff(T(A_2)) \ar[r]^{{\xi}_{2,3}^{\sharp}}\ar@/_/[u]_{\iota_2^{\sharp}}&
         \Aff(T(A_3)) \ar[r]\ar@/_/[u]_{\iota_3^{\sharp}}& \cd
         \Aff(T(B)), }
\end{displaymath}
where $\id_k: A_k\to A_k$ and $\iota_k :A_k\to A_k$ are {{both the}} identity maps{{---but}} we write {{them}} differently as they come from
two different systems.}}
{\blue{Recall that, for each $n,$  $\{{\cal F}_{n,k}: k\ge 1\}$ is an increasing sequence of finite subsets
of $A_n$ whose union is dense in {{the unit ball of}} $A_n.$  Recall also that ${\cal G}_1={\cal F}_{1,1},$
${\cal G}_n={\cal F}_{n,n}{{\cup}} {{\bigcup}}_{i=1}^{{n-1}}\phi_{i,n}({{{\cal F}_{i,n}\cup}}{\cal G}_i)\}
{{\cup\bigcup}}_{i=1}^{{n-1}}\xi_{i,n}({{{\cal F}_{i,n}\cup}}{\cal G}_i).$
For each $i$ and each $x\in A_i,$  and for each $n>i,$ there exist $j>i$ and  $y\in {\cal F}_{i,j}$ such that
$\|x-y\|<1/2^{n+1}.$ Then $\|\phi_{i,j}(x)-\phi_{i,j}(y)\|<1/2^{n+1}.$
Note $\phi_{i,j}(y) \in  {{\cal G}}_j,$ Let $n_0=\max\{j, n\}.$
Denote by  $\hat{x}\in \Aff(T(A_i))$ the function $\hat{x}(\tau)=\tau(x)$ for all $\tau\in T(A_i).$
  Put $z=\phi_{i,n_0}^{\sharp}(\hat{y}).$
By (ii)  of  \ref{range 0.29},  for any $m> n_0,$
\beq\label{intertw-1}
\| \xi_{n_0,m}^{\sharp}\circ \id_{n_0}^{\sharp}(z) -\id_m^{\sharp}\circ \phi_{n_0,m}^{\sharp}(z)\|<1/2^{{2n_0-2}}.
\eneq
It follows that
\beq
\|{{\xi}}_{n_0,m}^{\sharp}\circ \id_{n_0}^{\sharp}(\phi_{{{i}},n_0}^{\sharp}(\hat{x}))-\id_m^{\sharp}\circ \phi_{n_0,m}^{\sharp}(\phi_{i,n_0}(\hat{x}))\|<1/2^{n-1}.
\eneq
By exactly the same reason, for any $m>n_0,$
\beq\label{intertw-2}
\|\phi^{\sharp}_{{n_0,m}}\circ \iota_{n_0}^{\sharp}(\xi_{i,n_0}^{\sharp}(\hat{x}))-
  \iota_m^{\sharp}\circ \xi_{n_0,m}^{\sharp}(\xi_{i, n_0}^{\sharp}(\hat{x}))\|<1/2^{n-1}.
 \eneq
 Note also that $\id_k^{\sharp}$ and $\iota_k^{\sharp}$ are isometric isomorphisms.
It follows that the non-commutative diagram above is approximately intertwining, and {{the sequences of maps}}
$\{\id_k^{\sharp}\}$ and $\{\iota_k^{\sharp}\}$ induce
two isometric isomorphisms $j: \Aff(T(A))\to \Aff T(B))$ and $\iota^{\sharp}: \Aff(T(B))\to \Aff(T(A))$
between {{the}} Banach spaces $\Aff(T(A))$ and $\Aff(T(B))$ such that $j\circ \iota^{\sharp}=\id_{Aff(T(B))}$ and
$\iota^{\sharp}\circ j=\id_{\Aff(T(A))}.$
Moreover, since each $\id_k^{\sharp}$ and $\iota_k^{\sharp}$ are order preserving,
$j$ and $\iota^{\sharp}$ are order preserving. It follows from  \ref{Aq}
that they induce {{affine homeomorphisms}} $j_T: T(A)\to T(B)$ and $\iota_T: T(B)\to T(A)$
such that $j_T\circ \iota_T={\rm id}_{T(B)}$ and $\iota_T\circ j_T=\id_{T(A)}.$}}

{\blue{It remains to show, by identifying $K_0(B)$ with $K_0(A),$ that  $\iota^{\sharp}\circ \rho_B=\rho_A.$
Let $x\in K_0(B).$ We may assume that  there {{are}} an integer $i\ge 1$ and $y\in K_0(A_i)$
such that $(\xi_{i,\infty})_{*0}(y)=x.$   By the approximate intertwining diagram above, there  is an integer $k\ge i$ such that, for any $n>k,$
\beq
\|\iota_n^{\sharp}\circ \xi_{i,n}^{\sharp}\circ \rho_{A_i}(y)-\phi_{k,n}^{\sharp}\circ {{\iota_k^{\sharp}\circ}} \xi_{i,k}^{\sharp}\circ\rho_{A_i}(y)\|<\ep.
\eneq
Note that $\iota^{\sharp}\circ \xi_{i,\infty}^{\sharp}(\rho_{A_i}(y))=
\lim_{n\to\infty}\phi_{n,\infty}^{\sharp}\circ \iota_n^{\sharp}\circ\xi_{i,n}^{\sharp}(\rho_{A_i}(y)).$
It follows  that
\beq
\|\iota^{\sharp}\circ \xi_{i,\infty}^{\sharp}\circ \rho_{A_i}(y)-
\phi_{k, \infty}^{\sharp}\circ  \xi_{i,k}^{\sharp}\circ \rho_{A_i}(y)\|\le \ep{{,}}
\eneq
{{where we omit $\iota_k^{\sharp}$ since it is the identity map from $\Aff T(A_k)$ to itself.}} Since $(\xi_{i,k})_{*0}=(\phi_{i,k})_{*0},$ for all $k>i,$
\beq
&&\hspace{-0.4in}\phi_{k, \infty}^{\sharp}\circ \xi_{i,k}^{\sharp}(\rho_{A_i}(y))=\phi_{k,\infty}^{\sharp}\circ \rho_{A_k}\circ (\xi_{i,k})_{*0}(y)\\
&&=\phi_{k,\infty}^{\sharp}\circ \rho_{A_k}\circ (\phi_{i,k})_{*0}(y)=\rho_A\circ (\phi_{i,\infty})_{*0}(y)=\rho_A(x),\andeqn\\
&&\xi_{{{i}}, \infty}^{\sharp}\circ \rho_{A_i}(y)=\rho_B\circ (\xi_{i,\infty})_{*0}(y)=\rho_B(x)
\eneq
{{Therefore,}}
\beq
\|\iota^{\sharp}\circ \rho_B(x)-\rho_A(x)\|\le \ep
\eneq
for any given $\ep>0.$ {{This shows that}} $\iota^{\sharp}\circ \rho_B=\rho_A${{,
and}} completes the proof.
}}

\end{NN}

\begin{cor}\label{range 0.34}
For any $m>0$ and any $A_i$, there {{are}} an integer $n\geq i$ and a
projection $R\in M_m(A_{n+1})$ such that {{the following statements are {{true:}}}}
\begin{enumerate}
\item[(1)] $R$  commutes with $\LD{{\circ}}\xi_{{n,n+1}}(A_n)$, where $\LD:~
A_{n+1}\to M_m(A_{n+1})$ is the amplification map sending $a$ to an
$m\times m$ diagonal matrix:  $\LD(a)=\diag(a, ..., a)$;

\item[(2)] Recall  $
A_{{C},n+1}={{C([0,1],E_{n+1})\oplus_{\bt_{n+1,0}, \bt_{n+1,1}} F_{n+1}=}}A\left(  {{F}}_{n+1}, E_{n+1}, \bt_{{n+1,0}},\bt_{{n+1,1}}
\right)$, where $\bt_{{n+1,0}},\bt_{{n+1,1}}:~ {{F}}_{n+1}\to E_{n+1}$ {{are}} as in
the definition of $A_{n+1}$ (see (\ref{construction-1}).  {{There}} is a unital   injective \hm
$$\iota:~M_{m-1}(A_{{C},n+1}) \lr RM_m(A_{n+1})R$$
such that $R\LD(\xi_{{n,n+1}}(A_n))R \subset \iota(M_{m-1}(
A_{{C},n+1})).$
\end{enumerate}
\end{cor}

\begin{proof}
{\blue{Let $R_X:=R\in M_m(A_{X,n+1}^1)$  {{ be}} {{as}} obtained in \ref{range 0.25} (the definition is given by combining \ref{range 0.22} and \ref{range 0.24}) with the {{property described}}  in \ref{range 0.25}.  Let $\iota_X:=\iota: M_{m-1}(F_{n+1}^1)\to
R_XM_m(A_{X,n+1}^1)R_X$  be the unital injective \hm\, given by \ref{range 0.25}.
Let $\pi_X^1: A_{X,n+1}\to A_{X,n+1}^1.$
Then, since $\phi_{n,n+1}(I_n)\subset I_{n+1},$  one has $\pi_X^1(\phi_{n,n+1}(A_n))=\pi_X^1(\psi_{n,n+1}(A_{X,n})).$
Since $\xi_{{n,n+1}}|_{Sp ({{A_{X,n+1}^1}})}= \phi_{{n,n+1}}|_{Sp ({{A_{X,n+1}^1}})},$
one {{obtains}} \\
(i') $R_X$ commutes with $(\pi_X^1\otimes {\rm id}_m)\circ \Lambda\circ \xi_{n,n+1}(A_n).$

Moreover, one also {{has two additional properties}} }}\\
(ii') \beq\label{1346-1}R_X(\tht_{{n+1,{\blue{1}}}}) = \e_{_{{{F}}_{n+1}^1}}\otimes \left(
                                            \begin{array}{cc}
                                              \e_{m-1} & 0 \\
                                              0 & 0 \\
                                            \end{array}
                                          \right){{,}}
                                          \eneq
 {\blue{and {{consequently}} the {{point evaluation}} $\pi_{\theta_{n+1,1}}: {{M_m(A_{X, n+1}^1)}}\to {{M_{m(F_{n+1}^1)}}}$ {{takes $R_XM_m(A_{X,n+1}^1)R_X$ to $M_{m-1}(F_{n+1}^1)$. Below we will use the same notation $\pi_{\theta_{n+1,1}}$ to denote {{the}} restriction {{of this map}} to $R_XM_m(A_{X,n+1}^1)R_X$, whose codomain is $M_{m-1}(F_{n+1}^1)$.}}

\noindent
 (iii') $\pi_{\theta_{n+1,1}}\circ \iota_X={\rm id}_{{M_{m-1}(F_{n+1}^1)}},$ and $R_X(\pi_X^1\otimes {\rm id}_{M_m})(\Lambda\circ \xi_{n,n+1}(A_n))R_X\subset \iota_X(M_{m{{-1}}}(F_{n+1}^1))$.}}

 One  extends the definition of $R$ as follows. For each $x\in
Sp(A_{n+1})\setminus Sp({{A_{X,n+1}^1}})$, define
\beq\label{1346-2}R(x)=
\e_{_{A_{n+1}|_{x}}}\otimes \left(
                                            \begin{array}{cc}
                                              \e_{m-1} & 0 \\
                                              0 & 0 \\
                                            \end{array}
                                          \right).
                                          \eneq
                                          {\blue{Then define $R(x)=R_X(x)$ for $x\in Sp(A_{X,n+1}^1).$
By checking the boundary, one easily sees that $R\in A_{n+1}$ is a projection.
Then, by (i'), \eqref{1346-1},  and \eqref{1346-2}, $R$ commutes with $\Lambda(\xi_{n,n+1}(A_n)).$}}

{\blue{Define $\iota: M_{m-1}(A_{C,n+1})\to RM_m(A_{n+1})R$ by
\beq
\iota(f,a_1,a_2,...,a_{p_{n+1}})=(f, \iota_X(a_1),a_2,...,a_{p_{n+1}})
\eneq
for $f\in M_{m-1}(C([0,1], E_{n+1}))$ and $(a_1,a_2,...,a_{p_n})\in {{M_{m-1}(F_{n+1})}}={{\bigoplus}}_{i=1}^{p_{n+1}} {{M_{m-1}(F_{n+1}^i)}}.$
By (iii'), $\pi_{\theta_{n+1, 1}}\circ \iota_X(a_1)=a_1.$
Since $(f,a_1,a_2,...,a_{p_n})\in M_{m-1}({{A}}_{C,n+1}),$
we have
\beq
f(0)=\bt_{n+1,0}((a_1,a_2,...,a_{p_{n+1}})\andeqn f(1)=\bt_{n+1,1}((a_1,a_2,...,a_{p_{n+1}}).
\eneq
Note that
\beq\nonumber
\pi_{X, n+1}^e((\iota_X(a_1),a_2,...,a_{p_{n+1}})&=&(\pi_{\theta_{n+1,1}}(\iota_X(a_1)), a_2,...,a_{p_{n+1}})
= (a_1,a_2,...,a_{p_{n+1}}).
\eneq
Thus,
\beq\nonumber
f(0)=\bt_{n+1,0}\circ \pi_{X,n+1}^e((\iota(a_1), a_2,...,a_{p_{n+1}})\andeqn
f(1)=\bt_{n+1,1}\circ \pi_{X, n+1}^e(\iota(a_1), a_2,...,a_{p_{n+1}}).
\eneq
Since $\iota_X$ is unital and injective, one checks that $\iota$ just defined is also unital and injective.
In other words, $\iota$ maps $M_{m-1}(A_{{{C,}}n+1})$ to $M_m(A_{n+1}).$
Note that $\iota_X(M_{m-1}(F_{n+1}^1))\subset R_X(M_m(A_{X, n+1}^1)R_X.$ Then, by \eqref{1346-1} and
\eqref{1346-2}, $\iota(M_{m-1}(A_{{{C,}}n+1}))\subset RM_m(A_{n+1})R.$}}

By (iii'), $(\pi_{X, n+1}^1\otimes {\rm id}_{M_m})(R(\Lambda({{\xi}}_{n,n+1}(A_n)))R)\subset (\pi_{X, n+1}^1\otimes {\rm id}_{M_m})(\iota(M_{m-1}(A_{{C,n+1}}))).$  On the other hand, $(\pi_{C,n+1}^e\otimes {\rm id}_{M_m})\circ \iota={\rm id}_{M_{m-1}(A_{C, n+1})}.$
By \eqref{1346-2}, $(\pi_{J,n+1}\otimes {\rm id}_{M_m})(R)=1_{M_{m-1}(A_{C,n+1})}.$  It follows
that
$$
(\pi_{J,n+1}\otimes {\rm id}_{M_m})(R(\Lambda(\xi_{n,n+1}(A_n))R))\subset M_{m-1}(A_{C,n+1})
=(\pi_{J,n+1}\otimes {\rm id}_{M_m})(\iota(M_{m-1}(A_{C,n+1}))).
$$
Since $({\rm ker}\pi_{X, n+1}^1\otimes {\rm id}_{M_m})\cap (J_{n+1}\otimes {\rm id}_{A_m})=\{0\},$ it  follows  from
that  $R(\Lambda(\xi_{n,n+1}(A_n))R\subset \iota(M_{m-1}(A_{{C,n+1}})).$
This completes the proof.
\end{proof}



\begin{cor}\label{range 0.34a}
Let $B$ be as constructed above. Then $B\otimes U \in {\cal B}_0$  for every UHF-algebra $U$ of infinite dimension.
\end{cor}

\begin{proof}
{\blue{Note that, by \ref{range 0.5} (see lines below \eqref{1311k0}),}} $A_{C, n+1}\in {\cal C}_0.$
{\blue{Thus, in \ref{range 0.34},}}
$M_{m-1}(A_{{C},n+1})$ and $ \iota(M_{m-1}(A_{{C},n+1}))$ are in
${\cal C}_0$. Also, {\blue{for each $n$}} {{and each}}
$\tau \in
T(M_m(A_{n+1}))$, we have $\tau(\e-R)=1/m$.
Fix  an integer $k\ge 1$ and
{{a}}  finite subset ${\cal F}\subset B\otimes M_k,$ and
let ${\cal F}_1\subset B$ be a finite subset
such that $\{(f_{ij})_{k\times k}: f_{ij}\in {\cal F}_1\}\supset {\cal F}.$  Now, {\blue{by  applying  \eqref{range 0.34},
one shows that}} the
inductive limit algebra $B=\lim_{n\to\infty}(A_n, \xi_{n,m})$
{{has}} the following
property:
For any finite set ${\cal F}_1\subset B$, $\ep>0$, $\dt>0$, and any
$m>1/\dt$, there is a unital \SCA\,  $C\subset M_m(B)$ {{with different unit  $\e_C$}} which is
in ${\cal C}_0$ such that
\begin{enumerate}
\item[(i)] $\|[\e_C, \diag\{\underbrace{f,...,f}_m\}]\|<\ep/k^2$, for all
$f\in {\cal F}_1$,
\item[(ii)] $\dist(\e_C (\diag\{\underbrace{f,...,f}_m\})\e_C, C)<\ep/k^2$,
for all $f\in {\cal F}_1$, and
\item[(iii)] $\tau(\e_{{M_m(B)}}-\e_C)=1/m<\dt$ for all $\tau \in T(M_{m}(B))$.
\end{enumerate}
{{Consequently,
\begin{enumerate}
\item[(i')] $\|[\e_{M_k(C)}, \diag\{\underbrace{f,...,f}_m\}]\|<\ep$, for all
$f\in {\cal F}$,
\item[(ii')] $\dist(\e_{M_k(C)} (\diag\{\underbrace{f,...,f}_m\})\e_{M_k(C)}, M_k(C))<\ep$,
for all $f\in {\cal F}$, and
\item[(iii')] $\tau(\e_{M_{mk}(B)}-\e_{M_k(C)})=1/m<\dt$ for all $\tau \in T(M_{mk}(B))$.
\end{enumerate}}}

Now  $B\otimes U$ can be written as $\lim_{n\to\infty}(B\otimes M_{k_n}, \iota_{n,m})$ with $k_1|k_2|k_3\cd$ and $k_{n+1}/k_n \to \infty$, and ${\blue{\Lambda}_{n,{n+1}}}$ is the amplification {{map}}
by sending $f\in B\otimes M_{k_n}$ to  $\diag{(f,..., f)}\in B\otimes M_{k_{n+1}},$ where $f$ repeated $k_{n+1}/k_n$
times.

To show $B\otimes U\in {\cal  B}_0,$  let ${\cal F} \subset  B\otimes U$
be a finite subset and let $a \in (B\otimes U)_+\setminus  \{0\}.$ There is an integer $m_0>0$ such that $\tau(a)>1/m_0$ for all $\tau\in B\otimes U$.  Without loss of generality, we may assume that
 ${\cal F} \subset B\otimes M_{k_n}$ with ${{k_{n+1}/k_n}} >m_0$. Then by {{ (i'), (ii'), and (iii')}}
 for $B\otimes M_{k_n}$ {{(i.e., $k=k_n$)}} with $m=k_{n+1}/k_n$ (and recall that $\iota_{n,n+1}$ is the amplification), there is a unital \SCA\,
$D:=M_k(C) \subset B\otimes M_{k_{n+1}}$ with $C\in {\cal C}_0$ such that $\|[\e_D, \,\iota_{n,n+1}(f) ]\|<\ep$, for all
$f\in {\cal F}$, such that $\dist(\e_D (\iota_{n,n+1}(f))\e_D, D)<\ep$ for all $f\in {\cal F}$, and such that $\tau(\e-\e_D)=1/m<\dt$ for all $\tau \in T(M_{k_{n+1}}(B))$. Then ${\blue{\Lambda_{n+1,\infty}}}(D)$ is the desired subalgebra. (Note that $1-\e_D\lesssim a$ follows from {{the}} strict comparison property of $B\otimes U$ {\blue{(see 5.2 of \cite{RorUHF2})}}.
It follows that $B \otimes U\in {\cal B}_0.$
 \end{proof}


\begin{thm}\label{range 0.35}
For any simple weakly unperforated Elliott invariant $\big((G,G_+,u),~
K, \DT, r \big)$, there is a unital simple \CA\, $A\in {\cal
N}_0^{\cal Z}$ which is an inductive limit of $(A_n, \phi_{n,m})$ with $A_n$
{{as}} described {\blue{in}} \ref{range 0.19}, with $\phi_{n,m}$
injective,  such that
$$\big((K_0(A),K_0(A)_+, \e_A),~ K_1(A), T(A), r_{\!\!_A} \big)\cong
\big((G,G_+,u),~ K, \DT, r \big).$$
\end{thm}

\begin{proof} By \ref{range 0.34a}, $A\in {\cal N}_0$. Since $A$ is a unital simple inductive limit of subhomogeneous \CA s with no dimension growth{{, by}} Corollary 6.5 of \cite{Winter-Z-stable-01}, $A$ is ${\cal Z}$-stable.
\end{proof}

\begin{cor}\label{range 0.35a} For any simple weakly unperforated Elliott invariant $\big((G,G_+,u),~
K, \DT, r \big)$ with $K=\{0\}$ and $G$ torsion free, there is a unital ${\cal Z}$-stable simple \CA\, which is an inductive limit of $(A_n, \phi_{n,m})$ with $A_n$ in ${\cal C}_0$ {{as}}  described in \ref{DfC1}, with $\phi_{n,m}$
injective,  {{such}} that
$$\big((K_0(A),K_0(A)_+, \e_A),~ K_1(A), T(A), r_{\!\!_A} \big)\cong
\big((G,G_+,u),~ 0, \DT, r \big).$$

\end{cor}

\begin{proof} In the construction of $A_n$, just let all the spaces $X_n$ involved {{be}} the space {{consisting}} of a single point.

\end{proof}

\section{Models for \CA s in ${{\cal N}}_0$ with property (SP)}

{{Let us recall some notation {{concerning the}} classes of $C^*$-algebras used in this section. ${\cal C}$ is the class of Elliott-Thomsen building blocks defined in Definition \ref{DfC1}, and ${\cal C}_0$ consists of {{the  $C^*$-algebras}} in ${\cal C}$ with zero {{$K_1$-group}}. ${\cal D}_k$ is a class of recursive {{subhomogenous}} algebras defined in Definition \ref{8-N-3}.  ${\cal N}$ is the class of all separable amenable \CA s which satisfy the Universal Coefficient Theorem (UCT).  ${\cal B}_0$ and ${\cal B}_1$ are the classes of $C^*$-algebras defined in {{Definition}} \ref{DB1} (roughly speaking ${\cal B}_0$ (${\cal B}_1$, respectively) contains the $C^*$ algebras which can be approximated by $C^*$ algebras in the {{class}}  ${\cal C}_0$ (${\cal C}$, respectively) tracially). As in Definition \ref{Class0},
${\cal N}_0$ is  the class of  unital  simple \CA s
$A$  in ${\cal N}$ for which $A\otimes U\in {\cal
N}\cap {\cal B}_0$.}}

\begin{NN}\label{range 0.37} For technical  reasons, in the construction of our model algebras,
it is important for us to be able to decompose
$A_n$ into {{the}} direct sum of two parts: the homogeneous part
which stores the information of $\Inf K_0(A)$ and $K_1(A)$ and the
part of  {{the}} algebra in ${\cal C}_0$ which stores {{the}} information of
$K_0(A)/\Inf K_0(A)$, $T(A)$ and the pairing between {{these}}. {{This}} {{cannot}}
be done in general for the algebras in ${\cal N}_0$ (see \ref{Class0}), 
{{but}} we will prove that this can be done
if the Elliott invariant satisfies an extra condition{{, property (SP),}} described
below.


  Let $\big((G,G_+,u),~
K, \DT, r \big)$ be a weakly unperforated Elliott invariant as {{in}}
\ref{range 0.1}. We say that it has the {{property (SP)}} if for any
real number $s>0$, there is $g\in G_+\setminus \{0\}$ such that
$\tau (g) <s$ for any state $\tau$ on $G$, or equivalently,
$r(\tau)(g)<s$ for any $\tau \in \DT$. In this case, we will prove
that the algebra in \ref{range 0.35} can be chosen to be in {{the}} class
${\cal B}_0$ (rather than in the larger class ${\cal N}_0=\{A{:} \,
A\otimes M_{\mathfrak{p}}\in {\cal B}_0\}$). Roughly speaking, for
each $A_n$, we will separate the part of {{the}} homogeneous algebra which
will store all the information of {{the}} infinitesimal part of $K_0$  and
$K_1$, and it will be in the corner $P_nA_nP_n$ with $P_n$ small {{compared}} to
$\e_{A_n}$ in the limit algebra. In fact, the construction of this
case is much easier, since the homogeneous blocks can be {{separated}}
from the part {{in}} ${\cal C}_0$---we will first write the group
inclusion $G_n\hookrightarrow H_n$ as in \ref{range 0.5a}.

{{Let us point out}} that if $A\in {\cal N}_0$ then the Elliott invariant of $A\otimes {{U}}$ has
{{property}}  (SP) {{for any infinite dimensional UHF algebra $U$,}} even though the Elliott invariant of $A$ itself may not {{have the property}}. {{One can verify {{this}} fact as follows. As $U$ is a UHF algebra of infinite dimension, $(K_0(U), K_0(U)_+, [\one_U])=(\mathbb{P},\mathbb{P}\cap \R_+, 1)$, where $\mathbb{P}\subset \Q\subset \R$ is a dense {{subgroup of }}  $\R$. For any positive number $s$, by density of $\mathbb{P}$, we can choose a number $r\in (\mathbb{P}\cap \R_+)\setminus\{0\}$ such that $r<s$. Let $x\in U$ be a projection such that $[x]=r\in K_0(U)$. Then the projection $\one_A\otimes x\in A\otimes U$ {{satisfies}}  $\tau (\one_A\otimes x)=r<s$ for all $\tau \in T(A\otimes U)$.}}

\end{NN}

\begin{NN}\label{range 0.38}
Let $((G,G_+,u),K,\Delta,r)$ be {{as}} given in \ref{range 0.1} or
\ref{range 0.16}. As in \ref{range  0.3}, let $\rho:~ G\to \Aff\DT$
be dual to the map $r$. Denote {{ the}} kernel of the map $\rho$ by
$\Inf(G)$---the infinitesimal part of $G$, {{i.e.,}}
$$\Inf(G)=\{g\in G: \rho(g)(\tau)=0,\ {{\rforal}} \tau\in\DT\}.$$
Let $G^1\subset \Aff \DT$ be a countable dense subgroup which is $\Q$-linearly
independent with $\rho(G)$---that is, if $g\in \rho(G)\otimes\Q$ and
$g^1\in G^1\otimes\Q$ satisfy $g+g^1=0$, then both $g$ and $g^1$ are
zero. Note that such $G^1$ {{exists}}, since $\Q$  {\blue{is a}} vector space, and {{the}}
dimension of $\rho(G)\otimes \Q$ is countable, but the dimension of $\Aff
(\DT)$ is uncountable.  Again as in \ref{range 0.3}, let $H=G\oplus
G^1$ with $H_+\setminus \{0\}$ {{the }}  {set} of $(g,f)\in
G\oplus G^1$ with
$$\rho(g)(\tau)+f(\tau)>0  \rforal \tau\in \DT.$$
The scale $u\in G_+$ {{may be}} regarded as $(u,0)\in G\oplus G^1=H$
{{and so}} as the scale of $H_+$. Since $\rho(u)(\tau)>0,$ it follows  that $u$ is an order unit for $H.$ Since $G^1$ is {{$\Q$-linearly}}
independent {{of}} $\rho(G)$, we know $\Inf(G)=\Inf(H)$---that is, when
we embed $G$ into $H$, it does not create more elements in the
infinitesimal group. {{Evidently,}} ${\rm Tor}(G)={\rm Tor}(H)\subset \Inf(G)$. Let
$G'=G/\Inf(G)$ and $ H'=H/\Inf(H)${{. Then}} we have the following
diagram:

\begin{displaymath}
    \xymatrix{
       0 \ar[r] & \Inf(G) \ar[r]\ar@{=}[d]&G_{\,}
     \ar[r]\ar@{_{(}->}[d] &{ G'_{\,}}\ar[r]\ar@{_{(}->}[d]&0 \\
         0 \ar[r] & \Inf(H) \ar[r]&H
         \ar[r] &{ H'}\ar[r]&0.}
\end{displaymath}

Let ${ G'}_{{+}}$ (or $H'_{{+}}$), and $u'$  be the image of $G_+$ (or $H_+$)
and $u$ under the quotient map from $G$ to $G'$ (or from $H$ to
$H'$). Then $(G', G'_{{+}}, u')$ is {\blue{a}} weakly unperforated group without
infinitesimal elements.  Note that $G$ and $H$ share {\blue{the}} same unit $u$, and
therefore $G'$ and $H'$ share the same unit $u'$. Since
$r(\tau)|_{\Inf(G)}=0$ for any $\tau \in \DT$,  the map $r:\DT \to
S_u(G)$ induces a map $r': \DT \to S_{u'}(G')$. Hence $((G', G'_+, u'),
\{0\}, \DT, r')$ is a weakly unperforated Elliott invariant with trivial
$K_1$ group and no infinitesimal elements in the $K_0$-group.
\end{NN}

\begin{NN}\label{range 0.39}
   With the same argument as that of \ref{range 0.5}, we have the following diagram of
   inductive limit{s}:
\begin{displaymath}
    \xymatrix{{ G'}_1 \ar[r]^{\af'_{{12}}}\ar@{_{(}->}[d]_{\iota_1} & {G'}_2 \ar[r]^{\af'_{{23}}}\ar@{_{(}->}[d]_{\iota_2}&\cd \ar[r]&{ G'_{\,}} \ar@{_{(}->}[d]_{\iota}\\
         { H'}_1 \ar[r]^{\gm'_{{12}}}  & { H'}_2 \ar[r]^{\gm'_{{23}}}&\cd \ar[r]&{ H'},  }
\end{displaymath}
where each ${ H'}_n{{=\Z^{p_n}}}$ is a   direct sum of {{finitely  many}} copies of
{{the group}} $\Z$ {{with the positive cone $(H'_n)_+=(\Z_+)^{p_n}$}}, $\af'_{n, n+1}=\gm'_{n, n+1}|_{{G'}_n}$, and
${H'}_n/{G'}_n$ is a free abelian group. {{Note that not only {{is $(H'_n)_+$}}  finitely generated, but also $(G'_n)_+:=(H'_n)_+\cap G'_n$ is finitely generated{{, by  Theorem}} \ref{FG-Ratn}. As in \ref{range 0.16}, we may assume that
{{all}}
 $\gamma_{n,n+1}'$ are at least 2-large.}}

By \ref{range 0.15a}, we can construct an increasing sequence of
finitely generated subgroups
$$\Inf_1\subset \Inf_2\subset \Inf_3\subset\cd\subset \Inf(G),$$
with $\Inf(G) =\bigcup_{i=1}^{\infty} \Inf_n,$ and {{such that one has}} the inductive limit
\begin{displaymath}
\xymatrix{
\Inf_1\oplus H'_1 \ar[r]^-{\gm_{1,2}} &  \Inf_2\oplus H'_2 \ar[r]^-{\gm_{2,3}} & \Inf_3\oplus H'_3 \ar[r]^-{\gm_{3,4}} & \cd \ar[r] & {{H.}}
}
\end{displaymath}


Put $H_n:=\Inf_n\oplus
H'_n$ and $G_n:=\Inf_n\oplus G'_n.$ Since  $G'_n$ is a subgroup of
$H'_n$, the group   $G_n$ is also a subgroup of $H_n$. Define $\af_{n,n+1}: G_n\to
G_{n+1}$ by $\af_{n, n+1}=\gm_{n,n+1}|_{G_n}$, which is
compatible with $\af'_{n,n+1}$ in the sense that (ii) of \ref{range
0.15a} holds. Hence we obtain the following diagram of inductive limits:

\begin{displaymath}
    \xymatrix{{ G}_1 \ar[r]^{\af_{{12}}}\ar@{_{(}->}[d]_{\iota_1} & {G}_2 \ar[r]^{\af_{{23}}}\ar@{_{(}->}[d]_{\iota_2}&\cd \ar[r]&{ G_{\,}} \ar@{_{(}->}[d]_{\iota}\\
         { H}_1 \ar[r]^{\gm_{{12}}}  & { H}_2 \ar[r]^{\gm_{{23}}}&\cd \ar[r]&{ H},  }
\end{displaymath}
 with $\af_{n,n+1}(\Inf_n)\subset \Inf_{n+1}$ and $\af_{n,n+1}|_{\Inf_n}$ {{the}} inclusion map.
{\blue{\it We will fix, for each $n,$ a positive non-zero \hm\, $\lambda_n: H_n\to\Z$ such
that $\lambda_n(x)>0$ for any $x\in (H_n')_+\setminus \{0\}$}} {{and such that $\ld_n|_{\Inf_n}=0$.}}
Note that all notations discussed so far in this section will be used for the rest of this section.

\end{NN}

 \begin{lem}\label{range 0.40} Let $(G, G_+, u)=\lim_n((G_n,(G_n)_+,u_n), \af_{n,m})$
and $(H,H_+,u)=$ \linebreak $\lim_n((H_n,(H_n)_+,u_n), \gm_{n,m})$  be as above.
Suppose that $\big((G,G_+,u),~ K, \DT, r \big)$ has {{the}} {{property (SP)}}.
For any $n$ with
$G_n\stackrel{\iota_n}{\hookrightarrow}H_n,$  and for any $D=\Z^k$ (for any positive integer $k$), {{there are positive maps}} $(\kappa_n, \id): H_n \to D\oplus H_n$ {{and}} $(\kappa_n', \id): G_n \to D\oplus G_n$ {{such that, for}}
any integer
$L>0$, there are {{an}} integer $m(n)>n$
and positive maps
$\eta: D\oplus H_n \to H_{m(n)}$
{{and}} $\eta': D\oplus G_n \to G_{m(n)}$ such that the following  diagram
commutes:
\begin{displaymath}
\xymatrix{
G_n \ar[rr]^{\alpha_{n, m}} \ar[dr]_{({\kappa_n}', \id)} \ar[ddd]_{\iota_n} &  &  G_m \ar[ddd]^{\iota_m} \\
 &D \oplus G_n  \ar[ur]_{\eta'} \ar@{_{(}->}[d] _{(\id, \iota_n)}&   \\
 &D \oplus H_n \ar[dr]^{\eta}&  \\
H_n\ar[ur]^{({\kappa_n}, \id)} \ar[rr]_{\gm_{n, m}} &   & H_{m(n)}, }
\end{displaymath}
and such that the following {{statements}} are true:

{\rm (1)} {\blue{The map $\kappa_n: H_n\to D=\Z^k$ is defined by $\kappa_n(x)=(\lambda_n(x),\lambda_n(x), ..., \lambda_n(x))$
for $x\in H_n$ and $\kappa_n'={\kappa_n}|_{{G_n}},$ in particular,}}
each
component of $\kappa_n (u_n)=\kappa_n'(u_n)$ in $\Z^k$ is strictly
positive.

{\rm (2)} For any $\tau\in \DT$,
$$r(\tau)((\af_{m(n),\infty}\circ\eta')(\e_D))(=r(\tau)((\gm_{m(n),\infty}\circ\eta)(\e_D)))<1/L.$$

{\rm (3)} Each component of the map ${\blue{\pi'\circ \eta|_{D\oplus H_n'}:}}D\oplus H'_n
=\Z^k\oplus\Z^{p_n}\to H'_{m(n)}=\Z^{p_{m(n)}}$ {\blue{(where
$\pi': H_m\to H_m'$ is the projection map)}} is  strictly positive, and L-large---i.e., all
entries in the $(k+p_n)\times p_{m(n)}$ matrix corresponding to the map
are  larger than $L$. {\blue{Moreover, $\eta'=\eta\circ {{(\id,\iota_n)}}$
and  $\pi'\circ \eta'(D\oplus G_n')\subset G_{n+1}'.$}}

({{Note the maps $(\kappa_n, \id)$ and $(\kappa_n', \id)$}} are independent of $L$.)

\end{lem}

\begin{proof}
We will use the following  fact {{(which {{was}} pointed out in the first paragraph of \ref{range 0.39})}} several times: the positive cone of
$G_n'$ (and {{that}} of $H'_n$) is
 finitely generated (note that
 even though $G_n$ and $H_n$ are finitely generated, their positive {{cones}} may not be
 finitely generated).
Note that $H_n'$ is {{a}} subgroup of $H_n$ so we {{will  continue}} to use $\lambda_{{n}}$ for $\lambda_{{n}}|_{H_n'}.$

{\blue{We now fix $n$ and an integer $k>0.$}} {{Define $\kappa_n:H_n\to D=\Z^k$ and $\kappa_n:G_n\to D$ by}}
$$\kappa_n(a)=(\underbrace{\ld_n(a), ..., \ld_n(a)}_k)\in D
\andeqn\kappa_n'(a)=(\underbrace{\ld_n\circ \iota_n(a), ..., \ld_n\circ \iota_n(a)}_k)\in D.$$
 Since $G$ has the {{ property (SP)}} {\blue{and since $(H_n')_+$ is finitely generated,}}
 there is $p'\in G_+\setminus\{0\}$ such
 that for any $a\in (H_n')_+,$
 $$\gm_{n, \infty}(a)-k\cdot{\blue{\ld_n}}(a)\cdot p'\in H_+.$$
 {{Consequently,}}
 $$\alpha_{n, \infty}(a)-k\cdot{\blue{\ld_n\circ \iota_n(a)}}\cdot p'\in G_+\rforal a\in (G_n')_+,$$
 where
the {{maps}} $\alpha_{n, \infty}$ and $\gm_{n, \infty}$ are the
homomorphisms from $G_n$ to $G$ {{and}} from $H_n$ to $H$ respectively.
Moreover, {{for any positive integer $L$,}} one may require that
\begin{equation}\label{14star}
r(\tau)(\ld_n(u_n)\cdot p')<1/2kL  \rforal \tau\in
\DT.
\end{equation}
Since $(G_n')_+$ and ${{(H_n')_+}}$ are finitely generated, there {{are}} an integer  $m(n)\ge 1$
and $p\in {{(G_{m(n)})_+}}$ such that
\beq\nonumber
\af_{m(n),\infty}(p)=p',\,\,\,
\alpha_{n, m(n)}(a)-k\cdot\ld_n\circ \iota_n(a)\cdot p \in {{(G_{m(n)})_+}} \rforal a\in (G_n')_+,\andeqn\\\nonumber
\gm_{n, m(n)}(a)-k\cdot\ld_n(a)\cdot p \in {{(H_{m(n)})_+}}  \rforal a\in
{{(H_n')_+}}.
\eneq
Then define $\alpha_n'': G_n \to G_{m(n)}$ and $\gm_n'': H_n
\to H_{m(n)}$ by
\beq\nonumber
&&\alpha_n'': G_n\ni a\mapsto \alpha_{n, m(n)}(a)-k\cdot\ld_n\circ \iota_n(a)\cdot p\in G_{m(n)},\andeqn\\
&&\gm''_n: H_n\ni a\mapsto\gm_{n, m(n)}(a)-k\cdot\ld_n(a)\cdot p\in H_{m(n)}.
\eneq
By the choice of $p$, the maps $\alpha''_n$ and
$\gm_n''$ are positive. {\blue{(Note that
$\alpha_n''|_{\Inf_n}={\af_{n,m(n)}}|_{\Inf_n}$ and ${\gm''_n}|_{\Inf_n}={\gm_{n,m(n)}}|_{\Inf_n}.$)}}

 A direct calculation shows the following diagram commutes (where $D=\Z^k$):
\begin{displaymath}
\xymatrix{
G_n \ar[rr]^{\alpha_{n, m(n)}} \ar[dr]_{({\kappa_n}', \id)} \ar[ddd]_{\iota_n} &  &  G_{m(n)} \ar[ddd]^{\iota_{m(n)}} \\
 &D \oplus G_n  \ar[ur]_{\eta'} \ar@{_{(}->}[d] _{(\id, \iota_n)}&   \\
 &D \oplus H_n \ar[dr]^{\eta}&  \\
H_n\ar[ur]^{({\kappa_n}, \id)} \ar[rr]_{\gm_{n, m(n)}} &   & H_{m(n)}, }
\end{displaymath}
$$\eta'((m_1, ..., m_k, g))= (m_1+\cdots+m_k)p+\alpha''_n(g),\andeqn$$
$${{\eta}}((m_1, ..., m_k, g))= (m_1+\cdots+m_k)p+\gm''_n(g).$$
$$(\mbox{{Recall that}}~~\kappa_n'(a)=(\underbrace{\ld_n\circ \iota_n(a), ..., \ld_n\circ \iota_n(a)}_k)\in D\andeqn
\kappa_n(a)=(\underbrace{\ld_n(a), ..., \ld_n(a)}_k)\in D.)$$
 The order
 of $D\oplus G_n$  and $D\oplus H_n$ are the standard
  order on direct sums, i.e., $(a, b)\geq 0$ if and only if $a\geq 0$ and $b\geq 0$.
   Since the maps $\alpha''_n$ and $\gm''_n$ are positive, the
   maps $\eta'$ and $\eta$ are positive. Condition (1) follows from
   the {construction;} condition (2) follows from (\ref{14star}){{; and}} condition
   (3) follows from {\blue{the fact {{that }}
   $\gm_{k,k+1}'$ are $2$-large,}}
   if one passes to {{a}} further stage (choose larger
    $m(n)$).
\end{proof}

\begin{df}\label{AHblock}
A \CA\, is said to be in the class ${\bf H}$\index{${\bf H}$} if it is the direct sum
of  {{algebras}} of  the form $P(C(X)\otimes M_n)P$, where $X=\{pt\},
[0,1], $ $S^1, $ $S^2,$ $T_{2,k}${{, or $T_{3, k}$}} {{(see \ref{range 0.18} for the {{definitions}} of $T_{2,k}$ and $T_{3, k}$).}}
{{In addition, we assume that the rank of $P$ is at least ${{13}}$ when $X=T_{2,k}$
or $X=T_{3, k}.$}}
\end{df}

\begin{NN}\label{range 0.41}
Write $K$ (the $K_1$ part of the invariant) as the union of {{an}} increasing
sequence of finitely generated abelian subgroups:  $K_1\subset K_2 \subset
K_3\subset\cd\subset K$ with $K=\bigcup_{i=1}^{\infty} K_i.$
{\blue{Denote by $\chi_{n,n+1}: K_n\to K_{n+1}$ the embedding, $n=1,2,....$}}
For a finitely generated abelian group $G$, we use ${\rm rank}\,
G$ to denote {{the smallest cardinality of a}}   {{set of generators}} of
$G$---that is, $G$ can be written as a direct sum of
${\rm rank}(G)$ cyclic groups (e.g., $\Z$ or $\Z/m\Z$, $m\in \N$).

Let $d_n=\mbox{max}\{2,~1+\mathrm{rank}(\Inf_n)+\mathrm{rank}(K_n)$\}. Apply \ref{range
0.40} repeatedly, with $k=d_n$ and  $L> 13\cdot 2^{2n}$ and
$\td\gm_{n,n+1}=(\kappa_n, \id)\circ \eta_n,$ where $\kappa_n$  and $\eta_n$ are defined
in \ref{range 0.40},
for each $n${{.   Passing to a}} subsequence {{if necessary}}, we {\blue{obtain}} the following {{commutative}}
diagram of inductive limits:
\begin{displaymath}
\xymatrix{
\Z^{d_1} \oplus G_1\ar[r]^{\tilde{\alpha}_{1, 2}}\ar@{_{(}->}[d]_{(\id, \iota_1)} & \Z^{d_2} \oplus G_2 \ar@{_{(}->}[d]_{(\id, \iota_2)} \ar[r]^-{\tilde{\alpha}_{2, 3}} & \cdots \ar[r] & G_{\,} \ar@{_{(}->}[d]_{\iota} \\
\Z^{d_1} \oplus H_1\ar[r]_{\tilde{\gm}_{1, 2}} & \Z^{d_2} \oplus H_2
\ar[r]_-{\tilde{\gm}_{2, 3}} & \cdots \ar[r] & H,}
\end{displaymath}
{\blue{where ${\bar G}_n:=\Z^{d_n}\oplus G_n$ has the order unit $(\kappa_n'(u_n), u_n),$  $\Z^{d_n}\oplus H_n$ has
the order unit  $(\kappa_n(u_n), u_n),$
{{and}}  $\td\af_{n,n+1}={\td\gm_{n,n+1}}|_{\Z^{d_n}\oplus G_n}.$
Set $u_n^d:=\kappa_n'(u_n)\in \Z^{d_n}.$ Then
$u_n^d=\kappa_n(u_n)=\kappa_n'(u_n).$  {{Write}} ${\bar u}_n=(\kappa_n(u_n), u_n)$ and  $G_n''=\Z^d\oplus \Inf_n.$  So ${{\bar{G}_n=}}\Z^{d_n}\oplus G_n=G_n''\oplus G_n'.$  We may also write  ${\bar u}_n=(u_n'', u_n'),$ where $u_n''\in G_n''.$   Then $\td\af_{n,n+1}({\bar u}_n)={\bar u}_{n+1}.$
Let $\rho_n': \Z^{d_n}\oplus \Inf_n\oplus H_n'\to \Z^{d_n}\oplus H_n'$ be the quotient map, and
let $\rho_{G',n}':= {\rho_n'}|_{\Z^{d_n}\oplus {{G}}_n}$ be the map which maps ${{\Z^{d_n}\oplus G_n=}}\Z^{d_n}\oplus \Inf_n\oplus G_n'$ to $\Z^{d_n}\oplus G_n'.$
Put $\td\gm_{n,n+1}'=\rho_n'\circ \td \gm_{n,n+1}|_{\Z^{d_n}\oplus H_n'}$ and
$\td\af_{n,n+1}'=\rho_{G',n+1}'\circ {\td\af_{n,n+1}}|_{\Z^{d_n}\oplus G_n'}.$}}
{\blue{Note also  that, {{if we replace}} $G_n$ by $\Z^{d_n}\oplus G_n,$ and $H_n$ by $\Z^{d_n}\oplus H_n,$
respectively,  {{then}} the limit ordered groups
$G$ and $H$ do not change (see  \ref{range 0.40}).
In particular,  we still have $\Inf(G) =\bigcup_{i=1}^{\infty} \td\af_{n, \infty}(\Inf_n).$
Moreover,  by (2) of \ref{range 0.40},  $r(\tau)(\td\af_{n,\infty}(u_n''))<1/13 \cdot 2^{2n}$
for all $\tau\in \Delta.$
{{For exactly}} the same reason one has the following commutative diagram:
\beq\label{14afK0}
\xymatrix{
\Z^{d_1} \oplus G_1'\ar[r]^{\td \af_{1, 2}'}\ar@{_{(}->}[d]_{(\id, \iota_1)} & \Z^{d_2} \oplus G_2' \ar@{_{(}->}[d]_{(\id, \iota_2)} \ar[r]^-{\tilde{\alpha}_{2, 3}'} & \cdots \ar[r] & G'_{\,} \ar@{_{(}->}[d]_{\iota} \\
\Z^{d_1} \oplus H_1'\ar[r]_{\tilde{\gm}_{1, 2}'} & \Z^{d_2} \oplus H_2'
\ar[r]_-{\tilde{\gm}_{2, 3}'} & \cdots \ar[r] & H'.}
\eneq
}}\\
{\blue{Note that $\Z^{d_n}\oplus G_n'$ and $\Z^{d_n}\oplus H_n'$ share the order unit $(u_n^d, u_n').$}}
 {\blue{Let  ${\bar H}_n=Z^{d_n}\oplus H_n
=G_n''\oplus H_n'.$
Note that
each
$\td\gm_{k,k+1}'$ {{is}}
$13 \cdot 2^{2k}$-large{{; that is,}} {{ $\td\gm_{k,k+1}'=(c_{ij}^{k,k+1})\in M_{(d_{k+1}+p_{k+1})\times (d_{k}+p_{k})}(\Z_{+})$ with $c_{ij}^{k,k+1}\geq 13 \cdot 2^{2k}$}}.}}


{\blue{We now construct \CA s $\{C_n\},$ $\{B_n\},$ $\{F_n\}${{, and}} $\{A_n\},$
and unital injective {{\hm s}}\,
$\phi_{n, n+1}: A_n\to A_{n+1}${{,
inductively}}
as in Section 13.}}

As {{in}}  \ref{range 0.5a} (applied to $G_n'\subset H_n'$),
one can find  finite dimensional $C^*$-algebras $F_n={{\bigoplus}}_{i=1}^{p_n} F_n^i$ and $E_n={{\bigoplus}}_{j=1}^{l_n}E_n^j$,
unital homomorphisms ${\blue{\bt_{n,0},~\bt_{n,1}}}: F_n \to E_n$, and {{form}} the C*-algebra
\beq\nonumber
C_n=A(F_n,E_n,\bt_0,\bt_1):= \{(f,a)\in C([0,1], E_n)\oplus F_n;
f(0)=\bt_{n,0}(a), f(1)=\bt_{n,1}(a)\}\\ \nonumber
{\blue{:=C([0,1], E_n)\oplus_{\bt_{n,0}, \bt_{n,1}} F_n}},\qq\qq\qq\qq
\eneq
such that
$$(K_0(F_n),K_0(F_n)_+, [\e_{F_n}])=(H_n', {{(H_n')_+}},u_n'),$$
$$(K_0(C_n), K_0(C_n)_+,\e_{C_n})=(G_n', {\blue{(G_n')_+}},u_n'),\quad K_1(C_n)=\{0\},$$
and furthermore $K_0(C_n)$ is identified with
$$\ker((\bt_{n,1})_{*0}-(\bt_{n,0})_{*0})=\{x\in K_0(F_n); ((\bt_0)_{*0}-(\bt_1)_{*0})(x)=0\in
K_0(E_n)\}.$$

{\blue{ Write
$G_n''=\bigoplus_{i=1}^{d_n}(G_n'')^i$, with $(G_n'')^i=\Z$ for $i\leq
1+\mathrm{rank}(K_n),$ and $(G_n'')^i=\Z\oplus {{S_i,}}$ {{where $S_i$ is a cyclic group,}} for
$1+\mathrm{rank}(K_n)<i\leq d_n$, and {{${{\bigoplus}}_{i=2+{\rm rank}(K_n)} {{S}}_i=
\Inf_n$}}. Here the positive cone of $(G_n'')^i$ is given by the strict
positivity of {{the}} first coordinate for non-zero positive elements.  {{An}}
element in  $G_n''$ is positive if
each of its components in $(G_n'')^i$
is positive.}}
{\blue{Define a}}  unital \CA\,  $B_n \in {\bf H}$ such that (see \ref{1323proj})
\beq\label{0.43a}
(K_0(B_n), K_0(B_n)_+, \e_{B_n}, K_n)=(G_n'',{{(G_n'')_+}},u_n'', K_n).
\eneq
More precisely, we have that $B_n=\bigoplus_{i=1}^{d_n}B_n^i$,
with $K_0(B_n^i)=(G_n'')^i,$ and $K_1(B_n^i)$ is either  a cyclic group
for the case $2\leq i\leq 1+{\rm rank} (K_n)$ or zero for the other cases.
In particular, the algebra $B_n^1$ can be chosen to be a matrix
algebra over $\C,$ by the choice of $d_n.$   {{We may also}} assume that, for at least one block $B_n^2$,
the spectrum is not a single point (note that $d_n\geq2$);
otherwise, we will replace the single point spectrum by  the interval
$[0,1]$.  {\blue{We may also write $B_n^i=P_{X,n,i}M_{m_i}(X_{n,i})P_{X,n,i}$ as in \ref{AHblock}
where  $P_{X,n,i}$ has rank at least {{13}} and $X_{n,i}$ is connected.}}
{\blue{For each block $B_n^i$, choose a base point $x_{n,i}\in Sp(B_n^i)$.}}
{\blue{Denote by $\pi_{x_{n,i}}: B_n^i\to B_n^i(x_{n,i}):=F_{X,n}^i$
($\cong M_{{\rm rank} P_{X,n,i}(x_{n,i})}$)
the point evaluation at the point $\pi_{x_{n,i}}.$
Let $F_{X,n}=\bigoplus_{i=1}^{d_n} F_{X,n}^i$ and
let $\pi_{pts,n}: B_n\to F_{X,n}$ be
the quotient map. Put $I_{B,n}={\rm ker}\,\pi_{pts,n}.$
 Then  $(K_0(F_{X,n}), K_0(F_{X,n})_+, [1_{F_{X,n}}])=(\Z^{d_n}, \Z^{d_n}_+, [\pi_{pts,n}(1_{B_n}]).$
Note that $[\pi_{pts,n}(1_{B_n})]=u_n^d.$}}

{{Let us give some obvious properties of the homomorphism $\td\af_{n, n+1}: G_n''\oplus G_n'\to G_{n+1}''\oplus G_{n+1}'$ as below. Let $\af_{n,n+1}^{G_n'', G_{n+1}''}: G_n''\to G_{n+1}''$, $\af_{n,n+1}^{G_n'', G_{n+1}'}: G_n''\to G_{n+1}'$, $\af_{n,n+1}^{G_n', G_{n+1}''}: G_n'\to G_{n+1}''$ and $\af_{n,n+1}^{G_n', G_{n+1}'}: G_n'\to G_{n+1}'$ be the corresponding  partial maps.}}
Define  $\af_{n,n+1}^{{\bar G}_n, G_{n+1}'}:  {\bar G}_n\to G_{n+1}'$  by
$\af_{n,n+1}(g'',g')=\af_{n,n+1}^{G_n'', G_{n+1}'}(g'')+\af_{n,n+1}^{G_n', G_{n+1}'}(g').$
{{Since $\Inf_{n+1}\cap G_{n+1}'=0$ and $\td\af_{n,n+1}$ maps $\Inf_n$ to $\Inf_{n+1}$, we know $\af_{n,n+1}^{G_n'', G_{n+1}'}(\Inf_n)=0$. Hence $\af_{n,n+1}^{G_n'', G_{n+1}'}$  {{ factors}} through $G_n''/\Inf_n$ as
$$ \af_{n,n+1}^{G_n'', G_{n+1}'}: G_n'' \stackrel{\pi}{\longrightarrow} G_n''/\Inf_n\stackrel {\td\af''_{n,n+1}}{\longrightarrow} G_{n+1}'.$$}}
{\blue{In other words, we may write
$\af_{n,n+1}^{G_n'', G_{n+1}'}={\td\af''_{n,n+1}}\circ \pi, $ where  $\pi: G_n''\to G_n''/\Inf_n$
is the quotient map, and ${\td\af_{n,n+1}}'': G_n''/\Inf_n \to G_{n+1}'$ is the induced \hm.}}
Also, {\blue{since $\td\af_{n,n+1}={\td\gm_{n,n+1}}|_{G_n},$}}  $\af_{n,n+1}^{G_n', G_{n+1}''}$ {\blue{factors}} through $H_n'$ as
$$ \af_{n,n+1}^{G_n', G_{n+1}''}: G_n' \stackrel{\blue{{\iota_n}|_{_{G_n'}}}}
{\longrightarrow} H_n'\stackrel {\gm_{n,n+1}^{H_n', G_{n+1}''}}{\longrightarrow} G_{n+1}''.$$
{{That is,}} $\af_{n,n+1}^{G_n', G_{n+1}''}=\gm_{n,n+1}^{H_n', G_{n+1}''} \circ  {\blue{{\iota_n}|_{_{G_n'}}}},$
{\blue{where $\gm_{n,n+1}^{H_n', G_{n+1}''}: H_n'\to G_{n+1}''$ is the \hm\, induced by $\gm'_{n,n+1}.$
Note {{that}}
$\td\gm_{n,n+1}^{{{H}}_n', {{H}}_{n+1}'}$ has multiplicity at least $13 \cdot 2^{2n}.$
}}
\end{NN}

\begin{NN}\label{range 0.43}
We can extend the maps $\bt_{{n,0}}$ and $\bt_{{n,1}}$ to $\bt_{{n,0}},~\bt_{{n,1}}:
B_n\oplus F_n \to E_n$, by defining them to be zero on $B_n$.
Consider $A_n=B_n\oplus C_n$. Then the C*-algebra $A_n$ can be written as
$$A_n=\{(f,a)\in C([0,1],E_n)\oplus (B_n\oplus F_n): f(0)=\bt_{{n, 0}}(a), f(1)=\bt_{{n, 1}}(a)\}.$$
Then
$$
(K_0(A), K_0(A_n)_+, [1_{A_n}], K_1(A))=({\bar G}_n, ({\bar G}_n)_+, {\bar u}_n, K_n).
$$
$\\${\blue{Let $C_n'=F_{X,n}\oplus C_n$ (after $C_n$ is defined) and let  $\pi_{C_n',C_n}: C_n'\to C_n$ be the quotient map.
Note that $C_n'=C([0,1], E_n)\oplus_{\bt_{n,0}, \bt_{n,1}} F_{X,n}\oplus F_n,$
where $\bt_{n,0}$ and $\bt_{n,1}$ are extended to maps from $F_{X, n}\oplus F_n$
which are zero on $F_{X,n}.$
}}
{\blue{Let $I_n=C_0\big((0,1),E_n\big)$ be the ideal of $C_n$  (also an ideal of $A_n$).
Denote by $\pi_{I,n}: A_n\to B_n\oplus F_n$ and
$\pi_{_{A_n, C_n}}: A_n\to C_n$ the quotient maps. We also use $\pi_{_{A_n, B_n}}: A_n\to B_n$ for the quotient map.}}

{\blue{We will construct $B_m,$ $C_m$, $F_m$, and $A_m$
{{(in order for comparison
with \ref{condition2}, we use subscripts $m$ instead of $n$)}} as above together
with  a unital injective \hm\, $\phi_{m,m+1}: A_{{m}}\to A_{m+1},$ an injective \hm\, ${\bar \phi}_{m,m+1}': C_m'\to C_{m+1}',$
and  an injective \hm\, $\psi_{B, m,m+1}: B_m\to B_{m+1}$
satisfying the following conditions (similar to those of \ref{condition2})}}

{\blue{(a) $\phi_{m,m+1}(I_{{m}})\subset I_{m+1},$
$\phi_{m,m+1}(I_{B,m})\subset I_{B, m+1}$ {{(see \ref{range 0.41} for  $I_{B, m+1}$)}}. {{(Hence $\phi_{m,m+1}$ induces a homomorphism  $\phi^q_{m,m+1}: F_{X,m}\oplus F_m=A_m/(I_m\oplus I_{B,m})\to F_{X,m+1}\oplus F_{m+1}=A_{m+1}/(I_{m+1}\oplus I_{B,m+1})$.)}} {{ Furthermore }}
${\pi_{C_{m+1}', C_{m+1}}\circ {\bar \phi}_{{{m}},{{m}}+1}'}|_{C_{{m}}}=
{{\phi_{m,m+1}^{C_m,C_{m+1}},}}$ {{where $\phi_{m,m+1}^{C_m,C_{m+1}}:=\pi_{A_{m+1}, C_{m+1}}\circ \phi_{m,m+1}|_{C_m}$ is the partial map of $\phi_{m,m+1}$ from $C_m$ to $C_{m+1}$,}}
and $(\pi_{A_{m+1}, B_{m+1}}\circ \phi_{m,m+1})|_{B_m}=\psi_{B, {{m}},{{m}}+1};$\\

(b) $(\phi_{m,m+1})_{*0}={\td \af}_{m,m+1},$
{{$(\phi_{m,m+1})_{*1}=\chi_{m,m+1},$}} $({\bar \phi}_{m,m+1}')_{*0}=\td\af_{m,m+1}'$ and\\
 $(\pi_{C_{m+1}', C_{m+1}}\circ ({\bar \phi}_{m,m+1}')|_{C_{{m}}})_{*0}={{\af_{m,m+1}^{G_m', G_{m+1}'}}}$
 {{and}} $(\psi_{B, m,m+1})_{*0}=\af_{m,m+1}^{G_m'', G_{m+1}''};$\\

(c) {{({{compare}} (c) of \ref{condition2}),}} the map ${\bar \phi}_{m,m+1}'$ satisfies the {{conditions}} (1)--(8) {{of}} \ref{range 0.6}
(where  $C_m$ and $C_{m+1}$ are  replaced by $C_m'$ and $C_{m+1}',$ and
$H_m$ and $G_m$ are replaced by ${{\Z^{d_m}\oplus}}H_m'$ and ${{\Z^{d_m}\oplus}}G_m'$){{,
and}} satisfies  condition ${{\spdd_1}}$ in (c) of {{\ref{condition2} (with the number $L_m$ {{as}} described below). In particular,}} $${{(\phi^q_{m,m+1})_{*0}=(c_{ij}^{m,m+1}):K_0(F_{X,m}\oplus F_m)=\Z^{d_m}\oplus H'_m \to K_0(F_{X,m+1}\oplus F_{m+1})=\Z^{d_{m+1}}\oplus H'_{m+1};}}$$

(d) the matrices $\bt_{m+1,0}$ and $\bt_{m+1, 1}$ for $C_{m+1}'$ satisfy
the condition  $\spdd$ (where $(d_{ij}^{m,m+1}): (\Z^{d_m}\oplus H_m')/(\Z^{d_m}\oplus G_m')=H_m'/G_m'\to
H_{m+1}'/G_{m+1}'= (\Z^{d_{m+1}}\oplus H_{m+1}')/(\Z^{d_{m+1}}\oplus G_{m+1}')$ is the  map induced by $\td\gm_{m,m+1}'$).}}\\

{\blue{The number $L_m$ in ${{\spdd_1}}$ in (c) {{of \ref{condition2} (see (c) above)}} which  was to be chosen
at
the $m$-th step is described {{as follows}}: {{Set $\psi_k:=\psi_{B, m-1, m}\circ \psi_{B, m-2, m-1} \circ\cdots\circ \psi_{B, k, k+1}: B_k\to B_m.$ (Recall that $\psi_{B, i, i+1}=\pi_{A_{i+1},B_{i+1}}\circ \phi_{i,i+1}|_{B_i}$ is the partial map from $B_i\subset A_i$ to $B_{i+1}\subset A_{i+1}$ of the homomorphism $\phi_{i, i+1}$.}}
 There is a finite subset $Y_m\in Sp(B_m)$ such that, for each $k\leq {{m}}-1$,~ the set  ${{\bigcup}}_{x\in {{Y}}_m} Sp(\psi_k)_x$ is {{$1/m$}}-dense in $Sp(B_k)$, and $Y_m$ is also ${{1/m}}$-dense $Sp(B_{{m}})$. Let
$T_{{m}}:=\big\{ \big(\frac{k}{{m}}\big)_{n,j}\in Sp(E_{{m}}^j), ~~1\leq k\leq {{m}}-1, ~1\leq j\leq l_{{m}}\big\},$ and write $\Omega_m=Y_m\cup T_m$.}}
Choose
\beq\label{range 0.43-3}
L_{{m}} > 13\cdot
2^{2{{m}}}\cdot(\#(\Omega))\cdot(\max\{{{\mathrm{rank}
(\one_{{\blue{B^{l}_{{m}}}}}), \mathrm{rank}(\one_{F^i_{{m}}}), \mathrm{rank}(\one_{E^{\blue{j}}_{{m}}})}}\}),
\eneq
{{where by convention, if $C=M_k$, then ${\rm rank} (\one_C)=k$; and if $C=PM_{l}(C(X))P$, then ${\rm rank} (\one_C)={\rm rank}(P)$.}}
\end{NN}

\begin{NN}\label{14construct}
{\blue{~ We begin to construct $C_1$ and $F_1$ {{in}} exactly {{the}} same {{way}} as in {{Section 13}} (for
$G_1'$ and $H_1'$).
 Suppose that $C_1,C_2,...,C_n,$ $F_1,F_2,..., F_n,$ $B_1,B_2,...,B_n${{, $\phi_{k, k+1},$ ${\bar \phi}_{k,k+1}'$ and $\psi_{B, k, k+1}$ (for $k\leq n-1$) }}
have been constructed.
 We will choose $B_{n+1}$ as above.
 Since each $\phi_{B, k,k+1}$ is injective ($1\le k<n$){{, so also}} is $\psi_k$ above. Therefore
 $Y_n$ exists and $L_n$ can be defined.
 Since each {{map $\td \gm_{n,n+1}': \Z^{d_n}\oplus H_n' \to \Z^{d_{n+1}}\oplus H_{n+1}' $}} is strictly positive, exactly as in \ref{condition2},  (passing to larger $n+1$), we may assume that ${{\spdd_1}}$ holds.
 This determines the integer $n+1$ (which is originally denoted by some large integer greater than $n$).
 Then, as in  \ref{condition2}, using $\Lambda_n$ {{(defined in (\ref{Lambdanbt}))}}, we can choose $\bt_{n{{+1}},0}$ and $\bt_{n{{+1}},1}$ so that
 $\spdd$ holds.  Then we use these $\bt_{n{{+1}},0}$ and $\bt_{n{{+1}},1}$ (which are zero on $F_{X,n+1}$) to define $C_{n+1}'.$
 So  $C_{n+1}$ is also defined. Let $A_{n+1}=C_{n+1}\oplus B_{n+1}.$  Note that (d) has been verified.
 }}

{\blue{Let  ${\blue{P_n, Q_n}}\in A_{n+1}$ be projections such that $[P_n]=\td\af_{n,n+1}([\one_{C_n}])$, $[Q_n]=\td\af_{n,n+1}([\one_{B_n}]),$ and $P_n+Q_n=\one_{A_{n+1}}$. Let $P_{n,C}=\pi_{_{A_{n+1}, C_{n+1}}}(P_n),$
$P_{n,B}=\pi_{_{A_{n+1}, B_{n+1}}}(P_n),$  $Q_{n,C}=\pi_{_{A_{n+1}, C_{n+1}}}(Q_n),$ and $Q_{n,B}=\pi_{_{A_{n+1}, B_{n+1}}}(Q_n).$
Note that $\td\af_{n,n+1}'(u_n^d, u_n')=(u_{n+1}^d, u_{n+1}').$}}

{\blue{It follows by  Lemma \ref{range 0.6} that there is a unital injective \hm\, ${\bar \phi}_{n,n+1}': C_n'\to C_{n+1}'$
which satisfies the {{conditions (1)--(8)}} of \ref{range 0.6}.
In particular, $({\bar \phi}_{n,n+1}')_{*0}=\td\af_{n,n+1}'.$
Define a \hm\, $\phi_{n,n+1}^{A_n, C_{n+1}}: A_n\to C_{n+1}$ by
$\phi_{n,n+1}^{A_n, C_{n+1}}=\pi_{\small{C_{n+1}', C_{n+1}}}\circ {\bar \phi}_{n,n+1}'\circ \pi_{_{A_n, C_n'}}.$}}
{\blue{One checks that $(\phi_{n, n+1}^{A_n, C_{n+1}})_{*0}=\af_{n, n+1}^{{\bar G}_n, G_{n+1}'},$ and
$(\phi_{n, n+1}^{A_n, C_{n+1}})_{*1}=0.$   Moreover, $(\phi_{n,n+1}^{A_n, C_{n+1}})|_{I_{B,n}}=0,$ and \\
$(\pi_{C_{n+1}', C_{n+1}}\circ \phi_{n,n+1}^{A_n, C_{n+1}}|_{C_n})_{*0}=\af_{n,n+1}^{G_n',G_{n+1}'}.$}}
{\blue{Since $\td{{\gm}}_{n,n+1}'$ is $13\cdot 2^{2n}$-large, so {{also}} is
$\af_{n,n+1}^{G_n'',G_{n+1}''}.$ It follows from the second part of \ref{range 0.14} that there is an injective \hm\,
$\psi_{B,n, n+1}: B_n\to Q_{n,B}B_{n+1}Q_{n,B}$ which maps $I_{B,n}$ into $I_{B,n+1}$
such that $({{\psi_{B,n,n+1}}})_{*0}=\af_{n,n+1}^{G_n'', G_{n+1}''}$ and
$({{\psi_{B,n,n+1}}})_{*1}=\chi_{n,n+1}.$
$(\phi_{n,n+1}^{B_n, B_{n+1}})_{*1}=\chi_{n,n+1}.$
 Since  $\af_{n,n+1}^{H_n', G_{n+1}''}$ is  at least $13$-large,  by  the second part of \ref{range 0.14} (with each $X_i$  a point) again,
there is a unital injective \hm\, ${\bar \phi}_{n,n+1}^{F_n, B_{n+1}}: F_n\to  P_{n,B} B_{n+1}P_{n,B}$
such that $({\bar \phi}_{n,n+1}^{F_n, B_{n+1}})_{*0}=\af_{n, n+1}^{H_n', G_{n+1}''}.$
Define ${\bar \phi}_{n,n+1}^{C_n, B_{n+1}}={{({\bar \phi}_{n,n+1}^{F_n, B_{n+1}})\circ \pi_{I,n}}}|_{C_n}.$
Therefore $({\bar \phi}_{n,n+1}^{C_n, B_{n+1}})_{*0}=\af_{n,n+1}^{G_n', G_{n+1}''}.$
Then define $\phi_{n,n+1}: A_n\to A_{n+1}$
by, for all $(b,c)\in B_n\oplus C_n,$
\beq
\phi_{n,n+1}(b,c)= \psi_{B, n,n+1}(b)\oplus {\bar \phi}_{n,n+1}^{C_n, B_{n+1}}(c)\oplus \phi_{n,n+1}^{A_n,C_{n+1}}(b,c).
\eneq
One checks that $\phi_{n,n+1}$ ({{together with}} the induced maps) also satisfies (a), (b){{, and}} (c).
Let $A=\lim_{n\to\infty}(A_n, \phi_{n,n+1}).$}}
\end{NN}

\begin{NN}\label{14Atrace}
{\blue{{{For the algebra $A$ constructed above, we have}}
$$
(K_0({{A}}), K_0(A)_+, [1_A], K_1(A))=(G, G_+, u, K_1).
$$}}
By  (2) of \ref{range 0.40},
we also have
\beq\label{AlineC0}
\frac{\tau(\alpha_{n,\infty}(\e_{B_n}))}{\tau(\alpha_{n,\infty}(\e_{C_n}))}<\frac{1}{13\cdot 2^n-1},
\eneq
for all $\tau\in T(A).$
{\blue{Let $F_n'=F_{X,n}\oplus F_n,$ and define $\pi_{_{A_n, F_n'}}: A_n\to F_n'$ by
$\pi_{A_n, F_n'}(b, c)=(\pi_{pts,n}(b), c)$ for $(b,c)\in B_n\oplus F_n.$
By (a)  {{of}} \ref{14construct}, $\phi_{n,n+1}(I_n+I_{B,n})\subset I_{n+1}+{{I_{B,n+1}}}.$ 
{{Therefore,}} $\phi_{n,n+1}$ induces a unital \hm\, $\phi_{_{F',n,n+1}}: F_n'\to F_{n+1}'.$ {{(The map $\phi_{_{F',n,n+1}}$ is denoted by $\phi^q_{n,n+1}$ in \ref{range 0.43} {{to be consistence with section 13, see 13.31}}.)}}
Note $K_0(F_n')=\Z^{d_n}\oplus G_n'$ and $(\phi_{F',n,n+1})_{*0}=\td\gm_{n,n+1}'.$}}
{\blue{Define $F'=\lim_{n\to\infty}(F_n', \phi_{F,n,n+1}').$  Then $F'$ is a unital AF-algebra. By \eqref{14afK0},
$(K_0(F'), K_0(F')_+, [1_F'])=(H', H'_+, u').$ It follows that $F'\cong F$ by
Elliott's classification. Therefore $T(F')=\Delta.$
We also have the following commutative diagram:
\begin{displaymath}
\xymatrix{
A_1\ar[r]^{\phi_{1, 2}}\ar[d]_{\pi_{_{A_1, F_1'}}}& A_2 \ar[d]_{\pi_{_{A_2,F_2'}}} \ar[r]^-{\phi_{2, 3}} & \cdots \ar[r] & A \ar[d]
\\
F'_1\ar[r]_{\phi_{_{F',1,2}}} & F'_2
\ar[r]_{\phi_{_{F',2, 3}}} & \cdots \ar[r] & {{F'.}}}
\end{displaymath}}}
As in {{Section 13, using}} ${{\spdd_1}},$  \eqref{14Atrace}, and the above diagram,
one shows that
$$
(K_0(A), K_0(A)_+, [1_A], K_1(A), T(A), r_A)=(G, G_+, u, K, \Delta, r).
$$
Note that $A$ is not simple. So we make one {{more}} modification just as in {{Section}} 13 (but much simpler as $B_n$
is a direct summad of $A_n$). We describe this briefly as follows:
Let  $A_{n+1}^{-}={{\bigoplus}}_{i=2}^{d_{n+1}}B_n^i{{\oplus}} C_{n+1},$
let $\pi_{B_{n+1}^1}: A_{n+1}\to {{B_{n+1}^1,}}$ 
and    $\pi_{B_{n+1}}^-: A_{n+1}\to {{\bigoplus}}_{i\ge 2}B_{n+1}^i$
be the quotient maps. We also view $\pi_{B_{n+1}^1}$ as the quotient map from $B_{n+1}$ to
$B_{n+1}^1.$
Recall that $B_{n+1}^1\cong M_{R(n+1)}{{}}$ for integer $R(n+1)\ge 1.$
We may write
$\psi_{B,n, n+1}=\psi_{B,n,n+1}^-\oplus \psi_{B,n, n+1}^1,$ where
$
\psi_{B,n,n+1}^-=\pi_{B_{n+1}}^-\circ \psi_{B,n,n+1}\andeqn
\psi_{B,n,n+1}^1=\pi_{B_{n+1}^1}\circ \psi_{B,n,n+1}.
$
One may write, using {{the}} notation {{of}} Section 13, and using  ${{\spdd_1}},$
\beq
\psi_{B,n, n+1}^1(f)=(\pi_{x_{{n,1}}}^{\sim a_1}, \pi_{x_{{n,2}}}^{\sim a_2},...,\pi_{x_{{n,d_n}}}^{\sim a_{d_n}})(f)\rforal f\in B_n,
\eneq
{{where}}  $a_i{{>}} \#(Y_n)$ ($1\le i\le d_n$).   Let $q_{n+1}=\psi_{B,n,n+1}^1(1_{B_n}).$
There is a continuous path $\psi_{B,n,n+1}^{1,t}: B_n\to q_{n+1}B_{n+1}q_{n+1}$
($t\in [0,1]$)
such that $\psi_{B,n,n+1}^{1,0}={{\psi}}_{B,n,n+1}^1,$ and \\ {{$\psi_{B,n,n+1}^{1,1}:=\xi^{B_n,B_{n+1}^1}$ satisfies}}
\beq\label{13April2}
Sp({{\xi^{B_n,B_{n+1}^1}}})\supset Y_n{{\cup Sp(F_{X,n})}}.
\eneq
Define  ${{\xi^{B_n,B_{n+1}}}}=\psi_{B,n,n+1}^-\oplus {{\xi^{B_n,B_{n+1}^1}}}.$
  {{Then}} $Sp({{\xi^{B_n,B_{n+1}}}})\supset Y_n{{\cup Sp(F_{X,n})}}.$
  $q_{n+1}^\sim={{1_{B_{n+1}}^1}}-q_{n+1}={{\pi_{B_{n+1}^1}\circ {\bar \phi}^{C_n, B_{n+1}}(1_{C_n})}}$. 
 Consider {{the map}} $\pi_{B_{n+1}^1}\circ {\bar \phi}^{C_n, B_{n+1}}: C_n\to B_{n+1}^1.$
Exactly  as in \ref{range 0.31},  since ${{\spdd_1}}$ holds,  one has a continuous path of unital \hm s
$\Omega_{I,s}: C_n\to q_{n+1}^\sim B_{n+1}q_{n+1}^\sim$ {{with}}
$\Omega_{I,0}=\pi_{B_{n+1}^1}\circ {\bar \phi}^{C_n, B_{n+1}}$ and $\Omega_{I,1}:={{\xi}}^{C_n, B_{n+1}^1}$
such that
\beq\label{13April2-1}Sp({{\xi}}^{C_n, B_{n+1}^1})\supset  T_n\cup Sp(F_n).\eneq
Put ${\bar \phi}^{C_n, B_{n+1}^-}= \pi_{B_{n+1}}^-\circ {\bar \phi}^{C_n, B_{n+1}}.$ Define, for all $(b,c)\in A_n=B_n\oplus C_n,$
\beq
{{\xi_{n,n+1}(b,c)}}={{\xi^{B_n,B_{n+1}}}}(b)\oplus {{\xi}}^{C_n, B_{n+1}^1}(c)
\oplus {\bar \phi}^{C_n, B_{n+1}^-}(c)\oplus \phi_{n,n+1}^{A_n,C_{n+1}}(b,c).
\eneq
Then, {{by (\ref{13April2}) and (\ref{13April2-1}), we have}} $Sp({{\xi|_{x_{n+1,1}}}})\supset \Omega_n\cup Sp(F_n)${{. Hence,}}
 by (c) {{of}} \ref{range 0.43} ({{see}} (5)  and (8) of \ref{range 0.6}),
 $Sp(\xi_{m,n+2}|_x)$ is $1/n$-dense  in $Sp(A_{m})$ for all $m\le n.$
 Note that $\xi_{n,n+1}$ and $\phi_{n,n+1}$ are {{ homotopic}}.
 It follows that $(\xi_{n,n+1})_{*i}=(\phi_{n,n+1})_{*i},$ $i=0,1.$
 Let $B=\lim_{n\to\infty}(A_n, \xi_{n,n+1}).$
 Then $(K_0(B), K_0(B)_+, [1_B], K_1(B))=(K_0(A), K_0(A)_+, [1_B], K).$
 It follows from \ref{simplelimit} that  $B$ is a unital simple \CA.
 By  \eqref{AlineC0}, $\lim_{n\to\infty}\|\xi_{n,n+1}^{\sharp}-\phi_{n, n+1}^{\sharp}\|=0.$
Therefore, as in \ref{range 0.33} (but much {{more simply}}),
\beq
(K_0(B), K_0(B)_+, [1_B], K_1(B), T(B), r_B)=(G, G_+, u, K, \Delta, r).
\eneq
\end{NN}

{\it At this point, we would like to point out that,
if ${\rm ker}\rho_A=\{0\}$ and $K_1(A)=\{0\},$ then we do not need
the {{direct}} summand $B_n$ in the construction above.}

{\blue{We may summarize the {{ construction above}} as follows:}}

\begin{thm}\label{RangT}
Let $((G, G_+, u), K, \DT, r)$ be  a {{sextuple consisting}} of the following
objects: $(G, G_+, u)$ is a  {{countable}} weakly unperforated simple order-unit {{abelian}}
group, $K$ is a countable abelian group, $\DT$ is a
{{metrizable}}
Choquet simplex{{, and}} $r:  \DT\to \mathrm{S}_u(G) $ is a surjective
{{continuous}} affine map, where $\textrm{S}_u(G)$ {{is}} the compact convex set {{of}}
states {\blue{of}} $(G, G_+, u)$. Assume that $(G, G_+, u)$ has the
{{property (SP)}} in the sense (as above) that for any real number $s>0$, there is $g\in
G_+\setminus\{0\}$ such that $\tau(g)<s$ for any state $\tau$ on $G$.

Then there is a unital  simple C*-algebra $A\in {\cal B}_0$ which can be written
as $A={{\varinjlim}}(A_i, \psi_{i, i+1})$ with injective
$\psi_{i, i+1}$, where $A_i=B_i\oplus C_i,$ $B_i\in {\bf H},$ {{and}}
$C_i\in {\cal C}_0${{, in such a way}} 
that
\begin{enumerate}
\item $\lim_{i\to\infty}\sup\{\tau(\psi_{i, \infty}(1_{B_{{{i}}}})): \tau\in T(A)\}=0,$
\item $\ker \rho_A\subset \bigcup_{i=1}^{\infty} ({\psi}_{i, \infty})_{*0}({\rm ker}\rho_{B_i})$, and
\item ${\rm Ell}(A)\cong ((G, G_+, u), K, \DT, r).$
\end{enumerate}

Moreover, the inductive system $(A_i, \psi_i)$ can be chosen so that $\psi_{i, i+1}=\psi_{i, i+1}^{(0)}\oplus\psi_{i, i+1}^{(1)}$ with $\psi_{i, i+1}^{(0)}: A_i\to A_{i+1}^{(0)},$ {{$\psi_{i,i+1}^{(0)}$ is non-zero
on each summand of $A_i,$}} and $\psi_{i, i+1}^{(1)}: A_i\to A_{i+1}^{(1)}$ for C*-subalgebras $A_{i+1}^{(0)}$ and $A_{i+1}^{(1)}$ of $A_{i+1}$ with $1_{A_{i+1}^{(0)}}+1_{A_{i+1}^{(0)}}=1_{A_{i+1}}$ 
such that
$A_{i+1}^{(0)}$ is a non-zero finite dimensional C*-algebra,  
and $(\psi_{i, \infty})_{*1}$ is injective.

Furthermore, if $K=\{0\}$ and {{$\Inf(G)$=0, which implies that}} $G$ is torsion free, we can assume that $A_i=C_i,$
$i=1,2,....$

\end{thm}

\begin{proof} Condition (2) follows from  the fact that $\Inf(G) =\bigcup_{i=1}^{\infty} \td\af_{n, \infty}(\Inf_n)$
and (\ref{0.43a}), {{and (1) follows from (\ref{AlineC0})}}. {{Since}} $A$ is a unital simple inductive limit of subhomogeneous C*-algebras with no dimension growth,
it then follows from Corollary 6.5 of \cite{Winter-Z-stable-01} that $A$ is ${\cal Z}$-stable.
Hence $A$ has strict comparison for positive elements {{(and projections)}}.
By (1) above{{, the strict comparison property {{ mentioned}} above,}} and the fact that $C_i \in {\cal C}_0,$ {{we conclude}} that
$A\in {\cal B}_0$.  (This actually {{also}} follows from our construction immediately.)
\end{proof}



\begin{rem}\label{RRrangT}
Note that   $A_{i+1}^{(0)}$ can be chosen to be the first block
$B_{i+1}^1$, so we have
$$\lim_{i\to\infty}\tau(\psi_{i+1, \infty}(1_{A_{i+1}^{(0)}}))=0$$
uniformly for $\tau\in T(A)$. {{Moreover, for any $i,$ there exists $n\ge i$
such that
\beq\label{RRrangT-1}
\lim_{n\to\infty}\sup_{\tau\in T(A)}\tau(\psi_{n+1,\infty}\circ \psi_{n,n+1}^{(0)}\circ \psi_{i, n}(1_{A_i}))=0.
\eneq
}}
\end{rem}

\begin{rem}\label{range 0.46}
Let ${{\psi_{n,n+1}^B:=\psi_{n,n+1}^{B_n,B_{n+1}}:~B_n \to
B_{n+1}}}$ be the partial map of $\psi_{n,n+1}:{{A_n \to
A_{n+1}}}$, and let ${\blue{\psi_{n,m}^B}}: B_n \to B_m$ be the corresponding composition
${\blue{\xi_{m-1,m}^B\circ \xi_{m-2,m-1}^B\circ \cd \circ \xi_{n,n+1}^B}}$. Let
$e_n=\xi_{1,n}(\e_{B_1}).$ Then, from the construction, we know that
the algebra $B=\varinjlim(e_nB_n e_n, \psi_{n,m}^B)$ is simple, as we know
that $SP(\xi_{n,n+1}^B|_{{\blue{x_{n+1,1}}}})$ is dense enough in $Sp(B_n)$. Note
that the simplicity of $B$ does not follow from {{the}} simplicity of $A$
itself, since it is not a corner of $A$.

\end{rem}

\begin{cor}\label{CCRangT}
{{Let $B_1\in {\cal B}_0$ be a unital separable \CA\, and {{set}} $B=B_1\otimes U,$ where $U$ is a UHF-algebra of infinite type {{({{i.e.,}} $U\otimes U=U$)}}.
Then there exists a \CA\, $A$ with
all {{the}} properties  described in {{Theorem}} \ref{RangT} such
that ${\rm Ell}(A)={\rm Ell}(B).$ Moreover, $A$ may be chosen such that $A\otimes U\cong A$ and \ref{RRrangT} also is valid for $A.$}}
\end{cor}
{\begin{proof}
{{Let $((G, G_+, u), K, \DT, r)=(K_0(B), K_0(B)_+, [1_B], K_1(B), T(B), r_B).$
Since $A=A'\otimes U,$ it has the {{ property (SP)}}.
Then, by {{Theorem}} \ref{RangT}, there is a \CA\, $A'\in {\cal B}_0$ which {{has}} all the properties
of $A$ described in \ref{RangT} such that ${\rm Ell}(A')={\rm Ell}(B).$
Put $A=A'\otimes U.$
It is easy to check  ${\rm Ell}(A'\otimes U)={\rm Ell}(B\otimes U).$
Since $B\otimes U\cong B,$ we conclude that ${\rm Ell}(A)={\rm Ell}(B).$
Write $U=\lim_{n\to\infty}(M_{k(n)}, \imath_n),$ where $k(n+1)=r(n)k(n)$ and
$\imath_n: M_{k(n)}\to M_{k(n+1)}$ is defined by $\imath_n(a)=a\otimes 1_{M_{r(n)}}.$
Write $A'=\lim_{n\to\infty}(B_n{{\oplus}}C_n, \psi_{n,n+1}).$ Then one checks
$A\cong \lim_{n\to\infty} (B_n\otimes M_{k(n)}\oplus C_n\otimes M_{k(n)}, \psi_{n,n+1}\otimes \imath_n).$
It follows that $A$ has all {{the}} properties described in \ref{RangT} {{and also}} the one {{described}} in \ref{RRrangT}.
Since $A=A'\otimes U$ and $U$ is a UHF-algebra of infinite type, $A\otimes U\cong A.$}}
\end{proof}

The following result is an analog of Theorem 1.5 of \cite{Lnann}.

\begin{cor}\label{smallmap}
Let $A_1$ be a simple separable C*-algebra in $\mathcal B_1$,
 and let $A=A_1\otimes U$ for {an infinite dimensional}  UHF-algebra $U$. There exists an
  inductive limit algebra $B$ as constructed in Theorem \ref{RangT} such that $A$
  and $B$ have the same  Elliott invariant{.}
   Moreover, the C*-algebra $B$ {{may be chosen to have}} the following properties:

Let $G_0$ be a finitely generated subgroup of $K_0(B)$ with decomposition $G_0=G_{00}\oplus G_{01}$, where $G_{00}$  vanishes under all states of $K_0(A)$. Suppose $\mathcal P\subset \underline{K}(B)$ is a finite subset which generates a subgroup $G$ such that $G_0\subset G\cap K_0(B)$.

Then, for any $\epsilon>0$, any finite subset $\mathcal F\subset B$, any $1>r>0$, and any positive integer $K$, there is an
$\mathcal F$-$\epsilon$-multiplicative map $L:B\to B$ such that:
\begin{enumerate}
\item $[L]|_{\mathcal P}$ is well defined.
\item $[L]$ induces the identity maps on the infinitesimal part of  $G\cap K_0(B)$, $G\cap K_1(B)$,
      $G\cap K_0(B,\mathbb Z/k\mathbb Z)$ and $G\cap K_1(B, \mathbb Z/k\mathbb Z)$
      for $k=1, 2,{{\cdots.}}$
\item $\rho_B\circ[L](g)\leq r\rho_B(g)$ for all $g\in
      G\cap K_0(B)$, where $\rho_{{B}}$ is the canonical positive homomorphism from $K_0({{B}})$ to $\Aff(\mathrm{S}(K_0({{B}}), K_0({{B}})_+, [1_{{B}}]))$.
\item For any positive element $g\in G_{01}$, {{there exists $f\in K_0(B)_+$ with $g-[L](g)=Kf$.}}
\end{enumerate}
\end{cor}
\begin{proof}
{{As is pointed  out in the last paragraph of \ref{range 0.37}, the Elliott invariant of $A_1\otimes U$ has {{the}} property (SP).}} {{Replacing}} $A_1$ by $A_1\otimes U,$ we may assume that ${\rm Ell}(A_1)$ has  {{property (SP)}}.

 Consider ${\rm Ell}(A_1)$,
which satisfies the conditions of Theorem \ref{RangT}.
Therefore, by the first part of Theorem \ref{RangT}, there is {an} inductive
system $B_1=\varinjlim(T_i\oplus S_i, \psi_{i, i+1})$ such that

(i) $T_i\in\mathbf H$ and $S_i\in\mathcal C_0$
with $K_1(S_i)=\{0\}$,

(ii) $\lim \tau(\phi_{i, \infty}(1_{T_i}))= 0$ uniformly on $\tau\in T(B_1)$,

(iii) $\ker(\rho_{B_1})=\bigcup_{i=1}^\infty (\psi_{i, \infty})_{*0}(\ker(\rho_{T_i}))$, and

(iv) ${\rm Ell}(B_1)={\rm Ell}(A_1)$.

Put $B=B_1\otimes U$. Then ${\rm Ell}(A)={\rm Ell}(B)$. Let $\mathcal P {{\subset}} \underline{K}(B)$ be a finite subset, and let $G$ be the subgroup generated by {{$\mathcal P$,}} which {{we may assume}} contains $G_0$. Then there is a positive integer $M'$ such that $G \cap K_*(B, \mathbb Z/ k\mathbb Z)={\{0\}}$ if $k>M'$. Put $M=M'!$. Then $Mg=0$ for any $g\in G \cap K_*(B, \mathbb Z/ k\mathbb Z)$, $k=1, 2, ...$ .

Let $\ep>0$, {{a finite subset  $\mathcal F\subset B$,}} and $0< r<1$ be given. Choose a finite subset $\mathcal G\subset B$ and  $0<\ep'<\ep$
such that $\mathcal F\subset \mathcal G$  and for any $\mathcal G$-$\epsilon'$-multiplicative map $L: B\to B$,
the map $[L]_{\mathcal P}$ is well defined, and $[L]$ is
a homomorphism on $G$ (see \ref{KLtriple}).

{{{Choosing}} a sufficiently large $i_0,$ we may assume that
$[\psi_{i_0, \infty}](\underline{K}(T_{i_0}\oplus S_{i_0}))\supset G.$
In particular, we may assume,  {{by (iii) above, that}}
$G\cap {\rm ker}\rho_{B_1}\subset (\psi_{i_0, \infty})_{*0}({\rm ker}\rho_{T_{i_0}}).$ Let $G'\subset \underline{K}(T_{i_0}\oplus S_{i_0})$ be such that $[\psi_{i_0, \infty}]({ G}'){{\supset}}{ G}.$}
{Since $B=B_1\otimes U,$ we may {{write }}
$U=\varinjlim(M_{m(n)}, \imath_{n, n+1}),$
where $m(n)|m(n+1)$ and $\imath_{n, n+1}: M_{m(n)}\to M_{m(n+1)}$ is defined
by $a\mapsto a\otimes 1_{m(n+1)}.$ }
One may assume that for each $f\in\mathcal G$, {there exists $i> i_0$ such that}
\begin{equation}\label{diag-f}
f=(f_0\oplus f_1)\otimes 1_m
\in (T_i'\oplus S_i')\otimes M_m  
\end{equation}
for some $f_0\in T_i'$, ${{f_1}}\in S_i'$, and $m>2MK/r,$ {where
$m=m(i+1)m(i+2)\cdots m(n),$
$T_i'=\psi_{i, \infty}'(T_i\otimes M_{m(i)}),$  $S_i'=\psi_{i, \infty}'(S_i\otimes M_{m(i)}),$ and
where $\psi_{i, \infty}=\psi_{i, \infty}\otimes \imath_{i,\infty}.$}
Moreover, one may assume that $\tau(1_{T_i'})<r/2$ for all $\tau\in T(A_1)$. 

Choose a large $n$ {{such}} that $m={M_0}+l$ with ${M_0}$ divisible by $KM$ and $0\leq l<KM$. Then define the map $$L: (T_i'\oplus S_i')\otimes M_m \to (T_i'\oplus S_i')\otimes M_m$$ to be
$$
{L((f_{i, j} \oplus g_{i, j})_{m\times m})=(f_{i,j})_{m\times m}\oplus E_l(g_{i,j})_{m\times m}E_l,}
$$
{where $E_l={\rm diag}(\underbrace{1_{S_i'}, 1_{S_i'},...,1_{S_i'}}_l, \underbrace{0,0,...,0}_{M_0}),$}
which is a \morp\, from $(T_I'\oplus S_i')\otimes M_m$ to $B${{, where we identify $B$ with $B\otimes M_m$}}.
We then
 extend $L$ to a completely positive linear map $B\to {(1_B-{{(\one_{M_m(S_i')}-}}E_l{{)}})B(1_B-{{(\one_{M_m(S_i')}-}}E_l{{)}})}$.
 Also define $$R: {(T_i'\oplus S_i')\otimes M_m \to T_i'\oplus S_i'}$$ to be
 {{\begin{equation}
R
\left(
\begin{array}{cccc}
f_{1,1}\oplus g_{1,1} & f_{1,2}\oplus g_{1,2}&\cdots & f_{1,m}\oplus g_{1,m}\\
f_{2,1}\oplus g_{2,1} & f_{2,2}\oplus g_{2,2}&\cdots &f_{,m}\oplus g_{2,m}\\
 & \ddots & &\\
f_{m,1}\oplus g_{m,1} & f_{m,2}\oplus g_{m,2}&\cdots & f_{m,m}\oplus g_{m,m}
\end{array}
\right)=g_{1,1},
\end{equation}
where $f_{j,k}\in T'_i$ and $g_{j,k}\in S'_i$,}}
~and extend it to a {\morp}\, $B\to B$,
where ${T_i'\oplus S_i'}$ is regarded as a corner of ${(T_i'\oplus S_i')}\otimes M_m\subset B$. Then $L$ and $R$ are $\mathcal G$-$\epsilon'$-multiplicative. Hence $[L]|_\mathcal{P}$ is well defined. Moreover, $$\tau(L(1_A))<\tau(1_{T_1})+\frac{l}{m}<\frac{r}{2}+\frac{MK}{2MK/r}=r\rforal \tau\in T(A).$$

Note that for any $f$ in the form  \eqref{diag-f}, {{if $f$ is written in the form ${{(f_{jk}\oplus g_{jk})}}_{m\times m}$, then {{$g_{jj}=g_{11}$ and $g_{jk}=0$ for $j\not= k$}}. Hence }} one has
$$f=L(f)+{\overline{R}}(f),$$
{where $\overline{R}(f)$ may be written as}
$$
{\overline{R}(f)=\mathrm{diag}\{\underbrace{0,0,...,0}_l, \underbrace{(0\oplus g_{1,1}), ..., (0\oplus g_{1,1})}_{M_0}\}}.
$$
Hence for any $g\in G$,
$$g=[L](g)+ {M_0}[R](g).$$
Then, if $g\in\ (G_{0, 1})_+\subset (G_0)_+$, one {{has}}
$$g-[L](g)={M_0}[R](g)=K((\frac{M_0}{K})[R](g)).$$
And if $g\in G\cap K_i(B, \mathbb Z/ k\mathbb Z)$ ($i=0,1$), one also has
$$g-[L](g)={M_0}[R](g).
$$
Since $Mg=0$  and {$M|M_0$}, one has $g-[L](g)=0$.



{Since $L$ is {{the}} identity on $\psi_{i, \infty}'(T_i\otimes M_{m(i)})$ and $i>i_0,$ by (iii),
$L$ is the identity map on $G\cap \ker \rho_{B}.$  Since $K_1(S_i)=0$ for all $i,$
$L$ induces the identity map on $G\cap K_1(B).$ It follows that}
 $L$ is the desired map.
\end{proof}

Related to the {{considerations}} above we have the following decomposition:

{
\begin{prop}\label{lem-compress}
Let $A_1$ be a separable  amenable C*-algebra in ${\cal B}_1$ (or
${\cal B}_0$) and let $A=A_1\otimes U$ for some  infinite dimensional UHF-algebra $U$.
Let $\mathcal G\subset A,$
$\mathcal P\subset  {\underline{K}}(A)$ be finite subsets, ${\cal
P}_0\subset A\otimes {\cal K}$ be a finite subset of  projections,
and let $\epsilon>0$, $0<r_0<1$ and $M\in\mathbb N$ be arbitrary.
Then there {{are}} a projection $p\in A,$ a \SCA\, $B\in {\cal C}$ (or
${\cal C}_0$) with $p=1_B$ and $\mathcal
G$-$\epsilon$-multiplicative {{unital completely positive linear maps}} $L_1: A\to  (1-p)A(1-p)$
and $ L_2: A\to B$ such {{that:}}
\begin{enumerate}
\item
$\|L_1(x)+L_2(x)-x\|<\ep\tforal x\in {\cal G};$
\item $[L_i]|_{\cal P}$  is well defined, $i=1,2;$
\item $[L_1]|_{\cal P}+[\imath\circ L_2]|_{\cal P}=[{\rm id}]|_{\cal P};$
\item $\tau\circ[L_1](g)\leq r_0\tau(g)$ for all $g\in
      {\mathcal P}_0$ and $\tau\in T(A)$;
\item  For any $x\in {\cal P},$ there exists $y\in \underline{K}(B)$ such that
   $x-[L_1](x)=[\imath\circ L_2](x)=M[\imath](y)${{;
   and}}
\item for any $d\in {\cal P}_0,$ there exist positive element  $f\in {K_0}(B)_+$
      such that
      $$d -[L_1](d)=[\imath\circ L_2](d)=M[\imath](f),$$
      where $\imath: B\to A$ is the embedding.
      Moreover, we can require that $1-p\not=0.$
\end{enumerate}
\end{prop}
}
\begin{proof}
Since $A$ is in ${\cal B}_1$ (or ${\cal B}_0$), there is a sequence
of projections $p_n\in A$ and a sequence of \SCA s $B_n\in {\cal
B}_1$ ({{or}} ${\cal B}_0$) with $1_{B_n}=p_n$ such that
\beq\label{compress-1}
\lim_{n\to\infty}\|(1-p_n)a(1-p_n)+p_nap_n-a\|=0,\\
\lim_{n\to\infty} \dist(p_nap_n, B_n)=0{{,\andeqn}}\\
\lim_{n\to\infty}\max\{\tau(1-p_n): \tau\in T(A)\}=0 \eneq for all
$a\in A$.
 Since  each $B_n$ is
amenable, one obtains easily a sequence of {{unital completely positive linear maps}} $\Psi_n:
A\to B_n$ such that \beq\label{compress-4}
\lim_{n\to\infty}\|p_nap_n-\Psi_n(a)\|=0\rforal a\in A. \eneq In
particular, \beq\label{compress-5}
\lim_{n\to\infty}\|\Psi_n(ab)-\Psi_n(a)\Psi_n(b)\|=0\rforal a, b\in
A. \eneq Let  $j: A\to A\otimes U$ be defined by $j(a)=a\otimes
1_U.$ There is a unital \hm\, $s: A\otimes U\to A$ and  a sequence
of unitaries $u_n\in A\otimes U$ such that \beq\label{compress-6}
\lim_{n\to\infty}\|a-{\rm Ad}\, u_n\circ s\circ j(a)\|=0\rforal a\in
A. \eneq There are  non-zero projections $e_n'\in U$  and $e_n\in U$
such that \beq\label{compress-7} \lim_{n\to\infty}t(e_n)=0\andeqn
1-e_n={\rm diag}(\overbrace{e_n', e_n',...,e_n'}^M),
\eneq
{where $t\in T(U)$ is the unique tracial state on $U.$}
Choose
$N\ge 1$ such that
\beq\label{compress-8}
0<t(e_n)<r_0/2\andeqn
\max\{\tau(1-p_n):\tau\in T(A)\}<r_0/2.
\eneq
Define $\Phi_n: A\to
(1-p_n)A(1-p_n)$ by $\Phi_n(a)=(1-p_n)a(1-p_n)$ for all $a\in A.$
Define  $\Phi_n'(a)=\Phi_n(a) \oplus  {\rm Ad}\, u_n\circ s(a\otimes
e_n),$ and $\Psi_n'(a)={\rm Ad}\, u_n\circ s(\Psi(a)\otimes (1-e_n))$
for all $n\ge N.$ Note that $u_n^*s(B_n\otimes (1-e_n))u_n\in {\cal
C}_1$ (or
${\cal C}_0$). It is then easy to verify that, if we
choose a large $n,$ the maps $L_1=\Phi_n'$ and $L_2=\Psi_n'$ {meet}  the
requirements.
\end{proof}















\section{Positive maps from {{the}} $K_0$-group of a \CA\, in
$\mathcal C$ }

This section contains some technical lemmas about positive \hm s
from $K_0(C)$ for some $C\in {\cal C}.$

\begin{lem}[cf. 2.8 of \cite{Lnbirr}]\label{multiple-ext} Let
$G\subset \Z^l$ (for some $l>1$) be a subgroup. There is an integer $M>0$ satisfying the following condition:
Let $1>\sigma_1, \sigma_2>0$ { {be any given numbers.}}
{{Then, there}} is an integer ${{R}}>0$ {{such that}}:
{{if a set of  $l$ positive numbers $\af_i\in \R_+$ ($i=1,2,{...,} l$) satisfies
 $\af_i\ge \sigma_1$ for all $i$ and
satisfies}}
\beq\label{mule-0}
\sum_{i=1}^l \af_i m_i\in \Z \tforal (m_1, m_2,...,m_l)\in G,
\eneq
then for any integer $K\geq R$, there
{{exists}} a set of positive rational  numbers
$ \bt_i\in \frac{1}{KM}\mathbb {{Z_+}}$ ($i=1,2,..., l$)
such that
\beq\label{mule-1}
\sum_{i=1}^l|\af_i-\bt_i{{|}}<\sigma_2
\tand {\tilde \phi}|_G=\phi|_G,
\eneq
where $\phi((n_1, n_2,...,n_l))=\sum_{i=1}^l \af_in_i$ and ${\tilde \phi}((n_1,n_2,...,n_l))=\sum_{i=1}^l\bt_in_i$
for all
$(n_1,n_2,...,n_l)\in \Z^l.$
\end{lem}

\begin{proof}
Denote by $e_j\in \Z^l$ the element having 1 in the $j$-th coordinate and $0$ elsewhere.
First we consider the case that $\Z^l/G$ is finite.
In this case there is an integer $M\ge 1$ such that
$Me_j\in G$ for all $j=1,2,...,l.$  It follows that $\phi(Me_j)\in \phi(G){{\subset \Z}},$ $j=1,2,...,l$.
Hence $\af_j=\phi(e_j)\in \frac{1}{M}\mathbb {{Z_+}}$.
We choose $\bt_j=\af_j,$ $j=1,2,...,l,$ and $\tilde \phi=\phi$. The conclusion of the lemma  follows{{---that is, for any $\sigma_1, \sigma_2$, we can choose $R=1$.}}

Now we assume that $\Z^l/G$ is not finite.
Regard $\Z^l$ as a subset of $\Q^l$ and set
 $H_0$ to be the  vector subspace of $\Q^l$  spanned by elements in $G.$
{{Since  $G$  is a subgroup of $\Z^l,$ it must be finitely generated and hence
is isomorphic to $\Z^p,$ where $p$ is the (torsion--free) rank of $G.$  Since  the rank of
$\Z^l/G$ is strictly positive (otherwise $\Z^l/G$ would be a finitely generated torsion group, hence finite),
and using the fact that the rank is additive, one concludes that $0\le p<l.$}}

Let $g_1,g_2,...,g_p\in G$ be indepdent generators of $G.$ View them as elements
in $H_0\subset \Q^l$ and
write
\beq\label{mule-2}
g_i={{\left(g_{i, 1}, g_{i, 2},{...}, g_{i, l}\right)}},\,\,\,i=1,2,...,p.
\eneq

Define $L: \Q^p\to \Q^l$ to be $L=(f_{i,j})_{l\times p},$ where $f_{i,j}=g_{j,i},$ $i=1,2,...,l$ and
$j=1,2,...,p.$
Then $L^*=(g_{i,j})_{p\times l}.$  We also view $L^*:
\Q^l\to \Q^p.$
Define $T=L^*L: \Q^p\to \Q^p.$ {{This map}} is invertible.  Note that $T=T^*$ and $(T^{-1})^*=T^{-1}.$
Note  also that the matrix representation $(a_{i,j})_{{{l\times p}}}$ of $L\circ T^{-1}$ is {{an}} ${{l\times p}}$ matrix with entries
in $\Q.$ There is an integer $M_1\ge 1$ such that $a_{i,j}\in {1\over{M_1}}\Z,$ $i=1,2,...,{{l}}$ and
$j=1,2,...,{{p}}.$

Let $H_{00}={\rm ker} L^*.$ It has dimension $l-p>0.$ Let $P: \Q^l\to H_{00}$ be an orthogonal projection
which is a $\Q$-linear map. Represent $P$ as {{an}} $l\times l$ matrix.
Then its entries are in $\Q.$ There is an integer $M_2\ge 1$ such that
all entries are in
$\frac{1}{M_2}\Z.$
We will use the fact that $L^*=L^*(1-P).$

It is important to note that $M_1$ and $M_2$ depend on
$G$ only and are independent of $\{\af_j: 1\le j\le l\}.$  Let $M=M_1M_2.$

Suppose that $\sigma_1,\sigma_2 \in (0,1)$ are two positive numbers and the numbers $\alpha_i\geq \sigma_1$ ($i=1,2,{...},l$) satisfy
 condition (\ref{mule-0}).
The condition (\ref{mule-0}) is equivalent to {{the condition}} that ${{b_i\hspace{-0.05in}:=}}\sum_{j=1}^l \af_jg_{i,j}\in \Z,$ $i=1,2,...,p.$ Put
$b=(b_1,b_2,...,b_p)^T$ and $\af=(\af_1, \af_2,...,\af_l)^T.$   Then
$b=L^*\af.$


If we write
\beq\label{mule-3}
L (T^*)^{-1}b=c=\left(\begin{array}{c} c_1\\ c_2\\ \vdots\\ c_l\end{array}\right),
\eneq
then, since $b\in \Z^p,$ one has
$c_i\in \frac{1}{M_1}\Z.$
Choose an integer ${ {R}}\ge 1$ such that $1/{{R}}<\sigma_1\sigma_2/{{(4l^2)}}$. { {Let $K\geq R$ be any integer.}}
Note that
\beq\label{mule-3+}
L^*c=L^*L(T^*)^{-1}b=L^*LT^{-1}b=L^*\af.
\eneq
 Thus $\af-c\in{\rm ker} L^*$ as a subspace of $\R^l.$

For the space $\R^l$, we use $\|\cdot\|_1$ and $\|\cdot\|_2$ to denote the $l_1$ norm and $l_2$ norm on it, {respectively}. Then we have
 $$\|v\|_2 \leq \|v\|_1\leq l \|v\|_2~~~\rforal v\in \R^l.$$

 Since $H_{00}$ is dense in the real {{subspace}}
 ${\rm ker} L^*,$
there exists $\xi\in H_{00}$ such that
\beq\label{mule-4}
\|\af-c-\xi\|_2{{\leq \|\af-c-\xi\|_1}}<\sigma_1\sigma_2/{{(4l)}}.
\eneq
Pick $\eta\in \Q^l$ such that $\xi=P\eta.$
Then there is $\eta_0\in \Q^l$ such that $K\eta_0\in \Z^l$ and
\beq\label{mule-5}
{ {\|\eta_0-\eta\|_2\leq}}\|\eta_0-\eta\|_{{1}}<\sigma_1\sigma_2/{{(2l)}}.
\eneq
Since $P$ has norm
{{one}} with respect to {{the}} $l_2$ norm,
\beq\label{mule-6}
\|\af-c-P\eta_0\|_{{1}}{ {\leq l\|\af-c-P\eta_0\|_2 \leq l( \|\af-c-\xi\|_2+ \|P(\eta_0-\eta)\|_2) }}< \sigma_1\sigma_2.
\eneq
Put $\bt=c+P\eta_0=(\bt_1, \bt_2,...,\bt_l)^T.$ Note that $M_2K(P\eta_0)\in \Z^l$ {{, and that $M_1c\in \Z^l$.}}

{ {We have}}  ${ {K}}M\bt\in \Z^l,$ { {and}}
\beq\label{mule-7}
L^*\bt=L^*c=L^*\af\andeqn \|\af-\bt\|_{{1}}<\sigma_1\sigma_2.
\eneq
Moreover, since $\af_i\ge \sigma_1,$
\beq\label{mule-8}
\bt_i>0,\,\,\, i=1,2,...,l.
\eneq
Since $P\eta_0\in H_{00},$  one has that $L ^*\bt =L^*(1-P)\bt=L^*(1-P)c=L^*c=L^*\af=b.$
Define ${{\tilde \phi}}: \Q^l\to \Q$ by
\beq\label{mule-9}
{{\tilde \phi}}(x)=\langle x, \bt \rangle
\eneq
for all $x\in \Q^l.$
Note $L^*e_i=g_i,$ where
$e_i$ is the element in $\Z^p$ with  $i$-th coordinate  $1$ and $0$ elsewhere.  So
\beq\label{mule-10}
\phi(g_i)&=&\langle Le_i, \af\rangle =\langle e_i, L^*\af\rangle =\langle e_i, L^*\bt\rangle \nonumber \\
&=& \langle Le_i, \bt \rangle = \langle g_i, \bt \rangle ={{\tilde \phi}}(g_i),
   \eneq
$i=1,2,...,p.$
It follows that
${{\tilde \phi}}(g)=\phi(g)$ for all $g\in G.$  { {Hence ${\tilde \phi}|_G=\phi|_{G}.$}}
Note that
${\tilde \phi}(\Z^l)\subset {1\over{{{K}}M}}\Z$, since $\bt_i\in {1\over{{{K}}M}}
{{\Z_+}},$ $i=1,2,...,l.$
\end{proof}

If we do not need to approximate $\{\af_i: 1\le i\le l\},$ then ${{R}}$
can be chosen to be $1$, with $M={{M_1,}}$ which only depends on $G$ and $l$ (by replacing  $\bt$ by $c$ in the proof).

From the proof {{of}} \ref{multiple-ext}, since $L$ and $(T^*)^{-1}$ {depend} only on $g_1,g_2,...,g_p,$ we have the following {{result}}:

\begin{lem}\label{ext-norm}
Let $G\subset \Z^l$ be an ordered subgroup with order unit $e$, and let
$g_1,g_2,...,g_p$ ($p\le l$) be a set of free generators of $G.$  For any $\ep>0,$ there exists $\dt>0$ satisfying the
following condition:
if $\phi: G\to \R$ is a \hm\,
such that
$$
|\phi(g_i)|<\dt,\,\,\,i=1,2,...,p,
$$
then there is $\bt=(\bt_1, \bt_2,...,\bt_l)\in \R^l$ with
$|\bt_i|<\ep,$ $i=1,2,...,l,$ such that
$$
\phi(g)=\psi(g)\rforal g\in G,
$$
where $\psi: \Z^l\to \R$ {{is}} defined by $\psi((m_1,m_2,...,m_l))=\sum_{i=1}^l\bt_im_i$
for all $m_i\in \Z.$
\end{lem}


\begin{cor}\label{MextC}
Let $G\subset \Z^l$ be an {{ordered}} subgroup. Then, there exists  an integer $M\ge 1$
satisfying the following condition:
for any positive {{\hm}} $\kappa: G\to \Z^n$ (for any integer $n\ge 1$) with every element in $\kappa(G)$  divisible by $M$, there is  $R_0\ge 1$ such that, for any integer $K\ge R_0,$
there is a positive homomorphism ${\tilde \kappa}: \Z^l\to \Z^n$ such that
${\tilde \kappa}|_G=K\kappa.$
\end{cor}

\begin{proof}
We first prove the case that $n=1.$

Let $S\subset \{1,2,...,l\}$ be a subset and denote by $\Z^{(S)}$ the subset
$$
\Z^{(S)}=\{(m_1, m_2,...,m_l): m_i=0\,\,\,{\rm if}\,\,\, i\not\in S\}.
$$
Let $\Pi_S: \Z^l\to \Z^{{{(}}S{{)}}}$ be
the {{obvious}} projection and $G(S)=\Pi_S(G).$

{ {Let $M(S)$ be the  {{integer}} (in place of $M$) as in
 \ref{multiple-ext} associated with $G(S)\subset \Z^{(S)}.$}}
 Put $M=\prod_{S\subset \{1,2,...,l\}}{{M(S)}}.$

 Now assume that $\kappa: G\to \Z$ is a positive \hm\, with multiplicity $M$---that is, every element in $\kappa(G)$ is  divisible by $M$.

 By applying   2.8 of \cite{Lnbirr}, we obtain a
positive \hm\, $\bt: \Z^l\to \R$ such that $\bt|_G=\kappa.$
Define $f_i=\bt(e_i),$ where $e_i$ is the element in $\Z^l$ with $1$ at the $i$-th coordinate and
$0$ elsewhere, $i=1,2,...,l.$  Then $f_i\ge 0.$
Choose $S$ such that $f_i>0$ if $i\in S$ and $f_i=0$ if $i\not\in S.$

{{Evidently if $\xi_1,\xi_2\in \Z^l$ satisfy $\Pi_S(\xi_1)=\Pi_S(\xi_2)$, then $\bt(\xi_1)=\bt(\xi_2)$,  and if we  further
assume $\xi_1,\xi_2\in G$ then $\kappa(\xi_1)=\kappa(\xi_2)$. Hence  the maps $\bt$ and $\kappa$ induce  maps $\bt': \Z^{(S)} \to {{\R}}$ and $\kappa': G(S) \to \Z$ such that $\bt=\bt'\circ \Pi_S$ and $\kappa=\kappa'\circ
{{(\Pi_S)|_G}}$. In addition, we have $\bt'|_{G(S)}=\kappa'$.}}

Let $\sigma_1=2\sigma_2=\frac{\min \{f_i:~ i\in S\}}{2M}$. {{By applying}} Lemma \ref{multiple-ext} to $\af_i=f_i/M> \sigma_1$ for $i\in S$ and to $G(S)\subset \Z^{(S)},$ we obtain the number $R(\kappa)$ (depending on $\sigma_1$ and $\sigma_2$ and therefore depending on $\kappa$) as in the lemma. For any $K\geq R(\kappa)$, it follows from the lemma
{{that}} there are $\bt_i\in \frac{1}{KM}{{\Z_+}}$ (for $i\in S$) such that $\tilde \kappa'|_{G(S)}=\frac{1}{M}\kappa'$, where $\tilde \kappa': \Z^{(S)}\to \Q$ is defined by $\tilde \kappa'(\{m_i\}_{i\in S})=\sum_{i\in S}\bt_i m_i$. Evidently, $\tilde \kappa=KM(\tilde \kappa'\circ\Pi_S): \Z^l\to \Z$ is as desired for this case.



This {{proves}} the case $n=1.$

In general,
 let $s_i: \Z^n\to \Z$ be the projection
{{onto}} the $i$-th {{direct}} summand, $i=1,2,...,n.$  { {Apply the case $n=1$ to each of the the maps $\kappa_i:=s_i\circ \kappa$ {{(for $i=1,2,...,n$)}} to obtain $R(\kappa_i),$ {{and}} let $R_0=\max_{i} R(\kappa_i)$.}}  For any $K\ge R_{{0}},$
by what has been proved, we obtain ${\tilde\kappa}_i: \Z^l\to \Z$ such that
\beq\label{MextC-4}
{\tilde \kappa}_i|_G=Ks_i\circ \kappa|_G,\,\,\,i=1,2,...{{,n.}}
\eneq
Define ${\tilde \kappa}: \Z^l \to \Z^n$ by
${\tilde \kappa}(z)=({\tilde \kappa}_1(z),{\tilde \kappa}_2(z),...,{\tilde \kappa}_{{n}}(z)).$
The lemma follows.
\end{proof}




\begin{lem}[Lemma 3.2 of \cite{EN-Tapprox}]\label{annihilation2}
{{Let $G$ be a countable abelian unperforated ordered group  such that $G_+$ is finitely generated}}
and let $r: G\to \Z$ be a
strictly positive homomorphism, {{i.e., $r(G_+\setminus \{0\})\subset \Z_+\setminus \{0\}.$}} Then,  for any  order unit $u\in
G_+, $  there exists a natural number $m$ such that if the map
$\theta: G\to G$ is defined by $g\mapsto r(g)u$,
{{then  the}} positive homomorphism ${
{{\rm{id}}}}+m\theta: G\to G$ factors through $\bigoplus_{i=1}^n\Z$ positively for some $n$.
\end{lem}

\begin{proof}
Let $u$ be an order unit of $G,$ {{which exists as $G$ is finitely generated}}, and define the map $\phi: G\to G$ by $\phi(g)=g+ r(g)u$; that is, $\phi=\id+\theta$.
Define $G_n=G$ and $\phi_n: G\to G$ by $\phi_n(g)=\phi(g)$ for all $g$ and $n.$ Consider the inductive limit
$$
\xymatrix{
G\ar[r]^\phi & G \ar[r]^\phi & \cdots \ar[r] & \varinjlim G.
}
$$
Then the ordered group $\varinjlim G$ has the Riesz decomposition property.
In fact, let $a, b, c\in {{(\varinjlim G)_+}}$ such that
$$
a\le b+c.
$$
Without loss of generality, one may assume that $a\not=b+c.$

We may assume that there are $a', b', c'\in {{G_+}}$ for {{some}} $n$-th {{(finite stage)}} $G$ such that
$\phi_{n, \infty}(a')=a,$ $\phi_{n, \infty}(b')=b,$ and $\phi_{n, \infty}(c')=c,$ {{and furthermore }}
\beq\label{ann2-n1}
{{a' < b'+c'}}.
\eneq

A straightforward calculation shows that for each $k$,  there is $m(k)\in\mathbb N$ such that
$$
\phi_{n, n+k}(a')=a'+m(k)r(a')u,\,\, \phi_{n,n+k}(b')=b'+m(k)r(b')u, \andeqn \phi_{n,n+k}(c')={{c'}}+m(k)r(c')u.
$$
Moreover, the sequence $(m(k))$ is strictly increasing.
Since $r$ is strictly positive, {{combining}} with (\ref{ann2-n1}), we have that
$$
r(a{{'}})< r(b{{'}})+r(c{{'}})\,\,\,{\rm (in} \,\,\, {{\Z_+}}{\rm )}.
$$
There are, {{for $i=1,2,$}}
{{$l_i\in \Z_+$}} such that
$$
l_1+l_2=r(a{{'}}),\,\,l_1\le r(b{{'}}),\andeqn l_2\le r(c{{'}}).
$$
Without loss of generality, we may assume $d=r(b{{'}})-l_1
>0$ (otherwise
we let $d=r(c{{'}})-l_2$).
Since $u$ is an order unit, there is $m_1\in \Z_+$ such that
$$
m_1d u>a{{'}}.
$$
Choose $k\ge 1$ such that $m(k)>m_1.$
Let $a_1=a'+m(k){{l_1}}u$ and $a_2=m(k){{l_2}}u.$
Then
$$
a_2 = m(k)l_2u\le m(k)r(c')u\le c'+m(k)r(c')u=\phi_{n,n+k}(c').
$$
Moreover,
$$
a_1=a'+m(k)l_1u\le m(k)du+m(k)l_1u\le b'+m(k)r(b')u=\phi_{n, n+k}(b').
$$
Note
$$
\phi_{n, n+k}(a')=a_1+a_2\le \phi_{n, n+k}(b')+\phi_{n,n+k}(c').
$$
These {{inequalities}}  imply that
\beq
a=\phi_{n+k,\infty}(a_1)+\phi_{n+k{{,\infty}}}(a_2)\le b+c,\\
\phi_{n+k,\infty}(a_1)\le b\andeqn \phi_{n+k,\infty}(a_2)\le c.
\eneq

It follows  that the limit group $\varinjlim G$ has the Riesz decomposition property.
Since $G$ is unperforated, so is $\varinjlim G$. It then follows from the Effros-Handelman-Shen Theorem (\cite{EHS}) that  the ordered group $\varinjlim G$ is a dimension group. Therefore,
{{since $G$ and $G_+$ are finitely generated, for sufficiently large $k,$ the map $\phi_{1, k}: G\to G$ factors
positively through $\bigoplus_{i=1}^n\mathbb Z$ for some $n.$ As already pointed out, the map
$\phi_{1,k}$ is of the desired form  ${\rm id}+m\theta$  for some $m\ge 0.$
The lemma follows.}}
\end{proof}

{\Green{Let $G_1$ and $G_2$ be groups and $K\ge 1$ be an integer.  A \hm\, $\gamma: G_1\to G_2$
is said to have multiplicity $K,$ if there is, for each $g_1\in G_1,$ $\gamma(g_1)=Kg_2$ for some
$g_2\in G_2.$}}

\begin{lem}\label{Cuthm}
Let $(G,G_+,u)$ be  {{an ordered}} group
 with  order unit $u$ 
such that the positive cone $G_+$ is generated by finitely many positive elements which are {{strictly}} smaller than $u.$ Let  $\lambda: G\to K_0(A)$ be an order preserving map such that $\lambda(u)=[1_A]$  and
$\lambda(G_+\setminus \{0\})\subset K_0(A)_+\setminus \{0\},$ where $A\in {\cal B}_1$ ({{or}} $A\in {\cal B}_0$).
Let $a\in K_0(A)_+\setminus \{0\}$ with $a\le [1_A].$ 
Let ${\cal P}\subset G_+\setminus \{0\}$ be a finite subset.
Suppose that there exists an integer
$N\ge 1$ such that $N\lambda(x)>[1_A]$ for all $x\in {\cal P}.$

{{Then, there are two positive homomorphisms $\lambda_0: G\to
K_0(A)$ and $\gamma: G\to K_0(S'),$   with $S'\subset A$  and  $S'\in {\cal C}$ ({{or}} ${\cal
C}_0$)  such that  $\gamma(u)=[1_{S'}]$ and}}
\beq\label{Cuthm-1}
{{\lambda=\lambda_0+\lambda_1,\,\,\,
{{0\leq}}\lambda_0(u)<a}}\andeqn \gamma(g)>0 \eneq for all
$g\in G_+\setminus \{0\},$ {{where $\lambda_1=\imath_{*0}\circ
\gamma$}}
and $\imath: S'\to A$ is the
embedding.
Moreover, $N\gamma(x)\ge \gamma(u)$ for all $x\in
{\cal P}$. Furthermore, if $A=A_1\otimes U$, where $U$ is an
infinite dimensional UHF-algebra, and $A_1\in\mathcal B_1$ (or ${\cal
B}_0$),  then, for any integer $K\ge 1,$ we can require that
$S'=S\otimes M_K$ for some $S\in {\cal C}$ (or ${\cal C}_0$) and
$\gamma$ has multiplicity $K.$
\end{lem}

\begin{proof}
Let $\{g_1, g_2,..., g_m\}\subset G_+$ be {{a}} set of generators
of $G_+$ with $g_i < u$. Since $A$ has stable rank one, it is easy to check that
there are projections $q_1, q_2,...,q_m\in A$ such that
$\lambda(g_i)=[q_i],$ $i=1,2,...,m.$
To simplify
notation,
let us assume that ${\cal P}=\{g_1,g_2,...,g_m\}.$
 Define
$$
Q_i={\rm diag}(\overbrace{q_i, q_i,...,q_i}^{{{N}}}),\,\,\, i=1,2,...,m.
$$
By the assumption, there are $v_i\in M_N(A)$ such that
\beq\label{Cutm-14-n1} v_i^*v_i=1_A\andeqn v_iv_{{i}}^*\le
Q_i,\,\,\,i=1,2,...,m. \eneq
 Since $A\in {\cal B}_1$ (or ${\cal
B}_0$), there exist a sequence of projections $\{p_n\}$ {{in}} $A,$ a
sequence of \SCA s $S_n\in {\cal C}$ (${\cal C}_0$) with
$p_n=1_{S_n},$ and a sequence of unital
{{\cp s}} $L_n: A\to S_n$ such
that \beq\label{Cuthm-n1}
&&\hspace{-0.2in}\lim_{n\to\infty}\|a-((1-p_n)a(1-p_n)+p_nap_n)\|=0,\,\,\,\lim_{n\to\infty}\|L_n(a)-p_nap_n\|=0,\\\label{Cuthm-n1+}
&&
{{\andeqn \lim_{n\to\infty}\sup\{\tau(1-p_n): \tau\in T(A)\}=0,}}
\eneq
{{for all $a\in A.$}}
It is also standard to find, for each $i,$  a projection
$e_{i,n}'\in (1-p_n)A(1-p_n),$ a projection $e_{i,n}\in M_{{N}}(S_n),$
and {{a partial isometry}} $w_{i,n}\in {{M_N(}} S_n{{)}}$ such that
\beq\label{Cuthm-nn1}
&&\lim_{n\to\infty}\|(1-p_n)q_i(1-p_n)-e_{i,n}'\|=0,\\\label{Cuthm-nn2}
&&w_{i,n}^*w_{i,n}=p_n,\,\,\, w_{i,n}w_{i,n}^*\le
e_{i,n},\\\label{Cuthm-nn3}
&&\lim_{n\to\infty}\|{{(}}L_n{{\otimes {\rm id}_{M_N})}}(v_i)-w_{i,n}\|=0,\andeqn\\
&& \lim_{n\to\infty}\|
(L_n\otimes {\rm id}_{M_{{N}}})(Q_i)-e_{i,n}\|=0.
\eneq
Let
$\Psi_n: {{A}}\to (1-p_n)A(1-p_n)$ be defined by $\Psi_n(a)=(1-p_n)a(1-p_n)$ for
all $a\in A.$ We will use $[\Psi_n]\circ \lambda$  for $\lambda_0$
and ${{[L_n]\circ \lambda}}$ for $\gamma$ for  some large $n.$  {{Since $G$ is finitely generated,
choosing sufficiently large $n,$ we may assume that $\lambda_0$ and $\gamma$ are \hm s.}}
To see that $\lambda_0$ is  positive, we use
(\ref{Cuthm-nn1}) and the fact that $G_+$ is finitely generated. It
follows from (\ref{Cuthm-nn2}) and (\ref{Cuthm-nn3})  that
$N\gamma(x)\ge \gamma(u)$ for all $x\in {\cal P}.$ Since we assume
that the positive cone of $G_+$ is generated by ${\cal P},$ this
also shows that $\gamma(x)>0$ for all $x\in G_+\setminus \{0\}.$
By (\ref{Cuthm-n1+}), we can choose large $n$ so that
 ${{0\leq}}\lambda_0(u)<a.$

It should be noted {{that}} when $A$ does not have (SP), one can choose $\lambda=\lambda_1$  and $\lambda_0=0.$

If  $A=A_1\otimes U,$ then, \wilog, we may assume that $p_n\in A_1$ for all $n.$
Choose  a sequence of non-zero projections $e_n\in U$ such that $t(1-e_n)=r(n)/K,$
where $t$ is the unique tracial state on $U$ and $r(n)$ are positive rational numbers such that
$\lim_{n\to\infty}{{r(n)}}=0.$
Thus $S_n\otimes (1-e_n)\subset B_n$ where $B_n\cong S_n\otimes M_K$ and
$p_n\otimes (1-e_n)=1_{B_n}.$ We check that the lemma follows if we replace
$\Psi_n$ by $\Psi_n',$ where $\Psi_n'(a)=(1-p_n)a\otimes 1_U+p_na\otimes e_n$
\end{proof}

\begin{lem}[see Lemma 3.6 of \cite{EN-Tapprox} or Lemma 2.8 of \cite{Niu-TAS-II}]
\label{decomposition2}
Let $G={K_0}(S)$, where $S\in \mathcal C$. Let $H={K_0}(A)$ for  $A=A_1\otimes U,$ where $A_1\in\mathcal B_1$
(or ${\cal B}_0$), and $U$ is  {{an}} infinite dimensional UHF-algebra.
Let $M_1\ge 1$ be a given integer and $d\in K_0(A)_+\setminus \{0\}.$ Then for any strictly positive homomorphism $\theta: G\to H$ with multiplicity $M_1$, and any integers $M_2\ge 1$ and $K\ge 1$ such that
$K\theta(x)>[1_{A}]$ for all $x\in G_+\setminus\{0\},$ one has a decomposition $\theta=\theta_1+\theta_2$, where $\theta_1$ and $\theta_2$ are positive homomorphisms from $G$ to $H$ such that the following diagrams commute:
\begin{displaymath}
\xymatrix{
G \ar[rr]^{\theta_1} \ar[dr]_{\phi_1} &                      & H & & G \ar[rr]^{\theta_2} \ar[dr]_{\phi_2} &                     & H \\
                                      & G_1 \ar[ur]_{\psi_1} &   & &                                      & G_2 \ar[ur]_{\psi_2} & },
\end{displaymath}
where $\theta_1([1_{{S}}])\le d,$
$G_1\cong\bigoplus_n\Z:=\bigoplus_{i=1}^n\Z$ for some natural number $n$ and $G_2=K_0(S')$ for some
\SCA\, $S'$ of $A$ which is in the class $\mathcal C$ (or in ${\cal C}_0$),  $\phi_1, \psi_1$ are
positive \hm s
  and $\psi_2=\imath_{*0},$ where $\imath: S'\to A$  is the embedding.
Moreover, $\phi_1$ has  multiplicity $M_1,$ $\phi_2$ has multiplicity $M_1M_2$ {{and}}  $2K\phi_2(x)>\phi_2([1_{{S}}]){{=[1_{S'}]}}>0$ for all $x\in {{G_+\setminus\{0\}}}.$
\end{lem}

\begin{proof}
Let $u=[1_S].$
Suppose that $S={{A}}(F_1, F_2, \psi_0, \psi_1)$ with $F_1=M_{R(1)}\oplus M_{{R(2)}}\oplus \cdots \oplus M_{R(l)}.$
It is easy to find  a strictly positive  \hm\, $\eta_0: K_0( F_1)\to \Z.$  Define
$r: G\to \Z$ by $r(g)=\eta_0\circ (\pi_e)_{*0}.$
By replacing $S$ with $M_d(S)$  and $A$ by $M_{{d}}(A)$ for some integer ${{d}}\ge 1, $
\wilog{,} we may assume that $S$ {{contains}} a finite subset of projections
 ${\cal P}=\{p_1, p_2, ..., p_{{l}}\}$  such that
 every projection $q\in S$ is equivalent to one of {{the}} projections in ${\cal P}$ and
 $\{[p_i]:1\le i\le {{l}}\}$ generates $K_0(S)_+$ (see \ref{FG-Ratn}).
Let
 $$
 \sigma_0=\min\{\rho_A(d)(\tau): \tau\in T(A)\}.
 $$
 Note that since $A$ is simple, one has  $\sigma_0>0.$

Let
$$
\sigma_1=\inf\{\tau(\theta([p]){{)}}: p\in {\cal P},\,\,\, \tau\in T(A)\}>0.
$$
Since $A=A_1\otimes U,$ $A$  has the (SP) property, there is a projection $f_0\in A_+\setminus \{0\}$ such that
\beq\label{decomp6/7-1}
0<\tau(f_0)< \min\{\sigma_0, \sigma_1\}/8{{K}}r(u)\rforal \tau\in T(A).
\eneq

Since $A=A_1\otimes U,$ we may choose
$f_0$ so that ${{[f_0]}}=M_1\tilde{h}$ for some non-zero ${\tilde h}\in K_0(A)_+.$  {{Define}}
$\theta_0': G\to K_0(A)$ by $\theta_0'(g)=r(g){\tilde h}$ for all $g\in G.$ { {Set $\theta'=M_1\theta_0'$. Then
$2\theta'(x) <\theta(x)$ for all $x\in G_+\setminus \{0\}$.}}

Since  $\theta$ has multiplicity $M_{1}$, one has that $\theta(g)-\theta'(g)$ is divisible by $M_1$ for any $g\in G$. By the choice of $\sigma_0,$  one checks that $\theta-\theta'$ is strictly positive.
Moreover,
\beq\label{decom-nn1}
2K\rho_A((\theta(x)-\theta'(x))(\tau) &{{>}}& 2K\rho_A(\theta(x))(\tau)-{{K}}\rho_A(\theta(x))(\tau)\\
&=& K\rho_A(\theta(x))(\tau)
\ge  \rho_A([1_A])(\tau)\rforal \tau\in T(A).
\eneq

Applying \ref{Cuthm}, one obtains a \SCA\, $S'{{\subset}} A,$  a  {{positive}} \hm\, ${\tilde \theta}_1: G\to K_0(A)$ and
strictly positive \hm\, $\phi_2: G\to K_0(S')$ such that
\begin{equation}\label{decompo6/8-1}
\theta-\theta' = {\tilde \theta}_1+\imath_{*0}\circ \phi_2,
\end{equation}
\begin{equation}\label{decomp6/8-2}
0\le \tau({\tilde \theta}_1(u))< \frac{\tau({\tilde h})}{mM_1{{ M_2}}},\quad \tau\in T(A),\,\,
2K\phi_2(x)>\phi_2([1_S]),\,\,\,\phi_2([1_S])=[1_{S'}]{{,}}
\end{equation}
 where $m$ is from \ref{annihilation2} and $\phi_2$ has multiplicity  $M_1M_2,$
  {{and}} where $\imath: S'\to A$ is the embedding.
Put $$\theta_2=\imath_{*0}\circ \phi_2,\quad \textrm{and}\quad \psi_2=\imath_{*0}.$$

Since $\theta(g)-\theta'(g)$ is divisible by $M_1$ and any element in $\theta_2(G)$ is divisible by $M_1$, one has that any {{element}} in $\tilde{\theta}_1(G)$ is divisible by $M_1$. Therefore, the map $\tilde{\theta}_1$ can be decomposed further as $M_1\theta_1'$, and one has that $\theta-\theta'=M_1\theta_1'+\theta_2$. Therefore, there is a decomposition
$$\theta=\theta'+M_1\theta'_1+\theta_2=M_1\theta'_0+M_1\theta'_1+\theta_2.$$
Put $\theta_1=M_1(\theta'_0+\theta_1').$
Then,
$$
\rho_A(\theta_1([1_S]))(\tau)<
{{(1/2)\rho_A(d)(\tau)}}\rforal \tau\in T(A).
$$

We {{now}} show that $\theta_1$ has the desired factorization property.
For $\theta'_0+\theta'_1$, one has the following {{further}} decomposition: for any $g\in G$,
\begin{eqnarray*}
\theta'_0(g)+\theta'_1(g)&=&r(g){\tilde h}+\theta'_1(g)
=r(g)({\tilde h}-m\theta'_1(u))+r(g)m\theta'_1(u)+\theta_{{1}}'(g)\\
&=& r(g)({\tilde h}-m\theta'_1(u))+ \theta'_1(mr(g)u)+\theta'_1(g)\\
&=& r(g)(\tilde{h}-m\theta'_1(u))+\theta'_1(mr(g)u+g).
\end{eqnarray*}
By (\ref{decomp6/8-2}), ${\tilde{h}}-m\theta_1'(u)>0.$
By Lemma \ref{annihilation2}, $g\, {{\mapsto}}\, mr(g)u+g$ factors through $\bigoplus_n \Z(=\bigoplus_{i=1}^n\Z)$ positively for some $n$.
Therefore, the map $M_1(\theta'_0+\theta_1')$ factors though $\bigoplus_{(1+n)M_1}\Z$ positively.
So there are positive homomorphisms $\phi_1: G\to \bigoplus_{(1+n)M_1}\Z$ and $\psi_1: \bigoplus_{(1+n)M_1}\Z\to K_0(A)$
such that $\theta_1=\psi_1\circ \phi_1$ and $\phi_1$ has multiplicity  $M_1$.
\end{proof}








\section{Existence Theorems for affine maps
for building blocks}

{{The main purpose of this section is to present  Theorem \ref{ExtTraceC-D}.
We first consider the case that the target algebras are finite dimensional (Lemma \ref{ExtTraceMn}).
We then  replace
them by interval algebras via a path (Lemma \ref{ExtTraceI-D}).  Then we establish Theorem \ref{ExtTraceC-D}
for the target algebras in ${\cal C}.$
Methods used in this section  may also be found  later  in  \cite{EGLN}.}}

\begin{lem}\label{discretizeT0}
Let $A$ be a unital separable  \CA\,  with $T(A)\not=\emptyset$ and
let ${\cal H}\subset A$  be a finite  subset.
Then, for any $\sigma>0,$ there exist an integer $N>0$  and a finite subset
$E\subset \partial_e(T(A))$ with the following  property:
For any $\tau\in T(A)$ and any $k\ge N$, there are $t_1, t_2,...,t_k$  {{(which are not  necessarily distinct
points) in}} $E$ such that
\beq\label{161N-1}
|\tau(h)-\frac{1}{k}\sum_{i=1}^k t_i(h)|<\sigma \tforal h\in\mathcal H.
\eneq
 { {(If $\tau$ is a  (not necessarily normalized)
 trace on $A$ with $\|\tau\|\leq 1$, then there are $t_1, t_2,...,t_{k'}$
 {{in $E$}} with $k'\leq k$ such that
$$
|\tau(h)-\frac{1}{k}\sum_{i=1}^{k'} t_i(h)|<\sigma \tforal h\in\mathcal H.)
$$}}
  Suppose that $A,$ {{in addition,}}  is a subhomogeneous \CA. Then
there are $\pi_1, \pi_2, ..., \pi_k$ {{in}} ${\rm Irr}(A)$ such that
\beq\label{161N-2}
|\tau(h)-\frac{1}{k}\left({\rm tr}_1\circ \pi_1(h)+{\rm tr}_2\circ \pi_2(h)+\cdots+{\rm tr}_k\circ \pi_k(h)\right)|<\sigma \tforal h\in\mathcal H,
\eneq
where  { {$\pi_j\in E$ and}} ${\rm tr}_j$ is the  {{ unique tracial state}} of $\pi_j(A)$.
{{Moreover, if the space   $\hat{A}_{{l}}$ of irreducible representations of dimension exactly $l$  has no isolated points
for each $l,$  then
$\pi_1,\pi_2,...,\pi_k$ can be  chosen  to be distinct.}}
\end{lem}


{\bf Remark:} Note that in (\ref{161N-2}), {{the}} $\pi_i$ may not be distinct. But the subset $E$ of   irreducible
representations can chosen to be  independent of $\tau$ (but depending on $\sigma$ and
${\cal H}$).

\begin{proof}
Without loss of generality, one may assume that
$\|f\|\le 1$ for all $f\in {\cal H}.$
{{Since}} $T(A)$ is {weak*-compact},
{{there}} are $\tau_1, \tau_2,...,\tau_m\in T(A)$ such that, for any
$\tau\in T(A),$ there is $j\in\{1,2,...,m\}$ such that
\begin{equation}\label{12N1-1}
|\tau(f)-\tau_j(f)|<\sigma/4\tforal f\in {\cal H}.
\end{equation}
By the Krein-Milman Theorem, there are $t_1', t_2',...,t_n'\in \partial_e(T(A))$
and nonnegative numbers $\{\af_{i,j}\}$  such that, {{for each $j=1,2,...,m,$}}
\begin{equation}\label{12N1-2}
|\tau_j(f)-\sum_{i=1}^n \af_{i,j} t_i'(f)|<\sigma/8\andeqn\sum_{i=1}^n\af_{i,j}=1.
\end{equation}
Put $E=\{t_1',t_2',...,t_n'\}.$
Choose $N>32mn/\sigma.$ Let $\tau$ be a possibly unnormalized trace on $A$ with $0<\tau(1)\leq 1$. Suppose that $j$ is  chosen
so that $|\tau(f)/\tau(1)- \tau_j(f)|<\sigma/4 \tforal f\in {\cal H}$ as in \eqref{12N1-1}.
Then,  for any $k\ge N,$ there exist positive rational numbers $r_{i,j}$  and positive integers $p_{i,j}$ ($1\le i\le n$ and $1\le j\le m$)  such that
\begin{equation}\label{ratn-001}
{ {\sum_{i=1}^nr_{i,j}\leq 1, ~~{\mbox {or}}}}~~ \sum_{i=1}^nr_{i,j}=1~~ \mbox{if} ~~\tau(1)=1, ~{\blue{\mbox{for } ~1\leq j\leq m}}
\end{equation}
\begin{equation}\label{ratn-002}
r_{i,j}=\frac{p_{i,j}}{k},\quad\textrm{and}\quad |{{\tau(1)}}\alpha_{i,j}-r_{i,j}|<\frac{\sigma}{8n},\quad {\blue{\mbox{for}}}~1\le i\le n{\blue{,~\mbox{and}~1\leq j\leq m}}.\,\,
\end{equation}
Note that, {{by \eqref{ratn-001},}}
\begin{equation}\label{12N1-3}
{ {\sum_{i=1}^np_{i,j}\leq k,~~{\mbox {or}}}}~~ \sum_{i=1}^np_{i,j}=k~~{ {{\mbox {if}}~~ \tau(1)=1}}.
\end{equation}



Then, by (\ref{ratn-002}),
\beq\label{16nn}
|\tau_j(f)-\sum_{i=1}^n({p_{i,j}\over{k}})t_i'(f)|<\sigma/4+\sigma/8=3\sigma/8\rforal f\in {\cal H}.
\eneq
It is then clear that (\ref{161N-1}) holds on repeating each $t_i'$ $p_{i,j}$ times.

Now suppose that $A$ is subhomogeneous.
{{By Lemma 2.16 of \cite{Lncrell},}}
$t_i'$ has the form ${\rm tr}_i\circ \pi_i,$ where $\pi_1, \pi_2,..., \pi_n$ are in ${\rm Irr}(A).$
It follows that (\ref{161N-2}) holds.

{{Note that $\hat{A}_l$ is a locally compact Hausdorff space (see Proposition 4.4.10 of \cite{Pdbook},  for example).}}
{{Fix a metric on $\hat{A}_l.$ {{For each $\pi_i\in \hat{A}_l,$}}
there}} exists $\dt_i>0$ such that for any irreducible representation $x\in \hat{A}_l $ with ${\rm dist}(x, \pi_i)<\dt_i$, we have
\begin{equation*}
|{\rm tr}(x(f))-{\rm tr}(\pi_i(f))|<\sigma/64k\,\,\,
{\blue{\rforal f\in {\cal H},}}
\end{equation*}
{{where ${\rm tr}$ is the unique trace of $M_l.$}}

{{Suppose that, for each $l,$}} $\hat{A}_{{l}}$ has no isolated points.
{{Fix $j$ and,}}  for each $i$ {\blue{with ${\rm tr}_i\circ \pi_i\in E$,}} choose $p_{i,j}$ distinct points
 in a neighborhood $O({\blue{\pi_i}}, \dt_i)$ of ${\blue{\pi_i}}$ (in $\hat{A}_{{l}}$) with diameter
less than $\dt_i.$
Let $\{\pi_{1,j},\pi_{2,j},..., \pi_{k,j}\}$ be the resulting set of $k$ elements (see (\ref{12N1-3})).
Then, one has  (by \eqref{12N1-1} and \eqref{16nn})
\begin{equation}\label{12n1-5}
|\tau(f)-{1\over{k}}({\rm tr}_{1,j}(f(\pi_{1,j}))+{\rm tr}_{2,j}(f(\pi_{2,j}))+\cdots+ {\rm tr}_{k,j}(f(\pi_{k,j})))|<\sigma\rforal f\in {\cal H}
\end{equation}
{{(where ${\rm tr}_{i,j}$ is the tracial state of $\pi_{i,j}(A)$),}} as desired.
\end{proof}

{{
\begin{lem}\label{discretizeT}
Let $\mathcal H$ be a finite subset of $C([0, 1]\times\mathbb T)\otimes M_r$
(for some $r\ge 1$) and let $\sigma>0.$
Then there is an  integer  $N\ge 1$
such that for any finite  Borel   measure $\mu$ on $[a, b]\times \mathbb T$ with $\|\mu\|\le 1$
and any $k\ge N,$  there are $x_1, x_2, ..., x_{m}$ in $(0, 1)\times \mathbb T$
for some $m\le N$ such that
$$|\int_{(a,b)\times \T} {\rm tr}(h) d\mu-\frac{1}{k}({\rm tr}(h(x_1))+{\rm tr}(h(x_2))+\cdots+{\rm tr}(h(x_m)))|<\sigma
\tforal h\in\mathcal H,$$
where $[a,b]\subset [0,1]$ and ${\rm tr}$ is the tracial state of $M_r.$
\end{lem}}}
{{\begin{proof}
Let $A=C([0,1]\times \T)\otimes M_r.$
Note that, for each  $\mu$  {{as specified,}}
$\tau(f)=\int_{(a,b)\times \T} {\rm tr}(f)d\mu$ is a trace on $A$ with $\|\tau\|\le 1.$
Therefore,
the conclusion follows immediately from Lemma \ref{discretizeT0} if one allows $x_1, x_2,...,x_m$ {{to be}} in  $[0,1]\times \T.$
The equicontinuity of ${\cal H}$ {{then}}
allows us to require these points to be in $(0,1)\times \T.$
\end{proof}
}

The following fact is {{known to experts.}}

\begin{lem}\label{oldext}
Let $C=\bigoplus_{i=1}^kC(X_i)\otimes M_{r(i)},$  where each $X_i$ is a connected compact metric space. 
Let $\mathcal H\subset C$ be a finite subset and let $\sigma>0$. Then there is
an integer $N\ge 1$ satisfying the following {{condition}}:
for any positive homomorphism $\kappa: K_0(C)\to K_0(M_s)=\Z$   with $\kappa({{[1_{C(X_i)\otimes M_{r(i)}}]}})\ge N$
{{(for each $i$)}}
and any $\tau\in T(C)$ such that
\beq\label{16-n181}
{\rho_C(x)(\tau)}:={{(1/s)}}(\kappa(x))\tforal x\in K_0(C),
\eneq
{\blue{there}} is a homomorphism $\phi: {{C}}\to M_s$ such that
$\phi_{*_0}=\kappa$ and
$$|{\rm tr}\circ\phi(h)-\tau(h)|<\sigma \tforal h\in\mathcal H{\blue{,}}$$
{\blue{where ${\rm tr}$ is the tracial state on $M_s$.}}

\end{lem}	
\begin{proof}
{{ For  convenience, we  present a proof, using  Lemma \ref{discretizeT0}.
It is clear that the general case can be reduced to the case that
$C=C(X)\otimes M_r$ for some connected compact metric space $X$ and $r\ge 1.$
Let $\sigma$ and ${\cal H}$ be given.  \Wlog, we may assume that ${\cal H}$ is in the unit ball of $C.$
Let $N_1$ be as in \ref{discretizeT0} for $\sigma$ and ${\cal H}.$
Let $N=N_1r.$}}

{{Suppose now   $\tau$ and $\kappa$ are given such that $\kappa([1_C])\ge N$  {{and \eqref{16-n181} holds.}}
Let $D: K_0(C)\to \Z$ be defined by {{the}} rank function and let $e\in C$ be a projection of (constant)  rank one.
Then $\kappa([e])\ge N_1.$
The assumption \eqref{16-n181} means
that $\rho_C(x)(\tau)=(1/s)(\kappa([e])D(x))$ for all $x\in K_0(C).$
In particular, $1=(1/s)(\kappa([e])r)=(1/s)\kappa([1_C])$ and $s=\kappa([e])r.$
Let $k=\kappa([e])\ge N_1$ and $\pi_1, \pi_2,...,\pi_k$ be given by (the second part of) Lemma \ref{discretizeT0}.
Define $\phi(f)=\diag(\pi_1(f), \pi_2(f),...,\pi_k(f)).$
Then, by the choice of $N_1$ and {{by}} Lemma \ref{discretizeT0}, one has
$$
|{\rm {tr}}(\phi(h))-\tau(h)|<1/\sigma\tforal h\in {\cal H}.
$$
Moreover, $\phi_{0*}=\kappa.$}}

\end{proof}

\begin{lem}\label{appoldext}
Let $C=C(\T)\otimes F_1,$ where
$F_1=M_{R(1)}\oplus M_{R(2)}\oplus \cdots \oplus M_{R(l)},$ or $C=F_1.$
Let $\mathcal H\subset C$ be a finite subset, and let $\epsilon>0$.
There is $\dt>0$ satisfying the following condition:
For any $M_s$, any positive order-unit {{preserving}} map $\kappa: K_0(C)\to K_0(M_s),$
 and any tracial state $\tau\in T(C)$ such that
$$
 |{\rho_{ {M_s}}(\kappa(p))(\mathrm{tr})}-\tau(p)|<\dt
 $$
 for all projections {{$p$}} in $C,$
where $\mathrm{tr}$ is the tracial state  on $M_s$, there is a tracial state $\tilde{\tau}\in T(C)$ such that
$${{(1/s)}}(\kappa([p]))=\tilde{\tau}(p){{\tforal p\in C}}$$ and
$$|\tau(h)-\tilde{\tau}(h)|<\epsilon\rforal h\in \mathcal H.$$

\end{lem}

\begin{proof}
We may assume that ${\cal H}$ is in the unit ball of $C.$
Let $\dt=\ep/l.$
We may write  $\tau=\sum_{j=1}^l \af_j \tau_j,$ where $\tau_j$ is a tracial state
on $C(\T)\otimes M_{r(j)},$ $\af_j\in \R_+,$ and $\sum_{j=1}^k\af_j=1.$
Let $\bt_j={{(1/s)}}(\kappa([1_{{C(\T)\otimes M_{R(j)}}}]),$ $j=1,2,...,l.$
Put ${\tilde \tau}=\sum_{j=1}^l \bt_j\tau_j.$
Then
${\rm tr}(\kappa(p))={\tilde \tau}(p)$
for all projections $p\in C;$ and for any $h\in {\cal H},$
$$
|{\tilde \tau}(h)-\tau(h)|\le \sum_{j=1}^l|\bt_j-\af_j|<\ep,
$$
as desired.
\end{proof}

{\begin{rem}\label{botker}

{{(1) Let $C$ be a separable stably finite \CA\, and let $A=C\otimes C(\T)=C(\T, C).$
Recall that  $K_0(A)=K_0(C)\oplus \boldsymbol{\bt}(K_1(C))\cong K_0(C)\oplus K_1(C)$ (see \ref{Dbeta}).
If $K_i(C)$ is finitely generated, then  $K_i(A)$ is also finitely generated, $i=0,1$ (see \ref{2Lg13}).
Fix a point $t_0\in \T.$  Let $\pi_{t_0}: A\to C$ be the point evaluation  at $t_0$ defined by
$\pi_{t_0}(f)=f(t_0)$ for all $f\in C(\T, C).$ Define $\iota: C\to C(\T, C)=A$
by $\iota(c)(t)=c$ for all $t\in \T$ and $c\in C.$ Then $\pi_{t_0}\circ \iota={\rm id}_C.$
Thus, the  \hm\, $(\pi_{t_0})_{*0}: K_0(A)\to K_0(C)$ induced by $\pi_{t_0}$ maps
 $K_0(A)$ onto $K_0(C)$ and $((\pi_{t_0})_{*0})|_{K_0(C)}$ is an order
isomorphism.  In particular, we may write  ${\rm ker}(\pi_{t_0})_{*0}= \boldsymbol{\bt}(K_1(C)).$
In other words, if $p(t), q(t)\in M_N(C(\T, C))$ are projections and $[p]-[q]\in \boldsymbol{\bt}(K_1(C)),$
then $[p(t)]=[q(t)]$ in $K_0(C)$ for all $t\in \T.$
Note that, any tracial state $\tau\in T(A)$ may be written as
$\tau(f)=\int_{\T} s(f)(t)d\mu(t),$ where $s\in T(C)$ and $\mu$ is a probability Borel measure
on $\T.$  It follows that $\tau(p)-\tau(q)=0$ if $[p]-[q]\in {\rm ker}(\pi_{t_0})_{*0}$ for all $\tau\in T(A).$
This, in particular, implies that $\boldsymbol{\bt}(K_1(C))\subset {\rm ker}\rho_A.$}}

{{(2) Let $C=A(F_1, F_2, \phi_0, \phi_1)$ and
let $\pi_e: C\to F_1$ be as in
\ref{DfC1} which gives an order embedding
$(\pi_e)_{*0}: K_0(C)\to K_0(F_1)$ (see  \ref{2Lg13}).
Let $\phi_i': C(\T)\otimes F_1\to C(\T)\otimes F_2$ be defined
by $\phi_i'(f\otimes d)=f\otimes \phi_i(d)$ for all $f\in C(\T)$ and $d\in F_1,$ where $i=0,1.$}}
{{Then {{($A=C\otimes C(\T)$)}}
$$
A=\{(f,g)\in C([0,1]\times \T, F_2)\oplus (C(\T)\otimes F_1): f(0,t)=\phi_0'(g(t))\andeqn f(1,t)=\phi_1'(g(t))\}.
$$
Let $\pi_e': A\to C(\T, F_1)$ be defined by
$\pi_e'(c\otimes f)=\pi_e(c)\otimes f$ for all $c\in  C$ and $f\in C(\T).$
Note {{that}}, by \ref{2Lg13}, ${\rm ker}\rho_C=\{0\}.$ Then, by (1) above, ${\rm ker}\rho_A=\boldsymbol{\bt}(K_1(A)).$
One then verifies easily that {{the map}} $(\pi_e')_{*0}: K_0(A)=K_0(C)\otimes C(\T)=K_0(C)\oplus \boldsymbol{\bt}(K_1(C))\to F_1\otimes C(\T)
=K_0(F_1)$
is given by $(\pi_e')_{*0}(x\oplus y)=(\pi_e)_{*0}(x)$ for all $x\in K_0(A)$ and $y\in \boldsymbol{\bt}(K_1(C)).$ In other words,
$K_0(A)/{\rm ker}\rho_A$ is embedded into $K_0(C(\T, F_1).$
}}

\end{rem}
}}

\begin{lem}\label{ExtTraceMn}
Let $A=C$ for some $C\in\mathcal C$ or
$A=C\otimes C(\T)$ for some $C\in\mathcal C$.
Let $\Delta: A_+^{q, 1}\setminus\{0\}\to(0, 1)$ be an order preserving map. Let $\mathcal H\subset A$ be a finite subset, and let $\sigma>0$. Then there are a finite subset $\mathcal H_1\subset A_+$ and a positive integer $K$ such that for any  $\tau\in T(A)$ satisfying
\begin{equation}\label{extnn-0}
\tau(h) > \Delta(\hat{h}) \tforal h\in \mathcal H_{1}
\end{equation}
 and any positive homomorphism $\kappa: K_0(A)\to K_0(M_s)$
with $s=\kappa([1_A]){\blue{\in \N}}$
satisfying
\begin{equation}\label{extnn-0+}
\tau(x)=(1/s)(\kappa(x))\tforal x\in K_0(A),
\end{equation}
there is a unital homomorphism $\phi: A\to M_{sK}$ such that
$\phi_{*0}=K\kappa$ and $$|{\rm tr}'\circ\phi(h)-\tau(h)|<\sigma \tforal h\in\mathcal H,$$
where ${\rm tr}'$ is the tracial state on $M_{sK}.$
\end{lem}

\begin{proof}

We will {{consider}} the case that $A=C(\T)\otimes C.$ The case $A=C$ can be proved
in the same manner but simpler.

Let $C=\{(f,g)\in C([0,1], F_2)\otimes F_1: f(0)=\phi_0(g)\andeqn f(1)=\phi_1(g)\}.$
Write $F_1=M_{R_1}\oplus M_{R_2}\oplus \cdots \oplus M_{R_l}$ and
  $F_2=M_{r(1)}\oplus \cdots M_{r(k)}$ and $C([0,1]\times \T, F_2)=\bigoplus_{i=1}^k C([0,1]\times \T, M_{r(j)}).$
{{Then,  keeping the notation of part (2) of \ref{botker}, we write}}
\beq\nonumber
A=\{(f,g)\in C([0,1]\times \T, F_2)\oplus (C(\T)\otimes F_1): f(0,t)=\phi_0'(g(t))\andeqn f(1,t)=\phi_1'(g(t))\}.
\eneq
\noindent
{{In particular, $A\in {\cal D}_1.$}}
Let $\pi_i:  C(\T)\otimes F_1\to C(\T)\otimes M_{R_i}$ be the projection to the $i$th summand
  and $\ep_j: C([0,1]\times \T)\otimes F_2\to C([0,1]\times \T)\otimes M_{r(j)}$ the projection
  to the $j$-th summand, $1\le j\le k.$
  Set $\phi_{0,j}=\pi_j\circ \phi_0: F_1\to M_{r(j)}$
  and
  $\phi_{1,j}=\pi_j\circ \phi_1: F_1\to M_{r(i)}.$
  Let $\phi_{0,j}': C(\T)\otimes F_1\to C(\T)\otimes M_{r(j)}$ and
  $\phi_{1,j}' :C(\T)\otimes F_1\to C(\T)\otimes M_{r(j)}$ be the  \hm s induced
by $\phi_{0,j}$ and $\phi_{1,j},$ respectively.
   Let $\pi_e: C\to F_1$
   and $\pi_e': A\to C(\T)\otimes F_1$ be as in part (2) of \ref{botker}.

  Let $h=(h_f, h_g)\in A.$
In what follows,
{\blue{for}} $\xi\in (0,1)\times \T,$ we will use $\ep_i(h)(\xi)$ for $\ep_i(h_f)(\xi),$  and, we {{will}} use
$\ep_i(h)((0,t))$ for $\phi_{0,i}'(\pi_e'(h))(t)$ (for $t\in \T$) and $\ep_i(h)((1,t))$ for $\phi_{1,i}'(\pi_e'(h))(t)$
(for $t\in \T$), whenever it is convenient.

{{By}} \ref{FG-Ratn},  $K_0(C)$ is finitely generated by minimal projections
in $M_m(C).$ Replacing $C$ by $M_m(C),$ \wilog, we may assume that $K_0(A)$ is generated
by $\{p_1, p_2,...,p_c\},$ where $p_i\in C$ are minimal projections, $i=1,2,...,c.$   In what follows we will
identify $p_i$ with $p_i\otimes 1_{C(\T)}$ whenever it is convenient.

Note that $K_0(A)=K_0(C)\oplus \boldsymbol{\bt}(K_1(C))\cong K_0(C)\oplus K_1(C),$
  ${\rm ker}\rho_C=\{0\},$ and ${\rm ker}\rho_A=\boldsymbol{\bt}(K_1(C))$ (see \ref{botker}).
    Therefore, $\kappa|_{{\rm ker}\rho_A}=\{0\}.$  Let $\kappa_{00}: K_0(C)\to K_0(M_s)$ be the positive
  \hm\, induced by
  $\kappa.$

Note also $K_0(C(\T)\otimes F_1)\cong K_0(F_1)=\Z^l.$
Let $e_i$ be a  minimal projection of $M_{R_i},$ $i=1,2,...,l.$
Let $I={\rm ker}\pi_e'.$
Since $\pi_e$ is surjective (see  \ref{DfC1}), there are $h_i\in A_+$ such that $\|h_i\|\le 1$ and
$\pi_e'(h_i)=e_i,$ $i=1,2,...,l.$
We may  assume that
$\ep_i(h_j)(r,t)=\phi_{0,i}'(e_j)$ if $r\in [0,1/4]$ and $\ep_i(h_j)(r,t)=\phi_{1,i}'(e_j)$
if $r\in [3/4, 1].$

We may also assume that ${\cal H}$ is a subset of the unit ball of $A$ which contains $1_A$ as well as $\{p_1, p_2,...,p_c\}.$
Let $\overline{{\cal H}}=\{\pi_e'(h): h\in {\cal H}\}.$

Let $N_0$ (in place of $N$) be the integer for $\pi_e'(A)$ (in place of $C$), $\overline{{\cal H}}$
(in place of ${\cal H}$) and $\sigma/64$ (in place of $\sigma$) given by  {{Lemma}} \ref{oldext}. {{Let us  fix the}} metric $d(r\times t, r'\times t')=\sqrt{|r-r'|^2+|t-t'|^2}$ for all $r\times t, r'\times t'\in [0,1]\times \T.$
There exists $1/4>\dt_0>0$ such that, if ${\rm dist}(\xi, \xi')<\dt_0$ ($\xi, \xi'\in [0,1]\times \T$), or $|t-t'|<\dt_0$
($t,t'\in \T$),  or $0<r<\dt_0,$ then, for any  $h=(h_f, h_g)\in {\cal H}$, one has
\begin{equation}\label{extnn-1}
\|h_f(\xi)-h_f(\xi')\|<\sigma/64kN_0l,\quad \|h_g(t')-h_g(t)\| <\sigma/64N_0kl,
\end{equation}
and, for all $t\in \T,$
\begin{equation}\label{extnn-1+}
\|h_f(r, t)-\phi_0'(h_g)(t)\|<\sigma/64kN_0l\andeqn \|h_f(1-r,t)-\phi_1'(h_g)(t)\|<\sigma/64kN_0l.
\end{equation}

Choose ${{a_I}}\in I_+$ such that  $\|{{a_I}}\|\le 1,$ $a(r,t)=1_{F_2}$ and $a(1-r,t)=1_{F_2}$ if $1>r>\dt_0/2,$ and
${{a_I}}(r,t)={{a_I}}(1-r,t)=0$ if $0<r<\dt_0/4,$ for all $t\in \T.$

Now we choose ${\cal H}_1.$ For each $1\le j\le k,$ find a $g_j\in (C_0(0,1)\otimes \T \otimes M_{r(j)})_+\setminus\{0\}$ such that $g_j(r,t)=0$ if $r\not\in (\dt_0/2, 1-\dt_0/2)$ and $\|g_j\|\le 1.$
Find $g_j', g_j''\in (C_0(0,1)\otimes \T \otimes M_{r(j)})_+\setminus \{0\}$
such that $\|g_j'\|, \|g_j''\|\le 1,$
$g_j'(r,t)=0$ if $r\not\in (0,\dt_0/2)$ and $g_j''(r, t)=0$ {{if}}   $r\not\in (1-\dt_0/4, 1),$ $j=1,2,...,k.$
We {{will}} also use $g_j$ for the elements of  $I\subset A$ such that $\ep_j(g_j)(\xi)=g_j(\xi)$
for all $\xi\in (0,1)\times \T,$ $\ep_i(g_j)=0$ for all $i\not=j.$
We may also view $g_j'$ and $g_j''$ as elements of $I$ in exactly the same manner.
Let $h_i'=(1-{{a_I}})h_i(1-{{a_I}}),$ $i=1,2,...,l.$  Note that $h_i'(r, t)=h_i(r,t)$ if $r\in [0, \dt_0/4]\cup [1-\dt_0/4, 1]$
and $h_i'(r, t)=0$ if $r\in [\dt_0/2, 1-\dt_0/2].$

Recall that we  identify $p_i$ with $p_i\otimes 1_{C(\T)},$
$i=1,2,...,c.$
Put
$$
{\cal P}'=\{g_i'p_j: g_i'p_j\not=0, 1\le i\le k, 1\le j\le c\}\cup \{g_i''p_j'': g_i''p_j''\not=0,q\le i\le k, 1\le j\le c\}
$$
and put
$$
{\cal H}_1=\{1_A\}\cup \{h_i, h_i': 1\le i\le l\}\cup\{g_j, g_j', g_j'': 1\le j\le k\}\cup {\cal P}'.
$$
Put
\begin{equation}\label{extnn-2}
\sigma_1=\min\{\Delta(\hat{h}): h\in {\cal H}_1\}/2\quad\mathrm{and}\quad \sigma_2=\sigma_1\cdot  \sigma/64kl.
\end{equation}
Let $M$ be the integer in Lemma \ref{multiple-ext} associated with the pair
$K_0(A)$ (as $G$) and $\Z^l$ (see \ref{2Lg13}).
Let $K_1$ (in place of $R$) be the integer provided by \ref{multiple-ext} for $G=K_0(A),$
$\sigma_1$ and $\sigma_2.$
Let $N_j$  be the integer provided by {{Lemma}} \ref{discretizeT} for $\sigma_2/8k\prod_{j=1}^kr(j)$
(in place of $\sigma$) for $r=r(j),$ $j=1,2,...,l.$
Let ${\bar N}=\max\{N_j: 1\le j\le k\}.$
Put $d_0=\prod_{j=1}^k r(j).$
Let $K=K_1\cdot N_0\cdot {\bar N}\cdot M\cdot d_0.$

Now suppose that  $\kappa$ and $\tau\in T(A)$ are  given  satisfying
(\ref{extnn-0}) and (\ref{extnn-0+}).
We may write (see 2.14 of \cite{Lncrell})
\begin{equation}\label{extnn-trace}
\tau(a)=(\sum_{i=1}^k\int_{(0,1)\times \T} {\rm tr}_i(\ep_i(a)(\xi))d\mu_i(\xi))+ t_{{e}}\circ \pi_e'(a)\rforal  a\in A,
\end{equation}
where $\mu_i$ is a Borel measure on $(0,1)\times \T$ {{with $\|\mu_i\|\le 1,$}}
${\rm tr}_i$ is the normalized trace on $M_{r(i)},$
and
$t_e$ is a trace (with $\|t\|\le 1$) on $C(\T)\otimes F_1.$
Consider the finite set $\{\mu_1, \mu_2,...,\mu_k\}.$
By
\% applying
{{Lemma}} \ref{Lmeasdiv}, one can find two points $\dt_0''<\dt_0'$  in $(15\dt_0/16\dt_0, \dt_0]$ such that
\beq\label{16-18n-nn81}
\int_{[\dt_0'', \dt_0']\times \T} d\mu_i<\sigma_2/8k,\,\,\,i=1,2,...,k.
\eneq

{{Note, for each $j,$  $\ep_i(p_j)(\xi)$ is  constant on $(0,1)\times \T$ for each $i.$
If $\ep_i(p_j)\not=0,$ then $g_i'p_j\in {\cal P}'.$ We have  (as $0<\dt_0/2<\dt_0'\le \dt_0$ and $\ep_s(g_j')=0,$
if $s\not=j$)
\beq\label{16-18-n1-1}
\int_{(0,\dt_0')\times \T} {\rm tr}_i(\ep_i(p_{j}))d\mu_i\ge \int_{(0,1)\times \T} {\rm tr}_i(\ep_i(g_i'p_j))d \mu_i\\\label{16-18-n1-1+}
= \sum_{s=1}^k\int_{(0,1)\times \T} {\rm tr}_s(\ep_s(g_i'p_j))d \mu_s=\tau(g_i'p_j)\ge 2\sigma_1.
\eneq
}}
It follows from {{Lemma}} \ref{discretizeT} that, for each $i,$  there are
$t_{i,j}\in (0,1)\times \T,$ $j=1,2,...,m(i)\le {\bar N},$ such that
\beq\label{extnn-3}
|\int_{[\dt_0',1-\dt_0']\times \T}{\rm tr}_i(f)d\mu_i-(1/{\bar N})\sum_{j=1}^{m(i)} {\rm tr}_i(f(t_{i,j}))|<\sigma_2/8k
\eneq
for all $f\in {\cal H}.$
For each $i,$ define $\rho_i, \rho_i': A\to \C$ by
\beq\label{extnn-3+}
\rho_i(f)&=&\int_{(0,1)\times \T} {\rm tr}_i(\ep_i(f))d\mu_i-(1/{\bar N})\sum_{j=1}^{m(i)} {\rm tr}_i(\ep_i(f(t_{i,j}))),\\
\rho'_i(f)&=& \int_{(0,1-\dt_0']\times \T} {\rm tr}_i(\ep_i(f))d\mu_i-(1/{\bar N})\sum_{j=1}^{m(i)} {\rm tr}_i(\ep_i(f(t_{i,j})))
\eneq
for all $f\in A.$ Then,  since $\ep_i(p_j)(\xi)$ is constant on $[0,1]\times \T,$
\beq
\rho_i(p_j)&=&\int_{(0,1)\times \T} {\rm tr}_i(\ep_i(p_j))d\mu_i-(1/{\bar N})\sum_{s=1}^{m(i)} {\rm tr}_i(\ep_i(p_j(t_{i,s})))\\
\label{16-18n-n10}
&=& {\rm tr}_i(\ep_i(p_j))\rho_i(1_A).
\eneq
If $\ep_i(p_j)\not=0,$  {{by \eqref{extnn-3} and  \eqref{16-18-n1-1+},}}
\begin{eqnarray*}
\rho_i'(p_j)
&=&\int_{(0,\dt_0')\times \T} {\rm tr}_i(p_j)d\mu_i+\int_{[\dt_0', 1-\dt_0']\times \T} {\rm tr}_i(p_j)d\mu_i(t)-
(1/{\bar N})\sum_{s=1}^{m(i)} {\rm tr}_i(p_j(t_{i,s}))\\
&>&\int_{(0,\dt_0')\times \T} {\rm tr}_i(p_{j})d\mu_i-\sigma_2/2k\ge 2\sigma_1-\sigma_2/2k>0.
\end{eqnarray*}
Put
$\af_i'=\rho_i'(1_A),$ $i=1,2,..., k.$
Let $\nu_{0,i}$ and $\nu_{1,i}$ be the Borel measures on $\T$ given by
\beq\label{extnn-5+}
\int_{\T}{\rm tr}_i(f(t))d\nu_{0,i}(t)&=&\int_{(0,\dt_0')\times \T} {\rm tr}_i(1_C\otimes f)d\mu_i \andeqn\\\label{extnn-5++}
\int_{\T} {\rm tr}_i(f(t))d\nu_{1,i} &=& \int_{(1-\dt_0',1)\times \T} {\rm tr}_i(1_C\otimes f)d\mu_i
\eneq
for all $ f\in C(\T, M_{r(i)}),$ $1\le i\le k.$
Note that $\|v_{0,i}\|\ge 2\sigma_1,$ {{and by \eqref{extnn-3},
\beq\label{16-181208-n1}
|\rho_1'(1_A)-\|v_{0,i}\||=|\rho_1'(1_A)-\int_\T {\rm tr}_i(1_A)d\mu_i|<\sigma_2/8k,\,\,\,i=1,2,...,k.
\eneq}}
Define
$T_{0,i}, T_{1,i}: A\to \C$ by
\beq\label{extnn-6}
T_{0,i}(a)&=&{\af_i'\over{\|\nu_{0,i}\|}}\int_{\T} {\rm tr}_i\circ \phi'_{0,i}\circ \pi_e'(a)d\nu_{0,i}\andeqn\\
T_{1,i}(a)&=&\int_{\T} {\rm tr}_i\circ \phi_{1,i}'\circ \pi_e'(a)d\nu_{1,i}
\eneq
for all $a\in A.$
Note, for any $h\in A$ and $t\in \T,$
\beq\label{extnn-6n}
\hspace{-0.2in}{\rm tr}_i( (\phi_{0,i}'\circ \pi_e'(h))(t))={\rm tr}_i(\ep_i(h)((0,t)))\andeqn
{\rm tr}_i((\phi_{1,i}\circ \pi_e'(h))(t))={\rm tr}_i(\ep_i(h)((1,t))).
\eneq
Therefore, for $h\in {\cal H},$  by (\ref{extnn-1}), \eqref{extnn-5+}, \eqref{extnn-5++}, (\ref{extnn-1+}), and {{\eqref{16-181208-n1},}}
\begin{eqnarray}\label{extnn-6+}
&&|\int_{((0,\dt_0')\cup(1-\dt_0', 1))\times \T}{\rm tr}_i(\ep_i(h))d\mu_i-(T_{0,i}(h)+T_{1,i}(h))|\nonumber \\
&\le & \frac{2\sigma}{64kl}+|\int_{(0,\dt_0')\times \T}{\rm tr}_i(\ep_i(h)(\dt_0'/2,t))d\mu_i-T_{0,i}(h)| +|\int_{(1-\dt_0',1)\times \T}{\rm tr}_i(\ep_i(h)(1-\dt_0'/2,t))d\mu_i-T_{1,i}(h)| \nonumber \\
&=&|\int_{\T}{\rm tr}_i(\ep_i(h)(\dt_0'/2,t))d\nu_{0,i}-T_{0,i}(h)|+|\int_{\T}{\rm tr}_i(\ep_i(h)(1-\dt_0'/2,t))d\nu_{1,i}-T_{1,i}(h)|+\frac{\sigma}{32kl}\nonumber \\
&\le & |\int_{\T}{\rm tr}_i(\ep_i(h)(0,t))d\nu_{0,i}-T_{0,i}(h)|+|\int_{\T}{\rm tr}_i(\ep_i(h)(1,t))d\nu_{1,i}-T_{1,i}(h)|+\frac{\sigma}{32kl}+\frac{2\sigma}{64kl}\nonumber \\ \label{extnn-6+2}
&<& |1-\frac{\rho_i'(1_A)}{\|\nu_{0,i}\|}|\|\nu_{0,i}\| +0+\frac{\sigma}{16kl}<
{\sigma_2/8k}+\frac{\sigma}{16kl}.
\end{eqnarray}
Let ${\bar h}_j\in C([0,1]\times \T)\otimes F_2$ be
such that $\|{\bar h}_j\|\le 1,$ ${\bar h}_j(r,t)=\phi_0'(e_j)$ for $r\in [0, \dt_0''],$
${\bar h}_j(r, t)=0$ if $r\in [\dt_0', 1-\dt_0'],$ and ${\bar h}_j(r, t)=\phi_1'(e_j)$ if $r\in [1-\dt_0'', 1],$ $j=1,2,...,l.$
Then ${\bar h}_j\ge h_j',$ $j=1,2,...,l.$
Moreover {{(by \eqref{16-181208-n1}  again),}}
\begin{eqnarray}\label{16-18n-001}
&&|\int_{((0,\dt_0')\cup(1-\dt_0', 1))\times \T}{\rm tr}_i(\ep_i({\bar h}_j))d\mu_i-(T_{0,i}({\bar h}_j)+T_{1,i}({\bar h}_j))|\nonumber \\
&\le &\sigma_2/8k+|\int_{(0,\dt_0')\times \T}{\rm tr}_i(\phi_{0,i}'(e_j))d\mu_i-T_{0,i}({\bar h}_j)| +
|\int_{(1-\dt_0',1)\times \T}{\rm tr}_i(\phi_{1,i}'(e_j))d\mu_i-T_{1,i}({\bar h}_j)| \nonumber \\\label{16-18n-002}
&\le & \sigma_2/8k+|1-\frac{\rho_i'(1_A)}{\|\nu_{0,i}\|}|\|\nu_{0,i}\|+0<\sigma_2/8k+\sigma_2/8k=\sigma_2/4k.
\end{eqnarray}

Since ${\rm tr}_i(\ep_i(p_j))$ is constant on each open set $(0,1)\times \T,$
put $L_{i,j}={\rm tr}_i(\ep_i(p_j)).$ Then one checks (using \eqref{16-18n-n10} among other items)
that
\beq\label{extnn-7}
\hspace{-0.4in}T_{0,i}(p_j)+T_{1,i}(p_j)&=&{\af_i'\over{\|\nu_{0,i}\|}}\int_{\T} {\rm tr}_i\circ \phi_{0,i}\circ \pi_e'(p_j)d\nu_{0,i}
+\int_{\T} {\rm tr}_i\circ \phi_{1,i}\circ \pi_e'(p_j)d\nu_{1,i}\\ \label{extnn-7+1}
&=&\af_i'{\rm tr}_i(\ep_i(p_j))+\|\nu_{1,i}\|{\rm tr}_i(\ep_i(p_j))
=L_{i,j} (\rho_i'(1_A)+\int_{\T} d\nu_{1,i})\\\label{extnn-7+2}
&=&L_{i,j}(\rho_i(1_A))=\rho_i(p_j).
\eneq
Let
\begin{eqnarray}\label{extnn-8}
&&\hspace{-0.3in}T_1(a)=t_{{e}}\circ \pi_e'(a)+\sum_{i=1}^k (T_{0,i}(a)+T_{1,i}(a)),  T_2(a)=\sum_{i=1}^k[ (1/{\bar N})\sum_{j=1}^{m(i)} {\rm tr}_i(a(t_{i,j})]\\
&&\hspace{-0.3in}T(a)=T_1(a)+T_2(a) t_e\circ \pi_e'(a)+\sum_{i=1}^k (T_{0,i}(a)+T_{1,i}(a)+(1/{\bar N})\sum_{j=1}^{m(i)} {\rm tr}_i(a(t_{i,j})))
\end{eqnarray}
for all $a\in A.$
Thus, by \eqref{extnn-7+2} {{(see also \eqref{extnn-trace}),}}
\beq\label{16-18n-n11}
T(p_j)=t_e\circ \pi_e'(p_j)+\sum_{i=1}^k (\rho_i(p_j)+(1/{\bar N})\sum_{j=1}^{m(i)} {\rm tr}_i(p_j(t_{i,j})))=\tau(p_j)
\eneq
for all $j.$
Then $T_1$  and $ T_2$ are   traces on $A$ and $T$  is {{a}}  tracial state  on $A.$
Define
\begin{equation*}\label{16-18n-n31}
T_1'(b)=t_e(b)+\sum_{i=1}^k({\af_i'\over{\|\nu_{0,i}\|}}\int_{\T} {\rm tr}_i\circ \phi_{0,i}(b)d\nu_{0,i}+
\int_{\T} {\rm tr}_i\circ \phi_{1,i}(b)d\nu_{1,i})\rforal b\in C(\T)\otimes F_1.
\end{equation*}
In what follows we will also use $T_1'$ for the extension on $A\otimes M_m$ defined by
\beq\nonumber
T_1'(b\otimes x)&=&t_e(b){\rm Tr}_m(x)+
\\\label{16-18n-n21}
&&\hspace{-0.3in}\sum_{i=1}^k({\af_i'\over{\|\nu_{0,i}\|}}\int_{\T} {\rm tr}_i\circ \phi_{0,i}(b){\rm Tr}_m(x)d\nu_{0,i}
+\int_{\T} {\rm tr}_i\circ \phi_{1,i}(b){\rm Tr}_m(x)d\nu_{1,i})
\eneq
for all $b\in C(\T)\otimes F_1$ and $x\in M_m,$ where ${\rm{Tr}}_m: M_m\to \C$ is the non-normalized trace.

By  \eqref{16-18n-n11} and \eqref{extnn-0+}, and by \eqref{extnn-3}
and \eqref{extnn-6+2},
one has
\begin{eqnarray}\label{extnn-9}
&&(1/s)\circ \kappa(p)=T(p)\tforal p\in K_0(A)\andeqn\\\label{extnn-9++}
&&|\tau(h)-T(h)|<\sigma_2/8k+\sigma_2/8k+\sigma/16kl<\sigma/2\tforal h\in {\cal H}.
\end{eqnarray}
Put
$d_j=\prod_{i\not=j}r(i)$
and
$d=d_0(\sum_{i=1}^k m(i)).$  Note  $d_ir(i)=d_0$ for all $1\le i\le k.$
It follows
that
\beq\label{16-18n-n71}
d={\bar N} d_0T_2(1_A).
\eneq
Define
$\Psi: A\to M_d$ by $\Psi(a)=\bigoplus_{i=1}^k(\sum_{j=1}^{m(i)} {\bar \ep}_i(a)(t_{i,j}))$
for all $a\in A,$
where
\beq
{\bar \ep}_i(a)=\diag(\overbrace{\ep_i(a), \ep_i(a),...,\ep_i(a)}^{d_i{\bar N}})=\ep_i(a)\otimes 1_{d_i{\bar N}}\rforal a\in A.
\eneq
Denote by $t_d$ the tracial state of $M_d.$
Then
\beq\label{16-18n-n61}
T_2(1_A)t_d(\Psi(a))=T_2(a)\rforal a\in A.
\eneq
Let $\kappa_0: K_0(A)\to \Z$ be given by $\Psi.$
By {{hypothesis}},
$s\tau(p_i)\in \Z,$ $j=1,2,...,c.$
By {{\eqref{16-18n-n11}}},
\beq
&&\hspace{-0.6in}s{\bar N}d_0T_1(p_j)=s{\bar N}d_0(T(p_j)-T_2(p_j))=s{\bar N}d_0([\tau(p_j)-\sum_{i=1}^k((1/{\bar N})\sum_{j=1}^{m(i)} {\rm tr}_i(p(t_{i,j})))])\\\label{extnn-10}
\hspace{-0.3in}&&=s{\bar N}d_0\tau(p_j)-d_0\sum_{i=1}^k(\sum_{j=1}^{m(i)} {\rm tr}_i(p(t_{i,j})))\in \Z\,\,\,\hspace{0.6in}\, (s\tau(p_j)\in \Z).
\eneq

We have
\beq\nonumber
T_1'(e_j)&=&t_e(e_j)+\sum_{i=1}^k({\af_i'\over{\|\nu_{0,i}\|}}\int_{\T} {\rm tr}_i\circ \phi_{0,i}(e_j)d\nu_{0,i}+
\int_{\T} {\rm tr}_i\circ \phi_{1,i}(e_j)d\nu_{1,i})\\\nonumber
&=&t_e\circ \pi_e'({\bar h}_j')+\sum_{i=1}^k (T_{0,i}({\bar h}_j')+T_{1,i}({\bar h}_j'))\\\nonumber
&\ge & t_e\circ \pi_e'({\bar h}_j')+\sum_{i=1}^k(\int_{((0,\dt_0')\cup(1-\dt_0', 1))\times \T}{\rm tr}_i(\ep_i({\bar h}_j'))d\mu_i-\sigma_2/4k)
\hspace{0.4in} ({\rm{by}}\,\,\, \eqref{16-18n-002}) \\\nonumber
&=& t_e\circ \pi_e'({\bar h}_j')+\sum_{i=1}^k\int_{(0,1)\times \T}{\rm tr}_i(\ep_i({\bar h}_j'))d\mu_i-\sigma_2/4
\hspace{0.4in} (\ep_i({\bar h}_j)(\xi)=0,\,\,{\rm{for}}\,\, \xi\in [\dt_0', 1-\dt_0'])\\\nonumber
&=& \tau({\bar h}_j')-\sigma_2/4\ge \tau(h_j')-\sigma_2/4\\\label{16-18n-n41}
&\ge &  \Delta(h_j')-\sigma_2/4\ge 2\sigma_1-\sigma_2/4\ge \sigma_1\hspace{1in} ({\rm{by}}\,\,\eqref{extnn-0}).
\eneq

Let $\kappa_1: K_0(A)\to \Z$ be defined by $\kappa_1|_{{\rm ker}\rho_A}=0$ and by $\kappa_1([p_j])=s{\bar N}d_0T_1(p_j),$ $j=1,2,...,c.$
Let $\kappa_1': \Z^l\to \R$ {{be}} defined by $s{\bar N}d_0T_1'$ (see \eqref{16-18n-n21}).
As in  (2) of \ref{botker} (see also \ref{2Lg13} and (1) of \ref{botker}), we may view ${{K_0(A)/{\rm ker}\rho_A}}$
{{as}} $K_0(C)\subset \Z^l.$
Then  $\kappa_1'\circ (\pi_e')_{*0}=\kappa_1.$  In particular,
$\kappa_1'([\pi_e'(1_A)])=\kappa_1([1_A])=s{\bar N}d_0T_1([1_A]).$
Note that,  by  \eqref{16-18n-n41},
\beq\label{extnn-11}
\kappa_1'([e_j])=s{\bar N} d_0T_1'(e_j)\ge
s{\bar N}d_0\sigma_1\ge\sigma_1, \,\,\, j=1,2,...,l.
\eneq
By the choice of $K_1$  and $M,$ and by applying \ref{multiple-ext},
there is {{an order \hm\,}} $\kappa_2: K_0(C(\T)\otimes F_1)=\Z^l\to \Z$
such that
\begin{equation}\label{extnn-12}
\kappa_2|_{K_0(A)}=K_1M\kappa_1\quad\mathrm{and}\quad |\kappa_1'(e_j)-(1/K_1M)\kappa_2(e_j)|<\sigma_2.
\end{equation}
Write $T_1'=\sum_{j=1}^l \af_jt_j\circ \pi_j',$ where each $t_j$ is a tracial state on $C(\T)\otimes M_{R(j)},$
and $\af_j=\kappa_1'(e_j)/s{\bar N}d_0,$ $j=1,2,...,l.$
Write $\bt_j=(1/K_1Ms{\bar N}d_0)\kappa_2(e_j),$ $j=1,2,...,l.$ Then, by \eqref{extnn-12},
\beq\label{extnn-13}
|\af_j-\bt_j|<\sigma_2/s{\bar N}d_0,\,\,\,j=1,2,...,l.
\eneq
Put $T_1''=\sum_{j=1}^l\bt_j\circ \pi_j'$ and $T_1'''=T_1''/\|T_1''\|.$   Note  $T_1'''\in T(C(\T)\otimes F_1)$
and
\beq
\|T_1''\|&=&T''_1(\pi_e'(1_A))=(1/K_1Ms{\bar N}d_0)\kappa_2([\pi_e'(1_A)])\\
&=&(1/s{\bar N}d_0)\kappa_1([1_A])
=T_1([1_A]).
\eneq
We also have
$(1/K_1Ms{\bar N}d_0 T_1([1_A]))\kappa_2([p])=T_1'''(p)$ for all projections  $p$ in $C(\T)\otimes F_1.$
Put $K_2=K_1N_0M{\bar N}d_0T_1([1_A]).$
It follows from {{Lemma}} \ref{oldext} that there is a unital \hm\,
$\Phi: C(\T)\otimes F_1\to M_{sK_2}$ such that
\begin{equation}\label{extnn-14}
\Phi_{*0}=N_0\kappa_2 \quad\mathrm{and}\quad  |{\rm tr}\circ \Phi(h)-T_1'''(h)|<\sigma/64
\end{equation}
for all $h\in \overline{{\cal H}},$ where ${\rm tr}$ is the tracial state
on $M_{sK_2}.$
Recall (see \eqref{extnn-0+} and  \eqref{16-18n-n11}) that
$$
T_1(1_A)+T_2(1_A)=T(1_A)=(1/s)\kappa([1_A])=1.
$$
Thus,  $sK=sN_0{\bar N}K_1Md_0(T_1(1_A)+T_2(1_A)).$
Define $\phi: A\to M_{sK}$ by
\beq
\label{extnn-15}
\phi(a)=\Phi\circ \pi_e'(a)\oplus{\tilde \Psi}(a)\tforal a\in A,
\eneq
where ${\tilde \Psi}$ is the direct sum of $sN_0K_1M$ copies of $\Psi$ (recall \eqref{16-18n-n71}).
By  \eqref{extnn-14}, \eqref{extnn-12}, \eqref{16-18n-n61}, \eqref{extnn-8}, and \eqref{extnn-9}, for any projection $p\in A,$ one has
\begin{eqnarray*}
{(\phi)_{*0}(p)\over{sK}}&=&T_1(p) +\sum_{i=1}^k(1/{\bar N})\sum_{j=1}^{m(i)} {\rm tr}_j(p(t_{i,j}))
= T(p)=(1/s)\kappa,
\end{eqnarray*}
and hence
\begin{equation}
\phi_{*0}=K\kappa.
\end{equation}
By \eqref{extnn-14} and  \eqref{extnn-13},
\begin{equation}\label{extnn-18}
|(T_1(1_A){\rm tr}\circ \Phi(\pi_e'(h))-T_1(h)|<(T_1(1_A))(\sigma/64+l\sigma_2/s{\bar N}d_0)<\sigma/2
\end{equation}
for all $h\in {\cal H}.$
Note that  (recall the definition of $K$ and $K_2$, and \eqref{16-18n-n71})
\beq\label{16-181210-n1}
sK_2/sK=T_1(1_A)/(T_1(1_A)+T_2(1_A))=T_1(1_A)\andeqn sN_0K_1Md/sK=T_2(1_A).
\eneq
It follows from  \eqref{extnn-15}, the first part of  \eqref{16-181210-n1}, \eqref{extnn-9++}, \eqref{extnn-18},
the second part of \eqref{16-181210-n1} and \eqref{16-18n-n61} that
\beq
\hspace{-0.3in}|{\rm tr}'\circ \phi(h)-\tau(h)| &\le& |{\rm tr}'\circ \phi(h)-T(h)|+|T(h)-\tau(h)|\\
&<& |T_1(1_A){\rm tr}\circ \Phi(\pi_e'(h))-T_1(h)|+|{\rm tr}'\circ {\tilde \Psi}(h)-T_2(h) |+\sigma/2\\
&<&\sigma/2+|T_2(1_A)t_d(\Psi(h))-T_2(h)|+\sigma/2=\sigma
\eneq
for all $h\in {\cal H}.$
{{(Recall}} ${\rm tr}'$ is the tracial state of $M_{sK}$ {\blue{and ${\rm tr}$ is the tracial state on $M_{sK_2}$}}.)
\end{proof}

\begin{lem}\label{appextnn}
Let {{$A=C\otimes B$ be a unital \CA, where $C\in {\cal C}$ and $B=C(X),$ where $X$ is one point, or $X=\T.$}}
Let $\Delta: A_+^{q, {\bf 1}}\setminus\{0\}\to (0,1)$ be an order preserving map.  Let $\mathcal H\subset A$ be a finite subset and let $\epsilon>0$.
There exist a finite subset $\mathcal H_1\subset A_+^{\bf 1}\setminus \{0\}$  and
a finite subset of projections ${\cal P}\subset  M_n(A)$ (for some $n\ge 1$),
and
there is $\delta>0$ such that if a tracial state $\tau\in T(A)$ satisfies
\beq\label{extnn-0++-1}
\tau(h)>\Delta(\hat{h})\tforal h\in\mathcal H_1
\eneq
and
{{$\kappa: K_0(A)\to K_0(M_s)$ is any order preserving \hm\,  with $\kappa([1_A])=[1_{M_s}]$ satisfying}}
\begin{equation}\label{extnn-0++}
|{{\rho_{M_s}(\kappa([p]))(\mathrm{tr})}}-\tau(p)|<\delta
\end{equation}
for all projections $p\in  {\cal P},$
where $\mathrm{tr}$ is the tracial state on $M_s$, {{then}} there is a tracial state $\tilde{\tau}\in T(A)$ such that
$$
{{\rho_{M_s}(\kappa(x))(\mathrm{tr})}}=\tilde{\tau}({{x}}) \tforal {{x}}\in K_0(A)\tand
|\tau(h)-\tilde{\tau}(h)|<\epsilon\rforal h\in \mathcal H.$$
\end{lem}

\begin{proof}
{{Let
$C=A(F_1, F_2, \phi_1,\phi_2).$ Write $A=C\otimes C(X)=C(X, C),$ where $X$ is a single point or $X=\T.$ }}
{{Note that $K_0(A)=K_0(C)\oplus \boldsymbol{\bt}(K_1(C))\cong K_0(C)\oplus K_1(C),$
${\rm ker}\rho_C=\{0\}$ and ${\rm ker}\rho_A=\boldsymbol{\bt}(K_1(C))$ (see \ref{Dbeta}).}}
Without loss of generality, we may assume that {{the}}
projections in $A$ generate $K_0(A),$ by replacing $A$ by $M_N(A)$
for some integer $N \ge 1.$
{{Let ${\cal P}'$ be a finite subset of projections in $C$ which generates $K_0(C)_+$ (see \ref{FG-Ratn})
which we also assume to contain $[1_A].$
Let $\iota: C\to C(X, C)=A$ be defined by $\iota(c)(t)=c$ for all $t\in \T$ and $c\in C$ and
let ${\cal P}=\{\iota(p): p\in {\cal P}'\}.$ Then, by  {{Remark}} \ref{botker},}}
$\{\rho_A([p]):p\in {\cal P}\}$ generates
${ {\rho_{A}}}(K_0(A))_+.$ \Wlog, we may assume that ${\cal H}\subset A_+^{\bf 1}\setminus \{0\}.$
Let ${\cal H}_1={\cal H}\cup {\cal P}$ and let
$\sigma_0=\min\{\Delta(\hat{h}): h\in {\cal H}_1\}.$
{{\Wlog, we may assume that $0<\ep<1.$}}
Let $\dt=\ep\cdot\sigma_0/128.$

Now suppose that $\tau$ and $\kappa$ {{satisfy}} the {\blue{hypotheses}} for the above mentioned ${\cal H}_1,$
${\cal P},$ and $\dt.$
Note that, for any $s\in S_{[1_{M_s}]}(K_0(M_s)),$
$s\circ \kappa$ is a state on $K_0(A).$  By Corollary 3.4 of \cite{Blatrace},
$s\circ \kappa$ is induced by a tracial state of $A.$ It follows that
{{$\kappa({\rm ker}\rho_A)\subset {\rm ker}\rho_{M_s}=\{0\}.$
In what follows, identifying $K_0(M_s)$ with $\Z,$ we
may view $\kappa$ as an order preserving \hm\, from $K_0(A)$ to $\Z\subset \R.$}}
Define $\eta=(1-\ep/3)(\kappa/s-\tau): K_0(A)\to \R,$ {{where $\tau(x):=\rho_A(x)(\tau)$ for $x\in K_0(A).$}}
Let $d_A: K_0(A)\to K_0(A)/{\rm ker}\rho_A$   be {{the}} quotient map. {{Then $\kappa,$ $\tau,$ and $\eta$
factor through $K_0(A)/{\rm ker}\rho_A.$}}
{{Choose}} $\gamma: K_0(A)/{\rm ker}\rho_A\to \R$   such that
$\gamma\circ d_A=(\ep/3s)\kappa+\eta.$

{{Let $\pi_e': A\to F_1\otimes C(\T)$ be the \hm\,   defined in   Remark  \ref{botker}.
Let $\Psi_0=\pi_{t_0}\circ \pi_e': A\to F_1,$  where $\pi_{t_0}:  F_1\otimes C(\T)\to F_1$
is the point evaluation at $t_0\in \T.$}}
Note that, as in {{\ref{botker} and \ref{2Lg13}, }}
\beq\label{166-1}
K_0(A)/{\rm ker}\rho_A\cong K_0(C).
\eneq
{{As computed in \ref{2Lg13}, $K_0(A)/{\rm ker}\rho_A\cong K_0(C)=(\Psi_0)_{*0}(K_0(A))\subset K_0(F_1).$}}
For each $p\in {\cal P},$ from the assumption (\ref{extnn-0++}), one computes that
\beq\label{166-2}
|\eta([p])|<(1-\ep/3)\dt.
\eneq
Therefore, by \eqref{extnn-0++-1}, \eqref{extnn-0++}, \eqref{166-2},  and the choice of $\dt$ {\blue{and}} $\sigma_0,$
\beq\label{166-3}
\gamma(
d_A([p]))&=&(\ep/3)(\kappa([p])/s)+\eta([p])>
(\ep/3)(\Delta(\hat{p})-\dt)-(1-\ep/3)\dt\\
&\ge & (\ep/3)(1-1/128)\sigma_0-(1-\ep/3)\ep\sigma_0/128\\
&=& \ep\sigma_0[({1\over{3}}-{1\over{3\cdot 128}})-({1\over{128}}-{\ep\over{3\cdot 128}})]>0
\eneq
for all $p\in {\cal P}.$
In other words, $\gamma$ is positive. By 2.8 of \cite{Lnbirr}, there
is {{then}}  a  positive \hm\, $\gamma_1: K_0({{F_1}})\to \R$ such
that $\gamma_1\circ ({{\Psi_0}})_{*0}=\gamma {{\circ d_A}}.$
It is {{well known}} that there is a (non-normalized) trace  $T_0$ on ${{F_1}}$ such
that $\gamma_1([q])=T_0(q)$ for all projections $q\in {{F_1}}.$

Consider the trace $\tau'=(1-\ep/3)\tau+T_0\circ {{\Psi_0}}$ on $A.$
Then, for any projection $p\in A,$
\beq\label{166-4}
\hspace{-0.3in}\tau'(p)&=&(1-\ep/3)\tau(p)+T_0\circ {{\Psi}}_0(p)=(1-\ep/3)\tau(p)+(\ep/3s)\kappa([p])+\eta([p])\\
&=&(1-\ep/3)\tau(p) +(\ep/3s)\kappa([p]) +(1-\ep/3)(\kappa([p])/s-\tau(p))=(1/s)\kappa([p]).
\eneq
Since $(1/s)\kappa([1_A])=1,$ $\tau'\in T(A).$
We also compute that, by  (\ref{166-2}),
\beq\label{166-5}
|T_0\circ {{\Psi}}_0(1_A)|=|\gamma\circ d_A([1_A])|<\ep/2.
\eneq
Therefore, we also have
\beq\label{166-6}
|\tau'(h)-\tau(h)|<\ep\rforal h\in {\cal H}.
\eneq

\end{proof}

\begin{lem}\label{APPextnn}
Let  {{$C=C(X, C_0)$  for some $C_0\in {\cal C},$ where $X$ is a point or $X=\T.$}}
Let $\Delta: A_+^{q, {{\bf1}}}\setminus\{0\}\to(0, 1)$ be an order preserving map. Let $\mathcal H\subset A$ be a finite subset and let $\sigma>0$. Then there are  a finite subset $\mathcal H_1\subset A_+^{\bf 1}\setminus \{0\},$
$\dt>0,$ a finite subset ${\cal P}\subset K_0(A),$ and a positive integer $K$ such that for any  $\tau\in T(A)$ satisfying
\begin{equation*}
\tau(h) > \Delta(\hat{h}) \tforal h\in \mathcal H_{1}
\end{equation*}
 and any { {positive homomorphism} } $\kappa: K_0(A)\to K_0(M_s)=\Z$
with $s=\kappa([1_A])$
such that
\begin{equation*}
|\rho_A(x)(\tau)-(1/s)({\Green{\kappa(x)}})|<\dt
\end{equation*}
for all $x\in {\cal P},$
there is a unital homomorphism $\phi: A\to M_{sK}$ such that
$\phi_{*0}=K\kappa$ and $$|{\rm tr}'\circ\phi(h)-\tau(h)|<\sigma \tforal h\in\mathcal H,$$
where ${\rm tr}'$ is the tracial state on $M_{sK}.$
\end{lem}

\begin{proof}
Note that there is an integer $n\ge 1$ such that {\blue{the}} projections in $M_n(A)$ generate $K_0(A).$
Therefore this lemma  is a corollary of {{Lemma}} \ref{ExtTraceMn} and  {{Lemma}} \ref{appextnn}.
\end{proof}

\begin{lem}\label{ExtTraceI-D}
Let  {{$C=C(X, C_0)$  for some $C_0\in {\cal C},$ where $X$ is a point or $X=\T.$}}
Let $\Delta: C_+^{q, 1}\setminus\{0\}\to(0, 1)$ be an order preserving map. Let $\mathcal F, \mathcal H\subset C$ be  finite subsets, and let $\epsilon>0,\, \sigma>0$.
Then there are a finite subset $\mathcal H_1\subset C_+^{\bf 1}\setminus \{0\}$, $\delta>0$,
a finite subset ${\cal P}\subset K_0(C),$ and a positive integer $K$ such that for any continuous affine map $\gamma: T(C([0, 1]))\to T(C)$ satisfying
$$
\gamma(\tau)(h) > \Delta(\hat{h}) \tforal h\in \mathcal H_{1} \tand \tforal \tau\in T(C([0, 1])),
$$
and  any positive homomorphism $\kappa: K_0(C)\to K_0(M_s(C([0, 1])))$ with $\kappa([1_C])=s$
such that
$$
|\rho_A(x)(\gamma(\tau))-(1/s)\tau(\kappa(x))|<\delta  \tforal \tau\in {\blue{T(C([0, 1])),}}
$$
for all $x\in {\cal P},$
there is an 
$\mathcal F$-$\epsilon$-multiplicative completely positive linear map
$\phi: C\to M_{sK}(C([0, 1]))$ such that
$[\phi]|_{K_0(C)}=K\kappa$ {{(note that $K_0(C)$ is finitely generated; see \ref{botker} and the end of \ref{KLtriple})}} and
$$|\tau\circ\phi(h)-\gamma'(\tau)(h)|<\sigma\tforal h\in\mathcal H,$$
where $\gamma': T(M_{sK}(C([0,1])))\to T(C)$ is induced by $\gamma.$ Furthermore,
{{the map $\phi_0$ can be chosen such that}} $\phi_0=\pi_0\circ \phi$ and $\phi_1=\pi_1\circ \phi$ are  homomorphisms, {{where $\pi_t: M_{sK}(C([0,1])\to M_{sK}$
is the point evaluation at $t\in [0,1].$}}
In the case that $C\in\mathcal C$ (i.e., $X$ is a point), the map $\phi$ can be chosen to be a homomorphism.
\end{lem}

\begin{proof}
Since any
\CA s in $\mathcal C$ are semiprojective  {{(see the line above \ref{2Lg1}),}}
the second part of the statement follows directly from the first part of the statement. Thus, {{we need}} only show the first part of the statement. 
Without loss of generality, {{we}} may assume that  {{${\cal F}$ is in the unit ball of $A,$ $1_A\in {\cal F}$ and}}
${{\mathcal \{ab: a, b\in {\cal F}\}}}\subset \mathcal H$.
To simplify notation, without loss of generality, by replacing  $C$ by $M_r(C)$  for some
$r\ge 1,$ we may assume that the  set of projections in $C$ generates $\rho_C(K_0(C))$
{{(see \ref{FG-Ratn} and \ref{botker}).}}


Since the K-theory of $C$ is finitely generated, there is ${{m'}}\in\mathbb N$ such that, {{for any $x\in  {\rm Tor}(K_i(C))=0,$
$i=0,1,$}}
$${{mx=0\,\,\,{\rm for\,\, some\,\,\, integer}\,\,\, 0<m\le m'.}}$$
{{Put $m_1=(m')!.$}}
Let $\mathcal H_{1, 1}\subset C_+^{1}\setminus \{0\}$ (in place of $\mathcal H_{ {1}}$)
and $\sigma_{1}>0$ (in place of ${{\delta}}$) be the finite subsets and the positive constant {{provided  by}} Theorem \ref{UniqAtoM}
{{for}} $C$ (in place of ${{A}}$),
${{\ep_0=\min\{\sigma/3, \epsilon/3 \}}}$ (in place of $\epsilon$),  $\mathcal H$ ({\blue{in  place}} of $\mathcal F$), and $\Delta/2$.
{{(We will not need the finite set ${\cal P}$ {\blue{of}} Theorem \ref{UniqAtoM}, since $K_0(C)$ is finitely generated and when we apply Theorem \ref{UniqAtoM}, we will require that both maps induce {{the}} same $KL$ {\blue{map}}.)}}

Let $\mathcal H_{1, 2}\subset C$ (in place of $\mathcal H_1$)  be a  finite subset, $\delta>0$
be a positive number, ${\cal P}\subset K_0(C)$ be a finite subset,
and $K'$  be an integer {{as}} {\blue{provided}}  by Lemma \ref{APPextnn} {{for}} $C$, $\Delta/2$ (in place of $\Delta$),  $\mathcal H\cup \mathcal H_{1, 1}$
(in place of $\mathcal H$), and
{{the positive number}} $\min\{\sigma/16, \sigma_1/8,\{\Delta(\hat{h})/4: h\in\mathcal H_{1, 1}\}\}$ (in place of $\sigma$).
{{W}}e may assume that ${\cal P}$ is {{represented by $P,$ a finite}} set of projections in $C$
{{(see Remark \ref{ExtTraceMn}).}}



Put $\mathcal H_1=\mathcal H_{1, 1}\cup\mathcal H_{1, 2}$ and $K={{m_1}}K'$.
Then, let $\gamma: T(C([0, 1]))\to T(C)$ be a continuous affine map with $$\gamma(\tau)(h) > \Delta(\hat{h})\rforal h\in \mathcal H_{1},$$ and let $\kappa: K_0(C)\to K_0(M_s(C([0, 1])))$
with $\kappa([1_C])=s$ be
such that
$$
|\rho_\gamma(\tau)(x)-(1/s)\tau(\kappa(x))|<\delta\rforal  x\in {\cal P}
\rforal \tau\in T(C([0, 1])).
$$

Since $\gamma$ is continuous, there is a partition
$$0=x_0<x_1<\cdots<x_n=1$$
such that, for any $0\leq i\leq n-1$ and any $x\in[x_i, x_{i+1}]$, one has
\begin{equation}\label{eq-par}
|\gamma(\tau_x)(h)-\gamma(\tau_{x_i})(h)|< \min\{\sigma/8, \sigma_1/4\} \rforal h\in \mathcal H_1,
\end{equation}
where $\tau_x\in T(M_s(C([0, 1])))$ is the extremal trace concentrated at $x$.

For any $0\leq i\leq n$, consider the trace $\tilde{\tau}_i=\gamma(\tau_{x_i})\in T(C)$. It is clear that
\beq\nonumber
|\tilde{\tau}_i(x)-\mathrm{tr}({{(\pi_{x_i})_{*0}\circ}} \kappa(x))|<\delta\rforal x\in {\cal P}\andeqn
\tilde{\tau}_i(h)>\Delta(\hat{h})\rforal h\in\mathcal H_{1, 2},
\eneq
{{where ${\rm tr}$ is the tracial state of $M_s$
($\pi_{x}$ {{was}}  defined in the statement of the lemma).}}
By Lemma \ref{APPextnn}, there exists a unital homomorphism $\phi'_i: C\to M_{sK'}$ such that
$${{(\phi_i')_{*0}}}=K\kappa,$$
as we identify $K_0(C([0,1], M_s))$ with $\Z$ and
\begin{equation}\label{pt-exi}
|\mathrm{tr}\circ\phi'_i(h)-\tau_{x_i}(h)|<\min\{\sigma/16, \Delta(\hat{h})/4, \sigma_1/8,\ h\in\mathcal H_{1, 1}\}\rforal h\in\mathcal H\cup \mathcal H_{1, 1}.
\end{equation}
{{Note $(\phi_i')_{*1}=0$ since $K_1(M_{sK'})=\{0\}.$}}
In particular, by \eqref{pt-exi}, one has that, for any $0\leq i\leq n-1$,
$$
|\mathrm{tr}\circ\phi'_i(h)-\mathrm{tr}\circ\phi'_{i+1}(h)|<\sigma_1\rforal h\in {\cal H}_{1,1}.
$$
Note that $\gamma(\tau_{x_i})(h)>\Delta(\hat{h})$ for any $h\in\mathcal H_{1, 1}$ by  {{hypothesis.}}
It then also follows from \eqref{pt-exi} that,  for any $0\leq i\leq n$,
$$\mathrm{tr}\circ\phi_i'(h)>\Delta(\hat{h})/2\rforal h\in\mathcal H_{1, 1}.$$
Define {{the amplification $\phi_i''$ as}}
$$\phi''_i:=\phi_i'\otimes 1_{\mathrm{M}_{{{m_1}}}(\mathbb C)}: C\to M_{sK}{{=M_{{{m_1}}(sK')}}}.$$
{{Then}} {{$(\phi_i'')_{*j}=(\phi_{i+1}'')_{*j},$ $j=0,1.$
By the choice of ${{m_1}},$ $[\phi_i'']|_{K_j(C, \Z/k\Z)}=0=[\phi_{i+1}'']|_{K_j(C, \Z/k\Z)}$ on each non-zero $K_j(C, \Z/k\Z),$ $j=0,1,$ $k=2,3,....$}}
{{Therefore}}
$$[\phi''_i]=[\phi''_{i+1}]\quad \textrm{in $KL(C, M_{sK})$}.$$
It then follows from {{Theorem \ref{UniqAtoM}}} that there is a unitary $u_1\in M_s$ such that
$$\|\phi''_{0}(h)-\mathrm{Ad}u_1\circ\phi''_{1}(h)\|<
{{\ep_0}}\rforal h\in\mathcal H.$$
Consider the maps $\mathrm{Ad}u_1\circ\phi''_{1}$ and $\phi''_{2}$.
Applying {{Theorem \ref{UniqAtoM}}} again, one obtains a unitary $u_2\in M_{sK}$ such that
$$
\|\mathrm{Ad}u_1\circ\phi''_{1}(h)-\mathrm{Ad}u_2\circ\phi''_{2}(h)\|<{{\ep_0}}
\rforal h\in\mathcal H.
$$
{{Repeating}} this argument for all $i=1,...,n$, one obtains  unitaries  $u_i\in M_{sK}$ such that
$$
\|\mathrm{Ad}u_i\circ\phi''_{i}(h)-\mathrm{Ad}u_{i+1}\circ\phi''_{i+1}(h)\|<{{\ep_0}}
\rforal h\in\mathcal H.
$$
Then define $\phi_0=\phi_0''$ and $\phi_i=\mathrm{Ad}u_i\circ\phi''_i$, and one has
\begin{equation}\label{close-ev}
\|\phi_{i}(h)-\phi_{i+1}(h)\|<{{\ep_0}}
\rforal h\in\mathcal H.
\end{equation}

Define the linear map $\phi: C \to M_{sK}([0, 1])$ by
$$\phi(f)(t)=\frac{t-x_i}{x_{i+1}-x_i}\phi_i(f)+\frac{x_{i+1}-t}{x_{i+1}-x_i}\phi_{i+1}(f),\quad\textrm{if $t\in[x_i, x_{i+1}]$}.
$$
Since each $\phi_i$ is a homomorphism, by \eqref{close-ev},
the map $\phi$ is
$\mathcal F$-$\epsilon$-multiplicative.
It is clear that $[\phi]|_{K_0(C)}=K\kappa$. On the other hand, for any $x\in [x_i, x_{i+1}]$ for some $i=1, ..., n-1$, one has that for any $h\in\mathcal H$,
\begin{eqnarray*}
&&|\gamma(\tau_x)(h)-\tau_x\circ \phi(h)| \\
& = &|\gamma(\tau_x)(h)- (\frac{x-x_i}{x_{i+1}-x_i}\mathrm{tr}(\phi_i(h))+\frac{x_{i+1}-x}{x_{i+1}-x_i}\mathrm{tr}(\phi_{i+1}(h)))|\\
&<&|\gamma(\tau_x)(h)- (\frac{x-x_i}{x_{i+1}-x_i}\gamma(\tau_{x_i})(h)+\frac{x_{i+1}-x}{x_{i+1}-x_i}\gamma(\tau_{x_{i+1}})(h))|+\sigma/4 \quad\textrm{{(by \eqref{pt-exi})}}\\
&<&|\gamma(\tau_x)(h)- \gamma(\tau_{x_{i+1}})(h))|+3\sigma/8  \quad \quad \quad \quad \quad\quad \quad  \quad \quad \quad \quad \quad \quad \quad\,\,\,\textrm{{(}by \eqref{eq-par}{)}}\\
&<&\sigma/2 \quad  \quad \quad \quad \quad \quad\quad\quad \quad  \quad \quad \quad \quad \quad\quad\quad \quad  \quad \quad \quad \quad \quad\quad\quad\,\,\,\,\,\,\,\,\,\,\,\,\,\,\textrm{{(by \eqref{eq-par})}}.
\end{eqnarray*}
Hence for any $h\in\mathcal H$,
$$|\gamma(\tau)(h)-\tau\circ\phi(h)|<\sigma$$
for any $\tau\in T(M_{sK}(C([0, 1])))$.

Note that  {{$\pi_0\circ \phi=\phi_0$  and $\pi_1\circ \phi=\phi_n$}} which are \hm s.
Thus the map $\phi$ satisfies the {{conclusion}}
of the lemma.

\end{proof}

\begin{thm}\label{ExtTraceC-D}
 Let {{$C=C(X, C_0),$ where $C_0\in {\cal C}$ and $X$ is a point, or $X=\T.$}}
 Let $\Delta: C_+^{q, 1}\setminus\{0\}\to (0, 1)$ be an order preserving map. Let $\mathcal F, \mathcal H\subset C$ be finite subsets, and let $1>\sigma, \epsilon>0$.
There exist  a finite subset $\mathcal H_1\subset C_+^{\bf 1}\setminus \{0\}$, $\delta>0$,
a finite subset ${\cal P}\subset K_0(C),$ and a positive integer $K$ such that for any continuous affine map $\gamma: T(D)\to T(C)$ satisfying
$$\gamma(\tau)(h) > \Delta(\hat{h})\tforal h\in \mathcal H_1\tforal \tau\in T(D),$$
where $D$ is a \CA\, in $\mathcal C$,
any positive \hm\, $\kappa: K_0(C)\to K_0(D)$ {with $\kappa([1_C])=s[1_D]$} for some integer $s \ge 1$ satisfying
$$|\rho_C(x)(\gamma(\tau))-{(1/s)}\tau(\kappa(x))|<\delta\tforal \tau\in T(D)$$
and for all $x\in {\cal P},$
there is a $\mathcal F$-$\epsilon$-multiplicative positive linear map $\phi: C\to M_{{sK}}(D)$ such that
$${{[\phi]|_{K_0(A)}}}={K}\kappa$$ {{(see \ref{KLtriple})}} and
\beq\label{18-1610-nn}
|(1/(sK))\tau\circ\phi(h)-\gamma(\tau)(h)|<\sigma\tforal h\in\mathcal H\tand \tau\in T(D).
\eneq

In the case that $C\in\mathcal C$, the map $\phi$ can be chosen to be a homomorphism.
\end{thm}
\begin{proof}
{{Recall that if $\Phi: B_1\to B_2$ is a map from {\blue{a \CA\,}} $B_1$ to {\blue{a \CA\,}} $B_2,$ then
we will continue to use $\Phi$ for {{the amplification}} $\Phi\otimes {\rm id}_{M_n}: B_1\otimes M_n\to B_2\otimes M_n.$
We will use this practice  in this proof and  in the rest of this section.}}
As in the proof of {{Lemma}} \ref{ExtTraceI-D}, since \CA s in ${\cal C}$ are semiprojective, we  only need to prove the first part
of the statement.
Without loss of generality, one
 may assume that $\mathcal F$ {{is a subset of the unit ball
of $C,$ {\blue{$1_C\in {\cal F}$,}} and $\{ab: a, b\in {\cal F}\}\subset \mathcal H$.}}
{{Fix $1>\ep, \sigma>0.$}}
{{Replacing}} $C$ by $M_m(C)$ for some integer $m\ge 1,$  {{and applying \ref{FG-Ratn},}}
we may find a finite subset  {{$P$}}
of projections in $C$ {{such that ${\cal P}=\{[p]: p\in P\}$ generates $K_0(C).$}}
{{We}} may also assume that ${{P}}\subset \mathcal H$.

Since the K-groups of $C$ {{are}}  finitely generated (as abelian groups), there is ${{m'}}\in\mathbb N$ such that,
{{for any $x\in \mathrm{Tor}(K_i(C)),$ $i=0,1,$
$mx=0$ for some $0<m\le m'.$ Set $m_1=m!.$}}

Let $\mathcal H_{1, 1}\subset C_+^1\setminus \{0\}$ (in place of $\mathcal H_{1}$),
{{${\cal H}_{1,1}'\subset C_{s.a.}$ {{(in place of ${\cal H}_2$)}} be finite subsets
and ${\cal Q}\subset \underline{K}(C)$ (in place of ${\cal P}$) be another finite}} {\blue{subset,}}
and $\sigma_1>0$ (in place of ${{\delta}}$) be  a positive number {\blue{as provided}}  
by
{{Theorem \ref{UniqAtoM}}} with respect to $C$ (in  place of ${{A}}$), $\min\{\sigma/{4}, \epsilon/{{6}}\}$ (in  place of $\epsilon$), $\mathcal H$ (in place of $\mathcal F$), and $\Delta$.

{{Set $\sigma_0={\frac{1}{2}\min\{\sigma/16, \sigma_1/4, \min\{\Delta(\hat{h})/2: \ h\in\mathcal H_{1, 1}\}\}}.$}}
Let $\mathcal H_{1, 2}\subset C_+^{\bf 1}\setminus \{0\}$ (in place of $\mathcal H_1$) {be a finite subset}, let $\sigma_2$ (in place of $\delta$)
be a positive  number,
{and $K_1$ (in place of $K$) be an integer {{as provided}} by  Lemma \ref{APPextnn}} with respect to $\mathcal H\cup\mathcal H_{1, 1}\cup {\cal H}_{1,1}'$
and  $\sigma_0$
(in place of {$\sigma$}) and
$\Delta.$ {{\Wlog, we may assume that ${\cal H}_{1,2}\supset {\cal H}_{1,1}.$}}

Let $\mathcal H_{1, 3}\subset C_+^{\bf 1}\setminus \{0\}$ (in place of $\mathcal H_1$), $\sigma_3>0$ (in place of $\delta$){\blue{, and}} $K_2$
(in place of $K$) be  {\blue{finite subsets}} and a constant {\blue{as provided}} 
by  Lemma \ref{ExtTraceI-D}  with respect to $C$,
$\mathcal H\cup{\cal H}_{1,1}'\cup\mathcal H_{1, 2}$
(in place of $\mathcal H$), $\sigma_0$
(in place of $\sigma$), $\ep/{{12}}$ (in place of $\ep$), ${\cal H}$ (in place of
${\cal F}$){\blue{,  and}} $\Delta$ (with the same ${\cal P}$ {\blue{as above}}).

Put $\mathcal H_1=\mathcal H_{1, 1}\cup\mathcal H_{1, 2}\cup\mathcal H_{1, 3}\cup {{P}}$, $\delta=\min\{{\sigma_1/2}, \sigma_2,
\sigma_3,{1/4}\},$ and $K={ m_1}K_1K_2$.
Let
$$
D={{A}}(F_1, F_2, \psi_0, \psi_1)=\{(f,a)\in C([0,1], F_2)\oplus F_1: f(0)=\psi_0(a)\andeqn f(1)=\psi_1(a)\}
$$
be any
\CA\, in $\mathcal C$,
{{where $\psi_i: F_1\to F_2$ is a unital \hm, $i=0,1,$}} and let $\gamma: T(D)\to T(C)$ be a given continuous affine map satisfying $$\gamma(\tau)(h) > \Delta(\hat{h})\rforal h\in \mathcal H_{1} \rforal \tau\in T(D).$$

{{Write $F_1=M_{R(1)}\oplus M_{R(2)}\oplus\cdots \oplus M_{R(l)},$
$F_2=M_{r(1)}\oplus M_{r(2)}\oplus \cdots \oplus M_{r(k)}$ and
$I_j=C([0,1], M_{r(j)}),$ $j=1,2,...,k.$}}
{{Denote by ${{\pi_{e,j}}}: D\to M_{{R(j)}}$ the
\hm\, which is the composition of $\pi_e: D\to F_1$ (defined in \ref{DfC1}) and the projection from $F_1$ onto
$M_{R(j)}$ ($1\le j\le l$). Denote by $\pi^{I_j}: D\to I_j$  the restriction map
defined by $(f,a)\to f|_{[0,1]_j}$  (see \ref{2Rg15}).}}

{L}et $\kappa: K_0(C)\to K_0(M_s(D))$ be any positive map {with $s[1_D]=\kappa([1_C])$} satisfying
$$
|\rho_C(x)(\gamma(\tau))-{(1/s)}\tau(\kappa(x))|<\delta\rforal \tau\in T(D)
$$
{and for all $x\in {\cal P}.$}
Write $C([0, 1], F_2)=I_1\oplus I_2\oplus \cdots \oplus I_{{k}}$ with $I_i=C([0, 1], M_{r(i)})$, $i=1, ..., {{k}}$.
Note that $\gamma$ induces a continuous affine map $\gamma_i: T(I_i)\to T(C)$
 defined by $\gamma_i(\tau)=\gamma(\tau\circ{{\pi^{I_i}}})$ for each $1\leq i\leq {k}.$
 It is clear that
for any $1\leq i\leq {k}$,
one has that
\beq
\label{intv-dense}
\gamma_i(\tau)(h)>\Delta(\hat{h}) \rforal h\in\mathcal H_{1, 3}{\andeqn} \rforal \tau\in T(I_i),
\andeqn\\\label{intv-cpt}
|\rho_C(x)(\gamma_i(\tau))-\tau(({{\pi^{I_i}}})_{*0}\circ\kappa(x))|<\delta\leq \sigma_3\rforal  \tau\in T(M_s(I_i))
\eneq
{and for all $x\in {\cal P}$ and for any $1\le i\le k.$}
Since
\beq\label{1812-18n-1-1}
&&\gamma(\tau)(h)>\Delta(\hat{h})\rforal h\in\mathcal { {H}}_{1, 2}\,\,\, {\rm and}\rforal \tau\in T(D),
\andeqn\\
&&|\rho_C(x)(\gamma(\tau))-{(1/s)}\tau(\kappa(x))|<\delta  \rforal \tau\in T(D)
\eneq
{and for all  $x\in {\cal P}, $} one has that, for each $j,$
\beq\label{1611-18n-1}
\gamma\circ({{\pi_{e,j}}})^*(\mathrm{{{tr}}}_j')(h)>\Delta(\hat{h})\rforal h\in\mathcal H_{1, 2}\andeqn\\\label{1611-18n-2}
|\rho_C(x)(\gamma\circ({{\pi_{e,j}}})^*(\mathrm{tr}_j'))-\mathrm{{{tr}}}_j'(({{\pi_{e,j}}})_{*0}\circ\kappa(x))|<\delta\leq \sigma_2,
\eneq
where ${\rm tr}_j'$ is the tracial
state on $M_{R(j)},$ for all  $x\in K_0(C)$ and
where $\gamma\circ ({{\pi_{e,j}}})^*({\rm {{tr}}}_j')=\gamma({{{\rm tr}_j'\circ \pi_{e,j}}}).$

{{Using \eqref{1611-18n-1},  \eqref{1611-18n-2} and  applying \ref{APPextnn} to $(\pi_{e,j})_{*0}\circ K_1K_2\kappa,$
one obtains}} a homomorphism $\phi'_j: C\to M_{R(j)}\otimes M_{sK_1K_2}$ such that
\beq\label{C-D-nnn1}
&&(\phi'_j)_{*0}=({{\pi_{e,j}}})_{*0}\circ  K_1K_2\kappa\andeqn\\
\label{pt-pre}
&&\hspace{-0.3in}|\mathrm{tr}_j\circ \phi'_j(h)- (\gamma \circ ({{\pi_{e,j}}})^*(\mathrm{tr}_j'))(h) |
<{{\sigma_0}}\rforal h\in  {\cal H}\cup{\cal H}_{1,1}\cup {\cal H}_{1,1}',
\eneq
where
${\rm tr}_j$ is the tracial
state on $M_{R(j)}\otimes M_{s{{K_1K_2}}}.$
{{By \eqref{1611-18n-1},  \eqref{pt-pre} and the choice of $\sigma_0,$
\beq\label{1812-18n-1}
{\rm tr}_j\circ \phi_j'(h)\ge \Delta(\hat{h})/2\rforal h\in {\cal H}_{1,1}.
\eneq
}}
{{Set}} $\phi'=\bigoplus_{j=1}^l \phi'_j: C\to F_1\otimes\mathrm{M}_{sK_1K_2}.$
{{Then, for all $t\in T(F_1\otimes\mathrm{M}_{sK_1K_2}),$
\beq\label{1812-18n-2}
t\circ \phi'(h)\ge \Delta(\hat{h})/2\rforal h\in {\cal H}_{1,1}.
\eneq}}
{{By \eqref{pt-pre},  for all $h\in {\cal H}\cup {\cal H}_{1,1}{\cal H}_{1,1}'$ and for all $t\in T(F_1),$
\beq\label{1812-18n-1+1}
|(t\otimes {\rm tr}_{sK_1K_2})\circ \phi'(h)-\gamma(t\circ \pi_e)(h)|
<\sigma_0,
\eneq
where ${\rm tr}_{sK_1K_2}$ is the tracial state of $M_{sK_1K_2}.$}}
Applying Lemma \ref{ExtTraceI-D}, using \eqref{intv-dense} and \eqref{intv-cpt},  {{and
tensoring the resulting map with $1_{M_{K_1}},$}}
for any $1\leq i\leq {k}$,
{{one obtains}} an ${\cal H}$-$\ep/4$-multiplicative \morp\,
$\phi_i: C\to I_i\otimes M_{{sK_1K_2}}$ such that ${{[\phi_i]|_{K_0(C)}}}=({{\pi^{I_i}}})_{*0}\circ K_1K_2  \kappa$ and
\begin{equation}\label{int-pre}
|{\tilde \tau}\circ \phi_i(h)- ((\gamma\circ({{\pi^{I_i}}})^*(\tau))(h) |<\sigma_0
\end{equation}
for all $h\in \mathcal H\cup\mathcal H_{1, 1}\cup{\cal H}_{1,1}',$
and for all  $\tau\in T({I_i})$
{{and ${\tilde \tau}=\tau\otimes {\rm tr}_{sK_1K_2}.$}}
{{Furthermore, as
{{in the}}  conclusion of Lemma \ref{ExtTraceI-D}, the restrictions of $\phi_i$ to both boundaries (end points of the interval) are  homomorphisms.}}


For each $1\leq i\leq {k}$, denote by $\pi_{i, {{t}}}$
the evaluations of $I_i\otimes M_{{sK_1K_2}}$ at the point {{$t\in [0,1].$}}
{{Also, for each $1\le i\le k,$  define
$\psi_{0,i}=(q_i\circ \psi_0)\otimes M_{sK_1K_2}$ and $\psi_{1,i}=(q_i\circ \psi_1)\otimes {\rm id}_{M_{sK_1K_2}},$ where $q_i: F_2\to M_{r(i)}$
is the projection map.}}
Then one has
\begin{equation}\label{C-D-n1}
\psi_{0,i}\circ \pi_e=\pi_{i,0}\circ {{\pi^{I_i}}}.
\end{equation}
{{(Recall that {\blue{a map}} $\Phi$ {\blue{is}}   identified with $\Phi\otimes {\rm id}_{sK_1K_2}$).}}
It follows that
\begin{eqnarray}\label{C-Dn-2}
(\psi_{0, i}\circ\phi')_{*0} & = & (\psi_{0,i})_{*0}\circ (\sum_{j=1}^l({{\pi_{e,j}}})_{*0})\circ K_1K_2\kappa\\\label{C-Dn-2+}
&=&(\psi_{0,i})_{*0}\circ (\pi_e)_{*0}\circ K_1K_2\kappa=(\pi_{i,0}\circ {{\pi^{I_i}}})_{*0}\circ K_1K_2\kappa\\\label{18-CDn-2}
&=&(\pi_{i, 0})_{*0}\circ [\phi_i]|_{K_0(A)}.
\end{eqnarray}
{{Note that $\psi_{0,i}$ is unital.}} Therefore, by \eqref{1812-18n-2},  for any ${\rm tr}^{(i)}\in T(M_{r(i)}\otimes M_{sK_1K_2}),$
one has
\begin{equation}\label{1611-18n-10}
\mathrm{tr}^{(i)}\circ (\psi_{0, i}\circ\phi')(h)\geq \Delta(\hat{h})/2\rforal h\in \mathcal H_{1, 1},
\end{equation}
and, since $\pi_{i,0}$ is unital,  by \eqref{int-pre} and \eqref{1812-18n-1-1},
\begin{equation}\label{1611-18n-11}
\mathrm{tr}^{(i)}\circ (\pi_{i, 0}\circ \phi_i)(h)\geq \Delta(\hat{h})/2\rforal h\in \mathcal H_{1, 1}.
\end{equation}
Let  $t^{(i)}$ be the tracial  state of $M_{r(i)}.$
It follows from  \eqref{int-pre},
\eqref{C-D-n1} and \eqref{pt-pre}  that
\beq\label{1812-18n-10}
\mathrm{tr}^{(i)}\circ (\pi_{i, 0}\circ \phi_i)(h)\approx_{\sm_0} \gamma((\pi^{I_i})^*)(t^{(i)}\circ\pi_{i,0})(h)
=\gamma(t^{(i)}\circ \pi_{i,0}\circ \pi^{I_i})(h)\\\label{1812-18n-10+}
=\gamma(t^{(i)}\circ \psi_{i,0}\circ \pi_e)(h)\approx_{\sm_0}{\rm tr}^{(i)}\circ \phi'(h)
\rforal h\in \mathcal
H_{1,1}\cup {\cal H}_{1,1}'.
\eneq
Consider {{the amplifications}}
\beq\nonumber
\phi_i^\sim:&=&\phi_i\otimes 1_{\mathrm{M}_{{m_1}}(\mathbb C)}: C\to I_i\otimes M_{{sK}}\andeqn\\
\phi'':&=&\phi'\otimes 1_{\mathrm{M}_{{m_1}}(\mathbb C)}: C\to F_1 \otimes M_{{sK}}.
\eneq
{{Since $K_1(M_{r(i)})=\{0\},$ by \eqref{18-CDn-2}, $(\psi_{0, i}\circ\phi'')_{*j}=(\pi_{i,0}\circ \phi_i^\sim)_{*j},$ $j=0,1.$
By the choice of $m_1,$\linebreak
  $[\psi_{0, i}\circ\phi'']|_{K_j(C, \Z/k\Z)}=0=[\pi_{i,0}\circ \phi_i^\sim]|_{K_j(C, \Z/k\Z)},$ $j=0,1,$ and $i=0,1,$
for any $k$ {{such}} that $K_i(C,\Z/k\Z)\not=\{0\}.$}}
{{It follows that}}
$$[\psi_{0, i}\circ\phi'']=[\pi_{i, 0}\circ \phi^\sim_i]\quad\textrm{in $KL(C, \mathrm{M}_{r(i)sK})$}.$$
Therefore, by {{Theorem \ref{UniqAtoM}}} {{(by also \eqref{1611-18n-10}, \eqref{1611-18n-11},
and \eqref{1812-18n-10+}),}} there is a unitary $u_{i, 0}\in M_{r(i)} \otimes M_{{sK}}$ such that
$$
\|\mathrm{Ad}u_{i, 0}\circ \pi_{i, 0}\circ \phi^\sim_i(f) -  \psi_{0, i}\circ\phi''(f)\|<\min\{\sigma/{4}, \epsilon/{{6}}\}\rforal f\in\mathcal H.
$$
Exactly the same argument shows that
there is a unitary $u_{i, 1}\in M_{r(i)} \otimes M_{{sK}}$ such that
$$\|\mathrm{Ad}u_{i, 1}\circ \pi_{i, 1}\circ \phi_i^\sim(f) -  \psi_{1, i}\circ\phi''(f)\|<\min\{\sigma/4, \epsilon/{{6}}\}\rforal f\in\mathcal H.$$

Choose {two}  path{s} of {unitaries}  ${\{u_{i, 0}(t):t\in [0,1/2]\}\subset} M_{r(i)} \otimes M_{{sK}}$
such that $u_{i, 0}(0)=u_{i, 0}$ and {$u_{i, 0}(1/2)=1_{M_{r(i)} \otimes M_{sK}},$} {and
$\{u_{i,1}(t):t\in [1/2, 1]\}\subset M_{r(i)}\otimes M_{{sK}},$ such that
$u_{i,1}(1/2)=1_{M_{r(i)}\otimes M_{sK}}$ and $u_{i,1}(1)=u_{i,1}$ }
{Put $u_i(t)=u_{i,0}(t)$ if $t\in [0,1/2)$ and $u_i(t)=u_{i,1}(t)$ if $t\in [1/2,1].$
Define ${{\phi_{i,I}}}: C\to I_i\otimes M_{sK}$ by
\begin{equation*}
\pi_t\circ \phi_{i,I}={\rm Ad}\, u_i(t)\circ\pi_t\circ  \phi'_i,
\end{equation*}
where $\pi_t: I_i\otimes M_{sK}\to M_{r(i)}\otimes M_{sK}$ is the point evaluation at $t\in [0,1].$}

{One has that, for each $i,$
\beq\nonumber
&&\| \pi_{i, 0}\circ \phi_{i,I}(f) -  \psi_{0, i}\circ\phi''(f)\|<\min\{\sigma/4, \epsilon/{{6}}\}
\andeqn\\
&&\| \pi_{i, 1}\circ \phi_{i,I}(f) -  \psi_{1, i}\circ\phi''(f)\|<\min\{\sigma/4, \epsilon/{{6}}\}\rforal f\in\mathcal H.
\eneq
}
For each $1\leq i\leq {k}$, let $\epsilon_i<1/2$ be a positive number such that
\beq\label{1611-18n-21}
&&\hspace{-0.5in}\|\phi_{i,I}(f)(t) - \psi_{0, i}\circ\phi''(f)\|<\min\{\sigma/{4}, \epsilon/{{6}}\} \rforal f\in\mathcal H \rforal t\in [0, \epsilon_i],\andeqn\\\label{1611-18n-22}
&&\hspace{-0.5in}\|\phi_{i,I}(f)(t) -  \psi_{1, i}\circ\phi''(f)\|<\min\{\sigma/{4}, \epsilon/{{6}}\}\rforal f\in\mathcal H \rforal t\in [1-\epsilon_i, 1].
\eneq
Define ${\Phi}_i: C\to I_i\otimes M_{{sK}}$ to be
$${\Phi}_i(f)(t)=\left\{
\begin{array}{ll}
(\ep_i-t/\ep_i)(\psi_{0, i}\circ\phi'') + (t/\epsilon_i)\phi_{i.I}(f)(\epsilon_i), & \textrm{if $t\in [0, \epsilon_i]$},\\
\phi_{i,I}(f)(t), & \textrm{if $t\in [\epsilon_i, 1-\epsilon_i]$ },\\
((t-1+\epsilon_i)/\epsilon_i)(\psi_{1, i}\circ\phi'') + ((1-t)/\epsilon_i)\phi_{i,I}(f)(\epsilon_i), & \textrm{if $t\in [1- \epsilon_i, 1]$}.
\end{array}
 \right.$$
The map ${\Phi}_i$ is not necessarily a homomorphism, but it is
$\mathcal F$-$\epsilon$-multiplicative. Moreover, it satisfies the relations
\beq\label{C-D-n5}
\hspace{-0.2in}\pi_{i, 0}\circ {\Phi}_i(f) =  \psi_{0, i}\circ\phi''(f)
\andeqn
\pi_{i, 1}\circ {\Phi}_i(f) =  \psi_{1, i}\circ\phi''(f)\rforal f\in\mathcal H, i=1,..., {k}.
\eneq
{Define  $\Phi': C\to C([0,1], F_2)\otimes M_{sK}$ by $\pi_{i,t}\circ \Phi'={{\pi_t\circ}} \Phi_i,$ where
$\pi_{i,t}: C([0,1], F_2)\otimes M_{sK}\to M_{r(i)}\otimes M_{sK}$ {{is}} defined by the point evaluation
at $t\in [0,1]$ (on the $i$-th {\blue{direct}}  summand).
Define
$$\phi(f)=(\Phi'(f),{{\phi''}}(f)).$$
It follows from \eqref{C-D-n5}
that
$\phi$ is {\blue{an}} $\mathcal F$-$\epsilon$-multiplicative contractive completely positive linear map from $C$ to $D\otimes M_{sK}$}.
It follows from \eqref{C-D-nnn1} (see also \eqref{C-Dn-2+}) that
\beq\label{C-D-n9}
[\pi_e\circ \phi(p)]=[\phi{{''}}(p)]=(\pi_e)_{*0}\circ K\kappa([p]) \rforal p\in {\cal P}.
\eneq
Since $(\pi_e)_{*0}: K_0(D)\to \Z^l$ is injective {{(see \ref{2Lg13})}},  one has
\begin{equation}\label{C-D-n11}
\phi_{*0}=K\kappa.
\end{equation}
{{For any $\tau_0\in T(F_1\otimes M_{sK}),$ let $\tau=\tau_0\circ \pi_e.$
Note also $\pi_e\circ \phi=\phi''.$
 By \eqref{1812-18n-1+1},
\beq\label{1611-18n-31}
|\tau\circ  \phi(h)-\gamma(\tau)(h)|<\sigma/16\rforal h\in {\cal H}.
\eneq
}}
{{Let $\tau\in T(D)$ be defined by $\tau(f)=\sum_{j=1}^k \int_{(0,1)}
{\rm tr}_{r(j)sK}(\pi^{I_j}(f))d\mu_j$
for all $f\in D\otimes M_{sK},$
where  ${\rm tr}_{r(j)sK}$
is the tracial state of $M_{r(j)sK}$ and $\mu_j$ is a positive Borel
measure  on $(0,1)$ with $\sum_{j=1}^k\|\mu_j\|=1.$}}
{{Then, by the definition of $\phi,$  \eqref{1611-18n-21}, \eqref{1611-18n-22}, and \eqref{int-pre},
\beq\label{1611-18n-32}
|\tau\circ \phi(h)-\gamma(\tau)(h)|<\sigma/4+\sigma/16=5\sigma/16 \rforal h\in {\cal H}.
\eneq
In general, $\tau\in T(D\otimes M_{sK})$ has the form
\beq
\tau(f)=\tau_0\circ \pi_e(f)+\sum_{j=1}^k\int_{(0,1)}{\rm tr}_{r(j)sK}(\pi^{I_j}(f))d\mu_j\rforal f\in D\otimes M_{sK},
\eneq
where $\|\tau_0\|+\sum_{j=1}^k\|\mu_j\|=1$}} {{(see the proof of \ref{2Lg8})}}.
It follows from \eqref{1611-18n-31} and \eqref{1611-18n-32} that
$$
|(1/sK)\tau\circ \phi(h)-\gamma(\tau)(h)|<\sigma\rforal h\in {\cal H}
$$
and for all $\tau\in T(D).$
\end{proof}

\begin{lem}\label{TtoDelta}
Let $C\in {\cal C}.$ For any $\ep>0$ and any finite subset
${\cal H}\subset C_{s.a.},$ there exists a finite
subset of extremal traces ${ E}\subset T(C)$  and
a continuous affine map $\lambda: T(C)\to {{\triangle_E}},$
where $\triangle_{{E}}$ is the convex hull of of ${{E}},$ such that
\beq\label{TtoD-1}
|\lambda(\tau)(h)-\tau(h)|<\ep,\quad h\in\mathcal H,\ \tau\in T(C).
\eneq
\end{lem}

\begin{proof}
{{We}} may assume that
${\cal H}$ is in the unit ball of $C.$
Write $C={{A}}(F_1, F_2, \psi_0, \psi_1),$ where
$F_1=M_{R(1)}\oplus M_{R(2)}\oplus\cdots \oplus M_{R(l)}$ and
$F_2=M_{r(1)}\oplus M_{r(2)}\oplus \cdots \oplus M_{r(k)}.$
Let $\pi_{e,i}: C\to M_{R(i)}$ be the surjective \hm\, defined by the composition
of $\pi_e$ and the projection from $F_1$ onto $M_{R(i)},$ and
${{\pi^{I_j}}}: C\to C([0,1],M_{r(j)})$ the restriction which may also be
viewed as the restriction of the projection from $C([0,1], F_2)$ to $C([0,1], M_{r(i)}).$ Denote by $\pi_t\circ {{\pi^{I_j}}}$
the composition of ${{\pi^{I_j}}}$ and the point evaluation at $t\in [0,1].$
There is $\dt>0$ such that, for any $h\in {\cal H},$
\beq\label{TtoD-2}
\|{{\pi^{I_j}}}(h)(t)-{{\pi^{I_j}}}(h)(t')\|<\ep/16\rforal h\in {\cal H}
\eneq
and $|t-t'|<\dt,$ $t, t'\in [0,1].$

Let $g_1,g_2,..., g_n$ be a partition of unity over {{the}} interval $[\dt, 1-\dt]$  {{subordinate}} to
an open cover with order 2
such that each ${\rm supp}(g_i)$ has diameter $<\dt$ and
$g_{s}g_{s'}\not=0$ implies that $|s-s'|\le 1.$
Let $t_s\in {\rm supp}(g_s){{\cap}} [\dt, 1-\dt]$ be a point.  We may assume that
$t_s<t_{s+1}.$ We may further choose
$t_1=\dt$ and $t_n=1-\dt$ and assume that $g_1(\dt)=1$ and $g_n(1-\dt)=1,$
choosing an appropriate open cover of order 2.

Extend $g_s$ to $[0,1]$ by
defining $g_s(t)=0$ if $t\in [0, \dt)\cup (1-\dt,1]$ for $s=2,3,...,n-1$ and
\beq\label{ToD-3}
\hspace{-0.2in}g_1(t)=g_1(\dt)(t/\dt)\,\,\,{\rm for}\,\,\, t\in [0,\dt)\andeqn g_n(t)=g_n(1-\dt)(1-t)/\dt
\,\,\,{\rm for}\,\,\, t\in (1-\dt, 1].
\eneq
Define $g_0=1-\sum_{s=1}^n g_s.$ {{Then $g_0(t)=0 $ for all $t\in[\dt,1-\dt]$.}}
{{Put ${\bar g}_s=({g_s\cdot 1_{F_2}}, 0)\in C$ for $1\le s\le n,$ and
${\bar g}_0= ({g_0\cdot 1_{F_2}},1_{F_1}),$ so that ${{\psi_0\circ \pi_e({\bar g_0})}}=\psi_0(1_{F_1})$ and $
{{\psi_0(\pi_e(g_0))}}=\psi_1(1_{F_1}).$}}
 Let $g_{s,j}={{\pi^{I_j}(\bar g_s}}),$ $s=1,2,...,n,$ $j=1,2,...,k.$
Let $p_i\in F_1$ be the {{support}} projection corresponding to the summand
$M_{R(i)}.$
Choose $d_i\in C([0,1], F_2)$ such that $d_i(t)=\psi_0(p_i)$  for
$t\in [0,\dt]$ and $d_i(t)=\psi_1(p_i)$ for $t\in [1-\dt,1]$ and $0\le d_i(t)\le 1$
for $t\in (\dt, 1-\dt).$ Note that ${{{\bar d_i}=(d_i\cdot 1_{F_2}, p_i)}}\in C.$
{\blue{We}} may assume that
$\sum_{i=1}^l g_0d_i=g_0$ {{and $\sum_{i=1}^l {\bar g}_0{\bar d_i}={\bar g}_0.$}}
Denote by ${\rm tr}_i$ the  tracial state on $M_{R(i)}$ and
${\rm tr}_j'$ the tracial state on ${\blue{M_{r(j)}}},$ $i=1,2,...,l,$ and $j=1,2,...,k.$
Let
\beq\label{TtoD-4}
{{E}}=\{{\rm tr}_i\circ \pi_{e,i}:1\le i\le l\}\cup \bigcup_{s=1}^n\{{\rm tr}_j'\circ \pi_{t_s}\circ {{\pi^{I_j}}}: 1\le j\le k\}.
\eneq
Let ${{\triangle_E}}$ be the convex hull of ${{E}}.$
Define $\lambda: T(C)\to {{\triangle_E}}$ by
\beq\label{ToD-5}
\lambda(\tau)(f)=\sum_{j=1}^k \sum_{s=1}^n \tau(g_{s,j}){\rm tr}_j'\circ ({{\pi^{I_j}}}(f)(t_s))+
\sum_{i=1}^l \tau({{{\bar g_0}{\bar d_i}}}){\rm tr}_i\circ \pi_{e,i}(f){{,}}
\eneq
{{where {{we view}} $g_{s,j}
\in C_0((0,1), {\blue{M_{r(j)}}})\subset C$,}} for all $f\in C.$
It is clear that $\lambda$ is a continuous affine map.
Note that if $h\in C$  {{($\pi_{e,j}({\bar d}_i)=0,$ if {\blue{$i\not=j$),}}}}
\beq\label{ToD-6}
\lambda({\rm tr}_j\circ \pi_{e,j})(h)&=&\sum_{i=1}^l {\rm tr}_j\circ \pi_{e,j}({{{\bar g_0}{\bar d_i}}}){\rm tr}_i\circ \pi_{e,i}(h)\\\label{ToD-6+}
&=&{\rm tr}_j\circ \pi_{e,j}({{{\bar g_0}{\bar d_j}}}){\rm tr}_j\circ \pi_{e,j}(h)={\rm tr}_j\circ \pi_{e,j}(h).
\eneq
If  $\tau(f)={\rm tr}_j'\circ ( {{\pi^{I_j}}}(f)(t))$ with $t\in (\dt, 1-\dt),$ then, {{for}} $h\in {\cal H},$
\beq\label{ToD-7}
\tau(h)&=& {\rm tr}_j'\circ ({{\pi^{I_j}}}(h)(t))=(\sum_{s=1}^n{\rm tr}_j'\circ ({{\pi^{I_j}}}(h{{\bar g_s}})(t)))\\
&\approx_{2\ep/16} &\sum_{s=1}^n {\rm tr}_j'{{(\pi^{I_j}(h)(t_s)(\pi^{I_j})({\bar g}_s)(t))}}\\
&=&\sum_{s=1}^n g_{s,j}(t){\rm tr}_j'\circ {{\pi^{I_j}}}(h)(t_s)=
\sum_{s=1}^n \tau({{\bar g}}_{s,j}){\rm tr}_j'({{\pi^{I_j}}}(h))(t_s)\\\label{ToD-7+}
&=&\sum_{i=1}^k\sum_{s=1}^n \tau({{{\bar g}_{s,i}}}){\rm tr}_i'({{\pi^{I_i}}}(h))(t_s) +\sum_{i=1}^l 0 \cdot {\rm tr}_i\circ \pi_{e,i}(h)=\lambda(\tau)(h){{,}}
\eneq
where we {{note}} that $\tau({{\bar g}}_{s,i})=0$ if $i\not=j$. If  $\tau$ has the form $\tau(f)={\rm tr}_j'({{\pi^{I_j}}}(f)(t))$ for some fixed $t\in (0, \dt],$
then for $h\in {\cal H}$ with $h=(h_0, h_1),$ where
$h_0\in C([0,1], F_2)$ and $h_1\in F_1$
{{are}}
such that \nolinebreak
$\psi_0(h_1)=h_0(0)=h(0)$ and $\psi_1(h_1)=h_0(1)=h(1)$  (using
$\sum_{i=1}^lp_i=1_{F_1}$ and $p_i=1_{M_{r(i)}}$),
\beq
\tau(h)&=&{\rm tr}_j'({{\pi^{I_j}(h)}}(t))
={\rm tr}_j'({{\pi^{I_j}(h{\bar g}_1)}}(t))+{\rm tr}_j'({{\pi^{I_j}(h{\bar g}_0)}}(t)))\nonumber \\
&\approx_{\ep/8}&{\rm tr}_j'({{\pi^{I_j}(h)(\dt)\pi^{I_j}({\bar g}_1)}}(t))+{\rm tr}_j'{{(\pi^{I_j}(h)(0)\pi^{I_j}({\bar g}_0)}}(t))\nonumber \\
&=& g_1(t) {\rm tr}_j'{{(\pi^{I_j}(h)}}(\dt)) + g_0(t)
{\rm tr}_j'\circ
(\psi_0(h_1))\nonumber \\
&=&\tau({{\bar g}}_{1,j}){{{\rm tr}_j'}} {{(\pi^{I_j}(h)(t_1))}}+g_0(t)\sum_{i=1}^l{\rm tr}_j'
(\psi_0(h_1p_i))\nonumber\\
&=&\tau({{\bar g}}_{1,j}){\rm tr}_j'\circ {{(\pi^{I_j}(h)}}(t_1))+g_0(t)\big(\sum_{i=1}^l
{\rm tr}_j'( \psi_0(p_i
)){\rm tr}_i\circ\pi_{e,i}(h)\big)\nonumber
\nonumber\\
&=&\tau({{\bar g}}_{1,j}){\rm tr}_j'{{(\pi^{I_j}(h)}}(t_1))+g_0(t)\sum_{i=1}^l {\rm tr}_j'\circ ({{\pi^{I_j}({\bar{d}}_i)}}(t))) {\rm tr}_i(\pi_{e,i}(h))
%
\nonumber\\\nonumber
&=& \tau(g_{1,j}){\rm tr}_j'{{(\pi^{I_j}(h)}}(t_1))+\sum_{i=1}^l {\rm tr}_j'\circ ({{\pi^{I_j}({\bar g}_0{\bar d}_i)(t)}}) {\rm tr}_i(\pi_{e,i}(h))\\\label{1610-n1812}
&=&\tau(g_{1,j}){\rm tr}_j'{{(\pi^{I_j}(h)}}(t_1))+\sum_{i=1}^l \tau({{{\bar g}_0{\bar d}}}_i) {\rm tr}_i(\pi_{e,i}(h))
=\lambda(\tau)(h).
\eneq
The same argument as above shows that,
if
$$\tau(f)={\rm tr}_j'\circ ({{\pi^{I_j}}}(f)(t)),\quad t\in [1-\dt, 1),$$
then
\beq\label{ToD-8-n1}
\tau(h)\approx_{\ep/8} \lambda(\tau)(h)\rforal h\in {\cal H}.
\eneq
It follows {{from \eqref{ToD-6+}, \eqref{ToD-7+}, \eqref{1610-n1812}, and \eqref{ToD-8-n1}}} that
\begin{equation*}
|\tau(h)-\lambda(\tau)(h)|<\ep/8\tforal h\in {\cal H}
\end{equation*}
and for all extreme points of $\tau\in T(C).$
By Choquet's Theorem, for each $\tau\in T(C),$  there exist a Borel probability measure
$\mu_\tau$ on the extreme points $\partial_e{T(C)}$ of $T(C)$ such that
\begin{equation*}
\tau(f)=\int_{\partial_e{T(C)}}f(t)d\mu_\tau\rforal f\in \Aff(T(C)).
\end{equation*}
Therefore, for
each $h\in {\cal H},$
\begin{equation*}
\tau(h)=\int_{\partial_e(T(C)} {{\hat{h}}}(t)d\mu_{\tau}\approx_{\ep/8} \int_{\partial_e{T(C)} }{{\hat{h}}}(\lambda(t))d\mu_{\tau}
=\lambda(\tau)(h)\rforal
\tau\in T(C),
\end{equation*}
as desired.
\end{proof}

\begin{lem}\label{cut-trace}
Let $C$ be a unital stably finite \CA, and let $A\in \mathcal B_1$ (or $\mathcal B_0$). Let $\alpha: T(A)\to T(C)$ be a continuous affine map.

(1)\label{cut-trace-a} For any finite subset $\mathcal H\subset {{C_{s.a.}}},$ {{and}}
 any $\sigma>0$, there are a $C^*$-subalgebra $D\subset A$ and a continuous affine map $\gamma: T(D)\to T(C)$ such that $D\in\mathcal C$ (or $\mathcal C_0$), and
$$| {{\gamma(\imath(\tau))(h)}}- {{\alpha(\tau)(h)}}|<\sigma \tforal \tau\in T(A) \tforal h\in \mathcal H,$$
where $\imath: T(A)\ni \tau\to \frac{1}{\tau(p)}\tau|_D\in T(D)$,  $p=1_D$ {{and $\tau(1-p)<\sigma$
for all $\tau\in T(A).$}}

(2)\label{cut-trace-b} If  {\Green{there is a non-increasing map: $\Delta_0: C_+^{\bf 1}\setminus \{0\}\to (0,1),$
and}}
there are a finite subset $\mathcal H_1\subset C_+$ and $\sigma_1>0$ such that
$$\alpha(\tau)(g)\ge {\Green{\Delta_0(g)}}> \sigma_1
\tforal g\in\mathcal H_1{{\tand}}\tforal \tau\in T(A),$$
then the  affine map $\gamma$ can be chosen so that $$\gamma(\tau)(g)\ge {\Green{\Delta_0(g)}}
>\sigma_1
\tforal g\in\mathcal H_1
\tand\tforal \tau\in T(D).$$.

(3)\label{cut-trace-c} If the positive cone of $K_0(C)$ is generated by a finite subset ${\cal P}$  of projections and there is an order-unit map $\kappa: K_0(C) \to K_0(A)$ which is compatible {{with}} $\alpha$
{{and strictly positive}},  then, for any $\delta>0$, the C*-subalgebra $D$ and $\gamma$ 
{\Green{above}} can be chosen so that  there are {\Green{also}} {positive} homomorphisms $\kappa_0: K_0(C)\to
K_0((1-p)A(1-p))$ and $\kappa_1: K_0(C)\to K_0(D)$ such that  $\kappa_1$ is strictly positive,
{{$\kappa_1([1_C])=[1_D],$}}
$\kappa=\kappa_0+{{\imath_{*0}}}\circ \kappa_1,$ {where $\imath: D\to A$ is the embedding,},  and
\begin{equation}\label{june2-nn1}
|\gamma(\tau)({{q}})-{{\rho_D}}(\kappa_1([{{q}}])){{(\tau)}}|<\dt
{\blue{\tforal}} {{q}}\in {\cal P}{\tand\,\,\, \tau\in T(D)}.
\end{equation}


(4)\label{cut-trace-d} Moreover, in addition to  (3),
if $A\cong A\otimes U$ for some infinite dimensional UHF-algebra, for any given positive integer $K$, the \CA\, $D$ can be chosen so that  $D=M_K(D_1)$
for some $D_1\in {\cal C}$ ({\blue{or}} $D_1\in {\cal C}_0$) and $\kappa_1=K\kappa_1',$
where $\kappa_1': K_0(C)\to K_0(D_1)$ is a strictly positive homomorphism. {Furthermore,
$\kappa_0$ can also be chosen to be strictly positive.}

\end{lem}

\begin{proof}

Write $\mathcal H=\{h_1, h_2, ..., h_m\}.$  We may assume that $\|h_i\|\le 1,$ $i=1,2,...,m.$
Choose $f_1, f_2, ..., f_m\in A_{s.a.}$ such that that $\tau(f_i)=h_i(\alpha(\tau))$ for all $\tau\in T(A)$ and $\|f_i\|\le 2,$ $i=1,2,...,m$
(see 9.2 of \cite{LinTAI}).
Put ${\cal F} =\{1_A,  f_1, f_2,...,f_m\}.$

Let  $\dt>0$ and let ${\cal G}_1$ (in place of ${\cal G}$) be a finite subset
{{as given}} by Lemma 9.4 of \cite{LinTAI}  for $A,$ $\sigma/16$ (in place of $\ep$), and  ${\cal F}$.
Let $\sigma_1=\min\{\sigma/16, \dt/16, {{1/16}}\}.$
We may assume that ${\cal G}_1\supset {\cal F}$ {{and $\|g\|\le 1$ for $g\in {\cal G}_1.$}}
Put ${\cal G}=\{g,\,gh: g,\, h\in {\cal G}_1\}.$
Since $A\in {\cal B}_1$ (or ${\cal B}_0$), {{by the definition of ${\cal B}_1$ (or ${\cal B}_0$),}}
there is a $D\in {\cal C}$ ({{or}}  ${\cal C}_0$) {{with $p=1_D$ such that, for all $g\in {\cal G},$
$\|pg-gp\|<\sigma_1/32,$  $pgp\in_{\sigma_1/34} D,$  and
$\tau(1-p)<\sigma_1/4$ for all $\tau\in T(A).$}}
{{Define $L': A\to pAp$ by $L'(a)=pap$ for all $a\in A.$ Since $D$ is amenable, there
is a unital \cp\, $L'': pAp\to D$ such that $\|L''(pgp)-pgp\|<\sigma_1/32$ for all $g\in {\cal G}.$
Put $L=L''\circ L'.$  Then $L$ is ${\cal G}$-$\sigma_1/4$-multiplicative.}}
We estimate, for all $\tau\in T(A)$ and for all $g\in {\cal G}$ (since $\tau((1-p)gp)=\tau(pg(1-p))=0$),
\beq\label{june2-n1}
\tau(L(g)-g)\approx_{\sigma_1/32} \tau((1-p)g(1-p))\andeqn \tau((1-p)g(1-p))<\sigma_1/4
\eneq
{{Furthermore,
we may assume that $L_0: A\to (1-p)A(1-p)$ defined by
$L_0(a)=(1-p)a(1-p)$ is also ${\cal G}$-$\sigma_1/2$-multiplicative.}}
By the choice of $\dt$ and ${\cal G},$ it follows from Lemma 9.4 of \cite{LinTAI} that, for each $\tau\in  T(D),$ there is $\gamma'(\tau)\in T(A)$ such that
\begin{equation}\label{june2-n2}
|\tau(L(f))-\gamma'(\tau)(f)|<\sigma/16\rforal {{f}}\in {\cal F}.
\end{equation}
Applying {{Lemma}} \ref{TtoDelta},  one obtains $t_1, t_2,...,t_n\in \partial_e{T(D)}$ and
a continuous affine map $\lambda: T(D)\to {\triangle}$
such that
\begin{equation}\label{june6-n1}
|\tau({{d}})-\lambda(\tau)({{d}})|<{{\sigma}}/16\rforal \tau\in T(D) \andeqn {{d}}\in {{L}}( {\cal F}),
\end{equation}
where $\triangle$
is the convex hull of $\{t_1, t_2,...,t_n\}.$
Define an affine map $\lambda_1:
\triangle\to T(A)$ by
\begin{equation}\label{june6-n2}
\lambda_1(t_i)=\gamma'(t_i),\,\,\, i=1,2,...,m.
\end{equation}
Define $\gamma=\af\circ \lambda_1\circ \lambda.$
 Then, {{for $1\le j\le n,$  by \eqref{june2-n2}, \eqref{june6-n1}
 (recall $\tau(1-p)<\sigma_1/4$), }}
\begin{eqnarray*}
{{\gamma({{\imath}}(\tau))(h_j)}}&=&{{\alpha\circ \lambda_1\circ\lambda(\imath(\tau))(h_j)}}
= \lambda_1\circ\lambda({{\imath}}(\tau))(f_j)\\
&\approx_{\sigma/16}& \lambda(\imath(\tau))(L(f_j))
\approx_{\sigma/16}\imath(\tau)({{L(f_j)}})\\
&{{=}}&{{(1/\tau(p))\tau(L(f_j))\approx_{(9\sigma_1/32)/(1-\sigma_1/4)}
(1/\tau(p))\tau(f_j)}}\hspace{0.4in} {{{\rm ( by\,\, \eqref{june2-n1})}}}\\
&\approx_{(\sigma_1/2)/(1-\sigma_1/4)}&
\tau(f_j)= \alpha(\tau)(h_j),
\end{eqnarray*}
and this proves  (1).
Note that it follows from the construction that $\gamma(\tau)\in\alpha(T(A))$, and hence  (2)
also holds.

{{To show that (3) holds, let ${\cal P}=\{p_1, p_2,..,p_N\}$ be a finite subset of non-zero projections
such that $\{[p_1], [p_2],...,[p_N]\}$  generates $K_0(C)_+$ as {\blue{a positive cone}}.
Let $F_0=\kappa(K_0(C))$ and $(F_0)_+=\kappa(K_0(C)_+).$
Replacing $C$ by $M_N(C)$ and $A$ by $M_N(A),$
we
may assume that $\kappa({\cal P})$ is represented by a finite subset ${\cal Q}$ of projections in $A.$
Denote by $q_p$ a projection in $A$ such that $[q_p]=\kappa([p])$ for $p\in {\cal P}.$
In the proof of (1) above, we may choose that ${\cal G}_1\supset {\cal Q}.$
It is clear, in proof of  (1),  that, with sufficiently large ${\cal G}_1$ and small $\sigma_1/2,$ we may assume
that $[L_0]|_{\kappa(K_0(C))}$ and $[L]|_{\kappa(K_0(C))}$ are well defined, and $[L_0 (q)]\ge 0$  and $[L(q)]\ge 0$ for all
$q\in {\cal Q}.$ In other words,
$[L_0]|_{(F_0)_+}$ and $[L]|_{(F_0)_+}$ are a positive.  Moreover,
$[{\rm id}_A]|_{F_0}=[L_0]|_{F_0}+[\imath\circ L]|_{F_0}.$}}
{{Define $\kappa_0:=[L_0]\circ \kappa$ and $\kappa_1:=[L]\circ \kappa.$
Then $\kappa=\kappa_0+\imath_{*0}\circ \kappa_1.$  {{Since $\kappa([1_C])=[1_A]$ and
$L$ is unital, $\kappa_1([1_C])=[1_D].$}}
Thus, by case (1), it remains to show that $\kappa_1$ is strictly positive.
Let
$$
\sigma_0=\min\{\af(\tau)(p): p\in {\cal P}\andeqn \tau\in T(A) \}.
$$
Note that $\rho_A(\kappa([p]))(\tau)=\af(\tau)(p)$ for all $\tau\in T(A)$ and {{all}}
$p\in {\cal P}.$ Since $\kappa$ is strictly positive,
$\af(\tau)(p)>0$ for all $\tau\in T(A).$ Since $T(A)$ is compact and ${\cal P}$ is finite,  it follows that $\sigma_0>0.$}}

{{In the proof of  (1) above, choose $\sigma<\sigma_0/4$ and ${\cal H}\supset {\cal P}.$
Then,  by the proof of (1) (see \eqref{june2-n1}), for any $p\in {\cal P}\setminus \{0\},$
\beq
\tau(\imath\circ L(q_p))>3\sigma_0/4\rforal \tau\in T(A).
\eneq
There is also  a projection $q'\in A$ such that $\|\imath\circ L(q_p)-q'\|<\sigma_1$ and $[q']=[q_p]$
as
$L$ is contractive and ${\cal G}_1$-$\sigma_1/2$-multiplicative.}}
{{It follows that $\imath_{*0}\circ \kappa_1([p])>0$ for all $p\in {\cal P}.$
Let $x\in K_0(C)_+\setminus \{0\}.$
Since the positive cone $K_0(C)_+$ is finitely generated by $\{[p]: p\in {\cal P}\},$
one may write $x=\sum_{i=1}^l m_i [p_i]\in K_0(C)_+,$  where $m_i\in \Z_+,$ and
for some $i,$  $m_i>0.$   Therefore
$\imath_{*0}\circ \kappa_1(x)>0.$ It follows that $\kappa_1$ is strictly positive.}}


{To see that {{(4) holds,  we first assume that (3) holds.}} {{N}}ote {\blue{that
we}} may choose $D\subset A\otimes 1_U.$ Choose a projection
$e\in U$ such that
$$
0<t_0(e)<\dt_0<\dt-\max\{|\gamma(\tau)(p)-\tau(\kappa_1([p])|: p\in {\cal P} \tand\, \tau\in T(D)\},
$$
where $t_0$ is the unique tracial state of $U$
{{and $[1-e]$ is divisible by $K.$}} We then replace $\kappa_1$ by  $\kappa_2: K_0(A)\to K_0(D_2),$
where $D_2=D\otimes (1-e)$ and $\kappa_2([p])=\kappa_1([p])\otimes [1-e].$
Define $\kappa_3([p])=\kappa_1([p])\otimes [e].$
Then let $\kappa_4: K_0(C)\to K_0({{(1-(p\otimes (1-e)))A(1-(p\otimes (1-e))}}))$ be defined by
$\kappa_4=\kappa_0+{{\imath_{*0}}}\circ \kappa_3,$ where $\imath: D\otimes e\to A\otimes U\cong A$ is the embedding.
We then replace $\kappa_0$ by $\kappa_4$ {{and  $\kappa_1$ by $\kappa_2.$}} Note that, now, $\kappa_4$ is strictly positive.}
\end{proof}




\section{Maps from homogeneous \CA s to \CA s in ${\cal C}.$}

{The proof of the following lemma is similar to that of Theorem \ref{ExtTraceC-D}.}
%


\begin{lem}\label{ExtTraceH}
Let $X$ be a connected finite CW-complex.
Let $\mathcal H\subset C(X)$ be a finite subset, and let $\sigma>0$. There exists a finite subset $\mathcal H_{1, 1}\subset C(X)_+$ {{with}} the following property: for any $\sigma_{1, 1}>0$, there is a finite subset $\mathcal H_{1, 2}\subset C(X)_+$ {{satisfying }}  the following condition: for any $\sigma_{1, 2}>0$, there is a positive integer $M$ such that for any $D\in\mathcal C$ with the dimension of any irreducible representation of $D$ at least $M$, for any continuous affine map $\gamma: T(D)\to T(C(X))$ satisfying
\beq\nonumber
\gamma(\tau)(h) > \sigma_{1, 1}\tforal h\in \mathcal H_{1, 1}\tand\tau\in T(D),\tand\\
\gamma(\tau)(h) > \sigma_{1, 2} \tforal h\in \mathcal H_{1, 2}\tand\tau\in T(D),
\eneq
 there is a homomorphism $\phi: C(X)\to D$ such that
$$|\tau\circ\phi(h)-\gamma(\tau)(h)|<\sigma\tforal h\in\mathcal H {{\tand \tau\in T(D)}}.$$
Moreover, if $D\in {\cal C}_0,$ then there is a point evaluation $\Psi: C(X)\to D$ such that
$[\phi]=[\Psi].$
\end{lem}

\begin{proof}
{{We}} may assume that every element of $\mathcal H$ has norm at most $1$.

 Let $\eta>0$ {{be}} such that for any $f\in\mathcal H$ and any $x, x'\in X$ with $d(x, x')<\eta$, one has
$$|f(x)-f(x')|<\sigma/4.$$

 Since $X$ is compact, { {one can choose}} a finite subset $\mathcal H_{1, 1}\subset C(X)_+$ such that { {for any open ball $O_{\eta/24} \subset X$,  of radius $\eta/24$, there is {{a}} non-zero element $h\in  \mathcal H_{1, 1}$ with ${\rm supp}(h) \subset O_{\eta/24}$. We assume that $\|h\|\leq 1 $ for all $h\in \mathcal H_{1, 1}$. Consequently,} }
   if there is $\sigma_{1, 1}>0$ such that $$\tau(h)>\sigma_{1, 1}\rforal h\in\mathcal H_{1, 1},$$ then $$\mu_{\tau}(O_{\eta/24})>\sigma_{1, 1}$$ for any open ball $O_{\eta/24}$ with radius $\eta/24,$ where $\mu_\tau$ is the probability measure
induced by $\tau.$

{ { Fix $\sigma_{1,1} >0$. }}  Let $\delta$
and $\mathcal G\subset {{C(X)}}$ ({\blue{in  place}} of $\mathcal G$) be the constant and finite subset {{provided by}}  Lemma 6.2 of \cite{LnTAMS12} {{for}}  $\sigma/2$ ({\blue{in  place}} of $\epsilon$), $\mathcal H$ (in place of $\mathcal F$), and $\sigma_{1, 1}/\eta$ (in place of $\sigma$).

Let $\mathcal H_{1, 2}\subset {{C(X)_+}}$ (in  place of $\mathcal H_1$) be the finite subset provided by  Theorem \ref{uniCMn}
{{for}}  $\delta$ (in  place of $\epsilon$) and $\mathcal G$ (in place of $\mathcal F$).
{{We may assume that ${\cal H}_{1,1}\subset {\cal H}_{1,2}.$}}

Let $\sigma_{1, 2}>0$. Then let $\mathcal H_2\subset {{C(X)}}$ ({\blue{in place}} of $\mathcal H_2$) and {$\sigma_2$}
be the finite subset, and positive constant {{provided by}} Theorem \ref{uniCMn} {{for}} $\sigma_{1, 2}$ { {and $\mathcal H_{1, 2}$ }}(in place of $\sigma_1$ and $\mathcal H_1$).

Let $M$ (in place of $N$) be the constant {{provided by}} {{Corollary 2.5}} of \cite{Li-interval} {{for}}  $\mathcal H_2\cup \mathcal H_{1, 2}$
(in place of $F$) and $\min\{\sigma/4, \sigma_2/4, \sigma_{1, 2}/2, \sigma_{1, 1}/2\}$ (in  place of $\epsilon$).

Let $D=D(F_1, F_2, \psi_0, \psi_1)$ be a \CA\,  in $\mathcal{C}$ with the dimensions of {{its}}  irreducible representations at least $M$. {{Write $F_1=M_{R(1)}\oplus\cdots\oplus M_{R(l)}$ and $F_2=M_{r(1)}\oplus M_{r(2)}\oplus\cdots \oplus M_{r(k)}.$
Denote by $q_i: F_2\to M_{r(i)}$ the quotient map, and {{write}}
$\psi_{0,i}=q_i\circ \psi_0$ and $\psi_{1,i}=q_i\circ \psi_1.$}}

 Let $\gamma: T(D)\to T(C(X))$ be a map {{as in}} the lemma. Write $C([0, 1], F_2)=I_1\oplus\cdots\oplus I_k$ with $I_i=C([0, 1], M_{r(i)})$, $i=1, ..., k$. Then $\gamma$ induces a continuous {\blue{affine}} map $\gamma_i: T(I_i)\to T(C(X))$ by $\gamma_i(\tau)=\gamma(\tau\circ\pi^{{I_i}})$, where $\pi^{{I_i}}$ is the restriction map $D\to I_i$. It is then clear that
$$\gamma_i(\tau)(h)>\sigma_{1, 2}\rforal {{~~h\in \mathcal H_{1,2},~~}}\tau\in T(I_i).$$

{{D}}enote by $\pi_{{e,j}}: D\to M_{{R(j)}}$ the  composition of $\pi_e$ {{(see \ref{DfC1})}}  with the projection from $F_1$ onto $M_{R(j)}.$
By {{Corollary 2.5 }} of \cite{Li-interval}, for each $1\leq i\leq k$, there is a homomorphism $\phi_i:
C(X)\to I_i$ such that
\begin{equation}\label{app-tr-001}
|\tau\circ\phi_i(h)- \gamma_i(\tau)(h) | < \min\{\sigma/4, \sigma_2/4, \sigma_{1, 2}/2, \sigma_{1, 1}/2\}
\end{equation}
for all $h\in\mathcal H_2\cup\mathcal H_{1, 2}$
{{and for all $\tau\in T(I_i),$}} and for any $j$, there is also a homomorphism $\phi'_j: C(X)\to M_{R(j)}$ such that
\begin{equation}\label{app-tr-002}
|\mathrm{tr}_{R(j)}\circ\phi'_j(h)- \gamma\circ(\pi_{{e,j}})^*(\mathrm{tr}_{R(j)})(h) | < \min\{\sigma/4, \sigma_2/4, \sigma_{1, 2}/2, \sigma_{1, 1}/2\}
\end{equation}
for all $h\in\mathcal H_2\cup\mathcal H_{1, 2}\cup\mathcal H_{1, 1}$ {{(we use
${\rm tr}_{m}$ for  the tracial state of $M_{m}$).}}
{{Set}} $\phi'=\bigoplus_{j}\phi'_j$ and  {{let}}  $\pi_t: I_i\to M_{r(i)}$ be the point evaluation at $t\in [0,1].$
{{Then}}
\beq\nonumber
|\mathrm{tr}_{r(i)}\circ(\psi_{0, i}\circ \phi'){{(h)}}-\mathrm{tr}_{r(i)}\circ(\pi_0\circ\phi_i){{(h)}}|\leq \sigma_2/2\rforal h\in\mathcal H_2\andeqn\\\label{app-tr-002+}
\mathrm{tr}_{r(i)}\circ(\psi_{0, i}\circ\phi')(h)\geq \sigma_{1, 2}/2\quad\mathrm{and}\quad \mathrm{tr}_{r(i)}\circ(\pi_{0}\circ\phi_i)(h)\geq \sigma_{1, 2}/2\rforal h\in \mathcal H_{1, 2}.
\eneq
By Theorem \ref{uniCMn}, there is a unitary $u_{i, 0}\in M_{r(i)}$ such that
$$|| \mathrm{Ad}u_{i, 0}\circ\pi_{ 0}\circ\phi_i(f)-\psi_{0, i}\circ\phi'(f)||<\delta \rforal f\in\mathcal G.$$
Exactly the same argument shows that
there is a unitary $u_{i, 1}\in M_{r(i)}$ such that
$$\|\mathrm{Ad}u_{i, 1}\circ \pi_{1}\circ \phi_i(f) -  \psi_{1, i}\circ\phi'(f)\|<\dt\rforal f\in\mathcal {{G}}.$$

Choose two  paths of unitaries,  $\{u_{i, 0}(t):t\in [0,1/2]\}\subset M_{r(i)} $
such that $u_{i, 0}(0)=u_{i, 0}$ and $u_{i, 0}(1/2)=1_{M_{r(i)}},$ and
$\{u_{i,1}(t):t\in [1/2, 1]\}\subset M_{r(i)}$ such that
$u_{i,1}(1/2)=1_{M_{r(i)}}$ and $u_{i,1}(1)=u_{i,1}.$
Put $u_i(t)=u_{i,0}(t)$ if $t\in [0,1/2)$ and $u_i(t)=u_{i,1}(t)$ if $t\in [1/2,1].$
Define ${\tilde \phi}_i: C(X)\to I_i$ by
\begin{equation*}
\pi_t\circ {\tilde \phi_i}={\rm Ad}\, u_i(t)\circ \pi_t\circ \phi_i.
\end{equation*}

Then
\begin{equation}
\| \pi_{0}\circ {\tilde \phi}_i(f) -  \psi_{0, i}\circ\phi'(f)\|<\delta
\andeqn
\| \pi_{1}\circ {\tilde \phi}_i(f) -  \psi_{1, i}\circ\phi'(f)\|<\delta
\end{equation}
 for all $f\in\mathcal G, i=1,..., k.$

Note that it also follows from \eqref{app-tr-001} and \eqref{app-tr-002+} that
$$\mathrm{tr}_{r(i)}\circ(\psi_{0, i}\circ\phi')(h)\geq \sigma_{1, 1}/2\quad\mathrm{and}\quad \mathrm{tr}_{r(i)}\circ(\pi_{ 0}\circ{\tilde \phi}_i)(h)\geq \sigma_{1, 1}/2\rforal h\in \mathcal H_{1, 1}.$$ Hence,
$$
\mu_{{{{\rm tr}_{r(i)}}}\circ (\psi_{0, i}\circ\phi')}(O_{\eta/24})\geq\sigma_{1,1}\quad\mathrm{and}\quad \mu_{{{\rm tr}_{r(i)}}\circ (\pi_{0}\circ\tilde \phi_i)}(O_{\eta/24})\geq\sigma_{1,1}.
$$
Thus, by Lemma 6.2 of \cite{LnTAMS12}, for each $1\leq i\leq k$, there are two unital homomorphisms
$$\Phi_{0, i}, \Phi_{0, i}': C(X)\to C([0, 1], {\blue{M_{r(i)}}})$$ such that
$$\pi_0\circ\Phi_{0, i}=\psi_{0, i}\circ\phi',\quad \pi_0\circ\Phi_{0, i}'=\pi_{ 0}\circ\tilde\phi_i,$$
$$||\pi_t\circ\Phi_{0, i}(f)-\psi_{0, i}\circ\phi'(f)||<\sigma/2,\quad  ||\pi_t\circ\Phi_{0, i}'(f)-\pi_{0}\circ{\tilde \phi_i}{{(f)}} ||<\sigma/2$$
for all $f\in\mathcal H$ and $t\in[0, 1]$, and there is a unitary $w_{i, 0}\in {\blue{M_{r(i)}}}$ (in place of $u$) such that
$$\pi_1\circ\Phi_{0, i}=\mathrm{Ad}w_{i, 0}\circ \pi_1\circ \Phi_{0, i}'.$$

The same argument shows that, for each $1\le i\le k,$ there are two unital \hm s
$\Phi_{1,i}, \Phi_{1,i}': C(X)\to C([0,1], {\blue{M_{r(i)}}})$ such that
$$\pi_1\circ\Phi_{1, i}=\psi_{1, i}\circ\phi',\quad \pi_1\circ\Phi_{1, i}'=\pi_{1}\circ\tilde\phi_i,$$
$$||\pi_t\circ\Phi_{1, i}(f)-\psi_{1, i}\circ\phi'(f)||<\sigma/2,\quad  ||\pi_t\circ\Phi_{1, i}'(f)-\pi_{0}\circ {\tilde \phi}_i{{(f)}} ||<\sigma/2$$
for all $f\in\mathcal H$ and $t\in[0, 1]$, and there is a unitary $w_{i, 1}\in {\blue{M_{r(i)}}}$ (in place of $u$) such that
$$\pi_0\circ\Phi_{1, i}=\mathrm{Ad}w_{i, 1}\circ \pi_0\circ \Phi_{1, i}'.$$

Choose two continuous paths $\{w_{i,0}(t): t\in [0,1]\},\,\{w_{i,1}(t): t\in [0,1]\}$
in ${\blue{M_{r(i)}}}$ such that $w_{i,0}(0)=w_{i,0},\,\,\, w_{i,0}(1)=1_{\blue{M_{r(i)}}}$ and
$w_{i,1}(1)=1_{\blue{M_{r(i)}}}$ and $w_{i,1}(0)=w_{i,1}.$

For each $1\leq i\leq k$, by the continuity of $\gamma_i,$ there is  $1>\epsilon_i>0$ such that
$$|\gamma_i(\tau_x)(h)-\gamma_i(\tau_y)(h)|<\sigma/4 \rforal h\in\mathcal H,$$
provided that $|x-y|<\epsilon_i,$ where $\tau_x$ and $\tau_y$ are the extremal {\blue{traces}} of $I_i$ concentrated on {\blue{$x$ and $y$}}.

Define the map $\tilde{\tilde{\phi}}_i: C(X)\to I_i$ by
\begin{displaymath}
\pi_t\circ{\tilde{ \tilde{\phi}}}_i=
\left\{
\begin{array}{ll}
\pi_{\textstyle\frac{3t}{\epsilon_i}}\circ\Phi_{0, i}, & t\in[0, \epsilon_i/3),\\
\mathrm{Ad}(w_{i, 0}({\frac{\textstyle 3t}{\textstyle \epsilon_i}}-1))\circ \pi_1\circ\Phi_{0, i}', & t\in[\epsilon_i/3, 2\epsilon_i/3),\\
\pi_{3-{\textstyle \frac{3t}{\epsilon_i}}}\circ\Phi_{0, i}', & t\in[2\epsilon_i/3, \epsilon_i),\\
\pi_{\textstyle \frac{t-\epsilon_i}{1-2\epsilon_i}}
\circ{\tilde \phi}_i, & t\in[\epsilon_i, {{1-\ep_i}}),\\
\pi_{\textstyle{1 -2\ep_i/3-t\over{{{\ep}}_i/3}}}\circ \Phi_{1,i}' , &   t\in [{{1-\ep_i}}, 1-2\ep_i/3],\\
\mathrm{Ad}(w_{i, 1}({\frac{\textstyle (1-\ep_i/3)-t}{\textstyle \epsilon_i/3}}))\circ \pi_0\circ\Phi_{1, i}', & t\in[1-2\epsilon_i/3, 1-\epsilon_i/3],\\
\pi_{\textstyle \frac{t-1+\ep_i/3}{\epsilon_i/3}}\circ\Phi_{1, i}, & t\in[1-\epsilon_i/3,1].
\end{array}
\right.
\end{displaymath}
Then,
\begin{equation}\label{match-i}
\pi_0\circ\tilde{\tilde{\phi}}_i=\psi_{0, i}\circ\phi' \quad\mathrm{and}\quad \pi_1\circ\tilde{\tilde{\phi}}_i=\psi_{1, i}\circ
\phi'.
\end{equation}
One can also estimate, by the choice of $\ep_i$ and the definition of ${\tilde{\tilde{\phi}}}_i,$ that
\begin{equation}\label{app-tr-003}
|\tau_t\circ \tilde{\tilde{\phi}}_i(h)-\gamma_i(\tau_t)(h)|<\sigma\rforal t\in[0, 1] {{\andeqn h\in {\cal H}}},
\end{equation} where $\tau_t$ is the extremal tracial state of $I_i$ concentrated on $t\in[0, 1].$

Define $\Phi: {{C(X)}}\to C([0,1], F_2)$ by
$\Phi(f)=\bigoplus_{i=1}^k {\tilde{\tilde{\phi}}}_i(f)$ for all $f\in {{C(X)}}.$
Define $\phi: {{C(X)}}\to C([0,1],F_2)\oplus F_{ 1}$ by $\phi(f)=(\Phi(f),\phi'(f))$ for $f\in C(X).$
By \eqref{match-i}, $\phi$ is a \hm\, from ${{C(X)}}$ to $D.$
By \eqref{app-tr-003} and \eqref{app-tr-002}, one has that
$$|\tau\circ\phi({{h}})-\gamma(\tau)(h)|<\sigma\rforal h\in\mathcal H$$
and for all $\tau\in T(D),$
as desired.

To see the last part of the lemma, one assumes that  $D\in {\cal C}_0.$
Consider $\pi_e\circ \phi: C(X)\to F_1$ (where $D=D(F_1, F_2, \psi_0, \psi_1)$ as above).
 Since $\pi_e\circ \phi$ has finite dimensional range, it is a point evaluation.
 We may write $\pi_e\circ \phi(f)=\sum_{i=1}^m f(x_i)p_i$ for all $f\in C(X),$ where
 $\{x_1, x_2,...,x_m\}\subset X$ and $\{p_1, p_2,..., p_m\}\subset F_1$ is a set of mutually orthogonal projections.
Let $I=\{f\in C(X): f(x_1)=0\}$   and {{let}} $\imath: I\to C(X)$ be the embedding.
 It follows that $[\pi_e\circ \phi\circ \imath]=0.$ By
 \ref{2Lg13}, $[\pi_e]$ is injective on each $K_i({{D}})$ and {{on}} each $K_i({{D}}, \Z/k\Z)$ ($k\ge 2,$ $i=0,1$). Hence $[\phi\circ \imath]=0.$  Choose $\Psi(f)=f(x)\cdot 1_{{D}}.$ Since $X$ is connected,
$[\phi]=[\Psi]$ (see the end of Remark \ref{remark831KL}).
\end{proof}

\begin{cor}\label{ExtTraceH-D}
Let $X$ be a connected finite CW-complex.
Let $\Delta: C(X)_+^{q,  {\bf 1}}\setminus\{0\}\to (0, 1)$ be an order preserving map. Let $\mathcal H\subset C(X)$ be a finite subset and let $\sigma>0$. Then there {exist} a finite subset $\mathcal H_1\subset C(X)_+^{\bf 1}\setminus \{0\}$ and a positive integer $M$ such that  for any $D\in\mathcal C(X)$ with the dimension of any irreducible representation of $D$ at least $M,$ and for any continuous affine map $\gamma: T(D)\to T(C(X))$ satisfying $$\gamma(\tau)(h) > \Delta(\hat{h})\tforal h\in \mathcal H_{1}\tand \tforal \tau\in T(D),$$ there is a homomorphism $\phi: C(X)\to D$ such that
$$|\tau\circ\phi(h)-\gamma(\tau)(h)|<\sigma\tforal h\in\mathcal H.$$
\end{cor}
\begin{proof}
Let $\mathcal H_{1, 1}$ be the subset of Lemma \ref{ExtTraceH} with respect to $\mathcal H$ and $\sigma$. Then put
$$\sigma_{1, 1}=\min\{\Delta(\hat{h}):\ h\in\mathcal H_{1, 1}\}.$$
Let $\mathcal H_{1, 2}$ be the finite subset of Lemma \ref{ExtTraceH} with respect to $\sigma_{1, 1}$, and then put
$$\sigma_{1, 2}=\min\{\Delta(\hat{h}):\ h\in\mathcal H_{1, 2}\}.$$
Let $M$ be the positive integer of Lemma \ref{ExtTraceH} with respect to $\sigma_{1, 2}$. Then it follows from Lemma \ref{ExtTraceH} that the finite subset $$\mathcal H_1:=\mathcal H_{1, 1}\cup\mathcal H_{1, 2}$$ and the positive integer $M$ {\blue{are as desired in }} the corollary.
\end{proof}

\begin{thm}\label{istTr}
{{Let
$X$ be a  connected  finite CW complex,}}  and let $A\in {\cal B}_0$ be a unital separable simple \CA.
Suppose that $\gamma: T(A)\to T_f({{C(X)}}))$  {{(see Definition \ref{Aq})}} is a continuous affine map.
Then, for any $\sigma>0,$ {\blue{and}} any finite subset ${\cal H}\subset {{C(X)}}_{s.a.},$ there exists a
unital \hm\, $h: {{C(X)}}\to A$ such that
\beq\label{istTr-1}
[h]=[\Psi]{\blue{\in KL(C(X),A)}}\andeqn
|\tau\circ h(f)-\gamma(\tau)(f)|<\sigma\tforal f\in {\cal H} {\blue{\tand \tau\in T(A),}}
\eneq
where {\blue{$\Psi$ is a homomorphism with a finite dimensional image}}.
\end{thm}

\begin{proof}
{{We}} may assume that every element of $\mathcal H$ has norm at most one.
Let $\mathcal H_{1, 1}{{\subset C(X)_+\setminus \{0\}}}$  be the finite subset of Lemma \ref{ExtTraceH} with respect to $\mathcal H$ (in place of $\mathcal H$), $\sigma/4$ (in place of $\sigma$), and $C(X).$
Since $\gamma(T(A))\subset {\blue{T_{f}}}(C(X))$, there is $\sigma_{1, 1}>0$ such that
$$\gamma(\tau)(h)>\sigma_{1, 1}{,}\rforal h\in\mathcal H_{1, 1} {\textrm{ and}} \rforal \tau\in T(A).$$

Let $\mathcal H_{1, 2}\subset C(X)_+{{\setminus\{0\}}}$
 be the finite subset of Lemma \ref{ExtTraceH} with respect to $\sigma_{1, 1}$. Again, since $\gamma(T(A))\subset {\blue{T_{f}}}(C(X))$, there is $\sigma_{1, 2}>0$ such that
$$\gamma(\tau)(h)>\sigma_{1, 2}{,}\rforal h\in\mathcal H_{1, 2}\ {\textrm{and}} \rforal \tau\in T(A).$$

Let $M$
be the constant of Lemma \ref{ExtTraceH} with respect to $\sigma_{1, 2}.$
{ {Note}} that $A\in \mathcal B_0$. By (1) and (2)
of Lemma \ref{cut-trace}, there {{are}}  a $C^*$-subalgebra $D{{\subset}} A$ with $D\in\mathcal C_0$,
{{and}} a continuous affine map $\gamma': T(D)\to T(C(X))$ such that
\begin{equation}\label{istTr-eq1}
|\gamma'(\frac{1}{\tau(p)}\tau|_D)(f)-\gamma(\tau)(f)|<\sigma/4{,}\rforal \tau\in T(A){,} \rforal f\in\mathcal H,
\end{equation}
where $p=1_D$, $\tau(1-p)<\sigma/(4+\sigma)$ {{for all $\tau\in T(A),$}}
\beq\label{istTr-eq2}
&&\gamma'(\tau)(h)>\sigma_{1, 1}{,}\rforal \tau\in T(D)\rforal h\in\mathcal H_{1, 1}{\blue{,\andeqn}}\\
\label{istTr-eq3}
&&\gamma'(\tau)(h)>\sigma_{1, 2}{,}\rforal \tau\in T(D)\rforal h\in\mathcal H_{1, 2}.
\eneq

Moreover, since $A$ is simple, one may assume that the dimension of any irreducible representation of $D$ is at least $M$ (see \ref{Affon1}). Thus, by \eqref{istTr-eq2} and  \eqref{istTr-eq3}, one applies Lemma \ref{ExtTraceH} to $D$, $C(X)$, and $\gamma'$ (in  place of $\gamma$) to obtain a homomorphism $\phi: C(X)\to D$ such that
\begin{equation}\label{istTr-eq4}
|\tau\circ\phi(f)-\gamma'(\tau)(f)|<\sigma/4{,}\rforal f\in\mathcal H {\textrm{\ and}} \rforal \tau\in T(D).
\end{equation}

Moreover, we may assume that $[\phi]=[\Phi_0] \in KL(C(X),D)$ for {\blue{a homomorphism}} $\Phi_0: C(X)\to D$ {\blue{with a finite dimensional image}} since we assume that $D\in {\cal C}_0.$ Pick a point $x\in X$, and define $h: C(X)\to A$ by $$f \mapsto f(x)(1-p)\oplus\phi(f)\rforal f\in C(X).$$
For any $f\in\mathcal H$ {{and any $\tau\in T(A),$}}  one has
\begin{eqnarray*}
&&\hspace{-0.2in}|\tau\circ h(f)-\gamma(\tau)(f)\| \leq  |\tau\circ \phi(f)-\gamma(\tau)(f)|+\sigma/4
<|\tau\circ \phi(f)-\gamma'(\frac{1}{\tau(p)}\tau|_D)(f)|+\sigma/2\\
&&<|\tau\circ \phi(f)- \frac{1}{\tau(p)}\tau\circ\phi(f)|+3\sigma/4<\sigma.
\end{eqnarray*}
Define ${\blue{\Psi}}: C(X)\to A$ by ${\blue{\Psi}}(f)=f(x)(1-p)\oplus \Phi_0(f)$ for all $f\in C(X).$ Then $[h]=[{\blue{\Psi}}].$
\end{proof}




\section{KK-attainability
 of the building blocks}\label{KK-BBlock}

\begin{df}[9.1 of \cite{LinTAI}]\label{D181}
Let ${\cal D}$ be a class of unital \CA s.
A \CA\, $C$ is said to be KK-attainable with {{respect}}  to ${\cal D}$ if for any $A\in {\cal D}$ and any $\alpha\in {{ {\rm{Hom}}_{\Lambda}(\underline{K}(C), \underline{K}(A))}}^{++}$, there exists a sequence of
\morp s
$L_n: C\to {{A\otimes {\cal K}}}$ such that
\beq
&&\lim_{n\to\infty} \|L_n(ab)-L_n(a)L_n(b)\|=0\rforal a,\, b\in C\andeqn\\
&&{[}{{\{}}L_n{{\}}}{]}=\af.
\eneq
{{(The latter means that, for any finite subset ${\cal P}\subset \underline{K}(C),$  $[L_n]|_{\cal P}=\af|_{\cal P}$ for all large
$n.$)}}
{{If  $C$ satisfies the UCT, then ${\rm{Hom}}_{\Lambda}(\underline{K}(C), \underline{K}(A))^{++}$
may be replaced by $KL(C,A)^{++}.$
}}
In what follows, we will use ${\cal B}_{ui}$\index{${\cal B}_{u0}$}  {to denote} the class of those  {{separable}} 
 \CA s of the form\index{${\cal B}_{u1}$}
$A\otimes U,$ where $A\in {\cal B}_i$ and  $U$ is {a} UHF-algebra of infinite type, {{$i=0,1.$}}

\end{df}

\begin{thm}{\rm (Theorem 5.9 of \cite{LinTAF2}; see 6.1.11 of \cite{Lnbok})}\label{kkmaps}
Let $A$ be a separable \CA\, satisfying {{the}} UCT and let $B$ be a {{unital}}  separable \CA. Assume that $A$ is the closure of an increasing sequence $\{A_n\}$ of amenable residually finite dimensional \SCA s. Then for any $\alpha\in KL(A, B)$, there exist two sequences of completely positive contractions $\phi_n^{(i)}: A\rightarrow{B}\otimes\mathcal K\ (i=1, 2)$ satisfying the following
{{conditions}}:
\begin{enumerate}
\item $\|{\phi_n^{(i)}(ab)-\phi_n^{(i)}(a)\phi_n^{(i)}(b)}\|\rightarrow0$
      as $n\rightarrow\infty$;
\item {{for each finite subset ${\cal P},$ $[\phi_n^{(i)}]|_{\mathcal P}$ is well defined for
      sufficiently large $n$, $i=1,2,$}}
 and, for any $n$, the {{image}}  of $\phi_n^{(2)}$ {{is}}  contained in a finite
      dimensional sub-\CA\, of $B\otimes\mathcal K;$
      \item for each finite subset $\mathcal P \subset\underline{K}(A)$, there exists $m>0$
      such that $$[\phi_n^{(1)}]|_\mathcal{P}=(\alpha+[\phi_n^{(2)}])|_{\mathcal P} \tforal n\ge m;$$
\item for each $n$,
$\phi_n^{(2)}$ is a homomorphism
      on $A_n$.
\end{enumerate}
{\Green{Moreover, the condition that $A_n$ is amenable could be 
replaced by the condition that $B$ is amenable.}}
\end{thm}
{{Note that $\af$ does not need to be positive as stated in 6.1.11 of \cite{Lnbok}.  In fact, the proof
of Theorem 5.9 of \cite{LinTAF2} does not require that $\af$ {{be}}  positive, if one does not require the
part (5) of Theorem 5.9 of \cite{LinTAF2} (there was also a typo which is corrected above).}}

\begin{lem}\label{Nextension}
Let $C={{A}} (F_1, F_2, \phi_0, \phi_1)\in {\cal C}$ {{and let $N_1, N_2\ge 1$ be integers.
There exists $\sigma>0$ satisfying the following condition:}}
{{L}}et $A\in {\cal C}$  be another \CA\, {{and}}
let  $\kappa: K_0(C)\to K_0(A)$ be an order  preserving \hm\,
such that, for any non-zero element ${{x}}\in K_0(C)_+,$
$N_{{1}}\kappa({{x}})>[1_A]$ {{and $\kappa([1_C])\le N_2[1_A].$}}
{{Then, f}}or any $\tau\in T(A),$ there exists $t\in T(C)$ such that
\beq\label{Next-1}
t(h)\ge \sigma\int_{[0,1]} {\rm T}(\lambda(h)({\Green{s}}))d\mu({\Green{s}})\tforal h\in C_+\tand\\
{{{\rho_A(\kappa(x))(\tau)}}\over{{{\rho_A(\kappa([1_C]))(\tau)}}}}={{\rho_C(x)(t)}}\tforal x\in K_0(C)_+,
\eneq
where $\lambda: C\to C([0,1], F_2)$
is the \hm\,  given by \eqref{pull-back} and
${\rm T}(b)=\sum_{i=1}^k{\rm tr}_i({{\psi_i}}(b))$ for all $b\in F_2,$ where ${\rm tr}_i$ is
the normalized tracial state on the $i$-th simple {{direct}} summand {{$M_{r(i)}$}} of $F_2,$
{{$\psi_i: F_2\to M_{r(i)}$ is the projection map,}}
and
$\mu$ is Lebesgue measure on $[0,1].$
\end{lem}

\begin{proof}
{{In what follows, if $B$ is a \CA, $\tau\in T(B)$ and $x\in K_0(B),$
we will use the following convention:}}
$$
{{\tau(x):=\rho_B(x)(\tau).}}
$$
{{We also write $F_1=M_{R(1)}\bigoplus M_{R(2)}\bigoplus \cdots \bigoplus M_{R(l)}.$ Denote by
$q_j: F_1\to M_{R(j)}$ the projection map   and
${\rm {\bar{tr}}}_j$ the tracial state of $M_{R(j)}.$ Recall $\pi_e: C\to F_1$ is the quotient map defined in \ref{DfC1}.}}
{{Define $T: K_0(F_2)\to \R$ by $T(x)=\sum_{i=1}^k \rho_{M_{r(i)}}((\psi_i)_{*0}(x))({\rm tr}_i)$
for $x\in K_0(F_2),$  where $F_2=M_{r(1)}\bigoplus M_{r(2)}\bigoplus\cdots \bigoplus M_{r(k)}.$
{\Green{Note}} $T([\phi_0(\pi_e(1_C))])=k,$ where $\phi_0: F_1\to F_2$ is given by $C.$}}

{{Let $p_1,p_2,...,p_s$ be a set of minimal projections in $M_m(C)$ for some integer
$m$ such that $\{p_1, p_2,...,p_s\}$ generates $K_0(C)_+$ (see \ref{FG-Ratn}).
It follows that there is a $\sigma_{00}>0$ such that
\beq
\sigma_{00}mk<1/2.
\eneq	
{{Set}} $\sigma_0=\sigma_{00}/2N_1$ and $\sigma=\sigma_0/N_2.$}}
Write $$A={{A(F'_1,F'_2, \phi'_0,\phi'_1)=}}\{(g,c)\in C([0,1], F'_2)\oplus F'_1: g(0)=\phi_0{{'}}(c)\andeqn g(1)=\phi_1{{'}}(c),\, \, c\in F'_1\}.$$ {{Denote by $\pi_e': A\to F_1'$  the quotient map (see \ref{DfC1}).}}
{{Write $F_1'=M_{R(1)'}\oplus M_{R(2)'}\oplus \cdots\oplus  M_{R(l')'},$ and let
$q_i': F_1'\to M_{R(i)'}$ be the projection map and ${\rm tr}_i'$ the tracial state of $M_{R(i)'},$  $i=1,2,...,l'.$}}

{{Note
$K_0(F_1')=\Z^{l'}.$}}
View ${{(q_i'\circ \pi_e')_{*0}}}\circ \kappa$ as a positive \hm\, from $K_0(C)\to \R$
{{(as well as a \hm\, from $K_0(C)$ to $K_0(M_{R(i)'})=\Z$).}}
Since $$N_{{1}}\kappa(x)>[1_A]\rforal x\in K_0(C)_+\setminus\{0\},$$
one has that
\beq\label{18-3n}
{{N_1((q_i'\circ \pi_e')_{*0}}}\circ\kappa)(x)>R(i)'>0\rforal x\in K_0(C)_+\setminus\{0\}.
\eneq
{{Then, for $p\in \{p_1,p_2,...,p_s\},$
\beq
&&\hspace{-0.3in}N_1((q_i'\circ \pi_e')_{*0}\circ\kappa)([p])-\sigma_{00} R(i)'T\circ (\phi_0\circ \pi_e)(p)\\
&&=
N_1((q_i'\circ \pi_e')_{*0}\circ\kappa)([p])-\sigma_{00}R(i)'(\sum_{i=1}^k \rho_{M_{r(i)}}({\Green{\psi_i(p(0)}}))({\rm tr}_i))\\
&&\ge N_1((q_i'\circ \pi_e')_{*0}\circ\kappa)([p])-\sigma_{00}R(i)'km>0.
\eneq}}
{{Define $\Gamma_i: K_0(C)\to \R$  by
$$\Gamma_i(x)=((q_i'\circ \pi_e')_{*0}\circ\kappa)(x)- \sigma_0R(i)'\cdot T\circ (\phi_0\circ \pi)_{*0}(x)\rforal x\in K_0(C).$$
Then $\Gamma_i$ is positive (since $\sigma_0<\sigma_{00}/N_1$).}}
{{Note that   $(\pi_e)_{*0}: K_0(C)\to K_0(F_1)$
is an order embedding (see \ref{2Lg13}).}}
It follows from {{Theorem 3.2 of}} \cite{GH-RankFunction} (see also 2.8 of \cite{Lnbirr}) that there are positive homomorphisms ${\tilde \Gamma}_i:
K_0(F_1)
\to \R$ such that ${\tilde \Gamma}_i\circ (\pi_e)_{*0}=\Gamma_i$
and
$
{\tilde\Gamma}_i([e_j])={{\af}}_{i,j}{{\ge}}0,\,\,\,j=1,2,...,l,\,\,\, i=1,2,...,l',
$
where $e_j=q_j\circ  \pi_e(1_C),$ $j=1,2,...,l.$
{\Green{Using  the fact that the \hm s from $K_0(F_1)=\Z^l$ is determined by
their values on the canonical basis, one gets}}
\beq
&&\hspace{-0.4in}\sum_{j=1}^{l}{{\rm rank}q_j(\pi_e(p))\over{{\rm rank}q_j(\pi_e(1_C))}}\af_{i,j}+\sum_{j=1}^{k}{{\rm rank}\psi_j(\phi_0\circ\pi_e(p))\over{{\rm rank}\psi_j(\phi_0\circ \pi_e(1_C))}}\sigma_0R(i)'\\\label{Nextension-n1}
&&=
{{\Gamma_i([p])+\sigma_0 R(i)'\cdot T(\phi_0\circ \pi_e)_{*0}([p])}}
={{(q_i'\circ \pi_e')_{*0}}}\circ \kappa([p])
\eneq
for any projection $p\in M_m(C),$ where $m\ge 1$ is an integer. {{In particular,
$(q_i'\circ \pi_e')_{*0}\circ \kappa([1_C])=\sum_{j=1}^l\af_{i,j}+k\sigma_0 R(i)'.$}}
For any $\tau\in T(A),$ by 2.8 of \cite{Lnbirr}, {{since $K_0(A)$ is order embedded into $K_0(F_1')$ (see \ref{2Lg13}),}} there is $\tau'\in T(F_1')$ (see also Corollary 3.4 of \cite{Blatrace})
such that
\beq\label{18-18n-2}
\tau'\circ \pi_e'(x)=\tau(x)\rforal x\in K_0(A).
\eneq
Write
$
\tau'=\sum_{i=1}^{l'} \lambda_{i, \tau} {{{\rm tr}_i'}}{{\circ q_i'}},
$
where $0\le \lambda_i\le 1$, $\sum_{i=1}^{l'}\lambda_{i, \tau}=1.$
{{Put }}
\beq\label{18-18n-1}
{{\bt=\tau(\kappa([1_C]))=\tau'((\pi_e')_{*0}(\kappa([1_C])))=\sum_{i=1}^{l'}\lambda_{i,\tau} (1/R(i)')((q_i'\circ \pi_e')_{*0}\circ \kappa([1_C]))}}
\eneq
{{(with the identification $K_0(M_{R(i)'})=\Z$). Then $\beta\le N_2.$}}
For each $i,$ define, {{ for $(s,b)\in C,$}}
$$
t_i((s,b))={{{1\over{\bt R(i)'}}}}\left(\sum_{j=1}^{k}
{{\sigma_0 R(i)'}}\int_{[0,1]}{\rm tr}_j(\psi_j(s(t))) d\mu(t)+\sum_{j=1}^{l} \af_{i,j}{\rm {{{\bar{tr}}}}}_j(q_j(b))\right).
$$
{{Then, if $h=(s, b)\in C_+$ (recall $\bt \le N_2$),
\beq
t_i(h) &\ge &{\sigma_0 R(i)'\over{\bt R(i)'}}(\sum_{j=1}^{k}
\int_{[0,1]}{\rm tr}_j(\psi_j(s(t))) d\mu(t))\\
 &= & {\sigma_0\over{\bt }} \int_{[0,1]}T(\lambda(h)(t)) d\mu \ge \sigma \int_{[0,1]}T(\lambda(h)(t)) d\mu.
\eneq}}
For $\tau\in T(A)$ mentioned above,
define
$
t_{\tau}=\sum_{i=1}^{l'} \lambda_{i, \tau}t_i.
$
It is straightforward to verify
$$
t_{\tau}(h)\ge \sigma\int_{[0,1]} {\rm T}(\lambda(h){\Green{(t)}})d\mu(t)\rforal h\in C_+.
$$
Moreover, for each $i,$ by (\ref{Nextension-n1}), for {{every}} projection $p{{=(p,\pi_e(p))}}\in M_N(C)$ (for any $N\ge 1$),
\beq\nonumber
t_i(p)&=& {{{1\over{\bt R(i)'}}\left(\sum_{j=1}^{k} \sigma_0R(i)'\int_{[0,1]}{\rm tr}_j(\psi_i(p(t))) d\mu(t)+\sum_{j=1}^{l} \af_{i,j}{\rm {\bar{tr}}}_j(q_j(\pi_e(p)))
\right)}}\\\nonumber
&=& {{{1\over{\bt R(i)'}}}}\left(\sum_{j=1}^{k} \sigma_0R(i)'{{\rm rank}(\psi_j(\phi_0\circ \pi_e(p)))\over{{\rm rank}\psi_j(\phi_0\circ \pi_e(1_C))}}+\sum_{j=1}^{l} \af_{i,j}{{\rm rank}q_j(\pi_e(p))\over{{\rm rank}q_j(\pi_e(1_C))}}\right)\\\nonumber
&=& {{{1\over{\bt R(i)'}}\left(\sigma_0 R(i)'\cdot T\circ\phi_0\circ (\pi_e)_{*0}([p])+\Gamma_i([p])\right)}}\\\label{181818}
&=&{1\over{\bt}}(1/ R(i)')(q_i'\circ \pi_e')_{*0}(\kappa([p]))={1\over{\bt}}({\rm tr}_i'(q_i'\circ \pi_e')_{*0}(\kappa([p])).
\eneq
{{In particular,
\beq\label{18-18n-3}
t_i(1_C)={1\over{\bt}}({\rm tr_i'}(q_i'\circ \pi_e')_{*0}(\kappa([1_C]))={1\over{\bt}}(1/ R(i)')(q_i'\circ \pi_e')_{*0}(\kappa([1_C])).
\eneq
It follows that $t_{\tau}(1_C)=\bt/\bt=1$ (see \eqref{18-18n-1}) and  $t_\tau\in T(C).$}}
{{Finally, by  \eqref{18-18n-1}, \eqref{181818} and \eqref{18-18n-2},}}
\beq\nonumber
\hspace{-0.2in}t_{\tau}(x)&=&{{(1/\bt)\sum_{i=1}^{l'}\lambda_{i,\tau} {\rm tr}_i'( (q_i'\circ \pi_e')_{*0}(\kappa(x)))=(1/\tau(\kappa([1_C]))) (\tau'((\pi_e')_{*0}(\kappa(x))))}}\\
&=&\tau(\kappa(x))/\tau(\kappa([1_C]))\rforal x\in K_0(C).
\eneq
\end{proof}

\begin{prop}\label{pl-lifting}
Let $S\in\mathcal C$ and $N\ge 1.$  There exists an integer $K\ge 1$ satisfying the following {{condition}}:
 For any positive homomorphism  $\kappa: {K_0}(S)\to {K_0}(A)$
 which
 satisfies { {$\kappa([1_S])\leq [1_A]$ and}} $N\kappa([p])>[1_A]$ for any  ${\Green{[p]}}\in K_0(S)_+\setminus \{0\}, $ where $A\in {\cal C}$,
there exists a  homomorphism $\phi: S\to M_K(A)$ such that $\phi_{*0}=K\kappa.$ {{ If we further assume $\kappa([1_S])=[1_A]$, then $\phi$ can be chosen to be unital.}}

\end{prop}

\begin{proof}
Write
$$
S={ {A(F_1,F_2, \phi_0, \phi_1)=}}\{(f,g): (f,g)\in C([0,1], F_2)\oplus F_1: f(0)=\phi_0(g),\,\,\, f(1)=\phi_1(g){{\}}}.
$$
Denote by {{$\lambda: S\to C([0,1], F_2)$  the map given by \eqref{pull-back}.}}
{{Since $A$ has stable rank one {{(see \ref{2pg3})}}, the proposition is known for the case that $S$ is finite dimensional (see,
for example, Lemma 7.3.2 (ii) of  \cite{RLL-Ktheory}). So we may assume $S$ is not finite dimensional.
{\Green{Since $C$ has stabel rank one, by considering each summand of $S,$}}
we may reduce the general case to the case that $S$ has only one {{direct}} summand {{(or
is minimal---see \ref{DfC1})}}. In particular, we may assume that ${\rm ker}\phi_0\cap {\rm ker}\phi_1=\{0\}.$
It follows that that $\lambda$ is injective.}}

To simplify the notation, without loss of generality, replacing $S$ by $M_r(S)$ for some integer $r\ge 1$
{{(and replacing $A$ by $M_r(A)$),}} {{let us}} assume
that projections of $S$ generate $K_0(S)_{{+}}$ {{(see \ref{FG-Ratn}).
Note also, by \ref{2Lg13}, since there are only finitely many elements of $K_0(S)_+$ which {{are}} dominated by $[1_S],$
there are only finitely many elements of $K_0(S)_+$ which can be represented by projections in
$S.$ }}

Let $\sigma>0$ be {{as}} given by \ref{Nextension} (associated with {{the}}  integers $N_1=N$ and $N_2=1$).
Define { {$\Delta: S_+^{q,1}\setminus \{0\} \to (0,1)$ by}}

$$\Delta({{\hat h}})={{\sigma/2}}\int_{[0, 1]}\mathrm{T}(\lambda(h)(t))d\mu(t)$$
for all $h=(\lambda(h), \pi_e(h))\in S_{{+}},$ where
${\rm T}(c)=\sum_{i=1}^k {\rm tr}_i(\psi_i(c))$ for all $c\in  {{F_2}},$
where {{$\psi_i: F_2\to M_{r(i)}$ {{is}} the projection map onto}} the $i$-th simple direct summand of $F_2$ (so we assume
that $F_2$ has $k$ simple {{direct}} summands) {{and}}
${\rm tr}_i$
is the normalized tracial state on  {{$M_{r(i)}.$  Note that $\sigma$ depends only on $S$ and $N_1.$
So $\Delta$ depends only on $N$ and $S.$}}

Let $\mathcal H_1$, $\delta>0$, ${\cal P},$ and $K$  be the finite subsets and constants provided by {{Theorem}} \ref{ExtTraceC-D} with respect to $S$, $\Delta$, an arbitrarily chosen finite set $\mathcal H$ {\Green{(containing $1_S$),}} and an arbitrarily chosen $1>\sigma_1={\blue{1/2}}>0$
(in place of $\sigma$).
{{Note that, when the finite subset ${\cal P}$ is given, one can replace it by
a finite subset  ${\cal P}'$  {{such}}  that the subgroup generated by ${\cal P}'$ contains ${\cal P},$ as long as
$\dt$ is chosen to be sufficiently small. Therefore, since  we have assumed that projections in $S$ generate $K_0(S)_+,$
choosing a sufficiently small $\dt,$}}
we may assume that  ${\cal P}$ {{is represented by $P\subset S,$ where $P$}}
is a finite subset of projections such that
every projection $q\in S$ is equivalent to one of  {{the}} projections in ${\cal P}.$
{{(We will only apply a part of {{Theorem}} \ref{ExtTraceC-D} and will not use \eqref{18-1610-nn}).}}

{{Note that, by
{{hypothesis}}, $\kappa([p])$  {{is the class of}}
a full projection {{in}}  $M_N(A).$}}
{{W}}ithout loss of generality, {{applying \ref{cut-full-pj},}}
we may assume
that $[1_A]=\kappa([1_S]).$

{{Let ${\cal Q}\subset M_{r'}(A)$ for some $r'\ge 1$ be a finite set of projections such
that $\kappa({\cal P})$ can be represented by projections in ${\cal Q}.$}}
It follows from  {{Lemma}} \ref{TtoDelta} that there is a finite subset ${\cal T}$ of extreme points of $T({A})$
and there exists a continuous affine map $\gamma': T({A})\to C_{\cal T}$ such that
\beq\label{june7-1}
|\gamma'(\tau)(p)-\tau(p)|<\dt/2\rforal p\in  {{\cal Q} }\,\, {\Green{\text{and\,\,for\,\, all}\,\,\tau\in T(A),}}
\eneq
where $C_{\cal T}$ is the convex hull of ${\cal T}.$

Note that  any \CA\
in the class $\mathcal C$ is of type I, it is amenable and in particular it is exact. Therefore, by {{Corollary 3.4 of}} \cite{Blatrace},
for each $s\in C_{\cal T},$
 there is a tracial state $t_s\in T(S)$ such that
\beq\label{june7-2}
r_S(t_s)(x)=r_A(s)(\kappa(x))\rforal x\in K_0(S),
\eneq
where $r_S: T(S)\to S_{[1_S]}(K_0(S))$ and
$r_A: T(A)\to S_{[1_A]}(K_0(A))$ are the induced maps from
the tracial state spaces to the state spaces of the $K_0$-groups.
It follows from {{Lemma}} \ref{Nextension} that we may choose
$t_s$ such that
\beq\label{june7-3}
t_s(h){{>}} \Delta({{\hat{h}}})\tforal h\in S^{{{\bf 1}}}_+.
\eneq
For each ${{s}}\in {\cal T},$ define
$\lambda(s)=t_s$
{{(where $t_s$}} satisfies (\ref{june7-2}) and (\ref{june7-3})). This extends to a continuous affine map
$\lambda: C_{\cal T}\to T({{S}}).$
Put $\gamma=\lambda\circ \gamma'.$ Then, for any $\tau\in {\Green{T(A)}},$
\beq\label{june7-4}
\gamma(\tau)(h) {{>}} \Delta({{\hat{h}}}) \rforal h\in {\cal H}_1
\andeqn {\Green{(\text{by}\,\, \eqref{june7-2})}}
\eneq
\beq\label{june7-5}
\hspace{-0.5in}|\gamma(\tau)(q)-\tau(\kappa([q]))|&{{\le}}&|\lambda(\gamma'(\tau))(q)-\gamma'(\tau)(\kappa([q]))|
+|\gamma'(\tau)(\kappa([q]))-\tau(\kappa([q]))|\\
&=&|\gamma'{\Green{(\tau)}}(\kappa([q]))-{\Green{\tau(\kappa([q]))}}|<\dt/2
\eneq
for all projections $q\in S.$
One then applies {{Theorem}} \ref{ExtTraceC-D} to obtain a unital \hm\, $\phi: S\to M_K(A)$
such that $[\phi]=K\kappa.$
\end{proof}

\begin{lem}\label{liftingpl-M}
 Let $C\in {\cal C}.$ Then there is $M>0$ satisfying the following {{condition}}: Let $A_1\in\mathcal B_1$  and let
$A=A_1\otimes U$ for some UHF-algebra $U$ of {{infinite type}}
and let $\kappa:({K_0}(C), {K_0}^+(C))
\to({K_0}(A), {K_0}^+(A))$
 be a strictly positive homomorphism with multiplicity $M$ {\Green{(see the lines above \eqref{Cuthm}).}} Then there exists a  homomorphism $\phi: C\to
 {{M_m(A)}}$ {{(for some integer $m\ge 1$)}} such that $\phi_{*0}=\kappa$ and $\phi_{*1}=0.$

\end{lem}
\begin{proof}
Write $C={{A}}(F_1, F_2, \phi_0, \phi_{{1}}).$
Denote by $M$ the constant of {{Corollary}} \ref{MextC} for $G=K_0(C)\subset K_0(F_1)=\Z^l.$
Let $\kappa: K_0(C)\to K_0(A)$ be a
positive homomorphism satisfying the condition of the lemma.
{\Green{Let $e\in M_r(A)$ (for some integer $r\ge 1$) be a projection such that $\kappa([1_C])=e.$
Replacing $A$ by $eM_r(A)e,$ \wilog, we may assume that $\kappa$ is unital.}}
Since $K_0(C)_+$ is finitely generated,  {{$K_0(A)_+$}} is  simple and $\kappa$ is strictly positive, there is $N$ such that for any non-zero  element $x\in K_0(C)_+$, one has that $N\kappa(x)> 2[1_A]$. Let $K$ be the natural number of Proposition \ref{pl-lifting} with respect to $C$ and ${\Green{2N}}$.

We may also assume that $M_r(C)$ contains {{a set of}} minimal projections such that
every minimal element of $K_0(C)_+\setminus \{0\}$ is represented by {{a}} minimal projection from the set {{(see \ref{FG-Ratn})}}.

By Lemma \ref{decomposition2}, for any positive map $\kappa$ with multiplicity $M$, one has
$\kappa=\kappa_1+\kappa_2$ and there are  positive \hm s $\lambda_1: K_0(C)\to \Z^n,$
$\gamma_1: \Z^n \to K_0(A),$ {{and}}  $\lambda_2: K_0(C)\to K_0(C')$  such that
$\lambda_1$ has multiplicity $M,$ $\lambda_2$ has multiplicity $MK,$
$\kappa_1=\gamma_1\circ\lambda_1,$ $\kappa_2=\imath_{*0}\circ \lambda_2,$ and
$C'\subset A$ is a \SCA\, with $C'\in {\cal C},$
where
$\imath: C'\to A$ is the embedding. Moreover,
\begin{equation}\label{120514-ext-1}
\lambda_2([1_C])=[1_{C'}]\quad\textrm{and}\quad {{2}}N \lambda_2(x)>\lambda_2([1_C])>0
\end{equation}
for any $x\in K_0(C)_+\setminus \{0\}$ {{(here  $\lambda_1$ plays the role of $\phi_1$ and
$\lambda_2$ plays the role of $\phi_2$ in {{Lemma}} \ref{decomposition2}).}}
Let $\lambda_1([1_C])=(r_1, r_2,...,r_n),$ where $r_i\in \Z_+,$ $i=1,2,...,n.$
{{Note that we may assume that $r_i\not=0$ for all $i;$ otherwise,  we replace $\Z^n$ by $\Z^{n_1}$
for some $1\le n_1<n.$}}


{{Let $\imath_A: A\to A\otimes U$ be defined by $\imath_A(a)=a\otimes 1_U$ and
$j: A\otimes U\to A$ such that $j\circ \imath_A$ is approximately unitarily equivalent to ${\rm id}_A.$ In particular,
$[j\circ \imath_A]=[{\rm id}_A].$}}

{{
Let us  first  assume
that $\lambda_1\not=0.$ }}
Let $R_0$ be as in {\Green{Corollary}} \ref{MextC} associated with $K_0(C)=G\subset \Z^l$ and $\lambda_1: K_0(C)\to \Z^n$ {{, which has multiplicity $M$.}}
Put $F_3=M_{r_1}\oplus M_{r_2}\oplus \cdots \oplus M_{r_n}$ {{(recall $\lambda_1([1_C])=(r_1,r_2,...,r_n)$).}}
Since $A$ has stable rank one,  there is a \hm\, $\psi_0: F_3\to A$ such that
$(\psi_0)_{*0}=\gamma_1.$ Write $U=\lim_{n\to\infty}(M_{R(n)}, h_n),$ where $h_n: M_{R(n)}\to M_{R(n+1)}$ is a unital embedding. Choose $R(n)\ge R_0.$ Consider the unital \hm\, $j_{F_3}: F_3\to F_3\otimes M_{R(n)}$  defined
by $j_{F_3}(a)=a\otimes 1_{M_{R(n)}}$ for all $a\in F_3,$ and consider
the unital \hm\, $\psi_0\otimes h_{n, \infty}: F_{{3}}\otimes M_{R(n)}\to A\otimes U$ defined by
${(}\psi_0\otimes h_{n, \infty}{)}(a\otimes b)=\psi_0(a)\otimes h_{n, \infty}(b)$ for all $a\in F_3$ and $b\in M_{R(n)}.$
We {{have}}, for any projection $p\in F_3$ {\Green{(recall $A=A_1\otimes U$)}}
\beq\label{141204lift-1}
((\psi_0\otimes h_{n,\infty})\circ j_{F_3})_{*0}([p])=[\psi_0(p)\otimes 1_U]={{(\imath_A)_{*0}\circ}} (\psi_0)_{*0}([p])\in K_0(A).
\eneq
It follows that
\beq\label{141204-lift-2}
({{j_{*0}\circ}} (\psi_0\otimes h_{n, \infty})_{*0}\circ (j_{F_3})_{*0}{{={(j\circ \imath_A)_{*0}\circ (\psi_0)_{*0}}}}=(\psi_0)_{*0}.
\eneq
Now the map $(j_{F_3})_{*0}\circ {{\lambda}}_1: K_0(C)\to K_0({{F_3}})=\Z^n$ has multiplicity $MR(n).$
Applying  {\Green{Corollary}} \ref{MextC}, we obtain a positive \hm\, $\lambda_1': K_0(F_1)=\Z^l\to K_0(F_3)$
such that $(\lambda_1')_{*0}\circ (\pi_{{e}})_{*0}=(j_{F_3})_{*0}\circ \lambda_1.$
The {{construction above}} can be summarized by the following commutative diagram:
$$
\begin{array}{ccccc}
\Z^l & \stackrel{\lambda_1'}{\rightarrow}  & K_0(F_3\otimes M_{R(n)}{)} & \stackrel{(\psi_0\otimes h_{n, \infty})_{*0}}{\longrightarrow}& K_0(A\otimes U)\\
\hspace{0.3in}\uparrow_{(\pi_e)_{*0}} && \uparrow_{(j_{F_3})_{*0}} &&   {{{}_{\imath}\uparrow\downarrow^{j}}}\\
K_0(C) & \stackrel{\lambda_1}\longrightarrow  & K_0(F_3) & \stackrel{\gamma_1=(\psi_0)_{*0}}{\longrightarrow} & K_0(A).
\end{array}
$$
 We obtain a \hm\, $h_0: F_1\to F_3\otimes M_{R(n)}$ such that
 $(h_0)_{*0}=\lambda_1'.$
 Define ${{H}}_1=h_0\circ \pi_e: C\to F_3\otimes M_{R(n)}$ and
 ${{H}}_2=(j\circ \psi_0\otimes {{h}}_{n, \infty})\circ {\Green{H_1}}: C\to A\otimes U.
 $
Then, by the 
commutative diagram above,
\beq\label{120412-ext4}
({{H}}_2)_{*0}=\kappa_1.
\eneq
{{If $\lambda_1=0,$ then $\kappa_1=0,$ and we choose $H_2=0.$}}

Since $\lambda_2$ has multiplicity $K,$ there exists
$\lambda_2': K_0(C)\to K_0(C')$ such that
$K\lambda_2'=\lambda_2.$ Since $K_0(C')$ is weakly unperforated,
$\lambda_2'$ is positive. Moreover, by  (\ref{120514-ext-1}),
\beq\label{120414-ext4}
{{2}}KN\lambda_2'(x)>K\lambda_2'([1_C])=\lambda_2([1_C])=[1_{C'}]>0.
\eneq
Since $K_0(C')$ is weakly unperforated, we have
\beq\label{120414-ext5}
{\Green{2N}}\lambda_2'(x)>\lambda_2'([1_C])>0\rforal x\in K_0(C)_+\setminus \{0\}.
\eneq
There is a projection $e\in M_k(C')$ for some integer $k\ge 1$ such that
$\lambda_2'([1_C])=[e].$  Define $C''=eM_k(C')e.$  By (\ref{120414-ext4}), $e$ is full in $C'.$
In fact, $K[e]=[1_{C{{'}}}].$ In other words, $M_K(C'')\cong C'.$
By \ref{cut-full-pj}, $C''\in {\cal C}.$  Applying  Proposition \ref{pl-lifting}, we obtain a  unital \hm\,
${{H_3'}}: C\to C{{'}}$ such that $({{H_3'}})_{*0}=K\lambda_2'=\lambda_2.$
Put ${{H_3=\imath\circ H_3'}}.$
Note that $[1_A]=\kappa([1_C])=\kappa_1([1_C])+\kappa_2([1_C]).$  {{C}}onjugating by a unitary,
we may assume
that ${{H_2}}(1_C)+{{H_3}}(1_C)=1_A.$
Then it is easy to check that the map $\phi: C\to A\otimes U$ defined by $\phi(c)={{H_2}}(c)+{{H_3}}(c)$ for all $c\in C$ meets the requirements.
\end{proof}

\begin{lem}[cf.~Lemma 9.8 of \cite{LinTAI}]\label{smallkkk}
Let $A$ be a unital \CA\, and let $B_1$ be a unital separable simple \CA\, in ${\cal B}_0,$ and let $B=B_1\otimes U$
for some  UHF-algebra of infinite type and $C\in {\cal C}_0$ be a \SCA\, of $B.$
Let $G\subset \underline{K}(A)$ be a finitely generated subgroup.
Suppose that there exists an ${\cal F}$-$\dt$-multiplicative  \morp\, $\psi: A\to C\subset B$ such that {\Green{$\psi(1_A)=p$ is projection  and}}
$[\psi]|_G$ is well defined. Then, for any $\ep>0,$  there {{exist a}} \SCA\, $C_1\cong C$ of $B$ and an ${\cal F}$-$\dt$-multiplicative \morp\, $L: A\to C_1\subset B$ such that
\beq\label{smallkkk-1}
[L]|_{G\cap K_0(A, \Z/k\Z)}=[\psi]|_{G\cap K_0(A, \Z/k\Z)}\tand \tau(1_{C_1})<\ep
\eneq
for all $\tau\in T(B)$ and for all $k\ge 1$ {{and such}} that $G\cap K_0(A, \Z/k\Z)\not=\{0\},$
where $L$ and $\psi$ are viewed as maps to $B.$ Furthermore,
if $[\psi]|_{G\cap K_0(A)}$ is positive {{then}} so also is $[L]|_{G\cap K_0(A)}.$
\end{lem}

\begin{proof}
{{Let $\imath: U\to U\otimes U$ be defined by
$\imath(a)=a\otimes 1_U$ and {{let}} $j: U\otimes U\to U$ be  an isomorphism such that
$j\circ \imath$ is approximately unitarily equivalent to ${\rm id}_U.$
Put $j_B={\rm id}_{B_1}\otimes j: B\otimes U\,(=B_1\otimes U\otimes U)\to B\,(=B_1\otimes U).$}}
Let $1>\ep>0.$  Suppose that
$$
G\cap K_0(A,\Z/k\Z)=\{0\}\rforal k\ge K
$$
for some integer $K\ge 1.$
Find a projection $e_0\in U$ such that $\tau_0(e_0)<\ep$
and $1_U=e_0+\sum_{i=1}^m p_i,$ where $m=2l K!$ and $1/l<\ep$ and $p_1, p_2,...,p_m$ are mutually orthogonal
and mutually equivalent projections in $U.$
Choose $C_1'=C\otimes e_0\subset B\otimes U$ {{and $C_1=j_B(C_1').$}} Then $C_1\cong C.$
Let $\phi: C\to C_1$ be the isomorphism defined by $\phi(c)={{j_B(c\otimes e_0)}}$ for all $c\in C.$
Put $L=\phi\circ \psi.$ Note that $K_1(C)=K_1(C_1)=\{0\}.$  Both $[L]$ and $[\psi]$ map
${\Green{{\cal G}\cap}}K_0(A,\Z/k\Z)$ to $K_0(B)/kK_0(B){\Green{\subset K_0(B, \Z/k\Z)}}$ and factor through $K_0(C, \Z/k\Z).$ It follows that {\Green{(recall every element in $K_0(\cdot, \Z/k\Z)$ is $k$-torsion)}}
$$
[L]|_{G\cap K_0(A. \Z/k\Z)}=[\psi]|_{G\cap K_0(A, \Z/k\Z)}.
$$
In the case that $[\psi]|_{G\cap K_0(A)}$ is positive,  it follows from the definition of $L$
that $[L]|_{G\cap K_0(A)}$ is also positive.
\end{proof}

\begin{thm}\label{preBot2}
Let $C$ and $A$ be unital stably finite  \CA s and let  $\af\in {{KL}}_{e}(C, A)^{++}.$

{\rm  (i)}\,\, If  $C\in{\bf H}$, or $C\in \mathcal C,$ and $A_1$
is a unital simple \CA\, in ${\cal B}_{0}$ and
$A=A_1\otimes U$ for some UHF-algebra $U$ of infinite type,
 then there exists a sequence of completely positive linear maps  $L_n: C\to A$ such that
\beq\label{preBot2-1}
\lim_{n\to\infty} \|L_n(ab)-L_n(a)L_n(b)\|=0 \tforal a,\, b\in C\tand\,
{[}{\{}L_n{\}}{]}=\af;
\eneq

{\rm (ii)}\,\, if $C\in {\cal C}_0$ and $A_1\in {\cal B}_1,$ the conclusion  above also holds;

{\rm (iii)} if $C=M_n(C(S^2))$ for some integer $n\ge 1,$  $A={{A_1}}\otimes U$ and  $A_1\in {\cal B}_1,$
then there is a unital \hm\, $h: C\to A$ such that $[h]=\af;$

{\rm (iv)} if $C\in {\bf H}$ with $K_1(C)$ torsion, $M_n(C(S^2))$ {{is not a direct summand of $C,$}} and $A=A_1\otimes Q,$ where $A_1$ is unital and $A$ has stable rank one,
then there exists a unital \hm\, $h: C\to A$ such that $[h]=\af$;

{\rm (v)} if $C=M_n(C(\T))$ for some integer $n\ge 1,$ then
for any unital \CA\, $A$ with stable rank one, there is a unital \hm\, $h: C\to A$ such that
$[h]=\af.$
\end{thm}

\begin{proof}
Let us first consider (iii). This is a special case of Lemma 2.19 of \cite{Niu-TAS-II}. Let us provide a proof here.
In this case one has that $K_0(C)=\Z\oplus {\rm ker}\rho_C\cong \Z\oplus \Z$ is free {{abelian}} and $K_1(C)=\{0\}.$
Write $A=A_1\otimes U,$ where $K_0(U)=D\subset \Q$ is identified with a dense subgroup of $\Q$ and $1_U=1.$
Let $\af_0=\af|_{K_0(C)}.$ Then $\af_0([1_C])=[1_A]$ and $\af_0(x)\in {\rm ker}\rho_A$ for all $x\in {\rm ker}\rho_C$
{{(since $\af_0$ is order preserving).}}
Let $\xi\in {\rm ker}\rho_C=\Z$ be a generator and $\af_0(\xi)=\zeta\in {\rm ker}\rho_A.$
Let ${{B}}$ be a the unital simple AF-algebra with
$$
(K_0(B), K_0(B)_+, [1_{B_0}])=(D\oplus \Z, (D\oplus \Z)_+, (1,0)),
$$
where
$$
(D\oplus \Z)_+=\{(d,m): d>0,m\in \Z\}\cup \{(0,0)\}.
$$
{{Note that $K_0(\T^2)=\Z\oplus \Z$ which is order isomorphic to $K_0(C(S^2)).$
There is a standard continuous map which maps $s: \T^2\to S^2$ (with  $\T^2$ viewed as a 2-cell
attached to a figure 8, the map is defined by sending the figure 8 to a single point) such that
$\gamma: C(S^2)\to C(\T^2)$ defined by $\gamma(f)=f\circ s$ has the property
that $\gamma_{*0}={\rm id}_{\Z\oplus \Z}.$}}
It follows from \cite{EL-Emb}  that there is a unital \hm\, ${{h_0'}}: C(\T^2)\to B$ such that
$({{h'_0}})_{*0}(\xi)=(0, 1).$ {{Define $h_0=h_0'\circ \gamma.$}}
There is a positive order-unit preserving \hm\, $\lambda: D\to K_0(A)$ (given by the embedding
$a\to 1_A\otimes a$ from $U\to A_1\otimes U$).
Define a \hm\, $\kappa_0: K_0(B)\to K_0(A)$
by $\kappa_0((r,0))=\lambda(r)$ for all $r\in D$ and $\kappa_0((0,1))=\zeta.$
Since $A$ has stable rank one  {{and $B$ is AF}},   it is known {{that there is a}}
unital \hm\, $\phi: B\to A$ such that
\beq\label{12/24/14-1}
\phi_{*0}=\kappa_0.
\eneq
Define $L=\phi\circ h_0.$ Then, $[L]=\af.$  This proves (iii).

For ({{i}}v), we note that $K_i(A)$ is torsion free and divisible. {{One may reduce the general case to the case
that $C$ has only one direct summand}} $C=PM_n(C(X))P,$ where $X$ is a connected finite CW complex  and $P\in M_n(C(X))$ is a projection. Note that
we assume here $K_1(C)$ is torsion and $X\not=S^2.$ In this case {{(see \ref{AHblock}   and \ref{range 0.18}),}} $K_0(C)=\Z\oplus {\rm Tor}(K_0(C)),$ $K_1(C)=\{0\},$ or $K_0(C)=\Z$ and
$K_1(C)$ is finite.  Suppose that $P$ has rank $r\ge 1.$
Choose $x_0\in X.$ Let $\pi_{x_0}: C\to M_r$ be defined by $\pi_{x_0}(f)=f(x_0)$ for all $f\in C.$
Suppose that $e=(1,0)\in \Z\oplus {\rm Tor}(K_0(C))$ or $e=1\in \Z.$
Choose a projection $p\in A$ such that $[p]=\af_0(e)$ (this is possible since $A$ has stable rank one).
{{Note that $[P]=(r, x),$ where $x\in {\rm Tor}(K_0(C)),$ i.e.,  $[P]-re$ is a torsion element.}}
{{Since $K_i(A)$ is torsion free, $\af_0({\rm Tor}(K_0(C)))=0.$ Thus $r\af_0(e)=r[p]=\af_0([P])=[1_A].$
Moreover, since $K_0(A)$ is torsion free, $[1_A]$ is a sum of $r$ mutually equivalent projections which
are all  unitarily equivalent to $p.$}}
Thus there is a unital \hm\, $h_0: M_r\to A$ such that $h_0(e_{1,1})=p,$ where $e_{1,1}\in M_r$ is a rank one projection.
Define $h: C\to A$ by $h=h_0\circ \pi_{x_0}.$ One verifies that $[h]=\af.$

Now we prove (i) and (ii). {{First, if $C\in \mathbf H,$ as before,  we may assume that $C$ has a single
direct summand.
In this case, if $C=M_n(C(S^2)),$ this follows from case (iii).}}
If $C\in\mathbf H$ and $C\neq M_n(C(S^2))$, the statement follows from the same argument as that of Lemma 9.9 of \cite{LinTAI},
by replacing Lemma 9.8 of \cite{LinTAI} by \ref{smallkkk} above (and {{replacing}} $F$ by $C$ and $F_1$ by $C_1$) in
the proof of Lemma 9.9 of \cite{LinTAI}.

Assume that $C\in \mathcal C$. {{Considering}} ${\rm Ad}\, w\circ L_n|_C$ for suitable
unitary $w$ (in $M_r(A)$),  we may replace $C$ by $M_r(C)$ for some $r\ge 1$
so that $K_0(C)_+$ is generated by minimal projections $\{p_1,p_2,...,p_d\}\subset C$ (see \ref{FG-Ratn}).
{\Green{We may rearrange it so that $\{[p_1], [p_2],...,[p_{d'}]\}$  ($0<d'\le d$) forms
a base for $K_0(C).$}}
Since $A$ is simple and $\af({\Green{K_0(C)}}_+\setminus \{0\})\subset K_0(A)_+\setminus\{0\},$ there exists an integer
$N\ge 1$ such that
\begin{equation}\label{prebot2141204-n1}
N\af([p])>2[1_A]\rforal [p]\in K_0(C)_+\setminus\{0\}.
\end{equation}

Let $M\ge 1$
be the integer given by {{Lemma \ref{liftingpl-M}}}
associated
with
$C.$
Since $C$ has a separating family of finite dimensional representations, by Theorem \ref{kkmaps}, there exist two sequences of completely positive contractions $\phi_n^{(i)}: C\rightarrow{A}\otimes{\cal K}\ (i=0, 1)$ satisfying the following {{conditions}}:

\noindent
{\rm (a)} $\|{\phi_n^{(i)}(ab)-\phi_n^{(i)}(a)\phi_n^{(i)}(b)}\|\rightarrow0$, $\rforal a, b\in C$,
      as $n\rightarrow\infty$,

 \noindent
{\rm (b)} for any $n$, $\phi_n^{(1)}$  is a \hm\, with  finite
      dimensional  range
      and, consequently, for any finite subset
      ${\cal P}\subset \underline{K}(C)$, the map $[\phi_n^{(i)}]|_{\cal P}$ {{is}}  well defined for
     all sufficiently
     large $n$, {{and}}

\noindent
{\rm (c)} for each finite subset ${\cal P}\subset \underline{K}(C)$, there exists $m>0$
      such that $$[\phi_n^{(0)}]|_{\cal P}=\alpha+[\phi_n^{(1)}]|_{\cal P} \tforal n>m.$$


Since $C$ is semiprojective and the positive cone of the $K_0$-group is finitely generated {{(see the end of \ref{DfC1}),}}
there are homomorphisms $\phi_0$ and $\phi_1$ from $C \to A\otimes\cal K$ such that $$[\phi_0]=\alpha+[\phi_1].$$ Without los{{s}} of generality, {{we may}} assume that $\phi_0$ and $\phi_1$ are homomorphisms from $C$ to $M_r(A)$ for some $r$. Note that $M_r(A)\in \mathcal B_0$ (or $\in {\cal B}_1,$ when $C\in {\cal C}_0$).

Since $K_i(C)$ is finitely generated ($i=0,1$), there exists $n_0\ge 1$ such that
every element $\kappa\in KL(C,A)$ is determined by
$\kappa$ on $K_i(C)$ and $K_i(C, \Z/n\Z)$ for $2\le n\le n_0,$ $i={{0,1}}$  (see Corollary 2.11 of \cite{DL}).
Let ${\cal P}\subset \underline{K}(C)$ be a finite subset which generates {{the group}}
$$
\bigoplus_{i=0,1}(K_i(C)\oplus\bigoplus_{2\le n\le n_0}K_i(C, \Z/n\Z)).
$$
{{Let}} $K=n_0!.$
Let $\cal G$ be a finite subset of $M_r(A)$ which contains $\{\phi_0(p_i), \phi_1(p_i); i=1, ..., d\}$.
We may assume that ${{{\cal P}_0}}:=\{[p_i],i=1,2,{...,} d\} \subset {\cal P}.$
{\Green{Let $G({\cal P})$ be the subgroup generated by ${\cal P}.$ Note that $K_0(C)\subset G({\cal P}).$}}

Let
 $$
 T=\max\{\tau(\phi_0(p_i))+KM\tau(\phi_1(p_i)): 1\le i\le d; \tau\in T(A)\}.
 $$
 Choose $r_0>0$ such that
 \beq\label{prebot2-nn3}
 NTr_0<1/2.
 \eneq
 Let ${\cal Q}=[\phi_0]({\cal P})\cup [\phi_1]({\cal P})\cup\af({\cal P}).$
 Let $1>\ep>0.$
 By {\Green{Proposition}} \ref{lem-compress}, for $\ep$ and $r_0$ {{as}} above,
 there {{are}}  a non-zero projection $e\in M_r(A),$ a \SCA\, $B\in {\cal C}_0$ {{(or $B\in {\cal C}$ for case (ii))}} with
 $e=1_B,$  {{and}} $\mathcal{G}$-$\ep$-multiplicative \morp s $L_1 :M_r(A)\to(1-e)M_r(A)(1-e)$
 and $L_2: M_r(A)\to B$ with the following properties:
\begin{enumerate}
\item $\|L_1(a)+L_2(a)-a\|<\ep/2\tforal a\in {\cal G};$
\item $[L_i]|_{\cal Q}$  is well defined, $i=1,2;$
\item $[L_1]|_{\cal Q}+[\imath\circ L_2]|_{\cal Q}=[{\rm id}]|_{\cal Q};$
\item\label{08-10-lem-cond-02} $\tau\circ[L_1](g)\leq r_0\tau(g)$ for all $g\in {{\af}}({\mathcal P}_0)$ and $\tau\in T(A)$;
\item\label{08-10-lem-cond-03}  {{for}} any $x\in {\cal Q},$ there exists $y\in \underline{K}(B)$ such that
   $x-[L_1](x)=[\imath\circ L_2](x)=KM[\imath](y);$
   and
\item\label{cond3} {{there}} exist positive elements $\{f_i\}\subset{K_0}(B)_+$
      such that for  $i=1,...,{{d}},$
      $$\af([p_i])-[L_1](\af([p_i])=[\imath\circ L_2](\alpha([p_i]))=KM\imath_{*0}(f_i),$$
      \end{enumerate}
      where $\imath: B\to {\Green{M_r(A)}}$ is the embedding.
{\Green{Here we also write $[L_1]$ as a \hm\, on the subgroup generated ${\cal Q}.$}}
By (\ref{08-10-lem-cond-03}), since $K=n_0!,$
\beq\label{prebot2-nn4}
[\imath\circ L_2\circ \phi_0]|_{K_i(C, \Z/n\Z)\cap {\cal P}}=[\imath\circ L_2\circ \phi_1]|_{K_i(C, \Z/n\Z){{\cap}} {\cal P}}=0,\,\,\,i=0,1,
\eneq
and $n=1,2,...,n_0.$  It follows that
\beq\label{prebot2-n1}
&&{[}L_1\circ\phi_0{]}|_{K_i(C, \Z/n\Z)\cap {\cal P}}=[\phi_0]|_{K_i(C, \Z/n\Z)\cap {\cal P}}\quad\textrm{and}\quad\\
&& [L_1\circ\phi_1]|_{K_i(C, \Z/n\Z)\cap {\cal P}}=[\phi_1]|_{K_i(C, \Z/n\Z)\cap {\cal P}},
\eneq
$i=0,1$ and $n=1, 2, ...,n_0.$
Furthermore, { {for the case $B\in {\cal C}_0$ we have}} $K_1(B)=0,$ { {and consequently}}
\beq\label{prebot2-nn6}
[\imath\circ L_2]|_{K_1(C)\cap {\cal P}}=0.
\eneq

It follows that
\beq\label{prebot2-nn5}
[L_1\circ \phi_0]|_{K_1(C)\cap {\cal P}}=[\phi_0]|_{K_1(C)\cap{\cal P}},\,\,\, [L_1\circ \phi_1]|_{K_1(C)\cap {\cal P}}=[\phi_1]|_{K_1(C)\cap {\cal P}}.
\eneq
In the second case when  we assume that $C\in {\cal C}_0$ and $A_1\in {\cal B}_1,$ then $K_1(C)=0.$  Therefore
(\ref{prebot2-nn5}) above also holds.

Denote by $\Psi:=\phi_0\oplus\bigoplus_{KM-1}\phi_1$. One then has
\begin{eqnarray*}
[L_1\circ\Psi]_{K_i(C,\ \mathbb Z/n\mathbb Z)\cap {\cal P}} & = & [L_1\circ\phi_0]_{K_i(C,\ \mathbb Z/n\mathbb Z)\cap {\cal P}}+ (KM-1)[L_1\circ\phi_1]_{K_i(C,\ \mathbb Z/n\mathbb Z)\cap {\cal P}} \\
& = & [L_1\circ\phi_0]_{K_i(C,\ \mathbb Z/n\mathbb Z)\cap {\cal P}}-[L_1\circ\phi_1]_{K_i(C,\ \mathbb Z/n\mathbb Z)\cap {\cal P}}\\
&=&[\phi_0]_{K_i(C,\ \mathbb Z/n\mathbb Z)\cap {\cal P}}-[\phi_1]_{K_i(C,\ \mathbb Z/n\mathbb Z)\cap {\cal P}}\\
&=&\alpha|_{K_i(C,\ \mathbb Z/n\mathbb Z)\cap {\cal P}},
\end{eqnarray*}
where $i=0, 1$, $n=1, 2, ...,{{n_0}}.$  {\Green{Note that
we also assume $[L_1\circ \Psi]|_{{\cal G}({\cal P})}$ is a \hm.}}
%
By  (\ref{prebot2141204-n1}), (\ref{08-10-lem-cond-02}), and (\ref{prebot2-nn3}),
\beq\label{prebot1204-n2}
N(\tau(\alpha([p_i]))-[L_1\circ\Psi]([p_i]))\ge 2-Nr_0T\ge {{3/2}}  \rforal \tau\in T(A).
\eneq
Since the strict order on ${K_0}(A)$ is determined by traces {\Green{(see \ref{Comparison} and \ref{Ltrace})}}, one has that ${{N}}(\alpha([p_i])-[L_1\circ\Psi]([p_i]))>[1_A].$
Moreover, one also has
$$\begin{array}{ll}
  & \alpha([p_i])-[L_1\circ\Psi]([p_i]) \\
  =&\alpha([p_i])-([L_1\circ\alpha]([p_i])+KM[L_1\circ\phi_1]([p_i])) \\
  =&(\alpha([p_i])-[L_1\circ\alpha]([p_i]))-KM[L_1\circ\phi_1]([p_i]) \\
  =& KM(\imath_{*0}(f_j)-[L_1]\circ[\phi_1]([p_i])) \\
  =& KMf'_{{i}}, \hspace{0.8in}\quad \mbox{where}\ f'_{{i}}={\Green{\imath_{*0}(f_{{i}})}}-[L_1]\circ[\phi_1]([p_i]).
\end{array}$$
Note that $f_j'\in K_0(A)_+\setminus \{0\},$ $j=1,2,...,d.$
{\Green{Note also that $\af-[L_1\circ \Psi]$ defines a \hm\, on $K_0(C).$
Since $Mf_i'\in K_0(A),$ $i=1,2,...,d',$  the map $\bt: K_0(C)\to K_0(A)$
defined by $\bt=(1/K)(\af-[L_1\circ \Psi])|_{K_0(C)}$ which maps
$[p_i]$ to $Mf_i'$ ($1\le i\le d'$) is a \hm.  In fact, one has}}
{\Green{$\bt([p_j])=Mf_j'+z_j,$ $j=1,2,...,d,$
where $Kz_j=0.$   Therefore $\bt([p_j])\in K_0(A)_+\setminus\{0\}$ {\Green{(recall 
the order of $K_0(A)$ is determined by the traces),}}  $j=1,2,...,d.$}}

Since the {{semigroup}} $K_0(C)_+$  is generated by $[p_1], [p_2],...,[p_d],$
{{we have}}
$\bt(K_0(C)_+\setminus\{0\})\subset K_0(A)_+\setminus\{0\}.$
Since $\bt$ has  multiplicity $M$, {{by}} the choice of $M$ and by Lemma \ref{liftingpl-M}, there  exists a
homomorphism ${{H}}: C\to M_R(A)$ (it may not be unital)  {{for some integer $R\ge 1$}} such that
$$
{{H}}_{*0}=\beta\quad\mathrm{and}\quad {{H}}_{*1}=0.
$$

Consider the map $\phi':=L_{{1}}\circ\Psi\oplus (\bigoplus_{i=1}^K {{H}}): C \to A\otimes{\cal K}.$ {{One}} has that
$$[\phi']|_{K_0(C)\cap {\cal P}}=[L_{{1}}\circ\Psi]|_{K_0(C)\cap {\cal P}}+K\beta=\alpha|_{K_0(C)\cap {\cal P}}.$$
{\Green{Since every element in $K_i(\cdot, \Z/n\Z)$ is $n$-torsion,}}
$$K[H]|_{K_i(C,\ \mathbb Z/n\mathbb Z)\cap {\cal P}}=0,\quad i=0, 1,\ n=1, 2,...,n_0{,}$$
and therefore
$$[\phi']_{K_i(C,\ \mathbb Z/n\mathbb Z)\cap {\cal P}}=[L_{{1}}\circ\Psi]|_{K_i(C,\ \mathbb Z/n\mathbb Z)\cap {\cal P}}=\alpha|_{K_i(C,\ \mathbb Z/n\mathbb Z)\cap {\cal P}},$$
$i=0, 1$, $n=1, 2, ..., n_0${.}
We also have  that
$$[\phi']_{K_1(C)\cap {\cal P}}=[L_{{1}}\circ\Psi]|_{K_1(C)\cap {\cal P}}=\alpha|_{K_1(C)\cap {\cal P}}.$$
Therefore,
$$[\phi']|_{\cal P}=\alpha|_{\cal P}.$$

Since $[\phi'(1_C)]=[1_A]$ and $A$ has stable rank one, there is a unitary $u$ in a matrix algebra {{over}} $A$ such that the map $\phi=\mathrm{Ad}(u)\circ\phi'$ satisfies $\phi(1_C)=1_A$, as desired.

Case (v) is standard and is well known.
\end{proof}

\begin{cor}\label{est-m2a}
Any \CA\, $A$ given by  Theorem \ref{RangT} is KK-attainable with respect to ${\cal B}_{u0}$.
\end{cor}

\begin{proof}
Note that
{{$A=\lim_{n\to\infty}(A_n, \phi_n),$ where each $A_n$ is a}}  finite direct sum
of \CA s in $\mathbf H$ or  $\mathcal C_0$ {{ and  $\phi_n$ is injective.
Let $A_n'=\phi_{n, \infty}(A_n)\,(\cong A_n)$
and
$\imath_n: A_n'\to A$ be the embedding.}}
{{By Lemma 2.2 of \cite{DL},  for any \CA\, $B,$
${\rm Hom}_{\Lambda}(\underline{K}(A), \underline{K}(B))=
\lim_{n\to\infty}{\rm Hom}_{\Lambda}(\underline{K}(A_n'), \underline{K}(B)).$
Let $\af\in {\rm Hom}_{\Lambda}(\underline{K}(A),\underline{K}(B))^{++}.$ Then, for each $n,$
there is  $\af_n\in {\rm Hom}_{\Lambda}(\underline{K}(A_n'), \underline{K}(B))$
such that
$\af|_{[\imath_n](\underline{K}(A_n))}=\af_n$ and
$\af_{n+l}|_{[\phi_{n, n+l}](\underline{K}(A_n'))}=\af_n$ for all $l\ge 1.$}}
{{
Noting that $\imath_n(p)\not=0$ for any non-zero projection $p\in M_m(A_n')$ (for
any $m>0$), one has $\af_n\in {\rm Hom}_{\Lambda}(\underline{K}(A_n'), \underline{K}(B))^{++}.$}}
{{Let $q\in M_m(B)$  be a projection for some integer $m\ge 1$ such that $[q]=\af([1_A]).$
{\Green{Since $A$ is unital, we may assume that $\phi_n$ is unital (for large $n$).}}
Then $[q]=\af_n([1_{A_n'}]).$ Let $B_1=qM_m(B)q.$ Then $\af_n\in KL_e(A_n', B_1)^{++}.$}}
{\Green{Recall that $B$ has stable rank one. By considering each summand of $A_n',$
by}} Theorem \ref{preBot2} that, for each $n,$ there exists  a sequence of \morp s $L_{k,n}: A_n'\to B_1$
such that $\lim_{k\to\infty}\|L_{k,n}(ab)-L_{k,n}(a)L_{k,n}(b)\|=0$ for all $a, b\in A_k'$ and
$[\{L_{k,n}\}]=\af_n.$  Since {\Green{each $A_m'$}} is separable and amenable, there exists a sequence of \morp s
$\Psi_m: A\to A_m'$ such that
$\lim_{n\to\infty}\Psi_m(a)=a$ (see, for example, 2.3.13 of \cite{Lnbok}).
It is standard that,
{{after passing to}} two  suitable subsequences $\{m_n\}$ and $\{k_n\},$
{{the  sequence  $\{L_n=L_{k_n, m_n}\circ \Psi_{m_n}\}$  has the property that
\beq
\lim_{n\to\infty}\|L_n(ab)-L_n(a)L_n(b)\|=0\rforal a, b\in A\andeqn [\{L_n\}]=\af.
\eneq
Thus,}}
$A$ is KK-attainable with respect to ${\cal B}_{u0}.$
\end{proof}

\begin{cor}\label{C0ext}
Let $C\in {\cal C},$
$A\in {\cal B}_{u0},$ and $\af\in KK(C,A)^{++}$ be  such that
$\af([1_C])=[p]$ for some projection $p\in A.$
Then
there exists a \hm\, $\phi: C\to A$ such that,
$\phi_{*0}=\af.$
\end{cor}

\begin{proof}
This is a special case of Theorem \ref{preBot2} since $C$ is semiprojective.
\end{proof}

\begin{cor}\label{Cistbk}
Let $A\in\mathcal B_{u0}$.
Then there exist a unital simple \CA\, $B_1=\varinjlim(C_n, \phi_n),$ where
each $C_n$ is in ${\cal C}_0,$ and a UHF algebra $U$ of infinite type   such that, for $B=B_1\otimes {{U}},$  {{we have}}
\beq\label{Cistbk-n1}
(K_0(B), K_0(B)_+, [1_{B}], T(B), r_{B})
=(\rho_A(K_0(A)), \rho_A(K_0(A)_+), \rho_A([1_A]), T(A), r_A).
\eneq
Moreover, for each $n,$ there is a \hm\, $h_n: C_n\otimes U\to M_2(A)$ such that
\begin{equation}\label{Cistbk-1}
\rho_A\circ (h_n)_{*0}=(\phi_{n, \infty}\otimes\mathrm{id}_U)_{*0}.
\end{equation}
\end{cor}
\begin{proof}
{{Recall that $r_A: T(A)\to S_u(K_0(A))$ is defined
by $r_A(\tau)(x)=\rho_A(x)(\tau)$ for all $x\in K_0(A).$  Therefore $r_A$ also induces a continuous affine
map from $T(A)$ to $S_u(\rho_A(K_0(A)))$ by defining $r_A(\tau)(\rho_A(x))=\rho_A(x)(\tau)$ for all $x\in K_0(A).$
If $s\in S_u(\rho_A(K_0(A))),$ then $s\circ \rho_A\in S_u(K_0(A)).$
By Corollary 3.4 of \cite{Blatrace}, $r_A: T(A)\to S_u(K_0(A))$ is surjective.  It follows that $r_A$ also
maps $T(A)$ onto $S_u(\rho_A(K_0(A))).$ }}
Consider the tuple
$$(\rho_A(K_0(A)), \rho_A(K_0(A)_+), \rho_A([1_A]), T(A), r_A).$$
Since $A\cong A\otimes U_1,$ for a UHF algebra $U_1$ of infinite type,  {{$A$}} has the property (SP) (see \cite{BKR-ADiv}), and therefore the ordered group $(K_0(A), K_0(A)_+, [1_A])$ has the property (SP) in the sense of Theorem \ref{RangT}; that is, for any positive real number $0<s<1$, there is $g\in K_0(A)_+$ such that $\tau(g)<s$ for any $\tau\in T(A)$. Then it is clear that the scaled ordered group $(\rho_A(K_0(A)), (\rho_A(K_0(A)_+), \rho_A(1_A))$ also has the property (SP) in the sense of Theorem \ref{RangT}. Therefore, by Theorem \ref{RangT}, there is a simple unital \CA\, $B_1=\varinjlim(C_n, \phi_n)$, where each $C_n\in\mathcal C_0,$  such that
$$(K_0(B_1), K_0(B_1)_+, [1_{B_1}], T(B_1), r_{B_1})\cong(\rho_A(K_0(A)), \rho_A(K_0(A)_+), \rho_A(1_A), T(A), r_A).$$

{ {Let $U=U_1$ and}}  $B=B_1\otimes U$.
{{Then}}  $$(K_0(B), K_0(B)_+, [1_{B}], T(B), r_{B})
\cong(\rho_A(K_0(A\otimes U)), \rho_A(K_0(A\otimes U)_+), \rho_A(1_A), T(A\otimes U ), r_{A\otimes U}).$$
{{Let $\imath:  U\to  U\otimes U$ be defined by $\imath(a)=a\otimes 1_U$  for all $a\in U$ and
let $j: U\otimes U\to U$ be an isomorphism; {{recall}}
 that
$j^{-1}$ is approximately unitarily equivalent to $\imath.$
Define $\imath_A: A\to A\otimes U$ by $\imath_A={\rm id}_A\otimes \imath$ and
define $j_A={\rm id}_A\otimes j.$
It follows that $j_A$ induces an order isomorphism:
\beq
(\rho_A(K_0(A\otimes U)), \rho_A(K_0(A\otimes U)_+), \rho_A(1_{A\otimes U}), T(A\otimes U ), r_{A\otimes U})\\
\cong
(\rho_A(K_0(A)), \rho_A(K_0(A)_+), \rho_A(1_A), T(A ), r_{A}).
\eneq}}
{{Therefore \eqref{Cistbk-n1} holds.}}

Clearly, $B$ has the inductive limit decomposition $$B=\varinjlim(C_n\otimes U, \phi_n\otimes \mathrm{id}_U).$$
For each $n$, consider the positive homomorphism $(\phi_{n, \infty})_{*0}: K_0(C_n)\to K_0(B_1)\cong\rho_A(K_0(A))$. Since 
{{$\rho_A(K_0(A))$ is torsion free, $(\phi_{n, \infty})_{*0}(K_0(C_n))$ is a free abelian group.
There is a group}} \hm\, $\kappa_n: K_0(C_n)\to K_0(A)$ such that  {\Green{(recall that the order of $K_0(A)$ by traces)}}
is determined by
$$\rho_A\circ\kappa_n=(\phi_{n, \infty})_{*0} \andeqn \kappa_{{n}}([1_{C_n}]){{\le 2[1_A]}}.$$
{{Note that $\kappa_n$ is positive since $\rho_A\circ \kappa_n$ is.}}
By Corollary \ref{C0ext}, there is a  homomorphism $h'_n: C_n \to M_2(A)$ such that $(h'_n)_{*0}=\kappa_n.$ It is clear that $h_n:={{j_A\circ}}(h'_n\otimes \mathrm{id}_U)$ satisfies the desired condition.
\end{proof}

\begin{lem}\label{ExtTrace}
{Let $C\in {\cal C}$.
Let  $\sigma>0$  and let
${\cal H}\subset C_{s.a.}$ be {a} finite subset.
Let $A\in\mathcal B_{u0}$. Then for any $\kappa\in KL_e(C,A)^{++}$
and any continuous affine map $\gamma: T(A)\to T_f(C)$ {{(see \ref{Aq})}} which is compatible to $\kappa$, there is a unital homomorphism $\phi: C\to A$ such that
$$
[\phi]=\kappa\quad\textrm{and}\quad |\tau\circ \phi(h)-\gamma(\tau)(h)|<\sigma\tforal h\in {\cal H}\,\, {\Green{\text{and\,\,
for\,\, all}\,\,
\tau\in T(A).}}
$$
Moreover,
the conclusion above also holds if $C\in {\cal C}_0$ and $A\in {\cal B}_{u1}.$}
\end{lem}

\begin{proof}
We assume $C\in {\cal C}$ and $B\in {\cal B}_{u1}$ first. 
Let $\sigma>0.$
{{We}} may assume that every element of $\mathcal H$ has norm at most one.
Let  $\kappa$ and  $\gamma$ be given.  Define $\Delta: C_+^{q, 1}\setminus\{0\}\to (0, 1)$ by
$$\Delta(\hat{h})=\inf\{\gamma(\tau)(h)/2:\ \tau\in T({{A}})\}.$$

Let $\mathcal H_1\subset C_+{{\setminus \{0\}}}$, $\delta$, and $K$ be the finite subset, positive constant, and the positive integer of {\Green{Theorem}} \ref{ExtTraceC-D} with respect to $C$, $\Delta$, $\mathcal H,$ and $\sigma/{{4}}$ (in place of $\sigma$).
Let ${\cal P}\subset  {{M_m(C)}}$ {{(for some $m\ge 1$) be a finite subset of projections}}
 which generates $K_0(C)_+.$ {{We may assume
that ${\cal H}_1$ is in the unit ball of $C.$}}

{{Set
$$
\sigma_1=\min\{\sigma/2, \min\{\Delta({\hat{h}}): h\in {\cal H}_1\}\}>0.
$$}}
By Lemma \ref{cut-trace} {{(apply to $\sigma_1/(8+\sigma_1)$),}} there {{are}} a \SCA\, $D\subset  A$ with $D\in\mathcal C$
and with
$1_D=p\in A,$  {{and}} a continuous affine map $\gamma': T(D)\to T(C)$ such that
\begin{equation}\label{istTrC-eq1}
|\gamma'(\frac{1}{\tau(p)}\tau{\Green{|_D}})(f)-\gamma(\tau)(f)|<\sigma_{ 1}/(8+{\sigma_1})\rforal \tau\in T(A) \rforal f\in\mathcal H{{\cup {\cal H}_1}},
\end{equation}
where
$\tau(1-p)<\sigma_1/(8+\sigma_1)$ {{for all  $\tau\in T({\Green{A}}).$}}
{\Green{Moreover (by (2) of  Lemma \ref{cut-trace})}}
\begin{equation}\label{istTrC-eq2}
\gamma'(\tau)(h)>
\Delta(\hat{h})\rforal \tau\in T(D) \rforal h\in\mathcal H_{1}.
\end{equation}
Denote by $\imath: D\to pAp$  the embedding.
{\Green{Moreover, as by  (3)   and (4) of}}
 Lemma \ref{cut-trace}, there are {\Green{also}} positive homomorphisms $\kappa_{0,0}: K_0(C)\to K_0((1-p){{A}}(1-p))$ and $\kappa_{1,0}: K_0(C)\to K_0(D)$ such that  $\kappa_{1,0}$ is strictly positive, {{$\kappa_{1,0}([1_C])=[1_D],$}}  {{$\kappa_{1,0}$}} has multiplicity $K$, {{and}}
\beq\nonumber
\kappa|_{K_0({\Green{C}})}=\kappa_{0,0}+\imath_{*0}\circ \kappa_{1,0} \andeqn
|\gamma'(\tau)(q)-{{\rho_D}}(\kappa_{{1,0}}(q)){{(\tau)}}|<\delta,\quad q\in{\mathcal P},\ \tau\in T(D).
\eneq
{{Suppose that $A=A_1\otimes U$ for some $A_1\in {\cal B}_1$ and a  UHF-algebra $U$ of infinite type.}}
{{B}}y the last part of  \ref{cut-trace}, we may assume
that $\kappa_{0,0}$ is also strictly positive.
{{Let $x\in K_0(D)_+$ be such that $Kx=[1_D]=\kappa_{1,0}([1_C]).$ Then there is a projection
$e\in D$ such that $[e]=x$ as $D$ has stable rank one. Note  $e$ is full  and $D_1:=eDe$ {{belongs to}} ${\cal C}$ and
in case $D\in {\cal C}_0,$ $D_1\in {\cal C}_0$ (see \ref{cut-full-pj}).
We may write $D_1\otimes M_K=D.$}}
{{Define $\lambda: T(D_1)\to T(D)$ by $\lambda(t)(a\otimes b)=t(a){\rm tr}(b),$ where
${\rm tr}$ is the tracial state of $M_K.$  Define $\gamma'': T(D_1)\to T(C)$ by
$\gamma''(t)=\gamma'(\lambda(t)).$ Note that $K_0(D_1)=K_0(D).$}}
{{Then, for any $q\in {\cal P}$ and any $t\in T(D_1),$
\beq
&&\hspace{-0.4in}|\gamma''(t)({\Green{q}})-\rho_{D_1}((1/K)\kappa_{1,0})([q]))(t)|=|\gamma'(\lambda(t))({\Green{q}})-
\rho_{D}(\kappa_{1,0})([q])({\Green{\lambda(t)}})|<\dt\andeqn\\
&&\gamma''(t)(h)=\gamma'(\lambda(t))(h)
{\Green{\ge}} \Delta(\hat{h})\rforal h\in {\cal H}_1.\,\,\,\,\hspace{0.3in}
{\Green{{\text{(see\,\, also}\,\, \eqref{istTrC-eq2}).}}}
\eneq
}}
Therefore, by Lemma \ref{ExtTraceC-D},
there is a homomorphism $\phi_1: C\to {{M_K(D_1)=D\subset A}}$ {{such that
\beq
&&(\phi_1)_{*0}=K(1/K)\kappa_{1,0}=\kappa_{1,0}\andeqn\\\label{18-18n12-1}
&&|(1/K)(t\otimes {\rm Tr})( \phi_1(h))-\gamma''(t)(h)|<\sigma/4\rforal h\in {\cal H}\andeqn t\in T(D_1),
\eneq
where ${\rm Tr}$ is the unnormalized trace of $M_K.$}}
{{Note that $(1/K)(\tau\otimes {\rm Tr})=\lambda(\tau)$ for all $\tau\in T(D_1).$ By \eqref{18-18n12-1},}}
$$|\tau\circ\phi_1(h)-\gamma'(\tau)(h)|<\sigma/4\rforal h\in\mathcal H \andeqn \tau\in T(D).$$
Since $A$ is simple, $K_i((1-p)A(1-p))=K_i(A),$ $i=0,1.$
Let $\kappa_0=\kappa-[\imath\circ \phi_1]\in KL(C,A)=KL(C, (1-p)A(1-p)).$
Then
$\kappa_0|_{K_0(C)}=\kappa_{0,0}.$
Since $\kappa_0|_{K_0(C)}$ {{is equal to}}
$\kappa_{0{,}0},$ it is strictly positive
{{and $\kappa_0([1_C])=[1-p].$}} Note that $(1-p)A(1-p)\otimes U\cong (1-p)A(1-p).$ Therefore,
by Theorem \ref{preBot2}, since $C$ is semiprojective, there is a homomorphism $\phi_0: C\to (1-p)A(1-p)$ such that $[\phi_0]=\kappa_0.$
Note that this  holds for both the case that $A\in {\cal B}_{u0}$ {{(by  (i) of Theorem \ref{preBot2})}} and
the case that $C\in {\cal C}_0$ and $A\in {\cal B}_{u1}$  {{(apply (ii) of of Theorem \ref{preBot2}).}}

Consider the homomorphism $$\phi:=\phi_0\oplus\imath\circ \phi_1: C\to (1-p)A(1-p)\oplus D\subset A.$$ One has that $[\phi]=\kappa$ and, for all $h\in {\cal H}$ {\Green{and $\tau\in T(A),$}}
\begin{eqnarray*}
|\tau\circ \phi(h)-\gamma(\tau)(h)|&\leq & |\tau\circ \phi_1(h)-\gamma(\tau)(h)|+\sigma/4\\
&<&|\tau\circ \phi_1(h)-\gamma'(\frac{1}{\tau(p)}\tau|_D)({{h}})|+\sigma/2\\
&<&|\tau\circ \phi_1(h)- \frac{1}{\tau(p)}\tau\circ\phi_1(h)|+3\sigma/4<\sigma{,}
\end{eqnarray*}
as desired.
\end{proof}
It turns out that  KK-attainability implies the following existence theorem.
\begin{prop}\label{add-tr}
Let $A\in \mathcal B_0$, and assume that $A$ is KK-attainable with respect to ${\cal B}_{u0}$. Then for any $B\in \mathcal B_{u0}$, any $\alpha\in KL(A, B)^{{++}}$, and any continuous affine  map $\gamma: T(B)\to T(A)$ which is compatible with $\alpha$, there is a sequence of
{{\morp s}} $L_n: A\to B$ such that
 $$ \lim_{n\to\infty} ||L_n(ab)-L_n(a)L_n(b)||  =  0
 \tforal a, b \in A, $$
  $$ [{ \{}L_n{ \}}]  =  \alpha, \tand$$
$$ \lim_{n\to\infty}\sup\{|\tau\circ L_n(f)-\gamma(\tau)(f)|: \tau\in T(A)\} = 0 \tforal f\in A.$$
\end{prop}

\begin{proof}
The proof is the same as that of Proposition 9.7 of \cite{LinTAI}. Instead of using {{Corollary}} 9.6 of \cite{LinTAI}, one uses Lemma \ref{ExtTrace}.
\end{proof}


\section{The class ${\cal N}_1$}

Let $A$ be a unital \CA\, such that $A\otimes Q\in {\cal B}_1.$
In this section,  we will show that  $A\otimes Q\in {\cal B}_0.$ Note that this is proved without assuming
$A\otimes Q$ is nuclear.
However, it implies that ${\cal N}_1={\cal N}_0.$
If we assume that $A\otimes Q$ has finite nuclear dimension, {{then by}} using {{Lemma}} \ref{ExtTrace} and a characterization
of {{the class}}  $TAS$ by Winter, a much shorter proof of Theorem \ref{B0=B1} could {{be}}  given here.

\begin{lem}\label{K1inj}
Let $A\in {\cal B}_1$ be such that $A\cong A\otimes Q.$
Then the following property holds:
For any $\ep>0,$  any two non-zero mutually orthogonal elements $a_1, a_2\in A_+,$ and any finite subset ${\cal F}\subset A,$
there {{exist}} a projection $q\in A$ and a \SCA\, $C_1\in {\cal C}$ with $1_{C_1}=q$ such that
\begin{enumerate}
\item $\|[x, q]\|<\ep/16$, $x\in\mathcal F$,
\item ${{qxq}}\in_{\ep/16} C_1$, $x\in\mathcal F$, and
\item $1-q\lesssim a_1.$
\end{enumerate}

Suppose that $\Delta: (C_1)_{{+}}^{q,{\bf 1}}\setminus \{0\}\to (0,1)$ is an order preserving map
such that
$$
\tau(c)\ge \Delta(\hat{c})\rforal {{~~\tau \in T(A)~~ \mbox{and}~~}} c\in (C_1)^{\bf 1}_+\setminus \{0\}.
$$
{\rm(}By \ref{Ldet} such a $\Delta$  always exists.{\rm )}
Suppose also that
${\cal H}\subset (C_1)_+\setminus \{0\}$
and ${\cal F}_1\subset C_1$ are finite subsets. Then, there exist  another projection $p\in A$ with $p\le q,$
a \SCA\, $C_2\in {\cal C}$ {{of $A$}} with $1_{C_2}=p,$ and a  unital \hm\,
$H: C_1\to C_2$ such that

\begin{enumerate}\setcounter{enumi}{3}
\item $\|[x, p]\|<\ep/16$, $x\in {\cal F}$,
\item  $\|H(y)-{{p}}y{{p}}\|<\ep/16$, $y\in {\cal F}_1$, and

\item $1-p\lesssim a_1+a_2$.
\end{enumerate}
%
Moreover, {{$C_1$ and $H$ may be chosen {{in such a way}}  that
$K_1(C_1)$ may be written as}} $\Z^m\oplus G_0$ {{with}}
$ H_{*1}(G_0)=\{0\},$ and $H_{*1}|_{\Z^m}$ and $(j\circ H)_{*1}|_{\Z^m}$  both  injective,
where $j: C_2\to A$ is the embedding ($m$ could be zero, and in this case, $G_0=K_1(C_1)$).
Furthermore, we may {{choose $C_2$ and $H$ such that}}
\begin{equation*}
\tau(j\circ H(c))\ge  3\Delta(\hat{c})/4\tforal c\in {\cal H}\tand\tforal \tau\in T(A).
\end{equation*}

\end{lem}
\begin{proof}


Since $A\in {\cal B}_1,$
there exist a projection $q\in A$ and a \SCA\, $C_1\in {\cal C}$ {{of $A$}} with $1_{C_1}=q$ such that

(a) $\|[x, q]\|<\ep/16$, $x\in\mathcal F$,

(b) $qxq\in_{\ep/16} C_1$, $x\in\mathcal F$, and

(c) $1-q\lesssim a_1.$

There are two non-zero mutually orthogonal elements $a_2'$ and $a_3$ {{in $\overline{a_2Aa_2}.$}}
{{Since}} $A\cong A\otimes Q,$
$K_1(A)$ is torsion free. Denote by ${{\iota}}: C_1\to qAq$ the embedding.
Since $K_1(C_1)$ is finitely generated,
we may write $K_1(C_1)=G_1\oplus G_0,$ where $G_1\cong \Z^{m_1},$ ${{\iota}}_{*1}|_{G_1}$ is injective, and ${{\iota}}_{{*1}}|_{G_0}=0.$
Define
$$
\sigma={{\min}}\{\Delta(\hat{h})/{{16}}: h\in {\cal H}\}>0.
$$
Choose an element $a_2''\in \overline{a_2'Aa_2'}$ such that
$d_\tau(a_2'')<\sigma$ for all $\tau\in T(A)$
{{(recall $d_{\tau} (a) =\lim_{n\to \infty} \tau (a^{\frac1n})$)}}.

Suppose that $G_0$ is generated by $[v_1], [v_2],...,[v_l],$ where $v_i\in U(C_1)$ (note {{
that, by \ref{2pg3},}} $C_1$ has stable rank one).
Then, we may write ${{\iota(v_k)}}=\prod_{s=1}^{l(k)} \exp (i h_{s,k})$ (since $\imath_{*1}([v_k])=0$ {{in}} $K_1(A)$),
where $h_{s,k}\in (qAq)_{s.a.},$ $s=1,2,...,l(k),$ $k=1,2,...,l.$

Let ${\cal F}_1$ be a finite subset of $C_1$
{{with}} the following property:
if $x\in {\cal F},$ {{then}} there is $y\in {\cal F}_1$ such that $\|{{q}}x{{q}}-y\|<\ep/16.$

Let ${\cal F}_2$ be a finite subset  of $qAq,$ to be determined, which at least contains ${\cal F}\cup {\cal F}_1\cup{\cal H}$ and
{{the elements}} $h_{s,k},$ $\exp(i h_{s,k}),$ $s=1,2,...,l(k),$ $k=1,2,...,l.$

Let $0<\dt<\min\{\ep/64, \sigma/4\},$ to be determined.   Since, {{by \ref{B1hered},}} $qAq\in {\cal B}_1,$
one obtains a non-zero projection $q_1\in qAq$ and a \SCA\, $C_2{{\subset}} qAq$  such that

(d) $\|[x,\, q_1]\|<\dt$, $x\in {\cal F}_2$,

(e) $q_1xq_1\in_{\dt} C_2$, $x\in {\cal F}_2$, {{and}}

(f) $1-q_1\lesssim a_2{{''}}.$

With sufficiently small $\dt$ and large ${\cal F}_2,$ using the semiprojectivity of $C_1,$
we obtain unital \hm s $h_1: C_1\to C_2$
such that
\beq\label{K1inj-3}
\|h_1(a)-q_1aq_1\|<\min\{\ep/16, \sigma\}\rforal a\in {\cal F}_1\cup {\cal H}.
\eneq
By the choice of $\sigma$ and $a_2'',$ and by \eqref{K1inj-3},
one computes that, {{for all $\tau\in T(A),$}}
\beq\label{K1inj-4}
\tau(j\circ h_1(c))\ge {{7}}\Delta(\hat{c})/{{8}} \rforal c\in {\cal H},
\eneq
where  $j$ is  the embedding from $C_2$ into $qAq.$
Note that, when ${\cal F}_2$ is large enough, $(h_1)_{*1}(G_0)=\{0\}.$
{{Since $K_1(C)=\Z^{m_1}\oplus G_0,$ $G_0$ is a direct summand of ${\rm ker}\, (h_1)_{*1}.$}}

We may write $K_1(C_1)=G_2\oplus G_{2,0}\oplus G_{2,0,0},$ where $G_2\cong \Z^{m_2}$ with $m_2\le m_1,$
$G_2$ is a subgroup of $G_1,$ $G_{2,0,0}\supset G_0,$  $(h_1)_{*1}(G_{2,0,0})=\{0\},$
$(j\circ h_1)_{*1}|_{G_{2,0}}=0,$ and $(j\circ h_1)_{*1}|_{G_2}$  is injective.
Here we use the fact that $K_1(A)$ is torsion free {{(and $G_{2,0,0}$ could be chosen to be just $G_0$ itself).}}
If $G_{2,0}=\{0\}${{,}} we are done.
Otherwise, $m_2<m_1.$ We also note that $(j\circ h_1)_{*1}(G_{2,0}\oplus G_{2,0,0})=\{0\}.$

{{We will  repeat the process}} to construct $h_2$, and consider  $h_2\circ h_1.$
Then we may write $K_1(C_1)=G_3\oplus G_{3,0}\oplus G_{3,0,0}$ with
$G_3\cong \Z^{m_3}$ ($m_3\le m_2$), $G_3\subset G_2,$ $G_{3,0,0}\supset G_{2,0}\oplus G_{2,0,0},$
$(j\circ h_2\circ h_1)_{*1}(G_{3,0})=\{0\},$ and $(j\circ h_2\circ h_1)_{*1}|_{G_3}$  injective.
Again, if $G_{3,0}=\{0\},$ we are done (choose $H=h_2\circ h_1$).  Otherwise, $m_3<m_2<m_1.$
We continue this process.
This process stops at  a finite stage since $m_1$ is finite.
This proves the lemma.
%
\end{proof}

\begin{thm}\label{B0=B1}
Let $A_1$ be a unital simple separable
\CA\,  such that $A_1\otimes Q\in {\cal B}_1.$  Then $A_1\otimes Q\in {\cal B}_0.$
\end{thm}

\begin{proof}
Let $A=A_1\otimes Q.$ Suppose that $A\in {\cal B}_1.$
Let $\ep>0,$  let $a\in A_+\setminus\{0\},$ and let ${\cal F}\subset A$ {{be a finite subset.}}
Since $A$ has {{the}} property (SP) (see \ref{Psp}),  we obtain three non-zero and mutually orthogonal projections $e_0, e_1, e_2\in \overline{aAa}.$
There exist a projection $q_1\in A$ and a \SCA\, $C_1\in {\cal C}$ {{of $A$}} with $1_{C_1}=q_1$ such that
\beq\label{B0=B1-n1}
&&\|[x,\, q_1]\|<\ep/16\tand q_1xq_1\in_{\ep/16} C_1 \tforal x\in {\cal F}, \tand \\
&&
1-q_1\lesssim e_0.
\eneq
Let ${\cal F}_1\subset C_1\subset A$ be a finite subset
such that, for any $x\in {\cal F},$ there is $y\in {\cal F}_1$ such that
$\|q_1xq_1-y\|<\ep/16.$
For each $h\in (C_1)_+\setminus\{0\},$ define
$$
\Delta(\hat{h})=(1/2)\inf\{\tau(h): \tau\in T(A)\}.
$$

Then  the function $\Delta: (C_1)_+^{q, {\bf 1}}\setminus \{0\}\to (0,1)$ preserves the order.
Let ${\cal H}_1\subset (C_1)^{\bf 1}_+\setminus \{0\}$ be a finite subset, $\gamma_1, \gamma_2>0,$  $\dt>0,$ ${\cal G}\subset C_1$ be a finite subset,
${\cal P}\subset \underline{K}(C_1)$ be a finite subset, ${\cal H}_2\subset (C_1)_{s.a.}$ be a finite subset, and
${\cal U}\subset {{J_c(K_1(C_1)\subset U(C_1)/CU(C_1)}}$ {{(see \ref{Dcu})}}
{{be a}} finite subset
{{as provided}}   by  Theorem \ref{UniCtoA} (and Corollary \ref{RemUniCtoA}) for $C=C_1,$
$\ep/16$ (in place of $\ep$),  ${\cal F}_1$ (in place of ${\cal F}$), and $\Delta/2$ (in place of $\Delta$).

By {{Lemma}} \ref{K1inj}, there exist  another projection $q_2\in A$ such that
$q_2\le q_1,$
a \SCA\, $C_2\in {\cal C}$ {{of $A$}} with $q_2=1_{C_2},$
and a unital \hm\,
$H: C_1\to C_2$ such that
\beq\label{K1inj-2}
&&\|[x, q_2]\|<\ep/16\tforal x\in {\cal F},\\\label{K1inj-2+}
&&\|{{j\circ}}H(y)-q_2yq_2\|<\ep/16\tforal y\in {\cal F}_1,\\
&&\tau(j\circ H(c))\ge 3 \Delta(\hat{c})/4\tforal c\in {\cal H}\andeqn \rforal \tau\in T(A), \andeqn\\\label{K1inj-2++}
&&1-q_2\lesssim e_0+e_1.
\eneq
Moreover, we may write $K_1(C_1)=\Z^m\oplus G_{{00}},$ where
$H_{*1}(G_{00})=\{0\},$ {{and}}  $H_{*1}|_{\Z^m}$ and $( j\circ H)_{*1}|_{\Z^m}$ are injective, where
$j: C_2\to A$ is the embedding.
Let $A_2=q_2Aq_2$ and denote by $j_1: C_2\to A_2$ the embedding.
{{Put
\beq\label{1883-n+}
\sigma_0=\inf\{\tau(e_2): \tau\in T(A)\}>0.
\eneq}}
By Theorem \ref{RangT}, there exists a unital  simple \CA\, $B\cong B\otimes Q$  {{with}}
$B=\varinjlim
(B_n, \imath_n)$  {{in such a way}}
that each $B_n$ is equal to $B_{n,0}\oplus B_{n,1}$ with $B_{n,0}\in {\bf H}$
{{(see \ref{AHblock})}} and
$B_{n,1}\in {\cal C}_0,$  {{each}}  $\imath_n$ is injective, {{and}}
\beq\label{B0=B1-2}
\lim_{n\to\infty}\max\{\tau(1_{B_{n,0}}):\tau\in T(B)\}=0,\andeqn {\text{such that}}\,\,
{\rm Ell}(B)={\rm Ell}(A_2).
\eneq

We may assume {{(choosing a smaller $\dt$)}} that ${\cal U}={\cal U}_1\cup {\cal U}_0,$
with $\pi({\cal U}_1)$ generating $\Z^m\subset K_1(C_1)$ and $\pi({\cal U}_0)\subset G_{{00}},$
where $\pi: U(C_1)/CU(C_1)\to K_1(C_1)$ is the quotient map. {{Recall that}}
$J_c: K_1(C_1)\to U(C_1)/CU(C_1)$ is a fixed splitting map as defined in \ref{Dcu}.
Suppose that ${\bar v_1},{\bar v_2},...,{\bar v_m}$ form a set of {{independent}} generators for $J_c^{(1)}(\Z^m){{\cong \Z^m}}.$
{{We}} may assume
that ${\cal U}_1=\{{\bar v_1},{\bar v_2},...,{\bar v_m}\}.$
Put
$$
\gamma_3=\min\{\Delta(\hat{h})/2: h\in{\cal H}_1\}.
$$
Note $H^{\ddag}({\cal U}_0)\subset U_0(C_2)/CU(C_2).$ {{By Lemma \ref{HvsU-lem-2018}, we may}}
choose a finite subset ${\cal H}_3\subset (C_2)_{s.a.}$  and $\sigma>0$ {{with}} the following property:
for any two unital \hm s $h_1, h_2: C_2\to D,$ for any unital \CA\, $D$ of stable rank one,
if, {{for all $\tau\in T(D),$}}
\beq\label{B0=B1-n3}
|\tau\circ h_1(g)-\tau\circ h_2(g)|<\sigma\tforal g\in {\cal H}_3,
\eneq
then
\beq\label{B0=B1-n4}
{\rm dist}(h_1^{\ddag}({\bar v}), h_2^{\ddag}({\bar v}))<\gamma_2/8
\eneq
for all ${\bar v}\in H^{\ddag}( {\cal U}_0)\subset U_0(C_2)/CU(C_2).$
\Wlog, we may assume that $\|h\|\le 1$ for all $h\in {\cal H}_1\cup {\cal H}_2\cup {\cal H}_3.$

Let ${{\Gamma=(\kappa_0, \kappa_1, \kappa_T)}}: {\rm Ell}(A_2)\to {\rm Ell}(B)$  be the above identification,
where
$\kappa_0: K_0(A_2)\to K_0(B)$ is  an order isomorphism with $\kappa_0([q_2])=[1_B],$
$\kappa_1: K_1(A_2)\to K_1(B)$ is an isomorphism, and $\kappa_T: T(A_2)\to T(B)$ is
an affine homeomorphism such that
$r_{A_2}(\kappa_T^{-1}(\tau))(x)=r_B(\tau)(\kappa_0(x))$ for all $x\in K_0(A_2)$ and for all $\tau\in T(B)$ {{(or,
equivalently, $r_{A_2}(t)(\kappa_0^{-1}(y))=r_{B}(\kappa_T(t))(y)$ for all $y\in K_0(B)$ and for all $t\in T(A_2)$).}}
Since $B$ satisfies the UCT, there is an element $\kappa^{-1}\in KK_e(B, A_2)^{++}$ {{(see \ref{DKLtriple})}}
such that $\kappa^{-1}|_{K_i(B)}=\kappa_i^{-1},$ $i=0,1.$
{{Note that since $A_2\cong A_2\otimes Q,$ one has  ${\rm Hom}_{\Lambda}(\underline{K}(B), \underline{K}(A_2))=
{\rm Hom}(K_*(B), K_*(A_2)).$}}

Let $\kappa^{(2)}\in KK_e(C_2,B)^{++}$ be such that
$\kappa^{(2)}|_{K_i(C_2)}=\kappa_i\circ j_{*i},$ $i=0,1.$
{{In fact, since $K_i(B)=K_i(B)\otimes Q$ ($i=0,1$),  {{by the UCT,}} $\kappa^{(2)}$ is uniquely determined
by $\kappa_i\circ j_{*i}.$}}

It follows from Lemma \ref{ExtTrace} that there exists a unital \hm\, $\phi: C_2\to B$
such that
\beq\label{B0=B1-4}
&&[\phi]={{\kappa^{(2)}}}
\andeqn\\\label{B0=B1-4+}
&& |\tau(\phi(h))-\gamma(\tau)(h)|<\min\{\gamma_1, \gamma_2,\gamma_3,\sigma\}/8
\eneq
for all $h\in  H({\cal H}_1)\cup H({\cal H}_2)\cup {\cal H}_3,$ where $\gamma: T(B)\to T_f(C_{{2}})$ is induced
by $\kappa_T$ and the embedding $j.$
In particular, $\phi_{*1}$ is injective on $H_{*1}(\Z^m)$ {{(since $j_{*1}$ is injective
on $H_{*1}(\Z^m)$ and $\kappa_1$ is an isomorphism).}}


Since $C_2$ is semiprojective, \wilog, we may assume  $\phi(C_2)\subset {{j_{1, \infty}(B_1)\subset B}}.$
{{In what follows, we use  the notation $\phi'$ for the \hm\, from $C_2$ to $B_1$ such that
$\phi=j_{1, \infty}\circ \phi'.$}}
We may also assume
that, {{in the decomposition}}  $B_1=B_{1,0}\oplus B_{1,1},$
\beq\label{B0=B1-5}
\tau(1_{B_{1,0}})<\min\{\sigma_0/4, \gamma_1/8, \gamma_2/8, \gamma_3/8, \sigma/8\}\tforal \tau\in T(B).
\eneq
{{W}}e have  $(\imath_{1, \infty} \circ \phi')_{*1}=\phi_{*1}$, and $\phi_{*1}$  is injective on $H_{*1}(\Z^m)$. Consequently
$(\imath_{1, \infty})_{*1}$ is injective on $(\phi')_{*1}(H_{*1}(\Z^m)).$

Let $G_1=H^{\ddag}(J_c^{(1)}(\Z^m))\subset U(C_2)/CU(C_2)$ and
$G_0=H^{\ddag}(J_c^{(1)}(G_{00})).$ Note, by the construction, $G_0\subset U_0(C_2)/CU(C_2).$
Since $\kappa|_{K_1(B)}$ is an isomorphism and
$(\imath_{1, \infty})_{*1}$ is injective on $(\phi')_{*1}(H_{*1}(\Z^m)),$ $\phi^{\ddag}|_{G_1}$ is injective.
Let ${\cal H}_4=P(\phi({\cal H}_1\cup {\cal H}_2\cup {\cal H}_3)),$ where $P: B_1\to B_{1,1}$  is the projection.

{{Note that since $A\cong A\otimes Q,$ $KL(B_1, A)={\rm Hom}(K_*(B_1), K_*(A)).$}}
Let $e_3\in A_2$ be a projection such that $[e_3]=\kappa_0^{-1}([\imath_{1, \infty})(1_{B_{1,1}})]),$ {{and
let $A_3=e_3Ae_3.$}}
There is $\kappa_{B_{1,1}}'\in KK(B_{1,1}, A_2)^{++}$ such that $\kappa'_{B_{1,1}}$ induces
$(\kappa_*^{-1}\circ (\imath_{1,\infty})_*)|_{K_*(B_{1,1})}.$ {{As  $K_i(B_{1,1})$ is finitely generated and
$A\cong A\otimes Q$ {{(and $K_*(A_3)=K_*(A)$),  by the UCT,}} the map $\kappa_{B_{1,1}}'$ is uniquely determined by $(\kappa_*^{-1}\circ (\imath_{1,\infty})_*)|_{K_*(B_{1,1})}.$}}
It follows from the second part of  Lemma \ref{ExtTrace}
that there is a unital monomorphism
$\psi_1: B_{1,1}\to e_3Ae_3$ such that
\beq\label{B0=B1-6}
\hspace{-0.6in}&&[\psi_1]=\kappa_{B_{1,1}}'
\andeqn\\\label{B0=B1-6+}
\hspace{-0.6in}&&|\tau\circ \psi_1(g)-\gamma'(\tau)(\imath_{1,\infty}(g))|<\min\{\tau(([e_3]))/4, \gamma_1/8, \gamma_2/8, \gamma_3/8, \sigma/8\}\rforal g\in {\cal H}_4
\eneq
for all $\tau\in T(A_3),$ where $\gamma': T(A_3)\to T_f(B_{1,1})$ {{is the map}}  induced by $\kappa_T^{-1}$ and $\imath_{1, \infty}.$

Write $B_{1,0}=B_{1,0,1}\oplus B_{1,0,2},$ where $B_{1,0,1}$ is a finite direct sum of circle algebras and
$K_1(B_{1,0,2})$ is finite {{(according to the definition of ${\bf H}$---see \ref{AHblock}).}}
Since $B\cong B\otimes Q,$  we may assume that $(\imath_{1,\infty})_{*1}|_{K_1(B_{1,0.2})}=\{0\}.$

Since $A_2\cong A_2\otimes Q,$ by parts  (iii) and (iv) of Theorem \ref{preBot2},
there exists a unital \hm\, $\psi_2: B_{1,0,2}\to e_4A_2e_4$
such that $[\psi_{{2}}]$ is equal to $\kappa_{B_{1,0,2}}',$
which  induces
$\kappa_*^{-1}\circ (\imath_{1,\infty})_{*}|_{K_*(B_{1,0,{{2}}})},$ where $e_4\in A_2$ is a projection
orthogonal to $e_3$ and $[e_4]=\kappa_0^{-1}\circ (\imath_{1,\infty})_{*0}([1_{B_{1,0,2}}]).$
{{We have}} $(\psi_2)_{*1}=0,$ {{since $K_i(A)$ is torsion free.}}
Let $\psi_3: B_{1,1}\oplus B_{1,0,2}\to (e_3+e_4)A_2(e_3+e_4)$ be defined by $\psi_3{{(a,b)}}={{\psi}}_1{{(a)}}\oplus
\psi_2{{(b)}}$ {{for $a\in B_{1,1}$ and $b\in B_{1,0,2}.$}}

Let $P_1: B_1\to B_{1,0,1}.$ Then, since $(\imath_{1, \infty})_{*1}|_{K_1(B_{1,0,2})}=\{0\},$ {{and $K_1(B_{1,1})=\{0\},$
the restriction}}
 $(P_1)_{*1}|_{\phi_{*1}(H_{*1}(\Z^m))}$
is injective. Also, $P_1^{\ddag}$ is injective on $\phi^{\ddag}\circ H^{\ddag}(J_c^{(1)}(\Z^m)).$
Put $G_1'=P_1^{\ddag}\circ \phi^{\ddag}\circ H^{\ddag}(J_c^{(1)}(\Z^m)) { {\subset U(B_{1,0,1})/CU(B_{1,0,1})}}.$ Then $G_1'\cong \Z^m.$

It follows from Theorem \ref{preBot2} that there is a unital \hm\, $\psi_4': B_{1,0,1}\to (q_2-e_3-e_4)A_2(q_2-e_3-e_4)$
such that $[\psi_4']=\kappa_{B_{1,0,1}}'$ which
{{induces}} $(\kappa_*^{-1}\circ (\imath_{1, \infty})_*)|_{K_*(B_{1,0,1})}.$
Let
\beq\label{B0=B1-15}
z_i= P_1^{\ddag}\circ \phi'^{\ddag}\circ H^{\ddag}({\bar v}_i)\andeqn
\xi_i=\psi_3^{\ddag}\circ ({{\rm {\id}}}_{B_1}-P_1)^{\ddag}\circ \phi'^{\ddag}\circ H^{\ddag}({\bar v}_i),
\eneq
$i=1,2,...,m.$
It should be noted that, since $(\psi_1)_{*1}=(\psi_2)_{*1}=0,$ $\xi_i\in U_0(A_2)/CU(A_2),$ $i=1,2,...,m.$
Moreover, since
\beq\label{B0=B1-16}
(\psi_4')_{*1}\circ (P_1)_{*1}\circ \phi'_{*1}(x)=j_{*1}(x)\tforal x\in  H_{*1}(\Z^m)\subset K_1(C_2),
\eneq
\beq\label{B0=B1-17}
\pi((j\circ H)^{\ddag}({\bar v}_i))\pi((\psi_4')^{\ddag}(z_i))^{-1}=0\,\,\, {\rm in}\,\,\, K_1(A).
\eneq

Define a \hm\, $\lambda: G_1'\to U(A_2)/CU(A_2)$ by
\beq\label{B0=B1-18}
\lambda(z_i)=(j\circ H)^{\ddag}({\bar v}_i)((\psi_4')^{\ddag}(z_i)\xi_i)^{-1},\,\,\,i=1,2,...,m.
\eneq
Note that $\lambda(z_i)\in U_0(A)/CU(A),$  $i=1,2,...,m.$
By  Lemma \ref{UCUdiv},  {{the abelian group}} $U_0(A)/CU(A)$ is divisible. There exists  a \hm\, $\bar{\lambda}: U(B_{1,0,1})/CU(B_{1,0,1})
\to U(A_2)/CU(A_2)$
such that ${\bar \lambda}|_{G_1'}=\lambda.$
Define a \hm\, $\lambda_1:  U(B_{1,0,1})/CU(B_{1,0,1})\to U(A_2)/CU(A_2)$ by
$\lambda_1(x)=(\psi_4')^{\ddag}(x){\bar \lambda}(x)$ for all $x\in U(B_{1,0,1})/CU(B_{1,0,1}).$
By {{Theorem}} \ref{UCUiso}, the  \hm\,
$$
U((q_2-e_3-e_4)A_2(q_2-e_3-e_4))/CU((q_2-e_3-e_4)A_2(q_2-e_3-e_4))\to
U(A_2)/CU(A_2)
$$
given by sending $u$ to $\diag(u, e_3+e_4)$ is an isomorphism.  Since $B_{1,0,1}$ is a circle algebra, one easily obtains a unital \hm\,
$$
\psi_4: B_{1,0,1}\to (q_2-e_3-e_4)A_2(q_2-e_3-e_4)
$$
such that
\beq\label{B0=B1-19}
(\psi_4)_{*i}=(\psi_4')_{*i}=\kappa_i^{-1}\circ (\imath_{1, \infty})_{*i}|_{B_{1,0,1}} \,\,(i=0,1)\andeqn \psi_4^{\ddag}=\lambda_1.
\eneq
Define $\psi: B_1\to A_2$ by $\psi{{(a,b)}}=\psi_3{{(a)}}\oplus \psi_4{{(b)}}$
{{for $a\in B_{1,1}\oplus B_{1,0,2}$
and $b\in B_{1,0,1}.$}}
Then, we have
\beq\label{B0=B1-20}
\psi_{*i}=\kappa_i^{-1}\circ (\imath_{1, \infty})_{*i},\,\,\, i=0,1.
\eneq
Since $(\imath_{1, \infty}\circ \phi')_{*i}=\kappa_i\circ j_{*i},$ $i=0,1,$
we compute  that
\beq\label{B0=B1-21-}
(\psi\circ \phi'\circ H)_{*i}=\kappa_i^{-1}\circ (\imath_{1, \infty})_{*i}\circ (\phi'\circ H)_{*i}=(j\circ H)_{*i},\,\,\,i=0,1.
\eneq
{{Thus, since $A\cong A\otimes Q$ and $C_1$ satisfies the UCT,}}
\beq\label{B0=B1-21}
[\psi\circ \phi'\circ H]=[j\circ H]\,\,\,{{\mbox{in}~~}}KK(C_1,A).
\eneq
{{On the other hand,
\beq\label{1883-n1}
\tau(\psi(1_{B_{1,0}}))&=&r_{A_2}(\tau)(\kappa_{*0}^{-1}((\imath_{1, \infty})_{*0}([1_{B_{1,0}}]))\\
&=&
r_B(\kappa_T(\tau))((\imath_{1, \infty})_{*0}([1_{B_{1,0}}]))<\sigma_0\tforal \tau\in T(A).
\eneq
By  \eqref{1883-n+} and \eqref{B0=B1-5}, and by \ref{Comparison},
\beq\label{18-0803-n1}
[\psi(1_{B_{1,0}})]\le [e_2].
\eneq}}
By (\ref{B0=B1-4+}), (\ref{B0=B1-5}), and (\ref{B0=B1-6+}),
\beq\label{B0=B1-23}
|\tau\circ \psi\circ \phi'(g)-\tau(g)|<\min\{\sigma, \gamma_1, \gamma_2, \gamma_3,\sigma\}\rforal g\in H({\cal H}_1)\cup  H({\cal H}_2)\cup {\cal H}_3.
\eneq
{{Denote  by $\psi'$  the composition $\psi\circ \phi'.$  {{We}} have, for all $\tau\in T(A_2),$}}
\beq\label{B0=B1-24}
|\tau\circ \psi'(g)-\tau(j(g))|<\min\{ \gamma_1, \gamma_2, \gamma_4\}\rforal g\in H({\cal H}_1)\cup H({\cal H}_2).
\eneq
We then compute that
\beq\label{B0=B1-24+}
\tau\circ \psi'({{H(h)}})\ge \Delta(\hat{h})/2\rforal h\in {\cal H}_1,\andeqn \tau\in T(A_2).
\eneq

By  (\ref{B0=B1-18}) and the definition of $\lambda_1,$ we have
\beq\label{B0=B1-25}
(\psi'\circ H)^{\ddag}({\bar v}_i)&=&\psi^{\ddag}\circ (\phi'\circ H)^{\ddag}({\bar v}_i)\\
&=&\left(\psi_3^{\ddag}\circ ({\rm id}_{B_1}-P_1)^{\ddag}\circ (\phi'^{\ddag}\circ H^{\ddag})({\bar v}_i)\right)
\left(\psi_4^{\ddag}(P_1^{\ddag}(\phi'^{\ddag}\circ H^{\ddag}({\bar v}_i))\right)\\
&=&\xi_i\lambda_1(z_i)=(j\circ H)^{\ddag}({\bar v}_i),\,\,\, i=1,2,...,m.
\eneq
By the choice of ${\cal H}_3$ and $\sigma,$
we also have
\beq\label{B0=B1-26}
{\rm dist}((\psi')^{\ddag}(H({\bar v})), (j\circ H)^{\ddag}({\bar v}))<
\gamma_2\rforal {\bar v}\in {\cal U}_0.
\eneq
It follows from {{Theorem}} \ref{UniCtoA} and {{Corollary}} \ref{RemUniCtoA} ({{on using}}
\eqref{B0=B1-21}, \eqref{B0=B1-24},
\eqref{B0=B1-24+}, and \eqref{B0=B1-26}) that there is a unitary $U\in q_2Aq_2$ such that
\beq\label{B0=B1-27}
\|{\rm Ad}\, U\circ \psi'({{H(y)}})-{{j\circ H(y)}}\|<\ep/16\rforal {{y}}\in {\cal F}_1.
\eneq
Let $C_3={\rm Ad}\, U\circ \psi(B_{1,1})$ and  $p=1_{C_3}=
{\rm Ad}\, U\circ \psi(1_{B_{1,1}}).$
{{Since $1_{B_{1,1}}$ is in the center of $B_1,$}}
for any $y\in {\cal F}_1,$ $\phi'(y)1_{B_{1,1}}=1_{B_{1,1,}}\phi'(y).$
Therefore, ${\rm Ad}\, U\circ \psi'(y)p=p{\rm Ad}\, U\circ \psi\circ \phi'(y)=p{\rm Ad}\, U\circ \psi'(y).$
{{Note that, for any $y\in F_1$, we have
$$\|y-j\circ H(y)\oplus (1-q_2)y(1-q_2)\|<\ep/8,$$
and $p< q_2$. }}
Hence,
\beq\label{B0=B1-28}
\|py-yp\|\le \|py-p{\rm Ad}\,U\circ \psi'(y)\|+\|p{\rm Ad}\,U\circ \psi'(y)-yp\|<\ep/8+\ep/8=\ep/4
\eneq
for all $y\in {\cal F}_1.$
{{Combining this}} with
 (\ref{K1inj-2}), {{we have}}
\beq\label{B0=B1-29}
\|px-xp\|=\|pq_2x-xq_2p\|<2\ep/16+\|pq_2xq_2-q_2xq_2p\|<\ep\rforal x\in {\cal F}.
\eneq
Let $x\in {\cal F}.$ Choose $y\in {\cal F}_1$ such that $\|q_2xq_2-q_2yq_2\|
{{=\|q_2(q_1xq_1-y)q_2\|<\ep/16.}}$ Then, by
\eqref{K1inj-2+} and \eqref{B0=B1-27},
\beq\label{B0=B1-30}
\hspace{-0.5in}\|pxp-p({\rm Ad}\,U\circ \psi'(H(y)))p\| &\le & \|pxp-pq_2xq_2p\|
+\|pq_2xq_2p-pq_2yq_2p\|\\
&&\hspace{0.2in}+\|pq_2yq_2p-p(j\circ H(y))p\|\\
&&\hspace{0.4in}+\|p (j\circ H(y))p-p({\rm Ad}\, U\circ \psi'(H(y)))p\|\\
&&<{{0+\ep/16+}}\ep/16+\ep/16=3\ep/8.
\eneq	
However, $p{\rm Ad}\, U\circ \psi'(y)p={\rm Ad}\, U\circ(\psi({{1_{B_{1,1}}\phi'(y)1_{B_{1,1}}}}))\in C_3$ for all $y\in {\cal F}_1.$
Therefore,
\beq\label{B0=B1-31}
pxp\in_{\ep} C_3.
\eneq
We then estimate that, by \eqref{K1inj-2++} and by \eqref{18-0803-n1},
\beq\label{B0=B1-32}
[1-p]\le [1-q_2]+[\psi'(1_{B_{1,0}})]\le [e_0\oplus e_1\oplus e_2]\le [a].
\eneq
Since $B_{1,1}\in {\cal C}_0,$ by applying \ref{subapprox},  {{$C_3$ can be approximated by
\SCA s in ${\cal C}_0.$ It follows that}}
$A\in {\cal B}_0.$
\end{proof}

\begin{cor}\label{B0=B1tU}
If $A$ is a unital separable amenable simple \CA\, such that $A\otimes Q\in {\cal B}_1,$  then, for any infinite dimensional UHF-algebra $U,$ $A\otimes U\in {\cal B}_0.$
\end{cor}

\begin{proof}
It follows from {{Theorem}} \ref{B0=B1} that $A\otimes Q\in {\cal B}_0.$ Then, by  {{Lemma}} \ref{subapprox} and by
{{Theorem 3.4 of}} \cite{LS}, $A\otimes U\in {\cal B}_0$ for
every infinite dimensional UHF-algebra.
\end{proof}



\section{KK-attainability of the C*-algebras in $\mathcal B_0$}

{{The main purpose of this section is to establish Theorem \ref{MEST}. It is an} existence theorem for maps
from an algebra
in the {{sub-}}class $\mathcal B_0$ to a \CA\,  {{as}}  in Theorem \ref{RangT} {{(and elsewhere).}}
The {{construction}} is similar to that of Section 2 of \cite{Lnduke}, and, roughly, we will construct a map factoring through a \SCA\,
(in $\mathcal C_0$) of the given \CA\,  in $\mathcal B_0$, and also require this map to carry the given KL-element. But since the ordered $K_0$-group of a \CA\,  in $\mathcal C_0$ in general is not a Riesz group,
extra work has to be done to take care of this {{difficulty}} (similar work also appears in \cite{Niu-TAS-I}, \cite{Niu-TAS-II}, and \cite{Niu-TAS-III}).




\begin{NN}\label{201} Let us proceed as {{in}} Section 2 of \cite{Lnduke}. Let $A\in\mathcal B_0$
{{be a separable \CA\,}} and assume that $A$ has the  {{property}} (SP). By Lemma \ref{MF}, the \CA\, $A$ can be embedded as a C*-subalgebra of $\prod M_{n_k}/\bigoplus M_{n_k}$ for some {{sequence}} $\{n_k\};$  {{in other words,}}
$A$ is MF in the sense of Blackadar and Kirchberg (Theorem 3.2.2 of \cite{BK-inductive}). Since $A$ is assumed to be amenable, by Theorem 5.2.2 of \cite{BK-inductive}, the \CA\, $A$ is strong NF, and hence, by Proposition 6.1.6 of \cite{BK-inductive}, there is an increasing family of {{residually finite dimensional}} sub-\CA s $\{A_n\}$ {{with}} union dense in $A$. {\blue{Let us assume that $1\in A_1$.}}

{{Let us set up the initial stage:}}
{{Choose}} a dense sequence  $\{{{x_0,}} x_1, x_2,..., x_n,...\}$  of elements in the unit ball of $A$.
Let $\mathcal P_0\subset M_\infty(A)$ be a finite subset of projections. We assume that
$x_{{0}}{{\in A_{{1}}}}$ and $\mathcal P_0\subset M_\infty(A_1)$.
{{Choose}}  a finite subset
 ${\cal{F}}_{0}$ {{in the unit ball of}} $A_1\subset A$ with {\blue{$\{1, x_{{0}}\}\subset \mathcal F_{0}$,}}
  {{and $\eta_0>0$ such that
$[L']|_{{\cal P}_0}$ is well defined for any ${\cal F}_0$-$2\eta_0$-multiplicative \morp\, $L'$
from $A.$   Moreover, by  choosing even larger ${\cal F}_0$ and smaller $\eta_0,$
we may assume that
\beq
[L']|_{{\cal P}_0}=[L'']|_{{\cal P}_0},
\eneq
if  $\|L'(x)-L''(x)\|<2\eta_0$ for all $x\in {\cal F}_0$
and $L''$ is also an ${\cal F}_0$-$2\eta_0$-multiplicative \morp.}}
{\blue{Since $A_1$ is  a RFD algebra, one can  choose a homomorphism $h_0$ from $A_1$ to a finite-dimensional \CA\, $F_{0}$ which is non-zero on $\mathcal F_{0}.$}}

{{ Let $\dt_0<\min\{\eta_0/2,1/2\}$. }}

Since $A$ is assumed to have the  {{property}} (SP),
by Lemma 2.1 of \cite{Niu}, there is a non-zero homomorphism $h': F_0\to A$ {\blue {with $h'(1_{F_0})=e_0\in A$}} such that
\begin{enumerate}
\item[(1)] $||e_0x-xe_0||<\delta_0/256$\ \ and
\item[(2)] $|| h'\circ h_0(x)-e_0xe_0 ||<\delta_0/256$
for all $x\in {{{\cal F}_0}}$.
\end{enumerate}
Since $F_0$ has finite dimension, it follows from Arveson's Extension Theorem that the homomorphism $h_0: A_1 \to F_0$ can be extended to a
\morp\,  from $A$ to $F_0;$  let us still denote {{this}} by $h_0$.

   Put $H=h'\circ h_0: A\to A.$  Note that $e_0=H(1)$. Since the hereditary $C^*$-subalgebra $(1-e_0)A(1-e_0)$ is in the class $\mathcal B_0$ again {{(see Theorem \ref{B1hered})}}, there is a projection $q_1'\leq 1-e_0$ and a $C^*$-subalgebra $S_1'\in \mathcal C_0$ {\blue{(of $(1-e_0)A(1-e_0)$)}} with $1_{S_1'}=q_1'$ such that
\begin{enumerate}
\item[(3)] $|| q'_1x-xq'_1 ||<\delta_0/256$ for any $x\in (1-e_0)
{{{\cal F}_0}}(1-e_0),$
\item[(4)] $\textrm{dist}(q_1'xq'_1, S_1')<\delta_0/256$ for any $x\in (1-e_0)
    {{{\cal F}_0}}(1-e_0)$, and
\item[(5)] ${\tau(1-e_0-q_1')} < 1/16$ for any tracial state $\tau$ on $A$.
\end{enumerate}
Put $q_1=q_1'+e_0$ and $S_1=S'_1\oplus h'(F_0)$. One has
\begin{enumerate}
\item[(6)] $|| q_1x-xq_1 ||<\delta_0/{64}$ for any $x\in
{{{\cal F}_0}},$
\item[(7)] $\textrm{dist}(q_1xq_1, S_1)<\delta_0/64$ for any $x\in\ 
    {{{\cal F}_0}}$, and
\item[(8)] $\tau(1-q_1)=\tau(1-q'_1-e_0)  <{ 1/16}$ for any tracial state $\tau$ on $A$.
\end{enumerate}

{Let ${\bar{\cal F}}_0\subset S_1$ be a finite subset such that ${\rm dist}(q_1yq_1, {\bar{\cal F}}_0)<\dt_0/{{64}}$
for all $y\in 
{{{\cal F}_0}}.$ }
Let $\mathcal G_1$ be a finite generating subset of $S_1$ {{(see \ref{DfC1})}} which is in the unit ball.

Since $S_1$ is amenable, by {\blue{Theorem 2.3.13 of \cite{Lnbok}}}, there is a \morp\,  $L_0': q_1Aq_1\to S_1$ such that
\beq\label{June13-1}
\|L_0'(s)-s\|<\dt_0/256\tforal s\in   {\cal G}_1\cup {\bar{\cal F}}_0.
\eneq
Set $L_0(a)=L_0'(q_1aq_1)$ for any $a\in A$. Then $L_0$ is a completely positive contraction from $A$ to $S_1$
such that
\beq\label{June13-2}
\|L_0(s)-s\|<\dt_0/128 \rforal s\in {\cal G}_1\cup{\bar{\cal F}}_0.
\eneq
{{Consequently, {{by the}}  definition of ${\bar{\cal F}}_0$ we have
\beq\label{May-13-2019}
\|L_0(s)-s\|<\frac{\dt_0}{64}+\frac{\dt_0}{128}+\frac{\dt_0}{64}<\frac{\dt_0}{16} \rforal s\in {\cal G}_1\cup (q_1{\cal F}_0q_1).
\eneq}}
{{It follows}}
that $L_0$ is ${{{\cal F}_0\cup {\cal G}_1}}$-$\delta_0/16$-multiplicative.
{{This ends  the initial stage of the proof.}} {\blue{Let us {{now}} prove the following claim.}}

{\blue{\bf Claim:}} {\blue{There exist}} {{an increasing sequence}} of finite subsets
$\{{\cal F}_n\}_{{n=0}}^{{\infty}}$ {{of the unit ball of $A$, starting with ${\cal F}_0$ above,}}
with {\blue{${\cal F}_n\supset \{x_1,x_2,\cdots, x_n\}$}}
{{(consequently, $\bigcup_{k=0}^{\blue{\infty}}{\cal F}_k$ is dense in the unit ball of $A$)}},
{\blue{two}} decreasing sequence{\blue{s}} of positive numbers {\blue{$\{\eta_n\}_{{n=0}}^{{\infty}}$}}
 {{and}} $\{\delta_n\}_{{n=0}}^{{\infty}}$ {{starting with $\eta_0$, $\dt_0$ above,}}
 finite generating subsets ${\cal G}_n\subset S_n$ {{(in the unit ball of $S_n$),}} {{starting with ${\cal G}_1$}}
(see the end of \ref{DfC1}),  a sequence of homomorphisms $h_{\blue{n-1}}:S_{\blue{n-1}}\to S_{\blue{n}}$ {\blue{(for $n\geq 2$),}}
and a sequence of  {{${\cal F}_{n-1}\cup {\cal G}_n$-$\dt_{n-1}/2$}}-multiplicative \morp s $L_{\blue{n-1}}: A\to S_{\blue{n}}$
{{(for $n\ge 1$)}} such that:
\beq
&& || q_nx-xq_n ||<\delta_{{n-1}}/{\blue{64}} \rforal x\in\mathcal F_{{n-1}};\label{April28-2019}\\
&&{\rm dist}(q_n{\blue{x}}q_n,S_n)<\delta_{\blue{n-1}}/{\blue{64}},\,
\,{\blue{\rforal x\in\mathcal F_{{n-1}};}}\label{April28-2019-2}\\
&& \tau(1-q_n)<1/2^{n+1}\rforal  \tau\in T(A);\label{May16-2019-4}\\
 &&{{ \mathcal G_{n}}}\subset\mathcal F_{\blue{n}},~{\blue{\mbox{for}~ n\geq 1}}~
 ({\rm recall}\,\, \mathcal G_{n}\,\, \mbox{\blue{generates}}
 \,\, S_{n});\label{May12-2019}\\
&& || L_{\blue{n-1}}(a)-h_{\blue{n-1}}(a)||<{{\dt_{n-2}}}/{\blue{64}}\rforal {\blue{n\geq 2,}}~
      a\in L_{\blue{n-2}}({{\mathcal F_{n-2}}}) \cup \mathcal G_{\blue{n-1}};\label{April30-2019}\\
   && {{L_{n-1}(a)=L_{n-1}(q_{n}aq_{n})  \rforal a\in A;}}\label{May16-2019-1}\\ &&\|L_{\blue{n-1}}(a)-a\|<\dt_{{n-1}}/{\blue{16}} \rforal a\in {\cal G}_{\blue{n}}{{\cup q_n{\cal F}_{n-1}q_n;}}\label{April28-2019-1} \\
   && {{\dt_n<\eta_{n}/2,~~\dt_n<\dt_{n-1}/2}};\label{April30-2019-2}  \\
   && {\blue{\|\big((1-q_n)x(1-q_n)+ L_{n-1}(x)\big)-x\|<\dt_{{n-1}}/8, \rforal x\in {\cal F}_{{n-1}}}};\label{April30-2019-3}\\
   &&{{[L_{n-1}]|_{(\iota_{n-1})_{*0}(K_0(S_{n-1}))}~~\mbox{is well defined};\label{May16-2019-2} }}
      \eneq
{{if $\tilde{L}: A\to A$   is an ${\cal F}_{n-1}$-$2\eta_{n-1}$ multiplicative  \morp\,
 with $\|\tilde{L}(x)-L_{n-1}(x)\|<2 \eta_{n-1}$ for all $x\in {\cal F}_{n-1}$, then}}
 \beq\label{May16-2019-3}
{{[\iota_{n}\circ L_{n-1}]|_{(\iota_{n-1})_{*0}(K_0(S_{n-1}))}=[\tilde{L}]|_{(\iota_{n-1})_{*0}
(K_0(S_{n-1}))};}}
\eneq
{\blue{if $L: S_n\to B$ ($n\geq 1$, $B$ is any $C^*$-algebra,) is a ${\cal G} _n$-$2\eta_n$
-multiplicative \morp,  there exists
a \hm\, $h: S_{n}\to B$ such that}}
\beq\label{May1-2019}
\|L(a)-h(a)\|<\dt_{{n-1}}/64\rforal a\in {\cal G}_{n} {\blue{\cup L_{n-1}({\cal F}_{n-1})}};
\eneq
{\blue{if $L': A\to B$  is {{an}} ${{\cal F}}_n$-$2\eta_n$
-multiplicative \morp,
then
\beq\label{April29-2019}
[L']|_{(\iota_n)_{*0}(K_0(S_n))}~~\mbox{is well defined};
\eneq
and, if $L'': A\to B$   is another ${{\cal F}}_n$-$2\eta_n$ multiplicative  \morp\,
 with $\|L'(x)-L''(x)\|<2 \eta_n$ for all $x\in {{\cal F}}_n$, then
 \beq\label{April29-2019-1}
[L']|_{(\iota_n)_{*0}(K_0(S_n))}=[L'']|_{(\iota_n)_{*0}(K_0(S_n))}.
\eneq }}

{{We have constructed ${\cal F}_0$ $\eta_0$, $\dt_0$, $q_1$, $S_1$ (with $1_{S_1}=q_1$),  {{a}} finite generating set ${\cal G}_1\subset S_1$, and
 $L_0: A\to S_1$ in the initial stage.
 Note that (\ref{April28-2019})--\eqref{May16-2019-4}, \eqref{May16-2019-1}, \eqref{April28-2019-1}, and
 \eqref{April30-2019-3}
 for $n=1$
 follow from our construction (see (6), (7), (8), and (\ref{May-13-2019})).}}
 To begin the induction, let us
first define  $\eta_1, \dt_1$, and {{the}} subset ${\cal F}_1$ ($L_1$ and $h_1$ should be defined after we define $q_2$ and $S_2$, so it is not necessary to define them here).
{{Since $S_1$ is {{semiprojective}} (see the end of \ref{DfC1}), there exists a positive number $\eta_1<\eta_0/2$ satisfying the following condition: for any
${\cal G}_{1}$-$2\eta_{1}$-multiplicative \morp\, $L$ from $S_{1}$ to another \CA\, $B$ (recall that ${\cal G}_{1}$ is a generating
{{subset}} of $S_{1}$),  there exists
a \hm\, $h: S_{1}\to B$ such that
\beq
\|L(a)-h(a)\|<\dt_{0}/64\rforal a\in {\cal G}_{1}\cup L_{0}({\cal F}_{0}).
\eneq}}
{{Moreover,  since $K_0(S_1)$ is finitely generated,
we may assume
that there is a finite subset ${\cal F}_1'$ such that $[L']|_{(\iota_1)_{*0}(K_0(S_1))}$
is well defined, if $L'$ is an ${\cal F}_1'$-$2\eta_1$-multiplicative \morp\,
from $A$. Furthermore,
\beq
[L']|_{(\iota_1)_{*0}(K_0(S_1))}=[L'']|_{(\iota_1)_{*0}(K_0(S_1))}
\eneq
 if
$\|L'(x)-L''(x)\|<2\eta_1$ for all $x\in {\cal F}_1'$ and $L''$ is also
an ${\cal F}_1'$-$2\eta_1$-multiplicative \morp\, from $A.$}}
{{Let $\dt_1<\min\{\et_1/2, \dt_0/2\}$, and ${\cal F}_1={\cal F}_0\cup {\cal F}'_1\cup {\cal G}_1\cup L_0({\cal F}_0)\cup\{x_0,x_1\}$. Then (\ref{May12-2019}), (\ref{April30-2019-2}), (\ref{May1-2019}), (\ref{April29-2019}), and (\ref{April29-2019-1}) hold for $n=1$. (Note that (\ref{April30-2019}),  (\ref{May16-2019-2}), and (\ref{May16-2019-3}) only make sense for $n\geq 2$.)}}

{\blue{For $k\geq 2$, suppose {{that}} we have already constructed
${\cal F}_{k-1}\subset A$ (and all ${\cal F}_{i}$ with $1\leq i< k-1$), $S_{k-1}\in {\cal C}_0$ (and all $S_i$ with $1\leq i< k-1$),
{{a}}
  finite generating subset ${\cal G}_{k-1}$ of $S_{k-1}$ (and all ${\cal G}_i$ with $1\leq i< k-1$),
$\eta_{k-1}$, $\dt_{k-1}$ (and all $\eta_i, \dt_i$ with $1\leq i< k-1$),
and {{a}} homomorphism $h_{k-2}: S_{k-2} \to S_{k-1}$ if $k\geq 3$ (and all $h_i$ with $1\leq i< k-2$),  {{a}} completely positive linear contraction $L_{k-2}:A \to S_{k-1}$ (and all $L_i$ with $0\leq i< k-2$) as described  in (\ref{April28-2019})--(\ref{April29-2019-1}) with $k-1$ in place of $n$. }}

Since  $A\in {\cal B}_0,$ for {{the $k$  above,}} there {\blue{exist}} a projection $q_k\in A$ and a \SCA\, $S_k\in {\cal C}_0$
{\blue{with $1_{S_k}=q_k$}} such that
\beq\label{201-18-s-1}
\|q_kx-xq_k\|<\dt_{{k-1}}/64\rforal x\in {\cal F}_{{k-1}},\\\label{201-18-s-2}
{\rm dist}(q_kxq_k, S_k)<\dt_{{k-1}}/64\rforal x\in {\cal F}_{{k-1}},\andeqn\\\label{201-18-s-3}
\tau(1-q_k)<1/2^{k+1}\rforal \tau\in T(A).
\eneq
{{Hence (\ref{April28-2019}),(\ref{April28-2019-2}), and (\ref{May16-2019-4}) hold for $n=k$.}}

Let ${\cal G}_k\subset S_k$ be a finite generating {{subset}} {\blue{in the unit ball of ${\cal G}_k$}}.

{{Define}} $\Lambda_k: A\to q_kAq_k$ by $\Lambda_k(a)=q_kaq_k$ for all $a\in A.$
Then, by \eqref{201-18-s-1}, $\Lambda_k$ is ${\cal F}_{{k-1}}$-$\dt_{{k-1}}/{\blue{64}}$-multiplicative.
Let ${\bar{\cal F}}_{{k-1}}\subset S_k$ be a finite subset such that
${\rm dist}(q_kxq_k, {\bar{\cal F}_{{k-1}}})<\dt_k/{\blue{64}}$ for all $x\in {\cal F}_{{k-1}}.$

Since $S_k$ is amenable, by Theorem 2.3.12 of \cite{Lnbok}, there exists a unital \cp\, $L_k': q_kAq_k\to S_k$
such that
\beq\label{201-18-s-4}
\|L_k'(a)-a\|<\dt_{{k-1}}/64\rforal a\in {\cal G}_k\cup {\bar{\cal F}_{{k-1}}}.
\eneq

{\blue{ Let $L_{k-1}: A\to S_k$ be defined by $L_{k-1}=L'_k\circ \Lambda_k$.}} {{Evidently, (\ref{May16-2019-1}) holds for $k$ in place of $n$.}} Also,
\beq\label{April30-2019-1}
\|L_{k-1}(a) -a\|<\dt_{{k-1}}/64 ~~\mbox{for all}~~ a\in {\cal G}_k\cup {\bar{\cal F}_{{k-1}}}.
\eneq
{{Since ${\rm dist}(a, {\bar{\cal F}_{k-1}})<\dt_{k-1}/64$ for each $a\in q_k{\cal F}_{k-1}q_k$,  we have }}
\beq
{{\|L_{k-1}(a)-a\|<\dt_{k-1}/{16} \rforal a\in {\cal G}_{k}\cup {{q_k{\cal F}_{k-1}q_k,}}}}\label{May16-2019}
\eneq
which is (\ref{April28-2019-1}) for $k$ in place of $n$. {{Combining (\ref{201-18-s-1}) and (\ref{May16-2019}), we {{conclude}} that
$$\|\big((1-q_k)x(1-q_k)+ L_{k-1}(x)\big)-x\|<\dt_{k-1}/8, \rforal x\in {\cal F}_{k-1},$$
which is (\ref{April30-2019-3}) for $n=k$, and, further, $L_{k-1}$ is ${\cal F}_{k-1}$-$\dt_{k-1}$-multiplicative, and
{{hence}}
${\cal F}_{k-1}$-$\et_{k-1}$-multiplicative as $\dt_{k-1}<\et_{k-1}$. By (\ref{April29-2019}), and (\ref{April29-2019-1}) for $n=k-1$,
(\ref{May16-2019-2}) and (\ref{May16-2019-3}) hold for $n=k$. By (\ref{May12-2019}) for $n=k-1$, we have ${\cal G}_{k-1}\subset {\cal F}_{k-1}$. Then by (\ref{May1-2019}) for $n=k-1$, there is a homomorphism $h_{k-1}: S_{k-1}\to S_k$ such that
\beq
||L_{k-1}(a)-h_{k-1}(a)||<\dt_{k-2}/64\rforal ~
      a\in L_{k-2}(\mathcal F_{k-2}) \cup \mathcal G_{k-1},
      \eneq
which is (\ref{April30-2019}) for $n=k$.}}



Since $S_k$ is {{semiprojective}} (see the end of \ref{DfC1}), there exists {\blue{a positive number}} $\eta_k{\blue{< \eta_{k-1}/2<}} \eta_0/2$ {{satisfying}} the following condition: for any
${\cal G}_k$-$2\eta_k$-multiplicative \morp\, $L$ from $S_k$ to another \CA\, $B$ (recall that ${\cal G}_k$ is a generating
{{subset}} of $S_k$),  there exists
a \hm\, $h: S_k\to B$ such that
\beq\label{201-18-s-5}
\|L(a)-h(a)\|<\dt_{{k-1}}/64\rforal a\in {\cal G}_k\cup {{L_{k-1}({\cal F}_{k-1}),}}
\eneq
{{which is (\ref{May1-2019}) for $n=k$.}}

Moreover,  since $K_0(S_k)$ is finitely generated,
we may assume
that there is a finite subset ${\cal F}_k'$ such that $[L']|_{(\iota_k)_{*0}(K_0(S_k))}$
is well defined, if $L'$ is an ${\cal F}_k'$-$2\eta_k$-multiplicative \morp\,
from $A${\blue{. That is, (\ref{April29-2019}) holds for $k$ in place of $n$, provided ${\cal F}_k'\subset {\cal F}_k$.  {{Furthermore,}}}}
\beq\label{201-18-nnn1}
[L']|_{(\iota_k)_{*0}(K_0(S_k))}=[L'']|_{(\iota_k)_{*0}(K_0(S_k))}
\eneq
 if
$\|L'(x)-L''(x)\|<2\eta_k$ for all $x\in {\cal F}_k'$ and $L''$ is also
an ${\cal F}_k'$-$2\eta_k$-multiplicative \morp\, from $A.$
{{In other words,}} {\blue{(\ref{April29-2019-1}) holds for $k$ in place of $n$, provided ${\cal F}_k'\subset {\cal F}_k$.}} {{Let ${\cal F}_k={\cal F}_{k-1}\cup {\cal F}'_k\cup {\cal G}_k\cup L_{k-1}({\cal F}_{k-1})\cup\{x_0,x_1,\cdots,x_k\}$. Then (\ref{May12-2019}), (\ref{April29-2019}), and (\ref{April29-2019-1}) hold for $n=k$. Let  $\dt_k<\min\{\et_k/2, \dt_{k-1}/2\}$. Then (\ref{April30-2019-2}) holds for $n=k$.}}  {{This completes the proof of the  claim.}}

\end{NN}

\begin{NN}\label{June13-constr}
{\rm
Let $\Psi_n: A\to(1-q_{{n+1}})A(1-q_{{n+1}})$ denote the {{cutting-down}} map sending $a$ to $(1-q_{{n+1}})a(1-q_{{n+1}})$, and let $J_n:A\to A$ denote the map sending $a$ to $\Psi_n(a)\oplus L_n(a)$. {{By}} (\ref{April30-2019-3}),  {{the maps}}
$\Psi_n$ and $J_n$ are $\mathcal F_{{n}}$-$\delta_{{n}}/2$-multiplicative. Set $J_{m,n}=J_{n-1}\circ\cdots\circ J_m$ and
$h_{m,n}=h_{n-1}\circ\cdots\circ h_{m}: S_{m}\to S_{n}$. {\blue{Again by (\ref{April30-2019-3}),}}
$J_{m,n}$ is $\mathcal F_{{m}}$-$\delta_{{m}}$-multiplicative
and $\|J_{m,n}(x)-x\|<\dt_{{m}}$ for all $x\in {{{\cal F}_m}}.$ We {{shall}} also use $L_n,
\Psi_n, J_n, J_{m,n}, h_m, $ and $h_{m,n}$ for {{the extensions of these maps to}} a matrix algebra over $A$.}

{\blue{By (\ref{April30-2019-2}), (\ref{April29-2019}), and  (\ref{April29-2019-1}),}}
$J_{m,n}$ {{induces a well-defined map}}
$[J_{m,n}|_{{\cal P}_0\cup (\iota_k)_{*0}(K_0(S_m))}],$ which agrees with the identity map.

\end{NN}

Using the same argument as that of Lemma 2.7 of \cite{Lnduke}, one has the following lemma.

\begin{lem}[Lemma 2.7 of \cite{Lnduke}]\label{traces}
Let $\mathcal P\subset \text{M}_{ m}(A)$ be a finite set of projections. Assume that $\mathcal F_{\blue{0}}$ {\blue{(in \ref{201})}} is sufficiently large and $\delta_0$  {\blue{(in \ref{201})}}  is
sufficiently small  that $[L_{n}\circ J_{1,n}]|_{\mathcal P}$ and $[L_{n}\circ J_{1,n}]|_{G_0}$ are well defined, where $G_0$ is the subgroup {{of $K_0(A)$}} generated by $\mathcal P$. Then
$$\lim_{n\to\infty}\sup_{\tau\in T(A)}|\tau([\iota_{n+1}\circ L_{n}\circ J_{1,n}]([p]))-\tau([p])|=0$$
for any $p\in\mathcal P.$

Furthermore, we have
\beq
&&|\tau(h_{k, k+n+1}\circ [L_{k-1}]([p]))-\tau(h_{k, k+n}\circ [L_{k-1}]([p])|<(1/2)^{n+k}\tand\\
&& \tau(h_{k,k+n}\circ [L_{k-1}]([p]))\ge (1-\sum_{i=1}^n 1/2^{i+k})\tau([L_{k-1}]([p]))>0
 \eneq
 for all $p\in {\cal P}$ and  for all $\tau\in T(A),$  {{and}}
\beq\label{203-18-s1}
 \tau(h_{k,k+n}([q]))\ge (1-\sum_{i=1}^n 1/2^{i+k})\tau([q])>0 \tforal \tau\in T(A)
\eneq
{{and}} for all $q\in K_0(S_k)_{+}\setminus \{0\}.$
\end{lem}

\begin{rem}\label{Kstate}
Since $A$ is stably finite and assumed to be amenable, {{and}} therefore exact, any  state of $K_0(A)$ is the restriction of a tracial state of
$A$ (\cite{Blatrace} and \cite{Haagtrace}).  Thus, the lemma above still holds if one replaces the trace $\tau$ by
{{any}}
state $\tau_0$ on $K_0(A)$.
\end{rem}

\begin{NN}\label{Constr-1508}
Fix a finite subset ${\cal P}$ of projections of $M_r(A)$ (for some $r\ge 1$)
and an integer $N\ge 1$ such that $[L_{N+i}]|_{\cal P},$
$[J_{N+i}]|_{\cal P}$ and $[\Psi_{N+i}]|_{\cal P}$ are all well defined.
Keep {{the}} notation {{of}}  \ref{June13-constr}. Then, on ${\cal P},$
\beq\nonumber
[L_{N+1}\circ J_N]&=&[L_{N+1}\circ L_N]\oplus [L_{N+1}\circ \Psi_N]\\
&=& [h_{N+1}\circ L_N]\oplus [L_{N+1}\circ  \Psi_N], \andeqn
\eneq
\beq\nonumber
\hspace{-0.6in}[L_{N+2}\circ J_{N, N+2}]&=&[L_{N+2}\circ L_{N+1}\circ J_N]\oplus [L_{N+2}\circ \Psi_{N+1}\circ J_N]\\
&=&[L_{N+2}\circ L_{N+1}\circ L_N]\oplus [L_{N+2}\circ L_{N+1}\circ \Psi_N]\\
&&\oplus [L_{N+2}\circ  \Psi_{N+1}\circ J_N]\\\hspace{-0.1in}&=&[h_{N+1, N+3}]\circ [L_N]\oplus [L_{N+2}\circ L_{N+1}\circ \Psi_N]\oplus
[L_{N+2}\circ  \Psi_{N+1}\circ J_N].\\
\eneq
Moreover,  on ${\cal P},$
\beq\nonumber
[L_{N+n}\circ J_{N, {{N+n}}}] &=& [h_{N+1, N+n+1}]\circ [L_N]\oplus [L_{N+n}\circ \Psi_{N+n-1}\circ J_{N, N+n-1}]\\
&&\oplus [L_{N+n}\circ L_{N+n-1}\circ \Psi_{N+n-2}\circ J_{N, N+n-2}]\\
&&\oplus  [L_{N+n}\circ L_{N+n-1}\circ L_{N+n-2}\circ \Psi_{N+n-3}\circ J_{N, N+n-3}]\\
&&\oplus \cdots { \oplus[L_{N+n}\circ L_{N+n-1}\circ\cdots \circ L_{N+2}\circ \Psi_{N+1}\circ J_N}]\\
&& \oplus [L_{N+n}\circ L_{N+n-1}\circ\cdots \circ L_{N+1}\circ \Psi_N].
\eneq

Set $\psi^N_{N}=L_N,$ $\psi^N_{N+1}=L_{N+1}\circ \Psi_N,$
$\psi^N_{N+2}=L_{N+2}\circ \Psi_{N+1},..., {{\rm and}}\, \psi^N_{N+n}=L_{N+n}\circ \Psi_{N+n-1},$
$n=1,2,....$

{\blue{(Note that $\psi^N_{N+i}=\psi^{N+1}_{N+i}=\cdots=\psi^{N+i-1}_{N+i}$.  We insist {{on}} the notation $\psi^N_{N+i}$
{{in order}} to emphasize {{that}} our estimation {{starts}} with a fixed index $N$, and {{goes}} on to $N+1$, $N+2$,
{{and so on, and}} we only change the subscripts and keep the superscript {{the same}}
to emphasize our beginning index $N$ (see Corollary \ref{rho} below).)}}

\end{NN}

\begin{NN}\label{Constr-Sn} (a)
For each $S_n$, since the abelian group $K_0(S_n)$ is finitely generated and torsion free, there is a set of free generators $\{e^{n}_1, e^{n}_2, ..., e^{n}_{l_n}\}\subset K_0(S_n)$. By Theorem \ref{FG-Ratn}, the positive cone of the
{{ordered group}} $K_0(S_n)$ is finitely generated; denote {{the (semigroup)}}
generators {{(exactly the the minimal non-zero positive elements)}}
by $\{s^n_1, s^n_2, ..., s^n_{r_n}\}\subset K_0(S_n)_+{{\setminus\{0\}}}$. Then there is an
$r_n\times l_n$ integer-valued matrix $R'_n$ such that
$${{\vec{s}_n}}=R'_n\vec{e}_n,$$
where ${{\vec{s}_n}}=(s^n_1, s^n_2, ..., s^n_{r_n})^{\mathrm T}$ and $\vec{e}_n=(e^{n}_1, e^{n}_2, ..., e^{n}_{l_n})^{\mathrm T}$. In particular, for any ordered group $H$, and any elements $h_1, h_2, ..., h_{l_n}\in H$, the map $e^n_i\mapsto h_i$, $i=1, ..., l_n,$ induces an abelian-group homomorphism $\phi: K_0(S_n)\, {{\to}}\, H$, and the map $\phi$ is positive (or strictly positive) if and only if
$$R'_n\vec{h}\in H^{r_n}_+ \quad \textrm{(or $R_n'\vec{h}\in {(H_+\setminus \{0\})^{r_n}}$),}$$
where $\vec{h}=(h_1, h_2, ..., h_{l_n})^{{T}} \in H^{l_n}$.
Moreover, for each $e^n_i$, write it as $e^n_i=(e^n_i)_+-(e^n_i)_-$ {{with}} $(e^n_i)_+, (e^n_i)_-\in K_0(S_n)_+$ and fix this decomposition.
Define  a $r_n\times 2l_n$ matrix
\begin{displaymath}
R_n=
R'_n
\left(
\begin{array}{ccccccc}
1 & -1 & 0 & 0 &\cdots & 0 & 0 \\
0 & 0 & 1 & -1 & \cdots & 0 & 0 \\
\vdots & \vdots & \vdots & \vdots & \ddots & \vdots & \vdots \\
0 & 0 & 0 & 0 & \cdots & 1 & -1
\end{array}
\right)_{\blue{l_n\times l_{2n}}}.
\end{displaymath}
Then one has $${{\vec{s}_n}}=R_n\vec{e}_{n, \pm},$$ where $\vec{e}_{n, \pm}=((e^{n}_{1})_+, (e^{n}_{1})_-, ..., (e^{n}_{l_n})_+, (e^{n}_{l_n})_-)^{\mathrm T}$.  {{Thus,}}
for any ordered group $H$, and any elements $h_{1, +}, h_{1, -}, ..., h_{l_n, +}, h_{l_n, -}\in H$, the map $e^n_i\mapsto (h_{i, +}-h_{i, -})$, $i=1, ..., l_n,$ induces a positive (or strictly positive) homomorphism if and only if
$$R_n\vec{h}_{{\pm}}\in H^{r_n}_+\quad \textrm{ (or $R_n\vec{h}_{{\pm}}\in (H_+\setminus\{0\})^{r_n}$) },$$
where $\vec{h}_{\pm}=(h_{1, +}, h_{1, -}, ..., h_{l_n, +}, h_{l_n, -})^{{T}}   \in H^{l_n}$.

{ {(b) Let $A$ be as above. Let $B$ be an inductive limit $C^*$-algebra {{as}} in Theorem \ref{RangT} such that
 $$(K_0(A),{K_0(A)_+},[1_A],K_1(A))\cong(K_0(B),{K_0(B)_+},[1_B],K_1(B)).$$
Let $\alpha\in KL(A, B)$ be an element which implements the isomorphism above.

Let $e^n_i, e^n_{i,\pm} \in K_0(S_n)$ be as  above, $i=1,2,..., l_n.$
Let $s(0)=0, s(n)=\sum_{i=1}^n2l_i$.
Put
$$\alpha(\imath_n\circ  h_{j, n}(e^j_{i,+}))=g^{(n)}_{s(j-1)+2i-1}, ~~\alpha(\imath_n\circ h_{j, n}(e^j_{i,-}))=g^{(n)}_{s(j-1)+2i},~~\mbox{and}~~g^{(n)}_l=0~\mbox{if}~ l>s(n).$$
Let $a^{(n)}_j=\rho_B(g^{(n)}_j) \in \Aff({{T(B))_+}}$ for $j=1,2,{ ...}.$ Then{{,}} by Lemma \ref{traces},  $\lim_{n\to\infty} a^n_j=a_j >0$ uniformly on ${{T(B)}}$.

For $j\in \{1,2,{ ...}, n\}$,
let
$(s_1^j, s_2^j, { ...}, s_{r_j}^j) \in K_0(S_j)_{ +}\setminus \{0\}$
be
the generators of the positive cone {{$K_0(S_j)_+,$}}
and {{let $R_j$ be  the $r_j\times 2 l_j$}}  {{matrix}} as {{in}} part (a). Then
\beq
&&R_j (g^{(n)}_{s(j-1)+1}, g^{(n)}_{s(j-1)+2}, { ...}, g^{(n)}_{s(j)})^T\\
&=&(\alpha(\imath_n\circ h_{j,n}(s_1^j)), \alpha(\imath_n\circ h_{j,n}(s_2^j)), { ...}, \alpha(\imath_n\circ h_{j,n}(s_{r_j}^j)))^T
\andeqn\\
&&R_j (a_{s(j-1)+1}, a_{s(j-1)+2}, { ...}, a_{s(j)})^T\\
&=&\lim_{n\to \infty} (\rho_B(\alpha(\imath_n\circ h_{j,n}(s_1^j))), \rho_B(\alpha(\imath_n\circ h_{j,n}(s_2^j))), { ...}, \rho_B(\alpha(\imath_n\circ h_{j,n}(s_{r_j}^j))))^T.
\eneq
Note that  for each $\tau\in {{T(B)}}$, $\tau\circ \rho_B\circ \alpha$ {{is}} a state on $K_0(A),$ which
{{(by \cite{Blatrace} and \cite{Haagtrace})
can}} be extended to a trace on $A$. From Lemma \ref{traces} (\eqref{203-18-s1}), {{for all $1\le j\le n,$ and  for all $\tau\in T(B),$}}
$$
\tau (\rho_B(\alpha(\imath_n\circ h_{j,n}(s_i^j))))>(1-\sum_{k=1}^{\infty} 1/2^{j+k} )\tau(\rho_B{\blue{(\af(\imath_j (s_i^j)))}})>0,{{\rforal}} i \in\{1,2,{ ...}, r_j\}.$$
Hence each entry of $R_j (a_{s(j-1)+1}, a_{s(j-1)+2}, { ...}, a_{s(j)})^T$  is a strictly positive element of $\Aff(T(B))^{{++}}$.
Let ${\bar R}_n =\diag (R_1, R_2,..., R_n)$. Then
$$0\ll {\bar R}_n (a_1, a_2, ..., a_{s(n)})^T\in (\Aff({{T(B)}}))^{\sum_{i=1}^nr_k},$$
i.e.,  each coordinate  is strictly positive  on ${{T(B)}}$.

Furthermore, ${\bar R}_n(g^{(n)}_1, g^{(n)}_2, ..., g^{(n)}_{s(n)}) \in (K_0(B)_+\setminus \{0\})^{\sum_{i=1}^nr_k}.$ In particular, for each positive integer $N_0 <n$, we also have
$$\diag (R_{N_0+1},R_{N_0+1}, ..., R_n) (g^{(n)}_{s(N_0)+1}, g^{(n)}_{s(N_0)+2}, ..., g^{(n)}_{s(n)}) \in (K_0(B)_+\setminus \{0\})^{\sum_{i=N_0+1}^nr_k}.$$

}}

{ {(c)}} Since $\{e^{n}_1, e^{n}_2, ..., e^{n}_{l_n}\}$ is a set of {\blue{$\Z$-linearly independent}} generators of
{{the}} {\blue{free abelian group}} $K_0(S_n)$, for any projection $p$ in a matrix algebra {{over}} $S_n$, there {\blue{is a unique $l_n$-tuple of}} {{\em integers}} $m^n_1({{[p]}}), ..., m^n_{l_n}({{[p]}})$ such that
{{$[p]=\sum_{j=1}^{l_n}m_j^n([p])e_j^n$ and,}} for any homomorphism $\tau: K_0(S_n)\to\mathbb R$, one has
$$
\tau([p])=\langle\vec{m}_n({{[p]}}), \tau(\vec{e}_n)\rangle=\sum_{i=1}^{l_n} m_i^{{n}}({{[p]}})\tau(e^n_i)=\sum_{i=1}^{l_n} m^n_i({{[p]}})\tau((e^n_i)_+)-m^n_i({{[p]}})\tau((e^n_i)_-),$$
where $\vec{m}_n{{([p])}}=(m^n_1({{[p]}}), ..., m^n_{l_n}({{[p]}}))^{{T}}$ and $\vec{e_n}=(e^{n}_1, e^{n}_2, ..., e^{n}_{l_n}).$

For each $p\in M_m(A),$ for some integer $m\ge1,$ denote by $[{{\psi_{k+j}^k}}(p)]$ the element
of $K_0(S_{{k+j+1}})$ associated with ${{\psi^k_{k+j}}}(p).$  Let $\imath_n: S_n\to A$ be the embedding.
{{Consider the map}}
$${\overline{(\imath_n)_{*0}}}: \vec{e_n}\mapsto (((\imath_n)_{*0}(e_1^n), (\imath_n)_{*0}(e_2^n),...,(\imath_n)_{*0}(e_{l_n}^n)).$$

\end{NN}
Then, by Lemma \ref{traces} and Remark \ref{Kstate}, one has the following lemma.

\begin{lem}\label{rho-uniform}
With  {{notation}}
as
above, for any $p\in{{\mathcal P}},$  for each fixed $k$, one has that
\beq
\tau(p)=
\lim_{n\to\infty}\sum_{j=1}^n(\sum_{i=1}^{{{l_{k+j}}}} m_i^{{k+j}}([{{\psi^k_{k+j-1}}}(p)])\tau((\imath_{k+n}\circ h_{k+j,k+n})_{*0}(e_i^{{k+j}})_{{+}})\\
-m_i^{k+j}({{\psi^k_{k+j-1}}}(p)])\tau((\imath_{k+n}\circ h_{k+j,k+n})_{*0}(e_i^{{k+j}})_-))
\eneq
uniformly on $\mathrm{S}(K_0(A)).$  Moreover,  $\rho_A\circ (\imath_{{n+k}})_{*0}\circ h_{k+j,k+n}(e_{i,\pm}^{k+j})$
converge{{s}} {{uniformly}} to a strictly positive element {{of}} $\mathrm{Aff}(\mathrm{S}(K_0(A)))$ as $n\to\infty.$
\end{lem}

\begin{proof}
We first compute that, if $j>1,$ {{then}}
{{\beq\label{rho-comput-1508}
&&\sum_{i=1}^{l_{k+j}} m^{k+j}_i([\psi^k_{k+j-1}]([p]))\tau((\imath_{k+n}\circ h_{k+j,k+n})_{*0}(e^{k+j}_i))\\
&=& \tau([L_{k+n-1}\circ\cdots \circ L_{k+j+1}\circ L_{k+j}\circ L_{k+j-1}\circ \Psi_{k+j-2}]([p])),
\eneq}}
and, if $j=1,$ {{then}}
\beq\label{rho-compute-1508-1}
&&\sum_{i=1}^{l_{k+1}} m^{k+1}_i([{{\psi^k_{k}}}]([p]))\tau((\imath_{k+n}\circ h_{k+1,k+n})_{*0}(e^{k+1}_i))\\
&=& \tau([L_{{k+n-1}}\circ\cdots \circ L_{{k+1}}\circ L_k]([p])).
\eneq
Thus  (see \ref{Constr-1508}),
\beq\label{comput-1508-2}
&&\sum_{j=1}^n (\sum_{i=1}^{l_{k+j}} m^{k+j}_i([\psi^k_{{k+j-1}}(p)])\tau((\imath_n\circ h_{k+j,k+n})_{*0}(e^{k+j}_i)))\\
&=& \tau([L_{{k+n-1}}\circ\cdots \circ L_{{k+1}}\circ L_k]([p]))\\
&&+\sum_{j=2}^n\tau([L_{{k+n-1}}\circ\cdots \circ L_{k+j+1}\circ L_{k+j}{{\circ L_{k+j-1}}}\circ \Psi_{k+j-1}]([p]))\\
&=& \tau([L_{{k+n-1}}\circ\cdots \circ L_{{k+1}}\circ L_k]([p]))\\
&&+\sum_{j=2}^n\tau([L_{{k+n-1}}\circ\cdots \circ L_{k+j+1}\circ L_{k+j}{{\circ L_{k+j-1}}}\circ \Psi_{{k+j-2}}\circ J_{k,{{k+j-2}}}]([p]))\\
&=&\tau([L_{{k+n-1}}\circ J_{k, n+k-1}]([p])).
\eneq
{{The}} first part of the  {{conclusion}}
then follows from {{Lemma}} \ref{traces}. The second part also follows.

\end{proof}

One then has the following
\begin{cor}\label{rho}
{{With notation as above, in particular with}} $\mathcal P$ be a finite subset of projections in a matrix algebra over $A$,  let $G_0$ be the subgroup of $K_0(A)$ generated by $\mathcal P,$ {{ and let $k\ge 1$ be an integer}}.  Denote by $\tilde{\rho}: G_0\to{{\prod}}\mathbb Z$ the map defined  {{(see \ref{Constr-Sn} (c))}} by
\begin{eqnarray}
&&\hspace{-0.5in}[p] \mapsto \nonumber\\
&&\hspace{-0.5in} ({m}_1^{\blue{k+1}}({ g}_0), -m_1^{\blue{k+1}}({ g}_0), {m}_2^{\blue{k+1}}({ g}_0), -m_2^{\blue{k+1}}({ g}_0),\cdots m_{l_{k+1}}^{\blue{k+1}}({ g}_0), -m_{l_{k+1}}^{\blue{k+1}}({ g}_0), \nonumber\\
&&m^{k+{\blue{2}}}_1({ g}_1), -m^{k+{\blue{2}}}_1({ g}_1),m^{k+{\blue{2}}}_2({ g}_1), -m^{k+{\blue{2}}}_2({ g}_1), ...,m^{k+{\blue{2}}}_{l_{k+{\blue{2}}}}({ g}_1), -m^{{k+2}}_{l_{k+{\blue{2}}}}({ g}_1),\nonumber\\
&&\cdots\cdots \cdots \nonumber\\
\hspace{-0.2in}&&m^{k+{\blue{i+1}}}_1({ g}_2), -m^{k+{\blue{i+1}}}_1({ g}_{\blue{i}}),m^{k+{\blue{i+1}}}_2({ g}_{\blue{i}}), -m^{k+{\blue{i+1}}}_2({ g}_{\blue{i}}), ...,m^{k+{\blue{i+1}}}_{l_{k+{\blue{i+1}}}}({ g}_{\blue{i}}), -m^{k+{\blue{i+1}}}_{l_{k+{\blue{i+1}}}}({ g}_{\blue{i}}), \cdots),
\nonumber
\end{eqnarray}
where ${ g}_i=[\psi^{{{k}}}_{k+i}(p)],$ $i=0,1,2,....$
If {{$g\in {\cal G}_0$ and}} $\tilde{\rho}(g)=0$, then $\tau(g)=0$ for {{$\tau\in T(A).$}}
\end{cor}

By the definition of the {{maps}} $\tilde{\rho}$ and $H=h'\circ h_0: A\to F_0\to  A$, {\blue{where $h_0, h'$ are
{{as in}} \ref{201},}} using the same argument as that of Lemma 2.12 of \cite{Lnduke}, one has the following lemma.
\begin{lem}\label{ker}
Let $\mathcal P$ be a finite subset of projections in $M_k(A_1)\subset M_k(A)$. Then there {{are}} a finite subset
$\mathcal F_1\subset A_1$ and $\delta_0>0$ such that if the above construction starts with $\mathcal F_1$ and $\delta_0$, then
$$\ker\tilde{\rho}\subset\ker[h_0]\ \ \mbox{and}\ \ \ker\tilde{\rho}\subset\ker[H].$$
\end{lem}

The $K_0$-part of the existence theorem will  {{deal  with maps which} almost factor through the map $\tilde{\rho}$, and this lemma will help us to handle the elements of $K_0(A)$ which vanish under $\tilde\rho$. Moreover, to get  such  a $K_0$-homomorphism, one also needs to find a copy of the generating subset of the positive cone of $K_0(S)$ inside  the image of  ${{\tilde \rho}}$
as an ordered group
for a certain algebra $S\in \mathcal C_0$. In order to do so, one {needs} the following technical lemma, which is  {essentially} Lemma 3.4 of \cite{Lnduke}.

\begin{lem}[Lemma 3.4 of \cite{Lnduke}]\label{solveeq}
Part (a):
Let $S$ be a compact convex set, and let $\Aff(S)$ be the space of real continuous affine functions on $S$. Let $\mathbb{D}$ be a dense ordered subgroup of $\Aff(S)$ {{(with the strict pointwise order)}}, and let $G$
be an ordered group with the strict order determined by a surjective homomorphism $\rho:G\to\mathbb{D}$
{{(i.e., $g\ge 0,$ if and only if $g=0,$ or $\rho(g)>0).$}}  Let $\{x_{ij}\}_{1\leq i\leq r, 1\leq j<\infty}$ be an $r\times\infty$ matrix having rank $r$ and with $x_{ij}\in\mathbb Z$ for each $i,j$.  Let $g_j^{(n)}\in G$ {{be}} such that $\rho(g_j^{(n)})=a_j^{(n)}$, where $\{a_j^{(n)}\}$ is a sequence of positive elements {{of}}  $\mathbb{D}$ such that $a_j^{(n)}\to a_j \,(>0)$ uniformly on $S$ as $n\to\infty$.

{{S}}uppose that there is a sequence of integers $s(n)$ satisfying the following condition:

{Let $\widetilde{v_n}=(g_j^{(n)})_{s(n)\times1}$ be the {{initial}} part of $(g_j^{(n)})_{1\leq j<\infty}$
and let
$$\widetilde{y_n}=(x_{ij})_{r\times s(n)}\widetilde{v_n}.$$
{{Set}} $y_n=\rho^{(r)}(\widetilde{y_n})${\blue{, where $\rho^{(r)}=\diag(\underbrace{\rho,\rho,\cdots,\rho}_r)$}}.
{{Suppose that}} there exists $z=(z_j)_{r\times1}$ such that $y_n\to z$ on $S$ uniformly.}

With the condition above,  there exist $\delta>0${{,}}  a positive integer $K>0,$  {{and}} a positive integer $N$
{{with the following property}}:

{{ {If $n\geq N,$}}  $M$ is a positive integer,  and if $\tilde{z}'\in ({{ {K^3}}}G)^r$ (i.e.,  there is $\tilde{z}''\in G^r$ such that $K^3{{\tilde{z}''}}={{\tilde{z}'}}$) satisfies $|| z-M{z}' ||<\delta$, where $z'=(z'_1,z'_2,\cdots, z'_r){{^T}}$ with  $z'_j=\rho(\tilde{z}'_j)${ ,}
 then there is a $\tilde{u}=(\tilde{c}_j)_{s(n)\times1}\in {{G^{s(n)}}}$ such that
  $$(x_{ij})_{r\times s(n)}\tilde{u}=\tilde{z}'.$$}

Part (b): { W}ith the condition above if we further assume {{that}}  each $s(n)$ can be written as $s(n)=\sum_{k=1}^n l_k$, where {{the}} $l_k$ are positive integers,  and for each $k$,
{{there}} is {{an}}  $r_k\times l_k$ matrix  {{$R_k$}} with entries in $\Z$ {{(for some $r_k\in \N$)}} such that
{{the block diagonal matrix}}
$$
{\bar R}_n={\rm diag}(R_1, R_2,...,R_n)
$$
satisfies
\beq\label{June14-1}
{\bar R}_n{\bar g}_{n}{{ \in (G_{ +} \setminus \{0\})^{\sum_{k=1}^n r_k} ~~~~~\mbox{and} ~~~~{\bar R}_n{\tilde a}\in (\Aff(S)^{{++}})^{\sum_{k=1}^n r_k},}}
\eneq
where ${\bar g}_n={{(g_1^{(n)}, g_2^{(n)},...,g_{s(n)}^{(n)})}}^T$ {{and}} ${\tilde a}=(a_1, a_2, { ...}, a_{s(n)})^T$
(recall the definition of $\Aff(S)^{++}$ --see \ref{Aq}),
$n=1,2,{...}$,  { {then  there exist $\delta, K , N$ as described above but with  $\tilde u$ satisfying { the} extra condition}}  that
\beq\label{JUne14-2}
{\bar R}_{{n}}{\tilde u}>0.
\eneq
\end{lem}
\begin{proof}
{ {The part (a) is essentially Lemma 3.4 of \cite{Lnduke}. To prove part (b), we will repeat the argument {{of}} \cite{Lnduke} and show that if  (\ref{June14-1}) holds, then ${\tilde u}=(\tilde{c}_j)_{s(n)\times1}$ can be chosen to make (\ref{JUne14-2}) hold.

Without loss of generality, {{we may  suppose}} that $(x_{ij})_{r\times r}$ has rank $r$.  Choose $N_0$ such that $s(N_0)\geq r$. Write
$${\bar R}_{N_0} (a_1, a_2, {...,} a_{s(N_0)})^T=(b_1, b_2, {...,} b_l)^{\blue{T}}\in (\Aff(S)^{{++}})^l, $$
where $l=\sum_{i=1}^{N_0}r_i$.
 Let $\ep_0=\frac14\min_{1\leq j\leq l}\inf_{s\in S} \{b_j(s)\}>0$.  There is a positive number $\dt_0<\ep_0$ such that if
$$\|(a'_1, a'_2, {...,} a'_{s(N_0)})^T-(a_1, a_2, {...,} a_{s(N_0)})^T\|<\dt_0, $$
then
$$\|{\bar R}_{N_0}(a'_1, a'_2, {...,} a'_{s(N_0)})^T-(b_1, b_2, {...,} b_l)^{\blue{T}}\|<\ep_0.$$
We further assume that $\dt_0< \frac14\min_{1\leq j\leq s(N_0)}\inf_{s\in S}\{a_j(s)\}.$
Consequently, if $(h_1, h_2, {...,}, h_{s(N_0)})\in G^{{{s(N_0)}}}$ { satisfies}
$$\|(\rho(h_1), \rho({\blue{h_2}}), {...,}, \rho(h_{s(N_0)}))^T-(a_1, a_2, {...,} a_{s(N_0)})^T\|<\dt_0,$$
then
\beq\label{2017-may-9-1}
(h_1, h_2, {...,} h_{s(N_0)})\in (G_+\setminus \{0\})^{s(N_0)}, ~\mbox{and}~~{\bar R}_{N_0}(h_1, h_2, {...,} h_{s(N_0)})\in (G_+\setminus \{0\})^l.
\eneq

Since $A:=(x_{ij})_{r\times r}$ has rank $r$, there is an invertible matrix $B\in M_{r}(\mathbb Q)$ with $BA=I_r$. There is an integer $K>0$ such that all all entries of $KB$ and $K(B)^{-1}$ are integers.  Choose { a} positive number $\dt<\dt_0$ such that if $\|Mz'-z\|<\dt$, then
$\|B(Mz')-B(z)\|<\frac18\dt_0$.  (Note $\dt$ does not depend on $M$.)

Let $N\geq N_0$ {{be}} such that if $n\geq N$, then
$$\|y_n- z\|<\dt, ~~~~~\mbox{and}~~~\|a^{(n)}_j-a_j\| <\frac14 \dt_0, ~~\forall~~j\in\{1, 2, {...,}  s(N_0)\},$$
where {{we}} recall  {{that}} $y_n=(x_{ij})_{r\times s(n)} (a^{(n)}_1, a^{(n)}_2, {...,} a^{(n)}_{s(n)})^T \to z$ as $n \to \infty$. Hence
\beq\label{2017-may-9-2}
\|B(y_n)-B(z)\|<\frac18\dt_0.
\eneq

{{Let us show that}} $K$, $N$, and $\dt$ {{as defined above}} are as desired.

Suppose that $n\ge N\geq N_0$.
Set
$A_n=(x_{ij})_{r\times s(n)}.$  Then  $BA_n=C_n,$ where $C_n=(I_r, D_n')$ for some $r\times (s(n)-r)$ matrix $D_n'$. Since all entries of {{$A_n$ and}} $KB$  are integers, $KD_n'$ is also a matrix with integer entries.

{{Recall}}  that $\rho(g^{(n)}_j)=a^{(n)}_j$, {{and so}} from the first part of  (\ref{June14-1}), we have
$${\bar R}_n (a^{(n)}_1, a^{(n)}_2, {...,} a^{(n)}_{s(n)})^T\in (\Aff(S)^{{++}})^{\sum_{k=1}^n r_k}.$$
{{For each $n\ge N$ (by the continuity of the linear maps),}} there is
{{$0<\dt_1<\dt/4$}}
such that if $(x_1, x_2, { ...}, x_{s(n)})\in \Aff(S)^{s(n)}$ satisfies
$$\|(x_1, x_2, {...,} x_{s(n)})-(a^{(n)}_1, a^{(n)}_2, {...,} a^{(n)}_{s(n)})\|<\dt_1, $$
then
\beq\label{gg17-02}
{{\|(0_{r\times r}, D_n')(x_1, x_2, {...,} x_{s(n)})^T-(0_{r\times r}, D_n')(a^{(n)}_1, a^{(n)}_2, {...,} a^{(n)}_{s(n)})^T\|<\frac{\dt_0}4,}}
\eneq
and ${\bar R}_n (x_1, x_2, {...,} x_{s(n)})^T \in (\Aff(S)^{{++}})^{\sum_{k=1}^n r_k}$. In particular, we have
\beq\label{2017-may-9}
\diag(R_{N_0+1}, R_{N_0+2}, {...,} R_n) (x_{s(N_0)+1}, x_{s(N_0)+2}, {...,} x_{s(n)})^T \in (\Aff(S)^{{++}})^{\sum_{k=N_0+1}^n r_k}.
\eneq
Since $\mathbb D$ is dense in $\mathrm{Aff}(S)$, there {{is}} $\xi\in G_+^{s(n)}$ such that
$\xi=(\tilde{d}_j)_{s(n)\times 1}$   and
\begin{equation}\label{equ5001}
||MK^3\rho(\tilde{d}_j)-a_j^{(n)}||<\dt_1, ~~\rforal~~j=1,2,....,s(n).
\end{equation}

Suppose that $\tilde{z}'$ and $\tilde{z}''$ are {{as}} described in the lemma.  That is, $\tilde{z}''\in G^{r}$, $\tilde{z}'=K^3\tilde{z}'',$ and $\|M\rho^r(\tilde{z}')-z\|<\dt$.  Consequently, $\|MB\rho^r(\tilde{z}')-By_n\|<\frac14\dt_0$ and $\|M\rho^r(\tilde{z}')-y_n\|<2\dt.$

Let $D_n=(0_{r\times r}, D_n')$. {{Then
\beq\label{gg17-03}
C_n-D_n=(I_r, 0_{r\times s(n)-r}).
\eneq
}}
Since both $KB$ and $KD_n$ are matrices over $\mathbb Z$, we can define
{{
\beq\label{gg17-01}
u'=(KB)\tilde{z}''-(KD_n)\xi \in G^{r}.
\eneq
}}

(Warning: $G^r$ is a {{$\Z$-module}} {\blue{rather}}  than a $\Q$-vector space, so only an integer matrix can act on it. That is,  $B\tilde{z}''$ does not make sense, but  $(KB)\tilde{z}''$ makes sense.)

Let $\xi'=({\tilde d}_{r+1}, {\tilde d}_{r+1}, {...,}  {\tilde d}_{s(n)})^T\in  {{G}}^{s(n)-r}$
consist {{of}} the last $s(n)-r$ coordinates of $\xi$.  Then
$D_n \xi= D_n' \xi'$.
Put $u'=(\tilde{c}_1, \tilde{c_2}, ..., \tilde{c}_r)^T\in G_+^{r}, $   and
$$u=
\left(
\begin{array}{c}
u'  \\
K \xi'
\end{array}
\right)= (\tilde{c}_1, \tilde{c_2}, ..., \tilde{c}_r, K{\tilde d}_{r+1}, K{\tilde d}_{r+2}, ...,  K{\tilde d}_{s(n)})^T\in G^{s(n)}.$$
Then
\begin{eqnarray*}
A_n K^2 u&=&(KB^{-1}) (KB) A_n u\\
&=&(KB^{-1}(KI_r, KD_n') u\\
&=&(KB^{-1}) (Ku'+(K^2D_n') \xi')\\
&=&(KB^{-1})(K((KB){\tilde z}''-(KD_n)\xi)+(K^2D_n)\xi)\\
&=&(KB^{-1})(KB)(K{\tilde z}'')\\
&=&K^3{\tilde z}''={\tilde z}'.
\end{eqnarray*}

Let $a'=(a^{(n)}_1, a^{(n)}_2, ..., a^{(n)}_r)^T \in \Aff(S)^r$. Then we have  {{(using \eqref{gg17-01}, \eqref{gg17-02},
and \eqref{gg17-03}),}}
\begin{eqnarray*}
|| MK^2\rho^r(u')-a'||&=&|| MB\rho{{^r}}(\tilde{z}')-MK^3D_n\rho^{s(n)}(\xi)-a'||\\
&\leq&|| By_n-MK^3D_n\rho^{s(n)}(\xi)-a'{{\|}}+\frac{\dt_0}4\\
&\leq&|| By_n-D_n(a^{(n)}_1, a^{(n)}_2, ..., a^{(n)}_{s(n)})^T-a'||+\frac{\dt_0}4+\frac{\dt_0}4\\
&=&|| BA_n(a^{(n)}_1, a^{(n)}_2, ..., a^{(n)}_{s(n)})^T-D_n(a^{(n)}_1, a^{(n)}_2, ..., a^{(n)}_{s(n)})^T-a' ||+\frac{\dt_0}2\\
&=&|| (C_n-D_n)(a^{(n)}_1, a^{(n)}_2, ..., a^{(n)}_{s(n)})^T-a' ||+\frac{\dt_0}2
=|| 0  ||+\frac{\dt_0}2.
\end{eqnarray*}

Let  ${\tilde u}=K^2u=
\left(
\begin{array}{c}
{\tilde u}_1  \\
{\tilde u}_2
\end{array}
\right)$ with ${\tilde u}_1=(K^2\tilde{c}_1, K^2\tilde{c_2}, ..., K^2\tilde{c}_r, K^3{\tilde d}_{r+1}, K^3{\tilde d}_{r+2}, {...,}  K^3{\tilde d}_{s(N_0)})^T$ and ${\tilde u}_2=(K^3{\tilde d}_{s(N_0)+1}, K^3{\tilde d}_{s(N_0)+2}, ...,  K^3{\tilde d}_{s(n)})^T$.
{{Then}} we have $$A_n{\tilde u}= {\tilde z}'.$$
{{Combining}} (\ref{2017-may-9}) and (\ref{equ5001}), we have
$$\diag(R_{N_0+1}, R_{N_0+2}, ..., R_n) {\tilde u}_2 \in (G_+\setminus\{0\})^{\sum_{k=N_0+1}^n r_k}.$$
 Furthermore, {{combining}} (\ref{equ5001}) with the {{estimate above,}} we have
$$\|M\rho^{s(N_0)}({\tilde u}_1)-(a^{(n)}_1, a^{(n)}_2, {...,} a^{(n)}_{s(N_0)})^T\|<\frac{\dt_0}2.$$
Hence,
$$\|M\rho^{s(N_0)}({\tilde u}_1)-(a_1, a_2, {...,} a_{s(N_0)})^T\|<\frac{\dt_0}2+\frac{\dt_0}4<\dt_0.$$
By the choice of $\dt_0$ (see (\ref{2017-may-9-1})), we  have $M{\bar R}_{N_0} {\tilde u}_1\in (G_+\setminus\{0\})^{\sum_{k=1}^{N_0} r_k}$. {{In other words, ${\bar R}_{N_0}{\tilde u}_1\in (G_+\setminus\{0\})^{\sum_{k=1}^{N_0} r_k}.$}} Hence
${\bar R}_n {\tilde u}\in (G_+\setminus\{0\})^{\sum_{k=1}^{n} r_k}$. This ends the proof. }}

\end{proof}

\begin{df}
A unital stably finite \CA\, $A$ is said to have the {{property of}} $K_0$-density, if
$\rho_A(K_0(A))$ is dense in $\mathrm{Aff}(\mathrm{S}_{[1]}(K_0(A))$, where $\mathrm{S}_{[1]}{(K_0(A))}$ is the convex set of  states of $K_0(A)$ (i.e.,
the convex set of all
positive homomorphisms $r: K_0(A)\to \mathbb R$ satisfying
$r([1])=1$).
\end{df}

\begin{rem}
By Corollary 7.9 of \cite{Goodearl}, the linear space spanned by $\rho(K_0(A))$ is always dense in $\mathrm{Aff}(\mathrm{S}_{[1]}(K_0(A))$. Therefore, the unital stably finite \CA\, {{of}} the form of $A\otimes U$ for a UHF-algebra $U$ always has the {$K_0$-}density property. Moreover, any unital stably finite \CA\, $A$ which is tracially approximately divisible has the {$K_0$-}density property.
\end{rem}

\begin{rem}
Not all \CA s in $\mathcal B_1$ with the (SP) property satisfy the {{$K_0$-density}} property. The following  {{example
illustrates this}}:
Consider
\beq\nonumber
G=\{(a,b)\in {\mathbb Q}\oplus {\mathbb Q}:  2a-b \in \mathbb Z\}\andeqn\\
G_+\setminus \{0\}=({\mathbb Q}_+\setminus\{0\}\oplus \Q_+\setminus\{0\} )\cap G,
\eneq
and $1=(1,1)\in G_+$ as unit.
Then $(G, G_+, 1)$ is a weakly unperforated {{rationally}} Riesz simple ordered group (but not a  Riesz group---see \cite{LNjfa}).
Evidently $G$ has {{the property}} (SP); but the image of $G$ is not dense in $\mathrm{Aff}(\mathrm{S}_{[1]}(G)=\mathbb R\oplus \mathbb R$, as $(1/2,1/2)\in \mathbb R\oplus \mathbb R$ is not in the closure of the image of $G$. We leave the details to the readers.



\end{rem}

\begin{prop}\label{est-a2m}
Let $A\in\mathcal B_{0}$ {{ be an amenable \CA\, which has}} the {{$K_0$-}}density property {{and  satisfies the UCT,}} and let $B_1$ be an inductive limit \CA\, {{as}} in Theorem \ref{RangT} such that $$(K_0(A),{K_0(A)_+},[1_A],K_1(A))\cong(K_0(B),{K_0(B)_+},[1_B],K_1(B)),$$
{\blue{and such that $B$ satisfies the condition {{of}} Corollary \ref{smallmap} for $B$ there,}} where $B=B_1\otimes U$ for a UHF-algebra $U$ of infinite type.
 Let $\alpha\in KL(A, B)$ be an element which implements the isomorphism above. Then,
 there is a sequence of completely positive linear maps $L_n: A\to B$ such that
\beq\nonumber
\lim_{n\to\infty} ||L_n(ab)-L_n(a)L_n(b)||= 0\tforal a, b\in A\,\,
{{\tand [\{L_n\}]=\af.}}
\eneq
{\blue{(By $[\{L_n\}]=\af$, we mean for any finite set ${\cal P}\subset \underline{K}(A)$, there is a positive integer $N$ such that if $n\geq N$, then $[L_n]|_{\cal P}=\af|_{\cal P}$.)}}
\end{prop}

\begin{proof}
By Lemma \ref{MF}, $A$ is the closure of an increasing union of RFD \SCA s $\{A_n\}$.
{{Fix a finite subset ${\cal P}\subset \underline{K}(A).$}}
Let $\mathcal F_1$ be a finite subset of ${{A}} $ and let $\delta_0>0$ be such that
{{for}} any $\mathcal F_1$-$\delta_0$-multiplicative {{contractive}} linear map {{$L$ from $A,$}} the map $[L]|_{\mathcal P}$ is well defined.
We may assume {{
that ${\cal F}_1\subset A_1$ and there exists}}
${\mathcal P}'\subset \underline{K}(A_1)$ {{such that $[\iota]({\cal P}')={\cal P},$ where
$\iota: A_1\to A$ is the embedding.}}
 Let $G=G({\cal P})$ be the subgroup generated by ${\cal P}$ and let ${\mathcal P}_0{{'}}\subset{\mathcal P}{{'}}$ be such that ${\cal P}_0{{:=\iota_{*0}({\mathcal P}_0')}}$ generate{{s}} $G\cap{K_0}(A).$
{{We may suppose that}}
${\mathcal P}_0=\{{{[p_1],..., [p_l]}}\},$ where $p_1, ..., p_l$ are
projections in a matrix algebra over $A_{{1}},$ {{where
we identify $\iota(p_i)$ with $p_i$.}}  Let $G_0$ be the group generated by ${\cal P}_0.$
Moreover,
{{we may assume}} that $\mathcal F_1$ and
$\delta_0$ satisfy {{the conclusion of}} Lemma \ref{ker}. Let $k_0$ be an integer such that
$G({\mathcal P})\cap\mbox{K}_i(A,\mathbb Z/k\mathbb Z)=\{0\}$ for any $k\geq k_0$, $i=0,1.$

By {\Green{(the ``Moreover'' part of)}} Theorem \ref{kkmaps}, there are two $\mathcal F_1$-$\delta_0/2$ multiplicative \morp s
$\Phi_0, \Phi_1$ from $A$ to $B\otimes\mathcal {\mathcal K}$ ($B$ is amenable) such that
$$[\Phi_0]|_{\mathcal P}=\alpha|_{\mathcal P}+[\Phi_1]|_{\mathcal P}$$
and the image of $\Phi_1$ is {{contained}} in a finite dimensional \SCA.
Moreover, we may also  assume that $\Phi_1$ is a homomorphism when it is restricted {{to}} $A_1$, and the image is a finite dimensional \CA.
With $\Phi_1$ in the role of $h_0$, we can proceed with the construction as described at the beginning of this
section. We will keep the same notation.

Consider the map $\tilde\rho:G({\mathcal P})\cap\mbox{K}_0(A)\to l^\infty(\mathbb Z)$ defined in Corollary \ref{rho}. The
linear span of {{the subset}} $\{\tilde\rho(p_1), ..., \tilde\rho(p_l)\}$ over $\mathbb Q$ will have finite rank, say $r$. So,
we may assume that $\{\tilde\rho(p_1), ..., \tilde\rho(p_r)\}$ {{is}}  linearly independent and {{its}}  $\mathbb Q$-linear span
{{is}}
the whole subspace. Therefore, there is an integer $M$ such that for any $g\in\tilde\rho(G_0)$, the element $Mg$
is in the subgroup generated by $\{\tilde\rho(p_1), ..., \tilde\rho(p_r)\}$. Let
$x_{ij}=(\tilde\rho(p_i))_j$, and $z_i=\rho_B(\af([p_i])\in\mathbb{D}$, where $\mathbb{D}:=\rho_B(K_0(B))\subset
\text{Aff}(S_{[1]}(K_0(B))$. Since $A$ is assumed to have the  $K_0$-density property, so also is $B$. Therefore the image $\mathbb{D}$ is a dense subgroup of $\text{Aff}(\mathrm{S}_{[1]}({K_0}(B)))$.

Let  $\{S_j\}$ be the sequence of \SCA s in ${\cal C}_0$ in the construction at the beginning of this section.
Fix $k\ge 1.$
Let $e_i^{{{k+j}}}, e_{i,\pm}^{{{k+j}}}\in K_0(S_{k+j}), $ $i=1,2,..., l_{k+j},$  and {{the}}
$r_{k+j}\times 2l_{k+j}$ matrix $R_{k+j}$ {{be}} as described in \ref{Constr-Sn}.
{Let $s(j)=\sum_{i=1}^j2l_{k+i},$ $j=1,2,....$} Put
\beq\label{June14n2-1}
\af([ \imath_n\circ {\blue h}_{k+j, n}(e_{i,+}^{{k+j}})])=g^{(n)}_{{s(j-1)+2i-1}}{{\in K_0(B)_+}}, \\
\af([\imath_n\circ {\blue h}_{k+j,{{n}}}(e_{i,-}^{{k+j}})])=g^{(n)}_{s(j-1)+2i}{{\in K_0(B)_+}},
 \eneq
 $i=1,2,...,l_{\blue{k+j}},$
and $a_j^{(n)}=\rho_B(g_j^{(n)}),$ $j=1,2,..., s(n)=\sum_{j=1}^n l_{k+j},$ $n=1,2,....$ Note that $a_j^{(n)}\in\mathbb{D}^+\backslash\{0\}.$ It follows from Lemma \ref{traces} that
$\lim_{n\to\infty}a_{{s(j-1)+2i}}^{(n)}=a_{s(j-1)+2i}={{\lim_{n\to\infty}}}\rho_B(\af(g_{s(j-1)+2i}^{(n)}))\\>0$  and  $\lim_{n\to\infty}a_{s(j-1)+2i-1}^{(n)}=a_{s(j-1)+2i-1}
={{\lim_{n\to\infty}}}\rho_B(\af(g_{s(j-1)+2i-1}^{(n)}))>0$ uniformly.
 Moreover, by Lemma \ref{rho-uniform},
$\sum_{j=1}^nx_{ij}a_j^{(n)}\to z_i$ uniformly.  Furthermore, by \ref{Constr-Sn} {{(b), we know that $\{x_{ij}\}$, $g^n_j$,
{{and}} $R_n$ {{satisfy}} the conditions of Lemma \ref{solveeq}, with $K_0(B)$ in place of $G$ and ${{T(B)}}$ in place of $S$.}}
So Lemma \ref{solveeq} applies. Fix $K$ and $\delta$ obtained from Lemma \ref{solveeq}.

Let $\Psi:=\Phi_0\oplus(\overbrace{\Phi_1\oplus\cdots\oplus\Phi_1}^{MK^3(k_0+1)!-1})$. Since $\Phi_1$
factors through a finite dimensional \CA, it is zero when restricted to $K_1(A)\cap G$  and
${K_1}(A,\mathbb Z/k\mathbb Z)\cap G$ for
$2\le k\le k_0.$ Moreover,  the map $(\overbrace{\Phi_1\oplus\cdots\oplus\Phi_1}^{MK^3(k_0+1)!})$ vanishes on ${K_0}(A,\mathbb Z/k\mathbb Z)$ for $2\le k\le k_0.$ Therefore we have
$$[\Psi]|_{K_1(A)\cap G}=\alpha|_{K_1(A)\cap G},\quad [\Psi]|_{K_1(A,\mathbb Z/k{\blue{\mathbb Z}})\cap G}=\alpha|_{K_1(A,\mathbb Z/k\mathbb Z)\cap G}$$ and $[\Psi]|_{K_0(A,\mathbb Z/k\mathbb Z)\cap
G}=\alpha|_{K_0(A,\mathbb Z/k\mathbb Z)\cap G}$. We may assume {{that}} $\Psi(1_A)$ is a projection in M$_{ m}(B)$ for some integer ${ m}$.

We may also assume  that there exist projections $\{p'_1, ..., p'_l\}$ in  $B\otimes{\mathcal K}$ which are sufficiently close
to $\{\Psi(p_1), ..., \Psi(p_l)\},$ respectively, that $[p'_i]=[\Psi(p_i)]$. Note that $B\in\mathcal B_0$, and hence the strict order on the projections of $B$ is determined by traces. Thus there is a projection $q'_i\le p'_i$ such that $[q'_i]=MK^3(k_0+1)![\Phi_1(p_i)]$. Set $e'_i=p'_i-q'_i$, and let
${\mathcal P}_1=\Psi({\mathcal P})\cup\Phi_1({\mathcal P})\cup\{p'_i,q'_i,e'_i;i=1,...,l\}$. Denote by $G_1$ the group generated
by ${\mathcal P}_1$. Recall that $G_0=G({\cal P})\cap{K_0}(A)$ {{
and $\rho_A(G_0)$ is a  free abelian group}}, and decompose it as $G_{00}\oplus G_{01}$, where
$G_{00}$ is the infinitesimal part of $G_0$. Fix this decomposition and denote by $\{d_1, ..., d_t\}$ the positive elements which generate $G_{01}$.

{\blue{Note that $B$ belongs to  the class of  Corollary \ref{smallmap}, and so $M_r(B)$ does also.}} Applying Corollary \ref{smallmap}
to $M_r(B)$ with any finite subset $\mathcal{G}$, any ${\blue{\ep}}>0$ and any
$0<r_0<\delta<1$, one has a $\mathcal{G}$-$\epsilon$-multiplicative map $L:M_{ m}(B)\to M_{ m}(B)$ with the following properties:
\begin{enumerate}
\item $[L]|_{{\mathcal P}_1}$ and $[L]|_{G_1}$ are well defined;
\item $[L]$ induces the identity maps on the infinitesimal part
      of $G_1\cap$K$_0(B)$, $G_1\cap\mbox{K}_1(B)$,
      $G_1\cap\mbox{K}_0(B,\mathbb Z/k\mathbb Z),$ and $G_1\cap\mbox{K}_1(B, \mathbb Z/k\mathbb Z)$
      for the $k$ with $G_1\cap\mbox{K}_i(B, \mathbb Z/k\mathbb Z)\neq\{0\}$,
      $i=0,1$;
\item $\tau\circ[L](g)\leq r_0\tau(g)$ for all $g\in
      G_1\cap\mbox{K}_0(B)$ and $\tau\in T(B)$; {{and}}
\item there exist positive elements $\{f_i\}\subset{K_0}(B)_+$
      such that for  $i=1,...,t,$
      $$\alpha(d_i)-[L](\alpha(d_i))=MK^3(k_0+1)!f_i.$$
\end{enumerate}

 Using the compactness of $T(B)$ and strict comparison for positive elements {{of}}  $B$, the positive number $r_0$ can be chosen sufficiently small  that
 $\tau\circ[L]\circ[\Psi]([p_i])<\delta/2$ for all $\tau\in T(B)$,
 and $\alpha([p_i])-[L\circ\Psi]([p_i])>0,$  $i=1,2,...,l.$

Let $[p_i]=\sum_{j=1}^tm_jd_j+s_{{i}},$ where $m_j\in\Z$ and $s_{{i}}\in G_{00}.$
Note {{that}}, by (2) above, $(\af-[L]\circ \af)({s_i})=0.$  Then we have
\beq\nonumber
 & \hspace{-0.8in} \alpha([p_i])-[L\circ\Psi]([p_i])
  =\alpha([p_i])-([L\circ\alpha]([p_i])+MK^3(k_0+1)![L\circ\Phi_1]([p_i])) \\\nonumber
  &=(\alpha(\sum m_jd_j)-[L\circ\alpha](\sum m_jd_j))-MK^3(k_0+1)![L\circ\Phi_1]([p_i]) \\\label{june232019}
  &= MK^3(k_0+1)!(\sum m_jf_j-[L]\circ[\Phi_1]([p_i]))=MK^3(k_0+1)!f_i',
 \eneq
 {where $f_i'=\sum m_jf_j-[L]\circ[\Phi_1]([p_i]),$}
for $i=1,2,...,l.$
 {{Set}} $\bt([p_i])=K^3(k_0+1)!f_i',$ $i=1,2,...,l.$



Let us now construct a map $h': A \to B.$
It will be constructed by
factoring through the $K_0$-group of some \CA\, in the class $\mathcal C_0$ in the construction given at the beginning of this section.
Let $\tilde{z_i}'=\beta([p_i])$, and $z'_i=\rho_B(\tilde{z_i}')\in\Aff(S_{[1]}(K_0(B)))$.
Then we have
$$\begin{array}{lll}
    ||{Mz'-z}||_\infty & = & \max_i\{||\rho_B(\alpha([p_i])-[L\circ\Psi]([p_i]))-\rho(\af([p_i]))||\}\\
     & = & \max_i\{\sup_{\tau\in T(B)}\{\tau\circ[L]\circ[\Psi]([p_i])\}
     \le \delta/2,
  \end{array}$$
  where $z=(z_1,z_2,...,z_r)$ and $z'=(z_1',z_2',...,z_r').$
By Lemma \ref{solveeq},  for sufficiently large $n,$ one obtains
$\tilde{u}=(u_1,u_2,...,u_{s(n)})\in {{K_0(B)}}^{{s(n)}}$
such that

\begin{equation}\label{June-16-n1}
\sum_{j=1}^{s(n)}  x_{ij}u_j=\tilde{z}'_i.
\end{equation}
More importantly,
\beq\label{June16-3}
{\bar R}_n{\tilde u}>0.
\eneq
It follows from \ref{Constr-Sn}{{(a)}} {{(recall that $k$ is fixed) that, for each $1\le j\le n,$}}  the map
$$e_i^{k+j}\mapsto (u_{{s(j-1)}+2i-1}-u_{s(j-1)+2i}),\,\,
1\leq i\leq l_{k+j},$$
defines {{a}}
strictly positive homomorphism $\kappa_0^{(k+j)}$ from $K_0(S_{k+j})$ to $K_0(B).$
Since $B\in\mathcal B_{u0}$, by Corollary \ref{C0ext},  there is a
 \hm\, $h':D\to \mbox{M}_m(B)$ for some large $m$
 such that $h'_{*0}|_{S_k}=\kappa_0^{(k)},$ {{ where $D=S_{{{k+1}}}\oplus\cdots\oplus S_{{k+n}}.$}}
By \eqref{June-16-n1}, one has, keeping the notation {{of}} the construction at the beginning of this section,
{{for $[\psi_{k+j}^{k}(p_i)]=\sum_{l=1}^{l_{k+j}}m_i^{k+j}([\psi_{k+j}^k(p_i)])e_l^{k+j},$}}
{{\beq\nonumber
\kappa_0^{k+j}([\psi_{k+j}^{k}(p_i)])&=&\sum_{l=1}^{l_{k+j}}m_i^{k+j}([\psi_{k+j}^k(p_i)])\kappa_0^{k+j}(e_l^{k+j})\\\nonumber
&=&\sum_{l=1}^{l_{k+j}}m_i^{k+j}([\psi_{k+j}^k(p_i)])(u_{{s(j-1)}+2l-1}-u_{s(j-1)+2l})
=\sum_{l=1}^{l_{k+j}} x_{i,j+l}u_{s(j-1)+l}.
\eneq}}
{{Hence,}}
\beq
h'_{*0}({{[\psi^{{{k}}}_{k}(p_i)],}}[\psi^{{{k}}}_{k+1}(p_i)], [{{\psi^k_{k+2}}}(p_i)], ...,[{{\psi^k_{{n+k-1}}}}(p_i)])
{{=(x_{i,j})_{s(n)\times 1}{\tilde u}}}
=\beta([p_i]),
\eneq
$i=1,...,r.$
Now, define $h'': A\to  {{D \to}} \mathrm{M}_{{{m}}}(B)$ by $$h''=h'\circ ({{\psi^k_k\oplus}}{{{\psi^k_{k+1}}}}\oplus {{\psi^k_{k+2}}}{\oplus}\cdots\oplus {{\psi^k_{{k+n-1}}}}).$$
Then $h''$ is $\mathcal F$-$\delta$-multiplicative.

For any $x\in\ker\tilde{\rho}$, by Lemma \ref{ker}, $x\in\ker {{(\rho_B\circ\alpha)}}\cap\ker [H]$ and $x\in\ker[h_0]=\ker[\Phi_1]$. Therefore, we have $[\Phi_1](x)=0$ and $[\Psi](x)=\alpha(x).$ Note that $\alpha(x)$ also vanishes under any state of $({K_0}(B), {K_0}^+(B))$, {{and}} we have $[L]\circ\alpha(x)=\alpha(x)$. So, we get $$\alpha(x)-[L\circ\Psi](x)=0.$$
{{Thus,}} $(\alpha-[L\circ\Psi])|_{{\rm ker}{\tilde \rho}}=0.$
{{Hence}} we may view $\af-[L\circ \Psi]$ as a \hm\, {from} ${\tilde \rho}(G_0).$
{Recall} that
\beq\label{June16-11}
(\af-[L\circ \Psi])([p_i])=M\bt([p_i]),\,\,\,i=1,2,...,r.
\eneq

Set $h$ to be {{the direct sum of}} $M$ copies of $h''$. The map $h$ is $\mathcal F$-$\delta$-multiplicative, and $$[h]([p_i])=\alpha([p_i])-[L]\circ[\Psi]([p_i])\quad i=1,...,{r}.$$ Note that $[h]$  has  multiplicity $MK^3(k_0+1)!$ {{(see \eqref{june232019})}}, and $D
\in\mathcal C_0$ (the algebras in ${\cal C}$ with trivial ${K_1}$ group). One concludes that $h$ induces {{the}} zero map on $G\cap{K_1}(A)$, $G\cap{K_1}(A,\mathbb Z/m\mathbb Z),$ and $G\cap{K_1}(A,\mathbb Z/m\mathbb Z)$ for $m\leq k_0$. Therefore, we have $$[h]|_{\mathcal P}=\alpha|_{\mathcal P}-[L]\circ[\Psi]|_{\mathcal P}.$$

Set $L_1=(L\circ\Psi)\oplus h$. {{$L_1$}} is $\mathcal F$-$\delta$ multiplicative and $$[L_1]|_{\mathcal P}=[h]|_{\mathcal P}+[L]\circ[\Psi]|_{\mathcal P}=\alpha|_{\mathcal P}.$$
We may assume $L_1(1_A)=1_B$ by conjugating with a {{unitary in $M_m(B).$}}
Then $L_1$ is an $\mathcal F$-$\delta$-multiplicative map from $A$ to $B$, and $[L_1]|_{\mathcal P}=\alpha|_{\mathcal P}$.

{{Since $A$ is separable, $\underline{K}(A)$ is countable.
It follows that one obtains a sequence of \morp s from $A$ to $B$ such that
\beq\nonumber
\lim_{n\to\infty} ||L_n(ab)-L_n(a)L_n(b)||= 0\rforal a, b\in A\,\,
\tand [\{L_n\}]=\af.
\eneq
}}
\end{proof}

\begin{cor}\label{kk-attain}
Let $A$ be a {{separable}} amenable  \CA\, in the class $\mathcal B_{{0}}$ satisfying
 the UCT. Then $A$ is KK-attainable with respect to ${\cal B}_{u0}$ (see  \ref{D181}).
\end{cor}


\begin{proof}
Let $C$ be any \CA\, in $\mathcal B_{u0}$, and let $\alpha\in KL(A, C)^{++}$. We may write  $C=C_1\otimes U$ for some $C_1\in {\cal B}_0$ and for some UHF-algebra of infinite type.

{\Green{For any \CA\, $D,$ let us denote {{by}} $j_D: D\to D\otimes U$ the map defined by $j_D(d)=d\otimes 1_U$ for all $d
\in D.$}}

{\Green{
Let $\af\in  KL(A, C_1\otimes U)^{++}={\rm Hom}_{\Lambda}(\underline{K}(A), \underline{K}(C_1
\otimes U))^{++}.$
Let us point out an easy fact that there is ${\bar \af}\in KL(A\otimes U, C_1\otimes U)^{++}$
such that $\af={\bar \af}\circ [j_A].$ To see this, recall that, since $U$ is a UHF-algebra of infinite type,
$[j_U]=[{\rm id}_U].$   Thus $[j_C]\circ \af=\af$ (as we identify $C$ with $C\otimes U$).}}

{\Green{For each $k\ge 1,$ since $K_i(U)$ is torsion free, by K\"{u}nneth  formula,  for $i=0,1,$
\beq
K_i(D\otimes U)=K_i(D)\otimes K_0(U)\andeqn K_i(D\otimes U, \Z/k\Z)=K_i(D,\Z/k\Z)\otimes K_0(U),
\eneq
where $D=A,$ or $D=C,$
$k=2,3,....$
Define, for each $k\ge 1,$ if $y=x\otimes r,$ where $x\in K_i(A,\Z/k\Z)$ and $r\in K_0(U),$
${\bar \af}(y)=\af(x)\otimes r.$  One checks
that
\beq
{\bar \af}\circ [j_A](x)=\af(x\otimes 1)=[j_C]\circ \af(x)=\af(x)\rforal x\in \underline{K}(A).
\eneq
Therefore
\beq\label{alpha=bar}
\af={\bar \af}\circ [j_A].
\eneq
}}
By Theorem \ref{RangT}, there is a \CA\, $B$ which is an inductive limit of \CA s in the class $\mathcal C_0$ together with homogeneous \CA s in the class $\mathbf H$ {{(see \ref{AHblock})}}
such that
\beq\label{6232019-1}
(K_0(A_1), {K_0(A_1)_+}, [1_{A_1}]_0, K_1(A_1)) \cong {\Green{(K_0(B_1), {K_0(B_1)_+}, [1_{B_1}]_0, K_1(B_1))}},
\eneq
{\Green{where $A_1=A\otimes U,$}}
{\blue{and such that ${\Green{B_1:=}}B\otimes U$ satisfies the condition {{of}} Corollary \ref{smallmap} for $B$ there.}}
Since ${\Green{A_1}}$ satisfies the UCT, there is an invertible $\beta\in KL({\Green{A_1, B_1}})^{++}$ such that $\beta$ carries the isomorphism {{in
\eqref{6232019-1}.}}
{{Applying}} Corollary \ref{est-m2a}, one obtains a sequence of unital \cp s ${\blue{\{\Psi_n: {\Green{B_1}} \to C\}_{n=1}^{\infty}}}$ such that
\beq
\lim_{n\to\infty}\|\Psi_n(ab)-\Psi_{\blue{n}}(a)\Psi_{\blue{n}}(b)\|=0\rforal a, b\in {\Green{B_1}}\andeqn [\{\Psi_{\blue{n}}\}]={\Green{\bar \af}}\circ \bt^{-1}.
\eneq
{{Note that $B\otimes U$ can be chosen as in Proposition \ref{est-a2m}.}} Applying Proposition \ref{est-a2m} to
${\Green{\bt}}\in KL({\Green{A_1}}, B\otimes U),$  one obtains a sequence of unital \cp s $\{\Phi_n\}$ from
${\Green{A_1}}$ to $B{{\otimes U}}$  such that $[\{\Phi_n\}]={\Green{\bt}}$
and $\lim_{n\to\infty}\|\Phi_n(ab)-\Phi_n(a)\Phi_n(b)\|=0$
for all $a,\, b\in {\Green{A_1}}.$
Choosing a subsequence $\{{\Green{\Psi}}_{k(n)}\}$ and
{{defining}} $L_n={\Green{\Psi}}_{k(n)}\circ \Phi_n{\Green{\circ j_A}},$ one checks that {\Green{(see \eqref{alpha=bar})}}
$$[\{L_n\}]={\Green{[\{\Psi_{k(n)}\}]\circ[\{\Phi_n\}]{\Green{\circ [j_A]}}
={\bar \af}\circ \bt^{-1}\circ \beta\circ [j_A]={\bar \af}\circ [j_A]=}}\af\in KL(A,B)$$ and $\lim_{n\to\infty}\|L_n(ab)-L_n(a)L_n(b)\|=0$ for all $a,\, b\in A,$ as desired.
\end{proof}

\begin{thm}\label{MEST}
Let $A\in\mathcal B_{{0}}$ be {{a}}  amenable  \CA\,
satisfying the UCT, and let $B\in {\mathcal B}_{u0}$. Then for any $\alpha\in KL(A, B)^{++}$, and any affine
continuous map $\gamma: T(B)\to T(A)$ which is compatible {{with}}  $\alpha$, there is a sequence of completely positive linear maps $L_n: A\to B$ such that
$$\lim_{n\to\infty} ||L_n(ab)-L_n(a)L_n(b)||= 0\rforal a, b\in A,$$
$$[\{L_n\}]=\alpha,\quad\mathrm{and}\quad \lim_{n\to\infty}\sup_{\tau\in T(B)}|\tau\circ L_n(f)-\gamma(\tau)(f)|= 0\rforal f\in A.$$
\end{thm}

\begin{proof}
{{This}} follows from Corollary \ref{kk-attain} and Proposition \ref{add-tr} directly.
\end{proof}

\section{The isomorphism theorem}

\begin{df}\label{bdbk-k1}
{{Let $n$ and $s$  be integers. Let $Y$ be}} a connected finite CW complex  with  dimension no more than $3$
{{such that}}
$K_1(C(Y))$ is a finite group (no restriction on $K_0(C(Y))$),
and let $P$ be a projection in $\text{M}_n(\text{C}(Y))$ with rank $r\geq 6.$ {{Note that $Y$ could be a point.}}
{{Let $D'=\bigoplus_{i=1}^s E_i,$ where $E_i=M_{r_i}(C(\T))$ and where $r_i$ is an integer, $i=1,2,...,s.$
Hence}} $K_1(D')\cong \Z^s.$
{{Put $C'=D'\oplus PM_n(C(Y))P.$}}
 Then $K_1(C')=\mathrm{Tor}(K_1(C'))\oplus{{\Z^s}}.$

{\blue{Let}} $C$ {\blue{be}} a finite direct sum of  \CA s of the form $C'$ above
{and} \CA s in {{the class}} $\mathcal C_0$ (with trivial $K_1$-group).
{\blue{Write}} $C=D\oplus C_0\oplus C_1,$ where  {{$D$ is a finite direct sum of \CA s of the form
$M_{r_i}(C(\T)),$}}
$C_0$ is a direct sum of \CA s in {{the class}} ${\cal C}_0,$ and
{{$C_1$ is a finite direct sum
of \CA s of}} the form $PM_n(C(Y))P,$ where $Y$ is a connected finite CW complex with  dimension no more than 3 such that
$K_1(C(Y))$ is a finite abelian group and $P$  has rank at least $6.$
Then, one has {{(recall that $P$ has rank at least 6 on $Y$) $U(C)/U_0(C)=K_1(C)=K_1(D) \oplus \mathrm{Tor}(K_1(C))$ and}}
 \beq\label{UCrank}
 U(C)/CU(C)\cong U_0(C)/CU(C)\oplus K_1(D) \oplus \mathrm{Tor}(K_1(C)).
 \eneq
Here we identify $K_1(D)\oplus {\rm Tor}(K_1(C))$ with a subgroup of $U(C)/CU(C).$
{{Denote}} by $\pi_0, \pi_1, \pi_2$  the projection maps from $U(C)/CU(C)$ to each {{direct summand}}  according to the decomposition above {{(so $\pi_0(U(C)/CU(C))=U_0(C)/CU(C),$
$\pi_1(U(C)/CU(C))=K_1(D)$ and $\pi_2(U(C)/CU(C))=\mathrm{Tor}(K_1(C))$).}} {{
We note that $\pi_i(U(C)/CU(C))$ is a subgroup
of $U(C)/CU(C),$ $i=0,1,2.$}} 
{{We may also write $U(C)/CU(C)=U(D)/CU(D)\oplus U_0(C_0)/CU(C_0)\oplus U(C_1)/CU(C_1).$}}
 {{Denote by $\Pi: C\to D$  the projection map.  Viewing $D$ as a \SCA\, of  $C,$ {{one has}} $\Pi|_D={\rm id}_D.$
Since $\pi_1(U(C)/CU(C))$ is a subgroup
 of $U(C)/CU(C),$ $\Pi^{\ddag}\circ \pi_1$
 is a \hm\, from
$U(C)/CU(C)$ into $U(D)/CU(D)$ and
$\Pi^{\ddag}|_{\pi_1(U(C)/CU(C))}$ is injective (see \ref{DLddag} for the definition {{of}} $\Pi^{\ddag}$).}}

{{In what follows in this section, if $A$ is a unital \CA\, and $F\subset U(A)$ is a subset,
then $\overline{F}$ is the image of $F$ in $U(A)/CU(A)$ (see \ref{Dcu}).}}
We will frequently refer to the above notation later in  this section.
\end{df}

{{Definition \ref{bdbk-k1} plays a {{role similar to that which}} Definition 7.1 of \cite{LinTAI} plays in \cite{LinTAI}.
The only difference is the appearance of $C_0\in {\cal C}_0.$ Since $K_1(C_0)=\{0\},$
 as one will see in this section, {{this}} will
not cause a new problem. However, since \CA s in ${\cal C}_0$  have a different form, we will repeat
many of the same arguments for the {{sake of}} completeness.}}

As in \cite{LinTAI}, we have the following lemmas to control the maps from  $U(C)/CU(C)$ in the approximate intertwining argument in the proof of \ref{IST0}. The proofs are  repetitions of the corresponding arguments in \cite{LinTAI}.

\begin{lem}[see Lemma 7.2 of \cite{LinTAI}]\label{LinTAI-72}
Let $C$ be the \CA\, {{defined in \ref{bdbk-k1}, let $D,$ a finite direct sum of circle algebras,
be the direct summand of $C$  specified
in \ref{bdbk-k1},}}
let $\mathcal U\subset U(C)$ be a finite subset, and {{denote by}} $F$  the subgroup generated by $\mathcal U$.
{{Let}} $G$ {{be}} a subgroup of $U(C)/CU(C)$ which contains $\overline{F},$ the image of $F$ in $U(C)/CU(C),$ {{and}}  also contains $\pi_1(U(C)/CU(C))$ {{and}}\\ $\pi_2(U(C)/CU(C))$. Suppose that the
composed
map $\gamma: \overline{F}\to U(D)/CU(D)\to U(D)/U_0(D)$ is injective{{---i.e., if $x,y \in \overline{F}$ and $x\not=y$, then $[x]\not=[y]$ in $U(D)/U_0(D)$.}} Let $B$ be a unital \CA\, and $\Lambda: G\to U(B)/CU(B)$ be a homomorphism such that $\Lambda(G\cap(U_0(C)/CU(C)))\subset U_0(B)/CU(B)$.
Let $\theta: \pi_2(U(C)/CU(C))\to U(B)/CU(B)$ be defined by
$\theta(g)=\Lambda|_{\pi_2(U(C)/CU(C))}(g^{-1})$ for any $g\in\pi_2(U(C)/CU(C))$.
Then there is a homomorphism $\beta: U(D)/CU(D)\to U(B)/CU(B)$ with $$\beta(U_0(D)/CU(D))\subset U_0(B)/CU(B)$$  such that $$\beta\circ\Pi^{\ddagger}\circ\pi_1(\bar{w})=\Lambda(\bar{w})(\theta\circ\pi_2(\bar{w}))\tforal w\in F,$$
{{where $\Pi: C\to D$  is defined in   Definition {\blue{\ref{bdbk-k1}.}}}}
If, furthermore, $B\cong B_1\otimes {{V}}$ for a unital \CA\, $B_1\in\mathcal B_0$ and a UHF-algebra ${{V}}$ {\blue{of infinite type}}, and
$\Lambda(G)\subset U_0(B)/CU(B)$, then $\beta\circ\Pi^{\ddagger}\circ(\pi_1)|_{\bar{F}}=\Lambda|_{\bar{F}}$.
\end{lem}
The {{statement}} above
{{is}} summarized
in  the following commutative diagram:
\begin{displaymath}
\xymatrix{
{\bar F} \ar[d]^{\pi_1}  \ar[r]^{\mathrm{inclusion}} & G \ar[d]^{\Lambda+\theta\circ\pi_2} \\
\pi_1(\bar F) \ar[d]^{\Pi^\ddag} &  U(B)/CU(B) \\
U(D)/CU(D) \ar[ur]_{\beta}\,\,.&
}
\end{displaymath}
%

\begin{proof}
{{Note  that $U(D)/U_0(D)\cong K_1(D),$ and this is a free abelian group, as $D$ is a finite direct sum of \CA s of the form $M_n(C(\T)).$
Therefore $\gamma(\overline{F})$ as a subgroup of a free abelian group  is free abelian.}}
The proof is exactly the same as that of Lemma 7.2 of \cite{LinTAI} {{(since it has nothing to do with the
direct summands of $C_0,$ and, $K_1(C_0)=\{0\}$ and $U(C_0)/CU(C_0)=U_0(C_0)/CU(C_0)$). However,
we will repeat the proof.}}

Let $\kappa_1: U(D)/CU(D)\to K_1(D)\subset U(C)/CU(C)$ be  the quotient map and let\\
 $\eta: \pi_1(U(C)/CU(C))\to K_1(D)$ be the map defined by $\eta=\kappa_1\circ\Pi^{\ddagger}|_{\pi_1(U(C)/CU(C))}$.
{{Note that $\pi_1(U(C)/CU(C))$ is identified with $K_1(D)$ and  $\Pi^{\ddagger}|_{\pi(U(C)/CU(C)}$ is injective(see \ref{bdbk-k1}).}} {{Therefore}}
$\eta$ is an isomorphism.
Since {{(the composed map)}} $\gamma$ is injective and $\gamma(\overline{F})$ is free abelian, we conclude that
 $\kappa_1\circ\Pi^{\ddagger}\circ\pi_1$ is also injective on $\overline{F}$. Since $U_0(C)/CU(C)$ is divisible ({{see, for example, Lemma \ref{UCUdiv}}}),
there is a homomorphism $\lambda: K_1(D)\to U_0(C)/CU(C)$ such that
$$\lambda|_{\kappa_1\circ\Pi^\ddagger\circ\pi_1(\overline{F})}=\pi_0\circ((\kappa_1\circ\Pi^{\ddagger}\circ\pi_1)|_{\overline{F}})^{-1}{{,}}$$
{{where $((\kappa_1\circ\Pi^{\ddagger}\circ\pi_1)|_{\overline{F}})^{-1}: ~~\eta\circ \pi_1(\overline{F})\to \overline{F}$ is the inverse of the injective map $(\kappa_1\circ\Pi^{\ddagger}\circ\pi_1)|_{\overline{F}}$.}} This could be viewed as the following commutative diagram:
\begin{displaymath}
\xymatrix{
& K_1(D) \ar[dr]^\lambda& \\
(\eta\circ \pi_1)({\bar F}) \ar[rr]^{\pi_0\circ(\eta\circ\pi_1)^{-1}} \ar[ur] & & U_0(C)/CU(C). \\
}
\end{displaymath}

Define $$\beta=\Lambda((\eta^{-1}\circ\kappa_1)\oplus(\lambda\circ\kappa_1)),$$
{{a \hm\, from $U(D)/CU(D)$ to $U(B)/CU(B).$}}
Then, for any $\overline{w}\in\overline{F}$,
$$\beta(\Pi^{\ddagger}\circ\pi_1(\overline{w}))=\Lambda(\eta^{-1}(\kappa_1\circ\Pi^{\ddagger}(\pi_1(\overline{w})))\oplus \lambda\circ\kappa_1(\Pi^\ddagger(\pi_1(\overline{w}))))=\Lambda(\pi_1(\overline{w})\oplus\pi_0(\overline{w})).$$
{{Recall that}}
$\theta: \pi_2(U(C)/CU(C))\to U(B)/CU(B)$ {{is defined}}
by
$$\theta(x)=\Lambda(x^{-1})\rforal x\in\pi_2(U(C)/CU(C)).$$
Then
\beq\label{731-nn1}
\beta(\Pi^{\ddagger}(\pi_1(\overline{w})))=\Lambda(\overline{w})\theta(\pi_2(\overline{w}))\rforal w\in F.
\eneq

For the second part of the statement, one assume{{s}} that $\Lambda(G)\subset U_0(B)/CU(B)$. Then $\Lambda(\pi_2(U(C)/CU(C)))$ is a torsion subgroup of $U_0(B)/CU(B)$. But $U_0(B)/CU(B)$ is torsion free, {{as $B=B_1\otimes V,$}}  by Lemma \ref{UCUdiv}, and hence $\theta=0$. {{Thus, in the first diagram, with five groups,  above,
$(\Lambda+\theta\circ \pi_2)|_{\overline{F}}=\Lambda|_{\overline{F}},$ or with multiplicative notation,
$\theta(\pi_2(\overline{w}))={\bar 1}.$ The second part of the lemma follows from \eqref{731-nn1}.}}

\end{proof}


\begin{lem}[see Lemma 7.3 of \cite{LinTAI}]\label{exp-length3}
Let $B\in\mathcal B_1$ be a separable simple \CA, and let $C$ be  {{as defined in \ref{bdbk-k1}.}}
Let $\mathcal{U}\subset \text{U}(B)$ be a finite subset, and  {{suppose that, with}}
$F$ the subgroup generated by $\mathcal{U},$  $\kappa_1^{ B}(\bar{F})$ is free abelian, where $\kappa_1^{ B}: \text{U}(B)/\text{CU}(B)\to K_1(B)$ is the quotient map. Suppose that $\alpha :K_1(C)\to K_1(B)$ is an injective homomorphism and $L: \bar{F}\to \text{U}(C)/\text{CU}(C)$ is an injective homomorphism with ${L(\bar{F}\cap \text{U}_0(B)/\text{CU}(B))\subset \text{U}_0(C)/\text{CU}(C) }$ such that $\pi_1\circ L$ is {{injective}} and $$\alpha\circ\kappa_1^{  C}\circ L(g)=\kappa_1^{{B}}(g)\quad\mbox{for all}\quad g\in\bar{F},$$ where $\kappa_1^{  C}: \text{U}(C)/\text{CU}(C)\to K_1(C)$ is the quotient map. Then there exists a homomorphism $\beta: \text{U}(C)/\text{CU}(C)\to \text{U}(B)/\text{CU}(B)$ with $\beta(\text{U}_0(C)/\text{CU}(C))\subset \text{U}_0(B)/\text{CU}(B)$ such that
$$\beta\circ L(f)=f\tforal f\in\bar{F}.$$
\end{lem}
{{Part of the statement can be summarized by the following commutative diagram:
\begin{displaymath}
\xymatrix{
& {\bar F} \ar[dr]^L \ar[dl]^{inclusion}& \\
U(B)/CU(B) \ar[d]^{\kappa_1^B}
& & U(C)/CU(C) \ar@{-->}[ll]^\bt \ar[d]^{\kappa_1^C}\\
K_1(B) &&  K_1(C) \ar[ll]^\af.
}
\end{displaymath}
}}

\begin{proof}
{{We will repeat}} the proof of
Lemma 7.3 of \cite{LinTAI}.

Let $G$ be the preimage of  {{the subgroup}} $\alpha\circ\kappa_1^{  C}(U(C)/CU(C))$ {{of
$U(B)/CU(B)$}} under $\kappa_1^{ B}$. So we have the short exact sequence
$$0\to U_0(B)/CU(B) \to G \to \alpha\circ\kappa_1^{  C}(U(C)/CU(C)) \to 0.$$
Since $U_0(B)/CU(B)$ is divisible {{(see \ref{UCUdiv})}}, there is an injective homomorphism
$${{\lambda}}: \alpha\circ\kappa_1^{  C}(U(C)/CU(C)) \to G$$
such that $\kappa_1^{ B}\circ\lambda(g)=g$
for any $g\in\alpha\circ\kappa_1^{  C}(U(C)/CU(C))$. Since $\alpha\circ\kappa_1^{  C}\circ L(f)=\kappa_1^{ B}(f)$ for any $f\in\overline{F}$, we have $\overline{F}\subset G$.  Moreover, note that
$$({{\lambda}}\circ\alpha\circ\kappa_1^{  C}\circ L(f))^{-1} f\in U_0(B)/CU(B)\rforal f\in \overline{F}.$$
Define $\psi: L(\overline{F})\to U_0(B)/CU(B)$ by
$$
\psi(x)={{\lambda}}\circ\alpha\circ\kappa_1^{  C}(x^{-1})L^{-1}(x)
$$
for $x\in L(\overline{F})$. Since $U_0(B)/CU(B)$ is divisible, there is a homomorphism $\tilde{\psi}: U(C)/CU(C)\to U_0(B)/CU(B)$ such that $\tilde{\psi}|_{L(\overline{F})}=\psi$. Now define
$$\beta(x)={{\lambda}}\circ\alpha\circ\kappa_1^{  C}(x)\tilde{\psi}(x)\rforal x\in U(C)/CU(C).$$
{{Then}} $\beta(L(f))=f$ for $f\in \overline{F}$.
\end{proof}

\begin{lem}\label{FDCU}
Let $A$ be a unital separable \CA\, such that {{the subgroup generated by}}
$\{\rho_A([p]): p\in A \,\,\, a  projection\}$  is dense in {{its}}
real linear span.
Then, for any finite dimensional \SCA\, $B\subset A$ {{(with $1_B=1_A$),}}
$U(B)\subset CU(A).$
\end{lem}
\begin{proof}
Let $u\in U(B).$  Since $B$ is finite dimensional,  $u=\exp(i 2\pi h)$ for some $h\in B_{s.a.}.$ We may write
$h=\sum_{i=1}^n \lambda_i p_i,$ where $\lambda_i\in (0, 1]$ and $\{p_1, p_2,...,p_n\}$ is a set of mutually orthogonal
projections.  {{Let $h(t)=\sum_{j=1}^n t\lambda_ip_i$ and  $u(t)=\exp( i h(t))$ ($t\in [0,1]$).
Then $\hat{h}$ {{is}} in the real linear span of  $\{\rho_A([p]): p\in A\,\,\, {\text{a  projection}}\},$
where $\hat{h}(\tau)=\tau(h(t))$ for all $\tau\in T(A)$
and $t\in [0,1].$}}
By the assumption, {{this implies that $\Delta^1(u(t))=\hat{h}\in \overline{\rho_A^1(K_0(A))}$ (notation  as in 2.8 and 2.10 of \cite{GLX-ER})}}.
 {{Then,}}  applying Proposition 3.6 {{(2)}}  of \cite{GLX-ER} {{ with $k=1,$}} one has
$u\in CU(A).$
\end{proof}

\begin{lem}[see Lemma 7.4 of \cite{LinTAI}]\label{exp-length2}
Let $B\cong A\otimes {{V}}$, where $A\in\mathcal B_0$ and ${{V}}$ is an infinite dimensional  UHF-algebra. Let
$C={{D}}\oplus C_{{0}}\oplus C_1$ be
as defined in \ref{bdbk-k1}.
Let $F$ be a group generated by a finite subset $\mathcal{U}\subset \text{U}(C)$ such that $(\pi_1)|_{\overline{F}}$ is
{{injective, where $\overline{F}$ is the image of $F$ in  $U(C)/CU(C).$}}
Suppose that $\alpha: \text{U}(C)/\text{CU}(C)\to \text{U}(B)/\text{CU}(B)$ is a homomorphism  such that $\alpha(\text{U}_0(C)/\text{CU}(C))\subset \text{U}_0(B)/\text{CU}(B)$. Then, for any $\epsilon>0$, there are $\sigma>0,$  $\delta>0,$ and {{a}} finite subset $\mathcal G\subset C$ satisfying the following {{condition}}: if $\phi=\phi_0\oplus\phi_1: C\to B$ {{ (by such a decomposition, we mean there is a projection $e_0$ such that $\phi_0: C \to e_0Be_0$ and $\phi_1: C\to (1_B-e_0)B(1_B-e_0)$)}} {{are}} unital  $\mathcal G$-${\bf \delta}$-multiplicative completely positive linear {{maps}} such that
\begin{enumerate}
\item  $\phi_0$ {sends the} identity of each {direct} summand of $C$ to a projection
{{and is non-zero on each (non-zero) direct summand of $D,$}}
\item $\mathcal G$ is sufficiently large and ${\bf \delta}$ is sufficiently small, depending only on $F$ and $C,$  that $\phi^\ddagger,$ {{$\phi_0^{\ddag},$ and $\phi_1^{\ddag}$}} {{(associated with a finite subset containing ${\cal U}$ and $\ep/4$---see \ref{DLddag})}} are well defined on a subgroup of $\text{U}(C)/\text{CU}(C)$ containing all of $\bar{F}$, $\pi_0(\bar{F})$, $\pi_1(\text{U}(C)/\text{CU}(C))$, and  $\pi_2(\text{U}(C)/\text{CU}(C))$,
\item $[\phi]|_{K_1(C)}=\alpha_*$ where $\af_*: K_1(C)\to K_1(B)$ is {{the}} map {{ induced by $\af,$}}
and  $[\phi_0]={{[\phi_{00}]}}$  {{in $KK(C,B)$  (note
that $K_i(C)$ is finitely generated---see the end of \ref{KLtriple}),}} where $\phi_{00}$  is a \hm\, from $C$ to $B$ which  has a
finite dimensional image, {{and}}
\item $\tau(\phi_0(1_C))<{{\sigma}}$ for all $\tau\in T(B)$ (assume $e_0=\phi_0(1_C)$),
\end{enumerate}
then there is a homomorphism $\Phi: C\to e_0Be_0$ such that

{\rm (i)}
$[\Phi]=[\phi_{00}]$ {{in $KK(C,B)$}} and

{\rm (ii)}  $\alpha(\bar{w})^{-1}(\Phi\oplus \phi_1)^\ddagger(\bar{w})=\bar{g_w}$ where $g_w\in \text{U}_0(B)$ and
${\rm{cel}}(g_w)<\epsilon$ for any $w\in\mathcal{U}$ {{(see \ref{DLddag})}}.
\end{lem}
\begin{proof}

The argument is almost exactly the same as that of Lemma 7.4 of \cite{LinTAI} {{(see also \cite{EGL-AH} and \cite{NT})}}. Since the source algebra and target algebra
in this lemma are different from the ones in 7.4 of \cite{LinTAI}, we will repeat {{most}} of the argument here.
{{We would like to point out that, since $K_1(D)$ is free {{abelian}}, so also is $\pi_1(\overline{F}).$
Also $\gamma(\overline{F})$ is free abelian, as $U(D)/U_0(D)$ is free abelian, where $\gamma: \overline{F}\to U(D)/CU(D)\to U(D)/U_0(D)$
is the {{composed}} map.
Since $\pi_1$ is injective on ${\overline{F}},$
$\overline{F}\cap (U(C_0)/CU(C_0))=\{{\bar 1}\}.$  }}
We will retain the notation of \ref{bdbk-k1}.

{{We rewrite $D=\bigoplus_{i=1}^s E_i,$
where $E_i=M_{r_i}(C(\T)).$}}  {\blue{Since $E_i$ are semiprojective,
we {{may}} assume $\phi_0|_D$ and $\phi_1|_D$ are homomorphisms,
{{choosing}} ${\cal G}$ sufficiently large and $\dt$ sufficiently small.}} {{Denote by $\Pi_i: D\to E_i$ the quotient map.}}
{\blue{Let $z'_i\in M_{r_i}(C(\T))$ be given by $z'_i(z)=\diag(z,1,\cdots,1)$ for any $z\in \T$, and $z_i\in D$ be defined by $z_i=(1_{r_1},\cdots,1_{r_{i-1}}, z'_i, 1_{r_{i+1}},\cdots, 1_{r_{s}})$.}}
{{Let $S_1$ be the subgroup of $\pi_1(U(C)/CU(C))\cong K_1(D)$ (recall
that we may view $\pi_1(U(C)/CU(C))$ as a subgroup of $U(C)/CU(C);$ see  \eqref{UCrank}) generated by $\{\overline{z_1'}, \overline{z_2'},...,\overline{z_s'}\}.$
Then $S_1=\pi_1(U(C)/CU(C))$ and
$\pi_1(\overline{F})\subset S_1.$
Thus one obtains a finitely generated  subgroup
$F_1$ of $U(C)$ such that $F\subset F_1$ and $\overline{F_1}\supset S_1.$
Moreover, $\pi_1|_{\overline{F_1}}$ is injective and
$\pi_1(\overline{F_1})= S_1\supset \pi_1(\overline{F}).$  There is a finite subset
${\cal U}_1\subset U(C)$ which generates the subgroup $F_1.$ }}

{{Therefore, to simplify matters, from now on, we may assume
that ${\cal U}={\cal U}_1.$
and $\overline{F}=\overline{F_1}.$
In particular, $\kappa_1^C(\overline{F_1})=K_1(D),$
where $\kappa_1^C: U(C)/CU(C)\to K_1(C)$ is the quotient map.}}
{{Let $G=\overline{F_1}\oplus \pi_2(U(C)/CU(C))=\overline{F_1}\oplus {\rm Tor}(K_1(C)).$ $G$ is a finitely generated
subgroup of $U(C)/CU(C)$ which contains $\overline{F}$ and $\pi_1(U(C)/CU(C))\oplus \pi_2(U(C)/CU(C)).$}}

{{Let $w\in U(C).$  We write $w=(w_1, w_2, w_3),$ where
$w_1\in D$ and $w_2\in C_0$ and $w_3\in C_1$ and $\pi_1({\bar w})=\pi_1({\bar w}_1)=
(\overline{z_1}^{k(1,w)}, \overline{z_2}^{k(2,w)},...,\overline{z_s}^{k(s,w)}),$ where $k(i,w)$ is an integer, $1\le i\le s.$
Then $\Pi_i^{\ddag}(\pi_1(\overline{w_1}))=\overline{z_i}^{k(i,w)}$ and $\Pi_i^{\ddag}(\overline{w_1})=\overline{z_i^{k(i,w)}g_{i,w}},$ where $g_{i,w}\in U_0(E_i),$ $i=1,2,...,s.$}} {\blue{Note that one can choose a function $g\in C(\T)$ of {{the}} form $g(e^{2\pi it})=e^{2\pi if(t)}$ such that $\det (g_{i,w}(z))=\det(g(z)\cdot 1_{E_i})$. That is, $g_{i,w}=g\cdot 1_{E_i}$ in $U_0(E_i)/CU(E_i)$. Therefore, we can write
\beq\label{May5-2019}
\Pi_i^{\ddag}(\overline{w_1})=\overline{z_i^{k(i,w)}g_{i,w}\cdot 1_{E_i}}.
\eneq
}} {{Later, in this proof, we will   view $g_{{\blue{i}},w}$ as a complex-valued continuous function on $\T.$}}

{{Let $l=\max\{{\rm cel}(g_{i,w}): 1\le i\le s, w\in {\cal U}\}.$ Choose an integer $n_0\ge 1$ such that
$(2+l){{/n_0}}<\ep/4\pi.$  Choose $0<\sigma<1/(n_0+1).$}}

{{Since $K_i(C)$ is finitely generated ($i=0,1$), we may assume, with sufficiently small
$\dt$ and large ${\cal G},$  {{that}}  $[\psi]$ gives an element of $KK(C,B)$ (see \ref{KLtriple}) and
\hm s $\psi^{\ddag}$ can be defined on  $G$ {{(defined above)}}
for any \morp\, $\psi: C\to A'$
(for any unital \CA\, $A'$)
which maps $G\cap (U_0(C)/CU(C))$ into
$U_0(A')/CU(A'),$ and ${\rm dist}(\psi^{\ddag}(\overline{u}), \overline{\la \psi(u)\ra})<\ep/16$ for all $u\in {\cal U}$
(see \ref{DLddag}). }}
{{Fix such $\dt$ and ${\cal G}.$ Let $\phi,$ $\phi_0,$ and $\phi_1$ be as given.
Define $\phi_0^{\ddag}$ and $\phi_1^{\ddag}$ as in \ref{DLddag}.
Let $\phi^{\ddag}$ be chosen to be $\phi_0^{\ddag}+\phi_1^{\ddag}.$
We note that $\phi^{\ddag}, \phi_0^\ddag,$ and $\phi_1^\ddag$ are defined on $G$ and map
$G\cap (U_0(C)/CU(C))$ into $U_0(B)/CU(B).$
Define  $L': G\to  U(B)/CU(B)$ by
$L'(g)=\phi_1^{\ddag}(g^{-1})$ for all $g\in G.$}}
{{Applying Lemma \ref{LinTAI-72} twice (once  for the case that $\af$ plays the role of $L$ and  then again
for the case that $L'$ plays the role of $L$),}}
one obtains  \hm s $\bt_1, \bt_2: U(D)/CU(D)\to U(B)/CU(B)$
with $\bt_i(U_0(D)/CU(D))\subset U_0(B)/CU(B)$ ($i=1,2$)
such that
\beq\label{74-e1}
\bt_1\circ \Pi^{\ddag}(\pi_1({\bar w}))=\af({\bar w})\theta_1(\pi_2({\bar w}))\andeqn
\bt_2\circ \Pi^{\ddag}(\pi_1({\bar w}))=\phi_1^{\ddag}({\bar w}^*)\theta_2(\pi_2({\bar w}))
\eneq
for all ${\bar w}\in {\bar F},$  {{where}} $\theta_1, \theta_2: \pi_2(U(C)/CU(C))\to U(B)/CU(B)$ {{are defined by}}
$\theta_{{1}}(g)=\af(g^{{-1}})$ and $\theta_2(g)={{\phi_1}}^{\ddag}(g)=L'(g^{-1})$ for all $g\in \pi_2(U(C)/CU(C)).$
{{Denote by $\kappa_1: U(B)/CU(B)\to K_1(B)$  the quotient map.
Then, by (3) (as $[\phi_{00}]|_{K_1(C)}=0$) for any $g\in \pi_2({\overline F}),$
\beq\label{n74-e2n}
(\kappa_1\circ \theta_1)(g)=\kappa_1(\af(g^{-1}))=(\kappa_1\circ  \phi^{\ddag}(g))^{-1}=(\kappa_1\circ \theta_2(g))^{-1}.
\eneq
}}
{{Thus,}}
{{by \eqref{n74-e2n},}}
\beq\label{74-e2}
\theta_1(g)\theta_2(g)\in {{{\rm ker}\kappa_1}}=U_0(B)/CU(B)\rforal g\in \pi_2({\bar F}).
\eneq
Since $\pi_2(U(C)/CU(C))$ is {{torsion}}  and $U_0(B)/CU(B)$ is torsion free (see Lemma \ref{UCUdiv}), we have
\beq\label{72-e3}
\theta_1(g)\theta_2(g)={{\bar{1}}}\rforal g\in \pi_2({\bar F}).
\eneq



Let $e_0=\phi_0(1_C).$ {{Write $e_0=e_0^d\oplus e_0^0\oplus e_0^1,$
where $e_0^d=\phi_0(1_D),$ $e_0^0=\phi_0(1_{C_0}),$ and $e_0^1=\phi_0(1_{C_1}).$
Let $e_{i,1}$ be a rank one projection in $E_i{\blue{=M_{r_i}(C(\T))}},$ $i=1,2,...,s.$
Let $\{q_{i,j,1}:1\le j\le r_i\}$ be a set of mutually orthogonal and mutually
equivalent projections in $B$ such that $[q_{i,j,1}]=[\phi_0(e_{j,1})]=[\phi_{00}(e_{j,1})],$ $1\le i\le {\blue{s}},$ $1\le j\le {\blue{r_i}}${\blue{, and such that $\sum_{j=1}^{r_i}q_{i,j,1}=\phi_0(1_{E_i}).$}}
Let $q_{i,1,1,x},q_{i,1,1,y}$ be two non-zero mutually orthogonal projections in $B$ such
that $q_{i,1,1,x}+q_{i,1,1,y}=q_{i,1,1}.$}}

{{Choose $x_i'\in U(q_{i,1,1,x}Bq_{i,1,1,x})$ and
$y_i'\in U(q_{i,1,1,y}Bq_{i,1,1,y})$ such that $\overline{x_i'}=\bt_1(\overline{z_i})$ and
$\overline{y_i'}=\bt_2(\overline{z_i})$
and $\overline{y_i}'=\bt_2(\overline{z_i}),$ $i=1,2,...,s.$  This is possible because of Theorem \ref{UCUiso}.
{\blue{Here we use the simplified notation $\overline{u}$ to denote $\overline{u\oplus (1-p)}$ for any unitary $u$ in the {{cut-down}} algebra $pBp$ of $B$.}}
Put $x_i=x_i'\oplus q_{i,1,1,y}$ and $y_i=y_i'\oplus q_{i,1,1,x},$ $i=1,2,...,s.$
Note that $x_iy_i=y_ix_i.$ }}

{\blue{Let ${\tilde \phi}^i: C(\T)\to q_{i,1,1}Bq_{i,1,1}$ be defined by
${\tilde \phi}^i(f)= f(x_iy_i)$ for any $f\in C(\T)$ (defined as the  continuous function $f$ acting on the unitary element $x_iy_i$ by functional calculus). Identifying $\phi_0(1_{E_i})B\phi_0(1_{E_i})$ with $M_{r_i}(q_{i,1,1}Bq_{i,1,1})$, we can define $\Phi_1:  D=\bigoplus_{i=1}^s E_i=\bigoplus_{i=1}^sM_{r_i}(C(\T))\to \bigoplus_{i=1}^s M_{r_i}(q_{i,1,1}Bq_{i,1,1})$ by $\Phi_1=\bigoplus_{i=1}^s {\tilde \phi}^i\otimes {\id}_{r_i}$.  It is clear that for any $f\in (C(\T)$, $$(\Phi_1)^{\ddag}(f\cdot 1_{E_i})=\overline{f(x_iy_i)\phi_0(1_{E_i})}.$$ }}

{{Then $(\Phi_1)_{*0}=((\phi_{00})_{*0})|_{K_0(D)}.$
Define $h: C\to B$ by $h(c)=\Phi_1(\Pi(c))\oplus \phi_1(c)$ for all $c\in C$ (recall that $\Pi: C\to D$
is the quotient map). Note that $[h]|_{K_1(C)}=[\Phi_1\circ \Pi]|_{K_1(C)}+
[\phi_1]|_{K_1(C)}.$
Define $h^{\ddag}=(\Phi_1\circ \Pi)^{\ddag}+\phi_1^{\ddag}.$}}
{\blue{Recall that any $w\in U(C)$ can be written as $w=(w_1,w_2,w_3)\in D\oplus C_0\oplus C_1$ such that $
\Pi_i^{\ddag}(\overline{w_1})$ is {{as}} described in (\ref{May5-2019}). Then
$$(\Phi_1)^{\ddag}(\overline{w_1})=
\prod_{i=1}^s\overline{x_i^{k(i,w)}\big(g_{i,w}\cdot 1_{E_i}\big)(x_iy_i)y_i^{k(i,w)}}.$$}}
We compute (viewing $g_{i,w}$ as a continuous function on $\T$), for all $w\in {\cal U}$
(also using \eqref{74-e1} and \eqref{72-e3}),
{{\beq
h^{\ddag}(\overline{w})&=&
(\prod_{i=1}^s\overline{x_i^{k(i,w)}{\blue{\big(g_{i,w}\cdot 1_{E_i}\big)}}y_i^{k(i,w)}})\phi_1^{\ddag}(\overline{w})\\
&=&\bt_1((\Pi)^{\ddag}(\pi_1(\overline{w})))\Phi_1^{\ddag}
(\bigoplus_{i=1}^s\overline{g_{i,w}{\blue{\cdot 1_{E_i}}}})
\bt_2((\Pi)^{\ddag}(\pi_1(\overline{w})))\phi_1^{\ddag}(\overline{w})\\\label{n74-n3}
&=&\af(\overline{w})\theta_1(\pi_2(\overline{w}))\theta_2(\pi_2(\overline{w}))
\Phi_1^{\ddag}(\bigoplus_{i=1}^s\overline{g_{i,w}{\blue{\cdot 1_{E_i}}}})
=
\af(\overline{w})\Phi_1^{\ddag}(\bigoplus_{i=1}^s\overline{g_{i,w}{\blue{\cdot 1_{E_i}}}}).
\eneq
Put $g_w'=\Phi_1(\bigoplus_{i=1}^s g_{i,w}{\blue{\cdot 1_{E_i}}})\oplus (1_B-\phi_0(1_C)).$
{{Recalling}} $B=A\otimes V,$ choose
mutually orthogonal and mutually equivalent
projections $\{p_1,p_2,...,p_{n_0}\}$ in $1_A\otimes V\subset B$ such that $\tau(p_i)\ge 1/(n_0+1)> \sigma$ and
$\sum_{i=1}^{n_0}\tau(p_i)<1-\sigma$ for all $\tau\in T(A).$
Since $\tau(\phi_0(1_C))<\sigma$
for all $\tau\in T(B),$   and $B$ has strict comparison (see \ref{Comparison}),
$(1-\phi_0(1_C))B(1-\phi_0(C))$ contains $n_0$ mutually orthogonal and mutually equivalent
projections which are equivalent to $\phi_0(1_C).$
By Lemma 6.4 of \cite{LinTAI},  there exists $g_{w,-}\in CU(B)$ such that
${\rm cel}(g_w'g_{w,-})<(l/n_0)\pi<\ep/2.$}}

{{Since $[\phi_0]=[\phi_{00}]$ in $KK(C, B),$  by (3), $[\phi_1]|_{K_1(C)}=\af_*.$
It follows from \eqref{n74-n3} that $h_{*1}=\af_*.$ Therefore
$(\Phi_1)_{*1}=0.$ Thus $[\Phi_1]=[\phi_{00}|_{D}]$
in $KK(D, B).$ Define $\Phi_2: C_0\oplus C_1\to (e_0^0\oplus e_0^1)B(e_0^0\oplus e_0^1)$ by
$\Phi_2=(\phi_{00})|_{C_0\oplus C_1}.$
Define $\Phi: C\to e_0Be_0$ by $\Phi((f,g))=\Phi_1(f)\oplus \Phi_2(g)$ for $f\in D$ and $g\in C_0\oplus C_1.$
Thus (i) holds.}}

{{For $w\in {\cal U},$ put $w''=\Phi_2(\Pi^{c}(w))\oplus (1_B-e_0^d),$
where $\Pi^c: C\to C_0\oplus C_1$ is the quotient map.
Note, since $B=A\otimes V,$
 by \ref{UHFdense},  the subgroup generated $\{\rho_B(p): p\in B\}$
is dense in the real linear span of $\{\rho_B(p): p\in B\}.$ By \ref{FDCU}, $\overline{w''}\in CU(B),$
as $\Phi_2$ has finite dimensional range.
 Thus, for $w\in {\cal U},$ let $g_w=g_w'g_{w-},$  Then, ${\rm cel}(g_w)<\ep.$
 Moreover,
 by \eqref{n74-n3},
 \beq
 \af(\overline{w})^{-1}(\Phi\oplus \phi_1)^{\ddag}(\bar{w})=
 \af(\overline{w})^{-1}\overline{h(w)}{\overline{w''}}={\overline{g_w'}}={\overline{g_w'g_{w-}}}=\overline{g_w}.
 \eneq
Thus (ii) holds.}}
\end{proof}


\begin{lem}[see Lemma 7.5 of \cite{LinTAI}]\label{exp-length1}
Let $B\cong A\otimes V$, where $A \in\mathcal B_1$ and $V$ is an {{infinite dimensional}}  UHF-algebra. Let $\mathcal{U}\subset \text{U}(B)$ be a finite subset  such that
$\kappa_1^{{B}}(\bar{F})$ is free abelian, where $F$ {{is}} the subgroup generated by $\mathcal{U},$ and
$\kappa_1^{{B}}: \text{U}(B)/\text{CU}(B)\to K_1(B)$ is the quotient map. Let $C={{D\oplus {{C_0\oplus C_1}}}}$
{{be as defined in \ref{bdbk-k1}}}
 and let {\blue{$\gamma: U(C)/CU(C)\to U(B)/CU(B)$ be a continuous \hm\, such that
 $\gamma_*:=\gamma|_{K_1(C)}$ is
 {{injective}} (viewing $K_1(C)\subset U(C)/CU(C),$ see \eqref{UCrank}).}}
 Suppose that $j, L: \overline{F}\to \text{U}(C)/\text{CU}(C)$ are two injective homomorphisms with $j(\overline{F\cap U_0(B)})$, $L(\overline{F\cap U_0(B)})\subset \text{U}_0(C)/\text{CU}(C)$ such that
 $\kappa_1^{{B}}\circ\gamma\circ L=\kappa_1^{{B}}\circ\gamma\circ j=\kappa_1^{{B}}|_{\bar{F}}.$

Then, for any ${\blue{\ep}}>0$, there exists $\delta>0$ satisfying the following condition:
Suppose that there is a \hm\, $\phi: C\to B$ such that $\phi^{\ddagger}=\gamma,$  $\phi_{*1}=\gamma|_{K_1(C)},$ and
$\phi=\phi_0\oplus\phi_1: C\to B$, where $\phi_0$ and $\phi_1$ are homomorphisms  such that
\begin{enumerate}
\item $\tau(\phi_0(1_C))<\delta$ for all $\tau\in T(B)$ and 
\item $[\phi_0]{{=[\phi_{00}]}},$ {{where $\phi_{00}: C\to B$ is a \hm\,}}
 with finite dimensional image {{and}}\end{enumerate}
 $\phi_{00}$ is not zero on each summand of $D,$
%

\vspace{0.05in}
\noindent
then there is a homomorphism $\psi: C\to e_0Be_0$ ($e_0=\phi_0(1_C)$) such that
\begin{enumerate}\setcounter{enumi}{3}
\item  $[\psi]=[\phi_0]$ in $KL(C, B)$ and
\item  $(\phi^\ddagger\circ j(\bar{w}))^{-1}(\psi\oplus\phi_1)^\ddagger(L(\bar{w}))=\overline{g_w}$ where $g_w\in \text{U}_0(B)$ and ${\rm{cel}}(g_w)<\epsilon$ for any $w\in\mathcal{U}$.
\end{enumerate}
\end{lem}
\begin{proof}
First, if $z\in \overline{F}$
 and
$\pi_1(L(z))={\bar 1},$ then  $L(z)=\pi_0(L(z))\oplus \pi_2({{L}}(z)).$
Therefore, $\kappa_1^B\circ\phi^{\ddag}(L(z))=\kappa_1^B\circ \phi^{\ddag}(\pi_2(L(z))$ is a torsion {{element}}
as $\pi_2(U(C_{{1}})/CU(C_{{1}}))$ is torsion. Since $\kappa_1^B(\overline{F})$ is free
abelian {{and $\kappa_1^{B}\circ\gamma\circ L =\kappa_1^{B}|_{\bar{F}}$, we have}}
$\kappa_1^B\circ \phi^{\ddag}(L(z))={\bar 1}.$  Since $\kappa_1^B\circ \phi^{\ddag}\circ L$
is {{injective,}}
 $L(z)={\bar 1}.$  This implies that $\pi_1|_{(L(\overline{F}))}$ is
 {{injective.}} Note also
 $\pi_1(L(\overline{F}))$ is free .   Exactly the same reason implies that $\pi_1|_{j(\overline{F})}$
 is
 {{injective}} and $\pi_1(j(\overline{F}))$ is free {{abelian}}.


Now we repeat the proof of Lemma 7.5 of \cite{LinTAI}.
Let $\kappa_1^C: U(C)/CU(C)\to K_1(C)$ be the quotient map and
$G=(\kappa_1^B)^{-1}(\gamma_*\circ \kappa_1^C(U(C)/CU(C))).$
Consider the short exact sequence
$$
0\to U_0(B)/CU(B)\to G\,\,{\stackrel{\kappa_1^B}{\longrightarrow}}\,\, \gamma_*\circ \kappa_1^C(U(C)/CU(C))\to 0.
$$
Since $U_0(B)/CU(B)$ is divisible, there {{exists}} an injective \hm\,\\
 $\lambda: \gamma_*\circ \kappa_1^C(U(C)/CU(C))\to G$
such that $\kappa_1^B\circ \lambda(g)=g$ for all $g\in \gamma_*\circ \kappa_1^C(U(C)/CU(C)).$
Since $\kappa_1^B\circ\gamma\circ L(f)=\kappa_1^B(f)=\kappa_1^B\circ\gamma\circ j(f)\subset \phi_{*1}(K_1(C))
=\phi_{*1}(\kappa_1^C(C(U(C)/CU(C)))$ for all $f\in \overline{F},$ one obtains
$\overline{F}\subset G.$
Note that
$
\gamma_*\circ \kappa_1^C=\kappa_1^B\circ \gamma.
$
Thus, for any $f\in \overline{F},$
\beq
\hspace{-0.6in}\kappa_1^B((\lambda\circ \gamma_*\circ \kappa_1^C\circ L(f))^{-1}(\gamma\circ j(f)))
&=&(\gamma_*\circ \kappa_1^C\circ L(f))^{-1}\kappa_1^B(\gamma\circ j(f))\\
&=&(\kappa_1^B\circ \gamma\circ L(f))^{-1}(\kappa_1^B\circ \gamma\circ j(f))=[1_B].
\eneq
It follows that
\beq
(\lambda\circ \gamma_*\circ \kappa_1^C\circ L(f))^{-1}(\phi^{\ddag}\circ j(f))
\in U_0(B)/CU(B)\rforal f\in \overline{F}.
\eneq
Define $\zeta: L(\overline{F})\to U_0(B)/CU(B)$ by
\beq
\zeta(x)=\left((\lambda\circ \gamma_*\circ \kappa_1^C(x)\right)^{-1}(\gamma\circ j\circ L^{-1}(x))\rforal x\in L(\overline{F}).
\eneq
Since $U_0(B)/CU(B)$ is divisible,  there exists ${\tilde \zeta}: U(C)/UC(C)\to U_0(B)/CU(B)$ such that
${\tilde \zeta}|_{L(\overline{F})}=\zeta.$ Define $\af: U(C)/CU(C)\to U(B)/CU(B)$ by
$\af(x)=(\lambda\circ \gamma_*\circ \kappa_1^C(x)){\tilde \zeta}(x)$ for all $x\in U(C)/CU(C).$
{{Note that, for all $x\in U(C)/CU(C),$
\beq
\af(\kappa_1^C(x))=(\lambda\circ \gamma_*\circ \kappa_1^C(x))\tilde{\zeta}(\kappa_1^C(x))
=\lambda\circ \gamma_*\circ \kappa_1^C(x)=\gamma_*\circ \kappa_1^C(x).
\eneq
In other words, $\af_*:=\af|_{K_1(C)}=\gamma_*.$}}
Note {{also}} that
\beq
\af(L(f))=\gamma\circ j(f)\rforal f\in \overline{F}.
\eneq

Let ${\cal U}'\subset U(C)$ be a finite subset  {{such that the subgroup}}
$F'\subset U(C)$
generated by ${\cal U}'$
{{satisfies}} $\overline{F'}=L(\overline{F}).$ Let $\ep>0.$
Choose $\dt=\sigma$ as in \ref{exp-length2}
associated ${\cal U}'$ (in place of ${\cal U}$), $F'$ (in place of $F$),  $\af$ {{as}} mentioned above, and $L(\overline{F})$ (in place  of $\overline{F}$).
Suppose that $\phi=\phi_0\oplus \phi_1$  as described.
Applying \ref{exp-length2},
one obtains a \hm\, $\psi: C\to e_0Be_0$ such that
$[\psi]=[\phi_{00}]$ in $KK(C,B)$ and
\beq
(\af(z))^{-1}(\psi\oplus \phi)(z)={\overline{f_z}} \andeqn {\rm cel}(f_z)<\ep,
\eneq
where $f_z\in U(B),$ for all $z\in L(\overline{F}).$
{{Then,}} for any $w\in {\cal U},$ let $g_w=f_{L(\overline{w})},$
{{so}} that
\beq
(\phi^{\ddag}\circ j({\overline{w}}))^{-1}(\psi\oplus \phi)(\overline{w})=
\af(L(\overline{w})^{-1}(\psi\oplus \phi)(L(\overline{w})=\overline{g_w}\andeqn {\rm cel}(g_w)<\ep.
\eneq
{{This shows that $\psi$ has the desired properties.}}
\end{proof}

\begin{rem}
{The} roles {{that}}  Lemma \ref{exp-length2} and Lemma \ref{exp-length1} play  in the proof of the
isomorphism theorem, {{Theorem \ref{IST0} below,}} are the same as those {{played by Lemma 7.4 and 7.5 of
\cite{LinTAI}}}  in the proof of Theorem 10.4 of \cite{LinTAI}.
\end{rem}

{The following statement is well known. For the reader's convenience, we include a proof.
\begin{lem}\label{indlim-inv}
Let   $(A_n, \phi_{n, n+1})$ be a unital inductive sequence of separable \CA s, and
{{consider the inductive limit}}
$A=\varinjlim A_n$. Assume that $A$ is amenable. Let $\mathcal F\subset A$ be a finite subset, and let ${\blue{\ep}}>0$. Then there is  an integer $m\ge 1$ and a unital completely positive linear map $\Psi:A\to A_m$ such that
$$
\|\phi_{m, \infty}\circ\Psi(f)-f\|<\epsilon
\rforal f\in\mathcal F.
$$
\end{lem}
\begin{proof}

Regard $A$ as the $C^*$-subalgebra of $\prod A_n/\bigoplus A_n$  {{consisting of}}
the equivalence classes of the sequences $(x_1, x_2, ..., x_n, ...)$ {{such that}}
{{that there is $N$ with}}
$x_{n+1}=\phi_n(x_n),\quad n={{N, N+1}}, ...\ .$
Since $A$ is amenable, by the Choi-Effros lifting theorem, there is a unital completely positive linear map $\Phi: A\to \prod A_n$ such that
$\pi\circ\Phi=\id_A,$
where $\pi: \prod A_n\to \prod A_n/\bigoplus A_n$ is the quotient map. In particular, this implies that
\begin{equation}\label{coincidelim}
\lim_{k\to\infty}\|\pi_k\circ\Phi(a)-a_k\|=0,
\end{equation}
if $a=\pi((a_1, a_2, ..., a_k, ...))\in A$.

Write $\mathcal F=\{f_1, f_2, ..., f_l\}$, and for each $f_i$, fix a representative
$$f_i=\pi((f_{i, 1}, f_{i, 2}, ..., f_{i, k}, ...)).$$ In particular
\begin{equation}\label{convlim}
\lim_{k\to\infty} \phi_{k, \infty}(f_{i, k})=f_i.
\end{equation}
Then, for each $f_i$, one has
\begin{eqnarray*}
&&\limsup_{k\to\infty}\|\phi_{k, \infty}\circ\pi_k\circ\Phi(f_i)-f_i\| \\
& = & \limsup_{k\to\infty} \|\phi_{k, \infty}\circ\pi_k\circ\Phi(f_i)-\phi_{k, \infty}(f_{i, k})\| \quad(\textrm{by \eqref{convlim}}) \\
&\leq&\limsup_{k\to\infty}\|\pi_k\circ\Phi(f_i)-f_{i, k}\|
= 0\quad \textrm{(by \eqref{coincidelim}}).
\end{eqnarray*}
There then {{exists}} $m\in \N$ such that
$$\|\phi_{m, \infty}\circ\pi_m\circ\Phi(f_i)-f_i\|<\epsilon,\quad 1\leq i\leq l.$$
{{This shows that}} the unital completely positive linear map $\Psi:=\pi_m\circ\Phi$
satisfies the {{condition of the}} lemma.
\end{proof}
}


\begin{thm}\label{IST0}
Let $A_1\in \mathcal B_{{0}}$ be a unital separable simple {{amenable}} \CA\, satisfying the UCT, and {{let}} $A=A_1\otimes U$ for
{{some}}  UHF-algebra {$U$} of infinite type. Let $C$ be a \CA\,  {{constructed as}} in Theorem \ref{RangT}
{{(denoted by $A$ there)}}
{{such that}} ${\rm Ell}(C)\cong {\rm Ell}(A).$ Then there is an isomorphism $\phi: C\to A$ which carries the identification of ${\rm Ell}(C) \cong {\rm Ell}(A)$.
\end{thm}

\begin{proof}

Let  $C=\varinjlim (C_n, \iota_{n,n+1})$
be as  {{constructed}} in \ref{RangT} {{(as { {$A=\varinjlim (A_n, \psi_{n, n+1})$ with $\iota_{n,n+1}$ in place of $\psi_{n, n+1}$}} there),
where
$\iota_{n, n+1}$ is injective{\blue{,}} unital, {{and has the decomposition $\iota_{n,n+1}=\iota^{(0)}_{n,n+1}\oplus \iota^{(1)}_{n,n+1}$}}.  Put $\iota_n=\iota_{n, \infty},\, n=1,2,....$}}  {{By \ref{CCRangT}, we may assume
that $C\in  B_{u0}.$ Note also $A\in {\cal B}_{u0}.$ }}
{{The proof will use the fact that both $A$ and $C$  have stable rank 1 (see {{Theorem}} \ref{B1stablerk}) without further {{notice.}}
We will also use $\kappa_1^B: U(B)/CU(B)\to K_1(B)$
for the quotient map for any unital \CA\, $B$ with the property $U(B)/U_0(B)=K_1(B).$
We will  fix splitting maps: $J_c^C: K_1(C)\to U(C)/CU(C)$ and $J_c^A: K_1(A)\to U(A)/CU(A)$
such that $\kappa_1^C\circ J_c^C={\rm id}_{K_1(C)}$ and $\kappa_1^A\circ J_c^A={\rm id}_{K_1(A)}$ (see \ref{Dcu}).}}

Let
 $\gamma: T(A)\to T(C)$ be {{as}} given by the isomorphism
${\rm Ell}(C)\cong {\rm Ell}(A),$  and choose
$\alpha\in KL(C, A)$ {with $\alpha^{-1}\in KL(A, C)$}
lifting the $K_i$-group isomorphisms from the isomorphism ${\rm Ell}(C)\cong {\rm Ell}(A)$
{{(see \ref{DEll})}}.

Let $\mathcal G_1\subset \mathcal G_2\subset \cdots \subset C$ and $\mathcal F_1\subset \mathcal F_2\subset\cdots \subset A$ be increasing sequences of finite subsets with dense union. Let ${{1/2>}}\ep_1>\ep_2> \cdots >0$ be a decreasing sequence of positive numbers with finite sum.
{{Let ${\cal P}_{c,n}\subset \underline{K}(C)$ be finite subsets such that
${\cal P}_{c,n}\subset {\cal P}_{c,n+1}$ ($n\ge 1$) and $\bigcup_{n=1}^{\infty}{\cal P}_{c,n}=\underline{K}(C),$ let
${\cal Q}_{c,n}\subset \underline{K}(A)$ be finite subsets such that ${\cal Q}_{c,n}\subset {\cal Q}_{c,n+1}$ ($n\ge 1$) and
$\bigcup_{n=1}^{\infty}{\cal Q}_{c,n}=\underline{K}(A),$ and let ${\cal H}(c,n)\subset C_{s.a.}$
and
${\cal H}(a,n)\subset A_{s.a.}$ be finite subsets such that ${\cal H}(c,n)\subset {\cal H}(c,n+1)$ ($n\ge 1$) and $\bigcup_{n=1}^{\infty} {\cal H}(c,n)$ is dense in $C_{s.a.},$  and ${\cal H}(a,n)\subset {\cal H}(a, n+1)$ ($n\ge 1$)
and $\bigcup_{n=1}^{\infty}{\cal H}(a,n)$ is dense in $A_{s.a.}.$}}
{{
We may assume that  ${\cal P}_{c,n}$ is in the image
of ${{\underline{K}(C_{m(n)})}}$ under $\iota_{m(n){{*}}}$ for some $m(n)\ge 1.$}}

We will repeatedly apply {{the part (a) of Theorem \ref{MUN1}.}}
 Let $\delta_c^{(1)}>0$ (in place of $\delta$), $\mathcal G^{(1)}_c\subset C$ ({{${\cal G}_c^{(1)}$}} in place of $\mathcal G$ {{ and $C$ in place of $A$}}), $\sigma^{(1)}_{c, 1}, \sigma^{(1)}_{c, 2}>0$ (in place of $\sigma_1$ and $\sigma_2$),
 $\mathcal P^{(1)}_c\subset\underline{K}(C)$ (in place of $\mathcal P$), ${{\overline{\mathcal U^{(1)}_c}\subset U(C)/CU(C))}}$ (in place of $\mathcal U$),
and $\mathcal H^{(1)}_c\subset C_{s.a}$ (in place of ${{\mathcal H}} $) be as
{{provided}} by {{the part(a)}} of  \ref{MUN1} for $C$ (in place of $A$), $\ep_1$ (in  place of $\ep$), {{and}}  $\mathcal G_1$ (in  place of $\mathcal F$).
{{ Here, ${\cal U}^{(1)}_c$ is a finite subset of $U(C).$
As in  Remark \ref{ReMUN1}, we may assume that
$\overline{\mathcal U^{(1)}_c}\subset J_c^C(K_1(C)).$}}
{{We may also assume that the image of $\mathcal U^{(1)}_c$ in $K_1(C)$ is contained in ${\cal P}^{(1)}_c.$}}
As in  {{Remark  \ref{ReMUN1}}},
we may {{also}} assume that ${{{\cal U}_c^{(1)}}}$ is in the image of $U(C_n)$ for {{all $n\ge n_0$}}  for some
large  ${{n_0}}\ge m(1)$ {{under the map $\iota_n.$}}
 {{We may further assume that ${\cal P}_{c,1}\subset {\cal P}_c^{(1)}$
and ${\cal H}(c,1)\subset {\cal H}_c^{(1)}.$}}

Denote by $F_c^{{(1)}}\subset U(C)$ the subgroup generated by $\mathcal U^{(1)}_c$.
{{We may}} write $\overline{F_c}^{{(1)}}=(\overline{F_{c}^{(1)}})_0 \oplus \mathrm{Tor}(\overline{F_c^{{(1)}}})$
where $(\overline{F_{c}^{(1)}})_0$ is torsion free, {{as $F_c^{(1)}$ is finitely generated.
Note $\overline{F_c}^{{(1)}},$ $(\overline{F_{c}^{(1)}})_0,$ and $\mathrm{Tor}(\overline{F_c^{{(1)}}})$
are in $J_c(K_1(C)).$ }}
Choosing a smaller $\sigma_{c,2}^{(1)}$, {{we}} may assume that
$$\mathcal U^{(1)}_c=\mathcal U^{(1)}_{c, 0}\sqcup\mathcal U^{(1)}_{c, 1},$$
where $\overline{\mathcal U^{(1)}_{c, 0}}$ generates $(\overline{F_{c}^{(1)}})_0$ and $\overline{\mathcal U^{(1)}_{c, 1}}$ generates $\mathrm{Tor}(\overline{F_c^{(1)}})$ --namely, we can choose $\mathcal U^{(1)}_{c, 0}$ and $\mathcal U^{(1)}_{c, 1}$ so that $\mathcal U^{(1)}_c\subset\mathcal U^{(1)}_{c, 0}\cdot\mathcal U^{(1)}_{c, 1}$; then, choosing smaller $\sigma_{c,2}^{(1)},$ one can replace $\mathcal U^{(1)}_c$ by  $\mathcal U^{(1)}_{c, 0}\sqcup\mathcal U^{(1)}_{c, 1} .$ Note that for each $u\in \mathcal U^{(1)}_{c, 1}$, one has  $u^k \in CU(C)$, where $k$ is the order of $\overline{u}$.

{{Let  a finite subset ${\cal G}_{uc}^{(1)}\subset C$
and $\dt_{uc}^{(1)}>0$ satisfy
the following condition:
for any ${\cal G}_{uc}^{(1)}$-$\dt_{c,u}^{(1)}$-multiplicative unital completely positive linear map $L': C\to A'$ (for any unital \CA\, $A'$
with
$K_1(A')=U(A')/U_0(A')$),
$(L')^{\ddag}$ can be defined as  a \hm\, on $\overline{F_c^{(1)}},$
${\rm dist}((L')^{\ddag}(\overline{u}), \overline{\la L'(u)\ra})<\sigma_{c,2}^{(1)}/8 $ for all $u\in {\cal U}_c^{(1)},$
and $\kappa_1^{A'}\circ (L')^{\ddag}(\overline{u})=[L']\circ \kappa_1^C([u])$ for all $u\in {\cal U}_c^{(1)}$
(see \ref{DLddag}). {{Since $\overline{F_c^{(1)}}\subset J_c(K_1(C))$, and   $J_c(K_1(C))\cap U_0(C)/CU(C)$ only contains
{{the}}  unit  $\overline{1_C},$ and since  $(L')^{\ddag}$ is a homomorphism on $\overline{F_c^{(1)}},$ we have $(L')^{\ddag}(\overline{F_c^{(1)}}\cap U_0({{C}})/CU({{C}}))=\overline{1_{A'}}\in U_0(A')/CU(A').$}} Moreover, {{for any $u\in {\cal U}_{c,1}^{(1)}$ with $\overline{u}$ {{of}} order $k$, as $u^k\in CU(C)$,}} we may
assume (see \ref{DLddag}) that
\beq\label{21-18713-f-1}
{\rm dist}(\la L'(u^k)\ra, CU(A'))<\sigma_{c,2}^{(1)}/8.
\eneq
Put ${\cal G}_c^{(1)+}={\cal G}_c^{(1)}\cup {\cal G}_{uc}^{(1)}$ and
$\dt_c^{(1)+}={{\frac12}}\min\{\dt_c^{(1)}, \dt_{uc}^{(1)}\}.$}}

{{Recall that both $A$ and $C$ are in ${\cal B}_{u0}.$}}
By Theorem \ref{MEST},
there is a {{${\cal G}^{(1)+}_c$-$\dt^{(1)+}_c$}}
-multiplicative \morp\, $L_1: C\to A$ such that
\beq\label{eq-kk-001}
&&[L_1]|_{\mathcal P^{(1)}_c}=\alpha|_{\mathcal{P}^{(1)}_c}\andeqn\\\label{eq-tr-001}
&& |\tau\circ L_1(f)-\gamma(\tau)(f)|
<\sigma^{(1)}_{c, 1}/8\rforal f\in\mathcal H^{(1)}_c\rforal \tau\in T(A).
 \eneq
{{By choosing  ${\cal G}^{(1)+}_c$ large enough and $\dt^{(1)+}_c$ small enough, o}}ne may assume that $L_1^\ddagger$ is  a \hm\, defined  on $(\overline{F_c^{(1)}})$ {{(see \ref{DLddag}).
  We {{may further }}  assume
 that, ${\rm dist}(L_1^{\ddag}(\overline{u}), \overline{\la L_1(u)\ra})<\sigma_{c, 2}^{(1)}/8$ for all
 $u\in  {\cal U}_c^{(1)}$ and  $\kappa_1^A\circ L_1^{\ddag}=[L_1]\circ \kappa_1^C$ on $(\overline{F_c^{(1)}}).$}}
 {{Since $\af|_{K_1(C)}$ is an isomorphism, we may assume that $L_1^{\ddag}|_{\overline{F_c^{(1)}}}$ is injective
 as $\kappa_1^C$ is injective on $J_c(K_1(C)).$ In particular,}} $L_1^{\ddag}|_{(\overline{F_c^{(1)}})_0}$ is injective.
  Moreover, {{if $u\in {\cal U}_{c,1}^{(1)}$ and}}
  $\overline{u} $  is of order $k,$ one may also assume  {{(see \eqref{21-18713-f-1}) that}}
 $$\mathrm{dist}(L_1(u^k), CU(A))<\sigma_{c, 2}^{(1)}/{{8}}.$$

{{Applying}} {{the part (a) of }} Theorem \ref{MUN1} 
{{a}} second time, let $\delta^{(1)}_a>0$ (in  place of $\delta$), $\mathcal G^{(1)}_a\subset {{A}}$ (in  place of $\mathcal G$), $\sigma^{(1)}_{a, 1}, \sigma^{(1)}_{a, 2}>0$ (in  place of $\sigma_1$ and $\sigma_2$), $\mathcal P^{(1)}_a\subset\underline{K}(A)$ (in  place of $\mathcal P$), $\overline{\mathcal U^{(1)}_a}\subset U(A)/CU(A)$ (in place of $\mathcal U$), and $\mathcal H^{(1)}_a\subset A_{s.a}$ (in  place of $\mathcal {{H}}$) be as {{provided}} by {{the part (a) of}} Theorem \ref{MUN1} for
$A$ (in place of $A$), $\ep_2$ (in  place of $\ep$), and $\mathcal F_1$ (in place of $\mathcal F$). 
{{We may assume
that $\dt_a^{(1)}<\dt_c^{(1)+}/2,$ $L_1({\cal G}_c^{(1)+})\subset {\cal G}_a^{(1)},$
\beq\label{May31-2019}
{{J_c^A\circ \kappa_1^A(L_1^{\ddag}(\overline{{\cal U}_c^{(1)}}))\subset \overline{{\cal U}^{(1)}_a}}},
\eneq
and $\kappa_1^A(\overline{{\cal U}_a^{(1)}})\cup [L_1]({\cal P}_c^{(1)})\subset {\cal P}_a^{(1)}.$
Here we also assume that ${\cal U}_a^{(1)}\subset U(A)$ is a finite subset.
As in Remark \ref{ReMUN1}, we may further assume that
$\overline{\mathcal U^{(1)}_a}\subset J_c^A(K_1(A)).$
Moreover, we may assume that ${\cal P}_{a,1}\subset {\cal P}_c^{(2)}$ and
${\cal H}(a,1)\subset {\cal H}_a^{(1)}.$}}

Denote by $F_a^{{(1)}}\subset U(A)$ the subgroup generated by $\mathcal U_a^{(1)}$. Since ${{\mathcal U_a^{(1)}}}$ is finite, we can write $\overline{F_a^{{(1)}}}=(\overline{F_a^{(1)}})_0 \oplus \mathrm{Tor}(\overline{F_a^{(1)}})$, where $(\overline{F_{a}})_0$ is torsion free. Fix this decomposition. Without loss of generality (choosing a smaller $\sigma_{{{a}},2}^{(1)}$),
one may assume that
$$\mathcal U^{(1)}_a=\mathcal U^{(1)}_{a, 0}\sqcup\mathcal U^{(1)}_{a, 1},$$
where ${ {\mathcal U}}^{(1)}_{a, 0}$ generates $(\overline{F_a^{(1)}})_0$ and $\mathcal U^{(1)}_{a, 1}$ generates $\mathrm{Tor}(\overline{F_a^{(1)}})$ {{(then the condition (\ref{May31-2019}) should be changed to
{{the condition}} that $J_c^A\circ \kappa_1^A(L_1^{\ddag}(\overline{{\cal U}_c^{(1)}}))$ is in the subgroup generated by $ \overline{{\cal U}^{(1)}_a}$)}}.
{{Enlarging}} $\mathcal P^{(1)}_a$, we {{may}}  assume $\mathcal P^{(1)}_a\supset \kappa_1^A((\overline{F_a^{(1)}})_0)$ (in $ K_1(A)$).
Note that for each $u\in \mathcal U^{(1)}_{a, 1}$, $u^k \in CU(A)$, where $k$ is the order of $\overline{u}$.

{{Let ${\cal U}_{c,a}^{(1)'}\subset U(A)$ be a finite subset such that
$\overline{{\cal U}_{c,a}^{(1)'}}=L_1^{\ddag}(\overline{{\cal U}_c^{(1)}})$ and \\
${\cal U}_a^{(1)'}={\cal U}_a^{(1)}\cup {\cal U}_{c,a}^{(1)'}\cup\{\la L_1(u)\ra: u\in {\cal U}_c^{(1)}\}.$
Let $F_{a,1}$ be the subgroup of $U(A)$ generated by ${\cal U}_a^{(1)'}.$
Since $\kappa_1^A(L_1^{\ddag}(\overline{{\cal U}_c^{(1)}}))=\{\kappa_1^A(\overline{\la L_1(u)\ra}): u\in {\cal U}_c^{(1)}\}$ and
 since $J_c^A\circ \kappa_1^A(L_1^{\ddag}(\overline{{\cal U}_c^{(1)}}))$ { {is in the subgroup generated by}} $ \overline{{\cal U}_a^{(1)}},$
$\overline{F_{a,1}}=\overline{F_a^{(1)}}+(U_0(A)/CU(A))\cap{\overline{F_{a,1}}}.$ }}

{{Let  a finite subset ${\cal G}_{ua}^{(1)}\subset A$ and $\dt_{ua}^{(1)}>0$ {{satisfy}} the following
{{condition}}:
for any ${\cal G}_{ua}^{(1)}$-$\dt_{c,a}^{(1)}$-multiplicative \morp\, $L': A\to B'$ (for any unital \CA\, $B'$
with
$K_1(B')=U(B')/U_0(B')$),
$(L')^{\ddag}$ can be chosen to  be a \hm\, on $\overline{F_{a,1}},$
${\rm dist}((L')^{\ddag}(\overline{u}), \overline{\la L'(u)\ra})<\sigma_{a,2}^{(1)}/8 $ for all $u\in {\cal U}_a^{(1)'},$
$(L')^{\ddag}(\overline{F_{a,1}}\cap U_0({{A}})/CU({{A}}))\subset U_0(B')/CU(B'),$
and $\kappa_1^{B'}\circ (L')^{\ddag}(\overline{u})=[L']\circ \kappa_1^C([u])$ for all $u\in {\cal U}_a^{(1)'}$
(see \ref{DLddag}). We  assume that, for  $u\in {\cal U}_{a,1}^{(1)},$
\beq\label{21-18713-f-2}
{\rm dist} (\la L'(u^k)\ra, CU(B'))<\sigma_{a,2}^{(1)}/8,
\eneq
where $k$ is the order of $\overline{u}$ (see \ref{DLddag}).
Put ${\cal G}_a^{(1)+}={\cal G}_a^{(1)}\cup {\cal G}_{ua}^{(1)}$ and
$\dt_a^{(1)+}=\min\{\dt_a^{(1)}, \dt_{ua}^{(1)}\}.$}}

By Theorem \ref{MEST} and {{the}} amenability of $C,$  there are a finite subset ${\cal G}_a'\supset \mathcal G^{{(1)+}}_a,$ a positive number $\dt'_a<\dt_a^{{(1)+}},$ a sufficiently large integer $n\ge {{n_0}},$ and
a $\mathcal G_a' $-$\delta_a'$-multiplicative map $\Phi'_1: A\to C_n$ such that, {{with $\sigma_{a,c,1}^{(1)}=\min\{\sigma^{(1)}_{a, 1}, \sigma^{(1)}_{c, 1}\}/6,$}}
\beq\label{eq-kk-002}
&&\hspace{-0.5in}[\iota_n\circ \Phi'_1]_{\mathcal P^{(1)}_a\cup [L_1](\mathcal P^{(1)}_c)}=\alpha^{-1}|_{\mathcal P^{(1)}_a \cup [L_1](\mathcal P^{(1)}_c)}\andeqn\\
\label{eq-tr-002}
&&\hspace{-0.5in}|\tau\circ \iota_n\circ \Phi'_1(f)-\gamma^{-1}(\tau)(f)|
<\sigma_{a,c,1}^{(1)}/8
\rforal f\in\mathcal H^{(1)}_a\cup L_1(\mathcal H^{(1)}_c)\,\,\,{\rm and}\,\,\tau\in T(C).
\eneq
Moreover, one  may assume that $\Phi'_1\circ L_1$ is $\mathcal G^{(1)}_c$-$\delta^{(1)}_c$-multiplicative, $(\Phi_1')^\ddagger$ {{is a homomorphism}}  defined on {{$\overline{F_{a,1}},$
${\rm dist}((\Phi_1')^{\ddag}(\overline{u}), \overline{\la \Phi_1'(u)\ra})<\sigma_{a,c,1}^{(1)}/8,$
and $\kappa_1^{C_n}\circ (\Phi_1')^\ddag(\overline{u})=[\Phi_1'](\kappa_1^A(\overline{u}))$
for all $u\in {\cal U}_a^{(1)'}.$}}
{{Then, by \eqref{eq-kk-002}, since $\af^{-1}|_{K_1(A)}$
is injective, $(\iota_n)_{*1}\circ [\Phi_1']$
is  injective on $\kappa_1^A(\overline{F_{a,1}}),$  which implies that $[\Phi_1']$ is injective
on  $\kappa_1^A(\overline{F_{a,1}}).$ Since $\kappa_1^A$ is injective on $J_c^A(K_1(A)),$
we conclude that $\kappa_1^{C_n}\circ (\Phi_1')^{\ddag}$ is injective on $\overline{F_a^{(1)}},$
which implies that}} $(\Phi_1')^{\ddag}$ is injective on $\overline{F_a^{(1)}}.$
Write $(\Phi_1'\circ L_1)^{\ddag}{{:=}}(\Phi_1')^{\ddag}\circ L_1^{\ddag}.$
Then
$(\Phi'_1\circ L_1)^\ddagger$ is defined on $\overline{F_c^{(1)}}${{, since $L_1^{\ddag}(\overline{F_c^{(1)}})\subset \overline{F_{a,1}}$.}}

{{If $u\in F_c^{(1)},$    then $\overline{u}\in \overline{F_c^{(1)}}\subset J_c^C(K_1(C)).$
Suppose that $\overline{u}\not=0.$ {{Then}} $[u]\not=0,$ since $\kappa_1^C\circ J_c^C={\rm id}_{K_1(C)}.$
Then $\kappa_1^A(L_1^{\ddag}(\overline{u}))=
[L_1]\circ \kappa_1^C(\overline{u})=\af([u])\not=0.$ It follows that $J_c^A\circ \kappa_1^A(L_1^{\ddag}(\overline{u}))\not=0.$
Then
$(\Phi_1')^{\ddag}(J_c^A\circ \kappa_1^A(L_1^{\ddag}(\overline{u})))\not=0$
since $J_c^A\circ \kappa_1^A(L_1^{\ddag}(\overline{u}))\subset \overline{F_a^{(1)}}$ and
$(\Phi_1')^{\ddag}$ is injective on $\overline{F_a^{(1)}}.$}} {{ Moreover,
since $\kappa_1^C(\overline{{\cal U}_c^{(1)}})\subset {\cal P}_c^{(1)},$ by \eqref{eq-kk-002},
for any $u\,{{\in}}\,{\cal U}_c^{(1)},$
\beq\nonumber
\kappa_1^{C_n}((\Phi_1')^{\ddag}(J_c^A\circ \kappa_1^A(L_1^{\ddag}(\overline{u}))))
=[(\Phi_1')](\kappa_1^A(J_c^A\circ \kappa_1^A(L_1^{\ddag}(\overline{u}))))\\\label{21-n18-n1}
=[(\Phi_1')]( \kappa_1^A(L_1^{\ddag}(\overline{u})))=[\Phi_1']( [L_1](\kappa_1^C(\overline{u})))
\not=0.
\eneq}}
{{Put $z=L_1^{\ddag}(\overline{u})-J_c^A\circ \kappa_1^A(L_1^{\ddag}(\overline{u}).$ Then $z\in U_0(A)/CU(A).$
It follows that $(\Phi_1')^{\ddag}(z)\in U_0(C_n)/CU(C_n).$
If $(\Phi_1')^{\ddag}(L_1^{\ddag}(\overline{u}))=0,$ then
$\kappa_1^{C_n}((\Phi_1')^{\ddag}(J_c^A\circ \kappa_1^A(L_1^{\ddag}(\overline{u})))){{=
\kappa_1^{C_n}((\Phi_1')^{\ddag}(L_1^{\ddag}(\overline{u})))-
\kappa_1^{C_n}((\Phi_1')^{\ddag}(z))}} =0,$ which contradicts \eqref{21-n18-n1}. This implies
that $(\Phi_1')^{\ddag}(L_1^{\ddag}(\overline{u}))\not=0.$ In other words,}}
$(\Phi_1'\circ L_1)^{\ddag}$ is injective on $\overline{F_c^{(1)}}.$

Furthermore,  one may also assume {{(see \eqref{21-18713-f-2}) that}}, for $u\in \mathcal U_{a, 1}^{(1)},$
 \begin{equation}\label{almost-cd-01}
 \mathrm{dist}(\iota_n\circ \Phi_1'(u^k), CU(C))<\sigma_{a, 2}^{(1)}/8,
 \end{equation}
 where $k$ is the order of $\overline{u}$, and  {{(by \eqref{21-18713-f-1}, as $L_1({\cal G}_c^{(1)+})\subset {\cal G}_a^{(1)}$ and $\dt_a^{(1)+}<\dt_c^{(1)+}/2)$}}
 \begin{equation}\label{almost-cd-02}
 \mathrm{dist}((\iota_n\circ \Phi_1'\circ L_1)(v^{k'}), CU(A))<\sigma_{c, 2}^{(1)}/8,
 \end{equation}
 if $v\in \mathcal U_{c, 1}^{(1)}$ and  $k'$ is the order of $\overline{v}$.
It then follows from \eqref{eq-kk-001} and \eqref{eq-kk-002} that
\begin{equation}\label{comp-001}
[\iota_n\circ \Phi'_1\circ L_1]|_{\mathcal P^{(1)}_c}=[\id]|_{\mathcal P^{(1)}_c};
\end{equation}
and it follows from \eqref{eq-tr-001} and \eqref{eq-tr-002}  that
\begin{equation}\label{comp-002}
|\tau\circ\iota_n\circ \Phi'_1\circ L_1(f)-\tau(f)|<2\sigma^{(1)}_{c, 1}/3\rforal f\in\mathcal H^{(1)}_{c} \rforal \tau\in T(C).
\end{equation}

Recall that $(\overline{F_c^{(1)}})_0\subset U(C)/CU(C)$ is  {{a free abelian}} subgroup, generated by $\overline{\mathcal U^{(1)}_{c, 0}}$.  Since we have assumed that $\overline{\mathcal U^{(1)}_{c, 0}}$ is in the image of $U(C_n)/CU(C_n),$
there is an injective homomorphism
$j: (\overline{F_c^{(1)}})_0 \to U(C_n)/CU(C_n)$
such that
\begin{equation}\label{june-29-n-001}
\iota_{n}^\ddagger\circ j=\id|_{\overline{(F_c^{(1)})_0}}.
\end{equation}
Moreover, by (\ref{comp-001}),  $$\kappa_1^C\circ\iota_n^\ddagger\circ(\Phi_1'\circ L_1)^\ddagger|_{\overline{(F_c^{(1)})_0}}={{(\iota_n)_{*1}\circ [\Phi_1'\circ L_1]\circ \kappa_1^C|_{\overline{(F_c^{(1)})_0}}}}=\kappa_1^C|_{\overline{(F_c^{(1)})_0}}=
\kappa_1^C\circ\iota_n^\ddagger\circ j.$$

Let $\delta$ be the constant of Lemma \ref{exp-length1} with respect to $C_n$ (in place of $C$), $C$ (in place of $B$), $\sigma_{c, 2}^{(1)}/2$ (in place of $\ep$), {\blue$\iota_n^{\ddagger}$ (in place of $\gamma$),}
$j,$ and $(\Phi_1'\circ L_1)^\ddagger |_{\overline{(F_c^{(1)})_0}}$ (in place of $L$). {{For any $n'>n$,
{{we have}}
\beq\label{May11-2019}
\iota_n=(\iota_{n',\infty}\circ \iota^{(0)}_{n',n'+1}\circ \iota_{n,n'})
\oplus(\iota_{n',\infty}\circ \iota^{(1)}_{n',n'+1}\circ \iota_{n,n'}).
\eneq}}
{{Denote $\iota_{n',\infty}\circ \iota^{(0)}_{n',n'+1}\circ \iota_{n,n'}$ by $\imath_n^{(0)}$ and $\iota_{n',\infty}\circ \iota^{(1)}_{n',n'+1}\circ \iota_{n,n'}$ by $\imath_n^{(1)}.$}}  {{By}}
\eqref{RRrangT-1}, {{if $n'$ {{is}} large enough (in particular, depending on $\dt$ above), the}} decomposition ${ \iota_n=\imath_n^{(0)}\oplus\imath_n^{(1)}}$ {{satisfies}}
\begin{enumerate}
\item $\tau({\blue{\imath}}_n^{(0)}(1_{C_n}))<\min\{\delta, \sigma^{(1)}_{c, 1}/6, \sigma_{a,1}^{(1)}/6\}$ for all $\tau\in T(C)$, and
\item ${\blue{\imath}}_n^{(0)}$ has finite dimensional range, {{and
is non-zero on each direct summand of $C_n.$}}
\end{enumerate}
Then, by Lemma \ref{exp-length1}, there is a homomorphism $h: C_n\to e_0Ce_0$, where $e_0={\blue{\imath}}_n^{(0)}(1_{C_n})$, such that
\begin{enumerate}\setcounter{enumi}{2}
\item\label{purt-sm-cor} $[h]=[{\blue{\imath}}_n^{(0)}]$ in $KL(C_n, C)$, and
\item for each $u\in\mathcal U^{(1)}_{c, 0}$, one has that
\begin{equation}\label{eq-ak-001}
(\iota_n^\ddagger\circ j(\overline{u}))^{-1}(h\oplus{\blue{\imath}}_n^{(1)})^\ddagger((\Phi'_1\circ L_1)^\ddagger(\overline{u}))=\overline{g_u}\end{equation}
\end{enumerate}
for some $g_u\in U_0(C)$ with $\textrm{cel}(g_u)<\sigma^{(1)}_{c, 2}/2$.

Define  $\Phi_1=(h\oplus {\blue{\imath}}_n^{({{1}})})\circ\Phi'_1$. By \eqref{eq-kk-002} and \eqref{purt-sm-cor}, and, by \eqref{eq-tr-002}
and (1),
\beq\label{June20-n-002}
&& {[}\Phi_1{]}_{\mathcal P^{(1)}_a\cup [L_1](\mathcal P^{(1)}_c)}=\alpha^{-1}|_{\mathcal P^{(1)}_a \cup [L_1](\mathcal P^{(1)}_c)}\andeqn\\
&&|\tau(\Phi_1(f))-\gamma(\tau)(f)|<\sigma_{c,1}^{(2)}/3\rforal f\in {\cal H}_c^{(2)}.
\eneq
Note that $\Phi_1$ is still $\mathcal G_a^{(1)+}$-$\delta_a^{(1)+}$-multiplicative, and hence \eqref{almost-cd-01} and \eqref{almost-cd-02} still hold with $\Phi_1'$ replaced by $\Phi_1$. That is, for $u\in \mathcal U_{a, 1}^{(1)},$
\begin{equation}\label{almost-cd-001}
 \mathrm{dist}(\la \Phi_1(u)\ra^k, CU(C))<\sigma_{a, 2}^{(1)}/2,
 \end{equation}
 where $k$ is the order of $\overline{u}$, and
 \begin{equation}\label{almost-cd-002}
 \mathrm{dist}(\la (\Phi_1\circ L_1)(v)\ra^{k'}, CU(A))<\sigma_{c, 2}^{(1)},
 \end{equation}
if $v\in \mathcal U_{c, 1}^{(1)}$ and if  $k'$ is the order of $\overline{v}$.
By \eqref{comp-001}, (3),  and, \eqref{comp-002}, and (4),  one has
\beq\label{uniq-kk-001}
&&[\Phi_1\circ L_1]|_{\mathcal P^{(1)}_c}=[\id]|_{\mathcal P^{(1)}_c}\andeqn\\
\label{uniq-tr-001}
&&|\tau\circ\Phi_1\circ L_1(f)-\tau(f)|<\sigma^{(1)}_{c, 1}\rforal f\in\mathcal H^{(1)}_{c} \rforal \tau\in T(C).
\eneq
Moreover, for any $u\in \mathcal U^{(1)}_{c, 0}$, one has (by \eqref{june-29-n-001} and \eqref{eq-ak-001})
\begin{eqnarray}\label{alg-k1-c-tf}
(\Phi_1 \circ L_1)^\ddagger(\overline{u})&=&(\iota_n^\ddagger\circ j(\overline{u}))\cdot\overline{g_u}=\overline{u}\cdot\overline{g_u}\approx_{\sigma^{(1)}_{c, 2}}\overline{u}.
\end{eqnarray}
Let $\overline{u}\in \overline{\mathcal{U}_{c, 1}^{(1)}},$ with order $k$. By \eqref{almost-cd-002}, there is a self-adjoint element $b\in C$ with $||b||<\sigma_{c, 2}^{(1)}$ such that
\begin{equation*}
(u^*)^k(\la (\Phi_1 \circ L_1)(u)\ra) ^k \exp(2\pi i b)\in CU(C)
\end{equation*}
{{(where we notice that  $ (u^*)^k\in CU(C)$),}} and hence
\begin{equation*}
((u^*)(\la \Phi_1 \circ L_1(u)\ra)\exp(2\pi i b/k))^k\in CU(C).
\end{equation*}
Note that
$$(u^*)(\la \Phi_1 \circ L_1(u)\ra) \exp(2\pi i b/k)\in U_0(C)$$
and $U_0(C)/CU(C)$ is torsion free (Corollary \ref{Unotrosion}). One has
$$(u^*)(\la \Phi_1 \circ L_1(u)\ra)\exp(2\pi i b/k)\in CU(C).$$
In particular, this implies that
\begin{equation}\label{alg-k1-c-tor}
\mathrm{dist}((\Phi_1 \circ L_1)^{\ddagger}({\bar u}), {\bar u})<\sigma_{c, 2}^{(1)}/k \rforal  {\bar u}\in \mathcal{U}_{c, 1}^{(1)}.
\end{equation}
{{Combining this}} with \eqref{alg-k1-c-tf},  {{we have}}
\begin{equation}\label{alg-k1-c}
\mathrm{dist}((\Phi_1 \circ L_1)^{\ddagger}({\bar u}), {\bar u})<\sigma_{c, 2}^{(1)} \rforal u\in \mathcal U_{c}^{(1)}.
\end{equation}

Therefore, by \eqref{uniq-kk-001}, \eqref{uniq-tr-001}, and \eqref{alg-k1-c}, applying {{part (a) of}} Theorem \ref{MUN1},
{{we obtain}} a unitary $U_1$ such that
$$||U_1^*(\Phi_1\circ L_1 (f))U_1-f||<\ep_1 \rforal f\in\mathcal G_1.$$
{{Replacing}} $\Phi_1$  by $\mathrm{Ad}(U_1)\circ \Phi_1$,
{{we}} may assume that
 $$||\Phi_1\circ L_1 (f)-f||<\ep_1 \rforal f\in\mathcal G_1.$$ In other words, one has the  diagram
\begin{displaymath}
\xymatrix{
C \ar[r]^{\id} \ar[d]_{L_1} & C\\
A \ar[ur]_{\Phi_1},
}
\end{displaymath}
which is approximately commutative on the subset ${\cal G}_1$ {{to}} within $\ep_1.$
{{We also notice that, by \eqref{June20-n-002}, $[\Phi_1]$ is injective on $\kappa_1^A(\overline{F_a^{(1)}}),$
which implies that  $\kappa_1^C\circ \Phi_1^{\ddag}$ is injective on $\overline{F_a^{(1)}},$ since $\overline{F_a^{(1)}}
\subset J_c^A(K_1(A)).$
It follows that $\Phi_1^{\ddag}$ is injective on $\overline{F_a^{(1)}}.$}}

We will continue to  apply {{the part (a) of}} Theorem \ref{MUN1} 
Let $\delta^{(2)}_c>0$ (in  place of $\delta$), $\mathcal G^{(2)}_c\subset C$ (in place of $\mathcal G$), $\sigma^{(2)}_{c, 1}, \sigma^{(2)}_{c, 2}>0$ (in place of $\sigma_1$ and $\sigma_2$), $\mathcal P^{(2)}_c\subset\underline{K}(C)$ (in  place of $\mathcal P$), $\overline{\mathcal U^{(2)}_c}\subset U(C)/CU(C)$ (in  place of $\mathcal U$), and $\mathcal H^{(2)}_c\subset C_{s.a}$ (in place of $\mathcal {{H}}$)
be as
{{provided}} by {{the part (a) of}} Theorem \ref{MUN1}
for  $C$ (in place of $A$), 
$\ep_3$ (in place of $\ep$), and $\mathcal G_2$ (in  place of $\mathcal F$). {{
We may assume that ${\cal U}_c^{(2)}\subset U(C)$ is a finite subset. We may also assume, \wilog, that
$\dt_c^{(2)}<\dt_a^{(1)+}/2,$ $\Phi_1({\cal G}_a^{(1)+})\subset {\cal G}_c^{(2)},$
\beq\label{May31-2019-1}
J_c^C\circ \kappa_1^C(\Phi_1^{\ddag}(\overline{{\cal U}_a^{(1)}}))\subset \overline{{\cal U}^{(2)}_c}.\eneq
By \ref{ReMUN1}, we may assume that $\overline{{\cal{U}}_c^{(2)}}\subset J_c^C(K_1(C))$
and ${\cal U}_c^{(2)}$ is in the image of $U(C_m)$ under $\iota_m$ for all $m\ge n_1>n_0.$
We may also assume, \wilog, that
$\kappa_1^{{C}}(\overline{{\cal U}_c^{(2)}})\cup [\Phi_1]({\cal P}_a^{(1)})\subset {\cal P}^{(2)}_c.$
We may further assume that ${\cal P}_{c,2}\subset {\cal P}_c^{(2)}$ and
${\cal H}(c,2)\subset {\cal H}_c^{(2)}.$}}

Denote by $F^{(2)}_c\subset U(C)$ the subgroup generated by $\mathcal U^{(2)}_c$. Since ${\cal U}^{(2)}_c$
 is finite,  we {{may}} write $\overline{F^{(2)}_c}=(\overline{F^{(2)}_{c}})_0 \oplus \mathrm{Tor}(\overline{F^{(2)}_c})$, where $(\overline{F^{(2)}_{c}})_0$ is torsion free. Fix this decomposition.
 {{Choosing a}} smaller $\sigma_{c,2}^{(2)},$ {{we}} may assume that
$$
\mathcal U^{(2)}_c=\mathcal U^{(2)}_{c, 0}\sqcup\mathcal U^{(2)}_{c, 1},
$$
where $\overline{\mathcal U^{(2)}_{c, 0}}$ generates $(\overline{F^{(2)}_{c}})_0$ and $\overline{\mathcal U^{(2)}_{c, 1}}$ generates $\mathrm{Tor}(\overline{F^{(2)}_c})$ {{(then the condition (\ref{May31-2019-1}) should be changed to {{the condition}}
that $J_c^C\circ \kappa_1^C(\Phi_1^{\ddag}(\overline{{\cal U}_a^{(1)}}))  $ is in the subgroup generated by $ \overline{{\cal U}^{(2)}_c}$)}}.
Note that, for each $u\in \mathcal U^{(2)}_{c, 1}$, one has $u^k \in CU(C)$, where $k$ is the order of $\overline{u}$.

{{Let ${\cal U}_{a,c}^{(2)'}\subset U(C)$ be a finite subset such that
$\overline{{\cal U}_{a,c}^{(2)'}}=\Phi_1^{\ddag}(\overline{{\cal U}_a^{(2)}})$ and \\
${\cal U}_c^{(2)'}={\cal U}_c^{(2)}\cup {\cal U}_{a,c}^{(2)'}\cup\{\la \Phi_1(u)\ra: u\in {\cal U}_a^{(1)}\}.$
Let $F_{c,2}$ be the subgroup of $U(C)$ generated by ${\cal U}_c^{(2)'}.$
Since $\kappa_1^A(\Phi_1^{\ddag}(\overline{{\cal U}_a^{(1)}}))=\{\kappa_1^A(\overline{\la \Phi_1(u)\ra}): u\in {\cal U}_a^{(1)}\}$
and
 since $J_c^C\circ \kappa_1^C(\Phi_1^{\ddag}(\overline{{\cal U}_a^{(1)}}))$ {{is in the subgroup generated by $ \overline{{\cal U}^{(2)}_c}$}},
$\overline{F_{c,2}}=\overline{F_c^{(2)}}+(U_0(C)/CU(C))\cap{\overline{F_{c,2}}}.$ }}

{{Let a finite subset ${\cal G}_{uc}^{(2)}\subset C$  and $\dt_{uc}^{(2)}>0$ {{satisfy}} the following condition:
for any ${\cal G}_{uc}^{(2)}$-$\dt_{c,u}^{(2)}$-multiplicative unital \cp\, $L': C\to A'$ (for any unital \CA\, $A'$
with
$K_1(A')=U(A')/U_0(A')$),
$(L')^{\ddag}$ can be defined as  a \hm\, on $\overline{F_{c,2}},$
${\rm dist}((L')^{\ddag}(\overline{u}), \overline{\la L'(u)\ra})<\min\{\sigma_{a,1}^{(1)},\sigma_{c,2}^{(2)}\}/12 $ for all $u\in {\cal U}_c^{(2)'},$
$(L')^{\ddag}({{\overline{F_{c,2}}}}\cap U_0({{C}})/CU({{C}}))\subset U_0(A')/CU(A'),$
and $\kappa_1^{A'}\circ (L')^{\ddag}(\overline{u})=[L']\circ \kappa_1^C([u])$ for all $u\in {\cal U}_c^{(2)}$
(see \ref{DLddag}). Moreover, we may also assume that
\beq\label{21-18713-f-3}
{\rm dist}(\la L'(u^k)\ra, CU(A'))<\sigma_{c,2}^{(2)}/8,
\eneq
if $u\in {\cal U}_{c,1}^{(2)}$ and if $k$ is the order of $\overline{u}$ (see \ref{DLddag}).}}

There are a finite subset ${\cal G}_0\subset C$ and a positive number $\dt_0>0$  such that, for any two
${\cal G}_0$-$\dt_0$-multiplicative \morp s $L_1'', L_2'': C\to A,$ if $$\|L_1''(c)-L_2''(c)\|<\dt_0\rforal c\in {\cal G}_0,$$
then
\beq\nonumber
&&\hspace{-0.1in}[L_1'']|_{{\cal P}_c^{(2)}}=[L_2'']|_{{\cal P}_c^{(2)}}\andeqn\\\nonumber
&&\hspace{-0.2in}|\tau\circ L_1''(h)-\tau\circ L_2''(h)|< \min\{\sigma^{(2)}_{c, 1}, \sigma^{(1)}_{a, 1} \}/12 \rforal h\in \mathcal H^{(2)}_c\cup\Phi_1(\mathcal H^{(1)}_a)\andeqn \text{for all}\,\,\tau\in T(A){{.}}
\eneq
{{Put ${\cal G}_c^{(2)+}=
{\cal G}_c^{(2)}\cup {\cal G}_{uc}^{(2)}\cup {\cal G}_0$ and
$\dt_c^{(2)+}=\min\{\dt_c^{(2)}, \dt_{uc}^{(2)}, \dt_0\}/4.$}}

{{Let $M_c=\max\{\|g\|: g\in {\cal G}_c^{(2)+}\}.$}}
Note {{that}}, by Lemma \ref{indlim-inv},
there exist a
large $m\ge n_1$ and a unital \morp\, $L_{0,2}: C\to C_m$ such that
\begin{equation}\label{June20-2}
\|\iota_m\circ L_{0,2}(g)-g\|<\dt_c^{(2)+}/8(M_c+1)\rforal g\in {\cal G}_c^{(2)+}.
\end{equation}
{{Then $L_{0,2}$ is ${\cal G}_c^{(2)+}$-$\dt_c^{(2)+}/4$-multiplicative.}}
{{Let us  now fix such an $m\ge n_1.$}}

Let $\kappa_1^{C_m}: U(C_m)/CU(C_m)\to K_1(C_m)$
be the quotient map.
We may then assume that
$L_{0,2}^{\ddag}$ is defined {{on $\overline{F_{c,2}}$}} and injective on $\overline{F_c^{(2)}}$
{{(see the discussion regarding injectivity of $(\Phi_1')^{\ddag}$---with $\af^{-1}$ {{replaced}} by $[{\rm id}_C]$)
and {{that}} ${\rm dist}(L_{0,2}^{\ddag}(\overline{u}),\overline{u})<\min\{\sigma^{(2)}_{c, 1}, \sigma^{(1)}_{a, 1} \}/12$
for all $u\in {\cal U}_c^{(2)}.$}}
It follows that $(L_{0,2}\circ \Phi_1)^\ddagger{{:=L_{0,2}^{\ddag}\circ \Phi_1^{\ddag}}}$ is defined and injective on $\overline{(F_a^{(1)})}$ {{(see the discussion of the injectivity of $(\Phi_1'\circ L_1^{\ddag}$)).}}  {{M}}oreover,
\beq\label{21-18-716-n1}
\kappa_1^{C_m}\circ(L_{0, 2}\circ\Phi_1)^{\ddagger}(g)=[L_{0, 2}\circ\Phi_1](\kappa_1^A(g)), \quad g\in\overline{(F_a^{(1)})_0},
\eneq
and by \eqref{June20-n-002}, for any $g\in\overline{(F_a^{(1)})_0}$ {{(note that $\mathcal P^{(1)}_a\supset \kappa_{1, A}((\overline{F_{a}})_0)$(in $K_1(A)$)),}}
\beq\label{21-18-716-n2}
\alpha \circ [\iota_m] \circ [L_{0, 2}\circ\Phi_1](\kappa_1^A(g)) & = & \alpha \circ [\iota_m \circ L_{0, 2}]\circ[\Phi_1](\kappa_1^A(g))\\\label{21-18-716-n2+}
&=& \alpha \circ [\Phi_1](\kappa_1^A(g))
= \kappa_1^A(g).
\eneq
Hence,
$$\alpha\circ[\iota_m]\circ \kappa_1^{C_m}\circ(L_{0, 2}\circ\Phi_1)^{\ddagger}(g)=\alpha\circ[\iota_m]\circ[L_{0, 2}\circ\Phi_1](\kappa_1^A(g))=\kappa_1^A(g) \rforal g\in\overline{(F_a^{(1)})_0},$$
{{which  also implies that $\af\circ [\iota_m]$ is injective on
$[L_{0,2}\circ \Phi_1](\kappa_1^A(\overline{(F_a^{(1)})_0}))$ and
$[L_{0,2}\circ \Phi_1]$ is injective on $\kappa_1^A(\overline{(F_a^{(1)})_0}).$}}
{{By \eqref{21-18-716-n1},
$\kappa_1^{C_m}\circ(L_{0, 2}\circ\Phi_1)^{\ddagger}$ is injective on ${\overline{(F_a^{(1)})_0}}.$ Note that $\kappa_1^A({\overline{(F_a^{(1)})_0}})$ is free abelian. Therefore $\pi_1\circ (L_{0, 2}\circ\Phi_1)^{\ddagger}$
is {{injective}}
on ${\overline{(F_a^{(1)})_0}}$ (recall that $\pi_1$ is defined in \ref{bdbk-k1}).}}
It follows from Lemma \ref{exp-length3} {{(where $B$ is replaced by $A,$ $C$ is replaced  by $C_m,$
$F$ is replaced by $(F_a^{(1)})_0,$ $\af$ is replaced by $\af\circ (\iota_m)_{*0},$  and
$L$ {{is replaced}} by $(L_{0,2}\circ \Phi_1)^\ddag$)}} that there is a homomorphism $\beta: U(C_m)/CU(C_m)\to U(A)/CU(A)$ with $\beta(U_0(C_m)/CU(C_m))\subset U_0(A)/CU(A)$ such that
\begin{equation}\label{eq-ak-lift}
\beta\circ (L_{0,2}\circ \Phi_1)^\ddagger(f)=f\rforal f\in\overline{(F_a^{(1)})_0}.
\end{equation}

{{
Put $\overline{F_{a,c,m}^{(1)}}=(L_{0,2}\circ \Phi_1)^{\ddag}(\overline{(F_a^{(1)})_0}).$
Since $\overline{(F_a^{(1)})_0}\subset J_c^A(K_1(A)),$  the equations
\eqref{21-18-716-n1} and \eqref{21-18-716-n2+}
also imply that $\kappa_1^{C_m}\circ (L_{0,2}\circ \Phi_1)^{\ddag}$ is injective on $\overline{(F_a^{(1)})_0}.$
In particular, $\overline{F_{a,c,m}^{(1)}}$ is free {{abelian}} and $\kappa_1^{C_m}$ is injective on $\overline{F_{a,c,m}^{(1)}}.$
It follows that $\pi_1|_{\overline{F_{a,c,m}^{(1)}}}$ is {{injective}} (see \ref{bdbk-k1} for the definition of $\pi_1$).
Let ${\cal U}_{a,c,m}^{(1)}\subset U(C_m)$ be a finite subset whose image  in $U(C_m)/CU(C_m)$ generates
$\overline{F_{a,c,m}^{(1)}}.$
We may assume that,
\beq\label{21-18-716-n10}
(L_{0,2}\circ \Phi_1)^{\ddag}(\overline{{\cal U}_{a,0}^{(1)}})\subset  \overline{{\cal U}_{a,c,m}^{(1)}}.
\eneq
}}

{{Let $\sigma>0,$   $\dt'>0$ (in place of $\dt$), and the finite subset ${\cal G}_c\subset C_m$ (in place of ${\cal G}$) be
 as provided by Lemma \ref{exp-length2} with respect to
$\sigma_{a, 2}^{(2)}/4$ (in place of $\ep$) (and $C_m$ in place of $C,$ $A$ in place of $B,$
${\cal U}_{a,c,m}^{(1)}$ in place of ${\cal U},$
$\overline{F_{a,c,m}^{(1)}}$ in place of $\overline{F},$
and $\bt$
in place of $\af$).}}
{{By Theorem \ref{RangT} and Remark \ref{RRrangT}, {{{{just as in}} the decomposition of $\iota_n$ in  (\ref{May11-2019}),}} one may write $\iota_m=\imath_m^{(0)}\oplus \imath_m^{(1)},$ where
$\imath_m^{(i)}$ is  a \hm\, ($i=0,1$), $\imath_m^{(0)}$ has finite dimensional range, and
$\imath_m^{(0)}$ is non-zero on each direct summand of $C_m.$ Moreover,
\beq\label{21-18713h0}
\tau(\imath_m^{(0)}(1_{C_m}))<\min\{\sigma, \sigma_{c,1}^{(2)}, \sigma_{a,1}^{(1)}\}/12\rforal \tau\in T(C).
\eneq}}

{{Let $E'$ be a finite set of generators (in the unit ball) of the finite dimensional \SCA\, $\imath_m^{(0)}(C_m)$
containing $\imath_m^{(0)}(1_{C_m}).$
Since $C$ is simple,
\beq\label{719-nnn}
\sigma_{00}=\inf\{\tau(\iota_m(a^*a): a\in E',\,\tau\in T(C)\}>0.
\eneq
Put ${\cal H}_c^{(2)+}={\cal H}_c^{(2)}\cup \Phi_1({\cal H}_a^{(1)})\cup \{a^*a: a\in E'\}.$ }}

By Theorem \ref{MEST}, there is a
${\cal G}''$-$\dt''$-multiplicative map $L'_2: C\to A$ such that
\beq\label{21-18-nalpha1}
&&\hspace{-0.4in}[L'_2]|_{\mathcal P^{(2)}_c}=\alpha|_{\mathcal{P}^{(2)}_c}\andeqn\\
\label{210-18n713}
&&\hspace{-0.4in} |\tau\circ L'_2(f)-\gamma(\tau)(f)|<\min\{{{\sigma, \sigma_{00}}},\sigma^{(2)}_{c, 1}, \sigma^{(1)}_{a, 1} \}/12\rforal f\in\mathcal H_c^{(2)+}
\eneq
 and for all  $\tau\in T(A),$ where ${\cal G}''\subset C$ is a finite subset and $\dt''>0.$
 We may assume that
 \begin{equation*}
 {\cal G''}\supset {\cal G}_0\cup \mathcal G_c^{(2)+}\cup \iota_m({\cal G}_c)\andeqn  \dt''<\min\{\dt_0, \delta_c^{(2)+},\dt'\}/2.
 \end{equation*}

 {{Fix a finite subset ${\cal G}_{c,m}^{(2)}\subset C_m$ and  $0<\dt_0^{(2)}<\min\{\sigma_{00}, \dt_0\}/2.$
 We may assume that $E'\subset \imath_m^{(0)}({\cal G}_{c,m}^{(2)}).$
  Since every finite dimensional
\CA\, is semiprojective, {{since $\imath_m^{(0)}(C_m)$ is finite dimensional,}}  and since $L_2'$ is chosen after $C_m$ is
 chosen, with sufficiently large ${\cal G}''$ and small $\dt'',$ we may assume, \wilog, that    there exists a \hm\, $h_0: C_m\to A$ with finite dimensional range}}
 {{such that
 \beq\label{729-n1}
 &&\hspace{-0.6in}\|h_0(g)-L_2'\circ \imath_m^{(0)}(g)\|<\min\{\dt_0^{(2)}, \sigma/2,  \sigma_{a, 1}^{(1)}/6,\sigma_{c,1}^{(2)}/6\}\rforal g\in {\cal G}_{c,m}^{(2)},\\
 &&\hspace{-0.6in}\|(1-h_0(1_{C_m})L_2'\circ \imath_m^{(1)}(g)(1-h_0(1_{C_m}))-L_2'\circ \imath_m^{(1)}(g)\|<\dt_0^{(2)}\rforal g\in {\cal G}_{c,m}^{(2)},\andeqn\\\label{21-18713-trace}
 &&\hspace{-0.6in} \tau(h_0(1_{C_m}))<\min\{{{\sigma/2}}, \sigma_{a, 1}^{(1)}/6,\sigma_{c,1}^{(2)}/6\}\rforal
\tau\in T(A).
 \eneq
 Let  $l_m: C_m\to (1-h_0(1_{C_m}))A(1-h_0(1_{C_m}))$
 be defined by
 $$l_m(c)=(1-h_0(1_{C_m})L_2'\circ \imath_m^{(1)}(g)(1-h_0(1_{C_m}))\rforal c\in C_m.$$}}

 {{
 Since $m$ is now fixed and $\gamma$ is a homeomorphism,
 by \eqref{210-18n713}, \eqref{729-n1}
 and \eqref{719-nnn},
 we may assume that $L_2'$ is injective on  $\imath_m^{(0)}(C_m).$
 Since $\imath_m^{(0)}$ is non-zero on each summand of $C_m,$
 by \eqref{729-n1},  \eqref{210-18n713},
 and \eqref{719-nnn},
 we may also assume that $h_0$
 is non-zero on each direct summand of $C_m.$ }}
 {{Note that}}  $L_2'\circ \imath_m={{h'_0}}\oplus  {{l}}_m^{(1)}$,
 where
 {{$h'_0=L_2'\circ \imath_m^{(0)}$  and}}
 $l_m^{(1)}=L_2'\circ \imath_m^{(1)}.$  {{(Note that $h'_0$ is close to $h_0$ by (\ref{729-n1}).)}}

 {{Choosing}} a sufficiently large ${\cal G}''$ and small $\dt'',$ we may  {{assume that}}
$(L_2'\circ\iota_m)^\ddagger$
{{and $(l_m^{(1)})^{\ddag}$ are}} defined on a subgroup of $U(C_m)/CU(C_m)$ containing $(L_{0,2}\circ \Phi_1)^\ddagger(\overline{(F_a^{(1)})_0})$, $\pi_0((L_{0,2}\circ \Phi_1)^\ddagger(\overline{(F_a^{(1)})_0}))$, $\pi_1(U(C_m)/CU(C_m)),$ and $\pi_2(U(C_m)/CU(C_m))$.
{{Moreover, for all $u\in {\cal U}_{a,c,m}^{(1)},$
\beq
&&{\rm dist}((L_2'\circ \iota_m)^{\ddag}(\overline{u}), \overline{\la L_2'\circ \iota_m(u)\ra})
<\sigma_{a,2}^{(1)}/4
\andeqn\\\label{21-18-717-n1}
&&{\rm dist}(l_m^{\ddag}(\overline{u\ra}),
\overline{\la l_m(u)\ra})<\sigma_{a,2}^{(1)}/4.
\eneq}}

Then, by Lemma \ref{exp-length2}
{{(with $h_0\oplus  l_m^{(1)}$ in place
of $\phi,$ $h_0$ in place of $\phi_0,$ and ${{l_m^{(1)}=}}L_2'\circ \iota_m^{(1)}$ in place of $\phi_1$),}} there is a homomorphism $\psi_0: C_m\to e_0'Ae_0'$, where $e_0'={{h}}_0(1_{C_m})$, such that

(i) $[\psi_0]={{[h_0]}}$ {{in $KK(C_m, A),$}}   and

(ii) for any $u\in \mathcal U_{a, 0}^{(1)}$, one has
\begin{equation}\label{eq-ak-002}
\beta({{(L_{0,2})^{\ddag}}}\circ \Phi_1^\ddagger(\overline{u}))^{-1}(\psi_0\oplus {{l}}_m^{(1)})^\ddagger({{(L_{0,2})^{\ddag}}}\circ \Phi_1^\ddagger(\overline{u}))=\overline{g_u}
\end{equation}
 for some $g_u\in U_0(A)$ with $\mathrm{cel}(g_u)<\sigma_{a, 2}^{(1)}$.

Define $L_2=(\psi_0\oplus{{l}}_m^{1})\circ L_{0,2}: C_m\to A$ {{and
$L_2^{\ddag}=(\psi_0\oplus l_m^{(1)})^{\ddag}\circ (L_{0,2})^{\ddag}.$}}
{{Then  $[L_2]|_{{\cal P}_a^{(2)}}=[{{L_2'}}]|_{{\cal P}_a^{(2)}}$  and, by {{\eqref{210-18n713}}} and
\eqref{21-18713-trace},
\beq\label{729-nn2}
|\tau(L_2(f))-\gamma(\tau)(f)|<\min\{ \sigma_{c,1}^{(2)}, \sigma_{c,2}^{(2)}\}/6\rforal f\in {\cal H}_c^{(2)}\cup\Phi_1(\mathcal H^{(1)}_a).
\eneq
}}
Also, for any $u\in\mathcal U_{a,0}^{{(1)}}$, by  \eqref{eq-ak-002}
and \eqref{eq-ak-lift},  one  then
has
\begin{equation}\label{alg-k1-a-tf}
(L_2\circ\Phi_1)^\ddagger({{\overline{u}}})=\beta({{L_{2,0}^{\ddag}}}\circ \Phi_1^\ddagger(\overline{u}))\cdot\overline{g_u}=\overline{u}\cdot \overline{g_u}\approx_{\sigma_{a, 2}^{(1)}}\overline{u} \rforal u\in\mathcal U_{a, 0}^{(1)}.
\end{equation}
Moreover,  {{since $[\Phi_1]({\cal P}_a^{(1)})\subset P_c^{(2)},$ by  {{\eqref{21-18-nalpha1}}}, (i), \eqref{June20-n-002},
and by {{\eqref{210-18n713}}}, and by  \eqref{21-18713-trace},}}
\begin{equation}\label{uniq-kk-002}
[L_2\circ \Phi_1]|_{\mathcal P^{(1)}_a}=[\id]|_{\mathcal P^{(1)}_a},\andeqn
\end{equation}
\begin{equation}\label{uniq-tr-002}
|\tau\circ L_2\circ \Phi_1(f)-\tau(f)|<\sigma^{(1)}_{a, 1} \rforal f\in\mathcal H^{(1)}_{a} \andeqn\rforal \tau\in T(A).
\end{equation}
Note that $L_2$ is still $\mathcal G_c^{(2)+}$-$\delta_c^{(2)+}$-multiplicative {{ and
$L_2\circ \Phi_1$ is ${\cal G}_a^{(1)+}$-$\dt_a^{(1)+}$-multiplicative.}} One then has that for
{{any}} $u\in \mathcal U_{a, 1}^{(1)}$ (with $k$ the order of $\overline{u}$),
\begin{equation*}
 \mathrm{dist}(\overline{\la (L_2\circ\Phi_1)(u)\ra}^k, CU(A))<\sigma_{a, 2}^{(1)}
 \end{equation*}
 {{(see \eqref{21-18713-f-2}).}}
Therefore, there is a self-adjoint element $h\in A$ with $\|h\|<\sigma_{a, 2}^{(1)}$
such that
\begin{equation*}
(u^*)^k(\la L_2\circ \Phi_1 (u)\ra)^k \exp(2\pi i h)\in CU(A),
\end{equation*}
and hence
\begin{equation*}
((u^*)(\la L_2 \circ \Phi_1(u)\ra) \exp(2\pi i h/k))^k\in CU(A).
\end{equation*}
Note that
$$(u^*)(\la L_2\circ \Phi_1 (u)\ra)\exp(2\pi i h/k)\in U_0(A)$$
and $U_0(A)/CU(A)$ is torsion free (Corollary \ref{Unotrosion}). One has that
$$(u^*)(\la L_2\circ \Phi_1(u)\ra) \exp(2\pi i h/k)\in CU(A).$$
In particular, this implies that
\begin{equation*}
\mathrm{dist}((L_2\circ \Phi_1)^{\ddagger}(u), u)< \sigma_{a, 2}^{(1)}.
\end{equation*}
{{Combining this with}} \eqref{alg-k1-a-tf}, one has
\begin{equation}\label{alg-k1-a}
\mathrm{dist}((L_2\circ \Phi_1)^{\ddagger}(u), u)<\sigma_{a, 2}^{(1)} \rforal u\in \mathcal U_{a}^{(1)}.
\end{equation}

Then,  with \eqref{uniq-kk-002}, \eqref{uniq-tr-002}, and \eqref{alg-k1-a}, applying Theorem \ref{MUN1},   one obtains a unitary $W\in A$ such that
$$||W^*(L_2 \circ \Phi_1(f))W-f||<\ep_2 \rforal f\in\mathcal F_1.$$
{{Replacing}} $L_2$ {{by}} $\mathrm{Ad}(W)\circ L_2$,  {{one then}} has
$$||L_2 \circ \Phi_1(f)-f||<\ep_2 \rforal f\in\mathcal F_1.$$ That is, one has the following diagram
\begin{displaymath}
\xymatrix{
C \ar[r]^{\id} \ar[d]_{L_1} & C \ar[d]^{L_2}\\
A \ar[ur]_{\Phi_1}\ar[r]_\id & A,
}
\end{displaymath}
with the upper triangle approximately {{commuting}} on $\mathcal G_1$
{{ to within}}  $\ep_1$ and the lower triangle approximately commuting {{on}} $\mathcal F_1$
{{to}} within $\ep_2$.
{{Recall that $L_2$ is ${\cal G}_a^{(2)+}$-$\dt_a^{(2)+}$-multiplicative,
\beq
&&[L_2]|_{{\cal P}_a^{(2)}}=\af|_{{\cal P}_a^{(2)}},\andeqn\\
&&|\tau(L_2(f))-\gamma(\tau)(f)|<\sigma_{a,1}^{(2)}/6\rforal f\in {\cal H}_a^{(2)}\,\,\hspace{0.3in}({\rm{see}}\,\,\, \eqref{729-nn2}).
\eneq
 Note, by  the choice of ${\cal G}_a^{(2)+}$ and $\dt_a^{(2)+}$ and by {{\eqref{21-18-nalpha1},}}}}
\beq
\kappa_1^C\circ L_2^{\ddag}|_{\overline{(F_a^{(2)})_0}}=[L_2]\circ \kappa_1^A|_{\overline{(F_a^{(2)})_0}}=\af\circ
\kappa_1^A|_{\overline{(F_a^{(2)})_0}}.
\eneq
{{This implies that $L_2^{\ddag}$ is injective on $\overline{(F_a^{(2)})_0}.$}}

{{Just as}} $\dt_c^{(2)+},$ ${\cal G}_c^{(2)+},$
$\sigma_{c,1}^{(2)},$ $\sigma_{c,2}^{(2)},$  ${\cal P}_c^{(2)},$ and ${\cal H}_c^{(2)}$ {{were}}
chosen {{during the construction of
 $L_2,$}}  the construction can continue.
By repeating  this argument, one obtains the following approximate intertwining diagram
\begin{displaymath}
\xymatrix{
C \ar[r]^{\id} \ar[d]_{L_1} & C \ar[d]^{L_2}  \ar[r]^{\id} & C \ar[d]^{L_3}  \ar[r]^{\id}  & C \ar[d]^{L_4} \ar[r] & \cdots\\
A \ar[ur]_{\Phi_1}\ar[r]_\id & A \ar[ur]_{\Phi_2} \ar[r]_\id & A \ar[ur]_{\Phi_3} \ar[r]_\id & A \ar[r] \ar[ur]& \cdots ,
}
\end{displaymath}
where
\beq
&&||\Phi_{n} \circ L_n(g)-g||<\ep_{2n-1}\rforal g\in\mathcal G_n,\\
&&||L_{n+1} \circ \Phi_n(f)-f||<\ep_{2n}\rforal f\in\mathcal F_n,\,\,\, n=1,2,...,\\
&&{{{[}L_n{]}|_{{\cal P}_c^{(n)}}=\af|_{{\cal P}_c^{(n)}},\,\,\, [\Phi_n]|_{{\cal P}_a^{(n)}}=\af^{-1}|_{{\cal P}_c^{(n)}},}}\\
&&{{|\tau(L_n(f))-\gamma(\tau)(f)|<\sigma_{c,1}^{(n)}/3\rforal f\in {\cal H}_c^{(n)}, \rforal \tau\in T(A),\,\,\andeqn}}\\
&&{{|t(\Phi_n)(g)-\gamma^{-1}(t)(g)|<\sigma_{a,1}^{(n)}/3\rforal g\in {\cal H}_a^{(n)},\rforal t\in T(C).}}
\eneq
By the choices of ${\cal G}_n$ {{and}} ${\cal F}_n$ and the fact that $\sum_{n=1}^{\infty}\ep_n<\infty,$
the standard Elliott approximate intertwining argument  (Theorem 2.1 of \cite{Ell-RR0}) applies, {{and}}  shows that
{{there is an isomorphism $L: A\cong B$ with inverse $\Phi$ such that  $[L]=\af,$
and $L$ induces $\gamma$}}
as desired.

\end{proof}

\begin{thm}\label{IST1}
Let $A_1, B_1\in {\cal B}_0$ be two unital separable amenable simple \CA s
satisfying the UCT.
Let $A=A_1\otimes {{U_1}}$ and $B=B_1\otimes {{U_2}},$ where $U_1$ and $U_2$ are two
 UHF-algebras of infinite type.
Suppose that ${\rm Ell}(A)\cong{\rm Ell}(B).$ Then there exists an isomorphism $\phi: A\to B$
which carries the {{isomorphism}} ${\rm Ell}(A)\cong {\rm Ell}(B).$
\end{thm}

\begin{proof}
By Theorem \ref{RangT}, there is  a \CA\, $C,$  {{constructed}} as in {{Theorem}} \ref{RangT}, such that ${\rm Ell}(C)\cong {\rm Ell}(A)\cong {\rm Ell}(B)$. By  Theorem \ref{IST0}, one has that $C\cong A$ and $C\cong B$. In particular, $A\cong B$.
\end{proof}

\begin{cor}\label{CIST0}
Let $A$ and $B$ be as in \ref{IST1}.
If there is a homomorphism $\Gamma: {\rm Ell}(B)\to {\rm Ell}(A)$
(see \ref{DEll}), {{in particular,}} $\Gamma(K_0(B)_+\setminus \{0\})\subset K_0(A)_+\setminus \{0\},$ {{and $\Gamma([1_B])=[1_A]$,}}
{{t}}hen there is a unital homomorphism $\phi: B\to A$ such that $\phi$ induces $\Gamma.$

\end{cor}

\begin{proof}

By {{Theorem}} \ref{IST1},  we may assume that $B\cong C$ for some $C$ as constructed in {{Theorem}}
\ref{RangT}. So, \wilog, we may assume that $B=C.$

{{The}} proof  is basically the same  as that of {{Theorem}} \ref{IST0} but simpler
 since we only need to
have a one-sided  approximate intertwining.  In particular, we do not need to construct $\Phi_1.$
Thus, once $L_1$ is constructed, we can go on to construct $L_2.$

{{Nevertheless,}} we will repeat the argument here.
First we keep the first paragraph at the beginning of the proof of {{Theorem}} \ref{IST0}.

Now, since $C$ satisfies the UCT, by
{{hypothesis,}} there exist an element
$\af\in KL(C,A)^{++}$ such that $\af|_{K_i(C)}=\Gamma|_{K_i(C)},$ $i=0,1,$
and a continuous affine map $\gamma: T(A)\to T(C)$
such that
$$
r_A(\gamma(t))(x)=r_B(t)(\af(x))\rforal x\in K_0(C)\andeqn \rforal t\in T(A).
$$
{{Note that $\Gamma([1_C])=[1_A]$.}}

Let $\mathcal G_1\subset \mathcal G_2\subset \cdots \subset C$
be an increasing sequence of finite subsets with dense union. Let $1/2>\ep_1>\ep_2> \cdots >0$ be a decreasing sequence of positive numbers with finite sum.
Let ${\cal P}_{c,n}\subset \underline{K}(C)$ be finite subsets such that
${\cal P}_{c,n}\subset {\cal P}_{c,n+1}$ and $\bigcup_{n=1}^{\infty}{\cal P}_{c,n}=\underline{K}(C),$
and  ${\cal H}(c,n)\subset C_{s.a.}$ be finite subsets
such that ${\cal H}(c,n)\subset {\cal H}(c,n+1)$ and $\bigcup_{n=1}^{\infty} {\cal H}(c,n)$ is dense in $C_{s.a.}.$

We will repeatedly apply {{the part (a) of}} Theorem \ref{MUN1}. Let $\delta_c^{(1)}>0$ (in place of $\delta$), $\mathcal G^{(1)}_c\subset C$ (in place of $\mathcal G$ {{and in place of $A$}}), $\sigma^{(1)}_{c, 1}, \sigma^{(1)}_{c, 2}>0$ (in place of $\sigma_1$ and $\sigma_2$), $\mathcal P^{(1)}_c\subset\underline{K}(C)$ (in place of $\mathcal P$), ${{\overline{\mathcal U^{(1)}_c}\subset U(C)/CU(C))}}$ (in place of $\mathcal U$),
and $\mathcal H^{(1)}_c\subset C_{s.a}$ (in place of $\mathcal {{H}}$) be as provided by  {{the}} {{part (a) of}} \ref{MUN1} for $C$ (in place of $A$), $\ep_1$ (in place of $\ep$), {{and}} $\mathcal G_1$ (in  place of $\mathcal F$).
{{ Here ${\cal U}^{(1)}_c$ is a finite subset of $U(C).$
As in Remark \ref{ReMUN1}, we may assume that
$\overline{\mathcal U^{(1)}_c}\subset J_c^C(K_1(C)).$}}
{{\Wlog, we may assume that the image of $\mathcal U^{(1)}_c$ in $K_1(C)$ is contained in ${\cal P}^{(1)}_c.$}}
As in {{Remark}}  \ref{ReMUN1},
we may {{also}} assume that, {{for {{all $n\ge n_0$}}  for some  large  ${{n_0}}\ge 1,$}}
there is a finite subset ${\cal V}_{n,c}^{(1)}\subset U(C_n)$ such that
${\cal U}_c^{(1)}=\iota_n({\cal V}_{n,c}^{(1)})$
{{under the map $\iota_n.$}}
We {{may}} further assume that ${\cal P}_{c,1}\subset {\cal P}_c^{(1)}$ and ${\cal H}(c,1)\subset {\cal H}_c^{(1)}.$


Let $F_c^{(1)}$ be the subgroup generated by ${\cal U}_c^{(1)}.$  To simplify notation,
{{let us}} assume that $\kappa_1^C(\overline{F_c^{(1)}})$ is
equal to
the subgroup generated by ${\cal P}_c^{(1)}\cap K_1(C).$
Put $K_{1,c}^{(1)}=\kappa_1^C(\overline{F_c^{(1)}}).$
Write $\af(K_{1,c}^{(1)})=K_{1,c,a,f}\oplus {\rm Tor}(\af(K_{1,c}^{(1)})){{\subset K_1(A)}},$
where $K_{1,c,a,f}$ is free {{abelian}}.
There is an injective \hm\, $j_\af: K_{1,c,a,f}\to K_{1,c}^{(1)}$ such
that $\af\circ j_\af={\rm id}_{K_{1,c,a,f}}.$
We may write $K_{1,c}^{(1)}=j_\af(K_{1,c,a,f})\oplus K_{2,c}^{(1)},$
where $K_{2,c}^{(1)}$ is the preimage of ${\rm Tor}(\af(K_{1,c}^{(1)}))$ under $\af.$
Recall that we assumed that $\overline{F_c^{(1)}}\subset J_c^C(K_1(C)).$
Therefore, $\overline{F_c^{(1)}}=J_c^C({{K_{1,c}^{(1)}}})=J_c^C(j_\af(K_{1,c,a,f}))\oplus J_c^C(K_{2,c}^{(1)}).$
Put $\overline{(F_c^{(1)})_0}=J_c^C(j_\af(K_{1,c,a,f}))$ (this is one of the differences from  the proof
of \ref{IST0}),
which is free {{abelian}}. Note that $\af|_{\kappa_1^C(\overline{(F_c^{(1)})_0}}$ is injective.
\Wlog\, (choosing a smaller $\dt_c^{(1)}$), we may assume that
\beq
{\cal U}_c^{(1)}={\cal U}_{c,0}^{(1)}\sqcup {\cal U}_{c,1}^{(1)},
\eneq
where $\overline{{\cal U}_{c, 0}^{(1)}}$ generates $\overline{(F_c^{(1)})_0}$  and
${\cal U}_{c,1}^{(1)}$ generates $J_c^C(K_{2,c}^{(1)})$ (another difference from the proof of \ref{IST0}).
Let $k_0$ be an integer such that $x^{k_0}=0$ for all $x\in {\rm Tor}(\af(K_{1,c}^{(1)})).$
Let ${\cal V}_{c,0}^{(1)}\subset U(C_m)$ be a finite subset such that
$\iota_n({\cal V}_{c,0}^{(1)})={\cal U}_{c,0}^{(1)}.$

Since $\overline{(F_c^{(1)})_0}$ is free {{abelian}}, there exists an injective \hm\, $j_0: \overline{(F_c^{(1)})_0}\to J_c^{C_n}(K_1(C_n))$
such that $J_c^C\circ \iota_n^{\ddag}\circ j_0={\rm id}_{\overline{(F_c^{(1)})_0}}.$
In particular, $J_c^C\circ \iota_n^{\ddag}$ is injective on $j_0(\overline{(F_c^{(1)})_0}).$
It follows that $\pi_1|_{j_0(\overline{(F_c^{(1)})_0})}$ is injective.

Let {{a finite subset}} ${\cal G}_{uc}^{(1)}\subset C$ and $\dt_{uc}^{(1)}>0$ {{satisfy}} the following
{{condition}}:
for any ${\cal G}_{uc}^{(1)}$-$\dt_{c,u}^{(1)}$-multiplicative {{unital}} \cp\, $L': C\to A'$ (for any unital \CA\, $A'$
with
$K_1(A')=U(A')/U_0(A')$),
$(L')^{\ddag}$ can be defined as  a \hm\, on $\overline{F_c^{(1)}},$
${\rm dist}((L')^{\ddag}(\overline{u}), \overline{\la L'(u)\ra})<\sigma_{c,2}^{(1)}/4 $ for all $u\in {\cal U}_c^{(1)},$
and $\kappa_1^{A'}\circ (L')^{\ddag}(\overline{u})=[L']\circ \kappa_1^C([u])$ for all $u\in {\cal U}_c^{(1)}$
(see \ref{DLddag}). {{Since $\overline{F_c^{(1)}}\subset J_c(K_1(C))$, and   $J_c(K_1(C))\cap U_0(C)/CU(C)$ only contains
{{the}} unit  $\overline{1_C}$ and since  $(L')^{\ddag}$ is a homomorphism on $\overline{F_c^{(1)}},$ we have $(L')^{\ddag}(\overline{F_c^{(1)}}\cap U_0({{C}})/CU({{C}}))=\overline{1_{A'}}\in U_0(A')/CU(A').$}}
Moreover, we may also assume that,
if $[L']|_{\kappa_1^C(\overline{F_c}^{(1)})}=[\af]|_{\kappa_1^C(\overline{F_c}^{(1)})},$ and if $u\in {\cal U}_{c,1}^{(1)},$
then
\beq\label{21-18718-n1}
{\rm dist}(\la L'(u)\ra^k, CU(A))<\sigma_{c,2}^{(1)}/4,
\eneq
where $1\le k\le k_0$ is the order of $\af([u])$ (this is another difference from the  the proof of \ref{IST0}).
Put ${\cal G}_c^{(1)+}={\cal G}_c^{(1)}\cup {\cal G}_{uc}^{(1)}$ and
$\dt_c^{(1)+}=\min\{\dt_c^{(1)}, \dt_{uc}^{(1)}\}.$

Recall that both $A$ and $C$ are in ${\cal B}_{u0}.$
By Theorem \ref{MEST},  
there is a {{${\cal G}^{(1)+}_c$-$\dt^{(1)+}_c$}}
-multiplicative \morp\, $L_1: C\to A$ such that
\beq\label{eq-kk-001++}
&&[L_1]|_{\mathcal P^{(1)}_c}=\alpha|_{\mathcal{P}^{(1)}_c}\andeqn\\\label{eq-tr-00++}
&& |\tau\circ L_1(f)-\gamma(\tau)(f)|<\sigma^{(1)}_{c, 1}/4\rforal f\in\mathcal H^{(1)}_c\rforal \tau\in T(A).
 \eneq
{{Choosing ${\cal G}^{(1)+}_c$ large enough and $\dt^{(1)+}_c$ small enough o}}ne may assume that $L_1^\ddagger$ is  a \hm\, defined  on $(\overline{F_c^{(1)}})$ (see \ref{DLddag}).
  We {{may further}}
   assume
 that ${\rm dist}(L_1^{\ddag}(\overline{u}), \overline{\la L_1(u)\ra})<\sigma_{c, 2}^{(1)}/8$ for all
 $u\in  {\cal U}_c^{(1)}$ and  $\kappa_1^A\circ L_1^{\ddag}=[L_1]\circ \kappa_1^C$ on $(\overline{F_c^{(1)}}).$
 Since $\af|_{\kappa_1^C(\overline{(F_c^{(1)})_0})}$ is injective,
 we may assume that $L_1^{\ddag}|_{\overline{(F_c^{(1)})_0}}$ is injective
 as $\kappa_1^C$ is injective on $J_c(K_1(C)).$

 Moreover, if $1\le k (\le k_0)$ is the order of $\af([u]), $  {{we}} may also assume  (see \eqref{21-18713-f-1})
 {{that}}
 \beq\label{21-18-718-n10}
 \mathrm{dist}(\la L_1(u)\ra^k, CU(A))<\sigma_{c, 2}^{(1)}/4\rforal{\rm such}\,\, u\in\mathcal U_{c, 1}^{(1)}.
 \eneq
{{Furthermore,}}
we may assume that there exists a finite subset
${\cal G}_{c,c}^{(1)}\subset C_n$ such that
$\iota_n({\cal G}_{c,c,}^{(1)})= G_c^{(1)+}.$

We now construct $L_2.$
We will continue to  apply {{the part (a) of}} Theorem \ref{MUN1}.
Let $\delta^{(2)}_c>0$ (in place of $\delta$), $\mathcal G^{(2)}_c\subset C$ (in place of $\mathcal G$), $\sigma^{(2)}_{c, 1}, \sigma^{(2)}_{c, 2}>0$ (in place of $\sigma_1$ and $\sigma_2$), $\mathcal P^{(2)}_c\subset\underline{K}(C)$ (in place of $\mathcal P$), $\overline{\mathcal U^{(2)}_c}\subset U(C)/CU(C)$ (in place of $\mathcal U$), and $\mathcal H^{(2)}_c\subset C_{s.a}$ (in place of ${{\mathcal H}}$) be as provided by {{the}} {{part (a) of}} Theorem \ref{MUN1}
for  $C$ (in place of $A$), $\ep_2$ (in place of $\ep$), and $\mathcal G_2$ (in place of $\mathcal F$).
We {{may}} assume that ${\cal U}_{{c}}^{(2)}\subset {{U(C)}}$ is a finite subset. We may also assume, \wilog, that
$\dt_c^{(2)}<\dt_c^{(1)+}/2,$ ${\cal G}_c^{(1)+}\subset {\cal G}_c^{(2)}$ and
\beq\label{May31-2019-2}{{\overline{{\cal U}_c^{(1)}}\subset \overline{{\cal U}^{(2)}_c}.}}\eneq
By {{Remark}} \ref{ReMUN1}, we may assume that
$\overline{{\cal{U}}_c^{(2)}}\subset J_c^C(K_1(C))$ and there exists an integer $n_1>n_0$ such that
${\cal U}_c^{(2)}\subset \iota_m(U(C_m))$ for all   $m\ge \max \{n_1, n\}\ge n_0.$
We may also assume, \wilog, that
${\cal G}_c^{(1)}\subset {\cal G}_c^{(2)},$ ${\cal P}_c^{(1)}\cup {\cal P}_{c,2}\subset {\cal P}_c^{(2)},$
${\cal U}_c^{(1)}\subset {\cal U}_c^{(2)},$ ${\cal H}_c^{(1)}\cup {\cal H}(c,2)\subset {\cal H}_c^{(2)},$ and $\dt_c^{(2)}<\dt_c^{(1)}.$ We may further assume that
${\cal P}_c^{(2)}\subset \iota_m({\cal P}_{c,c}^{(2)})$ for some finite subset ${\cal P}_{c,c}^{(2)}$ of
$\underline{K}(C_m),$ and ${\cal H}_c^{(2)}\subset \iota_m({\cal H}_{c,c}^{(2)})$ for some
finite subset ${\cal H}_{c,c}^{(2)}$ of $C_m$ (for all $m\ge \max\{n_1, n\}$).

Let $F_c^{(2)}$ be the subgroup generated by ${\cal U}_c^{(2)}.$
\Wlog, to simplify notation, we may assume that $\kappa_1^C(\overline{F_c^{(2)}})$ is the same
as the subgroup generated by ${\cal P}_c^{(2)}\cap K_1(C).$
Put $K_{1,c}^{(2)}=\kappa_1^C(\overline{F_c^{(2)}}).$
Write $\af(K_{1,c}^{(2)})=K_{2,c,a,f}\oplus {\rm Tor}(\af(K_{1,c}^{(2)})),$
where $K_{2,c,a,f}$ is free {{abelian}}.
There is an injective \hm\, $j_\af^{(2)}: K_{2,c,a,f}\to K_{1,c}^{(2)}$ such
that $\af\circ j_\af^{(2)}={\rm id}_{K_{2,c,a,f}}.$
We may write $K_{1,c}^{(2)}=j_\af^{(2)}(K_{2,c,a,f})\oplus K_{2,c}^{(2)},$
where $K_{2,c}^{(2)}$ is the preimage of ${\rm Tor}(\af(K_{1,c}^{(2)}))$ under $\af.$
Recall that we assumed that $\overline{F_c^{(2)}}\subset J_c^C(K_1(C)).$
Therefore $\overline{F_c^{(2)}}=J_c^C({{K_{1,c}^{(2)}}})=J_c^C(j_\af(K_{2,c,a,f}))\oplus J_c^C(K_{2,c}^{(2)}).$
Put $\overline{(F_c^{(2)})_0}=J_c^C(j_\af(K_{2,c,a,f}));$ {{this group is free abelian.}}
\Wlog\,(choosing {{a}} smaller $\dt_c^{(2)}$), we may assume that
\beq
{\cal U}_c^{(2)}={\cal U}_{c,0}^{(2)}\sqcup {\cal U}_{c,1}^{(2)},
\eneq
where $\overline{{\cal U}_{c, 0}^{(2)}}$ generates $\overline{(F_c^{(2)})_0}$  and
${\cal U}_{c,1}^{(2)}$ generates $J_c^C(K_{2,c}^{(2)})$ {{(then the condition (\ref{May31-2019-2}) should be changed to
{{the condition}}  that $ \overline{{\cal U}^{(1)}_c}$  {{be contained}} in the subgroup generated by $ \overline{{\cal U}^{(2)}_c}$, which also implies $F_c^{(1)}\subset F_c^{(2)}$)}}.
Let $k_0$ be an integer such that $x^{k_0}=0$ for all $x\in {\rm Tor}(\af(K_{1,c}^{(1)})).$


Let {{a finite subset}} ${\cal G}_{uc}^{(2)}\subset C$ and $\dt_{uc}^{(2)}>0$ {{satisfy}} the following {{condition}}:
for any ${\cal G}_{uc}^{(2)}$-$\dt_{c,u}^{(2)}$-multiplicative {{unital}} \cp\, $L': C\to A'$ (for any unital \CA\, $A'$
with
$K_1(A')=U(A')/U_0(A')$),
$(L')^{\ddag}$ can be defined as  a \hm\, on ${{\overline{F_{c}^{(2)}}}},$
${\rm dist}((L')^{\ddag}(\overline{u}), \overline{\la L'(u)\ra})<\min\{\sigma_{a,1}^{(1)},\sigma_{c,2}^{(2)}\}/12 $ for all $u\in {\cal U}_c^{(2)'},$\\
$(L')^{\ddag}(\overline{F_c^{(2)}}\cap U_0({{C}})/CU({{C}})){{=\overline{1_{A'}}\in}}\, U_0(A')/CU(A'),$
and $\kappa_1^{A'}\circ (L')^{\ddag}(\overline{u})=[L']\circ \kappa_1^C([u])$ for all $u\in {\cal U}_c^{(2)}$
(see \ref{DLddag}).
Moreover, we may also assume that, if
$[L']|_{\kappa_1^C(\overline{F_c}^{(2)})}=[\af]|_{\kappa_1^C(\overline{F_c}^{(2)})},$
\beq\label{21-18-718-n11}
{\rm dist}(\la L'(u)\ra^k, CU(A'))<\sigma_{c,2}^{(2)}/4\rforal u\in {\cal U}_{c,1}^{(2)},
\eneq
where $k$ is the order of $\af([u])$ (see {{Definition}} \ref{DLddag}).

There are a finite subset ${\cal G}_0\subset C$ and a positive number $\dt_0>0$  such that, for any two
${\cal G}_0$-$\dt_0$-multiplicative unital \cp s $L_1'', L_2'': C\to A,$ if $$\|L_1''(c)-L_2''(c)\|<\dt_0\rforal c\in {\cal G}_0,$$
then
\beq\nonumber
&&[L_1'']|_{{\cal P}_c^{(2)}}=[L_2'']|_{{\cal P}_c^{(2)}}\andeqn\\\nonumber
&&\hspace{-0.1in}|\tau\circ L_1''(h)-\tau\circ L_2''(h)|<\min\{\sigma^{(2)}_{c, 1}, \sigma^{(1)}_{a, 1} \}/12 \rforal h\in \mathcal H^{(2)}_c
\rforal \tau\in T(A){{.}}
\eneq
Put ${\cal G}_c^{(2)+}={\cal G}_c^{(2)}\cup {\cal G}_{uc}^{(2)}\cup {\cal G}_0$ and
$\dt_c^{(2)+}=\min\{\dt_c^{(2)}, \dt_{uc}^{(2)}, \dt_0, \sigma_{c,1}^{(1)}/8\pi\}/8.$

Since $\overline{(F_c^{(2)})_0}$ is free {{abelian}}, there exists an injective \hm\, $j_1: \overline{(F_c^{(2)})_0}\to
J_c^{C_m}(K_1(C_m))$
such that $J_c^C\circ \iota_m^{\ddag}\circ j_1={\rm id}_{\overline{(F_c^{(2)})_0}}.$
In particular, $J_c^C\circ \iota_m^{\ddag}$ is injective on $j_1(\overline{(F_c^{(2)})_0}).$
It follows that $\pi_1|_{j_1(\overline{(F_c^{(2)})_0})}$ is injective.
Let ${\cal V}_c^{(2)}\subset U(C_m)$ be a finite subset such that  $\iota_m(V_c^{(2)})={\cal U}_{c,0}^{(2)}$
and $\overline{V_c^{(2)}}\subset j_1(\overline{(F_c^{(1)})_0}).$

To simplify notation, \wilog, {{let us}} assume that there exists {{a}} finite subset ${\cal G}_{c,c}^{(2)}\subset C_m$
for some $m\ge n_1$ such that $\iota_m({\cal G}_{c,c}^{(2)})={\cal G}_c^{(2)+}\cup {\cal H}_c^{(2)}.$

{{Let $M_c=\max\{\|g\|: g\in {\cal G}_c^{(2)+},\,\,{\rm or}\,\, g\in {\cal G}_{c,c}^{(2)}\}.$}}
Note,  since $C_m$ and $C$ are both amenable,
there exists a unital \cp\, $L_{0,2}: C\to C_m$ such that
\beq\label{June20-2++}
&&\|\iota_m\circ L_{0,2}(g)-g\|<\dt_c^{(2)+}/8(M_c+1)\rforal g\in {\cal G}_c^{(2)+}\andeqn\\\label{7-30-n1}
&&\|L_{0,2}\circ \iota_m(c)-c\|<\dt_c^{(2)+}/8(M_c+1)\rforal c\in {\cal G}_{c,c}^{(2)}\cup {\cal V}_c^{(2)}\cup \iota_{n,m}({\cal V}_c^{(1)}).
\eneq
{{Then $L_{0,2}$ is ${\cal G}_c^{(2)+}$-$\dt_c^{(2)+}/4$-multiplicative.}}
Let $\kappa_1^{C_m}: U(C_m)/CU(C_m)\to K_1(C_m)$
be the quotient map.
We may then assume that
$L_{0,2}^{\ddag}$ is defined on $\overline{F_c^{(2)}}$ {{and}} injective on $\overline{F_c^{(2)}},$ and
{{that}}
${\rm dist}(L_{0,2}^{\ddag}(\overline{u}),\overline{u})<\min\{\sigma^{(2)}_{c, 1}, \sigma^{(1)}_{c, 1} \}/12$
for all $u\in {\cal U}_c^{(2)}.$

Note that $J_c^C\circ \iota_m^{\ddag}\circ \iota_{n,m}^{\ddag}=J_c^C\circ \iota_n^{\ddag}$
is injective on $j_0(\overline{(F_c^{(1)})_0}).$ It follows that
$\iota_{n,m}^{\ddag}$ is injective on $j_0(\overline{(F_c^{(1)})_0}).$ Also,
$(\pi_1)|_{\iota_{n,m}^{\ddag}(j_0(\overline{(F_c^{(1)})_0}))}$  is injective.
Let $\bt':=(\iota_{n,m}^{\ddag})^{-1}$ be defined on $\iota_{n,m}^{\ddag}(j_0(\overline{(F_c^{(1)})_0})).$


Let $\sigma>0,$  $\dt_c'>0$ (in place of $\dt$) and  {{the}} finite subset ${\cal G}_c'\subset C_m$ (in place of ${\cal G}$) be
{{as provided by}}  Lemma \ref{exp-length2}, with respect to
$\sigma_{a, 2}^{(2)}/16\pi$ (in place of $\ep$), $C_m$ (in place of $C$),  $A$ (in place of $B$),
$V_c^{(2)}$ (in place of ${\cal U}$),
$\iota_{n,m}^{\ddag}(j_0(\overline{(F_c^{(1)})_0})$ (in place of $\overline{F}$),
and $L_1^{\ddag}\circ \iota_n^{\ddag}\circ \bt'$
(in place of $\af$).
By  {{Theorem}} \ref{RangT} and {{Remark}} \ref{RRrangT}, {{
as in  the decomposition in (\ref{May11-2019}) of $\iota_n$ in the proof of Theorem \ref{IST0},}} one may write $\iota_m=\imath_m^{(0)}\oplus \imath_m^{(1)},$ where
$\imath_m^{(i)}$ is {{a}}  \hm\, ($i=0,1$), $\imath_m^{(0)}$ has finite dimensional range, and
$\imath_m^{(0)}$ is non-zero on each {{non-zero}} direct summand of $C_m.$ Moreover,
\beq\label{21-18713h0++}
\tau(\imath_m^{(0)}(1_{C_m}))<\min\{\sigma, \sigma_{c,1}^{(2)}, \sigma_{c,1}^{(1)}\}/12\rforal \tau\in T(C).
\eneq
Let $E'$  be a
finite set of generators (in the unit ball) of {{the}} finite dimensional \SCA\, $\imath_m^{(0)}(C_m).$
Then
\beq\label{730-n10}
\sigma_{00}:=\inf\{\tau((e')^*e'): e'\in E',\,\tau\in T(C)\}>0.
\eneq
Put ${\cal H}_c^{(2)+}={\cal H}_c^{(2)}\cup \Phi_1({\cal H}_a^{(1)})\cup \{a^*a: a\in E'\}.$

Recall that both $A$ and $C$ are in {{class}} ${\cal B}_{u0}.$
By Theorem \ref{MEST},  {{for any finite subset ${{{\cal G}''}}\supset {\cal G}^{(2)+}_c$ and
$0<{{\dt''}}< \dt^{(2)+}_c,$}}
there is a ${{{\cal G}''}}$-${{\dt''}}$
-multiplicative map $L_2': C\to A$ such that
\beq\label{eq-kk-001++++}
&&\hspace{-0.5in}[L_2']|_{\mathcal P^{(1)}_c}=\alpha|_{\mathcal{P}^{(2)}_c}\andeqn\\\label{eq-tr-00++++}
&&\hspace{-0.5in} |\tau\circ L_2'(f)-\gamma(\tau)(f)|<\min\{\sigma, \sigma_{00},\sigma^{(2)}_{c, 1}, \sigma^{(1)}_{a, 1} \}/12\rforal f\in\mathcal H_c^{(2)+}
\eneq
 and for all  $\tau\in T(A),$ where ${\cal G}''\subset C$ is a finite subset and $\dt''>0.$
 We may choose that
 \begin{equation*}
 {\cal G''}\supset {\cal G}_0\cup \mathcal G_c^{(2)+}\cup \iota_m({\cal G}_c)\andeqn \dt''<\min\{\dt_0, \delta_c^{(2)+},\dt'\}/2.
 \end{equation*}


Fix any finite subset ${\cal G}_{c,m}^{(2)}\subset C_m$ and any $0<\dt_0^{(2)}<\min\{\sigma_{00}, \dt_0\}/2.$
 We {{may}} assume that $E'\subset \imath_m^{(0)}({\cal G}_{c,m}^{(2)}).$
  Since every finite dimensional
\CA\, is semiprojective, {{$\imath_m^{(0)}(C_m)$ is finite dimensional,}}  and since $L_2'$ is chosen after $C_m$ is
 chosen, with sufficiently large ${\cal G}''$ and small $\dt'',$ we may assume, \wilog, that    there exists a \hm\, $h_0: C_m\to A$ with finite dimensional range
 such that
 \beq\label{729-n1++}
 &&\hspace{-0.5in}\|h_0(g)-L_2'\circ \imath_m^{(0)}(g)\|<\dt_0^{(2)}\rforal g\in {\cal G}_{c,m}^{(2)},\\
 &&\hspace{-0.5in}\|(1-h_0(1_{C_m}))L_2'\circ \imath_m^{(1)}(g)(1-h_0(1_{C_m}))-L_2'\circ \imath_m^{(1)}(g)\|<\dt_0^{(2)}\rforal g\in {\cal G}_{c,m}^{(2)}\andeqn\\\label{21-18713-trace++++}
 &&\hspace{-0.5in} \tau(h_0(1_{C_m}))<\min\{\sigma/2, \sigma_{a, 1}^{(1)}/3,\sigma_{c,1}^{(2)}/3\}\rforal
\tau\in T(A).
 \eneq
 Let  $l_m: C_m\to (1-h_0(1_{C_m}))A(1-h_0(1_{C_m}))$  be defined by
 $$l_m(c)=(1-h_0(1_{C_m})L_2'\circ \imath_m^{(1)}(g)(1-h_0(1_{C_m}))\rforal c\in C_m.$$

One may assume that $(L_2')^\ddagger$ is  a \hm\, defined  on $(\overline{F_c^{(2)}})$ (see \ref{DLddag}).
  We {{may}} also  assume
 that ${\rm dist}((L_2')^{\ddag}(\overline{u}), \overline{\la (L_2')(u)\ra})<\sigma_{c, 2}^{(2)}/8$ for all
 $u\in  {\cal U}_c^{(2)},$ and  $\kappa_1^A\circ (L_2')^{\ddag}=[L_2']\circ \kappa_1^C$ on $(\overline{F_c^{(2)}}).$
 Since $\af|_{\kappa_1^C(\overline{(F_c^{(2)})_0})}$ is injective, we may assume that $(L_2')^{\ddag}|_{\overline{(F_c^{(2)})_0}}$ is injective
 as $\kappa_1^C$ is injective on $J_c(K_1(C)).$
 Moreover, if $1\le k (\le k_0)$ is the order of $\af([u]), $  {{we}}  may also assume  (see \eqref{21-18713-f-1})
 that
 \beq\label{21-18-718-n10++}
 \mathrm{dist}(\la L_2'(u)\ra^k, CU(A))<\sigma_{c, 2}^{(2)}/4\rforal u\in\mathcal U_{c, 1}^{(2)}.
 \eneq

  Since $m$ is now fixed and $\gamma$ is a homeomorphism,
 by \eqref{eq-tr-00++++}, \eqref{729-n1++}, and
 \eqref{730-n10},
 we may assume that $L_2'$ is injective on  $\imath_m^{(0)}(C_m).$
 Since $\imath_m^{(0)}$ is {{non-zero}} on each {{direct}} summand of $C_m,$
 by \eqref{729-n1++},  \eqref{729-n1++}, and
 \eqref{730-n10},
 we may also assume that $h_0$
 is {{non-zero}} on each {{direct}} summand of $C_m.$

 {{Choosing}} sufficiently large ${{{\cal G}''}}$ and small ${{\dt''}},$ we may  assume that\\
$(L_2'\circ\iota_m)^\ddagger$
and $(l_m^{(1)})^{\ddag}$ are defined on a subgroup of $U(C_m)/CU(C_m)$ containing
$\iota_{n,m}^{\ddag}(j_0(\overline{(F_c^{(2)})_0}))$,
$\pi_0(\iota_{n,m}^{\ddag}(j_0(\overline{(F_c^{(2)})_0})))$, $\pi_1(U(C_m)/CU(C_m)),$ and $\pi_2(U(C_m)/CU(C_m))$.
Moreover, for all $u\in \iota_{n,m}({\cal V}_{c,0}^{(1)}),$
\beq
&&{\rm dist}((L_2'\circ \iota_m)^{\ddag}(\overline{u}), \overline{\la L_2'\circ \iota_m(u)\ra})
<\sigma_{c,2}^{(1)}/4
\andeqn\\\label{21-18-717-n1}
&&{\rm dist}(l_m^{\ddag}(\overline{u\ra}),
\overline{\la l_m(u)\ra})<\sigma_{c,2}^{(1)}/4.
\eneq

 Then, by Lemma \ref{exp-length2}
({{with}} $L_2'\circ \imath_m$ in place
of $\phi,$
and $L_2'\circ \iota_m^{(1)}$ in place of $\phi_1$), there is  a homomorphism $\psi_0: C_m\to e_0'Ae_0'$, where $e_0'=h_0(1_{C_m})$, such that

(i) $[\psi_0]=[\phi_{0}]$ in $KK(C_m, A),$ and

(ii) for any $u\in  \iota_{n,m}^{\ddag}(\mathcal V_{c,0}^{(1)}),$
one has
\begin{equation}\label{eq-ak-002++++}
(L_1^{\ddag}\circ \iota_n^{\ddag}\circ \bt'(\overline{u}))^{-1}(\psi_0\oplus l_m^{(1)})^\ddagger(\overline{u}))=\overline{g_u}
\end{equation}
 for some $g_u\in U_0(C)$ with $\mathrm{cel}(g_u)<\sigma_{a, 2}^{(1)}/64\pi$.

Define $L_2''=(\psi_0\oplus l_m^{1})\circ L_{2,0}: C\to A$ and
$(L_2'')^{\ddag}=(\psi_0\oplus l_m^{(1)})^{\ddag}\circ L_{2,0}^{\ddag}.$ Then $[L_2'']|_{{\cal P}_a^{(2)}}=[L_2']|_{{\cal P}_a^{(2)}}.$
Therefore
\beq\label{21-18-718-n21}
[L_2'']|_{{\cal P}_c^{(2)}}=\af|_{{\cal P}_c^{(2)}}.
\eneq
We also have from \eqref{eq-tr-00++++} and \eqref{21-18713-trace++++}
\beq\label{21-18-718-n22}
|\tau(L_2''(f))-\gamma(\tau)(f)|<\min\{\sigma_{c,1}^{(1)}, \sigma_{c,1}^{(2)}\}/3\rforal f\in {\cal H}_c^{(2)}.
\eneq
For $u\in
\mathcal U_{c,0}^{(1)}$,   write $u=\iota_n(v)$ for some $v\in {\cal V}_c^{(1)},$
Then, by \eqref{7-30-n1},
\beq
&&L_1^{\ddag}(\overline{u})=L_1^{\ddag}(\iota_n(v))=L_1^{\ddag}(\iota_n^{\ddag}\circ \bt'(\overline{\iota_{n,m}(v)})\andeqn\\
&&L_2''(u)=(\psi_0\oplus l_m^{(1)})(L_{2,0}(u))=(\psi_0\oplus l_m^{(1)})(L_{2,0}(\iota_m\circ \iota_{n,m}(v)))\\
&&\hspace{0.5in}\approx_{\sigma_{c,2}^{(1)}/16\pi} (\psi_0\oplus l_m^{1})(\iota_{n,m}(v)).
\eneq
It follows from \eqref{eq-ak-002++++} that, {{on the one hand,}} for $u\in {\cal U}_{c,0}^{(1)},$
\beq\label{21-18-718-n12}
{\rm dist}(L_1^{\ddag}(\overline{u}), (L_2'')^{\ddag}(\overline{u}))<\sigma_{c,2}^{(1)}.
\eneq

On the other hand, for $u\in {\cal U}_{c,1}^{(1)},$ by \eqref{21-18-718-n10} and  \eqref{21-18-718-n11},
there exists $h\in A_{s.a.}$ with $\|h\|\le \sigma_{c,2}^{(1)}/2$ such that
$$
(\la L_1(u)\ra ^k)^*\la L_2''(u)\ra^k(\exp(2\pi h/k))^k\in CU(A).
$$
Note that $[L_1(u)]=[L_2''(u)]$ for $u\in {\cal U}_c^{(1)}.$
Then, since $U_0(A)/CU(A)$ is torsion free (see Corollary \ref{Unotrosion}),
$$
(\la L_1(u)\ra)^*\la L_2''(u)\ra(\exp(2\pi h/k))\in CU(A).
$$
Then,  for $u\in {\cal U}_{c,1}^{(1)},$
\beq\label{21-18-718-n13}
{\rm dist}(L_1^{\ddag}(\overline{u}), (L_2'')^{\ddag}(\overline{u}))<\sigma_{c,2}^{(1)}.
\eneq
Now, combining
\eqref{eq-kk-001++},
\eqref{eq-tr-00++},
\eqref{21-18-718-n21}, \eqref{21-18-718-n22},
\eqref{21-18-718-n12}, and \eqref{21-18-718-n13}, applying Theorem \ref{MUN1}, {{we have}}  a unitary $u\in
A$ such that
\beq
\|{\rm Ad}\, u\circ L_2''(c)-L_1(c)\|<\ep_1\rforal c\in {\cal F}_1.
\eneq
Put $L_2= {\rm Ad}\, u\circ L_2''.$ We obtain the diagram
\begin{displaymath}
\xymatrix{
C \ar[r]^{\id} \ar[d]_{L_1} & C \ar[d]^{L_2}\\
A \ar[r]_\id & A,
}
\end{displaymath}
which is approximately commutative on ${\cal G}_1$ to within $\ep_1.$
Note that, by \eqref{21-18-718-n21},
we have
\beq
&&[L_2]|_{{\cal P}_c^{(2)}}=\af|_{{\cal P}_c^{(2)}},\\
&&|\tau(L_2(f))-\gamma(\tau)(f)|<\sigma_{c,1}^{(2)}/3\rforal f\in {\cal H}_c^{(2)}.
\eneq
We also have that  the restriction $\pi_1|_{j_1(\overline{(F_c^{(2)})_0})}$ is injective.
By
the choice of ${\cal G}_c^{(2)+},$
$\dt_c^{(2)+},$  ${\cal U}_{c,0}^{(2)},$ and ${\cal U}_{c,1}^{(2)},$  this process can continue.

{{Repeating}} this argument, one obtains the following approximate intertwining diagram:
\begin{displaymath}
\xymatrix{
C \ar[r]^{\id} \ar[d]_{L_1} & C \ar[d]^{L_2}  \ar[r]^{\id} & C \ar[d]^{L_3}  \ar[r]^{\id}  & C \ar[d]^{L_4} \ar[r] & \cdots\\
A
\ar[r]_\id & A
\ar[r]_\id & A
\ar[r]_\id & A \ar[r]
& \cdots ,
}
\end{displaymath}
where
$$||{\rm id}_A \circ L_n(g)-L_{n+1}\circ {\rm id}_C(g)||<\ep_{2n-1}\rforal g\in\mathcal G_n,\,\, n=1,2,....$$
Moreover, $[L_n]|_{{\cal P}_c^{(n)}}=\af|_{{\cal P}_c^{(n)}}$ and
$|\tau(L_n(f))-\gamma(\tau)(f)|<\sigma_{c,1}^{(n)}/3\rforal f\in {\cal H}_c^{(n)}.$
By the choices of ${\cal G}_n$  and  ${\cal F}_n$ and the fact that $\sum_{n=1}^{\infty}\ep_n<\infty,$
the standard Elliott approximate intertwining argument (Theorem 2.1 of \cite{Ell-RR0})  applies, and shows that there is a \hm\,
$\phi:C\to A$ such that $[\phi]=\af$ and
$\tau(\phi(f))=\gamma(\tau)(f)$ for all $f\in A_{s.a.}$ and $\tau\in T(C).$

\end{proof}








%
%
%




{\small
\bibliographystyle{amsplain}

\providecommand{\bysame}{\leavevmode\hbox to3em{\hrulefill}\thinspace}
\providecommand{\MR}{\relax\ifhmode\unskip\space\fi MR }
\providecommand{\MRhref}[2]{%
  \href{http://www.ams.org/mathscinet-getitem?mr=#1}{#2}
}
\providecommand{\href}[2]{#2}

}

\noindent
G. Gong: Department of Mathematics, Hebei Normal University, Shijiazhuang, Hebei 050016, China and
     Department of mathematics, University of Puerto Rico, Rio Piedras, PR 00936, USA,\\
ghgong@gmail.com,

\noindent
H.Lin:  Department of Mathematic, East China Normal University, Shanghai 200062, China\\
   and Current\\
    Department of Mathematics, University of Oregon,
    Eugene, Oregon, 97402, USA. \\
hlin@uoregon.edu,

\noindent
Z. Niu: Department of Mathematics, University of Wyoming, Laramie, WY, USA, 82071,\\
zniu@uwyo.edu.

\pagebreak

\printindex
\end{document}